\pdfoutput=1

\documentclass[letterpaper,11pt]{report}
\usepackage{amsfonts,amssymb,amsmath,amsthm,amsopn,latexsym}
\usepackage{stmaryrd}
\usepackage[all]{xy}
\usepackage{fullpage}
\usepackage{graphicx}
\usepackage{epstopdf}
\usepackage{epsfig}
\usepackage{hyperref}
\usepackage{enumitem,todonotes}
\usepackage{verbatim}
\usepackage{mathdots}

\usepackage{algorithm}
\usepackage[noend]{algpseudocode}
\usepackage{graphicx}
\usepackage[nottoc,notlot,notlof]{tocbibind}

\usepackage{sistyle}
\usepackage{mathtools}
\usepackage{stmaryrd}
\SIthousandsep{,}

\usepackage{tikz}
\usetikzlibrary{topaths,calc,decorations.pathmorphing,quotes,positioning,fit,shapes}

\hypersetup{
    colorlinks=true,
    linkcolor=blue,
    filecolor=magenta,      
    urlcolor=cyan,      
    citecolor=blue,
}

\usepackage{appendix}

\DeclareFontFamily{U}{matha}{\hyphenchar\font45}
\DeclareFontShape{U}{matha}{m}{n}{
    <5> <6> <7> <8> <9> <10> gen * matha
    <10.95> matha10 <12> <14.4> <17.28> <20.74> <24.88> matha12
}{}
\DeclareSymbolFont{matha}{U}{matha}{m}{n}
\DeclareFontSubstitution{U}{matha}{m}{n}
\DeclareMathSymbol{\thinsubset}{3}{matha}{"80}
\DeclareMathSymbol{\thinsupset}{3}{matha}{"81}

\makeatletter
  \@addtoreset{chapter}{part}
  \@addtoreset{@ppsaveapp}{part}
\makeatother

\DeclareMathOperator{\Clo}{Clo}
\DeclareMathOperator{\Pol}{Pol}
\DeclareMathOperator{\Inv}{Inv}
\DeclareMathOperator{\lcm}{lcm}
\DeclareMathOperator{\CSP}{CSP}
\DeclareMathOperator{\GenSAT}{GenSAT}
\DeclareMathOperator{\Sg}{Sg}
\DeclareMathOperator{\Cg}{Cg}
\DeclareMathOperator{\Sig}{Sig}
\DeclareMathOperator{\Forks}{Forks}
\DeclareMathOperator{\Aut}{Aut}
\DeclareMathOperator{\End}{End}
\DeclareMathOperator{\Con}{Con}
\DeclareMathOperator{\Hom}{Hom}
\DeclareMathOperator{\Sn}{Sn}
\DeclareMathOperator{\typ}{typ}
\DeclareMathOperator{\arity}{ar}

\DeclareMathOperator{\supp}{supp}
\DeclareMathOperator{\maj}{maj}
\DeclareMathOperator{\aff}{aff}

\DeclareMathOperator{\Id}{Id}

\begin{document}

\makeatletter
\newtheorem*{rep@theorem}{\rep@title}
\newcommand{\newreptheorem}[2]{%
\newenvironment{rep#1}[1]{%
 \def\rep@title{#2 \ref{##1}}%
 \begin{rep@theorem}}%
 {\end{rep@theorem}}}
\makeatother

\newtheorem{thm}{Theorem}[section]
\newreptheorem{thm}{Theorem}
\newtheorem{prop}[thm]{Proposition}
\newtheorem{cor}[thm]{Corollary}
\newtheorem{lem}[thm]{Lemma}
\newreptheorem{lem}{Lemma}

\theoremstyle{definition}
\newtheorem{defn}[thm]{Definition}
\newtheorem{claim}{Claim}
\newtheorem{conj}{Conjecture}[section]
\newtheorem{prob}{Problem}[section]

\theoremstyle{remark}
\newtheorem{rem}{Remark}[section]
\newtheorem{ex}{Example}[section]
\newtheorem{exer}{Exercise}[section]

\newcommand{\llhd}{{\mathrlap{<}\;\lhd}}

\newcommand{\Rho}{\mathrm{P}}
\newcommand{\cS}{\mathcal{S}}
\newcommand{\cM}{\mathcal{M}}
\newcommand{\cN}{\mathcal{N}}
\newcommand{\gk}{\kappa}
\newcommand{\gS}{\Sigma}
\newcommand{\gl}{\lambda}
\newcommand{\gt}{\theta}


\newcommand{\cF}{\mathcal{F}}
\newcommand{\cE}{\mathcal{E}}
\newcommand{\cG}{\mathcal{G}}
\newcommand{\cP}{\mathcal{P}}
\newcommand{\cV}{\mathcal{V}}
\newcommand{\cW}{\mathcal{W}}
\newcommand{\cB}{\mathcal{B}}
\newcommand{\cA}{\mathcal{A}}
\newcommand{\cH}{\mathcal{H}}
\newcommand{\ZZ}{\mathbb{Z}}
\newcommand{\NN}{\mathbb{N}}
\newcommand{\QQ}{\mathbb{Q}}
\newcommand{\bA}{\mathbb{A}}
\newcommand{\bB}{\mathbb{B}}
\newcommand{\bC}{\mathbb{C}}
\newcommand{\bD}{\mathbb{D}}
\newcommand{\bE}{\mathbb{E}}
\newcommand{\bF}{\mathbb{F}}
\newcommand{\bG}{\mathbb{G}}
\newcommand{\bI}{\mathbb{I}}
\newcommand{\bL}{\mathbb{L}}
\newcommand{\bM}{\mathbb{M}}
\newcommand{\bN}{\mathbb{N}}
\newcommand{\bP}{\mathbb{P}}
\newcommand{\bQ}{\mathbb{Q}}
\newcommand{\bS}{\mathbb{S}}
\newcommand{\bX}{\mathbb{X}}
\newcommand{\fA}{\mathbf{A}}
\newcommand{\fB}{\mathbf{B}}
\newcommand{\fC}{\mathbf{C}}
\newcommand{\fD}{\mathbf{D}}
\newcommand{\fG}{\mathbf{G}}
\newcommand{\fP}{\mathbf{P}}
\newcommand{\fR}{\mathbf{R}}
\newcommand{\fS}{\mathbf{S}}
\newcommand{\fT}{\mathbf{T}}
\newcommand{\fX}{\mathbf{X}}
\newcommand{\fY}{\mathbf{Y}}
\newcommand{\fZ}{\mathbf{Z}}

\newcommand{\RR}{\mathbb{R}}
\newcommand{\CC}{\mathbb{C}}
\newcommand{\FF}{\mathbb{F}}
\newcommand{\HH}{\mathbb{H}}
\newcommand{\PP}{\mathbb{P}}
\newcommand{\EE}{\mathbb{E}}

\newcommand{\cC}{\mathcal{C}}
\newcommand{\cD}{\mathcal{D}}
\newcommand{\cK}{\mathcal{K}}
\newcommand{\cL}{\mathcal{L}}
\newcommand{\cO}{\mathcal{O}}
\newcommand{\cT}{\mathcal{T}}
\newcommand{\cU}{\mathcal{U}}

\newcommand{\fp}{\mathfrak{p}}

\newcommand{\dotcup}{\ensuremath{\mathaccent\cdot\cup}}

\newcommand{\rvline}{\hspace*{-\arraycolsep}\vline\hspace*{-\arraycolsep}}

\def\D{\mathrm{d}}

\title{Notes on CSPs and Polymorphisms}
\date{}
\author{Zarathustra Brady}
\maketitle

\tableofcontents


\section{General Outline}

These notes were born from a multi-year learning seminar at MIT\footnote{This material is based upon work supported by the NSF Mathematical Sciences Postdoctoral Research Fellowship under Grant No. (DMS-1705177). Any opinions, findings, and conclusions or recommendations expressed in this material are those of the author(s) and do not necessarily reflect the views of the National Science Foundation.}. Each of the early sections corresponds roughly to a one-hour talk from the seminar, with details filled in, while the later sections were written after the seminar had completed. The subsections that occur after some of the sections consist of optional extra material that wasn't covered in the learning seminar due to time constraints. The appendices consist of longer portions of optional material - summaries of famous universal algebraic theories that are useful to know about in order to navigate the literature.

In the next section we give a teaser for the sorts of results we'll try to prove here, mainly to convince the reader that there are highly nontrivial results in this area, and that it is not just abstract nonsense. The text proper begins in Chapter \ref{chapter-intuition}, which consists of the foundational abstractions we need later, together with several fundamental examples illustrating three different behaviors of CSPs that need to be understood in order to understand the general case.

In Chapter \ref{chapter-few-subpowers}, we go over the breakthrough theory of algebras few subpowers, which lead to the first truly nontrivial algorithmic result in this area. In Chapter \ref{chapter-bounded-width}, we go over the more technically challenging theory of absorbing subalgebras and its application to CSPs of ``bounded width'' - although the algorithms used to solve CSPs in this chapter are much simpler than the ones from the previous chapter, the algebraic machinery necessary to prove that these algorithms always succeed is much more difficult (but more broadly applicable). Chapters \ref{chapter-few-subpowers} and \ref{chapter-bounded-width} do not necessarily need to be read in order, as the algebraic approaches used are quite different. In Chapter \ref{chapter-taylor} we move to trying to understand the general case of finite Taylor algebras, starting with the simpler case of conservative Taylor algebras to introduce a few of the new ideas that will be necessary to handle the general case.

Currently these notes are in an unfinished state - maybe half way through the material needed for the CSP dichotomy for finite structures, with much more planned if that is ever finished.

\section{Introduction / Advertisement}

In this section we'll state many of the results and motivating questions that we'll try to understand in these notes. If you don't understand something written here right away, don't despair - we'll go over everything in more detail later. The impatient reader can safely skip ahead to Section \ref{section-inv-pol}.

The story starts with a result of Schaefer \cite{schaefer} on a problem he called ``Generalized Satisfiability''.

\begin{defn} If $\Gamma$ is a set of relations on $\{0,1\}$, then $\GenSAT(\Gamma)$ is the decision problem which takes as input a set of variables $V$ and a collection of \emph{constraints}, where each constraint is of the form ``the relation $R(v_1, ..., v_k)$ must be satisfied'' where $(v_1, ..., v_k)$ is a tuple of variables of $V$ and $R$ is a relation from $\Gamma$ of arity $k$, and where the desired output is whether or not it is possible to assign values in $\{0,1\}$ to the variables such that the assignment satifies all of the given constraints.
\end{defn}

\begin{thm}[Schaefer \cite{schaefer}] If $\GenSAT(\Gamma)$ is not NP-complete, then $\Gamma$ is contained in one of the following sets of relations:
\begin{itemize}
\item the set of relations containing the all-$0$s vector,

\item the set of relations containing the all-$1$s vector,

\item the set of relations which can be written as an intersection of Horn clauses, where a Horn clause is a disjunction of literals such that at most one variable appears positively,

\item the set of relations which can be written as an intersection of dual-Horn clauses, where a dual-Horn clause is a disjunction of literals such that at most one variable appears negatively,

\item the set of relations which can be written as an intersection of relations involving at most two variables,

\item the set of relations which can be written as solution sets to systems of linear equations over $\bF_2$.
\end{itemize}
In each of these cases, $\GenSAT(\Gamma)$ can be solved in polynomial time.
\end{thm}

The next result of this form is due to Hell and Ne\v{s}etril \cite{h-coloring}, on a generalization of $n$-coloring which they call ``$H$-coloring''.

\begin{defn} If $H$ is a graph, then $H$-coloring is the decision problem which takes a graph $G$ as input, and where the desired output is whether or not there is a graph homomorphism from $G$ to $H$.
\end{defn}

Note that if we take $H = K_n$, then $K_n$-coloring is equivalent to $n$-coloring.

\begin{thm}[Hell, Ne\v{s}etril \cite{h-coloring}] $H$-coloring is in P if $H$ is bipartite, and it is NP-complete otherwise.
\end{thm}

These two results led Feder and Vardi \cite{feder-vardi} to ask whether there is a general dichotomy between P and NP. However, any such dichotomy has to avoid Ladner's \cite{ladner} anti-dichotomy result.

\begin{thm}[Ladner \cite{ladner}] If P $\ne$ NP, then there are problems in NP which are neither in P nor NP-complete.
\end{thm}

In order to avoid Ladner's result, Feder and Vardi focused on a special type of problem: ``constraint satisfaction problems'' (abbreviated as CSPs) with a fixed ``template''.

\begin{defn} A CSP-\emph{template} $T$ consists of a finite set $D$ together with a finite collection $\Gamma = (R_1, ..., R_n)$ of relations on $D$ - equivalently, we can think of $T$ as a relational structure $(D, R_1, ..., R_n)$. The decision problem $\CSP(T)$ takes as input a list of variables $V$ and for each $i \le n$ a list of tuples $C_i$ of variables of $V$ which are required to satisfy the constraint $R_i$, and accepts if there exists an assignment of variables to values in the set $D$ satisfying the given constraints.
\end{defn}

\begin{ex}\label{ex-k-color} The problem $k$-coloring (given a graph, determine if it can be colored with $k$ colors) is equivalent to $\CSP(\{1, ..., k\}, \ne) = \CSP(K_k)$, where $K_k$ is the complete graph of $k$ vertices (considered as a relational structure).
\end{ex}

\begin{ex} The problem 2-SAT is equivalent to $\CSP(\{0,1\}, \le, \ne)$. This problem is in P - in fact, it is known to be NL-complete (NL stands for nondeterministic logspace), and it can be solved in linear time.
\end{ex}

\begin{ex} The problem 3-SAT can be thought of as $\CSP(\{0,1\}, R_{(0,0,0)}, ..., R_{(1,1,1)})$, where $R_{(i,j,k)} = \{0,1\}^3 \setminus \{(i,j,k)\}$. We can also simplify this to the equivalent problem $\CSP(\{0,1\}, \{0,1\}^3\setminus\{(0,0,0)\}, \ne)$.
\end{ex}

\begin{ex}\label{ex-nae} The problem NAE-SAT is $\CSP(\{0,1\}, NAE)$, where $NAE = \{0,1\}^3 \setminus \{(0,0,0), (1,1,1)\}$ is the relation that states that the three variables in question are not all equal. This CSP template is known to be NP-complete (see Theorem \ref{thm-nae-np-complete}).
\end{ex}

\begin{ex}\label{ex-1-in-3} The problem 1-IN-3 SAT is $\CSP(\{0,1\}, \{(0,0,1), (0,1,0), (1,0,0)\})$. This CSP template is known to be NP-complete (see Theorem \ref{thm-1-in-3-np-complete}).
\end{ex}

\begin{ex} The problem HORN-SAT is $\CSP(\{0,1\}, \{0\}, \{1\}, \{0,1\}^3\setminus\{(1,1,0)\})$ (the third constraint is $(x \wedge y) \implies z$). This problem is known to be P-complete, and it can be solved in linear time \cite{hornsat-linear}.
\end{ex}

\begin{ex} The problem XOR-SAT is $\CSP(\{0,1\}, \{(0,0,0), (0,1,1), (1,0,1), (1,1,0)\}, \ne)$. This problem is in P - in fact, it can be solved in deterministic time $n^{\log_2(7)}$ and randomized quadratic time \cite{sparse-systems} (whether it can be solved in linear time is unknown).
\end{ex}

Generalizing the XOR-SAT example to a larger domain, we have the following very general family of problems which can be thought of as the natural generalization of systems of linear equations, over a possibly noncommutative group.

\begin{ex} Let $G$ be a finite group, and consider the CSP template with domain $G$, and with a relation $gH$ for every subgroup $H \leq G^n$ and every element $g \in G^n$, for every $n$. Note that strictly speaking, this is not a CSP (as we have defined it) since the set of relations is infinite. Feder and Vardi \cite{feder-vardi} prove that this general subgroup problem is polynomially solvable.
\end{ex}

Based on the examples they knew at the time, Feder and Vardi guessed that tractable CSPs fall into two types: ``bounded width'' problems, which are solved by local propagation of information, and problems with ``the ability to count'' such as the subgroup problems above. They further divided the bounded width problems into two main subclasses: problems with ``width $1$'' (such as HORN-SAT) and problems with ``bounded strict width'' (such as 2-SAT).

The bounded width problems can be defined formally in terms of a logic programming language called \emph{Datalog} (a simple subset of the programming language Prolog), where a program consists of rules for updating a database of known facts about tuples of variables by adding new facts if certain preconditions are met. For instance, a program to determine whether a graph is connected might have two predicates, one for the edges of the graph and another for connectivity, and a rule that says ``if connected(a,b) and edge(b,c), then add connected(a,c) to the database''. This example program maintains facts about pairs of variables, but has rules that involve examining three variables at a time.

\begin{defn} A CSP has \emph{width} $(l,k)$, $k \ge l$ if it can be solved by a Datalog program which keeps track of facts about tuples of at most $l$ variables, and updates its database by using rules that examine at most $k$ variables at a time. We say that it has \emph{width} $l$ if there exists any $k$ such that it has width $(l,k)$.
\end{defn}

In some cases we want to consider CSPs with relations of arbitrarily large arities. In these cases, one uses the concept of \emph{relational width}, introduced by Bulatov \cite{2-semilattice}, where our Datalog program is also allowed to update its database of facts about $l$-tuples of variables by using rules that examine any set of variables which is contained in the scope of some constraint relation, and to shrink our constraint relations based on facts about $l$-tuples of variables.

As it turns out, there is a canonical Datalog program for solving problems of width $(l,k)$, which correctly solves every instance of a CSP if and only if the CSP has width $(l,k)$. This program just keeps track of the set of all possible assignments to each tuple of at most $l$ variables, and eliminates possibilities from these lists by brute-forcing the set of possible assignments to each $k$-tuple of variables in turn (checking for consistency with each subset of these variables of size $\le l$), until it can no longer eliminate any further possible assignments from its database. If there are $n$ variables, this runs in time $O(n^k)$ and space $O(n^l)$.

A slight weakening of the above canonical Datalog program with width $1$, in which we only consider one relation at a time in order to remove potential values for the variables, is called ``arc-consistency'', or sometimes ``generalized arc-consistency'' if the relations have arity greater than $2$. CSPs which can be solved by arc-consistency have a special property called ``tree-duality'', which says that an instance has a solution if and only its ``universal cover'' has a solution (the universal cover is an instance with variables and constraints forming an infinite tree that corresponds to the universal cover of the (hyper-)graph of variables and constraints of the original instance).

The width of a CSP can also be defined in terms of a two player game (see \cite{pebble-game-width}), in which one player (the Prover) tries to convince the other player (the Verifier) that an instance of the CSP has a solution. The game goes as follows: in each round of the game, the Prover has assigned values to a certain tuple of at most $l$ variables (at the beginning of the game, this tuple is empty). The Verifier then picks a superset of the previous tuple of size at most $k$, and challenges the Prover to extend their assignment to this larger collection of variables. After this the Verifier selects any subset of the variables of size at most $l$, restricting the assignment to that subset, and the next round begins. The Verifier wins if at any point the Prover's assignment fails to satisfy some constraint of the CSP. Then a CSP has width $(l,k)$ if the Prover has a winning strategy only when the problem has a valid global solution.

\begin{defn} A CSP has \emph{strict width} $l$ if, whenever a partial solution to an instance of the CSP has no extension to a full solution, there exists a subset of the partial solution of size at most $l$, such that this subset already has no extension to a full solution. Equivalently, for every instance of the CSP, the projection of the solution set onto any set of $k > l$ variables is completely determined by the projections of the solution set onto subsets of those variables of size $l$.
\end{defn}

As a consequence of the above definition, if a CSP has strict width $l$, then any constraint having arity greater than $l$ must be expressible as a conjunction of constraints involving at most $l$ variables. Feder and Vardi \cite{feder-vardi} prove that one can check whether a CSP has strict width $l$ in time polynomial in the size of the domain and the constraints (for a fixed $l$), and give a necessary and sufficient criterion in terms of the existence of a near-unanimity operation of arity $l+1$ which ``preserves'' the constraints of the CSP.

In trying to understand the set of CSPs which do \emph{not} have bounded width, Feder and Vardi \cite{feder-vardi} introduced the concept of the \emph{ability to count}. Their definition of this concept is quite technical, and it was later realized that it's enough to focus on a simpler case: the affine CSP over an abelian group.

\begin{defn} For every abelian group $A$, we define the associated \emph{affine CSP} to be the CSP with domain $A$, with the ternary relation $\{(x,y,z) \mid x+y+z = a\}$ and the unary singleton relation $\{a\}$ for each element $a \in A$.
\end{defn}

In case the reader wants to see the general definition of the ability to count, we have reproduced it below.

\begin{defn} A CSP has the \emph{ability to count} if there are elements $0,1$ in the domain and there are relations $C, Z$ in the library of constraints such that $C$ is ternary, $Z$ is unary, $(0,0,1),(0,1,0),(1,0,0) \in C$, $0 \in Z$, and any instance of the CSP which satisfies the following properties has no solution:
\begin{itemize}
\item the instance only uses the constraints $C, Z$,
\item the constraints of the instance can be partitioned into two parts $A, B$ such that each variable of the instance shows up in exactly one constraint from $A$ and exactly one constraint from $B$, and
\item $A$ contains exactly one more copy of the constraint $C$ than $B$ does.
\end{itemize}
\end{defn}

Following an argument of Razborov for bipartite matching, Feder and Vardi prove the following.

\begin{thm}[Feder, Vardi \cite{feder-vardi}] Any CSP with the ability to count can't be solved by polynomial size monotone circuits. A CSP with the ability to count can never have bounded width.
\end{thm}

They then make the following two outrageous conjectures.

\begin{conj} Any CSP which can't ``simulate'' a CSP which has the ability to count \emph{does} have bounded width.
\end{conj}

\begin{conj} Any CSP which can't ``simulate'' 1-IN-3 SAT can be solved in polynomial time.
\end{conj}

Shockingly, despite seeming hopelessly vague and intractable, both of these conjectures were recently proven to be \emph{correct}! In fact, the conjecture about the ability to count holds even if we only require that our CSP can't simulate any affine CSP.

The examples of subgroup problems given above together with the concept of the ability to count also prompt the following question.

\begin{prob} What is the largest possible generalization of the Gaussian elimination algorithm?
\end{prob}

Feder and Vardi \cite{feder-vardi} made a first attempt at answering this by introducing the concept of \emph{near-subgroups} of a group, and conjectured that they also lead to CSPs that could be solved in polynomial time. Using a result of Aschbacher \cite{near-subgroups-aschbacher}, Feder \cite{near-subgroups-feder} later succeeded in showing that near-subgroup problems can be solved in polynomial time.

In this case, however, they could have asked for more. Hubie Chen \cite{chen-few-subpowers} studied the ``expressive rate'' of a constraint language $\Gamma$, which is defined as the function that takes $n$ to the logarithm of the number of $n$-variable relations which can be defined as solutions sets to CSPs over $\Gamma$. He observed that on a two element domain, this expressive rate always either grows as a polynomial or as an exponential function, and that the cases where it grows polynomially are exactly the cases where the class of relations which can be defined from $\Gamma$ is ``polynomially learnable''. The same conjecture occurs in chapter 10 of V\'ictor Dalmau's thesis \cite{dalmau-thesis}, in an algebraic form.

\begin{conj} For any constraint language $\Gamma$, the logarithm of the number of distinct $n$-variable relations which can be defined by primitive positive formulas over $\Gamma$ always either grows as a polynomial or as an exponential function. In the case of polynomial growth this class of relations is efficiently learnable and the associated CSP can be solved in polynomial time.
\end{conj}

This conjecture was resolved via the theory of algebras with ``few subpowers'', which classifies CSPs such that the solution sets always have ``compact representations'', and gives general procedures for manipulating these compact representations.

In order to approach these questions, the key conceptual ingredient turned out to be a Galois duality from universal algebra, relating a family of relations to the set of operations which ``preserve'' the relations. This allows us to view CSPs as algebraic structures in disguise, and to use algebraic techniques to study the structure of their solution sets and to design algorithms. However, the algebraic structures we end up studying are much less structured than groups or lattices - they are in a sense the most basic algebraic structures that have any good properties at all.

The new framework was introduced by Jeavons \cite{jeavons}, who reinterpreted an instance of a CSP as a homomorphism problem between relational structures.

\begin{defn} An instance of the \emph{general combinatorial problem}, or GCP, is a pair of relational structures $\langle \mathbf{A},\mathbf{B}\rangle$ having the same signature (a \emph{relational structure} is a set together with a family of named relations on that set, and the \emph{signature} of a relational structure is a list of names of relations together with specifications of their arities). The question is whether there exists a homomorphism from $\mathbf{A}$ to $\mathbf{B}$.
\end{defn}

\begin{ex} Suppose that $\mathbf{T}$ is a CSP template (in the sense of Feder and Vardi above), interpreted as a relational structure $(D, R_1, ..., R_n)$. To any instance of the CSP, we can associate a relational structure $\mathbf{X} = (V, C_1, ..., C_n)$, where $V$ is the set of variables of the instance, and each $C_i$ is a list of those tuples of variables of $V$ which are required to satisfy the constraint $R_i$. Then a homomorphism of relational structures $\mathbf{X} \rightarrow \mathbf{T}$ is the same as an assignment of values in $D$ to each variable in $V$, such that each tuple of variables in each $C_i$ is mapped to an element of $R_i$. In other words, the GCP instance $\langle \mathbf{X}, \mathbf{T}\rangle$ is equivalent to the instance of $\CSP(\mathbf{T})$ corresponding to $\mathbf{X}$.
\end{ex}

Jeavons also gives a few ways for other well-known problems (not CSPs) to be realized as instances of his general combinatorial problem.

\begin{ex} If $G$ is a graph and $K_q$ is a clique with $q$ vertices, then the GCP instance $\langle K_q, G\rangle$ is the $q$-CLIQUE problem. Note that in this case, the \emph{target} of the homomorphism is the main variable, while the source stays fixed (aside from the parameter $q$).
\end{ex}

\begin{ex} Let $G = (V,E)$ be a graph on $n$ vertices, and let $C_n = (W,F)$ be a cycle on $n$ vertices. Then the GCP instance $\langle (W,F,\ne), (V,E,\ne)\rangle$ is the problem of determining whether $G$ has a Hamiltonian circuit.
\end{ex}

These other sorts of problems, where the target of the homomorphism varies arbitrarily and the source varies according to some parameter can be studied from the point of view of \emph{parametrized complexity} and \emph{fixed parameter tractability}. It turns out that hardness and easiness in this alternative setting is determined by the \emph{treewidth} of the source structures \cite{treewidth-homomorphism}. We won't discuss this research direction much.

After demonstrating the generality of the framework, Jeavons \cite{jeavons} restricts to studying homomorphism problems with a fixed target structure $\mathbf{T}$. He calls this $\operatorname{GCP}(\Gamma)$, where $\Gamma$ is the list of relations of $\mathbf{T}$, but we will call it $\CSP(\mathbf{T})$ in these notes. Note that this is the same problem as the CSP defined in the sense of Feder and Vardi above, but the instances are now treated as relational structures (which is useful notationally), and the new perspective in terms of homomorphisms gives a hint of a more algebraic approach. For instance, the homomorphism point of view prompts the following definition.

\begin{defn} Two relational structures $\mathbf{A}, \mathbf{B}$ with the same signature are \emph{homomorphically equivalent} if there exist homomorphisms $\mathbf{A} \rightarrow \mathbf{B}, \mathbf{B} \rightarrow \mathbf{A}$.
\end{defn}

The homomorphism point of view now makes it obvious that if $\mathbf{A}$ and $\mathbf{B}$ are homomorphically equivalent, then $\CSP(\mathbf{A})$ and $\CSP(\mathbf{B})$ are equivalent problems - that is, a ``yes'' instance of one will always be a ``yes'' instance of the other. For instance, every bipartite graph $H$ having at least one edge is homomorphically equivalent to the complete graph $K_2$ on two vertices, so if $H$ is bipartite then the $H$-coloring problem is equivalent to the $K_2$-coloring problem.

Jeavons \cite{jeavons} points out that for a given CSP template, one can introduce new relations without changing the complexity of the CSP so long as these new relations are built out of the old relations in certain ways. Specifically, Jeavons shows that up to logspace reductions, we may as well assume that the collection of relations $\Gamma$ contains the equality relation, and is closed under the following four operations:
\begin{itemize}
\item permutation of inputs,

\item adding dummy variables (extra variables which are ignored by the relation),

\item existential projection onto a subset of the variables, and

\item intersection.
\end{itemize}
Note that any new relation which can be built out of these four operations can be viewed as the solution set to some instance of $\CSP(\Gamma)$, projected onto some subset of the variables. We can also think of the new relation as being defined by a \emph{primitive positive formula}, that is, a formula built out of the existential quantifier $\exists$, the relations $R_i$ of $\Gamma$ (and equality), and conjunctions $\wedge$, but which does not involve negation, disjunction, implication, or universal quantification (such a formula is called a \emph{conjunctive query} in database theory).

\begin{ex} The template we gave for HORN-SAT did not contain all possible Horn clauses - it stopped at the $3$-ary Horn clause $x\wedge y \implies z$. The $4$-ary Horn clause $x \wedge y \wedge z \implies w$ can be represented by the following primitive positive formula over HORN-SAT:
\[
\exists u\ (x \wedge y \implies u) \wedge (u \wedge z \implies w).
\]
The Horn clauses of higher arity can be represented by primitive positive formulas over HORN-SAT in a similar way.
\end{ex}

\begin{defn} A set of relations $\Gamma$ on a fixed domain $D$ is called a \emph{relational clone} if it contains the equality relation, and is closed under permutations, adding dummy variables, projection, and intersections. Equivalently, a relational clone is a set of relations which is closed under defining new relations via primitive positive formulas.
\end{defn}

The connection to algebra comes from the following fundamental result.

\begin{thm} There is a Galois duality between relational clones and clones. In particular, a relational clone is completely determined by its set of ``polymorphisms'', that is, the set of functions that ``preserve'' all of the relations of $\Gamma$.
\end{thm}

In order to understand this result we must define clones, polymorphisms, and the concept of a function preserving a relation.

\begin{defn} A set of functions $D^k \rightarrow D, k \in \mathbb{N}$ is called a \emph{clone} if it contains the \emph{projections} $\pi_i^k : D^k \rightarrow D$ which satisfy $\pi_i^k(x_1, ..., x_k) = x_i$ (generally the superscript $k$ is omitted when it is clear), and is closed under \emph{composition}, the operation which takes a $k$-ary function $f$ and $k$ $l$-ary functions $g_1, ..., g_k$ to the function
\[
(f\circ (g_1, ..., g_k)) : (x_1, ..., x_{l}) \mapsto f(g_1(x_1, ..., x_l), ..., g_k(x_{1}, ..., x_{l})).
\]
\end{defn}

The reader should play with the above definition in order to convince themself that every natural method of building new functions from old functions can be described in terms of the composition operation given above together with the projections $\pi_i^k$. For instance, the function $f(x,g(y,x))$ can be built out of $f$ and $g$ as follows:
\[
(f\circ (\pi_1, g\circ (\pi_2, \pi_1)))(x,y) = f(x,g(y,x)).
\]

\begin{defn} A $k$-ary function $f$ is said to \emph{preserve} an $m$-ary relation $R$, written $f \rhd R$, if for every choice of $k$ $m$-tuples in $R$, applying $f$ componentwise produces a new $m$-tuple which is also in $R$. If we think of elements of $R$ as column vectors, we can write this as
\[
\begin{bmatrix} x_{11}\\ \vdots\\ x_{1m} \end{bmatrix}, ..., \begin{bmatrix} x_{k1}\\ \vdots\\ x_{km} \end{bmatrix} \in R \implies f\left( \begin{bmatrix} x_{11}\\ \vdots\\ x_{1m} \end{bmatrix}, ..., \begin{bmatrix} x_{k1}\\ \vdots\\ x_{km} \end{bmatrix}\right) = \begin{bmatrix} f(x_{11}, ..., x_{k1})\\ \vdots\\ f(x_{1m}, ..., x_{km}) \end{bmatrix} \in R.
\]
A function $f$ is a \emph{polymorphism} of a relational structure $(D, \Gamma)$ or of a relational clone $\Gamma$ if $f$ preserves $R_i$ for each relation $R_i \in \Gamma$.
\end{defn}

The concept of preservation can be understood in two different ways. From the relational point of view, we have $f \rhd R$ iff $f : D^k \rightarrow D$ is a homomorphism of relational structures $(D,R)^k \rightarrow (D,R)$, where $(D,R)^k$ is the categorical $k$th power of the relational structure $(D,R)$ (the $k$th power of $(D,R)$ has underlying set $D^k$ and relation $R^k$ given by listing all $m$-tuples of $k$-tuples such that the $m$-tuple of $i$th coordinates is in $R$ for each $i \le k$). From the algebraic point of view, we have $f \rhd R$ iff the subset $R \subseteq D^m$ is a subalgebra of the algebraic structure $(D,f)^m$, where $(D,f)^m$ is the categorical $m$th power of the algebraic structure $(D,f)$, where the basic operation is simply $f$ acting componentwise on $D^m$.

The Galois duality between relational clones and clones prompts a shift in ones way of thinking about CSPs. Instead of studying a CSP template, one studies an algebraic structure whose operations are the polymorphisms of the original CSP template. Constraints that can be expressed in terms of the original library of relations become \emph{subalgebras} of powers of this algebraic structure. Instances of a CSP become questions about whether intersections of various subalgebras of a power of the original algebra are empty or not.

\begin{ex} Suppose $\mathbb{A} = (\mathbb{Z}/p, f)$ is the algebraic structure with basic operation $f : (x,y,z) \mapsto x - y + z \pmod{p}$ for some prime $p$. Then a subalgebra of $\bA^n$ - that is, a subset which is closed under $f$ - is exactly the same as an \emph{affine linear subspace} of $(\mathbb{Z}/p)^n$ (recall that affine linear subspaces are like vector subspaces, but that they might not pass through the origin). Checking whether a collection of affine linear subspaces has a nonempty intersection is equivalent to solving a system of linear equations $\pmod{p}$.
\end{ex}

By using an old result classifying the minimal (nontrivial) clones on the domain $\{0,1\}$ (under the Galois duality, a minimal clone of functions corresponds to a maximal relational clone), Jeavons \cite{jeavons} was able to give a new and relatively simple proof of Schaefer's dichotomy theorem \cite{schaefer}. The algebraic structures corresponding to the basic polynomial time solvable problems are as follows.

\begin{ex} If $\Gamma = (\{0\}, \{1\}, \{0,1\}^3\setminus \{1,1,0\})$ is the template corresponding to HORN-SAT, then the clone of polymorphisms is generated by the function $\min : \{0,1\}^2 \rightarrow \{0,1\}$. This operation is an example of a \emph{semilattice} operation.
\end{ex}

\begin{ex} If $\Gamma = (\le, \ne)$ is the template corresponding to 2-SAT, then the clone of polymorphisms is generated by the \emph{majority} (or \emph{median}) function $\operatorname{maj} : \{0,1\}^3 \rightarrow \{0,1\}$.
\end{ex}

\begin{ex} If $\Gamma = (\{(0,0,0), (0,1,1), (1,0,1), (1,1,0)\}, \ne)$ is the template corresponding to XOR-SAT, then the clone of polymorphisms is generated by the \emph{affine linear} function $(x,y,z) \mapsto x-y+z \pmod{2}$ (this function is sometimes referred to as the \emph{minority} function).
\end{ex}

Early results focused on generalizing these basic examples, and developing the algebraic perspective further:
\begin{itemize}
\item If all polymorphisms of $\Gamma$ are unary, then $\CSP(\Gamma)$ is NP-hard by a gadget reduction from NAE-SAT (if the domain has size $2$) or $k$-coloring (if the domain has size $k \ge 3$).

\item Generalized arc-consistency solves any CSP which has an associative, commutative, idempotent polymorphism. These types of operations were called ACI operations at the time, but are now generally referred to as semilattice operations.

\item Later, Dalmau and Pearson \cite{dalmau-width-1} showed that generalized arc-consistency solves a CSP \emph{iff} it has a ``set'' polymorphism, also known as a family of ``totally symmetric'' polymorphisms, where the output depends only on the set of inputs and not on their order or multiplicity.

\item Already in Feder and Vardi's work \cite{feder-vardi}, it was shown that a CSP has strict width $l$ iff it has an $l+1$-ary ``near-unanimity'' polymorphism, that is, an operation such that whenever all but one of the inputs are equal, their common value is the output. This fact is closely connected to a result in universal algebra known as the Baker-Pixley theorem \cite{baker-pixley}.

\item Bulatov and Dalmau \cite{bulatov-dalmau-malcev} gave an algorithm generalizing Gaussian elimination as well as the algorithm for the general subgroup problem introduced by Feder and Vardi to the case of CSPs with a \emph{Mal'cev polymorphism} $p(x,y,z)$, that is, a polymorphism satisfying the identity $p(x,y,y) = x = p(y,y,x)$ for all $x,y$. In the case of groups such an operation is given by $(x,y,z) \mapsto xy^{-1}z$, but such operations also exist in quasigroups, making this a very wide generalization.

\item One can restrict to the case where $\Gamma$ contains a unary relation for every singleton subset of the domain. On the algebraic side, this corresponds to restricting to the case of idempotent algebras (that is, algebras where every singleton subset forms a subalgebra).

\item At this point multiple authors started to realize that whether a CSP is hard or not doesn't depend on the particular polymorphisms, but rather on the \emph{identities} that are satisfied by the polymorphisms. One of the first papers to point this out was a paper by Bulatov and Jeavons \cite{bulatov-jeavons-varieties} which also introduced a notion of polymorphisms for \emph{multisorted} relations, as well as the use of ``tame congruence theory'' from universal algebra.

\item In particular, it was shown that if no finite subset of the identities satisfied by the polymorphisms imply that the polymorphisms can't be unary, then $\CSP(\Gamma)$ is NP-hard by a gadget reduction from NAE-SAT or $k$-coloring. It was conjectured that this is an if and only if, that is, if a CSP has a family of polymorphisms that satisfy a nontrivial set of identities, then the CSP can always be solved in polynomial time.

\item Identities which do not involve composing functions with each other were singled out as special (such identities are called \emph{linear} identities, or identities of \emph{height at most} $1$). In \cite{indicator-instance}, a simple procedure was given to transform the search for polymorphisms of $\Gamma$ into an equivalent \emph{indicator instance} of $\CSP(\Gamma)$ - and a simple modification of this procedure can be used to search for the polymorphisms which satisfy a given collection of linear identities. In order to attack the \emph{meta-problem} (which takes the set of relations $\Gamma$ as input, and determines whether a given type of algorithm can solve $\CSP(\Gamma)$ as output), one then only has to find a way to solve the corresponding type of indicator instance.
\end{itemize}

The first big result was a comprehensive generalization of Gaussian elimination, generalizing the algorithm for Mal'cev operations as far as was reasonably possible. The basic idea here is to represent the solution space of all the constraints processed so far by giving a small generating set for that solution space, considered as a subalgebra of a power of the domain. In order for any algorithm along these lines to exist, there must first be a guarantee that every subalgebra of any power of the domain actually \emph{has} a small generating set.

\begin{thm}[Few Subpowers \cite{few-subpowers}, \cite{few-subpowers-algorithm}] The following are equivalent for an algebraic structure $\mathbb{A}$ on a finite domain:
\begin{itemize}
\item the number of subalgebras of $\bA^n$ grows like $2^{O(n^k)}$ for some fixed $k$,

\item every subalgebra of $\bA^n$ has a (nice) generating set of size $O(n^k)$ for some fixed $k$, called a \emph{compact representation},

\item $\bA$ has a $k$-\emph{edge term} for some $k$, that is, a term satisfying the ``shepherd's crook'' identity
\[
f\left(\begin{bmatrix}y & y & x & x & \cdots & x\\ y & x & y & x & \cdots & x\\ x & x & x & y & \cdots & x\\ \vdots & \vdots & \vdots & \vdots & \ddots & \vdots \\ x & x & x & x & \cdots & y \end{bmatrix}\right) = \begin{bmatrix}x\\ x\\ x\\ \vdots \\ x\end{bmatrix},
\]
where all but the first column has exactly one $y$ and $k-1$ $x$s, and

\item a Gaussian-elimination-like algorithm solves $\CSP(\bA)$ in polynomial time (the degree of the polynomial may depend on $k$).
\end{itemize}
\end{thm}

Note that when $k = 2$, a $k$-edge term is the same as a Mal'cev operation (up to permuting the inputs). Additionally, if a $k$-edge term ignores its first input, then it is a $k$-ary near-unanimity term. So few subpowers algebras generalize both subgroup CSPs and CSPs with bounded strict width. There is still an important open question connected to few subpowers algebras, related to the following algebraic problem.

\begin{prob}[Subpower Membership Problem] Given a finite subset $S \subseteq \bA^n$ and an element $x \in \bA^n$, determine if $x$ is in the subalgebra of $\bA^n$ generated by $S$.
\end{prob}

A result of Kozik \cite{kozik-subpower-exptime} shows that for general algebraic structures, this problem can be EXPTIME-complete. A recent result of Shriner \cite{subpower-hardness} showed that this problem remains hard even if we restrict ourselves to congruence distributive algebras.

\begin{conj} If $\bA$ has few subpowers, then the subpower membership problem can be solved in polynomial time.
\end{conj}

Peter Mayr \cite{subpower-supernilpotent} has shown that this conjecture holds for nilpotent Mal'cev algebras of prime power order (and for expansions of such algebras). In a different direction, a recent result of Bulatov, Mayr, and Szendrei \cite{subpower-residually-small} has proved the conjecture in the special case that the algebra $\bA$ is ``residually small'' (for those with a little universal algebra background, this means that every subdirectly irreducible algebra in the variety it generates has size bounded by some fixed cardinal - in the case of groups, it is equivalent to all Sylow subgroups being abelian). In the same paper, they also show that the subpower membership problem for algebras with few subpowers is always in NP. As far as I know, the above conjecture is still open even for the special case of quasigroups.

The second big result in this story was the classification of CSPs with bounded width, together with a surprising ``collapse'' of the bounded width hierarchy. The ideas used in the proof of this result - especially the theory of absorbing subalgebras - led to a number of breakthrough results in universal algebra.

\begin{thm}[Bounded Width Algebras \cite{bulatov-bounded}, \cite{barto}, \cite{sdp}, \cite{ability-to-count}, \cite{local-consistency}, \cite{maltsev}, \cite{slac}] For an idempotent algebra on a finite domain, the following are equivalent:
\begin{itemize}
\item $\CSP(\bA)$ has bounded width,

\item $\CSP(\bA)$ can't simulate any CSP which has the ability to count,

\item $\CSP(\bA)$ has relational width $(2,3)$,

\item $\CSP(\bA)$ can be solved by a ``cycle consistency'' algorithm, which checks arc-consistency and checks that every ``cycle'' of constraints has a valid solution extending each possible value of every variable in the cycle,

\item $\bA$ has terms $f,g$ of arity $3$ satisfying the identities
\[
g(x,x,y) = g(x,y,x) = g(y,x,x) = f(x,x,y) = f(x,y,x) = f(x,y,y),
\]

\item $\CSP(\bA)$ can be ``robustly'' solved by the basic semidefinite programming relaxation (i.e., if an $\epsilon$-portion of the constraints are garbled, then the basic SDP can be used to find an assignment that satisfies all but an $f(\epsilon)$-portion of the constraints, where $f(\epsilon) \rightarrow 0$ as $\epsilon \rightarrow 0$).
\end{itemize}
\end{thm}

Furthermore, there is a polynomial time algorithm for checking whether a relational structure (which has all unary singleton relations) has bounded width, and to find terms $f,g$ as in the above theorem if it does. This algorithm leverages the fact that the canonical width-$(2,3)$ Datalog program will correctly solve any instance of any bounded width CSP in polynomial time, and constructs a CSP whose solution corresponds to a pair of polymorphisms satisfying a nice set of identities. A similar algorithm is not known to exist for checking that a CSP has width $1$.

The third big result of the subject was a nice classification of the algebras which were conjectured to correspond to CSPs with polynomial time algorithms. This result was proved with the absorption theory that had been developed for the study of bounded width algebras.

\begin{thm}[Cyclic Terms for Taylor Algebras \cite{cyclic}] For an idempotent algebraic structure $\bA$ on a finite domain, the following are equivalent:
\begin{itemize}
\item there is a finite set of identities satisfied by the operations of $\bA$ which can't be satisfied by essentially unary functions,

\item for every prime $p > |\bA|$, $\bA$ has a ``cyclic term'' $f$ of arity $p$, that is, a term which satisfies the identity
\[
f(x_1, ..., x_{p-1}, x_p) = f(x_2, ..., x_p, x_1),
\]

\item $\bA$ has a $4$-ary term $t$ which satisfies the identity
\[
t(x,x,y,z) = t(y,z,z,x).
\]
\end{itemize}
\end{thm}

With this in hand, the main conjecture of the subject was finally possible to state simply: ``$\CSP(\bA)$ is in P iff $\bA$ has a cyclic term.''

The fourth big result of the subject concerns a generalization of CSPs to ``valued constraints''.

\begin{defn} A \emph{valued constraint} on $k$ variables is a cost function from $D^k$ to $(-\infty, \infty]$. An instance of a valued CSP (abbreviated VCSP) consists of a sum of valued constraints applied to various tuples of the variables, possibly with nonnegative coefficients. The goal is to minimize the sum of the cost functions. The associated CSP to a VCSP is the problem of finding an assignment that makes all of the costs finite.
\end{defn}

The Galois duality between clones and relational clones can be generalized to a duality between VCSP templates and ``fractional polymorphisms'' - essentially just formal convex combinations of ordinary polymorphisms, with the property that when they are applied to any tuple of elements of $D^k$, on average they decrease the cost assigned by any cost function from the VCSP template.

The standard example of a valued constraint with an interesting fractional polymorphism is a \emph{submodular function}, that is, a cost function $c$ defined on a lattice (or a power of a lattice) which satisfies the inequality
\[
\tfrac{1}{2}c(A) + \tfrac{1}{2}c(B) \ge \tfrac{1}{2}c(A\vee B) + \tfrac{1}{2}c(A\wedge B).
\]
It is well known that submodular cost functions can be minimized using a linear programming relaxation.

\begin{thm}[VCSP Dichotomy \cite{vcsp-hardness}, \cite{vcsp-algorithm}] If a VCSP has its associated CSP in P and has a cyclic fractional polymorphism, then by using the algorithm for the associated CSP as a black box to get the set of ``feasible'' values for each variable and applying the basic linear programming relaxation to the restriction of the VCSP to the feasible values we get the minimum cost solution. If the VCSP has no cyclic fractional polymorphisms, then it is NP-hard.
\end{thm}

Finally, the biggest result of all was recently proved independently by Bulatov and by Zhuk.

\begin{thm}[CSP Dichotomy \cite{bulatov-dichotomy}, \cite{zhuk-dichotomy}] A finite algebra $\bA$ has a cyclic term iff $\CSP(\bA)$ is in P (assuming P $\ne$ NP).
\end{thm}

A major open problem is whether one can test if $\CSP(\Gamma)$ is in P in time polynomial in the size of the description of the constraints of $\Gamma$ (given that $\Gamma$ contains all singleton unary relations). Zhuk tells me that he conjectures it to be NP-hard to test for the existence of cyclic terms. If so, perhaps this could lead to a new form of public key cryptography, where the private key is a cyclic term, and the public key is a CSP template which is preserved by that cyclic term...

The story has not ended with the proof of the main conjecture. There are at least six interesting research directions that are still being actively investigated: qualitative CSPs, counting complexity, promise problems, quantified CSPs, ``hybrid'' tractability (combining restrictions on both the source and the target relational structures), and planar CSPs.

Qualitative CSPs come from allowing the domain of the CSP to be infinite. Of course, this immediately leads to problems, for instance, how can one even specify a set of relations on an infinite domain? The idea, capturing good old fashioned AI intuition about ``qualitative'' reasoning, is to require that the specific values of the solutions are not important, just the qualitative relationships between them. To make this more precise, we require our domain (and the relations on it) to have a very large automorphism group.

\begin{defn} A permutation group $G$ acting on a set $S$ is \emph{oligomorphic} if $S^n$ has finitely many $G$-orbits for every $n \ge 1$.
\end{defn}

A standard example of an oligomorphic group is the group of order-preserving bijections on the rational numbers. Relations invariant under this group give rise to ``temporal'' CSPs, where the goal is to find some assignment of variables to times satisfying constraints about their relative ordering.

Bodirsky, in his habilitation thesis \cite{bodirsky-thesis} introduced qualitative CSPs and gave a number of classification results. Before beginning a classification, he first chooses an oligomorphic group $G$ acting on a countable set $S$. He then uses results from model theory (specifically, $\omega$-categorical theories and Fra\"iss\'e limits) as well as structural Ramsey theory (and the theory of extremely amenable groups) to understand the relations which are invariant under $G$ and their polymorphisms, and for several groups $G$ he succeeds in finding a complete classification of problems into ``easy'' and ``hard''. The main three cases considered by Bodirsky \cite{bodirsky-thesis} are the following:
\begin{itemize}
\item the automorphism group of $(\mathbb{Q},<)$, corresponding to temporal CSPs,

\item the automorphism group of the random graph, for which he proves ``Schaefer's Theorem for graphs'' (such CSPs can be interpreted as problems where the variables correspond to decisions about whether certain pairs of vertices of an unknown graph are connected by an edge or not), and

\item the automorphism group of an infinite branching tree structure $(L,|)$, where $|$ is a $3$-ary relation where $ab|c$ means that the youngest common ancestor of $a,b$ lies below the youngest common ancestor of $b,c$ - the invariant relations correspond to ``branching time constraints'', or ``phylogeny constraints'', and the associated CSPs could in principle be of interest to biologists.
\end{itemize}

Recent results on infinite CSPs indicate that the difficulty of the classification results seems to be related to the orbit growth function of the oligomorphic group $G$, which takes $n$ to the number of orbits of $n$-tuples under $G$ \cite{small-orbit-growth}. For sufficiently small orbit growth functions, a dichotomy result has been proven (using the finite case as a black box). The main conjecture in this field is the following somewhat technical statement.

\begin{conj}[\cite{topology-irrelevant}] Let $\fA$ be the core of a reduct of a finitely bounded homogeneous structure. Then $\CSP(\fA)$ is in P iff $\fA$ has a 6-ary polymorphism $s$ and unary polymorphisms $\alpha, \beta$ satisfying the ``pseudo-Siggers'' identity:
\[
\alpha\circ s(x,y,x,z,y,z) = \beta\circ s(y,x,z,x,z,y).
\]
Otherwise, $\CSP(\fA)$ is NP-complete.
\end{conj}

The most recent development in the study of CSPs is the study of promise problems. Promise problems are similar in spirit to approximation algorithms, but much more amenable to an algebraic approach. A promise problem is defined here to be a pair of problems, one more restrictive than the other, where the goal is to give an algorithm which correctly says ``yes'' if the less restrictive problem has a solution and says ``no'' if the more restrictive problem has no solution (if neither case holds, any output is allowable).

\begin{defn} If $\mathbf{A}, \mathbf{B}$ are relational structures such that a homomorphism $\mathbf{A} \rightarrow \mathbf{B}$ exists, then $\operatorname{PCSP}(\mathbf{A},\mathbf{B})$ is the following problem. The input is a relational structure $\mathbf{X}$ such that there exists a homomorphism $\mathbf{X} \rightarrow \mathbf{A}$ (the promise), although this map is not revealed to us. The desired output is a homomorphism from $\mathbf{X}$ to $\mathbf{B}$.
\end{defn}

A typical strategy for proving tractability of $\operatorname{PCSP}(\mathbf{A},\mathbf{B})$ is to find a relational structure $\mathbf{C}$ such that there exist homomorphisms $\mathbf{A} \rightarrow \mathbf{C} \rightarrow \mathbf{B}$ and such that $\CSP(\mathbf{C})$ is in P.

\begin{ex} Let $\mathbf{A}$ be 1-IN-3 SAT and let $\mathbf{B}$ be NAE-SAT (where the 1-IN-3 relation and the NAE relation have the same name in the signature). The identity map on the domain gives a homomorphism $\mathbf{A} \rightarrow \mathbf{B}$ since the 1-IN-3 relation is contained in the NAE relation. Although both problems are NP-complete, the PCSP associated to the pair is tractable: let $\mathbf{C} = (\mathbb{Z}, x+y+z=1)$, note that the inclusion map $\mathbf{A} \rightarrow \mathbf{C}$ is a homomorphism, and that the map $\operatorname{sgn}: \mathbb{Z} \rightarrow \{0,1\}$ given by
\[
\operatorname{sgn}(x) = \begin{cases} 0 & x \le 0\\ 1 & x \ge 1\end{cases}
\]
defines a homomorphism $\mathbf{C} \rightarrow \mathbf{B}$. The CSP associated to $\mathbf{C}$ is tractable (even though it is defined over an infinite domain), since it boils down to solving a (sparse) system of linear equations over the integers. It is not possible to find a \emph{finite} relational structure $\mathbf{C}$ with polynomial time CSP that fits between 1-IN-3 SAT and NAE-SAT (see \cite{pcsp-coloring-full} for a proof).
\end{ex}

The relevant algebraic object in this context is $\Pol(\mathbf{A},\mathbf{B})$, the set of homomorphisms $\mathbf{A}^k \rightarrow \mathbf{B}$. At first this structure doesn't seem algebraic at all, since there is no way to compose elements of $\Pol(\mathbf{A},\mathbf{B})$. However, one can still write down ``minor identities'' between the functions in $\Pol(\mathbf{A},\mathbf{B})$ such as $f(x,x,y) = g(y,x)$, and compare the set of minor identities obtained to the identities that occur in polymorphism algebras of tractable CSPs. This approach of studying minor identities has been surprisingly useful, and has led to the proposal to call sets of functions such as $\Pol(\mathbf{A},\mathbf{B})$ ``minions'' (a competing proposed name is ``clonoid'').

Unlike the situation for CSPs, it is still quite hard to prove hardness results for PCSPs. The following basic problem is still wide open.

\begin{conj} For any $k \ge l \ge 3$, the promise problem $\operatorname{PCSP}(K_l,K_k)$ is NP-hard (this problem is the problem of $k$-coloring a graph which is promised to be $l$-colorable).
\end{conj}

One of the first results in this direction concerned a PCSP called $(2+\epsilon)$-SAT, where one is given clauses of $2k+1$ variables and wants to satisfy the associated instance of SAT given the promise that it is possible to find an assignment in which every clause has at least $k$ satisfied literals. The $(2+\epsilon)$-SAT problem was proven to be NP-hard \cite{2+epsilon}, and this result was slightly generalized and put into the PCSP framework in \cite{pcsp-symmetric-boolean}.

Recent results in the study of PCSPs include a result of Barto, Bulin, Opr\v{s}al, and Krokhin \cite{pcsp-coloring-full} in which they used minion techniques to show that $\operatorname{PCSP}(K_d, K_{2d-1})$ is NP-hard for every $d \ge 3$, reducing from the hypergraph promise problem $\operatorname{PCSP}(\operatorname{NAE}_2,\operatorname{NAE}_k)$. The hypergraph coloring problem $\operatorname{PCSP}(\operatorname{NAE}_2,\operatorname{NAE}_k)$ was itself shown to be hard via a reduction from a variant of the PCP theorem \cite{hypergraph-promise-hardness}. The big result in PCSPs is the following result which connects computational complexity to height 1 identities satisfied by the minion of polymorphisms.

\begin{thm}[Barto, Bulin, Opr\v{s}al, Krokhin \cite{pcsp-coloring-full}] If there is a ``minion homomorphism'' from $\Pol(\fA_1,\fB_1)$ to $\Pol(\fA_2,\fB_2)$, then $\operatorname{PCSP}(\fA_2,\fB_2)$ has a logspace reduction to $\operatorname{PCSP}(\fA_1,\fB_1)$.
\end{thm}

\begin{rem} For those who like category theory, an abstract minion is just a covariant functor from the category of (finite) sets to the category of sets, and a minion homomorphism is just a natural transformation of functors. We could say that a ``representation'' of an abstract minion over $A,B$ is a natural transformation to the functor $I \mapsto \Hom(A^I,B)$.

$\operatorname{PCSP}(\fA,\fB)$ ends up being logspace equivalent to the problem of distinguishing between diagrams in the category of sets of size at most $N$ ($N$ any fixed large enough number) which have a nonempty limit (``yes'' instances), and diagrams such that the image under the minion $\Pol(\fA,\fB)$ has an empty limit (``no'' instances) (this is the ``promise satisfaction of a minor condition'' problem of \cite{pcsp-coloring-full}).
\end{rem}

\section{Incomplete list of Notation and Definitions}\label{s-definitions}

Most of the notation is either standard, or will be defined as it is introduced. In this section we record some of the definitions so that the reader can refer back to it if necessary. (Much more comprehensive background on structures can be found in \cite{hodges-model} or \cite{hodges-shorter}.)

\begin{defn} A (single sorted, first order) \emph{structure} $\fA$ is a tuple $(A, \{f_i\}, \{R_j\})$ such that $A$ is a set (called the \emph{underlying set} of $\fA$), each $f_i$ is a function of some arity $n_i$ on $A$, i.e. $f_i : A^{n_i} \rightarrow A$, and each $R_j$ is a relation of some arity $m_j$ on $A$, i.e. $R_j \subseteq A^{m_j}$. We sometimes write $\arity(f_i)$ for $n_i$ and $\arity(R_j)$ for $m_j$.

The \emph{signature} of a first order structure is the assignment of each function symbol $f_i$ to an arity $n_i$ together with an assignment of each relation symbol $R_j$ to some arity $m_j$. If two structures $\fA, \fB$ share the same signature, then we sometimes write $f_i^\fA$, $R_j^\fA$ to refer to the interpretations of the function and relation symbols in $\fA$, and $f_i^\fB, R_j^\fB$ to refer to the interpretations of the same function and relation symbols in $\fB$.

If $\fA, \fB$ are structures with the same signature, then a \emph{homomorphism} $\varphi : \fA \rightarrow \fB$ is a map $\varphi : A \rightarrow B$ of the underlying sets, such that for each function symbol $f_i$ we have
\[
\varphi(f_i^\fA(a_1, ..., a_{n_i})) = f_i^\fB(\varphi(a_1), ..., \varphi(a_{n_i}))
\]
for all $a_1, ..., a_{n_i} \in A$, and such that for each relation symbol $R_j$ we have
\[
(a_1, ..., a_{m_j}) \in R_j^\fA \;\; \implies \;\; (\varphi(a_1), ..., \varphi(a_{m_j})) \in R_j^\fB
\]
for all $a_1, ..., a_{m_j} \in A$. A homomorphism $\varphi : \fA \rightarrow \fB$ is called an \emph{isomorphism} if there is a homomorphism $\psi : \fB \rightarrow \fA$ such that $\psi \circ \varphi = \operatorname{id}_{\fA}$ and $\varphi \circ \psi = \operatorname{id}_\fB$.

A subset $B \subseteq A$ is called a \emph{subuniverse} of the structure $\fA = (A, \{f_i\}, \{R_j\})$ if $B$ is closed under each function $f_i$. If $B$ is a subuniverse of $\fA$, then the corresponding \emph{induced substructure} $\fB$ is defined to be $(B, \{f_i|_B\}, \{R_j \cap B^{\arity(R_j)}\})$ - note that the inclusion map $B \hookrightarrow A$ automatically defines a homomorphism $\fB \hookrightarrow \fA$.

A structure is called a \emph{relational structure} if it has no functions in its signature, and a structure is called an \emph{algebraic structure} or an \emph{algebra} if it is nonempty and has no relations in its signature. We usually write a relational structure with the bold font, i.e. $\fA$, while we usually write an algebraic structure with the blackboard bold font, i.e. $\bA$ (note that many authors reverse this convention).

We define the \emph{total size} of a relational structure $\fA = (A, \{R_j\})$, written $\|\fA\|$, to be
\[
\|\fA\| \coloneqq \sum_j m_j|R_j|.
\]
Note that the number of bits needed to describe $\fA$ is larger than $\|\fA\|$ by a factor of about $\log_2 |A|$, ignoring the overhead needed to describe the signature of $\fA$.
\end{defn}

When convenient, we often abuse notation to treat a structure $\fA$ like its underlying set $A$: we write $a \in \fA$ to mean that $a \in A$, we write $|\fA|$ for $|A|$, etc.

We will also need to consider a certain restricted type of multisorted relational structures.

\begin{defn} A \emph{multisorted relational structure} $\fA$ is defined to be a tuple $(\{A_i\}_{i \in I}, \{R_j\}_{j \in J})$, where for each $j \in J$ there is a tuple $i^j = (i^j_1, ..., i^j_{m_j}) \in I^{m_j}$ such that
\[
R_j \subseteq A_{i^j_1} \times \cdots \times A_{i^j_{m_j}}.
\]
The sets $A_i$ are called the \emph{sorts} of the structure $\fA$, and the $R_j$s are the \emph{multisorted relations} of $\fA$.

The \emph{signature} of a multisorted relational structure $\fA$ consists of the set $I$ of indices of sorts of $\fA$, together with the map $J \mapsto \bigcup_{m \ge 0} (\{m\} \times I^m)$ which sends each relation symbol to its arity $m_j$ and tuple $i^j \in I^{m_j}$ of sorts (which is sometimes called the \emph{type} of the relation symbol).

If $\fA, \fB$ are multisorted relational structures of the same signature, then a \emph{homomorphism} $\varphi : \fA \rightarrow \fB$ is a collection of maps $\varphi_i : A_i \rightarrow B_i$, such that for each relation symbol $R_j$ we have
\[
(a_1, ..., a_{m_j}) \in R_j^\fA \;\; \implies \;\; (\varphi_{i^j_1}(a_1), ..., \varphi_{i^j_{m_j}}(a_{m_j})) \in R_j^\fB. 
\]
\end{defn}

Note that it is possible to make a much more general definition of a multisorted first order structure: for instance, it would be natural to consider a vector space as a multisorted structure with scalars as one sort of element and vectors as another sort of element, along with a binary multiplication operation which has a signature specifying that multiplication must be defined for pairs of scalars (producing a scalar) as well as for scalars and vectors (producing a vector) but not for vectors and vectors (unless one of the vectors happens to also be a scalar). We won't need such a general concept in these notes: multisorted relational structures and collections of algebras sharing a common signature will suffice for most of what we do.

\begin{defn} A \emph{subalgebra} of an algebra $\bA$ is an algebraic structure $\bB$ with the same signature such that the underlying set $B$ of $\bB$ is a subset of the underlying set $A$ of $\bA$, and such that the inclusion map $\iota : B \hookrightarrow A$ defines a homomorphism $\iota : \bB \rightarrow \bA$ (in which case $B$ is a subuniverse of $\bA$, and $\bB$ is the induced substructure). In this case we write $\bB \le \bA$. If $S \subseteq A$, then we define $\Sg_\bA(S)$, the \emph{subalgebra generated by} $S$, to be the smallest subalgebra of $\bA$ whose underlying set contains $S$.

If $\bA, \bB$ are structures with the same signature, then we define their \emph{product} $\bA \times \bB$ to be the structure with underlying set $A\times B$ where $A,B$ are the underlying sets of $\bA,\bB$, with each function symbol $f_i$ interpreted by
\[
f_i^{\bA\times\bB}\Big(\begin{bmatrix}a_1\\b_1\end{bmatrix},...,\begin{bmatrix}a_{n_i}\\b_{n_i}\end{bmatrix}\Big) = \begin{bmatrix}f_i^\bA(a_1, ..., a_{n_i})\\f_i^\bB(b_1, ..., b_{n_i})\end{bmatrix},
\]
and with each relation symbol $R_j$ interpreted by
\[
\Big(\begin{bmatrix}a_1\\b_1\end{bmatrix},...,\begin{bmatrix}a_{m_j}\\b_{m_j}\end{bmatrix}\Big) \in R_j^{\bA\times\bB} \;\; \iff \;\; (a_1, ..., a_{m_j}) \in R_j^\bA \wedge (b_1, ..., b_{m_j}) \in R_j^\bB.
\]
Arbitrarily large products $\prod_{i \in I} \bA_i$ of structures $\bA_i$ all having the same signature are defined similarly. (This definition matches with the category-theoretic definition of the product, in the category of structures with a fixed signature.)

A \emph{homomorphic image} of $\bA$ is defined to be any $\bB$ of the same signature such that there is a surjective homomorphism $\varphi : \bA \twoheadrightarrow \bB$. An algebraic structure $\bA$ is called \emph{simple} if every surjective homomorphism $\bA \twoheadrightarrow \bB$ is either an isomorphism, or has $|\bB| = 1$.
\end{defn}

\begin{defn} A \emph{congruence} on an algebraic structure $\bA$ is an equivalence relation $\theta$ on $\bA$ which is also a subalgebra: $\theta \le \bA^2$. The set of all congruences of $\bA$ is written as $\Con(\bA)$. If $S$ is some collection of ordered pairs $(a,b) \in \bA^2$, then the \emph{congruence generated by} $S$, written $\Cg_\bA(S)$, is the least congruence of $\bA$ which contains $S$. We generally use greek letters for congruences.

If $\theta \in \Con(\bA)$ and $a \in \bA$, then we write $a/\theta$ for the \emph{congruence class} of $a$ with respect to $\theta$, that is, $a/\theta = \{b \mid (a,b) \in \theta\}$. We say that $a,b$ are \emph{congruent} with respect to $\theta$ if $(a,b) \in \theta$, and we may write this in symbols in several different ways:
\[
(a,b) \in \theta \;\; \iff \;\; b \in a/\theta \;\; \iff \;\; a \equiv b \pmod{\theta} \;\; \iff \;\; a \equiv_\theta b \;\; \iff \;\; a\ \theta\ b.
\]

If $\theta \in \Con(\bA)$, then we write $\bA/\theta$ for the set of equivalence classes of $\theta$ considered as an algebraic structure, with
\[
f_i^{\bA/\theta}(a_1/\theta, ..., a_{n_i}/\theta) = f_i^\bA(a_1, ..., a_{n_i})/\theta.
\]
The algebra $\bA/\theta$ is called a \emph{quotient} of $\bA$, and there is a canonical \emph{quotient map} $\bA \twoheadrightarrow \bA/\theta$ which takes each element of $\bA$ to its equivalence class in $\theta$.

If $\varphi : \bA \rightarrow \bB$ is a homomorphism, then we define the \emph{kernel} of $\varphi$, written $\ker \varphi$, as the equivalence relation on $\bA$ given by $(x,y) \in \ker \varphi$ iff $\varphi(x) = \varphi(y)$. The kernel of a homomorphism is always a congruence, and if $\varphi$ is surjective then $\bA/\ker \varphi$ is isomorphic to $\bB$. In particular, every homomorphic image of $\bA$ is isomorphic to a quotient of $\bA$.

For $\alpha, \beta \in \Con(\bA)$, we define their \emph{meet} $\alpha \wedge \beta$ to be their intersection, and we define their \emph{join} $\alpha \vee \beta$ to be the congruence generated by $\alpha \cup \beta$. We define the least congruence $0_\bA$ to be the equivalence relation on $\bA$ where each equivalence class is a singleton, and we define the greatest congruence $1_\bA$ to be the equivalence relation on $\bA$ consisting of just one equivalence class.
\end{defn}

\begin{defn} If $R \subseteq A_1 \times \cdots \times A_n$ is a (possibly multisorted) relation, then for any $I \subseteq \{1, ..., n\}$ we define the projection map $\pi_I : R \rightarrow \prod_{i \in I} A_i$ in the obvious way. A relation $R \subseteq A_1 \times \cdots \times A_n$ is \emph{subdirect} if $\pi_i(R) = A_i$ for all $i$.

If $\bA_1, ..., \bA_n$ are algebraic structures with the same signature, then any subalgebra $\RR \le \bA_1 \times \cdots \times \bA_n$ is called a (multisorted) \emph{relation} on the $\bA_i$s. If $\pi_i(\RR) = \bA_i$ for all $i$ then we say that $\RR$ is subdirect, and we write this in symbols as $\RR \le_{sd} \bA_1 \times \cdots \times \bA_n$.
\end{defn}

\begin{defn} A \emph{tolerance} $\bS$ on an algebraic structure $\bA$ is a binary relation $\bS \le \bA \times \bA$ which is symmetric and contains the diagonal. Note that the transitive closure of any tolerance is automatically a congruence. We say that a tolerance is \emph{connected} if its transitive closure is the full congruence $1_\bA$.

If $\RR \le_{sd} \bA_1 \times \cdots \times \bA_n$, then we define the $i$th \emph{link tolerance} of $\RR$ to be the binary relation
\[
\{(b,c) \mid \exists a_j \text{ s.t. } (a_1, ..., a_{i-1}, b, a_{i+1}, ..., a_n), (a_1, ..., a_{i-1}, c, a_{i+1}, ..., a_n) \in \RR\}
\]
on $\bA_i$. A binary relation is called \emph{linked} if its link tolerances are connected.
\end{defn}

\begin{defn}
If $R \subseteq A \times B$ is a binary relation, then we define its \emph{reverse} $R^- \subseteq B\times A$ by
\[
R^- = \{(b,a) \mid (a,b) \in R\}.
\]
If $R \subseteq A \times B$ and $S \subseteq B \times C$, then we define their \emph{relational composition} $R\circ S \subseteq A\times C$ by
\[
R \circ S = \{(a,c) \mid \exists b\in B\text{ s.t. }(a,b) \in R \wedge (b,c) \in S\}.
\]
If $R \subseteq A \times B$ and $U \subseteq A$, then we define the \emph{sum} $U+R \subseteq B$ by
\[
U + R = \{b \mid \exists a \in U\text{ s.t. }(a,b)\in R\},
\]
and for $V \subseteq B$ we define the difference $V - R \subseteq A$ by $V - R = V + R^-$. Note that we have $(U + R) + S = U + (R \circ S)$ for $U \subseteq A$, $R \subseteq A \times B$, $S \subseteq B \times C$.
\end{defn}

\begin{defn} If $\{R_j\}$ is some collection of relations, then a \emph{primitive positive formula} over the $R_j$s is defined to be a formula of the form
\[
\exists y_1 \cdots \exists y_n\text{ s.t. }\bigwedge_{i \in I}\varphi_i(x_1, ..., x_m, y_1, ..., y_n),
\]
where $I$ is a finite set, and each $\varphi_i$ is either some relation $R_j$ applied to some of the variables $x_1, ..., x_m, y_1, ..., y_n$, or is the equality relation applied to some pair of variables.

A relation $R$ is \emph{primitively positively definable} over the $R_j$s if there is a primitive positive formula $\varphi$ such that the elements of $R$ are exactly the tuples of values $(x_1, ..., x_m)$ that satisfy $\varphi$. For instance, the relational composition $R \circ S$ is primitively positively definable over $R$ and $S$.

A \emph{relational clone} is a collection of relations (on a common domain) which is closed under primitive positive definitions. The smallest relational clone which contains $\{R_i\}$ is written as $\langle\{R_i\}\rangle$.
\end{defn}

The algebraic analogue of primitive positive formulas is the concept of \emph{terms}.

\begin{defn} If $\{f_i\}$ is a collection of function symbols in some fixed signature, then we define a \emph{term} inductively as follows:
\begin{itemize}
\item each variable $x_j$ is a term, corresponding to the operation $\pi_j$, and

\item if $f_i$ is function symbol of arity $k$ and $t_1, ..., t_k$ are terms, then $f_i(t_1, ..., t_k)$ is also a term, corresponding to the operation $f_i \circ (t_1, ..., t_k)$.
\end{itemize}
For instance, if $f$ is binary and $g$ is ternary, then the expression
\[
f(g(x,y,f(x,z)), f(u,v))
\]
is a term. Every term can be visualized as a labeled ordered tree, where every leaf is labeled by a variable and where every internal vertex is labeled by a function whose arity is equal to the number of children of that vertex. The \emph{height} of a term is the largest distance between the root of this tree and any leaf.

In some cases, it may be more efficient to visualize terms via directed acyclic graphs, to avoid repeatedly drawing many copies of a common subtree - to distinguish these points of view, a tree representation of a term is called a \emph{formula}, while a directed acyclic graph representation of a term is called a \emph{circuit}. In written math, a circuit representation of a term corresponds to a sequence of definitions of subterms in terms of function symbols applied to previously defined subterms. For instance, the formula
\[
t(x,y) \coloneqq f(f(f(f(x,y), f(y,x)), f(f(y,x), f(x,y))), f(f(f(y,x), f(x,y)), f(f(x,y), f(y,x))))
\]
corresponds to the circuit
\begin{align*}
g(x,y) &\coloneqq f(f(x,y), f(y,x)),\\
g(y,x) &\coloneqq f(f(y,x), f(x,y)),\\
h(x,y) &\coloneqq f(g(x,y), g(y,x)),\\
h(y,x) &\coloneqq f(g(y,x), g(x,y)),\\
t(x,y) &\coloneqq f(h(x,y), h(y,x)),
\end{align*}
where we have included the redundant definitions of $g(y,x)$ and $h(y,x)$ since they would correspond to additional nodes in a computation graph used to compute $t$.

The collection of all terms in a given signature $\sigma$ defines the \emph{term algebra} $\cF_\sigma(\{x_j\})$, with each function symbol $f_i$ of arity $k$ interpreted as the operation
\[
f_i^{\cF_{\sigma}(\{x_j\})} : (t_1, ..., t_k) \mapsto f_i(t_1, ..., t_k),
\]
where the right hand side is interpreted as an abstract term. The term algebra is also called the \emph{absolutely free algebra} in the signature $\sigma$.
\end{defn}

\begin{defn} An \emph{identity} is just a pair of terms $s, t$ with the symbol $\approx$ in between them:
\[
s \approx t.
\]
The identity $s\approx t$ has \emph{height 1} if both terms $s$ and $t$ have height 1, and the identity $s \approx t$ is called \emph{linear} if both $s$ and $t$ have height at most 1. Height 1 identities are also called \emph{minor identities}.

An algebraic structure \emph{satisfies} the identity $s \approx t$, written $\bA \models s \approx t$, if
\[
\forall x_1, ..., x_k \in \bA, \;\; s^\bA(x_1, ..., x_k) = t^\bA(x_1, ..., x_k),
\]
where $x_1, ..., x_k$ is a list of all of the variables which occur in $s$ or $t$.
\end{defn}

\begin{defn} If $\cT$ is a set of identities in the signature $\sigma$, then we define the \emph{variety} $\cV(\cT)$ to be the collection of all algebraic structures $\bA$ with signature $\sigma$ such that $\bA \models \cT$.

If $\{\bA_i\}$ is a collection of algebraic structures in the signature $\sigma$, then the \emph{variety generated by} $\{\bA_i\}$, written $\cV(\{\bA_i\})$, is the variety $\cV(\cT)$ where $\cT$ is the collection of all identities $s \approx t$ which are satisfied in every single $\bA_i$.
\end{defn}

\begin{defn} If $\cV = \cV(\cT)$ is a variety with signature $\sigma$ and defining identities $\cT$, then we define the congruence $\approx_\cV$ on the absolutely free algebra $\cF_\sigma(\{x_j\})$ to be the congruence generated by the set of pairs of terms
\[
(s \circ (u_1, ..., u_k), t \circ (u_1, ..., u_k))
\]
such that $s \approx t \in \cT$, where $u_1, ..., u_k$ is an arbitrary list of terms which we use to replace the variables $x_1, ..., x_k$ which occur in $s$ and $t$. The algebraic structure
\[
\cF_\sigma(\{x_j\}) /\!\! \approx_\cV
\]
is called the \emph{free algebra on the generators} $\{x_j\}$ in the variety $\cV$, and is written as $\cF_\cV(\{x_j\})$. By construction, the free algebra $\cF_\cV(\{x_j\})$ satisfies every identity in $\cT$. In the language of category theory, the functor
\[
\cF_\cV : S \mapsto \cF_{\cV}(\{x_j\}_{j \in S})
\]
is adjoint to the forgetful functor from algebras in $\cV$ to their underlying sets:
\[
\Hom_{\text{Set}}(S, A) = \Hom_{\cV}(\cF_\cV(S), \bA)
\]
when $\bA \in \cV$ is an algebraic structure with underlying set $A$.

We define a \emph{term operation} of a variety $\cV$ to be an equivalence class $t/\!\!\approx_{\cV}$ together with an ordered sequence of variables containing all variables which occur in $t$. If $\bA$ is an algebraic structure, then we define the \emph{term operations} of $\bA$ to be the collection of interpretations $t^\bA$ of terms $t$ as operations on $\bA$ - these are easily seen to correspond with the term operations of the variety $\cV(\bA)$ generated by $\bA$.

Viewing each $k$-ary term operation $t^\bA : \bA^k \rightarrow \bA$ as an element of $\bA^{\bA^k}$, we get an alternative construction of the free algebra $\cF_{\cV(\bA)}(x_1, ..., x_k)$:
\begin{align*}
\cF_{\cV(\bA)}(x_1, ..., x_k) &= \{k\text{-ary term operations }t\text{ of }\cV(\bA)\}\\
&\cong \{k\text{-ary term operations }t^\bA\text{ of }\bA\}\\
&\cong \Sg_{\bA^{\bA^k}}\{\pi_1, ..., \pi_k\}.
\end{align*}
More generally, if $\cV = \cV(\{\bA_i\}_{i \in I})$, then $\cF_\cV(x_1, ..., x_k)$ is isomorphic to a subalgebra of $\prod_{i \in I} \bA_i^{\bA_i^k}$.
\end{defn}

\begin{defn} If $\bA$ is an algebraic structure, then the \emph{clone} of $\bA$, written as $\Clo(\bA)$, is the collection of all term operations of $\bA$. If $\{f_i\}$ is a collection of operations on the domain $A$, then we write $\langle \{f_i\} \rangle$ for the clone of the algebraic structure $(A, \{f_i\})$. Alternatively, a \emph{clone} on a domain $A$ is just a collection of operations of $A$ which is closed under composition.

If $\cV$ is a variety, then the \emph{clone} of $\cV$, written as $\Clo(\cV)$, is defined to be the collection of free algebras
\[
\cF_\cV(x_1, ..., x_k),
\]
together with the rules for composing a $k$-ary term operation $t \in \cF_\cV(x_1, ..., x_k)$ with $k$-tuples of $m$-ary term operations $u_1, ..., u_k \in \cF_\cV(x_1, ..., x_m)$ to produce the $m$-ary term operation $t \circ (u_1, ..., u_k)$.

An \emph{abstract clone} is a collection of function symbols $\{f_i\}$ in a signature $\sigma$ which is closed under a composition law $\circ$ which takes a function symbol $f$ of arity $k$ and a $k$-tuple of function symbols $g_1, ..., g_k$ of arity $m$ as input and produces a function symbol $f\circ (g_1, ..., g_k)$ of arity $m$ as output, satisfying a generalized associativity law:
\[
(f\circ(g_1, ..., g_k)) \circ (h_1, ..., h_m) = f \circ (g_1 \circ (h_1, ..., h_m), ..., g_k \circ (h_1, ..., h_m)),
\]
with special ``projection'' function symbols $\pi_i^k$ of every arity $k$ which satisfy
\[
\pi_i^k \circ (g_1, ..., g_k) = g_i
\]
and
\[
f \circ (\pi_1^k, ..., \pi_k^k) = f.
\]
Note that in any abstract clone, the collection of function symbols of arity $1$ always forms a semigroup under $\circ$ with identity element $\pi_1^1$.

A \emph{clone homomorphism} is a map $\xi$ from the function symbols of one abstract clone to the function symbols of another abstract clone which preserves arities, sends the projections $\pi_i^k$ to themselves, and respects composition:
\[
\xi(f \circ (g_1, ..., g_k)) = \xi(f) \circ (\xi(g_1), ..., \xi(g_k)).
\]
A \emph{height 1 clone homomorphism}, also known as a \emph{minion homomorphism}, is a map $\xi$ which preserves arities and respects composition with projections:
\[
\xi(f \circ (\pi_{i_1}^m, ..., \pi_{i_k}^m)) = \xi(f) \circ (\pi_{i_1}^m, ..., \pi_{i_k}^m).
\]
\end{defn}

\begin{defn} A \emph{poset} is a relational structure $(P, \le)$ such that $\le$ is a \emph{partial order} - that is, a reflexive and transitive relation which satisfies $a \le b \wedge b \le a \implies a = b$. More generally, a \emph{quasiorder} (also called a \emph{preorder}) is any reflexive and transitive relation. We usually use $\preceq$ to denote a quasiorder and $\le$ to denote a partial order. If $\preceq$ is a quasiorder, then we define an \emph{associated equivalence relation} $\sim$ by
\[
a \sim b \;\; \iff \;\; a \preceq b \wedge b \preceq a.
\]

For any $a \le b \in P$, we define the \emph{interval} between $a$ and $b$, written $\llbracket a, b \rrbracket$, to be the set of $c \in P$ such that $a \le c \le b$. (We use this notation to avoid confusion with the commutator, which is written as $[\cdot, \cdot]$.)

We say that $b$ is a \emph{cover} of $a$, written $a \prec b$ (if there is no danger of confusion with a strict quasiorder), if $a < b$ and there is no $c \in P$ such that $a < c < b$, that is, if $\llbracket a, b \rrbracket = \{a,b\}$. If $a$ has exactly one cover, then we will write $a^*$ to refer to the unique cover of $a$.
\end{defn}

\begin{defn} A \emph{lattice} is either an algebraic structure $\cL = (L, \wedge, \vee)$ which satisfies the identities
\begin{align*}
&x \wedge x \approx x, &x \vee x \approx x,\\
&x \wedge y \approx y \wedge x, &x \vee y \approx y \vee x,\\
&x \wedge (y \wedge z) \approx (x \wedge y) \wedge z, &x \vee (y \vee z) \approx (x \vee y) \vee z,\\
&x \wedge (x \vee y) \approx x, &x \vee (x \wedge y) \approx x,
\end{align*}
or it is a poset $\cL = (L, \le)$ such that every pair of elements $\{a,b\} \subseteq \cL$ has a least upper bound $a \vee b$ and a greatest lower bound $a \wedge b$, or it is a first-order structure $\cL = (L, \wedge, \vee, \le)$ which satisfies
\[
a \le b \iff a = a \wedge b \iff a \vee b = b
\]
as well as the algebraic identities above. Note that the operations $\wedge, \vee$ are determined by $\le$, and the partial order $\le$ is determined by either of the operations $\wedge$ or $\vee$. A \emph{$0,1$-lattice} is a lattice with named constants $0, 1$ which are respectively the least and greatest elements of $\cL$.

A \emph{lattice homomorphism} is a homomorphism of the algebraic structure $(L, \wedge, \vee)$. Note that there may be some homomorphisms of the relational structure $(L, \le)$ which do not count as lattice homomorphisms - a map which respects the partial order $\le$ is just called a \emph{monotone} map.
\end{defn}

\section{Crash course on NP-completeness}

The class of problems NP can be viewed as a formalization of the intuitive concept of a puzzle. Informally, a problem is a \emph{puzzle} if it is easy to tell whether a guessed solution actually solves it or not. We formalize this by requiring that there is a mechanical, algorithmic procedure which takes as input a purported solution to the puzzle, and efficiently (i.e. in polynomial time) tells us whether or not the input counts as an actual solution. Before getting too deep into the formalization, let's go through a few examples of problems in NP.

\begin{ex} The \emph{subset sum problem}: given a set $U$ of integers (expressed in binary) and a target sum $t$ (also expressed in binary), determine whether or not there is a subset $S$ of $U$ whose elements sum to exactly $t$. To see that this is in NP, note that if you want to convince me that the answer is Yes you can simply tell me what the set $S$ is, and then I can easily check that $S$ is a subset of $U$ and that the sum of the elements of $S$ is $t$.
\end{ex}

\begin{ex} The \emph{graph isomorphism problem}: given two graphs $\fG_1, \fG_2$, determine whether or not they are isomorphic. To see that this is in NP, note that if you want to convince me that the answer is Yes you can simply tell me what the isomorphism $f : \fG_1 \rightarrow \fG_2$ is, and then I can easily verify that $f$ is indeed an isomorphism.
\end{ex}

\begin{ex} The \emph{boolean satisfiability problem}: given a formula of propositional logic, determine whether or not it is possible to assign True/False values to the propositional variables to make the formula evaluate to True. To see that this is in NP, note that if you want to convince me that the answer is Yes you can simply tell me an assignment of True/False values to the propositional variables, and then I can easily evaluate the formula and check that the result is True.
\end{ex}

\begin{ex} An example of a problem not believed to be in NP is determining whether or not a given player has a winning strategy in a complex board game such as Chess or Go (generalized to have an arbitrarily large board). The reason for this is that there is no obvious way for you to compactly express a winning strategy into something I could understand in a reasonable amount of time, and even if you somehow did have a compact way to describe your winning strategy, there is still no obvious way for me to efficiently test that the other player can't find a clever way to defeat it.
\end{ex}

For the sake of formalizing the concept of an efficient, mechanical, and algorithmic procedure, first note that for most purposes we can ignore the internals of a computational machine $M$ and just think about the function $f_M$ implemented by it, which takes a problem description as input and produces an answer as output:
\[
f_M : \text{Questions} \rightarrow \{\text{Yes}, \text{No}\}.
\]
Of course, not every function $f$ corresponds to a physical machine $M$ which can be built and which runs efficiently. The set of functions (with Yes/No outputs) which can be implemented by \emph{Turing machines} which finish their computation after a number of steps which is at most polynomial in the number of bits needed to describe the input is called ``P'' (an abbreviation for \emph{polynomial time}). Formally, a Turing machine can be defined as follows (although the details won't be particularly important to us).

\begin{defn}[Turing \cite{turing-computable}] A (vanilla) \emph{Turing machine} consists of the following:
\begin{itemize}
\item a finite set of \emph{symbols} $\Sigma$, one of which is singled out as the \emph{blank} symbol,
\item a finite set of \emph{states} $S$, one of which is singled out as the \emph{initial} state $s_0$ and another of which is singled out as the \emph{halt} state,
\item a \emph{transition function} $\delta$, which takes as input a pair consisting of a current symbol and a current state (other than the halt state), and produces as output a triple consisting of a new symbol and a new state (to overwrite the old symbol and state) and an instruction to move either one step to the left or one step to the right,
\item an infinite \emph{tape} storing a sequence of symbols, which can be modeled as a function $\mathbb{N} \rightarrow \Sigma$, such that all but finitely many of the symbols stored in the tape are blank, together with a \emph{head} $p \in \mathbb{N}$ pointing to a particular location in the tape.
\end{itemize}
The input to the Turing machine is encoded into the tape as a finite sequence of non-blank symbols (followed by an infinite sequence of blank symbols). When the Turing machine starts running, the head is placed at the beginning of the tape (i.e. $p = 0$) and the state $s$ is set to the initial state (i.e. $s = s_0$). Then in each step of operation, the Turing machine updates the tape as well as its current state as follows:
\begin{itemize}
\item first the machine reads the symbol $a \in \Sigma$ stored on the tape in the location pointed to by $p$,
\item then the machine consults the transition function $\delta$ for the pair $(a,s)$, where $s \in S$ is the current state, and obtains the triple $(a',s',\pm 1) = \delta(a,s)$,
\item then the machine overwrites the symbol $a$ at location $p$ of the tape with the symbol $a'$, changes its current state to $s'$, and increases or decreases the location $p$ by $\pm 1$ based on the third output of $\delta(a,s)$. (If you try to decrease $p$ when $p$ is $0$, then the machine either crashes or silently fails to change the value of $p$, depending on the author.)
\end{itemize}
If at any point the Turing machine reaches the halt state, then it immediately stops running, and at that point the contents of the tape are interpreted as the output of the machine.

There are many variations on the vanilla Turing machine described above:
\begin{itemize}
\item there can be multiple halting states, for instance to record whether the program \emph{succeeded} or \emph{failed},
\item the tape can be infinite in both directions, so that it is modeled by a function $\mathbb{Z} \rightarrow \Sigma$,
\item there can be multiple heads pointing to unrelated positions in the tape - in this case, the transition function $\delta$ will take as input the ordered tuple of symbols pointed to by the sequence of different heads (as well as the current state), and returns instructions for how to change the values pointed to by all of the heads simultaneously (as well as which direction each head should move in, and how the state should update), and
\item there can be multiple tapes (with one or more heads per tape), some of which might be read-only and some of which might be write-only (in order to cleanly separate the input from the output, and to keep track of the amount of scratch space used by the machine throughout the computation).
\end{itemize}
\end{defn}

\begin{exer}[Two heads are better than one] Consider the task of reversing a finite sequence of bits (i.e., if the input is 00101, then the desired output is 10100). Prove that any vanilla Turing machine which correctly executes this task must run for at least $\frac{\lfloor n^2/2 \rfloor}{\log_2(|S|)} - 2(n-1)$ steps on some input sequence of length $n$ (hint: imagine placing a tollbooth between two consecutive locations on the memory tape, and find a lower bound on the amount of information that must travel through the tollbooth). Find a \emph{two-headed} Turing machine which can accomplish the same task using just $2(n+1)$ steps, where $n$ is the length of the input.
\end{exer}

\begin{rem} In these notes, when we make precise claims about running times of algorithms, we implicitly use a slightly more powerful model of computation than vanilla Turing machines, where we also allow the manipulation and storage of \emph{pointers} (i.e. numbers describing locations on the memory tape) which can be used to instantly jump to various locations of the memory tape. This is essentially the \emph{random access machine} model of computation, where the cost of pointer arithmetic is ignored. The \emph{word RAM} model (which we will not define precisely here) formalizes what we mean by ignoring the cost of pointer arithmetic - the main thing to know about it is that it prevents us from abusing pointer arithmetic too much by requiring pointers to be described by \emph{machine words}, which are required to be representable with only logarithmically many bits. If the word RAM model is taken too seriously, then it has some counter-intuitive consequences: for instance, it allows you to sort a list of $n$ integers in time $O(n\log(\log(n)))$ \cite{sorting-word-RAM-deterministic}, even though the number of bits needed to describe an ordering on $n$ numbers is $\log_2(n!) \gg n \log(n)$.
\end{rem}

The choice to define the class P in terms of Turing machines is due to the empirical observation (known as the \emph{extended Church-Turing thesis}) that whenever a serious attempt is made to formalize the intuitive concept of an ``algorithmic procedure'' from first principles (e.g. total recursive functions, cellular automata, lambda calculus, or random-access machines), the resulting type of machine can be simulated by a Turing machine with at most a polynomial slowdown to the running time, and conversely can simulate a Turing machine (again with at most a polynomial slowdown in the running time). The fact that running time can vary from one concept of ``algorithmic procedure'' to the next, but only by a polynomial amount, leads naturally to singling out the polynomial-time computable functions as worthy of special consideration.

\begin{rem} The intuitive concept of an ``algorithmic procedure'' is meant to capture what a trained human mathematician is capable of mechanically computing using a pencil and paper. It is possible to imagine that physical devices could exist which can reliably answer questions that no human could ever resolve with pencil and paper, regardless of the amount of time they devote to the problem - such a hypothetical device is known as a \emph{hypercomputer}. As an example, one model of a hypercomputer is a Turing machine simulator that (somehow) executes each step of its simulation twice as quickly as the previous step until it halts, allowing us to determine if a given Turing machine will ever halt on a given input - the function this would implement is known as a \emph{halting oracle}, and a well-known self-reference argument shows that no Turing machine which always halts can implement a halting oracle (if such a Turing machine existed, then we could use it to construct a Turing machine which halts iff it does not halt).

Our current understanding of physics supports the \emph{physical Church-Turing thesis}, which can be summarized as the belief that hypercomputers (and, in particular, halting oracles) can't be built in the real world. More precisely, the physical Church-Turing thesis states that every physical process can be simulated by a probabilistic Turing machine (that is, a Turing machine augmented with an extra tape containing an infinite sequence of random bits), and that the statistics of every finite physical process can be computed to any desired accuracy with a deterministic Turing machine. Note that if the physical Church-Turing thesis turns out to be false, then the job of the theoretical physicist (i.e. making predictions about the statistics of lab experiments based on pencil and paper computations) would appear to be fundamentally impossible.

The \emph{extended physical Church-Turing thesis}, also known as the \emph{classical complexity-theoretic Church-Turing thesis}, is the belief that physical processes can be simulated by (probabilistic) Turing machines \emph{efficiently}, i.e. with at most polynomial slowdown. The collection of functions which can be reliably and efficiently computed by a probabilistic Turing machine is known as ``BPP'' (an abbreviation for \emph{bounded-error probabilistic polynomial time}) - it is defined as the collection of functions (with Yes/No outputs) such that there exists a probabilistic Turing machine which produces the correct output at least $2/3$ of the time (any fixed fraction $> 1/2$ could be used instead of $2/3$, since the success probability can be amplified by repeatedly executing the probabilistic Turing machine and taking a majority vote of its outputs). A widely believed conjecture states that BPP = P, i.e. that randomness does not fundamentally expand the set of efficiently computable functions - the intuition is that we appear to be able to design cryptographically secure pseudorandom number generators which can't be efficiently distinguished from true randomness without knowing their seeds (this intuition is made precise in a paper of Ipagliazzo and Wigderson \cite{P-vs-BPP-exponential-circuits} which shows that BPP = P follows from the existence of problems which can be solved in exponential time and which require exponentially large circuits for their solution). An example of a problem which is in BPP but which hasn't yet been proven to be in P is the \emph{polynomial identity testing} problem: given a multivariate polynomial over a finite field, expressed via a circuit built out of addition and multiplication gates, determine whether the polynomial is equal to $0$ as a polynomial (this can be solved by plugging in random inputs from a big enough extension field and checking if we get $0$ when we evaluate the circuit).

Surprisingly, our current understanding of quantum physics (and the level of difficulty of \emph{ab initio} quantum chemistry) appears to imply that the classical complexity-theoretic physical Church-Turing thesis is \emph{false}. Instead, it should be replaced with the \emph{quantum complexity-theoretic Church-Turing thesis}, which states that any physical process can be efficiently simulated by a ``quantum Turing machine'' (see \cite{deutsch-quantum-turing-machine} and \cite{quantum-computation} for details on how this can be defined). The collection of functions (with Yes/No outputs) which can be reliably and efficiently computed by a quantum Turing machine is known as ``BQP'' (an abbreviation for \emph{bounded-error quantum polynomial time}). BQP is believed to be strictly larger than P, but not ``too much larger'' - for instance, all constraint satisfaction problems which are known to be contained in BQP are also known to be contained in P. The main example of a problem (\emph{not} a CSP) which is known to be contained in BQP but which is conjectured not to be in P is \emph{factoring} (we phrase factoring as a Yes/No decision problem as follows: given a large integer $n$ written out in binary, together with integers $u, \ell$ also written out in binary, determine whether or not $n$ has a prime factor $p$ satisfying $\ell \le p \le u$).
\end{rem}

Now that we have formalized the concept of an efficient algorithmic decision procedure as the class P of functions (with Yes/No outputs) which can be implemented by Turing machines which always halt within polynomially many steps, we can finally move on to formally defining NP. The name ``NP'' is an abbreviation for \emph{non-deterministic polynomial time}: the idea is that we have a ``non-deterministic'' variant of an ordinary Turing machine which sometimes has more than one choice about how to proceed, and whenever faced with such a choice it (somehow) always chooses the option which is more likely to result in the final output being Yes.

Instead of making the non-deterministic choices all throughout the computation, we can simulate the operation of a general non-deterministic Turing machine with a two phase process: first we (non-deterministically) write out the perfect sequence of choices onto a separate tape, and then we run an ordinary (multi-tape) Turing machine which consults the next entry of the choice tape every time the original non-deterministic machine would have needed to make a choice. Letting $x$ stand for the original input to the non-deterministic machine and letting $y$ be the perfect sequence of choices, we are lead to the following definition.

\begin{defn} The class NP is defined as the collection of functions $g$ (with Yes/No outputs) which have the form
\[
g(x) = \begin{cases}\text{Yes} & \text{if }\exists y\text{ s.t. } f(x,y) = \text{Yes},\\\text{No} & \text{if }\forall y, f(x,y) = \text{No},\end{cases}
\]
where the length of $y$ is required to be bounded by a polynomial of the length of $x$, and the function $f$ is in P.

The class coNP is defined the same way, but with the roles of ``Yes'' and ``No'' interchanged.
\end{defn}

This leads to our first example of an \emph{NP-complete} problem.

\begin{ex}\label{ex-np-complete} The problem we are about to describe is the ``generic'' NP problem. The input is a triple $(M, x, k)$ where $M$ is a description of a Turing machine, $x$ is a partial input to the Turing machine $M$, and $k$ is a natural number expressed in unary. The desired output is ``Yes'' if there is another partial input $y$ such that when we run the machine $M$ with the pair $(x,y)$ as input, $M$ halts within $k$ steps and outputs ``Yes'', and the desired output is ``No'' if no such $y$ exists.
\end{ex}

The reason we say that this problem is NP-complete is that there is an efficiently computable transformation from any problem in NP into the problem above: if we can solve the problem from Example \ref{ex-np-complete} efficiently, then we can solve any problem in NP efficiently (and, additionally, the problem from Example \ref{ex-np-complete} is clearly in NP).

In fact, the transformation from any problem in NP to the problem from Example \ref{ex-np-complete} is not just efficiently computable - it's computable even in computational models which have extremely limited resources. We might attempt to formalize this by introducing a new complexity class ``L'', for problems which can be solved using only a logarithmically large workspace.

\begin{defn} The class L is defined as the collection of functions (with Yes/No outputs) which can be solved by a two-tape Turing machine, where one tape is a read-only input tape, and the other (readable and writable) tape only has space for $O(\log(n))$ symbols, where $n$ is the number of symbols in the input.
\end{defn}

\begin{rem} Recall that in the word RAM model, each pointer is stored in a machine word consisting of $O(\log(n))$ bits (in order to represent positions in the tape of size $n^{O(1)}$), so a logspace computation can simulate pointer arithmetic with a constant number of pointers - this is pretty much the bare minimum needed to navigate from one arbitrary location in the input tape to another when the input tape is read-only.
\end{rem}

\begin{prop} $\mathrm{L} \subseteq \mathrm{P}$.
\end{prop}
\begin{proof} If a Turing machine halts, then it can never be in the exact same configuration twice during its computation, so we just need to check that the total number of different configurations that a logspace Turing machine can take on during a single computation is bounded by $n^{O(1)}$. The number of possible states is $|S| = O(1)$, the number of possible locations for the read head is $O(n)$, the number of possible locations for the read/write head is $O(\log(n)) = n^{o(1)}$, and the number of possible symbol sequences on the read/write tape is $|\Sigma|^{O(\log(n))} = n^{O(1)}$, and multiplying these all together we see that the total number of configurations is bounded by $n^{O(1)}$.
\end{proof}

However, the purpose of introducing the logspace model of computation was to describe \emph{transformations} from inputs to one problem to inputs to another problem. This leads to the following definition.

\begin{defn} A \emph{logspace transducer} is a function which can be computed by a three-tape Turing machine, where one tape is a read-only input tape, one tape is readable and writeable but only has space for $O(\log(n))$ symbols, and the final tape is a write-only, write-once output tape where the write head can only move in one direction.
\end{defn}

\begin{prop} A function $f$ from finite sequences $x$ of (non-blank) symbols in $\Sigma$ to finite sequences of symbols $f(x)$ can be computed by a logspace transducer if and only if the following function is in L: the input is a triple $(x, i, a)$, where $x$ is a sequence of non-blank symbols, $i$ is a number encoded in binary, $a$ is a symbol, and desired output is ``Yes'' iff the $i$th symbol of $f(x)$ is $a$.
\end{prop}
\begin{proof}[Sketch] To turn a logspace transducer into the second type of machine, the main idea is to simulate the transducer while maintaining an extra counter with $O(\log(n))$ bits at the end of the read/write tape, and to increment the counter every time the transducer attempts to write a symbol to the transducer's write tape. As the counter reaches $i$, we check if the transducer was attempting to write the symbol $a$.

For the other direction, given a Yes/No machine in L, we make a transducer which uses a binary counter to keep track of the number $i$ of symbols it has already output, and then cycles through each possible next symbol $a$ in order, simulating the Yes/No machine on the input $(x,i+1,a)$ to figure out which symbol it needs to output next (and halting if none of the possible next symbols $a$ results in a ``Yes'').
\end{proof}

Both L and P have nice closure properties with respect to subroutines. For P, this can be phrased as follows: if we call a subroutine which runs in time $O(n^k)$ at most $O(n^m)$ times while running our program, then the total time spent within the subroutine will be at most $O(n^{k+m})$ - so we can take any particular subroutine which is already known to run in polynomial time and pretend that it magically runs within a single step, without changing our assessment of whether or not the overall algorithm runs within polynomial time. For L, we phrase the closure property as the ability to compose logspace transducers.

\begin{prop} If $f$ and $g$ are logspace transducers, then so is their composition $f \circ g$. More generally, if $f, g_1, ..., g_m$ are all logspace transducers, then so is the map $x \mapsto f(g_1(x), ..., g_m(x))$, where by ``$(g_1(x), ..., g_m(x))$'' we mean a concatenation of the outputs of $g_1, ..., g_m$ (possibly with a separator symbol between the outputs of $g_i$ and $g_{i+1}$).
\end{prop}
\begin{proof}[Sketch] Given logspace machines computing $f$ and $g$, we need to construct a logspace machine which computes $f \circ g$. We do this by simulating the computation of $f$, while keeping track of a counter $i$ (using $O(\log(n))$ bits, so that the maximum value of $i$ is $n^{O(1)}$) describing the position of $f$'s read head in an imaginary input tape containing the hypothetical output of $g$ on the actual input to our machine. For each step of the simulation of the computation of $f$, when we need to know the symbol in the $i$th location of $f$'s imaginary input tape, we use the logspace machine for $g$ (using a partitioned off section of the read/write memory tape of size $O(\log(n))$, which we re-fill with blank symbols between consecutive calls to $g$) to compute from scratch the $i$th output symbol of $g$ on the actual input to our machine.

For the more general statement, note that if we have logspace machines computing $g_1, ..., g_m$ then we can easily build a logspace transducer $g$ which simulates each of $g_1, ..., g_m$ in turn and concatenates their outputs.
\end{proof}

Using this fact, we can say that two problems are logspace-equivalent if we can transform the inputs of either one into the inputs of the other using a logspace transducer, so that the desired output of one on the initial input matches the desired output on the other on the transformed input. We can also define a partial order on problem difficulty based on logspace transducers, which we use as the basis for the definition of NP-completeness.

\begin{defn} We say that a function $g$ has a \emph{logspace reduction} to a function $f$ if there is a logspace transducer $t_{g \rightarrow f}$ such that
\[
g(x) = f(t_{g \rightarrow f}(x))
\]
for all possible inputs $x$ to $g$.

We say that a function $f$ (with Yes/No outputs) is \emph{NP-hard under logspace reductions} if for every problem $g$ which is in NP, there is a logspace reduction from $g$ to $f$. We say that $f$ is \emph{NP-complete under logspace reductions} if $f$ is both contained in NP and NP-hard under logspace reductions.
\end{defn}

\begin{prop} The ``generic NP problem'' from Example \ref{ex-np-complete} is NP-complete under logspace reductions.
\end{prop}
\begin{proof} For any NP problem $g$, by the definition of NP there is a machine $M$ such that the output $g(x)$ is ``Yes'' iff there is some $y$ such that $M$ returns ``Yes'' on the input $(x,y)$, along with a polynomial $p$ such that $M$ always halts within $p(|x|)$ steps, where $|x|$ is the length of $x$. To define a logspace transducer from inputs $x$ of $g$ to inputs of the problem from Example \ref{ex-np-complete}, we simply concatenate a description of $M$ with the original input $x$ and the value $k = p(|x|)$ expressed in unary to get the triple $(M,x,k)$. (Finding a logspace transducer which computes $k = p(|x|)$ expressed in unary is left as an easy exercise to the reader.)
\end{proof}

Now that we have found one NP-complete problem, we would like to find some simpler NP-complete problems. The first simplification is based on the observation that while Turing machines are a good framework for open-ended computations with arbitrarily large inputs, they are overkill for discussing computations where the inputs have a fixed size and which are guaranteed to finish in a fixed amount of time. A simpler framework for computations on fixed-size inputs is the framework of \emph{circuits}.

\begin{defn} Fix a domain $D$ and a collection of basic operations $\cO$ on $D$ (i.e. each operation $f \in \cO$ is a function $D^k \rightarrow D$, for some arity $k \ge 0$ depending on $f$). A \emph{circuit} $C = (V, E, s, t, \ell, <)$ over $(D, \cO)$ with $n$ inputs and $m$ outputs is defined to be the following type of labeled directed multigraph:
\begin{itemize}
\item $V$ is the set of vertices of $C$, and $E$ the set of directed edges, with $s : E \rightarrow V$ describing the source of each edge and $t : E \rightarrow V$ describing the target of each edge (so $e \in E$ is interpreted as a directed edge from $s(e)$ to $t(e)$),
\item the directed multigraph $(V,E,s,t)$ is required to be finite and to have no directed cycles,
\item the collection of incoming edges to any vertex $v \in V$ (that is, the collection of $e \in E$ such that $t(e) = v$) is totally ordered by $<$,
\item $\ell$ is a labeling $\ell : V \rightarrow \cO \sqcup \{\mathrm{in}_1, ..., \mathrm{in}_n\} \sqcup \{\mathrm{out}_1, ..., \mathrm{out}_m\}$,
\item $\ell$ labels each vertex of in-degree $0$ by either an element of $\{\mathrm{in}_1, ..., \mathrm{in}_n\}$ or an operation in $\cO$ of arity $0$,
\item $\ell$ labels each vertex of in-degree $k > 0$ and positive out-degree by an operation in $\cO$ with arity $k$, and
\item $\ell$ labels each vertex of out-degree $0$ by an element of $\{\mathrm{out}_1, ..., \mathrm{out}_m\}$, with each of the labels $\mathrm{out}_1, ..., \mathrm{out}_m$ occurring exactly once, and each vertex of out-degree $0$ having in-degree $1$.
\end{itemize}
If each basic operation in $\cO$ is a symmetric function of its inputs, then we will leave the edge-ordering $<$ out of our descriptions of circuits over $(D, \cO)$.

To each circuit $C$ with $n$ inputs and $m$ outputs, we associate a function $\mathrm{ev}_C : D^n \rightarrow D^m$ as follows. For each input tuple $x = (x_1, ..., x_n) \in D^n$, we inductively assign values $\mathrm{ev}_{x,C}(v)$ to the vertices $v$ of the circuit $C$:
\begin{itemize}
\item to each vertex $v \in V$ of in-degree $0$ which is labeled by $\ell(v) = \mathrm{in}_i$, we assign the value $\mathrm{ev}_{x,C}(v) = x_i$,
\item to each vertex $v \in V$ which is labeled by an operation $\ell(v) = g \in \cO$ of arity $k$, if $e_1 < \cdots < e_k$ are the incoming edges to $v$ (i.e. $t(e_i) = v$ for each $i$) in sorted order and $v_1 = s(e_1), ..., v_k = s(e_k)$ are the source vertices of these incoming edges, and if we have already inductively assigned the values $\mathrm{ev}_{x,C}(v_1), ..., \mathrm{ev}_{x,C}(v_k)$, then we set
\[
\mathrm{ev}_{x,C}(v) = g(\mathrm{ev}_{x,C}(v_1), ..., \mathrm{ev}_{x,C}(v_k)),
\]
\item finally, to each vertex $v \in V$ of out-degree $0$ labeled by $\ell(v) = \mathrm{out}_i$, if $e$ is the unique incoming edge with $t(e) = v$, and $v' = s(e)$ is the corresponding source vertex of the edge $e$, then we assign the same value to $v$ as we assigned to $v'$, that is, $\mathrm{ev}_{x,C}(v) = \mathrm{ev}_{x,C}(v')$.
\end{itemize}
After we have finished assigning values to the vertices of the circuit, we read off the output $\mathrm{ev}_{C}(x_1, ..., x_n) = (\mathrm{ev}_{x,C}(v_1), ..., \mathrm{ev}_{x,C}(v_m)) \in D^m$, where $v_1, ..., v_m$ are the vertices satisfying $\ell(v_i) = \mathrm{out}_i$. (Equivalently, instead of assigning input-dependent values in $D$ to the vertices of the circuit, we could have assigned operations $\mathrm{ev}_C(v) \in \langle\cO\rangle$, where $\langle\cO\rangle$ is the clone generated by $\cO$, to each vertex $v \in V$, starting by assigning the projection $\pi_i$ to each vertex $v$ labeled by $\ell(v) = \mathrm{in}_i$.)

We often refer to the vertices $v \in V$ as the \emph{gates} of the circuit $C$ - if $\ell(v) \in \cO$, then we say that $v$ is an \emph{$\ell(v)$-gate} of $C$, and otherwise $v$ is an \emph{input gate} if $\ell(v) \in \{\mathrm{in}_1, ..., \mathrm{in}_n\}$ or an \emph{output gate} if $\ell(v) \in \{\mathrm{out}_1, ..., \mathrm{out}_m\}$.

We define the \emph{size} of a circuit $C$ as the number $|V| - n - m$ of vertices in the directed multigraph $(V,E,s,t)$ which are not labeled as input or output gates. We define the \emph{depth} of $C$ to be the number of non-input, non-output gates occurring along the longest directed path in the directed multigraph $(V,E,s,t)$.
\end{defn}

\begin{rem} Note that in our definition of a circuit, we place no restrictions on the out-degrees of gates, and to keep things simple we also assume that each gate produces just a single element of the domain $D$ as output (which is copied to all of the outgoing wires from this gate), rather than an ordered tuple of outputs. We can get away with this because copying information is assumed to be a free operation in classical computing. However, when studying \emph{reversible} circuits (where every function $D^n \rightarrow D^n$ computed by a reversible circuit must be a permutation, with the same number of inputs as outputs), or when studying \emph{quantum} circuits (where quantum information must obey the \emph{no-cloning} theorem), it is no longer appropriate to make this assumption, and the definition of a circuit must be modified in those contexts.
\end{rem}

The next result shows that the particular choice of basic operations $\cO$ doesn't affect circuit complexity by more than a constant factor, as long as $\cO$ is finite and the clone $\langle\cO\rangle$ generated by $\cO$ is fixed.

\begin{prop} If $\cO_1, \cO_2$ are finite sets of basic operations on a set $D$ such that $\cO_1$ and $\cO_2$ generate the same clone, then there is a constant $K$ such that for any circuit $C_1$ over $(D, \cO_1)$ of size $s$ and depth $d$ which computes a function $\mathrm{ev}_{C_1} : D^n \rightarrow D^m$, there is a circuit $C_2$ over $(D, \cO_2)$ of size at most $K\cdot s$ and depth at most $K\cdot d$ which satisfies $\mathrm{ev}_{C_2} = \mathrm{ev}_{C_1}$.
\end{prop}
\begin{proof} Since $\cO_1$ and $\cO_2$ generate the same clone, every basic operation in $\cO_1$ can be computed by some circuit over $(D, \cO_2)$. Since $\cO_1$ is finite, we can find a finite upper bound $K$ on the size of the circuits over $(D, \cO_2)$ needed to compute the basic operations of $\cO_1$. Now we replace each gate in the circuit $C_1$ which is labeled by a basic operation $g \in \cO_1$ by a circuit over $(D, \cO_2)$ of size (and depth) at most $K$ which computes $g$, and check that the resulting circuit $C_2$ over $(D, \cO_2)$ has size at most $K\cdot s$ and depth at most $K \cdot d$.
\end{proof}

Since we can encode the elements of any particular finite set as binary sequences, we will focus on circuits over the domain $D = \{0,1\}$ - and we identify $0$ with ``False'' and $1$ with ``True'' to interpret the basic logical operations $\wedge, \vee, \neg$ as operations on the set $\{0,1\}$.

\begin{defn} We define a \emph{Boolean circuit} to be a circuit over $(\{0,1\}, \wedge, \vee, \neg, 0, 1)$. Note that since all of the basic operations in a Boolean circuit are symmetric functions of their inputs, we don't need to specify an ordering on the incoming edges to any of the gates.
\end{defn}

Eventually we will wish to simulate Turing machines by Boolean circuits. The most complex part of a Turing machine is the transition rule $\delta : \Sigma \times S \rightarrow \Sigma \times S \times \{\pm 1\}$ used to decide the next action and state, given the current state and the current symbol on the tape. The standard way to represent such a transition rule is to simply list out the desired output tuples for each possible input tuple. Representing the elements of $\Sigma$ and $S$ as binary sequences, we see that we need to find a way to convert a truth-table description of an arbitrary function $\delta : \{0,1\}^n \rightarrow \{0,1\}^m$ into a Boolean circuit which computes $\delta$.

\begin{prop} There is a logspace transducer which takes as input a truth-table description of an arbitrary function $\delta : \{0,1\}^n \rightarrow \{0,1\}^m$ and produces as output a Boolean circuit of size $O(mn\cdot 2^n)$ and depth $O(n)$ which computes $\delta$.
\end{prop}
\begin{proof} Note that the size of the input is $m\cdot 2^n$, since for each of the $2^n$ possible input tuples to $\delta$ we need $m$ bits to describe the corresponding output tuple of $\delta$, so by ``logarithmic space'' we mean that our read/write work tape can have size $O(\log(m\cdot 2^n)) = O(\log(m)) + O(n)$. We will just describe the circuit which will be produced as output - finding a particular logspace Turing machine which constructs this circuit will be left as an exercise to the reader.

The circuit we construct will actually compute each of the $m$ output bits completely separately, so we may as well assume that $m = 1$ to keep things simple. The circuit we construct will be based on the formula
\[
\delta(x_1, ..., x_n) = \bigvee_{y \in \{0,1\}^n} (x_1 = y_1 \wedge \cdots \wedge x_n = y_n \wedge \delta(y_1, ..., y_n)).
\]
For a fixed tuple $y_1, ..., y_n$, each formula ``$x_i = y_i$'' can be replaced by ``$x_i$'' if $y_i = 1$ or by ``$\neg x_i$'' if $y_i = 0$, while the value of $\delta(y_1, ..., y_n)$ can be replaced by either $0$ or $1$ depending on what is written in the truth-table for $\delta(y)$. Thus for fixed values of $y_1, ..., y_n$, the value of the subformula
\[
(x_1 = y_1 \wedge \cdots \wedge x_n = y_n \wedge \delta(y_1, ..., y_n))
\]
can be computed by a circuit with $n$ input gates, at most $n$ $\neg$-gates, one constant gate (either $0$ or $1$, depending on the value of $\delta(y_1, ..., y_n)$), and $n$ $\wedge$-gates, which has size (and depth) $O(n)$.

To finish the construction, we will use $2^n-1$ $\vee$-gates to combine the results of these subformula circuits computing $x = y \wedge \delta(y)$ for various values of $y \in \{0,1\}^n$ to get a circuit computing $\delta(x)$. In order to keep the depth of the resulting circuit bounded by $O(n)$, we can arrange the computation as a binary tree corresponding to the following formula:
\[
\bigvee_{y \in \{0,1\}^n} (x = y \wedge \delta(y)) = \bigvee_{y_1 \in \{0,1\}} \Big(\bigvee_{y_2 \in \{0,1\}} \Big( \cdots \Big(\bigvee_{y_n \in \{0,1\}} (x = y \wedge \delta(y)) \Big)\cdots\Big)\Big).\qedhere
\]
\end{proof}

\begin{cor} The clone generated by the basic Boolean functions $\wedge,\vee,\neg,0,1$ contains all operations $\delta : \{0,1\}^n \rightarrow \{0,1\}$.
\end{cor}

\begin{ex} Some authors prefer to include the \emph{exclusive-or} operation $\oplus : \{0,1\}^2 \rightarrow \{0,1\}$ as a basic operation in their circuits. This operation may be described by its truth table:
\begin{center}
\begin{tabular}{cc|c}
$x$ & $y$ & $x \oplus y$\\
\hline
0 & 0 & 0\\
0 & 1 & 1\\
1 & 0 & 1\\
1 & 1 & 0
\end{tabular}
\end{center}
The recipe above shows us how to mindlessly convert this truth table into a Boolean formula:
\begin{align*}
x \oplus y &= \bigvee_{u \in \{0,1\}}\Big(\bigvee_{v \in \{0,1\}} \big((x = u) \wedge (y = v) \wedge (u\oplus v)\big)\Big)\\
&= \Big(\big((\neg x \wedge \neg y \wedge 0)\vee (\neg x \wedge y \wedge 1)\big)\vee\big((x \wedge \neg y \wedge 1)\vee (x \wedge y \wedge 0)\big)\Big).
\end{align*}
With a little bit of post-processing, we can clean this up to get a tidier Boolean formula:
\begin{align*}
x \oplus y &= (\neg x \wedge y) \vee (x \wedge \neg y).
\end{align*}
This formula can be directly converted to a Boolean circuit of size $5$ and depth $3$, which has two $\neg$-gates, two $\wedge$-gates, and one $\vee$-gate.
\end{ex}

\begin{ex} Some authors prefer to include the \emph{multiplexer} $\mathrm{mux} : \{0,1\}^3 \rightarrow \{0,1\}$ as a basic operation in their circuits. This operation is defined by
\[
\mathrm{mux}(s,x,y) \coloneqq \begin{cases} x & \text{if }s = 0,\\ y & \text{if }s = 1,\end{cases}
\]
and in some programming languages this is represented with the compact if-then-else notation $s\ ?\ y : x$, and is referred to as the \emph{ternary operator}.

Mindlessly converting the truth table for $\mathrm{mux}$ into a Boolean formula, we get
\begin{align*}
\mathrm{mux}(s,x,y) &= \Big(\big((\neg s \wedge \neg x \wedge \neg y \wedge 0)\vee (\neg s \wedge \neg x \wedge y \wedge 0)\big)\vee\big((\neg s \wedge x \wedge \neg y \wedge 1)\vee (\neg s \wedge x \wedge y \wedge 1)\big)\Big)\\
&\;\;\;\;\; \vee \Big(\big((s \wedge \neg x \wedge \neg y \wedge 0)\vee (s \wedge \neg x \wedge y \wedge 1)\big)\vee\big((s \wedge x \wedge \neg y \wedge 0)\vee (s \wedge x \wedge y \wedge 1)\big)\Big).
\end{align*}
With a tiny bit of post-processing, this simplifies to the Boolean formula
\begin{align*}
\mathrm{mux}(s,x,y) &= \big((\neg s \wedge x \wedge \neg y)\vee (\neg s \wedge x \wedge y)\big) \vee \big((s \wedge \neg x \wedge y) \vee (s \wedge x \wedge y)\big),
\end{align*}
and using slightly more thought we can simplify it further to the Boolean formula
\[
\mathrm{mux}(s,x,y) = (\neg s \wedge x) \vee (s \wedge y).
\]
This formula can be directly converted to a Boolean circuit of size $4$ and depth $3$, which has one $\neg$-gate, two $\wedge$-gates, and one $\vee$-gate.
\end{ex}

Now we can describe how we will simulate Turing machines by Boolean circuits, when both the input size and the number of computation steps are fixed.

\begin{prop} There is a logspace transducer which takes as input a tuple $(M,n,k)$, where $M$ is a description of a Turing machine and $n,k$ are natural numbers written in unary, and produces as output a description of a Boolean circuit $C$ which takes $n$ inputs $x_1, ..., x_n$ and produces as output a single bit, which evaluates to $1$ if and only if running the Turing machine $M$ on the input $x = (x_1, ..., x_n)$ halts within $k$ steps and outputs ``Yes''.
\end{prop}
\begin{proof}[Sketch] Note that if we only run the Turing machine $M$ for $k$ steps, then we can never examine any entry of the tape which is more than $k$ steps away from the starting position $0$. So we can represent the configuration of the tape at each step of the computation by a sequence of $O(k)$ symbols from $\Sigma$.

Now, there are two reasonable approaches for representing the current location of the Turing machine's head $p$ in the tape. One approach is to set aside $O(\log_2(k))$ bits to record the value of $p$ in binary. A simpler approach is to store a bit next to each tape location, which keeps track of the answer to the question ``is this location of the tape the location which the head $p$ is pointing to?''. If we follow the second approach, then we should also set aside an extra sequence of bits next to each tape location which are meant keep track of what state the Turing machine is in, if the head $p$ is pointing to that tape location (and which can contain meaningless garbage if the head $p$ is not pointing to that location).

So the complete description of a configuration of the Turing machine $M$ at any particular step can be represented as a sequence of $O(k\cdot (\log(|\Sigma|) + 1 + \log(|S|)))$ binary bits. Given a binary representation the configuration of the machine at the $i$th step, we can compute the binary representation of the configuration at the $(i+1)$th step with a circuit built out of $k$ copies of the circuit computing the transition function $\delta$ of the Turing machine $M$ (one copy for each location in the tape where we might currently find the head $p$).

Stringing these circuits together, we get a circuit with $O(k \cdot k \cdot (\log(|\Sigma|) + 1 + \log(|S|))^2 \cdot |\Sigma|\cdot |S|)$ gates which computes a binary representation of the configuration of the Turing machine $M$ after $k$ steps, and we can then easily read off whether the machine has halted and whether the output is ``Yes'' from that binary representation.
\end{proof}

We are now ready to describe our second NP-complete problem.

\begin{ex} The (Boolean) \emph{Circuit-SAT} problem takes as input a Boolean circuit $C$ with $n$ input bits and one output bit, and the desired output is ``Yes'' if there exists a tuple of $n$ bits $x = (x_1, ..., x_n) \in \{0,1\}^n$ such that $\mathrm{ev}_C(x) = 1$ (i.e., when we run the circuit $C$ on the input $x$, the output bit is $1$).
\end{ex}

The previous proposition shows that Circuit-SAT is NP-complete under logspace reductions. While circuits are much simpler objects than Turing machines, they are still fairly complex. The next trick is to atomize a circuit $C$ by introducing new variables corresponding to the outputs of each gate occurring within $C$ - since each gate in a Boolean circuit has at most two inputs and one output value, we will get a collection of simple constraints which each involve at most $2+1 = 3$ of the expanded collection of variables. This leads us to an extremely simple NP-complete problem called 3-SAT.

\begin{defn} If $x_1, ..., x_n$ are Boolean variables, then a \emph{literal} is defined to be an expression which is either of the form $x_i$ or $\neg x_i$ for some variable $x_i$. A \emph{clause} of \emph{length} $k$ is an expression of the form $(\ell_1 \vee \cdots \vee \ell_k)$, where each $\ell_i$ is a literal.

An \emph{$k$-CNF formula} (where CNF is an abbreviation for \emph{conjunctive normal form}) is defined to be a formula of the form $C_1 \wedge \cdots \wedge C_m$, where each $C_i$ is a clause of length at most $k$.
\end{defn}

\begin{ex} For any $k$, we define \emph{$k$-SAT} to be the problem where the input is a $k$-CNF formula, and the desired output is ``Yes'' if there is an assignment of values in $\{0,1\}$ to the variables which makes the formula true.
\end{ex}

\begin{thm}\label{thm-3-sat-np-complete} 3-SAT is NP-complete under logspace reductions.
\end{thm}
\begin{proof} By the previous proposition, there is a logspace reduction from any problem in NP to Circuit-SAT, so we just need to find a logspace reduction from Circuit-SAT to 3-SAT. As usual, we will describe the desired output of the logspace transducer for a given circuit $\cC$, but we will leave the detailed construction of such a transducer as an exercise to the reader.

Let $\cC = (V,E,s,t,\ell)$ be a Boolean circuit with $n$ input bits and one output bit. Suppose for simplicity that each input $\mathrm{in}_i$ occurs as a label for an input gate of $\cC$ exactly once. We will produce a 3-CNF formula $\Phi$ with $|V|$ variables, written $x_v$ for $v \in V$. The intended meaning of the variable $x_v$ is the unknown value of $\mathrm{ev}_{\mathrm{in},\cC}(v)$ (i.e. the unknown output of gate $v$, which depends on the unknown values of the inputs $\mathrm{in}_1, ..., \mathrm{in}_n$).

For each gate $v$ which is labeled by an operation $\ell(v) \in \{\wedge,\vee,\neg,0,1\}$, we will introduce a collection of at most four clauses $C_{v,j}$ of length at most $3$ which will be used to guarantee that the value of $x_v$ is computed correctly, given the values of $x_{v_i}$ for the incoming edges $e_i$ which satisfy $s(e_i) = v_i$, $t(e_i) = v$.

If $v_\wedge$ is an $\wedge$-gate, and if $e_1, e_2$ are two the incoming edges to $v_\wedge$ and $v_1 = s(e_1), v_2 = s(e_2)$ are the corresponding source vertices of these edges, then we proceed as follows. By examining a truth table for $\neg(x_{v_1} \wedge x_{v_2})$, we see that
\begin{align*}
(x_{v_\wedge} \ne x_{v_1} \wedge x_{v_2}) \;\;\; \iff \;\;\; &(\neg x_{v_1} \wedge \neg x_{v_2} \wedge x_{v_\wedge}) \vee (\neg x_{v_1} \wedge x_{v_2} \wedge x_{v_\wedge})\\
&\vee (x_{v_1} \wedge \neg x_{v_2} \wedge x_{v_\wedge}) \vee (x_{v_1} \wedge x_{v_2} \wedge \neg x_{v_\wedge}),
\end{align*}
so by De Morgan's law we have
\begin{align*}
(x_{v_\wedge} = x_{v_1} \wedge x_{v_2}) \;\;\; \iff \;\;\; &(x_{v_1} \vee x_{v_2} \vee \neg x_{v_\wedge}) \wedge (x_{v_1} \vee \neg x_{v_2} \vee \neg x_{v_\wedge})\\
&\wedge (\neg x_{v_1} \vee x_{v_2} \vee \neg x_{v_\wedge}) \wedge (\neg x_{v_1} \vee \neg x_{v_2} \vee x_{v_\wedge}).
\end{align*}
Note that the right hand side is a conjunction of four clauses of length $3$ - we will name these clauses $C_{v_\wedge,0}, C_{v_\wedge,1}$, $C_{v_\wedge,2}, C_{v_\wedge,3}$.

We handle the other gates which are labeled by operations in $\{\wedge,\vee,\neg,0,1\}$ similarly. For the output gate $v_{\mathrm{out}}$, if we let $e_{\mathrm{out}}$ be the unique incoming edge and $v_0$ be the source of the edge $e_{\mathrm{out}}$, then we will introduce the clause $C_{v_{\mathrm{out}}, 0} = (x_{v_0})$ of length one to ensure that the output of the circuit evaluates to $1$. (We don't need to introduce any clauses for the input gates, since the $n$ inputs can have arbitrary values.)

To finish the construction, we define the formula $\Phi$ to be
\[
\Phi = \bigwedge_{v \in V} \bigwedge_j C_{v,j}.
\]
The formula $\Phi$ is a 3-CNF formula by construction. It is straightforward to check that $\Phi$ is satisfiable if and only if the circuit $\cC$ is satisfiable (and in fact that there is a bijection between satisfying assignments to $\Phi$ and $\cC$), which completes the proof.
\end{proof}

Simple as it is, it can be argued that 3-SAT is still a rather complex problem - clauses can have anywhere between $0$ and $3$ negated variables, so we have four distinct types of clauses on three variables. The 1-IN-3 SAT problem from Example \ref{ex-1-in-3} has just one type of three-variable constraint, which is additionally symmetric in its arguments.

\begin{thm}\label{thm-1-in-3-np-complete} 1-IN-3 SAT is NP-complete under logspace reductions.
\end{thm}
\begin{proof} We just need to find a logspace reduction from 3-SAT to 1-IN-3 SAT. We will do this by finding a mechanical transformation from any instance of 3-SAT to an instance of 1-IN-3 SAT which will have a solution if and only if the original 3-SAT instance had a solution. We will use the abbreviation
\[
\text{1-IN-3} = \{(0,0,1), (0,1,0), (1,0,0)\}.
\]

First, note that we can introduce two new variables $z, w$ in the 1-IN-3 SAT instance which we can force to take the values $0$ and $1$ (respectively) by imposing the constraint
\[
(z, z, w) \in \text{1-IN-3}.
\]
Next, we can introduce a negated copy $y_i$ of each variable $x_i$ from the original 3-SAT instance, and we can ensure that $y_i = \neg x_i$ by throwing in the constraint
\[
(x_i, y_i, z) \in \text{1-IN-3},
\]
using the fact that $z$ is forced to take the value $0$. This way every literal of the original 3-SAT instance corresponds to some variable of the 1-IN-3 SAT instance which we are constructing.

To finish, we just need to construct a primitive positive formula over 1-IN-3 which is equivalent to the clause $(\neg x_1 \vee \neg x_2 \vee \neg x_3)$. Here is one way to do it:
\begin{align*}
(\neg x_1 \vee \neg x_2 \vee \neg x_3) \;\;\; \iff \;\;\; \exists &a_1, a_2, a_3, b_1, b_2, b_3 \text{ s.t. } (a_1, a_2, a_3) \in \text{1-IN-3}\ \wedge\\
&(x_1, a_1, b_1) \in \text{1-IN-3} \wedge (x_2, a_2, b_2) \in \text{1-IN-3} \wedge (x_3, a_3, b_3) \in \text{1-IN-3}.
\end{align*}
By replacing some $x_i$s with their corresponding negated copies $y_i$, we can simulate the other types of clauses as well. Thus we can simulate any clause of the original 3-SAT instance by introducing $6$ new auxiliary variables (i.e. the $a_i$s and $b_i$s) and $4$ new 1-IN-3 constraints. Putting it all together, if the original 3-SAT instance had $n$ variables and $m$ constraints, then the 1-IN-3 SAT instance we construct will have $2 + 2n + 6m$ variables and $1 + n + 4m$ constraints, and it will have a solution if and only if the original 3-SAT instance had a solution.
\end{proof}

The NAE-SAT (where NAE is short for ``Not All Equal'') problem from Example \ref{ex-nae} is arguably even more symmetric that 1-IN-3 SAT, since the NAE constraint is also preserved by simultaneously replacing all three inputs with their negations. This extra symmetry makes proving that NAE-SAT is NP-complete slightly more tricky.

\begin{thm}\label{thm-nae-np-complete} NAE-SAT is NP-complete under logspace reductions.
\end{thm}
\begin{proof} It's enough to find a logspace reduction from 1-IN-3 SAT to NAE-SAT. We will do this by finding a mechanical transformation from any instance of 1-IN-3 SAT to an instance of NAE-SAT which will have a solution if and only if the original 1-IN-3 SAT instance had a solution. We will use the abbreviation
\[
\text{NAE} = \{0,1\}^3 \setminus \{(0,0,0), (1,1,1)\},
\]
as well as the abbreviation 1-IN-3 introduced in the previous proof.

The main difficulty is that there is no way to guarantee that any particular variable in our NAE-SAT instance is forced to take the value $1$, due to the symmetry between $0$ and $1$. We get around this by exploiting that very symmetry: we introduce two new variables $z, w$ in the NAE-SAT instance, which we force to take opposite values by imposing the constraint
\[
(z,z,w) \in \text{NAE}.
\]
Now, from here on we will simply pretend that $z = 0$ and that $w = 1$ - if this is not the case, then we will swap $0$ with $1$ throughout any hypothetical solution to our NAE-SAT instance to replace it with a solution where we do have $z = 0$ and $w = 1$.

To finish, we just need to construct a primitive positive formula over NAE and the unary constraint $\{1\}$ which is equivalent to $(x_1, x_2, x_3) \in \text{1-IN-3}$. Here is one way to do it:
\begin{align*}
(x_1, x_2, x_3) \in \text{1-IN-3} \;\;\; \iff \;\;\; &(x_1, x_2, x_3) \in \text{NAE} \wedge (x_1, x_2, 1) \in \text{NAE}\ \wedge\\
&(x_1, 1, x_3) \in \text{NAE} \wedge (1, x_2, x_3) \in \text{NAE}.
\end{align*}
Replacing the $1$s with the variable $w$ which we introduced earlier, we get a NAE formula which doesn't involve the unary constraint $\{1\}$. If the original 1-IN-3 SAT instance had $n$ variables and $m$ constraints, then the NAE-SAT instance we construct this way will have $2 + n$ variables and $4m$ constraints, and will have exactly twice as many solutions as the original 1-IN-3 SAT instance.
\end{proof}


If we want an NP-complete constraint satisfaction problem which only involves binary constraints, we will have to consider a larger domain. Arguably the simplest binary relation on the domain $[k]$ is the $\ne$ relation, which is fully symmetric. This corresponds to the $k$-coloring problem from Example \ref{ex-k-color}.

\begin{thm}\label{thm-k-color-np-complete} For every $k \ge 3$, the $k$-coloring problem is NP-complete under logspace reductions.
\end{thm}
\begin{proof} It's enough to find a logspace reduction from NAE-SAT to $k$-coloring. To make things a little bit more notationally convenient, let's replace the domain of $k$-coloring with the set $\{0, 1, ..., k-1\}$. We use the same abbreviation NAE for $\{0,1\}^3 \setminus \{(0,0,0), (1,1,1)\}$ as in the previous proof.

Just as in the previous argument, we have to deal with the fact that any permutation of the colors $0, ..., k-1$ leaves any instance of $k$-coloring unchanged. We deal with this by introducing $k$ new variables $c_0, ..., c_{k-1}$, which we force to all have different values by imposing the $\binom{k}{2}$ constraints
\[
c_i \ne c_j
\]
for each pair $i,j < k$ with $i \ne j$. Now, from here on we will simply pretend that $c_i = i$ for each $i < k$ - if this is not the case, then we will permute the colors throughout any hypothetical solution to our $k$-coloring instance to replace it with a solution where $c_i = i$ for all $i < k$.

To finish, we just need to construct a primitive positive formula over $\ne$ together with the collection of unary constraints $\{i\}$ for $i < k$ which is equivalent to $(x_1, x_2, x_3) \in \text{NAE}$ (note that $(x_1, x_2, x_3) \in \text{NAE}$ implies that each of $x_1, x_2, x_3$ is in $\{0,1\}$). Note that the collection of unary constraints $\{i\}$ together with the binary constraint $\ne$ allow us to primitively positively define \emph{every} unary relation on $\{0, ..., k-1\}$, by the formula
\[
x \in U \;\;\; \iff \;\;\; \bigwedge_{i \not\in U} x \ne i.
\]
Now we can use the formula
\begin{align*}
(x_1, x_2, x_3) \in \text{NAE} \;\;\; \iff \;\;\; &(x_1, x_2, x_3 \in \{0,1\})\ \wedge \\
&\big(\exists y_1, y_2, y_3 \in \{0,1,2\} \text{ s.t. } (y_1\ne y_2) \wedge (y_1 \ne y_3) \wedge (y_2 \ne y_3)\\
&\ \wedge (x_1 \ne y_1) \wedge (x_2 \ne y_2) \wedge (x_3 \ne y_3)\big).
\end{align*}
If the original NAE-SAT instance had $n$ variables and $m$ constraints, then the $k$-coloring instance we construct this way will have $k + n + 3m$ variables and $\binom{k}{2} + (3(k-2) + 3(k-3) + 6)m$ constraints, and will have exactly $k!$ times as many solutions as the original NAE-SAT instance.
\end{proof}

\chapter{Initial Intuition}\label{chapter-intuition}

\section{The $\Inv$-$\Pol$ Galois connection}\label{section-inv-pol}


We begin by recalling some definitions from the introduction.

\begin{defn} A set of relations $\Gamma$ on a fixed domain $D$ is called a \emph{relational clone} if it contains the equality relation, and is closed under permutations, adding dummy variables, existential projection, and intersections. Equivalently, a relational clone is a set of relations which is closed under defining new relations via primitive positive formulas.
\end{defn}

\begin{defn} A set of functions $D^k \rightarrow D, k \in \mathbb{N}$ is called a \emph{clone} if it contains the \emph{projections} $\pi_i^k : D^k \rightarrow D$ which satisfy $\pi_i^k(x_1, ..., x_k) = x_i$ (generally the superscript $k$ is omitted when it is clear), and is closed under \emph{composition}, the operation which takes a $k$-ary function $f$ and $k$ $l$-ary functions $g_1, ..., g_k$ to the function
\[
(f\circ (g_1, ..., g_k)) : (x_1, ..., x_{l}) \mapsto f(g_1(x_1, ..., x_l), ..., g_k(x_{1}, ..., x_{l})).
\]
\end{defn}

\begin{defn} A $k$-ary function $f$ is said to \emph{preserve} an $m$-ary relation $R$, written $f \rhd R$, if for every choice of $k$ $m$-tuples in $R$, applying $f$ componentwise produces a new $m$-tuple which is also in $R$. If we think of elements of $R$ as column vectors, we can write this as
\[
\begin{bmatrix} x_{11}\\ \vdots\\ x_{1m} \end{bmatrix}, ..., \begin{bmatrix} x_{k1}\\ \vdots\\ x_{km} \end{bmatrix} \in R \implies f\left( \begin{bmatrix} x_{11}\\ \vdots\\ x_{1m} \end{bmatrix}, ..., \begin{bmatrix} x_{k1}\\ \vdots\\ x_{km} \end{bmatrix}\right) = \begin{bmatrix} f(x_{11}, ..., x_{k1})\\ \vdots\\ f(x_{1m}, ..., x_{km}) \end{bmatrix} \in R.
\]
A function $f$ is a \emph{polymorphism} of a relational structure $(D, \Gamma)$ or of a relational clone $\Gamma$ if $f$ preserves $R_i$ for each relation $R_i \in \Gamma$.
\end{defn}

We can write the condition for $f \rhd R$ more compactly as $M \in R^k \implies f(M) \in R$, where $M \in R^k$ means that $M$ is a matrix with $k$ columns, each of which belongs to $R$, and $f(M)$ is the column vector obtained by applying $f$ to the rows of $M$.

In order to state the Galois connection, we need a few additional definitions.

\begin{defn} If $\Gamma$ is any set of relations on a domain $D$, then we define $\langle \Gamma \rangle$ to be the relational clone generated by $\Gamma$ (that is, $\langle \Gamma \rangle$ is the smallest relational clone which contains $\Gamma$). Similarly, if $\cO$ is any set of operations on $D$, we define $\langle \cO \rangle$ to be the clone generated by $\cO$. If $\bA = (D, \cO)$ is an algebraic structure, we let $\Clo(\bA)$ be the clone generated by the basic operations of $\bA$, so $\Clo(\bA) = \langle \cO\rangle$.
\end{defn}

\begin{defn} If $\Gamma$ is any set of relations on a domain $D$, then we define $\Pol(\Gamma)$ to be the the set of operations on $D$ that preserve every relation of $\Gamma$. If $\cO$ is any set of operations on $D$, we define $\Inv(\cO)$ to be the set of relations which are preserved by every operation in $\cO$. If we want to restrict to operations or relations of a specific arity, we use the notations
\begin{align*}
\Pol_k(\Gamma) &= \{f : D^k \rightarrow D \mid \forall R \in \Gamma,\ f \rhd R\},\\
\Inv_m(\cO) &= \{R \subseteq D^m \mid \forall f \in \cO,\ f \rhd R\}.
\end{align*}
\end{defn}

It is worth thinking about what sort of information about an algebraic structure $(D, \cO)$ can be found in $\Inv(\cO)$.

\begin{ex} If $\bA = (D, \cO)$ is an algebraic structure, then $\Inv_2(\cO)$ determines (among other things)
\begin{itemize}
\item the lattice of subalgebras of $\bA$,

\item $\Aut(\bA)$, the automorphism group of $\bA$,

\item $\End(\bA)$, the semigroup of endomorphisms of $\bA$,

\item $\Con(\bA)$, the lattice of congruences on $\bA$,

\item the set of partial orders on $D$ which are compatible with the operations of $\bA$, and

\item $\Inv_2(\bB)$ for any subalgebra $\bB \subset \bA$ or quotient $\bB = \bA/\!\!\sim$.
\end{itemize}
\end{ex}

It is easy to see that for all $\Gamma$, $\Pol(\Gamma)$ will be a clone, and that for all $\cO$, $\Inv(\cO)$ will be a relational clone. As a consequence, we have $\langle \Gamma \rangle \subseteq \Inv(\Pol(\Gamma))$ and $\langle \cO \rangle \subseteq \Pol(\Inv(\cO))$. The next two results show that these inclusions are actually equalities.

Before diving into the proof, the following concrete example will be useful for understanding the notation. Consider what it means for a ternary function $f$ to preserve the binary relation $\le$ (functions which preserve $\le$ are often called \emph{monotone}). Since $0 \le 0$, $0 \le 1$, and $1 \le 1$, we have
\[
\begin{bmatrix} 0\\ 0 \end{bmatrix}, \begin{bmatrix} 0\\ 1 \end{bmatrix}, \begin{bmatrix} 1\\ 1 \end{bmatrix} \in \; \le \;\;\; \implies \;\;\; \begin{bmatrix} f(0,0,1)\\ f(0,1,1) \end{bmatrix} \in \; \le,
\]
that is, $f(0,0,1) \le f(0,1,1)$. It's convenient to abbreviate the above as follows:
\[
\begin{bmatrix} 0 & 0 & 1\\ 0 & 1 & 1 \end{bmatrix} \in \; \le^3 \;\;\; \implies \;\;\; f\left(\begin{bmatrix} 0 & 0 & 1\\ 0 & 1 & 1 \end{bmatrix}\right) \in \; \le.
\]

\begin{thm} If $\Gamma$ is a set of relations on a finite domain $D$, then $\Inv(\Pol(\Gamma)) = \langle \Gamma \rangle$. In fact, if a relation $S \subseteq D^m$ is preserved by $\Pol(\Gamma)$ and can be generated by $k$ elements of $D^m$ (using operations of $\Pol(\Gamma)$), then $S$ can be defined by a primitive positive formula over $\Gamma$ which involves at most $|D|^k$ auxiliary variables.
\end{thm}
\begin{proof} Suppose that $S$ is generated by elements $x_1, ..., x_k \in D^m$, and let $X$ be the matrix having the $x_i$s as columns. Then $S = \{f(X) \mid f \in \Pol_k(\Gamma)\}$, so as a starting point we will construct a primitive positive formula $\Phi$ that describes $\Pol_k(\Gamma)$.

Note that $D^{D^k}$ is naturally interpreted as the set of functions $f : D^k \rightarrow D$: if $f \in D^{D^k}$, then the $(a_1, ..., a_k)$-coordinate of $f$ is $f(a_1, ..., a_k)$. We can now give a positive primitive formula for $\Pol_k(\Gamma) \subseteq D^{D^k}$:
\[
\Phi(f) \coloneqq \bigwedge_{R \in \Gamma} \bigwedge_{M \in R^k} f(M) \in R.
\]
If $\Gamma$ is infinite, the outer $\bigwedge$ will be an infinite conjunction. However, since there are only finitely many possible subsets of $D^{D^k}$, some finite subset $\Phi'$ of the inner conjunctions will define the same subset of $D^{D^k}$.

Finally, to define $S$ we use the primitive positive formula
\[
S(a) \coloneqq \exists f \in D^{D^k}\ \Phi'(f) \wedge (f(X) = a).\qedhere
\]
\end{proof}

\begin{thm} If $\cO$ is a set of operations on a finite domain $D$, then $\Pol(\Inv(\cO)) = \langle \cO \rangle$.
\end{thm}
\begin{proof} Suppose that $f \in \Pol(\Inv(\cO))$ is a $k$-ary function. Let $\cF(k) \subseteq D^{D^k}$ be the subalgebra of the algebraic structure $(D,\cO)^{D^k}$ generated by the functions $\pi_i : D^k \rightarrow D$, $\pi_i(x_1, ..., x_k) = x_i$. Then $\cF(k)$, interpreted as a set of functions from $D^k$ to $D$, is exactly the set of $k$-ary functions in $\langle \cO \rangle$.

Since $f \in \Pol(\Inv(\cO))$ and $\cF(k) \in \Inv(\cO)$, we must have $f \rhd \cF(k)$, so in particular we must have $f(\pi_1, ..., \pi_k) \in \cF(k)$. But $f(\pi_1, ..., \pi_k)$ is exactly $f$ thought of as an element of $D^{D^k}$, so this means that $f \in \langle \cO \rangle$.
\end{proof}

\begin{cor} There is an order reversing bijection between clones and relational clones, given by the operations $\Inv$ and $\Pol$.
\end{cor}

\begin{rem} The map $\{1, ..., k\} \rightarrow D^{D^k}$ given by $i \mapsto \pi_i$, where $\pi_i : D^k \rightarrow D$ is given by $\pi_i(x_1, ..., x_k) = x_i$, shows up in the theory of approximation algorithms as the \emph{long code}, which is the longest way of encoding $\{1, ..., k\}$ over the alphabet $D$ which doesn't have any redundant coordinates.
\end{rem}

\begin{ex} In the next section we will prove the following three correspondences between clones and relational clones on the domain $\{0,1\}$:
\begin{itemize}
\item $\langle \text{2-SAT} \rangle = \langle \le, \ne \rangle$ corresponds to $\langle \operatorname{maj} \rangle$ (the majority function on three inputs),

\item $\langle \text{HORN-SAT} \rangle = \langle \{0\}, \{1\}, x\wedge y\implies z \rangle$ corresponds to $\langle \min \rangle$ (the minimum function on two inputs), and

\item $\langle \text{XOR-SAT} \rangle = \langle \{1\}, x+y+z \equiv 0 \pmod{2} \rangle$ corresponds to $\langle x-y+z \pmod{2} \rangle$.
\end{itemize}
\end{ex}

\begin{defn} If $\bA, \bA'$ are two algebraic structures on the same domain such that every basic operation of $\bA'$ is in $\Clo(\bA)$, then we say that $\bA'$ is a \emph{reduct} of $\bA$ and that $\bA$ is an \emph{expansion} of $\bA'$. If $\Clo(\bA) = \Clo(\bA')$, then $\bA$ and $\bA'$ are called \emph{term equivalent}.
\end{defn}

The lattice of clones on the domain $\{0,1\}$ has been completely described - it has countably many elements, and is known as Post's lattice \cite{post-lattice} (see also chapter II.3 of \cite{lau-clone-theory}). It is known that on a domain of size $\ge 3$, there are uncountably many clones \cite{uncountable-clones}, \cite{zhuk-selfdual} (see also chapter II.8 of \cite{lau-clone-theory}). In particular, we see that most clones and relational clones can't be generated by finitely many functions or relations.

\begin{defn} A clone $\cO$ is said to be \emph{finitely generated} if there is a finite set $S$ of operations such that $\cO = \langle S \rangle$. It is said to be \emph{finitely related} if there is a finite set of relations $\Gamma$ such that $\cO = \Pol(\Gamma)$.
\end{defn}

\begin{ex}\label{ex-non-finitely-related} The clone on $\{0,1\}$ generated by the binary implication function $\to$, given by $\operatorname{\to}(x,y) = \neg x \vee y$, is finitely generated but not finitely related. One quick way to prove this is to show that for every $n \ge 3$, the $n$-ary threshold function $t_2^n$ defined by
\[
t_2^n(x_1, ..., x_n) = \begin{cases} 1 & \sum_i x_i \ge 2\\ 0 & \sum_i x_i \le 1\end{cases}
\]
is not in $\langle \to \rangle$, but every way of identifying two coordinates of $t_2^n$ gives a function which is in $\langle \to \rangle$ (exercise: why does this prove that $\langle \to \rangle$ can not be finitely related?). $\Inv(\to)$ is generated by the infinite sequence of relations $R_1, R_2, ...$ given by $R_n = \{0,1\}^n \setminus \{(0, ..., 0)\}$, and $\langle\to\rangle$ consists of all functions of the form $f(x_1, ..., x_n)\vee x_i$.
\end{ex}

Matthew Moore \cite{finitely-related-undecidable} has shown that determining whether a given finitely generated clone is finitely related is a Turing-complete problem, and therefore undecidable in general. It is conjectured that determining whether a given finitely related clone is finitely generated is also undecidable in general.

\begin{rem} The Galois connection between relational clones and clones on a finite set was originally discovered by Geiger in 1968 \cite{geiger-galois}, and Reinhard P\"oschel investigated the general case (where the domain may be infinite) in \cite{poschel-galois} - in the infinite case, the main difference is that clones must also be taken to be closed in the topology of pointwise convergence. (Jeavons reproved one direction of the connection - that relational clones on a finite domain are determined by their polymorphisms - in \cite{jeavons}.)
\end{rem}

\begin{rem} The Galois connection presented here, between operations and relations, can be viewed as being induced by the two-sorted preservation relation $\rhd$. In general, whenever one has a two-sorted binary relation $R$ on a pair of sets $A, B$, one can define operations $F,G$ on the power sets of $A, B$ respectively by
\begin{align*}
F(S) &= \{b \in B \mid \forall a \in S\ aRb\},\\
G(T) &= \{a \in A \mid \forall b \in T\ aRb\}.
\end{align*}
The abstract order-theoretic properties of such a pair $F,G$ are
\begin{itemize}
\item $F$ and $G$ are antitone: $S \subseteq S' \implies F(S) \supseteq F(S')$, and similarly for $G$, and

\item for any $S \in \cP(A)$ and $T \in \cP(B)$, we have
\[
T \subseteq F(S)\ \iff\ S \subseteq G(T).
\]
\end{itemize}
Actually the first property listed is redundant, as we have
\[
(S\subseteq S') \wedge (F(S') \subseteq F(S'))\ \implies\ S \subseteq S' \subseteq G(F(S'))\ \implies\ F(S') \subseteq F(S),
\]
and either of $F,G$ is determined by the other together with the second property: $F(S) = \bigcup_{S \subseteq G(T)} T$. Additionally, the second property follows from the first property together with $S \subseteq G(F(S))$ and $T \subseteq F(G(T))$: for one direction, we have
\[
T \subseteq F(S)\ \implies S \subseteq G(F(S)) \subseteq G(T).
\]
Any such pair $F,G$ determines the binary relation $R$, since
\[
(a,b) \in R\ \iff\ b \in F(\{a\})\ \iff\ a \in G(\{b\}),
\]
and $F,G$ are both determined by the second property and their restrictions to singletons, since
\[
b \in F(S)\ \iff\ \forall a \in S,\ a \in G(\{b\})\ \iff\ \forall a \in S,\ b \in F(\{a\}).
\]

Then one can define closure operators $G\circ F, F\circ G$ on subsets of $A$ and $B$. When we say these are ``closure operators'', we mean that the images of these operators form collections of ``closed'' sets, such that any intersection of closed sets is closed, and for $S \subseteq A$, $G\circ F(S)$ is equal to the smallest closed set which contains $S$. All of these properties are easy to show directly in terms of the binary relation $R$, but they can also be proved order theoretically.

For the order theoretic proof, note that
\[
F(S) \subseteq F(S)\ \implies\ S \subseteq G\circ F(S),
\]
and similarly for $F\circ G$, and so we have
\[
F(S) \subseteq F\circ G(F(S)) = F(G\circ F(S)) \subseteq F(S),
\]
and we see that a set in $X \subseteq B$ is closed iff it is of the form $F(S)$ for some $S \subseteq A$. For the intersection property, note that
\[
S \subseteq G(X) \cap G(Y)\ \iff\ X \cup Y \subseteq F(S)\ \iff\ S \subseteq G(X\cup Y),
\]
and for the characterization of the closure of $S$ we have
\[
G\circ F(S) \subseteq G(Y)\ \iff\ Y \subseteq F\circ G\circ F(S) = F(S)\ \iff\ S \subseteq G(Y).
\]
Then $F$ and $G$ will provide a Galois correspondence between the closed subsets of $A$ and the closed subsets of $B$. The nontrivial thing to do is to describe the closure operators explicitly.

In our case, the relation $R$ was given by $\rhd$, and the sets $A,B$ were the sets of operations and relations on a given domain. Our main difficulty was in proving that the closure operators $G\circ F = \Pol \circ \Inv$ and $F \circ G = \Inv \circ \Pol$ were concretely described by the closure operators $\cO \mapsto \langle \cO \rangle$ for clones and $\Gamma \mapsto \langle \Gamma \rangle$ for relational clones, respectively. In ordinary Galois theory, the sets $A,B$ are taken to be a field and a group of automorphisms of the field, and the relation $R$ determines whether a given element of the field is fixed by a given automorphism (exercise: find the corresponding closure operations).
\end{rem}

\subsection{Galois connection for multisorted relational clones}

While single-sorted CSPs might seem cleaner from a theoretical point of view, any attempt to describe algorithms for CSPs on larger domains ends up needing to consider the multisorted case. The main reason for this is that most algorithms make progress by showing that some of the variables must take their values within a smaller subset of their original domain - as soon as we restrict some of the variables to smaller domains and not others, we end up with a problem that is most naturally described in the framework of multisorted CSPs.

A multisorted CSP template consists of a list of domains $D_1, D_2, ...$ for the sorts, together with a library of relations
\[
R_j \subseteq D_{i^j_1} \times \cdots \times D_{i^j_{m_j}}.
\]
For each relation $R_j$, we keep track of both its arity $m_j$ as well as its \emph{type}: the tuple of sorts $(i^j_1, ..., i^j_{m_j})$ which tells us what sorts of variables we are allowed to apply this relation to. Note that a multisorted CSP template is exactly the same thing as a multisorted relational structure (as defined in Section \ref{s-definitions}).

In order to describe an \emph{instance} of the multisorted CSP over the template $(\{D_i\}, \{R_j\})$, we first need a list of sets of variables $V_1, V_2, ...$, where $V_i$ is the set of variables whose values must be taken in the domain $D_i$. If $x \in V_i$, then we say that $D_i$ is the \emph{variable domain} for the variable $x$. Next, we give a list of sets of constraints $C_j$, where for each $j$,
\[
C_j \subseteq V_{i^j_1} \times \cdots \times V_{i^j_{m_j}}
\]
is a list of tuples $(x_1, ..., x_{m_j})$ of variables which we will require to be contained in the corresponding constraint relation $R_j$ from the template.

The careful reader may have noticed that an instance of the multisorted CSP over the template $(\{D_i\}, \{R_j\})$ is exactly the same thing as a multisorted relational structure $(\{V_i\}, \{C_j\})$ sharing the same signature as the template. A \emph{solution} to the instance is then the same thing as a homomorphism from the instance to the template, that is, a mapping from variables to elements of their variable domains so that each tuple of variables which is related by a constraint gets mapped to a tuple which lies in the corresponding constraint relation.

Just as in the single-sorted case, we are interested in the whole collection of all multisorted relations which can be expressed as projections (onto a subset of the variables) of solution sets to instances of the multisorted CSP with a given template $(\{D_i\}, \{R_j\})$. These correspond to formulas built out of bounded existential quantifiers (i.e. $\exists x \in D_i$), conjunctions, and relations $R_j$, where each existential quantifier must have some sort $D_i$ associated with the variable it binds, and each relation must be applied to a tuple of variables of the correct sorts for that relation. For instance, if $R_1$ was a binary relation
\[
R_1 \subseteq D_1 \times D_2
\]
and $R_2$ was a ternary relation
\[
R_2 \subseteq D_2 \times D_1 \times D_2,
\]
then the formula
\[
\exists x \in D_1\; \exists y \in D_2 \text{ s.t. } (x,y) \in R_1 \wedge (y,x,z) \in R_2
\]
describes a unary relation on $z$ which can be applied to variables of sort $2$. On the other hand, the formula
\[
(x,x,x) \in R_2
\]
would be considered invalid for this CSP template, since $x$ would be forced to have sort $1$ and sort $2$ simultaneously.

In the singly sorted case, there was no harm in also allowing the binary equality relation $=$ into our library of constraints. In the multisorted case, it doesn't make sense to do this since we haven't specified the type of the equality relation. Instead, we introduce a different equality relation $=_i$ for each sort $D_i$ of the domain:
\[
x =_i y \;\; \iff \;\; x,y \in D_i \wedge x = y.
\]
Now adding the equality relations $=_i$ into our constraint library causes no harm. The \emph{multisorted primitive positive formulas} over $(\{D_i\}, \{R_j\})$ can now be defined as formulas built out of bounded existential quantifiers (with associated sorts), conjunctions, relations $R_j$ (from the template) applied to tuples of variables of the correct sorts, and equality relations $=_i$ applied to pairs of variables of the same sort $i$.

Taking the closure of a multisorted library of constraint relations leads us to the definition of a multisorted relational clone.

\begin{defn} A collection $\Gamma$ of multisorted relations with sorts $\{D_i\}$ is called a \emph{multisorted relational clone} if it contains an equality relation $=_i$ for each sort $D_i$, and is closed under permuting variables (in tandem with their sorts), adding dummy variables (of any sorts), existential projection onto subsets of the variables, and intersection (of relations of the same arity and type).

Equivalently, a multisorted relational clone is a collection of multisorted relations which is closed under defining new relations via multisorted primitive positive formulas.
\end{defn}

What is the algebraic counterpart of a multisorted relational clone? Our first guess might be a multisorted algebraic structure, but the usual definition of a multisorted algebraic structure will not do what we want. The purpose of the algebraic operations is for us to be able to apply them coordinatewise to solutions to instances of CSPs in order to produce new solutions. This means that each of our algebraic operations must make sense when applied to values of any given sort, as long as all of its arguments have the same sort. Following Romov \cite{romov-multisorted-galois}, we will call these ``vectors of operations''.

\begin{defn} A \emph{vector of operations} on $\{D_i\}$ of arity $k$ is a tuple $f = (f^{D_1}, f^{D_2}, ...)$ where each $f^{D_i}$ is an operation on $D_i$ of arity $k$:
\[
f^{D_i} : D_i^k \rightarrow D_i.
\]
For any $j \le k$, we define the vector of operations $\pi_j^k$ by
\[
(\pi_j^k)^{D_i}(x_1, ..., x_k) = x_j
\]
for $x_1, ..., x_k \in D_i$, and we call $\pi_j^k$ a \emph{projection}.

For $f$ a vector of operations of arity $k$, and for $g_1, ..., g_k$ a collection of $k$ vector operations of arity $l$, we define their \emph{composition} by
\[
(f \circ (g_1, ..., g_k))^{D_i} : (x_1, ..., x_l) \mapsto f^{D_i}(g_1^{D_i}(x_1, ..., x_l), ..., g_l^{D_i}(x_1, ..., x_l))
\]
for $x_1, ..., x_l \in D_i$.

A \emph{multi-clone} is defined to be a collection of vectors of operations on a collection of sorts $\{D_i\}$ which contains the projections and is closed under composition.
\end{defn}

\begin{defn} If $f$ is a vector of operations on the sorts $\{D_i\}$ of arity $k$, and if
\[
R \subseteq D_{i_1} \times \cdots \times D_{i_m}
\]
is a multisorted relation of arity $m$ and type $(i_1, ..., i_m)$, then we say that $f$ \emph{preserves} $R$, written $f \rhd R$, if for every choice of $k$ $m$-tuples of $R$
\[
\begin{bmatrix} x_{11}\\ \vdots \\ x_{1m} \end{bmatrix}, ..., \begin{bmatrix} x_{k1}\\ \vdots \\ x_{km} \end{bmatrix} \in R,
\]
we have
\[
f\left(\begin{bmatrix} x_{11}\\ \vdots \\ x_{1m} \end{bmatrix}, ..., \begin{bmatrix} x_{k1}\\ \vdots \\ x_{km} \end{bmatrix}\right) \coloneqq \begin{bmatrix} f^{D_{i_1}}(x_{11}, ..., x_{k1})\\ \vdots\\ f^{D_{i_m}}(x_{1m}, ..., x_{km})\end{bmatrix} \in R.
\]
A vector of operations $f$ on $\{D_i\}$ is a \emph{polymorphism} of a multisorted relational structure $(\{D_i\}, \Gamma)$ or of a multisorted relational clone $\Gamma$ on $\{D_i\}$ if $f$ preserves $R_j$ for each relation $R_j \in \Gamma$.
\end{defn}

Note that a collection of vectors of operations $f_n$ on $\{D_i\}$ defines a collection of algebraic structures $(D_i, \{f_n^{D_i}\})$ which all share the same signature. So the ``correct'' algebraic counterpart to a multisorted relational structure ends up being a collection of algebraic structures which share a common signature.

\begin{defn} For $\Gamma$ a collection of multisorted relations on $\{D_i\}$, we define $\langle \Gamma \rangle$ to be the multisorted relational clone generated by $\Gamma$. Similarly, if $\cO$ is a collection of vectors of operations on $\{D_i\}$, then we define $\langle \cO \rangle$ to be the multi-clone generated by $\cO$.

If $\{\bD_i\} = \{(D_i, \{f_n^{D_i}\})\}$ is a collection of algebraic structures which share a common signature, then we define $\Clo(\{\bD_i\})$ to be the multi-clone $\langle \{f_n\} \rangle$ generated by the vectors $f_n$ of basic operations of $\{\bD_i\}$.
\end{defn}

We extend the definitions of $\Inv$ and $\Pol$ to the multisorted case in the obvious way.

\begin{defn} If $\Gamma$ is any collection of multisorted relations on $\{D_i\}$, then we define $\Pol(\Gamma)$ to be the the collection of vectors of operations on $\{D_i\}$ that preserve every relation of $\Gamma$. If $\cO$ is any collection of vectors of operations on $\{D_i\}$, we define $\Inv(\cO)$ to be the collection of multisorted relations which are preserved by every operation in $\cO$.
\end{defn}

Now we can state and prove the Galois connection for the multisorted case.

\begin{thm}[Romov \cite{romov-multisorted-galois}] If $\Gamma$ is a collection of multisorted relations on a (possibly infinite) collection of finite domains $\{D_i\}$, then $\Inv(\Pol(\Gamma)) = \langle \Gamma \rangle$.
\end{thm}
\begin{proof} Suppose that
\[
S \subseteq D_{i_1} \times \cdots \times D_{i_m}
\]
is generated by elements $x_1, ..., x_k \in \prod_j D_{i_j}$, and let $X$ be the matrix having the $x_j$s as columns. Then $S = \{f(X) \mid f \in \Pol_k(\Gamma)\}$, so as a starting point we will construct a primitive positive formula $\Phi$ that describes $\Pol_k(\Gamma)$.

Note that $\prod_i D_i^{D_i^k}$ is naturally interpreted as the set of vectors of functions $f^{D_i} : D_i^k \rightarrow D_i$: if $f \in \prod_i D_i^{D_i^k}$, then the $(i, (a_1, ..., a_k))$-coordinate of $f$ is $f^{D_i}(a_1, ..., a_k)$. We can now give an infinitary positive primitive formula for $\Pol_k(\Gamma) \subseteq \prod_i D_i^{D_i^k}$:
\[
\Phi(f) \coloneqq \bigwedge_{R \in \Gamma} \bigwedge_{M \in R^k} f(M) \in R.
\]
If $\Gamma$ is infinite, the outer $\bigwedge$ will be an infinite conjunction, and if there are infinitely many sorts $D_i$ then the number of variables of $\Phi(f)$ will also be infinite. However, since there are only finitely many possible subsets of $\prod_{j \le m} D_{i_j}^{D_{i_j}^k}$, and since $D_{i'}^{D_{i'}^k}$ is also finite for each $i' \not\in i_1, ..., i_m$, by K\"onig's Lemma some finite subset $\Phi'$ of the inner conjunctions will have
\[
\exists \{f^{D_{i'}} \mid i' \not\in i_1, ..., i_m\}\ \Phi(f) = \exists \{f^{D_{i'}} \mid i' \not\in i_1, ..., i_m\}\ \Phi'(f).
\]
Note that only finitely many operations $f^{D_{i'}}$ actually show up in $\Phi'(f)$, so the right hand side can be interpreted as a finite primitive positive formula.

Finally, to define $S$ we use the primitive positive formula
\[
S(a) \coloneqq \exists f \in \prod_i D_i^{D_i^k}\ \Phi'(f) \wedge (f(X) = a).\qedhere
\]
\end{proof}

\begin{rem} The same proof can be used to show that if we allow the domains $D_i$ to be infinite, $\Inv(\Pol(\Gamma))$ will be equal to the collection of multisorted relations which can be written as increasing unions of multisorted relations which can be defined by infinitary primitive positive formulas over $\Gamma$ which have finitely many free variables.
\end{rem}

\begin{thm}[Romov \cite{romov-multisorted-galois}] If $\cO$ is a collection of vectors of operations on a finite collection of finite domains $(D_1, ..., D_n)$, then $\Pol(\Inv(\cO)) = \langle \cO \rangle$.
\end{thm}
\begin{proof} Suppose that $f \in \Pol(\Inv(\cO))$ is a vector of $k$-ary operations. Let $\cF(k) \subseteq \prod_{i \le n} D_i^{D_i^k}$ be the subalgebra of the algebraic structure $\prod_i (D_i,\cO^{D_i})^{D_i^k}$ generated by the projections $\pi_j$, which are given by
\[
\pi_j^{D_i} : D_i^k \rightarrow D_i, \;\;\; \pi_j^{D_i}(x_1, ..., x_k) = x_j.
\]
Then $\cF(k)$, interpreted as a set of vectors of functions from $D_i^k$ to $D_i$, is exactly the set of $k$-ary vectors of operations in $\langle \cO \rangle$.

Since $f \in \Pol(\Inv(\cO))$ and $\cF(k) \in \Inv(\cO)$, we must have $f \rhd \cF(k)$, so in particular we must have $f(\pi_1, ..., \pi_k) \in \cF(k)$. But $f(\pi_1, ..., \pi_k)$ is exactly $f$ thought of as an element of $\prod_i D_i^{D_i^k}$, so this means that $f \in \langle \cO \rangle$.
\end{proof}

\begin{rem} The same argument can be used to give a multisorted generalization of Reinhard P{\"o}schel's result from \cite{poschel-galois} characterizing polymorphism clones on infinite domains. We say that a vector of operations $f$ of arity $k$ is in $\overline\cO$, the \emph{closure of $\cO$ under the topology of pointwise convergence}, if for every finite collection of tuples $\{(i^j, (a^j_1, ..., a^j_k))\} \subseteq \bigsqcup_i D_i^k$, there is some $g \in \cO$ which satisfies
\[
\forall j\ g^{D_{i^j}}(a^j_1, ..., a^j_k) = f^{D_{i^j}}(a^j_1, ..., a^j_k).
\]
Then the general result is that $\Pol(\Inv(\cO)) = \overline{\langle \cO \rangle}$, even if we have infinitely many infinite domains.
\end{rem}

\begin{rem} There are several ways to try to force the multisorted setting into the single sorted straitjacket. One approach is to start from a multisorted relational structure $(\{D_i\}, \{R_j\})$, and to make a new single sorted relational structure with domain
\[
D = \bigsqcup_i D_i
\]
out of it, with each of the original relations $R_j$ interpreted as relations on $D$ in the obvious way, and with each sort $D_i$ interpreted as a unary relation on $D$. Although this approach works, it is inelegant for the following reasons:
\begin{itemize}
\item the equality relation $=$ on $D$ doesn't correspond to anything natural in the original multisorted structure, and
\item a polymorphism $f$ of the new relational structure is forced to assign values to tuples $(a_1, ..., a_k) \in D^k$ which don't necessarily all come from the same sort, but the value of $f(a_1, ..., a_k)$ is irrelevant to every natural question we might ask unless there is some $i$ with $\{a_1, ..., a_k\} \subseteq D_i$.
\end{itemize}

If the number of domains $D_i$ is small, a more elegant approach was described in Bulatov and Jeavons \cite{bulatov-jeavons-varieties}. They define a single-sorted relational structure with domain
\[
D = \prod_i D_i,
\]
with each relation $R_j \subseteq D_{i^j_1} \times \cdots \times D_{i^j_{m_j}}$ interpreted as the $m_j$-ary relation on $D$ given by
\[
(x_1, ..., x_{m_j}) \in R_j^D \;\;\; \iff \;\;\; (\pi_{i^j_1}(x_1), ..., \pi_{i^j_{m_j}}(x_{m_j})) \in R_j,
\]
and with extra binary equivalence relations $\sim_i\ \subseteq D^2$ given by
\[
x \sim_i y \;\;\; \iff \;\;\; \pi_i(x) = \pi_i(y).
\]
On the algebraic side of the picture, this corresponds to replacing the collection $\{\bD_i\}$ of algebras (all sharing a common signature) with their product
\[
\bD = \prod_i \bD_i.
\]
Let's verify that the polymorphisms of the single-sorted relational structure $(D, \{R_j^D\} \cup \{\sim_i\})$ really correspond to the vectors of operations in $\Pol(\{D_i\}, \{R_j\})$.
\begin{itemize}
\item If $f \in \Pol(\{D_i\}, \{R_j\})$, then we define $f^D$ by
\[
\pi_i(f^D(a_1, ..., a_k)) = f^{D_i}(\pi_i(a_1), ..., \pi_i(a_k)).
\]
That this $f^D$ preserves each relation $R_j^D$ and each equivalence relation $\sim_i$ can be checked by unwinding the definitions.

\item If $f \in \Pol(D, \{R_j^D\} \cup \{\sim_i\})$, then since $f$ preserves $\sim_i$ we can define $f^{D_i}$ by arbitrarily picking reference elements $b_1, ..., b_k \in \prod_{j \ne i} D_j$, and setting
\[
f^{D_i}(a_1, ..., a_k) = \pi_i(f((a_1, b_1), ..., (a_k, b_k)))
\]
for $a_1, ..., a_k \in D_i$, and the resulting functions $f^{D_i}$ will not depend on the choice of $b_1, ..., b_k$. Then we can check that the vector of operations $i \mapsto f^{D_i}$ preserves $R_j$ iff $f$ preserves $R_j^D$ by unwinding the definitions.
\end{itemize}
\end{rem}

Throughout the rest of these notes, we will drop the distinction between operations and vectors of operations, as well as the distinction between multi-clones and clones.

\section{Three basic examples}

We start with the correspondence between 2-SAT and majority.

\begin{thm} Suppose that a relation $R \subseteq \{0,1\}^m$ is preserved by the majority function $\operatorname{maj} : \{0,1\}^3 \rightarrow \{0,1\}$. Then $R$ is bijunctive, that is, $R$ can be written as a conjunction of binary and unary relations.
\end{thm}
\begin{proof} We prove this by induction on $m$. If $m \le 2$ then there is nothing to prove. Otherwise, for each $i \le 3$ let $R_i$ be the existential projection of $R$ onto all variables except for the $i$th. We will show that $R$ is equivalent to
\[
\Phi(x_1, ..., x_m) \coloneqq R_1(x_2, x_3, ..., x_m) \wedge R_2(x_1,x_3, ..., x_m) \wedge R_3(x_1,x_2,...,x_m),
\]
and the result will then follow by the induction hypothesis. It is clear that $R \subseteq \Phi$, so suppose $(x_1, ..., x_m) \in \Phi$. Then by the definitions of $R_1, R_2, R_3$, there exist $x_1', x_2', x_3'$ such that $(x_1',x_2,x_3,...,x_m), (x_1,x_2',x_3,...,x_m), (x_1,x_2,x_3',...,x_m) \in R$, and applying $\maj$ to these three tuples we see that $(x_1, ..., x_m) \in R$ as well.
\end{proof}

\begin{defn} An operation $f : \{0,1\}^k \rightarrow \{0,1\}$ is called \emph{monotone} if it preserves the relation $\le$. It is called \emph{self-dual} if it preserves the relation $\ne$.
\end{defn}

\begin{thm} Suppose that a function $f : \{0,1\}^k \rightarrow \{0,1\}$ is monotone and self-dual. Then $f \in \langle \maj \rangle$.
\end{thm}
\begin{proof} We prove this by induction on $k$. It's easy to check that there are no monotone self-dual functions of arity $\le 2$ other than the projections, so assume that $k \ge 3$. By the induction hypothesis, any function we can make by identifying two variables of $f$ is in $\langle \maj \rangle$. We claim that we have
\[
f(x,y,z,...) = \maj(f(x,y,y,...), f(z,y,z,...), f(x,x,z,...)),
\]
where the $...$ always represent the remaining $k-3$ variables. To see this, note that the formula is trivially true when $x=y=z$, so we only need to check it when one of the variables is different from the other two. We will check it in the case $(x,y,z) = (0,1,0)$, since every other case is analogous (via cyclically permuting $x,y,z$ or swapping $0$s and $1$s throughout). In this case, since $f$ is monotone we have
\[
f(0,1,1,...) \ge f(0,1,0, ...) \ge f(0,0,0,...),
\]
so the median of these three values will be $f(0,1,0,...) = f(x,y,z,...)$, and the majority is equal to the median on $\{0,1\}$.
\end{proof}

Examining the proof, we see that every $n$-ary monotone self-dual function $f$ can be written in terms of $\maj$ as a term of depth at most $n-2$, such that every subterm is obtained by identifying some of the variables of $f$.

\begin{cor} For any odd $n$, the $n$-ary function $m_n$ given by
\[
m_n(x_1, ..., x_n) \coloneqq \begin{cases}1 & \sum_i x_i > \frac{n}{2},\\ 0 & \sum_i x_i < \frac{n}{2}\end{cases}
\]
is in the clone generated by $\maj$. In fact, we may write $m_n$ as a term of depth at most $n-2$, such that every subterm is also a linear threshhold function, where for $a \in \NN^n$ with $\sum_i a_i = n$ we define the $n$-ary linear threshhold function $t_a$ by
\[
t_a(x_1, ..., x_n) \coloneqq \begin{cases}1 & \sum_i a_ix_i > \frac{n}{2},\\ 0 & \sum_i a_ix_i < \frac{n}{2}.\end{cases}
\]
\end{cor}

For the majority function $m_n$, we can actually find a substantially smaller term using a probabilistic construction. (A deterministic construction, based on sorting networks, can be found in \cite{majority-sorting-networks}.)

\begin{prop}[Valiant \cite{valiant-majority}] For any odd $n$, the majority function $m_n$ can be represented by a term of depth $O(\log(n))$ and size $O(n^{4.3})$.
\end{prop}
\begin{proof} We'll follow Goldreich's exposition \cite{goldreich-majority}. Consider the completely generic formula $f_\ell(y_1, ..., y_{3^\ell})$ of depth $\ell$, defined recursively by $f_0 = \pi_1$, $f_1 = \maj$, and
\[
f_{\ell+1}(y) \coloneqq \maj(f_\ell(y_1, ..., y_{3^\ell}), f_\ell(y_{3^\ell+1}, ..., y_{2\cdot 3^\ell}), f_\ell(y_{2\cdot 3^\ell+1}, ..., y_{3^{\ell+1}})).
\]
Then define a random function $g_{\ell}(x_1, ..., x_n)$ by replacing each $y_i$ in $f_\ell$ with a random choice of $x_{j_i}$, where the $j_i$ are independently and uniformly randomly chosen from the set $\{1, ..., n\}$. For any particular input $x \in \{0,1\}^n$, if $p_i$ is the probability that $g_i(x) = m_n(x)$, then we have
\[
p_0 \ge \frac{1}{2} + \frac{1}{2n}
\]
and
\begin{align*}
p_{i+1} &= 3(1-p_i)p_i^2 + p_i^3\\
&= 0.5 + (1.5 - 2(p_i-0.5)^2)(p_i-0.5)\\
&= 1 - (3 - 2(1-p_i))(1-p_i)^2.
\end{align*}
A little computation then shows that for $\ell \approx (1+1/\log_2(1.5))\log_2(n) \approx 2.71\log_2(n)$ we have $1-p_\ell < 2^{-n}$, so a union bound shows that for this choice of $\ell$ at least one assignment to the $y_i$s has $g_\ell(x) = m_n(x)$ for all $x \in \{0,1\}^n$.
\end{proof}

Monotone self-dual functions can be interpreted as voting functions. They also have a combinatorial interpretation in terms of maximal ``intersecting families'' of sets.

\begin{defn} Let $S$ be a set. A family $\cF\subseteq \cP(S)$ is called an \emph{intersecting family} of subsets of $S$ if $A,B\in \cF$ implies $A\cap B\ne 0$.
\end{defn}

\begin{prop} An intersecting family of subsets of a set $S$ is maximal (with respect to containment) if and only if for every set $A\subseteq S$ we have either $A\in\cF$ or $(S\setminus A)\in\cF$. For every $n\ge 1$ there is a bijection between maximal intersecting families $\cF$ of subsets of $\{1,...,n\}$ and monotone self-dual boolean functions $f:\{0,1\}^n\rightarrow \{0,1\}$.
\end{prop}

We can describe a maximal intersecting family of subsets of a set more compactly by describing its collection of minimal elements. We can mutate an intersecting family by taking one of its minimal elements $A$, deleting it, and replacing it with its complement - this is called ``switching'' the subset $A$ with its complement.

\begin{defn} For every $n$, we define an undirected graph $\mathcal{M}_n$ whose vertices are the maximal intersecting families of subsets of $\{1,...,n\}$, and whose edges are the pairs of families $\cF, \mathcal{G}$ such that $|\cF \setminus \mathcal{G}| = 1$.
\end{defn}

The graph $\cM_4$ is depicted below, with vertices labeled by the minimal elements of the corresponding intersecting families as well as the corresponding monotone self-dual functions (written in terms of the majority function, which we abbreviate as $m$).

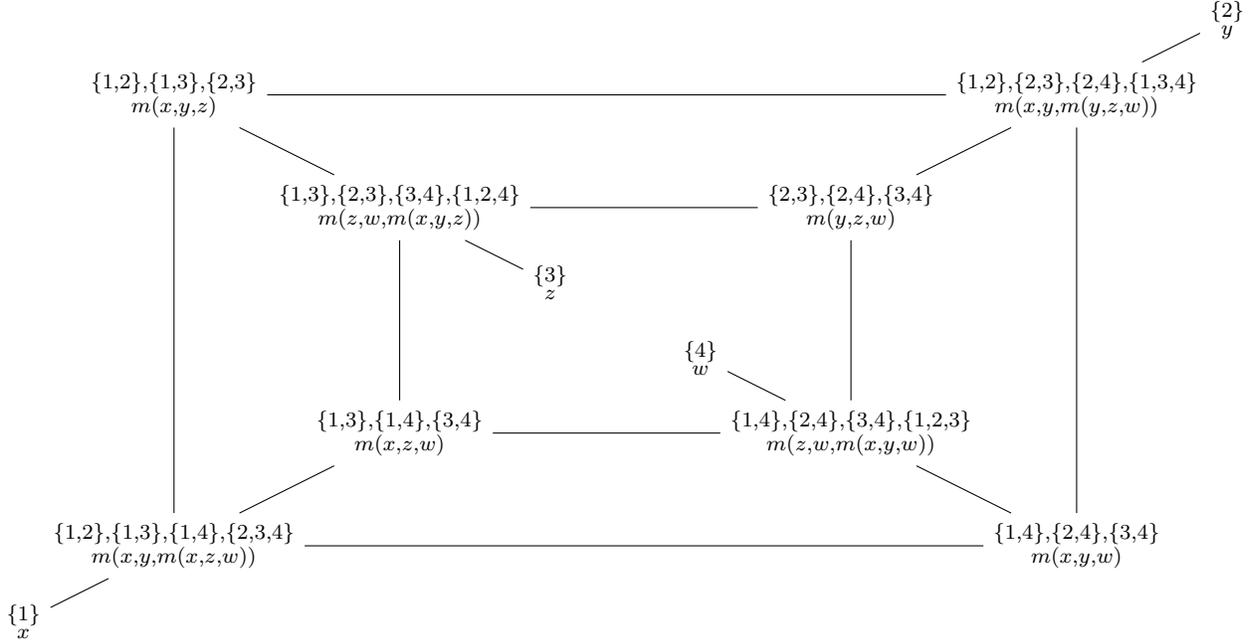
\begin{center}
\begin{tikzpicture}[scale=1]
  \node (xyz) at (-6,3) {$\substack{\{1,2\},\{1,3\},\{2,3\}\\ m(x,y,z)}$};
  \node (zwxyz) at (-3,1.5) {$\substack{\{1,3\},\{2,3\},\{3,4\},\{1,2,4\}\\ m(z,w,m(x,y,z))}$};
  \node (z) at (-1,0.5) {$\substack{\{3\}\\ z}$};
  \node (w) at (1,-0.5) {$\substack{\{4\}\\ w}$};
  \node (zwxyw) at (3,-1.5) {$\substack{\{1,4\},\{2,4\},\{3,4\},\{1,2,3\}\\ m(z,w,m(x,y,w))}$};
  \node (xyw) at (6,-3) {$\substack{\{1,4\},\{2,4\},\{3,4\}\\ m(x,y,w)}$};
  \node (x) at (-8,-4) {$\substack{\{1\}\\ x}$};
  \node (xyxzw) at (-6,-3) {$\substack{\{1,2\},\{1,3\},\{1,4\},\{2,3,4\}\\ m(x,y,m(x,z,w))}$};
  \node (xzw) at (-3,-1.5) {$\substack{\{1,3\},\{1,4\},\{3,4\}\\ m(x,z,w)}$};
  \node (yzw) at (3,1.5) {$\substack{\{2,3\},\{2,4\},\{3,4\}\\ m(y,z,w)}$};
  \node (xyyzw) at (6,3) {$\substack{\{1,2\},\{2,3\},\{2,4\},\{1,3,4\}\\ m(x,y,m(y,z,w))}$};
  \node (y) at (8,4) {$\substack{\{2\}\\ y}$};
  \draw (x) -- (xyxzw) -- (xyz) -- (xyyzw) -- (y);
  \draw (xyz) -- (zwxyz) -- (xzw) -- (xyxzw) -- (xyw) -- (xyyzw) -- (yzw) -- (zwxyz) -- (z);
  \draw (w) -- (zwxyw) -- (xyw);
  \draw (xzw) -- (zwxyw) -- (yzw);
\end{tikzpicture}
\end{center}

The graph $\mathcal{M}_n$ is always connected: given two maximal intersecting families $\mathcal{F}, \mathcal{G}$, there will always be some minimal element of $\cF$ which is not contained in $\mathcal{G}$, and switching this set with its complement gives us a maximal intersecting family $\cF'$ which is adjacent to $\cF$ and has one more element in common with $\mathcal{G}$ than $\cF$ does. For more about maximal intersecting families of sets, see \cite{maximal-intersecting}.

Next we move to the correspondence between HORN-SAT and the minimum operation.

\begin{thm} Suppose that a relation $R \subseteq \{0,1\}^m$ is preserved by the minimum function $\min : \{0,1\}^2 \rightarrow \{0,1\}$. Then $R$ can be written as a conjunction of Horn clauses.
\end{thm}
\begin{proof} Write $R = \bigwedge_i C_i$ in conjunctive normal form, such that each clause $C_i$ is minimal. Note that this means that for each literal $l$ in $C_i$, there is some assignment to the variables that satisfies $R$, and has the rest of the literals in $C_i$ other than $l$ set to $0$.

Suppose, for a contradiction, that some clause $C_i$ has at least two non-negated variables in it, and assume without loss of generality that $C_i = x_1 \vee \cdots \vee x_p \vee \bar{x}_{p+1} \vee \cdots \vee \bar{x}_{p+q}$, $p \ge 2$. By the minimality of $C_i$, there are assignments $a, b$ which satisfy $R$ and such that $a_2 = \cdots = a_p = \bar{a}_{p+1} = \cdots = \bar{a}_{p+q} = 0$ and $b_1 = \cdots = b_{p-1} = \bar{b}_{p+1} = \cdots = \bar{b}_{p+q} = 0$. But then $\min(a,b)$ fails to satisfy $C_i$, and hence fails to satisfy $R$.
\end{proof}

\begin{thm} Suppose that a function $f : \{0,1\}^k \rightarrow \{0,1\}$ preserves the relations $\{0\}, \{1\},$ and $x\wedge y \implies z$. Then there is a nonempty subset $I \subseteq \{1, ..., k\}$ such that $f(x_1, ..., x_k) = \min_{i \in I} x_i$.
\end{thm}
\begin{proof} Since $f$ preserves $\{0\}$ and $\{1\}$, we have $f(0,...,0) = 0$ and $f(1, ..., 1) = 1$. Since $\le$ is in the relational clone generated by $x \wedge y \implies z$, $f$ must be monotone.

For each subset $I \subseteq \{1, ..., k\}$, let $\chi_I$ be the indicator vector of $I$. Suppose that $I,J$ have $f(\chi_I) = f(\chi_J) = 1$, then from $\chi_I \wedge \chi_J \implies \chi_{I\cap J}$ (coordinatewise) we see that we must have $f(\chi_{I\cap J}) = 1$ as well. Thus, there is a unique minimum subset $I^*$ satisfying $f(\chi_{I^*}) = 1$. Since $f$ is monotone, we have $f(\chi_J) = 1 \iff J \supseteq I^*$.
\end{proof}

\begin{rem} The fact that $\min$-closed relations on the domain $\{0,1\}$ can always be written as intersections of Horn clauses has the following useful consequence in logic.

Suppose that $P_1, ..., P_m$ is a list of logical statements about some type of structure $M$ in some collection of structures $\cM$. Suppose that for every pair of structures $M_1, M_2 \in \cM$ there is a structure $M' \in \cM$ such that for each $i$, $P_i$ holds in $M'$ iff $P_i$ holds in both $M_1$ and $M_2$. Then there is a collection of Horn clauses $\phi_1, ..., \phi_n$ in the propositions $P_1, ..., P_m$ such that an assignment of true/false values to the $P_i$s can be realized by some $M \in \cM$ iff the assignment satisfies the collection of Horn-clauses $\phi_1, ..., \phi_n$.
\end{rem}

Finally, we come to the affine linear case. We leave the proofs of the following two results to the reader.

\begin{thm} Suppose that a relation $R \subseteq (\ZZ/p)^m$ is preserved by the ternary operation $x - y + z \pmod{p}$. Then $R$ is an affine linear subspace of $(\ZZ/p)^m$ - that is, a vector subspace of $(\ZZ/p)^m$ offset by a fixed vector - and $R \in \langle \{1\}, x+y\equiv z\pmod{p}\rangle$.
\end{thm}

\begin{thm} Suppose that a function $f : (\ZZ/p)^k \rightarrow \ZZ/p$ preserves the relations $\{1\}$ and $x+y \equiv z \pmod{p}$. Then $f$ is an affine linear function - that is, a linear function such that the sum of the coefficients is $1$ - and $f \in \langle x-y+z \pmod{p} \rangle$. If $p$ is odd, we have $\langle x-y+z \pmod{p} \rangle = \langle \frac{x+y}{2} \pmod{p} \rangle$.
\end{thm}

\section{Varieties, Birkhoff's HSP theorem, and the hardness proof}

From here on we switch over to the algebraic language. To a relational structure $\fA = (D, \Gamma)$ we associate an algebraic structure $\bA = (D, \mathcal{O})$ with $\langle \cO \rangle = \Pol(\Gamma)$. We let $\CSP(\bA)$ be another name for $\CSP(\fA) = \CSP(\Inv(\bA))$.

\begin{rem} Suppose $\bA, \bB$ are two algebraic structures with associated relational structures $\fA, \fB$. It is tempting to think that a homomorphism $\bA \rightarrow \bB$ will correspond to a homomorphism $\fA \rightarrow \fB$, or vice versa. Unfortunately, this is total nonsense - if the (functional) signatures of $\bA$ and $\bB$ match, the (relational) signatures of $\fA$ and $\fB$ will likely have nothing to do with each other!

In a similar vein, the automorphism groups $\Aut(\bA)$ and $\Aut(\fA)$ have almost nothing to do with each other. A trivial but illuminating example is the case where $\bA$ has no functions at all, so that $\Aut(\bA)$ is the full symmetric group - in this case, $\fA$ has every possible relation in its signature, including named singleton unary relations for every element of the domain. Thus, if $\bA$ is trivial, then $\fA$ is \emph{rigid}, with $\Aut(\fA) = \{\Id\}$.
\end{rem}

We will now use the algebraic language to relate the complexities of CSPs with different domains. This will finally clarify what we meant by one CSP ``simulating'' another CSP in the introduction (well, there is one more method of simulation that will be introduced in the next section).

\begin{thm} If $\bA$ is an algebraic structure, and $\bB$ is either
\begin{itemize}
\item a subalgebra of $\bA$,

\item a power of $\bA$, or

\item a quotient of $\bA$,
\end{itemize}
then there is a logspace reduction from $\CSP(\bB)$ to $\CSP(\bA)$.
\end{thm}
\begin{proof} If $\bB$ is a subalgebra of $\bA$, we can convert any instance of $\CSP(\bB)$ into an instance of $\CSP(\bA)$ by simply adding an extra unary constraint for each variable corresponding to the relation $\bB \subseteq \bA^1$.

If $\bB = \bA^n$ for some $n$, then we can convert an instance of $\CSP(\bB)$ to an instance of $\CSP(\bA)$ by replacing each variable with an $n$-tuple of variables, and using the fact that every subalgebra of $(\bA^n)^m$ is a subalgebra of $\bA^{mn}$.

If $\bB = \bA/\!\!\sim$ for some congruence $\sim\ \subseteq \bA^2$ on $\bA$, then every relation $R \subseteq \bB^m$ lifts to a relation $\tilde{R} \subseteq \bA^m$ by the rule $x \in \tilde{R} \iff x/\!\!\sim\ \in R$.
\end{proof}

More generally, if we have several algebras $\bA_1, \bA_2, ...$ in the same (functional) signature, we can define $\CSP(\{\bA_1, \bA_2, ...\})$ to be the problem where each variable comes with a \emph{sort} - that is, a specific algebra $\bA_i$ that it lives in - and each relation is \emph{multisorted}, where a multisorted relation is ``allowed'' if it cuts out a subalgebra of the relevant product of the $\bA_i$s. This sort of multisorted relation was considered by Bulatov and Jeavons \cite{bulatov-jeavons-varieties}. In this framework, there is a logspace equivalence between $\CSP(\bA_1 \times \bA_2)$ and $\CSP(\{\bA_1, \bA_2\})$.

So we see that we are naturally led to study families of finite algebras (all sharing a signature) which are closed under taking finite products, subalgebras, and quotients. This leads us to the concept of a \emph{variety} (or \emph{pseudovariety}, if the family of finite algebras is not finitely generated). Lurking in the background here is a new Galois connection, this time between families of \emph{identities} and families of algebras.

\begin{defn} A \emph{term} (in a given functional signature) is either a variable name or a $k$-ary function symbol applied to a $k$-tuple of previously constructed terms. An \emph{identity} is a formal expression $s \approx t$, where $s$ and $t$ are terms. An algebra $\bA$ \emph{satisfies} an identity $s \approx t$, written $\bA \models s \approx t$, if
\[
\forall x_1, ..., x_n \in \bA\  s(x_1, ..., x_n) = t(x_1, ..., x_n)
\]
(here we are assuming that the variables of $s$ and $t$ are drawn from $x_1, ..., x_n$).
\end{defn}

The $\approx$ notation is confusing at first, since in the context of universal algebra it is viewed as a statement which is \emph{stronger} than ordinary equality. The idea here is that approximate equality is never considered in universal algebra, so there should be no confusion in repurposing the symbol $\approx$ into an abbreviation for universal quantifiers. For instance, the intended meaning of the expression ``$f(x,y) \approx f(y,x)$'' is ``$\forall x,y\ f(x,y) = f(y,x)$''. An alternate point of view is that $\approx$ refers to the congruence on the absolutely free algebra corresponding to the identities which are satisfied by the algebras we are interested in.

\begin{defn} The \emph{variety} $\cV(\mathcal{T})$ determined by a set of identities $\mathcal{T}$ is the set of algebras that satisfy all of the identities in $\mathcal{T}$. If $\bA_1, \bA_2, ...$ is a collection of algebras, then $\cV(\bA_1, \bA_2, ...)$ is the variety associated to the set of all identities that hold simultaneously in all of the algebras $\bA_i$.
\end{defn}

Birkhoff \cite{birkhoff} introduced a convenient notation for manipulating sets (strictly speaking these are classes, not sets) of algebras: if $\cS$ is a set of algebras, then $P\cS$ is the set of products of algebras from $\cS$ (possibly infinite - $P_{fin}$ is the notation if one restricts to finite products), $S\cS$ is the set of subalgebras of algebras from $\cS$, and $H\cS$ is the set of quotients (homomorphic images) of algebras from $\cS$.

\begin{thm}[Birkhoff's HSP Theorem \cite{birkhoff}] For any collection $\cS$ of algebras, we have $\cV(\cS) = HSP(\cS)$, that is, an algebra $\bA$ satisfies every identity which is satisfied in every element of $\cS$ if and only if it is the homomorphic image of a subalgebra of a product of elements of $\cS$. Furthermore, if $\cS$ is a finite collection of finite algebras, then the set of finite algebras in $\cV(\cS)$ is equal to $HSP_{fin}(\cS)$.
\end{thm}
\begin{proof} It is easy to check that if $\bA, \bB \models s \approx t$, then $\bA \times \bB \models s \approx t$, and similarly every subalgebra and quotient of $\bA$ also satisfies $s \approx t$. Thus $HSP(\cS) \subseteq \cV(\cS)$.

For the other containment, suppose that $\bA \in \cV(\cS)$, and suppose that $\bA$ is generated by a subset $I \subseteq \bA$. We let $\mathbb{P}$ be the product of all the algebras of $\cS$, and define the ``free algebra'' $\cF(I)$ to be the subalgebra of $\mathbb{P}^{\mathbb{P}^I}$ which is generated by the projection functions $\pi_i$ for $i \in I$, given by $\pi_i(x) = x_i$. We claim that there is a surjective homomorphism $h : \cF(I) \twoheadrightarrow \bA$ with $h(\pi_i) = i$.

Suppose not. Then there are two terms $s,t$ with $s(\pi_{i_1}, ..., \pi_{i_n}) = t(\pi_{i_1}, ..., \pi_{i_n})$ in $\mathbb{P}^{\mathbb{P}^I}$, but $s(i_1, ..., i_n) \ne t(i_1, ..., i_n)$ in $\bA$. But then $s \approx t$ is satisfied by $\mathbb{P}$, and hence by every algebra in $\cS$, and is not satisfied in $\bA$, contradicting our assumption that $\bA \in \cV(\cS)$.

For the last claim, note that if $\bA, \cS$, and every element of $\cS$ are finite, then so are $I, \mathbb{P}, \mathbb{P}^I$, and $\mathbb{P}^{\mathbb{P}^I}$.
\end{proof}

Birkhoff's HSP theorem gives one half of the Galois connection between identities and algebras. The other half is a result from model theory, which explains why elementary results in algebra can always be proved by writing down a long string of equalities.

\begin{thm} If $\mathcal{T}$ is a family of identities, then the set of identities which hold in $\cV(\mathcal{T})$ is equal to the closure of $\mathcal{T} \cup \{x \approx x\}$ under:
\begin{itemize}
\item substituting a term for a variable in an identity,
\item applying a $k$-ary function to both sides of a $k$-tuple of identities,
\item deducing $s\approx t$ from $t \approx s$, and
\item deducing $s \approx u$ from $s \approx t$ and $t \approx u$.
\end{itemize}
\end{thm}
\begin{proof} Define the free algebra $\cF_{\mathcal{T}}(x_1, ...)$ on countably many variables by taking the set of all terms on these variables, and then taking the quotient of this term algebra by the congruence generated by the images under all possible substitutions of the identities in $\mathcal{T}$. The result will be an algebra satisfying all of the identities of $\mathcal{T}$, and one can check directly from the definition of a congruence that the identities that hold in this free algebra are exactly the ones described in the theorem statement.
\end{proof}

Using Birkhoff's theorem, we can give a criterion for NP-completeness.

\begin{thm} If $\CSP(\bA)$ is \emph{not} NP-complete, then there is a finite set of identities $s_i \approx t_i$ which are satisfied by $\bA$, which can't be satisfied by assigning each function symbol to a projection of the same arity.
\end{thm}
\begin{proof} If $\cV(\bA)$ contains an algebra $\bB$ of size at least $2$ where each function symbol acts as a projection, then $\CSP(\bB)$ is NP-complete and has a logspace reduction to $\CSP(\bA)$. Such an algebra $\bB$ will exist if there is a way to assign the function symbols to projections that satisfies \emph{every} identity satisfied by $\bA$. To see that we only have to consider finite sets of identities, we use a compactness argument: each function symbol has only a finite number of projections it can be assigned to, so we can apply K\H{o}nig's Lemma.
\end{proof}

\begin{ex} Consider the algebra $\bA = (\{0,1\}, \min)$, and use the binary function symbol $s$ to abbreviate $\min$. Then it turns out that
\[
\cV(\bA) = HSP(\bA) = SP(\bA) = \cV(\mathcal{T}_{semi}),
\]
where $\mathcal{T}_{semi}$ is the following set of identities:
\[
s(x,x) \approx x, \;\;\; s(x,y) \approx s(y,x), \;\;\; s(x,s(y,z)) \approx s(s(x,y),z).
\]
The second identity above, $s(x,y) \approx s(y,x)$, can't be satisfied by assigning $s$ to either of the projections $\pi_1, \pi_2$.

An algebra in $\cV(\mathcal{T}_{semi})$ is called a \emph{semilattice}, and can be visualized as a poset where every nonempty finite subset has a greatest lower bound (if we visualize it this way, we often call it a \emph{meet} semilattice).

Any finite meet semilattice which has a greatest element can be extended to a lattice, since every finite subset will also have a least upper bound (just take the greatest lower bound of the collection of all upper bounds, which is nonempty by the assumption that there is a greatest element). Since we can always adjoin a new ``top'' element to any finite meet semilattice, we see that every finite semilattice is isomorphic to a subalgebra of the meet semilattice reduct of some lattice (in fact, this is also true for infinite semilattices).

Alternatively, a semilattice can be thought of as a poset where every nonempty finite subset has a least upper bound, if we are thinking in terms of an operation like $\max$ - if we are visualizing it in this way, we call it a \emph{join} semilattice. (I generally prefer to visualize semilattices as join semilattices, but most authors prefer to visualize semilattices as meet semilattices.)

Since it is often confusing when people who think of semilattices as meet semilattices try to talk to people who think of them as join semilattices (i.e. minimal elements in one language become maximal elements in the other language), it is useful to have some vocabulary which is agnostic to the meet/join distinction. We say that an element $a$ is \emph{absorbing} with respect to $s$ if it satisfies
\[
s(a,x) = s(x,a) = a
\]
for all $x$, and we say that an element $b$ is \emph{neutral} with respect to $s$ if it satisfies
\[
s(b,x) = s(x,b) = x
\]
for all $x$. Every two-element semilattice has a neutral element and an absorbing element, and knowing which is which determines the semilattice operation. In general, every finite semilattice has an absorbing element, but might not have a neutral element (for instance, the free semilattice on two generators has absorbing element $s(x,y)$ and has no neutral element). In a meet semilattice, the absorbing element will be the bottom and the neutral element (if it exists) will be the top, while in a join semilattice, the absorbing element will be the top and the neutral element (if it exists) will be the bottom.
\end{ex}

\begin{ex} Consider the algebra $\bA = (\{0,1\}, \maj)$, and use the ternary function symbol $m$ to abbreviate $\maj$. Then it turns out that
\[
\cV(\bA) = HSP(\bA) = SP(\bA) = \cV(\mathcal{T}_{med}),
\]
where $\mathcal{T}_{med}$ is the following set of identities:
\begin{align*}
m(x,y,z) &\approx m(y,z,x) \approx m(x,z,y),\\
m(x,x,y) &\approx x,\\
m(m(x,y,z),u,v) &\approx m(x,m(y,u,v),m(z,u,v)).
\end{align*}
The identity $m(x,y,z) \approx m(y,z,x)$ can't be satisfied by assigning $m$ to one of the projections $\pi_1, \pi_2, \pi_3$.

An algebra in $\cV(\mathcal{T}_{med})$ is called a \emph{median algebra}. A finite median algebra corresponds to a \emph{median graph}, that is, a graph with the property that for every three vertices $x,y,z$ there exists a unique vertex which lies on a shortest path connecting every pair of $x,y,z$ (to recover the graph structure, we draw an edge from $x$ to $y$ whenever $m(x,y,z) \in \{x,y\}$ for all $z$). Examples of median graphs are paths, trees, planar grids, ``squaregraphs'', hypercubes, Hasse diagrams of distributive lattices, and the graph $\cM_n$ of maximal intersecting families from the last section. For more about the theory of median algebras, see \cite{median-poc} or \cite{bowditch-median}.

Median algebras are very closely connected to distributive lattices. It isn't hard to show that in any distributive lattice, the following identity holds:
\[
(x \wedge y) \vee (y \wedge z) \vee (z \wedge x) \approx (x \vee y) \wedge (y \vee z) \wedge (z \vee x),
\]
and in fact this identity is \emph{equivalent} to the lattice being distributive. The common value of both sides is called the median operation $m(x,y,z)$ on the lattice - the reader can easily check that it satisfies the identities $\mathcal{T}_{med}$. In fact, if a median algebra has two elements $0,1$ with $m(0,x,1) = x$ for all $x$, then it forms a distributive lattice under the operations $x \wedge y = m(0,x,y)$ and $x\vee y = m(x,y,1)$, and the median operation $m$ can be recovered from $\wedge, \vee$ via the above formula \cite{median-distributive}.
\end{ex}

\begin{ex} The operation $m(x,y,z) = x - y + z \pmod{p}$ satisfies the identity
\[
m(x,y,y) \approx x \approx m(y,y,x),
\]
and this identity can't be satisfied by assigning $m$ to one of the projections $\pi_1, \pi_2, \pi_3$.

Similarly, if $p$ is odd, the operation $m(x,y) = \frac{x+y}{2} \pmod{p}$ satisfies the identity
\[
m(x,y) \approx m(y,x),
\]
which can't be satisfied by projections.
\end{ex}

As with clones and relational clones, there are several natural finiteness questions that come up with varieties.

\begin{defn} A variety $\cV$ is \emph{finitely generated} if there is a finite list of finite algebras $\bA_1, ..., \bA_n$ such that $\cV = \cV(\bA_1, ..., \bA_n)$. A variety $\cV$ is \emph{locally finite} if the free algebra on $n$ generators $\cF_{\cV}(x_1, ..., x_n)$ is finite for every $n$. A variety $\cV$ is \emph{finitely based} if there is a finite set of equations $\mathcal{T}$ such that $\cV = \cV(\mathcal{T})$.
\end{defn}

A variety $\cV$ is locally finite iff for all $\bA \in \cV$ and for all finite subsets $\{a_1, ..., a_n\} \subseteq \bA$, the subalgebra of $\bA$ generated by $a_1, ..., a_n$ is finite. Every finitely generated variety is locally finite (by the proof of the HSP Theorem). In general, determining whether a given finitely generated variety is finitely based, or vice versa, is a very difficult problem. For instance, the famous Burnside problem is the problem of determining whether the variety of groups satisfying the identity $x^n \approx e$ is locally finite.

\begin{rem} Sometimes we want to consider infinite families of finite algebras with a finite functional signature, closed under \emph{finite} products, subalgebras, and homomorphisms. Such a family of algebras is called a \emph{pseudovariety}. There are two different ways to describe pseudovarieties in terms of identities.

Eilenberg and Sch{\"u}tzenb{\'e}rger \cite{pseudovarieties-ultimate} show that a pseudovariety is determined by an infinite sequence of identities, such that a finite algebra is contained in the pseudovariety iff it satisfies \emph{all but finitely many} of the identities in the sequence. The trick is to sort the isomorphism classes of finite algebras by their sizes, and for each size $k$ write down a finite set of identities in $k$ variables which characterizes the free algebra on $k$ generators in the subvariety generated by the set of algebras of size at most $k$.

Reiterman \cite{pseudovarieties-implicit} shows that a pseudovariety is determined by identities between ``implicit operations'': operations which aren't defined from terms directly, but which can still be defined on any particular finite algebra in a way that is compatible with homomorphisms. Examples of implicit operations in the language of a unary function $f$ are
\[
f^\infty = \lim_{n \rightarrow \infty} f^{\circ n!}, \;\;\; f^{\infty-1} = \lim_{n \rightarrow \infty} f^{\circ (n!-1)},
\]
where the limits are taken pointwise (note that the functions $f^{\circ n!}$ stabilize once $n$ exceeds the size of the domain). For any function $f$ on a finite domain, $f^\infty$ will always satisfy the identity $f^\infty(f^\infty(x)) \approx f^\infty(x)$, while the pseudovariety of \emph{invertible} functions on finite sets is cut out by the identities
\[
f(f^{\infty-1}(x)) \approx f^{\infty-1}(f(x)) \approx x.
\]

For those who like category theory, a $k$-ary implicit operation of a pseudovariety $\cV$ with underlying set functor $S : \cV \rightarrow \text{Set}$ is a natural transformation from $S^k$ to $S$. If a free algebra on $k$ elements exists in $\cV$, then a standard argument shows that every $k$-ary implicit operation of $\cV$ is actually \emph{explicit}, that is, a term operation of $\cV$. In general, every finite subset of $\cV$ will generate a locally finite subvariety of $\cV$, which shows that the restriction of any implicit operation to this subset agrees with some term operation of $\cV$. Reiterman \cite{pseudovarieties-implicit} puts a metric structure on the set of implicit operations of a pseudovariety such that the collection of implicit operations becomes the completion of the collection of explicit operations.
\end{rem}


\section{Cores and Idempotent Reducts}

In this section we briefly return to the relational point of view, and the concept of homomorphic equivalence, to provide one last algebraic ingredient: the restriction to \emph{idempotent} algebraic operations.

\begin{defn} Two relational structures $\mathbf{A}, \mathbf{B}$ with the same signature are \emph{homomorphically equivalent} if there exist homomorphisms $\mathbf{A} \rightarrow \mathbf{B}, \mathbf{B} \rightarrow \mathbf{A}$.
\end{defn}

The prototypical example of homomorphic equivalence is a (non-surjective) endomorphism from a relational structure to itself, providing a homomorphic equivalence between the original relational structure and the restriction of the relational structure to a proper subset of its domain. On the algebraic side, this manifests as a unary operation which is not invertible. The algebraic implications of such unary operations in the polynomial clone of an algebra are at the heart of the subject called ``tame congruence theory'', which was introduced to give the first structure theory for finite algebras in the book by Hobby and McKenzie \cite{hobby-mckenzie}.

\begin{ex} Consider the relational structure $\fA$ corresponding to the binary implication algebra $\bA = (\{0,1\}, \to)$. This relational structure has as basic relations $R_n = \{0,1\}^n \setminus \{(0, ..., 0)\} = x_1 \vee \cdots \vee x_n$. The unary algebraic operation $\operatorname{\to}(x,x)$ of $\bA$ takes every element to $1$, and defines an endomorphism of relational structures $\fA \rightarrow \fA$ whose image is $\{1\}$. Together with the inclusion relation, we get a homomorphic equivalence between $\fA$ and the one-element relational structure with domain $\{1\}$ and relations $R_n \mid_{\{1\}^n} = \{1\}^n$, whose CSP is clearly trivial.
\end{ex}

As the example shows, non-surjective endomorphisms provide trivial ways to simplify CSPs.

\begin{defn} A relational structure on a finite domain $\fA$ is called a \emph{core} if every endomorphism of $\fA$ is also an automorphism of $\fA$. If $\fA$ is not a core, then $\fB$ is called a \emph{core of} $\fA$ if $\fB$ is a core and $\fB$ is homomorphically equivalent to $\fA$.
\end{defn}

\begin{ex} Every complete graph $K_n = ([n], \ne)$ is a core.
\end{ex}

\begin{ex} If $G$ is a bipartite graph with at least one edge, then the core of $G$ is $K_2$, the complete graph on two vertices.
\end{ex}

\begin{rem} In the infinite case, the definition of a core must be modified: an infinite relational structure is called a core if every endomorphism is an \emph{embedding}, i.e. an injective map that is an isomorphism onto the restriction of the target relational structure to its image. An example of an infinite core is $(\mathbb{Q},<)$. See section 3.6 of \cite{bodirsky-thesis} for more information about cores of infinite structures.
\end{rem}

\begin{prop} Every relational structure on a finite domain $\fA$ has a core. Any two cores of $\fA$ are isomorphic.
\end{prop}
\begin{proof} The first statement follows directly from induction on the size of $\fA$: if $\fA$ is not a core, then it is homomorphically equivalent to its restriction to some proper subset of itself. For the second statement, note that if $\fB, \fB'$ are two cores of $\fA$ then they are homomorphically equivalent, and composing the maps $\fB \rightarrow \fB'$, $\fB' \rightarrow \fB$ gives us endomorphisms of $\fB, \fB'$ which must both be invertible by the definition of a core.
\end{proof}

Note that although restricting our attention to cores seems like a trivial step, we are sweeping the following problem under the rug.

\begin{prob} Given a finite relational structure $\fA$ as input, determine whether or not $\fA$ is a core.
\end{prob}

Obviously there is a brute-force approach to checking if $\fA$ is a core: simply write down every possible endomorphism, and go through them one by one. Since we only have to do this brute force once for a given CSP template, this is not as bad as it sounds, but it is still far from ideal. Unfortunately, as it turns out, a brute force approach is pretty much the best one can do.

\begin{thm} [Hell, Ne\v{s}et\v{r}il \cite{core-graph}] The problem of determining whether or not a given undirected graph is not a core is NP-complete, even if the graph is assumed to be $3$-colorable (with a given $3$-coloring).
\end{thm}

The next main idea comes from ``self-reducibility'': often, when solving a CSP, one makes a guess (or deduces) that a certain variable should have a certain value. We would like to be able to express a CSP together with some constraints stating that certain variables have certain values using the language of the original CSP. If this is possible, then an algorithm for deciding whether the CSP has a solution can be directly converted into an algorithm for \emph{finding} a solution to the CSP.

\begin{defn} A relational structure $\fA$ is a \emph{rigid core} if it has no endomorphisms other than the identity. (In general, a structure is called \emph{rigid} if it has no automorphisms.)
\end{defn}

\begin{thm} A relational structure $\fA$ on a finite domain $D$ is a rigid core if and only if it has the following property: for every element $a \in D$, the unary relation $\{a\}$ is contained in the relational clone generated by the relations of $\fA$.
\end{thm}
\begin{proof} This follows directly from the $\Inv$-$\Pol$ Galois connection: $\{a\} \in \langle \fA \rangle$ iff $\{a\}$ is closed under $\Pol(\fA)$, and since $\{a\}$ is generated by a single element, we only need to check that it is closed under $\Pol_1(\fA)$, which is exactly the set of endomorphisms of $\fA$.

We can also give a direct proof, by unraveling the proof of the $\Inv$-$\Pol$ connection in this special case, as follows. Define a CSP with a variable $f_a$ for each $a \in D$. For every relation $R \subseteq D^m$ of $\fA$ and every tuple $(a_1, ..., a_m) \in R$, we impose the constraint $(f_{a_1}, ..., f_{a_m}) \in R$ on our CSP. Now the solution-set to our CSP exactly corresponds to the set of endomorphisms of $\fA$, and if $\fA$ is a rigid core then existentially projecting onto the variable $f_a$ produces the unary relation $\{a\}$.
\end{proof}

So it is very desirable to restrict our attention to rigid cores. Most of the example CSPs from the introduction were rigid cores, with the notable exceptions of $k$-coloring and NAE-SAT. The $k$-coloring problem is an excellent toy example: the reader may be already be aware of the fact that $\CSP(\{1, ..., k\}, \ne)$ (the $k$-coloring problem) is logspace equivalent to $\CSP(\{1,...,k\}, \ne, \{1\}, ..., \{k\})$ - the rigid core obtained by adjoining the unary singleton relations to $k$-coloring. It is worth examining the proof of that equivalence and understanding how the next result generalizes it.

\begin{thm}\label{rigid-core-reduction} Suppose that $\fA = (D, \Gamma)$ is a core on a finite domain $D$, and let $\fA^{rig}$ be the rigid core obtained by adjoining all singleton unary relations to $\fA$. Then $\CSP(\fA)$ is equivalent to $\CSP(\fA^{rig})$ under logspace reductions.
\end{thm}
\begin{proof} We need to find a way to convert an instance of $\CSP(\fA^{rig})$ to an instance of $\CSP(\fA)$ without changing whether it has a solution. As in the previous result, introduce a set of variables $f_a$ for each element $a \in D$, and define a primitive positive formula $\Phi$ by
\[
\Phi(f) \coloneqq \bigwedge_{R \in \Gamma} \bigwedge_{(a_1, ..., a_m) \in R} R(f_{a_1}, ..., f_{a_m}).
\]
Suppose that our instance of $\CSP(\fA^{rig})$ has the form
\[
\Psi(x) = \exists x_{n+1}, ..., x_{n+m}\ \Psi_0(x) \wedge \bigwedge_{(i,a) \in E} x_i \in \{a\},
\]
where $\Psi_0$ is a primitive positive formula using the relations of $\Gamma$, and $E$ is a set describing the additional unary singleton constraints of $\Psi$. Let $\Psi'$ be the following formula:
\[
\Psi'(x) \coloneqq \exists f\ \exists x_{n+1}, ..., x_{n+m}\ \Phi(f) \wedge \Psi_0(x) \wedge \bigwedge_{(i,a) \in E} x_i = f_{a}.
\]
We claim that the instance $\Psi'$ of $\CSP(\fA)$ has a solution iff the instance $\Psi$ of $\CSP(\fA^{rig})$ has a solution. Suppose that $f,x$ solves $\Psi'$, then by the construction of $\Phi(f)$ $f$ describes an endomorphism $f : \fA \rightarrow \fA$, and since $\fA$ is a core this endomorphism must have an inverse $f^{-1}$. Then $f^{-1}(x)$ satisfies $\Psi_0$ (since $f^{-1}$ is an endomorphism of $\fA$), and for $(i,a) \in E$ we have $f^{-1}(x_i) = f^{-1}(f_{a}) = a$, so $f^{-1}(x)$ is a solution to $\Psi(x)$.
\end{proof}

Now we look at what the restriction to rigid cores means on the algebraic side.

\begin{defn} A function $f : D^k \rightarrow D$ is \emph{idempotent} if it satisfies the identity $f(x, x, ..., x) \approx x$. An algebraic structure $\bA = (D, \cO)$ is \emph{idempotent} if every $f \in \cO$ is idempotent. Equivalently, $\bA$ is idempotent if every singleton subset of $D$ is a subalgebra of $\bA$.
\end{defn}

\begin{defn} If $\bA = (D, \cO)$ is an algebraic structure, then the \emph{idempotent reduct} $\bA^{id}$ of $\bA$ has the same domain, and has as its operations the set of all idempotent functions $f \in \langle \cO \rangle$ (or, alternatively, some smaller generating set of idempotent functions).
\end{defn}

\begin{ex} If $\bA = (\ZZ/p, +, 0, 1)$, then $\bA^{id}$ has as its operations the set of all affine linear functions on $\ZZ/p$, and one can take $\{x-y+z \pmod{p}\}$ as a generating set of basic operations (or, if $p$ is odd, one can alternatively take $\{\frac{x+y}{2} \pmod{p}\}$ as a generating set of basic operations).
\end{ex}

\begin{prop} If $\fA$ is a core corresponding to the algebraic structure $\bA$, then the rigid core $\fA^{rig}$ corresponds to the idempotent reduct $\bA^{id}$. In particular, every CSP is equivalent up to logspace reductions to $\CSP(\bA)$ for some idempotent algebra $\bA$.
\end{prop}

The reader might be worried that there is no obvious way to generate the collection of all idempotent operations contained in a given clone. For core structures this is not an issue: the polymorphisms of a core structure always decompose neatly into idempotent parts and invertible unary parts.

\begin{prop} Suppose that $\cO$ is a clone such that all of the unary operations in $\cO$ are invertible. Then for every $k$-ary function $f \in \cO$, if we define the unary function $f_{un}$ by
\[
f_{un}(x) \coloneqq f(x,...,x)
\]
and the $k$-ary function $f_{id}$ by
\[
f_{id}(x_1,...,x_k) \coloneqq f_{un}^{-1}(f(x_1, ..., x_k)),
\]
then $f_{id}$ is idempotent and
\[
f = f_{un} \circ f_{id}.
\]
In particular, if $G$ is the group of unary operations in $\cO$, then for every $k$ there are precisely $|G|$ times as many $k$-ary operations in $\cO$ as there are $k$-ary idempotent operations in $\cO$.

If $\cO$ is generated by the functions $f_1, ..., f_m$ of arities $k_1, ..., k_m$, then the set of idempotent operations of $\cO$ is generated by the functions
\[
(f_i \circ (g_{1}, ..., g_{k_i}))_{id},
\]
over all choices of $i$ and all choices of $g_1, ..., g_{k_i} \in G$. In particular, the set of idempotent operations of $\cO$ is finitely generated if and only if the full clone $\cO$ is finitely generated.
\end{prop}

\begin{ex} There is an example of a core structure $\fA$ which has polymorphisms satisfying a nontrivial system of identities, but such that its rigidification $\fA^{rig}$ has no such polymorphisms and is therefore NP-complete. This example is due to Ross Willard and can be found in \cite{topology-irrelevant}.

The underlying set of $\fA$ is the set of expressions $a_i$ with $a \in \{1,2,3\}$ and $i \in \{0,1\}$. The relations of $\fA$ are given by
\begin{align*}
R(a_i,b_j) &\coloneqq (i = j) \wedge (a \ne b),\\
S(a_i,b_j) &\coloneqq i \ne j.
\end{align*}
It is easy to check that this structure is a core.

Polymorphisms of $\fA$ include the unary automorphism $\alpha(a_i) = a_{1-i}$ and the ternary function $s$ given by
\[
s(a_i,b_j,c_k) = \begin{cases}c_k & i = j,\\ a_i & i \ne j.\end{cases}
\]
These polymorphisms satisfy the identity
\[
s(x,x,y) \approx s(y,\alpha(y),x) \approx y,
\]
which can't be satisfied by projections.

Since the unary relation $\{a_i \mid i = 0\}$ is definable in $\fA^{rig}$, we see that polymorphisms of $\fA^{rig}$ restrict to idempotent polymorphisms of the triangle $K_3$. We will show that $K_3$ has no nontrivial idempotent polymorphisms: in fact, we'll show that every polymorphism of $K_3$ is the composition of a projection with an automorphism of $\{1,2,3\}$.

To see that all polymorphisms of $K_3$ are essentially unary, suppose that $f : K_3^n \rightarrow K_3$ depends nontrivially on its first coordinate, that is, that there are $x,y \in K_3^n$ with $x_i = y_i$ for all $i > 1$ with $f(x) \ne f(y)$. By composing with automorphisms of $\{1,2,3\}$, we may assume without loss of generality that
\[
f(1,1,...,1) = 1, f(2,1,...,1) = 2.
\]
Since $f$ preserves the $\ne$ relation, we must then have
\[
f(3,2,...,2) = f(3,3,...,3) = 3.
\]
These imply that
\[
f(\{1,2\}^n) \subseteq \{1,2\}, f(\{1,2\}\times \{1,3\}^{n-1}) \subseteq \{1,2\}.
\]
For any $z_2, ..., z_n$, we can find $x_2, ..., x_n \in \{1,2\}$ and $y_2, ..., y_n \in \{1,3\}$ with $x_i,y_i,z_i$ all distinct. Thus we must have
\[
f(3,z_2, ..., z_n) = 3
\]
for all $z_2, ..., z_n$, and in particular $f(3,1,...,1) = 3$. Now we can repeat the argument with $1$ or $2$ in place of $3$ to see that $f(x_1, ..., x_n) = x_1$ for all $x_1, ..., x_n$, that is, $f = \pi_1$.

Alternatively, we could have shown that $\Pol(K_3^{rig})$ is trivial by instead showing that every relation on $\{1,2,3\}$ is primitively positively definable from the singleton relations together with $\ne$. We leave this as an exercise for the reader (hint: once you have all ternary relations of the form $(x=a)\wedge(y=b) \implies (z=c)$, it's easy to construct the rest).
\end{ex}

\subsection{Reflections and Height 1 Identities}

Let's recap the various methods we have used to reduce different CSPs to each other:
\begin{itemize}
\item Reduce the set of basic relations $\Gamma$ of a relational structure $\fA = (A, \Gamma)$ to some collection of relations $\Gamma' \subseteq \langle \Gamma \rangle$. Equivalently, expand the collection of basic operations $\cO$ in an algebraic structure $\bA = (A, \cO)$ to a collection of operations $\cO'$ with $\cO \subseteq \langle \cO' \rangle$. The collection of algebraic structures $\bA' = (A,\cO')$ with $\cO \subseteq \cO'$ is called the collection of \emph{expansions} of $\bA$, and we use the notation $E(\{\bA\})$ for it in analogy with Birkhoff's HSP operations.

\item Each of Birkhoff's algebraic HSP operations, on the algebraic side: we can replace an algebraic structure $\bA$ by any power $\bA^n$, any subalgebra $\bB \le \bA$, or any quotient $\bA/\!\!\sim$ to get a CSP which is no harder than $\CSP(\bA)$.

\item Homomorphic equivalence of relational structures: when there are homomorphisms $\fA \rightarrow \fB$ and $\fB \rightarrow \fA$, then $\CSP(\fA) = \CSP(\fB)$, since a relational structure $\fX$ will have a homomorphism to $\fA$ iff $\fX$ has a homomorphism to $\fB$. We mainly use this to reduce the general case to the case where $\fA$ is a core relational structure.

\item Starting from a core relational structure $\fA$, we showed that the rigidification $\fA^{rig}$ which we get by adding each singleton unary relation to the basic relations of $\fA$ has $\CSP(\fA^{rig})$ no harder than $\CSP(\fA)$. On the algebraic side, this lets us reduce the general case to the case where every basic operation of $\bA$ is idempotent.
\end{itemize}

Barto, Opr{\v{s}}al, and Pinsker \cite{barto-reflections} find it unsatisfactory to have so many unrelated methods of proving reductions between CSPs, and looked for a single framework which could encompass all known techniques for proving reductions. They show that every single method of proving a reduction between $\CSP(\fA)$ and $\CSP(\fB)$ introduced so far can be described by combining just two basic cases:
\begin{itemize}
\item if $\fB$ is a ``pp-power'' (defined below) of $\fA$, then $\CSP(\fB)$ has a logspace reduction to $\CSP(\fA)$, and

\item if $\fB$ is homomorphically equivalent to $\fA$ then $\CSP(\fB) = \CSP(\fA)$.
\end{itemize}
Furthermore, they show that even if we chain several such reductions together, we can always find an equivalent reduction where the pp-power step is taken before the homomorphic equivalence step.

\begin{defn} A \emph{pp-power} of a relational structure $\fA$ is a relational structure $\fB$ with domain $\fA^n$ for some $n$, such that every relation of $\fB$ can be defined by a primitive positive formula using the relations of $\fA$ (note that the signatures of $\fA$ and $\fB$ will generally be different).
\end{defn}

\begin{prop} If $\fB$ is homomorphically equivalent to a pp-power of $\fA$, then there is a reduction from $\CSP(\fB)$ to $\CSP(\fA)$ which can be computed in linear time and logarithmic space.
\end{prop}

\begin{defn}\label{defn-pp-construct} We say that $\fA$ \emph{pp-constructs} $\fB$ if $\fB$ is homomorphically equivalent to some pp-power of $\fA$.
\end{defn}

For most of the reductions we have described so far, it is easy to see how we can express them in the pp-constructability framework. For instance, in order to simulate the relational structure corresponding to $\bA/\!\!\sim$, we first construct the relational structure $\fA'$ whose basic relations consist of all subpowers $\RR \le \bA^m$ which are compatible with the congruence $\sim$ (using the special case of the pp-power construction where the power $n$ is equal to $1$), and then we follow up with a homomorphic equivalence from $\fA'$ to a relational structure where each equivalence class of $\sim$ is collapsed to a single element. The only really tricky reduction is the reduction to rigid cores, i.e. adding singleton unary relations to a core structure.

Here is how we can go about adding a singleton unary relation $\{a\}$ to a core $\fA$ in the pp-constructability framework. Let $\fB$ be the relational structure which has the new unary relation $\{a\}$ (along with all of the original relations which $\fA$ had). We will define a relational structure $\mathbf{C}$ which will be a pp-power of $\fA$ having domain $\fA^2$, and show that $\mathbf{C}$ is homomorphically equivalent to $\fB$.

Let $O$ be the orbit of $a$ under $\Aut(\fA)$ - note that $O$ is a unary relation contained in the relational clone defined by $\fA$ - and for every $m$-ary relation $R$ of $\fA$, make a corresponding relation $\tilde{R}$ of $\mathbf{C}$ by
\[
((x_1,y_1), ..., (x_m,y_m)) \in \tilde{R} \iff (x_1, ..., x_m) \in R \wedge y_1 = \cdots = y_m \in O.
\]
For the relation $\{a\}$ of $\fB$, we make a corresponding relation $S$ of $\mathbf{C}$ given by
\[
(x,y) \in S \iff x = y \in O.
\]
To show that $\fB$ and $\mathbf{C}$ are homomorphically equivalent, we just need to exhibit a pair of homomorphisms between them. The homomorphism $\fB \rightarrow \mathbf{C}$ is given by $x \mapsto (x,a)$. To define the homomorphism from $\mathbf{C}$ to $\fB$, we need to choose an automorphism $g_y$ of $\fA$ with $g_y(y) = a$ for every $y \in O$. Then the homomorphism $\mathbf{C} \rightarrow \fB$ is given by $(x,y) \mapsto g_y(x)$ if $y \in O$ (and $(x,y)$ maps to an arbitrary element if $y \not\in O$).

Now let's check that pp-constructability is transitively closed.

\begin{prop}[From \cite{barto-reflections}] If $\fA$ pp-constructs $\fB$ and $\fB$ pp-constructs $\fC$, then $\fA$ pp-constructs $\fC$.
\end{prop}
\begin{proof} It's easy to check that homomorphic equivalence is an equivalence relation, and that a pp-power of a pp-power is a pp-power of the original structure. We just need to check that if $\fA$ is homomorphically equivalent to $\fB$ and $\fC$ is a pp-power of $\fB$, then there is some $\fC'$ such that $\fC'$ is a pp-power of $\fA$ and $\fC'$ is homomorphically equivalent to $\fC$.

Suppose that $\fC$ is a pp-power of $\fB$ with power $n$, so that the underlying set $C$ of $\fC$ is equal to $B^n$. We will construct the pp-power $\fC'$ of $\fA$ using the same power $n$, with the same relational signature as $\fC$, as follows. For each $m$-ary relation symbol $R$ of $\fC$, by the definition of a pp-power there is some primitive positive formula in terms of relations $S_i$ of $\fB$ which defines $R$ as an $mn$-ary relation on $B$:
\[
x = (x_1, ..., x_m) = ((x_{11}, ..., x_{1,n}), ..., (x_{m1}, ..., x_{mn})) \in R \iff \exists y_1, ..., y_k \in B \text{ s.t. } \bigwedge_i \pi_{I_i}(x, y) \in S_i.
\]
Then since $\fA$ and $\fB$ have the same relational signature (this is part of our assumption that they are homomorphically equivalent), we can interpret each relation symbol $S_i$ in $\fA$ to define an $mn$-ary relation on $A$, which we will use to give an interpretation of the relation symbol $R$ in $\fC'$:
\[
x = (x_1, ..., x_m) = ((x_{11}, ..., x_{1,n}), ..., (x_{m1}, ..., x_{mn})) \in R^{\fC'} \iff \exists y_1, ..., y_k \in A \text{ s.t. } \bigwedge_i \pi_{I_i}(x, y) \in S_i^{\fA}.
\]
It is now easy to check that any homomorphism $\varphi: \fA \rightarrow \fB$ defines a homomorphism $\varphi^n : \fC' \rightarrow \fC$ by simply letting $\varphi$ act componentwise on elements of $C' = A^n$. Similarly, any homomorphism $\fB \rightarrow \fA$ defines a homomorphism $\fC \rightarrow \fC'$, so $\fC$ and $\fC'$ are homomorphically equivalent.
\end{proof}

Barto, Opr{\v{s}}al, and Pinsker \cite{barto-reflections} also characterize what happens on the algebraic side of the picture when one relates two relational structures by a pp-power or a homomorphic equivalence. The new thing here is really the homomorphic equivalence: if $g : \fA \rightarrow \fB$ and $h : \fB \rightarrow \fA$, then there is a relationship between $\Pol(\fA)$ and $\Pol(\fB)$ which they call a \emph{reflection}, which takes a function $f \in \Pol_k(\fA)$ to the operation
\[
\xi(f) : (x_1, ..., x_k) \mapsto g(f(h(x_1), h(x_2), ..., h(x_k)))
\]
in $\Pol_k(\fB)$. Note that $\xi$ does not respect composition: $\xi(f_0 \circ (f_1, ..., f_k))$ is not in general equal to $\xi(f_0)\circ (\xi(f_1), ..., \xi(f_k))$. However, $\xi$ \emph{does} preserve \emph{height 1 identities}.

\begin{defn} An identity is called a \emph{height 1 identity}, or a \emph{minor identity}, if it has the form $f(x_1, ..., x_k) \approx g(y_1, ..., y_l)$, where the $x_i$s and $y_j$s are (not necessarily distinct) variables. A map $\Pol(\fA) \rightarrow \Pol(\fB)$ (taking functions to functions) which respects height 1 identities is called a \emph{height 1 clone homomorphism} or a \emph{minion homomorphism}.
\end{defn}

\begin{defn} If $\bA = (A, \cO)$ is an algebraic structure and $B$ is a set, and maps $g : A \rightarrow B$, $h : B \rightarrow A$ are given, then the \emph{reflection} of $\bA$ induced by $g,h$ is defined to be the algebraic structure $\bB$ with domain $B$ and the same signature as $\bA$, with the operation $g\circ f \circ (h, ..., h)$ on $B$ corresponding to the operation $f \in \cO$.
\end{defn}

\begin{prop}\label{erp-prop} $\fB$ is homomorphically equivalent to a pp-power of $\fA$ iff $\Pol(\fB)$ contains a reflection of $\Pol(\fA)^n$ for some $n$ (by $\Pol(\fA)^n$ we mean the clone of operations of $\Pol(\fA)$ acting on a power of the \emph{domain}).
\end{prop}
\begin{proof} We prove the non-obvious direction. Let $A,B$ be the underlying sets of $\fA, \fB$, and suppose that $g: A^n \rightarrow B$ and $h: B \rightarrow A^n$ induce a reflection $\xi : \Pol(\fA)^n \rightarrow \Pol(\fB)$. We will construct a pp-power $\mathbf{C}$ of $\fA$ with underlying set $A^n$ which is homomorphically equivalent to $\fB$. For every relation $R$ of $\fB$, let $\tilde{R}$ be the relation
\[
\tilde{R} \coloneqq \{f(h(r_1), ..., h(r_k)) \mid f \in \Pol_k(\fA), r_1, ..., r_k \in R\}.
\]
By definition, $\tilde{R}$ is the closure of $h(R)$ under $\Pol(\fA)$, so $\tilde{R}$ is defined by a primitive positive formula over $\fA$. We use $\tilde{R}$ as the relation corresponding to $R$ in $\mathbf{C}$. Finally, we just need to check that $g : \mathbf{C} \rightarrow \fB$ and $h : \fB \rightarrow \mathbf{C}$ are homomorphisms. That $h$ is a homomorphism follows from $h(R) \subseteq \tilde{R}$. To check that $g$ is a homomorphism, note that if $x = f(h(r_1), ..., h(r_k)) \in \tilde{R}$ with $r_1, ..., r_k \in R$, then $g(x) = \xi(f)(r_1, ..., r_k)$ is an element of $R$ since $\xi(f) \in \Pol(\fB)$ by assumption.
\end{proof}

\begin{thm}[ERP Theorem \cite{barto-reflections}]\label{erp-thm} $\Pol(\fB)$ contains a reflection of $\Pol(\fA)^n$ for some $n$ iff there is a height 1 clone homomorphism $\Pol(\fA) \rightarrow \Pol(\fB)$.
\end{thm}
\begin{proof} We prove the non-obvious direction. Let $\bA = (A,\Pol(\fA))$ be the algebraic structure corresponding to $\fA$, and suppose $\xi : \Pol(\fA) \rightarrow \Pol(\fB)$ is a height 1 clone homomorphism. Let $\cF$ be the subalgebra of $\bA^{\bA^B}$ generated by the operations $\pi_b : \bA^B \rightarrow \bA$ given by $\pi_b : x \mapsto x_b$. Note that $\cF$ is secretly the free algebra over $\bA$ on $|B|$ generators.

Define maps $g : \bA^{\bA^B} \rightarrow B$ and $h : B \rightarrow \bA^{\bA^B}$ by $h(b) = \pi_b$ and
\[
g(f(\pi_{b_1}, ..., \pi_{b_k})) = \xi(f)(b_1, ..., b_k)
\]
for $f \in \Pol_k(\fA)$, and define $g(x)$ arbitrarily for $x \not\in \cF$. To see that $g$ is well-defined, note that if $f_0(\pi_{b_1}, ..., \pi_{b_k}) = f_1(\pi_{c_1}, ..., \pi_{c_l})$, then $f_0, f_1$ are related by a height 1 identity in $\bA$ which implies that
\[
\xi(f_0)(b_1, ..., b_k) = \xi(f_1)(c_1, ..., c_l).
\]
Finally, we see that $g,h$ induce $\xi$ as a reflection from $\bA^{\bA^B}$:
\[
\xi(f)(b_1, ..., b_k) = g(f(\pi_{b_1}, ..., \pi_{b_k})) = g(f(h(b_1), ..., h(b_k))).\qedhere
\]
\end{proof}

As a consequence, we see that the complexity of a CSP only depends on the set of height 1 identities satisfied by its polymorphisms, and that identities involving composition of functions are in a sense superfluous. We also have the following result.

\begin{cor} Let $\fA$ be a relational structure with core $\fB$, and let $\fB^{rig}$ be $\fB$ together with any finite collection of singleton unary relations. Then a system of height 1 identities can be satisfied in $\Pol(\fA)$ iff it can be satisfied in $\Pol(\fB^{rig})$.
\end{cor}

\begin{rem} A height 1 clone homomorphism $\Pol(\fA) \rightarrow \Pol(\fB)$ is completely determined by its restriction to polymorphisms of $\fA$ of arity at most $|B|$, since every operation $f : B^k \rightarrow B$ is determined by its $|B|$-ary minors. So there are only finitely many candidates for height 1 clone homomorphisms from $\Pol(\fA)$ to $\Pol(\fB)$: if the underlying sets are $A,B$, then there are at most
\[
\left|\Pol_{|B|}(\fB)^{\Pol_{|B|}(\fA)}\right| \le |B|^{|B|^{|B|} \cdot |A|^{|A|^{|B|}}}
\]
candidates. Less obviously, requiring that these candidates respect minors of arity $\le |B|$ brings the number of candidates down to just
\[
\left|B^{\Pol_{|B|}(\fA)}\right| \le |B|^{|A|^{|A|^{|B|}}}.
\]
Unwinding the proofs of the ERP Theorem \ref{erp-thm} and Proposition \ref{erp-prop}, each of these corresponds to a candidate pp-construction of $\fB$.

More explicitly, suppose that
\[
\xi_{|B|} : \Pol_{|B|}(\fA) \rightarrow \Pol_{|B|}(\fB)
\]
defines our candidate height 1 homomorphism and respects minors. Then we define maps $g : \bA^{\bA^B} \rightarrow B$ and $h : B \rightarrow \bA^{\bA^B}$ as in the proof of the ERP Theorem \ref{erp-thm}, but with the definition of $g$ slightly modified: we set
\[
g(f(\pi_{b_1}, ..., \pi_{b_{|B|}})) = \xi_{|B|}(f)(b_1, ..., b_{|B|})
\]
for $|B|$-ary polymorphisms $f \in \Pol_{|B|}(\fA)$, with $b_1, ..., b_{|B|}$ a fixed enumeration of the elements of $B$. For each basic $m$-ary relation $R$ of $\fB$, we define the corresponding relation
\begin{align*}
\tilde{R} \subseteq \Big(\bA^{\bA^B}\Big)^m &= \{f(h(r_1), ..., h(r_k)) \mid f \in \Pol_k(\fA), r_1, ..., r_k \in R\}\\
&= \Sg_{(\bA^{\bA^B})^m}\{h(r) \mid r \in R\},
\end{align*}
exactly as in the proof of Proposition \ref{erp-prop}. Each $\tilde{R}$ is pp-definable from the basic relations of $\fA$, so the relational structure $\fC$ with underlying set $\bA^{\bA^B}$ and basic relations given by the $\tilde{R}$s is a pp-power of $\fA$ with the same relational signature as $\fB$. By construction, the map $h : \fB \rightarrow \fC$ will be a homomorphism, so the only challenge is to check whether or not $g : \fC \rightarrow \fB$ is actually a homomorphism, that is, to check whether or not
\[
g(\tilde{R}) \stackrel{?}{\subseteq} R
\]
for each basic relation $R$ of $\fB$. In fact, we can forget about the height 1 homomorphism $\xi_{|B|}$ and just search for $g : \fC \rightarrow \fB$ - in other words, we treat $\fC$ as an instance of $\CSP(\fB)$ with $|\Pol_{|B|}(\fA)|$ variables that actually participate in any constraints, and $\fA$ will pp-construct $\fB$ if and only if the instance $\fC$ has a solution.

If $\fB$ has only finitely many basic relations, then we can test each candidate pp-construction in finite time, so in this case there is an effective procedure to decide whether or not $\fA$ pp-constructs $\fB$. (Note that each $\tilde{R}$ is determined by the collection of $|R|$-ary polymorphisms of $\fA$, so this argument also shows that in order to check that $\xi_{|B|}$ extends to a height 1 clone homomorphism from $\Pol(\fA)$ to $\Pol(\fB)$, we just need to check that it extends to a height 1 clone homomorphism from $\Pol_k(\fA)$ to $\Pol_k(\fB)$ for $k = \max |R|$ over the basic relations $R$ of $\fB$.)
\end{rem}

\section{Taylor Algebras}

Once we restrict to idempotent algebras, we can start playing games with identities involving nesting functions to simplify our criterion for NP-completeness.

\begin{defn} An algebra $\bA$ is a \emph{Taylor algebra} if it has an idempotent term $t$ that satisfies a system of identities of the form
\[
t\left(\begin{bmatrix} x & ? & \cdots & ? \\ ? & x & \cdots & ? \\ \vdots & \vdots & \ddots & \vdots \\ ? & ? & \cdots & x \end{bmatrix}\right) \approx t\left(\begin{bmatrix} y & ? & \cdots & ? \\ ? & y & \cdots & ? \\ \vdots & \vdots & \ddots & \vdots \\ ? & ? & \cdots & y \end{bmatrix}\right),
\]
where the $?$s are filled in somehow with $x$s and $y$s. Such an operation $t$ is called a \emph{Taylor term}, and a variety with a Taylor term is called a \emph{Taylor variety}.
\end{defn}

Note that by the defining identities of any Taylor term $t$, $t$ can't be any projection (unless the algebra in question has only one element).

\begin{thm}[Taylor \cite{taylor-varieties}] If an idempotent algebra $\bA$ satisfies any set of identities that can't be satisfied by projections, then it has a Taylor term. Equivalently, an idempotent variety is Taylor iff it does not contain a two element algebra having no nontrivial operations.
\end{thm}

Before we prove Taylor's theorem, we will work through an example.

\begin{ex} Let $f$ be an idempotent ternary term satisfying the identity
\[
f(f(y,x,z),x,f(z,y,y)) \approx f(x,y,z).
\]
Then
\[
t(x_1, ..., x_9) \coloneqq f(f(x_1,x_2,x_3), f(x_4,x_5,x_6), f(x_7,x_8,x_9))
\]
is a Taylor term, since it satisfies the identities
\[
t(y,x,z,x,x,x,z,y,y) \approx t(x,x,x,y,y,y,z,z,z) \approx t(x,y,z,x,y,z,x,y,z),
\]
and by specializing these identities (substituting $x=y$, $y=z$, or $z=x$) we can get a system of Taylor identities for $t$.
\end{ex}

\begin{defn}\label{star-comp-defn} If $f : D^k \rightarrow D$ and $g : D^l \rightarrow D$, we define the \emph{star composition} $f * g : D^{kl} \rightarrow D$ to be $f\circ (g, g, ..., g)$.
\end{defn}

\begin{prop} If $f,g$ are idempotent, then $f,g \in \langle f*g \rangle$.
\end{prop}
\begin{proof} $f(x_1, ..., x_k) \approx f(g(x_1, ..., x_1), ..., g(x_k, ...., x_k))$ and $g(x_1, ..., x_l) \approx f(g(x_1, ..., x_l), ..., g(x_1, ..., x_l))$.
\end{proof}

\begin{defn} An identity is called a \emph{height 1 identity}, or a \emph{minor identity}, if it has the form $f(x_1, ..., x_k) \approx g(y_1, ..., y_l)$, where the $x_i$s and $y_j$s are (not necessarily distinct) variables.
\end{defn}

\begin{prop} If an idempotent term satisfies a system of height 1 identities which can't be satisfied by projections, then it is a Taylor term.
\end{prop}

\begin{proof}[Proof of Taylor's Theorem.] By a compactness argument, there is a finite set $\mathcal{T}$ of identities satisfied by a finite set of operations $f_1, ..., f_n$ of $\bA$ which can't be satisfied by projections. Let $s = f_1 * \cdots * f_n$. Then each $f_i \in \langle s \rangle$, so we can convert $\mathcal{T}$ into a collection $\cT'$ of identities in $s$ which can't be satisfied by projections either.

The identities of $\cT'$ might involve some amount of nesting of $s$ within itself, that is, they may not be height 1 identities. Let $m$ be the greatest nesting depth occuring in $\cT'$, and let $t = s * \cdots * s$, with $m$ copies of $s$. Let $\cT''$ be the set of height 1 identities involving $t$ only which are satisfied by $\bA$. We claim that $\cT''$ can't be satisfied by projections.

To see this, note first that for every $k \le m$, if we index the variables of $t$ by $m$-tuples $(i_1, ..., i_m)$ of indices for coordinates of $s$, and if we let $x^k$ be the tuple of variables given by
\[
x^k_{(i_1, ..., i_m)} = y_{i_k}
\]
for all $(i_1, ..., i_m)$, then we have
\[
t(x^1) \approx \cdots \approx t(x^m) \;\;\; (\approx s(y)).
\]
If this system of height 1 identities in $t$ is satisfied by a projection $\pi_{(i_1, ..., i_m)}$, then we see that we must have $i_1 = \cdots = i_m = i$, say, for some index $i$ of the variables of $s$. But then there is some identity of $\cT'$ which is incompatible with $s = \pi_i$, and this identity of $\cT'$ can be modified by replacing variables $z$ by expressions $s(z,...,z)$ repeatedly until it becomes a height 1 identity involving only $t$, which will then be incompatible with $t = \pi_{(i,...,i)}$.
\end{proof}

\begin{cor} If $\bA$ is an idempotent algebra and $\CSP(\bA)$ is not NP-complete, then $\bA$ has a Taylor term.
\end{cor}

When $\bA$ is not Taylor, the above result lets us conclude that there is some two element algebra $\bB \in HSP(\bA)$ with no nontrivial operations, but it doesn't give us a good bound on how large a power of $\bA$ we will need to take to find $\bB$. It turns out that, in fact, if such a $\bB$ exists then it already can be found inside $HS(\bA)$. We will prove a slight generalization of this fact, which applies to \emph{strictly simple} algebras.

\begin{defn} An algebra $\bA$ is called \emph{strictly simple} if $\bA$ is simple and every subalgebra of $\bA$ either has size $1$ or is equal to $\bA$.
\end{defn}

\begin{lem}\label{strictly-simple-hs} If $\bA$ is an idempotent algebra and $\bB \in HSP_{fin}(\bA)$ is strictly simple, then $\bB \in HS(\bA)$. Note that if $\bA, \bB$ are both finite, then $\bB \in HSP_{fin}(\bA)$ iff $\bB \in HSP(\bA)$.

More generally, if $\bA_i$ are idempotent for all $i \in I$ and if $\bB \in HSP_{fin}(\{\bA_i\}_{i \in I})$ is strictly simple, then $\bB \in HS(\bA_i)$ for some $i \in I$.
\end{lem}
\begin{proof} (Following Zhuk \cite{zhuk-strong}) We prove the first statement - the proof of the general case is nearly identical, just slightly more notationally involved. Pick $n$ minimal such that there is some $\bS \le \bA^n$ and some $\sigma \in \Con(\bS)$ with $\bS/\sigma \cong \bB$. If there is any pair $r,s \in \bS$ such that $\pi_1(r) = \pi_1(s)$ but $r/\sigma \ne s/\sigma$, then if we set $a = \pi_1(r)$ and
\[
\bS' = \pi_{\{2, ..., n\}}(\bS \cap (\{a\} \times \bA^{n-1})) \le \bA^{n-1},
\]
and define $\sigma' \in \Con(\bS')$ by restricting $\sigma$ in the obvious way, then $r' = \pi_{\{2, ..., n\}}(r)$ and $s' = \pi_{\{2, ..., n\}}(s)$ have $r'/\sigma' \ne s'/\sigma'$. Thus $\bS'/\sigma'$ is isomorphic to a subalgebra of a quotient of $\bB$ of size at least $2$, so $\bS'/\sigma' \cong \bB$, contradicting the minimality of $n$.

Otherwise, if there is no such pair $r,s$, then there is a congruence $\sigma_1 \in \Con(\pi_1(\bS))$ such that for all $r \in \bS$, the congruence class $r/\sigma$ is completely determined by $\pi_1(r)/\sigma_1$. But then we have
\[
\bB \cong \bS/\sigma \cong \pi_1(\bS)/\sigma_1 \in HS(\bA).\qedhere
\]
\end{proof}

\begin{cor} If $\bA$ is finite and idempotent, then either $\bA$ has a Taylor term, or there is some two element algebra $\bB \in HS(\bA)$ with no nontrivial operations.
\end{cor}

\begin{cor} If $\bA$ is finite, idempotent, and has no Taylor term, then there are nonempty subalgebras $\bB, \bC \le \bA$ such that $\bB \cap \bC = \emptyset$ and $(\bB \cup \bC)^3 \setminus (\bB^3 \cup \bC^3) \le \bA^3$. In particular, $\CSP(\bA)$ can simulate NAE-SAT in a trivial way.
\end{cor}

\begin{rem} A recent result of Ol{\v{s}}{\'a}k simplifies the identities we need to consider even further. Ol{\v{s}}{\'a}k \cite{olsak-weak} proves that in any Taylor algebra, whether finite or infinite, there is always a $6$-ary \emph{weak 3-cube} term $t$, that is, an idempotent term satisfying the identity
\[
t(x,y,y,y,x,x) \approx t(y,x,y,x,y,x) \approx t(y,y,x,x,x,y).
\]
The weak 3-cube term may be understood as saying that the ternary relation on the free algebra $\cF_{\cV}(x,y)$ which is generated by the ternary Not-All-Equal relation on $\{x,y\}$ has a diagonal element.

Ol{\v{s}}{\'a}k's proof that such a term exists first uses the theory of absorbing subalgebras to produce a $12$-ary term which he calls a double loop term, and then simplifies it down to a weak $3$-cube term by using an intricate collection of identities which are satisfied by binary idempotent operations.
\end{rem}

\begin{rem} There is a curious connection between systems of two-variable height 1 identities on ternary functions and the problem 1-IN-3 SAT. Suppose that you are given such a system of identities $\cT$ on ternary functions $f_1, ..., f_n$, and that you want to determine whether these identities rule out projections.
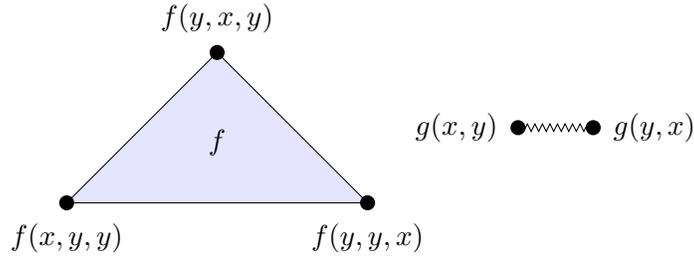
\begin{center}
\begin{tikzpicture}
    \node[label=above:{$f(y,x,y)$}] (v1) at (0,2) {};
    \node[label=below:{$f(x,y,y)$}] (v2) at (-2,0) {};
    \node[label=below:{$f(y,y,x)$}] (v3) at (2,0) {};
    \node[label=left:{$g(x,y)$}] (v4) at (4,1) {};
    \node[label=right:{$g(y,x)$}] (v5) at (5,1) {};

    \begin{scope}[fill opacity=0.2]
    \filldraw[fill=blue!50] ($(v1)$) 
        to[out=225,in=45] ($(v2)$) 
        to[out=0,in=180] ($(v3)$)
        to[out=135,in=-45] ($(v1)$);
    \end{scope}
    
    \draw[decorate, decoration={zigzag,segment length=1mm,amplitude=0.5mm}] ($(v4)$) to[out=0,in=180] ($(v5)$); 

    \foreach \v in {1,2,...,5} {
        \fill (v\v) circle (0.1);
    }
    
    \node at (0,0.8) {$f$};
\end{tikzpicture}
\end{center}
Define a set of binary functions $f_i^j(x,y)$, $j \le 3$, by
\begin{align*}
f_i^1(x,y) &= f_i(x,y,y),\\
f_i^2(x,y) &= f_i(y,x,y),\\
f_i^3(x,y) &= f_i(y,y,x),
\end{align*}
and identify any pair of $f_i^j$s which are indentified by $\cT$. Make a drawing of a hypergraph with a vertex for every equivalence class of $f_i^j$s, with a zigzag edge connecting any pair of vertices $g,h$ with $g(x,y) \approx h(y,x)$ under $\cT$, and with a hyperedge for each $f_i$ connecting it to $f_i^1, f_i^2, f_i^3$. An assignment of projections $\pi_j$ to the functions $f_i$ is the same as a choice $j$ of 1-IN-3 of the vertices on the hyperedge $f_i$ to be granted the value $\pi_1$, while every zigzag edge of the hypergraph corresponds to a $\ne$ constraint. (Ol{\v{s}}{\'a}k's paper \cite{olsak-weak} has one such picture, and I've found the technique enormously helpful for visualizing large systems of identities on ternary functions.)

As a concrete example, take following collection of four ternary terms $p,q,r,s$ defined in terms of Ol\v{s}\'ak's weak 3-cube term:
\begin{align*}
p(x,y,z) &= t(y,z,z,x,x,x),\\
q(x,y,z) &= t(x,z,y,z,y,z),\\
r(x,y,z) &= t(y,x,x,z,y,y),\\
s(x,y,z) &= t(x,x,y,z,y,x).
\end{align*}
Then the definitions together with the weak 3-cube identities imply the following system of identities on $p,q,r,s$:
\begin{align*}
p(x,y,x) &\approx q(y,x,x),\\
p(x,x,y) &\approx r(y,x,x),\\
q(x,y,x) &\approx s(x,y,x),\\
r(x,x,y) &\approx s(x,x,y),\\
s(x,y,y) &\approx q(x,x,y) \approx r(x,y,x). 
\end{align*}

It may not be apparent, at a glance, whether or not this system of identities can be satisfied by projections. If we draw the associated 1-IN-3 SAT instance, we find that it has 7 vertices (corresponding to binary terms), 4 occurences of the 1-IN-3 SAT constraint (for the four ternary terms $p,q,r,s$), and one occurence of the $\ne$ constraint (coming from the fact that the last identity above relates $s(x,y,y)$ to $q(x,x,y)$):
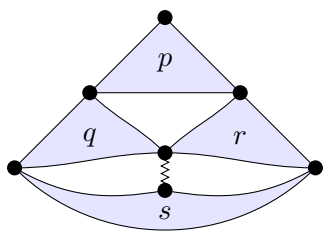
\begin{center}
\begin{tikzpicture}
    \node (v1) at (0,2) {};
    \node (v2) at (-1,1) {};
    \node (v3) at (1,1) {};
    \node (v4) at (-2,0) {};
    \node (v5) at (0,0.2) {};
    \node (v6) at (2,0) {};
    \node (v7) at (0,-0.3) {};

    \begin{scope}[fill opacity=0.2]
    \filldraw[fill=blue!50] ($(v1)$) 
        to[out=225,in=45] ($(v2)$) 
        to[out=0,in=180] ($(v3)$)
        to[out=135,in=-45] ($(v1)$);
    \filldraw[fill=blue!50] ($(v2)$) 
        to[out=225,in=45] ($(v4)$) 
        to[out=0,in=180] ($(v5)$)
        to[out=135,in=-45] ($(v2)$);
    \filldraw[fill=blue!50] ($(v3)$) 
        to[out=225,in=45] ($(v5)$) 
        to[out=0,in=180] ($(v6)$)
        to[out=135,in=-45] ($(v3)$);
    \filldraw[fill=blue!50] ($(v7)$) 
        to[out=190,in=-30] ($(v4)$) 
        to[out=-45,in=-135] ($(v6)$)
        to[out=-150,in=-10] ($(v7)$);
    \end{scope}
    
    \draw[decorate, decoration={zigzag,segment length=1mm,amplitude=0.5mm}] ($(v5)$) to[out=270,in=90] ($(v7)$); 

    \foreach \v in {1,2,...,7} {
        \fill (v\v) circle (0.1);
    }
    
    \node at (0,1.4) {$p$};
    \node at (-1,0.4) {$q$};
    \node at (1,0.4) {$r$};
    \node at (0,-0.6) {$s$};
\end{tikzpicture}
\end{center}
It is now easy (well, as easy as solving a small instance of 1-IN-3 SAT) to verify that the associated 1-IN-3 SAT instance has no solution, so this system of identities can't be satisfied by projections.
\end{rem}

\begin{rem}\label{gen-eckmann-hilton} Taylor's original reason for studying Taylor algebras was to try to deeply understand the reason that $\pi_1$ of a topological group is always abelian. Taylor \cite{taylor-varieties} considers, for any variety $\cV$, the category of topological $\cV$-objects, that is, topological algebraic structures satisfying the identities of $\cV$. Taylor showed that the $\pi_1$s of topological $\cV$-objects will share a nontrivial property iff $\cV$ has a Taylor term, and that this occurs iff $\pi_1$ is always abelian. The fact that a Taylor term must be taken to be idempotent is related to the fact that the fundamental group is really a groupoid (in the sense of category theory), rather than a group, so only the idempotent operations of $\cV$ can constrain its structure (I'm slightly fuzzy on the details).

Aside from the topological details, this can be viewed as an analogue of the Eckmann-Hilton principle \cite{eckmann-hilton} which states that a unital magma object in the category of unital magmas is necessarily commutative and associative. In fact, the following result holds for Taylor algebras: if $\bA$ is a Taylor algebra, and $m : \bA^2 \rightarrow \bA$ is a homomorphism such that there exists an element $0 \in \bA$ with $m(0,x) = m(x,0) = x$ for all $x$, then $m$ is commutative and associative.

Note that our assumption on $m$ implies that $m*m$ satisfies the identities
\[
m(x,y) \approx m*m(x,y,0,0) \approx m*m(x,0,y,0) \approx \cdots \approx m*m(0,0,x,y),
\]
where in each $m*m$ we always have the $x$ occuring to the left of the $y$. Additionally, since $m*m : \bA^4 \rightarrow \bA$ is a homomorphism, for any $n$-ary operation $t$ of $\bA$ we can evaluate $(m*m)*t$ on a $4\times n$ matrix of variables in two different ways: we may either start by applying $t$ to the rows and then apply $m*m$ to the resulting column vector, \emph{or} we may first apply $m*m$ to the columns and then apply $t$ to the resulting row vector - either way gives the same result.

Using these two observations together with the Taylor identities for an $n$-ary Taylor term $t$, we prove that $m$ is commutative by writing $m(x,y)$ as $(m*m)*t$ applied to a $4\times n$ matrix of $0$s, $x$s, and $y$s where every column has an $x$ above a $y$, and manipulate this expression until every $y$ is above an $x$. The strategy is to always keep the $x$s in the middle two rows and the $y$s in the top or bottom, and to move a $y$ up a column whenever that column is free of $x$s. To temporarily move $x$s out of the way, we apply the Taylor identities for $t$ to swap them with $0$s, possibly shifting the $x$s up and down between the middle two rows to get to a configuration where the Taylor identities will apply. A similar argument with $m*m*m$ in the place of $m*m$ can be used to prove associativity. If $t$ is a $6$-ary weak 3-cube term, for instance, then a portion of the proof of the commutativity of $m$ goes as follows:
\begin{align*}
m\left(\begin{bmatrix} x\\ y \end{bmatrix}\right) &= (m*m)*t\left(\begin{bmatrix} 0 & 0 & 0 & 0 & 0 & 0 \\ x & 0 & 0 & 0 & x & x \\ 0 & x & x & x & 0 & 0 \\ y & y & y & y & y & y\end{bmatrix}\right) = (m*m)*t\left(\begin{bmatrix} 0 & 0 & 0 & 0 & 0 & 0 \\ 0 & 0 & x & x & x & 0 \\ 0 & x & x & x & 0 & 0 \\ y & y & y & y & y & y\end{bmatrix}\right)\\
&= (m*m)*t\left(\begin{bmatrix} y & 0 & 0 & 0 & 0 & y \\ 0 & 0 & x & x & x & 0 \\ 0 & x & x & x & 0 & 0 \\ 0 & y & y & y & y & 0\end{bmatrix}\right) = (m*m)*t\left(\begin{bmatrix} y & 0 & 0 & 0 & 0 & y \\ x & 0 & 0 & 0 & x & x \\ 0 & x & x & x & 0 & 0 \\ 0 & y & y & y & y & 0\end{bmatrix}\right)\\
&= \cdots = (m*m)*t\left(\begin{bmatrix} y & y & y & y & y & y \\ 0 & 0 & x & x & x & 0 \\ x & x & 0 & 0 & 0 & x \\ 0 & 0 & 0 & 0 & 0 & 0\end{bmatrix}\right) = m\left(\begin{bmatrix} y\\ x \end{bmatrix}\right),
\end{align*}
where we have used the Taylor identity $t(x,0,0,0,x,x) \approx t(0,0,x,x,x,0)$ satisfied by a weak 3-cube term to temporarily move the first and last $x$ out of the way.
\end{rem}

\section{Two simple algorithms (width 1 and bounded strict width)}

\begin{defn} A CSP template $\fA = (D, \Gamma)$ has \emph{relational width 1} if the relational width $(1,k)$ algorithm below solves it for some $k$.
\end{defn}

\begin{algorithm}
\caption{Relational width $(1,k)$ algorithm}
\begin{algorithmic}[1]
\State Set $S_v \gets D$ for each variable $v$.
\Repeat
\ForAll{$v_1, ..., v_k$}
\State Let $X$ be the set of solutions to the restriction of the CSP to the variables $v_1, ..., v_k$ (projecting each constraint onto this subset of variables).
\State Set $S_{v_i} \gets \pi_i(X \cap (S_{v_1} \times \cdots \times S_{v_k}))$ for each $i \le k$.
\State For each constraint $R$ which involves some $v_i$, remove all tuples of $R$ which are incompatible with $S_{v_i}$.
\EndFor
\Until{the sets $S_v$ stop changing.}
\State If any $S_v = \emptyset$, there is no solution.
\end{algorithmic}
\end{algorithm}

Compare this to the generalized arc-consistency algorithm, which is more popular (and more efficient!) in practice. (After this section, I'll usually refer to generalized arc-consistency as just ``arc-consistency'' to save space.)

\begin{algorithm}
\caption{Generalized arc-consistency algorithm}
\begin{algorithmic}[1]
\State Set $S_v \gets D$ for each variable $v$.
\While{some constraint $R$ on variables $(v_1, ..., v_m)$ has $\pi_j (R \cap (S_{v_1} \times \cdots \times S_{v_m})) \ne S_{v_j}$}
\State Set $S_{v_j} \gets \pi_j (R \cap (S_{v_1} \times \cdots \times S_{v_m}))$.
\EndWhile
\State If any $S_v = \emptyset$, there is no solution.
\end{algorithmic}
\end{algorithm}

\begin{thm}[Feder, Vardi \cite{feder-vardi}]\label{thm-width-1-arc} A CSP template $\fA$ has relational width 1 iff it is solved by the generalized arc-consistency algorithm.
\end{thm}
\begin{proof}[Sketch] Suppose $\fA$ has width $(1,k)$, and let $\fB$ be an instance of $\CSP(\fA)$. By a generalization of the randomized construction of graphs with large girth and large chromatic number, there is a relational structure $\fB'$ which has a map to $\fB$, has girth larger than $k$, and which has a map to $\fA$ iff $\fB$ has a map to $\fA$ (alternatively, if $\Gamma$ contains the equality relation, we can cheat by adding long chains of equalities). Since $\fB'$ locally looks like a tree, the width $(1,k)$ algorithm and the generalized arc-consistency algorithm give the same results for $\fB'$, so if there is no homomorphism from $\fB$ to $\fA$, then generalized arc-consistency applied to $\fB'$ will correctly find that there is no solution. But then generalized arc-consistency applied to $\fB$ will also find that there is no solution, since every deduction on $\fB'$ can be mimicked on $\fB$.
\end{proof}

\begin{defn} A connected relational structure is a \emph{tree} if every collection of occurences of relations with arities $r_1, ..., r_k$ involves at least $1 + \sum_i (r_i-1)$ distinct elements. A relational structure $\fA$ has \emph{tree duality} if for every $\fB$, there is a map $\fB \rightarrow \fA$ iff every tree which maps to $\fB$ has a map to $\fA$.
\end{defn}

\begin{prop} $\fA$ has width 1 iff it has tree duality.
\end{prop}
\begin{proof} If generalized arc-consistency shows that there is no homomorphism $\fB \rightarrow \fA$, then we can make a proof tree that shows that some set $S_v$ eventually becomes empty. Each node of the proof tree corresponds to the fact that some variable $w$ of $\fB$ takes values from a set $S_w$, and the hyperedges of the proof tree are labeled by relations of $\fB$. So the proof tree is actually a relational structure with a map to $\fB$, and the same sequence of generalized arc-consistency deductions apply to the proof tree to show that it has no map to $\fA$.
\end{proof}

\begin{rem} Essentially the same arguments apply for any width $(l,k)$, with ``trees'' replaced by ``$(l,k)$-trees'' (definition left as an exercise to the reader). Note that $(l,k)$-trees have tree-width $k-1$. When studying \emph{relational} width, we replace ``trees'' by ``$(l,k)$-reltrees'' (defined in \cite{no-pure-width-2} for $k = l$).
\end{rem}

\begin{rem} Dalmau has shown that any CSP with relational width $(2,2)$ is also solved by generalized arc-consistency \cite{no-pure-width-2}. 2-SAT is an example of a CSP with width $(2,3)$ which is \emph{not} solved by arc-consistency, so Dalmau's result is best possible.
\end{rem}

Generalized arc-consistency has a close connection with the algebraic concept of a ``subdirect product''.

\begin{defn} A subalgebra $\RR\le \bA_1\times \cdots \times \bA_n$ is called a \emph{subdirect product}, written $\RR \le_{sd} \bA_1\times \cdots \times \bA_n$, if $\pi_i(\RR) = \bA_i$ for all $i$.
\end{defn}

So an algebraic way of thinking of arc-consistency is that we shrink the domains of the variables until we get to a situation where every relation is a subdirect product. It's worth noting that as we shrink our domains and relations, the new domains and relations we obtain will always be preserved by any polymorphisms which preserved the original relations, since the new domains and relations can be defined by primitive positive formulas from the original ones.

We now find an algebraic characterization of CSP templates with width 1. The main idea is to consider the ``most generic'' problem which arc-consistency requires to have a solution, and to ask what such a solution must look like. This most generic problem will have a different variable for each possible nonempty set $S \subseteq D$, and will have all relations which are consistent with these sets imposed.

\begin{defn} For $\fA = (D, \Gamma)$ a relational structure, define $\cP_\emptyset(\fA)$ to be the structure with ground set $\cP(D)\setminus \{\emptyset\}$, and for every $m$-ary relation $R \in \Gamma$ let the corresponding relation $\cP_\emptyset(R)$ be the set of all $m$-tuples $S_1, ..., S_m \in \cP(D)\setminus \{\emptyset\}$ such that there is some nonempty $X \subseteq R$ with $\pi_i(X) = S_i$ for each $i$.
\end{defn}

Note that $\cP_\emptyset(R)$ can be equivalently defined as the set of $m$-tuples $(S_1, ..., S_m)$ such that $\pi_i(R \cap (S_1 \times \cdots \times S_m)) = S_i$ for each $i$.

\begin{defn} A homomorphism $\cP_\emptyset(\fA) \rightarrow \fA$ is called a \emph{set polymorphism} of $\fA$.
\end{defn}

\begin{defn} A function $f : D^k \rightarrow D$ is called \emph{totally symmetric} if the value of $f(a_1, ..., a_k)$ only depends on $\{a_1, ..., a_k\}$. Note that this is stronger than being symmetric, since the multiplicity of the $a_i$s is also ignored.
\end{defn}

\begin{thm}\label{thm-set-polymorphism} The following are equivalent:
\begin{itemize}
\item $\fA$ has width 1,

\item $\fA$ has a set polymorphism, and

\item $\fA$ has totally symmetric polymorphisms of every arity.
\end{itemize}
\end{thm}
\begin{proof} If $\fA$ has width 1, then generalized arc-consistency applied to $\cP_\emptyset(\fA)$ shows that there is a homomorphism $f: \cP_\emptyset(\fA) \rightarrow \fA$, since at every step the set associated to the variable $S \subseteq D$ will contain $S$ (by induction on the number of steps and the definition of $\cP_\emptyset(R)$). So suppose that $f$ is a set polymorphism, and for every $k \ge 1$, let $f_k$ be the totally symmetric function
\[
f_k(a_1, ..., a_k) = f(\{a_1, ..., a_k\}).
\]
We need to check that $f_k$ is a polymorphism of $\fA$. Suppose that $x_1, ..., x_k \in R$, then if $X = \{x_1, ..., x_k\}$, $f_k(x_1, ..., x_k)$ has $i$th coordinate equal to $f(\pi_i(X))$. Since $(\pi_1(X), ..., \pi_m(X)) \in \cP_\emptyset(R)$ by the definition of $\cP_\emptyset(R)$, we see that $f_k(x_1, ..., x_k) = (f(\pi_1(X)), ..., f(\pi_m(X))) \in R$.

Finally, suppose that $\fA$ has totally symmetric polymorphisms $f_k$ of every arity, and let $\fB$ be a (finite) instance such that generalized arc-consistency stops after finding nonempty sets $S_v$ for every variable $v \in \fB$. Choose $k$ at least as large as the largest number of tuples in any relation that shows up in $\fB$, and let $f$ be the function on sets of size $\le k$ associated to $f_k$. We claim that the map $v \mapsto f(S_v)$ defines a homomorphism from $\fB$ to $\fA$. To see this, let $(v_1, ..., v_m)$ be a tuple with the constraint $R$ imposed, and let $X = R \cap (S_{v_1}, ..., S_{v_m}) = \{x_1, ..., x_k\}$ (possibly with repeated $x_i$s if $|X| < k$). Then $f_k(x_1, ..., x_k) = (f(\pi_1(X)), ..., f(\pi_m(X))) = (f(S_{v_1}), ... f(S_{v_m})) \in R$ since $f_k$ is a polymorphism.
\end{proof}

\begin{cor} A relational structure $\fA$ has width 1 iff it is homomorphically equivalent to a pp-power of HORN-SAT.
\end{cor}
\begin{proof} Let $f$ be a set polymorphism of $\fA$, and let $f_k$ be the associated totally symmetric polymorphism of arity $k$. We define a height 1 clone homomorphism from $\langle \min \rangle \rightarrow \Pol(\fA)$ by sending $\min(x_1, ..., x_k)$ to $f_k(x_1, ..., x_k)$. Now apply the ERP Theorem \ref{erp-thm} and Proposition \ref{erp-prop} from the subsection on reflections.
\end{proof}

\begin{ex} Suppose that $\fA$ has a binary polymorphism $s$ which is associative, commutative, and idempotent (such an $s$ is called a \emph{semilattice operation}). Then we can define $n$-ary polymorphisms $s_n$ inductively by $s_n(x_1, ..., x_n) = s(s_{n-1}(x_1, ..., x_{n-1}), x_n)$, and $s_n$ will be totally symmetric for every $n$. Thus, every relational structure with a semilattice polymorphism has width 1.
\end{ex}

\begin{ex} We give an example of a width 1 algebra which is not a semilattice. Let $f$ be the idempotent set operation on $\{a,b,c\}$ given by
\[
f(\{a,b\}) = b,\ f(\{b,c\}) = c,\ f(\{c,a\}) = a,\ f(\{a,b,c\}) = a,
\]
and let $f_k$ be the associated totally symmetric polymorphism of arity $k$. We have $f_k \in \langle f_3 \rangle$ for every $k$, and in fact a $k$-ary function $g$ which depends on all its inputs is in $\langle f_3 \rangle$ iff its restriction to every two element subset of $\{a,b,c\}$ is equal to the corresponding restriction of $f_k$ (tricky exercise). The relational clone $\Inv(f_3)$ is generated by the ternary relations $R_{ab}, R_{bc}, R_{ca},$ where $R_{ab}$ is defined by
\[
R_{ab}(x,y,z) \coloneqq (x \in \{a,b\}) \wedge (x=a \implies y=z),
\]
and $R_{bc}, R_{ca}$ are defined similarly.
\end{ex}

\begin{ex}\label{width-1-not-finitely-gen} Here we give a more surprising example, of a width 1 clone such that no finitely generated subclone has width 1. Let $f$ be the idempotent set operation on $\{-1,0,1\}$ (which we stylize as $\{-,0,+\}$) given by
\[
f(\{0,-\}) = -,\ f(\{0,+\}) = +,\ f(\{-,+\}) = f(\{-,0,+\}) = 0,
\]
and let $f_k$ be the associated totally symmetric polymorphism of arity $k$. The clone $\cO$ generated by the collection of all $f_k$ then has width 1. Every finitely generated subclone of $\cO$ is contained in $\langle f_k\rangle$ for some $k$. To see that $\cO \ne \langle f_k \rangle$, consider the $k+1$-ary relation $R_k$ given by
\[
R_k(x_0, ..., x_k) \coloneqq \bigwedge_{i < j} (x_i + x_j \ge 0) \wedge (x_0, ..., x_k) \ne (0, ..., 0).
\]
Then it is easy to check that $R_k$ is preserved by $f_k$, but is not preserved by $f_{k+1}$. To see that $\langle f_k \rangle$ does not have width 1, define $R_k^-$ similarly to $R_k$, but with each $x_i + x_j \ge 0$ replaced by $x_i + x_j \le 0$. Then for $k\ge 2$ the instance
\[
R_k(x_0, ..., x_k) \wedge R_k^-(x_0, ..., x_k)
\]
of $\CSP(\Inv(\langle f_k \rangle))$ is arc-consistent (since both $R_k$ and $R_k^-$ are subdirect) but has no solution.

The relational clone $\Inv(\cO)$ corresponding to this example is generated by the unary relation $\{+\}$, the binary relations $x = -y$ and $x \le y$, and the ternary relation $(x \ge 0) \wedge (x = 0 \implies y = z)$. The clone $\cO$ is an example of a clone which is finitely related but not finitely generated.
\end{ex}

Note that one doesn't need to \emph{know} what the set polymorphism of $\fA$ is to apply the arc-consistency algorithm. If $\fA$ is a rigid core, we can use the self-reducibility of $\CSP(\fA)$ to find a solution to every solvable instance $\fB$ of $\CSP(\fA)$ in polynomial time. By applying this to $\cP_\emptyset(\fA)$, we can then \emph{find} a set polymorphism of $\fA$ - in time polynomial in the size of $\cP_\emptyset(\fA)$, which is sadly exponential in the size of $\fA$. The following problem is currently open.

\begin{prob} Given a rigid core $\fA$, can we determine whether it has width 1 in time polynomial in the size of the description of $\fA$?
\end{prob}

Now we move to the case of bounded strict width. This has a connection to an intriguing paper of Dechter \cite{dechter} which predates the algebraic approach to the CSP. The next definition follows Dechter \cite{dechter}.

\begin{defn} A partial assignment of values to variables is \emph{locally consistent} if it satisfies every constraint which only involves those variables. A CSP instance is \emph{strong $i$-consistent} if every locally consistent partial assignment to less than $i$ variables can always be extended to a locally consistent partial assignment of any containing set of $i$ variables. An instance is \emph{globally consistent} if every locally consistent partial assignment extends to a global solution.
\end{defn}

There is a straightforward polynomial time algorithm to enforce strong $i$-consistency for any fixed $i$, introducing new constraints of arity $< i$ by intersecting and existentially projecting old constraints until no changes occur. It is desirable to have globally consistent problems, because then a solution may be found greedily. Can we check if a given problem is globally consistent?

\begin{thm}[Dechter \cite{dechter}] If a CSP with domain sizes bounded by $n$ and all constraint arities bounded by $m$ is strong $(n(m-1)+1)$-consistent, then it is globally consistent.
\end{thm}
\begin{proof} Suppose for contradiction that some locally consistent partial assignment $a_1, ..., a_k$ to $v_1, ..., v_k$ can't be extended to $v_{k+1}$, $k \ge n(m-1)+1$. Then for every possible value $a$ of $v_{k+1}$, there is some constraint $C_a$ involving at most $m-1$ of the variables $v_1, ..., v_k$ which is inconsistent with this choice of $a$ and whichever of the $a_i$s are relevant. Thus, there is a collection of at most $n$ constraints $C_a$ involving at most $n(m-1)$ of the variables from $v_1, ..., v_k$ together with the variable $v_{k+1}$, for which a locally consistent partial assignment of all but one of the variables can't be extended. But this contradicts the assumption of strong $(n(m-1)+1)$-consistency.
\end{proof}

The trouble with applying Dechter's result is that as we enforce strong consistency, we may need to add constraints of higher and higher arities (aside from the case where $n(m-1) \le m$, which neatly explains why 2-SAT is easy). To avoid this, we want to find situations in which the newly introduced constraints can always be written as intersections of constraints of low arity.

\begin{defn} A CSP template $\fA = (D, \Gamma)$ has \emph{strict width} $l$ if every strong $(l+1)$-consistent instance of $\CSP(D, \langle \Gamma \rangle)$ which contains the projections of its relations onto subsets of size at most $l$ is globally consistent, and has its solution-set determined by the collection of relations of arity at most $l$.
\end{defn}

Note that the definition of strict width only makes sense in terms of the whole relational clone generated by $\Gamma$, a hint that it is properly viewed as an algebraic condition. Algebraically, the relevant result is the Baker-Pixley theorem \cite{baker-pixley}.

\begin{thm}[Baker, Pixley \cite{baker-pixley}] The following are equivalent for an algebraic structure $\bA$:
\begin{itemize}
\item every subalgebra of $\bA^n$ is equal to the intersection of its projections onto sets of at most $l$ coordinates, and

\item $\bA$ has an $(l+1)$-ary \emph{near-unanimity} term, that is, a term $t$ satisfying the identities
\[
x \approx t(y,x,...,x) \approx t(x,y,...,x) \approx t(x,x,...,y),
\]
where in each case all but one of the inputs to $t$ is $x$.
\end{itemize}
If $\bA$ is idempotent, then these are both equivalent to every subalgebra of $\bA^{l+1}$ being equal to the intersection of its projections onto sets of $l$ coordinates.
\end{thm}
\begin{proof} Note that if $|\bA| \ge 2$, then either condition implies $l > 1$ (consider the equality relation as a subalgebra of $\bA^2$). Suppose that the first condition holds, and consider the free algebra on $l+1$ generators $\cF_\bA(l+1) \subseteq \bA^{\bA^{l+1}}$ which is generated by the projections $\pi_i : \bA^{l+1} \rightarrow \bA$. Let $A^{l+1}_{nu}$ be the set of tuples of elements in $\bA^{l+1}$ which have all but at most one entry equal to each other, and let $X \subseteq \bA^{A^{l+1}_{nu}}$ be the projection of $\cF_\bA(l+1)$ onto these coordinate tuples.

We claim that $X$ contains the tuple $t$ of near-unanimous votes of the entries of the coordinate tuples. By assumption, we just have to check that for every projection $\pi_{x_1, ..., x_l}(X)$ onto at most $l$ coordinates $x_1, ..., x_l \in A^{l+1}_{nu}$, there is some element $f \in \cF_\bA(l+1)$ with $\pi_{x_1, ..., x_l}(f) = \pi_{x_1, ..., x_l}(t)$. But each tuple $x_i$ has at most one dissenting coordinate, so there must be some coordinate $j \le l+1$ such that each $(x_i)_j$ is equal to the vote $t(x_i)$. Thus we can take $f = \pi_j$ to see that $\pi_{x_1, ..., x_l}(\pi_j) = \pi_{x_1, ..., x_l}(t)$.

Now suppose that $t$ is an $(l+1)$-ary near-unanimity term, and suppose that $\bB \subseteq \bA^n$. Let $b \in \bA^n$ be such that $\pi_I(b) \in \pi_I(\bB)$ for every $I \subseteq \{1, ..., n\}$ with $|I| \le l$, we will show by induction on $|J|$ that $\pi_J(b) \in \pi_J(\bB)$ for \emph{every} subset $J \subseteq \{1, ..., n\}$. For the inductive step, if $|J| \ge l+1$ then we may set $J_1, ..., J_{l+1}$ to be subsets of $J$ formed by deleting different elements of $J$, and for each $J_i$ there is some $b_{J_i} \in \bB$ with $\pi_{J_i}(b_{J_i}) = \pi_{J_i}(b)$ by induction. But then $b_J = t(b_{J_1}, ..., b_{J_{l+1}}) \in \bB$ and has $\pi_J(b_J) = \pi_J(b)$ by the near-unanimity equations.

For the last claim, if $\bA$ is idempotent and $\bB \subseteq \bA^n$ with $n > l+1$ and $b \in \bigcap_{|I| = l} \pi_I(\bB)$, then $\bB' = \pi_{\{1, ..., n-1\}}(\bB \cap (\bA^{n-1}\times \{b_n\}))$ is a subalgebra of $\bA^{n-1}$, and we may induct on $n$ to see that $\bB' = \bigcap_{|I| = l} \pi_I(\bB')$, while the assumption on subalgebras of $\bA^{l+1}$ gives $\pi_I(b) \in \pi_I(\bB')$ for every $I$ with $|I| = l$.
\end{proof}

\begin{thm} A relational structure $\fA$ has strict width $l$ iff it has an $(l+1)$-ary near-unanimity polymorphism.
\end{thm}
\begin{proof} Let $\bA$ be the associated algebraic structure. For any $n$ and any $\bB \subseteq \bA^n$, the strong $(l+1)$-consistent instance formed via the relations $\bB$ and $\pi_I(\bB)$ for all $I \subseteq \{1, ..., n\}$ with $|I| \le l$ together with the definition of strict width $l$ imply that $\bB = \bigcap_{|I| \le l} \pi_I(\bB)$, so by the Baker-Pixley Theorem $\bA$ has an $(l+1)$-ary near unanimity term.

For the other direction, suppose that $t$ is an $(l+1)$-ary near-unanimity term operation of $\bA$ and that we have a strong $(l+1)$-consistent instance of $\CSP(\bA)$, which we may assume by the Baker-Pixley Theorem to only involve relations of arity at most $l$. Suppose that we have a locally consistent partial solution which assigns the values $a_1, ..., a_k$ to the variables $v_1, ..., v_k$ which we want to extend to the variable $v_{k+1}$. By strong $(l+1)$-consistency, we can assume that $k \ge l+1$. By induction on $k$, we can assume that for each $i \le l+1$ there is some value $a_{k+1}^i$ that we can assign the the variable $v_{k+1}$ such that if we ignore $v_i$, we get a locally consistent partial solution.

We claim that assigning the value $a_{k+1} = t(a_{k+1}^1, ..., a_{k+1}^{l+1})$ to $v_{k+1}$ gives a locally consistent partial solution. To see this, consider some constraint $C$ which involves the variable $v_{k+1}$ and some variables from $v_1, ..., v_k$. For each $i \le l+1$, by $l$-consistency and the fact that $C$ has arity at most $l$ we can find a value $a_i'$ such that $(a_1, ..., a_i', ..., a_{k+1}^i)$ satisfies the constraint $C$. Applying $t$ to these $l+1$ tuples, we see that the tuple $(a_1, ..., a_{l+1}, ..., t(a_{k+1}^1, ..., a_{k+1}^{l+1}))$ also satisfies $C$, by the near-unanimity identities and the fact that $t$ is a polymorphism of $C$.
\end{proof}

\begin{algorithm}
\caption{Strict width $l$ algorithm}
\begin{algorithmic}[1]
\State Replace each constraint with its projections onto all subsets of at most $l$ variables.
\Repeat
\ForAll{sets $\{v_1, ..., v_k\}$ of variables with $k \le l+1$}
\State Let $X$ be the set of solutions to the restriction of the CSP to the variables $v_1, ..., v_k$.
\State If $\pi_I(X)$ is not implied by the restriction of the CSP to the variables in $I$ for some $I \subset \{v_1, ..., v_k\}$, add it as a new constraint.
\EndFor
\Until{no new constraints are added.}
\State Greedily assign values to variables until we find a global solution.
\end{algorithmic}
\end{algorithm}

\begin{ex} 2-SAT has the ternary polymorphism $\maj$, which is a near-unanimity operation. Therefore 2-SAT has strict width $2$, a fact which also follows from Dechter's result above \cite{dechter}.
\end{ex}

\begin{ex}\label{ex-dual-discriminator} Generalizing 2-SAT, let $D$ be any domain, and let $d : D^3 \rightarrow D$ be given by
\[
d(x,y,z) = \begin{cases}x & \text{if } y \ne z,\\ y & \text{if } y = z.\end{cases}
\]
This function $d$ is known as the \emph{dual discriminator}, and for $|D| \ne 4$ it is the only majority function (up to permuting inputs) on $D$ which preserves the graph of every bijection from $D$ to itself.

A binary relation $R \subseteq D^2$ is preserved by the dual discriminator iff it is a ``0/1/all constraint'', that is, a constraint such that when viewed as a bipartite graph on the disjoint union $D \sqcup D$, every vertex which doesn't have degree $0$ or $1$ connects to all vertices on the other side which have positive degree. Typical 0/1/all constraints are displayed below.

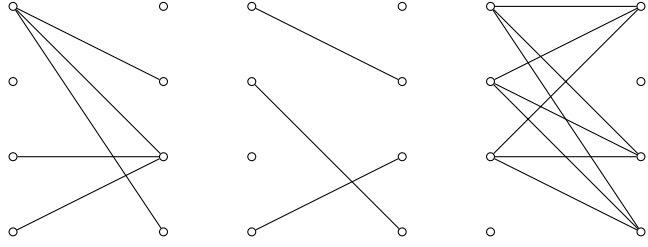
\begin{center}
\begin{tabular}{c@{\hskip 30pt}c@{\hskip 30pt}c}
\begin{tikzpicture}[scale=1]
  \node [circle, minimum width=3pt, draw, inner sep=0pt] (a1) at (0,0) {};
  \node [circle, minimum width=3pt, draw, inner sep=0pt] (a2) at (0,1) {};
  \node [circle, minimum width=3pt, draw, inner sep=0pt] (a3) at (0,2) {};
  \node [circle, minimum width=3pt, draw, inner sep=0pt] (a4) at (0,3) {};
  \node [circle, minimum width=3pt, draw, inner sep=0pt] (b1) at (2,0) {};
  \node [circle, minimum width=3pt, draw, inner sep=0pt] (b2) at (2,1) {};
  \node [circle, minimum width=3pt, draw, inner sep=0pt] (b3) at (2,2) {};
  \node [circle, minimum width=3pt, draw, inner sep=0pt] (b4) at (2,3) {};
  \draw (a1) -- (b2) -- (a2);
  \draw (b2) -- (a4) --(b3);
  \draw (b4) -- (a4) -- (b1);
\end{tikzpicture} &
\begin{tikzpicture}[scale=1]
  \node [circle, minimum width=3pt, draw, inner sep=0pt] (a1) at (0,0) {};
  \node [circle, minimum width=3pt, draw, inner sep=0pt] (a2) at (0,1) {};
  \node [circle, minimum width=3pt, draw, inner sep=0pt] (a3) at (0,2) {};
  \node [circle, minimum width=3pt, draw, inner sep=0pt] (a4) at (0,3) {};
  \node [circle, minimum width=3pt, draw, inner sep=0pt] (b1) at (2,0) {};
  \node [circle, minimum width=3pt, draw, inner sep=0pt] (b2) at (2,1) {};
  \node [circle, minimum width=3pt, draw, inner sep=0pt] (b3) at (2,2) {};
  \node [circle, minimum width=3pt, draw, inner sep=0pt] (b4) at (2,3) {};
  \draw (a1) -- (b2);
  \draw (a3) -- (b1);
  \draw (a4) -- (b3);
\end{tikzpicture} &
\begin{tikzpicture}[scale=1]
  \node [circle, minimum width=3pt, draw, inner sep=0pt] (a1) at (0,0) {};
  \node [circle, minimum width=3pt, draw, inner sep=0pt] (a2) at (0,1) {};
  \node [circle, minimum width=3pt, draw, inner sep=0pt] (a3) at (0,2) {};
  \node [circle, minimum width=3pt, draw, inner sep=0pt] (a4) at (0,3) {};
  \node [circle, minimum width=3pt, draw, inner sep=0pt] (b1) at (2,0) {};
  \node [circle, minimum width=3pt, draw, inner sep=0pt] (b2) at (2,1) {};
  \node [circle, minimum width=3pt, draw, inner sep=0pt] (b3) at (2,2) {};
  \node [circle, minimum width=3pt, draw, inner sep=0pt] (b4) at (2,3) {};
  \draw (a2) -- (b1) -- (a3) -- (b2) -- (a4) -- (b4) -- (a2);
  \draw (a2) -- (b2);
  \draw (a3) -- (b4);
  \draw (a4) -- (b1);
\end{tikzpicture}
\end{tabular}
\end{center}

For any $a$ in $D$, a generating set of binary relations for $\Inv(d)$ is given by the graphs of a pair of bijections which generate the symmetric group on $|D|$ elements, the unary relation $D \setminus\{a\}$, and the binary relation $x=a \vee y=a$.
\end{ex}


\begin{ex}\label{ex-strict-width-n} For every $n$, the relational structure $(\{0,1\}, \{0\}, \le, \{0,1\}^n\setminus\{(0,...,0)\})$ has strict width exactly $n$. A near-unanimity term for it is given by the threshold function
\[
t_2^{n+1}(x_1, ..., x_{n+1}) = \begin{cases} 1 & \sum_i x_i \ge 2,\\ 0 & \sum_i x_i \le 1.\end{cases}
\]
To see that it doesn't have strict width less than $n$, note that the relation $\{0,1\}^n\setminus\{(0,...,0)\}$ is not the intersection of its projections onto $n-1$ coordinates. Note that this template also has width 1 (it is preserved by the semilattice operation $\max$), so the strict width algorithm is far from being the best way to solve it for $n$ large.
\end{ex}

Note that the existence of an $(l+1)$-ary near-unanimity operation in $\Pol(\fA)$ is equivalent to the solvability of the CSP instance $\Phi$ (of $\fA$ together with singleton unary relations) with variables indexed by elements of $\fA^{l+1}$ described by the primitive positive formula
\[
\Phi(t) \coloneqq \bigwedge_{R \in \Gamma} \bigwedge_{M \in R^{l+1}} t(M) \in R \wedge \bigwedge_{a,b \in \fA} t(b,a, ..., a) = t(a,b,...,a) = \cdots = t(a,a,...,b) \in \{a\}.
\]
This instance may be solved in polynomial time by the strict width $l$ algorithm, giving us an $(l+1)$-ary near-unanimity term $t$ as output. Note, however, that the number of variables is exponential in $l$ - what if we just want to know whether the structure $\fA$ has bounded strict width, allowing $l$ to be arbitrarily large?

\begin{prob} Given a relational structure $\fA$, determine whether it has bounded strict width.
\end{prob}

The good news is that whether the structure is given as a finite relational structure or a finite algebraic structure, the existence of a near unanimity term is at least \emph{decidable} \cite{near-unanimity-maroti}, \cite{near-unanimity-congruence-distributive}, \cite{near-unanimity-zhuk}. The bad news is that the minimal arity of a near-unanimity term may be very large.

\begin{thm}[Zhuk \cite{near-unanimity-zhuk}, Barto, Draganov \cite{near-unanimity-minimal}, Corollary \ref{cor-finitely-related-near-unanimity}, Examples \ref{ex-tree-formula-high-arity}, \ref{ex-tree-formula-binary}] For any relational structure $\fA$ with $|\fA| = n$ such that every basic relation of $\fA$ has arity at most $m$, if $\fA$ has bounded strict width, then $\fA$ has strict width at most
\[
\max(m-1,2)^{3^n - 2^{n+1}} + 1.
\]

Conversely, for each $m \ge 3$ and $n \ge 2$, there is an example of a relational structure with bounded strict width such that every basic relation of $\fA$ has arity at most $m$, which has no near-unanimity polymorphism of arity at most
\[
(m-1)^{2^{n-2}},
\]
and for $m = 2, n \ge 3$ there is an example with no near-unanimity polymorphism of arity at most
\[
2^{2^{n-3}}.
\]
\end{thm}

Luckily, it is possible to determine whether a relational structure has bounded strict width without actually exhibiting a near-unanimity polymorphism. For instance, in \cite{deciding-absorption-relational} a nondeterministic polynomial time algorithm, which only tests for the existence of certain chains of ternary polymorphisms of $\fA$, is given for deciding whether a given subset of $\fA$ is an absorbing subalgebra (defined later). Using the fact that cycle consistency solves CSPs which have bounded width (which we will prove later), this can be converted into a polynomial time algorithm for testing whether $\fA$ has bounded strict width.

\subsection{The Basic Linear Programming relaxation of a CSP}\label{subsection-lp}

Another simple algorithm for solving CSPs, which is closely related to generalized arc-consistency, is the basic LP relaxation. If the domain of each variable $v$ is $D_v$, we replace the set of potential values $D_v$ with its formal convex hull, which we can think of as the set of \emph{probability distributions} on $D_v$. We represent the probability distribution corresponding to a variable $v$ as a tuple of real numbers $p_{v,a}$, one for each $a \in D_v$, satisfying
\[
0 \le p_{v,a} \le 1, \sum_{a \in D_v} p_{v,a} = 1.
\]
We also replace each constraint with its convex hull. That is, if the constraint $C$ imposes the relation $R = R_C$ on the variables $v_1, ..., v_m$, then we require the existence of a probability distribution $p_{C,r}$, on the tuples $r$ of $R$ such that
\[
0 \le p_{C,r} \le 1, \sum_r p_{C,r} = 1,
\]
and which is compatible with the probability distributions on the individual variables in the sense that
\[
p_{v_i,a} = \sum_{r_i = a} p_{C,r}.
\]

If a problem is known not to be fully satisfiable, we can relax it further by extending the probability distributions over relations $R \subseteq D_{v_1} \times \cdots \times D_{v_m}$ to probability distributions over all of $D_{v_1} \times \cdots \times D_{v_m}$, and then try to maximize the sum of the probabilities that tuples which are supposed to be in $R$ are actually in $r$:
\[
\tfrac{1}{\#C}\sum_C \sum_{r \in R_C} p_{C,r}.
\]
This system of linear equations and inequalities, with the optimization target above, is known as the \emph{basic LP} relaxation of a given CSP instance.

\begin{thm}[Kun, O'Donnell, Tamaki, Yoshida, Zhou \cite{lp-width-1}]\label{lp-robust} For any relational structure $\fA$, the following are equivalent:
\begin{itemize}
\item the basic LP relaxation correctly solves every instance of $\CSP(\fA)$,

\item $\fA$ has symmetric polymorphisms of every arity.
\end{itemize}
Furthermore, if $\fA$ has width 1 then the basic LP relaxation can be used to \emph{robustly} solve $\CSP(\fA)$, that is, if we are given an instance which is $1 - \epsilon$ satisfiable, then we can find a solution which satisfies a $1 - O(1/\log(1/\epsilon))$ fraction of the constraints.
\end{thm}
\begin{proof} Suppose first that the basic LP solves $\CSP(\fA)$, and consider the (by now standard) instance $\Phi$ that describes the existance of a symmetric polymorphism of arity $n$:
\[
\Phi(s) \coloneqq \bigwedge_{R \in \Gamma} \bigwedge_{M \in R^n} s(M) \in R \wedge \bigwedge_{a_1, ..., a_n \in \bA} \bigwedge_{\sigma \in S_n} s(a_1, ..., a_n) = s(a_{\sigma(1)}, ..., a_{\sigma(n)}).
\]
By the assumption that the basic LP decides $\CSP(\fA)$, we just need to exhibit a fractional solution to this CSP. This is achieved by taking $s = \frac{1}{n}\pi_1 + \cdots + \frac{1}{n}\pi_n$: as a convex combination of polymorphisms, it satisfies the relaxation of the first collection of constraints, and since it is a symmetric convex combination of its inputs it satisfies the second collection of constraints.

For the other direction, suppose that an instance of the CSP has a fractional solution to its basic LP relaxation, with probability distributions $p_{v,a}$ for each variable/value and $p_{C,r}$ for each constraint/tuple. We may assume that these probabilities are all rational (since the defining system of linear equations and inequalities had rational coefficients), and that they have a common denominator $n$. By assumption $\fA$ has a symmetric polymorphism $s$ of arity $n$, which we can think of as a function from probability distributions with denominator $n$ over the domain of $\fA$ to elements of $\fA$.

Applying $s$ to each $p_{v,\cdot}$ gives an element $a_v \in \fA$, and applying it to each probability distribution $p_{C,\cdot}$ gives a tuple $r_C$ in the associated relation $R$ (since $s$ is a polymorphism). Furthermore, the compatibility equations between the distributions $p_{v_i,\cdot}$ and $p_{C,\cdot}$ that we get when $v_i$ is the $i$th coordinate of the constraint $C$, together with the symmetry of $s$, imply that $a_{v_i} = (r_C)_i$ for each $i$, so $(a_{v_1}, ..., a_{v_m}) = r_C \in R$. Thus the $a_v$s form a valid solution to the CSP instance.

Finally, assume that $\fA$ has width 1, with set polymorphism $f$, and suppose that our original instance was $1-\epsilon$ satisfiable. Then the basic LP finds a fractional solution with value $\ge 1 - \epsilon$. We will use the polymorphism $f$ to make a randomized rounding scheme. First, we immediately give up on any constraints $C$ that the LP only satisfies with value $\le 1 - \sqrt{\epsilon}$ - these can form at most a $\sqrt{\epsilon}$ fraction of the constraints by Markov's inequality. Second, we will choose a threshold $\theta \le \frac{1}{|\fA|}$, and for each variable $v$ we assign the value
\[
a_v = f(\{a \in \fA \mid p_{v,a} \ge \theta\}).
\]
Note that the restriction $\theta \le \frac{1}{|\fA|}$ ensures that the sets on the right hand side are nonempty. We will show that if $\theta$ is chosen from a certain probability distribution, then on average we will obtain a good solution to the CSP, and deduce from this that some specific choice of $\theta$ works at least as well. For this we need the following claim.

{\bf Claim.} If $C$ is the constraint $(v_1, ..., v_m) \in R$ which is satisfied with value $\ge 1 - \sqrt{\epsilon}$, and if $2\sqrt{\epsilon} \le \theta \le \frac{1}{|\fA|}$ is such that
\[
\theta \not\in (p_{v_i,a}/(2|R|), p_{v_i,a}]
\]
for any pair $i \le m, a \in \fA$, then $(a_{v_1}, ..., a_{v_m})$ satisfies $C$.

{\bf Proof of Claim.} For each $v$, let $S_v = \{a \mid p_{v,a} \ge \theta\}$, so $a_v = f(S_v)$. In order to show that $(a_{v_1}, ..., a_{v_m})$ satisfies $R$, we just need to check that this colection of sets $S_{v_i}$ together with $R$ form a generalized arc-consistent instance. Let $a \in S_{v_i}$ for some $i$, then we have $p_{v_i,a} \ge \theta \ge 2\sqrt{\epsilon}$ by the definition of $S_{v_i}$. From
\[
\sum_{r \in R, r_i = a} p_{C,r} \ge p_{v_i,a} - \sqrt{\epsilon} \ge p_{v_i,a}/2,
\]
we see that there must be some $r \in R$ with $r_i = a$ and $p_{C,r} \ge p_{v_i,a}/(2|R|)$. Since $p_{v_i,a} \ge \theta$, by the assumption on $\theta$ we have $p_{v_i,a}/(2|R|) \ge \theta$, so $p_{C,r} \ge \theta$. But then $p_{v_j,r_j} \ge p_{C,r} \ge \theta$ for all $j$, so $r_j \in S_{v_j}$ for all $j$, and we see that $a$ extends to a solution of $R \cap (S_{v_1} \times \cdots \times S_{v_m})$.

To finish the proof, we choose $\theta$ uniformly at random from the set $\{\frac{1}{|\fA|}, \frac{1}{|\fA|T}, ... \frac{1}{|\fA|T^b}\}$, where $T$ is twice the maximum number of tuples in any relation $R$ and $b = \lfloor \log(1/2|A|\sqrt{\epsilon})/\log(T) \rfloor$. Note that $b$ grows like $\log(1/\epsilon)$, that's the only important thing to keep track of in the mess. Then every constraint of arity $m$ which we hadn't given up on is satisfied with probability at least $1 - m|A|/b$ (since there are at most $m|A|$ bad choices of $\theta$ where the claim doesn't apply), and asymptotically that looks like $1 - O(1/\log(1/\epsilon))$.
\end{proof}

\begin{rem} The dependence of the error in $1 - O(1/\log(1/\epsilon))$ on $\epsilon$ in the previous theorem is best possible in the case of HORN-SAT: Guruswami and Zhou \cite{robust-horn-gap} show that there are integrality gap instances even for the SDP relaxation (see Example \ref{ex-simul-ind}), and by a fundamental result of Raghavendra \cite{raghavendra-optimal} they deduce that under the Unique Games conjecture it is NP-hard to find an assignment satisfying a $1-o(1/\log(1/\epsilon))$ fraction of the constraints of a HORN-SAT instance which is promised to be $1-\epsilon$ solvable.
\end{rem}

\begin{rem} In \cite{lp-width-1}, it is also claimed that the basic LP solves every instance of $\CSP(\fA)$ if and \emph{only if} $\fA$ has width 1. The proof has a subtle error, however. The following counterexample, due to Kun, can be found in \cite{dalmau-approximation}.
\end{rem}

\begin{ex} Let $\fA = (\{-1,0,1\}, R_+, R_-)$, where $R_+ = \{(a,b,c) \mid a+b+c \ge 1\}$ and $R_- = \{(a,b,c) \mid a+b+c \le -1\}$. Then for every $h,n$ with $h < \frac{n}{3}$, the function
\[
s_{h,n}(x_1, ..., x_n) = \begin{cases}1 & \sum_i x_i > h\\ 0 & -h \le \sum_i x_i \le h\\ -1 & \sum_i x_i < -h\end{cases}
\]
is a symmetric polymorphism of $\fA$. Thus $\CSP(\fA)$ is solved by the basic LP relaxation. However, $\fA$ has no totally symmetric polymorphism of arity $3$, since such a polymorphism would necessarily map the matrices
\[
\begin{bmatrix} -1 & 1 & 1\\ 1 & -1 & 1\\ 1 & 1 & -1\end{bmatrix} \in R_+^3, \begin{bmatrix} 1 & -1 & -1\\ -1 & 1 & -1\\ -1 & -1 & 1\end{bmatrix} \in R_-^3
\]
to the same diagonal tuple, so $\fA$ does not have width 1.
\end{ex}

\begin{ex}\label{lp-not-width-1} The previous example can be generalized to a much larger relational structure on $\{-1,0,1\}$ as follows. Set $s_n = s_{0,n}$, then it isn't hard to show that $s_n \in \Clo(s_2)$ for all $n$ (hint: start by defining $t_n(x_1, ..., x_n) = s_2(x_1,s_{n-1}(x_2, ..., x_n))$), so $\Inv(s_2)$ also defines a CSP template which is solved by the basic LP relaxation.

\begin{center}
\begin{tabular}{c|ccc} $s_2$ & $-$ & $0$ & $+$\\ \hline $-$ & $-$ & $-$ & $0$\\ $0$ & $-$ & $0$ & $+$\\ $+$ & $0$ & $+$ & $+$\end{tabular}
\end{center}

$\Inv(s_2)$ is generated by the unary relation $\{1\}$, the binary relation $x = -y$, and the set of \emph{odd cycle relations}, where the $m$-th odd cycle relation $R_m$ is defined by
\[
R_m(x_1, ..., x_{2m-1},y,z) \coloneqq (x_1 + x_2 \ge 0) \wedge \cdots \wedge (x_{2m-1} + x_1 \ge 0) \wedge (x_1 = \cdots = x_{2m-1} = 0 \implies y = z).
\]
(I found this set of generating relations by a technique I learned from Zhuk \cite{zhuk-key}, in which we search for ``key'' relations $R$, for which there is some ``key tuple'' $x \not\in R$ such that the relation $R$ is maximal among those relations of $\Inv(s_2)$ which do not contain $x$. It isn't hard to show that any key tuple must consist mostly of $0$s, and using the negation symmetry we can assume that $R$ contains all tuples in $\{0,1\}^n$ aside from the key tuple. Then we look at the set of pairs of coordinates that can't simultaneously be set to $-1$, and prove that the resulting graph can't be bipartite...)

The clone $\langle s_2 \rangle$ is not finitely related. To see this, define an operation $s_n'$ for $n$ odd by the rule
\[
s_n'(x_1, ..., x_n) = \begin{cases}s_{0,n}(x_1, ..., x_n) & \text{if some }x_i = 0,\\ s_{1,n}(x_1, ..., x_n) & \text{if all } x_i \in \{-1,1\}.\end{cases}
\]
For every odd $n = 2m-1$, the operation $s_n' \not\in \langle s_2 \rangle$ - since it does not preserve the relation $R_m$ - but the function $s_n'(x,x,y_3,...,y_n)$ we get by identifying two of its inputs \emph{is} in $\langle s_2 \rangle$ (exercise for the reader), so it preserves every relation in $\Inv(s_2)$ which contains strictly less than $n$ tuples.

The clone $\langle s_2 \rangle$ is strictly contained in the width 1 clone from Example \ref{width-1-not-finitely-gen}, and corresponds to a strictly larger relational clone with a tractable CSP. Later we will see that this relational clone can be enlarged further, such that the CSP remains solvable by bounded width reasoning.
\end{ex}

Currently it is unknown if the following problem can be solved using at most an exponential amount of time - we will show that it is at least decidable later in these notes (see Corollary \ref{cor-lp-decidable}).

\begin{prob} Given a finite relational structure $\fA$, determine if it has symmetric polymorphisms of every arity.
\end{prob}

An interesting result in this direction is proved in \cite{symmetric-polymorphisms}: an algebraic structure $\bA$ has symmetric term operations of all arities iff there is no $\bB \in HSP(\bA)$ which has a pair of automorphisms in $\Aut(\bB)$ having no common fixed point (in fact, if $\bA$ has no symmetric term operation of arity $n$, we can take $\bB$ to be the free algebra on $n$ variables in the variety generated by $\bA$). If $HSP(\bA)$ could be replaced by $HS(\bA)$ in their result, then this would imply that it is enough to check for the existence of symmetric polymorphisms of arities up to $|\fA|$.

Later we will prove that any Taylor algebra has cyclic term operations of all arities which have no small prime factors, so we might hope that we could use these to help construct symmetric polymorphisms of higher arity. More ingredients are likely needed for such an argument, however: in \cite{symmetric-polymorphisms}, an example is given of a relational structure which has cyclic polymorphisms of every arity, but which has no symmetric polymorphism of arity $5$.

The mysterious relational structure found in \cite{symmetric-polymorphisms} was eventually explained by the paper \cite{symmetric-terms-simple-groups}, which constructed an analogous relational structure for every finite non-abelian simple group, along with a corresponding type of polymorphism generalizing cyclic operations. The example from \cite{symmetric-polymorphisms} comes from the smallest non-abelian simple group, the alternating group $A_5$ of even permutations on five symbols. That example is (up to pp-constructability) the simplest relational structure which fails to have a $21$-ary polymorphism $f$ satisfying the identities
\begin{align*}
&f(x_1,x_2,x_3,x_4,x_5,y_1,y_2,y_3,y_4,y_5,y_6,z_{1,2},z_{1,3},z_{1,4},z_{1,5},z_{2,3},z_{2,4},z_{2,5},z_{3,4},z_{3,5},z_{4,5})\\
\approx\ &f(x_1,x_2,x_4,x_5,x_3,y_2,y_3,y_1,y_5,y_6,y_4,z_{1,2},z_{1,4},z_{1,5},z_{1,3},z_{2,4},z_{2,5},z_{2,3},z_{4,5},z_{3,4},z_{3,5})
\end{align*}
and
\begin{align*}
&f(x_1,x_2,x_3,x_4,x_5,y_1,y_2,y_3,y_4,y_5,y_6,z_{1,2},z_{1,3},z_{1,4},z_{1,5},z_{2,3},z_{2,4},z_{2,5},z_{3,4},z_{3,5},z_{4,5})\\
\approx\ &f(x_3,x_4,x_1,x_2,x_5,y_1,y_2,y_4,y_3,y_6,y_5,z_{3,4},z_{1,3},z_{2,3},z_{3,5},z_{1,4},z_{2,4},z_{4,5},z_{1,2},z_{1,5},z_{2,5}),
\end{align*}
and consists of a $21$-element domain together with a pair of binary relations, each binary relation corresponding to one of the identities above. The three sets of variables $\{x_1, ..., x_5\}, \{y_1, ..., y_6\}$, and $\{z_{1,2}, ..., z_{4,5}\}$ correspond to the three primitive group actions of $A_5$ (up to isomorphism of sets with an $A_5$-action) - or, equivalently, to the three conjugacy classes of maximal subgroups of $A_5$ (these maximal subgroups are isomorphic to $A_4, D_5$, and $S_3$).

\section{Mal'cev algebras}

The goal in this section and the next is to generalize group theoretic algorithms (such as the algorithm for solving XOR-SAT) by isolating the special feature of groups which makes them so nice. First we should connect groups to CSPs, by defining the correct analogue of ``affine spaces'' for general groups.

\begin{prop} If $G$ is a group, then a nonempty subset $H \subseteq G^n$ is preserved by the ternary operation $(x,y,z) \mapsto xy^{-1}z$ iff $H$ is a coset of a subgroup of $G^n$.
\end{prop}
\begin{proof} Let $U$ be the subgroup of $G^n$ generated by expressions of the form $y^{-1}z$ for $y,z \in H$. Then $H$ is preserved under $(x,y,z) \mapsto xy^{-1}z$ iff $H$ is closed under the right action of $U$, so $H$ is a union of left cosets of $U$. To see that $H$ is just a single coset, note that for $x,y \in H$, we have $x^{-1}y \in U$ and $x(x^{-1}y) = y$.

Conversely, if $H = hU$ for some subgroup $U$ of $G^n$, then
\begin{align*}
HH^{-1}H &= hU(hU)^{-1}hU\\
&= hUUh^{-1}hU\\
&= hUUU = hU = H.\qedhere
\end{align*}
\end{proof}

\begin{rem} In any group $G$, the idempotent operation $(x,y,z) \mapsto xy^{-1}z$ generates the clone of all operations of the form
\[
(x_1, ..., x_n) \mapsto x_{i_1}x_{i_2}^{-1}x_{i_3}x_{i_4}^{-1}x_{i_5}\cdots x_{i_{2k}}^{-1}x_{i_{2k+1}}.
\]
This clone will never contain the idempotent operation $(x,y) \mapsto xy^{-2}x^2$, unless the group $G$ is abelian. To see this, note that $xy^{-1}z$ preserves the binary relation
\[
R_a \coloneqq \Big\{\begin{bmatrix} g\\ ag \end{bmatrix} \mid g \in G \Big\}
\]
for every $a \in G$, so the relation $R_a$ must be preserved by every operation in the clone generated by $xy^{-1}z$. Applying the operation $xy^{-2}x^2$ to the pair of pairs $(1, a), (b^{-1}, ab^{-1}) \in R_a$, we get
\[
\begin{bmatrix} 1 \\ a \end{bmatrix}\begin{bmatrix} b^{-1}\\ ab^{-1} \end{bmatrix}^{-2}\begin{bmatrix} 1\\ a \end{bmatrix}^2 = \begin{bmatrix} b^2\\ aba^{-1}ba \end{bmatrix},
\]
and the resulting pair will be contained in $R_a$ if and only if
\[
ab^2 = aba^{-1}ba,
\]
or equivalently if and only if
\[
ab = ba.
\]
\end{rem}

The idempotent operation $(x,y,z) \mapsto xy^{-1}z$ was isolated by universal algebraists who wanted to understand the underlying reason for the fact that normal subgroups commute: if $K,N\lhd G$ are normal subgroups of a group $G$, then
\[
KN = NK
\]
and $KN$ is also a normal subgroup of $G$. Of course this is easy to verify in the context of groups, but from the point of view of universal algebra it is really saying something interesting about \emph{congruences} of groups. If $K,N$ correspond to congruences $\alpha, \beta$ on $G$, then we can view this equality as the statement that
\[
\alpha \circ \beta = \beta \circ \alpha = \alpha\vee\beta,
\]
where composition of binary relations is defined as follows.

\begin{defn} Let $R,S$ be binary relations $R \subseteq A\times B, S\subseteq B\times C$. Then we define their \emph{composition} $R\circ S$ to be the subset of $A\times C$ consisting of pairs $(a,c)$ such that there exists a $b \in B$ with $aRb$ and $bSc$. As a primitive positive formula, we can write this as
\[
R\circ S(a,c) \coloneqq \exists b\in B\ R(a,b) \wedge S(b,c).
\]
\end{defn}

In general, it is not the case that congruences commute. In order to find the smallest congruence containing a pair of congruences in a general algebraic structure, one uses the following fact.

\begin{prop} If $\alpha,\beta$ are congruences on an algebraic structure $\bA$, then their least upper bound $\alpha \vee \beta$ is the transitive closure of $\alpha \circ \beta$, that is,
\[
\alpha \vee \beta = \bigcup_{n \ge 0} (\alpha \circ \beta)^{\circ n}.
\]
\end{prop}

If $\alpha,\beta$ do commute, then the above formula simplifies to $\alpha \vee \beta = \alpha \circ \beta$. So it is natural to try to understand the collection of all algebraic structures with commuting congruences. Of course, a structure with no congruences at all has this property - but we want to understand algebraic structures that have a \emph{reason} for their congruences to commute, so rather than studying algebras in isolation we study varieties with this property.

\begin{defn} We say that a variety $\cV$ is \emph{congruence permutable} if for all $\bA \in \cV$ and all $\alpha, \beta \in \Con(\bA)$ we have $\alpha \circ \beta = \beta \circ \alpha$.
\end{defn}

\begin{thm} A variety $\cV$ is congruence permutable iff $\cV$ has a ternary term $p$ which satisfies the identity
\[
p(x,y,y) \approx p(y,y,x) \approx x.
\]
\end{thm}
\begin{proof} Suppose first that $\cV$ is congruence permutable. Let $\cF = \cF_\cV(x,y,z)$ be the free algebra on three generators in $\cV$. Define a congruence $\alpha$ on $\cF$ to be the least congruence with $x/\alpha = y/\alpha$, that is, $\alpha$ is the kernel of the homomorphism $\cF_\cV(x,y,z) \rightarrow \cF_\cV(x,z)$ given by $x,y \mapsto x, z \mapsto z$. Similarly, let $\beta$ be the least congruence on $\cF$ with $y/\beta = z/\beta$.

Then $(x,z) \in \alpha \circ \beta$, so if $\cV$ has commuting congruences, then there must be some $p(x,y,z) \in \cF$ such that $x/\beta = p(x,y,z)/\beta$ and $p(x,y,z)/\alpha = z/\alpha$. But this is equivalent to the pair of identities $x \approx p(x,y,y)$, $p(x,x,z) \approx z$.

Conversely, suppose such a term $p$ exists, and let $\bA \in \cV$ and $\alpha, \beta \in \Con(\bA)$. Then for any $a,b,c$ with $a/\alpha = b/\alpha$ and $b/\beta = c/\beta$ we have
\[
p(a,b,c)/\beta = p(a,b,b)/\beta = a/\beta
\]
and
\[
p(a,b,c)/\alpha = p(a,a,c)/\alpha = c/\alpha,
\]
so $(a,c) \in \beta\circ\alpha$.
\end{proof}

\begin{defn} A ternary term $p$ is called a \emph{Mal'cev term} if it satisfies the identity $p(x,y,y) \approx p(y,y,x) \approx x$. An algebra with a Mal'cev term is called a \emph{Mal'cev algebra}, and a variety with a Mal'cev term is called a \emph{Mal'cev variety}.
\end{defn}

One reason universal algebraists like congruence permutability is that it implies that the congruence lattice is \emph{modular}, a property first isolated by Dedekind in his investigation of the lattice of submodules of a module over a ring.

\begin{defn}\label{modular-defn} A lattice $\mathcal{L}$ is \emph{modular} if for all $\alpha, \beta, \gamma \in \mathcal{L}$, we have
\[
\alpha \le \beta \implies \alpha \vee (\gamma \wedge \beta) = (\alpha \vee \gamma) \wedge \beta.
\]
Equivalently, a lattice $\mathcal{L}$ is modular if it has no five element sublattice isomorphic to the lattice $\mathcal{N}_5$ whose Hasse diagram is a pentagon: consider the sublattice generated by
\begin{align*}
\alpha' &= \alpha \vee (\gamma \wedge \beta),\\
\beta' &= (\alpha \vee \gamma) \wedge \beta,
\end{align*}
and $\gamma$, with top element $\alpha \vee \gamma = \alpha' \vee \gamma$ and bottom element $\gamma \wedge \beta = \gamma \wedge \beta'$, and note that we always have $\alpha' \le \beta'$.
\end{defn}

\begin{center}
\begin{tikzpicture}[scale=2]
  \node (ag) at (0,1) {$\alpha\vee\gamma$};
  \node (a) at (1.4,-0.7) {$\alpha$};
  \node (b) at (1.4,0.7) {$\beta$};
  \node (g) at (-1,0) {$\gamma$};
  \node (bg) at (0,-1) {$\gamma\wedge\beta$};
  \node (a-gb) at (0.9,-0.3) {$\alpha\vee(\gamma\wedge\beta)$};
  \node (ag-b) at (0.9,0.3) {$(\alpha\vee\gamma)\wedge\beta$};
  \draw (ag-b) -- (ag) -- (g) -- (bg) -- (a-gb);
  \draw (a) -- (a-gb) -- (ag-b) -- (b);
\end{tikzpicture}
\end{center}

\begin{prop}\label{prop-malcev-modular} If $\bA$ has permuting congruences, then $\Con(\bA)$ is a modular lattice.
\end{prop}
\begin{proof} We just have to check that if $\alpha \le \beta$, then $\alpha \circ (\gamma \wedge \beta) \ge (\alpha \circ \gamma) \wedge \beta$. Suppose that
\[
(x,z) \in (\alpha \circ \gamma) \wedge \beta,
\]
and choose $y$ such that $(x,y) \in \alpha, (y,z) \in \gamma$. Then $(y,x) \in \alpha \subseteq \beta$ and $(x,z) \in \beta$, so
\[
(y,z) \in \beta \circ \beta = \beta,
\]
so $(y,z) \in \gamma\wedge \beta$, so $(x,z) \in \alpha \circ (\gamma \wedge \beta)$.
\end{proof}


An unexpectedly large example of a Mal'cev variety is the variety of \emph{quasigroups}.

\begin{defn} A binary operation on a finite set is called a \emph{quasigroup} if its multiplication table is a Latin square (i.e. each element appears exactly once in each row and column). The \emph{variety of quasigroups} has three basic operations $\cdot, /, \backslash$, which satisfy the following identities:
\[
(a\cdot b)/b \approx a, \;\;\; (a/b)\cdot b \approx a, \;\;\; b\backslash(b\cdot a) \approx a, \;\;\; b\cdot (b\backslash a) \approx a.
\]
\end{defn}

Note that in the finite case, if $\cdot$ is a quasigroup operation, then the operations $/, \backslash$ can be defined in terms of $\cdot$ by an iteration argument (for any invertible unary function $f$ on an $n$ element set, $f^{-1} = f^{\circ (n!-1)}$). For infinite quasigroups, they have to be introduced into the language explicitly.

\begin{prop}\label{prop-quasigroup-malcev} If $\bA = (A, \cdot, /, \backslash)$ is a quasigroup, then
\[
p: (x,y,z) \mapsto (x/y)\cdot((y/y)\backslash z)
\]
is a Mal'cev term.
\end{prop}
\begin{proof} Plugging in $x=y$ we get
\[
p(y,y,z) = (y/y)\cdot((y/y)\backslash z) \approx z,
\]
and plugging in $z=y$ we get
\[
p(x,y,y) = (x/y)\cdot ((y/y)\backslash y) \approx (x/y)\cdot ((y/y)\backslash ((y/y)\cdot y)) \approx (x/y)\cdot y \approx x.\qedhere
\]
\end{proof}


The corresponding property on the CSP side of the picture is something known as the \emph{parallelogram property} (some authors call this \emph{rectangularity}, although the definition of rectangularity is often slightly weaker in the case of relations of higher arity).

\begin{defn} A binary relation $R \subseteq A\times B$ has the \emph{parallelogram property} if whenever $(a,b), (c,b), (c,d) \in R$, we also have $(a,d) \in R$. A relation of higher arity is said to have the parallelogram property if every way of grouping its coordinates into two groups gives a binary relation with the parallelogram property.
\end{defn}

\begin{thm} A finite algebraic structure $\bA$ has a Mal'cev term $p$ iff every relation $\RR \in \Inv(\bA)$ has the parallelogram property.
\end{thm}
\begin{proof} Suppose first that $\bA$ has a Mal'cev term $p$, let $\bB, \bC \in \cV(\bA)$ and let $\RR \le \bB\times \bC$ be a subalgebra of their product. Suppose that $(a,b), (c,b), (c,d) \in \RR$. Then
\[
\begin{bmatrix} a\\ d\end{bmatrix} = p\left(\begin{bmatrix} a\\ b\end{bmatrix}, \begin{bmatrix} c\\ b\end{bmatrix}, \begin{bmatrix} c\\ d\end{bmatrix}\right) \in \RR,
\]
so $\RR$ has the parallelogram property.

Conversely, suppose that every relation in $\Inv(\bA)$ has the parallelogram property. Let $\pi_1, \pi_2 \in \bA^{\bA^2}$ be the elements corresponding to the functions $\pi_i : (a_1,a_2) \mapsto a_i$. Let $\RR \le (\bA^{\bA^2})^2$ be the subalgebra generated by the three pairs $(\pi_1,\pi_1), (\pi_2,\pi_1), (\pi_2,\pi_2)$. Then since $\RR$ has the parallelogram property, we must have $(\pi_1,\pi_2) \in \RR$, so there must be a ternary term $p$ such that
\[
\begin{bmatrix} \pi_1\\ \pi_2\end{bmatrix} = p\left(\begin{bmatrix} \pi_1\\ \pi_1\end{bmatrix}, \begin{bmatrix} \pi_2\\ \pi_1\end{bmatrix}, \begin{bmatrix} \pi_2\\ \pi_2\end{bmatrix}\right).
\]
But then this $p$ is a Mal'cev term for $\bA$.
\end{proof}

If we want to test whether an algebra has a Mal'cev term, then the above result would make it seem like we need to test whether relations of arbitrarily large arity have the parallelogram property. As it turns out, for idempotent algebras we only need to test whether all binary relations have the parallelogram property.

\begin{thm}[\cite{horowitz-complexity-malcev-conditions}, \cite{deciding-some-malcev-conditions}, \cite{deciding-n-permutable}, \cite{deciding-minority}, \cite{zhuk-strong}] A finite idempotent algebra $\bA$ has a Mal'cev term if and only if every binary relation $\RR \in \Inv_2(\bA)$ has the parallelogram property. More explicitly, this occurs if and only if we have
\[
\begin{bmatrix}a\\ d\end{bmatrix} \in \Sg_{\bA^2}\left\{\begin{bmatrix}a\\ b\end{bmatrix}, \begin{bmatrix}c\\ b\end{bmatrix}, \begin{bmatrix}c\\ d\end{bmatrix}\right\}
\]
for all $a,b,c,d \in \bA$.
\end{thm}
\begin{proof} (Following \cite{zhuk-strong}) Suppose that $\bA$ is not Mal'cev, and consider a relation $\RR \in \Inv(\bA)$ of minimal arity $n$ among those which do not have the parallelogram property. If $n=2$, we are done. Otherwise, we will try to use $\RR$ to define a relation of lower arity which also fails to have the parallelogram property.

Suppose that $\RR$ fails to have the parallelogram property when considered as a binary relation on $\bA^k \times \bA^{n-k}$, and assume without loss of generality that $k \ge 2$. Since $\RR$ fails to have the parallelogram property, there are tuples $a,c \in \bA^k$ and tuples $b,d \in \bA^{n-k}$ such that $(a,b), (c,b), (c,d) \in \RR$ but $(a,d) \not\in \RR$. Write $a_1 = \pi_1(a)$ and $a' = \pi_{2, ..., k}(a)$, and define $c_1, c'$ similarly. Define $\RR' \le \bA^{k-1}\times \bA^{n-k}$ by
\[
(x',y) \in \RR' \;\; \iff \;\; \exists x_1 \in \bA \; ((x_1,x'),y) \in \RR \wedge ((x_1,x'),b) \in \RR.
\]
Then we have $(a',b), (c',b), (c',d) \in \RR'$, so if $\RR'$ has the parallelogram property then we must have $(a',d) \in \RR'$. Thus there is some $e_1 \in \bA$ such that $((e_1,a'),b), ((e_1,a'),d) \in \RR$. Now define $\RR'' \le \bA \times \bA^{n-k}$ by
\[
(x_1,y) \in \RR'' \;\; \iff \;\; ((x_1,a'),y) \in \RR.
\]
Then we have $(a_1,b), (e_1,b), (e_1,d) \in \RR''$, so if $\RR''$ has the parallelogram property then we must have $(a_1,d) \in \RR''$, which means that $(a,d) \in \RR$, a contradiction.
\end{proof}

\begin{rem}\label{rem-malcev-construction} In \cite{deciding-minority}, the authors give an explicit polynomial time procedure to construct a Mal'cev term out of a collection of idempotent ``local Mal'cev terms'' $t_{abcd}$ satisfying
\[
t_{abcd}(a,b,b) = a, \;\; t_{abcd}(c,c,d) = d.
\]
The construction consists of two stages. In the first stage we construct, for each $a,b$, a term $t_{ab}$ which satisfies
\[
t_{ab}(a,b,b) = a, \;\; t_{ab}(y,y,x) \approx x.
\]
To do this, we pick an ordering $(c_i,d_i)$ of the set of ordered pairs $(c,d)$, and inductively define terms $t_{ab}^i$ by $t_{ab}^0(x,y,z) \coloneqq x$ and
\[
t_{ab}^{i+1}(x,y,z) \coloneqq t_{abu_id_i}(t_{ab}^i(x,y,z), t_{ab}^i(y,y,z), z),
\]
where $u_i = t_{ab}^i(c_i,c_i,d_i)$. These terms will satisfy $t_{ab}^i(a,b,b) = a$ and $t_{ab}^i(c_j,c_j,d_j) = d_j$ for all $j < i$. We finish the first stage by taking $t_{ab} \coloneqq t_{ab}^{n^2}$, where $n$ is the number of elements in our algebra.

The second stage of the construction is similar. We first pick an ordering $(a_i,b_i)$ of the set of ordered pairs $(a,b)$, and then inductively define terms $t_i$ by $t_0(x,y,z) \coloneqq z$ and
\[
t_{i+1}(x,y,z) \coloneqq t_{a_iv_i}(x,t_i(x,y,y),t_i(x,y,z)),
\]
where $v_i = t_i(a_i,b_i,b_i)$. These terms will satisfy $t_i(y,y,x) \approx x$ and $t_i(a_j,b_j,b_j) = a_j$ for all $j < i$. The term $p(x,y,z) \coloneqq t_{n^2}(x,y,z)$ will then be a Mal'cev term.
\end{rem}

\begin{defn} If $\RR \le_{sd} \bA\times \bB$ is a subdirect binary relation, then the \emph{linking congruence} of $\RR$ can refer to any of the following three congruences: the congruence $\ker \pi_1 \vee \ker \pi_2$ on $\RR$, the congruence $\alpha$ on $\bA$ generated by pairs $a,a' \in \bA$ such that there exists a $b \in \bB$ with $(a,b), (a',b) \in \RR$, or the similar congruence $\beta$ defined on $\bB$. The relation $\RR$ is called \emph{linked} if these congruences are full.
\end{defn}

Note that in the above definition, we have $\alpha = \pi_1(\ker \pi_1 \vee \ker \pi_2), \beta = \pi_2(\ker \pi_1 \vee \ker \pi_2)$, and
\[
\bA/\alpha \cong \RR/(\ker \pi_1 \vee \ker \pi_2) \cong \bB/\beta.
\]
A more \emph{visual} way to understand the linking congruence is to think of the relation $\RR$ as a bipartite graph on $\bA \sqcup \bB$, and to define the congruence classes to be the connected components of this graph. In particular, $\RR$ is linked iff this bipartite graph is connected.

\begin{prop}[Goursat's Lemma]\label{goursat} A subdirect binary relation $\RR \le_{sd} \bA\times \bB$ has the parallelogram property iff there are congruences $\alpha,\beta$ on $\bA,\bB$ respectively and an isomorphism $f : \bA/\alpha \rightarrow \bB/\beta$ such that, writing $\pi_{\alpha},\pi_{\beta}$ for the quotient maps, we have $\RR = \pi_{\alpha}^{-1} \circ f^{-1} \circ \pi_{\beta}$ (treating $\pi_\alpha, \pi_\beta, f$ as binary relations with inputs on the right and outputs on the left).
\end{prop}
\begin{proof} Thinking of $\RR$ as a bipartite graph on $\bA \sqcup \bB$, we just have to prove that every connected component of $\RR$ is a complete bipartite graph. Suppose $a \in \bA$ and $b \in \bB$ are in the same connected component of $\RR$, and let $a = a_1, b_1, ..., a_k, b_k = b$ be a path from $a$ to $b$ with $(a_i,b_i) \in \RR$ and $(a_{i+1},b_i) \in \RR$ for each $i$. We will show that $(a_1,b_i) \in \RR$ by induction on $i$:
\[
\begin{bmatrix} a_1\\ b_i \end{bmatrix}, \begin{bmatrix} a_{i+1}\\ b_i\end{bmatrix}, \begin{bmatrix} a_{i+1}\\ b_{i+1}\end{bmatrix} \in \RR \implies \begin{bmatrix} a_1\\ b_{i+1}\end{bmatrix} \in \RR.\qedhere
\]
\end{proof}

Despite the trivial nature of binary relations with the parallelogram property, higher arity relations can encode more complicated global information.

\begin{ex} Consider the affine algebra $\bA = (\ZZ/p, x-y+z)$, and let $\RR \le_{sd} \bA^n$ be the relation $x_1 + \cdots + x_n \equiv 0\pmod{p}$. Then if we think of $\RR$ as a (subdirect) binary relation on $\bA\times \bA^{n-1}$, it is the graph of the homomorphism $\bA^{n-1} \rightarrow \bA$ given by $(x_2, ..., x_n) \mapsto -x_2-\cdots -x_n\pmod{p}$.

More generally, for any $i$, if we think of $\RR$ as a subdirect binary relation on $\bA^i\times \bA^{n-i}$, then the linking congruence gives homomorphisms $\bA^i \rightarrow \bA \leftarrow \bA^{n-i}$: $(x_1, ..., x_i)\mapsto x_1 + \cdots + x_i\pmod{p}$ and $(x_{i+1}, ..., x_n) \mapsto -x_{i+1}-\cdots - x_n\pmod{p}$.
\end{ex}

Ternary relations on simple Mal'cev algebras have a particularly interesting structure.

\begin{prop} Let $\bA_1, \bA_2, \bA_3$ be simple idempotent Mal'cev algebras, and suppose that $\RR\le_{sd} \bA_1\times\bA_2\times\bA_3$ has $\pi_{i,j}(\RR) = \bA_i\times\bA_j$ for each $i\ne j$ but that $\RR \ne \bA_1\times\bA_2\times\bA_3$. Then for each $a \in \bA_1$, the relation
\[
\RR_a = \pi_{2,3}(\RR \cap (\{a\}\times \bA_2\times \bA_3))
\]
is the graph of an isomorphism between $\bA_2$ and $\bA_3$, and for every $b \in \bA_2, c \in \bA_3$ there is a unique $a \in \bA_1$ such that $(b,c) \in \RR_a$.
\end{prop}
\begin{proof} Consider $\RR$ as a subdirect relation on $\bA_1\times (\bA_2\times \bA_3)$. Since the linking congruence on $\bA_1$ is not full (else $\RR$ would be the full relation by the parallelogram property), it must be trivial (since $\bA_1$ is simple), so $\RR$ is the graph of a homomorphism from $\bA_2\times \bA_3$ to $\bA_1$, which proves the last assertion.

Similarly, $\RR$ may be viewed as the graph of a homomorphism from $\bA_1\times \bA_2$ to $\bA_3$, so $\RR_a$ is the graph of a surjective homomorphism from $\bA_2$ to $\bA_3$ for any $a \in \bA_1$ (surjective because $\pi_{1,3}(\RR) = \bA_1 \times \bA_3$), and by simplicity of $\bA_2$ this homomorphism must be an isomorphism.
\end{proof}

If we fix an isomorphism $\bA_1 \cong \bA_2 \cong \bA_3 \cong \bA$ coming from the above proposition, then the kernel of the associated homomorphism $\bA \times \bA \cong \bA_2\times \bA_3 \rightarrow \bA_1$ contains the diagonal of $\bA\times \bA$ as a congruence class. In this case - that is, the case where $\bA\times\bA$ has the diagonal as a congruence class  of some congruence - $\bA$ is called an \emph{abelian} algebra.


\begin{ex}\label{ex-simple-nonab-malcev} Let $\bA_n = (\{0, ..., n-1\}, p)$, where $p$ is the ternary Mal'cev operation defined by
\[
p(x,y,z) = \begin{cases}z & \text{if }x=y,\\ y & \text{if }x=z,\\ x & \text{if } x\not\in\{y,z\}.\end{cases}
\]
For $n \ge 3$, $\bA_n$ is simple and non-abelian (i.e.\! the diagonal is not a congruence class of any congruence on $\bA_n^2$). $\Inv(\bA_n)$ is generated by a pair of graphs of permutations of $\{0,...,n-1\}$ which generate the full symmetric group, the unary relation $x \ne 0$, and the ternary relation
\[
(x,y,z \in \{0,1\}) \wedge (x+y+z \equiv 0\!\!\!\pmod{2}).
\]
It is a good exercise to prove that the above relations generate $\Inv(\bA_n)$.
\end{ex}

\begin{ex}\label{ex-solvable-nonab-malcev} Here we describe an example of a three element Mal'cev algebra which is ``solvable'', but which is not abelian. Let $\bA = (\{0,1,*\}, p)$, where $p$ is the ternary Mal'cev operation defined by
\[
p(x,y,z) = \begin{cases}x & \text{if }y=z,\\ y & \text{if }x=z,\\ z & \text{if } x=y,\\ * & \text{if } \{x,y,z\} = \{0,1,*\}.\end{cases}
\]
Every two element subset of $\bA$ is a subalgebra isomorphic to the idempotent reduct of $\ZZ/2$, and $\bA$ has a congruence $\theta$ corresponding to the partition $\{0,1\}, \{*\}$ such that $\bA/\theta$ is also isomorphic to the idempotent reduct of $\ZZ/2$.

Along with the binary relations on $\bA$, which can be described by applying Goursat's Lemma \ref{goursat} to $\bA$ and its subalgebras, there is also the ternary relation
\[
(x = y = z = *) \vee (x,y,z \in \{0,1\} \wedge x+y+z \equiv 0\!\!\!\pmod{2}),
\]
whose elements correspond to the columns of the matrix
\[
\begin{bmatrix} * & 0 & 0 & 1 & 1\\ * & 0 & 1 & 0 & 1\\ * & 0 & 1 & 1 & 0\end{bmatrix}.
\]
That this relation forms a subalgebra of $\bA^3$ is related to the fact that $\theta$ can be considered to be an ``abelian congruence'' (in a sense we will define later). Finding primitive positive definitions of the other ternary relations on $\bA$ from the relations described above, and showing that they generate all of $\Inv(\bA)$, are left as exercises (if these exercises are too difficult, it may be helpful to come back to them after reaching Theorem \ref{rectangular-structure}).
\end{ex}

\section{Mal'cev algorithm and compact representations}\label{s-malcev-algorithm}

The algorithm for solving CSPs invariant under a Mal'cev operation, due to Bulatov and Dalmau \cite{bulatov-dalmau-malcev}, is based on the fact that any Mal'cev constraint has a small generating set. More specifically, we will show that any subset of a relation $\RR$ which has the same projection to each factor and contains representatives of all of the ``forks'' of $\RR$ actually generates $\RR$.

\begin{defn} If $R \subseteq \bA_1 \times \cdots \times \bA_n$, then we define the \emph{signature} of $R$, written $\operatorname{Sig}(R)$, to be the set of triples $(i,a,b)$ with $i \in \{1, ..., n\}, a,b, \in \bA_i$ such that there are some $t_a,t_b \in R$ with $\pi_{1, ..., i-1}(t_a) = \pi_{1, ..., i-1}(t_b)$ and $\pi_i(t_a) = a, \pi_i(t_b) = b$. In this case we say that the pair $t_a,t_b$ \emph{witnesses} the triple $(i,a,b)$.
\end{defn}

\begin{thm}\label{malcev-compact} Suppose that a relation $\RR \le \bA_1 \times \cdots \times \bA_n$ is preserved by a Mal'cev term $p$, and that $S \subseteq \RR$ is a subset with $\operatorname{Sig}(S) = \operatorname{Sig}(\RR)$. Then $\RR$ is generated by $S$ (using only $p$).
\end{thm}
\begin{proof} Let $\bS$ be the subset of $\RR$ generated by $S$ using $p$. We will prove by induction on $i$ that $\pi_{1, ..., i}(\bS) = \pi_{1, ..., i}(\RR)$.

Suppose that $t \in \RR$. By the induction hypothesis, there is some $t' \in \bS$ with $\pi_{1, ...,i-1}(t) = \pi_{1, ..., i-1}(t')$. Let $a = \pi_i(t'), b = \pi_i(t)$. Since $\bS \subseteq \RR$, we have $(i,a,b) \in \operatorname{Sig}(\RR) = \operatorname{Sig}(S)$, so there must be a pair $t_a, t_b \in S$ witnessing the triple $(i,a,b)$. Define $t'' \in \bS$ by
\[
t'' = p(t',t_a,t_b).
\]
Then from $\pi_{1, ...,i-1}(t_a) = \pi_{1, ..., i-1}(t_b)$ and the fact that $p$ is Mal'cev, we have
\[
\pi_{1, ...,i-1}(t'') = \pi_{1, ...,i-1}(p(t',t_a,t_a)) = \pi_{1, ...,i-1}(t') = \pi_{1, ...,i-1}(t).
\]
Additionally, from $\pi_i(t') = \pi_i(t_a) = a$ and the fact that $p$ is Mal'cev, we have
\[
\pi_i(t'') = p(a,a,b) = b = \pi_i(t),
\]
so $\pi_{1, ..., i}(t'') = \pi_{1, ..., i}(t)$.
\end{proof}

\begin{defn} A subset $S \subseteq \RR$ is called a \emph{compact representation} of a Mal'cev relation $\RR$ if $\Sig(S) = \Sig(\RR)$ and $|S| \le 2|\Sig(\RR)|$.
\end{defn}

\begin{prop} Every Mal'cev relation $\RR \le \bA_1 \times \cdots \times \bA_n$ has a compact representation $S$. We always have $|S| \le 2n\cdot \max_i |A_i|^2$.
\end{prop}

Now we need some subroutines for manipulating compact representations. The first such procedure is called \texttt{Nonempty}: it takes as input a compact representation $R$ of a relation $\RR \le \bA_1 \times \cdots \times \bA_n$ and any description of a relation $\bS \le \bA_{i_1} \times \cdots \times \bA_{i_k}$ on a small subset $\{i_1, ..., i_k\}$ of the indices, and it tells us whether $\RR\cap \bS \ne \emptyset$. In the case $\RR \cap \bS \ne \emptyset$, \texttt{Nonempty} returns an element of the intersection.

\begin{algorithm}
\caption{\texttt{Nonempty}$(R, i_1, ..., i_k, \bS)$, $p$ a Mal'cev term, $R$ a compact representation of $\RR \le \bA_1\times \cdots\times \bA_n$, $\bS \le \bA_{i_1} \times \cdots \times \bA_{i_k}$.}
\begin{algorithmic}[1]
\State Set $R' \gets R$.
\While{$\pi_{i_1, ..., i_k}(R')$ is not closed under $p$ and $R' \cap \bS = \emptyset$}
\State Pick $t_1, t_2, t_3 \in R'$ with $\pi_{i_1, ..., i_k}(p(t_1,t_2,t_3)) \not\in \pi_{i_1, ..., i_k}(R')$.
\State Set $R' \gets R' \cup \{p(t_1,t_2,t_3)\}$.
\EndWhile
\If{$R' \cap \bS \ne \emptyset$}
\State \Return{any element of $R' \cap \bS$.}
\Else
\State \Return{$\emptyset$.}
\EndIf
\end{algorithmic}
\end{algorithm}

\begin{prop} \texttt{Nonempty} correctly determines whether $\RR \cap \bS \ne \emptyset$ in time polynomial in $n$, $|R|$, and $|\pi_{i_1, ..., i_k}(\RR)| \le \prod_{j \le k} |\bA_{i_j}|$.
\end{prop}
\begin{proof} Since $\RR$ is generated by $R$ using $p$, we also have $\pi_{i_1,...,i_k}(\RR)$ generated by $\pi_{i_1, ..., i_k}(R)$ using $p$. To see the bound on the running time, note that in each iteration of the while loop, the set $\pi_{i_1, ..., i_k}(R')$ gains a new element, and its size is clearly bounded by $|\pi_{i_1, ..., i_k}(\RR)|$.
\end{proof}

The next subroutine for manipulating compact representations is \texttt{Fix-values}. \texttt{Fix-values} converts a compact representation $R$ of $\RR \le \bA_1\times \cdots\times \bA_n$ to a compact representation of
\[
\RR \wedge (x_1 = a_1) \wedge \cdots \wedge (x_m = a_m),
\]
for any choice of $m \le n$ and $a_i \in \bA_i$ for all $i$. \texttt{Fix-values} is really the core of the algorithm, the other steps are mostly formal (in fact, \texttt{Nonempty} and \texttt{Fix-values} are the only two subroutines which use the Mal'cev term $p$).

\begin{algorithm}
\caption{\texttt{Fix-values}$(R, a_1, ..., a_m)$, $p$ a Mal'cev term, $R$ a compact representation of $\RR \le \bA_1\times \cdots\times \bA_n$.}
\begin{algorithmic}[1]
\State Set $R_0 \gets R$.
\For{$j$ from $1$ to $m$}
\If{$(j,a_j,a_j) \not\in \Sig(R_{j-1})$}
\State \Return $\emptyset$.
\Else
\State Set $R_{j} \gets \{t\}$, where $t \in R_{j-1}$ and the pair $t,t$ witnesses the triple $(j,a_j,a_j)$.
\EndIf
\ForAll{$(i,a,b) \in \Sig(R_{j-1})$ with $i > j$}
\State Let $t_a,t_b \in R_{j-1}$ witness the triple $(i,a,b)$.
\State Let $t \gets \texttt{Nonempty}(R_{j-1},j,i,\{(a_{j},a)\})$.
\If{$t \ne \emptyset$}
\State Set $R_{j} \gets R_{j} \cup \{t, p(t,t_a,t_b)\}$.
\EndIf
\EndFor
\EndFor
\State \Return $R_m$.
\end{algorithmic}
\end{algorithm}

\begin{prop} \texttt{Fix-values} correctly returns a compact representation of $\RR_m = \RR \wedge (x_1 = a_1) \wedge \cdots \wedge (x_m = a_m)$ in polynomial time.
\end{prop}
\begin{proof} We prove by induction on $j$ that $R_j$ is a compact representation of $\RR_j$ for each $j \le m$. Note that for any $(i,a,b) \in \Sig(R_j)$, if $a \ne b$ then we must have $i > j$. For $i \le j$, we have $(i,a_i,a_i) \in \Sig(R_j)$ iff $\RR_j \ne \emptyset$ by how we initialize $R_j$.

If $i > j$, then $(i,a,b) \in \Sig(\RR_j)$ implies $(i,a,b) \in \Sig(\RR_{j-1})$, witnessed by some pair $t_a,t_b \in R_{j-1}$. Additionally, if $(i,a,b) \in \Sig(\RR_j)$, then there is certainly some $t \in \RR_j$ with $\pi_i(t) = a$, so the call to \texttt{Nonempty} inside the loop will succeed. Then
\[
\pi_{1,...,i-1}(p(t,t_a,t_b)) = \pi_{1,...,i-1}(p(t,t_a,t_a)) = \pi_{1, ..., i-1}(t),
\]
so from $i > j$ we have $p(t,t_a,t_b) \in \RR_j$. From $\pi_i(t) = a, \pi_i(p(t,t_a,t_b)) = p(a,a,b) = b$, we see that the pair $t,p(t,t_a,t_b)$ witnesses the triple $(i,a,b)$.

To see that \texttt{Fix-values} runs in polynomial time, note that every call to \texttt{Nonempty} involves a constraint on two variables.
\end{proof}

\begin{cor} Given a compact representation $R$ of a relation $\RR \le \bA_1 \times \cdots \times \bA_n$ which is preserved by a given Mal'cev operation $p$, and given a tuple $t \in \bA_1 \times \cdots \times \bA_n$, we can determine whether $t \in \RR$ in polynomial time.
\end{cor}

The next subroutine will give a compact representation for the intersection of a relation $\RR$ given by a compact representation $R$ and a relation $\bS$ of small arity. In \cite{bulatov-dalmau-malcev} this subroutine was called \texttt{Next-beta}, so we will copy that notation here.

\begin{algorithm}
\caption{\texttt{Next-beta}$(R, i_1, ..., i_k, \bS)$, $R$ a compact representation of $\RR \le \bA_1\times \cdots\times \bA_n$, $\bS \le \bA_{i_1} \times \cdots \times \bA_{i_k}$.}
\begin{algorithmic}[1]
\State Set $R' \gets \emptyset$.
\ForAll{$(i,a,b) \in \Sig(R)$}
\State Set $t_a \gets \texttt{Nonempty}(R,i_1, ..., i_k,i,\bS\times \{a\})$.
\If{$t_a \ne \emptyset$}
\State Set $t_b \gets \texttt{Nonempty}(\texttt{Fix-values}(R, \pi_1(t_a), ..., \pi_{i-1}(t_a)), i_1, ..., i_k, i, \bS\times \{b\})$.
\If{$t_b \ne \emptyset$}
\State Set $R' \gets R' \cup \{t_a,t_b\}$.
\EndIf
\EndIf
\EndFor
\State \Return $R'$.
\end{algorithmic}
\end{algorithm}

\begin{prop} \texttt{Next-beta} correctly finds a compact representation of $\RR \cap \bS$ in time polynomial in $n$, $|R|$, and $|\pi_{i_1, ..., i_k}(\RR)|\cdot \max_i |\bA_i| \le \prod_{j \le k} |\bA_{i_j}|\cdot \max_i |\bA_i|$.
\end{prop}

Bulatov and Dalmau \cite{bulatov-dalmau-malcev} then go on to define a subroutine \texttt{Next} which calls \texttt{Next-beta} on larger and larger projections of $\bS$, ensuring that $|\pi_{i_1, ..., i_k}(\RR)| \le |\bS| \cdot \max_i |\bA_i|$ every time that \texttt{Next-beta} is called. A better approach, leading to a more powerful algorithm, was found by Mar\'oti \cite{malcev-on-top}. The subroutine \texttt{Intersect} takes two compact representations $R,S$ of relations $\RR, \bS$ as input and returns a compact representation of $\RR \cap \bS$ as output.

\begin{algorithm}
\caption{\texttt{Intersect}$(R, i_1, ..., i_k, S)$, $R$ a compact representation of $\RR \le \bA_1\times \cdots\times \bA_n$, $S$ a compact representation of $\bS \le \bA_{i_1} \times \cdots \times \bA_{i_k}$.}\label{alg-intersect}
\begin{algorithmic}[1]
\State Let $t_R \in R$ and $t_S \in S$ be any tuples.
\State Set $R' \gets (R \times \{t_S\}) \cup (\{t_R\}\times S) \subseteq \bA_1 \times \cdots \bA_n \times \bA_{i_1} \times \cdots \times \bA_{i_k}$.
\For{$j \le k$}
\State Set $R' \gets \texttt{Next-beta}(R', i_j, n+j, =_{\bA_{i_j}})$.
\EndFor
\State \Return a minimal subset of $\pi_{1, ..., n}(R')$ which witnesses every triple $(i,a,b) \in \Sig(\pi_{1, ..., n}(R'))$.
\end{algorithmic}
\end{algorithm}

\begin{thm} Any CSP which is preserved by a Mal'cev operation, where the relations are given by their compact representations, can be solved in time polynomial in the number of variables, the number of relations, and the size of the largest domain. In fact, we can find a compact representation of the solution set in polynomial time.
\end{thm}
\begin{proof} We start with any compact representation of $\bA_1 \times \cdots \times \bA_n$, and simply apply the subroutine \texttt{Intersect} repeatedly to find a compact representation of the intersection of all the constraint relations. To see that \texttt{Intersect} works correctly and efficiently, note that $R'$ is initialized as a compact representation of $\RR \times \bS$ and ends as a compact representation of $\RR\cap \bS$ followed by $k$ repeated coordinates. To see that \texttt{Intersect} runs in polynomial time, note that each call of \texttt{Next-beta} involves a relation of arity $2$.
\end{proof}

\begin{cor} For any primitive positive formula $\varphi$ in a collection of relations $\RR_i$, if we are given compact representations of each $\RR_i$ then we can efficiently find a compact representation of the relation described by $\varphi$.
\end{cor}
\begin{proof} If we are given a compact representation of a relation and we permute its variables, we can efficiently find a compact representation for the permuted relation by using the \texttt{Intersect} subroutine with $\RR$ equal to a full relation. To handle projections, note that we can project onto any initial segment of the variables by just projecting our compact representation and pruning it.
\end{proof}

While this might appear to be a fully satisfactory theory, there is still one big question remaining: what happens if instead of having relations described by compact representations, we have relations which are instead described by an arbitrary set of generators? It's clear that we just need to find a way to compute a compact representation of $\Sg_{\bA^n}(S)$ for any small set $S \subseteq \bA^n$, and a little thought shows that this can be reduced to the following problem.

\begin{prob} Let $\bA$ be a fixed Mal'cev algebra. Given a subset $S \subseteq \bA^n$, and given a tuple $t \in \bA^n$, can we determine whether $t \in \Sg_{\bA^n}(S)$ in time polynomial in $|S|$ and $n$?
\end{prob}

This is a special case of the Subpower Membership Problem \ref{subpower-membership}. Even this special case is open (the answer is conjectured to be yes). In the case of groups, the famous Schreier-Sims algorithm gives a positive solution (see \cite{group-algorithms} for a straightforward exposition).

\begin{rem} The proof of correctness of the subroutine \texttt{Nonempty} and the algorithm for \texttt{Fix-values} are both directly connected to the proof of Theorem \ref{malcev-compact}. The subroutines \texttt{Next-beta} and \texttt{Intersect} use the subroutines \texttt{Nonempty} and \texttt{Fix-values} as black boxes and don't involve the algebraic structure at all. Thus, in order to generalize the Mal'cev algorithm to more general algebraic structures, the only new ingredient needed is a proof of a generalization of Theorem \ref{malcev-compact}.
\end{rem}

\subsection{Near-subgroups}

In this subsection, we will describe the maximal polynomial-time solvable extension $\fG^*$ of the relational clone $\fG$ of cosets of subgroups of $\bG^m$, where $\bG$ is a finite group. The relational clone $\fG^*$ will turn out to have a Mal'cev polymorphism, so the algorithm for Mal'cev algebras can be used to prove the dichotomy for extensions of $\fG$.

First, consider the simple case where $\bG = \ZZ/n$ is cyclic of order $n$ at least $3$. It's easy to see that if we add the unary relation $\{0,1\}$ to $\ZZ/n$, then we can simulate 1-IN-3 SAT via the primitive positive formula
\[
x+y+z = 1\ \wedge\ x,y,z \in \{0,1\}.
\]
Using an inductive argument with this as the base case, Feder and Vardi \cite{feder-vardi} show that if we adjoin any unary relation to $\ZZ/n$ which isn't a coset of a subgroup, then we can simulate 1-IN-3 SAT as well.

\begin{prop}[Feder, Vardi \cite{feder-vardi}] If we adjoin any unary relation $K$ to the relational structure $\fG = (\ZZ/n, \{1\}, x+y=z)$, then the resulting CSP is NP-complete unless $K$ is a coset of a subgroup of $\ZZ/n$.
\end{prop}
\begin{proof} Using the binary relation $y = x+i$ for constants $i \in \ZZ/n$, we see that $K-i$ is in the relational clone generated by $K$ and $\fG$. Thus we may assume without loss of generality that $0 \in K$, and by possibly restricting to a subgroup we may assume that $\langle K\rangle = \ZZ/n$. By applying an automorphism of $\ZZ/n$, we may also assume that $1 \in K$.

We induct on $|K|, n$. If there is an $i \ne 0$ with $i,i+1 \in K$, then $K \cap (K-i)$ also contains $0,1$, and will be strictly smaller than $K$ unless $K = K-i$, in which case we may take the quotient by $\langle i\rangle$. Thus we may assume that $i,i+1$ are not both in $K$ for any $i \ne 0$.

If $K$ contains some $i$ with neither of $i,i-1$ relatively prime to $n$, then by induction $K\cap \langle i\rangle$ and $(K-1)\cap \langle i-1\rangle$ are subgroups, so $2i, 2i-1 \in K$ and we must have $2i-1 \equiv 0 \pmod{n}$, contradicting the assumption that $i$ has a common factor with $n$.

If $K$ contains $i\ne 1$ with $i$ relatively prime to $n$, then $K\cap (i-K)$ contains $0,i$ but not $1$, and we may apply the induction hypothesis to get a contradiction. Similarly, if $K$ contains $i\ne 0$ with $i-1$ relatively prime to $n$, then $(K-1)\cap (i-K)$ contains $0$ and $i-1$ but not $-1$, and we may apply the induction hypothesis.

Thus the only case to consider is the case $K = \{0,1\}$, and we have already seen that in this case we can simulate 1-IN-3 SAT (unless $n=2$, in which case $K = \ZZ/n$).
\end{proof}

Next, consider the case where $\bG$ is the Klein four-group $(\ZZ/2)^2$. The only unary relations which aren't already cosets of subgroups of $(\ZZ/2)^2$ are the relations with three elements. If we adjoin any three element unary relation to $(\ZZ/2)^2$, then we can again simulate 1-IN-3 SAT: if we adjoin the relation $K = \{(0,0),(0,1),(1,0)\}$, for instance, then we can use the primitive positive formula
\[
\exists t\ (x,y,z \in \{(0,0),(0,1)\}\ \wedge\ x+y+z = (0,1)\ \wedge\ (x,t) \in \{((0,0),(0,0)), ((0,1),(1,0))\}\ \wedge\ y+t \in K),
\]
which is satisfied iff exactly one of $x,y,z$ is $(0,1)$ and the other two are $(0,0)$.

Now consider the case where $\bG$ is any finite abelian group, and $K \subseteq \bG$ is a unary relation which can be added without creating NP-completeness. Then if any $a,a+b \in K$, we must have $a+ib \in K$ for all $i \in \ZZ$ by the cyclic case. By the Klein four-group case, if we have subgroups $\bN \le \bM \le \bG$ with $\bM/\bN \cong (\ZZ/2)^2$, then if $K$ meets any three elements of $\bM/\bN$ it must also meet the fourth.

\begin{prop} Suppose that $\bG$ is an abelian group and that $K \subseteq \bG$ has $0 \in K$, has the property that if $a,a+b \in K$ then $a+\langle b\rangle \subseteq K$, and the property that for any subgroups $\bN \le \bM \le \bG$ with $\bM/\bN \cong (\ZZ/2)^2$, we have $|(K\cap \bM)/\bN| \ne 3$. Then $K$ must be a subgroup of $\bG$.
\end{prop}
\begin{proof} From the first assumption, for any $a,b \in K$ we must have $-ia, jb \in K$, so
\[
ia+2jb = jb-(-ia-jb) \in jb+\langle -ia - jb\rangle \subseteq K,
\]
and similarly $2ja+ib \in K$ for all $i,j \in \ZZ$.

Thus, if we take $\bM = \langle a,b\rangle$ and $\bN = \langle 2a,2b \rangle$, we see that either $|\bM/\bN| < 4$ in which case $a+b \in \langle a,b \rangle \subseteq K$, or $\bM/\bN \cong (\ZZ/2)^2$ and $|(K\cap \bM)/\bN| \ge 3$. In the latter case, the second assumption implies that there are $i,j$ such that $(2i+1)a+(2j+1)b \in K$. Then $a+b = (2i+1)a+(2j+1)b-2ia-2jb \in K$ by repeated application of the first assumption. Either way, $a+b \in K$, so $K$ is closed under addition.
\end{proof}

\begin{cor} If $\bG$ is a finite abelian group and $\fG$ the associated relational structure, then for any $m$-ary relation $K$ which is not a coset of a subgroup of $\bG^m$, the CSP we get by adding $K$ to $\fG$ is NP-complete.
\end{cor}
\begin{proof} Apply the previous proposition to the abelian group $\bG^m$.
\end{proof}

In the case of nonabelian groups, however, we may be able to adjoin interesting new constraints. Note that if we adjoin any constraint, then we automatically adjoin all of its cosets, since for any constant $b \in \bG$ the relation $y = bx$ is a left coset of the diagonal subgroup of $\bG^2$. So by the abelian case, the only possibilities for new relations are those described by the following definition.

\begin{defn} A subset $K \subseteq \bG$ is a \emph{near subgroup} of $\bG$ if it contains $1$, and for any $b \in K^{-1}$, any $\mathbb{M} \le \bG$ and any $\mathbb{N} \lhd \mathbb{M}$ with $\mathbb{M}/\mathbb{N}$ abelian, the quotient set $(bK\cap \mathbb{M})/\mathbb{N}$ is a subgroup of $\mathbb{M}/\mathbb{N}$.
\end{defn}

\begin{prop} If $\mathbb{H} \le \bG$ is a subgroup and $K \subseteq \mathbb{H}$ is a near subgroup of $\mathbb{H}$, then $K$ is a near subgroup of $\bG$. Similarly, if $\varphi : \bG \twoheadrightarrow \mathbb{H}$ is a surjective group homomorphism and $K \subseteq \mathbb{H}$ is a near subgroup of $\mathbb{H}$, then $\varphi^{-1}(K)$ is a near subgroup of $\mathbb{H}$.
\end{prop}

\begin{thm}[Aschbacher \cite{near-subgroups-aschbacher}]\label{near-subgroups} The intersection of two near subgroups of a finite group is a near subgroup.
\end{thm}

\begin{cor}[Feder \cite{near-subgroups-feder}]\label{near-subgroups-malcev} Let $\bG$ be a finite group, and let $\mathbf{G}^*$ be the relational structure on the underlying set of $\bG$ having as relations all cosets of all near subgroups of $\bG^n$. Then $\mathbf{G}^*$ has a Mal'cev polymorphism.
\end{cor}
\begin{proof} Consider the ``free near subgroup generated by two elements'', that is, the smallest near subgroup $K$ of $\bG^{\bG^2}$ which contains $\pi_1, \pi_2$ (a smallest such near subgroup exists since the intersection of all of them is guaranteed to be a near subgroup as well). Let $\mathbb{N}$ be the commutator subgroup of the group generated by $\pi_1, \pi_2$. Since $\langle \pi_1,\pi_2\rangle / \mathbb{N}$ is abelian, there must be some $c \in \mathbb{N}$ with $\pi_1\pi_2c \in K$ by the definition of a near subgroup.

We define a binary operation $g$ by $g = \pi_1\pi_2c$, that is, $g(x,y) = xyc(x,y)$, where $c \in \bG^{\bG^2}$ is interpreted as a function $c : \bG^2 \rightarrow \bG$. Since for all $x,y$ we know that $c(x,y)$ is contained in the commutator subgroup of $\langle x, y\rangle$, we have $c(x,1) = c(1,x) = 1$ for all $x$, so $g(x,1) = g(1,x) = x$. Now we define a Mal'cev operation $p$ by
\[
p(x,y,z) = yg(y^{-1}x,y^{-1}z) = xy^{-1}zc(y^{-1}x,y^{-1}z).
\]
That $p$ is Mal'cev follows directly from the fact that $g$ satisfies the identities $g(1,x) \approx g(x,1) \approx x$.

To see that $p$ is really a polymorphism of $\mathbf{G}^*$, let $X$ be any coset of any near subgroup of $\bG^n$, and let $x,y,z \in X$. Then $y^{-1}X$ is a near subgroup of $\bG^n$. Since $g = \pi_1\pi_2c \in K$, $g$ preserves every near subgroup of $\bG^n$ (since for any $a,b \in \bG^n$, $K$ is contained in the near subgroup of $\bG^{\bG^2}$ obtained by taking the preimage of the map $\varphi$ from the subgroup generated by $\pi_1, \pi_2$ to $\bG^n$ which sends $\pi_1 \mapsto a, \pi_2 \mapsto b$). Thus from $y^{-1}x, y^{-1}z \in y^{-1}X$ we have $g(y^{-1}x,y^{-1}z) \in y^{-1}X$, and $p(x,y,z) = yg(y^{-1}x,y^{-1}z) \in X$, so $p$ does indeed preserve $X$.
\end{proof}

In order to prove Aschbacher's Theorem \ref{near-subgroups}, we first need a more convenient characterization of near-subgroups.

\begin{defn} A subset $K$ of a finite group $\bG$ is a \emph{twisted subgroup} if $1 \in K$ and $x,y \in K \implies xyx \in K$.
\end{defn}

\begin{prop} If $K$ is a twisted subgroup and $x \in K$, then $\langle x \rangle \subseteq K$, so in particular $K = K^{-1}$. If $b \in K$, then $bK$ is also a twisted subgroup.
\end{prop}
\begin{proof} For the first statement, for any $x \in K$ we have $x^k\cdot 1\cdot x^k, x^k\cdot x\cdot x^k \in K$ for all $k \ge 0$, so $\langle x \rangle \subseteq K$. For the second statement, if $x,y \in bK$ and $b \in K$, then $b^{-1}(x\cdot y\cdot x) = (b^{-1}x) \cdot (b\cdot b^{-1}y \cdot b) \cdot (b^{-1}x) \in K$, so $xyx \in bK$.
\end{proof}

\begin{prop} A subset $K \subseteq \bG$ is a near-subgroup iff it is a twisted subgroup such that for any $b \in K^{-1}$, any $\mathbb{M} \le \bG$ and any $\mathbb{N} \lhd \mathbb{M}$ with $\mathbb{M}/\mathbb{N}$ isomorphic to the Klein four-group, $|(bK\cap \mathbb{M})/\mathbb{N}| \ne 3$.
\end{prop}
\begin{proof} We just need to check that $K$ being a twisted subgroup is equivalent to $\langle x\rangle \subseteq bK$ for all $x,b$ with $x \in bK, b \in K^{-1}$. The previous proposition proves one direction of the equivalence. For the other direction, if $x,y \in K$, then $yx \in yK$ and $y^{-1} \in \langle y \rangle \subseteq K$, so $(yx)^2 \in \langle yx\rangle \subseteq yK$, which is equivalent to $xyx = y^{-1}(yx)^2 \in K$.
\end{proof}

\begin{ex}\label{ex-heisenberg-near-subgroup} An explicit example of a near-subgroup which is not a subgroup is given in \cite{feder-vardi}. Let $\bG$ be the Heisenberg group of order $p^3$ ($p$ odd):
\[
\bG = \left\{\begin{bmatrix} 1 & a & c\\ 0 & 1 & b\\ 0 & 0 & 1\end{bmatrix} \text{ s.t. } a,b,c \in \ZZ/p\right\} \le \operatorname{SL}_3(\ZZ/p).
\]
Let $K \subseteq \bG$ be given by
\[
K = \left\{\begin{bmatrix} 1 & a & \frac{ab}{2}\\ 0 & 1 & b\\ 0 & 0 & 1\end{bmatrix} \text{ s.t. } a,b \in \ZZ/p\right\}.
\]
Since $\bG$ has odd order, to check that $K$ is a near subgroup we just need to check that it is a twisted subgroup, i.e. that it contains the identity and is closed under the binary operation $x,y \mapsto xyx$. This can be checked by direct calculation: for any $a,b,c,d \in \ZZ/p$ we have
\[
\begin{bmatrix} 1 & a & \frac{ab}{2}\\ 0 & 1 & b\\ 0 & 0 & 1\end{bmatrix}\begin{bmatrix} 1 & c & \frac{cd}{2}\\ 0 & 1 & d\\ 0 & 0 & 1\end{bmatrix}\begin{bmatrix} 1 & a & \frac{ab}{2}\\ 0 & 1 & b\\ 0 & 0 & 1\end{bmatrix} = \begin{bmatrix} 1 & 2a+c & \frac{(2a+c)(2b+d)}{2}\\ 0 & 1 & 2b+d\\ 0 & 0 & 1\end{bmatrix}.
\]
That $K$ is not a subgroup follows from
\[
\begin{bmatrix} 1 & 0 & 0\\ 0 & 1 & 1\\ 0 & 0 & 1\end{bmatrix}\begin{bmatrix} 1 & 1 & 0\\ 0 & 1 & 0\\ 0 & 0 & 1\end{bmatrix} = \begin{bmatrix} 1 & 1 & 0\\ 0 & 1 & 1\\ 0 & 0 & 1\end{bmatrix} \not\in K.
\]
\end{ex}

That we needed to take $p$ odd in the above example is no coincidence, as the next proposition shows.

\begin{prop} If $\bG$ is a $2$-group, then any near-subgroup of $\bG$ is a subgroup of $\bG$.
\end{prop}
\begin{proof} We prove this by induction on the order of $\bG$. Let $K$ be a near-subgroup of $\bG$, and assume without loss of generality that $\langle K\rangle = \bG$.

Let $z \in Z(\bG)$ be a nontrivial involution in the center of $\bG$, which must exist since every $2$-group has a nontrivial center, and every nontrivial element of the center has a power which is a nontrivial involution. By induction we have $K/\langle z\rangle = \bG/\langle z\rangle$.

If $z \in K$ then for any $g \in \bG\setminus K$, we have $gz \in K$, and from $\langle gz,z\rangle$ abelian and $gz,z \in K$ we have $\langle gz,z\rangle \subseteq K$, so in particular $g = gz\cdot z \in K$. Thus if $z \in K$ we have $\bG = K$.

Thus we may assume that $z \not\in K$, and in fact that $Z(\bG) \cap K = 1$. Let $\bM$ be a maximal subgroup of $\bG$ containing $\langle z\rangle$, then by induction we have $\bM \cap K$ a subgroup of $\bM$. Since $(\bM \cap K)/\langle z\rangle = \bM/\langle z\rangle$, we have $\bM \cong (\bM\cap K)\times \langle z\rangle$.

Let $\Phi(\bM)$ be the Frattini subgroup of $\bM$, which for $2$-groups is given by $\Phi(\bM) = \bM^2[\bM,\bM]$ (where by $\bM^2$ we mean the collection of all squares $a^2$ for $a \in \bM$). Then from $\Phi(\langle z\rangle) = 1$ we have $\Phi(\bM) = \Phi(\bM\cap K)$, and since $\bM \lhd \bG$ we have $\Phi(\bM) \lhd \bG$. Thus if $\Phi(\bM) \ne 1$ then by considering parities of the sizes of the orbits of elements of $\Phi(\bM)$ under conjugation we see that $\Phi(\bM \cap K) = \Phi(\bM)$ contains a nontrivial element of $Z(\bG)$, contradicting $K \cap Z(\bG) = 1$. Thus $\Phi(\bM) = 1$, so $\bM$ has exponent $2$. Since this holds for every maximal subgroup of $\bG$ which contains $\langle z\rangle$, we see that $\bG$ has exponent $2$, so $\bG$ is abelian.
\end{proof}

Next we show that we can reduce to the situation where $\langle K\rangle$ has an automorphism of order two which sends $k$ to $k^{-1}$ for all $k \in K$.

\begin{defn} If $K$ is a twisted subgroup, we define the $K$-\emph{radical} $\Xi_K$ to be the set of elements of the form $k_1 \cdots k_n$ with $k_i \in K$ such that $k_1^{-1}\cdots k_n^{-1} = 1$.
\end{defn}

\begin{prop} If $K$ is a twisted subgroup and $\Xi_K$ is the $K$-radical, then $\Xi_K$ is a normal subgroup of $\langle K\rangle$, and for any $x \in K$ we have $x\Xi_K \subseteq K$.
\end{prop}
\begin{proof} To see that $\Xi_K$ is normal in $\langle K\rangle$, just note that for any $b \in K$ we have
\[
k_1^{-1}\cdots k_n^{-1} = 1 \iff b^{-1}k_1^{-1}\cdots k_n^{-1}b = 1 \implies bk_1\cdots k_nb^{-1} \in \Xi_K,
\]
so $b\Xi_Kb^{-1} \subseteq \Xi_k$.

For the second statement, note that $k_1^{-1} \cdots k_n^{-1} = 1 \iff k_n\cdots k_1 = 1$, so if $x \in K$ then we have
\[
x(k_1\cdots k_n) = (k_n\cdots k_1)x(k_1\cdots k_n) = k_n(\cdots(k_1 x k_1)\cdots)k_n \in K.\qedhere
\]
\end{proof}

\begin{prop} If $K$ is a twisted subgroup with $\Xi_K = 1$, and if $\tau$ satisfies $\tau^2 = 1$, $\tau k\tau = k^{-1}$ for $k \in K$, then $\tau K$ is preserved under conjugation by elements of $\langle K,\tau\rangle$.
\end{prop}
\begin{proof} If $x,y \in K$, then $x^{-1} \tau y x = \tau xyx \in \tau K$, and $\tau^{-1} \tau y\tau = \tau y^{-1} \in \tau K$.
\end{proof}

\begin{prop} A twisted subgroup $K \subseteq \bG$ is a near-subgroup of $\bG$ iff the intersection $bK \cap \bS$ is a subgroup of $\bS$ for every $2$-Sylow subgroup $\bS$ of $\bG$ and every $b \in K^{-1}$.
\end{prop}
\begin{proof} Suppose for contradiction that $\bM \le \bG$, $\bN\lhd \bM$ with $\bM/\bN$ isomorphic to the Klein four-group, and $b \in K^{-1}$ with $|(bK \cap \bM)/\bN| = 3$. We may assume without loss of generality that $\bM = \langle K\rangle \cap \bM$ and $\bN = \langle K\rangle \cap \bN$, that $\bG = \langle K\rangle$, that $b = 1$, and that $\Xi_K = 1$. From $\Xi_K = 1$, we see that there is an order $2$ automorphism $\tau$ of $\bG = \langle K\rangle$ with $k^\tau = k^{-1}$ for all $k \in K$, so we work in the semidirect product of $\bG$ and $\langle \tau \rangle$, with $\tau^2 = 1$ and $\tau g\tau = g^\tau$ for $g \in \bG$.

Let $x,y \in K$ be representatives of the nontrivial elements of $(K \cap \bM)/\bN$. We may assume without loss of generality that $x,y$ have orders equal to powers of $2$, since otherwise we may replace them with odd powers of themselves. Let $\bS_x,\bS_y$ be $2$-Sylow subgroups of $\bM\langle\tau\rangle$ containing $\langle x,\tau\rangle,\langle y,\tau\rangle$, respectively, then by the Sylow theorems there is some $g \in \bM\langle\tau\rangle$ with $g^{-1}\bS_yg = \bS_x$. Then $x,\tau,g^{-1}yg,g^{-1}\tau g \in \bS_x$, and our strategy is to show that $x,g^{-1}yg \in K\cap \bS_x$.

We have $g^{-1}\tau yg \in \tau K$ by the previous proposition, so $\tau g^{-1}\tau y g \in K \cap \bS_x$, and similarly $\tau g^{-1}\tau g \in K\cap \bS_x$. Since $K\cap \bS_x$ is assumed to be a subgroup, we have
\[
xg^{-1}yg = x(\tau g^{-1} \tau g)^{-1} (\tau g^{-1} \tau yg) \in K\cap \bS_x.
\]
Then since $\bM/\bN$ is abelian and $\tau y\tau = y^{-1} \equiv_\bN y$, we have $xg^{-1}yg \equiv_\bN xy$, contradicting the assumption $|(K\cap \bM)/\bN| = 3$.
\end{proof}

\begin{proof}[Proof of Theorem \ref{near-subgroups}] If $K,K'$ are near-subgroups of $\bG$, then they are both twisted subgroups and so their intersection $K \cap K'$ is also a twisted subgroup. Now let $\bS$ be any $2$-group contained in $\bG$, then for any $b \in K\cap K'$ we see that $bK\cap bK' \cap \bS = (bK \cap \bS) \cap (bK' \cap \bS)$ is an intersection of subgroups of $\bS$, so it is a subgroup of $\bS$, and the previous proposition shows that this implies that $K\cap K'$ is a near-subgroup of $\bG$.
\end{proof}

\begin{ex} By Corollary \ref{near-subgroups-malcev} to Theorem \ref{near-subgroups}, the Heisenberg group $\bG$ of order $p^3$ from Example \ref{ex-heisenberg-near-subgroup} should have a Mal'cev term operation which preserves all near-subgroups of $\bG^n$ for all $n$, and in particular which preserves the near-subgroup $K \subseteq \bG$. When $p = 3$, we can use the identity $(x^{-1}z)^3 \approx \operatorname{Id}$ to see that
\[
x^{-1}zx^{-1} \approx z^{-1}xz^{-1}
\]
in $\bG$, so the term operation
\[
q(x,y,z) \coloneqq yx^{-1}zx^{-1}y \approx yz^{-1}xz^{-1}y
\]
defines a Mal'cev term operation, contained in the clone generated by the standard Mal'cev operation $(x,y,z) \mapsto xy^{-1}z$, which satisfies the additional identity
\[
q(x,y,z) \approx q(z,y,x).
\]
This Mal'cev term preserves $K$, and $q$ is given explicitly by the formula
\begin{align*}
q\left(\begin{bmatrix} 1 & a_1 & c_1 + \frac{a_1b_1}{2}\\ 0 & 1 & b_1\\ 0 & 0 & 1\end{bmatrix}, \begin{bmatrix} 1 & a_2 & c_2 + \frac{a_2b_2}{2}\\ 0 & 1 & b_2\\ 0 & 0 & 1\end{bmatrix}, \begin{bmatrix} 1 & a_3 & c_3 + \frac{a_3b_3}{2}\\ 0 & 1 & b_3\\ 0 & 0 & 1\end{bmatrix}\right)\\
= \begin{bmatrix} 1 & a_1 - a_2 + a_3 & c_1 - c_2 + c_3 + \frac{(a_1 - a_2 + a_3)(b_1 - b_2 + b_3)}{2}\\ 0 & 1 & b_1 - b_2 + b_3\\ 0 & 0 & 1\end{bmatrix},
\end{align*}
so we see that in fact we have an isomorphism between the reduct $(\bG, q(x,y,z))$ and the power $(\ZZ/3, x - y + z)^3$ (this isomorphism is the exponential map from the Heisenberg Lie algebra to the Heisenberg group, taken modulo $3$).

We get similar isomorphisms from reducts of the Heisenberg group to $(\ZZ/p, x-y+z)^3$ for every odd prime $p$, but not for $p = 2$. For larger primes $p$, the Mal'cev operation $q(x,y,z)$ to use is given by
\[
q(x,y,z) \coloneqq (xy^{-1}z)\big((zy^{-1}x)^{-1}(xy^{-1}z)\big)^{\frac{p-1}{2}} \approx (zy^{-1}x)\big((xy^{-1}z)^{-1}(zy^{-1}x)\big)^{\frac{p-1}{2}}.
\]
When $p = 2$, the Heisenberg group is isomorphic to the dihedral group $D_4$.
\end{ex}


\section{Abelian Mal'cev algebras are affine}\label{s-abelian-malcev}

In this section we will prove that abelian Mal'cev algebras are affine. This is an important step in the proof that problems which do not have the ``ability to count'' have bounded width. First we will carefully define what an affine algebra is, starting with the more basic concept of a quasi-affine algebra.

\begin{defn} An algebra $\bA$ is called \emph{quasi-affine} if there is an abelian group $\bG = (G,0,+,-)$ with underlying set $G$ containing the underlying set of $\bA$, such that the restriction of the $4$-ary relation $x+y = z+w$ to $\bA$ is preserved by all the operations of $\bA$.
\end{defn}

We want to relate this to the more familiar concept of a module over a ring.

\begin{defn} If $\RR$ is a ring and $\bM$ is a module over $\RR$ with underlying group $(M,0,+,-)$, then we consider $\bM$ to be a universal algebraic object $(M,0,+,-,\{\phi_r\}_{r \in \RR})$, where for each $r\in\RR$ the unary operation $\phi_r : \bM \rightarrow \bM$ is given by $\phi_r : m \mapsto rm$.

In general, a universal algebraic object is called a \emph{module} if it is an expansion of an abelian group by any collection of unary operations that distribute over addition.
\end{defn}

The way these concepts are related is a coarser notion than term equivalence, known as \emph{polynomial equivalence} (warning: in some older references, ``polynomial equivalence'' means term equivalence and ``functional equivalence''/``algebraic equivalence'' means polynomial equivalence).

\begin{defn} If $\cO$ is any set of operations, then the \emph{polynomial clone} generated by $\cO$ is the clone generated by $\cO$ together with the constant functions (one for each element of the underlying set). Two algebras or clones on the same underlying set are called \emph{polynomially equivalent} if they have the same polynomial clones.
\end{defn}

\begin{prop} An algebra $\bA$ is quasi-affine iff it is a subalgebra of a reduct of the polynomial clone of a module.
\end{prop}
\begin{proof} Let $\bA$ be a quasi-affine algebra, and let $\bG = (G,0,+,-)$ be the corresponding group. We may assume without loss of generality that $0 \in \bA$, and that $\bG$ is the abelian group with the following presentation: the generators are the elements of $\bA\setminus\{0\}$, and the relations are given by $x+y-z-w = 0$ for every quadruple of elements $x,y,z,w \in \bA$ such that $x+y = z+w$ in $\bG$.

Suppose that $f$ is any $n$-ary operation of $\bA$, and for each $i \le n$ let $\phi_i : \bA \rightarrow G$ be the unary operation given by
\[
\phi_i(x) = f(0,...,0,x,0,...,0) - f(0,...,0),
\]
with the $x$ in the $i$th position. Since $f$ preserves the relation $x+y = z+w$ on $\bA$, we have
\[
\phi_i(x) + \phi_i(y) = \phi_i(z) + \phi_i(w)
\]
for any $x,y,z,w \in \bA$ such that $x+y = z+w$ in $\bG$. Thus $\phi_i$ is compatible with the defining relations of $\bG$, so the map $\phi_i : \bA \rightarrow G$ extends to a unique homomorphism $\phi_i : \bG \rightarrow \bG$.

To finish, we just need to prove that
\[
f(x_1, ..., x_n) = \phi_1(x_1) + \cdots + \phi_n(x_n) + f(0,...,0)
\]
for all $x_1, ..., x_n \in \bA$, since $f(0,...,0)$ is a constant operation.

We prove this by induction on the number $k$ of nonzero values among $x_1, ..., x_n$. The base cases $k = 0,1$ follow from the definition of the $\phi_i$. For the inductive step, assume without loss of generality that the nonzero values of the $x_i$s are $x_1, ..., x_{k+1}$. Since $f$ preserves the relation $x+y = z+w$, we have
\[
f(x_1, ..., x_{k+1}, 0, ..., 0) + f(0,...,0) = f(x_1, ..., x_k, 0, 0, ..., 0) + f(0, ..., 0, x_{k+1}, 0, ..., 0),
\]
so by the inductive hypothesis and the definition of $\phi_{k+1}$ we have
\[
f(x_1, ..., x_{k+1}, 0, ..., 0) = \phi_1(x_1) + \cdots + \phi_k(x_k) + f(0, ..., 0) + \phi_{k+1}(x_{k+1}).\qedhere
\]
\end{proof}

\begin{defn} An algebra $\bA$ is called \emph{affine} if it is polynomially equivalent to a module.
\end{defn}

\begin{prop}\label{affine-malcev} An algebra is affine iff it is quasi-affine and has a Mal'cev term.
\end{prop}
\begin{proof} The hardest step is showing that every affine algebra $\bA$ has a Mal'cev term. Since $\bA$ is polynomially equivalent to a module, there must be some $n+3$-ary term $t$ and some constants $a_1, ..., a_n \in \bA$ such that
\[
t(x,y,z,a_1, ..., a_n) = x-y+z
\]
for all $x,y,z$. Since any affine algebra is quasi-affine, we can write $t$ in the form
\[
t(x,y,z,u_1,...,u_n) = x-y+z + \sum_i \phi_i(u_i) + c
\]
for some unary $\phi_i$ and some constant $c$. Define $p(x,y,z)$ by
\[
p(x,y,z) = t(x,t(y,y,y,x,...,x),z,x,...,x).
\]
Then $p$ is a term operation of $\bA$, and we have
\[
p(x,y,z) = x-\big(y-y+y + \sum_i \phi_i(x) + c\big)+z + \sum_i \phi_i(x) + c = x-y+z,
\]
so $p$ is Mal'cev.

For the converse, if $\bA$ is quasi-affine and has a Mal'cev term $p$, then $p(x,y,y) \approx p(y,y,x) \approx x$ imply that $p(x,0,0) = x, p(0,0,z) = z$, and $p(y,y,0) = y + p(0,y,0) = 0$, so we must have $p(x,y,z) = x-y+z$. Thus $x+z = p(x,0,z)$ and $x-y = p(x,y,0)$ are polynomial operations of $\bA$, and therefore for each term $f$ of $\bA$ the unary function $\phi(x) = f(x,0,...,0) - f(0,...,0)$ is a polynomial operation of $\bA$ as well.
\end{proof}


It is less trivial to give a universal algebraic definition of what it means to be \emph{abelian}. We will give several different definitions and prove that they are equivalent to each other, and that they restrict to the right concept in the special case of groups.

\begin{defn} An algebraic structure $\bA$ is called \emph{abelian} if there is a congruence $\Theta$ on $\bA\times \bA$ such that the diagonal $\Delta_\bA = \{(a,a) \mid a \in \bA\}$ is one of the congruence classes of $\Theta$.
\end{defn}

\begin{prop} A group is abelian iff it is commutative.
\end{prop}
\begin{proof} A group $\bG$ is abelian iff the diagonal $\Delta_\bG$ is a normal subgroup of $\bG\times \bG$. To check that $\Delta_\bG$ is normal, we just need to check that it is closed under conjugation by elements of the form $(1,b)$ for all $b \in \bG$. Since
\[
(1,b)(a,a)(1,b)^{-1} = (a,bab^{-1}),
\]
the normality of $\Delta_\bG$ is equivalent to the identity $a \approx bab^{-1}$, which is equivalent to $ab \approx ba$.

Alternatively, we can argue as follows. The group $\bG$ is commutative iff the map $\bG \rightarrow \bG$ given by $x \mapsto x^{-1}$ is a homomorphism, and if this occurs then there is a homomorphism $\bG \times \bG \rightarrow \bG$ such that the restriction $\bG\times \{1\} \rightarrow \bG$ is the identity, and such that the diagonal maps to $\{1\}$. Conversely, if the diagonal is a normal subgroup, then every coset intersects $\bG\times \{1\}$ and $\{1\}\times \bG$ exactly once, so the quotient $\bG\times \bG/\Delta_\bG$ is isomorphic to $\bG$ in two different ways, and composing these isomorphisms we obtain the map $x \mapsto x^{-1}$, so $\bG$ is commutative.
\end{proof}

Now we give a second definition of abelian, which is phrased in a way which is closely related to the concept of a ``commutator'' of congruences in a general algebraic structure.

\begin{defn} We say that an algebraic structure $\bA$ satisfies the \emph{term condition} if for all terms $t \in \Clo_{n+1}(\bA)$ and all $u,v \in \bA$, $a_i,b_i \in \bA$ for $i \le n$, we have
\[
t(u,a_1, ..., a_n) = t(u,b_1, ..., b_n) \iff t(v,a_1, ..., a_n) = t(v,b_1, ..., b_n).
\]
\end{defn}

\begin{prop} An algebra $\bA$ is abelian iff it satisfies the term condition.
\end{prop}
\begin{proof} We think of congruences on $\bA^2$ as subalgebras of $\bA^{2\times 2}$, the set of $2\times 2$ matrices with entries in $\bA$ (here elements of $\bA^2$ are visualized as column vectors, and an element of $\bA^{2\times 2}$ is viewed as a row vector of column vectors). To understand the smallest congruence on $\bA^2$ with $\Delta_\bA$ contained in a congruence class, we consider the relation $\bM \le \bA^{2\times 2}$ generated by matrices of the form
\[
\begin{bmatrix} u & v\\ u & v\end{bmatrix}, \;\;\; \begin{bmatrix} a & a\\ b & b\end{bmatrix},
\]
where the first type of matrix corresponds to the fact that any two elements of $\Delta_\bA$ are congruent, while the second type of matrix corresponds to the fact that every element of $\bA^2$ is congruent to itself. Then considering $\bM$ as a binary relation on $\bA^2$, the transitive closure of $\bM$ is a congruence $\Theta$ on $\bA^2$, and it is clearly as small as possible given that $\Delta_\bA$ is contained in a congruence class of $\Theta$.

To understand whether $\Delta_\bA$ is a congruence class of $\Theta$, it's enough to check whether $\Delta_\bA$ meets any element of $\bA^2 \setminus \Delta_\bA$ in $\bM$. This occurs (that is, $\bA$ is nonabelian) iff there is some term $t \in \Pol_{m+n}(\bA)$ and some $u_i,v_i \in \bA$ for $i \le m$, $a_i,b_i \in \bA$ for $i \le n$ such that
\[
t(u_1, ..., u_m, a_1, ..., a_n) = t(u_1, ..., u_m, b_1, ..., b_n)
\]
but
\[
t(v_1, ..., v_m, a_1, ..., a_n) \ne t(v_1, ..., v_m, b_1, ..., b_n).
\]
So if $\bA$ is abelian, then it certainly satisfies the term condition (just take $m=1$ in the above). Conversely, if $\bA$ satisfies the term condition, then we will show that the above situation can't happen by induction on $m$. We just note that by the induction hypothesis, we have
\[
t(u_1, ..., u_m, a_1, ..., a_n) = t(u_1, ..., u_m, b_1, ..., b_n) \implies t(v_1, ..., v_{m-1}, u_m, a_1, ..., a_n) = t(v_1, ..., v_{m-1}, u_m, b_1, ..., b_n),
\]
and then by the term condition applied to a version of $t$ with variables permuted so that the $m$th variable becomes the first, this implies that
\[
t(v_1, ..., v_m, a_1, ..., a_n) = t(v_1, ..., v_m, b_1, ..., b_n).\qedhere
\]
\end{proof}

\begin{prop} Every quasi-affine algebra satisfies the term condition and is therefore abelian.
\end{prop}
\begin{proof} If $t$ is an $n+1$-ary term operation of a quasi-affine algebra, then we can write $t$ in the form
\[
t(x_0, ..., x_n) = \phi_0(x_0) + \cdots + \phi_n(x_n) + c,
\]
where the $\phi_i$ are unary and $c$ is a constant. Then for any $u \in \bA, a_i, b_i \in \bA$, we have
\[
t(u,a_1, ..., a_n) = t(u,b_1,...,b_n) \iff \phi_1(a_1) + \cdots + \phi_n(a_n) = \phi_1(b_1) + \cdots + \phi_n(b_n),
\]
and this is a condition which does not depend on the value of $u$.
\end{proof}

\begin{ex} If a group is commutative, then it is affine, so it satisfies the term condition. Conversely, if a group satisfies the term condition for the binary term $t(x,y) = yxy^{-1}$, then the group is commutative, since we have $t(1,1) = t(1,y) \iff t(x,1) = t(x,y)$, that is, $1 = yy^{-1} \iff x = yxy^{-1}$.
\end{ex}

\begin{ex} A ring is abelian in the sense of universal algebra iff it is a \emph{zero ring}, that is, a ring satisfying the identity $xy \approx 0$. To see the necessity, we apply the term condition with the term $t(x,y) = xy$ and the pairs $(u,v) = (0,x)$ and $(a,b) = (0,y)$, to see that $0\cdot 0 = 0 \cdot y \iff x\cdot 0 = x\cdot y$. To see the sufficiency, note that every zero ring is affine.
\end{ex}

\begin{ex} The quasigroup with multiplication table
\begin{center}
\begin{tabular}{c|cccc} $\cdot$ & $0$ & $1$ & $2$ & $3$\\ \hline $0$ & $3$ & $2$ & $0$ & $1$\\ $1$ & $2$ & $3$ & $1$ & $0$\\ $2$ & $1$ & $0$ & $2$ & $3$\\ $3$ & $0$ & $1$ & $3$ & $2$\end{tabular}
\end{center}
is abelian, but is neither commutative nor associative. In fact it is affine, with underlying group equal to the Klein four-group: the multiplication can be written as $x\cdot y = x \oplus \phi(y) \oplus 3$, where $\phi$ is the transposition $(2\ 3)$. This example is from \cite{commutator-theory}.
\end{ex}

In terms of congruence lattices, the main important feature of an affine algebra $\bA$ is that $\Con(\bA\times\bA)$ contains the following five element sublattice.
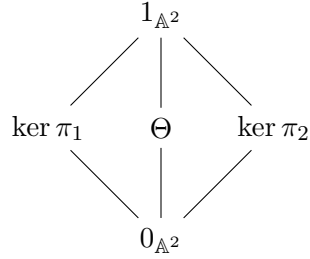
\begin{center}
\begin{tikzpicture}[scale=1.5]
  \node (1) at (0,1) {$1_{\bA^2}$};
  \node (p1) at (-1,0) {$\ker \pi_1$};
  \node (t) at (0,0) {$\Theta$};
  \node (p2) at (1,0) {$\ker \pi_2$};
  \node (0) at (0,-1) {$0_{\bA^2}$};
  \draw (1) -- (p1) -- (0) -- (p2) -- (1);
  \draw (1) -- (t) -- (0);
\end{tikzpicture}
\end{center}
The abstract five element lattice corresponding to this picture is known as the diamond lattice $\cM_3$. The lattice $\cM_3$ has a special role in lattice theory: every modular lattice which isn't distributive contains a sublattice which is isomorphic to $\cM_3$ (see Proposition \ref{distributive-m3} in the appendix).

\begin{thm}\label{abelian-M3} If $\bA$ is an abelian Mal'cev algebra, and if $\Theta$ is any congruence of $\bA^2$ which contains the diagonal $\Delta_\bA$ as a congruence class, then the congruences $\Theta, \ker \pi_1, \ker \pi_2$ generate a five element sublattice of $\Con(\bA^2)$ isomorphic to $\cM_3$.
\end{thm}
\begin{proof} In general, we always have $\ker \pi_1 \vee \ker \pi_2 = 1_{\bA^2}$ and $\ker \pi_1 \wedge \ker \pi_2 = 0_{\bA^2}$. Since every element of $\bA$ is congruent under $\ker \pi_1$ to an element of the diagonal $\Delta_\bA$, we have $\ker \pi_1 \vee \Theta = 1_{\bA^2}$, and similarly $\ker \pi_2 \vee \Theta = 1_{\bA^2}$.

All that remains is to check that $\Theta \wedge \ker \pi_1 = \Theta\wedge \ker \pi_2 = 0_{\bA^2}$, and this is where we will use the assumption that $\bA$ has a Mal'cev term $p$. If $(a,b)$ is congruent to $(c,d)$ modulo $\Theta \wedge \ker \pi_1$, then we must have $a = c$. Then 
\[
\begin{bmatrix} b\\ d\end{bmatrix} = p\left(\begin{bmatrix} b\\ b\end{bmatrix}, \begin{bmatrix} a\\ b\end{bmatrix}, \begin{bmatrix} a\\ d\end{bmatrix}\right) \equiv_\Theta p\left(\begin{bmatrix} b\\ b\end{bmatrix}, \begin{bmatrix} a\\ b\end{bmatrix}, \begin{bmatrix} a\\ b\end{bmatrix}\right) = \begin{bmatrix} b\\ b\end{bmatrix} \in \Delta_\bA,
\]
so $(b,d) \in \Delta_\bA$, that is, $b = d$. So from $(a,b) \equiv_{\Theta \wedge \ker \pi_1} (c,d)$ we have shown $(a,b) = (c,d)$, that is, we have $\Theta \wedge \ker \pi_1 = 0_{\bA^2}$.
\end{proof}


The idea now is to study the \emph{equivalence class geometry} on $\bA^2$, where points are elements of $\bA^2$, lines correspond to congruence classes of congruences, and two lines are considered \emph{parallel} if they are both congruence classes of the same congruence. The three congruences $\ker \pi_1, \Theta, \ker \pi_2$ on an abelian Mal'cev algebra give us a particularly nice type of combinatorial geometry.

\begin{defn} An \emph{S-3-system} is a set of points $S$ together with three parallel classes of lines $\Theta_1, \Theta_2, \Theta_3$ on $S$, which satisfy the following properties:
\begin{itemize}
\item for any point $p \in S$ and any $i \le 3$, there is exactly one line $l_i$ of $\Theta_i$ which contains $p$, and

\item if $l_i, l_j$ are lines of $\Theta_i, \Theta_j$, respectively, with $i \ne j$, then their intersection $l_i \cap l_j$ contains exactly one point $p \in S$.
\end{itemize}
Equivalently, an S-3-system is a relational structure $(S,\Theta_1,\Theta_2,\Theta_3)$ such that:
\begin{itemize}
\item each $\Theta_i$ is an equivalence relation on $S$,

\item for $i \ne j$ we have $\Theta_i \wedge \Theta_j = 0_S$, and

\item for $i \ne j$ we have $\Theta_i \circ \Theta_j = 1_S$.
\end{itemize}
The assumption $\Theta_i \wedge \Theta_j = 0_S$ says that any pair of non-parallel lines intersect in \emph{at most} one point, while the assumption $\Theta_i \circ \Theta_j = 1_S$ says that any pair of non-parallel lines intersect in \emph{at least} one point.
\end{defn}

\begin{cor} If $\bA$ is an abelian Mal'cev algebra and $\Theta$ is any congruence of $\bA^2$ with the diagonal as a congruence class, then $(\bA^2, \ker \pi_1, \ker \pi_2, \Theta)$ is an S-3-system with a Mal'cev polymorphism.
\end{cor}

From here on we will classify S-3-systems which have Mal'cev polymorphisms, following Gumm's approach \cite{abelian-malcev-affine}. As a preliminary result, we will show that every S-3-system has a coordinate system which describes the three parallel classes of lines in terms of a \emph{loop} (recall that a loop is just a quasigroup which has an identity).

\begin{lem}\label{S3-coords} If $(S,\Theta_1,\Theta_2,\Theta_3)$ is an S-3-system, and $e$ is any point of $S$, then there is a loop $\bL = (L,\cdot,1)$ and a bijection $L\times L \rightarrow S$ with $(1,1) \mapsto e$, such that for any $x,y,x',y' \in L$ we have
\begin{align*}
(x,y) \equiv_{\Theta_1} (x',y') \;\; &\iff \;\; x = x',\\
(x,y) \equiv_{\Theta_2} (x',y') \;\; &\iff \;\; y = y',\\
(x,y) \equiv_{\Theta_3} (x',y') \;\; &\iff \;\; x\cdot y = x'\cdot y',
\end{align*}
where we have implicitly identified $S$ with $L\times L$.
\end{lem}
\begin{proof} Take $L$ to be the line $l_1$ through $e$ in the parallel class $\Theta_1$, and take $1 = e$. Let $l_2$ be the line through $e$ in the parallel class $\Theta_2$. Then there is a bijection between elements of $l_1$ and elements of $l_2$, taking $x \in l_1$ to $y \in l_2$ when $x,y$ are on a line $l_3$ in the parallel class $\Theta_3$: each $x$ is in a unique such line $l_3$, and each $l_3$ intersects $l_2$ in a unique $y$. Using this bijection, we identify the elements of $l_2$ with $L$ as well.

Now we note that for any point $p \in S$, there is a unique pair of lines $l_1' \in \Theta_1, l_2' \in \Theta_2$ with $l_1' \cap l_2' = \{p\}$. So we can uniquely identify the point $p$ by describing the point $x \in l_1 \cap l_2'$ and the point $y \in l_2 \cap l_1'$ - this gives us the desired bijection between $L\times L$ and $S$.

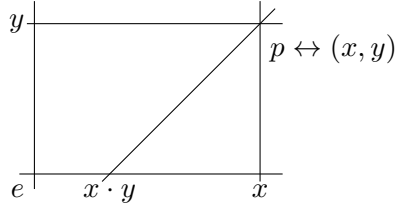
\begin{center}
\begin{tikzpicture}[scale=1.0] 
\draw (-0.2,0) -- (3.3,0);
\draw (0,-0.2) -- (0,2.3);
\draw (3,-0.1) -- (3,2.3);
\draw (-0.1,2) -- (3.3,2);
\draw (3.2,2.2) -- (0.9,-0.1);
\node[align=left, below left] (e) at (0,0){$e$};
\node[align=left, below] (x) at (3,0){$x$};
\node[align=left, below] (xy) at (1,0){$x\cdot y$};
\node[align=left, left] (y) at (0,2){$y$};
\node[align=left, below right] (p) at (3,2){$p \leftrightarrow (x,y)$};
\end{tikzpicture}
\end{center}

Finally, to define the multiplication $\cdot$ on $L$, note that for every $x,y \in L$ there is a point $p \in S$ corresponding to $(x,y)$, and this point $p$ is in a unique line $l_3 \in \Theta_3$. We then define $x\cdot y$ to be the element of $L$ corresponding to the point $l_3 \cap l_1$, or alternatively to the point $l_3 \cap l_2$ (which corresponds to the same element of $L$ by the way we identified points of $l_2$ with points of $l_1$).
\end{proof}

The key observation is that the Mal'cev operation is completely determined by the geometry of the configuration.

\begin{lem} If an S-3-system $\mathbf{S} = (S,\Theta_1,\Theta_2,\Theta_3)$ has a Mal'cev polymorphism $p$, then $p$ is completely determined by $\mathbf{S}$. In fact, $p(x,y,z)$ can be ``geometrically constructed'' from the points $x,y,z$.
\end{lem}
\begin{proof} First consider the special case where $x,y$ lie on a line $l_1$ and $y,z$ lie on a different line $l_2$. Suppose that $l_1 \in \Theta_1$ and $l_2 \in \Theta_2$. Then $p(x,y,z) \equiv_{\Theta_1} p(y,y,z) = z$ and $p(x,y,z) \equiv_{\Theta_2} p(x,y,y) = x$, so if we draw the line $l_2' \in \Theta_2$ through $x$ and the line $l_1' \in \Theta_1$ through $z$, we see that $p(x,y,z)$ is the intersection point $l_1' \cap l_2'$.

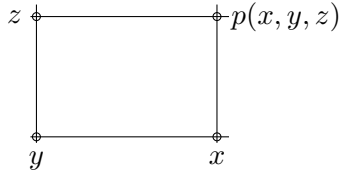
\begin{center}
\begin{tikzpicture}[scale=0.8] 
\draw (-0.1,0) -- (3.2,0);
\draw (0,-0.1) -- (0,2.2);
\draw (3,-0.1) -- (3,2.2);
\draw (-0.1,2) -- (3.2,2);
\node[circle, minimum width=3pt, draw, inner sep=0pt, label=below:$y$] (y) at (0,0){};
\node[circle, minimum width=3pt, draw, inner sep=0pt, label=below:$x$] (x) at (3,0){};
\node[circle, minimum width=3pt, draw, inner sep=0pt, label=left:$z$] (z) at (0,2){};
\node[circle, minimum width=3pt, draw, inner sep=0pt, label=right:{$p(x,y,z)$}] (p) at (3,2){};
\end{tikzpicture}
\end{center}

Next consider the special case where $x,y,z$ lie on a line $l_1$, and suppose $l_1 \in \Theta_1$. Draw the line $l_2 \in \Theta_2$ through $y$ and the line $l_3 \in \Theta_3$ through $x$, and let $y' \in l_2 \cap l_3$ be their point of intersection. Draw the line $l_1'$ through $y'$ parallel to $l_1$, draw the line $l_2'$ through $z$ parallel to $l_2$, and let $z' \in l_1' \cap l_2'$ be their point of intersection. Finally, draw the line $l_3'$ through $z'$ parallel to the line $l_3$, and let $p$ be the intersection point of $l_1$ and $l_3'$.

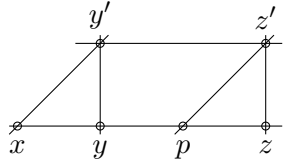
\begin{center}
\begin{tikzpicture}[scale=1.1] 
\draw (-0.1,0) -- (3.2,0);
\draw (1,-0.1) -- (1,1.1);
\draw (3,-0.1) -- (3,1.1);
\draw (0.7,1) -- (3.2,1);
\draw (-0.1,-0.1) -- (1.1,1.1);
\draw (1.9,-0.1) -- (3.1,1.1);
\node[circle, minimum width=3pt, draw, inner sep=0pt, label=below:$x$] (x) at (0,0){};
\node[circle, minimum width=3pt, draw, inner sep=0pt, label=below:$y$] (y) at (1,0){};
\node[circle, minimum width=3pt, draw, inner sep=0pt, label=below:$z$] (z) at (3,0){};
\node[circle, minimum width=3pt, draw, inner sep=0pt, label=above:{$y'$}] (y') at (1,1){};
\node[circle, minimum width=3pt, draw, inner sep=0pt, label=above:{$z'$}] (z') at (3,1){};
\node[circle, minimum width=3pt, draw, inner sep=0pt, label=below:$p$] (p) at (2,0){};
\end{tikzpicture}
\end{center}

We claim that $p = p(x,y,z)$. To see this, note that $x \equiv_{\Theta_3} y'$, so $p(x,y,z) \equiv_{\Theta_3} p(y',y,z)$, and $p(y',y,z) = z'$ by the first case we considered. Thus $p(x,y,z) \equiv_{\Theta_3} z'$, i.e. $p(x,y,z) \in l_3'$, and since $x \equiv_{\Theta_1} y \equiv_{\Theta_1} z$, we have $p(x,y,z) \equiv_{\Theta_1} p(x,x,x) = x$, i.e. $p(x,y,z) \in l_1$. Thus $p(x,y,z) \in l_1 \cap l_3'$, so $p(x,y,z) = p$. (Alternatively, we could have used $p(x,y,z) \equiv_{\Theta_2} p(x,y',z') = p$, by the first case.)

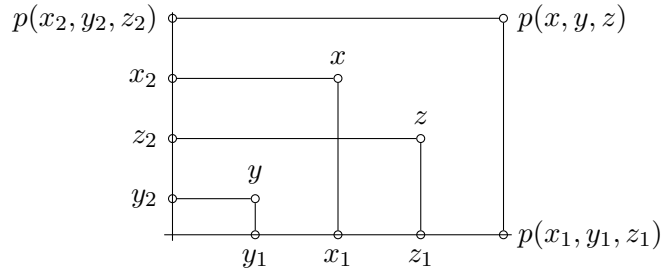
\begin{center}
\begin{tikzpicture}[xscale=1.1,yscale=0.8]
\node[circle, minimum width=3pt, draw, inner sep=0pt, label=above:$x$] (x) at (2,3){};
\node[circle, minimum width=3pt, draw, inner sep=0pt, label=above:$y$] (y) at (1,1){};
\node[circle, minimum width=3pt, draw, inner sep=0pt, label=above:$z$] (z) at (3,2){};
\node[circle, minimum width=3pt, draw, inner sep=0pt, label=below:{$x_1$}] (x1) at (2,0.4){};
\node[circle, minimum width=3pt, draw, inner sep=0pt, label=below:{$y_1$}] (y1) at (1,0.4){};
\node[circle, minimum width=3pt, draw, inner sep=0pt, label=below:{$z_1$}] (z1) at (3,0.4){};
\node[circle, minimum width=3pt, draw, inner sep=0pt, label=left:{$x_2$}] (x2) at (0,3){};
\node[circle, minimum width=3pt, draw, inner sep=0pt, label=left:{$y_2$}] (y2) at (0,1){};
\node[circle, minimum width=3pt, draw, inner sep=0pt, label=left:{$z_2$}] (z2) at (0,2){};
\node[circle, minimum width=3pt, draw, inner sep=0pt, label=right:{$p(x_1,y_1,z_1)$}] (p1) at (4,0.4){};
\node[circle, minimum width=3pt, draw, inner sep=0pt, label=left:{$p(x_2,y_2,z_2)$}] (p2) at (0,4){};
\node[circle, minimum width=3pt, draw, inner sep=0pt, label=right:{$p(x,y,z)$}] (p) at (4,4){};
\draw (x2) -- (x) -- (x1);
\draw (y2) -- (y) -- (y1);
\draw (z2) -- (z) -- (z1);
\draw (p2) -- (p) -- (p1);
\draw (-0.1,0.4) -- (4.1,0.4);
\draw (0,0.3) -- (0,4.1);
\end{tikzpicture}
\end{center}

For the general case, we can pick any lines $l_1 \in \Theta_1, l_2 \in \Theta_2$, set $x_1, y_1, z_1$ to be the projections of $x,y,z$ onto $l_1$ via lines in $\Theta_2$ and define $x_2,y_2,z_2 \in l_2$ similarly, and note that $p(x,y,z) \equiv_{\Theta_2} p(x_1,y_1,z_1)$ and $p(x,y,z) \equiv_{\Theta_1} p(x_2,y_2,z_2)$, and we can construct $p(x_1,y_1,z_1), p(x_2,y_2,z_2)$ using the second case considered.
\end{proof}

\begin{cor} If $p$ is a Mal'cev polymorphism of an S-3-system, then $p(x,y,z) \approx p(z,y,x)$.
\end{cor}
\begin{proof} The term $p(z,y,x)$ is also a Mal'cev polymorphism, so by the Lemma it must be identical to $p(x,y,z)$.
\end{proof}

\begin{cor}\label{S3-malcev-graph} If $p$ is a Mal'cev polymorphism of an S-3-system $(S,\Theta_1,\Theta_2,\Theta_3)$, then the graph $\Gamma_p$ of $p$, considered as a $4$-ary relation on $S$, is primitively positively definable from $\Theta_1,\Theta_2,\Theta_3$.
\end{cor}

Corollary \ref{S3-malcev-graph} can also be interpreted as saying that the map $p : \bS^3 \rightarrow \bS$ is a homomorphism of the \emph{algebraic} structure $\bS$ whose basic operations consist of all polymorphisms of the relational structure $\mathbf{S}$. In particular, $p$ ``commutes with itself'', that is, the two ways of computing $p*p$ on a $3\times 3$ grid of variables (columns first or rows first) agree with each other. We can summarize this fact by saying that the Mal'cev operation $p$ is \emph{central}.

\begin{defn} An $n$-ary term $t$ of an algebraic structure $\bA$ is called \emph{central} if the map $t : \bA^n \rightarrow \bA$ is a homomorphism.
\end{defn}

Now we relate the Mal'cev polymorphism to the coordinate loop $\bL$. First we will show that $\bL$ is associative.

\begin{lem}\label{S3-loop-associative} If $\mathbf{S} = (S,\Theta_1,\Theta_2,\Theta_3)$ is an S-3-system with a Mal'cev polymorphism $p$, and if $\bL$ is a coordinate loop of $\mathbf{S}$, then $\bL$ satisfies
\[
(x_1\cdot y_1 = x_2\cdot y_2) \wedge (x_1\cdot y_3 = x_2\cdot y_4) \wedge (x_3\cdot y_1 = x_4\cdot y_2) \implies (x_3\cdot y_3 = x_4 \cdot y_4).
\]
In particular, $\bL$ is associative, that is, $\bL$ is a group.
\end{lem}
\begin{proof} For those who prefer a purely algebraic proof, this follows from
\[
\begin{bmatrix} x_3\\ y_3\end{bmatrix} = p\left(\begin{bmatrix} x_1\\ y_3\end{bmatrix}, \begin{bmatrix} x_1\\ y_1\end{bmatrix}, \begin{bmatrix} x_3\\ y_1\end{bmatrix}\right) \equiv_{\Theta_3} p\left(\begin{bmatrix} x_2\\ y_4\end{bmatrix}, \begin{bmatrix} x_2\\ y_2\end{bmatrix}, \begin{bmatrix} x_4\\ y_2\end{bmatrix}\right) = \begin{bmatrix} x_4\\ y_4\end{bmatrix}.
\]

To see that this implies the associativity of $\bL$, let $x,y,z$ be any elements of $L$, and plug in $(x_1,x_2,x_3,x_4) = (1, y, x, x\cdot y), (y_1,y_2,y_3,y_4) = (y, 1, y\cdot z, z)$. Then we get
\[
(1\cdot y = y\cdot 1) \wedge (1\cdot (y\cdot z) = y\cdot z) \wedge (x\cdot y = (x\cdot y)\cdot 1) \implies (x\cdot (y\cdot z) = (x\cdot y) \cdot z).
\]

For a geometric way to visualize the proof, note that the stated property of $\bL$ corresponds to the existence of the dashed line in the following picture.

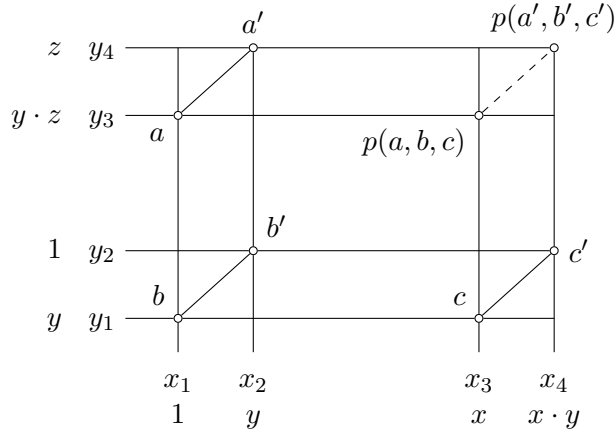
\begin{center}
\begin{tikzpicture}[xscale=1,yscale=0.9]
\node[circle, minimum width=3pt, draw, inner sep=0pt, label=above left:{$b$}] (11) at (1,1){};
\node[circle, minimum width=3pt, draw, inner sep=0pt, label=above right:{$b'$}] (22) at (2,2){};
\node[circle, minimum width=3pt, draw, inner sep=0pt, label=below left:{$a$}] (13) at (1,4){};
\node[circle, minimum width=3pt, draw, inner sep=0pt, label=above:{$a'$}] (24) at (2,5){};
\node[circle, minimum width=3pt, draw, inner sep=0pt, label=above left:{$c$}] (31) at (5,1){};
\node[circle, minimum width=3pt, draw, inner sep=0pt, label=right:{$c'$}] (42) at (6,2){};
\node[circle, minimum width=3pt, draw, inner sep=0pt, label=below left:{$p(a,b,c)$}] (33) at (5,4){};
\node[circle, minimum width=3pt, draw, inner sep=0pt, label=above:{$p(a',b',c')$}] (44) at (6,5){};
\node[align=left, below] (x1) at (1,0.5){$\begin{matrix} x_1\\ 1\end{matrix}$};
\node[align=left, below] (x2) at (2,0.5){$\begin{matrix} x_2\\ y\end{matrix}$};
\node[align=left, below] (x3) at (5,0.5){$\begin{matrix} x_3\\ x\end{matrix}$};
\node[align=left, below] (x4) at (6,0.5){$\begin{matrix} x_4\\ x\cdot y\end{matrix}$};
\node[align=left, left] (y1) at (0.3,1){$\begin{matrix} y & y_1\end{matrix}$};
\node[align=left, left] (y2) at (0.3,2){$\begin{matrix} 1 & y_2\end{matrix}$};
\node[align=left, left] (y3) at (0.3,4){$\begin{matrix} y\cdot z & y_3\end{matrix}$};
\node[align=left, left] (y4) at (0.3,5){$\begin{matrix} z & y_4\end{matrix}$};
\draw (y1) -- (11) -- (31) -- (6,1);
\draw (y2) -- (22) -- (42);
\draw (y3) -- (13) -- (33) -- (6,4);
\draw (y4) -- (24) -- (44);
\draw (x1) -- (11) -- (13) -- (1,5);
\draw (x2) -- (22) -- (24);
\draw (x3) -- (31) -- (33) -- (5,5);
\draw (x4) -- (42) -- (44);
\draw (11) -- (22);
\draw (13) -- (24);
\draw (31) -- (42);
\draw[dashed] (33) -- (44);
\end{tikzpicture}
\end{center}

If we set $a = (x_1,y_3)$, etc. as in the picture, then the existence of the dashed line follows from the fact that $p$ preserves the congruence $\Theta_3$ and the fact that $p(a,b,c)$ completes the parallelogram through $a,b,c$ and $p(a',b',c')$ completes the parallelogram through $a',b',c'$.
\end{proof}

\begin{lem}\label{S3-malcev-loop} If $\mathbf{S} = (S,\Theta_1,\Theta_2,\Theta_3)$ is an S-3-system with a Mal'cev polymorphism $p$, and if $\bL$ is a coordinate group of $\mathbf{S}$, then for $x = (x_1,x_2), y = (y_1,y_2), z = (z_1,z_2) \in S$, $p(x,y,z)$ is given by
\[
p\left(\begin{bmatrix} x_1\\ x_2\end{bmatrix}, \begin{bmatrix} y_1\\ y_2\end{bmatrix}, \begin{bmatrix} z_1\\ z_2\end{bmatrix}\right) = \begin{bmatrix} x_1\cdot y_1^{-1}\cdot z_1\\ x_2\cdot y_2^{-1}\cdot z_2\end{bmatrix}.
\]
\end{lem}
\begin{proof} It's enough to consider the case where $x,y,z$ are along the line $l_1 \in \Theta_1$ with $\Theta_1$-coordinate $1$. Consider the diagram
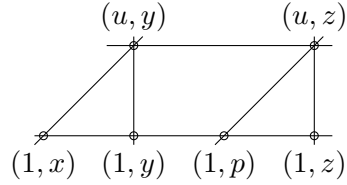
\begin{center}
\begin{tikzpicture}[scale=1.2] 
\draw (-0.1,0) -- (3.2,0);
\draw (1,-0.1) -- (1,1.1);
\draw (3,-0.1) -- (3,1.1);
\draw (0.7,1) -- (3.2,1);
\draw (-0.1,-0.1) -- (1.1,1.1);
\draw (1.9,-0.1) -- (3.1,1.1);
\node[circle, minimum width=3pt, draw, inner sep=0pt, label=below:{$(1,x)$}] (x) at (0,0){};
\node[circle, minimum width=3pt, draw, inner sep=0pt, label=below:{$(1,y)$}] (y) at (1,0){};
\node[circle, minimum width=3pt, draw, inner sep=0pt, label=below:{$(1,z)$}] (z) at (3,0){};
\node[circle, minimum width=3pt, draw, inner sep=0pt, label=above:{$(u,y)$}] (y') at (1,1){};
\node[circle, minimum width=3pt, draw, inner sep=0pt, label=above:{$(u,z)$}] (z') at (3,1){};
\node[circle, minimum width=3pt, draw, inner sep=0pt, label=below:{$(1,p)$}] (p) at (2,0){};
\end{tikzpicture}
\end{center}
which we used to construct $p(x,y,z)$. Then from $(1,x) \equiv_{\Theta_3} (u,y)$ we have $1\cdot x = u\cdot y$, and from $(1,p) \equiv_{\Theta_3} (u,z)$ we have $1 \cdot p = u \cdot z$. Solving for $u$ we get $u = xy^{-1}$, and solving for $p$ we get $p = xy^{-1}z$.
\end{proof}

\begin{cor} If $\mathbf{S} = (S,\Theta_1,\Theta_2,\Theta_3)$ is an S-3-system with a Mal'cev polymorphism $p$, and if $\bL$ is a coordinate group of $\mathbf{S}$, then $\bL$ is commutative.
\end{cor}
\begin{proof} From $p(x,y,z) \approx p(z,y,x)$ we get $xy^{-1}z \approx zy^{-1}x$ in $\bL$, and plugging in $y = 1$ gives $xz \approx zx$, so $\bL$ is commutative.
\end{proof}

Putting all of this together, we have the main result of this section.

\begin{thm}\label{abelian-malcev} Any abelian Mal'cev algebra $\bA$ is affine.
\end{thm}
\begin{proof} By Theorem \ref{abelian-M3} and its corollary, $\mathbf{S} = (\bA^2, \ker \pi_1, \ker \pi_2, \Theta)$ is an S-3-system with Mal'cev polymorphism $p$, where $p$ is the Mal'cev term operation of $\bA$ and $\Theta$ is a congruence on $\bA^2$ with the diagonal as a congruence class. By Lemma \ref{S3-coords}, there is a loop structure $\bL$ on the underlying set of $\bA$ which describes $\mathbf{S}$. By Lemma \ref{S3-loop-associative}, Lemma \ref{S3-malcev-loop}, and its corollary, $\bL$ is an abelian group and $p$ is given by $p(x,y,z) = x-y+z$ (writing the abelian group operation additively).

By Corollary \ref{S3-malcev-graph}, the relation $x-y+z=p$ is primitively positively definable from $\ker \pi_1, \ker \pi_2, \Theta$, so the relation $x+z = y+p$ is preserved by all operations of $\bA$, that is, $\bA$ is quasi-affine. Since $\bA$ was assumed to be Mal'cev, this means that $\bA$ is affine.
\end{proof}

We have proved the hardest part of the Fundamental Theorem of Abelian Algebras. For the sake of completeness, we include the rest of it.

\begin{thm}[Fundamental Theorem of Abelian Algebras] For an algebraic structure $\bA$, the following are equivalent:
\begin{itemize}
\item[(1)] $\bA$ is affine,

\item[(2)] $\bA$ is abelian and has a Mal'cev polynomial,

\item[(3)] $\bA$ has a central Mal'cev polynomial.
\end{itemize}
\end{thm}
\begin{proof} That (1) implies (2) and (3) is clear. For (2) implies (1) and (3), note that any polynomial of $\bA$ preserves every congruence of $\bA^2$, so the polynomial clone of $\bA$ is also abelian and we may apply the previous theorem. For (3) $\implies$ (1), we just need to show that any Mal'cev operation $p$ which commutes with itself comes from an abelian group, since then the fact that $p(x,y,z) = x-y+z$ is central will imply that $\bA$ is quasi-affine.

So suppose that $p$ is a Mal'cev operation which commutes with itself, and pick any element to call $0$ in $\bA$. We define addition and negation on $\bA$ by
\[
x+y \coloneqq p(x,0,y), \;\;\; -x \coloneqq p(0,x,0).
\]
That $0$ is an identity element for $+$ follows from the Mal'cev identities $p(x,0,0) = p(0,0,x) = x$.

To see that $+$ is associative, we evaluate the expression
\[
p*p\left(\begin{bmatrix} x & 0 & y\\ 0 & 0 & 0\\ 0 & 0 & z\end{bmatrix}\right)
\]
in two ways: evaluating it by rows first, we get $(x+y)+z$, and evaluating it by columns first, we get $x+(y+z)$.

To see that $-$ computes the inverse, we evaluate the expression
\[
p*p\left(\begin{bmatrix} x & 0 & 0\\ 0 & 0 & x\\ 0 & 0 & 0\end{bmatrix}\right)
\]
in two ways: by rows we get $p(x,x,0) = 0$, and by columns we get $x+(-x)$. A similar argument shows that $(-x)+x = 0$.

For commutativity of $+$, we evaluate the expression
\[
p*p\left(\begin{bmatrix} y & 0 & x\\ y & y & x\\ x & y & y\end{bmatrix}\right)
\]
in two ways: by rows we get $p(y+x,x,x) = y+x$, and by columns we get $p(x,0,y) = x+y$.

Finally, to express $p$ in terms of the group operations $+,-$, we evaluate the expression
\[
p*p\left(\begin{bmatrix} x & y & z\\ 0 & y & 0\\ -y & 0 & 0\end{bmatrix}\right)
\]
in two ways: by rows we get $p(p(x,y,z),-y,-y) = p(x,y,z)$, and by columns we get $p(x-y,0,z) = x-y+z$.
\end{proof}

The method of visualizing algebraic arguments via the geometry of equivalence classes was extended to congruence modular varieties by Gumm in his book ``Geometrical methods in congruence modular algebras'' \cite{gumm-geometric}, where he used it to show that any abelian algebra in a congruence modular variety is affine. This was extended further by Hobby and McKenzie \cite{hobby-mckenzie}, who used tame congruence theory to show that any finite abelian algebra in a Taylor variety is affine (in the infinite case, Kearnes and Szendrei \cite{kearnes-taylor-affine} show that any abelian Taylor algebra is quasi-affine - the example $(\RR, \frac{x+y}{2})$ shows that an additional assumption is needed for it to be affine). Later we will go over a simpler proof of the fact that finite abelian Taylor algebras are affine, from \cite{pointing-no-absorption}.

\begin{rem} If we leave the context of Taylor varieties, we can no longer expect abelian algebras to be affine, since they could fail to have any interesting operations at all. But we can still ask whether abelian algebras are quasi-affine. The following problem is open.
\end{rem}

\begin{prob} Under what conditions are abelian algebras quasi-affine? Is it true that every idempotent abelian algebra is quasi-affine?
\end{prob}

It is known that if we drop idempotence, then some extra condition is needed: Quackenbush \cite{quasi-affine-quackenbush} gives an example of an infinite, non-idempotent algebra which is abelian but not quasi-affine. Quackenbush's example is a slight modification of the completely free algebra on $8$ elements with a single binary operation, where the modification is that $x_1 \cdot x_2 = x_5 \cdot x_6, x_3 \cdot x_4 = x_7 \cdot x_8, x_1\cdot x_4 = x_5 \cdot x_8$, but $x_3 \cdot x_2 \ne x_7 \cdot x_6$. Another example with just five elements is given in Example \ref{ex-strongly-abelian-not-quasiaffine}.

Kearnes \cite{kearnes-simple} has shown that any \emph{simple} idempotent abelian algebra is quasi-affine - in fact, he shows that any simple idempotent algebra which has a skew congruence (that is, a congruence on some power $\bA^n$ which is not the kernel of some projection) either has a strongly absorbing element (that is, an element $a$ such that every term $t$ which depends on its first variable has $t(a,...) = a$) or is a subalgebra of a simple reduct of a module.

There are a few other contexts in which it is known that abelian implies quasi-affine. In \cite{kearnes-quasi-affine}, Kearnes shows that any abelian algebra with a central binary polynomial which is cancellative is quasi-affine, and in \cite{stronkowski-embedding} this is extended to the result that any abelian algebra with a commutative cancellative polynomial is quasi-affine. In \cite{affine-quandles}, it is shown that abelian quandles are quasi-affine.

\subsection{Commutators}\label{ss-commutators}

In this subsection we define an extension of the commutator from group theory to a commutator on congruences of general algebraic structures. The purpose of the commutator is to detect situations where the operations of an algebraic structure behave linearly. The theory of the commutator works best in congruence modular varieties, but it still has some use in general Taylor varieties, although slight differences in the technical details of the definition become important outside the world of congruence modular varieties. The commutator we will be discussing is called the \emph{term condition} commutator.

\begin{defn} If $\alpha, \beta, \delta \in \Con(\bA)$, we say that $\alpha$ \emph{centralizes} $\beta$ modulo $\delta$, written $C(\alpha,\beta;\delta)$ (or $C(\alpha,\beta)$ if $\delta = 0_\bA$), if for every $n+1$-ary term $t \in \Clo_{n+1}(\bA)$, for any $(u,v) \in \alpha$, and for any $(a_1,b_1), ..., (a_n,b_n) \in \beta$, we have
\[
t(u,a_1, ..., a_n) \equiv_\delta t(u,b_1, ..., b_n) \iff t(v,a_1,...,a_n) \equiv_\delta t(v,b_1,...,b_n).
\]
The smallest $\delta$ which satisfies $C(\alpha,\beta;\delta)$ is called the \emph{commutator} of $\alpha,\beta$, and is written as $[\alpha,\beta]$. If $\theta \le \alpha,\beta$, then we also define the \emph{relative commutator} $[\alpha,\beta]_\theta$ to be the least $\delta \ge \theta$ which satisfies $C(\alpha,\beta;\delta)$.
\end{defn}

As with the criterion for abelianness, the term condition implies a seemingly stronger version where more variables change at once.

\begin{prop} If $\alpha$ centralizes $\beta$ modulo $\delta$, then for every $m+n$-ary term $t \in \Clo_{m+n}(\bA)$, for any $(u_1,v_1), ..., (u_m,v_m) \in \alpha$, and for any $(a_1,b_1), ..., (a_n,b_n) \in \beta$, we have
\[
t(u_1,...,u_m,a_1, ..., a_n) \equiv_\delta t(u_1,...,u_m,b_1, ..., b_n) \iff t(v_1,...v_m,a_1,...,a_n) \equiv_\delta t(v_1,...,v_m,b_1,...,b_n).
\]
\end{prop}

Before we go on, let's check that this matches the usual commutator from group theory.

\begin{prop} If $\bM,\bN$ are normal subgroups of a group $\bG$, $[\bM,\bN]$ is the (normal) subgroup generated by commutators $[m,n] = mnm^{-1}n^{-1}$ for $m \in \bM, n \in \bN$, and $\theta_\bM, \theta_\bN, \theta_{[\bM,\bN]}$ are the associated congruences, then $\theta_{[\bM,\bN]} = [\theta_\bM,\theta_\bN]$.
\end{prop}
\begin{proof} We will show that $\theta_\bM$ centralizes $\theta_\bN$ iff every element of $\bM$ commutes with every element of $\bN$ - this will finish the proof, since $[\bM,\bN]$ is the smallest normal subgroup $\mathbb{K}$ of $\bG$ such that every element of $\bM/\mathbb{K}$ commutes with every element of $\bN/\mathbb{K}$ in $\bG/\mathbb{K}$.

First suppose that $\theta_\bM$ centralizes $\theta_\bN$. Let $t$ be the binary term $t(x,y) = yxy^{-1}$, then for any $m\in\bM, n\in\bN$, by the term condition applied to $(1,m) \in \theta_\bM, (1,n) \in \theta_\bN$, we have
\[
1 = nn^{-1} \iff m = nmn^{-1},
\]
so $m$ and $n$ commute.

Now suppose that every element of $\bM$ commutes with every element of $\bN$, and consider an arbitrary $n+1$-ary term $t \in \Clo_{n+1}(\bG)$ and any $(u,v) \in \theta_\bM, (a_1,b_1), ..., (a_n,b_n) \in \theta_\bN$ with
\[
t(u,a_1,...,a_n) = t(u,b_1,...,b_n).
\]
Thinking of $t(ux,a_1y_1,....,a_ny_n)$ as a function of $x,y_1,...,y_n$ with parameters $u,a_1,...,a_n$, we may rearrange it into the form
\[
t(ux,a_1y_1,....,a_ny_n) = t'(x,y_1,...,y_n)t(u,a_1,...,a_n)
\]
for some $t'$ in the clone generated by the group operations together with the unary conjugation operations $\phi_c: x \mapsto cxc^{-1}$, so we may rewrite our assumption as
\[
t'(1,1,...,1) = t'(1,a_1^{-1}b_1,...,a_n^{-1}b_n).
\]
To show that
\[
t(v,a_1,...,a_n) = t(v,b_1,...,b_n),
\]
we just need to show that
\[
t'(u^{-1}v,1,...,1) = t'(u^{-1}v,a_1^{-1}b_1,...,a_n^{-1}b_n),
\]
which follows from the assumed equality together with the fact that for each $c,d \in \bG$ and each $i$, $\phi_c(u^{-1}v) \in \bM$ commutes with $\phi_d(a_i^{-1}b_i) \in \bN$.
\end{proof}

\begin{ex} In the case of rings, the term condition commutator applied to a pair of ideals $I,J$ gives $[I,J] = IJ + JI$. Note that this is a bit different from what we might have expected (it has nothing to do with the Lie bracket), but it makes more sense when we remember that we only consider a ring to be abelian if it is a zero ring.
\end{ex}

\begin{ex} In a majority algebra, the commutator is given by $[\alpha,\beta] = \alpha \wedge \beta$. To see this, suppose $(a,b) \in \alpha \wedge \beta$, and apply the term condition to the majority operation to see that
\[
m(\boxed{a},a,a) = m(\boxed{a},a,b) \implies m(\boxed{b},a,a)\ [\alpha,\beta]\ m(\boxed{b},a,b),
\]
so $(a,b) \in [\alpha,\beta]$. A similar argument shows that the commutator is given by intersection in any variety with a near-unanimity term.
\end{ex}

\begin{ex}\label{semi-sd-meet} In a semilattice, the commutator is given by $[\alpha,\beta] = \alpha\wedge\beta$. Let $s$ be the semilattice operation, and let $s_3$ be the term given by $s_3(x,y,z) = s(x,s(y,z))$. Then for $(a,b) \in \alpha\wedge\beta$, we have
\[
s_3(\boxed{a},a,b) = s_3(\boxed{a},b,b) \implies s_3(\boxed{b},a,b)\ [\alpha,\beta]\ s_3(\boxed{b},b,b),
\]
so $s(a,b)\ [\alpha,\beta]\ b$, and similarly $s(a,b)\ [\alpha,\beta]\ a$, so $(a,b) \in [\alpha,\beta]$.
\end{ex}

Sometimes it is helpful to visualize the term condition via $2\times 2$ matrices.

\begin{defn}\label{commutator-matrix} For $\alpha, \beta \in \Con(\bA)$, we define the algebra $\bM(\alpha,\beta) \le \bA^{2\times 2}$ to be the subalgebra of $2\times 2$ matrices which is generated by the matrices of the form
\[
\begin{bmatrix} u & u\\ v & v\end{bmatrix} \text{ with } (u,v) \in \alpha, \;\;\; \begin{bmatrix} a & b\\ a & b\end{bmatrix} \text{ with } (a,b) \in \beta.
\]
\end{defn}

\begin{prop} If $\alpha, \beta, \delta \in \Con(\bA)$, then $\alpha$ centralizes $\beta$ modulo $\delta$ iff for all
\[
\begin{bmatrix} a & b\\ c & d\end{bmatrix} \in \bM(\alpha,\beta)
\]
we have
\[
a \equiv_\delta b \iff c \equiv_\delta d.
\]
\end{prop}

The usual picture which is drawn to represent the term condition for $C(\alpha,\beta;\delta)$ is this:
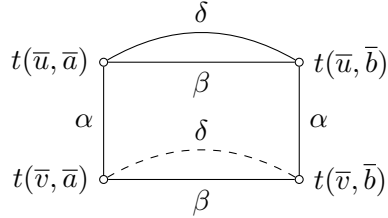
\begin{center}
\begin{tikzpicture}[scale=1.3]
\node[circle, minimum width=3pt, draw, inner sep=0pt, label=left:{$t(\overline{u},\overline{a})$}] (ua) at (0,1.2){};
\node[circle, minimum width=3pt, draw, inner sep=0pt, label=right:{$t(\overline{u},\overline{b})$}] (ub) at (2.0,1.2){};
\node[circle, minimum width=3pt, draw, inner sep=0pt, label=left:{$t(\overline{v},\overline{a})$}] (va) at (0,0){};
\node[circle, minimum width=3pt, draw, inner sep=0pt, label=right:{$t(\overline{v},\overline{b})$}] (vb) at (2.0,0){};
\draw (ua) to ["$\beta$"'] (ub) to ["$\alpha$"] (vb) to ["$\beta$"] (va) to ["$\alpha$"] (ua) [bend left] to ["$\delta$"] (ub);
\draw [bend right, dashed] (vb) to ["$\delta$"'] (va);
\end{tikzpicture}
\end{center}
where the positioning of the four corners matches with the way we have laid out the $2\times 2$ matrices in $\bM(\alpha,\beta)$. A mnemonic for remembering where the $\delta$ edges go is that in the term condition $C(\alpha,\beta;\delta)$, ``$\delta$ is next to $\beta$''.

We now list a few elementary properties of the commutator which hold in general, which are given as exercises in Hobby and McKenzie's book \cite{hobby-mckenzie}.

\begin{prop}\label{gen-commutator} For $\alpha,\beta,\delta \in \Con(\bA)$, we have
\begin{itemize}
\item[(a)] if $C(\alpha,\beta;\delta_i)$ for $i \in I$, then $C(\alpha,\beta;\bigwedge_{i\in I} \delta_i)$, so $[\alpha,\beta]$ and $[\alpha,\beta]_\theta$ are well-defined,

\item[(b)] if $(\alpha \vee (\beta \wedge \delta)) \wedge \beta \le \delta$ then $C(\alpha, \beta; \delta)$ holds, so $[\alpha,\beta] \le \alpha \wedge \beta$,

\item[(c)] if $\alpha' \le \alpha, \beta' \le \beta$, then $C(\alpha,\beta;\delta) \implies C(\alpha',\beta';\delta)$, so $[\alpha',\beta'] \le [\alpha, \beta]$,

\item[(d)] for any $\gamma$ we have $C(\alpha,\beta;\delta) \implies C(\alpha\wedge \gamma,\beta;\delta \wedge \gamma)$,

\item[(e)] if $C(\alpha_i,\beta;\delta)$ holds for all $i\in I$ then $C(\bigvee_{i \in I}\alpha_i,\beta;\delta)$ holds,

\item[(f)] if $\theta \le \alpha, \beta, \delta$ then $C(\alpha,\beta;\delta)$ holds iff $C(\alpha/\theta,\beta/\theta;\delta/\theta)$ holds in $\bA/\theta$, so $[\alpha/\theta,\beta/\theta] = [\alpha,\beta]_\theta/\theta$,

\item[(g)] if $\bB \le \bA$ then $C(\alpha,\beta;\delta) \implies C(\alpha|_\bB, \beta|_\bB; \delta|_\bB)$, so $[\alpha|_\bB,\beta|_\bB] \le [\alpha,\beta]|_\bB$,

\item[(h)] if $[\alpha,\alpha] = 0_\bA$, then any congruence class of $\alpha$ which is also a subalgebra of $\bA$ is an abelian subalgebra.
\end{itemize}
\end{prop}
\begin{proof} Parts (a), (c), (d), (f), (g), (h) follow immediately from the definitions. For (b), note that for any $\begin{bmatrix} a & b\\ c & d\end{bmatrix} \in \bM(\alpha,\beta)$ with $a \equiv_\delta b$, we have $c \equiv_\alpha a \equiv_{\beta \wedge \delta} b \equiv_\alpha d$ and $c \equiv_\beta d$, so $(c,d) \in (\alpha \circ (\beta \wedge \delta) \circ \alpha) \wedge \beta$, which is a subset of $\delta$ by assumption.

For (e), we string together several instances of the term condition: if $(u,v) \in \bigvee_i \alpha_i, (a_i,b_i) \in \beta$, and $t(u,\overline{a}) \equiv_\delta t(u,\overline{b})$, then if we let $u = u_0, u_1, ..., u_n = v$ be a sequence of elements of $\bA$ with $(u_i, u_{i+1}) \in \alpha_{j_i}$ for some $j_i \in I$, then by the term condition $C(\alpha_{j_i},\beta;\delta)$ we have
\[
t(u_i,\overline{a})\equiv_\delta t(u_i,\overline{b}) \implies t(u_{i+1},\overline{a})\equiv_\delta t(u_{i+1},\overline{b}),
\]
so by inducting on $i$ we get $t(v,\overline{a}) \equiv_\delta t(v,\overline{b})$.
\end{proof}

\begin{cor} If an idempotent algebra $\bA$ has any congruences $\alpha, \beta \in \Con(\bA)$ with $[\alpha,\beta] \ne \alpha \wedge \beta$, then some subalgebra of some quotient of $\bA$ is a nontrivial abelian algebra.
\end{cor}
\begin{proof} Let $\delta = \alpha \wedge \beta$, then from $\delta \le \alpha, \beta$ we have $[\delta,\delta] \le [\alpha,\beta] < \alpha\wedge \beta = \delta$. Thus $\delta' = \delta/[\delta,\delta]$ is a nontrivial congruence on $\bA/[\delta,\delta]$ with $[\delta',\delta'] = [\delta,\delta]/[\delta,\delta] = 0_{\bA/[\delta,\delta]}$, so there is some nontrivial congruence class $\bB$ of $\delta'$ and $\bB$ is an abelian subalgebra of $\bA/[\delta,\delta]$.
\end{proof}

\begin{prop}\label{sd-meet-commutator} If $[\alpha,\beta] = \alpha \wedge \beta$ for all $\alpha, \beta \in \Con(\bA)$, then $\Con(\bA)$ satisfies the meet-semidistributive law:
\[
\alpha \wedge \beta = \alpha\wedge \gamma \implies \alpha \wedge (\beta \vee \gamma) = \alpha \wedge \beta.
\]
\end{prop}
\begin{proof} If $\alpha,\beta,\gamma \in \Con(\bA)$ satisfy $\alpha \wedge \beta = \alpha \wedge \gamma$, then $C(\beta,\alpha;\alpha \wedge \beta)$ and $C(\gamma,\alpha;\alpha\wedge \beta)$ hold, so $C(\beta\vee \gamma,\alpha;\alpha\wedge \beta)$ holds, so $\alpha \wedge (\beta\vee \gamma) = [\beta\vee \gamma,\alpha] \le \alpha\wedge \beta$.
\end{proof}

\begin{defn} An algebra $\bA$ is \emph{congruence meet-semidistributive}, written SD($\wedge$) for short, if for all $\alpha,\beta,\gamma \in \Con(\bA)$ with $\alpha \wedge \beta = \alpha\wedge \gamma$, we have $\alpha \wedge (\beta \vee \gamma) = \alpha \wedge \beta$. A variety $\cV$ is SD($\wedge$) if every algebra $\bA \in \cV$ is SD($\wedge$).
\end{defn}

The next corollary is the key to classifying CSPs which do not have the ``ability to count'' - as we will see later, a finite idempotent algebra generates an SD($\wedge$) variety if and only if the associated CSP has bounded width.

\begin{cor}\label{cor-sd-meet} If an idempotent variety does not contain any nontrivial abelian algebras, then it is congruence meet-semidistributive. Conversely, a congruence meet-semidistributive variety does not contain any nontrivial affine algebra.
\end{cor}
\begin{proof} For the converse statement, note that if $\bA$ is affine, then $\Con(\bA^2)$ contains a copy of the diamond lattice $\cM_3$, and $\cM_3$ doesn't satisfy the meet-semidistributive law.
\end{proof}

Now we consider some definitions which are useful in the case where the commutator is not trivial (i.e., not given by $[\alpha,\beta] = \alpha \wedge \beta$).

\begin{defn}\label{defn-centralizer} Suppose that $\alpha \le \beta \in \Con(\bA)$. We say that $\beta$ is \emph{abelian} over $\alpha$ if the term condition $C(\beta,\beta;\alpha)$ holds. We say that $\beta$ is \emph{solvable} over $\alpha$ if there is a chain of congruences $\alpha = \alpha_0 \le \cdots \le \alpha_n = \beta$ such that $\alpha_{i+1}$ is abelian over $\alpha_i$ for each $i$.

A congruence $\alpha$ is called abelian if it is abelian over $0_\bA$ (equivalently $[\alpha,\alpha] = 0_\bA$), and similarly $\alpha$ is called solvable if $\alpha$ is solvable over $0_\bA$. An algebra $\bA$ is called solvable if $1_\bA$ is solvable.

The \emph{center} of an algebra $\bA$ is defined to be the largest $\zeta$ such that $C(\zeta,1_\bA)$ holds (equivalently, the largest $\zeta$ with $[\zeta, 1_\bA] = 0_\bA$). For $\beta$ a congruence, we define the \emph{centralizer} of $\beta$, written $(0:\beta)$, to be the largest congruence $\alpha$ such that $[\alpha,\beta] = 0$, and more generally for any $\delta$ we define the \emph{relative centralizer} $(\delta:\beta)$ to be the largest $\alpha$ such that $C(\alpha,\beta;\delta)$ holds.
\end{defn}

\begin{prop} For congruences on $\bA$, we have the following:
\begin{itemize}
\item[(a)] for any $\beta,\delta$ there exists a largest $\alpha$ such that $C(\alpha,\beta;\delta)$ holds, so $(\delta:\beta)$ (and, in particular, the center of $\bA$) is well-defined,

\item[(b)] if $\gamma$ is solvable over $\beta$ and $\beta$ is solvable over $\alpha$, then $\gamma$ is solvable over $\alpha$,

\item[(c)] if $\beta$ is solvable (abelian) over $\alpha$, then $\beta \wedge \gamma$ is solvable (abelian) over $\alpha \wedge \gamma$ for any $\gamma$,

\item[(d)] if $\theta \le \alpha \le \beta$, then $\beta$ is solvable (abelian) over $\alpha$ iff $\beta/\theta$ is solvable (abelian) over $\alpha/\theta$,

\item[(e)] $\bA/\theta$ is solvable (abelian) iff $1_\bA$ is solvable (abelian) over $\theta$.
\end{itemize}
\end{prop}
\begin{proof} Part (a) follows from Proposition \ref{gen-commutator}(e), part (b) is obvious, part(c) follows from Proposition \ref{gen-commutator}(d), part (d) follows from Proposition \ref{gen-commutator}(f), and part (e) is part (d) specialized to the case $\beta = 1_\bA, \alpha = \theta$.
\end{proof}

If our algebra is finite, then solvability has a surprisingly simple alternative characterization based on tame congruence theory, which is described in Appendix \ref{a-sec-snags}. To take the general theory further, we need to make an additional assumption on our variety, such as congruence modularity. The interested reader can find the (surprisingly deep) theory of commutators in congruence modular varieties in Appendix \ref{a-commutator}.

A weaker assumption which is still good enough to prove most of the basic properties of commutators is the existence of a ternary term known as a \emph{difference term}, generalizing the Gumm difference term found in congruence modular varieties, which acts like a Mal'cev term on abelian algebras.

\begin{defn} A ternary term $p$ is called a \emph{difference term} for a variety, if it satisfies the identity $p(y,y,x) \approx x$, and for every $(x,y) \in \theta$ for $\theta$ a congruence, we always have $p(x,y,y) \equiv_{[\theta,\theta]} x$.
\end{defn}

\begin{ex} Any SD($\wedge$) variety has a difference term: just take $p(x,y,z) = z$. That this works relies on the fact that $[\alpha,\beta] = \alpha\wedge\beta$ in SD($\wedge$) varieties, which we haven't proved - this can be found in \cite{kearnes-taylor-affine}.
\end{ex}

One property of a difference term is that it forces several alternative commutators to match with the term condition commutator, and one of these commutators is clearly symmetric.

\begin{defn} For any $n \ge 1$, we define the $n$-cycle commutator $[\alpha,\beta]_n$ to be the least congruence $\delta$ such that for any cycle of $n$ matrices
\[
\begin{bmatrix} a_1 & b_1\\ c_1 & d_1\end{bmatrix}, \begin{bmatrix} a_2 & b_2\\ c_2 & d_2\end{bmatrix}, ..., \begin{bmatrix} a_n & b_n\\ c_n & d_n\end{bmatrix} \in \bM(\alpha,\beta)
\]
such that $b_i \equiv_\delta a_{i+1}$ for all $i < n$, $b_n \equiv_\delta a_1$, and $d_i \equiv_\delta c_{i+1}$ for all $i < n$, we have additionally that $d_n \equiv_\delta c_1$.
\end{defn}

If $\bA$ is affine, then it is easy to check that $[1_\bA,1_\bA]_n = 0_\bA$ for every $n$. Note that for $n = 1$, we have $[\alpha,\beta]_1 = [\alpha,\beta]$. Additionally, since we can take the $n$th matrix in the cycle to have a pair of equal columns, we have $[\alpha,\beta]_i \le [\alpha,\beta]_{i+1}$ for all $i$.

Quackenbush's famous example of an abelian algebra which is not quasi-affine from \cite{quasi-affine-quackenbush} is an example of an algebra where $[1_\bA,1_\bA]_1 = 0_\bA$ but $[1_\bA,1_\bA]_2 \ne 0_\bA$.

For $n = 2$, the $2$-cycle commutator is clearly symmetric: $[\alpha,\beta]_2 = [\beta,\alpha]_2$. Since it is defined via two matrices in $\bM(\alpha,\beta)$, and since each matrix comes from some term, the commutator $[\alpha,\beta]_2$ is also called the \emph{two term commutator}. The two term condition is illustrated in the following diagram.
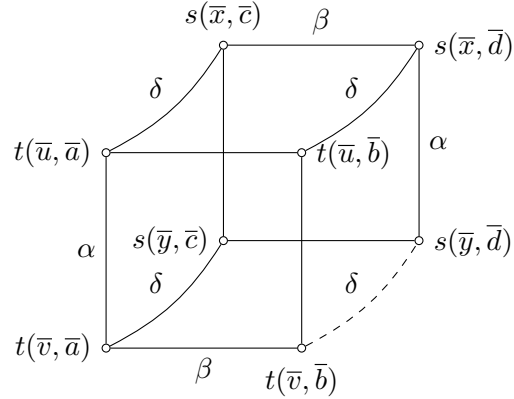
\begin{center}
\begin{tikzpicture}[scale=1.3]
\node[circle, minimum width=3pt, draw, inner sep=0pt, label=left:{$t(\overline{u},\overline{a})$}] (ua) at (0,2.0){};
\node[circle, minimum width=3pt, draw, inner sep=0pt, label=right:{$t(\overline{u},\overline{b})$}] (ub) at (2.0,2.0){};
\node[circle, minimum width=3pt, draw, inner sep=0pt, label=left:{$t(\overline{v},\overline{a})$}] (va) at (0,0){};
\node[circle, minimum width=3pt, draw, inner sep=0pt, label=below:{$t(\overline{v},\overline{b})$}] (vb) at (2.0,0){};
\node[circle, minimum width=3pt, draw, inner sep=0pt, label=above:{$s(\overline{x},\overline{c})$}] (xc) at (1.2,3.1){};
\node[circle, minimum width=3pt, draw, inner sep=0pt, label=right:{$s(\overline{x},\overline{d})$}] (xd) at (3.2,3.1){};
\node[circle, minimum width=3pt, draw, inner sep=0pt, label=left:{$s(\overline{y},\overline{c})$}] (yc) at (1.2,1.1){};
\node[circle, minimum width=3pt, draw, inner sep=0pt, label=right:{$s(\overline{y},\overline{d})$}] (yd) at (3.2,1.1){};
\draw (ua) to (ub) to (vb) to ["$\beta$"] (va) to ["$\alpha$"] (ua);
\draw (xc) to ["$\beta$"] (xd) to ["$\alpha$"] (yd) to (yc) to (xc);
\draw[bend right=15] (ua) to ["$\delta$"] (xc) (ub) to ["$\delta$"] (xd) (va) to ["$\delta$"] (yc);
\draw[bend right=15, dashed] (vb) to ["$\delta$"] (yd);
\end{tikzpicture}
\end{center}

If we have a difference term, then all of the $n$-cycle commutators turn out to be equal.

\begin{thm}[Lipparini \cite{lipparini-difference}] In a variety with a difference term, we have $[\alpha,\beta]_n = [\alpha,\beta]$ for all $n$. In particular, we have $[\alpha,\beta] = [\beta,\alpha]$.
\end{thm}
\begin{proof} Suppose that $p$ is a difference term. We will show that $[\alpha,\beta]$ satisfies the $n$-cycle term condition by induction on $n$. Suppose that matrices $\begin{bmatrix} a_i & b_i\\ c_i & d_i\end{bmatrix} \in \bM(\alpha,\beta)$ for $i \le n$ are as in the definition of the $n$-cycle condition for $\delta = [\alpha,\beta]$. Applying the difference term, we have
\[
p\left(\begin{bmatrix} a_i & b_i\\ c_i & d_i\end{bmatrix}, \begin{bmatrix} b_1 & b_1\\ d_1 & d_1\end{bmatrix}, \begin{bmatrix} a_1 & a_1\\ c_1 & c_1\end{bmatrix}\right) = \begin{bmatrix} p(a_i,b_1,a_1) & p(b_i, b_1, a_1)\\ p(c_i, d_1, c_1) & p(d_i, d_1, c_1)\end{bmatrix} \in \bM(\alpha,\beta)
\]
for $2 \le i \le n-1$, and
\[
p\left(\begin{bmatrix} a_n & b_n\\ c_n & d_n\end{bmatrix}, \begin{bmatrix} b_1 & a_1\\ d_1 & c_1\end{bmatrix}, \begin{bmatrix} a_1 & a_1\\ c_1 & c_1\end{bmatrix}\right) = \begin{bmatrix} p(a_n,b_1,a_1) & p(b_n, a_1, a_1)\\ p(c_n, d_1, c_1) & p(d_n, c_1, c_1)\end{bmatrix} \in \bM(\alpha,\beta).
\]
The reader can check that these form a system of matrices as in the definition of the $n-1$-cycle condition for $\delta = [\alpha,\beta]$, so by the inductive hypothesis we have
\[
p(c_2,d_1,c_1) \equiv_{[\alpha,\beta]} p(d_n,c_1,c_1).
\]
From $c_2 \equiv_{[\alpha,\beta]} d_1$ and the fact that $p$ is a difference term, the left hand side is congruenct to $c_1$ modulo $[\alpha,\beta]$. From the fact that $c_1 \equiv_{\beta} d_n$ and $(c_1,d_n) \in \alpha \circ [\alpha,\beta] \circ \alpha = \alpha$, we have $(c_1,d_n) \in \alpha \wedge \beta$, so from the fact that $p$ is a difference term we have $p(d_n,c_1,c_1) \equiv_{[\alpha\wedge\beta,\alpha\wedge\beta]} d_n$, and from $[\alpha\wedge\beta,\alpha\wedge\beta] \le [\alpha,\beta]$ we get $c_1 \equiv_{[\alpha,\beta]} d_n$.
\end{proof}

In fact, substantially more is true in varieties with a difference term. Kearnes \cite{kearnes-difference} shows that almost all properties of the commutator which hold in congruence modular varieties generalize to varieties with a difference term, other than $[\alpha_1 \vee \alpha_2,\beta] = [\alpha_1,\beta] \vee [\alpha_2,\beta]$. This property must be weakened, but it is at least true that if $[\alpha_1,\beta] = [\alpha_2,\beta]$ then $[\alpha_1\vee\alpha_2,\beta] = [\alpha_1,\beta]$ in varieties with difference terms.

If we go beyond varieties with a difference term, the commutator may no longer be symmetric. For instance, in the algebra $\bA = (\{0,1,2,*\},\cdot)$ with $\cdot$ given by
\begin{center}
\begin{tabular}{c|cccc} $\cdot$ & $0$ & $1$ & $2$ & $*$\\ \hline $0$ & $0$ & $2$ & $1$ & $*$\\ $1$ & $2$ & $1$ & $0$ & $*$\\ $2$ & $1$ & $0$ & $2$ & $*$\\ $*$ & $*$ & $*$ & $*$ & $*$\end{tabular},
\end{center}
if we let $\theta \in \Con(\bA)$ be the congruence corresponding to the partition $\{0,1,2\},\{*\}$, then we have
\[
[\theta,1_\bA] = 0_\bA, \;\;\; [1_\bA,\theta] = \theta.
\]
The simplest way to fix this asymmetry is to make the following definition from \cite{kearnes-taylor-affine}.

\begin{defn} If $\alpha, \beta \in \Con(\bA)$, then we define their \emph{symmetric commutator}, written $[\alpha,\beta]_s$, to be the least congruence $\delta$ such that both $C(\alpha,\beta;\delta)$ and $C(\beta,\alpha;\delta)$ hold.
\end{defn}

To see that $[\alpha,\beta]_s$ is well-defined, we use Proposition \ref{gen-commutator}(a) to see that the intersection of any collection of congruences that simultaneously satisfy $C(\alpha,\beta;\delta)$ and $C(\beta,\alpha;\delta)$ will also satisfy this pair of term conditions. Since the two-term commutator satisfies $C(\alpha,\beta;[\alpha,\beta]_2)$ and $C(\beta,\alpha;[\alpha,\beta]_2)$, we always have $[\alpha,\beta]_s \le [\alpha,\beta]_2$.

We can also go in the other direction, and define a more general commutator by trying to directly think about what an algebra needs to satisfy to be quasi-affine. This leads to the following definition.

\begin{defn} If $\alpha, \beta \in \Con(\bA)$, then we define their \emph{linear commutator}, written $[\alpha,\beta]_\ell$, as follows. Define a group $\bG$ with the following presentation: the generators of $\bG$ are the elements of $\bA$, and the relations are given by
\[
\begin{bmatrix} a & b\\ c & d\end{bmatrix} \in \bM(\alpha, \beta) \;\; \implies \;\; a+d = b+c \text{ in } \bG.
\]
Then we define the equivalence relation $[\alpha,\beta]_\ell$ to be the kernel of the natural map $\bA \rightarrow \bG$.
\end{defn}

\begin{prop} The linear commutator $[\alpha,\beta]_\ell$ always defines a congruence of $\bA$, we have $[\alpha,\beta]_n \le [\alpha,\beta]_\ell$ for each $n$ and $[\alpha,\beta]_\ell \le \alpha \wedge \beta$, and $[1_\bA,1_\bA]_\ell = 0_\bA$ iff $\bA$ is quasi-affine.
\end{prop}
\begin{proof} The linear commutator can be defined combinatorially as follows. For each matrix $M = \begin{bmatrix} a & b\\ c & d\end{bmatrix} \in \bM(\alpha, \beta)$, call $a,d$ the ``positive'' corners of the matrix $M$, and $b,c$ the ``negative'' corners of $M$. Then $(x,y) \in [\alpha,\beta]_\ell$ iff there is some collection of matrices $M_i \in \bM(\alpha, \beta)$ and a way to pair off values in the positive corners of the matrices $M_i$ to equal values in negative corners of the matrices, so that the only unpaired values are $x$ and $y$, with one occuring in a positive corner of some matrix and the other occuring in a negative corner of some matrix. This defines the linear commutator $[\alpha,\beta]_\ell$ as the union of a directed limit of relations defined by primitive positive formulas in $\bM(\alpha,\beta)$, so the equivalence relation $[\alpha, \beta]_\ell$ is compatible with the operations of $\bA$.

The inequality $[\alpha,\beta]_n \le [\alpha,\beta]_\ell$ follows from the combinatorial description of $[\alpha,\beta]_n$ above (the $n$-cycle condition is a special case of the general setup of matching corners of matrices together). To prove that $[\alpha,\beta]_\ell \le \alpha \wedge \beta$, we just need to check that $[\alpha,\beta]_\ell \le \alpha$ by symmetry, and this follows by chasing equalities and congruences through the matrices in the combinatorial description of $[\alpha,\beta]_\ell$.

If $\bA$ is quasi-affine, then it is easy to see that $[1_\bA,1_\bA]_\ell = 0_\bA$. Finally, if $[1_\bA,1_\bA]_\ell = 0_\bA$, then $\bA$ embeds injectively into the group $\bG$ from the definition of $[1_\bA,1_\bA]_\ell$, and we can generalize the combinatorial description of $[\alpha,\beta]_\ell$ to get a combinatorial description of the restriction of the $4$-ary relation $a+d = b+c$ to $\bA$, so this $4$-ary relation is preserved by the operations of $\bA$.
\end{proof}

We have the following relationship between the various commutators which have been defined so far:
\[
[\alpha,\beta] \le [\alpha,\beta]_s \le [\alpha,\beta]_2 \le [\alpha,\beta]_3 \le \cdots \le [\alpha,\beta]_\ell \le \alpha \wedge \beta,
\]
and among these, the commutators $[\alpha,\beta]_s, [\alpha, \beta]_2$, and $[\alpha,\beta]_\ell$ are symmetric by construction. In \cite{kearnes-taylor-affine}, Kearnes and Szendrei prove that in every Taylor variety we always have $[\alpha,\beta]_s = [\alpha,\beta]_\ell$, so almost all of the commutators collapse into a single concept in Taylor varieties (and they \emph{all} collapse in varieties with difference terms). They also give an alternative characterization of the linear commutator by showing that it is equivalent to the commutator obtained by first ``freely'' extending your variety to make the basic operations multinear inside some larger abelian group, and then computing commutators in the multilinear setting, which has a Mal'cev operation $x - y + z$.

\chapter{Compact Representations and algebras with Few Subpowers}\label{chapter-few-subpowers}

\section{Generalized Majority-Minority operations (motivating Few Subpowers)}

The Few Subpowers algorithm was heavily influenced by Dalmau's paper on generalized majority-minority operations \cite{dalmau-gmm}. Dalmau's motivation was that in both near-unanimity algebras and Mal'cev algebras, every subalgebra of $\bA^n$ has a nice generating set: in the Mal'cev case, we can use a compact representation, while in the near-unanimity case, if the arity is $l+1$, we can use any set of elements which has the same projection onto every subset of the coordinates of size at most $l$. The goal was to unify these two cases.

\begin{defn} An operation $\varphi$ is a \emph{generalized majority-minority} operation (abbreviated as \emph{gmm} operation) if for each pair $a,b$ we either have
\[
\varphi(x,y,...,y) = \varphi(y,x,...,y) = \cdots = \varphi(y,y,...,x) = y \;\;\; \text{ for all }x,y \in \{a,b\},
\]
or
\[
\varphi(x,y,...,y) = \varphi(y,y,...,x) = x \;\;\; \text{ for all }x,y \in \{a,b\}.
\]
In the second case we say that $a,b$ is a \emph{minority pair} for $\varphi$.
\end{defn}

\begin{defn} If $R \subseteq \bA_1 \times \cdots \times \bA_n$, then we define the \emph{signature} of $R$, written $\operatorname{Sig}(R)$, to be the set of triples $(i,a,b)$ with $i \in \{1, ..., n\}$, $a,b$ a minority pair in $\bA_i$, such that there are some $t_a,t_b \in R$ with $\pi_{1, ..., i-1}(t_a) = \pi_{1, ..., i-1}(t_b)$ and $\pi_i(t_a) = a, \pi_i(t_b) = b$. In this case we say that the pair $t_a,t_b$ \emph{witnesses} the triple $(i,a,b)$.
\end{defn}

\begin{thm}\label{gmm-compact} If $\RR \le \bA_1 \times \cdots \times \bA_n$ is preserved by an $l+1$-ary gmm operation $\varphi$ and $S \subseteq \RR$ has $\Sig(S) = \Sig(\RR)$ and $\pi_I(S) = \pi_I(\RR)$ for all $I \subseteq \{1, ..., n\}$ with $|I| \le l$, then $\RR$ is generated by $S$ (using only $\varphi$).
\end{thm}
\begin{proof} We prove this by induction on the arity $n$ of $\RR$. Suppose that $a = (a_1, ..., a_n) \in \RR$, by the induction hypothesis there is some $b_n$ with $(a_1, ..., a_{n-1}, b_n)$ in the subalgebra generated by $S$. We have two cases, based on whether $a_n,b_n$ is a majority pair or a minority pair.

{\bf Case 1: }$a_n,b_n$ is a majority pair. In this case we show that for every $I \subseteq \{1, ..., n\}$, we have $\pi_I a$ in the subalgebra generated by $\pi_I S$, by induction on $|I|$. We already know it for $|I| \le l$ and for $n \not\in I$. Suppose $I = \{i_1, ..., i_m\}$ with $i_1 < \cdots < i_m = n$ and $m \ge l+1$. By the inductive hypothesis, there are elements $b_{i_1}, ...$ such that
\[
(b_{i_1},a_{i_2},...,a_n), (a_{i_1}, b_{i_2}, ..., a_n), ..., (a_{i_1}, a_{i_2}, ..., b_n) \in \Sg_\varphi(S).
\]
If some $b_i = a_i$ then we are done. If some pair $a_i,b_i$ is minority then - assuming WLOG that $a_{i_1}, b_{i_1}$ is minority - we have
\[
\varphi\left(\begin{bmatrix} b_{i_1} & \cdots & b_{i_1} & a_{i_1}\\ a_{i_2} & \cdots & a_{i_2} & a_{i_2}\\ \vdots & \ddots & \vdots & \vdots\\ a_n & \cdots & a_n & b_n\end{bmatrix}\right) = \begin{bmatrix} a_{i_1}\\ a_{i_2}\\ \vdots\\ a_n\end{bmatrix} \in \Sg_\varphi(\pi_I S),
\]
where all but the last column of the displayed matrix are equal. Otherwise, if all pairs $a_i,b_i$ are majority, then we have
\[
\varphi\left(\begin{bmatrix} b_{i_1} & a_{i_1} & \cdots & a_{i_1}\\ a_{i_2} & b_{i_2} & \cdots & a_{i_2}\\ \vdots & \vdots & \ddots & \vdots\\ a_n & a_n & \cdots & b_n\end{bmatrix}\right) = \begin{bmatrix} a_{i_1}\\ a_{i_2}\\ \vdots\\ a_n\end{bmatrix} \in \Sg_\varphi(\pi_I S),
\]
where all of the columns of the displayed matrix are distinct, which is possible because $m \ge l+1$.

{\bf Case 2: }$a_n,b_n$ is a minority pair. In this case, by the assumption $\Sig(S) = \Sig(\RR)$, there are $c,d \in S$ witnessing the triple $(n,a_n,b_n)$. Set $b = (a_1, ..., a_{n-1}, b_n)$, then we claim that
\[
a = \varphi(b,b,...,b,\varphi(b,d,...,d,c)).
\]
First consider the last coordinate: since $a_n,b_n$ is a minority pair and $c_n = a_n, d_n = b_n$, we have
\[
\varphi(b_n, ..., b_n,\varphi(b_n,d_n,...,d_n,c_n)) = \varphi(b_n, ..., b_n, \varphi(b_n, ..., b_n, a_n)) = a_n,
\]
so the last coordinates agree. For $i < n$, we have $a_i = b_i$ and $c_i = d_i$, so
\[
\varphi(b_i, ..., b_i,\varphi(b_i,d_i,...,d_i,c_i)) = \varphi(a_i, ..., a_i,\varphi(a_i,c_i,...,c_i,c_i)) = a_i,
\]
where the last equality holds regardless of whether $a_i,c_i$ is a majority pair or a minority pair.
\end{proof}

\begin{defn} A subset $S \subseteq \RR$ is called a \emph{compact representation} of a relation $\RR$ preserved by an $l+1$-ary gmm operation if $\Sig(S) = \Sig(\RR)$, $\pi_I(S) = \pi_I(\RR)$ for every $I$ with $|I| \le l$, and $|S| \le 2|\Sig(\RR)| + \sum_{|I| \le l} |\pi_I(\RR)|$.
\end{defn}

In order to manipulate compact representations of relations, we again define subroutines \texttt{Nonempty}, \texttt{Fix-values}, \texttt{Next-beta}, and \texttt{Intersect}:
\begin{itemize}
\item \texttt{Nonempty}$(R, i_1, ..., i_k, \bS)$ takes $R$ a compact representation of $\RR \le \bA_1\times \cdots\times \bA_n$, $\bS \le \bA_{i_1} \times \cdots \times \bA_{i_k}$, computes the subalgebra generated by $\pi_{i_1, ..., i_k}(R)$ under $\varphi$, and if this intersects with $\bS$, then it returns an element of $\RR$ which maps to an element of the intersection,

\item \texttt{Fix-values}$(R, a_1, ..., a_m)$ takes $R$ a compact representation of $\RR \le \bA_1\times \cdots\times \bA_n$ and returns a compact representation of the relation $x \in \RR \wedge (x_1 = a_1) \wedge \cdots \wedge (x_m = a_m)$ by inductively fixing one coordinate $x_i$ to $a_i$ at a time, and for each new coordinate that is fixed we compute a new compact representation by computing projections onto at most $l$ coordinates using \texttt{Nonempty} and computing witnesses for triples in the signature using the proof of Case 2 of Theorem \ref{gmm-compact},

\item \texttt{Next-beta}$(R, i_1, ..., i_k, \bS)$ takes $R$ a compact representation of $\RR \le \bA_1\times \cdots\times \bA_n$, $\bS \le \bA_{i_1} \times \cdots \times \bA_{i_k}$, and returns a compact representation of $\RR \cap \bS$ by computing all projections onto at most $l$ coordinates using \texttt{Nonempty} and computing witnesses for triples in the signature using \texttt{Fix-values} and \texttt{Nonempty}, and

\item \texttt{Intersect}$(R, i_1, ..., i_k, S)$ takes $R$ a compact representation of $\RR \le \bA_1\times \cdots\times \bA_n$, $S$ a compact representation of $\bS \le \bA_{i_1} \times \cdots \times \bA_{i_k}$, and computes a compact representation for $\RR\cap \bS$ by first making a compact representation of $\RR\times \bS$ and then repeatedly calling \texttt{Next-beta} to intersect this with the equality relation on the pair of coordinates $i_j,n+j$.
\end{itemize}

The only subroutine which has changed substantially from the Mal'cev case is the \texttt{Fix-values} subroutine.

\begin{algorithm}
\caption{\texttt{Fix-values}$(R, a_1, ..., a_m)$, $\varphi$ an $l+1$-ary gmm term, $R$ a compact representation of $\RR \le \bA_1\times \cdots\times \bA_n$.}
\begin{algorithmic}[1]
\State Set $R_0 \gets R$.
\For{$j$ from $1$ to $m$}
\State Let $R_j \gets \emptyset$.
\ForAll{$I = \{i_1, ...\} \subseteq \{1,...,n\}$ with $|I| \le l$ and $(b_{i_1}, ...) \in \pi_I(R_{j-1})$}
\State Set $R_j \gets R_j \cup \texttt{Nonempty}(R_{j-1},j,i_1,...,i_{|I|},\{(a_j,b_{i_1},...,b_{i_{|I|}})\})$.
\EndFor
\ForAll{$(i,a,b) \in \Sig(R_{j-1})$ with $i > j$ and $a,b$ a minority pair}
\State Let $t_a,t_b \in R_{j-1}$ witness the triple $(i,a,b)$.
\State Let $t \gets \texttt{Nonempty}(R_{j-1},j,i,\{(a_{j},a)\})$.
\If{$t \ne \emptyset$}
\State Set $R_{j} \gets R_{j} \cup \{t, \varphi(t,t,...,t,\varphi(t,t_a,...,t_a,t_b))\}$.
\EndIf
\EndFor
\EndFor
\State \Return $R_m$.
\end{algorithmic}
\end{algorithm}

Reviewing what we've done, we have a procedure for converting proofs that compact representations generate relations into algorithms for computing compact representations of intersections for relations. The most critical step of the algorithm is the step of the \texttt{Fix-values} subroutine in which we convert a pair that witnesses a triple $(i,a,b)$ in $R_{j-1}$ to a pair that witnesses a triple $(i,a,b)$ in $R_j$.

Before we go on, we can use this algorithm to settle the dichotomy conjecture for constraint languages which contain ``swap'' relations $\{(a,b),(b,a)\}$ for every pair of elements $a,b$.

\begin{thm}\label{thm-swap-gmm} Suppose that $\fA = (A,\Gamma)$ is a relational structure where $\Gamma$ is a set of relations which contains the swap relation $S_{ab} = \{(a,b),(b,a)\}$ for every pair $a,b \in \fA$. Then either $\CSP(\Gamma)$ is NP-complete, or $\fA$ has a ternary generalized majority-minority polymorphism. In the second case, $\CSP(\Gamma)$ can be solved in polynomial time by Dalmau's algorithm.
\end{thm}
\begin{proof} Note that $\Gamma$ is automatically core, since any unary polymorphism of $S_{ab}$ must send $a,b$ to distinct values in $\{a,b\}$. Thus if $\CSP(\Gamma)$ is not NP-complete, then it must have a Taylor polymorphism $t$.

First we will show that this implies that for all $a,b \in \fA$ there is a ternary polymorphism $f_{ab}$ such that the restriction of $f_{ab}$ to $\{a,b\}$ is either the majority operation or the minority operation. Since $\pi_1(S_{ab}) = \{a,b\}$, the set $\{a,b\}$ is closed under $t$. Let $t' \in \Clo(t)$ have minimal arity such that the restriction of $t'$ to $\{a,b\}$ is not a projection. An elementary combinatorial argument known as {\'S}wierczkowski's Lemma \cite{semiprojection-lemma} shows that if $t'$ has arity at least four, then there is some way of identifying two variables of $t'$ to get a term $t''$ of smaller arity such that the restriction of $t''$ to $\{a,b\}$ is also not a projection. Thus the arity of $t'$ is at most three. The arity of $t'$ can't be one or two since $t'$ is idempotent and preserves $S_{ab}$.

Since every way of identifying two variables of $t'|_{\{a,b\}}$ gives a projection, up to reordering the variables of $t'$ there are just three cases. In two of these cases, $t'$ already restricts to a majority or minority operation on $\{a,b\}$. In the remaining case, after reordering the variables we may assume that $t'(x,y,y) = t'(y,y,x) = t'(x,y,x) = x$ for $x,y \in \{a,b\}$, and taking $f_{ab}(x,y,z) = t'(x,t'(x,y,z),z)$ gives a function $f_{ab}$ which restricts to a majority operation on $\{a,b\}$.

Now we choose any ordering of the collection of pairs $\{a,b\}$, with the $i$th pair given by $\{a_i,b_i\}$. We inductively define functions $f_i \in \Clo(t)$ by $f_0 = \pi_1$, and for $i \ge 0$ we set
\[
f_{i+1}(x,y,z) = f_{a_ib_i}(f_i(x,y,z),f_i(y,z,x),f_i(z,x,y)).
\]
We claim that the final function $f_n$ (with $n = \binom{|A|}{2}$) is a generalized majority-minority polymorphism of $\fA$. Since each $f_{ab}$ is idempotent, it's enough to check that the restriction of $f_{i+1}$ to $\{a_i,b_i\}$ is either a pure majority or pure minority function.

From the fact that $f_i$ preserves the unary relation $\pi_1(S_{a_ib_i}) = \{a_i,b_i\}$ and the fact that the restriction of $f_{a_ib_i}$ to $\{a_i,b_i\}$ is invariant under cyclically permuting its input variables, we see that $f_{i+1}$ also restricts to a cyclic term on $\{a_i,b_i\}$. Since $f_{i+1}$ preserves $S_{ab}$, it must therefore either restrict to the pure majority or pure minority function on $\{a_i,b_i\}$.
\end{proof}

There are two examples of generalized majority-minority algebras on a three element domain which do not come from majority or Mal'cev operations, and correspond to maximal tractable constraint languages.

\begin{ex}\label{gmm-ex-1} The first example is $\bA_1 = (\{a,b,c\}, \varphi_1)$, where $\varphi_1$ is a ternary gmm such that $\{a,x\}$ is a pure minority subalgebra of $\bA_1$ for all $x$, $\{b,c\}$ is a majority subalgebra of $\bA_1$, and the equivalence relation corresponding to the partition $\{a\}, \{b,c\}$ is a congruence $\alpha$ on $\bA_1$ such that the quotient $\bA_1/\alpha$ is a pure minority algebra. Explicitly, $\varphi_1$ is the symmetric idempotent function of its inputs which is given by
\[
\varphi_1(a,a,x) = x, \varphi_1(a,x,x) = a, \varphi_1(b,b,c) = b, \varphi_1(b,c,c) = c, \varphi_1(a,b,c) = a.
\]
The corresponding relational clone is generated by the partial order $\{(a,a),(b,b),(b,c),(c,c)\}$, the order two automorphism $\{(a,a),(b,c),(c,b)\}$, and the affine ternary relation $\{(a,a,b),(a,b,a),(b,a,a),(b,b,b)\}$.
\end{ex}

\begin{ex}\label{gmm-ex-2} The second example is $\bA_2 = (\{a,b,c\}, \varphi_2)$, where $\varphi_2$ is a ternary gmm such that $\{a,x\}$ is a majority subalgebra of $\bA_2$ for all $x$, $\{b,c\}$ is a pure minority subalgebra of $\bA_2$, the equivalence relation corresponding to the partition $\{a\}, \{b,c\}$ is a congruence $\alpha$ on $\bA_2$ such that the quotient $\bA_2/\alpha$ is a majority algebra, and the permutation $(b\; c)$ is an automorphism of $\bA_2$. Explicitly, $\varphi_2$ is the cyclically symmetric idempotent function of its inputs which is given by
\[
\varphi_2(a,a,x) = a, \varphi_2(a,x,x) = x, \varphi_2(b,b,c) = c, \varphi_2(b,c,c) = b, \varphi_2(a,b,c) = b, \varphi_2(a,c,b) = c.
\]
The corresponding relational clone is generated by the binary relations $\{(a,b),(b,a)\}$, $\{(a,a),(a,b),(b,b)\}$, $\{(a,a),(b,c),(c,b)\}$ and the ternary relation $\{(a,a,a),(b,b,b),(b,c,c),(c,b,c),(c,c,b)\}$.
\end{ex}

The reader might notice that generalized majority-minority operations are \emph{not} defined in terms of satisfying a system of identities. So we should be able to immediately generalize Dalmau's result to the variety of algebras generated by algebras with a gmm operation, by finding the identities which are satisfied by a gmm operation that were critical to the correctness of the algorithm. How did we apply the operation $\varphi$, throughout the algorithm \texttt{Fix-values} and the proof of Theorem \ref{gmm-compact}?

The first thing to note is that we often set almost all of the entries of $\varphi$ to the same value. So define auxiliary binary and ternary terms $p,d$ by
\begin{align*}
d(x,y) &= \varphi(x,y,...,y,y),\\
p(x,y,z) &= d(\varphi(x,y,...,y,z),z).
\end{align*}
The important property of $d$ is that we have $d(a,b) = a$ when $a,b$ are a minority pair. For $p$, the important property is that when $a,b$ are a minority pair, then we have $p(a,b,b) = a$, and in every case we always have
\[
p(y,y,z) = z.
\]
We can express the fact that $p(a,b,b) = a$ when $a,b$ are a minority pair by the equation
\[
p(x,y,y) = d(x,y),
\]
which also holds for majority pairs.

Where did we actually use the function $\varphi$? It is only called directly in the subroutine \texttt{Nonempty}. It is crucial that it is actually used there, because the full function $\varphi$ was necessary for Case 1 of Theorem \ref{gmm-compact}. The proof of that case does not immediately appear to generalize, as there was substantial casework within it, based on whether there was a minority pair $a_i,b_i$ or not. However, clever use of the function $d(x,y)$ can mimic the casework that appeared there. For each $a_i, b_i$, the expression $d(a_i,b_i)$ has the nice property that $a_i, d(a_i,b_i)$ automatically forms a majority pair (or an equal pair, which we can think of as a degenerate case of a majority pair). So if we define a function $s(x_0,x_1, ..., x_l)$ by
\[
s(x_0, x_1,..., x_l) = \varphi(x_0,d(x_0,x_1), ..., d(x_0,x_l)),
\]
then we find that
\begin{align*}
s(y,x,x,...,x) &= \varphi(y,d(y,x), ..., d(y,x)) = d(y,x),\\
s(x,y,x,...,x) &= \varphi(x, d(x,y), x, ..., x) = x,\\
s(x,x,y,...,x) &= \varphi(x, x, d(x,y), ..., x) = x,\\
&\;\; \vdots\\
s(x,x,x,...,y) &= \varphi(x, x, x, ..., d(x,y)) = x.
\end{align*}
This function $s$ lets us generalize Case 1 of Theorem \ref{gmm-compact}, the case where $d(a_n,b_n) = b_n$, while the function $p$ was necessary to generalize Case 2. To unify them, we should slightly modify our construction of $s$ to create the following term $e$:
\[
e(u,v,x_1, ..., x_l) = \varphi(v, d(u,x_1), ..., d(u,x_{l-1}),d(x_1,x_l)).
\]
Then $s$ is related to $e$ by
\[
s(x_0, x_1, ..., x_l) = e(x_1,x_0,x_1,...,x_l)\text{ if all but one of the }x_i\text{ are equal},
\]
$p$ is related to $e$ by
\[
p(x,y,z) = e(y,x,z,...,z)\text{ if }x = y\text{ or }y = z,
\]
and $e$ satisfies the identities
\begin{align*}
e(y,y,x,x,...,x) &= \varphi(y,d(y,x),...,d(y,x),x) = x,\\
e(y,x,y,x,...,x) &= \varphi(x, y, d(y,x), ..., d(y,x)) = x,\\
e(x,x,x,y,...,x) &= \varphi(x, x, d(x,y), ..., x) = x,\\
&\;\; \vdots\\
e(x,x,x,x,...,y) &= \varphi(x, x, x, ..., d(x,y)) = x.
\end{align*}
Can we use this system of identities to prove an analogue of Theorem \ref{gmm-compact}? Yes! The trick is to plug things back into $e$, to make the following term $t$:
\[
t(u,v,w,x_1,...,x_l) = e(p(v,u,x_1),s(w,x_1,...,x_l),x_1,...,x_l).
\]
Now if we have a tuple $a = (a_1,...,a_n)$ which we want to prove is in the subalgebra generated by $S$, and if this subalgebra already contains $(a_1,...,a_{i-1},b_i,a_{i+1}, ..., a_n)$ for each $i$, as well as a pair $(c_1, ..., c_{n-1}, a_n), (c_1, ..., c_{n-1}, b_n)$ which witnesses the triple $(n,a_n,b_n)$, then we have
\[
t\left(\begin{bmatrix} c_1 & c_1 & a_1 & b_1 & a_1 & \cdots\\ c_2 & c_2 & a_2 & a_2 & b_2 & \cdots\\ \vdots & \vdots & \vdots & \vdots & \vdots & \ddots\\ a_n & b_n & b_n & a_n & a_n & \cdots\end{bmatrix}\right) = e\left(\begin{bmatrix} b_1 & a_1 & b_1 & a_1 & \cdots\\ a_2 & a_2 & a_2 & b_2 & \cdots\\ \vdots &\vdots & \vdots & \vdots & \ddots\\ d_n & d_n & a_n & a_n & \cdots\end{bmatrix}\right) = \begin{bmatrix} a_1\\ a_2\\ \vdots\\ a_n\end{bmatrix},
\]
where $d_n = d(b_n,a_n)$.

While playing these sorts of games with identities may yield more and more general examples of algebraic structures where relations have compact representations, we are not being very systematic here. So perhaps we should work backwards: what absolutely needs to be true for something like compact representations to exist?

\begin{prop} If every subpower $\RR \le \bA^n$ has a compact representation $S$ consisting of at most $p(n)$ tuples, then the number of different subalgebras of $\bA^n$ is at most $|\bA^n|^{p(n)} = |\bA|^{np(n)}$.
\end{prop}

\begin{cor} No analogue of compact representations can exist for subpowers of a nontrivial semilattice.
\end{cor}
\begin{proof} It's enough to consider the case $\bA = (\{0,1\},\max)$, since every semilattice contains a subalgebra isomorphic to it. The number of subpowers of $\bA^n$ is at least the number of subsets on $\{0,1\}^n$ which are generated by subsets $S \subseteq \{x \in \{0,1\}^n, \sum_i x_i = n/2\}$ (suppose $n$ is even). Any two distinct subsets $S, S'$ of the set of tuples with weight $n/2$ will generate different subalgebras of $\bA^n$, so the number of subalgebras of $\bA^n$ is at least
\[
2^{\binom{n}{n/2}} \ge 2^{2^n/n},
\]
and $2^n/n$ clearly grows faster than any polynomial.
\end{proof}

What makes the semilattice case so different from the Mal'cev case and the near-unanimity case? The main difference is that the identities satisfied by a semilattice do not allow us to get back to $x$ once we start combining it with other values, while the identities for Mal'cev and near-unanimity terms all have $x$s on the right hand sides.

So we should start by trying to prove that having few subpowers implies that there are terms satisfying a nontrivial system of identities which have $x$s on the right hand sides of each identity, such as the system of identities satisfied by the term $e$ constructed earlier. The trick, as we will see, is to apply the existence of compact representations to the case of a power of the free algebra on two generators, considered as a subalgebra of $(\bA^{\bA^2})^n$.

\section{Algebras with Few Subpowers}

First we define an invariant of an algebraic structure and the variety it generates, which is slightly more well-behaved than the function that takes $n$ to the number of subalgebras of $\bA^n$.

\begin{defn} If $\bA$ is an algebraic structure and $a_1, ..., a_k \in \bA$, we say that $a_1, ..., a_k$ are \emph{independent} if no $a_i$ is in the subalgebra generated by the rest of the $a_j$s. For every $n$, we define $i_\bA(n)$ to be the size of the largest independent set in $\bA^n$.
\end{defn}

\begin{prop} If $\bA$ is a finite algebra, then any subalgebra of $\bA^n$ can be generated by at most $i_\bA(n)$ elements, so the number of subalgebras of $\bA^n$ is bounded above by $|\bA^n|^{i_\bA(n)} = 2^{n\lg(|\bA|)i_\bA(n)}$. The number of subalgebras of $\bA^n$ is also bounded below by $2^{i_\bA(n)}$.
\end{prop}
\begin{proof} Since $\bA$ is finite, every subalgebra of $\bA^n$ has a minimal generating set, and this minimal generating set is necessarily independent. The upper bound on the number of subalgebras follows from counting the number of possible minimal generating sets.

For the lower bound on the number of subalgebras, suppose that $a_1, ..., a_k$ are independent in $\bA^n$. Then every subset $S$ of $\{a_1, ..., a_k\}$ generates a distinct subalgebra of $\bA^n$, since $\Sg_{\bA^n}(S) \cap \{a_1, ..., a_k\} = S$ by the definition of independence. Thus $\bA^n$ has at least $2^k$ distinct subalgebras.
\end{proof}

\begin{prop} If $\bB \in HSP(\bA)$ is also finite, then $i_\bB(n) \le i_\bA(cn)$ for some constant $c$ depending only on $\bB$.
\end{prop}
\begin{proof} If $\bA,\bB$ are both finite, then there is some finite number $c$ such that $\bB \in HS(\bA^c)$, that is, there is a subalgebra $\bC \le \bA^c$ and a surjective homomorphism $f : \bC \rightarrow \bB$. Then every independent set in $\bB^n$ lifts to an independent set in $(\bA^c)^n = \bA^{cn}$ by choosing any section of $f$ and applying it coordinate-wise.
\end{proof}

We will apply the above result to the free algebra on two generators $\cF_{\cV(\bA)}(x,y) \le \bA^{\bA^2}$ to prove that if an algebra has few subpowers, then it has a \emph{cube term}. Since cube terms have exponentially high arity, it's necessary to develop some notation to define them properly.

\begin{defn} For every subset $S \subseteq \{1, ..., k\}$, we define the $k$-dimensional column vector $v^S$ by
\[
v^S_i = \begin{cases} y & i \in S,\\ x & i \not\in S.\end{cases}
\]
A $k$-\emph{cube term} is a term $t$ with variables indexed by nonempty subsets of $\{1, ..., k\}$, such that if we fix an enumeration $S_1, ..., S_{2^k-1}$ of these subsets, we have the identity
\[
t(v^{S_1}, ..., v^{S_{2^k-1}}) \approx v^{\emptyset}.
\]
\end{defn}

For instance, if $k = 3$ then (with one possible choice of variable ordering) a $3$-cube term is a $7$-ary term $t$ satisfying the identity
\[
t\left(\begin{bmatrix} y & y & y & x & y & x & x\\ y & y & x & y & x & y & x\\ y & x & y & y & x & x & y\end{bmatrix}\right) \approx \begin{bmatrix} x\\ x\\ x\end{bmatrix}.
\]
Note that a Mal'cev term is the same as a $2$-cube term (up to reordering variables).

\begin{thm}[Few subpowers implies cube term \cite{few-subpowers}] Let $\bF = \cF_{\cV(\bA)}(x,y) \le \bA^{\bA^2}$ be the free algebra on two generators in the variety generated by $\bA$.
\begin{itemize}
\item If $i_{\bF}(k) < 2^k$ for any $k$, then $\bA$ has a $k$-cube term.

\item If $i_{\bF}(m) < \binom{m}{k}$ for any $m,k$, then $\bA$ has a $k$-cube term.
\end{itemize}
In particular, if $i_\bA(n) = o(n^k)$ then $\bA$ has a $k$-cube term, and if $i_\bA(n) = 2^{o(n)}$ then there exists some $k$ such that $\bA$ has a $k$-cube term.
\end{thm}
\begin{proof} For the first statement, if $i_{\bF}(k) < 2^k$, then the vectors $v^S$ for $S \subseteq \{1, ..., k\}$ can't be independent, so some $v^S$ is in the subalgebra generated by the others. By applying an automorphism of $\bF^k$ which swaps $x$s and $y$s in the coordinates belonging to $S$, we may assume without loss of generality that $S = \emptyset$. From $v^\emptyset \in \Sg_{\bF^k}\{v^S \mid S \ne \emptyset\}$, we see that there is a term $t$ such that $t(v^{S_1}, ...) = v^\emptyset$, and since $\bF$ is the free algebra on two generators, this implies the $k$-cube term identities.

For the second statement, consider the set of vectors $v^S$ with $S \in \binom{\{1, ..., m\}}{k}$. By assumption, these are not independent, so some $v^S$ is in the subalgebra generated by the others. Then if we project onto the coordinates of $S$ and use the fact that for $S \ne T$ with $|S| = |T|$ we never have $S \subseteq T$, we get the situation of the previous paragraph inside $\bF^S \cong \bF^k$.
\end{proof}

Next, we upgrade the $k$-cube term by repeatedly plugging it into itself to produce simpler terms, finally arriving at the $k$-edge term.

\begin{defn} If $\Delta \subseteq \cP(\{1, ..., k\}) \setminus \{\emptyset\}$, then we say that $t$ is a $\Delta$-\emph{cube} term if it has variables indexed by elements of $\Delta$ and satisfies the identity $t(v^{S_1}, ...) = v^\emptyset$, where $S_1, ...$ is an enumeration of the elements of $\Delta$.

If we set $\Delta^e = \{\{1,2\}, \{1\}, \{2\}, ..., \{k\}\}$, then a $\Delta^e$-cube term is called a $k$-\emph{edge} term.
\end{defn}

A $k$-edge term is simple enough that we can write out the identities it satisfies explicitly: a $k+1$-ary term $e$ is a $k$-edge term iff it satisfies
\[
e\left(\begin{bmatrix}y & y & x & x & \cdots & x\\ y & x & y & x & \cdots & x\\ x & x & x & y & \cdots & x\\ \vdots & \vdots & \vdots & \vdots & \ddots & \vdots \\ x & x & x & x & \cdots & y \end{bmatrix}\right) \approx \begin{bmatrix}x\\ x\\ x\\ \vdots \\ x\end{bmatrix}.
\]

\begin{thm}[Cube term implies edge term \cite{few-subpowers}]\label{cube-edge} If $\bA$ has a $k$-cube term, then it also has a $k$-edge term.
\end{thm}
\begin{proof} Since it is hard to deal with terms having exponentially many variables, we will do the last step of the proof first, and show that if $\bA$ has a $\Delta^*$-cube term $t^*$ then it has a $k$-edge term, where
\[
\Delta^* = \{\{1,2\}, ..., \{1,k\}, \{1\}, \{2\}, ..., \{k\}\}
\]
only has $2k-1$ elements. The $\Delta^*$-cube term identities for $t^*$ state that
\[
t^*\left(\begin{bmatrix} y & y & \cdots & y & y & x & x & \cdots & x\\ y & x & \cdots & x & x & y & x & \cdots & x\\ x & y & \cdots & x & x & x & y & \cdots & x\\ \vdots & \vdots & \ddots & \vdots & \vdots & \vdots & \vdots & \ddots & \vdots \\ x & x & \cdots & y & x & x & x & \cdots & y \end{bmatrix}\right) \approx \begin{bmatrix}x\\ x\\ x\\ \vdots \\ x\end{bmatrix}.
\]
In order to show that there is a $k$-edge term, we just need to show that $v^\emptyset$ can be generated from $\{v^S \mid S \in \Delta^e\}$ using the $\Delta^*$-cube term $t^*$.

Let $a = t^*(x,...,x,y,x,...,x)$, where the only $y$ occurs at the index corresponding to $\{1\}$ (this is the middle index if we order the variables of $t$ as in the displayed indentities above). First we will use $t$ to generate vectors $v^{S,a}$ for $S \in \Delta^*$ which look just like the vectors $v^S$, except $y$s in the first coordinate are replaced by $a$s. If $S \in \Delta^*$ and $1 \not\in S$, then $S$ is already in $\Delta^e$ and $v^{S,a} = v^S$, so we don't have to worry about these. If $S = \{1\}$, then we use
\[
t^*\left(\begin{bmatrix} x & x & \cdots & x & y & x & x & \cdots & x\\ y & x & \cdots & x & x & y & x & \cdots & x\\ x & y & \cdots & x & x & x & y & \cdots & x\\ \vdots & \vdots & \ddots & \vdots & \vdots & \vdots & \vdots & \ddots & \vdots \\ x & x & \cdots & y & x & x & x & \cdots & y \end{bmatrix}\right) = \begin{bmatrix}a\\ x\\ x\\ \vdots \\ x\end{bmatrix},
\]
and note that every column of the matrix on the left hand side is $v^S$ for some $S \in \Delta^e$. If $S = \{1,2\}$, then we use
\[
t^*\left(\begin{bmatrix} x & x & \cdots & x & y & x & x & \cdots & x\\ y & y & \cdots & y & y & y & y & \cdots & y\\ x & x & \cdots & x & x & x & x & \cdots & x\\ \vdots & \vdots & \ddots & \vdots & \vdots & \vdots & \vdots & \ddots & \vdots \\ x & x & \cdots & x & x & x & x & \cdots & x \end{bmatrix}\right) = \begin{bmatrix}a\\ y\\ x\\ \vdots \\ x\end{bmatrix},
\]
again noting that every column corresponds to an element of $\Delta^e$. Finally, if $S = \{1,i\}$, say $S = \{1,3\}$ without loss of generality, then we use
\[
t^*\left(\begin{bmatrix} x & x & \cdots & x & y & x & x & \cdots & x\\ y & y & \cdots & y & y & x & x & \cdots & x\\ x & x & \cdots & x & x & y & y & \cdots & y\\ x & x & \cdots & x & x & x & x & \cdots & x\\ \vdots & \vdots & \ddots & \vdots & \vdots & \vdots & \vdots & \ddots & \vdots \\ x & x & \cdots & x & x & x & x & \cdots & x \end{bmatrix}\right) = \begin{bmatrix}a\\ x\\ y\\ x\\ \vdots \\ x\end{bmatrix},
\]
where every row other than the first three (or other than the first, second, and $i$th in the general case) is all $x$s, and again every column belongs to $\Delta^e$.

Now that we've constructed the $v^{S,a}$s for all $S \in \Delta^*$, we use $t^*$ to put them all together:
\[
t^*\left(\begin{bmatrix} a & a & \cdots & a & a & x & x & \cdots & x\\ y & x & \cdots & x & x & y & x & \cdots & x\\ x & y & \cdots & x & x & x & y & \cdots & x\\ \vdots & \vdots & \ddots & \vdots & \vdots & \vdots & \vdots & \ddots & \vdots \\ x & x & \cdots & y & x & x & x & \cdots & y \end{bmatrix}\right) = \begin{bmatrix}x\\ x\\ x\\ \vdots \\ x\end{bmatrix}.
\]
Thus if $\bA$ has a $\Delta^*$-cube term, then it has a $k$-edge term. Explicitly, the construction we just worked through corresponds to the formula
\begin{align*}
e(x_0,x_1,...,x_k) = t^*(t^*&(x_2, ..., x_2, x_0, x_2, ..., x_2), t^*(x_2, ..., x_2, x_0, x_3, ..., x_3), ...,\\
&t^*(x_2, ..., x_2, x_0, x_k, ..., x_k), t^*(x_2, ..., x_k, x_1, x_2, ..., x_k), x_2, ..., x_k).
\end{align*}

Now that we have the general idea down, we work through the inductive argument needed to prove that if we have a $k$-cube term, then we have a $\Delta^*$-cube term. Let $\Delta^{\ell *} = \Delta^* \cup \cP(\{1,..., \ell\})\setminus \emptyset$. Note that a $k$-cube term is the same as a $\Delta^{k *}$-cube term, and a $\Delta^*$-cube term is the same as a $\Delta^{0*}$-cube term.

{\bf Claim:} If $\bA$ has a $\Delta^{\ell *}$-cube term $t^\ell$, then it also has a $\Delta^{(\ell -1)*}$-cube term.

{\bf Proof of Claim:} We argue as before, this time taking $a = t^\ell(x,...,x,y,x,...,x)$, where the lone $y$ occurs in the index corresponding to $\{\ell\}$. For $S \in \Delta^{\ell *}$, we let $v^{S,a}$ be the vector similar to $v^S$, but with any $y$ in the $\ell$th coordinate replaced with an $a$. We just need to generate each $v^{S,a}$ for $S \in \Delta^{\ell *}$ using the vectors coming from $\Delta^{(\ell-1)*}$. Again, if $\ell \not\in S$ then $v^{S,a} = v^S$ and $S \in \Delta^{(\ell-1)*}$ already.

If $S = \{\ell\}$, then we plug in the matrix $M$ to $t^\ell$ which looks just like the matrix which gives the defining identities for $t^\ell$, but has the $\ell$th row replaced by the sequence of $x$s and $y$s we used to define $a$. Explicitly, $M$ is given by
\[
\begin{array}{c|cc} M_{i,T} & T \ne \{\ell\} & T = \{\ell\}\\ \hline i \ne \ell & v^T_i & v^T_i=x\\ i = \ell & x & y.\end{array}
\]
Then $t^\ell(M) = v^{\{\ell\},a}$, and the $T$th column of $M$ is $v^{T\setminus\{\ell\}}$ if $T \ne \{\ell\}$ and is $v^{\{\ell\}}$ if $T = \{\ell\}$.

If $\ell \in S$ but $S \ne \{\ell\}$, then we plug in a matrix $M^S$ such that each of its columns is equal to one of $v^{S\setminus \{\ell\}}, v^{\{1\}}, v^{\{1,\ell\}}$: if $\ell \not\in T$, then the $T$th column of $M^S$ is $v^{S\setminus \{\ell\}}$, if $\ell \in T$ but $T \ne \{\ell\}$ then the $T$th column is $v^{\{1\}}$, and if $T = \{\ell\}$ then the $T$th column is $v^{\{1,\ell\}}$. Explicitly, $M^S$ is given by
\[
\begin{array}{c|ccc} M^S_{i,T} & \ell \not\in T & \ell \in T \ne \{\ell\} & T = \{\ell\}\\ \hline i = 1 \in S & y & y & y \\ i = 1 \not\in S & x & y & y \\ i \ne 1,\ell,\ i \in S & y & x & x\\ i \ne 1,\ell,\ i \not\in S & x & x & x\\ i = \ell & x & x & y.\end{array}
\]
These choices ensure that $t^\ell(M^S) = v^{S,a}$.

To finish, we apply $t^\ell$ to the set of vectors $v^{S,a}$ for $S \in \Delta^{\ell *}$, and see that the defining identities for $t^\ell$ imply that the resulting vector is $v^\emptyset$. Thus there is a $\Delta^{(\ell-1)*}$-cube term $t^{\ell-1}$ which can in principle be written explicitly by plugging in variables to the star composition $t^\ell * t^\ell$.
\end{proof}

From a $k$-edge term $e$, we can now construct terms $s,p$ that act like near-unanimity and Mal'cev terms which have been ``glued together'' by a binary term $d$. I've rearranged the variables of these terms from the notation used in \cite{few-subpowers}, for the sake of readability and for consistency with the notation used in Appendix \ref{a-commutator}.

\begin{thm}[Edge terms imply terms $s,p,d$ \cite{few-subpowers}]\label{edge-spd} If $e$ is a $k$-edge term on a finite algebra $\bA$, then there are terms $s,p,d \in \Clo(e)$ with $s$ $k$-ary which satisfy the system of identities
\begin{align*}
s(y,x,x,...,x) &\approx d(y,x),\\
s(x,y,x,...,x) &\approx x,\\
&\;\; \vdots\\
s(x,x,x,...,y) &\approx x,\\
p(y,y,x) &\approx x,\\
p(x,y,y) &\approx d(x,y),\\
d(d(x,y),y) &\approx d(x,y).
\end{align*}
Furthermore, these terms can be computed from $e$ in time $O(|\bA|^k)$. If $\bA$ is infinite, then we can find terms $s,p,d \in \Clo(e)$ satisfying all but the last displayed identity.
\end{thm}
\begin{proof} If we ignore the last identity involving $d$, we can find terms $s_1,p_1,d_1$ satisfying the other identities as follows:
\begin{align*}
s_1(x_1, x_2, ..., x_k) &= e(x_2,x_1,x_2, ..., x_k),\\
p_1(x,y,z) &= e(y,x,z,...,z),\\
d_1(x,y) &= e(y,x,y,...,y).
\end{align*}

We can get the last identity by an iteration argument. For each $i$, we set
\begin{align*}
s_{i+1}(x_1,x_2, ..., x_k) &= s_1(s_i(x_1,x_2,....,x_k),x_2,...,x_k),\\
p_{i+1}(x,y,z) &= p_1(d_i(x,y),y,z),\\
d_{i+1}(x,y) &= d_1(d_i(x,y),y).
\end{align*}
Then for each $i$, the terms $s_i,p_i,d_i$ satisfy the desired identities aside from the last one. Since $\bA$ is finite, we can take $i = |\bA|!$ to find that
\[
d_{|\bA|!}(d_{|\bA|!}(x,y),y) = d_{|\bA|!}(x,y)
\]
for all $x,y \in \bA$.

To compute $s_{|\bA|!}$ efficiently from $e$, first we compute $s_1$, and then for each choice of $a_2, ..., a_k \in \bA$ we find the induced unary polynomial $f_{a_2,...,a_k} : x_1 \mapsto s_1(x_1, a_2, ..., a_k)$. To finish, we note that for every unary function $f : \bA \rightarrow \bA$ we can compute $f^{\infty} \coloneqq \lim_{n \rightarrow \infty} f^{\circ n!}$ in time $O(|\bA|)$ using a clever algorithm which we will go over later, but which the reader may enjoy trying to discover now as an exercise.
\end{proof}

Now we can use the binary term $d$ to define minority pairs and signatures.

\begin{defn} If $s,p,d$ are terms as in the Theorem \ref{edge-spd}, then we say that $a,b \in \bA$ are a \emph{minority pair} if $d(b,a) = b$. If $R \subseteq \bA_1 \times \cdots \times \bA_n$, then we say that $(i,a,b)$ is a \emph{minority index} of $R$ which is \emph{witnessed} by a pair $t_a, t_b \in R$ if:
\begin{itemize}
\item $a,b$ are a minority pair, i.e. $d(b,a) = b$,

\item the pair $t_a, t_b$ agree up to coordinate $i$: $\pi_{1, ..., i-1}(t_a) = \pi_{1, ..., i-1}(t_b)$, and

\item we have $\pi_i(t_a) = a, \pi_i(t_b) = b$.
\end{itemize}
We define the \emph{signature} of $R$, written $\Sig(R)$, to be the set of minority indices which are witnessed by pairs in $R$.
\end{defn}

\begin{defn} If $\RR \le \bA_1 \times \cdots \times \bA_n$ and the $\bA_i$ are in a variety with a $k$-edge term, then we say that a set $S \subseteq \RR$ is a \emph{compact representation} of $\RR$ if:
\begin{itemize}
\item $\Sig(S) = \Sig(\RR)$,

\item for every $I \subseteq \{1, ..., n\}$ with $|I| \le k-1$ we have $\pi_I(S) = \pi_I(\RR)$, and

\item $|S| \le 2|\Sig(\RR)| + \sum_{I \subseteq \{1,...,n\}, |I| \le k-1} |\pi_I(\RR)|$.
\end{itemize}
\end{defn}

\begin{thm}[Subpowers with edge terms are generated by compact representations \cite{few-subpowers}]\label{edge-gen} If $\RR \le \bA_1 \times \cdots \times \bA_n$ and the $\bA_i$ are finite algebras in a variety with a $k$-edge term $e$, then for any compact representation $S$ of $\RR$, we have $\RR = \Sg_{e}(S)$.
\end{thm}
\begin{proof} Let $s,p,d$ be terms as in Theorem \ref{edge-spd}. We induct on $n$. Suppose $a = (a_1, ..., a_n) \in \RR$, then by the induction hypothesis there is $b_n \in \bA_n$ with $(a_1, ..., a_{n-1}, b_n) \in \Sg_e(S)$. Then if we let $d_n = d(b_n,a_n)$ then we see that $a_n,d_n$ is a minority pair and $(a_1, ..., a_n, d_n) \in \RR$, so $(n,a_n,d_n) \in \Sig(\RR)$, and from the definition of a compact representation we see that there must be some $c_1, ..., c_{n-1}$ such that
\[
(c_1, ..., c_{n-1}, a_n), (c_1, ..., c_{n-1}, d_n) \in S.
\]

We show by an inner induction on subsets $I \subseteq \{1, ..., n\}$ that for each $I$, we have $\pi_I(a) \in \pi_I(\Sg_e(S))$. If $|I| \le k-1$ this follows from the definition of a compact representation, while if $n \not\in I$ then this follows from the outer inductive hypothesis. For the sake of notational simplicity we will assume that $I = \{1, ..., n\}$. Then by the inductive hypothesis, there are $b_1, ..., b_{n-1}$ such that for each $i$, we have
\[
(a_1, ..., a_{i-1}, b_i, a_{i+1}, ..., a_n) \in \Sg_e(S).
\]
Then we have
\[
s\left(\begin{bmatrix} a_1 & b_1 & a_1 & \cdots\\ a_2 & a_2 & b_2 & \cdots\\ \vdots & \vdots & \vdots & \ddots\\ b_n & a_n & a_n & \cdots\end{bmatrix}\right) = \begin{bmatrix} a_1\\ a_2\\ \vdots\\ d_n\end{bmatrix} \in \Sg_e(S).
\]
Additionally, we have
\[
p\left(\begin{bmatrix} c_1 & c_1 & b_1\\ c_2 & c_2 & a_2\\ \vdots & \vdots & \vdots\\ d_n & a_n & a_n\end{bmatrix}\right) = \begin{bmatrix} b_1\\ a_2\\ \vdots\\ d_n\end{bmatrix} \in \Sg_e(S).
\]
Now we can apply the $k$-edge term $e$ to see that
\[
e\left(\begin{bmatrix} b_1 & a_1 & b_1 & a_1 & \cdots\\ a_2 & a_2 & a_2 & b_2 & \cdots\\ \vdots & \vdots & \vdots & \vdots & \ddots\\ d_n & d_n & a_n & a_n & \cdots\end{bmatrix}\right) = \begin{bmatrix} a_1\\ a_2\\ \vdots\\ a_n\end{bmatrix} \in \Sg_e(S).\qedhere
\]
\end{proof}

\begin{cor} For a fixed finite algebra $\bA$:
\begin{itemize}
\item $\bA$ has a $k$-edge term but no $k-1$-edge term iff $i_\bA(n) = \Theta(n^{k-1})$, and

\item $\bA$ has no $k$-edge term for any $k$ iff $i_\bA(n) = 2^{\Theta(n)}$.
\end{itemize}
\end{cor}
\begin{proof} We only need to check that if $\bA$ has a $k$-edge term, then $i_\bA(n) = O(n^{k-1})$. Suppose that $a_1, ..., a_m \in \bA^n$ are independent, and consider the relations $\RR_i = \Sg_{\bA^n}\{a_1, ..., a_i\}$. We can easily find a sequence of compact representations $S_1, ..., S_m$ of $\RR_1, ..., \RR_m$ with $S_i \subseteq S_{i+1}$ for each $i$. From the independence of the $a_i$s, we have $\RR_i \ne \RR_{i+1}$ for all $i$, so by induction we see that $|S_i| \ge i$ for all $i$. Then from the fact that $S_m$ is a compact representation, we have
\[
m \le |S_m| \le 2n|\bA|^2 + \sum_{I \subseteq \{1, ..., n\}, |I| \le k-1} |\bA|^{k-1} = O(n^{k-1}).\qedhere
\]
\end{proof}

We can now generalize Dalmau's generalized majority-minority algorithm to an algorithm for computing compact representations of intersections of two relations which are both described by compact representations. The only changes we need to make are to use the edge term $e$ in the \texttt{Nonempty} subroutine in the place of the gmm term $\varphi$, and to modify the \texttt{Fix-values} subroutine to use the ternary term $p$ from Theorem \ref{edge-spd}.

\begin{algorithm}
\caption{\texttt{Fix-values}$(R, a_1, ..., a_m)$, $p,d$ terms as in Theorem \ref{edge-spd}, $R$ a compact representation of $\RR \le \bA_1\times \cdots\times \bA_n$.}
\begin{algorithmic}[1]
\State Set $R_0 \gets R$.
\For{$j$ from $1$ to $m$}
\State Let $R_j \gets \emptyset$.
\ForAll{$I = \{i_1, ...\} \subseteq \{1,...,n\}$ with $|I| < k$ and $(b_{i_1}, ...) \in \pi_I(R_{j-1})$}
\State Set $R_j \gets R_j \cup \texttt{Nonempty}(R_{j-1},j,i_1,...,i_{|I|},\{(a_j,b_{i_1},...,b_{i_{|I|}})\})$.
\EndFor
\ForAll{$(i,a,b) \in \Sig(R_{j-1})$ with $i > j$ and $a,b$ a minority pair (i.e. $d(b,a) = b$)}
\State Let $t_a,t_b \in R_{j-1}$ witness the triple $(i,a,b)$.
\State Let $t \gets \texttt{Nonempty}(R_{j-1},j,i,\{(a_{j},a)\})$.
\If{$t \ne \emptyset$}
\State Set $R_{j} \gets R_{j} \cup \{t, p(t_b,t_a,t)\}$.
\EndIf
\EndFor
\EndFor
\State \Return $R_m$.
\end{algorithmic}
\end{algorithm}

That the modified \texttt{Fix-values} subroutine works follows from the following Proposition.

\begin{prop} If the pair of tuples $t_a,t_b$ witness the minority index $(i,a,b)$, then for any $t$ with $\pi_i(t) = a$ the pair of tuples $t,p(t_b,t_a,t)$ also witnesses the minority index $(i,a,b)$.
\end{prop}
\begin{proof} From the identity $p(y,y,x) \approx x$ we have
\[
\pi_{<i}(p(t_b,t_a,t)) = \pi_{<i}(p(t_a,t_a,t)) = \pi_{<i}(t),
\]
and since $(a,b)$ is a minority pair, we have
\[
\pi_i(p(t_b,t_a,t)) = p(b,a,a) = d(b,a) = b.\qedhere
\]
\end{proof}

\begin{ex}\label{ex-few-subpowers} There is an example of an algebra $\bA = (\{a,b,c\}, g)$ with $g$ a ternary operation such that $\bA$ has a $3$-edge term, but is not in the variety generated by generalized majority-minority algebras of any arity (up to term equivalence). The ternary operation $g$ is the idempotent symmetric function given by
\[
g(a,b,b) = b, g(a,a,b) = a, g(a,c,c) = a, g(a,a,c) = c, g(b,c,c) = a, g(b,b,c) = c, g(a,b,c) = c.
\]
You can understand this as follows: the subset $\{a,b\}$ is a majority subalgebra, the subset $\{a,c\}$ is a pure minority subalgebra, and there is a congruence with equivalence classes $\{a,b\}, \{c\}$ so that the quotient is a pure minority algebra. Also, the only way to get $b$ out of an application of $g$ is if at least two of the inputs are $b$s (this property is called ``absorption'': the subalgebra $\{a,c\}$ absorbs $\{a,b,c\}$ with respect to $g$).

To see that this isn't in the variety generated by generalized majority-minority algebras, recall that in any gmm algebra there are functions $s,p,d$ as in Theorem \ref{edge-spd}, where $d$ satisfies the additional identity $d(x,d(y,x)) \approx x$ since $d$ either acts as first or second projection for any particular pair $x,y$. Since the quotient corresponding to $\{a,b\}, \{c\}$ is a pure minority algebra, we must have $d(c,b) = c$, so by the extra identity we have $d(b,c) = d(b,d(c,b)) = b$. Then the function $p$ would satisfy
\[
p\left(\begin{bmatrix} b & c & c\\ c & c & b\end{bmatrix}\right) = \begin{bmatrix} d(b,c)\\ b\end{bmatrix} \stackrel{?}{=} \begin{bmatrix} b\\ b\end{bmatrix}.
\]
But this is impossible: the subalgebra of $\bA^2$ generated by $(b,c), (c,c), (c,b)$ doesn't contain $(b,b)$, because of the absorption property of $\{a,c\}$ with respect to $g$.

To see that $\bA$ has a $3$-edge term, we define an auxiliary $4$-ary term $f$ by
\[
f(u,x,y,z) = g(g(u,x,z),g(u,y,z),g(u,z,z)),
\]
and then define our $3$-edge term by
\[
e(u,x,y,z) = g(g(f(u,x,y,z),x,x),g(f(u,x,y,z),y,y),g(f(u,x,y,z),z,z)).
\]
If we define functions $s,p,d$ from the $3$-edge term $e$ as in Theorem \ref{edge-spd}, then $d$ is given by
\begin{center}
\begin{tabular}{c|ccc}
$d(x,y)$ & $a$ & $b$ & $c$\\ \hline $a$ & $a$ & $b$ & $a$\\ $b$ & $a$ & $b$ & $a$\\ $c$ & $c$ & $c$ & $c$
\end{tabular}
\end{center}
and the minority pairs are $(a,c), (c,a), (b,c)$. The fact that $(c,b)$ is \emph{not} a minority pair is witnessed by the fact that the relation $\Sg_{\bA^2}\{(b,c),(c,c),(c,b)\}$ contains $(b,c)$ but does not contain $(b,b)$, even though it has $(2,b,c)$ in its signature.

The associated relational clone is generated by the order two automorphism $\{(a,b),(b,a)\}$ of $\{a,b\}$, the partial order $\{(a,a),(a,b),(b,b),(c,c)\}$, the binary relation $\{(a,a),(a,b),(a,c),(b,a),(b,c)\}$ which witnesses the fact that $\{a,c\}$ is a ``central'' subalgebra in Zhuk's terminology \cite{zhuk-dichotomy} (which is closely related to $\{a,c\}$ being a ternary absorbing subalgebra), and the affine ternary relation $\{(a,a,c),(a,c,a),(c,a,a),(c,c,c)\}$.
\end{ex}

For an idempotent algebra $\bA$ with a nontrivial congruence $\theta \in \Con(\bA)$, such as the previous example, we can test whether $\bA$ has few subpowers by checking that $\bA/\theta$ has few subpowers and that each congruence class of $\theta$ has few subpowers separately. This follows from the following easy results from \cite{cube-term-blockers}.

\begin{prop}\label{prop-few-subpowers-composition} Suppose $\bA$ is an idempotent algebra, $\theta \in \Con(\bA)$, and that there are terms $t_1,t_2$ such that $t_1$ acts as a $\Delta_1$-cube term on $\bA/\theta$ and $t_2$ acts as a $\Delta_2$-cube term on each congruence class of $\theta$. Then $t_2*t_1$ is a $\Delta$-cube term for $\bA$, where $\Delta = \{S\times T \mid S \in \Delta_2, T \in \Delta_1\}$.
\end{prop}
\begin{proof}[Proof by example] Suppose that $t_1(x,y,z)$ is a Mal'cev term on $\bA/\theta$ which is idempotent on $\bA$ and that $t_2(x,y,z)$ is a Mal'cev term on each congruence class of $\theta$. Then for any $a,b \in \bA$, if we write
\begin{align*}
c &= t_1(a,a,b) \in b/\theta,\\
d &= t_1(b,a,a) \in b/\theta,
\end{align*}
then we have
\begin{align*}
(t_2*t_1)\left(\begin{bmatrix} a & a & b & a & a & b & b & b & b\\ b & a & a & b & a & a & b & b & b\\ b & b & b & a & a & b & a & a & b\\ b & b & b & b & a & a & b & a & a\end{bmatrix}\right) &\coloneqq t_2\left(\begin{bmatrix} t_1(a,a,b) & t_1(a,a,b) & t_1(b,b,b)\\ t_1(b,a,a) & t_1(b,a,a) & t_1(b,b,b)\\ t_1(b,b,b) & t_1(a,a,b) & t_1(a,a,b)\\ t_1(b,b,b) & t_1(b,a,a) & t_1(b,a,a)\end{bmatrix}\right)\\
&= t_2\left(\begin{bmatrix} c & c & b\\ d & d & b\\ b & c & c\\ b & d & d\end{bmatrix}\right) = \begin{bmatrix} b\\ b\\ b\\ b \end{bmatrix}.\qedhere
\end{align*}
\end{proof}

\begin{cor}\label{cor-few-subpowers-product} If $\bA_1, ..., \bA_n$ are idempotent algebras with the same signature such that each $\bA_i$ has a $\Delta_i$-cube term $t_i$, then $t_1*\cdots*t_n$ is a $\Delta$-cube term for $\bA_1 \times \cdots \times \bA_n$, where $\Delta = \{S_1\times \cdots \times S_n \mid S_i \in \Delta_i\}$.
\end{cor}

\begin{cor} Suppose $\bA$ is a finite idempotent algebra and $\theta \in \Con(\bA)$. Then $\bA$ has few subpowers iff $\bA/\theta$ has few subpowers and each congruence class of $\theta$ has few subpowers.
\end{cor}

\subsection{Some connections with congruence modularity}

\begin{thm} If an algebra has an edge term, then it generates a congruence modular variety.
\end{thm}
\begin{proof} By Theorem \ref{directed-gumm-terms} from Appendix \ref{a-commutator}, we just need to check that an algebra with an edge term has directed Gumm terms, that is, terms $f_1, ..., f_k, p$ satisfying the system of identities
\begin{align*}
f_1(x,x,y) &\approx x,\\
f_i(x,y,x) &\approx x\text{ for all }i,\\
f_i(x,y,y) &\approx f_{i+1}(x,x,y)\text{ for all }i,\\
f_k(x,y,y) &\approx p(x,y,y),\\
p(x,x,y) &\approx y.
\end{align*}
If the reader wants to understand why this system of identities implies congruence modularity \emph{without} reading all of Appendix \ref{a-commutator}, then they can take the following path: first, read the discussion before Theorem \ref{directed-gumm-terms} to see why the existence of directed Gumm terms implies the existence of Gumm terms, then read part of the proof of Theorem \ref{gumm-terms} to see how to construct Day terms from Gumm terms, and finally, read Section \ref{s-shifting} of Appendix \ref{a-commutator} to see why the existence of Day terms is equivalent to congruence modularity.

Suppose that $e$ is a $k$-edge term. Define terms $f_i(x,y,z)$ for $i < k$ by
\[
f_i(x,y,z) = e(x,...,x,y,z,...,z),
\]
such that there are $i-1$ $z$s, a single $y$, and $k+1-i$ $x$s. Then we have
\[
f_1(x,x,y) = e(x,...,x,x) = x,
\]
and for $i < k$ we have
\[
f_i(x,y,x) = e(x,...,x,y,x,...,x) = x.
\]
From the construction of the $f_i$s we have
\[
f_i(x,y,y) = e(x,...,x,y,y,...,y) = f_{i+1}(x,x,y)
\]
for $i+1 < k$. Finally, if we define $f_k(x,y,z)$ by
\[
f_k(x,y,z) = e(y,x,y,z,...,z)
\]
and $p(x,y,z)$ by
\[
p(x,y,z) = e(y,x,z,z,...,z),
\]
then
\[
f_{k-1}(x,y,y) = e(x,x,x,y,...,y) = f_k(x,x,y)
\]
and
\[
f_k(x,y,x) = e(y,x,y,x,...,x) = x
\]
by the $k$-edge identities, while
\[
f_k(x,y,y) = e(y,x,y,y,...,y) = p(x,y,y)
\]
and
\[
p(x,x,y) = e(x,x,y,y,...,y) = y
\]
by the $k$-edge identities again. Thus $f_1, ..., f_k, p$ are a sequence of directed Gumm terms.
\end{proof}

\begin{thm}\label{cd-few-subpowers-nu} For $k \ge 3$, an algebra has a $k$-edge term and generates a congruence distributive variety iff it has a $k$-ary near-unanimity term.
\end{thm}
\begin{proof} First the easy direction. If an algebra $\bA$ has a $k$-ary near-unanimity term $t$, then adding an extra variable at the beginning of $t$ produces a $k$-edge term. Additionally, the discussion before Theorem \ref{directed-gumm-terms} shows that we can construct a sequence of J\'onsson terms from $t$, and then Theorem \ref{jonsson-terms} shows that $\bA$ generates a congruence distributive variety.

Now the harder direction: assume that $\bA$ generates a congruence distributive variety and has a $k$-edge term $e$. By Theorem \ref{directed-gumm-terms}, there is a sequence of directed J\'onsson terms $f_1, ..., f_m$, that is, a sequence satisfying the system of identities
\begin{align*}
f_1(x,x,y) &\approx x,\\
f_i(x,y,x) &\approx x\text{ for all }i,\\
f_i(x,y,y) &\approx f_{i+1}(x,x,y)\text{ for all }i,\\
f_m(x,y,y) &\approx y.
\end{align*}

Let $\cF = \cF_{\cV(\bA)}(x,y) \le \bA^{\bA^2}$ be the free algebra on two generators in the variety generated by $\bA$. Let $\bS \le \cF^k$ be generated by the vectors $(x,...,x,y,x,...,x)$ with all but one entry equal to $x$ and the remaining entry equal to $y$. Note that $\bS$ is symmetric under permuting its coordinates. We just need to prove that $(x,...,x) \in \bS$.

{\bf Claim:} For all $i$, we have $(f_i(y,x,x),x,...,x) \in \bS$.

{\bf Proof of Claim:} We induct on $i$, taking $(y,x,...,x) \in \bS$ as our base case. By the induction hypothesis, we have
\[
(f_i(y,y,x),x,...,x) = (f_{i-1}(y,x,x),x,...,x) \in \bS.
\]
Additionally, the tuples
\[
\begin{bmatrix} f_i(y,x,x)\\ f_i(x,x,y)\\ x\\ \vdots \\ x\end{bmatrix} = f_i\left(\begin{bmatrix} y & x & x\\ x & x & y\\ x & y & x\\ \vdots & \vdots & \vdots \\ x & x & x\end{bmatrix}\right)
\]
and
\[
\begin{bmatrix} f_i(y,y,x)\\ f_i(x,x,y)\\ x\\ \vdots \\ x\end{bmatrix} = f_i\left(\begin{bmatrix} y & y & x\\ x & x & y\\ x & x & x\\ \vdots & \vdots & \vdots \\ x & x & x\end{bmatrix}\right)
\]
are both in $\bS$. Now we apply the $k$-edge term $e$:
\[
e\left(\begin{bmatrix} f_i(y,y,x) & f_i(y,y,x) & f_i(y,x,x) & f_i(y,x,x) & \cdots & f_i(y,x,x)\\ f_i(x,x,y) & x & f_i(x,x,y) & x & \cdots & x\\ x & x & x & f_i(x,x,y) & \cdots & x\\ \vdots & \vdots & \vdots & \vdots & \ddots & \vdots \\ x & x & x & x & \cdots & f_i(x,x,y)\end{bmatrix}\right) = \begin{bmatrix} f_i(y,x,x)\\ x\\ x\\ \vdots \\ x\end{bmatrix}.
\]

To finish the proof, we apply the Claim with $i = m$ to see that $(x,x,...,x) = (f_m(y,x,x),x,...,x) \in \bS$.
\end{proof}

\begin{ex} We give an example of a congruence distributive algebra without few subpowers. Recall from Example \ref{ex-strict-width-n} that for each $n$, the relational structure $(\{0,1\}, \{0\}, \le, \{0,1\}^n\setminus\{(0,...,0)\})$ has strict width exactly $n$. The limiting relational clone on $\{0,1\}$ generated by the relations $\{0\}, \le$, and $\{0,1\}^n\setminus\{(0,...,0)\}$ for all $n \in \NN$ corresponds to the clone generated by the ternary operation
\[
f(x,y,z) = x \vee (y \wedge z).
\]
Since the $n$-ary critical relation $\{0,1\}^n\setminus\{(0,...,0)\}$ doesn't have the parallelogram property and is preserved by $f$ for all $n$, the clone generated by $f$ can't have few subpowers by Theorem \ref{parallelogram-critical}.

To check that the algebra $\bA = (\{0,1\}, x \vee (y \wedge z))$ generates a congruence distributive variety, consider the sequence of ternary terms given by
\[
f_1(x,y,z) = x\vee (y\wedge z), \;\;\; f_2(x,y,z) = (x\wedge y) \vee z.
\]
To see that this is a sequence of directed J\'onsson terms, note that they satisfy $f_i(x,y,x) = x \vee (y \wedge x) = x$, are connected by
\[
f_1(x,y,y) = x \vee (y \wedge y) = x \vee y = (x \wedge x) \vee y = f_2(x,x,y),
\]
and have $f_1(x,x,y) = x, f_2(x,y,y) = y$. By Theorem \ref{jonsson-terms} and the discussion before Theorem \ref{directed-gumm-terms}, this implies that $\bA$ is congruence distributive.
\end{ex}

\begin{ex}\label{ex-con-semi-square} We've seen earlier that the two-element semilattice $\bA = (\{0,1\}, \max)$ does not have few subpowers. Here we will check that the two-element semilattice does not generate a congruence modular variety. In fact, the congruence lattice $\Con(\bA^2)$ already fails to be modular. It turns out that every congruence on $\bA^2$ is generated (as a congruence) by just one pair of elements $a,b$ of $\bA^2$, so we can label the nontrivial congruences on $\bA^2$ by pairs of elements $a,b \in \bA^2$, yielding the following congruence lattice.
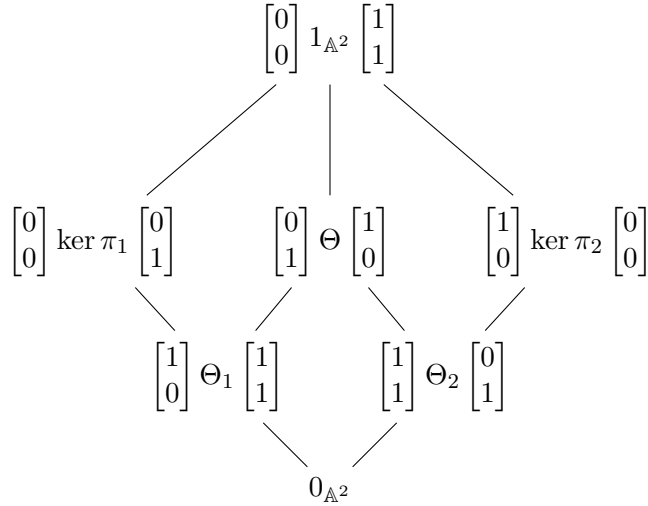
\begin{center}
\begin{tikzpicture}[scale=1.5]
  \node (1) at (0,2) {$\begin{bmatrix}0\\ 0\end{bmatrix} 1_{\bA^2} \begin{bmatrix}1\\ 1\end{bmatrix}$};
  \node (p1) at (-2.1,0.2) {$\begin{bmatrix}0\\ 0\end{bmatrix} \ker \pi_1 \begin{bmatrix}0\\ 1\end{bmatrix}$};
  \node (t) at (0,0.2) {$\begin{bmatrix}0\\ 1\end{bmatrix} \Theta \begin{bmatrix}1\\ 0\end{bmatrix}$};
  \node (p2) at (2.1,0.2) {$\begin{bmatrix}1\\ 0\end{bmatrix} \ker \pi_2 \begin{bmatrix}0\\ 0\end{bmatrix}$};
  \node (t1) at (-1,-1) {$\begin{bmatrix}1\\ 0\end{bmatrix} \Theta_1 \begin{bmatrix}1\\ 1\end{bmatrix}$};
  \node (t2) at (1,-1) {$\begin{bmatrix}1\\ 1\end{bmatrix} \Theta_2 \begin{bmatrix}0\\ 1\end{bmatrix}$};
  \node (0) at (0,-2) {$0_{\bA^2}$};
  \draw (1) -- (p1) -- (t1) -- (0) -- (t2) -- (p2) -- (1);
  \draw (1) -- (t);
  \draw (t1) -- (t) -- (t2);
\end{tikzpicture}
\end{center}
To see that this isn't modular, note that the sublattice generated by $\ker \pi_1, \ker \pi_2, \Theta_2$ is isomorphic to the pentagon lattice $\cN_5$. Considered as an abstract lattice, $\Con(\bA^2)$ is the standard example of a lattice which is meet-semidistributive (recall from Example \ref{semi-sd-meet} and Proposition \ref{sd-meet-commutator} that the variety of semilattices is SD($\wedge$)) but not join-semidistributive (we have $\Theta \vee \ker \pi_1 = \Theta \vee \ker \pi_2 = 1_{\bA^2}$, but $\Theta \vee (\ker \pi_1 \wedge \ker \pi_2) = \Theta \ne 1_{\bA^2}$).
\end{ex}

Although congruence modularity is slightly weaker than having few subpowers, the concepts are quite close. One hint at the connection between them comes from counting \emph{congruences} on subpowers of $\bA$.

\begin{defn} If $\bA$ is an algebra, then we define the function $c_\bA(n)$ to be the base-$2$ logarithm of the maximum size of $\Con(\RR)$ over all $\RR \le \bA^n$.
\end{defn}

\begin{prop} A variety $\cV$ is congruence distributive iff for all subdirect products $\RR \le_{sd} \bA_1 \times \cdots \times \bA_n$ in $\cV$, every congruence on $\RR$ can be written as a product of congruences on the $\bA_i$s.
\end{prop}
\begin{proof} Suppose first that $\cV$ is congruence distributive. Then for any congruence $\theta$ on $\RR$, by distributivity and $\bigwedge_i \ker \pi_i = 0_\RR$ we have
\[
\bigwedge_i (\theta \vee \ker \pi_i) = \theta \vee \bigwedge_i \ker \pi_i = \theta \vee 0_\RR = \theta,
\]
so $\theta$ is the product of the congruences $\pi_i(\theta \vee \ker \pi_i) \in \Con(\bA_i)$.

Conversely, suppose that $\bA \in \cV$, and suppose that $\alpha, \beta, \gamma \in \Con(\bA)$. Then $\bA/(\beta \wedge \gamma)$ is a subdirect product of $\bA/\beta$ and $\bA/\gamma$, so the congruence
\[
\alpha \vee (\beta \wedge \gamma),
\]
considered as a congruence on $\bA/(\beta \wedge \gamma)$, is a product congruence iff it is equal to
\[
(\alpha \vee \beta) \wedge (\alpha \vee \gamma).\qedhere
\]
\end{proof}

\begin{cor} If $\cV(\bA)$ is congruence distributive, then $c_\bA(n) = nc_\bA(1)$.
\end{cor}

If a variety is congruence modular but \emph{not} congruence distributive, then it necessarily contains a (finitely generated) nontrivial affine algebra. So we need to understand $c_\bA(n)$ for $\bA$ a finite affine algebra, and since the congruence lattice only depends on the polynomial clone, we may assume that $\bA$ is a module over a ring. In this case, there is a bijection between congruences on $\bA^n$ and submodules of $\bA^n$.

\begin{prop} If $\bA$ is a nontrivial finite module over a ring, then $c_\bA(n) \ge \frac{n^2-1}{4}$.
\end{prop}
\begin{proof} We may as well assume that $\bA$ is simple. Let $c$ be any nonzero element of $\bA$. For $n = 2m$, the span of the columns of the $n\times m$ matrix $\begin{bmatrix} cI\\ M\end{bmatrix}$ completely determines the $m \times m$ matrix $M$, so $c_\bA(2m) \ge m^2\log_2(|\bA|) \ge m^2$.
\end{proof}

\begin{cor} If $\bA$ is finite and $\cV(\bA)$ is congruence modular but not congruence distributive, then $c_{\bA}(n) = \Omega(n^2)$.
\end{cor}

How can we get an upper bound on $c_\bA(n)$ when $\bA$ is congruence modular? The trick is to use the fact that in modular lattices, the \emph{height} of the lattice is well-behaved. We can relate the height of a congruence lattice to its size using the following elementary bound.

\begin{prop} If $\bA$ is a finite algebra such that $\Con(\bA)$ has height $h$, then
\[
|\Con(\bA)| \le \sum_{i = 0}^h \tbinom{|\bA|}{2}^i \le |\bA|^{2h}.
\]
\end{prop}
\begin{proof} Consider any congruence $\alpha \in \Con(\bA)$. Since every cover of $\alpha$ is generated (as a congruence) by $\alpha$ together with some pair $(a,b) \not\in \alpha$, the number of covers of $\alpha$ is bounded by $\binom{|\bA|}{2}$. Since every element of $\Con(\bA)$ can be reached from $0_\bA$ by repeatedly choosing covers at most $h$ times, we get the stated bound on $|\Con(\bA)|$.

We can get a slightly better bound as follows: the above argument shows that every congruence can be generated (as a congruence) by at most $h$ pairs in $\binom{\bA}{2}$. Additionally, there is only one congruence at height $h$, since $\Con(\bA)$ has a top element $1_\bA$. So we have
\[
|\Con(\bA)| \le 1 + \sum_{i = 0}^{h-1} \binom{\binom{|\bA|}{2}}{i}.\qedhere
\]
\end{proof}

\begin{cor} If $\bA$ is finite and generates a congruence modular variety, then $c_{\bA}(n) \le n^2 \cdot 2|\bA|\log_2(|\bA|)$.
\end{cor}
\begin{proof} Let $c$ be the maximum height of $\Con(\bB)$ over all subalgebras $\bB \le \bA$ ($c$ is automatically bounded by $|\bA|$). We claim that for any $\RR \le \bA^n$, the height of $\Con(\RR)$ is bounded by $cn$. Since $\Con(\RR)$ is modular, we can compute its height by looking at the size of \emph{any} maximal chain in $\Con(\RR)$.

We will choose our maximal chain to be any maximal extension of the chain
\[
0_\RR \le \ker \pi_{[n-1]} \le \cdots \le \ker \pi_{[2]} \le \ker \pi_1 \le 1_\RR.
\]
By the Diamond Isomorphism Theorem \ref{diamond-isom}, the interval $\llbracket \ker \pi_{[i]}, \ker \pi_{[i-1]} \rrbracket$ is isomorphic to the interval $\llbracket \ker \pi_i, \ker \pi_{[i-1]}\vee \ker \pi_i\rrbracket$, so its height is bounded by the height of the interval $\llbracket \ker \pi_i, 1_\RR \rrbracket$, which is isomorphic to $\Con(\RR/\ker \pi_i)$. Since $\RR/\ker \pi_i \cong \pi_i(\RR) \le \bA$, the height of $\Con(\RR/\ker \pi_i)$ is bounded by $c$, and putting these intervals together we see that the height of $\Con(\RR)$ is bounded by $cn$.

Using the previous bound, we get
\[
\log_2(|\Con(\RR)|) \le \log_2(|\RR|^{2cn}) \le 2cn \log_2(|\bA|^n) = 2cn^2 \log_2(|\bA|).\qedhere
\]
\end{proof}

\begin{thm}[Few congruences on subpowers iff congruence modular \cite{few-subpowers}] Let $\bA$ be a finite algebra with at least two elements, and let $\cV(\bA)$ be the variety it generates.
\begin{itemize}
\item If $\cV(\bA)$ is congruence distributive, then $c_\bA(n) = \Theta(n)$.
\item If $\cV(\bA)$ is congruence modular but not congruence distributive, then $c_\bA(n) = \Theta(n^2)$.
\item If $\cV(\bA)$ is not congruence modular, then $c_\bA(n) = 2^{\Theta(n)}$.
\end{itemize}
\end{thm}
\begin{proof} By the previous results, all we need to check is that if $\cV(\bA)$ is not congruence modular, then $c_\bA(n) = 2^{\Omega(n)}$. Let $\bF = \cF_{\cV(\bA)}(x,y,z,w) \le \bA^{\bA^4}$ be the free algebra on four generators. We will show that if $c_\bF(2n) < 2^n$ for any $n$, then $\bA$ has Day terms, and is therefore congruence modular by Appendix \ref{s-shifting}.

Define congruences on $\bF$ as in Corollary \ref{day-terms}: let $\theta_{ab}$ be the congruence generated by the pair $(a,b)$ for any pair of variables $a,b$, set $\alpha = \theta_{xy} \vee \theta_{zw}, \beta = \theta_{xz} \vee \theta_{yw}$, and $\gamma = (\alpha \wedge \beta) \vee \theta_{zw}$. This is the generic Shifting Lemma configuration:
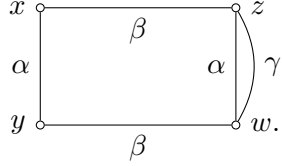
\begin{center}
\begin{tikzpicture}[scale=1.3]
\node[circle, minimum width=3pt, draw, inner sep=0pt, label=left:{$x$}] (a) at (0,1.2){};
\node[circle, minimum width=3pt, draw, inner sep=0pt, label=right:{$z$}] (c) at (2.0,1.2){};
\node[circle, minimum width=3pt, draw, inner sep=0pt, label=left:{$y$}] (b) at (0,0){};
\node[circle, minimum width=3pt, draw, inner sep=0pt, label=right:{$w.$}] (d) at (2.0,0){};
\draw (a) to ["$\beta$"'] (c) to ["$\alpha$"'] (d) to ["$\beta$"] (b) to ["$\alpha$"] (a);
\draw [bend left] (c) to ["$\gamma$"] (d);
\end{tikzpicture}
\end{center}
To show the existence of Day terms, we just need to show that $(x,y) \in \gamma$.

Pick an $n$ such that $c_{\bF}(2n) < 2^n$, and consider the subalgebra $\RR \le \bF^{2n}$ consisting of tuples such that every pair of coordinates are related by $\beta$ (it helps to imagine elements of $\RR$ written out horizontally as row vectors, following the convention that variables which are related by $\beta$ are laid out on horizontal lines). We will define a family of $2^n$ pairs of elements of $\RR$ as follows.

First, we define elements $x^0, x^1, y^0, y^1 \in \bF^2$ by $x^0 = (x,z), x^1 = (z,x)$ and similarly $y^0 = (y,w), y^1 = (w,y)$. Then, for any $i = (i_1, ..., i_n) \in \{0,1\}^n$, we define $f_i, g_i \in \RR$ by
\begin{align*}
f_i &= (x^{i_1}, ..., x^{i_n}),\\
g_i &= (y^{i_1}, ..., y^{i_n}).
\end{align*}
For each $i \in \{0,1\}^n$, we define a congruence $\Theta(i)$ to be the congruence of $\RR$ generated by the pair $(f_i, g_i)$. Since $c_{\bF}(2n) < 2^n$, there must be some $i \in \{0,1\}^n$ such that
\[
\Theta(i) \le \bigvee_{j \ne i} \Theta(j),
\]
and by permuting the coordinates of $\RR$, we see that in fact this must hold for every $i$, and in particular for $i = (0, ..., 0)$. By dropping half of the coordinates of $\RR$ to get a similar algebra $\RR' \le \bF^n$ such that $f_0$ becomes the vector $f_0' = (x, ..., x)$ and $g_0$ becomes the vector $g_0' = (y, ..., y)$, and defining elements $f_j', g_j'$ by dropping half the coordinates of $f_j, g_j$, we see that
\[
(f_0', g_0') \in \bigvee_{j \ne (0, ..., 0)} \Theta'(j),
\]
where $\Theta'(j)$ is the congruence of $\RR'$ generated by the pair $(f_j', g_j')$.

Each $\Theta'(j)$ has the following property: if $(a,b) \in \Theta'(j)$ and every pair of coordinates of $a$ are related by $\alpha \wedge \beta$, then every pair of coordinates of $b$ are also related by $\alpha \wedge \beta$. To see this, just note that for each coordinate $i \le n$ we have $(a_i,b_i) \in \alpha$, since this holds in the case where $(a,b) = (f_j', g_j')$.

For $j \ne (0, ..., 0)$, $\Theta'(j)$ has the following additional property: there exists some coordinate $i \le n$ such that if $(a,b) \in \Theta'(j)$, then $(a_i,b_i) \in \gamma$. In fact, we can take the coordinate $i$ to be the first coordinate of $j$ such that $j_i = 1$, and note that the $i$th coordinates of $f_j', g_j'$ are $z, w$ respectively, with $(z,w) \in \gamma$ by the definition of $\gamma$.

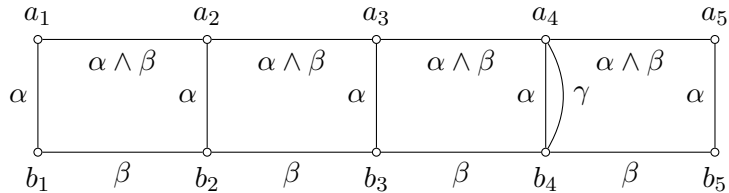
\begin{center}
\begin{tikzpicture}[scale=1.5]
\node[circle, minimum width=3pt, draw, inner sep=0pt, label=above:{$a_1$}] (a1) at (0,1){};
\node[circle, minimum width=3pt, draw, inner sep=0pt, label=above:{$a_2$}] (a2) at (1.5,1){};
\node[circle, minimum width=3pt, draw, inner sep=0pt, label=above:{$a_3$}] (a3) at (3,1){};
\node[circle, minimum width=3pt, draw, inner sep=0pt, label=above:{$a_4$}] (a4) at (4.5,1){};
\node[circle, minimum width=3pt, draw, inner sep=0pt, label=above:{$a_5$}] (a5) at (6,1){};
\node[circle, minimum width=3pt, draw, inner sep=0pt, label=below:{$b_1$}] (b1) at (0,0){};
\node[circle, minimum width=3pt, draw, inner sep=0pt, label=below:{$b_2$}] (b2) at (1.5,0){};
\node[circle, minimum width=3pt, draw, inner sep=0pt, label=below:{$b_3$}] (b3) at (3,0){};
\node[circle, minimum width=3pt, draw, inner sep=0pt, label=below:{$b_4$}] (b4) at (4.5,0){};
\node[circle, minimum width=3pt, draw, inner sep=0pt, label=below:{$b_5$}] (b5) at (6,0){};
\draw (a1) to ["$\alpha\wedge \beta$"'] (a2) to ["$\alpha\wedge \beta$"'] (a3) to ["$\alpha\wedge \beta$"'] (a4) to ["$\alpha\wedge \beta$"'] (a5);
\draw (b1) to ["$\beta$"'] (b2) to ["$\beta$"'] (b3) to ["$\beta$"'] (b4) to ["$\beta$"'] (b5);
\draw (a1) to ["$\alpha$"'] (b1);
\draw (a2) to ["$\alpha$"'] (b2);
\draw (a3) to ["$\alpha$"'] (b3);
\draw (a4) to ["$\alpha$"'] (b4);
\draw (a5) to ["$\alpha$"'] (b5);
\draw [bend left] (a4) to ["$\gamma$"] (b4);
\end{tikzpicture}
\end{center}

Putting the above properties together, and using $\alpha \wedge \beta \le \gamma$, we see that $(f_0', b) \in \bigvee_{j \ne (0, ..., 0)} \Theta'(j)$ implies that every coordinate of $f_0'$ is congruent modulo $\gamma$ to every coordinate of $b$, and taking $b = g_0'$ we see that $(x,y) \in \gamma$, which completes the proof.
\end{proof}

\begin{ex} Consider the two-element semilattice $\bA = (\{0,1\},\max)$ once again. In this case, we can check directly that $c_{\bA}(n) \ge \binom{n}{n/2}$. To see this, note that for every nonempty upwards closed subset $U \le \bA^n$, there is a congruence $\theta_U$ which collapses all elements of $U$ into a single top element of $\bA^n/\theta_U$, and which does not identify any pair of elements $a \ne b$ such that $\{a,b\} \not\subseteq U$. In other words, $\theta_U = U^2 \cup \Delta_{\bA^n}$.

We just need to check that the number of distinct nonempty upwards closed subsets $U$ of $\{0,1\}^n$ is at least $2^{\binom{n}{n/2}}$: for this, note that upwards closed sets $U$ are in a one-to-one correspondence with antichains (every upwards closed set $U$ is determined by its antichain of minimal elements), and every set of elements of $\bA^n$ which each have exactly $n/2$ coordinates equal to $1$ forms an antichain.
\end{ex}

\section{Parallelogram terms}

Examining the proof of Theorem \ref{edge-gen}, we can extract useful terms known as \emph{parallelogram terms}, which we can use to give a better description of the relational clone corresponding to an algebra with few subpowers.

\begin{defn} If $k = m+n$, then an $m,n$-\emph{parallelogram term} is a $k+3$-ary term $r$ which satisfies the identities
\[
r\left(\begin{bmatrix} y & y & x & z & \cdots & x & x & \cdots & x\\ \vdots & \vdots & \vdots & \vdots & \ddots & \vdots & \vdots & \ddots & \vdots\\ y & y & x & x & \cdots & z & x & \cdots & x\\ x & y & y & x & \cdots & x & z & \cdots & x\\ \vdots & \vdots & \vdots & \vdots & \ddots & \vdots & \vdots & \ddots & \vdots\\ x & y & y & x & \cdots & x & x & \cdots & z \end{bmatrix}\right) = \begin{bmatrix} x\\ \vdots\\ x\\ x\\ \vdots\\ x \end{bmatrix},
\]
where the upper left $m\times 3$ block has all rows given by $y, y, x$, the lower left $n \times 3$ block has all rows given by $x,y,y$, and the right $k\times k$ block has $z$s on the diagonal and $x$s elsewhere.
\end{defn}

\begin{thm}[Edge term implies parallelogram terms \cite{parallelogram-terms}] For any $m,n > 0$ with $m+n = k$, a variety has a $k$-edge term $e$ iff it has an $m,n$-parallelogram term $r$.
\end{thm}
\begin{proof} It's clear that every $m,n$-parallelogram term is a $\Delta$-cube term for
\[
\Delta = \{\{1, ..., m\}, \{1, ..., k\}, \{m+1, ..., m+n\}, \{1\}, ..., \{k\}\},
\]
so by Theorem \ref{cube-edge} if $\bA$ has a parallelogram term then it has an edge term.

Now suppose that $e$ is a $k$-edge term. We will build $m,n$-parallelogram terms $r_m$ by induction on $m$. For $m = 1$, we need to show that the vector in $\cF(x,y,z)^k$ of all $x$s is in the subalgebra generated by the columns of the matrix defining a $1,k-1$-parallelogram term. These vectors are the vectors where all entries other than one are $x$s and the last is a $z$, the vector of all $y$s, and the vectors $(x,y,...,y), (y,x,...,x)$.

Letting $d = d(y,x) = e(x,y,x,...,x)$, we have
\[
e\left(\begin{bmatrix} x & y & x & x & \cdots & x\\ z & x & z & x & \cdots & x\\ x & x & x & z & \cdots & x\\ \vdots & \vdots & \vdots & \vdots & \ddots & \vdots\\ x & x & x & x & \cdots & z \end{bmatrix}\right) = \begin{bmatrix} d\\ x\\ x\\ \vdots \\ x \end{bmatrix}
\]
and
\[
e\left(\begin{bmatrix} x & y & x & x & \cdots & x\\ y & y & z & z & \cdots & z\\ y & y & x & x & \cdots & x\\ \vdots & \vdots & \vdots & \vdots & \ddots & \vdots\\ y & y & x & x & \cdots & x \end{bmatrix}\right) = \begin{bmatrix} d\\ z\\ x\\ \vdots \\ x \end{bmatrix},
\]
so the vectors $(d,x,x,...,x)$ and $(d,z,x,...,x)$ are in the subalgebra of $\cF(x,y,z)^k$ generated by the columns of the matrix defining a $1,k-1$-parallelogram term. Note that the previous two applications of the edge term $e$ correspond to applications of the terms
\[
s(x_1, ..., x_k) = e(x_2, x_1, x_2, ..., x_k) \text{ and } p(x,y,z) = e(y,x,z,...,z)
\]
which act like near-unanimity and Mal'cev terms, respectively. To get the vector of all $x$s, we apply $e$ one more time:
\[
e\left(\begin{bmatrix} d & d & x & x & \cdots & x\\ z & x & z & x & \cdots & x\\ x & x & x & z & \cdots & x\\ \vdots & \vdots & \vdots & \vdots & \ddots & \vdots\\ x & x & x & x & \cdots & z \end{bmatrix}\right) = \begin{bmatrix} x\\ x\\ x\\ \vdots \\ x \end{bmatrix}.
\]
Explicitly, our $1,k-1$-parallelogram term $r_1$ is defined from the edge term $e$ by
\begin{align*}
r_1(x,y,z,u_1,...,u_k) &= e(p(y,z,u_2),s(x,u_2,...,u_k),u_2,...,u_k)\\
&= e(e(z,y,u_2,...,u_2),e(u_2,x,u_2,...,u_k),u_2,...,u_k).
\end{align*}

For $m > 1$, we construct the $m,k-m$-parallelogram term $r_m$ using the previous term $r_{m-1}$. Here we focus on the $m$th rows of our matrices. Let
\[
a = r_{m-1}(y,y,x,x,...,x,z,x,...,x),
\]
where the $z$ occurs in the $m+3$rd entry. We want to construct tuples $(x,...,x,a,x,...,x)$ and $(x,...,x,a,y,...,y)$ from the columns of the defining matrix for an $m,k-m$-parallelogram term. We construct these tuples via
\[
r_{m-1}\left(\begin{bmatrix} y & y & x & z & \cdots & x & x & x & \cdots & x\\ \vdots & \vdots & \vdots & \vdots & \ddots & \vdots & \vdots & \vdots & \ddots & \vdots \\ y & y & x & x & \cdots & z & x & x & \cdots & x\\ y & y & x & x & \cdots & x & z & x & \cdots & x\\ x & y & y & x & \cdots & x & x & z & \cdots & x\\ \vdots & \vdots & \vdots & \vdots & \ddots & \vdots & \vdots & \vdots & \ddots & \vdots \\ x & y & y & x & \cdots & x & x & x & \cdots & z\end{bmatrix}\right) = \begin{bmatrix} x\\ \vdots \\ x\\ a\\ x\\ \vdots \\ x\end{bmatrix}
\]
and
\[
r_{m-1}\left(\begin{bmatrix} y & y & x & x & \cdots & x & x & x & \cdots & x\\ \vdots & \vdots & \vdots & \vdots & \ddots & \vdots & \vdots & \vdots & \ddots & \vdots \\ y & y & x & x & \cdots & x & x & x & \cdots & x\\ y & y & x & x & \cdots & x & z & x & \cdots & x\\ y & y & y & y & \cdots & y & x & y & \cdots & y\\ \vdots & \vdots & \vdots & \vdots & \ddots & \vdots & \vdots & \vdots & \ddots & \vdots \\ y & y & y & y & \cdots & y & x & y & \cdots & y\end{bmatrix}\right) = \begin{bmatrix} x\\ \vdots \\ x\\ a\\ y\\ \vdots \\ y\end{bmatrix}.
\]
To get to the vector of all $x$s, we use
\[
r_{m-1}\left(\begin{bmatrix} x & x & x & x & x & \cdots & z & x & \cdots & x\\ \vdots & \vdots & \vdots & \vdots & \vdots & \iddots & \vdots & \vdots & \ddots & \vdots \\ x & x & x & x & z & \cdots & x & x & \cdots & x\\ a & a & x & z & x & \cdots & x & x & \cdots & x\\ x & y & y & x & x & \cdots & x & z & \cdots & x\\ \vdots & \vdots & \vdots & \vdots & \vdots & \iddots & \vdots & \vdots & \ddots & \vdots \\ x & y & y & x & x & \cdots & x & x & \cdots & z\end{bmatrix}\right) = \begin{bmatrix} x\\ \vdots \\ x\\ x\\ x\\ \vdots \\ x\end{bmatrix},
\]
where the middle row works out because $m > 1$. Explicitly, $r_{m}$ is defined in terms of $r_{m-1}$ by
\[
r_{m}(x,y,z,u_1,...,u_k) = r_{m-1}(r_{m-1}(x,y,z,u_1,...,u_k),r_{m-1}(y,y,z,z,...,u_{m},...,z),z,u_{m},...,u_1,u_{m+1},...,u_k).\qedhere
\]
\end{proof}

To understand what parallelogram terms tell us, it is necessary to restrict to certain special relations, known as \emph{critical} relations.

\begin{defn} A subalgebra $\RR \le \bA_1 \times \cdots \times \bA_n$ is \emph{critical} if it is $\cap$-irreducible, that is, if it can't be written as an intersection of strictly larger subalgebras, and if furthermore the relation $\RR$ has no dummy variables (that is, it depends on all of its inputs).
\end{defn}

A standard result in the theory of algebraic lattices (Proposition \ref{meet-irreducible-rep} from Appendix \ref{a-commutator}) shows that every relation can be written as an intersection of critical relations (possibly of lower arity). The following result shows that every relation in an algebra with $k$-parallelogram terms can be written as an intersection of relations of arity less than $k$ and relations with the parallelogram property.

\begin{thm}[Parallelogram terms constrain critical relations \cite{parallelogram-terms}]\label{parallelogram-critical} A variety $\cV$ has $k$-parallelogram terms iff for all critical $\RR \le \bA_1 \times \cdots \times \bA_n$ with $\bA_i \in \cV$, either $n < k$ or $\RR$ has the parallelogram property.
\end{thm}
\begin{proof} First suppose that $\cV$ has $k$-parallelogram terms, and let $\RR \le \bA_1 \times \cdots \times \bA_n$ be a critical relation. Let $\RR^*$ be the cover of $\RR$, i.e., $\RR^*$ is the intersection of all relations which properly contain $\RR$, and let $a = (a_1, ..., a_n) \in \RR^*\setminus \RR$. Then a relation $\bS$ which contains $\RR$ will properly contain $\RR$ iff $\bS$ contains $a$. Following Zhuk \cite{zhuk-key}, we call $a$ a \emph{key tuple} for the critical relation $\RR$.

Since $\RR$ is critical, $\RR$ is properly contained in its existential projections onto any proper subset of the coordinates $1, ..., n$. Thus, there must exist elements $b_1, ..., b_n$ such that the tuples $(b_1,a_2,...,a_n), (a_1,b_2,...,a_n), ..., (a_1, a_2, ..., b_n)$ are all in $\RR$.

Now suppose, for contradiction, that $n \ge k$ and that $\RR$ does not have the parallelogram property when considered as a binary relation on $(\bA_1 \times \cdots \times \bA_i)\times (\bA_{i+1} \times \cdots \times \bA_n)$. Then the are $x_1, ..., x_n, y_1, ..., y_n$ such that the three tuples $(y_1, ..., y_n), (y_1, ..., y_i, x_{i+1}, ..., x_n), (x_1, ..., x_i, y_{i+1}, ..., y_n)$ are in $\RR$, but $(x_1, ..., x_n)$ is not in $\RR$. Since $x = (x_1, ..., x_n)$ is not in $\RR$, the subalgebra generated by $\RR \cup \{x\}$ must properly contain $\RR$, so
\[
a \in \Sg(\RR \cup \{x\}).
\]
Thus there are tuples $c^1, ..., c^m \in \RR$ and an $m+1$-ary term $t$ such that
\[
t(x,c^1, ..., c^m) = a.
\]
Defining a tuple $d$ by
\[
t(y,c^1,..., c^m) = d,
\]
we see that the three tuples $(d_1, ..., d_n), (d_1, ..., d_i, a_{i+1}, ..., a_n), (a_1, ..., a_i, d_{i+1}, ..., d_n)$ are all in $\RR$. But then we can use an $i, n-i$-parallelogram term $r$ (which exists because $n \ge k$) to see that
\[
\begin{bmatrix} a_1\\ \vdots\\ a_i\\ a_{i+1}\\ \vdots\\ a_n \end{bmatrix} = r\left(\begin{bmatrix} d_1 & d_1 & a_1 & b_1 & \cdots & a_1 & a_1 & \cdots & a_1\\ \vdots & \vdots & \vdots & \vdots & \ddots & \vdots & \vdots & \ddots & \vdots\\ d_i & d_i & a_i & a_i & \cdots & b_i & a_i & \cdots & a_i\\ a_{i+1} & d_{i+1} & d_{i+1} & a_{i+1} & \cdots & a_{i+1} & b_{i+1} & \cdots & a_{i+1}\\ \vdots & \vdots & \vdots & \vdots & \ddots & \vdots & \vdots & \ddots & \vdots\\ a_n & d_n & d_n & a_n & \cdots & a_n & a_n & \cdots & b_n \end{bmatrix}\right) \in \RR,
\]
contradicting the assumption that $a \not\in \RR$.

For the converse direction, suppose that $\cV$ is a variety such that every critical $k$-ary relation has the parallelogram property, and suppose that $m+n = k$. Let $\cF = \cF_\cV(x,y,z)$ be the free algebra on three generators in $\cV$. Suppose for contradiction that $\cV$ doesn't have an $m,n$-parallelogram term. Then
\[
\begin{bmatrix} x\\ \vdots\\ x\\ x\\ \vdots\\ x \end{bmatrix} \not\in \Sg_{\cF^k}\left\{\begin{bmatrix} y & y & x & z & \cdots & x & x & \cdots & x\\ \vdots & \vdots & \vdots & \vdots & \ddots & \vdots & \vdots & \ddots & \vdots\\ y & y & x & x & \cdots & z & x & \cdots & x\\ x & y & y & x & \cdots & x & z & \cdots & x\\ \vdots & \vdots & \vdots & \vdots & \ddots & \vdots & \vdots & \ddots & \vdots\\ x & y & y & x & \cdots & x & x & \cdots & z \end{bmatrix}\right\},
\]
so by Zorn's Lemma there exists a maximal $k$-ary relation $\RR$ on $\cF$ which contains the right hand side but does not contain the tuple $(x,...,x)$. The relation $\RR$ is then a critical $k$-ary relation on $\cF$, since every relation which properly contains $\RR$ must contain $\RR^* = \Sg(\RR \cup \{(x,...,x)\})$ and since every existential projection of $\RR$ onto a proper subset of the coordinates contains a vector of all $x$s (by the last $k$ columns of the matrix of generators above). However, $\RR$ does not have the parallelogram property when considered as a binary relation on $\cF^m\times \cF^n$, by the first three columns of the matrix of generators above, contradicting our assumption on $\cV$.
\end{proof}

\begin{cor} A variety $\cV$ has $k$-parallelogram terms iff for every relation $\RR \le \bA_1 \times \cdots \times \bA_n$ with $\bA_i \in \cV$, there exists a relation $\RR' \le \bA_1 \times \cdots \times \bA_n$ such that $\RR'$ has the parallelogram property and
\[
\RR = \RR' \cap \bigcap_{I \subseteq [n], |I| < k} \pi_{I}(\RR).
\]
\end{cor}

The relation $\RR'$ from the corollary need not be so mysterious: we can take it to be the \emph{least} relation $\RR'$ which contains $\RR$ and has the parallelogram property, since any intersection of relations which have the parallelogram property also has the parallelogram property. This choice of $\RR'$ can also be ``generated'' from $\RR$, by repeatedly adjoining tuples which are required to be inside in order for the parallelogram property to hold.

More explicitly, for any $I \subseteq [n]$, we can find the least relation $\RR^I$ which contains $\RR$ and has the (binary) parallelogram property when considered as a subalgebra of
\[
\Big(\prod_{i \in I} \bA_i\Big) \times \Big(\prod_{j \not\in I} \bA_{j}\Big),
\]
by finding the linking congruence of $\RR$ when considered as a subalgebra of the above, which restricts to a congruence $\alpha_I \in \Con(\pi_I(\RR))$, and taking $\RR^I$ to be the relation $\alpha_I \circ \RR$. We can then take
\[
\RR' = \bigcup_{I_1, I_2, ... \subseteq [n]} \RR^{I_1I_2\cdots}.
\]
In particular, if all of the algebras $\bA_i$ are finite, then $\RR'$ is contained in the (multisorted) relational clone generated by $\RR$.

\subsection{Critical rectangular relations in congruence modular varieties}

Using the commutator theory for congruence modular varieties, we can give a more detailed structure theory for the high-arity critical relations preserved by algebras with few subpowers. In fact, this structure theory applies more generally in congruence modular varieties, so long as we restrict our attention to critical relations with a weak form of the parallelogram property.

\begin{defn} A relation $\RR \le \bA_1 \times \cdots \times \bA_k$ is said to have the $1,k-1$-\emph{parallelogram property}, or alternatively is called \emph{rectangular}, if for any $i \le k$, when we regard $\RR$ as a binary relation on
\[
(\bA_1 \times \cdots \times \bA_{i-1} \times \bA_{i+1} \times \cdots \times \bA_k) \times \bA_i,
\]
it has the (binary) parallelogram property.
\end{defn}

The main property of subdirect rectangular relations which we need - and which holds in complete generality, not just in the context of congruence modularity - is that if we define a congruence $\theta_i$ on $\bA_i$ from the linking congruence of $\RR$ (considered as a binary relation on $(\cdots)\times \bA_i$), then we have $x \in \RR$ iff $x/\prod_i \theta_i \in \RR/\prod_i \theta_i$. Thus we may as well study the relation
\[
\RR/\prod_i \theta_i \le_{sd} \bA_1/\theta_1 \times \cdots \times \bA_k/\theta_k
\]
instead of studying $\RR$. The reduced relation is critical if the original $\RR$ is critical, is still rectangular, and has trivial linking congruences on each $\bA_i/\theta_i$, so it can be viewed as the graph of a surjective homomorphism
\[
\pi_{[k]\setminus\{i\}}\Big(\RR/\prod_i \theta_i\Big) \twoheadrightarrow \bA_i/\theta_i
\]
for each $i$.

\begin{defn} A subdirect rectangular relation $\RR \le_{sd} \bA_1 \times \cdots \times \bA_k$ is called \emph{reduced} if for each $i \le k$, $\RR$ is the graph of a surjective homomorphism
\[
\pi_{[k]\setminus\{i\}}(\RR) \twoheadrightarrow \bA_i,
\]
or equivalently, for each $i$ the map
\[
\pi_{[k]\setminus \{i\}} : \RR \rightarrow \pi_{[k]\setminus\{i\}}(\RR)
\]
is an isomorphism, i.e. $\ker \pi_{[k]\setminus \{i\}} = 0_\RR$.
\end{defn}

\begin{prop} If $\RR \le_{sd} \bA_1 \times \cdots \times \bA_k$ is a reduced subdirect critical rectangular relation, then each $\bA_i$ is subdirectly irreducible.
\end{prop}
\begin{proof} Let $\RR^*$ be the cover of $\RR$ in the lattice of subalgebras of $\bA_1 \times \cdots \times \bA_k$, and let $a = (a_1, ..., a_k)$ be a key tuple for $\RR$, that is, an element of $\RR^*\setminus \RR$. Since $\RR$ is critical, for every $i$ there is some $b_i \in \bA_i$ such that $(a_1, ..., a_{i-1}, b_i, a_{i+1}, ..., a_k) \in \RR$ (and this $b_i$ is unique, since $\RR$ is reduced). The claim is that for each $i$, every nontrivial congruence on $\bA_i$ contains the pair $(a_i,b_i)$ - that is, each $\bA_i$ is subdirectly irreducible with monolith equal to the congruence generated by the pair $(a_i,b_i)$.

Let $\psi_i \in \Con(\bA_i)$ be any nontrivial congruence. Then the relation
\[
\exists y_i\ ((x_1, ..., x_{i-1}, y_i, x_{i+1}, ..., x_k) \in \RR) \wedge (x_i \equiv_{\psi_i} y_i)
\]
strictly contains $\RR$ (since $\RR$ is reduced), so it contains $\RR^*$, and in particular contains the key tuple $a$. Using the fact that $\RR$ is reduced again, we see that the pair $(a_i,b_i)$ must be contained in $\psi_i$.
\end{proof}

As it turns out, reduced critical rectangular relations are closely related to the concept of \emph{similarity} between subdirectly irreducible algebras (see Appendix \ref{ss-similarity}). We won't need the full theory of similarity, just the following definition.

\begin{defn} If $\bA_1, ..., \bA_k$ are subdirectly irreducible algebras, then we say that an algebra $\RR$ is the \emph{graph of a joint similarity} between the $\bA_i$s if for each $i$, $\RR$ has a (critical) congruence $\alpha_i$ with $\RR/\alpha_i \cong \bA_i$, and for each pair $i,j$ there are congruences $\gamma_{ij}, \delta_{ij} \in \Con(\RR)$ such that
\[
\llbracket \alpha_i, \alpha_i^* \rrbracket \searrow \llbracket \gamma_{ij}, \delta_{ij} \rrbracket \nearrow \llbracket \alpha_j, \alpha_j^*\rrbracket.
\]
More explicitly, this means that $\alpha_i \vee \delta_{ij} = \alpha_i^*$, $\alpha_j \vee \delta_{ij} = \alpha_j^*$, and $\alpha_i \wedge \delta_{ij} = \alpha_j \wedge \delta_{ij}$.

Note that by Proposition \ref{similarity-graph}, $\RR/(\alpha_1 \wedge \cdots \wedge \alpha_k)$ is also a graph of a joint similarity, so there is no real loss in restricting to the case where $\RR$ is a subdirect product of the $\bA_i$s, with $\alpha_i = \ker \pi_i$.
\end{defn}

\begin{thm}[Kearnes, Szendrei \cite{parallelogram-terms}]\label{rectangular-structure} If $\RR \le_{sd} \bA_1 \times \cdots \times \bA_k$ is a reduced subdirect critical rectangular relation of arity $k \ge 3$ in a congruence modular variety, then
\begin{itemize}
\item[(a)] $\RR$ is the graph of a joint similarity between the $\bA_i$s,

\item[(b)] for each $i,j$, the image of $\pi_{i,j}(\RR)$ in $\bA_i/(0_{\bA_i}:0_{\bA_i}^*) \times \bA_j/(0_{\bA_j}:0_{\bA_j}^*)$ is the graph of an isomorphism
\[
\bA_i/(0_{\bA_i}:0_{\bA_i}^*) \xrightarrow{\sim} \bA_j/(0_{\bA_j}:0_{\bA_j}^*),
\]

\item[(c)] each monolith $0_{\bA_i}^*$ is abelian, and

\item[(d)] the cover $\RR^*$ is also rectangular, and the linking congruence of $\RR^*$ on $\bA_i$ is the monolith $0_{\bA_i}^*$.
\end{itemize}
If $\RR$ has the parallelogram property, then so does its cover $\RR^*$.
\end{thm}
\begin{proof} (a) Let $a = (a_1, ..., a_k) \in \RR^*\setminus \RR$ be a key tuple for $\RR$, and for each $i$ let $b_i \in \bA_i$ such that $(a_1, ..., a_{i-1}, b_i, a_{i+1}, ..., a_k) \in \RR$. Let $a^i = (a_1, ..., a_{i-1}, b_i, a_{i+1}, ..., a_k)$. Then for any $i \ne j$, if we let $\delta_{ij}$ be the congruence generated by the pair $(a^i, a^j)$, we claim that
\[
\llbracket \ker \pi_i, (\ker \pi_i)^* \rrbracket \searrow \llbracket 0_\RR, \delta_{ij} \rrbracket.
\]
The equality $\ker \pi_i \vee \delta_{ij} = (\ker \pi_i)^*$ was proved in the previous proposition. For the equality $\ker \pi_i \wedge \delta_{ij} = 0_\RR$, note that
\[
\ker \pi_i \wedge \delta_{ij} \le \ker \pi_i \wedge \ker \pi_{[k]\setminus\{i,j\}} = \ker \pi_{[k]\setminus \{j\}} = 0_\RR,
\]
where the last equality follows from the fact that $\RR$ is reduced.

(b) This follows directly from (a) and the Diamond Isomorphism Theorem \ref{diamond-isom} - for details, see Proposition \ref{similarity-graph}.

(c) By Proposition \ref{similarity-graph} again, if $\pi_{i,j}(\RR)$ is not the graph of an isomorphism for any pair $i,j$, then each monolith $0_{\bA_i}^*$ must be abelian.

(d) Suppose that $u,v,w \in \RR$ with $\pi_{[k]\setminus\{i\}}(u) = \pi_{[k]\setminus \{i\}}(v)$ and $v_i = w_i$. We need to show that there is some element $t \in \RR$ with $\pi_{[k]\setminus\{i\}}(t) = \pi_{[k]\setminus\{i\}}(w)$ and $t_i = u_i$.

Since $\RR^*$ is contained in the relation
\[
\exists y_i\ ((x_1, ..., x_{i-1}, y_i, x_{i+1}, ..., x_k) \in \RR) \wedge (x_i \equiv_{0_{\bA_i}^*} y_i)
\]
and $\RR$ is reduced, we have $(u_i,v_i) \in 0_{\bA_i}^*$. Let $p(x,y,z)$ be a Gumm difference term as in Theorem \ref{gumm-difference}, i.e. a term such that $p(y,y,x) \approx x$, and such that for $(x,y) \in \theta$ and $\theta$ any congruence we have $p(x,y,y)\ [\theta, \theta]\ x$. Then taking $\theta = 0_{\bA_i}^*$, we have $p(u_i,v_i,v_i) = u_i$ by part (c), so we can take $t = p(u,v,w)$.

For the last claim, suppose that we view $\RR^*$ as a binary relation on $\bA_I \times \bA_{[n]\setminus I}$, where we set $\bA_I = \prod_{i \in I} \bA_i$, and that we have $(a,b), (c,b), (c,d) \in \RR^*$. Pick some $i \in I$ and $j \not\in I$. Then there is some $a'$ such that $\pi_{I\setminus\{i\}}(a') = \pi_{I\setminus\{i\}}(a)$, $a_i' \equiv_{0^*_{\bA_i}} a_i$, and $(a',b) \in \RR$. Similarly find $c'$ which only differs from $c$ in the $i$th coordinate, has $c_i' \equiv_{0^*_{\bA_i}} c_i$, and has $(c',b) \in \RR$. Then $(c',d) \in \RR^*$ by part (d), so we can find $d'$ which only differs from $d$ in the $j$th coordinate, has $d_i' \equiv_{0^*_{\bA_j}} d_i$, and has $(c',d') \in \RR$. Then by the parallelogram property for $\RR$, we have $(a',d') \in \RR$, so by part (d) we have $(a,d) \in \RR^*$.
\end{proof}

\begin{ex}\label{gmm-ex-2-critical} Consider the generalized majority-minority algebra $\bA = (\{a,b,c\}, \varphi_2)$ from Example \ref{gmm-ex-2}, which is subdirectly irreducible with abelian monolith $0_{\bA}^*$ corresponding to the partition $\{a\}, \{b,c\}$ of its elements, and has $\bA/0_\bA^*$ isomorphic to a two element majority algebra. We can check that the monolith $0_{\bA}^*$ of $\bA$ is equal to its own centralizer by verifying that $[1_{\bA}:0_\bA^*] = 0_\bA^*$ and $[0_\bA^*,0_\bA^*] = 0_\bA$: to see this, note that
\[
\varphi_2\left(\begin{bmatrix} a & a\\ b & b\end{bmatrix}, \begin{bmatrix} a & a\\ b & b\end{bmatrix}, \begin{bmatrix} b & c\\ b & c\end{bmatrix}\right) = \begin{bmatrix} a & a\\ b & c\end{bmatrix} \in \bM(1_\bA,0_\bA^*),
\]
so $(b,c) \in [1_\bA,0_\bA^*]$, while every element of $\bM(0_\bA^*,0_\bA^*)$ either has all entries equal to $a$, or has all entries in $\{b,c\}$ with an even number of $b$s and an even number of $c$s.

The ternary relation $\RR \le_{sd} \bA^3$ corresponding to the columns of the matrix
\[
\begin{bmatrix} a & b & b & c & c\\ a & b & c & b & c\\ a & b & c & c & b\end{bmatrix}
\]
is a reduced subdirect critical rectangular relation of arity $3$ (with key tuple $(c,c,c)$), so by the structure theorem it is the graph of a joint similarity between three copies of $\bA$. Every two-coordinate projection $\pi_{i,j}(\RR)$ is equal to the congruence $0_{\bA}^* = (0_\bA : 0_\bA^*)$, and the cover $\RR^*$ of $\RR$ in $\Inv_3(\bA)$ is the relation $x\ 0_\bA^*\ y\ 0_\bA^*\ z$.

More generally, for any $k$ we can define a relation $\RR_k \le_{sd} \bA^k$ which contains the tuple $(a,...,a)$ together with the $2^{k-1}$ tuples in $\{b,c\}^k$ such that the total number of $c$s is even, and we see that $\RR_k$ is a reduced critical rectangular relation for each $k$. We claim that for every $k \ge 3$, there are exactly four critical relations in $\Inv_k(\bA)$: $\RR_k$, $\RR_k \setminus\{(a,...,a)\}$, and the two relations we get from these by swapping $b$s and $c$s in the last coordinate.

To prove the claim, we first note that the only algebras in $HS(\bA)$ which have abelian monoliths are $\bA$ and $\{b,c\}$, and that these two algebras are not similar to each other (since $\bA/0_\bA^*$ is not isomorphic to any quotient of $\{b,c\}$). Thus by the structure theorem, we only need to consider relations which are either subdirect in $\bA^k$ or subdirect in $\{b,c\}^k$. The interesting case is the case of relations which are subdirect in $\bA^k$.

The next thing we need to check is that no graph of a similarity $\bC \le_{sd} \bA^2$ from $\bA$ to $\bA$ induces the isomorphism $\bA/0_\bA^* \rightarrow \bA/0_\bA^*$ which corresponds to swapping the equivalence classes $\{a\}$ and $\{b,c\}$ of $0_\bA^*$. Note that the only candidate for $\bC$ is the relation $\{(a,b),(a,c),(b,a),(c,a)\}$, and for this choice of $\bC$ the congruence lattice $\Con(\bC)$ is given by the following picture.
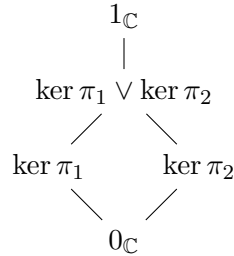
\begin{center}
\begin{tikzpicture}[scale=1]
  \node (1) at (0,3) {$1_{\bC}$};
  \node (t) at (0,2) {$\ker \pi_1 \vee \ker \pi_2$};
  \node (p1) at (-1,1) {$\ker \pi_1$};
  \node (p2) at (1,1) {$\ker \pi_2$};
  \node (0) at (0,0) {$0_{\bC}$};
  \draw (1) -- (t) -- (p1) -- (0) -- (p2) -- (t);
\end{tikzpicture}
\end{center}
As the reader can see, there is no pair $\gamma,\delta \in \Con(\bC)$ such that $\llbracket \ker \pi_1, \ker \pi_1 \vee \ker \pi_2 \rrbracket \searrow \llbracket \gamma, \delta \rrbracket \nearrow \llbracket \ker \pi_2, \ker \pi_1 \vee \ker \pi_2 \rrbracket$, so $\bC$ is not the graph of a similarity. Alternatively, we can see that $\bC$ can't be the graph of a similarity using the characterization in Corollary \ref{similar-detail}, since the corresponding congruence classes of $0_\bA^*$ which are linked by $\bC$ do not have the same sizes.

Thus, in any subdirect critical relation $\RR \le_{sd} \bA^k$ of arity $k > 2$, each $\pi_{i,j}(\RR)$ must be the congruence $0_\bA^*$, so $\RR$ will consist of the tuple $(a,...,a)$ together with some subalgebra of $\{b,c\}^k$. Since for any $\bS \le \{b,c\}^k$ the set $\bS \cup \{(a,...,a)\}$ will always be closed under $\varphi_2$, if $\RR$ is critical then so is $\RR \setminus \{(a,...,a)\}$, and it's easy to check that there are only two critical relations $\bS \le_{sd} \{b,c\}^k$. This completes the classification of critical relations in $\Inv_k(\bA)$ for $k > 2$.
\end{ex}

\begin{rem} Using the structure theorem \ref{rectangular-structure} and the fact that the centralizer of the monolith $(0:0^*)$ is automatically abelian for subdirectly irreducible algebras in residually small congruence modular varieties (Corollary \ref{residual-crit}), one can easily reduce the subpower membership problem \ref{subpower-membership} for residually small congruence modular varieties to the subpower membership problem for abelian groups by taking advantage of the properties of the Gumm difference term (see Corollary \ref{difference-graph}). For details of the reduction, see \cite{subpower-residually-small}.
\end{rem}

\begin{ex} We give an example of a minimal algebra with few subpowers which does not generate a residually small variety. Let $\bA = (\{a,b,c,d\}, g)$, where $g$ is the idempotent ternary symmetric operation which is determined by that fact that it commutes with the cyclic permutation $\sigma = (a\ b\ c\ d)$ and satisfies
\begin{align*}
g(a,a,b) &= a,\\
g(a,a,c) &= c,\\
g(a,a,d) &= c,\\
g(a,b,c) &= c.
\end{align*}
Then $\bA$ has a unique nontrivial congruence $0^*_\bA$ corresponding to the partition $\{a,c\},\{b,d\}$, and $\bA/0^*_\bA$ is isomorphic to a two element majority algebra. The congruence classes of $0^*_\bA$ are affine over $\ZZ/2$, and the algebra $\bS = \Sg_{\bA^2}\{(a,b),(b,a)\}$ has a congruence $\psi$ corresponding to the partition
\[
\left\{\begin{bmatrix} a\\ b\end{bmatrix},\begin{bmatrix} b\\ c\end{bmatrix},\begin{bmatrix} c\\ d\end{bmatrix}, \begin{bmatrix} d\\ a\end{bmatrix}\right\}, \left\{\begin{bmatrix} a\\ d\end{bmatrix},\begin{bmatrix} b\\ a\end{bmatrix},\begin{bmatrix} c\\ b\end{bmatrix}, \begin{bmatrix} d\\ c\end{bmatrix}\right\},
\]
such that $\bS/\psi$ is isomorphic to a two element affine algebra over $\ZZ/2$ (which is isomorphic to $\{a,c\}$). In fact, we have an isomorphism $\bS \cong \bA\times \{a,c\}$.

To see that $\bA$ has few subpowers, let $e$ be the term
\[
e(u,x,y,z) = g(x,g(u,y,y),g(y,g(x,y,z),g(x,y,z))).
\]
Then $e$ acts like the majority operation $g(x,y,z)$ on $\bA/0^*_\bA$, acts like the minority operation $g(x,u,y)$ on $\{a,c\}$, and has
\begin{align*}
e\left(\begin{bmatrix} b & b & a & a\\ b & a & b & a\\ a & a & a & b\end{bmatrix}\right) &= g\left(\begin{bmatrix} b & a & a\\ a & b & a\\ a & a & a\end{bmatrix}\right) = \begin{bmatrix} a\\ a\\ a\end{bmatrix},\\
e\left(\begin{bmatrix} d & d & a & a\\ d & a & d & a\\ a & a & a & d\end{bmatrix}\right) &= g\left(\begin{bmatrix} d & c & a\\ a & d & c\\ a & a & a\end{bmatrix}\right) = \begin{bmatrix} a\\ a\\ a\end{bmatrix}.
\end{align*}
Thus $e$ is a $3$-edge term.

Note that applying $\sigma$ to the second coordinate of $\bS$ turns it into $0^*_\bA$, and under the isomorphism $(1,\sigma) : \bS \xrightarrow{\sim} 0^*_\bA$, one of the congruence classes of $\psi$ becomes the diagonal $\{(x,x) \mid x \in \bA\}$. Thus $0^*_\bA$ is the center of $\bA$, and $\bA$ is similar to the idempotent reduct of $\ZZ/2$. Since $1_\bA = (0_\bA : 0^*_\bA)$ is not abelian, we see that $\bA$ can't generate a residually small variety.

We can check that
\[
\begin{bmatrix} a & b\\ a & b\\ a & b\end{bmatrix} \not\in \Sg_{\bA^{3\times 2}}\left\{\begin{bmatrix} a & b\\ a & b\\ b & a\end{bmatrix}, \begin{bmatrix} a & b\\ b & a\\ a & b\end{bmatrix}, \begin{bmatrix} b & a\\ a & b\\ a & b\end{bmatrix}\right\}
\]
by taking the rows modulo $\psi$. Thus none of the subsets $\{a,b\}, \{b,c\}, \{c,d\}, \{d,a\}$ (which are taken to each other by powers of the automorphism $\sigma$) are closed under any term which acts nontrivially on $\bA/0^*_\bA$. Using this, one can show that $\Clo(g)$ does not contain any proper Taylor subclones.

What do critical relations on $\bA$ look like? Suppose that $\RR \le_{sd} \bA^m\times \{a,c\}^n$ is critical and subdirect for some $m,n$ with $m+n \ge 3$. By Theorem \ref{parallelogram-critical}, $\RR$ has the parallelogram property. All we can conclude from Theorem \ref{rectangular-structure} is that $\RR^*$ has the parallelogram property and has linking congruence $(0^*_\bA)^m\times 1_{\{a,c\}}^n$, so the reduction $\RR^*_{red}$ of $\RR^*$ is a subdirect $m$-ary relation on the two element majority algebra $\bA/0^*_\bA$ which has the parallelogram property.

Luckily, it turns out that any such $\RR^*_{red}$ has $\pi_{ij}(\RR^*_{red})$ either a full relation or the graph of an automorphism of $\bA/0^*_\bA$ for any $i,j \in [m]$. This can be proved directly by reasoning about globally consistent instances of 2-SAT whose solution sets have the parallelogram property, or it can be proved more abstractly by using the fact that the two element majority algebra is subdirectly irreducible and generates a congruence distributive variety.

However we prove the claim about $\RR^*_{red}$, we see that if we assume without loss of generality that $(a,...,a) \in \RR$ (by applying powers of $\sigma$ to coordinates of $\RR$), then we can group the coordinates of $\RR$ into groups of size $m_1, ..., m_k$,
\[
\RR \le_{sd} \bA^{m_1} \times \cdots \times \bA^{m_k} \times \{a,c\}^n,
\]
such that $\pi_{ij}(\RR)$ is full for coordinates $i,j$ coming from separate groups, and $\pi_{ij}(\RR) = 0^*_{\bA}$ for coordinates $i,j$ coming from the same group.

Since we have assumed $(a,...,a) \in \RR$, $\RR$ must be closed under the unary polynomial $\phi : x \mapsto g(a,x,x)$. Since $\phi(a) = \phi(c) = a$ and $\phi(b) = \phi(d) = d$, we see that any vector of $a$s and $d$s which is constant on each group of coordinates will be contained in $\RR$. From this we see that in fact, any piecewise-constant vector
\[
((x_1,...,x_1),(x_2,...,x_2), ..., (x_k, ..., x_k), (a,...,a)) \in \bA^{m_1} \times \cdots \times \bA^{m_k} \times \{a,c\}^n
\]
must be contained in $\RR$. If we now consider the restriction $\RR \cap \{a,c\}^{m+n}$, then we find that it is an affine space defined by a system of linear equations over $\ZZ/2$, where the number of coordinates from any single group which show up in any equation must be even, since we may swap $(a,...,a),(c,...,c) \in \bA^{m_i}$ in any element of $\RR$. Thus we see that $\RR$ can be written as an intersection of relations $\RR'$ where the coordinates pair up in groups $\{i,j\}$ of size two, such that $\pi_{ij}(\RR') = 0^*_\bA$ and the relation $\RR'$ factors through the map $0^*_\bA \twoheadrightarrow \{a,c\}$ for each such pair of coordinates.

Using the above analysis, we see that the relational clone corresponding to $\bA$ is generated by the graph of the automorphism $\sigma$, which is $\Sg_{\bA^2}\{(a,b),(d,a)\}$, the critical binary relation $\Sg_{\bA^2}\{(a,a),(a,b),(b,b)\}$, which corresponds to a partial order on the majority algebra $\bA/0^*_\bA$, and the ternary relation $\Sg_{\bA^3}\{(a,a,a),(a,c,c),(b,b,a)\}$, which is the graph of the homomorphism $0^*_\bA \twoheadrightarrow \{a,c\}$.
\end{ex}

\section{Learnability of relations encoded by compact representations}

We'll start off by reviewing some of the standard definitions of learning theory.

\begin{defn} Fix a universe $U$. We call a collection $\cC$ of subsets of $U$, together with a rule for encoding the elements of $\cC$, a \emph{concept class}. An encoding of an element $C \in \cC$ is called a \emph{concept} (from $\cC$). The encoding scheme is called \emph{polynomially evaluable} if there is an algorithm which takes an encoding of a concept $C \in \cC$ and an element $x \in U$, and determines whether $x \in C$ in polynomial time.
\end{defn}

Generally we imagine a situation in which a teacher knows a target concept $C \in \cC$, and a student tries to learn the target concept $C$ from the teacher, either by seeing (random) examples of elements in $U$ and being told whether or not they are in the target concept $C$, or by asking the teacher certain types of questions. The teacher is modeled as an oracle which can be queried by the learner.

The main model which we will be examining in this section is the model of \emph{exact learning with (improper) equivalence queries} from \cite{angluin-learning}. Learnability results in the equivalence query model can be converted directly into learnability results in the \emph{probably approximately correct} model (which is often abbreviated as PAC-learning).

\begin{defn} Let $\cC'$ be a concept class which contains $\cC$, and call $\cC'$ the \emph{hypothesis} class. We define an \emph{equivalence oracle} $O_C$ with \emph{target concept} $C \in \cC$ to be the function which takes as input a hypothesis $C' \in \cC'$, returns ``true'' if $C = C'$, and otherwise returns an (arbitrary) element of the symmetric difference $C\Delta C'$.
\end{defn}

\begin{defn} An algorithm which makes calls to an oracle $O$ is said to \emph{learn} the concept class $\cC$ in the exact model with equivalence queries if, when the oracle $O$ is the equivalence oracle $O_C$ with target concept $C \in \cC$, the algorithm makes finitely many calls to the oracle $O$ with encodings of hypotheses $C' \in \cC'$ before finally discovering the concept $C$. The learning algorithm is called \emph{proper} if $\cC' = \cC$, and \emph{improper} otherwise. If there is an algorithm which learns $\cC$ in time polynomial in $\log |U|$, then we say that $\cC$ is \emph{polynomially learnable}.
\end{defn}

We are interested in the case where the universe $U$ is $\bA^n$ for $n$ large and $\bA$ a fixed algebraic structure, and where the concept class $\cC$ consists of the set of subalgebras of $\bA^n$, i.e. $\cC = \Inv_n(\bA)$ (recall $\Inv_n(\bA)$ is the set of $n$-ary relations which are preserved by the basic operations of $\bA$). In order for polynomial (in $n$) length encodings of the concepts in $\cC$ to exist, we need $\log |\cC|$ to be bounded by a polynomial in $n$, that is, we need $\bA$ to have few subpowers.

Suppose that $\bA$ has a $k$-edge term, and fix a particular $k$-edge term $e$. In this case, $n$-ary relations on $\bA$ are naturally encoded by compact representations, so we will use compact representations as our encoding scheme for the concept class $\cC = \Inv_n(\bA)$.

For the sake of definiteness, we will slightly modify the definition of a compact representation $R$ by requiring that for each element $x_I$ of $\pi_I(R)$ (where $|I| < k$), a specific element $x \in R$ with $\pi_I(x) = x_I$ has been marked (by $x_I$), and similarly for each minority index $(i,a,b)$ of $R$, a particular ordered pair $(u_a,u_b) \in R^2$ witnessing this index has been marked (by $(i,a,b)$). We will also require that each element $x$ of the compact representation $R$ is marked at least once (i.e., either $x$ is part of a marked witness to a minority index of $R$, or $x$ is a marked witness for some element of a projection of $R$ onto a small set of coordinates).

One feature which we would like this encoding scheme to satisfy is that there should be a polynomial time procedure to check whether an element $a \in \bA^n$ is contained in the relation $\RR$ encoded by the compact representation $R$. In other words, we want our encoding scheme to be \emph{polynomially evaluable}. The next lemma can be used to show that our encoding scheme is polynomially evaluable. We use the notation $[i]$ for the set $\{1, ..., i\}$.

\begin{lem}\label{compact-rep-filter} Suppose that $R \subseteq \bA^n$ is a compact representation of $\RR \le \bA^n$, $i \le n$, $a \in \bA^n$, $b \in \RR$ with $\pi_{[i-1]}(a) = \pi_{[i-1]}(b)$, and set $c_i = d(b_i, a_i)$. Suppose that
\begin{itemize}
\item for each $I \subseteq [i]$ with $|I| < k$ and $i \in I$, the element $x^I \in R$ is the marked element of $R$ witnessing $\pi_I(x^I) = \pi_I(a)$, and

\item the pair $(u_a,u_c) \in R^2$ is the marked witness of the minority index $(i,a_i,c_i)$.
\end{itemize}
Then there is a term $t^{[i]}$ of $\bA$ which can be built out of the terms $e,s,p,d$ of Theorem \ref{edge-spd} in time polynomial in $n$, such that $b^{[i]} = t^{[i]}(b,u_a,u_c,x^{I_1}, ...) \in \RR$ satisfies $\pi_{[i]}(a) = \pi_{[i]}(b^{[i]})$.
\end{lem}
\begin{proof} The proof is a modification of the proof of Theorem \ref{edge-gen}, with the induction over subsets of $[i]$ modified to only involve polynomially many subsets of $[i]$. The trick is to consider sets $I$ of the form $[j] \cup J$, where $j \le i$, $|J| = k-1$, and $i \in J$. There are only polynomially many such sets $I$, and we can induct on $j$ to handle them.

So we will show by induction on $j$ that for every set $I = [j] \cup J$ with $|J| = k-1$ and $i \in J$, there is a term $t^I$ such that $b^I = t^I(b,u_a,u_c,x^{I_1}, ...)$ satisfies $\pi_I(b^I) = \pi_I(a)$. The base case $j = 0$ is handled by taking $t^I = x^I$ for $|I| = k-1$.

For the inductive step, note that if $I = [j] \cup J$, then we can also write $I = [j-1] \cup (\{j\} \cup J)$. Let $\{j\} \cup J = \{l_1, ..., l_{k-1}, i\}$, and define sets $I_1, ..., I_{k-1}$ by deleting $l_1, ..., l_{k-1}$, respectively, from $I$, and note that each of the sets $I_m$ has the form $I_m = [j-1] \cup J_m$, where $J_m = (\{j\} \cup J)\setminus \{l_m\}$ and $i \in J_m$. By the induction hypothesis, we have already constructed terms $t^{I_m}$ and corresponding elements $b^{I_m} \in \RR$ with $\pi_{I_m}(b) = \pi_{I_m}(a)$. Then if we consider
\[
s(b,b^{I_1}, ..., b^{I_{k-1}}),
\]
we see that if we restrict to the coordinates in $\{j\} \cup J$, we have
\[
s\left(\begin{bmatrix} a_{l_1} & b^{I_1}_{l_1} & a_{l_1} & \cdots & a_{l_1}\\ a_{l_2} & a_{l_2} & b^{I_2}_{l_2} & \cdots & a_{l_2}\\ \vdots & \vdots & \vdots & \ddots & \vdots\\ a_{l_{k-1}} & a_{l_{k-1}} & a_{l_{k-1}} & \cdots & b^{I_{k-1}}_{l_{k-1}}\\ b_i & a_i & a_i & \cdots & a_i\end{bmatrix}\right) = \begin{bmatrix} a_{l_1}\\ a_{l_2}\\ \vdots\\ a_{l_{k-1}}\\ c_i\end{bmatrix}.
\]
Additionally, if we consider
\[
p(u_c,u_a,b^{I_1}),
\]
then if we restrict to the coordinates in $\{j\} \cup J$, we have
\[
p\left(\begin{bmatrix} u_{l_1} & u_{l_1} & b^{I_1}_{l_1}\\ u_{l_2} & u_{l_2} & a_{l_2}\\ \vdots & \vdots & \vdots\\ u_{l_{k-1}} & u_{l_{k-1}} & a_{l_{k-1}}\\ c_i & a_i & a_i\end{bmatrix}\right) = \begin{bmatrix} b^{I_1}_{l_1}\\ a_{l_2}\\ \vdots\\ a_{l_{k-1}}\\ c_i\end{bmatrix}.
\]
Thus, if we take $t^I$ to be given by
\[
t^I = e(p(u_c,u_a,t^{I_1}), s(b, t^{I_1}, ..., t^{I_{k-1}}), t^{I_1}, ..., t^{I_{k-1}}),
\]
then when we restrict to the coordinates in $\{j\} \cup J$, we get
\[
e\left(\begin{bmatrix} b^{I_1}_{l_1} & a_{l_1} & b^{I_1}_{l_1} & a_{l_1} & \cdots & a_{l_1}\\ a_{l_2} & a_{l_2} & a_{l_2} & b^{I_2}_{l_2} & \cdots & a_{l_2}\\ \vdots & \vdots & \vdots & \vdots & \ddots & \vdots\\ a_{l_{k-1}} & a_{l_{k-1}} & a_{l_{k-1}} & a_{l_{k-1}} & \cdots & b^{I_{k-1}}_{l_{k-1}}\\ c_i & c_i & a_i & a_i & \cdots & a_i\end{bmatrix}\right) = \begin{bmatrix} a_{l_1}\\ a_{l_2}\\ \vdots\\ a_{l_{k-1}}\\ a_i\end{bmatrix},
\]
which completes the induction step.
\end{proof}

We can now check if an element $a \in \bA^n$ is in the relation encoded by the compact representation $R$ as follows. First we check that $\pi_I(a) \in \pi_I(S)$ for each $I$ with $|I| < k$, and let $x^I$ be the marked element of $R$ with $\pi_I(x^I) = \pi_I(a)$. Then we try to construct elements $b^{[i]} \in \Sg_{\bA^n}(R)$ iteratively with $\pi_{[i]}(b^{[i]})= \pi_{[i]}(a)$. We start by taking $b^{[k-1]} = x^{[k-1]}$, and repeatedly invoke Lemma \ref{compact-rep-filter} to see that if $(i,a_i,c_i) \in \Sig(R)$, where $c_i = d(b^{[i-1]}_i,a_i)$, then we can construct $b^{[i]} \in \Sg_{\bA^n}(R)$ in polynomial time. We formalize this procedure as a subroutine which I will call \texttt{Approximate}$(R,a)$ (this is almost the same as the combination of the subroutines \texttt{Interpolate} and \texttt{New-Fix-values} from \cite{few-subpowers-algorithm}).


\begin{algorithm}
\caption{\texttt{Approximate}$(R, a)$, $e,s,p,d$ terms as in Theorem \ref{edge-spd}, $R \subseteq \bA^n$ a compact representation such that $\pi_I(a) \in \pi_I(R)$ for all $I$ with $|I| < k$.}
\begin{algorithmic}[1]
\ForAll{$I \subseteq [n]$ with $|I| < k$}
\State Let $x^I$ be the marked element of $R$ with $\pi_I(x^I) = \pi_I(a)$.
\EndFor
\State Set $b^{[k-1]\cap[n]} = x^{[k-1]\cap[n]}$.
\For{$i$ from $k$ to $n$}
\State Set $c_i \gets d(b^{[i-1]}_i,a_i)$.
\If{$(i,a_i,c_i) \not\in \Sig(R)$}
\State \Return $b^{[i-1]}$.
\Else
\State Let $(u^i_a,u^i_c)$ be the marked witness of the minority index $(i,a_i,c_i)$ in $R$.
\EndIf
\For{$j$ from $1$ to $i-k+1$}
\ForAll{$l_1, ..., l_{k-1}$ with $j = l_1 < l_2 < \cdots < l_{k-1} < i$}
\State Set $I \gets [j] \cup \{l_2, ..., l_{k-1}, i\}$.
\For{$m$ from $1$ to $k-1$}
\State Set $I_m \gets I\setminus \{l_m\}$.
\EndFor
\State Set $b^I \gets e(p(u^i_c,u^i_a,b^{I_1}),s(b^{[i-1]},b^{I_1},...,b^{I_{k-1}}),b^{I_1},...,b^{I_{k-1}})$.
\EndFor
\EndFor
\EndFor
\State \Return $b^{[n]}$.
\end{algorithmic}
\end{algorithm}

The running time of \texttt{Approximate} is $O(n^{k+1})$: there are less than $n^k$ choices for $j = l_1 < \cdots < l_{k-1} < i$, and for each choice, computing the new $b^I$ takes $O(n)$ steps (since $b^I$ has $n$ coordinates). By only maintaining the values of $b^{[i-1]}$ and $b^I$ with $i \in I$ in the $i$th step through the outer loop, the memory required is reduced to $O(n^k)$, which is the same as the space required to store a typical compact representation $R$.

\begin{prop} If $R \subseteq \bA^n$ is a compact representation and $a \in \bA^n$ with $\pi_I(a) \in \pi_I(R)$ for all $I$ with $|I| < k$, then either \texttt{Approximate}$(R,a)$ returns $a$ and $a \in \Sg_{\bA^n}(R)$, or \texttt{Approximate}$(R,a)$ returns $b \ne a$ such that $b \in \Sg_{\bA^n}(R)$, and such that if $i$ is minimal with $b_i \ne a_i$, then the minority index $(i,a_i,d(b_i,a_i))$ is not witnessed in $R$.
\end{prop}

At this point everything seems wonderful, but there is one major wrinkle: we have no idea how to (efficiently) test whether a given ``compact representation'' $R \subseteq \bA^n$ is actually a compact representation of the subalgebra $\RR = \Sg_{\bA^n}(R)$ it generates - in other words, we don't know how to test whether $R$ is a valid encoding of a concept from the concept class $\cC = \Inv_n(\bA)$. While it's easy to test whether $R$ and $\RR$ have the same projections onto small subsets of the coordinates (just check whether $\pi_I(R)$ is closed under the operations of $\bA$ for all $I$ with $|I| < k$), what is missing is a way to test whether $R$ witnesses every minority index which is witnessed in $\RR$.

Let's think for a moment about the problem of checking whether $R$ and $\Sg(R)$ witness the same minority indices. Since there are only $n|\bA|^2$ possible minority indices, we may as well focus on one particular minority index $(i,a,b)$. By replacing $R$ with $\pi_{[i]}(R)$ and $n$ with $i$, we may reduce to the case $i = n$.

\begin{prop} Suppose that the minority index $(i,a,b)$ is witnessed by some pair $(u_a,u_b)$ in a relation $\RR \le \bA^n$. Then for any tuple $t_a \in \RR$ with $\pi_i(t_a) = a$, there is a tuple $t_b \in \RR$ such that the pair $(t_a,t_b)$ also witnesses the minority index $(i,a,b)$. If $i = n$, then $t_b$ is uniquely determined by $t_a$.
\end{prop}
\begin{proof} Take $t_b = p(u_b,u_a,t_a)$. Then the identity $p(y,y,x)\approx x$ implies that $\pi_{[i-1]}(t_b) = \pi_{[i-1]}(t_a)$, and the fact that $a,b$ are a minority pair (that is, that $d(b,a) = b$) and the identity $p(x,y,y) \approx d(x,y)$ imply that $\pi_i(t_b) = p(b,a,a) = b$.
\end{proof}

So we can check whether a minority index $(i,a,b)$ is witnessed by $\Sg(R)$ as follows. First we pick any tuple $t_a \in R$ with $\pi_i(t_a) = a$. Then we modify it to make a tuple $t_b$, by replacing the $i$th coordinate with a $b$. Finally, we check whether $\pi_{[i]}(t_b) \in \Sg(\pi_{[i]}(R))$. By the above results, we have $\pi_{[i]}(t_b) \in \Sg(\pi_{[i]}(R))$ if and only if $(i,a,b) \in \Sig(\Sg(R))$. We find ourselves naturally led to consider the \emph{subpower membership problem}.

\begin{prob}[Subpower Membership Problem]\label{subpower-membership} Given a finite subset $S \subseteq \bA^n$ and an element $x \in \bA^n$, determine if $x$ is in the subalgebra of $\bA^n$ generated by $S$.
\end{prob}

\begin{thm}[Bulatov, Mayr, Szendrei \cite{subpower-residually-small}] For a fixed finite algebra $\bA$ with few subpowers, the following problems are polynomial time reducible to each other:
\begin{itemize}
\item The subpower membership problem for $\bA$: determine if $x \in \Sg_{\bA^n}(S)$, given $x \in \bA^n$ and $S \subseteq \bA^n$.

\item Find a compact representation for $\Sg_{\bA^n}(S)$, given a subset $S \subseteq \bA^n$.

\item The subpower intersection problem for $\bA$: given subsets $R,S \subseteq \bA^n$, find a set of generators for $\Sg_{\bA^n}(R) \cap \Sg_{\bA^n}(S)$.
\end{itemize}
If $\bA$ has few subpowers and has a finite number of basic operations, then the subpower membership problem for $\bA$ is in NP.
\end{thm}
\begin{proof} Left as an exercise to the reader. The hardest bit is the claim that the subpower membership problem is in NP: for this, imagine that we have a set $R$ which looks like a compact representation, and consider the set $C$ of all $a \in \bA^n$ such that \texttt{Approximate}$(R,a)$ returns $a$. If $C$ is not closed under the basic operations of $\bA$, then there should be a \emph{witness} to the fact that $C$ is not closed, and a nondeterministic algorithm can guess such a witness, verify that it works, and use it to enlarge $R$.
\end{proof}

Unfortunately, whether the subpower membership problem is in P for algebras with few subpowers is currently an open problem (even in the special case of quasigroups). So we need to find a workaround for this issue.

The workaround is to enlarge the concept class $\cC = \Inv_n(\bA)$ to a larger concept class $\cC'$, where concepts in $\cC'$ are encoded by ``compact representations'' $R \subseteq \bA^n$, where we allow sets $R$ which are not compact representations of the subalgebra $\RR = \Sg_{\bA^n}(R)$ which they generate. In order to be precise about exactly what concept $C$ is encoded by $R$, we use the \texttt{Approximate} subroutine.

\begin{defn} If $R \subseteq \bA^n$ is a ``compact representation'', then the corresponding concept $C \subseteq \bA^n$ encoded by $R$ is defined by the following rule. An element $a \in \bA^n$ is in $C$ iff the following two conditions are satisfied:
\begin{itemize}
\item for every $I \subseteq [n]$ with $|I| < k$, we have $\pi_I(a) \in \pi_I(R)$, and

\item the subroutine \texttt{Approximate}$(R,a)$ returns $a$.
\end{itemize}
\end{defn}

The penalty we will pay for this workaround is that since the new concept class $\cC'$ is larger than $\cC = \Inv_n(\bA)$, our learning algorithm will now be making \emph{improper} equivalence queries. If the subpower membership problem for $\bA$ can be proved to be in P, then we will be able to upgrade to a learning algorithm which makes only proper equivalence queries.

Now we can finally describe the learning algorithm, which is remarkably simple.

\begin{algorithm}
\caption{\texttt{Learn}$(O)$, $O$ an equivalence oracle for an unknown target concept $C \in \Inv_n(\bA)$.}
\begin{algorithmic}[1]
\State Set $R \gets \emptyset$.
\While{$O(R)$ does not return ``true''}
\State Set $a \gets O(R)$.
\ForAll{$I \subseteq [n]$ with $|I| < k$ such that $\pi_I(a)$ has no designated witness in $R$}
\State Set $R \gets R\cup\{a\}$.
\State Mark $a$ as the designated witness for $\pi_I(a)$.
\EndFor
\While{$\texttt{Approximate}(R,a)$ does not return $a$}
\State Set $b \gets \texttt{Approximate}(R,a)$.
\State Let $i$ be minimal such that $a_i \ne b_i$.
\State Set $R \gets R\cup\{a,d(b,a)\}$.
\State Mark the pair $(a,d(b,a))$ as the designated witness for the minority index $(i,a_i,d(b_i,a_i))$.
\EndWhile
\State Optionally, enlarge $R$ further to make it closer to a compact representation of $\Sg(R)$.
\EndWhile
\end{algorithmic}
\end{algorithm}

\begin{prop} For a fixed algebra $\bA$ with few subpowers, the algorithm \texttt{Learn}$(O)$ takes time polynomial in $n$ to find an encoding $R$ for the target concept $C \in \Inv_n(\bA)$.
\end{prop}
\begin{proof} At every step of the algorithm, we have $\Sg(R) \subseteq C$: this is true at the beginning, and if it is true before we call $O(R)$, then since the concept $C'$ encoded by $R$ has $C' \subseteq \Sg(R) \subseteq C$, the value $a$ returned by $O(R)$ will be contained in $C\Delta C' = C \setminus C' \subseteq C$, so $\Sg(R \cup \{a\}) \subseteq C$.

Furthermore, every time we process a new $a \in C\setminus C'$, we strictly enlarge $R$ to make the new concept encoded by $R$ contain $a$, either by adding $a$ as a designated witness to $\pi_I(a) \in \pi_I(R)$ for some $I$, or by adding new minority indices which were not present in the original $R$. Since $R$ can only increase in size polynomially many times (as a compact representation has size bounded by a polynomial in $n$), we can only call the oracle polynomially many times before the process must terminate.
\end{proof}

\begin{rem} If we did not insist on polynomial evaluability of the encoding scheme (or if we could solve the subpower membership problem), then we could instead encode relations via generating sets. The learning algorithm would then be even simpler: at every step, the learner guesses that the target concept is the relation generated by all the examples it has seen so far. This learning algorithm is known as the \emph{closure algorithm}. The issue is that now the equivalence oracle becomes hard to implement, as the teacher is forced to determine whether a given set generates the target relation they have in mind.
\end{rem}

Now we will explain how all of this relates to Valiant's PAC-learning model \cite{valiant-pac}. In the PAC-learning model, the teacher (oracle) has access to both a target concept $C \in \cC$ and a probability distribution $\mu$ over the universe $U$, both of which are unknown to the learner. The learner is allowed to request random classified examples, sampled from the distribution $\mu$ (by a ``classified'' example, I mean that the learner is given an example and told whether or not it is in the target concept $C$).

\begin{defn} If $C \in \cC$ is a target concept and $\mu$ is a probability distribution on the universe $U$, then the \emph{sampling oracle} for the pair $C,\mu$ is a randomized oracle which samples a random element $a \in U$ drawn from the distribution $\mu$, and returns the ordered pair $(a,a \in C)$, where by ``$a \in C$'' we mean either ``true'' or ``false'' based on whether $a$ is in the target concept $C$.
\end{defn}

In the PAC-learning model, the goal of a learning algorithm is to output an encoding of a concept $C'$ in the hypothesis class $\cC'$, such that the $\mu$-measure $\mu(C \Delta C')$ of the symmetric difference between $C$ and $C'$ is small. We can't hope to do better than this, since the chance of seeing an example which lets us distinguish between $C$ and $C'$ is at most $\mu(C \Delta C')$ times the number of classified examples we request.

\begin{defn} We say that an algorithm with access to a sampling oracle \emph{learns} a concept class $\cC$ in the \emph{probably approximately correct} model with \emph{error} $\epsilon$ and \emph{confidence} $1-\delta$ if for any target concept $C \in \cC$ and any probability distribution $\mu$ over the universe, the algorithm eventually returns a hypothesis $C' \in \cC'$ such that
\[
\PP[\mu(C \Delta C') \le \epsilon] \ge 1 - \delta.
\]
The probability here is taken over the random choices made by the oracle (and possibly the learning algorithm) - the target concept $C$ is \emph{not} being randomized here, we require this for \emph{all} $C \in \cC$ and \emph{all} $\mu$.

We say that a concept class $\cC$ is \emph{efficiently PAC-learnable} if there is an algorithm which learns $\cC$ in the PAC-model and takes time polynomial in $\log(|U|)$, $\frac{1}{\epsilon}$, and $\log(\frac{1}{\delta})$ for $\epsilon, \delta > 0$.
\end{defn}

The standard learning algorithm in the PAC model is to request a large number of classified examples, and then choose \emph{any} hypothesis $C' \in \cC'$ which is consistent with all of the classified examples we have seen so far. For this to work, it is necessary that the hypothesis class $\cC'$ is in some sense ``small'', and we also need to have a way to efficiently find at least one hypothesis which is consistent with the examples. First we will define a measure of the ``size'' a the concept class $\cC$, known as the VC-dimension.

\begin{defn} If $\cC$ is a collection of subsets of some universe $U$, then we say that a set $S$ is \emph{shattered} by $\cC$ if for all $X \subseteq S$, there is some $C \in \cC$ with $C \cap S = X$. We define the \emph{Vapnik-Chervonenkis dimension} of $\cC$, written $\operatorname{VC}(\cC)$, to be the size of the largest set $S$ which is shattered by $\cC$.
\end{defn}

To see that the VC-dimension is a good measure of the complexity of a concept class, we recall the Sauer-Shelah Lemma.

\begin{lem}[Sauer-Shelah Lemma]\label{sauer-shelah} If $\cC$ is a collection of subsets of $U$ with VC-dimension $d$, then
\[
|\cC| \le \sum_{i = 0}^d \binom{|U|}{i}.
\]
In fact, we have the stronger result that the number of sets $S \subseteq U$ which are shattered by $\cC$ is at least $|\cC|$.
\end{lem}
\begin{proof} We show that $\cC$ shatters at least $|\cC|$ sets by induction on $|\cC|$. For the base case, note that the empty set is shattered by $\cC$ as long as $|\cC| \ge 1$. For the inductive step, let $x \in U$ be an element which is in some of the sets in $\cC$ but not all of them, and let $\cC_x$ be the collection of $C \in \cC$ with $x \in C$ and $\cC_x' = \cC \setminus \cC_x$. Inductively, $\cC_x$ shatters at least $|\cC_x|$ sets and $\cC_x'$ shatters at least $|\cC_x'|$ sets, and any set shattered by $\cC_x$ or $\cC_x'$ must not contain $x$.

To finish the induction, we just need to check that for any set $S$ which is shattered by both $\cC_x$ and $\cC_x'$, the set $S\cup \{x\}$ is shattered by $\cC$.
\end{proof}

If the set $S$ is shattered by $\cC$, then the sampling oracle could sample from a uniform distribution on $S$, and in this case the learner is faced with the problem of learning an arbitrary subset $X = C \cap S$ of $S$ given an oracle which returns uniformly random classified examples. If the learner examines $o(|S|)$ classified examples, then clearly they can't hope to succeed. The following result makes this precise.

\begin{prop} If an algorithm learns a concept class $\cC$ with error $\epsilon$ and confidence $1-\delta$ after requesting at most $m$ classified examples, then
\[
m \ge (2(1-\epsilon)(1-\delta) - 1)\operatorname{VC}(\cC).
\]
\end{prop}
\begin{proof} Let $S$ be a set with $|S| = \operatorname{VC}(\cC)$ which is shattered by $\cC$, and for each $X \subseteq S$ consider the sampling oracle $O_X$ which samples from the uniform distribution $\mu$ on $S$, and has target concept some $C_X \in \cC$ with $C_X \cap S = X$. If we average the performance of the learning algorithm over the sampling oracles $O_X$ (with $X$ chosen as a uniformly random subset of $S$), we see that if it outputs a hypothesis $C'$, then
\[
\EE[\mu(C_X \Delta C')] \ge \frac{1}{2}\Big(1 - \frac{m}{|S|}\Big).
\]
By Markov's inequality, this implies that
\[
(1-\epsilon)\PP[\mu(C_X\Delta C') \le \epsilon] \le \frac{1}{2} + \frac{m}{2|S|},
\]
so
\[
\frac{1}{2} + \frac{m}{2|S|} \ge (1-\epsilon)(1-\delta).\qedhere
\]
\end{proof}

Conversely, if the VC-dimension of $\cC$ is small, then the standard learning algorithm in the PAC model performs well, so long as it can be implemented.

\begin{thm}[VC-dimension determines sample-complexity \cite{vc-dim-learnability}] If $\operatorname{VC}(\cC') = d$, then any algorithm which takes
\[
m \ge \max\left(\frac{4}{\epsilon}\log\Big(\frac{2}{\delta}\Big),\ \frac{8d}{\epsilon}\log\Big(\frac{13}{\epsilon}\Big)\right).
\]
samples from a sampling oracle and outputs any hypothesis $C' \in \cC'$ consistent with the data will learn $\cC$ with error $\epsilon$ and confidence $1-\delta$.
\end{thm}
\begin{proof}[Sketch] Consider the following process: pick $2m$ samples from the probability distribution $\mu$, permute them randomly, feed the first $m$ samples (after permuting) to the learning algorithm, and count how many of the last $m$ samples are classified incorrectly by the hypothesis $C'$ chosen by the learning algorithm.

If the learning algorithm fails to learn $\cC$ with error $\epsilon$ and confidence $1-\delta$, then for some choice of target concept $C$ and distribution $\mu$, the process described will incorrectly classify at least $\frac{\epsilon m}{2}$ of the last $m$ samples with probability at least $\frac{\delta}{2}$, by Chebyshev's inequality (at least for $m \ge \frac{8}{\epsilon}$). Thus there will be some specific set $X$ of size $2m$, such that at least a $\frac{\delta}{2}$ fraction of its permutations lead to an incorrect classification of at least $\frac{\epsilon m}{2}$ of its last $m$ elements.

By the Sauer-Shelah Lemma \ref{sauer-shelah}, the number of distinct subsets of $X$ which can be written as $C' \cap X$ for some $C' \in \cC'$ is bounded by $\sum_{i \le d} \binom{2m}{i}$. For each possible intersection $C' \cap X$, the chance of the first $m$ samples from $X$ being consistent with $C'$ and the last $m$ samples from $X$ having at least $\frac{\epsilon m}{2}$ inconsistencies with $C'$ is at most $2^{-\epsilon m/2}$. Thus if the learning algorithm fails, then by the union bound we must have
\[
2^{-\epsilon m/2}\sum_{i \le d} \binom{2m}{i} \ge \frac{\delta}{2},
\]
and plugging in the assumed bounds on $m$ and chugging through the inequalities gives a contradiction.
\end{proof}

Note that if $\bA$ is an algebraic structure and we take $U = \bA^n$, $\cC = \Inv_n(\bA)$, then a set $S$ is shattered by $\Inv_n(\bA)$ iff $S$ is an \emph{independent} subset of $\bA^n$. Thus the VC-dimension of $\Inv_n(\bA)$ is exactly the same thing as the number $i_\bA(n)$, so if the concept classes $\Inv_n(\bA)$ are efficiently PAC-learnable as $n$ varies, then $\bA$ must have few subpowers.

We can convert learning algorithms in the equivalence query model into learning algorithms in the PAC model by using the sampling oracle to simulate an equivalence oracle.

\begin{prop}[Angluin \cite{angluin-learning}] If a concept class $\cC$ is efficiently learnable in the (improper) equivalence query model using a hypothesis class $\cC'$ which has a polynomially evaluable encoding scheme, then $\cC$ is also efficiently learnable in the PAC model.
\end{prop}
\begin{proof} Given a sampling oracle $O$, we simulate an equivalence oracle as follows. The $i$th time the equivalence oracle is called by the learner, say to determine whether the target concept $C$ is equivalent to a hypothesis $C' \in \cC'$, we call the sampling oracle $O$ some number $q_i$ times to get $q_i$ random classified examples, and we check whether the way they are classified agrees with the hypothesis $C'$ (here is where we are using polynomial evaluability). If their classifications do agree with $C'$, then we pretend that the equivalence oracle returned ``true'', and otherwise we pick one of the examples $a$ whose classification does not agree with $C'$ and return $a$ as the counterexample in $C \Delta C'$.

By the union bound, the probability that the simulated equivalence oracle \emph{ever} returns ``true'' for a hypothesis $C'$ with $\mu(C\Delta C') \ge \epsilon$ is at most
\[
\sum_i (1-\mu(C\Delta C'))^{q_i} \le \sum_i (1-\epsilon)^{q_i}.
\]
If we take
\[
q_i \ge \frac{1}{\epsilon}(\ln(1/\delta) + i\ln(2)),
\]
for instance, then we get
\[
\sum_i (1-\epsilon)^{q_i} \le \sum_i e^{-\epsilon q_i} \le \sum_i e^{\ln \delta - i\ln(2)} = \sum_i \frac{\delta}{2^i} \le \delta.\qedhere
\]
\end{proof}

\begin{rem} Another learning model is the on-line learning model described by Littlestone \cite{littlestone-online-learning}. In this model, the learner is repeatedly presented with examples, and for each example must guess its classification before being told whether its guess is correct. The goal of the learner is to have an upper bound on the number of incorrect guesses it makes, even if the sequence of examples is chosen adversarially. It is easy to convert a learnability result in the (improper) equivalence query model into an algorithm for on-line learning.
\end{rem}

\begin{rem} There is a variant of the PAC learning model in which the learner is also allowed to use \emph{membership queries}: in a membership query, the learner picks an element $x \in U$, and asks the teacher (oracle) whether $x$ is in the target concept.

In \cite{pac-membership}, several situations are given where the addition of membership queries can be shown not to help with learning, under some standard cryptographic assumptions. In \cite{learnable-gmm}, there is a claim that some of the impossibility results for PAC learning of $\Inv_n(\bA)$ when $\bA$ doesn't have few subpowers can be generalized to impossibility results in the model of PAC learning with membership queries (under cryptographic assumptions), but the exact statement and the proof are left to a ``full version'' of the paper which I have been unable to track down. The more recent paper \cite{membership-query-hardness} by Chen and Valeriote proves such a hardness result for algebraic structures which are not congruence modular, and for finitely related structures congruence modularity is equivalent to few subpowers by the main result of \cite{barto-valeriote-conjecture}.
\end{rem}


\section{Algebras with few subpowers are finitely related}

Suppose a clone $\cO$ on a finite domain $A$ has a $k$-edge term $e$. We want to show that there exists some finite set of relations $R_1, ..., R_m$ which generate the relational clone which is dual to $\cO$. This is equivalent to $\cO$ being exactly the set of operations $\Pol(R_1, ..., R_m)$ which preserve the relations $R_1, ..., R_m$. If $R_1, ..., R_m$ are all preserved by $\cO$, then the clone $\Pol(R_1, ..., R_m)$ will certainly contain $\cO$, but might end up being too large. In this case, $\Pol(R_1, ..., R_m)$ will still contain the $k$-edge term $e$, and we can use this to our advantage.

To understand the structure of a clone $\cO$ with a $k$-edge term, we go back to the explicit representation of the set $\cO_n$ of $n$-ary operations of $\cO$ as the free algebra over $\bA = (A,\cO)$ on $n$ generators, which is concretely given by the subalgebra
\[
\cO_n = \cF_\bA(x_1, ..., x_n) \le \bA^{\bA^n}
\]
generated by the elements $\pi_i : \bA^n \rightarrow \bA$ given by $\pi_i(a_1, ..., a_n) = a_i$, where $x_i \in \cF_\bA(x_1, ..., x_n)$ is identified with the element $\pi_i \in \bA^{\bA^n}$. Similarly, recall that the set of $n$-ary operations $f \in \Pol_n(R_1, ..., R_m)$, considered as a subalgebra of $\bA^{\bA^n}$, is given by the primitive positive formula
\[
f \in \Pol_n(R_1, ..., R_m) \iff \bigwedge_{i \le m} \bigwedge_{M \in R_i^n} f(M) \in R_i.
\]
To check that these two subalgebras of $\bA^{\bA^n}$ are equal, by Theorem \ref{edge-gen} and the fact that one is contained within the other, it suffices to check that they have the same projections onto subsets $I \subseteq \bA^n$ of the coordinates with $|I| < k$, and to check that they have the same forks. If $R_1, ..., R_m$ generate all relations of $\Inv(\bA)$ with arity less than $k$, then the first condition will be satisfied. The hard part is dealing with the forks.

In order to make precise statements about the set of forks in $\bA^{\bA^n}$, we first need to choose an ordering on the coordinates of $\bA^{\bA^n}$, that is, an ordering on the elements of $A^n$. A natural choice is to first fix any total order $\le$ on the set $A$, and to extend this to the \emph{lexicographic order} on $A^n$.

\begin{defn} If $(A, \le)$ is a set with a total order, then we define the \emph{lexicographic order} $\le_{lex}$ on $A^n$ by $a \le_{lex} b$ iff either $a = b$ or there is some $i \le n$ such that $a_j = b_j$ for $j < i$ and $a_i < b_i$. In other words, $a <_{lex} b$ if $a_i < b_i$ at the first coordinate $i$ where $a$ and $b$ differ.
\end{defn}

\begin{defn} If $(I,\le)$ is a totally ordered set and $R \subseteq A^I$ is a relation on $A$, then for $i \in I$ we define the set of \emph{forks} of $R$ at the $i$th coordinate to be the set of pairs $(a,b) \in A^2$ given by
\[
\Forks(R,i) \coloneqq \{(a,b) \mid \exists t_a, t_b \in R,\ \pi_{<i}(t_a) = \pi_{<i}(t_b), \pi_i(t_a) = a, \pi_i(t_b) = b\}.
\]
\end{defn}

So in order to understand a clone $\cO$ with a $k$-edge term, we need to understand the relations of arity less than $k$, together with the set of forks $\Forks(\cO_n,a)$ for all $a \in A^n$ and all $n$. The issue is that while each set $\Forks(\cO_n,a)$ is given by a finite collection of pairs of elements, there are infinitely many elements $a \in A^n, n \in \bN^+$ to consider. So we need a way to relate $\Forks(\cO_n,a)$ to $\Forks(\cO_m,b)$ for some choices of $a \in A^n, b \in A^m$.

\begin{prop} Suppose $a \in A^n, b \in A^m$. If there is a map $\phi : [m] \rightarrow [n]$ such that the associated function $\phi^* : A^n \rightarrow A^m$ given by $\phi^*(x_1, ..., x_n) = (x_{\phi(1)}, ..., x_{\phi(m)})$ satisfies the conditions
\begin{itemize}
\item $\phi^*(a) = b$ and
\item for all $c <_{lex} a$ we have $\phi^*(c) <_{lex} b$,
\end{itemize}
then for any clone $\cO$ on $A$, we have $\Forks(\cO_m,b) \subseteq \Forks(\cO_n,a)$.
\end{prop}
\begin{proof} Letting $\bA = (A,\cO)$, $\phi$ induces a natural map of free algebras $\cO_m \rightarrow \cO_n$ given by $x_i \mapsto x_{\phi(i)}$. We will write this natural map as $f \mapsto f_\phi$. Considering $\cO_m, \cO_n$ as subalgebras of $\bA^{A^m}, \bA^{A^n}$, respectively, we see that for $c \in A^n$ and $f \in \cO_m$, the $c$th coordinate of the image $f_\phi$ of $f$ under this map is given by
\[
f_\phi(c) = f(\phi^*(c)).
\]
In particular, if $t, t' \in \cO_m$ with $\pi_{<_{lex} b}(t) = \pi_{<_{lex} b}(t')$, then
\[
\pi_{<_{lex} a}(t_\phi) = \pi_{<_{lex} a}(t'_\phi),
\]
and
\[
\begin{bmatrix} t_\phi(a)\\ t'_\phi(a)\end{bmatrix} = \begin{bmatrix} t(\phi^*(a))\\ t'(\phi^*(a))\end{bmatrix} = \begin{bmatrix} t(b)\\ t'(b)\end{bmatrix},
\]
so every fork in $\Forks(\cO_m,b)$ is also a fork in $\Forks(\cO_n,a)$.
\end{proof}

\begin{prop} A map $\phi$ as in the previous proposition exists if there is a strictly increasing function $h : [n] \rightarrow [m]$ such that
\begin{itemize}
\item the same elements of $A$ occur in both $a$ and $b$,
\item $h^*(b) = a$, that is, $a_i = b_{h(i)}$ for all $i \in [n]$, and
\item for all $s \in A$, if the index of the first occurence of $s$ in $a$ is $i$, then $h(i)$ is the index of the first occurence of $s$ in $b$.
\end{itemize}
If no coordinate $a_i$ of $a$ is minimal or maximal with respect to the order $<$ on $\bA$, then the converse is true: such a $\phi$ exists iff such an $h$ exists.
\end{prop}
\begin{proof} Given such an $h$, we define $\phi$ as follows. We set $\phi(h(i)) = i$, and for $j$ not in the image of $h$ let $\phi(j)$ be the first index $i$ such that $a_i = b_j$, so that $h(\phi(j)) \le j$ for all $j$. Then for any $c <_{lex} a$, if $i$ is the first index where $a_i \ne c_i$, and if $j$ is the first coordinate where $\phi^*(a)$ and $\phi^*(c)$ differ, then we have $a_{\phi(j)} \ne c_{\phi(j)}$, so $i \le \phi(j)$, so $h(i) \le h(\phi(j)) \le j$, so we must have $h(i) = j$ since $\phi^*(a)$ and $\phi^*(c)$ also differ at $h(i)$. Thus $\phi^*(c) <_{lex} \phi^*(a) = b$.

Now suppose that no coordinate $a_i$ of $a$ is minimal or maximal with respect to the order $<$ on $\bA$. Then the map $\phi$ in the previous proposition must be surjective: if $i$ is not in the image of $\phi$, then we can define $c <_{lex} a$ which only differs from $a$ on the $i$th coordinate, and $\phi^*(c) = \phi^*(a) \not<_{lex} b$, contradicting the choice of $\phi$. Thus we can define $h : [n] \rightarrow [m]$ by
\[
h(i) = \min\{j \in [m] \mid \phi(j) = i\},
\]
so $h^*(b) = a$ and we see that $a$ and $b$ have the same set of symbols.

For any $i$, if we define $c <_{lex} a$ which matches $a$ up to the $i$th coordinate, has $c_i < a_i$, and $c_j > a_j$ for all $j > i$, then from $\phi^*(c) < b = \phi^*(a)$, we see that $h(i) < h(j)$ for all $j > i$. Thus $h$ must be strictly increasing. Finally, from the definition of $h$, we see that if $i$ is the index of the first occurence of $s$ in $a$, then $h(i)$ must be the index of the first occurence of $s$ in $b$.
\end{proof}

\begin{defn} Let $A^+ = \bigcup_{n \ge 1} A^n$, and define the partial order $\le_E$ on $A^+$ by $a \le_E b$ iff there exists a map $h$ as in the previous proposition. Equivalently, $a \le_E b$ iff the same set of elements of $A$ occur in $a$ and $b$, and $b$ can be formed from $a$ by inserting elements $s \in A$ after their first occurences in $a$.
\end{defn}

Note that the partial ordering $\le_E$ on $A^+$ has no dependence on the arbitrary choice of ordering $<$ we introduced on the elements of $A$ (of course, the set $\Forks(\cO,a)$ still depends on the choice of $<$). The partial order $\le_E$ is a refinement of the embeddability partial ordering that occurs in Higman's Lemma \cite{higmans-lemma}. We can now simplify the description of the sets $\Forks(\cO,a)$ using the ordering $\le_E$.

\begin{defn} For any pair $(c,d) \in A^2$, we define the set $\lambda(\cO,(c,d)) \subseteq A^+$ to be the set of $a \in A^+$ such that $(c,d) \not\in \Forks(\cO,a)$.
\end{defn}

\begin{cor} For any clone $\cO$ on a set $A$, the set $\lambda(\cO,(c,d))$ is upwards closed in $A^+$ with respect to $\le_E$, that is, if $a \in \lambda(\cO,(c,d))$ and $a \le_E b$, then $b \in \lambda(\cO,(c,d))$.
\end{cor}

To describe an upwards closed subset (also called an \emph{upset}) of a \emph{finite} poset, it is enough to describe its minimal elements. We want to show that $\lambda(\cO,(c,d))$ can be described in terms of its minimal elements, but for this to work, it's necessary to show that $(A^+,\le_E)$ is a \emph{well partial order}.

\begin{defn} A partial order $(X,\le)$ is a \emph{well partial order} if for every infinite sequence $x_1, x_2, ...$ of elements of $X$, there exists an infinite increasing subsequence $i_1 < i_2 < \cdots$ such that $x_{i_1} \le x_{i_2} \le \cdots$. 
\end{defn}

\begin{prop} A partial order $(X,\le)$ is a well partial order iff it has no infinite descending chains and no infinite antichains.
\end{prop}
\begin{proof} Let $x_1, x_2, ...$ be any infinite sequence of elements of $X$. Color the edges of the complete graph on $\NN^+$ with three colors, as follows: for $i < j$, the edge $\{i,j\}$ is colored red if $x_i > x_j$, colored blue if $x_i, x_j$ are incomparable, and colored green if $x_i \le x_j$. By Ramsey's Theorem, there must be some infinite monochromatic clique in this graph, so either there is an infinite descending chain, an infinite antichain, or an infinite subsequence $i_1 < i_2 < \cdots$ with $x_{i_1} \le x_{i_2} \le \cdots$.
\end{proof}

\begin{prop} A partial order $(X,\le)$ is a well partial order iff for all upsets $U \subseteq X$, every element of $U$ is $\ge$ some minimal element of $U$ and $U$ has finitely many minimal elements, that is, there exists a finite set of elements $u_1, ..., u_k \in U$ such that $U = \{x \mid x \ge u_i\text{ for some }i\}$.
\end{prop}
\begin{proof} Suppose first that $(X,\le)$ is a well partial order. If some element $u \in U$ is not above any minimal element of $U$, then we can find an infinite descending chain in $U$. Since any pair of distinct minimal elements of $U$ are incomparable, the number of minimal elements of $U$ must be finite. The converse follows from the previous proposition.
\end{proof}

So the last ingredient of the argument will be the proof that $\le_E$ is a well partial order. While the tools we have available are capable of proving this directly, it is useful to reduce this to the fact that the simpler (and more well-known) embeddability partial ordering $\le_e$, due to Higman, is a well partial order - this allows us to transfer other results about Higman's ordering to the partial order $\le_E$.

\begin{defn} Define the partial order $\le_e$ on $A^+$ by $a \le_e b$ if $b$ can be formed from $a$ by inserting elements of $A$.
\end{defn}

\begin{prop} If $B$ is the disjoint union of $A$ with the two-element set of symbols $\{\#,'\}$, then there is an embedding of partial orders
\[
F : (A^+, \le_E) \hookrightarrow (B^+, \le_e),
\]
i.e. a function $F$ such that for $x,y \in A^+$, we have $x \le_E y$ iff $F(x) \le_e F(y)$.
\end{prop}
\begin{proof} We define $F : A^+ \rightarrow B^+$ to be the function which modifies $x \in A^+$ by inserting a $'$ after the first occurence of each symbol within $x$, inserting a $\#$ at the end of $x$, and then following that with a $'$ for each symbol which doesn't occur within $x$. For instance, if $A = \{a,b,c,d,e\}$, then
\[
F(adaadca) = a'd'aadc'a\#'',
\]
where the two $'$s at the end keep track of the fact that $b,e$ did not occur within the word $adaadca$.

Note that $F(x)$ is always formed from $x$ by inserting exactly $1$ copy of $\#$ and exactly $|A|$ copies of the symbol $'$. Thus, if $F(x) \le_e F(y)$, then $F(y)$ must be obtained by inserting only symbols from $A$ into the word $F(x)$. If any symbol $s \in A$ is inserted before its first occurence in $F(x)$, or inserted directly in front of a $'$, then we can see that the resulting word can't be of the form $F(y)$, by considering the first location with an invalid insertion.
\end{proof}

\begin{thm} If $A$ is a finite set, then the partial order $\le_e$ on $A^+$ is a well partial order.
\end{thm}
\begin{proof} We prove this by induction on $|A|$. Since $\le_e$ clearly has no infinite descending chains (as $a <_e b$ implies $|a| < |b|$), we just need to prove that $\le_e$ has no infinite antichains. Suppose for contradiction that $\le_e$ has an infinite antichain, and let $x_1, x_2, ...$ be a lexicographically minimal infinite antichain, that is, suppose that $x_1$ is minimal such that there exists an infinite antichain containing $x_1$, that $x_2$ is minimal such that there exists an infinite antichain containing $\{x_1, x_2\}$, etc.

By the infinite pigeonhole principle, we see that there is an infinite subsequence $i_1 < i_2 < \cdots$ such that every element $x_{i_j}$ ends in the same element of $A$, say $a$. Let $x_{i_j}'$ be the element of $A^+$ we obtain by deleting the $a$ in the last coordinate of $x_{i_j}$, then from the definition of $\le_e$ we see that $x_{i_j} \le_e x_{i_k} \iff x_{i_j}' \le_e x_{i_k}'$. Let $j$ be minimal such that $x_j >_e x_{i_k}'$ for some $k$, and note that $j \le i_1$ so $j$ is well-defined. Then the sequence
\[
x_1, x_2, ..., x_{j-1}, x_{i_k}', x_{i_{k+1}}', ...
\]
is also an infinite antichain, and is lexicographically smaller than $x_1, x_2, ..., x_{j-1}, x_j, ...$, a contradiction.
\end{proof}

\begin{cor} If $A$ is a finite set, then the partial order $\le_E$ on $A^+$ is a well partial order.
\end{cor}

\begin{thm}[Few subpowers implies inherently finitely related \cite{few-subpowers-finitely-related}]\label{few-finite} If a clone $\cO$ contains a $k$-edge term, then it is finitely related. In fact, a set $\Gamma \subseteq \Inv(\cO)$ generates $\Inv(\cO)$ iff the following two conditions are satisfied:
\begin{itemize}
\item every relation of arity strictly less than $k$ in $\Inv(\cO)$ is contained in $\langle \Gamma \rangle$, and
\item for each minority pair $(c,d) \in A^2$ and each minimal element $a \in \lambda(\cO,(c,d))$, if we set $n = |a|$, then the relation $\Pol_n(\Gamma)$ on $\bA^{\bA^n}$ defined by the primitive positive formula
\[
\bigwedge_{R \in \Gamma} \bigwedge_{M \in R^n} f(M) \in R
\]
has $(c,d) \not\in \Forks(\Pol_n(\Gamma),a)$.
\end{itemize}
\end{thm}
\begin{proof} By the fact that $\le_E$ is a well partial order, we see that there is a finite set $\Gamma \subseteq \Inv(\cO)$ which satisfies the conditions given: for instance, we may take $\Gamma$ to consist of the collection of all relations in $\Inv(\cO)$ of arity less than $k$, together with the relations $\cO_n \le \bA^{\bA^n}$ for every $n$ such that some minimal element $a \in \lambda(\cO,(c,d))$ has $|a| = n$ for some $(c,d) \in A^2$.

Now suppose that $\Gamma$ satisfies the given conditions. Then for any $(c,d) \in A^2$ and any $b \in \lambda(\cO,(c,d))$, there exists some minimal $a \in \lambda(\cO,(c,d))$ with $a \le_E b$. Thus if $|a| = n, |b| = m$, then $\Forks(\Pol_m(\Gamma),b) \subseteq \Forks(\Pol_n(\Gamma),a)$, and by the second condition on $\Gamma$ we have $(c,d) \not\in \Forks(\Pol_n(\Gamma),a)$. Thus for any $b \in A^+$ with $|b| = m$, we have
\[
(c,d) \not\in \Forks(\cO_m,b) \implies (c,d) \not\in \Forks(\Pol_m(\Gamma),b),
\]
so $\Forks(\Pol_m(\Gamma),b) \subseteq \Forks(\cO_m,b)$. Since $\Pol_m(\Gamma)$ contains $\cO_m$ (by $\Gamma \subseteq \Inv(\cO)$), and every projection of $\Pol_m(\Gamma)$ onto fewer than $k$ coordinates of $\bA^{\bA^m}$ is contained in the corresponding projection of $\cO_m$ (by the first condition on $\Gamma$), we can apply Theorem \ref{edge-gen} to see that $\Pol_m(\Gamma) = \cO_m$.
\end{proof}

\begin{cor} The number of clones on a finite set which contain an edge term is countable.
\end{cor}

\begin{rem} There is a converse to Theorem \ref{few-finite}: if $\cO$ is a clone on a finite set such that every clone $\cO'$ with $\cO' \supseteq \cO$ is finitely related, then $\cO$ has an edge term. The proof of this relies on the theory of \emph{cube term blockers} (aka \emph{projective subalgebras}) from  \cite{cube-term-blockers}, which roughly states that a clone $\cO$ fails to contain a cube term iff there is an infinite sequence of invariant relations which look like the relations $\{0,1\}^n\setminus\{(0,...,0)\}$ - recall that the clone corresponding to this sequence of relations on $\{0,1\}$ was our basic example of a clone which was not finitely related (Example \ref{ex-non-finitely-related}).
\end{rem}

\begin{ex} Consider the algebra $\bA = (\{a,b,c\},g)$ from Example \ref{ex-few-subpowers}, which has $\{a,b\}$ a majority subalgebra and $\{a,c\}$ an absorbing minority subalgebra. Recall that the minority pairs of $\bA$ were $(a,c), (c,a), (b,c)$. Since $a \in \Sg_{\bA}\{b,c\}$, for any $s \in A^+$ we have
\[
(b,c) \in \Forks(\langle g\rangle, s) \implies (a,c) \in \Forks(\langle g\rangle, s).
\]
Take the standard alphabetical ordering $<$ on $\{a,b,c\}$. It's easy to check that $\lambda(\langle g\rangle,(a,c))$ contains $a, b, c, ab, ba, bc, ca, cb, abc, acb, acc, bac, bca, cab, cba$ and that $\lambda(\langle g \rangle, (b,c)) = A^+$: for the strings of length $2$, the free algebra $\cF_{\bA}(x,y)$ only has six elements so we may compute the forks directly, for permutations of $abc$ we note that a corresponding permutation of $aac$ comes before it and $g$ preserves the congruence corresponding to the partition $\{a,b\},\{c\}$, and for $acc$ we note that $aac$ and $aca$ come before it and that $g$ preserves the affine ternary relation $\{(a,a,a),(a,c,c),(c,a,c),(c,c,a)\}$.

To complete the description of $\lambda(\langle g\rangle,(a,c))$, we just need to check that for all $2 \le i \le n$, the word $s_{in} = a\cdots aca\cdots a \in A^+$ of length $n$ with a $c$ in the $i$th position and $a$s elsewhere has $(a,c) \in \Forks(\langle g\rangle, s_{in})$. For this, we take the terms $x_1$ and $g(x_1,x_1,x_i)$ in the free algebra, and check that they make a fork at $s_{in}$. For $s' <_{lex} s_{in}$, we have $s_1' = a$ and $s_i' < c$, so
\[
g(s_1',s_1',s_i') = g(a,a,a)\text{ or }g(a,a,b) = a = s_1',
\]
so $x_1$ and $g(x_1,x_1,x_i)$ agree on tuples which come lexicographically before $s_{in}$. At $s_{in}$, we get the fork $(a,g(a,a,c)) = (a,c)$.
\end{ex}

\begin{ex} Consider the gmm algebra $\bA_2 = (\{a,b,c\},\varphi_2)$ from Example \ref{gmm-ex-2}, which had majority subalgebras $\{a,b\}, \{a,c\}$ and minority subalgebra $\{b,c\}$. The only minority pair to worry about is $(b,c)$, and under the standard alphabetical ordering $<$ on $\{a,b,c\}$, we find that $\lambda(\langle \varphi_2 \rangle, (b,c))$ contains the following $16$ elements of $A^+$:
\[
a, b, c, ab, ac, ba, ca, cb, acb, bcc, cab, cba, abcc, bacc, bcac, bcca.
\]
Again, it is easy to check the strings of length $2$ as $\cF_{\bA_2}(x,y)$ only has $4$ elements, strings which have a $c$ preceding the first $b$ such as $acb$ don't work because the corresponding word with $b$s and $c$s swapped (i.e. $abc$ in this case) comes before it and $\varphi_2$ preserves order two the automorphism swapping $b$ and $c$, and strings containing $bcc$ such as $abcc$ don't work because the two strings where one of the $c$s is replaced by a $b$ (i.e. $abbc$ and $abcb$ in this case) come before it and $\varphi_2$ preserves the ternary relation corresponding to the columns of the matrix
\[
\begin{bmatrix} a & b & b & c & c\\ a & b & c & b & c\\ a & b & c & c & b\end{bmatrix}.
\]

It's much harder to show that the remaining elements $s$ which are not $\ge_E$ to one of the $16$ strings displayed above all have $(b,c) \in \Forks(\langle \varphi_2 \rangle, s)$. Each such $s$ has at least one $b$, exactly one $c$, and has its first $b$ before its $c$. We may assume without loss of generality that $s$ begins with a $b$, and suppose $s$ has its only $c$ at the $i$th position for some $i \ge 2$. We need to show that there is some term $t \in \cF_{\bA_2}(x_1, ..., x_n)$ such that the pair $(x_1, t)$ gives us a fork at $s$. In other words, we need to show that we can find a term $t$ such that for each $s' <_{lex} s$ we have $t(s') = s_1'$ and $t(s) = c$.

The only way I know to show the existence of such a term $t$ is to use the analysis of critical relations in $\Inv_k(\bA_2)$ carried out in Example \ref{gmm-ex-2-critical}. By that analysis, we see that every relation $\RR \le \bA_2^k$ is the intersection of some family of binary relations and some family of relations $\RR_I \le \bA^I$ such that for each $I$ and each $i,j \in I$, we have $\pi_{i,j}(\RR_I) \subseteq 0_{\bA_2}^*$, where $0_{\bA_2}^*$ is the congruence corresponding to the partition $\{a\}, \{b,c\}$. Thus, if the term $t$ we are looking for does not exist, then either there is some $s' <_{lex} s$ such that
\[
\begin{bmatrix} s_1'\\ c\end{bmatrix} \not\in \Sg_{\bA_2^2}\begin{bmatrix}s' \\ s\end{bmatrix},
\]
or there is some family $s^1, ..., s^k <_{lex} s$ such that for each $j,l$, we have $(s_j^l, s_j) \in 0_{\bA_2}^*$ but
\[
\begin{bmatrix} s_1^1\\ \vdots \\ s_1^k\\ c\end{bmatrix} \not\in \Sg_{\bA_2^2}\begin{bmatrix} s^1 \\ \vdots \\ s^k\\ s\end{bmatrix}.
\]
To rule out the first possibility, we note that if $s' <_{lex} s$ then $s_1' \in \{a,b\}$, and if $(s_1',s_i') \ne (b,c)$, then $(s_1',s_i')$ is a majority pair and taking $\varphi_2(x_1,x_1,x_i)$ does the trick, while if $(s_1',s_i') = (b,c)$, then at the first location $j$ where $s'$ and $s$ differ we must have $s_j' = a, s_j = b$, so taking $\varphi_2(x_j, x_1, x_i)$ does the trick:
\[
\varphi_2\left(\begin{bmatrix} a & b & c\\ b & b & c\end{bmatrix}\right) = \begin{bmatrix} b\\ c\end{bmatrix}.
\]
To rule out the second possibility, note that if $s^l <_{lex} s$ and $(s_j^l, s_j) \in 0_{\bA_2}^*$ for all $j$, then the first coordinate where $s^l$ and $s$ can differ is at the coordinate $i$ with $s_i = c$, so we must have $s_i^l = b$, $s_i = c$ and $s_1^l = s_1 = b$. Thus the term $x_i$ rules out the second possibility.
\end{ex}

\section{Cube term blockers: projective subalgebras and crosses}

The goal of this section is to describe a criterion for the existence of an edge term operation of low arity, based on the nonexistence of certain families of relations or special types of subalgebras. The forbidden type of subalgebra was originally called a cube term blocker in \cite{cube-term-blockers}, but now it is slightly more common to see it referred to as a projective subalgebra.

\begin{defn}[\cite{cube-term-blockers}]\label{defn-projective} If $\PP \le \bA$ is a subalgebra such that for each term $t \in \Clo(\bA)$ there is some coordinate $i(t) \le \arity(t)$ such that
\[
\forall x_1, .., x_{\arity(t)}\in \bA, \;\; x_{i(t)} \in \PP \;\;\; \implies \;\;\; t(x_1, ..., x_{\arity(t)}) \in \PP,
\]
then we say that $\PP$ is a \emph{projective subalgebra} of $\bA$. Following \cite{cube-terms-crosses}, we say that $t$ is \emph{$\PP$-absorbing in the $i$th coordinate} if the above implication holds for $t$ and $i = i(t)$.

We say that a pair $(\PP,\bB)$ of subalgebras with $\PP < \bB \le \bA$ is a \emph{cube term blocker} for $\bA$ if $\PP$ is a proper projective subalgebra of $\bB$.
\end{defn}

\begin{rem} The name ``projective subalgebra'' refers to the way every term acts somewhat like a projection, when it comes to deciding whether or not its output will be contained in $\PP$. Note that it has nothing to do with the concept of a projective module from ring theory (or to the category-theoretic generalization).
\end{rem}

Although the definition of a projective subalgebra refers to all terms $t \in \Clo(\bA)$, it's actually enough to check the condition for the basic operations of $\bA$.

\begin{prop}\label{prop-projective-basic} If $\bA = (A, \{f_j\})$, then $\PP \le \bA$ is a projective subalgebra of $\bA$ iff for each $j$, there is some $i(f_j) \le \arity(f_j)$ such that
\[
x_{i(f_j)} \in \PP \;\;\; \implies \;\;\; f_j(x_1, ..., x_{\arity(f_j)}) \in \PP.
\]
\end{prop}
\begin{proof} We just need to check that the collection of operations $t$ which have a $\PP$-absorbing coordinate $i(t)$ is closed under composition. If we have already found $\PP$-absorbing coordinates $i(f)$ and $i(g_j)$ for $j = 1, ..., \arity(f)$, and if we define a new operation $h$ by
\[
h = f \circ (g_1, ..., g_{\arity(f)}),
\]
then we have
\begin{align*}
x_{i(g_{i(f)})} \in \PP \;\;\; &\implies \;\;\; g_{i(f)}(x_1, ..., x_{\arity(h)}) \in \PP\\
&\implies \;\;\; f(g_1(x_1, ..., x_{\arity(h)}), ..., g_{i(f)}(x_1, ..., x_{\arity(h)}), ..., g_{\arity(f)}(x_1, ..., x_{\arity(h)})) \in \PP,
\end{align*}
so we can take $i(h) = i(g_{i(f)})$.
\end{proof}

Projective subalgebras can be described in terms of the existence of an infinite family of high-arity relations which look like cubes with a corner removed (such relations are referred to as ``symmetric crosses'' in \cite{cube-terms-crosses}, or as ``chipped cubes'' in \cite{cube-terms-chipped-cubes}).

\begin{prop}\label{prop-cube-term-blocker-relational} For any $\PP < \bB \le \bA$, the pair $(\PP,\bB)$ is a cube term blocker for $\bA$ (that is, $\PP$ is a projective subalgebra of $\bB$) if and only if for every $n \ge 1$, the $n$-ary relation $R_n$ defined by
\[
(x_1, ..., x_n) \in R_n \;\;\; \iff \;\;\; \bigwedge_{i \in [n]} x_i \in \bB \wedge \bigvee_{j \in [n]} x_j \in \PP,
\]
or equivalently by
\begin{align*}
R_n &= (\PP \times \bB^{n-1}) \cup (\bB \times \PP \times \bB^{n-2}) \cup \cdots \cup (\bB^{n-1} \times \PP)\\
&= \bB^n \setminus (\bB \setminus \PP)^n,
\end{align*}
is preserved by all of the term operations of $\bA$ (that is, $R_n$ is a subuniverse of $\bA^n$).
\end{prop}
\begin{proof} Let $t$ be an $m$-ary term operation of $\bA$. If $\PP$ is a projective subalgebra of $\bB$, then there must be some $\PP$-absorbing (with respect to $\bB$) coordinate $i(t)$, i.e.
\[
t(\bB, ..., \PP, ..., \bB) \subseteq \PP,
\]
where the $\PP$ occured in the $i(t)$th coordinate. Then for any $x^1, ..., x^m \in R_n$, if $x^{i(t)}$ has $x_j^{i(t)} \in \PP$, then we have
\[
t(x_j^1, ..., x_j^{i(t)}, ..., x_j^m) \in \PP
\]
and
\[
t(x_i^1, ..., x_i^m) \in \bB
\]
for all $i \ne j$, so
\[
t(x^1, ..., x^m) \in \bB^{j-1}\times\PP\times \bB^{n-j} \subseteq R_n.
\]
Thus $t$ preserves $R_n$ as long as $\PP$ is a projective subalgebra of $\bB$.

Now suppose that $\PP$ is not a projective subalgebra of $\bB$, with $t$ some $m$-ary term operation having no corresponding $\PP$-absorbing (with respect to $\bB$) coordinate $i(t)$. Then for each $i \le m$, there must be some tuple
\[
x_i^1, ..., x_i^m \in \bB, \;\;\; x_i^i \in \PP
\]
such that
\[
t(x_i^1, ..., x_i^m) \not\in \PP.
\]
Then we have
\[
\begin{bmatrix} x_1^i\\ \vdots\\ x_m^i\end{bmatrix} \in \bB^{i-1}\times\PP\times \bB^{n-i} \subseteq R_m
\]
for each $i$, but
\[
t\left(\begin{bmatrix} x_1^1\\ \vdots\\ x_m^1\end{bmatrix}, \cdots, \begin{bmatrix} x_1^m\\ \vdots\\ x_m^m\end{bmatrix}\right) \in (\bB\setminus\PP)^m,
\]
so $t$ does not preserve $R_m$.
\end{proof}

\begin{cor} A subalgebra $\PP \le \bA$ is projective iff for each $n \ge 1$, the $n$-ary relation
\[
\bigvee_{i \in [n]} x_i \in \PP
\]
is preserved by the basic operations of $\bA$.
\end{cor}

\begin{rem} If $R_n = B^n \setminus (B \setminus P)^n$ is a subuniverse of $\bA^n$ for some $n \ge 2$, then (as long as $P \subseteq B$) it will automatically be the case that $B$ and $P$ are also subuniverses of $\bA$, by
\[
B = \pi_1(R_n)
\]
and
\[
x \in P \;\;\; \iff \;\;\; (x,x,...,x) \in R_n.
\]
As a consequence, there was no loss of generality in restricting our definition of cube term blockers to the case where $P$ and $B$ are subuniverses of $\bA$.
\end{rem}

Now we finally explain the reason for the name ``cube term blocker''.

\begin{prop}[\cite{cube-term-blockers}]\label{prop-cube-term-blocked} If an algebra $\bA$ has a cube term blocker $\PP < \bB \le \bA$, then $\bA$ has no cube term.
\end{prop}
\begin{proof} Pick any $a \in \PP$ and any $b \in \bB \setminus \PP$. If $\bA$ had a $k$-cube term, then we would have
\[
(b, ..., b) \in \Sg_{\bA^k}\big(\{a,b\}^k \setminus \{(b, ..., b)\}\big),
\]
but by Proposition \ref{prop-cube-term-blocker-relational} the right hand side is contained in $\bB^k \setminus (\bB \setminus \PP)^k$, which does not contain $(b, ..., b)$.
\end{proof}

\begin{prop} If an algebra $\bA$ has a cube term blocker $\PP < \bB \le \bA$, then there is an expansion $\hat{\bA}$ of $\bA$ (i.e., an algebra with the same underlying set whose basic operations are a superset of the basic operations of $\bA$) such that $\hat{\bA}$ is \emph{not} finitely related. If $\bA$ is idempotent then we can take $\hat{\bA}$ to be idempotent as well.
\end{prop}
\begin{proof} Let the basic operations of $\hat{\bA}$ be the collection of all operations which preserve the relations
\[
R_n = \bB^n \setminus (\bB \setminus \PP)^n
\]
for all $n$, and which also preserve all unary relations which are preserved by $\bA$. Then $\hat{\bA}$ is an expansion of $\bA$, which is idempotent iff $\bA$ is idempotent, so we just need to check that $\hat{\bA}$ is not finitely related.

By the Galois correspondence between $\Inv$ and $\Pol$, the relational clone $\Inv(\hat{\bA})$ is generated by the infinite sequence of relations $R_n$ together with some collection of unary relations. Thus, using
\[
(x_1, ..., x_k) \in R_k \;\;\; \iff \;\;\; (x_1, ..., x_k, \underbrace{x_k, ..., x_k}_{m-k}) \in R_m
\]
for $k \le m$, we see that every finite collection of relations in $\Inv(\hat{\bA})$ is in the relational clone generated by some $R_m$ together with some unary relations. Therefore we just need to check that for each $m$, $R_{m+1}$ is \emph{not} contained in the relational clone $\Gamma_m$ generated by $R_m$ together with all unary relations.

In order to show that $R_{m+1}$ is not contained in the relational clone $\Gamma_m$, we just need to find an operation $t$ which preserves (a generating set of) the relations in $\Gamma_m$ but which does not preserve $R_{m+1}$. One such operation is the $m+1$-ary operation given by
\[
t(x_1, ..., x_{m+1}) = \begin{cases} x_1 & \text{if }x_1 \not\in \PP\text{ and at most one } x_i \in \PP,\\ x_2 & \text{if }x_1 \in \PP\text{ and no other }x_i \in \PP,\\ x_j & \text{if }i < j\text{ are the first two indices  s.t. }x_i, x_j \in \PP.\end{cases}
\]
Since the output of $t$ is always one of the inputs to $t$, every unary relation $S \subseteq \bA$ will be preserved by $t$. That $t$ preserves $R_m$ follows from the pigeonhole principle: among any $m+1$ tuples $x^1, ..., x^{m+1} \in R_m$, there must be some pair $i < j \le m+1$ and some coordinate $k \le m$ such that $x^i_k, x^j_k \in \PP$, in which case we have $t(x_k^1, ..., x_k^{m+1}) \in \PP$ as well.

To see that $t$ doesn't preserve $R_{m+1}$, pick any $m+1$ tuples $x^1, ..., x^{m+1}$ such that
\[
x^i \in (\bB \setminus \PP)^{i-1} \times \PP \times (\bB \setminus \PP)^{m+1-i} \subset R_{m+1}
\]
for each $i$, and note that we then have $t(x^1, ..., x^{m+1}) \in (\bB \setminus \PP)^{m+1}$, so $t(x^1, ..., x^{m+1}) \not\in R_{m+1}$.
\end{proof}

One of the main results of \cite{cube-term-blockers} is the converse of Proposition \ref{prop-cube-term-blocked} (for idempotent algebras). Since Corollary \ref{cor-few-subpowers-product} shows that any product of idempotent algebras with few subpowers has few subpowers, the same should be true of idempotent algebras which have no cube term blockers. This will follow from a few basic results about projective subalgebras.

\begin{prop} If $\PP \le \bB \le \bA$ and $\PP$ is a projective subalgebra of $\bB$, then
\begin{itemize}
\item[(a)] for any $n$, $\PP^n$ is a projective subalgebra of $\bB^n$,
\item[(b)] for any $\theta \in \Con(\bA)$, the quotient $\PP/\theta$ is a projective subalgebra of $\bB/\theta$,
\item[(c)] for any $\bC \le \bA$ with $\bC \cap \PP \ne \emptyset$, the intersection $\PP \cap \bC$ is a projective subalgebra of $\bB \cap \bC$, and
\item[(d)] for any binary relation $\RR \le \bA \times \bD$, if we define $\PP + \RR$ by
\[
\PP + \RR = \pi_2(\RR \cap (\PP \times \bD))
\]
and similarly for $\bB + \RR$, then $\PP + \RR$ is a projective subalgebra of $\bB + \RR$.
\end{itemize}
Note that as a special case of (d), if there is a homomorphism $\pi : \bD \rightarrow \bA$, then $\pi^{-1}(\PP)$ is a projective subalgebra of $\pi^{-1}(\bB)$.
\end{prop}
\begin{proof} In each case, if we are trying to show that $\PP'$ is a projective subalgebra of $\bB'$, the strategy is to show that if the restriction of some term operation $t$ to $\bB$ is $\PP$-absorbing in its $i$th coordinate, then the restriction of $t$ to $\bB'$ is $\PP'$-absorbing in the (same) $i$th coordinate.

We will prove (d), since the other cases are easier. Supposing that $t^\bB$ is $\PP$-absorbing in its first coordinate, we will prove that $t^{\bB+\RR}$ is also $\PP+\RR$-absorbing in its first coordinate. For this, note that for any $q \in \PP+\RR$ we can (by the definition of $\PP+\RR$) find some $p \in \PP$ with $(p,q) \in \RR$, and similarly for $d_2, ... \in \bB+\RR$ there are $b_i \in \bB$ with $(b_i, d_i) \in \RR$, so
\[
t\left(\begin{bmatrix} p & b_2 & \cdots\\ q & d_2 & \cdots \end{bmatrix}\right) \in \RR,
\]
so
\begin{align*}
t(q, d_2, ...) &\in t(p, b_2, ...) + \RR \subseteq \PP + \RR.\qedhere
\end{align*}
\end{proof}

\begin{cor} Suppose that $\PP \le \bB \le \bA$ and $\PP$ is a projective subalgebra of $\bB$. If a relation $\RR \le \bA^m$ is defined by a primitive positive formula $\Phi(x_1, ..., x_m)$ where some of the (free or bound) variable domains are $\bB$, and if $\bS \le \RR$ is defined by modifying the formula $\Phi$ by replacing some of those variable domains with $\PP$, then $\bS$ is a projective subalgebra of $\RR$.
\end{cor}

\begin{rem} It is \emph{not} always the case that if $\PP$ is a projective subalgebra of $\bB$ and $\bQ$ is a projective subalgebra of $\bC$, then $\PP \times \bQ$ is a projective subalgebra of $\bB \times \bC$. The reason is that $t^\bB$ might be $\PP$-absorbing in its $i$th coordinate (only) and $t^\bC$ might be $\bQ$-absorbing in its $j$th coordinate (only) for $i \ne j$. For instance, consider the case where $\bA = \bB \times \bC$ has just one basic operation $t(x,y)$ such that $t^\bB$ is first projection and $t^\bC$ is second projection (such an algebra $\bA$ is known as a \emph{rectangular band}).

Similarly, it is easy to construct an example of $\bQ \le \PP \le \bB$ such that $\bQ$ is a projective subalgebra of $\PP$, $\PP$ is a projective subalgebra of $\bB$, but $\bQ$ is not a projective subalgebra of $\bB$.
\end{rem}

\begin{prop}[\cite{cube-terms-crosses}] If $\{\bA_i\}$ is a collection of idempotent algebras and $\bA \in HSP_{fin}(\{\bA_i\})$ has a cube term blocker, then some $\bA_i$ has a cube term blocker.
\end{prop}
\begin{proof} The previous results immediately show that any subalgebra or quotient of an algebra with no cube term blockers also has no cube term blocker, so we just need to consider the case where $\bA = \bA_1 \times \cdots \times \bA_n$. Suppose that
\[
\PP < \bB \le \bA_1 \times \cdots \times \bA_n
\]
and that $\PP$ is a projective subalgebra of $\bB$. If $\pi_n(\PP) \ne \pi_n(\bB)$, then $\pi_n(\PP)$ is a proper projective subalgebra of $\pi_n(\bB)$, so in this case $\bA_n$ has a cube term blocker.

Otherwise, if $\pi_n(\PP) = \pi_n(\bB)$, then there must be some $b_n \in \pi_n(\bB)$ such that
\[
\pi_{[n-1]}(\PP \cap (\bA_1 \times \cdots \times \bA_{n-1} \times \{b_n\})) \ne \pi_{[n-1]}(\bB \cap (\bA_1 \times \cdots \times \bA_{n-1} \times \{b_n\})),
\]
in which case (using the fact that $\{b_n\}$ is a subalgebra of $\bA_n$) the left hand side is a proper projective subalgebra of the right hand side, so we see by induction on $n$ that one of $\bA_1, ..., \bA_{n-1}$ has a cube term blocker.
\end{proof}

\begin{cor}[\cite{cube-terms-crosses}]\label{cor-free-cube-term-blocker} If $\bA$ is a finite idempotent algebra with no cube term blocker, then the free algebra on two generators $\cF_{\cV(\bA)}(x,y)$ in the variety $\cV(\bA)$ generated by $\bA$ also has no cube term blocker.
\end{cor}

Putting the previous corollary together with the next result, we can show that every finite idempotent algebra with no cube term blocker has a cube term.

\begin{lem}[\cite{cube-term-blockers}]\label{lem-partial-cube-term} If $\bA$ is a finite idempotent algebra with no cube term blocker, then for every pair $a \ne b \in \bA$ there is some $n$ such that
\[
(b, ..., b) \in \Sg_{\bA^n}\big(\{a,b\}^n\setminus\{(b, ..., b)\}\big).
\]
\end{lem}
\begin{proof} Suppose for the sake of contradiction that for each $n$ we have $(b, ..., b) \not\in \Sg_{\bA^n}(\{a,b\}^n \setminus \{(b, ..., b)\})$. Let $n$ be large, and pick a tuple $(b_1, ..., b_n) \in (\bA\setminus\{a\})^n$ with
\[
\sum_{i \le n} |\Sg_\bA\{a,b_i\}|
\]
minimal among $n$-tuples satisfying
\[
(b_1, ..., b_n) \not\in \Sg_{\bA^n}\Big(\prod_{i\le n} \{a,b_i\} \setminus \{(b_1, ..., b_n)\}\Big).
\]
For each $i \le n$, define $P_i \subset \bB_i \le \bA$ by
\begin{align*}
\bB_i &= \Sg_\bA\{a,b_i\},\\
P_i &= \{p \in \bB_i \mid \Sg_{\bA}\{a,p\} \ne \bB_i\}.
\end{align*}
Writing
\[
\RR = \Sg_{\bA^n}\Big(\prod_{i\le n} \{a,b_i\} \setminus \{(b_1, ..., b_n)\}\Big),
\]
we claim that
\[
(x_1, ..., x_n) \in \RR \;\;\; \stackrel{?}{\iff} \;\;\; \bigwedge_{i \le n} x_i \in \bB_i \wedge \bigvee_{i \le n} x_i \in P_i.
\]
To see this, first note that for any $c_i \in \bB_i$ we have
\[
\Big(\prod_{i\le n}\{a,c_i\}\Big) \setminus \{(c_1, ..., c_n)\} \subseteq \RR,
\]
since for each $j \le n$ we have
\[
\prod_{i < j} \{a,c_i\} \times \{a\} \times \prod_{k > j} \{a,c_k\} \subseteq \Sg_{\bA^n}\Big(\prod_{i < j} \{a,b_i\} \times \{a\} \times \prod_{k \ge j}\{a,b_k\}\Big) \subseteq \RR.
\]
Thus by the choice of $(b_1, ..., b_n)$, for $c_i \in \bB_i$ we have
\[
(c_1, ..., c_n) \not\in \RR \;\;\; \implies \;\;\; \bigwedge_{i \le n} \Sg_\bA\{a,c_i\} = \bB_i,
\]
so
\[
\RR \supseteq \Big(\prod_{i \le n} \bB_i\Big) \setminus \Big(\prod_{i \le n} (\bB_i \setminus P_i)\Big).
\]
For the other containment, suppose for the sake of contradiction that there is a tuple $(c_1, ..., c_n) \in \RR$ with $c_i \in \bB_i \setminus P_i$ for all $i$. In this case, we will show by induction on $i$ that we have
\[
(b_1, ..., b_i, c_{i+1}, ..., c_n) \in \RR
\]
for $i = 0, ..., n$. For the inductive step, from $(b_1, ..., b_{i-1}, c_i, ..., c_n) \in \RR$ and
\[
(b_1, ..., b_{i-1}, a, c_{i+1}, ..., c_n) \in \RR,
\]
we see that
\[
(b_1, ..., b_{i-1}, \Sg_{\bA}\{a, c_i\}, c_{i+1}, ..., c_n) \subseteq \RR,
\]
and by the definition of $P_i$ we have
\[
c_i \in \bB_i \setminus P_i \;\; \implies \;\; b_i \in \Sg_{\bA}\{a,c_i\}.
\]
This completes the inductive step, and taking $i = n$ gives us a contradiction to the assumption $(b_1, ..., b_n) \not\in \RR$, so we must in fact have
\[
\RR = \Big(\prod_{i \le n} \bB_i\Big) \setminus \Big(\prod_{i \le n} (\bB_i \setminus P_i)\Big).
\]
Now assuming that $n > (m-1)(|\bA|-1)$, we see by the pigeonhole principle that there must be some $b'$ such that $b_i = b'$ for at least $m$ different choices of $i$. Rearranging the coordinates, we may assume without loss of generality that $b_1 = \cdots = b_m = b'$, so
\[
\RR_m = \pi_{[m]} \big(\RR \cap \big(\bB_1^m \times \{(b_{m+1}, ..., b_n)\}\big)\big)
\]
is given by
\[
\RR_m = \bB_1^m \setminus (\bB_1 \setminus P_1)^m,
\]
and is a subalgebra of $\bA^m$ (note that $P_1$ is automatically a subalgebra of $\bA$ as a consequence). Since we can take $m$ arbitrarily large, there must be some pair $P \subset \bB$ which shows up infinitely often as $m$ grows, and this pair is a cube term blocker for $\bA$.
\end{proof}

\begin{thm}[Markovi{\'c}, Mar{\'o}ti, McKenzie \cite{cube-term-blockers}]\label{thm-cube-term-blockers} If a finite idempotent algebra $\bA$ has no cube term blocker, then $\bA$ has few subpowers.
\end{thm}
\begin{proof} By Corollary \ref{cor-free-cube-term-blocker}, the free algebra $\cF_{\cV(\bA)}(x,y)$ has no cube term blocker. Thus by Lemma \ref{lem-partial-cube-term}, there is some $n \ge 1$ such that
\[
(x, ..., x) \in \Sg_{\cF_{\cV(\bA)}(x,y)^n}\big(\{x,y\}^n \setminus \{(x, ..., x)\}\big),
\]
and any term $t$ witnessing this is an $n$-cube term.
\end{proof}

\begin{cor}[Markovi{\'c}, Mar{\'o}ti, McKenzie \cite{cube-term-blockers}]\label{cor-inherently-finitely-related} A finite idempotent algebra $\bA$ has few subpowers if and only if every (idempotent) expansion of $\bA$ is finitely related.
\end{cor}

One of the reasons for the interest in Corollary \ref{cor-inherently-finitely-related} was that it was seen as indirect support for the following conjecture due to Matthew Valeriote, which was eventually settled by Barto \cite{barto-valeriote-conjecture}.

\begin{conj}[Valeriote's Edinburgh conjecture \cite{csp-comparing-solution-sets}] A finite algebra $\bA$ has few subpowers if and only if $\bA$ is finitely related and generates a congruence modular variety.
\end{conj}

Using Theorem \ref{thm-cube-term-blockers}, Kazda and Zhuk \cite{cube-terms-chipped-cubes} were able to design a simple polynomial time algorithm (Algorithm \ref{alg-projective}) for checking whether a given idempotent algebraic structure (described by listing out the tables of the basic operations) has few subpowers.

\begin{algorithm}
\caption{Algorithm for finding a cube term blocker $\PP < \bB \le \bA$ with $a \in \PP$, from \cite{cube-terms-chipped-cubes}.}\label{alg-projective}
\begin{algorithmic}[1]
\State Set $S \gets \{a\}$.
\While{$S \ne \bA$}
\State Find $b \in \bA \setminus S$ minimizing $|\Sg_\bA\{a,b\}|$.
\If{$S \cap \Sg_\bA\{a,b\}$ is a projective subalgebra of $\Sg_\bA\{a,b\}$} \Comment{Proposition \ref{prop-projective-basic}}
\State Set $\PP \gets S \cap \Sg_\bA\{a,b\}$ and $\bB \gets \Sg_\bA\{a,b\}$.
\State \Return $(\PP,\bB)$.
\Else
\State Set $S \gets S \cup \Sg_\bA\{a,b\}$.
\EndIf
\EndWhile
\State \Return ``No cube term blockers have $a \in \PP$.''
\end{algorithmic}
\end{algorithm}

\begin{prop} Algorithm \ref{alg-projective} correctly decides whether there is a cube term blocker $\PP < \bB \le \bA$ with $a \in \PP$.
\end{prop}
\begin{proof} We just need to prove that at every step of the algorithm, for every element $c \in S$ the subalgebra $\Sg_\bA\{a,c\}$ has no proper projective subalgebra containing $a$.

Suppose that $b \in \bA \setminus S$ minimizes $|\Sg_\bA\{a,b\}|$. Inductively, we see that any projective subalgebra $\PP$ of $\Sg_\bA\{a,b\}$ which contains $a$ must contain every element $c \in S \cap \Sg_\bA\{a,b\}$, since $\PP \cap \Sg_\bA\{a,c\}$ is a projective subalgebra of $\Sg_\bA\{a,c\}$. By the choice of $b$, every $c \in \Sg_\bA\{a,b\} \setminus S$ must have $b \in \Sg_\bA\{a,c\}$, so we also see that any proper subalgebra of $\Sg_\bA\{a,b\}$ must contain $S \cap \Sg_\bA\{a,b\}$, so $S \cap \Sg_\bA\{a,b\}$ is the only possible proper projective subalgebra of $\Sg_\bA\{a,b\}$.
\end{proof}

Now let's turn our attention to trying to find a $k$-cube term with $k$ as small as possible. The key is the following local-to-global result from Horowitz's thesis \cite{horowitz-thesis}.

\begin{thm}[Horowitz \cite{horowitz-thesis}, \cite{horowitz-complexity-malcev-conditions}]\label{thm-local-cube-terms} A finite idempotent algebra $\bA$ has a $k$-cube term iff for every sequence of $k$ ordered pairs $a_i \ne b_i \in \bA$, we have
\[
(b_1, ..., b_k) \in \Sg_{\bA^k}\Big(\prod_{i \le k} \{a_i,b_i\} \setminus \{(b_1, ..., b_k)\}\Big).
\]
\end{thm}
\begin{proof} This is a generalization of the construction of a Mal'cev term out of local Mal'cev terms from Remark \ref{rem-malcev-construction}. Define $v_i^S(x,y)$ by
\[
v_i^S(x,y) \coloneqq \begin{cases} y & i \in S,\\ x & i \not\in S, \end{cases}
\]
and fix some ordering $S_1, ..., S_{2^k-1}$ of the nonempty subsets of $[k]$. If $(b_1, ..., b_k) \in \Sg(\prod_i \{a_i,b_i\} \setminus \{(b_1, ..., b_k)\})$, then there must be some (idempotent) $(2^k-1)$-ary term $t_{a_1b_1\cdots a_kb_k}$ such that
\[
t_{a_1b_1\cdots a_kb_k}(v_i^{S_1}(b_i,a_i), ..., v_i^{S_{2^k-1}}(b_i,a_i)) = b_i.
\]
Our goal is to build our way up to a $k$-edge term by composing the terms $t_{a_1b_1\cdots a_kb_k}$ with each other.

The strategy is to construct, by backward induction on $0 \le j < k$, terms $t_{a_1b_1\cdots a_jb_j}$ such that
\begin{align*}
i \le j \;\;\; &\implies \;\;\; t_{a_1b_1\cdots a_jb_j}(v_i^{S_1}(b_i,a_i), ..., v_i^{S_{2^k-1}}(b_i,a_i)) = b_i,\\
i > j \;\;\; &\implies \;\;\; t_{a_1b_1\cdots a_jb_j}(v_i^{S_1}(x,y), ..., v_i^{S_{2^k-1}}(x,y)) \approx x,
\end{align*}
so that once $j$ reaches $0$ we will have constructed a $k$-edge term. Fixing an enumeration $(c_1, d_1), ...$ of the ordered pairs of elements of $\bA$, we will construct $t_{a_1b_1\cdots a_jb_j}$ by inductively constructing a sequence of (idempotent) terms $t_{a_1b_1\cdots a_jb_j}^\ell$ satisfying
\begin{align*}
i \le j \;\;\; &\implies \;\;\; t_{a_1b_1\cdots a_jb_j}^\ell(v_i^{S_1}(b_i,a_i), ...) = b_i,\\
i > j+1 \;\;\; &\implies \;\;\; t_{a_1b_1\cdots a_jb_j}^\ell(v_i^{S_1}(x,y), ...) \approx x,\\
n \le \ell \;\;\; &\implies \;\;\; t_{a_1b_1\cdots a_jb_j}^\ell(v_{j+1}^{S_1}(d_n,c_n), ...) = d_n,
\end{align*}
starting with
\[
t_{a_1b_1\cdots a_jb_j}^0(x_1, ...) \coloneqq x_{\#\{j+1\}},
\]
where $\#\{j+1\}$ is defined by $S_{\#\{j+1\}} = \{j+1\}$. To see that this works for $\ell = 0$, note that
\[
i \ne j+1 \;\;\; \implies \;\;\; v_i^{\{j+1\}}(x,y) \approx x.
\]

For the inductive step, assume that we have already constructed $t_{a_1b_1\cdots a_jb_j}^{\ell-1}$, and let
\[
e_\ell = t_{a_1b_1\cdots a_jb_j}^{\ell-1}(v_{j+1}^{S_1}(d_\ell,c_\ell), ...).
\]
If we already have $e_\ell = d_\ell$, then we can take $t_{a_1b_1\cdots a_jb_j}^\ell \coloneqq t_{a_1b_1\cdots a_jb_j}^{\ell-1}$. Otherwise, we define $(2^k-1)$-ary terms $u^{S_m}$ by
\[
u^{S_m}(x_1, ...) \coloneqq \begin{cases} x_m & \text{if } j+1 \not\in S_m,\\ t_{a_1b_1\cdots a_jb_j}^{\ell-1}(x_1, ...) & \text{if } S_m = \{j+1\},\\ t_{a_1b_1\cdots a_jb_j}^{\ell-1}(v_{j+1}^{S_1}(x_w, x_m), ...) & \text{if } S_m = \{j+1\} \sqcup S_w, \end{cases}
\]
and we take
\[
t_{a_1b_1\cdots a_jb_j}^\ell \coloneqq t_{a_1b_1\cdots a_jb_je_\ell d_\ell} \circ (u^{S_1}, ...).
\]
To verify that this works, it's enough to check that we have
\begin{align*}
i \le j \;\;\; &\implies \;\;\; u^{S_m}(v_i^{S_1}(b_i,a_i), ...) = v_i^{S_m}(b_i,a_i),\\
i > j+1 \;\;\; &\implies \;\;\; u^{S_m}(v_i^{S_1}(x,y), ...) \approx v_i^{S_m}(x,y),\\
n \le \ell-1 \;\;\; &\implies \;\;\; u^{S_m}(v_{j+1}^{S_1}(d_n,c_n), ...) = d_n,
\end{align*}
and
\[
u^{S_m}(v_{j+1}^{S_1}(d_\ell,c_\ell), ...) = v_{j+1}^{S_m}(d_\ell, e_\ell).
\]
These are all automatic when $j+1 \not\in S_m$, and when $S_m = \{j+1\}$ they follow from the inductive hypothesis and the definition of $e_\ell$. The cases where $i \ne j+1$ and $S_m = \{j+1\} \sqcup S_w$ follow from
\[
i \ne j+1,\;\; S_m = \{j+1\} \sqcup S_w \;\;\; \implies \;\;\; v_i^{S_w}(x,y) \approx v_i^{S_m}(x,y)
\]
together with idempotence of $t_{a_1b_1\cdots a_jb_j}^{\ell-1}$. Finally, the cases where $n \le \ell$ and $S_m = \{j+1\} \sqcup S_w$ follow from
\[
S_m = \{j+1\} \sqcup S_w \;\;\; \implies \;\;\; v_{j+1}^{S_w}(d_n, c_n) = d_n, \;\; v_{j+1}^{S_m}(d_n,c_n) = c_n,
\]
and $v_{j+1}^{S_m}(d_\ell,e_\ell) = e_\ell = t_{a_1b_1\cdots a_jb_j}^{\ell-1}(v_{j+1}^{S_1}(d_\ell,c_\ell), ...)$.
\end{proof}

\begin{defn}[\cite{cube-terms-chipped-cubes}, \cite{cube-terms-crosses}] A relation $\RR \le \bA_1 \times \cdots \times \bA_k$ is called a $k$-dimensional \emph{chipped cube} if there are subalgebras $\PP_i < \bB_i \le \bA_i$ such that
\[
(x_1, ..., x_k) \in \RR \;\;\; \iff \;\;\; \bigwedge_{i \le k} x_i \in \bB_i \wedge \bigvee_{j \le k} x_j \in \PP_j,
\]
that is, if
\[
\RR = \Big(\prod_{i \le k} \bB_i\Big) \setminus \Big(\prod_{i \le k} \bB_i \setminus \PP_i\Big).
\]
Following \cite{cube-terms-crosses}, we say that a chipped cube is a \emph{cross} if we have $\bB_i = \bA_i$ for all $i \le k$, that is, if
\[
(x_1, ..., x_k) \in \RR \;\;\; \iff \;\;\; \bigvee_{i \le k} x_i \in \PP_i.
\]
\end{defn}

\begin{thm}[Kazda, Zhuk \cite{cube-terms-chipped-cubes}]\label{thm-chipped-cubes} A finite idempotent algebra $\bA$ has a $k$-cube term iff there is no $k$-dimensional chipped cube $\RR \le \bA^k$.
\end{thm}
\begin{proof} By Horowitz's Theorem \ref{thm-local-cube-terms}, we just need to check that if there is some sequence of $k$ ordered pairs $a_i \ne b_i \in \bA$ such that
\[
(b_1, ..., b_k) \not\in \Sg_{\bA^k}\Big(\prod_{i\le k} \{a_i,b_i\} \setminus \{(b_1, ..., b_k)\}\Big),
\]
then $\bA$ has a $k$-dimensional chipped cube. For this, we fix $(a_1, ..., a_k)$ and choose such a $(b_1, ..., b_k)$ such that
\[
\sum_{i \le k} |\Sg_\bA\{a_i,b_i\}|
\]
is minimized. Then by the same argument as the one in Lemma \ref{lem-partial-cube-term}, if we define $P_i, \bB_i$ by
\begin{align*}
\bB_i &= \Sg_\bA\{a_i,b_i\},\\
P_i &= \{p \in \bB_i \mid \Sg_\bA\{a_i,p\} \ne \bB_i\},
\end{align*}
then
\[
\RR \coloneqq \Sg_{\bA^k}\Big(\prod_{i\le k} \{a_i,b_i\} \setminus \{(b_1, ..., b_k)\}\Big) = \Big(\prod_{i \le k} \bB_i\Big) \setminus \Big(\prod_{i \le k} \bB_i \setminus P_i\Big).
\]
To see that $P_i$ is actually a subuniverse of $\bA$, note that we have
\[
p \in P_i \;\;\; \iff \;\;\; (b_1, ..., b_{i-1}, p, b_{i+1}, ..., b_k) \in \RR.\qedhere
\]
\end{proof}

Surprisingly, Kearnes and Szendrei \cite{cube-terms-crosses} showed that a version of this result holds in the infinite case, by using some special properties of free algebras.

\begin{thm}[Kearnes, Szendrei \cite{cube-terms-crosses}]\label{thm-infinite-cross} Let $\bA$ be a (possibly infinite) idempotent algebra, and let $\bF = \cF_{\cV(\bA)}(x,y)$ be the free algebra on two generators in the variety generated by $\bA$. Then
\begin{itemize}
\item[(a)] $\bA$ has a $k$-cube term iff $\bF$ has no $k$-dimensional cross $\RR < \bF^k$, and
\item[(b)] $\bA$ has a cube term iff $\bF$ has no proper projective subalgebra $\PP < \bF$.
\end{itemize}
\end{thm}
\begin{proof} We start with (a). Since $\bF$ is the free algebra on two generators, $\bA$ has a $k$-cube term iff
\[
(x, ..., x) \in \Sg_{\bF^k}\big(\{x,y\}^k \setminus \{(x, ..., x)\}\big).
\]
Defining $\operatorname{Cross}(P_1, ..., P_k)$ by
\[
(x_1, ..., x_k) \in \operatorname{Cross}(P_1, ..., P_k) \;\;\; \iff \;\;\; \bigvee_{i \le k} x_i \in P_i,
\]
we see that if $\bA$ does not have a $k$-cube term then
\[
(x, ..., x) \not\in \Sg_{\bF^k}\big(\operatorname{Cross}(\{y\}, ..., \{y\})\big).
\]
We now inductively define a sequence of proper subalgebras $\PP_i < \bF$ with $y \in \PP_i$ and $x \not\in \PP_i$ as follows. Assuming that we have already constructed $\PP_1, ..., \PP_{i-1}$ such that
\[
(x, ..., x) \not\in \Sg_{\bF^k}\big(\operatorname{Cross}(\PP_1, ..., \PP_{i-1}, \{y\}, \{y\}, ..., \{y\})\big),
\]
we use Zorn's Lemma to find a maximal subalgebra $\PP_i \supseteq \{y\}$ among those satisfying
\[
(x, ..., x) \not\in \Sg_{\bF^k}\big(\operatorname{Cross}(\PP_1, ..., \PP_{i-1}, \PP_i, \{y\}, ..., \{y\})\big).
\]
Note that even the initial choice of $\PP_1$ might need to depend on the particular value of $k$.

We claim that for this particular choice of $\PP_1, ..., \PP_k$, the cross $\operatorname{Cross}(\PP_1, ..., \PP_k)$ forms a subuniverse of $\bF^k$. Suppose for the sake of contradiction that this is not true, i.e. suppose that there are $u_i \in \bF \setminus \PP_i$ such that
\[
(u_1, ..., u_k) \in \Sg_{\bF^k}\big(\operatorname{Cross}(\PP_1, ..., \PP_k)\big).
\]
We will prove, by induction on $i$, that in this case we have
\[
(x, ..., x, u_{i+1}, ..., u_k) \in \Sg_{\bF^k}\big(\operatorname{Cross}(\PP_1, ..., \PP_k)\big)
\]
for each $0 \le i \le k$. Since $u_i \not\in \PP_i$, the choice of $\PP_i$ implies that we have
\begin{align*}
(x, ..., x) &\in \Sg_{\bF^k}\big(\operatorname{Cross}(\PP_1, ..., \PP_{i-1}, \Sg(\PP_i \cup \{u_i\}), \{y\} ..., \{y\})\big)\\
&= \Sg_{\bF^k}\big(\operatorname{Cross}(\PP_1, ..., \PP_{i-1}, \PP_i, \{y\}, ..., \{y\}) \cup \{(x, ..., x, u_i, x, ..., x)\}\big),
\end{align*}
since $\bF$ is generated by $\{x,y\}$ and every element of $\{x,y\}^{i-1}\times \{u_i\} \times \{x,y\}^{k-i}$ other than $(x, ..., x, u_i, x, ..., x)$ is already contained in the cross. Since $\bF$ is the free algebra on two generators, for each $j > i$ there is a homomorphism $\varphi_j : \bF \rightarrow \bF$ such that
\begin{align*}
\varphi_j(x) &= u_j,\\
\varphi_j(y) &= y.
\end{align*}
Since homomorphisms commute with all term operations, this implies that we have
\begin{align*}
(x, ..., x, x, \varphi_{i+1}(x), ..., \varphi_k(x)) \in \Sg_{\bF^k}\big(&\operatorname{Cross}(\PP_1, ..., \PP_{i-1}, \PP_i, \{\varphi_{i+1}(y)\}, ..., \{\varphi_k(y)\})\\
& \cup \{(x, ..., x, u_i, \varphi_{i+1}(x), ..., \varphi_k(x))\}\big),
\end{align*}
that is,
\[
(x, ..., x, x, u_{i+1}, ..., u_k) \in \Sg_{\bF^k}\big(\operatorname{Cross}(\PP_1, ..., \PP_{i-1}, \PP_i, \{y\}, ..., \{y\}) \cup \{(x, ..., x, u_i, u_{i+1}, ..., u_k)\}\big),
\]
which completes the inductive step. Taking $i = k$, we get the desired contradiction.

For part (b), we first use a variation of the argument of (a) to find an infinite sequence of proper subalgebras $\PP_i < \bF$ with $y \in \PP_i$ and $x \not\in \PP_i$ such that $\operatorname{Cross}(\PP_1, ..., \PP_k)$ is a subuniverse of $\bF^k$ for all $k \ge 1$. As in (a), we inductively apply Zorn's Lemma to pick a maximal $\PP_i$ such that
\[
(x, ..., x) \not\in \Sg_{\bF^k}\big(\operatorname{Cross}(\PP_1, ..., \PP_i, \{y\}, ..., \{y\})\big)
\]
for all $k \ge i$. We then inductively show that if there are $u_i \in \bF \setminus \PP_i$ such that
\[
(u_1, ..., u_k) \in \Sg_{\bF^k}\big(\operatorname{Cross}(\PP_1, ..., \PP_k)\big),
\]
then for each $i \le k$ there is some $m \ge k$ such that
\[
(x, ..., x, u_{i+1}, ..., u_k, x, ..., x) \in \Sg_{\bF^m}\big(\operatorname{Cross}(\PP_1, ..., \PP_k, \{y\}, ..., \{y\})\big),
\]
reaching a contradiction once $i = k$.

Once we have found the infinite sequence $\PP_1, \PP_2, ...$, we define $\PP < \bF$ by
\[
\PP = \bigcup_{i \ge 0} \Big(\bigcap_{j \ge i} \PP_j\Big),
\]
that is, $\PP$ consists of all elements of $\bF$ which are contained in all but finitely many of the $\PP_i$s. $\PP$ is a subalgebra of $\bF$ since it is an increasing union of the subalgebras $\bigcap_{j\ge i} \PP_i$. Also, we have $y \in \PP$ and $x \not\in \PP$, so $\PP$ is a proper subalgebra of $\bF$.

To show that $\PP$ is a projective subalgebra of $\bF$, suppose for the sake of contradiction that there was some tuple $u = (u_1, ..., u_k)$ with $u_i \in \bF \setminus \PP$ and
\[
(u_1, ..., u_k) \in \Sg_{\bF^k}\big(\operatorname{Cross}(\PP, ..., \PP)\big).
\]
Then there is a finite collection of tuples $v^1, ..., v^n \in \operatorname{Cross}(\PP, ..., \PP)$ such that $u \in \Sg_{\bF^k}\{v^1, ..., v^n\}$. Pick $N$ large enough that
\[
j_1, ..., j_k \ge N \;\;\; \implies \;\;\; v^m \in \operatorname{Cross}(\PP_{j_1}, ..., \PP_{j_k})
\]
for all $m \le n$. Since each $u_i$ is not contained in $\PP$, there must be an infinite sequence of $j$s such that $u_i \not\in \PP_j$, so we can find $N \le j_1 < j_2 < \cdots < j_k$ such that
\[
u_i \not\in \PP_{j_i}
\]
for each $i \le k$. But then we have
\[
(u_1, ..., u_k) \in \Sg_{\bF^k}\{v^1, ..., v^n\} \subseteq \Sg_{\bF^k}\big(\operatorname{Cross}(\PP_{j_1}, ..., \PP_{j_k})\big),
\]
contradicting the choice of the sequence $\PP_1, \PP_2, ...$.
\end{proof}

Kearnes and Szendrei \cite{cube-terms-crosses} used this to prove the following striking result, which can be viewed as a partial converse to the fact that a semilattice does not have few subpowers.

\begin{cor}[\cite{cube-terms-crosses}]\label{cor-cyclic-cube-term} Suppose that $\bA = (A, c(x_1, ..., x_m))$, where the $m$-ary basic operation $c$ is idempotent and \emph{cyclic}, i.e. $c$ satisfies the identity
\[
c(x_1, x_2, ..., x_m) \approx c(x_2, ..., x_m, x_1).
\]
Then $\bA$ fails to have a cube term iff the two-element algebra
\[
(\{0,1\}, \max(x_1, ..., x_m)),
\]
which is term-equivalent to a semilattice, is contained in the variety generated by $\bA$.
\end{cor}
\begin{proof} By Theorem \ref{thm-infinite-cross}, if $\bA$ has no cube term then the free algebra $\bF = \cF_{\cV(\bA)}(x,y)$ has a proper projective subalgebra $\PP < \bF$. Note that if $c$ is $\PP$-absorbing in any coordinate, then $c$ must be $\PP$-absorbing in \emph{all} coordinates, since $c$ is cyclic.

Pick any element $b \not\in \PP$, and consider the subalgebra $\bB \le \bF$ generated by $\PP \cup \{b\}$. The fact that $c$ is $\PP$-absorbing in all coordinates implies that we have
\[
(x_1, ..., x_m) \in (\PP \cup \{b\})^m \setminus \{(b, ..., b)\} \;\;\; \implies \;\;\; c(x_1, ..., x_m) \in \PP,
\]
and since $c$ is idempotent we have
\[
c(b, ..., b) = b,
\]
so in particular $\bB = \PP \cup \{b\}$. Then if we define an equivalence relation $\theta$ on $\bB$ with blocks $\PP$ and $\{b\}$, we see that $\theta$ is a congruence of $\bB$, and the quotient $\bB/\theta$ is isomorphic to $(\{0,1\}, \max(x_1, ..., x_m))$ by identifying $b/\theta$ with $0$ and $\PP/\theta$ with $1$.
\end{proof}

Our final goal for this section is to go over the proofs from \cite{cube-terms-crosses} and \cite{cube-terms-chipped-cubes} of the fact that if the basic operations of an idempotent algebra have low arity, and if the algebra has no cube term blocker, then it must have a $k$-cube term with $k$ small.

\begin{lem}[\cite{cube-terms-crosses}, \cite{cube-terms-chipped-cubes}] If an $m$-ary idempotent operation $t(x_1, ..., x_m)$ preserves the $k$-dimensional chipped cube $\RR$ given by
\[
(x_1, ..., x_k) \in \RR \;\;\; \iff \;\;\; \bigwedge_{i \le k} x_i \in \bB_i \wedge \bigvee_{j \le k} x_j \in \PP_j,
\]
then there are at most $m - 1$ values of $i \le k$ such that $t^{\bB_i}$ has no $\PP_i$-absorbing coordinate.
\end{lem}
\begin{proof} Suppose for the sake of contradiction that there are at least $m$ values of $i \le k$ such that $t^{\bB_i}$ has no $\PP_i$-absorbing coordinate. We may assume without loss of generality that these include all $i \le m$.

For each $i \le m$, by the assumption that $t^{\bB_i}$ is not $\PP_i$-absorbing in its $i$th coordinate we can pick an $m$-tuple $(b^i_1, ..., b^i_{i-1}, p^i_i, b^i_{i+1}, ..., b^i_m) \in \bB_i^{i-1} \times \PP_i \times \bB_i^{m-i}$ such that
\[
t(b^i_1, ..., b^i_{i-1}, p^i_i, b^i_{i+1}, ..., b^i_m) \not\in \PP_i.
\]
For $m < i \le k$, pick any element $b^i \in \bB_i \setminus \PP_i$ and set $b^i_j = b^i$ for all $j \le m$, so that by idempotence of $t$ we have
\[
t(b^i_1, ..., b^i_m) = t(b^i, ..., b^i) = b^i \not\in \PP_i.
\]
For each $j \le m$, define a $k$-dimensional column vector
\[
v_j = (b^1_j, ..., b^{j-1}_j, p^j_j, b^{j+1}_j, ..., b^k_j)^T \in \RR,
\]
and note that by our choices we have
\[
t(v_1, ..., v_m) \in \prod_{i \le k} \bB_i \setminus \PP_i,
\]
contradicting the assumption that $t$ preserves $\RR$.
\end{proof}

\begin{thm}[Kearnes, Szendrei \cite{cube-terms-crosses}]\label{thm-cube-term-arity-crosses} If an idempotent algebra $\bA = (A, f_1, ..., f_n)$ has a cube term, then it has a $k$-cube term with
\[
k = 1 + \sum_{i \le n} \big(\arity(f_i) - 1\big),
\]
where $\arity(f_i)$ is the arity of $f_i$.
\end{thm}
\begin{proof} Supposing that $\bA$ has no $k$-cube term, Theorem \ref{thm-infinite-cross} shows that somewhere in the variety generated by $\bA$ there is a $k$-dimensional chipped cube $\RR$ given by
\[
(x_1, ..., x_k) \in \RR \;\;\; \iff \;\;\; \bigwedge_{i \le k} x_i \in \bB_i \wedge \bigvee_{j \le k} x_j \in \PP_j.
\]
Our goal is to prove that in this case, at least one of the pairs $(\PP_i, \bB_i)$ is a cube term blocker.

By the lemma, for each basic operation $f_j$ there are at most $\arity(f_j)-1$ values of $i \le k$ such that $f_j^{\bB_i}$ has no $\PP_i$-absorbing coordinate. By the choice of $k$, there must then be some $i^* \le k$ such that each $f_j^{\bB_{i^*}}$ has some $\PP_{i^*}$-absorbing coordinate. Then by Proposition \ref{prop-projective-basic} $\PP_{i^*}$ is a projective subalgebra of $\bB_{i^*}$, so $\bA$ has no cube term.
\end{proof}

\begin{ex} We can now outline an easier way to verify that the three-element algebra $\bA = (\{a,b,c\}, g)$ from Example \ref{ex-few-subpowers} has a $3$-edge term. First, we use Corollary \ref{cor-cyclic-cube-term} to see that since $g$ is symmetric, in order to check whether $\bA$ has few subpowers we just need to check whether $\bA$ has a two-element semilattice as a subquotient (we only need to check subquotients by Lemma \ref{strictly-simple-hs}, since two-element algebras are strictly simple). Once this is done, we note that since $\arity(g) = 3$, Theorem \ref{thm-cube-term-arity-crosses} implies that if $\bA$ has few subpowers then it must have a $3$-cube term, and therefore also a $3$-edge term by Theorem \ref{cube-edge}.
\end{ex}

A better bound when $|\bA|$ is small was given in \cite{cube-terms-chipped-cubes}. The key idea is to show that there is a chipped cube where the number of distinct pairs $(\PP_i,\bB_i)$ is bounded by $\binom{|\bA|}{2}$.

\begin{lem}[\cite{cube-terms-chipped-cubes}] If $\bA$ is a finite idempotent algebra which has a $2(k-1)$-cube term but no $k$-cube term, and if we pick a sequence of $k$ ordered pairs $a_i \ne b_i \in \bA$ which minimize
\[
\sum_{i \le k} |\Sg_{\bA}\{a_i,b_i\}|
\]
among all such sequences of pairs which satisfy
\[
(b_1, ..., b_k) \not\in \Sg_{\bA^k}\Big(\prod_{i \le k} \{a_i,b_i\} \setminus \{(b_1, ..., b_k)\}\Big),
\]
then $\RR = \Sg_{\bA^k}\Big(\prod_{i \le k} \{a_i,b_i\} \setminus \{(b_1, ..., b_k)\}\Big)$ is a $k$-dimensional chipped cube given by
\[
(x_1, ..., x_k) \in \RR \;\;\; \iff \;\;\; \bigwedge_{i \le k} x_i \in \bB_i \wedge \bigvee_{j \le k} x_j \in \PP_j.
\]
with
\begin{align*}
\bB_i &= \Sg_\bA\{a_i,b_i\},\\
\PP_i &= \{p \in \bB_i \mid \Sg_{\bA}\{a_i,p\} \ne \bB_i\},
\end{align*}
and the number of distinct ordered pairs $(\PP_i, \bB_i)$ is at most $\binom{|\bA|}{2}$.
\end{lem}
\begin{proof} Everything other than the bound on the number of pairs $(\PP_i, \bB_i)$ follows from same argument as the one used in Lemma \ref{lem-partial-cube-term} and Theorem \ref{thm-chipped-cubes}, without using the fact that the $a_i$ can also be varied. Since the pair $(\PP_i,\bB_i)$ is determined by the ordered pair $(a_i,b_i)$, we just need to prove that
\[
(a_i,b_i) \ne (b_j,a_j)
\]
for all $i,j$.

Suppose for the sake of contradiction that $(a_i,b_i) = (b_j, a_j)$. Since $a_i \ne b_i$, we can't have $i = j$, so we may assume without loss of generality that $i = 1$ and $j = k$, i.e. that
\[
(a_1, b_1) = (b_k, a_k).
\]
Note that if there is any $c \in \PP_1 \cap \PP_k$, then $\Sg_\bA\{c,b_1\} = \Sg_\bA\{c,a_k\} \ne \bB_k$ and $c \ne b_1$, so replacing $a_1$ with $c$ gives
\[
(b_1, ..., b_k) \not\in \Sg_{\bA^k}\Big(\Big(\{c,b_1\}\times \prod_{2 \le i \le k} \{a_i,b_i\}\Big) \setminus \{(b_1, ..., b_k)\}\Big) \subset \RR,
\]
contradicting the minimality of our initial choice of the sequence of pairs $a_i \ne b_i$. Thus we must have
\[
\PP_1 \cap \PP_k = \emptyset.
\]
Now if we define a relation $\bS \le \bA^{2(k-1)}$ by the primitive positive formula
\[
(x_1, ..., x_{k-1}, y_2, ..., y_k) \in \bS \;\;\; \iff \;\;\; \exists z \text{ s.t. } (x_1, ..., x_{k-1}, z) \in \RR \wedge (z, y_2, ..., y_k) \in \RR,
\]
then
\[
(b_1, ..., b_{k-1}, b_2, ..., b_k) \not\in \bS
\]
since $z$ would need to be in $\PP_1 \cap \PP_k$, but for each $i \le k-1$ and $j \ge 2$ we have
\[
(b_1, ..., a_i, ..., b_{k-1}, b_2, ..., b_k), (b_1, ..., b_{k-1}, b_2, ..., a_j, ..., b_k) \in \bS,
\]
by taking $z = b_k = a_1$ or $z = a_k = b_1$, so $\bA$ can't have a $2(k-1)$-cube term.
\end{proof}

\begin{thm}[Kazda, Zhuk \cite{cube-terms-chipped-cubes}] If a finite idempotent algebra $\bA = (A, f_1, ..., f_n)$ with $|\bA| > 2$ has a cube term, and if
\[
\arity(f_1) \ge \cdots \ge \arity(f_n),
\]
and if we define $m$ by
\[
m = \min\Big(n, \binom{|\bA|}{2}\Big),
\]
then $\bA$ has a $k$-cube term with
\[
k = 1 + \sum_{i \le m} \big(\arity(f_i) - 1\big).
\]
\end{thm}
\begin{proof} We only have to handle the case where $m = \binom{|\bA|}{2} < n$, since the other case follows from Theorem \ref{thm-cube-term-arity-crosses}. Suppose that $K$ is the largest arity of any chipped cube in $\bA$, and note that if $K < 3$ then there is nothing to prove, since $m \ge \binom{3}{2} \ge 3$ and we may assume that none of the $f_i$ are unary. Then since $2(K-1) > K$ we can apply the previous lemma to see that there is a $K$-dimensional chipped cube where the number of distinct pairs $(\PP_i, \bB_i)$ is at most $m$.

If $\bA$ has a cube term then for each distinct pair $(\PP_i, \bB_i)$ there must be some $j$ such that $f_j^{\bB_i}$ has no $\PP_i$-absorbing coordinate. Since there are at most $m$ distinct pairs $(\PP_i, \bB_i)$, at most $m$ functions $f_{j_1}, ..., f_{j_m}$ are needed to witness that none of the pairs $(\PP_i, \bB_i)$ is a cube term blocker. Now we argue as in Theorem \ref{thm-cube-term-arity-crosses} to see that we must have
\[
K \le \sum_{i \le m} \big(\arity(f_{j_i}) - 1\big) \le \sum_{i \le m} \big(\arity(f_i) - 1\big).\qedhere
\]
\end{proof}

Kazda and Zhuk \cite{cube-terms-chipped-cubes} also prove a polynomial bound on the least dimension of a cube term in the non-idempotent case (specifically, they show that we can take $k \le |\bA|^3\arity(\bA)$, where $\arity(\bA)$ is the largest arity of any basic operation of $\bA$), and use it to find an EXPTIME algorithm for determining whether a non-idempotent algebra has few subpowers.

\begin{prob} If a (non-idempotent) algebra $\bA$ is described to us by listing out the tables of its basic operations, is it possible to determine whether or not $\bA$ has few subpowers in time polynomial in the number of bits needed to describe $\bA$?
\end{prob}

\chapter{Absorption and Bounded Width}\label{chapter-bounded-width}

\section{Fourth basic example: the Rock-Paper-Scissors algebra}\label{s-rps}

We're going to start building intuition for the bounded width case with a detailed investigation of a fourth basic algebra on three elements, which is sometimes called the ``rock-paper-scissors'' algebra. This algebra is $\bA_{rps} = (\{a,b,c\}, \cdot)$, where $\cdot$ is the binary, commutative, idempotent operation described by the following table.
\begin{center}
\begin{tabular}{cc}
\begin{tabular}{c|ccc} $\cdot$ & $a$ & $b$ & $c$\\ \hline $a$ & $a$ & $b$ & $a$\\ $b$ & $b$ & $b$ & $c$\\ $c$ & $a$ & $c$ & $c$\end{tabular} & 
\begin{tikzpicture}[scale=1.5,baseline=0.5cm]
  \node (a) at (-0.5,0) {$a$};
  \node (b) at (0,0.8) {$b$};
  \node (c) at (0.5,0) {$c$};
  \draw [->] (a) edge (b) (b) edge (c) (c) edge (a);
\end{tikzpicture}
\end{tabular}
\end{center}
The algebra $\bA_{rps}$ is not a semilattice, but every two-element subset of $\bA_{rps}$ is a semilattice. Thus, the binary operation $\cdot$ satisfies the following identities:
\[
xx \approx x, \;\;\; xy \approx yx, \;\;\; x(xy) \approx xy.
\]
Any binary operation satisfying the above identities is known as a 2-\emph{semilattice} operation, and the algebra $\bA_{rps}$ is the smallest 2-semilattice which is not a semilattice.

As we will see, the corresponding relational clone is generated by the binary relation $\{(a,b), (b,c), (c,a)\}$ (which corresponds to the fact that the algebra has a cyclic automorphism) and the ternary relation $R_{a,b}$ given by the formula
\[
R_{a,b}(x,y,z) \coloneqq (x \in \{a,b\}) \wedge (x = a \implies y = z).
\]
The ternary relation $R_{a,b}$ has a special role, which is closely connected to the fact that $\{a,b\}$ is a semilattice subalgebra of $\bA_{rps}$.

\begin{prop} If $\RR, \bS \subseteq \bA^n$ are any $n$-ary relations with $\bS \subseteq \RR$, then the $n+1$-ary relation
\[
((x,y) \in \RR\times \{a,b\}) \wedge (y = a \implies x \in \bS)
\]
can be defined by a primitive positive formula in $\RR,\bS$, and $R_{a,b}$.
\end{prop}
\begin{proof} Just use the following primitive positive formula:
\[
\exists z \in \bA^n\ \ x \in \RR \wedge z \in \bS \wedge \bigwedge_{i \le n} R_{a,b}(y,x_i,z_i).\qedhere
\]
\end{proof}

\begin{prop}\label{rps-decompose} If $\RR \le_{sd} \bA_{rps}^n$ is a subdirect $n$-ary relation, then $\RR$ is the intersection of its two-variable projections, each of which is either a full relation or the graph of an automorphism of $\bA_{rps}$ which is either the identity or is cyclic. In particular, there is some subset of the coordinates $I \subseteq \{1, ..., n\}$ such that the projection $\pi_I$ is an isomorphism from $\RR$ to $\bA_{rps}^I$.
\end{prop}
\begin{proof} We prove this by induction on $n$. The base case, $n = 2$, is easily verified: since $\bA_{rps}$ is simple, every subdirect binary relation on $\bA_{rps}$ is either the graph of an automorphism or is linked, and we can check that every connected subgraph of the complete bipartite graph $K_{3,3}$ either contains a bipartite matching or is a tree with two leaves on both parts (e.g. using Hall's Marriage Lemma). Therefore up to automorphisms of $\bA_{rps}$ we just need to consider relations which contain $\{(a,a),(a,b),(b,b),(c,c)\}$, $\{(a,a),(b,c),(c,b)\}$, or $\{(a,b),(a,c),(b,a),(c,a)\}$, and all three of these generate $\bA_{rps}^2$.

Now consider the case $n > 2$. By the induction hypothesis, we may assume without loss of generality that $\pi_{[n]\setminus \{i\}}(\RR) = \bA_{rps}^{n-1}$ for every $i \le n$. Suppose for contradiction that $\RR \ne \bA_{rps}^n$.

Since the automorphism group of $\bA_{rps}$ is transitive, we may assume without loss of generality that $(a,...,a) \not\in \RR$. Since $\bA_{rps}$ is idempotent, the set $\RR'$ of triples $(x,y,z)$ such that $(x,y,z,a,...,a) \in \RR$ is a subalgebra of $\bA_{rps}^3$, and since $\pi_{[n]\setminus \{i\}}(\RR) = \bA_{rps}^{n-1}$ for $i = 1,2,3$ we see that every projection of $\RR'$ onto any pair of coordinates is full. So we can reduce to the case $n = 3$.

If any two of $(a,a,c), (a,c,a), (c,a,a)$ are in $\RR$, then we can combine them to obtain $(a,a,a)$. So we may suppose that $(a,a,b) \in \RR$. If we consider the binary relation consisting of pairs $(y,z)$ with $(a,y,z) \in \RR$, then by the $n = 2$ case, we must have $(a,c,a) \in \RR$. Similar reasoning with the roles of the first and second coordinates reversed then shows that we must have $(c,a,a) \in \RR$, a contradiction.
\end{proof}

\begin{prop}\label{rps-absorb} If $\RR \le_{sd} \bA_{rps}^n\times \{a,b\}^k$ has full projection onto $\bA_{rps}^n$, then we have $\bA_{rps}^n\times \{(b,...,b)\} \subseteq \RR$.
\end{prop}
\begin{proof} For any $x \in \bA_{rps}^n$, let $x^-$ be the tuple obtained from $x$ by applying the cyclic permutation $(a\ c\ b)$ componentwise. Then it's easy to check that for any $x,y \in \bA_{rps}^n$, we have
\[
(xy^-)y = y.
\]

By multiplying all of the elements of $\RR$ together in any order (with parentheses placed arbitrarily), we see that there is some $x \in \bA_{rps}^n$ such that $(x_1, ..., x_n, b, ..., b) \in \RR$. For any $y \in \bA_{rps}^n$, there are tuples $c,d \in \{a,b\}^n$ such that $(y,c), (y^-,d) \in \RR$ by the assumption $\pi_{[n]}(\RR) = \bA_{rps}^n$. Thus
\[
(y,b) = ((x,b)\cdot(y^{-},d))\cdot (y,c) \in \RR.\qedhere
\]
\end{proof}

The previous two propositions are enough to describe an algorithm which solves $\CSP(\bA_{rps})$. The algorithm first establishes arc-consistency, reducing some of the domains of the variables until every constraint relation becomes subdirect. Then for each variable with a two element domain, the last proposition shows that we may as well take that variable equal to the top/absorbing element of that domain. After this restriction, if we consider the remaining variables, each relation decomposes into binary relations, each of which is either an equality relation or the graph of a cyclic automorphism. This final problem can be solved by checking that no cycle of binary relations implies that any variable is related to itself by a nontrivial cyclic automorphism.

\begin{defn} An instance of a CSP is \emph{cycle-consistent} if for every sequence of variables $v_1, ..., v_n$ and relations $R_1, ..., R_n$ and pairs of coordinates $(i_k,j_k)$ such that $v_k, v_{k+1}$ are related by $\pi_{(i_k,j_k)}(R_k)$ for each $k$ (indices taken modulo $n$), the composition
\[
\pi_{(i_1,j_1)}(R_1) \circ \cdots \circ \pi_{(i_n,j_n)}(R_n)
\]
contains the equality relation on the domain of the variable $v_1$.
\end{defn}

\begin{cor} Any cycle-consistent instance of $\CSP(\bA_{rps})$ has a solution.
\end{cor}

If we want to understand the complete structure of a general relation $\RR \le \bA_{rps}^n$, things become more complicated. Typical relations we need to consider have the form
\begin{align*}
x_1 \in &\{a,b\} \wedge (x_1 = a \implies x_2 \in \{a,b\}) \wedge (x_1 = x_2 = a \implies x_3 \in \{a,b\})\\
&\wedge \cdots \wedge (x_1 = \cdots = x_k = a \implies y=z)
\end{align*}
or
\begin{align*}
x_1 \in &\{a,b\} \wedge (x_1 = a \implies x_2 \in \{a,b\}) \wedge (x_1 = x_2 = a \implies x_3 \in \{a,b\})\\
&\wedge \cdots \wedge (x_1 = \cdots = x_{k-1} = a \implies x_k \in \{a,b\}).
\end{align*}
We can slightly modify these relations by applying cyclic automorphisms of $\bA_{rps}$ to some of the variables, or by renaming the variables. We call any relation obtained by making such a modification to the two types of relation above a \emph{basic relation} on $\bA_{rps}$.

\begin{thm} Suppose $\RR \le \bA_{rps}^n$. Then $x \in \RR$ iff $x$ satisfies every basic relation on $\bA_{rps}$ which contains $\RR$. In particular, $\RR$ is contained in the relational clone generated by $\{(a,b), (b,c), (c,a)\}$ and $R_{a,b}$.
\end{thm}
\begin{proof} Suppose $x$ satisfies every basic relation which contains $\RR$. Let $I = \{i_1, ..., i_k\} \subseteq [n]$ be maximal such that, after applying cyclic automorphisms to coordinates in $I$, we have $x_{i_j} = a$ for all $j \le k$, and such that the basic relation
\begin{align*}
y_{i_1} \in &\{a,b\} \wedge (y_{i_1} = a \implies y_{i_2} \in \{a,b\}) \wedge (y_{i_1} = y_{i_2} = a \implies y_{i_3} \in \{a,b\})\\
&\wedge \cdots \wedge (y_{i_1} = \cdots = y_{i_{k-1}} = a \implies y_{i_k} \in \{a,b\})
\end{align*}
contains $\RR$. Assume without loss of generality that the coordinates are ordered such that $I = \{n-k+1, ..., n\}$ and such that the $n-k$-ary relation $\RR'$ defined by
\[
(y_1, ..., y_{n-k}) \in \RR' \iff (y_1, ..., y_{n-k}, a, ..., a) \in \RR
\]
has $\RR' \le_{sd} \bA_{rps}^m\times \{a,b\}^{n-m-k}$ for some $m$ (possibly after further applications of cyclic automorphisms). Then by the maximality of $I$, we have $x = (x_1, ..., x_m, b, ..., b, a, ..., a)$. By Propositions \ref{rps-decompose}, \ref{rps-absorb}, and our assumption that $x$ satisfies all basic relations containing $\RR$, we have $(x_1, ..., x_m) \in \pi_{[m]}(\RR')$ and $\pi_{[m]}(\RR')\times \{(b,...,b)\} \subseteq \RR'$, so $(x_1, ..., x_m, b, ..., b) \in \RR'$, so $x \in \RR$.
\end{proof}

\begin{rem} The intricate yet understandable structure of the basic relations considered above is at the heart of the uncountable region found by Zhuk \cite{zhuk-selfdual} in the lattice of clones on a three-element domain. Each of the clones in Zhuk's uncountable region properly contains the clone of the rock-paper-scissors algebra, so the generating relations for the corresponding relational clones can be written in terms of the basic relations considered above.
\end{rem}

\begin{prop} Suppose that $f : \bA^n \rightarrow \bA$ is any idempotent operation which depends on all of its inputs and preserves the relation $R_{a,b}$. Then the restriction of $f$ to $\{a,b\}$ must be the $n$-ary semilattice operation on $\{a,b\}$, that is, for any $(x_1, ..., x_n) \in \{a,b\}^n \setminus \{(a,...,a)\}$, we have $f(x_1, ..., x_n) = b$.
\end{prop}
\begin{proof} Suppose for contradiction that there is some $(x_1, ..., x_n) \in \{a,b\}^n \setminus \{(a,...,a)\}$ with $f(x_1, ..., x_n) \ne b$. Since $\{a,b\} = \pi_1(R_{a,b})$ is preserved by $f$, we must then have $f(x_1, ..., x_n) = a$. We will show that for all $i$ with $x_i = b$, $f$ does not depend on its $i$th input.

Let $y, z \in \bA^n$ be any pair of tuples with $y_i = z_i$ whenever $x_i = a$. Then each $(x_i,y_i,z_i) \in R_{a,b}$, so
\[
\begin{bmatrix} a\\ f(y)\\ f(z)\end{bmatrix} = f\left(\begin{bmatrix} x_1 & x_2 & \cdots & x_n\\ y_1 & y_2 & \cdots & y_n\\ z_1 & z_2 & \cdots & z_n\end{bmatrix}\right) \in R_{a,b},
\]
so $f(y) = f(z)$.
\end{proof}

\begin{thm} An $n$-ary operation $f$ is contained in $\Clo_n(\bA_{rps})$ iff it preserves the relations $\{(a,b),(b,c),(c,a)\}$ and $R_{a,b}$. If $f$ depends on all its inputs, this occurs iff $f$ preserves the cyclic automorphism of $\bA_{rps}$ and $f|_{\{a,b\}}$ is the $n$-ary semilattice operation on $\{a,b\}$.
\end{thm}
\begin{proof} We just need to check this in the case where $f$ depends on all of its inputs. Let $\cF = \cF_{\cV(\bA_{rps})}(x_1, ..., x_n) \le \bA_{rps}^{\bA_{rps}^n}$ be the subalgebra generated by $\pi_1, ..., \pi_n : \bA_{rps}^n \rightarrow \bA_{rps}$. The projection $\pi_x(\cF)$ of $\cF$ onto the coordinate of $\bA_{rps}^{\bA_{rps}^n}$ corresponding to $x \in \bA_{rps}^n$ is the subalgebra of $\bA_{rps}$ generated by $\{\pi_1(x), ..., \pi_n(x)\} = \{x_1, ..., x_n\}$.

If $x$ is a diagonal tuple, say $x = (a,...,a)$, then $\pi_x(\cF) = \{a\}$, corresponding to the fact that any $f \in \cF$ must be idempotent, with $f(a,...,a) = a$. If exactly two elements of $\bA_{rps}$ occur in $x$, say $x \in \{a,b\}^n$, then $\pi_x(\cF) = \{a,b\}$, and if $f$ depends on all its inputs and preserves $R_{a,b}$, this implies that we must have $f(x) = b$, i.e. $\pi_x(f) = b$. Thus, if $I \subseteq \bA_{rps}^n$ is the set of $x$ such that all three of $a,b,c$ show up in the coordinates of $x$, we see that $\pi_I(\cF) \le_{sd} \bA_{rps}^I$, and by Proposition \ref{rps-absorb} we have $f \in \cF \iff \pi_I(f) \in \pi_I(\cF)$.

By Proposition \ref{rps-decompose}, $\pi_I(\cF)$ is the intersection of its two-variable projections, each of which is either full or the graph of a cyclic automorphism of $\bA_{rps}$. A two variable projection $\pi_{x,y}(\cF)$ will only be the graph of a cyclic automorphism $\sigma \in \Aut(\bA_{rps})$ if $(\pi_i(x),\pi_i(y))$ is in the graph of $\sigma$ for all $i$, that is, if $y_i = \sigma(x_i)$ for all $i$. Thus, $\pi_I(f) \in \pi_I(\cF)$ iff whenever $y = \sigma(x)$, we have $f(y) = \sigma(f(x))$.
\end{proof}

Note that one of the key steps behind the analysis of the rock-paper-scissors algebra was Proposition \ref{rps-decompose} which classified the subdirect powers of the algebra, and that the method of proof depended only on checking properties of subdirect binary and ternary relations on $\bA_{rps}$. The general pattern behind this is best understood in terms of a property of the polynomial clone known as \emph{polynomial completeness}.

\begin{defn} An algebra is \emph{polynomially complete} if its polynomial clone is the clone of all operations, that is, if every operation on the underlying set can be expressed using the basic operations of the algebra together with the constant operations.
\end{defn}

\begin{thm} A finite idempotent algebra $\bA$ is polynomially complete if every binary relation on $\bA$ which contains the diagonal is either the equality relation or the full relation, and every ternary relation $\RR \le_{sd} \bA^3$ such that every two variable projection of $\RR$ is full is equal to the full relation $\bA^3$.
\end{thm}
\begin{proof} We will show by induction on $n$ that every $n$-ary relation $\RR \le \bA^n$ which contains the subalgebra of diagonal tuples $(x,...,x)$, $x \in \bA$ is given by a conjunction of equalities between pairs of coordinates. The base case $n = 2$ follows from our assumption on $\bA$. By the inductive hypothesis, we may assume without loss of generality that $\pi_{[n]\setminus\{i\}} \RR = \bA^{n-1}$ for each $i$.

If $n=3$, then our assumption on $\bA$ implies that $\RR = \bA^3$. Otherwise, suppose for contradiction that $(x_1, ..., x_n) \not\in \bA^n$, and consider the ternary relation $\RR'$ consisting of triples $(u,v,w)$ such that $(u,v,w,x_4, ..., x_n) \in \RR$. Since $\bA$ is idempotent, $\RR'$ is a subalgebra of $\bA^3$, and every two-variable projection of $\RR'$ is full, so by the $n=3$ case we must have $(x_1,x_2,x_3) \in \RR'$, a contradiction.

Note that we have shown that the relational clone corresponding to the polynomial clone of $\bA$ is generated by the equality relation. The general $\Inv-\Pol$ Galois duality now shows that $\bA$ is polynomially complete. To see this concretely, consider the subalgebra of $\bA^{\bA^n}$ generated by the functions $\pi_i$ and the constant (diagonal) tuples. Then this subalgebra is described by a conjunction of equalities between pairs of coordinates. But no two-variable projection of this subalgebra can be an equality relation: if $x \ne y \in \bA^n$, then there is always some $i$ such that $\pi_i(x) \ne \pi_i(y)$. Thus this subalgebra of $\bA^{\bA^n}$ must be the full set of operations $\bA^n \rightarrow \bA$.
\end{proof}

\begin{cor} The rock-paper-scissors algebra is polynomially complete.
\end{cor}

As far as relations go, the main impact of polynomial completeness is that it strongly constrains subdirect relations where each factor is polynomially complete. As we have seen, if some factors are not polynomially complete, then the structure of an arbitrary relation can be quite intricate. In the case of the rock-paper-scissors algebra, we are able to side-step this intricacy by restricting each factor which is a proper subalgebra of $\bA_{rps}$ to its top/absorbing element. This is a general strategy that can be used in the study of bounded width algebras, as well as finite Taylor algebras.

We conclude this section with a few classical results about polynomial completeness.

\begin{defn} The ternary \emph{discriminator} function is the function $t$ defined by
\[
t(x,y,z) = \begin{cases} z & x = y,\\ x & x \ne y.\end{cases}
\]
\end{defn}

\begin{prop} A finite algebra is polynomially complete iff it has the ternary discriminator as a polynomial operation.
\end{prop}
\begin{proof} One direction is obvious. For the other direction, it's enough to show that the idempotent algebra $\bA = (A,t)$ whose only basic operation is the ternary discriminator $t$ is polynomially complete. We may assume that the underlying set $A$ contains at least two distinct elements $a,b$. Suppose first that $\RR \le_{sd} \bA^2$ is a relation properly containing the diagonal of $\bA^2$, and assume without loss of generality that $(a,b) \in \RR$ with $a \ne b$. Then for any $c \in \bA$, we have
\[
\begin{bmatrix} a\\ c\end{bmatrix} = t\left(\begin{bmatrix} a & b & c\\ b & b & c \end{bmatrix}\right) \in \RR,
\]
and similarly $(d,b) \in \RR$ for any $d \in \bA$. Then for any $c,d \in \bA$ we have
\[
\begin{bmatrix} d\\ c\end{bmatrix} = t\left(\begin{bmatrix} a & a & d\\ c & b & b \end{bmatrix}\right) \in \RR,
\]
so $\RR = \bA^2$.

To finish, we just need to show that any ternary relation $\RR \le_{sd} \bA^3$ such that every two variable projection is full must be the full relation $\bA^3$. Since $\bA$ has full automorphism group, if $\RR \ne \bA^3$ then we may assume without loss of generality that $(a,a,a) \not\in \RR$, while all three of $(a,a,b),(a,b,a),(b,a,a)$ are in $\RR$. Then we have
\[
\begin{bmatrix} b\\ a\\ b\end{bmatrix} = t\left(\begin{bmatrix} a & a & b\\ a & b & a\\ b & a & a \end{bmatrix}\right) \in \RR,
\]
so
\[
\begin{bmatrix} a\\ a\\ a\end{bmatrix} = t\left(\begin{bmatrix} a & b & b\\ a & a & a\\ b & b & a \end{bmatrix}\right) \in \RR,
\]
contradicting the assumption $(a,a,a) \not\in \RR$.
\end{proof}

\begin{ex} We can give an alternative proof of the fact that the rock-paper-scissors algebra is polynomially complete by expressing the ternary discriminator as a polynomial. First, we can define the unary polynomial $x^+$ corresponding to the cyclic permutation $(a\ b\ c)$ by
\[
x^+ = ((xa)c)(xb),
\]
and we can define the inverse of this by $x^- = (x^+)^+$.
Note that we now have
\[
xy^+ = \begin{cases} x^+ & x = y,\\ x & x\ne y.\end{cases}
\]
Thus if we set $u(x,y,z) = (z(xy^+)^-)x$, then we have
\[
u(x,y,z) = (z(xy^+)^-)x = \begin{cases} xz & x = y,\\ x & x \ne y,\end{cases}
\]
so we may take
\[
t(x,y,z) = ((u(x,y,z)u(x,y,z^+)^-)u(x,y,z^-)^+)^-.
\]
To see that this works, note that if $x=y$, then two of $xz, (xz^+)^-, (xz^-)^+$ are equal to $z$ while the third is equal to $z^+$, so since $\{z,z^+\}$ is a semilattice we see that in this case $t(x,y,z)$ is given by
\[
(((xz)(xz^+)^-)(xz^-)^+)^- = (zzz^+)^- = (z^+)^- = z,
\]
while if $x \ne y$ then $u(x,y,?) = x$, so $t(x,y,z)$ is given by
\[
((xx^-)x^+)^- = (xx^+)^- = (x^+)^- = x.
\]
\end{ex}

The ternary discriminator $t$ satisfies the system of identities
\begin{align*}
t(x,y,y) &\approx x,\\
t(x,y,x) &\approx x,\\
t(y,y,x) &\approx x.
\end{align*}
Any ternary term satisfying this system of identities is known as a \emph{Pixley term}. Note that any Pixley term is automatically a Mal'cev term, and that the term $d$ defined from $t$ by
\[
d(x,y,z) = t(x,t(x,y,z),z)
\]
is automatically a majority term. In the case where $t$ is the ternary discriminator, $d$ becomes the dual discriminator of Example \ref{ex-dual-discriminator}.

\begin{thm}[Pixley \cite{Pixley-term}]\label{pixley-poly} An algebra $\bA$ generates a variety which is both congruence permutable and congruence distributive iff it has a Pixley term. If $\bA$ is also simple, then it is polynomially complete.
\end{thm}
\begin{proof} If $\bA$ has a Pixley term, then it has both a Mal'cev term and a majority term, so it generates a congruence permutable and congruence distributive variety. Conversely, suppose that $\bA$ generates a congruence permutable and congruence distributive variety. Let $\cF = \cF_{\cV(\bA)}(x,y,z)$ be the free algebra on three generators in this variety, and for $a,b \in \{x,y,z\}$ let $\theta_{ab}$ be the smallest congruence with $a \equiv_{\theta_{ab}} b$. Then $(x,z) \in \theta_{xz} \wedge (\theta_{xy} \circ \theta_{yz})$, so by congruence distributivity and permutability, we have
\[
(x,z) \in \theta_{xz} \wedge (\theta_{xy} \vee \theta_{yz}) = (\theta_{xz} \wedge \theta_{xy}) \vee (\theta_{xz} \wedge \theta_{yz}) = (\theta_{xz} \wedge \theta_{yz}) \circ (\theta_{xz} \wedge \theta_{xy}).
\]
Thus there is some $t \in \cF$ such that
\[
x\ (\theta_{xz} \wedge \theta_{yz})\ t\ (\theta_{xz} \wedge \theta_{xy})\ z.
\]
Thus $t$ is a ternary term which satisfies the Pixley identities.

Now suppose that $\bA$ is simple. Since $\bA$ is Mal'cev, every binary relation on $\bA$ is the graph of an isomorphism modulo the linking congruence, and the linking congruence is necessarily either $0_\bA$ or $1_\bA$. Thus every binary relation on $\bA$ which contains the diagonal is either full or equal to the diagonal. Since $\bA$ has a majority term, every ternary relation on $\bA$ whose two variable projections are all full must itself be a full relation. Thus $\bA$ is polynomially complete.
\end{proof}

Varieties which are both congruence distributive and congruence permutable are known as \emph{arithmetical} varieties. The name arithmetical comes from the theory of arithmetical rings, which are rings where the ``Chinese remainder condition'' holds: for any ideals $I_1, ..., I_n$ and elements $a_1, ..., a_n$ with $a_i \equiv a_j \pmod{I_i + I_j}$ for all $i,j$, there exists some $x$ with $x \equiv a_i \pmod{I_i}$ for all $i$.

\section{Partial semilattice operations and the digraph of semilattice subalgebras}\label{s-partial-semi}

In this section we will go over a binary analogue of a standard result about iterating unary functions to make (compositionally) idempotent functions, that is, functions satisfying $e \circ e = e$. First we review the case of unary iteration.

\begin{defn} If $f : A \rightarrow A$ is a unary function, we define $f^{\circ n}$ to be $f \circ \cdots \circ f$, with $n$ copies of $f$. If $(A,f)$ is either finite or profinite, we define $f^{\infty}$ by
\[
f^{\infty}(x) \coloneqq \lim_{n \rightarrow \infty} f^{\circ n!}(x).
\]
Alternatively, we can define $f^{\infty}$ as the limit of $f^{\circ n}$ over the net of positive integers $n$, ordered by divisibility. Similarly, we define $f^{\infty - 1}$ by
\[
f^{\infty - 1}(x) \coloneqq \lim_{n \rightarrow \infty} f^{\circ (n! - 1)}(x).
\]
\end{defn}

\begin{prop} If $(A,f)$ is profinite, then the limit defining $f^{\infty}$ exists, and $f^{\infty}$ satisfies the identity
\[
f^{\infty}(f^{\infty}(x)) \approx f^{\infty}(x).
\]
Furthermore, if $A$ is finite, then
\[
f^{\infty} = f^{\circ \lcm\{1, ..., |A|\}},
\]
and the graph of $f^{\infty}$ can be computed from the graph of $f$ in time linear in $|A|$.
\end{prop}
\begin{proof} It's enough to prove this in the case where $A$ is finite. Let $m, m'$ be any positive multiples of $\lcm\{1, ..., |A|\}$, we will show that $f^{\circ m} = f^{\circ m'}$: this will show that the limit is equal to $f^{\circ m}$, and taking $m' = 2m$ will show that $f^\infty \circ f^\infty = f^\infty$. To see that $f^{\circ m} = f^{\circ m'}$, note that for any $x$, the sequence $x, f(x), f(f(x)), ..., f^{\circ k}(x), ...$ must be eventually periodic with period $p$ at most $|A|$, and the periodic behavior must begin within the first $|A|$ steps, so for any $k \ge |A|$ we have $f^{\circ k}(x) = f^{\circ (k+p)}(x)$. Since $|m-m'|$ is a multiple of $p$ and $m,m' \ge |A|$, this implies that $f^{\circ m} = f^{\circ m'}$.

In order to compute the graph of $f^{\infty}$ efficiently, we will also compute the function $f^{\infty - 1}$ simultaneously. First, make a list of elements of $A$, and mark all of them as ``unprocessed''. In each round, we pick the next unprocessed element $x$ from the list, and compute the sequence of iterates $x, f(x), f(f(x)), ...$, marking each one as ``processed'' as we compute it, until the first time we compute $f^{\circ k}(x)$ and find that it has already been marked as ``processed''. There are two cases: either $f^{\circ k}(x)$ is equal to $f^{\circ i}(x)$ for some $i < k$, or $f^{\circ k}(x)$ was processed in some previous round. We can distinguish between the two cases by checking whether the value of $f^\infty(f^{\circ k}(x))$ has already been computed.

In the case where $f^{\circ k}(x) = f^{\circ i}(x)$ for some $i < k$, we first set $f^\infty(f^{\circ j}(x)) \coloneqq f^{\circ j}(x)$ and $f^{\infty-1}(f^{\circ j}(x)) \coloneqq f^{\circ (j-1)}(x)$ for $i < j \le k$. For $j < i$, we iterate downwards, setting
\[
f^\infty(f^{\circ j}(x)) \coloneqq f^{\infty - 1}(f^{\circ (j+1)}(x))
\]
and
\[
f^{\infty-1}(f^{\circ j}(x)) \coloneqq f^{\infty - 1}(f^\infty(f^{\circ j}(x))).
\]
In the case where $f^{\circ k}(x)$ was processed in a previous round, we iterate downwards using the above rules to handle all $j < k$.

Since the number of steps needed for each round is linear in the number of elements which are marked as processed in that round, and since each element of $A$ is marked as processed at most once, the entire procedure for computing $f^\infty$ and $f^{\infty-1}$ runs in time linear in $|A|$.
\end{proof}

In the context of CSPs, the reduction to the case of core structures was based on the observation than any non-surjective unary polymorphism $f : \fA \rightarrow \fA$ allows us to replace the underlying set $A$ by the smaller set $f(A)$ to obtain a homomorphically equivalent CSP on a smaller domain. In this case, the map $f^\infty : \fA \rightarrow \fA$ will also be non-surjective, and in fact we have the guarantee that
\[
f^\infty(A) \subseteq f^{\circ n}(A)
\]
for all $n \ge 0$. So whenever we shrink the domain of a non-core CSP using a unary polymorphism, we may as well assume that the unary polymorphism in question is (compositionally) idempotent.

On the algebraic side, if $e \circ e = e$ and $e \in \Clo_1(\bA)$, we can define a reduct $\bA_e$ of $\bA$ as follows. For every $n$-ary operation $f \in \Clo_n(\bA)$, we define the corresponding operation $f_e : A^n \rightarrow A$ by
\[
f_e(x_1, ..., x_n) = e(f(e(x_1), ..., e(x_n))).
\]
Then we define $\bA_e$ to be the algebraic structure $(A, \{f_e \mid f \in \Clo(\bA)\})$ having a basic operation $f_e$ for each term $f$ of $\bA$.

Each operation $f_e$ only depends on the restriction of $f$ to $e(A)$, and takes values in $e(A)$. Also, if $f$ preserves $e(A)$, then $f_e$ and $f$ agree when they are restricted to $e(A)$. The reduct $\bA_e$ has $e(A)$ as a subalgebra, and is completely determined by its restriction to the subalgebra $e(A)$ together with the description of the map $e : A \rightarrow e(A)$. So the reduct $\bA_e$ and its subalgebra $e(A)$ are essentially interchangeable, and the subalgebra $e(A)$ of $\bA_e$ has as its basic operations the terms of $\bA$ which preserve $e(A)$.

As a special case of the general result relating reflections to height $1$ identities, we have the following basic result.

\begin{prop} If a system of height $1$ identities is satisfied by terms $f^1, ..., f^k$ of $\bA$, then the same system of height $1$ identities is satisfied by the corresponding operations $f^1_e, ..., f^k_e$ of $\bA_e$ (defined as above).
\end{prop}

Note that identities which involve nesting functions may not survive the process of passing from $\bA$ to the reduct $\bA_e$.

Now we return to the world of idempotent operations, and describe a surprisingly powerful binary analogue of unary iteration. Rather than (compositionally) idempotent operations, we will produce a type of binary operation which I call a \emph{partial semilattice} operation.

\begin{defn}\label{defn-partial-semilattice} We say that an idempotent binary operation $s$ is a \emph{partial semilattice} if it satisfies the identity
\[
s(x,s(x,y)) \approx s(s(x,y),x) \approx s(x,y).
\]
Equivalently, $s$ is a partial semilattice if for all $x,y$, the set $\{x,s(x,y)\}$ is closed under $s$, and acts like a semilattice subalgebra with absorbing element $s(x,y)$ under $s$.
\end{defn}

Note that unlike semilattices and $2$-semilattices, partial semilattices are \emph{not necessarily} Taylor operations. The binary projection $\pi_1$ is an extreme example of a partial semilattice operation which is not Taylor. This is a necessary feature of the definition, since we will show that \emph{any} idempotent binary operation can be used to produce a partial semilattice operation (in a nontrivial way).

In order to produce partial semilattice operations, we will start by treating our binary operation as a unary function of the second variable, with the first variable treated as a (constant) parameter.

\begin{defn} If $t : \bA^2 \rightarrow \bA$ is a binary function and $\bA$ is finite (or profinite), then we define $t^\infty$ to be the pointwise limit
\[
t^\infty(x,y) \coloneqq \lim_{n \rightarrow \infty} t^{n!}(x,y),
\]
where $t^1 \coloneqq t$ and $t^{n+1}(x,y) \coloneqq t(x,t^n(x,y))$.
\end{defn}

\begin{prop} For any binary term $t$, we have
\[
t^\infty(x,t^\infty(x,y)) \approx t^\infty(x,y).
\]
If $t$ is idempotent, then so is $t^\infty$.
\end{prop}

The function $t^\infty$ now satisfies one of the two defining identities for a partial semilattice. Note that $t^\infty$ can be computed from $t$ in time linear in $|A|^2$. To find a term $u$ which satisfies the second identity $u(u(x,y),x) \approx u(x,y)$, we plug $t^\infty$ into itself in a surprisingly counterintuitive way.

\begin{prop} If $f$ is an idempotent binary term which satisfies the identity
\[
f(x,f(x,y)) \approx f(x,y),
\]
and if we define a term $u$ by
\[
u(x,y) \coloneqq f(x,f(y,x)),
\]
then $u$ satisfies the identity
\[
u(u(x,y),x) \approx u(x,y).
\]
\end{prop}
\begin{proof}
We have
\[
f(x,u(x,y)) \approx f(x,f(x,f(y,x))) \approx f(x,f(y,x)) \approx u(x,y),
\]
so
\[
u(u(x,y),x) \approx f(u(x,y),f(x,u(x,y))) \approx f(u(x,y),u(x,y)) \approx u(x,y).\qedhere
\]
\end{proof}

Finally, to get a term which satisfies \emph{both} defining identities of a partial semilattice, we iterate the function $u$ on its second variable.

\begin{prop} If $u$ is an idempotent binary term which satisfies the identity
\[
u(u(x,y),x) \approx u(x,y),
\]
then $s \coloneqq u^\infty$ satisfies the identity
\[
s(x,s(x,y)) \approx s(s(x,y),x) \approx s(x,y).
\]
\end{prop}
\begin{proof} Define $u^n$ as in the definition of $u^\infty$. Then for any $m$ we have
\[
u^m(u(x,y),x) \approx u(x,y),
\]
and on replacing $y$ by $u^{n-1}(x,y)$, we get
\[
u^m(u^n(x,y),x) \approx u^n(x,y)
\]
for any $m,n$.
\end{proof}

The full process, going from $t$ to $f = t^\infty$ to $u(x,y) = f(x,f(y,x))$ to $s = u^\infty$, is functorial, and the final function $s : A^2 \rightarrow A$ can be computed from $t$ in time linear in $|A|^2$. Since $s$ was defined from $t$ in a nontrivial way, we get the following result.

\begin{prop}\label{prop-semilattice-iteration} If $t$ is a binary idempotent term and $a,b$ are such that $t(a,b) = t(b,a) = b$, then the partial semilattice term $s \in \Clo(t)$ defined by the above process also satisfies $s(a,b) = s(b,a) = b$.

More generally, if $B,C$ are subsets of $\bA$ such that for any $x \in B\cup C$ and any $y \in C$ we have $t(x,y), t(y,x) \in C$, then the same holds for $s$.
\end{prop}

\begin{cor} If $(b,b) \in \Sg_{\bA^2}\{(a,b), (b,a)\}$, then there is a partial semilattice term $s \in \Clo(\bA)$ such that $s(a,b) = s(b,a) = b$.
\end{cor}

Once we have a partial semilattice term $s$ with $s(a,b) = s(b,a) = b$, we can use it to preprocess the inputs to other $n$-ary functions to force them to preserve the subset $\{a,b\}$ and act like the $n$-ary semilattice operation on this subset. To do this, we first need to find terms $s_n \in \Clo(s)$ which act like the $n$-ary semilattice operation.

\begin{prop}\label{higher-semilattice} If $s$ is a partial semilattice operation, then for all $n$ there are terms $s_n \in \Clo(s)$ of arity $n$ such that if $\{x,x_2, ..., x_n\} = \{x,y\}$, then
\[
s_n(x,x_2, ..., x_n) \approx s(x,y).
\]
\end{prop}
\begin{proof} If $\{x,x_2, ..., x_n\} = \{x,y\}$, then the expressions $s(x,x_2), ..., s(x,x_n)$ are all equal to either $x$ or $s(x,y)$, and at least one of them is equal to $s(x,y)$, so since $\{x,s(x,y)\}$ acts like a semilattice oriented from $x$ to $s(x,y)$ under $s$, we can combine these expressions in any order to produce such a term $s_n$.

For concreteness, we define $s_n$ inductively, as follows: $s_1(x) = x, s_2(x,y) = s(x,y)$ and
\[
s_n(x_1, ..., x_n) = s(s_{n-1}(x_1, ..., x_{n-1}), s(x_1,x_n)).\qedhere
\]
\end{proof}

Now we can use the terms $s_n$ to preprocess the inputs to $n$-ary functions. If $f$ is an $n$-ary term of $\bA$, define the term $f_s$ by
\[
f_s(x_1, ..., x_n) = f(s_n(x_1, ..., x_n), s_n(x_2, ..., x_n, x_1), ..., s_n(x_n, x_1, ..., x_{n-1})).
\]
As in the case of unary operations, we will consider the reduct $\bA_s$ with basic operations $f_s$ for every term $f$ of $\bA$. This reduct will be simpler in the sense that for any $a,b$ with $s(a,b) = s(b,a) = b$, each term $f_s$ will act like the $n$-ary semilattice operation on $\{a,b\}$. Additionally, every two-variable height $1$ identity which holds in $\bA$ will also hold in $\bA_s$.

\begin{prop}\label{semilattice-preparation} Let $\bA = (A, (f^i)_{i \in I})$ be a finite idempotent algebra, and let $\Sigma$ be the set of all two-variable height $1$ identities which involve both variables on each side and are satisfied in $\bA$. Then the operations $(f^i_s)_{i \in I}$ of $\bA_s$ will also satisfy the identities in $\Sigma$.

Additionally, if $\bB, \bC$ are subalgebras of $\bA$ such that for any $x \in \bB$ and any $y \in \bC$ we have $s(x,y), s(y,x) \in \bC$, then for any $n$-ary term $f$ of $\bA$ and any $x_1, ..., x_n \in \bB \cup \bC$ such that at least one $x_i \in \bC$, we have $f_s(x_1, ..., x_n) \in \bC$.
\end{prop}
\begin{proof} Suppose we have an identity
\[
f^i(a_1, ..., a_m) \approx f^j(b_1, ..., b_n),
\]
with $\{a_1, ..., a_m\} = \{b_1, ..., b_n\} = \{x,y\}$. Define $a_1', ..., a_m'$ by $a_k' = s(x,y)$ if $a_k = x$ and $a_k' = s(y,x)$ if $a_k = y$, and define $b_1', ..., b_n'$ similarly. Then for each $k$, we have
\[
s_m(a_k, ..., a_m, a_1, ..., a_{k-1}) \approx a_k',
\]
and similarly for the $b_l'$s, so
\[
f^i_s(a_1, ..., a_m) \approx f^i(a_1', ..., a_m') \approx f^j(b_1', ..., b_n') \approx f^j_s(b_1, ..., b_n).
\]

For the last statement, we just need to check that for any $x_1, ..., x_n \in \bB\cup \bC$ with at least one of the $x_i$s in $\bC$ we have $s_n(x_1, ..., x_n) \in \bC$ (since $\bC$ is closed under each term $f$ of $\bA$). This follows from the fact that $s_n$ is defined from $s$ in a way that involves all of its variables.
\end{proof}

Since an algebra $\bA$ is Taylor iff it satisfies a nontrivial system of two-variable height $1$ identities, if $\bA$ is Taylor then $\bA_s$ will also be Taylor. Later, we will see that algebras with bounded width are also characterized by two-variable height $1$ identities, so the same sort of implication (i.e. $\bA$ has bounded width implies $\bA_s$ has bounded width) will hold in that case as well. Algebras with few subpowers are \emph{not} characterized by height $1$ identities, essentially because no semilattice can have few subpowers, so such an implication fails in that case.

There are two other interesting cases which are not characterized by two-variable height $1$ identities: algebras of width $1$, and algebras such that the associated CSP is solved by the linear programming relaxation. It turns out that we can still prove a similar result in these cases.

\begin{prop} If $\bA$ has symmetric terms $f_n$ of every arity, then it has symmetric terms $f_n^s$ which act like the semilattice operation on each set $\{a,b\}$ with $s(a,b) = s(b,a) = b$.
\end{prop}
\begin{proof} Let $f_n$ be a symmetric term operation of arity $n$, for each $n$. Then for any $n$, let $\sigma_1, ..., \sigma_{n!}$ be an enumeration of the permutations of $\{1, ..., n\}$, and define $f_n^s$ by
\[
f_n^s(x_1, ..., x_n) \coloneqq f_{n!}(s_n(x_{\sigma_1(1)}, ..., x_{\sigma_1(n)}), ..., s_n(x_{\sigma_{n!}(1)}, ..., x_{\sigma_{n!}(n)})).
\]
Then $f_n^s$ is a symmetric term operation of arity $n$.
\end{proof}

\begin{prop}\label{prop-totally-symmetric-preparation} If $\bA$ has totally symmetric terms $f_n$ of every arity, then it has totally symmetric terms $f_n^s$ which act like the semilattice operation on each set $\{a,b\}$ with $s(a,b) = s(b,a) = b$.
\end{prop}
\begin{proof} Fix $n$. For every $m \ge 1$, let $S^n_m$ be the set of $n$-ary terms $t$ of $\bA$ such that there exist variables $y_1, ..., y_l$ with $\{y_1, ..., y_l\} = \{x_1, ..., x_n\}$ and such that for each $i$, the number of $j$ with $y_j = x_i$ is at least $m$, and
\[
t(x_1, ..., x_n) = s_l(y_1, ..., y_l).
\]
Note that for $m' > m$ we have $S^n_{m'} \subseteq S^n_m$, and each $S^n_m$ is finite and nonempty, so the intersection $S^n = \bigcap_m S^n_m$ is also finite and nonempty. Furthermore, for any $a_1, ..., a_n \in \bA$, the set of values
\[
\{t(a_1, ..., a_n) \mid t \in S^n\}
\]
depends only on the set $\{a_1, ..., a_n\}$. Thus we can take
\[
f_n^s(x_1, ..., x_n) \coloneqq f_{|S^n|}(\{t(x_1, ..., x_n) \mid t \in S^n\}).\qedhere
\]
\end{proof}

\begin{rem} The previous two propositions only used the fact that the restrictions of the $s_n$s to $\{a,b\}$ are symmetric and totally symmetric, respectively. So they can be generalized to show that if an algebra $\bA$ has symmetric/totally symmetric operations of each arity, then for every subset $X$ of $\bA$ such that some collection of terms $t_n$ of $\bA$ preserve $X$ and have symmetric/totally symmetric restrictions to $X$, we can find symmetric/totally symmetric operations of $\bA$ which preserve $X$ and such that their restrictions to $X$ agree with the restrictions of the terms $t_n$. It turns out that a similar general result holds for Taylor clones and clones of bounded width, but the proof of that will need to wait until we show that Taylor algebras have cyclic terms.
\end{rem}

Recall that for any $a,b$, the set $\{a,b\}$ is a semilattice subalgebra of $\bA$ iff the ternary relation $(x \in \{a,b\}) \wedge (x = a \implies y = z)$ defines a subalgebra of $\bA^3$. We can generalize this somewhat.

\begin{prop} If $B,C$ are subsets of $\bA$, then the ternary relation
\[
(x \in B\cup C) \wedge (x \not\in C \implies y=z)
\]
defines a subalgebra of $\bA^3$ iff $B \cup C$ is a subalgebra of $\bA$, and for any $n$, any $n$-ary term $f \in \Clo_n(\bA)$ which depends on all of its inputs, and any $x_1, ..., x_n \in B\cup C$ such that at least one $x_i \in C$, we have $f(x_1, ..., x_n) \in C$.
\end{prop}

\begin{defn} If $\bC \le \bB$ are subalgebras of $\bA$ such that there exists a term $t$ with $t(\bB,\bC), t(\bC,\bB) \subseteq \bC$, then we say that $\bC$ \emph{binary absorbs} $\bB$, and write $\bC \lhd_{bin} \bB$. If for any $n$, any $n$-ary term $f \in \Clo_n(\bA)$ which depends on all of its inputs, and any $x_1, ..., x_n \in \bB$ such that at least one $x_i \in \bC$, we have $f(x_1, ..., x_n) \in \bC$, then we say that $\bC$ \emph{strongly absorbs} $\bB$, and write $\bC \lhd_{str} \bB$.
\end{defn}

We can summarize the previous results in the following proposition, which shows that binary absorption and strong absorption are very nearly the same thing.

\begin{prop}\label{bin-abs-semi} If $\bC \lhd_{bin} \bB$, then there is a partial semilattice term $s$ with $s(\bB,\bC), s(\bC,\bB) \subseteq \bC$, and in the reduct $\bA_s$ the subalgebras $\bB_s, \bC_s$ satisfy $\bC_s \lhd_{str} \bB_s$. Furthermore, $\bC \lhd_{str} \bB$ iff the ternary relation $(x \in \bB) \wedge (x \not\in \bC \implies y = z)$ defines a subalgebra of $\bA^3$ (and $\bC \le \bB$).
\end{prop}

In general, a binary absorbing subalgebra of a binary absorbing subalgebra might not be binary absorbing (consider the $4$ element lattice $(\{0,1\}^2, \wedge, \vee)$ and the sequence $\{(0,1)\} \lhd_{bin} \{(0,0),(0,1)\} \lhd_{bin} \{0,1\}^2$), and similarly for strongly absorbing subalgebras (consider the idempotent commutative groupoid $(\{a,b,c\},\cdot)$ given by $ab = ac = b, bc = c$ and the sequence $\{c\} \lhd_{str} \{b,c\} \lhd_{str} \{a,b,c\}$). However, we can always chain together binary and strong absorption in one particular order.

\begin{prop} If $\bC \lhd_{bin} \bB \lhd_{str} \bA$, then $\bC \lhd_{bin} \bA$. Applying this repeatedly, we see that if
\[
\bC \lhd_{bin} \bB_n \lhd_{str} \cdots \lhd_{str} \bB_1 \lhd_{str} \bA,
\]
then $\bC \lhd_{bin} \bA$.
\end{prop}
\begin{proof} Suppose that $\bC$ absorbs $\bB$ with respect to the binary term $t$. Define a term $u$ by
\[
u(x,y) \coloneqq t(t(x,t(x,y)),t(y,t(y,x))).
\]
Then for any $a \in \bA, c \in \bC$, we have $t(a,c) \in \bB$ and $t(a,t(a,c)) \in \bB$ since $c \in \bB \lhd_{str} \bA$, so
\[
u(a,c) \in t(\bB,t(c,\bB)) \subseteq t(\bB,\bC) \subseteq \bC,
\]
and similarly $u(c,a) \in \bC$.
\end{proof}

By iteratively replacing $\bA$ with reducts $\bA_s$ for partial semilattice terms $s$ quadratically many times, we can reduce to the case where for all $a,b$, we have $(b,b) \in \Sg_{\bA^2}\{(a,b),(b,a)\}$ iff $\{a,b\}$ is a semilattice subalgebra of $\bA$ with absorbing element $b$.

\begin{defn}\label{defn-prepared} We say that an idempotent algebra $\bA$ has been \emph{prepared} if for every pair $a,b$ such that
\[
\begin{bmatrix} b\\ b\end{bmatrix} \in \Sg_{\bA^2}\Big\{\begin{bmatrix} a\\ b\end{bmatrix},\begin{bmatrix} b\\ a\end{bmatrix}\Big\},
\]
the set $\{a,b\}$ is a semilattice subalgebra of $\bA$ (in the sense that the restriction of any $k$-ary basic operation of $\bA$ to $\{a,b\}$ is the $k$-ary semilattice operation on $\{a,b\}$ with absorbing element $b$).
\end{defn}

For algebras which have been prepared, it makes sense to define a digraph whose edges correspond to semilattice subalgebras of $\bA$.

\begin{defn} If $s$ is a partial semilattice operation and $a,b$ have $s(a,b) = b$, then we write $a \rightarrow_s b$, or just $a \rightarrow b$ if $s$ is understood (or if the algebra has been prepared).
\end{defn}

\begin{thm} Let $s$ be a fixed nontrivial partial semilattice term operation of an idempotent algebra $\bA$. If $\bA$ is prepared, then the following are equivalent.
\begin{itemize}
\item[(a)] $s(a,b) = b$, that is, $a \rightarrow b$,

\item[(b)] the restriction of $s$ to $\{a,b\}$ acts like the semilattice operation on $\{a,b\}$ with absorbing element $b$,

\item[(c)] there exists $c$ such that $s(a,c) = b$,

\item[(d)] $\begin{bmatrix}b \\ b\end{bmatrix} \in \Sg_{\bA^2}\left\{\begin{bmatrix}a \\ b\end{bmatrix}, \begin{bmatrix}b \\ a\end{bmatrix}\right\}$

\item[(e)] there is a binary term $t$ of $\bA$ with $t(a,b) = t(b,a) = b$,

\item[(f)] there is a partial semilattice term $s'$ of $\bA$ with $s'(a,b) = b$,

\item[(g)] for every $n$ and every $n$-ary term $f \in \Clo_n(\bA)$ which depends on all its inputs, the restriction of $f$ to $\{a,b\}$ acts like the $n$-ary semilattice operation on $\{a,b\}$ with absorbing element $b$,

\item[(h)] the ternary relation $(x \in \{a,b\}) \wedge (x = a \implies y = z)$ defines a subalgebra of $\bA^3$.
\end{itemize}
If $\bA$ has not been prepared, then (a), (b), (c) are equivalent to each other, (d), (e), (f) are equivalent to each other, (g), (h) are equivalent to each other, and (g) implies (a) implies (d).
\end{thm}

\begin{prop}\label{prepared-bin-abs} If $\bA$ is prepared, then the following hold:
\begin{itemize}
\item[(a)] for $\bB \lhd_{bin} \bA$ and any $a \in \bA$, there is some $b \in \bB$ such that $a \rightarrow b$,
\item[(b)] if $\bB \lhd_{bin} \bA$ and $a \in \bA$, $b \in \bB$ have $b \rightarrow a$, then $a \in \bB$,
\item[(c)] if $\bC \lhd_{bin} \bB \lhd_{bin} \bA$, then $\bC \lhd_{bin} \bA$,
\item[(d)] if $\bB_1, \bB_2 \lhd_{bin} \bA$, then $\bB_1 \cap \bB_2 \ne \emptyset$ and $\bB_1 \cap \bB_2 \lhd_{bin} \bA$.
\end{itemize}
In particular, there is a unique minimal binary absorbing subalgebra $\bB \lhd_{bin} \bA$, and this $\bB$ has no proper binary absorbing subalgebra.
\end{prop}
\begin{proof} Part (a) follows from the existence of a partial semilattice term $s$ with $s(\bA,\bB) \subseteq \bB$ and part (b) follows from part (g) of the previous proposition.

For part (c), choose a partial semilattice term $s$ with $s(\bB,\bC), s(\bC,\bB) \subseteq \bC$, and choose any binary term $t$ with $t(\bA,\bB), t(\bB,\bA) \subseteq \bB$. Define a binary term $u$ by
\[
u(x,y) \coloneqq s(s(t(x,y),y),s(t(y,x),x)).
\]
Then for $a \in \bA, c \in \bC$ we have $t(a,c), t(c,a) \in \bB$, and we have $t(c,a) \rightarrow s(t(c,a),a)$, so by part (b) we have $s(t(c,a),a) \in \bB$. Thus
\[
u(a,c) \in s(s(\bB,c),\bB) \subseteq s(\bC,\bB) \subseteq \bC,
\]
and similarly $u(c,a) \in \bC$.

For part (d), pick $b_1 \in \bB_1$, then by part (a) there is some $b_2 \in \bB_2$ with $b_1 \rightarrow b_2$, and then by part (b) we have $b_2 \in \bB_1$, so $b_2 \in \bB_1 \cap \bB_2$. Then from $\bB_2 \lhd_{bin} \bA$ we have $\bB_1 \cap \bB_2 \lhd_{bin} \bB_1$, and we can apply part (c) to finish.
\end{proof}

\begin{prop}\label{arrow-partial} If $\bA$ has been prepared and $a,b,c \in \bA$ have $c \in \Sg\{a,b\}$ with $a \rightarrow c$, then $\bA$ has a partial semilattice term $s$ with $s(a,b) = c$.
\end{prop}
\begin{proof} Let $s'$ be an arbitrary nontrivial partial semilattice term operation of $\bA$, and choose $p$ a binary term operation of $\bA$ with $p(a,b) = c$. Then take $s(x,y) = s'(x,p(x,y))$. We clearly have $s(a,b) = s'(a,p(a,b)) = s'(a,c) = c$, so we just have to check that $s$ is a partial semilattice.

If $p$ is second projection then $s = s'$ and we are done. Otherwise, since $\bA$ has been prepared, $p$ and $s'$ both act as the semilattice operation on $\{x,s'(x,p(x,y))\} = \{x,s(x,y)\}$, so $s$ also acts as the semilattice operation on $\{x,s(x,y)\}$.
\end{proof}

In any digraph, the strongly connected components have a natural partial order.

\begin{defn} We say that $b$ is \emph{reachable} from $a$ if there is a sequence $a = a_0, a_1, ..., a_k = b$ such that $a_i \rightarrow a_{i+1}$ for $i = 0, ..., k-1$.
\end{defn}

\begin{prop} If $\bA$ is prepared and $s^1, ..., s^k$ are partial semilattice terms of $\bA$, then for any $n$-ary term $f \in \langle s^1, ..., s^k\rangle$, $f(x_1, ..., x_n)$ is always reachable from at least one of the variables $x_1, ..., x_n$.
\end{prop}

\begin{defn} We say that a subset $S$ of an algebra $\bA$ which has a partial semilattice operation $s$ is \emph{upwards closed} if whenever $a\in S$ and $a' \in \bA$ have $a \rightarrow_s a'$, we also have $a' \in S$.
\end{defn}

\begin{defn} We say that a set $A$ is \emph{strongly connected} if for every subset $S \subset A$ with $S\ne \emptyset, A$ there is an $a \in S$ and a $b \in A\setminus S$ such that $a \rightarrow b$. We say that a set $A$ is a \emph{maximal strongly connected component} of an algebra $\bA$ if $A$ is a strongly connected subset which is upwards closed (note that every finite upwards closed set contains at least one maximal strongly connected component). Finally, we call an element of an algebra $\bA$ \emph{maximal} if it is contained in any maximal strongly connected component of $\bA$.
\end{defn}

The main application of partial semilattice terms to CSPs is the following general idea: if a solvable instance of a CSP is arc-consistent (i.e. all relations are subdirect), then it probably has a solution where each variable is assigned a value in a maximal strongly connected component of the corresponding domain. So a basic case to try to understand is the case where every domain is a strongly connected algebra.

\begin{rem} The digraph considered in this section is the same as the set of ``thin red edges'' of Andrei Bulatov's colored graph \cite{colored-graph} attached to any Taylor algebra. Bulatov has a different construction of a partial semilattice operation $s$ from a binary term $t$, which is still based on the couterintuitive idea of taking a function $f$ which satisfies $f(x,f(x,y)) \approx f(x,y)$ and plugging in $f(x,f(y,x))$.
\end{rem}

\section{Maximal strongly connected components and polynomial completeness}

In this section we prove a few results of Andrei Bulatov \cite{bulatov-bounded} about the way maximal strongly connected components of partial semilattice algebras interact with binary and ternary relations. A consequence of the results of this section is that simple, strongly connected algebras are always polynomially complete. Throughout this section, we will always fix a partial semilattice operation $s$.

\begin{thm}\label{strong-binary} Fix a partial semilattice operation $s$. Suppose $\RR \le_{sd} \bA \times \bB$ is subdirect and $A,B$ are maximal strongly connected subsets of $\bA, \bB$, respectively.
\begin{itemize}
\item[(a)] The set of $a$ such that $(\{a\} \times B) \cap \RR \ne \emptyset$ is upwards closed. In particular, if $(A\times B)\cap \RR$ is nonempty, then it is subdirect in $A\times B$.

\item[(b)] The set of $a$ such that $\{a\} \times B \subseteq \RR$ is upwards closed.

\item[(c)] If $A$ is contained in a linked component of $\RR$ (that is, a connected component of $\RR$ considered as a bipartite graph on $\bA \sqcup \bB$), $(A\times B) \cap \RR \ne \emptyset$, and $A,B$ are finite, then $A \times B \subseteq \RR$.
\end{itemize}
Additionally, the product $A\times B$ is a maximal strongly connected subset of $\bA \times \bB$.
\end{thm}
\begin{proof} For part (a), suppose that $(a,b) \in \RR$ and $b \in B$, and let $a \rightarrow a'$. Since $\RR$ is subdirect, there is some $b'$ with $(a',b') \in \RR$. Then
\[
\begin{bmatrix} a'\\ s(b,b')\end{bmatrix} = s\left(\begin{bmatrix} a\\ b\end{bmatrix}, \begin{bmatrix} a'\\ b'\end{bmatrix}\right) \in \RR,
\]
and $b \rightarrow s(b,b')$, so $s(b,b') \in B$.

For part (b), suppose that $\{a\}\times B \subseteq \RR$ and $a \rightarrow a'$. Let $S$ be the set of $b \in B$ such that $(a',b) \in \RR$, that is, $S = \pi_2((\{a'\}\times B) \cap \RR)$. By part (a), $S$ is nonempty. To finish, we just have to show that $S$ is upwards closed. Suppose $b \in S$ and $b \rightarrow b'$. Then by assumption we have $(a,b') \in \RR$, so
\[
\begin{bmatrix} a'\\ b'\end{bmatrix} = s\left(\begin{bmatrix} a'\\ b\end{bmatrix}, \begin{bmatrix} a\\ b'\end{bmatrix}\right) \in \RR.
\]

For part (c), suppose first that $A \times A \subseteq \RR\circ \RR^-$, where $\RR^- \le \bB \times \bA$ is the reverse of $\RR$ (we will later reduce the general case to this case). Let $a$ be any element of $A$, and let $X$ be the set of $b \in \bB$ such that $(a,b) \in \RR$, that is, $X = \pi_2((\{a\}\times \bB) \cap \RR)$. By part (a), $X \cap B \ne \emptyset$, and by the finiteness of $B$, the intersection $X \cap B$ has a maximal strongly connected component $S$. Since $B$ is a maximal strongly connected component of $\bB$, $S$ is a maximal strongly connected component of $X$.

By the assumption $A \times A \subseteq \RR\circ \RR^-$ and the definition of $X$, we see that $(A\times X)\cap \RR$ is subdirect in $A\times X$. Thus by part (b) and the fact that $\{a\} \times S \subseteq (A\times X)\cap \RR$, we see that $A \times S \subseteq (A\times X)\cap \RR$, so $A \times S \subseteq \RR$. Then by part (b) applied to $\RR^-$, we see that $A \times B \subseteq \RR$.

Now suppose that $A \times A \not\subseteq \RR\circ \RR^-$. From the finiteness of $A$ we see that there is some $k$ such that $A \times A \subseteq (\RR\circ \RR^-)^{\circ k}$. Choose $k$ minimal, and let $\RR' = (\RR \circ \RR^-)^{\circ (k-1)} \le_{sd} \bA^2$. Then $\RR'$ is equal to its own reverse $\RR'^-$, and $A \times A \subseteq \RR' \circ \RR'$ since $2(k-1) \ge k$ for $k \ge 2$. Thus the previous paragraphs applied to $\RR'$ (using $\RR' = \RR'^-$) show that $A \times A \subseteq \RR'$, contradicting the minimality of $k$.
\end{proof}


\begin{cor}\label{maximal-projection} If $\pi : \bA \twoheadrightarrow \bB$ is a surjective homomorphism of finite algebras, then the subalgebra of $\bA$ generated by the maximal elements of $\bA$ maps surjectively onto the subalgebra of $\bB$ generated by the maximal elements of $\bB$.
\end{cor}

\begin{cor} If we start with any arc-consistent instance of $\CSP(\bA_1, ..., \bA_n)$ and replace every domain and every relation by the subalgebra generated by its maximal elements, then the resulting instance will still be arc-consistent.
\end{cor}

\begin{cor}\label{strong-simple} Fix a partial semilattice operation $s$. Suppose that $\RR \le_{sd} \bA\times \bB$ is a subdirect product of finite algebras $\bA, \bB$, and that $\bB$ is simple and $\bB = \Sg(B)$, with $B$ a maximal strongly connected component of $\bB$. Then:
\begin{itemize}
\item[(a)] if $\bA$ is also simple and $\bA = \Sg(A)$ with $A$ a maximal strongly connected component of $\bA$, and if $\RR \cap (A\times B) \ne \emptyset$, then $\RR$ is either the graph of an isomorphism or is $\bA\times \bB$, and

\item[(b)] if $\bA$ is arbitrary and $\RR$ is not the graph of a homomorphism from $\bA$ to $\bB$, then there is an $a \in \bA$ with $\{a\}\times \bB \subseteq \RR$.
\end{itemize}
\end{cor}
\begin{proof} If $\bB$ is simple, then the linking congruence of $\RR$ on $\bB$ must either be the trivial congruence $0_\bB$, in which case $\RR$ is the graph of a homomorphism from $\bA$ to $\bB$, or the full congruence $1_\bB$, in which case $\RR$ is linked. In the second case, the results follow from Theorem \ref{strong-binary}(c).
\end{proof}

\begin{thm}\label{strong-ternary} Fix a partial semilattice operation $s$. Suppose $R \subseteq A \times B \times C$ is closed under $s$, $A$ is strongly connected, $\pi_{23}(R)$ is strongly connected, $\pi_{12}(R) = A \times B$, $\pi_{13}(R) = A\times C$, and $A,B,C$ are finite. Then $R = A \times \pi_{23}(R)$.
\end{thm}
\begin{proof} By Theorem \ref{strong-binary}(c), we just need to show that $R$ is linked as a subset of $A \times \pi_{23}(R)$. We will do this by showing that for any $a \rightarrow a'$ in $A$, some fork of $R$ links $a$ to $a'$ in one step.

Since $\pi_1(R) = A$, there exist $b \in B, c \in C$ such that $(a,b,c) \in R$. Since $\pi_{13}(R) = A\times C$, there exists some $b' \in B$ such that $(a',b',c) \in R$. Since
\[
\begin{bmatrix} a'\\ s(b,b')\\ c\end{bmatrix} = s\left(\begin{bmatrix} a\\ b\\ c\end{bmatrix}, \begin{bmatrix} a'\\ b'\\ c\end{bmatrix}\right) \in R,
\]
we may assume without loss of generality that $b' = s(b,b')$, that is, that $b \rightarrow b'$.

Since $\pi_{12}(R) = A\times B$, there exists some $c' \in C$ such that $(a,b',c') \in R$. Since
\[
\begin{bmatrix} a\\ b'\\ s(c,c')\end{bmatrix} = s\left(\begin{bmatrix} a\\ b\\ c\end{bmatrix}, \begin{bmatrix} a\\ b'\\ c'\end{bmatrix}\right) \in R,
\]
we may assume without loss of generality that $c' = s(c,c')$, that is, that $c \rightarrow c'$.

Since $(a',b',c)$ and $(a,b',c')$ are in $R$, we have
\[
\begin{bmatrix} a'\\ b'\\ c'\end{bmatrix} = s\left(\begin{bmatrix} a'\\ b'\\ c\end{bmatrix}, \begin{bmatrix} a\\ b'\\ c'\end{bmatrix}\right) \in R.
\]
Thus both $a$ and $a'$ meet $(b',c') \in \pi_{23}(R)$.
\end{proof}

\begin{rem} The proof of Theorem \ref{strong-ternary} actually proves something slightly more general: if $R \subseteq A\times B \times C$ is closed under $s$, $\pi_{12}(R) = A\times B, \pi_{13}(R) = A\times C$, and $A$ is weakly connected, then $R$ is linked when considered as a subalgebra of $A \times \pi_{23}(R)$.
\end{rem}

\begin{cor} Fix a partial semilattice operation $s$. Suppose $R \subseteq A_1 \times \cdots \times A_n$ is closed under $s$, $A_1$ is strongly connected, $\pi_{[2,n]}(R)$ is strongly connected, $\pi_{1i}(R) = A_1 \times A_i$ for $i \in [2,n]$, and $A_i$ are finite for all $i$. Then $R = A_1 \times \pi_{[2,n]}(R)$.
\end{cor}

\begin{cor}\label{strong-product} Fix a partial semilattice operation $s$. Suppose $R \subseteq A_1 \times \cdots \times A_n$ is closed under $s$, all $A_i$ are strongly connected, $\pi_{ij}(R) = A_i \times A_j$ for all $i \ne j$, and $A_i$ are finite for all $i$. Then $R = A_1 \times \cdots \times A_n$.
\end{cor}

\begin{cor}\label{simple-strong-poly} Fix a partial semilattice operation $s$. If $\bA$ is simple and is generated by a finite maximal strongly connected component $A$, then $\bA$ is polynomially complete.
\end{cor}
\begin{proof} We just need to show that every relation $\RR \le \bA^n$ which contains the set of constant tuples $\Delta_n = \{(a, ..., a) \mid a \in \bA\}$ is an intersection of equality relations. First consider the case $n = 2$. From the assumption that $\bA$ is simple we see that either $\RR$ is the equality relation, or $\RR$ is linked. If $\RR$ is linked, then Theorem \ref{strong-binary}(c) and the fact that $\RR$ contains $\Delta_2$ implies that $A \times A \subseteq \RR$, and from the assumption $\bA = \Sg(A)$ we see that $\RR = \bA \times \bA$.

Now consider the case $n \ge 3$. If any two-variable projection $\pi_{ij}(\RR)$ is the equality relation, then we can ignore one of the coordinates $i,j$, so we may assume without loss of generality that $\pi_{i,j}(\RR) = \bA\times \bA$ for all $i,j$. Let $a$ be any element of $A$, and let $R$ be a maximal strongly connected component of $\RR$ which is reachable from $(a, ..., a)$. Then for any $i,j$ we must have $\pi_{i,j}(R) = A\times A$, so by the previous corollary we have $R = A^n$. Thus $A^n \subseteq \RR$, and from the assumption $\bA = \Sg(A)$ we see that $\RR = \bA^n$.
\end{proof}

\begin{ex}\label{ex-subdirect-2-semi-not-absorption-free} The reader may be wondering whether we can weaken the assumption that $\pi_{23}(R)$ is strongly connected from Theorem \ref{strong-ternary} to the assumption that $B,C$ are strongly connected. It seems plausible that if $B,C$ are both strongly connected and $\pi_{23}(R)$ is a subdirect product of $B$ and $C$, $\pi_{23}(R)$ might automatically be strongly connected.

However, there is an example of a strongly connected $2$-semilattice $\bA$ and a subdirect product $\RR \le_{sd} \bA^2$ which is \emph{not} strongly connected. The $2$-semilattice $\bA$ is pictured below.
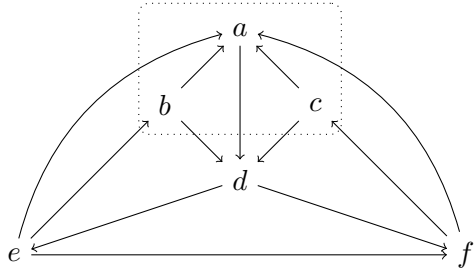
\begin{center}
\begin{tikzpicture}[scale=1]
  \node (a) at (0,3) {$a$};
  \node (b) at (-1,2) {$b$};
  \node (c) at (1,2) {$c$};
  \node (d) at (0,1) {$d$};
  \node (e) at (-3,0) {$e$};
  \node (f) at (3,0) {$f$};
  \draw [->] (b) edge (a) (c) edge (a);
  \draw [->] (d) edge (e) (d) edge (f) (e) edge (f);
  \draw [->] (a) edge (d) (b) edge (d) (c) edge (d);
  \draw [->] (e) edge (b) (f) edge (c);
  \draw [->, bend left] (e) to (a);
  \draw [->, bend right] (f) to (a);
  \node [fit=(a)(b)(c),draw,dotted,rounded corners] {};
\end{tikzpicture}
\end{center}
The missing values are given by $s(b,c) = s(b,f) = s(e,c) = a$.

If we let $\theta \le_{sd} \bA^2$ be the smallest congruence containing $(b,c)$, then $\theta$ corresponds to the partition $\{a,b,c\},\{d\},\{e\},\{f\}$, and $\bA/\theta$ is a four element tournament. Considering $\theta$ as an algebra, we find that $\theta$ is \emph{not} strongly connected: $(b,c)$ and $(c,b)$ are incomparable minimal elements of $\theta$, and the remaining elements of $\theta$ form a maximal strongly connected component.
\end{ex}

\begin{ex} Here we will give an example of a subdirect product of strongly connected algebras which has two maximal strongly connected components (such an example is necessarily not a $2$-semilattice, since every $2$-semilattice has a unique maximal strongly connected component).

As in the previous example, we will consider a congruence $\theta$ on a six-element algebra $\bA$. This time $\bA/\theta$ will be the three-element rock-paper-scissors algebra, and every congruence class of $\bA$ will have two elements, with $s$ acting as $\pi_1$ on the congruence class. As a digraph, $\bA$ is just a directed six-cycle, pictured below.
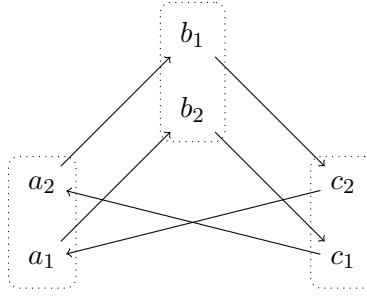
\begin{center}
\begin{tikzpicture}[scale=1]
  \node (a1) at (-2,0) {$a_1$};
  \node (a2) at (-2,1) {$a_2$};
  \node (b1) at (0,3) {$b_1$};
  \node (b2) at (0,2) {$b_2$};
  \node (c1) at (2,0) {$c_1$};
  \node (c2) at (2,1) {$c_2$};
  \draw [->] (a1) edge (b2) (a2) edge (b1);
  \draw [->] (b1) edge (c2) (b2) edge (c1);
  \draw [->] (c1) edge (a2) (c2) edge (a1);
  \node [fit=(a1)(a2),draw,dotted,rounded corners] {};
  \node [fit=(b1)(b2),draw,dotted,rounded corners] {};
  \node [fit=(c1)(c2),draw,dotted,rounded corners] {};
\end{tikzpicture}
\end{center}
Given the above digraph structure and the assumption that there is a congruence $\theta$ corresponding to the partition $\{a_1,a_2\},\{b_1,b_2\},\{c_1,c_2\}$, there is only one way to fill in the values of the partial semilattice operation $s$. The reader can check that the congruence $\theta$, considered as a subalgebra of $\bA^2$, has two maximal strongly connected components which are both isomorphic to $\bA$.
\end{ex}

Despite the above examples, we do at least have the following result, which is important for understanding how restricting to maximal strongly connected components interacts with cycle-consistency.

\begin{thm}\label{strong-diagonal} Fix a partial semilattice operation $s$. Suppose that $R \subseteq A \times A$ is closed under $s$, $A$ is finite and strongly connected, and $R$ contains the diagonal $\Delta_A = \{(a,a) \mid a \in A\}$. Then $R$ has a maximal strongly connected component which contains $\Delta_A$.
\end{thm}
\begin{proof} Since $\Delta_A$ is strongly connected, it's enough to show that if $(a,b)$ is reachable from $(a,a)$ in $R$, then some element $(c,c)$ of $\Delta_A$ is reachable from $(a,b)$ in $R$. We will define a unary polynomial $\phi$ of $R$ such that $\phi((a,a)) = (a,b)$ and such that $\phi(x)$ is reachable from $x$ in $R$ for all $x \in R$.

To construct $\phi$, choose some sequence $(a_i,b_i) \in R$ such that $(a,a) = (a_0,b_0)$, $(a_i,b_i) \rightarrow (a_{i+1},b_{i+1})$ for all $i$, and $(a_k,b_k) = (a,b)$ for some $k$. Then define $\phi$ by
\[
\phi(x) = s\left(s\left(\cdots s\left(s\left(x,\begin{bmatrix}a_1\\b_1\end{bmatrix}\right), \begin{bmatrix}a_2\\b_2\end{bmatrix}\right),\cdots\right), \begin{bmatrix}a_k\\b_k\end{bmatrix}\right).
\]
Note that since $\phi((a,a)) = (a,b)$, we have $\pi_1(\phi((a,x))) = a$ for all $x \in A$.

Since $A$ is finite, we can find $m \ge 1$ such that $\phi^{\circ 2m} = \phi^{\circ m}$. Define another unary polynomial $\phi_\Delta$ of $R$ by
\[
\phi_\Delta(x) = s\left(s\left(\cdots s\left(s\left(x,\begin{bmatrix}b_1\\b_1\end{bmatrix}\right), \begin{bmatrix}b_2\\b_2\end{bmatrix}\right),\cdots\right), \begin{bmatrix}b_k\\b_k\end{bmatrix}\right),
\]
that is, by replacing each $(a_i,b_i)$ in the definition of $\phi$ by $(b_i,b_i)$. Then if $\phi^{\circ m}((a,a)) = (a,c)$, we have
\[
\phi_\Delta^{\circ m}\left(\phi^{\circ (m-1)}\left(\begin{bmatrix} a\\ b\end{bmatrix}\right)\right) = \phi_\Delta^{\circ m}\left(\begin{bmatrix} a\\ c\end{bmatrix}\right) = \begin{bmatrix} c\\ c\end{bmatrix}.
\]
Thus $(c,c)$ is reachable from $(a,b)$ in $R$.
\end{proof}

\begin{cor}\label{cycle-consistency-maximal} If we start with any cycle-consistent instance of $\CSP(\bA_1, ..., \bA_n)$ and replace every domain and every relation by the subalgebra generated by its maximal elements, then the resulting instance will still be cycle-consistent.
\end{cor}
\begin{proof} By Theorem \ref{strong-binary}(a), we just need to check this in the special case where our cycle-consistent instance is a cycle of binary relations $\RR_i \le_{sd} \bA_i\times \bA_{i+1}$ with indices taken modulo $n$. Let $\RR \le_{sd} \bA_1 \times \cdots \times \bA_n \times \bA_1$ be the relation given by the formula
\[
(x_1,x_2) \in \RR_1 \wedge \cdots \wedge (x_n,x_{n+1}) \in \RR_n.
\]
The assumption that the instance is cycle-consistent implies that $\Delta_{\bA_1} \subseteq \pi_{1,n+1} \RR$. Set $\RR_\Delta = \pi_{1,n+1} \RR$.

For any algebra $\bA$, let $\bA^{\max}$ denote the subalgebra of $\bA$ generated by the maximal elements of $\bA$. We see from Theorem \ref{strong-binary}(a) that $\pi_{i,i+1}(\RR^{\max}) = \RR_i^{\max}$ for each $i$ and that $\pi_{1,n+1}(\RR^{\max}) = \RR_\Delta^{\max}$. By Theorem \ref{strong-diagonal} we have $\Delta_{\bA_1^{\max}} \subseteq \RR_\Delta^{\max}$, so the new instance is cycle-consistent at the first variable.
\end{proof}

\section{$2$-semilattices, spirals, and ancestral algebras}

In this section we'll discuss a pretty general class of partial semilattice algebras which are nice enough for the associated CSP to have bounded width, due to Bulatov \cite{colored-graph-prelim}. Following the strategy of replacing domains of variables with the subalgebras generated by their maximal elements, and noting that many of the structural results proved in the preceeding section apply best to strongly connected algebras, we see that it would be quite convenient if every domain of every variable in our CSP has a unique maximal strongly connected component. The most straightforward examples of algebras with this property are $2$-semilattices.

\begin{defn} A binary operation $s$ is a $2$-\emph{semilattice} operation if it satisfies the identities
\[
s(x,y) \approx s(y,x), \;\;\; s(x,s(x,y)) \approx s(x,y), \;\;\; s(x,x) \approx x.
\]
In other words, a $2$-semilattice is a partial semilattice operation which is also commutative.
\end{defn}

\begin{prop} An algebra $\bA = (A,s)$ is a $2$-semilattice iff for all $a,b \in \bA$, the subalgebra $\Sg_{\bA}\{a,b\}$ is a semilattice under $s$.
\end{prop}

\begin{prop} If $\bA = (A,s)$ is a finite $2$-semilattice, then $\bA$ has a unique maximal strongly connected component.
\end{prop}
\begin{proof} If $a, b$ are any two maximal elements of $\bA$, then $s(a,b) = s(b,a)$ is reachable from both $a$ and $b$, so $a$ and $b$ must be in the same maximal strongly connected component.
\end{proof}

The first difficult results about bounded width CSPs were proved for $2$-semilattices. However, the proofs only depended on the fact that every $2$-semilattice has a unique maximal strongly connected component. Bulatov \cite{colored-graph-prelim} calls this the ``maximal red component condition''. I've chosen to call such algebras ``ancestral'' instead, because they can be equivalently defined as follows.

\begin{defn} An idempotent algebra $\bA$ with a fixed partial semilattice operation $s$ is called \emph{ancestral} if for all $a,b \in \bA$, there is some $c \in \Sg_{\bA}\{a,b\}$ which is reachable from both $a$ and $b$. We call any such $c$ a \emph{common ancestor} of $a$ and $b$.
\end{defn}

\begin{prop}\label{ancestral-maximal} A finite idempotent algebra $\bA$ is ancestral iff every proper subalgebra of $\bA$ has a unique maximal strongly connected component.
\end{prop}

A nice generalization of $2$-semilattices is the collection of algebras which I call ``spirals''. Spirals are defined in terms of a single commutative binary operation, so they can be described more rapidly than general ancestral algebras. As we will see later, a minimal Taylor clone is ancestral if and only if it is a minimal spiral, so we would not lose too much generality by restricting the study of ancestral algebras to the study of spirals.

\begin{defn} An algebra $\bA = (A,f)$ is a \emph{spiral} if $f$ is a commutative idempotent binary operation and every subalgebra of $\bA$ which is generated by two elements either has size two or has a surjective homomorphism to the free semilattice on two generators.
\end{defn}

\begin{ex} Here we give an example of a minimal spiral $\bA_6$ which is not a $2$-semilattice.
\begin{center}
\begin{tabular}{cc}
\begin{tabular}{c | c c c c c c} $\bA_6$ & $a$ & $b$ & $c$ & $d$ & $e$ & $f$\\ \hline $a$ & $a$ & $c$ & $e$ & $d$ & $e$ & $d$\\ $b$ & $c$ & $b$ & $c$ & $c$ & $f$ & $f$\\ $c$ & $e$ & $c$ & $c$ & $c$ & $e$ & $c$\\ $d$ & $d$ & $c$ & $c$ & $d$ & $d$ & $d$\\ $e$ & $e$ & $f$ & $e$ & $d$ & $e$ & $f$\\ $f$ & $d$ & $f$ & $c$ & $d$ & $f$ & $f$ \end{tabular} &
\begin{tikzpicture}[scale=1.5,baseline=-0.4cm]
  \node (a) at (-1,-1) {$a$};
  \node (b) at (1,-1) {$b$};
  \node (c) at (-0.4,-0.1) {$c$};
  \node (d) at (0.4,-0.1) {$d$};
  \node (e) at (-0.4,0.6) {$e$};
  \node (f) at (0.4,0.6) {$f$};
  \draw [->] (a) edge (d) (a) edge (e) (b) edge (c) (b) edge (f);
  \draw [->] (c) edge (e) (e) edge (f) (f) edge (d) (d) edge (c) (f) edge (c) (e) edge (d);
\end{tikzpicture}
\end{tabular}
\end{center}
Every proper subalgebra of $\bA_6$ is a $2$-semilattice - in fact, every pair of elements other than $\{a,b\}$ generates a two or three element semilattice subalgebra of $\bA_6$. The pair $\{a,b\}$ generates $\bA_6$, and $\bA_6$ has a congruence $\theta$ corresponding to the partition $\{a\}, \{b\}, \{c,d,e,f\}$ such that $\bA_6/\theta$ is isomorphic to the free semilattice on two generators.

The reader may check that any nonempty subset $S$ of $\bA_6$ which is closed under multiplication by $a$ and by $b$ must necessarily contain all four of $c,d,e,f$ - using this observation, it is easy to check that $\Clo(\bA_6)$ contains no nontrivial proper subclones.
\end{ex}

\begin{thm} If $\bA = (A,f)$ is a spiral, then for any partial semilattice term $s \in \Clo(f)$ which is defined nontrivially in terms of $f$, the reduct $\bA_s = (A,s)$ is ancestral.
\end{thm}
\begin{proof} We prove this by induction on the size of $A$. Let $a,b$ be any two elements of $A$. If $\Sg_\bA\{a,b\}$ has size two, then since $f$ is commutative we must either have $a \rightarrow b$ or $b \rightarrow a$, so one of $a,b$ is a common ancestor of $a$ and $b$.

Otherwise, by the definition of a spiral, there is a surjective homomorphism $\alpha$ from $\Sg_\bA\{a,b\}$ to the free semilattice on two generators. Clearly $a$ and $b$ must be sent to the two generators of the free semilattice by $\alpha$, say $\alpha(a) = x$ and $\alpha(b) = y$, and every nontrivial binary term $t \in \Clo(f)$ must have $\alpha(t(a,b)) = t(x,y) = f(x,y)$. Thus the kernel of $\alpha$ has congruence classes $\{a\}, \{b\},$ and $S = \Sg_\bA\{a,b\}\setminus\{a,b\}$, and $S$ is a binary absorbing subalgebra of $\bA$ with respect to $f$.

Since $S$ is a binary absorbing subalgebra of $\bA$ with respect to $f$ and $s \in \Clo(f)$ is defined nontrivially, we must have $s(a,b), s(b,a) \in S$. Since $|S| \le |A|-2$, we can apply the inductive hypothesis to see that $s(a,b),s(b,a)$ have a common ancestor in $\Sg_{(S,s)}\{s(a,b),s(b,a)\} \subseteq \Sg_{\bA_s}\{a,b\}$.
\end{proof}

\begin{ex}\label{ex-ancestral-strong} An example of an ancestral algebra which is not a $2$-semilattice or a spiral is the algebra $\bA_4 = (\{a,b,c,d\},s)$, where $s$ is the partial semilattice operation described below.
\begin{center}
\begin{tabular}{cc}
\begin{tabular}{c|cccc} $s$ & $a$ & $b$ & $c$ & $d$\\ \hline $a$ & $a$ & $b$ & $b$ & $a$\\ $b$ & $b$ & $b$ & $c$ & $c$\\ $c$ & $d$ & $c$ & $c$ & $d$\\ $d$ & $a$ & $a$ & $d$ & $d$\end{tabular} & 
\begin{tikzpicture}[scale=1.5,baseline=0.5cm]
  \node (a) at (-0.5,0.9) {$a$};
  \node (b) at (0.5,0.9) {$b$};
  \node (c) at (0.5,-0.1) {$c$};
  \node (d) at (-0.5,-0.1) {$d$};
  \draw [->] (a) edge (b) (b) edge (c) (c) edge (d) (d) edge (a);
\end{tikzpicture}
\end{tabular}
\end{center}
The algebra $\bA_4$ has the cyclic automorphism $(a\ b\ c\ d)$, and is generated by the pair $a,c$, since $s(a,c) = b, s(c,a) = d$. The binary term $s'$ given by
\[
s'(x,y) \coloneqq s(x,s(y,x))
\]
is another (nontrivial) partial semilattice term operation of $\bA_4$, such that $s'(a,c) = a, s'(c,a) = c$. So the reduct $(\{a,b,c,d\},s')$ of $\bA_4$ is \emph{not} an ancestral algebra, as it has the subalgebra $(\{a,c\},s')$ which has the two maximal strongly components $\{a\}$ and $\{c\}$.

It is easy to check that $\bA_4$ is simple, and every proper subalgebra of $\bA_4$ is a two element semilattice. By Corollary \ref{simple-strong-poly}, $\bA_4$ is polynomially complete, and in fact Theorem \ref{strong-binary} and Theorem \ref{strong-ternary} imply that every subdirect relation $\RR \le_{sd} \bA_4^n$ can be written as an intersection of two variable relations, each of which is the graph of an automorphism of $\bA_4$. In particular, if we consider the ternary relation
\[
\RR_{ac} = \Sg_{\bA^3}\left\{\begin{bmatrix}a\\ a\\ c\end{bmatrix}, \begin{bmatrix}a\\ c\\ a\end{bmatrix}, \begin{bmatrix}c\\ a\\ a\end{bmatrix}\right\},
\]
we find that $\RR_{ac} = \bA_4^3$. Since there is an automorphism of $\bA_4$ which interchanges $a$ and $c$, we see that there are ternary terms $g, g' \in \Clo(\bA_4)$ such that $\{a,c\}$ is closed under $g$ and $g'$, with $(\{a,c\},g)$ a two element majority algebra and $(\{a,c\}, g')$ a two element affine algebra. Either of the reducts $(\{a,b,c,d\},g)$ or $(\{a,b,c,d\},g')$ defines a Taylor algebra, since $g$ satisfies the identity
\[
g(x,x,y) \approx g(x,y,x) \approx g(y,x,x) \approx s'(x,y),
\]
and $g'$ satisfies the similar identity
\[
g'(x,x,y) \approx g'(x,y,x) \approx g'(y,x,x) \approx s'(y,x).
\]
\end{ex}

\begin{prop} Every quotient of an ancestral algebra is ancestral.
\end{prop}

\begin{thm} If $\bA_1, ..., \bA_n$ are ancestral algebras with partial semilattice operation $s$, then so is $\bA_1 \times \cdots \times \bA_n$.
\end{thm}
\begin{proof} We prove this by induction on $n$. Let $a,b \in \bA_1 \times \cdots \times \bA_n$. Since $\bA_1$ is ancestral, there is some $c_1 \in \Sg_{\bA_1}\{a_1,b_1\}$ which is reachable from both $a_1$ and $b_1$. Lifting the path from $a_1$ to $c_1$ to a path from $a$ to some element $c' \in \Sg\{a,b\}$ with $c_1' = c_1$, and lifting the path from $b_1$ to $c_1$ to a path from $b$ to some $c'' \in \Sg\{a,b\}$ with $c_1'' = c_1$, we see that we just need to find a common ancestor of $c'$ and $c''$. Since $c_1' = c_1''$ and $\bA_1$ is idempotent, we see that $c',c''$ have a common ancestor so long as $\bA_2 \times \cdots \times \bA_n$ is ancestral, which follows from the inductive hypothesis.
\end{proof}

\begin{cor} If $\bA$ is ancestral and $\bB \in HSP_{fin}(\bA)$, then $\bB$ is also ancestral.
\end{cor}

It turns out that ancestral algebras can be defined entirely in terms of collections of partial semilattice operations.

\begin{thm}\label{ancestral-terms} A finite idempotent algebra $\bA$ with a fixed partial semilattice operation $s$ is ancestral iff for some $m \ge n \ge 0$ it has a sequence of partial semilattice terms $p_1, p_2, ..., p_m$ such that
\begin{itemize}
\item $a \rightarrow_s p_i(a,b)$ for all $a,b \in \bA$ and all $i$,
\item $a \rightarrow_s b$ implies $p_i(a,b) = b$ for all $i$, and
\item if we define binary operations $f_i$ recursively by $f_0(x,y) \coloneqq s(x,y)$ and
\[
f_i(x,y) \coloneqq p_i(f_{i-1}(x,y),f_{i-1}(y,x))
\]
for $i \ge 1$, then $f_m(x,y) \approx f_n(y,x)$.
\end{itemize}
\end{thm}
\begin{proof} That the existence of such a sequence implies $\bA$ is ancestral follows from the fact that for any $a,b$, each $f_i(a,b)$ is reachable from $a$ and each $f_j(b,a)$ is reachable from $b$.

For the converse direction, let $\bF = \cF_{\bA}(x,y) \le \bA^{\bA^2}$ be the free algebra on two generators in the variety generated by $\bA$. Since $\bF \in SP_{fin}(\bA)$, $\bF$ is ancestral, so there is some sequence of elements $f_0, ..., f_n \in \bF$ with $f_0(x,y) = s(x,y)$, such that each $f_{i-1} \rightarrow_s f_i$, each $f_i \in \Sg_{\bF}\{f_{i-1}(x,y),f_{i-1}(x,y)\}$, and such that the subset $S$ of elements of the subalgebra $\bS = \Sg_{\bF}\{f_n(x,y),f_n(y,x)\}$ which are reachable from $f_n$ in $\bS$ is minimal given these constraints. Then $S$ must be strongly connected, and for every $g \in S$ we must have $\bS = \Sg_{\bF}\{g(x,y), g(y,x)\}$. Thus we can extend our sequence $f_0, ..., f_n$ by $f_{n+1}, ..., f_m$ such that each $f_{i-1} \rightarrow_s f_i$, and $f_m(x,y) \approx f_n(y,x)$, and we will automatically have $f_i \in \Sg_{\bF}\{f_{i-1}(x,y),f_{i-1}(x,y)\}$ for each $i$.

Note that $f_{i-1} \rightarrow_s f_i$ and $f_i \in \Sg_{\bF}\{f_{i-1}(x,y),f_{i-1}(x,y)\}$ implies the existence of a binary term $p_i$ such that $x \rightarrow_s p_i(x,y)$ and $f_i(x,y) = p_i(f_{i-1}(x,y), f_{i-1}(y,x))$, by the argument of Proposition \ref{arrow-partial}. Note that the reduct with basic operations $s, f_i$ is ancestral, and has the property that $a \rightarrow_s b$ implies $f_i(a,b) = f_i(b,a) = b$ for all $i$, so $\{a,b\}$ is a semilattice subalgebra with respect to any nontrivial binary term in $\Clo(f_0, ..., f_m)$. Thus we may assume without loss of generality that $a \rightarrow_s b$ implies $p_i(a,b) = b$ for all $i$, and then the argument of Proposition \ref{arrow-partial} implies that each $p_i$ is a partial semilattice term.
\end{proof}

In fact, we can go further: every ancestral algebra has an ancestral reduct which is prepared. Recall that $\bA$ is \emph{prepared} if for all $a,b \in \bA$, we have $(b,b) \in \Sg_{\bA^2}\{(a,b),(b,a)\}$ iff $\{a,b\}$ is a semilattice subalgebra of $\bA$ with $a \rightarrow b$.

\begin{thm} Every finite ancestral algebra $\bA$ has a reduct which is prepared and ancestral.
\end{thm}
\begin{proof} Let $s, f_i$ be as in Theorem \ref{ancestral-terms}, and assume without loss of generality that these are the basic operations of $\bA$. Suppose there is a pair $a,b \in \bA$ with $(b,b) \in \Sg_{\bA^2}\{(a,b),(b,a)\}$ but $s(a,b) \ne b$. Let $s'$ be a partial semilattice term with $s'(a,b) = b$. Then $c \rightarrow_{s} d$ implies $c \rightarrow_{s'} d$, and if we define
\[
f_0'(x,y) \coloneqq s'(x,y)
\]
and
\[
f_i'(x,y) \coloneqq f_{i-1}(s'(x,y),s'(y,x))
\]
for $i \ge 1$, then the reduct with basic operations $s', f_i'$ is an ancestral algebra (with respect to $s'$) with strictly more semilattice subalgebras than $\bA$.
\end{proof}

Due to the structural simplifications we can obtain by passing to reducts, it makes sense to focus on ancestral algebras such that no proper reduct is also ancestral.

\begin{defn} A finite algebra $\bA$ is called a \emph{minimal ancestral} algebra if $\bA$ is ancestral, and no proper reduct of $\bA$ is ancestral.
\end{defn}

Since every minimal ancestral algebra is automatically prepared, we don't need to specify a particular choice of partial semilattice operation to define the digraph of semilattice subalgebras.

\begin{prop} Every finite ancestral algebra has a reduct which is a minimal ancestral algebra.
\end{prop}
\begin{proof} Whether an algebra is ancestral only depends on the collection of partial semilattice operations in its clone. Since there are only finitely many partial semilattice operations on a given finite set, we don't need to worry about infinite descending chains of smaller and smaller ancestral reducts.
\end{proof}

\begin{prop} If $\bA$ is a minimal ancestral algebra and $\bB \in HSP_{fin}(\bA)$, then $\bB$ is also a minimal ancestral algebra.
\end{prop}
\begin{proof} Let $f_i$ be terms for $\bA$ as in Theorem \ref{ancestral-terms}. If we can find a proper reduct of $\bB$ which is ancestral, then there is a sequence of terms $f_i'$ of this reduct such that $f_0'(x,y) \approx x$, $f_i'(x,y) \rightarrow f_{i+1}'(x,y)$, and $f_m'(a,b) = f_n'(b,a)$ holds for all $a,b \in \bB$. Then if we define additional terms $f_{m+i}'$ by
\[
f_{m+i}'(x,y) \coloneqq f_i(f_m'(x,y), f_n'(y,x)),
\]
we see that these terms $f_0', ..., f_m', f_{m+1}', ...$ generate the same reduct on $\bB$ as $f_0', ..., f_m'$, and generate an ancestral reduct of $\bA$.
\end{proof}

\begin{thm}\label{ancestral-absorption} If $\bA$ is a minimal ancestral algebra, then for any $a,b \in \bA$, if $S$ is the maximal strongly connected component of $\Sg_\bA\{a,b\}$, then we have $\Sg_\bA\{a,b\} = S \cup \{a,b\}$. If $\{a,b\} \not\subseteq S$, then $\Sg_\bA\{a,b\}$ has a semilattice quotient with $S$ as a congruence class which acts as the top element.
\end{thm}
\begin{proof} Choose terms $f_i$ as in Theorem \ref{ancestral-terms}. Let $\bF = \cF_\bA(x,y)$ be the free algebra on two generators in the variety generated by $\bA$. Pick any element $g(x,y)$ in the maximal strongly connected component of $\bF$, and note that since $g(x,y)$ is reachable from both $x$ and $y$ in $\bF$, every term $t(x_1, ..., x_k) \in \Clo(g)$ which depends on all its inputs has the property that $t(x_1, ..., x_k)$ is reachable from each $x_i$ in $\cF_\bA(x_1, ..., x_k)$.

Applying the semilattice iteration argument, we get a partial semilattice term $s'(x,y) \in \Clo(g)$, which is reachable from each of $x$, $y$, and $g(x,y)$ in $\bF$. In particular, we see that $s'(x,y)$ is contained in the maximal strongly connected component of $\bF$, and if we define terms $f_i'$ by
\[
f_0'(x,y) \coloneqq s'(x,y)
\]
and
\[
f_i'(x,y) \coloneqq f_{i-1}(s'(x,y),s'(y,x))
\]
for $i \ge 1$, then the reduct with basic operations $f_i'$ is an ancestral algebra, and each $f_i'(x,y)$ is contained in the maximal strongly connected component of $\bF$. Thus the clone generated by the $f_i'$s must be equal to the clone of $\bA$, and we see that every element of $\bF$ is either equal to one of $x,y$ or is contained in the maximal strongly connected component of $\bF$.
\end{proof}

\begin{cor}\label{minimal-ancestral-maximal} If $\bA$ is a minimal ancestral algebra, then the maximal strongly connected component of $\bA$ is a strongly absorbing subalgebra of $\bA$.
\end{cor}

There is a sense in which even the class of minimal ancestral algebras is unnecessarily large: it contains algebras such as the algebra $\bA_4$ from Example \ref{ex-ancestral-strong} which have proper Taylor reducts with two element majority or affine subalgebras.

\begin{thm} Suppose $\bA$ is a minimal ancestral algebra which is generated by $a$ and $b$, is strongly connected, and is simple. Then there are ternary terms $g,g' \in \Clo(\bA)$ such that $\{a,b\}$ is closed under $g$ and $g'$, $(\{a,b\},g)$ is a two element majority algebra, and $(\{a,b\},g')$ is a two element affine algebra.
\end{thm}
\begin{proof} Let $\bS = \Sg_{\bA^2}\{(a,b),(b,a)\}$. If $\bS$ is linked, then by Theorem \ref{strong-binary}(c) we must have $(b,b) \in \bS$, so $a \rightarrow b$, a contradiction. Otherwise, $\bS$ is the graph of an automorphism swapping $a$ and $b$. In this case, the ternary relation $\RR = \Sg_{\bA^3}\{(a,a,b),(a,b,a),(b,a,a)\}$ has $(a,a),(a,b),(b,a) \in \pi_{i,j}(\RR)$ for each $i,j$, so by Theorem \ref{strong-binary}(c) we have $\pi_{i,j}(\RR) = \bA^2$, and then by Theorem \ref{strong-ternary} we have $\RR = \bA^3$. Thus $(a,a,a) \in \RR$ and $(b,b,b) \in \RR$, and we can take $g,g'$ to be ternary terms of $\bA$ which witness these facts.
\end{proof}

Later we will see that the above result implies that a minimal ancestral algebra which is both strongly connected and generated by two elements has a proper Taylor reduct (and, in fact, has a proper bounded width reduct). For now we will show that minimal ancestral algebras which avoid this situation are actually spirals.

\begin{thm} If $\bA$ is a minimal ancestral algebra such that for all $a,b$ the subalgebra $\Sg_{\bA}\{a,b\}$ has no strongly connected quotient, then $\bA$ is term equivalent to a spiral.
\end{thm}
\begin{proof} Let $s$ be a nontrivial partial semilattice operation on $\bA$. Define a sequence of terms $f_i$ inductively by $f_0 \coloneqq s$ and
\[
f_{i+1}(x) \coloneqq f_i(s(x,y),s(y,x)).
\]
We will show by induction on $|\bA|$ that for each $a,b \in \bA$, there is an $n$ such that $f_n(a,b) = f_n(b,a)$. To see this, note that by Theorem \ref{ancestral-absorption}, for any $a,b$ the subalgebra generated by $s(a,b), s(b,a)$ is contained in the maximal strongly connected component $S$ of $\Sg_\bA\{a,b\}$, so as long as $S \ne \bA$ we can apply the induction hypothesis to see that there is some $i$ such that
\[
f_i(s(a,b),s(b,a)) = f_i(s(b,a),s(a,b)),
\]
and for this $i$ we then have $f_{i+1}(a,b) = f_{i+1}(b,a)$.

Thus there is some $n$ such that $f = f_n$ is commutative (in fact, we can take $n = |\bA|$). To finish, we need to show that if $\bA$ is generated by two elements $a,b$ with $|\bA| > 2$, then the maximal strongly connected component $S$ of $\bA$ does not contain either of $a,b$. To this end, suppose for a contradiction that $S$ contains $b$. Let $\bS = \Sg_{\bA^2}\{(a,b),(b,a)\}$. If $S$ is contained in a linked component of $\bS$, then by Theorem \ref{strong-binary}(c) we must have $(b,b) \in \bS$, so $a \rightarrow b$, a contradiction. Otherwise, the linking congruence $\theta \in \Con(\bA)$ of $\bS$ has $|S/\theta| > 1$ and $b/\theta \in S/\theta$, and so we may assume without loss of generality that $\theta$ is trivial. But if $\theta$ is trivial, then $\bA$ has an automorphism which interchanges $a$ and $b$, so $S$ contains both $a$ and $b$, so $\bA$ is both strongly connected and generated by two elements, a contradiction.
\end{proof}

\section{Cycle-consistency solves ancestral CSPs}

In this section we will prove that any cycle-consistent instance of an ancestral CSP has a solution. This proof is a simple case of Kozik's proof \cite{slac} of the fact that cycle-consistency solves CSPs over templates with bounded width: the main purpose of presenting the argument in this special case is to allow the reader to focus on the overall proof strategy before getting into the technical algebraic details.

The ingredients which we will need for the proof are the following facts about ancestral algebras.
\begin{itemize}
\item Every ancestral algebra $\bA$ has a unique maximal strongly connected component $\bA^{\max}$ (Proposition \ref{ancestral-maximal}).

\item If $\pi : \bA \twoheadrightarrow \bB$ is a surjective homomorphism, then $\pi(\bA^{\max}) = \bB^{\max}$ (Corollary \ref{maximal-projection} to Theorem \ref{strong-binary}(a)).

\item If $\RR \le_{sd} \bA \times \bB$ and $\bA^{\max}$ is contained in a linked component of $\RR$, then $\RR^{\max} = \bA^{\max}\times \bB^{\max}$ (Theorem \ref{strong-binary}(c)).

\item In particular, if $\RR \le_{sd} \bA \times \bB$, $\bA$ is generated by $\bA^{\max}$, $\bB$ is generated by $\bB^{\max}$, and $\bB$ is simple, then $\RR$ is either the graph of a homomorphism $\bA \twoheadrightarrow \bB$ or $\RR = \bA \times \bB$ (Corollary \ref{strong-simple}).

\item If $\RR \le_{sd} \bA \times \bB \times \bC$ has $\pi_{12}(\RR) = \bA \times \bB$ and $\pi_{13}(\RR) = \bA \times \bC$, then $\RR^{\max} = \bA^{\max} \times \pi_{23}(\RR)^{\max}$ (Theorem \ref{strong-ternary}).

\item Applying the above inductively, if $\RR \le_{sd} \bA_1 \times \cdots \times \bA_n$ has $\pi_{ij}(\RR) = \bA_i \times \bA_j$ for all $i \ne j$, then $\RR^{\max} = \bA_1^{\max} \times \cdots \times \bA_n^{\max}$ (Corollary \ref{strong-product}).

\item If $\RR \le_{sd} \bA\times \bA$ and $\RR$ contains the diagonal $\Delta_\bA$, then $\Delta_{\bA^{\max}} \subseteq \RR^{\max}$ (Theorem \ref{strong-diagonal}).

\item If we start with any cycle-consistent instance of $\CSP(\bA_1, ..., \bA_n)$ and replace every domain and every relation by the subalgebra generated by its maximal elements, then the resulting instance will still be cycle-consistent (Corollary \ref{cycle-consistency-maximal}).
\end{itemize}
If we assume that our algebras are minimal ancestral (rather than just ancestral), then each $\bA^{\max}$ becomes a subalgebra (Corollary \ref{minimal-ancestral-maximal}), which slightly simplifies the arguments. We won't use this simplification, but the reader should keep it in mind.

The general strategy is to start with a cycle-consistent instance, and to find a way to shrink some of the variable domains and relations to get a strictly smaller cycle-consistent instance. Eventually, we reach a situation where all the variable domains have size $1$ and the instance is still cycle-consistent - at this point, there is obviously a solution to the CSP. We have already seen that by shrinking variable domains, we can reach a situation where each variable domain $\bA_x$ is generated by $\bA_x^{\max}$ (the last bullet point above).

To finish the argument, we need to find another strategy for reducing the variable domains when each $\bA_x = \Sg(\bA_x^{\max})$. The intuition is that if $\bA_x = \Sg(\bA_x^{\max})$, then there is some congruence $\theta_x \in \Con(\bA_x)$ such that $\bA_x/\theta_x$ is simple, and in fact $\bA_x/\theta_x$ will be polynomially complete by Corollary \ref{simple-strong-poly}. Since polynomially complete algebras should have few interesting subdirect relations, it's plausible that we can replace the domain $\bA_x$ with an arbitrary congruence class of $\theta_x$, and always obtain a cycle-consistent instance.

So fix a variable $x$ with $|\bA_x| > 1$, a maximal congruence $\theta_x$ in $\Con(\bA_x)$, and a congruence class $\bA_x'$ of $\theta_x$. We now have to restrict the other variable domains in order to, at the very least, get an arc-consistent sub-instance. We will show that a very minimalistic sort of reduction strategy suffices: instead of worrying about all possible issues with ensuring arc-consistency, we will only consider paths from variables $y$ to $x$ through the instance.

\begin{defn}\label{path-defn} If $\mathbf{X}$ is an instance of a CSP and $x,y$ are variables of $\fX$, then a \emph{path} $p$ from $x$ to $y$ is defined as a sequence $x = v_0, (\RR_1,i_1,j_1), v_1, ..., v_{n-1}, (\RR_n,i_n,j_n), v_n = y$ such that each $v_k$ is a variable, and each $\RR_k$ is a relation such that one of the constraints of the instance $\fX$ imposes the relation $\RR_k$ on a tuple $u = (u_1, ...)$ of variables with $u_{i_k} = v_{k-1}$ and $u_{j_k} = v_k$.

To every path $p$ from $x$ to $y$, we associate the binary relation $\PP_p \le \bA_x \times \bA_y$ which is given by
\[
\PP_p \coloneqq \pi_{i_1j_1}(\RR_1) \circ \cdots \circ \pi_{i_nj_n}(\RR_n).
\]
In other words, $\PP_p$ is the set of pairs of values in $\bA_x \times \bA_y$ which are consistent with the path $p$.

We define addition and negation of paths in the natural way, so that if $p$ is a path from $x$ to $y$ and $q$ is a path from $y$ to $z$, then $p+q$ is a path from $x$ to $z$ with $\PP_{p+q} = \PP_p \circ \PP_q$, and $-p$ is a path from $y$ to $x$ with $\PP_{-p} = \PP_p^-$.
\end{defn}

In particular, we see that an instance is arc-consistent iff for all paths $p$ the associated binary relations $\PP_p$ are subdirect, and it is cycle-consistent iff we additionally have $\Delta_{\bA_v} \subseteq \PP_p$ for every path $p$ from a variable $v$ back to itself.

\begin{defn}\label{defn-proper-red} Suppose that $\fX$ is a cycle-consistent instance such that for all variable domains we have $\bA_v = \Sg(\bA_v^{\max})$, that $x$ is any variable with $|\bA_x| > 1$, that $\theta_x$ is any maximal congruence on $\bA_x$, and that $\bA_x'$ is any congruence class of $\bA_x/\theta_x$.

For each variable $y$, we say that $y$ is \emph{proper} if there is a path $p$ from $y$ to $x$ such that $\PP_p/\theta_x \le \bA_y \times \bA_x/\theta_x$ is the graph of a homomorphism $\iota_y : \bA_y \twoheadrightarrow \bA_x/\theta_x$. In this case, we define the congruence $\theta_y \in \Con(\bA_y)$ to be the kernel of $\iota_y$, and we define $\bA_y'$ to be the preimage of $\bA_x'$ under $\iota_y$. If $y$ is not proper, then we define $\bA_y'$ to be $\bA_y$.

We define the reduced instance $\fX'$ by replacing the domain of each variable $v$ by $\bA_v'$, and replacing each constraint relation $\RR \le \bA_{v_1} \times \cdots \times \bA_{v_m}$ of $\fX$ by $\RR' = \RR \cap (\bA_{v_1}' \times \cdots \times \bA_{v_m}')$.
\end{defn}

The reason for the name ``proper'' is that a variable $v$ is proper iff the reduced domain $\bA_v'$ is a proper subalgebra of $\bA_v$. First we need to check that the maps $\iota_y$ for the proper variables $y$ are well-defined.

\begin{lem} If $y$ is a proper variable and $p,q$ are two paths from $y$ to $x$ such that $\PP_p/\theta_x, \PP_q/\theta_x$ are graphs of homomorphisms $\iota_p, \iota_q : \bA_y \twoheadrightarrow \bA_x/\theta_x$, then in fact we have $\iota_p = \iota_q$. Thus $\iota_y, \theta_y$, and $\bA_y'$ are all well-defined.
\end{lem}
\begin{proof} The path $p-q$ connects $y$ to itself, so by cycle-consistency we must have $\Delta_{\bA_y} \subseteq \PP_{p-q} = \PP_p \circ \PP_q^-$. Taking the quotient by $\theta_x$, we see that $\Delta_{\bA_y} \subseteq (\PP_p/\theta_x) \circ (\PP_q/\theta_x)^-$, so for every element $a \in \bA_y$ we must have $\iota_p(a) = \iota_q(a)$.
\end{proof}

We sometimes abuse notation, and think of $\iota_y$ as an isomorphism from $\bA_y/\theta_y$ to $\bA_x/\theta_x$.

\begin{lem}\label{ancestral-proper-path} Suppose $p$ is a path from $y$ to a proper variable $z$. Then one of the following is true:
\begin{itemize}
\item $\PP_p/\theta_z = \bA_y \times \bA_z/\theta_z$, or

\item $y$ is also proper, and $\PP_p/(\theta_y\times\theta_z)$ is the graph of an isomorphism $\iota_p : \bA_y/\theta_y \xrightarrow{\sim} \bA_z/\theta_z$ such that $\iota_y = \iota_z \circ \iota_p$.
\end{itemize}
\end{lem}
\begin{proof} This follows from Corollary \ref{strong-simple} and cycle-consistency (note that $\bA_z/\theta_z$ is simple, since it is isomorphic to $\bA_x/\theta_x$).
\end{proof}

We have the ingredients necessary to check that the reduced instance $\fX'$ is cycle-consistent. We start with arc-consistency.

\begin{lem}\label{ancestral-red-arc} Suppose $\RR \le_{sd} \bA_{v_1} \times \cdots \times \bA_{v_n}$ is a constraint of $\fX$. Then the reduced constraint $\RR' = \RR \cap (\bA_{v_1}' \times \cdots \times \bA_{v_n}')$ is subdirect inside $\bA_{v_1}' \times \cdots \times \bA_{v_n}'$, that is, $\pi_i(\RR') = \bA_{v_i}'$ for each $i$.
\end{lem}
\begin{proof} By symmetry, it's enough to prove that $\pi_1(\RR') = \bA_{v_1}'$. In other words, for each element $a \in \bA_{v_1}'$, we want to find a tuple $s \in \RR$ such that $s_i \in \bA_{v_i}'$ for all $i$. We may ignore variables $v_i$ such that $i \ne 1$ and $v_i$ is not proper, since for such $i$ the restriction from $\bA_{v_i}$ to $\bA_{v_i}' = \bA_{v_i}$ has no effect. Similarly, for any two proper variables $v_i, v_j$ such that $\pi_{ij}(\RR)$ induces an isomorphism between $\bA_{v_i}/\theta_{v_i}$ and $\bA_{v_j}/\theta_{v_j}$, we may ignore one of the two variables $v_i, v_j$, since any element $s \in \RR$ which satisfies $s_i \in \bA_{v_i}'$ will automatically also satisfy $s_j \in \bA_{v_j}'$.

To formalize the process of ignoring variables, we define an equivalence relation $\sim$ on the set of indices of proper variables of $\RR$, with $i \sim j$ when $\pi_{ij}(\RR)$ induces an isomorphism between $\bA_{v_i}/\theta_{v_i}$ and $\bA_{v_j}/\theta_{v_j}$ (that $\sim$ is an equivalence relation is easy to check). Then we let $I \subseteq [n]$ be a set of variable indices such that each $\sim$-class has exactly one representative in $I$, $1 \in I$, and no index of any non-proper variable other than possibly $1$ is in $I$. We then define a relation $\bS \le \bA_{v_1} \times \prod_{i \in I \setminus\{1\}} \bA_{v_i}/\theta_{v_i}$ by
\[
\bS \coloneqq \pi_I(\RR)\Big/\prod_{i \in I \setminus\{1\}} \theta_{v_i}.
\]

We just need to show that for every $a \in \bA_{v_1}'$ there is some $s \in \bS$ with $s_1 = a$ and $s_i = \bA_{v_i}'/\theta_{v_i}$ for each $i \in I \setminus\{1\}$. Note that by Lemma \ref{ancestral-proper-path} and the construction of $I$, for every pair $i,j \in I$ the projection $\pi_{ij}(\bS)$ is full. Thus by Corollary \ref{strong-product}, we in fact have
\[
\bS^{\max} = \bA_{v_1}^{\max} \times \prod_{i \in I \setminus\{1\}} \bA_{v_i}^{\max}/\theta_{v_i},
\]
and since each $\bA_{v_i}$ is generated by $\bA_{v_i}^{\max}$, we have
\[
\bS = \bA_{v_1} \times \prod_{i \in I \setminus\{1\}} \bA_{v_i}/\theta_{v_i}.\qedhere
\]
\end{proof}

Now we can check that cycle-consistency also holds for the reduced instance.

\begin{lem}\label{ancestral-red-cycle} Suppose $p$ is a path from $v$ to $v$ in the instance $\fX$, and let $p'$ be the corresponding path in $\fX'$. If $\Delta_{\bA_{v_1}} \subseteq \PP_p$, then $\Delta_{\bA_{v_1}'} \subseteq \PP_{p'}$.
\end{lem}
\begin{proof} Suppose that $p$ is the path $v = v_0, (\RR_1,i_1,j_1), v_1, ..., v_{n-1}, (\RR_n,i_n,j_n), v_n = v$. Note that in the corresponding path $p'$, we must replace each $\RR_i$ with $\RR_i'$, so we must also worry about the proper variables which occur in $\RR_i$ but do not lie along the path $p$. In order to do this cleanly, we consider the relation $\RR$ defined by
\[
\RR \coloneqq \Big\{(v_0, u^1, ..., u^n, v_n) \in \bA_v \times \prod_{i \le n} \RR_i \times \bA_v \; \Big| \;
v_0 = u^1_{i_1}, u^1_{j_1} = u^2_{i_2}, ..., u^{n-1}_{j_{n-1}} = u^n_{i_n}, u^n_{j_n} = v_n\Big\}.
\]
If each $\RR_i$ has arity $m_i$, then $\RR$ is thought of as a relation of arity $m = 2 + \sum_i m_i$, and the indices of $\RR$ might contain several copies of variables of the instance $\fX$. Let the $i$th index of $\RR$ correspond to the variable $y_i$ in $\fX$, with $y_1 = v_0 = v$ and $y_m = v_n = v$, so
\[
\RR \le_{sd} \bA_{y_1} \times \cdots \times \bA_{y_m}.
\]
Note that by the arc-consistency of the instance $\fX$, for any two indices $i,j$ of the relation $\RR$, the projection $\pi_{ij}(\RR)$ is the same as $\PP_q$ for some path $q$ from $y_i$ to $y_j$ formed out of the relations $\RR_i$, and that $\pi_{1m}(\RR) = \PP_p$, so $\pi_{1m}(\RR) \supseteq \Delta_{\bA_v}$.

As in the argument for arc-consistency, we define an equivalence relation $\sim$ on the proper indices of $\RR$ defined by $i \sim j$ when $\pi_{ij}(\RR)$ induces an isomorphism between $\bA_{y_i}/\theta_{y_i}$ and $\bA_{y_j}/\theta_{y_j}$. We let $I \subseteq [m]$ to be a set of indices of $\RR$ with $1, m \in I$, such that $I$ contains no indices of non-proper variables of $\RR$ other than possibly $1$ and $m$, such that $I \setminus \{m\}$ contains one representative from each $\sim$ class of $\{1, ..., m-1\}$, and such that $I \setminus \{1\}$ contains one representative from each $\sim$ class of $\{2, ..., m\}$. As before, we define a relation $\bS$ by
\[
\bS \coloneqq \pi_I(\RR)\Big/\prod_{i \in I \setminus\{1,m\}} \theta_{y_i}.
\]

We just need to show that for every $a \in \bA_v'$, there is some $s \in \bS$ with $s_1 = s_m = a$ and $s_i = \bA_{v_i}'/\theta_{v_i}$ for each $i \in I \setminus\{1,m\}$. By Lemma \ref{ancestral-proper-path} and the construction of $I$, for every pair $i,j \in I$ with $\{i,j\} \ne \{1,m\}$ the projection $\pi_{ij}(\bS)$ is full. Thus by Corollary \ref{strong-product}, we have
\[
\pi_{I\setminus\{m\}}(\bS) = \bA_{y_1} \times \prod_{i \in I \setminus\{1,m\}} \bA_{y_i}/\theta_{y_i}
\]
and
\[
\pi_{I\setminus\{1\}}(\bS) = \bA_{y_m} \times \prod_{i \in I \setminus\{1,m\}} \bA_{y_i}/\theta_{y_i}.
\]
Thus by Theorem \ref{strong-ternary}, we have
\[
\bS^{\max} = \pi_{1m}(\bS)^{\max} \times \prod_{i \in I \setminus\{1,m\}} \bA_{y_i}^{\max}/\theta_{y_i},
\]
and by Theorem \ref{strong-diagonal} and the assumption $\pi_{1m}(\bS) = \pi_{1m}(\RR) \supseteq \Delta_{\bA_v}$, we have $\pi_{1m}(\bS)^{\max} \supseteq \Delta_{\bA_v^{\max}}$. Since each $\bA_y$ is generated by $\bA_y^{\max}$, we have
\[
\bS \supseteq \Delta_{\bA_v} \times \prod_{i \in I \setminus\{1,m\}} \bA_{y_i}/\theta_{y_i},
\]
so in particular for every $a \in \bA_v'$ we have $\{a\} \times \prod_{i \in I \setminus\{1,m\}} \bA_{y_i}'/\theta_{y_i} \times \{a\} \subseteq \bS$, so $(a,a) \in \pi_{1m}(\RR') = \PP_{p'}$.
\end{proof}

Thus the reduced instance $\fX'$ is cycle-consistent. Since we can iteratively shrink our instance whenever some variable $x$ has $\bA_x \ne \Sg(\bA_x^{\max})$ or has $\bA_x = \Sg(\bA_x^{\max})$ but $|\bA_x| > 1$, we see that we eventually reach a situation where each $\bA_x$ consists of a single element, and then arc-consistency proves that this collection of single elements gives a solution to the original instance. We have proved our main result.

\begin{thm} If $\fX$ is a cycle-consistent instance of an ancestral CSP, then $\fX$ has a solution.

In fact, for any variable $x$ of $\fX$, and for any element $a \in \bA_x$ such that there is a sequence of subalgebras $\bA_x \supseteq \bA_0 \supseteq \cdots \supseteq \bA_n = \{a\}$ with $\bA_0 = \Sg(\bA_x^{\max})$ and such that for each $i$, there is a maximal congruence $\theta_i \in \Con(\bA_i)$ and a congruence class $\bA_i'$ of $\theta_i$ with $\bA_{i+1} = \Sg(\bA_i'^{\max})$, there is a solution to the instance $\fX$ in which $x$ is assigned the value $a$.
\end{thm}

The simple construction of the reduced instance $\fX'$ can be used to show that we can find a solution to any cycle-consistent instance of an ancestral CSP in linear time.

\section{Cycle-consistency solves majority CSPs}

The paper which prompted the study of cycle-consistency was a preliminary investigation by Chen, Dalmau, and Gru{\ss}ien \cite{arc}, which studied a slightly stronger consistency notion: singleton arc-consistency. Singleton arc-consistency refers to the strategy of fixing a particular value for some variable, and checking if applying arc-consistency to the remaining variables produces a contradiction. Singleton arc-consistency is clearly at least as powerful as cycle-consistency, and can be established in quadratic time using an algorithm from \cite{singleton-arc-consistency-quadratic}. One of the main results of \cite{arc} showed that singleton arc consistency solves majority CSPs, but in fact their proof strategy was to show that cycle-consistent instances of majority CSPs always have solutions.

The argument for majority algebras is simpler than the argument for ancestral algebras, essentially because the analogue of the case where all the variables domains are strongly connected doesn't need to be considered. Instead, we are always in the situation where some variable domain $\bA_x$ has a proper absorbing subalgebra (every singleton is an absorbing subalgebra of a majority algebra), although we need to work slightly harder than we did in the absorbing case of ancestral CSPs since the absorption is no longer binary absorption. Rather than working with absorbing subalgebras, \cite{arc} used the closely related concept of an \emph{ideal} of a majority algebra.

\begin{defn} If $\bA = (A,m)$ is a majority algebra, then $\bB \le \bA$ is called an \emph{ideal} of $\bA$ if $m(\bB,\bA,\bB) \subseteq \bB$.
\end{defn}

The word ``ideal'' comes from the theory of median algebras - a subset $\bB$ is an ideal of a median algebra $\bA$ iff there is a congruence $\theta$ of $\bA$ such that $\bB$ is a congruence class of $\theta$. The corresponding statement is not true of majority algebras in general: every subset of the dual discriminator algebra from Example \ref{ex-dual-discriminator} is an ideal, but the dual discriminator algebra on $n$ elements is simple (and polynomially complete) for $n \ge 3$.

The next result shows that ideals interact with standard algebraic constructions (products, quotients, intersections) nicely. A similar result holds for absorbing subalgebras, with the same proof.

\begin{prop} Suppose that a relation $\RR$ is defined by a primitive positive formula $\Phi$ involving the relations $\RR_1, ..., \RR_k$. If we replace each $\RR_i$ with an ideal $\RR_i'$ of $\RR_i$ to make a primitive positive formula $\Phi'$, then the relation $\RR'$ which is defined by $\Phi'$ is an ideal of $\RR$.
\end{prop}
\begin{proof} Let $\Phi(x) = \exists y \Psi(x,y)$, with $\Psi$ quantifier-free, and let $\Psi'$ be the corresponding formula with $\RR_i$s replaced by $\RR_i'$s. Then for any $a,b,c$ with $a,c \in \RR'$ and $b \in \RR$, there exist $d,e,f$ such that $\Psi'(a,d), \Psi(b,e), \Psi'(c,f)$ hold, so $\Psi'(m(a,b,c),m(d,e,f))$ holds since each $\RR_i'$ is an ideal, so $\Phi'(m(a,b,c))$ holds.
\end{proof}

Recall the definition of a path in an instance (Definition \ref{path-defn}). It's notationally convenient to allow paths to act on subsets of the variable domains.

\begin{defn}\label{path-action-defn} If $p$ is a path connecting variables $x,y$ of an instance $\fX$, and if $B$ is a subset of the variable domain $\bA_x$, then we define $B + p$ to be the subset of $\bA_y$ given by
\[
B + p \coloneqq \{c \in \bA_y \mid \exists b \in B \text{ s.t. } (b,c) \in \PP_p\} = \pi_2(\PP_p\cap (B \times \bA_y)).
\]
\end{defn}

\begin{prop} If $\bB \le \bA_x$ and $p$ is a path from $x$ to $y$, then $\bB + p$ is a subalgebra of $\bA_y$. If $\bB$ is an ideal of $\bA_x$ and the instance is arc-consistent, then $\bB + p$ is an ideal of $\bA_y$.
\end{prop}

Our overall strategy will be to start with a cycle-consistent instance $\fX$, and find a collection of ideals $\bA_x'$ of the variable domains $\bA_x$ such that reducing each domain to $\bA_x'$ produces an arc-consistent instance $\fX'$. Then we will show that any such $\fX'$ is automatically cycle-consistent.

In order to find an arc-consistent family of ideal subdomains, we consider the set $\mathcal{I}$ of pairs $(x,\bB)$ where $x$ is a variable and $\bB$ is a proper ideal of $\bA_x$. Note that $\mathcal{I}$ is nonempty as long as some $x$ has $|\bA_x| > 1$, since every singleton is an ideal.

\begin{defn} Let $\mathcal{I}$ be the set of pairs $(x,\bB)$ where $x$ is a variable and $\bB$ is a proper ideal of $\bA_x$. We define a quasiorder $\preceq$ on $\mathcal{I}$ by $(x,\bB) \preceq (y, \bB+p)$ for every path $p$ from $x$ to $y$ with $\bB + p \ne \bA_y$.
\end{defn}

\begin{prop} If $\fX$ is a cycle-consistent instance, $x$ is a variable, and $(x,\bB) \preceq (x, \bC)$, then $\bB \le \bC$.
\end{prop}
\begin{proof} Suppose $p$ is a path from $x$ to itself with $\bB + p = \bC$. By cycle-consistency we must have $\Delta_{\bA_x} \subseteq \PP_p$, so $\bB \subseteq \bB + p$.
\end{proof}

\begin{defn} Suppose $\fX$ is a cycle-consistent instance of a majority CSP, and assume without loss of generality that each constraint of $\fX$ is binary. Fix a maximal element $(x,\bA_x')$ of $\mathcal{I}$ under the quasiorder $\preceq$.

Call a variable $y$ \emph{proper} if there is a path $p$ from $x$ to $y$ such that $\bA_x' + p \ne \bA_y$, and in this case set $\bA_y' = \bA_x' + p$. If $y$ is not proper, then set $\bA_y' = \bA_y$.

Define the reduced instance $\fX'$ by replacing the domain of each variable $v$ by $\bA_v'$, and by replacing each constraint $\RR \le \bA_u \times \bA_v$ with $\RR' = \RR \cap (\bA_u' \times \bA_v')$.
\end{defn}

First we need to check that the sets $\bA_y'$ are well-defined.

\begin{lem} If there are paths $p,q$ from $x$ to $y$ such that $\bA_x' + p \ne \bA_y$ and $\bA_x' + q \ne \bA_y$, then $\bA_x' + p = \bA_x' + q$.
\end{lem}
\begin{proof} Since $(x,\bA_x')$ is maximal and $(x,\bA_x') \preceq (y, \bA_x' + p)$, we must have $(y, \bA_x' + p) \preceq (x,\bA_x') \preceq (y, \bA_x' + q)$, so $\bA_x' + p \le \bA_x' + q$. Similarly we have $\bA_x' + q \le \bA_x' + p$, so $\bA_x' + p = \bA_x' + q$.
\end{proof}

Next we check arc-consistency.

\begin{lem} If $p$ is a path from $y$ to $z$ and $p'$ is the corresponding path in $\fX'$, then $\bA_y' + p' = \bA_z'$.
\end{lem}
\begin{proof} We just need to check this in the case when $p$ has length $1$, corresponding to a binary relation $\RR \le_{sd} \bA_y \times \bA_z$. If $\bA_y' + p \ne \bA_z$, then $y,z$ must both be proper with $\bA_y' + p = \bA_z'$. Either way we see that $\bA_y' + p \supseteq \bA_z'$, and since $\RR' = \RR \cap (\bA_y' \times \bA_z')$ we have $\bA_y' + p' = \bA_z'$ in the reduced instance.
\end{proof}

Finally, we check that arc-consistency of $\fX'$ and cycle-consistency of $\fX$ implies cycle-consistency of $\fX'$. For this, we note that if $p$ is a path from $v$ back to itself in $\fX$, and if $p'$ is the corresponding path in $\fX'$, then $\PP_{p'}$ is an ideal of $\PP_p$. Since $\PP_p \supseteq \Delta_{\bA_v'}$ we have
\[
m(\PP_{p'}, \Delta_{\bA_v'}, \PP_{p'}) \subseteq \PP_{p'},
\]
so the cycle-consistency of $\fX'$ follows from the following result.

\begin{thm}\label{ideal-diagonal} Suppose that $\RR \le_{sd} \bA\times \bA$ is subdirect with $m(\RR,\Delta_{\bA},\RR) \subseteq \RR$, where $m$ is a majority operation. Then $\Delta_{\bA} \subseteq \RR$.

In fact, if $\RR \le_{sd} \bA_1 \times \cdots \times \bA_n$ is subdirect and satisfies $m(\RR,S,\RR) \subseteq \RR$, where $S$ is any subset of $\bA_1 \times \cdots \times \bA_n$, then $S \subseteq \RR$.
\end{thm}
\begin{proof} First we prove the statement about binary relations, since this is all we will need. Let $a$ be any element of $\bA$. Since $\RR$ is subdirect, there are $b, c \in \bA$ such that $(a,b) \in \RR$ and $(c,a) \in \RR$. Then since $(a,a) \in \Delta_{\bA}$, we have
\[
\begin{bmatrix} a\\ a \end{bmatrix} = m\left(\begin{bmatrix} a\\ b \end{bmatrix}, \begin{bmatrix} a\\ a \end{bmatrix}, \begin{bmatrix} c\\ a \end{bmatrix}\right) \in m(\RR,\Delta_{\bA},\RR) \subseteq \RR.
\]

For the more general statement, we show by induction on $k$ that $\pi_{[k]}(S) \subseteq \pi_{[k]}(\RR)$ for each $k \le n$. The base case $k = 1$ follows from the assumption that $\RR$ is subdirect, and for the inductive step we may as well assume that we have already proven this for $k = n-1$, and wish to show it for $n$. Let $(a_1, ..., a_n)$ be any element of $S$. Then by the inductive hypothesis there is some $b$ such that $(a_1,...,a_{n-1},b) \in \RR$, and by the assumption that $\RR$ is subdirect there are $c_1, ..., c_{n-1}$ such that $(c_1, ..., c_{n-1},a_n) \in \RR$. Then we have
\[
\begin{bmatrix} a_1\\ \vdots \\ a_{n-1} \\ a_n \end{bmatrix} = m\left(\begin{bmatrix} a_1 \\ \vdots \\ a_{n-1} \\ b \end{bmatrix}, \begin{bmatrix} a_1\\ \vdots \\ a_{n-1} \\ a_n \end{bmatrix}, \begin{bmatrix} c_1\\ \vdots \\ c_{n-1} \\ a_n \end{bmatrix}\right) \in m(\RR,S,\RR) \subseteq \RR.\qedhere
\]
\end{proof}

\begin{cor} The reduced instance $\fX'$ is cycle-consistent.
\end{cor}

We have proved the main result of this section.

\begin{thm} Every cycle-consistent instance $\fX$ of a majority CSP has a solution.

In fact, for any variable $v$ of $\fX$ and any value $a \in \bA_v$, the instance $\fX$ has a solution in which the variable $v$ is assigned the value $a$.
\end{thm}
\begin{proof} For the second statement, we note that if $|\bA_v| > 1$, then $(v,\{a\}) \in \mathcal{I}$, so there is some maximal element $(x,\bA_x') \in \mathcal{I}$ such that $(v,\{a\}) \preceq (x,\bA_x')$, and we define the reduction $\fX'$ in terms of the maximal element $(x,\bA_x')$. If $v$ is proper, then from $(v,\{a\}) \preceq (x,\bA_x') \preceq (v,\bA_v')$ we must have $a \in \bA_v'$, and if $v$ is not proper then we have $a \in \bA_v = \bA_v'$. Either way, we see by induction that the reduced instance $\fX'$ has a solution in which the variable $v$ is assigned the value $a$.
\end{proof}

\begin{cor} Suppose $\bA$ is an algebra with a partial semilattice term $s$ and a ternary term $g$ such that for any subalgebra $\bB \le \bA$, the restriction of $g$ to $\Sg(\bB^{\max})$ is a majority operation. Then every cycle-consistent instance of $\CSP(\bA)$ has a solution.
\end{cor}
\begin{proof} By Corollary \ref{cycle-consistency-maximal}, if we start with a cycle-consistent instance $\fX$ and restrict all the variable domains $\bA_i$ to $\Sg(\bA_i^{\max})$ to create a new instance $\fX'$, then $\fX'$ will still be cycle-consistent, and by assumption $\fX'$ will be preserved by the majority operation $g$. Then by the previous theorem, $\fX'$ will have a solution.
\end{proof}

\begin{ex}\label{ex-slippery} Consider $\bA = (\{-,0,+\},g)$, where $g$ is the idempotent cyclic ternary operaton with
\begin{align*}
g(0,0,-) &= g(0,-,-) = -,\\
g(0,-,+) &= g(-,-,+) = -,\\
g(0,0,+) &= g(0,+,+) = +,\\
g(0,+,-) &= g(-,+,+) = +.
\end{align*}
This can be described more succinctly as follows: the permutation $(-\ +)$ is an automorphism of $\bA$, $\{-,+\}$ is a majority subalgebra of $\bA$, and $\{0,-\}, \{0,+\}$ are semilattice subalgebras of $\bA$ with $0 \rightarrow -,+$. The term $s(x,y) \coloneqq g(x,x,y)$ is a partial semilattice, and $s,g$ satisfy the assumptions of the Corollary above, so every cycle-consistent instance of $\CSP(\bA)$ has a solution. We give a table for $s$ and draw the graph of two element subalgebras of $\bA$ (with undirected edges for majority subalgebras and directed edges for semilattice subalgebras) below.
\begin{center}
\begin{tabular}{cc}
\begin{tabular}{c|ccc} $s$ & $-$ & $0$ & $+$\\ \hline $-$ & $-$ & $-$ & $-$\\ $0$ & $-$ & $0$ & $+$\\ $+$ & $+$ & $+$ & $+$\end{tabular} & 
\begin{tikzpicture}[scale=1,baseline=0.5cm]
  \node (0) at (0,0) {$0$};
  \node (-) at (-0.6,1) {$-$};
  \node (+) at (0.6,1) {$+$};
  \draw [->] (0) edge (-) (0) edge (+);
  \draw (-) -- (+);
\end{tikzpicture}
\end{tabular}
\end{center}

The relational clone $\Inv(g)$ is generated by the unary relation $x \ne 0$, the binary relation $x = -y$, the binary relation $x \le y$, and the ternary relation $x = 0 \implies y = z$.

The clone $\langle g \rangle$ is properly contained in the clone $\langle s_2 \rangle$ from Example \ref{lp-not-width-1}, and it does not contain any proper subclone with a Taylor operation. In some sense the algebra considered in this example is the prototypical example of a bounded width algebra: Bulatov \cite{bulatov-bounded} has shown that in every minimal bounded width clone, the maximal strongly connected components behave as if there is a majority operation preserving them, and for every pair of maximal strongly connected components there is a two-element majority subalgebra which connects them.
\end{ex}

\begin{rem} It's tempting to try to generalize Theorem \ref{ideal-diagonal} to near-unanimity operations. We say that a subalgebra $\bB$ \emph{absorbs} $\bA$ with respect to a near-unanimity operation $t$ if
\[
t(\bB, ..., \bB, \bA, \bB, ..., \bB) \subseteq \bB
\]
for each possible location of $\bA$. Suppose that $\RR \le_{sd} \bA \times \bA$ absorbs $\Delta_{\bA}$ with respect to $t$ - can we conclude that $\RR$ contains the diagonal?

Unfortunately the answer is no: even if $\RR$ is subdirect and absorbs $\bA^2$ with respect to a near-unanimity term, we might not have $\Delta_{\bA} \subseteq \RR$. Consider the threshold function $t_2^n$ from Example \ref{ex-non-finitely-related} defined by
\[
t_2^n(x_1, ..., x_n) = \begin{cases} 1 & \sum_i x_i \ge 2,\\ 0 & \sum_i x_i \le 1.\end{cases}
\]
For $n \ge 4$, the relation
\[
\RR = \left\{\begin{bmatrix}0 \\ 1\end{bmatrix}, \begin{bmatrix}1 \\ 0\end{bmatrix}, \begin{bmatrix}1 \\ 1\end{bmatrix}\right\}
\]
absorbs $\{0,1\}^2$ with respect to $t_2^n$, but does not contain the diagonal element $(0,0)$. However, $\RR$ \emph{does} intersect the diagonal at $(1,1)$. In the next section we will see that this weaker claim generalizes: if $\RR \le_{sd} \bA \times \bA$ absorbs $\Delta_{\bA}$, then $\RR \cap \Delta_{\bA} \ne \emptyset$.
\end{rem}

\section{Absorption, J\'onsson absorption, and connectivity}

Absorption is a common generalization of ideals of majority algebras and maximal strongly connected components of minimal ancestral algebras, and a lot of the theory of absorbing subalgebras applies to general (finite, idempotent) algebras, without assuming the existence of a Taylor term. After introducing absorption, we will show that absorbing subalgebras $\RR'$ of binary relations $\RR$ retain some of the connectivity properties of the original relations $\RR$.

\begin{defn} A subalgebra $\bB \le \bA$ \emph{absorbs} $\bA$ with respect to an idempotent term $t$ if
\[
t(\bB, ..., \bB, \bA, \bB, ..., \bB) \subseteq \bB
\]
for each possible location of $\bA$. We just say that $\bB$ absorbs $\bA$, written $\bB \lhd \bA$, if there exists some idempotent term $t$ such that $\bB$ absorbs $\bA$ with respect to $t$.

More generally, we sometimes say that a set $B$ absorbs a set $A$ with respect to an idempotent term $t$ if
\[
t(B, ..., B, A, B, ..., B) \subseteq B
\]
for each possible location of $A$. Note that if $B \subseteq A$, then $B$ must be closed under $t$.
\end{defn}

The reason we avoid specifying the idempotent term $t$ in the notation $\bB \lhd \bA$ is that there exists a common term $t$ which witnesses all absorption within any finite collection of pairs $\bB_i \lhd \bA_i$.

\begin{prop} If $\bB_1 \lhd \bA_1$ with respect to $t_1$ and $\bB_2 \lhd \bA_2$ with respect to $t_2$, then each $\bB_i \lhd \bA_i$ with respect to the star composition $t_1*t_2$ (see Definition \ref{star-comp-defn}). If $\bA_1 = \bB_2$, then $\bB_1 \lhd \bA_2$ with respect to $t_1*t_2$.
\end{prop}

\begin{cor}\label{cor-near-unanimity-absorption} A finite algebra $\bA$ has a near-unanimity term iff for all $a \in \bA$, the singleton $\{a\}$ absorbs $\bA$.
\end{cor}

A common strategy in arguments involving absorbing operations $t$ of high arity $n$ is to consider expressions of the form
\[
t(x,...,x,y,z,...,z),
\]
where just a single $y$ occurs, and iteratively march the location of the $y$ one step to the left at a time. We can make such arguments more transparent by phrasing them in terms of the sequence of ternary terms
\[
d_i(x,y,z) \coloneqq t(\underbrace{x,...,x}_{n-i},y,\underbrace{z,...,z}_{i-1}),
\]
with $d_0(x,y,z) \coloneqq x$ and $d_{n+1}(x,y,z) \coloneqq z$, so that the $d_i$ satisfy the system of identities
\begin{align*}
d_0(x,y,z) &\approx x,\\
d_i(x,y,y) &\approx d_{i+1}(x,x,y),\\
d_{n+1}(x,y,z) &\approx z.
\end{align*}
If $\bB$ absorbs $\bA$ with respect to the term $t$, then we will additionally have
\[
d_i(\bB,\bA,\bB) \subseteq \bB
\]
for all $i$.

\begin{defn} A \emph{J\'onsson absorption chain} is a sequence of ternary terms $d_1, ..., d_n$ which satisfy the identities
\begin{align*}
d_1(x,x,y) &\approx x,\\
d_i(x,y,y) &\approx d_{i+1}(x,x,y),\\
d_n(x,y,y) &\approx y.
\end{align*}
We say that $\bB$ \emph{J\'onsson absorbs} $\bA$ with respect to the J\'onsson chain $d_1, ..., d_n$ if for each $i \in [n]$ we have
\[
d_i(\bB,\bA,\bB) \subseteq \bB.
\]
If $\bB$ J\'onsson absorbs $\bA$ with respect to some J\'onsson chain, then we write $\bB \lhd_J \bA$.
\end{defn}

\begin{prop} If $\bB \lhd \bA$, then $\bB \lhd_J \bA$.
\end{prop}

As with absorption, we can witness several instances of J\'onsson absorption simultaneously with a single J\'onsson absorption chain $d_1, ..., d_n$.

\begin{prop}\label{simultaneous-jonsson-absorption} If $\bB_1 \lhd_J \bA_1$ with respect to $d_1, ..., d_m$ and $\bB_2 \lhd_J \bA_2$ with respect to $e_1, ..., e_n$, then the sequence of terms $f_1, ..., f_{mn}$ defined by
\[
f_{n(i-1)+j}(x,y,z) \coloneqq d_i(x,e_j(x,y,z),z)
\]
is a J\'onsson absorption chain which witnesses both $\bB_1 \lhd_J \bA_1$ and $\bB_2 \lhd_J \bA_2$. If $\bA_1 = \bB_2$, then $\bB_1 \lhd_J \bA_2$ with respect to $f_1, ..., f_{mn}$.
\end{prop}

\begin{cor}\label{cor-congruence-distributive-jonsson-absorption} A finite algebra $\bA$ generates a congruence distributive variety iff for all $a \in \bA$, the singleton $\{a\}$ J\'onsson absorbs $\bA$.
\end{cor}
\begin{proof} A J\'onsson absorbing chain which witnesses $\{a\} \lhd_J \bA$ for all $a \in \bA$ is the same as a sequence of terms $d_1, ..., d_m$ which satisfy the system of identities
\begin{align*}
d_1(x,x,y) &\approx x,\\
d_i(x,y,x) &\approx x\text{ for all }i,\\
d_i(x,y,y) &\approx d_{i+1}(x,x,y)\text{ for all }i,\\
d_m(x,y,y) &\approx y,
\end{align*}
that is, $d_1, ..., d_m$ are a sequence of directed J\'onsson terms. By Theorem \ref{directed-gumm-terms}, a variety is congruence distributive iff it has directed J\'onsson terms.
\end{proof}

\begin{ex} If $\bA = (A,s)$ is a 2-semilattice, then $\bB \lhd_J \bA$ iff $s(\bA,\bB) = s(\bB,\bA) \subseteq \bB$, that is, iff $\bB \lhd_{str} \bA$.
\end{ex}

\begin{ex} If $\bA$ is abelian, then $\bA$ has no J\'onsson absorbing singleton subalgebras. To see this, note that if $\bA$ is abelian, then for any J\'onsson chain $d_1, ..., d_n$ witnessing $\{b\} \lhd_J \bA$ and any $a \ne b$, we have $d_1(b,b,a) = b$, and then by induction we have
\[
d_i(b,\boxed{b},a) = b = d_i(b,\boxed{b},b) \implies d_i(b,\boxed{a},a) = d_i(b,\boxed{a},b) = b \implies d_{i+1}(b,b,a) = d_i(b,a,a) = b,
\]
so $a = d_n(b,a,a) = b$, a contradiction.

In particular, no affine algebra $\bA$ has any proper J\'onsson absorbing subalgebra $\bB$, because we can apply the above argument to the quotient $\bA/\theta_\bB$, where $\theta_\bB$ is the congruence of $\bA$ which has $\bB$ as a congruence class.
\end{ex}

\begin{ex} Suppose $\bB$ is an ideal of a majority algebra $\bA = (A,m)$. Then $\bB \lhd_J \bA$ with respect to the J\'onsson absorption chain $d_1(x,y,z) = m(x,y,z)$ (of length $1$):
\begin{align*}
m(x,x,y) &\approx x,\\
m(x,y,y) &\approx y,\\
m(\bB,\bA,\bB) &\subseteq \bB.
\end{align*}
In fact, the converse holds: if $\bB \lhd_J \bA$, then there must be a majority term $m' \in \Clo(m)$ such that $m'(\bB,\bA,\bB) \subseteq \bB$. This follows from the fact that every ternary term in a majority algebra is either a projection or another majority operation.

If $\bA$ generates a locally finite variety, then by applying the construction of Proposition \ref{simultaneous-jonsson-absorption} iteratively to all the majority operations in $\Clo(m)$, we can find a single majority term $\hat{m} \in \Clo(m)$ such that for any $\bC \le \bB \in HSP(\bA)$ we have
\[
\bC \lhd_J \bB \;\; \iff \;\; \hat{m}(\bC,\bB,\bC) \subseteq \bC.
\]

As we will see later, for finite majority algebras $\bB \lhd_J \bA$ implies that $\bB \lhd \bA$ - possibly with respect to a term of very high arity (for instance, in the case where $\bA$ is the dual discriminator algebra from Example \ref{ex-dual-discriminator} and $|\bB| = |\bA|-1$, the minimal arity of a term $t$ which witnesses $\bB \lhd \bA$ is $|\bA| + 1$). So ideals of finite majority algebras are actually the same thing as absorbing subalgebras!
\end{ex}

\begin{rem} If we define a concept called \emph{ideal absorption} by $\bB \lhd_I \bA$ when there is a ternary term $d$ such that $d(x,x,y) \approx x \approx d(y,x,x)$ and $d(\bB,\bA,\bB) \subseteq \bB$, then all of the results about ideals of majority algebras generalize. I don't know any applications of this idea outside the context of majority algebras.
\end{rem}

Like ideals of majority algebras, absorbing subalgebras play nice with primitive positive formulas.

\begin{prop} Suppose that a relation $\RR$ is defined by a primitive positive formula $\Phi$ involving the relations $\RR_1, ..., \RR_k$. If we replace each $\RR_i$ with an absorbing subalgebra $\RR_i' \lhd \RR_i$ to make a primitive positive formula $\Phi'$, then the relation $\RR'$ which is defined by $\Phi'$ is an absorbing subalgebra of $\RR$. The same is true with ``absorbing'' replaced by ``J\'onsson absorbing''.
\end{prop}
\begin{proof} Since only finitely many relations $\RR_i$ show up in $\Phi$, we can find a single absorbing term (or J\'onsson chain) which witnesses all absorptions $\RR_i' \lhd \RR_i$ (or $\RR_i' \lhd_J \RR_i$) simultaneously. From here the proof is similar to the proof in the case of ideals of majority algebras.
\end{proof}

Now we will illustrate how J\'onsson absorption is used, by proving a few connectivity results. Recall that every binary relation $\RR \le \bA\times \bA$ can be visualized as a graph in two different ways: we can either think of $\RR$ as a bipartite graph on the disjoint union $\bA \sqcup \bA$, or we can think of $\RR$ as a directed graph on $\bA$. The next result is perhaps the most crucial.

\begin{thm}[Absorbing directed paths \cite{deciding-absorption}]\label{absorbing-directed-path} If $\bS, \RR \le \bA \times \bA$ are binary relations with $\bS \lhd_J \RR$, and $a,b \in \bA$ satisfy
\begin{itemize}
\item $(a,a), (b,b) \in \bS$, and

\item $(a,b) \in \RR$,
\end{itemize}
then if we think of $\bS$ as a directed graph on $\bA$, there is a directed path from $a$ to $b$ in $\bS$, that is, $(a,b) \in \bS^{\circ n}$ for some $n$.
\end{thm}
\begin{proof} Suppose $\bS \lhd_J \RR$ with respect to the J\'onsson chain $d_1, ..., d_n$. Then for each $i$ we have
\[
\begin{bmatrix} d_i(a,a,b)\\ d_{i+1}(a,a,b) \end{bmatrix} = \begin{bmatrix} d_i(a,a,b)\\ d_i(a,b,b) \end{bmatrix} = d_i\left(\begin{bmatrix} a\\ a \end{bmatrix}, \begin{bmatrix} a\\ b \end{bmatrix}, \begin{bmatrix} b\\ b \end{bmatrix}\right) \in d_i(\bS, \RR, \bS) \subseteq \bS.
\]
Stringing these together, we get a directed path from $d_1(a,a,b) = a$ to $d_n(a,b,b) = b$ of length $n$, so in fact $(a,b) \in \bS^{\circ n}$.
\end{proof}

Applying the above to $\bS^{\circ m} \lhd_J \RR^{\circ m}$ for a sufficiently large $m$, we get the following stronger-looking corollary.

\begin{cor}\label{cor-absorbing-directed-path} If $\bS, \RR \le \bA \times \bA$ have $\bS \lhd_J \RR$, and $a,b \in \bA$ satisfy
\begin{itemize}
\item each of $a,b$ is contained in a directed cycle of the digraph $\bS$, and

\item there is a directed path from $a$ to $b$ in the digraph $\RR$,
\end{itemize}
then there is a directed path from $a$ to $b$ in the digraph $\bS$.
\end{cor}

For the sake of applying the previous result, it is useful to keep in mind the following basic fact about finite directed graphs.

\begin{prop}\label{prop-directed-cycle} If $(A,R)$ is a finite directed graph such that each vertex of $A$ has in-degree at least $1$ (in other words, such that $\pi_2(R) = A$), then for every vertex $a \in A$ there is some $a' \in A$ and some $n$ such that $a'$ is contained in a directed cycle of length $n$ and such that there is a directed path from $a'$ to $a$ of length $n$ (that is, $(a',a'), (a',a) \in R^{\circ n}$).
\end{prop}
\begin{proof} Define a function $\varphi : A \rightarrow A$ such that for each $a \in A$ we have $(\varphi(a),a) \in R$. Then there is some $n$ such that $\varphi^{\circ 2n} = \varphi^{\circ n}$ by the finiteness of $A$: in fact, we may take $n = \lcm\{1, ..., |A|\}$.
\end{proof}

In the next result, we think of binary relations as bipartite graphs. Recall that the \emph{linked components} of a binary relation $\RR \le \bA \times \bB$ are the connected components of $\RR$ considered as a bipartite graph on $\bA \sqcup \bB$, and that the linked components of size greater than $1$ are the same as the congruence classes of the linking congruence $\ker \pi_1 \vee \ker \pi_2$ on $\RR$.

\begin{thm}[Absorbing linked components \cite{deciding-absorption}]\label{absorbing-linked} If $\bS, \RR \le \bA \times \bB$ are binary relations with $\bS \lhd_J \RR$, and $a,b \in \pi_1(\bS)$ are in the same linked component of $\RR$, then $a,b$ are in the same linked component of $\bS$.
\end{thm}
\begin{proof} If $a,b$ are linked in $\RR$, then there is some $m$ such that $(a,b) \in (\RR \circ \RR^{-})^{\circ m}$. Since $(\bS \circ \bS^{-})^{\circ m} \lhd_J (\RR \circ \RR^{-})^{\circ m}$ and $(a,a), (b,b) \in \bS \circ \bS^-$ by $a,b \in \pi_1(\bS)$, we can apply Theorem \ref{absorbing-directed-path} to see that there is some $n$ such that $(a,b) \in (\bS \circ \bS^{-})^{\circ mn}$. Thus $a,b$ are in the same linked component of $\bS$.
\end{proof}

The next result is an analogue of Theorem \ref{strong-diagonal} and Theorem \ref{ideal-diagonal} for J\'onsson absorption.

\begin{thm}[Loop Lemma, finite absorbing case \cite{deciding-absorption-relational}]\label{absorbing-diagonal} If $\RR \le_{sd} \bA\times \bA$ is subdirect, $\bA$ is finite, and $\RR$ J\'onsson absorbs the diagonal $\Delta_\bA$, then $\RR \cap \Delta_\bA \ne \emptyset$.
\end{thm}
\begin{proof} We may assume without loss of generality that $\bA$ is idempotent. As long as $|\bA| > 1$, we will try to find a proper subalgebra $\bB \le \bA$ with $\RR \cap (\bB \times \bB)$ subdirect. Then $\RR\cap (\bB\times \bB)$ will J\'onsson absorb $\Delta_{\bB}$, and we can show by induction that $\RR \cap \Delta_{\bB} \ne \emptyset$.

Let $b$ be any element of $\bA$, and define a sequence of subalgebras $\bB_i$ by $\bB_0 = \{b\}$,
\[
\bB_{i+1} = \bB_i + \RR,
\]
i.e. $\bB_{i+1} = \pi_2(\RR \cap (\bB_i \times \bA))$. If there is any $i$ such that $\bB_i \ne \bA$ but $\bB_{i+1} = \bA$, then for every $\bC \le \bA$ we have $(\bC + \RR^-) \cap \bB_i \ne \emptyset$, so by the finiteness of $\bB_i$ we may take
\[
\bB = \bigcap_{k \ge 0} \bB_i+(\RR^-)^{\circ k} = \{a \mid \exists a_0, a_1, ... \in \bB_i \text{ s.t. } a_0 = a \text{ and } \forall j \; (a_j,a_{j+1}) \in \RR\}.
\]
Otherwise, each $\bB_i \ne \bA$, and by the finiteness of $\bA$ there must be some $m,n$ such that $\bB_m = \bB_{m+n}$. We will show that in this case we have $\bB_{m+i} = \bB_m$ for each $i$, so we may take $\bB = \bB_m$.

Consider any directed cycle $a_0, ..., a_{kn} = a_0$ of $\RR$ (considered as a digraph) with $a_0 \in \bB_m$. We will show that each $a_i \in \bB_m$. Note that
\[
(a_0,a_0), (a_i,a_i) \in \RR^{\circ kn},
\]
that $\RR^{\circ kn}$ J\'onsson absorbs $(\RR\cup\Delta_{\bA})^{\circ kn}$, and that
\[
(a_0,a_i) \in \RR^{\circ i} \subseteq (\RR\cup\Delta_{\bA})^{\circ kn}.
\]
Thus by Theorem \ref{absorbing-directed-path} there is some $l$ such that $(a_0,a_i) \in \RR^{\circ ln}$, and since $\bB_m + \RR^{\circ ln} = \bB_{m+ln} = \bB_m$, we see that $a_i \in \bB_m$.

Since $\bB_{m+i} + \RR^{\circ n} = \bB_{m+i}$ and $\bB_{m+i}$ is finite, for each element $a$ of $\bB_{m+i}$ there is an $a_i$ contained in a directed cycle of $\RR^{\circ n}$ and a directed path of $\RR^{\circ n}$ from $a_i$ to $a$, so in fact we have $a \in \bB_m$ as well, and we see that $\bB_{m+i} \subseteq \bB_m$. Similarly we have $\bB_m \subseteq \bB_{m+i}$, so $\bB_m = \bB_{m+i}$.
\end{proof}

\begin{cor}\label{cor-absorbing-diagonal} If $\RR \le_{sd} \bA\times \bA$ is subdirect, $\bA$ is finite and has no proper absorbing subalgebra, and $\RR$ absorbs the diagonal $\Delta_\bA$, then $\Delta_\bA \subseteq \RR$.
\end{cor}
\begin{proof} Since $\RR \cap \Delta_{\bA} \ne \emptyset$ is an absorbing subalgebra of $\Delta_{\bA}$ and $\Delta_{\bA} \cong \bA$ has no proper absorbing subalgebra, we must have $\RR \cap \Delta_\bA = \Delta_\bA$.
\end{proof}

\begin{defn} We say that $\bB$ is a \emph{minimal absorbing subalgebra} of $\bA$, written $\bB \llhd \bA$, if $\bB \lhd \bA$ and $\bB$ has no proper absorbing subalgebra.
\end{defn}

\begin{prop} Every finite absorbing subalgebra of $\bA$ contains a minimal absorbing subalgebra of $\bA$, and any pair of distinct minimal absorbing subalgebras of $\bA$ are disjoint.
\end{prop}
\begin{proof} This follows from the fact that $\bC \lhd \bB \lhd \bA$ implies $\bC \lhd \bA$, and the fact that the intersection of any pair of absorbing subalgebras is an absorbing subalgebra.
\end{proof}

\begin{thm}\label{arc-consistent-cycles} Suppose that $\fX$ is an arc-consistent instance of a CSP, and suppose that for each variable domain $\bA_v$ there is a minimal absorbing subalgebra $\bA_v' \llhd \bA_v$ such that the reduced instance $\fX'$ with variable domains replaced by $\bA_v'$ and relations $\RR \le_{sd} \bA_{v_1} \times \cdots \times \bA_{v_n}$ replaced by $\RR' = \RR \cap (\bA_{v_1}' \times \cdots \times \bA_{v_n}')$ is arc-consistent.

Then for any path $p$ in $\fX$ from a variable $v$ to itself such that $\PP_p \supseteq \Delta_{\bA_v}$, the corresponding path $p'$ of $\fX'$ has $\PP_{p'} \supseteq \Delta_{\bA_v'}$. In particular, if $\fX$ is cycle-consistent then so is $\fX'$.
\end{thm}
\begin{proof} Note that $\PP_{p'}$ is an absorbing subalgebra of $\PP_p$, so $\PP_{p'}$ absorbs $\Delta_{\bA_v'}$. Since $\PP_{p'}$ is subdirect in $\bA_v' \times \bA_v'$ by the arc-consistency of $\fX'$ and $\bA_v'$ has no proper absorbing subalgebra, we may apply Corollary \ref{cor-absorbing-diagonal} to see that $\PP_{p'} \supseteq \Delta_{\bA_v'}$.
\end{proof}

Later we will show that any cycle-consistent instance $\fX$ has an arc-consistent reduction $\fX'$ where all variable domains are replaced by minimal absorbing subalgebras, which will set us up to apply Theorem \ref{arc-consistent-cycles}. The argument strategy will be fairly generic, not using any specific properties of absorbing subalgebras other than Theorem \ref{absorbing-directed-path} and the fact that absorption is compatible with primitive positive formulas. Additionally, we will be able to weaken cycle-consistency to a property known as $pq$-consistency, which says that for any pair of paths $p,q$ from a variable $v$ to itself, there is some $j \ge 0$ such that $\PP_{j(p+q)+p} \supseteq \Delta_{\bA_v}$.

\subsection{Local criterion for J\'onsson absorption}

Since a finite algebra $\bA$ has bounded strict width iff every singleton is an absorbing subalgebra of $\bA$, we'd like to have a way to test whether a given subalgebra $\bB \le \bA$ is an absorbing subalgebra. Since the arity of a potential absorbing term is unbounded, we'll start with the easier problem of testing whether $\bB$ is a \emph{J\'onsson} absorbing subalgebra, since in this case there is an obvious algorithm which will at least eventually halt: list out every possible ternary term operation of $\bA$ by brute force, and make a digraph of possible J\'onsson chains.

The idea behind finding a better way to test whether $\bB \lhd_J \bA$ is to try to find a converse to the fundamental digraph connectivity result characterizing J\'onsson absorption (Theorem \ref{absorbing-directed-path}). In order to formulate the converse, we need to consider generic pairs of digraphs $\bS \le \RR \le \bC\times \bC$ such that
\[
\bB \lhd_J \bA \implies \bS \lhd_J \RR.
\]
One natural way to do this is to write $\RR$ as the projection to the last two coordinates of a ternary relation $\bX \le \bA\times\bC\times\bC$, and to take $\bS$ to be the corresponding projection of $\bX \cap (\bB\times\bC\times\bC)$.

\begin{defn} For $\bB \le \bA$ and $\bC \in \cV(\bA)$ all idempotent, we say that $\bA,\bB,\bC$ satisfy the condition $J(\bA,\bB;\bC)$ if for every $a \in \bA,\ b, b' \in \bB$, and $c,d \in \bC$, if we set
\[
\bS = \pi_{23}\left(\Sg_{\bA\times\bC\times\bC}\left\{\begin{bmatrix}b\\ c\\ c\end{bmatrix},\begin{bmatrix}a\\ c\\ d\end{bmatrix}, \begin{bmatrix}b'\\ d\\ d\end{bmatrix}\right\} \cap \begin{bmatrix} \bB\\ \bC\\ \bC \end{bmatrix}\right),
\]
then there is some $n$ such that $(c,d) \in \bS^{\circ n}$.
\end{defn}

We will show that the condition $J(\bA,\bB;\bA)$ is equivalent to $\bB \lhd_J \bA$, following the strategy of \cite{deciding-absorption}. Note that Theorem \ref{absorbing-directed-path} proves one direction of the equivalence, so we just need to prove that $J(\bA,\bB;\bA) \implies \bB \lhd_J \bA$. The strategy will be to use induction to show that $J(\bA^m,\bB^m;\bA^n)$ holds for all $m,n$, and then to take $m = |\bA||\bB|^2, n = |\bA|^2$ to show that a certain directed path exists between binary terms in the free algebra on two generators, which will correspond to a J\'onsson absorption chain. Before diving into the details, we will outline how this criterion could be used to test whether $\bB \lhd_J \bA$.

Note that if $\bA$ is given in terms of tables for its basic operations, then the condition $J(\bA,\bB;\bA)$ can be tested in time polynomial in $|\bA|, |\bB|$ (with the degree of the polynomial depending on the arities of the basic operations), since the total number of tuples $a,b,b',c,d$ is $|\bA|^3|\bB|^2$, computing $\bS$ requires us to compute a ternary relation of size at most $|\bA|^3$, and we only need to check whether $(c,d) \in \bS^{\circ n}$ for $n \le |\bA|$.

If $\bA$ is instead given in terms of a list of basic relations, then testing the condition $J(\bA,\bB;\bA)$ can be reduced to solving polynomially many polynomially large constraint satisfaction problems over the domain $\bA$ - so in particular if $\CSP(\bA)$ can be solved in polynomial time, then we can test $J(\bA,\bB;\bA)$ in polynomial time. To see this, note that in order to test whether a given edge $(e,f)$ is an element of $\bS$, we just need to test whether $\bA$ has a ternary polymorphism $t$ such that
\begin{align*}
t(b,a,b') &\in \bB,\\
t(c,c,d) &= e,\\
t(c,d,d) &= f,
\end{align*}
and the set of ternary polymorphisms $t \in \cF_{\bA}(x,y,z) \le \bA^{\bA^3}$ can be described by a primitive positive formula involving only $|\bA|^3$ variables.

\begin{lem} If $J(\bA_1,\bB_1;\bC)$ and $J(\bA_2,\bB_2;\bC)$ both hold, then so does $J(\bA_1\times \bA_2,\bB_1\times\bB_2;\bC)$.
\end{lem}
\begin{proof} Suppose $a = (a_1,a_2) \in \bA_1 \times \bA_2$ and $b = (b_1, b_2), b' = (b_1', b_2') \in \bB_1\times \bB_2$, $c,d \in \bC$. Define $\bS, \RR \le \bC \times \bC$ as usual, and define an intermediate digraph $\bS_1$, where instead of restricting to $\bB_1 \times \bB_2$, we restrict to $\bB_1 \times \bA_2$ instead - so for the purposes of computing $\bS_1$, we can ignore the $\bA_2$ components. Then by $J(\bA_1,\bB_1;\bC)$, from $(c,d) \in \RR$ we see that there is a directed path from $c$ to $d$ in $\bS_1$.

To finish, we just need to check that for each $(e,f) \in \bS_1$, there is a directed path from $e$ to $f$ in $\bS$. Note that $(e,f) \in \bS_1$ means that there are some $b_1'' \in \bB_1, a_2'' \in \bA_2$ such that
\[
\begin{bmatrix} (b_1'', a_2'')\\ e\\ f \end{bmatrix} \in \Sg\left\{\begin{bmatrix}(b_1,b_2)\\ c\\ c\end{bmatrix},\begin{bmatrix}(a_1,a_2)\\ c\\ d\end{bmatrix}, \begin{bmatrix}(b_1',b_2')\\ d\\ d\end{bmatrix}\right\}.
\]
Then from $e,f \in \Sg\{c,d\}$, there are some $(b_1''',b_2''') \in \bB_1 \times \bB_2$ with
\[
\begin{bmatrix} (b_1''', b_2''')\\ e\\ e \end{bmatrix} \in \Sg\left\{\begin{bmatrix}(b_1,b_2)\\ c\\ c\end{bmatrix}, \begin{bmatrix}(b_1',b_2')\\ d\\ d\end{bmatrix}\right\},
\]
and similarly for $(f,f)$, so we just need to check that
\[
\pi_{23}\left(\Sg\left\{\begin{bmatrix}(b_1''',b_2''')\\ e\\ e\end{bmatrix},\begin{bmatrix}(b_1'',a_2'')\\ e\\ f\end{bmatrix}, \begin{bmatrix}(b_1'''',b_2'''')\\ f\\ f\end{bmatrix}\right\} \cap \begin{bmatrix}\bB_1\times \bB_2\\ \bC\\ \bC \end{bmatrix}\right)
\]
contains a directed path from $e$ to $f$. But now we can ignore the $\bB_1$ component, so this follows from $J(\bA_2,\bB_2;\bC)$.
\end{proof}

\begin{lem} If $J(\bA,\bB;\bC_1)$ and $J(\bA,\bB;\bC_2)$ both hold and $\bC_1, \bC_2$ are finite and idempotent, then $J(\bA,\bB;\bC_1\times\bC_2)$ holds as well.
\end{lem}
\begin{proof} Suppose not. Choose $c = (c_1,c_2), d = (d_1,d_2) \in \bC_1 \times \bC_2$ such that $\Sg\{c,d\}$ is minimal among all pairs such that there exist $a \in \bA, b, b' \in \bB$ so that the associated digraph $\bS$ has no directed path from $c$ to $d$.

Ignoring the $\bC_2$ components, we can apply $J(\bA,\bB;\bC_1)$ to find a sequence of edges $(e^i,f^{i+1}) \in \bS$ such that $f^i_1 = e^i_1$ for each $i \le n$, $c = f^1$, and $e^n = d$. Since we assumed that there is no directed path from $c$ to $e^n = d$, we can consider the first $i$ such that there is no directed path from $c$ to $e^i$.

Since $e^i \in \Sg\{c,d\}$, we have
\[
\begin{bmatrix}c\\ e^i\end{bmatrix} \in \Sg\left\{\begin{bmatrix}c\\ c\end{bmatrix},\begin{bmatrix}c\\ d\end{bmatrix},\begin{bmatrix}d\\ d\end{bmatrix}\right\} = \RR,
\]
and since there is no directed path from $c$ to $e^i$ in $\bS$, we see that we must have $\Sg\{c,e^i\} = \Sg\{c,d\}$ by our minimality assumption, so in particular we have $f^i \in \Sg\{c,e^i\}$. Thus we have
\[
\begin{bmatrix}f^i\\ e^i\end{bmatrix} \in \Sg\left\{\begin{bmatrix}c\\ e^i\end{bmatrix},\begin{bmatrix}e^i\\ e^i\end{bmatrix}\right\} \subseteq \RR.
\]
By the choice of $i$ there is a path from $c$ to $f^i$ in $\bS$ (passing through $e^{i-1}$ if $i > 1$). To get a contradiction, we just need to show that there is a directed path from $f^i$ to $e^i$ in $\bS$. Since $(e^i,e^i), (f^i,f^i) \in \bS$, there are $a' \in \bA, b'',b''' \in \bB$ such that
\[
\pi_{23}\left(\Sg\left\{\begin{bmatrix}b''\\ f^i\\ f^i\end{bmatrix}, \begin{bmatrix} a'\\ f^i\\ e^i\end{bmatrix}, \begin{bmatrix}b'''\\ e^i\\ e^i\end{bmatrix}\right\}\cap \begin{bmatrix}\bB\\ \bC_1\times\bC_2\\ \bC_1\times\bC_2\end{bmatrix}\right) \subseteq \bS.
\]
Since $f^i_1 = e^i_1$, we can ignore the $\bC_1$ components in the above, so by $J(\bA,\bB;\bC_2)$ there is a directed path from $f^i$ to $e^i$ in $\bS$.
\end{proof}

\begin{rem} The same argument can be used to show that if $\bC$ has a congruence $\theta$ such that $J(\bA,\bB;\bC/\theta)$ holds and such that for each $c \in \bC$ the condition $J(\bA,\bB;c/\theta)$ holds (here $c/\theta$ is treated as a subalgebra of $\bC$), then $J(\bA,\bB;\bC)$ holds.
\end{rem}

\begin{thm}[Local criterion for J\'onsson absorption \cite{deciding-absorption}]\label{local-jonsson} If $\bB \le \bA$ are finite and idempotent, then $\bB \lhd_J \bA$ if and only if $J(\bA,\bB;\bA)$ holds.
\end{thm}
\begin{proof} By the previous two lemmas, $J(\bA^m,\bB^m;\bA^n)$ holds for $m = \bB\times\bA\times\bB$ and $n = \bA\times\bA$. There is a natural map $\Phi: \cF_\bA(x,y,z) \rightarrow \bA^{\bB\times\bA\times\bB}$ and a pair of natural maps $\Psi_1, \Psi_2: \cF_\bA(x,y,z) \rightarrow \bA^{\bA\times\bA}$: the first takes $f$ to the restriction $f|_{\bB\times\bA\times\bB}$, the other two take $f$ to the functions $f(x,x,y)$, $f(x,y,y)$.

Then we can apply $J(\bA^m,\bB^m;\bA^n)$ with $a = \Phi(\pi_2)$, $b = \Phi(\pi_2), b' = \Phi(\pi_3)$, $c = \Psi_i(\pi_1) = \Psi_1(\pi_2), d = \Psi_i(\pi_3) = \Psi_2(\pi_2)$. If we set
\[
\bS = \pi_{23}\left(\Sg\left\{\begin{bmatrix}x|_{x,z\in\bB}\\ x\\ x\end{bmatrix}, \begin{bmatrix} y|_{x,z\in\bB}\\ x\\ y\end{bmatrix}, \begin{bmatrix}z|_{x,z\in\bB}\\ y\\ y\end{bmatrix}\right\}\cap \begin{bmatrix}\bB^m\\ \bA^n\\ \bA^n\end{bmatrix}\right),
\]
then the inner ternary subalgebra is exactly $\operatorname{Im}(\Phi,\Psi_1,\Psi_2)$, so $\bS$ is exactly the digraph of pairs of binary terms $g(x,y),h(x,y)$ such that there is some ternary term $f(x,y,z)$ satisfying
\begin{align*}
f(\bB,\bA,\bB) &\subseteq \bB,\\
f(x,x,y) &\approx g(x,y),\\
f(x,y,y) &\approx h(x,y).
\end{align*}
The condition $J(\bA^m,\bB^m;\bA^n)$ says that this digraph contains a path from the term $x$ to the term $y$, which is the same as a J\'onsson absorption chain for $\bB \lhd_J \bA$.
\end{proof}

Note that the same argument shows that it is enough to check $J(\bA,\bB;\bC_i)$ for any collection of algebras $\bC_1, ..., \bC_n$ generating a variety $\cV$ such that $\cF_\bA(x,y) = \cF_{\cV}(x,y)$. In cases where $\cF_\bA(x,y)$ is small the criterion becomes especially nice.


\begin{cor} If $\bA = (A,m)$ is a majority algebra, then $\bB \lhd_J \bA$ iff there do \emph{not} exist $a\in \bA$ and $b,c \in \bB$ such that
\begin{itemize}
\item $a,b,c$ are distinct,
\item $\Sg_\bA\{a,b,c\} \cap \bB = \{b,c\}$,
\item the partitions $\{\{b\},\Sg\{a,b,c\}\setminus \{b\}\}$ and $\{\{c\},\Sg\{a,b,c\}\setminus \{c\}\}$ of $\Sg\{a,b,c\}$ correspond to congruences $\theta_b, \theta_c$ on $\Sg\{a,b,c\}$.
\end{itemize}
The third bullet point can also be stated in the equivalent form: $\Sg\{a,b,c\}/(\theta_b\wedge\theta_c)$ is isomorphic to the three element median algebra, with median element $a/(\theta_b \wedge \theta_c) = \Sg\{a,b,c\}\setminus\{b,c\}$.
\end{cor}


\section{Absorption and $\bB$-essential relations}

In this section we'll give a relational description of absorption, as well as a first simplification via Ramsey theory. The relational description is a generalization of the way relations over near-unanimity algebras decompose.

\begin{defn} Suppose $\bB \le \bA$. We say that a relation $\RR \le \bA^m$ is $\bB$-\emph{essential} if for every $1 \le i \le n$ we have
\[
\RR \cap (\bB^{i-1} \times \bA \times \bB^{n-i}) \ne \emptyset,
\]
but
\[
\RR \cap \bB^n = \emptyset.
\]

More generally, if $\bB_i \le \bA_i$ for all $i$, then we say that $\RR \le \bA_1 \times \cdots \times \bA_m$ is $(\bB_1, ..., \bB_m)$-\emph{essential} if
\[
\RR \cap (\bB_1 \times \cdots \times \bB_{i-1} \times \bA_i \times \bB_{i+1} \times \cdots \times \bB_m) \ne \emptyset
\]
for each $i$, but
\[
\RR \cap (\bB_1 \times \cdots \times \bB_m) = \emptyset.
\]
\end{defn}

\begin{prop}\label{prop-essential-down} If $\RR \le \bA^m$ is $\bB$-essential, then so is
\[
\pi_{[m-1]}(\RR \cap (\bA^{m-1}\times\bB)).
\]
In particular, if there is a $\bB$-essential relation of some arity, then there are $\bB$-essential relations of all smaller arities.
\end{prop}

\begin{prop}\label{prop-absorption-essential} If $\bB$ absorbs $\bA$ with respect to a term $t$ of arity $m$, then there are no $\bB$-essential relations $\RR \le \bA^m$ of arity $m$.
\end{prop}
\begin{proof} Suppose for contradiction that $\RR \le \bA^m$ is $\bB$-essential, and let $b_{ij} \in \bB, a_i \in \bA$ be such that
\[
\begin{bmatrix}a_1\\ b_{21}\\ \vdots\\ b_{m1}\end{bmatrix}, \begin{bmatrix}b_{12}\\ a_2\\ \vdots\\ b_{m2}\end{bmatrix}, ..., \begin{bmatrix}b_{1m}\\ b_{2m}\\ \vdots\\ a_m\end{bmatrix} \in \RR.
\]
Then if we apply $t$, we have
\[
t\left(\begin{bmatrix}a_1 & b_{12} & \cdots & b_{1m}\\ b_{21} & a_2 & \cdots & b_{2m}\\ \vdots & \vdots & \ddots & \vdots\\ b_{m1} & b_{m2} & \cdots & a_m\end{bmatrix}\right) \in \RR \cap \bB^m
\]
since $\bB$ absorbs $\bA$ with respect to $t$, which is a contradiction.
\end{proof}

\begin{cor}\label{cor-absorption-essential} If $\bB$ absorbs $\bA$ with respect to a term $t$ of arity $m$, then for any $n \ge m-1$ and any relation $\RR \le \bA^n$ such that
\[
\pi_I(\RR) \cap \bB^{m-1} \ne \emptyset
\]
for all $I \subseteq [n]$ with $|I| = m-1$, we have
\[
\RR \cap \bB^n \ne \emptyset.
\]
\end{cor}
\begin{proof} We prove this by induction on $n \ge m-1$. The base case $n = m-1$ follows by taking $I = [n]$. For the inductive step, note that by the inductive hypothesis we have
\[
\pi_{[n]\setminus\{i\}}(\RR) \cap \bB^{n-1} \ne \emptyset
\]
for all $i \in [n]$, and we must have $\RR \cap \bB^n \ne \emptyset$ since there are no $\bB$-essential relations of arity $n \ge m$ by Propositions \ref{prop-essential-down} and \ref{prop-absorption-essential}.
\end{proof}

Our main result is the converse to Proposition \ref{prop-absorption-essential}.

\begin{thm}[Relational description of absorption \cite{deciding-absorption}]\label{absorption-essential} If $\bA$ is finite and idempotent, then $\bB$ absorbs $\bA$ with respect to a term of arity $m$ if and only if there are no $\bB$-essential relations of arity $m$. In particular, we have $\bB \lhd \bA$ if and only if there is a bound on the arity of $\bB$-essential relations.
\end{thm}

The strategy of the proof is to show that if there are no $m$-ary terms $t$ which absorb $\bB$, then the projection of the free algebra $\cF_{\bA}(x_1, ..., x_m) \le \bA^{\bA^m}$ onto the coordinates where all but one input $x_i$ are in $\bB$ looks like a $\bB$-essential relation. The arity of this projection will be much higher than $m$, but the set of coordinates can be naturally grouped into $m$ parts.

\begin{lem}\label{essential-groups} If $n_1, ..., n_m \ge 1$ and $\RR \le \bA^{n_1}\times \cdots \times \bA^{n_m}$ is $(\bB^{n_1}, ..., \bB^{n_m})$-essential, then there is a $\bB$-essential relation $\RR' \le \bA^m$ of arity $m$. In fact, $\RR'$ can be chosen to have the form
\[
\RR' = \pi_I\Big(\RR \cap \Big(\prod_i \bC_i\Big)\Big)
\]
for some $I \subseteq [n_1 + \cdots + n_m]$ with $|I| = m$ such that
\[
\big|I \cap (n_1 + \cdots + n_{j-1}, n_1 + \cdots + n_j]\big| = 1
\]
for each $j \le m$, and for some choice of $\bC_i \in \{\bA, \bB\}$ for each $i$.
\end{lem}
\begin{proof} We prove this by induction on $n = n_1 + \cdots + n_m$. If all $n_i = 1$, then $\RR$ is an $m$-ary $\bB$-essential relation already. Otherwise, we may assume $n_m > 1$ without loss of generality. First consider the relation
\[
\RR_1 = \pi_{[n-1]}(\RR \cap (\bA^{n-1}\times \bB)) \le \bA^{n_1}\times \cdots \times \bA^{n_m-1}.
\]
We have
\[
\RR_1 \cap (\bB^{n_1}\times \cdots \times \bA^{n_i} \times \cdots \times \bB^{n_{m-1}} \times \bB^{n_m-1}) \ne \emptyset
\]
for each $i \ne m$, and
\[
\RR_1 \cap \bB^{n-1} = \emptyset,
\]
so the only way for $\RR_1$ to fail to be $(\bB^{n_1}, ..., \bB^{n_{m-1}}, \bB^{n_m-1})$-essential is if
\[
\pi_{[n-n_m]\cup\{n\}}(\RR) \cap \bB^{n-n_m+1} = \emptyset.
\]
In this case, we see that
\[
\RR_2 = \pi_{[n-n_m]\cup\{n\}}(\RR)
\]
is a $(\bB^{n_1}, ..., \bB^{n_{m-1}}, \bB)$-essential relation.
\end{proof}

\begin{proof}[Proof of Theorem \ref{absorption-essential}] We just need to prove that if there is no $m$-ary $\bB$-essential relation, then $\bB$ absorbs $\bA$ with respect to some $m$-ary term $t$. For each $i$, let $X_i$ be the set of tuples $(x_1, ..., x_m) \in \bA^m$ such that $x_j \in \bB$ for $j \ne i$, and $x_i \in \bA\setminus \bB$. Consider the relation
\[
\RR = \pi_{X_1 \cup \cdots \cup X_m} (\cF_{\bA}(x_1, ..., x_m)) \le \bA^{X_1} \times \cdots \times \bA^{X_m}.
\]
Since $\cF_{\bA}(x_1, ..., x_m)$ contains the projection functions $\pi_i : \bA^m \rightarrow \bA$, by the definition of the sets $X_i$ we have
\[
\RR \cap (\bB^{X_1} \times \cdots \times \bA^{X_i} \times \cdots \times \bB^{X_m}) \ne \emptyset
\]
for all $i$. Since there is no $\bB$-essential relation of arity $m$, we see that $\RR$ can't be $(\bB^{X_1}, ..., \bB^{X_m})$-essential by Lemma \ref{essential-groups}, so we must have
\[
\RR \cap (\bB^{X_1} \times \cdots \times \bB^{X_m}) \ne \emptyset
\]
as well. Then by the definition of $\RR$, we see that there is a term $t \in \cF_{\bA}(x_1, ..., x_m)$ which absorbs $\bB$.
\end{proof}

We can simplify this slightly as follows.

\begin{cor} We have $\bB \lhd \bA$ with respect to an $m$-ary term $t$ iff for all $b_{ij} \in \bB, a_i \in \bA$ we have
\[
\Sg_{\bA^m}\left\{\begin{bmatrix}a_1\\ b_{21}\\ \vdots\\ b_{m1}\end{bmatrix}, \begin{bmatrix}b_{12}\\ a_2\\ \vdots\\ b_{m2}\end{bmatrix}, ..., \begin{bmatrix}b_{1m}\\ b_{2m}\\ \vdots\\ a_m\end{bmatrix}\right\} \cap \bB^m \ne \emptyset.
\]
\end{cor}

We leave the following generalization as an exercise to the reader.

\begin{thm}\label{absorption-essential-gen} If $\bA_1, ..., \bA_k$ are finite and idempotent, and $\bB_i \le \bA_i$ for each $i$ are such that there is no $(\bB_{i_1}, ..., \bB_{i_m})$-essential relation $\RR \le \bA_{i_1} \times \cdots \times \bA_{i_m}$ for any choice of $i_1, ..., i_m \in [k]$, then there is an $m$-ary term $t$ such that each $\bB_i$ absorbs $\bA_i$ with respect to $t$.
\end{thm}

\begin{cor} A finite idempotent algebra $\bA$ has a near-unanimity term operation of arity $m$ iff for each choice of $a_i, b_i \in \bA$, we have
\[
\begin{bmatrix}b_1\\ b_2\\ \vdots\\ b_m\end{bmatrix} \in \Sg_{\bA^m}\left\{\begin{bmatrix}a_1\\ b_2\\ \vdots\\ b_m\end{bmatrix}, \begin{bmatrix}b_1\\ a_2\\ \vdots\\ b_m\end{bmatrix}, ..., \begin{bmatrix}b_1\\ b_2\\ \vdots\\ a_m\end{bmatrix}\right\}.
\]
\end{cor}

Now we move our focus to finding a simpler characterization of $\bB \lhd \bA$, without restricting to terms of a particular arity. We'll use the notation $r_k(m)$ for the multicolored Ramsey number $R(m,...,m)$ (with $k$ copies of $m$), defined as the least number $n$ such that any edge coloring of $K_n$ with $k$ colors must have a monochromatic copy of $K_m$.

\begin{thm} If $\bA$ is finite and idempotent, then $\bB \lhd \bA$ iff there do \emph{not} exist $a \in \bA$ and $b,c \in \bB$ such that for every $m$, we have
\[
\Sg_{\bA^m}\left\{\begin{bmatrix} a & b & b & \cdots & b\\ c & a & b & \cdots & b\\ c & c & a & \cdots & b\\ \vdots & \vdots & \vdots & \ddots & \vdots\\ c & c & c & \cdots & a\end{bmatrix}\right\} \cap \bB^m = \emptyset.
\]
\end{thm}
\begin{proof} We just need to show that if $\bB$ does not absorb $\bA$, then such $a,b,c$ exist for every $m$. Let $n = |\bA|(r_{|\bB|^2}(m)-1) + 1$. Then since $\bB$ doesn't absorb $\bA$ with respect to any term of arity $n$, there is some collection of $a_i \in \bA, b_{ij} \in \bB$ such that
\[
\Sg_{\bA^n}\left\{\begin{bmatrix} a_1 & b_{12} & \cdots & b_{1n}\\ b_{21} & a_2 & \cdots & b_{2n}\\ \vdots & \vdots & \ddots & \vdots\\ b_{n1} & b_{n2} & \cdots & a_n\end{bmatrix}\right\} \cap \bB^n = \emptyset.
\]

By the pigeonhole principle, there is some $a$ which occurs at least $n' = r_{|\bB|^2}(m)$ times among $a_1, ..., a_n$. Suppose without loss of generality that $a_1, ..., a_{n'}$ are all equal to $a$. If we restrict to the rows and columns with $a_i = a$, we find that
\[
\Sg_{\bA^{n'}}\left\{\begin{bmatrix} a & b_{12} & \cdots & b_{1n'}\\ b_{21} & a & \cdots & b_{2n'}\\ \vdots & \vdots & \ddots & \vdots\\ b_{n'1} & b_{n'2} & \cdots & a\end{bmatrix}\right\} \cap \bB^{n'} = \emptyset.
\]

Now we color the complete graph $K_{n'}$ with $|\bB|^2$ colors, coloring the edge $\{i,j\}$ (with $i < j$) with the color corresponding to the ordered pair $(b_{ij},b_{ji})$. Then by the definition of the Ramsey number $r_{|\bB|^2}(m)$, there is a monochromatic copy of $K_m$, with all edges colored by the color corresponding to some pair $(b,c) \in \bB^2$. By restricting to the rows and columns corresponding to the vertices of this monochromatic $K_m$, we see that
\[
\Sg_{\bA^{m}}\left\{\begin{bmatrix} a & b & \cdots & b\\ c & a & \cdots & b\\ \vdots & \vdots & \ddots & \vdots\\ c & c & \cdots & a\end{bmatrix}\right\} \cap \bB^{m} = \emptyset.\qedhere
\]
\end{proof}

\begin{cor} If $\bA = (A,m)$ is a finite majority algebra, then $\bB \lhd \bA$ iff there is a majority term $m' \in \Clo(m)$ such that $m'(\bB,\bA,\bB) \subseteq \bB$. Equivalently, we have $\bB \lhd \bA \iff \bB \lhd_J \bA$.

More precisely, if $\bB \lhd_J \bA$, then $\bB$ absorbs $\bA$ with respect to a term of arity at most $\lceil e\cdot |\bB|!\rceil$, where $e$ is Euler's constant $\sum_{n \ge 0} \frac{1}{n!} \approx 2.718$.
\end{cor}
\begin{proof} The weaker bound $\lceil e|\bA|\cdot |\bB|^2!\rceil$ on the arity of an absorbing term follows from the estimate $r_k(3) \le \lceil e\cdot k!\rceil$ and the fact that
\[
\begin{bmatrix}b\\ m'(c,a,b)\\ c\end{bmatrix} \in \Sg_{\bA^3}\left\{\begin{bmatrix} a & b & b\\ c & a & b\\ c & c & a\end{bmatrix}\right\}.
\]
However, we don't need the exact setup above. It's enough to find $n$ sufficiently large that for every $n\times n$ matrix
\[
\begin{bmatrix} a_1 & b_{12} & \cdots & b_{1n}\\ b_{21} & a_2 & \cdots & b_{2n}\\ \vdots & \vdots & \ddots & \vdots\\ b_{n1} & b_{n2} & \cdots & a_n\end{bmatrix}
\]
with off-diagonal entries in $\bB$, we can find $i,j,k$ distinct such that $b_{ij} = b_{ik}$ and $b_{ki} = b_{kj}$. If we set $b = b_{ij} = b_{ik}$ and $c = b_{ki} = b_{kj}$, then this will give us a submatrix of the form
\[
\begin{bmatrix} a_i & b & b\\ b_{ji} & a_j & b_{jk}\\ c & c & a_k\end{bmatrix},
\]
and applying $m'$ will give us an element of $\bB^3$.

Taking $n = \lceil e\cdot |\bB|!\rceil$ is good enough to find $i,j,k$ with $b_{ij} = b_{ik}$ and $b_{ki} = b_{kj}$. The proof is a minor adaptation of the proof of the upper bound on $r_k(3)$, and is left as an exercise to the reader.
\end{proof}

\begin{rem} It's intriguing that in the case of majority algebras, the bound on the arity of the absorbing operation only depends on the size of $|\bB|$.
\end{rem}

\begin{prob} Define $m(k)$ to be the least number such that whenever $\bA$ is a finite majority algebra, $\bB \lhd_J \bA$, and $|\bB| \le k$, we can always find a term $t$ of arity at most $m(k)$ such that $\bB$ absorbs $\bA$ with respect to $t$. How quickly does $m(k)$ grow?
\end{prob}

The dual discriminator algebra from Example \ref{ex-dual-discriminator} shows that we always have $m(k) \ge k+2$, while the previous result shows that $m(k) \le \lceil e\cdot k!\rceil$. For $k = 1,2$ we have $m(1) = 3, m(2) = 4$. Could it be that we have $m(k) = k+2$ for every $k$?

\subsection{Finitely related implies J\'onsson absorption is equivalent to absorption}\label{ss-finitely-related-jonsson}

The story of this subsection began with a paper of Larose, Loten, and Z{\'a}dori \cite{zadori-near-unanimity-graphs}, which gave a concrete characterization of undirected (reflexive) graphs which are preserved by a near-unanimity operation. Surprisingly, it turned out that the corresponding algebraic structure has a near-unanimity term iff it generates a congruence distributive variety (i.e., iff it has J\'onsson terms). This coincidence inspired what came to be known as the \emph{Z\'adori conjecture}: could it be the case that for every finitely related algebra, the existence of J\'onsson terms is equivalent to the existence of near-unanimity terms?

The Z\'adori conjecture was settled affirmatively by Barto \cite{near-unanimity-congruence-distributive}, showing as a consequence that it is possible to decide whether or not a given finite relational structure has a near-unanimity polymorphism. Zhuk \cite{near-unanimity-zhuk} independently proved the decidability result, with a better bound on the smallest arity of a near-unanimity polymorphism. Barto and Bul\'in \cite{deciding-absorption-relational} then extended Zhuk's approach to show that in a finitely related algebra, J\'onsson absorption becomes equivalent to absorption (note that by Corollaries \ref{cor-near-unanimity-absorption} and \ref{cor-congruence-distributive-jonsson-absorption}, this is a generalization of the Z\'adori conjecture).

The basic idea of (Barto and Bul\'in's extension of) Zhuk's approach is to suppose that $\bB \lhd_J \bA$, and that we are able to construct a $\bB$-essential relation $\RR \le \bA^n$ for some huge $n$ via a primitive positive formula $\Phi$ built out of relations $\RR_i \le \bA^m$ for small $m$. The goal is to prove a contradiction if $n$ is sufficiently large compared to $m$ and $|\bA|$, which will show that $\bB \lhd \bA$ by Theorem \ref{absorption-essential}. The core of the strategy is to reduce to the case where the primitive positive formula which defines the huge relation $\RR$ has a simple structure.

We'll start with a particularly simple case, where $\RR \le \bA^n$ is defined by a ``comb formula'' $\Phi$ over a collection of ternary relations $\RR_i \le \bA^3$ in the following way:
\[
(x_1, ..., x_n) \in \RR \;\;\; \iff \;\;\; \exists y_0, ..., y_n \in \bA \text{ s.t. } \bigwedge_{i \in [n]} (y_{i-1}, x_i, y_i) \in \RR_i.
\]
We call this a comb formula because the diagram we get by drawing the projection maps $\pi_j : \RR_i \rightarrow \bA$ looks like a comb:
\begin{center}
\begin{tikzpicture}[scale=2]
  \node (R1) at (-1.5,0) {$\RR_1$};
  \node (R2) at (-0.5,0) {$\RR_2$};
  \node (R3) at (0.5,0) {$\RR_3$};
  \node (R4) at (1.5,0) {$\RR_4$};
  \node (ay0) at (-2,0.1) {$\bA$};
  \node (ay1) at (-1,0.1) {$\bA$};
  \node (ay2) at (0,0.1) {$\bA$};
  \node (ay3) at (1,0.1) {$\bA$};
  \node (ay4) at (2,0.1) {$\bA$};
  \node (ax1) at (-1.5,-0.5) {$\bA$};
  \node (ax2) at (-0.5,-0.5) {$\bA$};
  \node (ax3) at (0.5,-0.5) {$\bA$};
  \node (ax4) at (1.5,-0.5) {$\bA$};
  \draw (R1) [->] -- (ay0);
  \draw (R1) [->] -- (ax1);
  \draw (R1) [->] -- (ay1);
  \draw (R2) [->] -- (ay1);
  \draw (R2) [->] -- (ax2);
  \draw (R2) [->] -- (ay2);
  \draw (R3) [->] -- (ay2);
  \draw (R3) [->] -- (ax3);
  \draw (R3) [->] -- (ay3);
  \draw (R4) [->] -- (ay3);
  \draw (R4) [->] -- (ax4);
  \draw (R4) [->] -- (ay4);
\end{tikzpicture}
\end{center}
(Note that $y_0, y_n$ do not actually add anything to the expressive power of a comb formula: in fact, there is no loss of generality in assuming that $\RR_1$ is the constraint $x_1 = y_1$ and that $\RR_n$ is the constraint $x_n = y_{n-1}$.) Our goal is to show that if $n$ is sufficiently large compared to $|\bA|$, then such a relation $\RR$ can never be $\bB$-essential. We will need to use the following consequence of Corollary \ref{cor-absorbing-directed-path} and Proposition \ref{prop-directed-cycle}.

\begin{lem} Suppose that $\bS \lhd_J \PP \le \bA\times \bA$, and that there are subsets $C, D \subseteq \bA$ such that
\begin{align*}
C + \bS &= C,\\
D - \bS &= D,
\end{align*}
and
\[
\PP \cap (C\times D) \ne \emptyset.
\]
Then $C \cap D \ne \emptyset$.
\end{lem}
\begin{proof} Pick any $(c,d) \in \RR \cap (C \times D)$. Since $C + \bS = C$ is equivalent to each vertex having in-degree at least $1$ in the digraph on $C$ with edge set $\bS \cap (C\times C)$, we can apply Proposition \ref{prop-directed-cycle} to see that there is some $c' \in C$ and some $k \ge 1$ such that $(c', c') \in \bS^{\circ k}$ and $(c', c) \in \bS^{\circ k}$. Similarly, there is some $d' \in D$ such that $(d,d'), (d',d') \in \bS^{\circ l}$ for some $l \ge 1$.

Since $\bS \subseteq \PP$, we see that
\[
(c',d') \in \bS^{\circ k} \circ \PP \circ \bS^{\circ l} \subseteq \PP^{\circ (k + 1 + l)},
\]
so we can apply Corollay \ref{cor-absorbing-directed-path} to see that there is a directed path from $c'$ to $d'$ in $\bS$. Since $c' \in C$ and $C + \bS = C$, we see that in fact we must have $d' \in C$ as well, so $d' \in C \cap D$.
\end{proof}

\begin{thm}[Barto and Bul\'in \cite{deciding-absorption-relational}]\label{thm-comb-essential} Suppose that $\RR \le \bA^n$ is defined by the formula
\[
(x_1, ..., x_n) \in \RR \;\;\; \iff \;\;\; \exists y_0, ..., y_n \in \bA \text{ s.t. } \bigwedge_{i \in [n]} (y_{i-1}, x_i, y_i) \in \RR_i,
\]
for some collection of ternary relations $\RR_i \le \bA^3$. If $\bB \lhd_J \bA$ and
\[
n \ge 3^{|\bA|} - 2^{|\bA|+1} + 3,
\]
then $\RR$ is not $\bB$-essential.
\end{thm}
\begin{proof} Suppose for the sake of contradiction that $\RR$ is $\bB$-essential. For each $i \in [n]$, define binary relations $\bS_i \lhd_J \PP_i \le \bA \times \bA$ by
\begin{align*}
\PP_i &= \pi_{13}(\RR_i),\\
\bS_i &= \pi_{13}\big(\RR_i \cap (\bA\times\bB\times\bA)\big).
\end{align*}
Since $\RR$ was assumed to be $\bB$-essential, for each $i \in [n]$ we have
\[
\bS_1 \circ \cdots \circ \bS_{i-1} \circ \PP_i \circ \bS_{i+1} \circ \cdots \circ \bS_n \ne \emptyset,
\]
but
\[
\bS_1 \circ \cdots \circ \bS_n = \emptyset.
\]
For each $0 \le i \le n$, if we define $\bC_i, \bD_i \le \bA$ by
\begin{align*}
\bC_i &= \bA + \bS_1 \circ \cdots \circ \bS_i,\\
\bD_i &= \bA - \bS_{i+1} \circ \cdots \circ \bS_n,
\end{align*}
then the assumption that $\RR$ is $\bB$-essential implies that
\[
\bC_i \cap \bD_i = \emptyset
\]
for all $0 \le i \le n$, and that for any $0 \le i < j \le n$ we have
\[
(\bC_i\times \bD_j) \cap (\PP_{i+1} \circ \cdots \circ \PP_j) \ne \emptyset,
\]
so in particular $\bC_i, \bD_j \ne \emptyset$ for $i,j \in [n-1]$. Since the number of ordered pairs $(C,D)$ of disjoint, nonempty subsets of $\bA$ is given by
\[
\#\big\{(C,D) \mid \emptyset \ne C,D \subseteq \bA, C\cap D = \emptyset\big\} = 3^{|\bA|} - 2^{|\bA|+1} + 1,
\]
by the pigeonhole principle there must be some pair $i < j \in [n-1]$ such that $(\bC_i, \bD_i) = (\bC_j, \bD_j)$. Now we can apply the previous lemma to $\PP = \PP_{i+1}\circ\cdots\circ\PP_j$, $\bS = \bS_{i+1}\circ\cdots\circ\bS_j$ and $C = \bC_i = \bC_j$, $D = \bD_i = \bD_j$ to see that $\bC_i \cap \bD_i \ne \emptyset$, contradicting the assumption that $\RR$ was $\bB$-essential.
\end{proof}

The next step of the proof is to consider relations $\RR$ defined by ``tree formulas'' $\Phi$, i.e. formulas such that the bipartite graph between variables and relations imposed on those variables forms a tree. The goal is to show that if such an $\RR$ is $\bB$-essential, then we can define a relation $\RR'$ whose arity is not too small compared to the arity of $\RR$, such that $\RR'$ is defined by a comb formula $\Phi'$ and $\RR'$ is also $\bB$-essential. The main idea is to use the fact that if a tree has bounded degree and has many leaves, then it must contain a long path.

\begin{thm}\label{thm-tree-to-comb} Suppose that the primitive positive formula $\Phi$ is a tree formula over relations $\RR_i \le \bA^{m_i}$, with all $m_i \le m$, which defines a $\bB$-essential relation $\RR$ of arity at least
\[
\max(m-1,2)^{n-3} + 2.
\]
Then there is a collection of ternary relations $\RR_i' \le \bA^3$, for $1 \le i \le n$, such that the relation $\RR'$ defined by the comb formula
\[
(x_1, ..., x_n) \in \RR' \;\;\; \iff \;\;\; \exists y_0, ..., y_n \in \bA \text{ s.t. } \bigwedge_{i \in [n]} (y_{i-1}, x_i, y_i) \in \RR_i'
\]
is also $\bB$-essential.
\end{thm}
\begin{proof} We may as well assume that $m \ge 3$. We can also make a series of simplifying assumptions about the structure of the tree formula $\Phi$.

First, we may assume that every free variable of $\Phi$ is a leaf. To see this, note that if some free variable $x$ of $\Phi$ is not a leaf, then we can modify $\Phi$ by replacing it with the tree formula
\[
\exists x \text{ s.t. } \Phi \wedge x = x',
\]
where $x'$ is a new free variable. Furthermore, we may assume that every bound variable of $\Phi$ is \emph{not} a leaf (by replacing any relation involving a bound leaf variable with its existential projection onto the other variables), so that the free variables of $\Phi$ are exactly the leaves.

Next we may assume that each bound variable has degree at most three in $\Phi$ (that is, that each variable is involved in at most three constraints). To do this, suppose that some bound variable $x$ is involved in four or more constraints, so that $\Phi$ has the form
\[
\Phi = \exists x \text{ s.t. } (x,...) \in \RR_1 \wedge \cdots \wedge (x, ...) \in \RR_k \wedge \Psi,
\]
with $k \ge 4$. Then we replace $\Phi$ with the equivalent tree formula
\[
\exists x, x' \text{ s.t. } x = x' \wedge (x,...) \in \RR_1 \wedge (x, ...) \in \RR_2 \wedge (x', ...) \in \RR_3 \wedge \cdots \wedge (x', ...) \in \RR_k \wedge \Psi,
\]
where $x$ now occurs in just three constraints, and $x'$ now occurs in $k-1$ constraints.

In fact, as long as $m \ge 3$, we may even assume that every bound variable has degree exactly two in $\Phi$. For this, we introduce the ternary relation $\RR_=$ defined by
\[
(x,y,z) \in \RR_= \;\;\; \iff \;\;\; x = y = z \in \bA.
\]
Now if $\Phi$ has a bound variable $x$ of degree three, so that $\Phi$ has the form
\[
\Phi = \exists x \text{ s.t. } (x,...) \in \RR_1 \wedge (x, ...) \in \RR_2 \wedge (x, ...) \in \RR_3 \wedge \Psi,
\]
then we replace $\Phi$ with the equivalent tree formula
\[
\exists x,x',x'' \text{ s.t. } (x,x',x'') \in \RR_= \wedge (x,...) \in \RR_1 \wedge (x', ...) \in \RR_2 \wedge (x'', ...) \in \RR_3 \wedge \Psi,
\]
where $x, x', x''$ now each occur in just two constraints.

Now that the simplifications are complete, let $E$ be the set of variables of $\Phi$ and let $V$ be the set of constraints of $\Phi$. Letting $F \subseteq E$ be the set of free variables of $\Phi$, we can interpret $(V \sqcup F, E)$ as a graph, where each bound variable $x \in E$ is interpreted as an edge connecting the two constraints it is involved in, and each free variable $y \in E$ is interpreted as an edge connecting the constraint it is involved in to the leaf $y$ considered as an element of $F$. Since $\Phi$ is a tree formula, the graph $(V\sqcup F, E)$ is also a tree, and the maximum degree of any vertex is bounded by $m$ by assumption.

Now pick any $v_1 \in F$, and inductively define a non-backtracking path $v_1, v_2, ... \in V \sqcup F$ through the tree $(V\sqcup F, E)$ as follows:
\begin{itemize}
\item we pick $v_2 \in V$ to be the unique constraint involving the free variable $v_1$,

\item if we have already picked $v_{i-1}, v_i$ and $v_i \in V$, then we pick $v_{i+1} \ne v_{i-1}$ among the neighbors of $v_i$ so that the number of free variables which are reachable from $v_{i+1}$ in the subgraph formed by deleting $v_i$ is maximized, and

\item if $v_i \in F$ for $i > 1$, then we stop.
\end{itemize}
For each $i \ge 1$, let $F_i$ be the set of free variables which are reachable from $v_{i+1}$ in the subgraph formed by deleting $v_i$. By construction, we have $F_1 = F \setminus \{v_1\}$ and
\[
|F_i| > 1 \;\;\; \implies \;\;\; |F_i| \ge |F_{i+1}| \ge \Big\lceil\frac{|F_i|}{m-1}\Big\rceil.
\]
Setting $F_0 = F$, let $0 = i_1 < i_2 < \cdots < i_k$ be the collection of $i$s such that $F_i \ne F_{i+1}$. For each $i_j$, we pick some free variable $x_j \in F_{i_j}\setminus F_{i_j + 1}$, so that $x_1 = v_1$ and $x_k$ is the final free variable at the end of the path $v_1, v_2, ...$. Also, let $y_0 = y_1 = x_1$, $y_{k-1} = y_k = x_k$, and for $2 \le j \le k-2$ we let $y_j$ be the bound variable which connects the constraint $v_{i_j}$ to the constraint $v_{i_j+1}$.

Now we define a pp-formula $\Phi'$ by starting from $\Phi$, imposing the unary constraint $z \in \bB$ for each unused free variable $z \in F \setminus \{x_1, x_2, ..., x_k\}$, and existentially quantifying all such $z$s, so that the remaining free variables are just $x_1, x_2, ..., x_k$ and the resulting constraint is still $\bB$-essential. Chopping the formula $\Phi'$ into pieces by removing the variables $x_1, ..., x_k$ and $y_2, ..., y_{k-2}$, we see that each piece has boundary equal to some triple $\{y_{j-1}, x_j, y_j\}$, and by existentially projecting the remaining variables in each piece we get an equivalent constraint
\[
(y_{j-1}, x_j, y_j) \in \RR_j' \le \bA^3,
\]
with $\RR_1' = \RR_k' = \RR_=$, so we see that $\Phi'$ is equivalent to the comb formula
\[
\exists y_0, ..., y_k \in \bA \text{ s.t. } \bigwedge_{i \in [k]} (y_{i-1}, x_i, y_i) \in \RR_i'.
\]

To finish the proof, we just need to check that $k \ge n$ (and we may as well assume that $n \ge 3$, since comb formulas are silly when $n \le 2$). For this, we note that by our assumption on the arity of the original $\bB$-essential relation $\RR$ we have
\[
|F_{i_2}| = |F_1| = |F_0| - 1 \ge \max(m-1,2)^{n-3} + 1 > \max(m-1,2)^{n-3},
\]
so inductively we have
\[
|F_{i_j}| > \max(m-1,2)^{n-j-1}
\]
for $2 \le j \le n-1$, so in particular
\[
|F_{i_{n-1}}| > 1,
\]
which shows that $k > n-1$.
\end{proof}

The last ingredient of the proof is to show that if some primitive positive formula $\Phi$ defines a $\bB$-essential relation of large arity, then we can repeatedly unravel cycles in $\Phi$ to produce a tree formula $\Phi'$ which defines a $\bB$-essential relation of the same arity. This is tricky in two different ways: first, we need to apply Theorem \ref{absorbing-diagonal} (a nontrivial application of J\'onsson absorption) in order to unravel a cycle, and second, it is not obvious that repeatedly unravelling cycles will eventually end up producing a tree formula. We start by showing how we will apply Theorem \ref{absorbing-diagonal}.

\begin{lem}\label{lem-jonsson-unravelling} If $\bS \lhd_J \PP \le \bA \times \bA$, then we have
\[
\exists x \text{ s.t. } (x,x) \in \bS \;\;\; \iff \;\;\; \exists x_0, ..., x_{|\bD|} \in \bD \text{ s.t. } (x_0,x_1) \in \bS \wedge \cdots \wedge (x_{|\bD|-1}, x_{|\bD|}) \in \bS,
\]
where $\bD = \{x \mid (x,x) \in \PP\}$.
\end{lem}
\begin{proof} Left to right is clear: if $(x,x) \in \bS$, then $x \in \bD$ and we can take $x_0 = \cdots = x_{|\bD|} = x$.

For the other direction, suppose that there are $x_0, ..., x_{|\bD|} \in \bD$ such that $(x_i, x_{i+1}) \in \bS$ for each $i$. By the pigeonhole principle, there must be some $i < j$ such that $x_i = x_j$, so the digraph $(\bD, \bS\cap \bD^2)$ contains a cycle of length $j-i$. Letting $\bC \le \bD$ be the subalgebra consisting of elements $x \in \bD$ which are contained in a cycle of length $j-i$ in this digraph, we see that $\bC$ is nonempty and that $\bS \cap \bC^2$ is subdirect in $\bC \times \bC$.

Since $\bC \subseteq \bD$, we see that $\PP \cap \bC^2$ contains the diagonal $\Delta_\bC$, so $\bS \cap \bC^2$ is subdirect in $\bC \times \bC$ and J\'onsson absorbs the diagonal $\Delta_\bC$. Thus by Theorem \ref{absorbing-diagonal} there is some $x \in \bC$ such that $(x,x) \in \bS$.
\end{proof}

\begin{thm}[Barto, Bul\'in \cite{deciding-absorption-relational}, generalizing Zhuk \cite{near-unanimity-zhuk}] If $\bB \lhd_J \bA$, and if $\Phi$ is a primitive positive formula over relations $\RR_i \le \bA^{m_i}$ which defines a $\bB$-essential relation $\RR$, then there is a tree formula $\Phi'$ over the same relations $\RR_i$ together with some unary relations $\bD_j \le \bA$ which defines a $\bB$-essential relation $\RR'$ of the same arity as $\RR$.
\end{thm}
\begin{proof} Let $n$ be the arity of $\RR$. As in the proof of Theorem \ref{thm-tree-to-comb}, we may as well assume that each free variable of $\Phi$ occurs in just one constraint of $\Phi$. Furthermore, by splitting variables into copies and adding equality constraints between the copies, we may as well assume that $\Phi$ is defined by starting with a (quantifier-free) tree formula $\Psi$, adding a bunch of equality constraints between variables of $\Psi$, and then existentially quantifying, i.e.
\[
\Phi(x_1, ..., x_n) = \exists y_1, ..., y_k \text{ s.t. } \Psi(x,y) \wedge \bigwedge_{(i,j) \in E} y_i = y_j
\]
for some $E \subseteq [k] \times [k]$.

For any $(i_0,j_0) \in E$, we can unravel a cycle of $\Phi$ to produce the formula $\Phi_{(i_0,j_0)}$ with $n\cdot (|\bA|+1)$ free variables, which is defined by
\begin{align*}
\Phi_{(i_0,j_0)}(x_1^0, ..., x_n^{|\bA|}) = \exists y_1^0, ..., y_k^{|\bA|} \text{ s.t. } &\bigwedge_{0 \le l \le |\bA|} \Big(\Psi(x^l,y^l) \wedge \bigwedge_{(i,j) \in E\setminus (i_0,j_0)} y_i^l = y_j^l\Big)\\
&\wedge \bigwedge_{0 \le l < |\bA|} y_{i_0}^{l+1} = y_{j_0}^l \in \bD_{i_0},
\end{align*}
where
\[
\bD_{i_0} = \Big\{y_{i_0} \mid \exists x_1, ..., x_n, y_1, ..., y_k \text{ s.t. } \Psi(x,y) \wedge \bigwedge_{(i,j) \in E} y_i = y_j\Big\}.
\]

{\bf Claim.} For any $(i_0, j_0) \in E$, the formula $\Phi_{(i_0,j_0)}$ defines a relation $\RR_{(i_0,j_0)} \le \big(\bA^{|\bA|+1}\big)^n$ which is $\bB^{|\bA|+1}$-essential.

{\bf Proof of claim.} By the assumption that the original relation $\RR$ was $\bB$-essential, for each $i \le n$ there is some tuple $(x_1, ..., x_n) \in \RR$ such that $x_j \in \bB$ for all $j \ne i$. Thus there are $y_1, ..., y_k$ such that
\[
\Psi(x,y) \wedge \bigwedge_{(i,j) \in E} y_i = y_j,
\]
and taking $x^l_j = x_j$ and $y^l_j = y_j$ for all $j$ shows that $\RR_{(i_0,j_0)}$ contains a tuple in $\big(\bB^{|\bA|+1}\big)^{i-1}\times\bA^{|\bA|+1}\times \big(\bB^{|\bA|+1}\big)^{n-i}$.

To finish the proof of the claim, we just need to check that $\RR_{(i_0,j_0)}$ doesn't contain any tuple in $\big(\bB^{|\bA|+1}\big)^n$. Define binary relations $\bS \lhd_J \PP \le \bA\times \bA$ by
\begin{align*}
(y_{i_0},y_{j_0}) \in \PP \;\;\; &\iff \;\;\; \exists x_1, ..., x_n \in \bA, \exists y_1, ..., y_k \text{ s.t. } \Psi(x,y) \wedge \bigwedge_{(i,j) \in E\setminus (i_0,j_0)} y_i = y_j,\\
(y_{i_0},y_{j_0}) \in \bS \;\;\; &\iff \;\;\; \exists x_1, ..., x_n \in \bB, \exists y_1, ..., y_k \text{ s.t. } \Psi(x,y) \wedge \bigwedge_{(i,j) \in E\setminus (i_0,j_0)} y_i = y_j.
\end{align*}
Then we have
\[
\bD_{i_0} = \{y_{i_0} \mid (y_{i_0}, y_{i_0}) \in \PP\},
\]
and $\RR_{(i_0,j_0)}$ contains a tuple in $\big(\bB^{|\bA|+1}\big)^n$ iff
\[
\exists y_{i_0}^0, ..., y_{i_0}^{|\bA|+1} \in \bD_{i_0} \text{ s.t. } \bigwedge_{0 \le l < |\bA|} (y_{i_0}^l, y_{i_0}^{l+1}) \in \bS.
\]
By Lemma \ref{lem-jonsson-unravelling}, this is equivalent to the existence of a loop $(y_{i_0}, y_{i_0}) \in \bS$, which is equivalent to $\RR$ containing a tuple from $\bB^n$ - which can't happen by our assumption that $\RR$ is $\bB$-essential.

Now by Lemma \ref{essential-groups}, we can find $\bD_i^l \in \{\bA,\bB\}$ and a set $I \subseteq [n]\times[|\bA|+1]$ with $|I| = n$ and
\[
\big|I \cap \{i\} \times [|\bA|+1]\big| = 1
\]
for all $i$, such that
\[
\RR_{(i_0,j_0)}' = \pi_I\Big(\RR_{(i_0,j_0)} \cap \Big(\prod_{i,l} \bD_i^l\Big)\Big) \le \bA^n
\]
is $\bB$-essential. This corresponds to a formula $\Phi_{(i_0,j_0)}'$ which is derived from $\Phi_{(i_0,j_0)}$ by restricting some of the free variables to lie in $\bB$ and then existentially quantifying all free variables aside from those labeled by $I$.

Now for the tricky part: we need to show that there is some way of choosing $(i_0,j_0) \in E$ to guarantee that the new formula $\Phi_{(i_0,j_0)}'$ is, in some sense, closer to being a tree formula than the original $\Phi$ was. We will need to use the fact that if there is some bound variable $y$ such that removing it splits $\Phi_{(i_0,j_0)}'$ into several components, one of which contains no free variables, then we can prune the component which contains no free variables and replace it with a unary constraint on $y$ without changing the relation $\RR_{(i_0,j_0)}'$.

Let $G_\Phi$ be the bipartite graph of variables and constraints of the formula $\Phi$, with edges between variables and the constraints which involve them. A vertex of $G_\Phi$ is called a \emph{cut-vertex} if deleting it from $G_\Phi$ splits $G_\Phi$ into more than one connected component. An induced subgraph $\mathcal{C}$ of $G_\Phi$ is called a \emph{biconnected component} of $G_\Phi$ if it is maximal with respect to the property of $\mathcal{C}$ having no cut-vertices of its own. The following are standard results from graph theory:
\begin{itemize}
\item every cycle of $G_\Phi$ (with no repeated vertices) is contained in some biconnected component of $G_\Phi$,
\item each edge of $G_\Phi$ belongs to exactly one biconnected component of $G_\Phi$,
\item the vertices which appear in more than one biconnected component are exactly the cut-vertices of $G_\Phi$, and
\item the bipartite graph of cut-vertices and biconnected components of $G_\Phi$ forms a tree.
\end{itemize}

Following Zhuk \cite{near-unanimity-zhuk}, we call a biconnected component of $G_\Phi$ ``trivial'' if it consists of just a single edge (between a variable and a constraint of $\Phi$), and we mainly focus on how large a given nontrivial biconnected component of $G_\Phi$ is, with a secondary focus on how close we are to being able to prune a given biconnected component. We can consider ourselves very close to pruning a biconnected component $\mathcal{C}$ if it contains a cut-vertex $y$ such that when we delete $y$, very many of the free variables of $\Phi$ are not reachable from $\mathcal{C}$. Thus we define
\[
\operatorname{MaxCut}(\mathcal{C}) \coloneqq \max_{y \in \mathcal{C}} \#\{\text{free variables of }\Phi\text{ which can't be reached from }\mathcal{C}\text{ without passing through }y\},
\]
and we pick a nontrivial biconnected component $\mathcal{C}_0$ with the following properties:
\begin{itemize}
\item the size $|\cC_0|$ should be maximal among biconnected components, and
\item $\operatorname{MaxCut}(\cC_0)$ should be minimal among biconnected components $\cC$ satisfying $|\cC| = |\cC_0|$.
\end{itemize}
In other words, we pick the biconnected component which we are currently farthest from being able to prune, among the components which are as large as possible. Once we have picked such a component $\mathcal{C}_0$, we pick any of the extra equality constraints $y_i = y_j$, $(i,j) \in E$, which occurs within it and take $(i_0,j_0) = (i,j)$.

To finish the proof, we just need to show that the number of nontrivial biconnected components $\cC'$ of $G_{\Phi_{(i_0,j_0)}'}$ with $|\cC'| = |\cC_0|$ and $\operatorname{MaxCut}(\cC') = \operatorname{MaxCut}(\cC_0)$ is decreased, and that no $\cC'$ has $|\cC'| > |\cC_0|$ or $|\cC'| = |\cC_0|$ and $\operatorname{MaxCut}(\cC') < \operatorname{MaxCut}(\cC_0)$. It's easy to see that if $\cC'$ involves any copy of any edge of $\cC_0$, then $|\cC'| < |\cC_0|$, since $\cC'$ can't contain any copy of the constraint relation $y_{i_0} = y_{j_0}$. Otherwise, $\cC'$ must be a copy of some biconnected component $\cC \ne \cC_0$ of $G_\Phi$, so at least we have $|\cC'| = |\cC| \le |\cC_0|$.

We may assume that $\Phi$ has already been pruned in the following two ways:
\begin{itemize}
\item for any nontrivial biconnected component $\mathcal{C}$ of $G_\Phi$ with $\operatorname{MaxCut}(\cC) = n$, we can prune $\cC$ and replace it by a unary constraint on the cut-vertex separating $\cC$ from the free variables of $\Phi$, and
\item for any pair of cut-vertices $y_i, y_j$ such that there are no free variables $x$ which can reach both of $y_i, y_j$ without passing through the other, we can prune any nontrivial biconnected components which can reach both of $y_i, y_j$ without passing through the other, replacing them with a path of binary constraints connecting $y_i$ to $y_j$.
\end{itemize}
The second form of pruning allows us to make the following statement about $\cC_0$: for any other nontrivial biconnected component $\cC \ne \cC_0$ with $|\cC| = |\cC_0|$, the cut-vertex $y$ which achieves the maximum in the definition of $\operatorname{MaxCut}(\cC)$ is the one which separates $\cC$ from $\cC_0$. To see this, let $y_0 \in \cC_0$ be the cut-vertex which separates $\cC_0$ from $\cC$, and note that if $y$ does not separate $\cC$ from $\cC_0$, then every free variable which is separated from $\cC$ by $y$ is also separated from $\cC_0$ by $y_0$, so the assumption $\operatorname{MaxCut}(\cC_0) \le \operatorname{MaxCut}(\cC)$ implies that we should have pruned $\cC$ and replaced it by a path of binary constraints from $y_0$ to $y$.

Now consider any nontrivial biconnected component $\cC'$ of $G_{\Phi_{(i_0,j_0)}'}$ which is a copy of some nontrivial biconnected component $\cC \ne \cC_0$ of $G_\Phi$ with $|\cC| = |\cC_0|$. If $y \in \cC$ is the cut-vertex separating $\cC$ from $\cC_0$, and if $y' \in \cC'$ is the corresponding vertex of $G_{\Phi_{(i_0,j_0)}'}$, then $y'$ separates $\cC'$ from every one of the equality relations $y_{i_0}^{l+1} = y_{j_0}^l$, from every copy of any free variable $x$ which was separated from $\cC$ by $y$, and from every other copy of $\cC$ in $G_{\Phi_{(i_0,j_0)}'}$. Therefore we have
\[
\operatorname{MaxCut}(\cC') \ge \operatorname{MaxCut}(\cC),
\]
with equality iff every single free variable $x_i$ which is not separated from $\cC$ by $y$ has the corresponding free variable $x_i^l$ (with $(i,l) \in I$) of $\Phi_{(i_0,j_0)}'$ not separated from $\cC'$ by $y'$. In particular, there is at most one copy $\cC'$ of $\cC$ which satisfies $\operatorname{MaxCut}(\cC') = \operatorname{MaxCut}(\cC)$, and if such a copy exists then every other copy of $\cC$ is ready to be pruned from $\Phi_{(i_0,j_0)}'$.
\end{proof}

We can now put the pieces together to get the main result of this section. Note that although we have stated all of the results so far as if there was only a single J\'onsson absorbing subalgebra $\bB \lhd_J \bA$ of interest, they can all be straightforwardly generalized to the case where there are several J\'onssons absorbing subalgebras $\bB_i \lhd_J \bA_i \le \bA$, and where we are trying to put bounds on the potential arity of $(\bB_{i_1}, ..., \bB_{i_n})$-essential relations $\RR \le \bA_{i_1} \times \cdots \times \bA_{i_n} \le \bA^n$. Thus we can use Theorem \ref{absorption-essential-gen} to get the following slightly stronger result.

\begin{thm}[Barto, Bul\'in \cite{deciding-absorption-relational}]\label{thm-finitely-related-jonsson-absorption} If $\bA$ is finitely related, with $\Inv(\bA)$ generated by relations of arity at most $m$, then there is a term $t$ of arity at most
\[
\max(m-1,2)^{3^{|\bA|} - 2^{|\bA|+1}} + 2,
\]
such that for every J\'onsson absorption $\bB_i \lhd_J \bA_i \le \bA$, $\bB_i$ absorbs $\bA_i$ with respect to $t$.
\end{thm}

\begin{cor}\label{cor-finitely-related-near-unanimity} If $\bA$ is finitely related, then $\bA$ generates a congruence distributive variety iff $\bA$ has a near-unanimity term. In this case, if $\Inv(\bA)$ is generated by relations of arity at most $m$, then $\bA$ has a near-unanimity term operation of arity at most
\[
\max(m-1,2)^{3^{|\bA|} - 2^{|\bA|+1}} + 2.
\]
\end{cor}

It's natural to wonder if these bounds can be improved. As it turns out, Zhuk \cite{near-unanimity-zhuk} found examples showing that sometimes the least arity of a near-unanimity term really is doubly exponential in $|\bA|$. To understand these examples, it is helpful to start by constructing a comb formula which defines a $\bB$-essential relation of singly exponential arity, for an absorbing subalgebra $\bB$.

Since ternary relations are tricky to visualize, we will actually start by trying to construct an exponentially long sequence of pairs of binary relations $\bS_i \lhd \PP_i$ such that
\[
\bS_1 \circ \cdots \circ \bS_{i-1} \circ \PP_i \circ \bS_{i+1} \circ \cdots \circ \bS_n \ne \emptyset
\]
for all $i$, but such that $\bS_1 \circ \cdots \circ \bS_n = \emptyset$. In order to compactly represent both relations $\bS_i, \PP_i$ simultaneously, we will draw a bipartite graph where solid edges represent elements which are contained in both $\bS_i$ and $\PP_i$, while dashed lines represent elements which are contained in $\PP_i$ but not in $\bS_i$. For instance, if we had
\begin{align*}
\PP &= \Big\{\begin{bmatrix} a\\ a\end{bmatrix}, \begin{bmatrix} a\\ b\end{bmatrix}, \begin{bmatrix} b\\ b\end{bmatrix}\Big\},\\
\bS &= \Big\{\begin{bmatrix} a\\ a\end{bmatrix}, \begin{bmatrix} b\\ b\end{bmatrix}\Big\},
\end{align*}
then we would represent this situation by drawing the bipartite graph below.
\begin{center}
\begin{tikzpicture}[scale=1,baseline=0.5cm]
  \node (a1) at (-1,1) {$a$};
  \node (b1) at (-1,0) {$b$};
  \node (a2) at (0.5,1) {$a$};
  \node (b2) at (0.5,0) {$b$};
  \draw (a1) edge (a2) (b1) edge (b2);
  \draw [dashed] (a1) edge (b2);
\end{tikzpicture}
\end{center}
Note that for this particular pair $\bS$, $\PP$, it is \emph{not possible} to have $\bS \lhd_J \PP$ by Theorem \ref{absorbing-directed-path}.

\begin{ex}[Zhuk \cite{near-unanimity-zhuk}]\label{ex-comb-formula} Let $A = \{*\}\cup \{0, 1, ..., n\}$, and for each $k \le n$ define a pair of binary relations $S_k \subset P_k$ by
\begin{align*}
x + P_k &= \begin{cases} \{x\} & x > k,\\ \{*, 0, ..., k-1\} & x \in \{k, *\},\\ \emptyset & x \in \{0, ..., k-1\},\end{cases}\\
S_k &= P_k \setminus \{(k,*)\},
\end{align*}
so that when $n = 4, k = 2$ we can represent the situation visually by the bipartite graph below.
\begin{center}
\begin{tikzpicture}[scale=1,baseline=0.5cm]
  \node (41) at (-1,4) {$4$};
  \node (31) at (-1,3) {$3$};
  \node (21) at (-1,2) {$2$};
  \node (11) at (-1,1) {$1$};
  \node (01) at (-1,0) {$0$};
  \node (*1) at (-1,-1) {$*$};
  \node (42) at (1,4) {$4$};
  \node (32) at (1,3) {$3$};
  \node (22) at (1,2) {$2$};
  \node (12) at (1,1) {$1$};
  \node (02) at (1,0) {$0$};
  \node (*2) at (1,-1) {$*$};
  \draw (41) edge (42) (31) edge (32);
  \draw (21) edge (12) (21) edge (02);
  \draw [dashed] (21) edge (*2);
  \draw (*1) edge (12) (*1) edge (02) (*1) edge (*2);
\end{tikzpicture}
\end{center}
We will find an exponentially long sequence of binary relations $S_{k_0}, S_{k_1}, S_{k_2}, ...$ which compose to the empty set, such that for each $i$, if we replace $S_{k_i}$ by $P_{k_i}$ then the composition becomes nonempty. When $n = 2$, we can represent these all-but-one paths as the rows of a visually self-similar matrix:
\[
\begin{bmatrix} * & * & * & * & * & * & * & *\\ 0 & * & * & * & * & * & * & *\\ 1 & 1 & * & * & * & * & * & *\\ 1 & 1 & 0 & * & * & * & * & *\\ 2 & 2 & 2 & 2 & * & * & * & *\\ 2 & 2 & 2 & 2 & 0 & * & * & *\\ 2 & 2 & 2 & 2 & 1 & 1 & * & *\\ 2 & 2 & 2 & 2 & 1 & 1 & 0 & *\end{bmatrix},
\]
and by examining this matrix we can see that the desired sequence of relations is given by
\[
S_* \circ S_0 \circ S_1 \circ S_0 \circ S_2 \circ S_0 \circ S_1 \circ S_0 = \emptyset,
\]
where $S_* = \Delta_{A \setminus \{*\}}$, $P_* = \Delta_A$ (the first step of the paths has been left out of the rows of the matrix, since it obscures the pattern). The pattern generalizes: for $i > 0$ we will take
\[
k_i = \max\{j \text{ s.t. } 2^j \text{ divides } i\},
\]
and for $i = 0$ we use $S_* \subset P_*$. In particular, the sequence will have a total of $2^{n+1} = 2^{|A|-1}$ binary relations composed together to produce the empty set.

Now we will show that there is a near-unanimity operation $t$ of arity $2^{n+1} + 1$ which preserves the relations $S_k, P_k$ and witnesses the absorptions $S_k \lhd P_k$ for all $k \in A$. Such a $t$ is given by
\[
t(x_0, ..., x_{2^{n+1}}) = \begin{cases} n & \text{if }\#\{i \mid x_i = n\} > 2^n,\\ n-1 & \text{otherwise, if }\#\{i \mid x_i = n-1\} > 2^{n-1},\\ \vdots & \vdots\\ 0 & \text{otherwise, if }\#\{i \mid x_i = 0\} > 1,\\ * & \text{otherwise, i.e. if }\#\{i \mid x_i = k\} \le 2^k\text{ for all }k.\end{cases}
\]
To see that $t$ witnesses $S_* \lhd P_*$, or equivalently that $t$ witnesses $A\setminus\{*\} \lhd A$, note that if all but one of the $x_i$ are not $*$, then there must be some $k$ such that $\#\{i \mid x_i = k\} > 2^k$, since
\[
1 + 2 + \cdots + 2^n < 2^{n+1}.
\]
The argument to show that $t$ witnesses $S_k \lhd P_k$ is similar: if $(x_i,y_i) \in P_k$ and all but one are in $S_k$, and if $t(x_0, ..., x_{2^{n+1}}) = k$, then by the definition of $t$ and $P_k$ at least $2^k$ of the $y_i$s must be in the set $\{0, ..., k-1\}$, so there must be some $j \in \{0, ..., k-1\}$ such that $\#\{i \mid y_i = j\} > 2^j$, so $t(y_0, ..., y_{2^{n+1}}) \ne *$.

In order to turn this example into a comb formula, we just define ternary relations $R_{i+1}$ by
\[
(y,x,y') \in R_{i+1} \;\;\; \iff \;\;\; x \in \{*, 0\} \; \wedge \; (y,y') \in P_{k_i} \; \wedge \; (x = * \implies (y,y') \in S_{k_i}).
\]
To see that $R_{i+1}$ is preserved by $t$, note that for $x_i \in \{*,0\}$ we have $t(x_0, ...) = 0$ unless all but at most one $x_i$ is $*$, so the fact that $R_{i+1}$ is preserved by $t$ follows from the fact that $t$ witnesses $S_{k_i} \lhd P_{k_i}$. The comb formula
\[
(x_1, ..., x_{2^{n+1}}) \in R \;\;\; \iff \;\;\; \exists y_0, ..., y_{2^{n+1}} \text{ s.t. } \bigwedge_{i \in [2^{n+1}]} (y_{i-1}, x_i, y_i) \in R_i
\]
then defines a $\{*\}$-essential relation $R$ of arity $2^{n+1} = 2^{|A|-1}$ which is preserved by the near-unanimity operation $t$.
\end{ex}

In order to construct an example demonstrating doubly exponential behavior, we reverse-engineer the proof of Theorem \ref{thm-tree-to-comb} to produce a tree formula of exponential radius, so that paths through it look like the comb formula constructed in Example \ref{ex-comb-formula}. However, we will need to introduce new relations of higher arity which are not preserved by the near-unanimity operation $t$ from Example \ref{ex-comb-formula}.

\begin{ex}[Zhuk \cite{near-unanimity-zhuk}, Barto, Draganov \cite{near-unanimity-minimal}]\label{ex-tree-formula-high-arity} Let $A, S_k, P_k$ be defined as in Example \ref{ex-comb-formula}, and let $m \ge 3$. For each $k \in \{0, ..., n\}$, define an $m$-ary relation $R_{m,k}$ by
\[
(x_1, ..., x_{m-1}, y) \in R_{m,k} \;\;\; \iff \;\;\; \bigwedge_{i \in [m-1]} (x_i,y) \in P_k \; \wedge \bigvee_{j \in [m-1]} (x_j,y) \in S_k.
\]
Note that the formula above is \emph{not} a primitive positive formula, and that the $m$-ary relation $R_{m,k}$ is not preserved by the near-unanimity operation $t$ from Example \ref{ex-comb-formula} for $m \ge 3$.

To illustrate how we will make use of the relations $R_{m,k}$, note that for any $k_1, k_2$ we have
\begin{align*}
&\exists y_1, ..., y_{m-1} \text{ s.t. } \bigwedge_{i \in [m-1]} (x_{i,1}, ..., x_{i,m-1}, y_i) \in R_{m,k_1} \wedge \; (y_1, ..., y_{m-1}, z) \in R_{m,k_2}\\
&\implies \;\;\; \bigwedge_{(i,j) \in [m-1]^2} (x_{i,j},z) \in P_{k_1}\circ P_{k_2} \; \wedge \bigvee_{(i,j) \in [m-1]^2} (x_{i,j},z) \in S_{k_1}\circ S_{k_2}.
\end{align*}
We now inductively pp-define relations $T_{m,j}$ of arity $(m-1)^j + 1$ by $T_{m,1} = R_{m,n}$ and
\begin{align*}
(x_1, ..., x_{(m-1)^{j+1}}, z) \in T_{m,j+1} \;\;\; \iff \;\;\; \exists y_1, ... y_{m-1} \text{ s.t. } &\bigwedge_{i \in [m-1]} (x_{(i-1)(m-1)^j+1}, ..., x_{i(m-1)^j}, y_i) \in T_{m,j}\\
&\wedge \; (y_1, ..., y_{m-1}, z) \in R_{m,k_{2^n+j}},
\end{align*}
where $n = k_{2^n}, k_{2^n+1}, ..., k_{2^{n+1}-1}$ is the second half of the sequence from Example \ref{ex-comb-formula}. Then we define a relation $T_m$ of arity $(m-1)^{2^n} = (m-1)^{2^{|A|-2}}$ by
\[
(x_1, ..., x_{(m-1)^{2^n}}) \in T_m \;\;\; \iff \;\;\; \exists z \text{ s.t. } (x_1, ..., x_{(m-1)^{2^n}}, z) \in T_{m,2^n}.
\]
It's now fairly straightforward to verify that $T_m$ is an $\{n\}$-essential relation, so there is no near-unanimity operation of arity $(m-1)^{2^n}$ which preserves all of the $m$-ary relations $R_{m,k}$. (For a direct proof that no such near-unanimity operation exists, see \cite{near-unanimity-minimal}.)

To finish up this example, we need to find a near-unanimity operation $u$ of arity $(m-1)^{2^n}+1$ which preserves the $m$-ary relations $R_{m,k}$ for all $k \in \{0, ..., n\}$. If we make the convention that $* < k$ for every $k \in \{0, ..., n\}$, then such an operation is given by
\[
u(x_0, ..., x_{(m-1)^{2^n}}) = \begin{cases} n & \text{if }\;\#\{i \mid x_i \le n\} \; > \; (m-1)^{2^n}\cdot \#\{i \mid x_i < n\},\\ n-1 & \text{otherwise, if }\;\#\{i \mid x_i \le n-1\} \; > \; (m-1)^{2^{n-1}}\cdot \#\{i \mid x_i < n-1\},\\ \vdots & \vdots\\ 0 & \text{otherwise, if }\;\#\{i \mid x_i \le 0\} \; > \; (m-1)\cdot \#\{i \mid x_i < 0\},\\ * & \text{otherwise, i.e. if }\;\#\{i \mid x_i \le k\} \; \le \; (m-1)^{2^k}\cdot \#\{i \mid x_i < k\}\text{ for all }k.\end{cases}
\]
In order to check that $u$ preserves $R_{m,n}$, the important point is that if none of the inputs $x_i$ to $u$ were equal to $n$ and the output of $u$ ended up being $*$, then we must have had
\begin{align*}
(m-1)^{2^n} + 1 &= \#\{i \mid x_i < n\}\\
&\le (m-1)^{2^{n-1}}\cdot\#\{i \mid x_i < n-1\}\\
&\le \cdots\\
&\le (m-1)^{2^{n-1} + 2^{n-2} + \cdots + 1}\cdot\#\{i \mid x_i < 0\}\\
&= (m-1)^{2^n-1}\cdot\#\{i \mid x_i = *\},
\end{align*}
so strictly more than $m-1$ of the inputs $x_i$ to $u$ must have been equal to $*$. Thus, if we have
\[
\begin{bmatrix} x_0^1\\ \vdots\\ x_0^{m-1}\\ y_0\end{bmatrix}, \begin{bmatrix} x_1^1\\ \vdots\\ x_1^{m-1}\\ y_1\end{bmatrix}, \cdots \in R_{m,n}
\]
but somehow
\[
\begin{bmatrix} u(x_0^1, x_1^1, ...)\\ \vdots\\ u(x_0^{m-1}, ...)\\ u(y_0, y_1, ...)\end{bmatrix} = \begin{bmatrix} n\\ \vdots\\ n\\ *\end{bmatrix} \not\in R_{m,n},
\]
then each $y_i < n$ since $(x_i^1, y_i) \in P_n$, so by what we just showed we see that
\[
\#\{i \mid y_i = *\} \ge m,
\]
and for each $i$ with $y_i = *$ we see that there must be some $j_i \in [m-1]$ with
\[
(x_i^{j_i},y_i) = (x_i^{j_i}, *) \in S_n \;\;\; \implies \;\;\; x_i^{j_i} = *,
\]
so by the pigeonhole principle there must be some $j \in [m-1]$ such that $x_i^j = * \ne n$ for at least two $i$s, in which case we have
\[
u(x_0^j, ..., x_{(m-1)^{2^n}}^j) \ne n,
\]
contrary to our assumption. The verification that $u$ preserves $R_{m,k}$ for the remaining values of $k$ is similar, and can be found in \cite{near-unanimity-minimal}.
\end{ex}

Example \ref{ex-tree-formula-high-arity} shows that the bound from Theorem \ref{thm-tree-to-comb} is nearly best possible when $m \ge 3$, so if we want to significantly improve the bound on the arity of a near-unanimity operation in Corollary \ref{cor-finitely-related-near-unanimity} then our only hope is to improve the comb-formula estimate in Theorem \ref{thm-comb-essential} from $3^{|\bA|}$ to something closer to $2^{|\bA|}$.

What about the case $m = 2$? We can't build a doubly exponential essential relation out of the binary relations $S_k, P_k$ from Example \ref{ex-comb-formula}, since we already know that they are preserved by a near-unanimity operation of singly exponential arity. The trick is to introduce one more element into the domain, after which we can use slight variations on the binary relations $S_k, P_k$ to pp-define the ternary relations $R_{3,k}$ from Example \ref{ex-tree-formula-high-arity}.

\begin{ex}[Zhuk \cite{near-unanimity-zhuk}, Barto, Draganov \cite{near-unanimity-minimal}]\label{ex-tree-formula-binary} Let $A^\pm = \{*_-, *_+, 0, 1, ..., n\}$. For each $k \in \{0, ..., n\}$, define the binary relation $Q_k^\pm$ by
\begin{align*}
x + Q_k^\pm &= \begin{cases} \{x\} & x > k,\\ \{*_+, 0, ..., k-1\} & x = k,\\ \emptyset & x \in \{0, ..., k-1\},\\ \{*-, *_+, 0, ..., k-1\} & x \in \{*_-, *_+\},\end{cases}
\end{align*}
so that when $n = 4, k = 2$ we represent $Q_2^\pm$ by the bipartite graph below.
\begin{center}
\begin{tikzpicture}[scale=1,baseline=0.5cm]
  \node (41) at (-1,4) {$4$};
  \node (31) at (-1,3) {$3$};
  \node (21) at (-1,2) {$2$};
  \node (11) at (-1,1) {$1$};
  \node (01) at (-1,0) {$0$};
  \node (*p1) at (-1,-1) {$*_+$};
  \node (*m1) at (-1,-2) {$*_-$};
  \node (42) at (1,4) {$4$};
  \node (32) at (1,3) {$3$};
  \node (22) at (1,2) {$2$};
  \node (12) at (1,1) {$1$};
  \node (02) at (1,0) {$0$};
  \node (*p2) at (1,-1) {$*_+$};
  \node (*m2) at (1,-2) {$*_-$};
  \draw (41) edge (42) (31) edge (32);
  \draw (21) edge (12) (21) edge (02);
  \draw (21) edge (*p2);
  \draw (*p1) edge (12) (*p1) edge (02) (*p1) edge (*p2) (*p1) edge (*m2);
  \draw (*m1) edge (12) (*m1) edge (02) (*m1) edge (*p2) (*m1) edge (*m2);
\end{tikzpicture}
\end{center}
Similarly, we define $Q_k^\mp$ by swapping the roles of $*_+$ and $*_-$ in $Q_k^\pm$. Then the relation $R_{3,k}^\pm$ defined by the primitive positive formula
\[
(x_1, x_2, y) \in R_{3,k}^\pm \;\;\; \iff \;\;\; (x_1,y) \in Q_k^\pm \; \wedge \; (x_2,y) \in Q_k^\mp
\]
is essentially the same as the relation $R_{3,k}$ from Example \ref{ex-tree-formula-high-arity}, so we can use the same construction to show that there is no near-unanimity operation of arity $2^{2^n} = 2^{2^{|A^\pm| - 3}}$ which preserves $Q_k^\pm$ and $Q_k^\mp$ for all $k$. A near-unanimity operation of arity $2^{2^n} + 1$ which does preserve $Q_k^\pm, Q_k^\mp$ is the operation $u_\pm$ given by
\[
u_\pm(x_0, ..., x_{(m-1)^{2^n}}) = \begin{cases} n & \text{if }\;\#\{i \mid x_i \le n\} \; > \; 2^{2^n}\cdot \#\{i \mid x_i < n\},\\ n-1 & \text{otherwise, if }\;\#\{i \mid x_i \le n-1\} \; > \; 2^{2^{n-1}}\cdot \#\{i \mid x_i < n-1\},\\ \vdots & \vdots\\ 0 & \text{otherwise, if }\;\#\{i \mid x_i \le 0\} \; > \; 2\cdot \#\{i \mid x_i < 0\},\\ *_+ & \text{otherwise, if }\;\#\{i \mid x_i = *_+\} \; > \; \#\{i \mid x_i = *_-\},\\ *_- & \text{otherwise.}\end{cases}
\]
Note that $u_\pm$ preserves the equivalence relation $\theta$ which identifies $*_+$ with $*_-$, and that the quotient operation $u_\pm / \theta$ is exactly $u$. The proof of the fact that $u_\pm$ preserves $Q_k^\pm, Q_k^\mp$ is similar to the proof of the fact that $u$ preserves $R_{3,k}$, and can be found in \cite{near-unanimity-minimal}.
\end{ex}

\section{Finding an arc-consistent absorbing subinstance}\label{s-absorbing-arc-consistent}

In this section we'll go over Marcin Kozik's proof from \cite{pq-consistency} (which refined the argument from \cite{slac}) of the fact that every cycle-consistent instance has a cycle-consistent subinstance such that every domain is absorption-free. In fact, Kozik proves something stronger, involving a weaker consistency notion known as $pq$-consistency. The technique for the proof can be viewed as a generalization of the argument for the case of majority algebras, but it is much more difficult because we can't assume that all the relations involved are binary. The main idea of the proof was originally developed in \cite{nu-pathwidth}, for the sake of proving a technical lemma about absorption generalizing Theorem \ref{absorbing-diagonal} and Theorem \ref{ideal-diagonal}, which was needed to show that near-unanimity CSPs can be solved in NL (nondeterministic logspace).

First we define the weaker consistency notion known as $pq$-consistency (Kozik names it $jpq$-consistency in \cite{pq-consistency}). The basic idea behind this definition is that it is a consistency check which (aside from assuming arc-consistency) only involves pairwise projections of constraints, only computes compositions of these binary relations along cycles, and is strong enough to rule out the existence of a cycle such that each binary relation along it is the graph of a permutation and the composition of all these permutations is not the identity permutation.

\begin{defn}\label{defn-pq-consistent} A CSP instance $\fX$ with domains $\bA_{v_i}$ corresponding to variables $v_i$ is called $pq$\emph{-consistent} if
\begin{itemize}
\item it is arc-consistent, i.e. each relation $\RR \le \bA_{v_1}\times\cdots \times \bA_{v_k}$ imposed on the variables is subdirect, and
\item for each variable $v$ and each pair of cycles $p,q$ of $\fX$ which begin and end at $v$, there exists some $j \ge 0$ such that the binary relation $\PP_{j(p+q)+p}$ corresponding to the path $j(p+q)+p$ (see Definition \ref{path-defn}) contains the diagonal $\Delta_{\bA_v}$, i.e. for each $a \in \bA_v$, we have $a \in \{a\} + j(p+q)+p$ (see Definition \ref{path-action-defn}).
\end{itemize}
\end{defn}

The reader may find it interesting to check that in the proofs that we have already given for the fact that cycle-consistency solves ancestral CSPs and majority CSPs, we may substitute $pq$-consistency for cycle-consistency everywhere without significantly complicating the arguments. The reason for introducing the slightly more technical notion of $pq$-consistency is that the ``standard semidefinite relaxation'' of a CSP naturally produces a $pq$-consistent instance, but doesn't always produce a cycle-consistent instance - and the semidefinite relaxation is the tool used to ``robustly'' solve bounded width CSPs in \cite{sdp}.

The main step of the argument is the following technical result.

\begin{thm}[Kozik \cite{pq-consistency}]\label{absorbing-reduction} If $\fX$ is a $pq$-consistent instance and $\fY$ is an arc-consistent subinstance of $\fX$ defined by restricting each domain and each relation of $\fX$ to an absorbing subalgebra, and if any domain of $\fY$ has a proper absorbing subalgebra, then there is a proper arc-consistent subinstance $\fZ$ of $\fY$ defined by restricting each domain to an absorbing subalgebra.
\end{thm}

Before diving into the proof of this result, we'll show how it can be used.

\begin{thm}\label{absorbing-pq} If $\fX$ is a $pq$-consistent instance with domains $\bA_v$, then there is a $pq$-consistent subinstance $\fX'$ of $\fX$ defined by restricting each domain $\bA_v$ to a minimal absorbing subalgebra $\bA_v'$. If $\fX$ is cycle-consistent, then so is $\fX'$.
\end{thm}
\begin{proof} By repeatedly applying Theorem \ref{absorbing-reduction}, we may find an arc-consistent subinstance $\fX'$ of $\fX$ such that each domain has no proper absorbing subalgebra. Then by Theorem \ref{arc-consistent-cycles}, for every cycle $r$ from $v$ to $v$ of $\fX$ such that $\PP_r \supseteq \Delta_{\bA_v}$, if $r'$ is the corresponding cycle in $\fX'$, then we have $\PP_{r'} \supseteq \Delta_{\bA_v'}$. If $\fX$ is $pq$-consistent, then for any cycles $p,q$ from $v$ to $v$ there is some $j$ such that $\PP_{j(p+q)+p} \supseteq \Delta_{\bA_v}$, so on taking $r = j(p+q)+p$ we see that the corresponding cycles $p',q'$ have $\PP_{j(p'+q')+p'} \supseteq \Delta_{\bA_v'}$.
\end{proof}

The proof of Theorem \ref{absorbing-reduction} will only rely on three properties of absorption. Since there are several absorption-like concepts that have proven useful, and most of them satisfy these properties, we will consider an arbitrary ``absorption concept'' $\lhd_X$ which applies to certain pairs $\bB \le \bA$, and which satisfies the following three properties.
\begin{itemize}
\item {\bf Compatibility with pp-formulas.} If $\bS_i \lhd_X \RR_i$ are relations, and if a relation $\RR$ is defined by a pp-formula $\Phi$ involving the relations $\RR_1, ..., \RR_k$ (and possibly some other relations), then if we define a relation $\bS$ by the pp-formula $\Phi'$ defined by replacing each $\RR_i$ by $\bS_i$ in $\Phi$, we have $\bS \lhd_X \RR$.

\item {\bf Transitive closure.} If $\bC \lhd_X \bB \lhd_X \bA$, then $\bC \lhd_X \bA$.

\item {\bf Connectivity transfers.} If $\bS \lhd_X \RR$ and $\RR \le \bA \times \bA$, and if $a,b \in \bA$ are such that $(a,a), (b,b) \in \bS$ and $(a,b) \in \RR$, then there is some $k$ such that $(a,b) \in \bS^{\circ k}$.
\end{itemize}
Note that by the local criterion for J\'onsson absorption (Theorem \ref{local-jonsson}), if $\lhd_X$ is compatible with pp-formulas, then the connectivity transfer property of $\lhd_X$ is equivalent to the implication $\bB \lhd_X \bA \implies \bB \lhd_J \bA$. Also, a trivial case of compatibility with pp-formulas implies that for all $\bA$, we have $\bA \lhd_X \bA$.

Throughout most of the proof, we will be focusing on the arc-consistent instance $\fY$. Therefore, for each variable $v$ of $\fY$, we let $\bA_v$ be the corresponding domain in $\fY$, and let $\bA_v^{\fX}$ be the corresponding domain in the original $pq$-consistent instance $\fX$, so we have $\bA_v \lhd_X \bA_v^{\fX}$. Similarly, if $\RR$ refers to a relation in $\fY$, then we let $\RR^{\fX}$ refer to the corresponding relation in $\fX$, with $\RR \lhd_X \RR^{\fX}$.

The argument strategy generalizes the strategy used for majority algebras. We will consider the set $\cB$ of ordered pairs $(x,\bB)$ such that $x$ is a variable of $\fY$, $\bB \lhd_X \bA_x$, and $\bB \ne \emptyset, \bA_x$. We want to define a quasiorder $\preceq$ on $\cB$, such that if restricting the domain of the variable $x$ to $\bB$ and imposing arc-consistency forces another variable $y$ to have its domain restricted to $\bC$, then we have $(x,\bB) \preceq (y,\bC)$. Unfortunately, it is not enough to consider paths alone to define this partial order: general deductions involving arc-consistency involve reasoning about \emph{trees}.

\begin{defn} To every relational structure $\fA = (A, R_1, ...)$ we associate the bipartite graph $\cG_\fA$ with vertex sets $A$ and $R_1 \sqcup \cdots$, and edge set consisting of pairs $(a,r)$ for every $a \in A$ and $r \in R_i$ such that some coordinate of $r$ is equal to $a$ (if $a$ occurs as a coordinate of $r$ multiple times, then we make multiple copies of the edge $(a,r)$).

We say that $\fA$ is a \emph{tree} if the associated bipartite graph $\cG_\fA$ is a tree (so in particular, no tuple $r$ in any relation $R_i$ can have any repeated coordinates).
\end{defn}

Kozik \cite{pq-consistency} extends the concepts of paths and addition of paths to trees in order to define the partial order $\preceq$ on $\cB$ properly.

\begin{defn} If $\fY$ is a CSP instance, viewed as a relational structure, then we define a \emph{tree pattern} $p$ from $x$ to $y$ to consist of the following information:
\begin{itemize}
\item a relational structure $\fA = (A,R_1, ...)$ which is a tree, with each relation of $\fA$ corresponding to a relation of $\fY$,
\item a homomorphism of relational structures $h : \fA \rightarrow \fY$,
\item a subset $I \subseteq A$ of the elements of $A$ which we call the set of \emph{inputs} to the pattern, such that for all $i \in I$ we have $h(i) = x$, and
\item an element $o \in A$ which we call the \emph{output} of the pattern, such that $h(o) = y$.
\end{itemize}
If $p$ is a tree pattern from $x$ to $y$, then we may view it as a CSP instance via the homomorphism $h : \fA \rightarrow \fY$. If $\bB \le \bA_x$, then we define $\bB + p$ to be the subalgebra of values $b \in \bA_y$ such that the instance $\fA$ has a solution with the variables from $I$ assigned to values in $\bB$, and with the variable $o$ assigned to the value $b$.

If $p$ is a tree pattern from $x$ to $y$, and if $q$ is a tree pattern from $y$ to $z$, then we define the tree pattern $p+q$ by attaching a copy of $p$ to each input of $q$, combining the output of each copy of $p$ to the corresponding input of $q$. This definition is set up to ensure that $\bB+(p+q) = (\bB+p)+q$ for any $\bB \le \bA_x$.
\end{defn}

\begin{prop} If $p$ is a tree pattern from $x$ to $y$ in an arc-consistent instance $\fY$ and $\bB \lhd_X \bA_x$, then $\bB + p \lhd_X \bA_y$.
\end{prop}
\begin{proof} This follows from the fact that $\lhd_X$ is compatible with pp-formulas: we have $\bA_x + p = \bA_y$ if $\fY$ is arc-consistent, and so $\bB+p \lhd_X \bA_x+p = \bA_y$.
\end{proof}

Note that unlike the situation for path patterns, arc-consistency of the instance $\fY$ is no longer enough to ensure that $\bB \ne \emptyset \implies \bB+p \ne \emptyset$ for all tree patterns $p$. So we can no longer take as given that the subalgebras we construct will always be nonempty.

\begin{defn} Define the quasiordered set $(\cB, \preceq)$ to be the set of ordered pairs $(x,\bB)$ such that $x$ is a variable of the instance $\fY$, $\bB \lhd_X \bA_x$, and $\bB \ne \emptyset,\bA_x$, with the quasiorder defined by $(x,\bB) \preceq (y,\bC)$ if there exists a tree pattern $p$ from $x$ to $y$ with $\bB + p = \bC$.
\end{defn}

As in the argument for majority algebras, we now pick a maximal component $\cC$ of the quasiordered set $(\cB,\preceq)$ (since $\cB$ is nonempty by assumption and is finite, such a maximal component exists). We would like to use $\cC$ to define our reduced instance $\fZ$, but we no longer have a guarantee that there is at most one set $\bB$ with $(x,\bB) \in \cC$ for a given variable $x$.

A worst case scenario would be that there exist $\bB_1, \bB_2$ with $(x,\bB_i) \in \cC$ such that $\bB_1 \cap \bB_2 = \emptyset$: in this case, we would have no hope of using $\cC$ to define an arc-consistent reduction, because no matter which $(y,\bC) \in \cC$ we pick, there exist tree patterns $p_1, p_2$ from $y$ to $x$ with $\bC + p_i = \bB_i$, so reducing the domain $\bA_y$ to $\bC$ and imposing arc-consistency would make it impossible to assign any value to $x$. The main step of the proof is ruling out this scenario.

\begin{lem}\label{lem-absorbing-reduction} If $\cC$ is a maximal component of $(\cB,\preceq)$, and if $(x,\bB), (x,\bC) \in \cC$, then $\bB \cap \bC \ne \emptyset$.
\end{lem}

Before proving the lemma, we'll show how we can use it to finish the proof of Theorem \ref{absorbing-reduction}. This step won't use the fact that the instance $\fX$ is $pq$-consistent, or the fact that $\lhd_X$ transfers connectivity: the lemma is where these crucial facts are used.

\begin{proof}[Proof of Theorem \ref{absorbing-reduction}, assuming the lemma] Note that if $(x,\bB), (x,\bC) \in \cC$, then we can splice together tree patterns to show that $(x,\bB\cap\bC) \in \cC$ as well (so long as $\bB \cap \bC \ne \emptyset$, which follows from the lemma). So for every $x$, we can define a subalgebra $\bB_x \lhd_X \bA_x$ by taking $\bB_x$ to be the intersection of all $\bB$ such that $(x,\bB) \in \cC$ (or taking $\bB_x = \bA_x$ if no such $\bB$ exist). We define the absorbing subinstance $\fZ$ by reducing the domains of $\fY$ from $\bA_x$ to $\bB_x$. We need to check that $\fZ$ is arc-consistent.

Consider a single relation $\RR \le_{sd} \bA_{x_1} \times \cdots \times \bA_{x_k}$ of $\fY$. We wish to show that $\RR\cap \prod_i \bB_{x_i}$ is subdirect in $\prod_i \bB_{x_i}$. We will show by induction on $i$ that
\[
\pi_i(\RR \cap \prod_{j\le i} \bB_{x_j} \times \prod_{l > i} \bA_{x_l}) = \bB_{x_i}.
\]
The base case $i=1$ follows from the fact that $\fY$ is arc-consistent. For the inductive step, we pick any $(y,\bC) \in \cC$ and splice together tree patterns $p_j$ from $y$ to $x_j$ with $\bC + p_j = \bB_{x_j}$ for $j < i$ such that $\bB_{x_j} \ne \bA_{x_j}$ together with the relation $\RR$ to make a tree pattern $p$ from $y$ to $x_i$ with
\[
\bC + p = \pi_i(\RR \cap \prod_{j \le i-1} \bB_{x_j} \times \prod_{l > i-1} \bA_{x_l}),
\]
and note that by the induction hypothesis the right hand side is nonempty. Thus we either have $\bC + p = \bA_{x_i}$ or $(x_i,\bC+p) \in \cC$, and in either case we have $\bB_{x_i} \subseteq \bC+p$ (by the lemma), which completes the proof.
\end{proof}

Now we finally prove the crucial lemma.

\begin{proof}[Proof of the lemma] Suppose for contradiction that the lemma is not true, and choose $\bC$ maximal such that $(x,\bC) \in \cC$ and such that there exists $(x,\bB) \in \cC$ with $\bB \cap \bC = \emptyset$. Let $\bB_1, ..., \bB_k$ be the set of minimal $\bB$s such that $(x,\bB) \in \cC$ and $\bB \cap \bC = \emptyset$. Note that since the set of $\bB$s with $(x,\bB) \in \cC$ is closed under nonempty intersection, we must have $\bB_i \cap \bB_j = \emptyset$ for all $i \ne j$. Additionally, any $\bB$ with $(x,\bB) \in \cC$ and $\bB \cap \bC = \emptyset$ must contain at least one $\bB_i$.

Choose tree patterns $p_i, q, r$ from $x$ to $x$ such that $\bB_i + p_i = \bB_{i+1}$, $\bC + q = \bB_1$, $\bB_k + r = \bC$. Define the tree pattern $p$ by $p = q + p_1 + \cdots + p_{k-1} + r$, and note that $\bC + p = \bC$. We will mainly work inside the instance $\fA$ corresponding to the tree pattern $p$.

First we prune the inputs of the tree pattern $p$ a little bit to make a new tree pattern $p'$ (with the same instance $\fA$), removing variables of $\fA$ from the input set one at a time as long as we can remove one while keeping $\bC + p' = \bC$. Now pick any remaining input variable $s \in \fA$ of $p'$ (at least one input variable remains at the end of the pruning process, by the arc-consistency of $\fY$), and let $t$ be the output variable of $p'$ (note that $s,t$ are both mapped to $x$ in $\fY$). Let $p''$ be $p'$ with $s$ removed from its input set. Consider the binary relation $\bS \le \bA_x \times \bA_x$ consisting of pairs $(a,b)$ such that some solution of the instance $\fA$ assigns the value $a$ to $s$, assigns the value $b$ to $t$, and assigns all input variables of $p''$ to values in $\bC$.

Since $\bC + p' = \bC$, we have
\[
\bC + \bS = \bC + p' = \bC,
\]
and because of the pruning process we have
\[
\pi_2(\bS) = \bC + p'' \ne \bC,
\]
so by the maximal choice of $\bC$ we have $\pi_2(\bS) \cap \bB_i \ne \emptyset$ for all $i$. By splicing $p''$ together with a tree pattern $q_i$ with $\bC + q_i = \bB_i$ (merging their outputs together), we see that $(x,\pi_2(\bS) \cap \bB_i) \in \cC$, so by the minimality of $\bB_i$ we have
\[
\pi_2(\bS) \supseteq \bB_i
\]
for all $i$. Thus the subalgebra
\[
\bB_i - \bS = \pi_1(\bS \cap \bA_x\times \bB_i)
\]
is nonempty, has $(\bB_i-\bS) \cap \bC = \emptyset$ since $(\bC+\bS) \cap \bB_i = \emptyset$, and by splicing $p''$ with the same $q_i$ and changing the output to $s$, we see that $(x,\bB_i-\bS) \in \cC$. Thus there is some $j_i$ such that $\bB_i-\bS \supseteq \bB_{j_i}$. Then we have
\[
(\bB_{j_i} + \bS) \cap \bB_i \ne \emptyset,
\]
and by another tree splice (this time splicing $q_{j_i}$ into $p''$ by merging the output of $q_{j_i}$ with $s$) we see that either $\bB_{j_i} + \bS = \bA_x$ or $(x,\bB_{j_i} + \bS) \in \cC$, so by the minimality of $\bB_i$ we have
\[
\bB_{j_i} + \bS \supseteq \bB_i.
\]

Thus we have
\[
\cup_i \bB_i + \bS \supseteq \cup_i\bB_i,
\]
so if we consider $\bS$ as a digraph on $\bA_x$, we see that there is some directed cycle of $\bS$ which is entirely contained in $\cup_i \bB_i$. From $\bC + \bS = \bC$, we also see that there is some directed cycle of $\bS$ which is entirely contained in $\bC$. The plan is to apply Corollary \ref{cor-absorbing-directed-path} to produce a directed path in $\bS$ from an element of $\bC$ to an element of $\cup_i \bB_i$, which will give us a contradiction since any directed path in $\bS$ which starts in $\bC$ must end up in $\bC$.

In order to apply Corollary \ref{cor-absorbing-directed-path}, we need to construct a binary relation $\RR$ such that $\bS \lhd_X \RR$ and such that there is a directed path from $\bC$ to $\cup_i \bB_i$ in $\RR$. This is where we will finally use the assumption that $\fY$ absorbs a bigger instance $\fX$ which is $pq$-consistent. We define an instance $\fA^\fX$ similarly to $\fA$, but with each domain replaced with the corresponding domain in $\fX$ and similarly for the relations, and define $\RR$ to be the projection of the solution set to $\fA^\fX$ onto the variables $s,t$. Then since $\lhd_X$ is compatible with pp-formulas and since every domain/relation restriction in sight is absorbing, we have $\bS \lhd_X \RR$.

Now pick any $b \in \cup_i \bB_i$ which is contained in a directed cycle of $\bS$. Suppose $b \in \bB_i$. Consider the path from $s$ to the output variable of $q + p_1 + \cdots + p_{i-1}$ in $\fA$, call this path $\alpha$, and let $\beta$ be the path from that output variable to $t$ in $\fA$. The images of these paths in $\fX$ are cycles $\alpha_\fX, \beta_\fX$ from $x$ to $x$, so by the $pq$-consistency of $\fX$ there must exist some $j \ge 0$ such that $b \in \{b\} + j(\beta_\fX+\alpha_\fX)+\beta_\fX$. Note that by the arc-consistency of $\fX$, $\RR$ is the binary relation corresponding to the cycle $\alpha_\fX+\beta_\fX$. Additionally, since
\[
\bC + q + p_1 + \cdots + p_{i-1} = \bB_i,
\]
there is some $c \in \bC$ such that $b \in \{c\} + \alpha_\fX$. Thus we have
\[
b \in \{c\} + \alpha_\fX + j(\beta_\fX + \alpha_\fX) + \beta_\fX = \{c\} + (j+1)(\alpha_\fX + \beta_\fX) = \{c\} + \RR^{\circ (j+1)}.
\]
Additionally, by following paths of $\bS$ backwards sufficiently many times, we see that $c$ is reachable from a directed cycle of $\bS$ which is entirely contained in $\bC$. Thus there is some $m$ such that for some $a \in \bC$, we have $(a,a), (b,b) \in \bS^{\circ m}$ and $(a,b) \in \RR^{\circ m}$, and since $\bS^{\circ m} \lhd_X \RR^{\circ m}$ we may apply the transfer of connectivity property to see that for some $n$ we have $(a,b) \in \bS^{\circ n}$, which gives us our contradiction.
\end{proof}


To finish the analysis of bounded width algebras, we just need to understand the case where all the domains are absorption free. For this we need two main ingredients: first is that binary relations are forced to be boring unless some absorption occurs, and second is that if a simple algebra has an exciting ternary relation whose binary projections are boring, then the algebra must be affine and therefore does not have bounded width.

\subsection{Absorption constants}

In this subsection, we'll go over the proof of an interesting generalization of the Theorem \ref{absorbing-diagonal} to higher arity relations, from \cite{nu-pathwidth}, which we mentioned at the beginning of the last section. Since the proof of this result is so similar to the proof of Kozik's result from the last section, this seems like an appropriate place to cover the argument.

We will need some notation for the diagonal of a power $\bA^n$. One option is to use the notation $\Delta_\bA^n$, but this looks very similar to the notation we use in the appendix on commutator theory (Appendix \ref{a-commutator}). Another notation some authors use is $0_\bA^n$, so that when $n = 2$ we get the least congruence $0_\bA$, but I am not a big fan of this notation either. Yet another possibility is $=_\bA^n$. I decided on a fourth option, which allows me to refer to specific elements of the diagonal without too much ugliness, and which emphasizes the fact that the diagonal is isomorphic to $\bA$.

\begin{defn} For any $n$, define the \emph{diagonal subalgebra} of arity $n$ to be the subalgebra $\bA^{(n)} \le_{sd} \bA^n$, given by
\[
\bA^{(n)} = \{(a, ..., a) \mid a \in \bA\}.
\]
Additionally, for each $a \in \bA$, we define the corresponding \emph{constant tuple} to be
\[
a^{(n)} = (a, ..., a) \in \bA^{(n)}.
\]
\end{defn}

\begin{thm}[Theorem 6 of \cite{nu-pathwidth}] If $\bA$ is finite, $\RR \le_{sd} \bA^n$ is subdirect, and $\RR$ J\'onsson absorbs the diagonal $\bA^{(n)}$, then $\RR \cap \bA^{(n)} \ne \emptyset$.
\end{thm}
\begin{proof} The proof strategy is similar to the proof of Theorem \ref{absorbing-diagonal}. We assume without loss of generality that $\bA$ is idempotent, and we induct on $|\bA|$. It's enough to show that there is some proper subalgebra $\bB \le \bA$ such that
\[
\RR \cap \bB^n \le_{sd} \bB^n,
\]
since $\RR \cap \bB^n$ will automatically J\'onsson absorb $\bB^{(n)}$.

Similarly to the proof of Theorem \ref{absorbing-reduction}, we define the a quasiordered set $(\cB, \preceq)$ to be the set of subalgebras $\bB \le \bA$ with $\bB \ne \emptyset, \bA$, with the quasiorder defined by $\bB \preceq \bC$ if there is a tree pattern $p$ built out of copies of the relation $\RR$ such that $\bB + p = \bC$.

Pick some maximal component $\cC$ of $(\cB, \preceq)$. Note that if $\bB, \bC \in \cC$ have $\bB \cap \bC \ne \emptyset$, then we can splice together tree patterns to see that $\bB \cap \bC \in \cC$ as well.

First suppose that every pair $\bB, \bC \in \cC$ have $\bB \cap \bC \ne \emptyset$. Then there is some $\bB \in \cC$ which is contained in all other elements of $\cC$. We claim that in this case, we have
\[
\RR \cap \bB^n \le_{sd} \bB^n,
\]
which will allow us to complete the proof. To check this, we prove by induction on $i$ that
\[
\bB \subseteq \pi_i(\RR \cap \bB^{i-1}\times\bA^{n-i+1}).
\]
For $i = 1$ this follows from the assumption that $\RR$ is subdirect. For the induction step, note that the induction hypothesis implies that
\[
\pi_i(\RR \cap \bB^{i-1}\times\bA^{n-i+1}) \ne \emptyset,
\]
so
\[
\pi_i(\RR \cap \bB^{i-1}\times\bA^{n-i+1}) \in \cC \cup \{\bA\},
\]
and either way it contains $\bB$.

Now suppose, for the sake of contradiction, that there are $\bB, \bC \in \cC$ with $\bB \cap \bC = \emptyset$. As in the proof of Lemma \ref{lem-absorbing-reduction}, we take $\bC \in \cC$ maximal under inclusion such that there exists some $\bB \in \cC$ with $\bB \cap \bC = \emptyset$, and we let $\bB_1, ..., \bB_k$ be the set of minimal (under inclusion) $\bB$s such that $\bB \in \cC$ and $\bB \cap \bC = \emptyset$.

We continue following the proof of Lemma \ref{lem-absorbing-reduction}, defining tree patterns $p_i, q, r$ built out of copies of the relation $\RR$ such that $\bB_i + p_i = \bB_{i+1}$, $\bC + q = \bB_1$, $\bB_k + r = \bC$, and defining $p$ by
\[
p = q + p_1 + \cdots + p_k + r.
\]
Then we prune the inputs of the pattern $p$ to make a pattern $p'$ with as few inputs as possible such that
\[
\bC + p' = \bC,
\]
and let $p''$ be the pattern we get from $p'$ by removing one additional input $s$ from the input set, and define
\[
\bS \le \bA \times \bA
\]
as the set of possible pairs of values for the pruned input $s$ and the output of the pattern $p''$ which extend to assignments where every remaining input of $p''$ is given a value in $\bC$, as in the proof of Lemma \ref{lem-absorbing-reduction}.

By the exact same argument from Lemma \ref{lem-absorbing-reduction}, we have
\[
\bC + \bS = \bC + p' = \bC
\]
and
\[
\cup_i \bB_i + \bS \supseteq \cup_i \bB_i,
\]
so every element $c \in \bC$ is reachable from a directed cycle of $\bS$ which is entirely contained in $\bC$, and $\cup_i\bB_i$ contains a directed cycle of $\bS$.

Now we finally deviate slightly from the proof of Lemma \ref{lem-absorbing-reduction}. Define a pattern $p_=$ by replacing each occurence of $\RR$ by $\RR \cup \bA^{(n)}$ in the pattern $p$, and similarly define $p_=', p_=''$, and define
\[
S_= \subseteq \bA \times \bA
\]
as the set of possible pairs of values for the pruned input $s$ and the output of the pattern $p_=''$ which extend to assignments where every remaining input of $p_=''$ is given a value in $\bC$. The $\bS \subseteq S_=$ and $\bS$ J\'onsson absorbs $S_=$. Additionally, for each $i$ we have
\[
\bC + p_= \supseteq \bC + q + p_1 + \cdots + p_{i-1} = \bB_i,
\]
since we can simply feed a bunch of equal copies of an element $b \in \bB_i$ to each of the remaining levels of the tree pattern. Thus we have
\[
\bC + S_= \supseteq \cup_i \bB_i.
\]
Thus we can find a directed path in $S_=$ from some element $c \in \bC$ which is contained in a directed cycle of $\bS$ to some element of $\cup_i \bB_i$ which is contained in a directed cycle of $\bS$. This allows us to apply Corollary \ref{cor-absorbing-directed-path} to conclude that there is some directed path from $\bC$ to $\cup_i\bB_i$ in $\bS$, which gives us our contradiction.
\end{proof}

Ross Willard points out the following consequence of this result.

\begin{cor} If $\bA$ is a finite algebra, then there is at least one element $a \in \bA$ such that, for all subdirect relations $\RR \le_{sd} \bA^n$, we have
\[
\RR \lhd_J \Sg_{\bA^n}(\RR \cup \bA^{(n)}) \;\; \implies \;\; a^{(n)} \in \RR.
\]
The set of such elements $a$ forms a J\'onsson absorbing subalgebra of $\bA$.
\end{cor}

\begin{defn} Say that $a$ is an \emph{absorption constant} of $\bA$ with respect to the absorption concept $\lhd_X$ if
\[
\RR \lhd_X \Sg_{\bA^n}(\RR \cup \bA^{(n)}) \;\; \implies \;\; a^{(n)} \in \RR
\]
for all subdirect relations $\RR \le_{sd} \bA^n$. Let
\[
\operatorname{Abs}_X(\bA) \lhd_X \bA
\]
be the set of absorption constants of $\bA$ with respect to $\lhd_X$.
\end{defn}

\begin{prob}[Ross Willard] Can we give an independent characterization of the canonical absorbing subalgebra $\operatorname{Abs}(\bA)$? What can we do with it?
\end{prob}

\section{Zhuk's centers and ternary absorption}

In this section we'll go over a very strong technique introduced by Zhuk in his proof of the dichotomy conjecture \cite{zhuk-dichotomy}, which produces ternary absorption as soon as we have a certain type of binary relation on a pair of Taylor algebras. This technique allows us to both simplify and strengthen one of the key results needed for the study of general Taylor algebras, known as the ``absorption theorem''.

First, we'll go over the history of this idea, so the reader can understand where the definition comes from and why it is (somewhat) natural.

The main idea behind Zhuk's approach in \cite{zhuk-dichotomy} is to note that if an algebra is not polynomially complete, then its polynomial clone must be contained in a maximal proper subclone of the clone of all functions (that every proper subclone is contained in a \emph{maximal} proper subclone follows from the fact that the clone of all functions is finitely generated: in fact, it's generated by the set of functions of arity $2$). A maximal clone corresponds under the $\Inv-\Pol$ Galois connection to a minimal relational clone, and every minimal relational clone can be generated by a single relation, of one of several special forms. Zhuk is very familiar with the theory of relational clones, so he was aware of Rosenberg's Completeness Theorem \cite{rosenberg-completeness} (see \cite{pinsker-rosenberg} or chapter II.6 of \cite{lau-clone-theory} for alternate expositions), which completely classifies the special relations which correspond to maximal clones into six different types.

Zhuk then considered each of the types of relations from Rosenberg's classification, and investigated which of them might be preserved by the polynomial clone of a Taylor algebra, and what the existence of such a relation implies about the structure of the Taylor algebra. The most interesting case is the case of the relations known as \emph{central relations}.

\begin{defn} A relation $\RR \le \bA^n$ is \emph{central} if it has the following properties:
\begin{itemize}
\item $\RR$ is symmetric under permuting its coordinates,
\item $\RR$ contains every tuple which has any pair of equal coordinates, and
\item the set $\bC \le \bA$ defined by
\[
\bC = \{c \in \bA \mid \forall a_2, ..., a_n \in \bA, (c,a_2,...,a_n) \in \RR\}
\]
is not empty and is not equal to $\bA$.
\end{itemize}
The set $\bC$ is known as the \emph{center} of the central relation $\RR$.
\end{defn}

Since relations of high arity are hard to think about, Zhuk simplifies this to a special type of binary relation on $\bA \times \bB$, where $\bB$ is secretly taken to be $\bA^{n-1}$. To see that this step doesn't lose anything essential, we use the following fact about absorbing subalgebras of powers.

\begin{prop}\label{prop-absorption-in-powers} Suppose that $\bA$ is idempotent and that $\bA^k$ has a proper absorbing subalgebra for some $k$. Then $\bA$ has a proper absorbing subalgebra.

In fact, this holds for any absorption concept $\lhd_X$ which is compatible with pp-formulas.
\end{prop}
\begin{proof} We induct on $k$. Suppose that $\bB \lhd_X \bA^k$. If $\pi_1(\bB) \ne \bA$ then $\pi_1(\bB) \lhd_X \bA$ and we are done, otherwise since $\bB \ne \bA^k$ there must exist some $a \in \bA$ such that $\pi_{[k]\setminus\{1\}}(\bB \cap \{a\} \times \bA^{k-1}) \ne \bA^{k-1}$. Since $\lhd_X$ is compatible with pp-formulas and $\{a\} \le \bA$ by the idempotence of $\bA$, we have
\[
\pi_{[k]\setminus\{1\}}(\bB \cap \{a\} \times \bA^{k-1}) \lhd_X \bA^{k-1},
\]
so we can apply the induction hypothesis.
\end{proof}

With this in mind, it's natural to restrict our attention to binary relations $\RR \le \bA \times \bB$ which have a nontrivial proper ``left center'', and to try to use them to produce an absorbing subalgebra inside either $\bA$ or $\bB$.

\begin{defn} If $\RR \le_{sd} \bA \times \bB$ is subdirect and $\bB$ is finite and idempotent, then the \emph{left center} of $\RR$ is the subalgebra $\bC \le \bA$ defined by
\[
\bC = \{c \in \bA \mid \forall b \in \bB, (c,b) \in \RR\}.
\]
The \emph{right center} of a subdirect binary relation is defined similarly (so the right center of $\RR$ is the left center of $\RR^-$, and is a subalgebra of $\bB$).
\end{defn}

To see that the left center $\bC$ is automatically a subalgebra of $\bA$, note that it can be defined by the following pp-formula:
\[
c \in \bC \iff \bigwedge_{b \in \bB} \exists x (x \in \{b\} \wedge (c,x) \in \RR).
\]

In order to do anything useful with such a binary relation, we will need to assume that $\bB$ is Taylor. We will attempt to exploit the Taylor term to produce binary absorption on $\bB$, using the following lemma.

\begin{lem}\label{lem-projective-bin-absorbing} Suppose $\bB \le \bA$ and that there is an idempotent term $t \in \Clo_k(\bA)$ with the following two properties:
\begin{itemize}
\item $t$ satisfies an identity of the form $t(x,u_2,...,u_k) \approx t(y,v_2,...,v_k)$, where each $u_i,v_i \in \{x,y\}$, and
\item $t(\bB, \bA, ..., \bA) \subseteq \bB$.
\end{itemize}
Then $\bB$ absorbs $\bA$ with respect to some idempotent binary operation $f$.
\end{lem}
\begin{proof} To make the notation more clear, we treat each $u_i,v_i$ as a binary function, with $u_i = u_i(x,y)$ and $v_i = v_i(x,y)$. Define $f(x,y)$ by
\[
f(x,y) \coloneqq t(x,u_2(x,y),...,u_k(x,y)) \approx t(y,v_2(x,y),...,v_k(x,y)).
\]
Then for any $a \in \bA$ and $b \in \bB$, we have
\[
f(a,b) = t(b,v_2(a,b),...,v_k(a,b)) \in t(\bB,\bA,...,\bA) \subseteq \bB,
\]
and
\[
f(b,a) = t(b,u_2(b,a),...,u_k(b,a)) \in t(\bB,\bA,...,\bA) \subseteq \bB.\qedhere
\]
\end{proof}

\begin{thm}[Zhuk \cite{zhuk-dichotomy}]\label{thm-center-implies-absorbing} Suppose that $\bA, \bB$ are finite idempotent algebras, and that there is a term $t$ which is Taylor on $\bB$. If $\RR \le_{sd} \bA \times \bB$ is subdirect and has a nontrivial left center $\bC$, then either $\bB$ has a proper binary absorbing subalgebra, or $\bC$ absorbs $\bA$ with respect to the term $t*\cdots *t$, with $|\bB|-1$ copies of $t$.
\end{thm}
\begin{proof} Suppose $t$ has arity $k$. We will show that if $\bB$ has no proper binary absorbing subalgebra, then for any $a \in \bA \setminus \bC$ and for any $c_1, ..., c_k \in \bC$ and any $i \le k$, the value
\[
t(c_1, ..., c_{i-1}, a, c_{i+1}, ..., c_k)
\]
is ``closer'' to being in $\bC$ than $a$ is. To make this precise, we measure how close an element $a$ is to being in $\bC$ by looking at the size of the set
\[
a+\RR = \pi_2(\RR \cap \{a\}\times \bB).
\]
By the definition of $\bC$, we have $|a + \RR| = |\bB|$ if and only if $a \in \bC$.

Since $\RR$ is preserved by $t$, we have
\[
t(c_1, ..., a, ..., c_k) + \RR \supseteq t(c_1+\RR, ..., a+\RR, ..., c_k+\RR) = t(\bB, ..., a + \RR, ..., \bB).
\]
Since $t$ is idempotent, the right hand side of the above must contain $a+\RR$, and if it is equal to $a+\RR$ then we can apply the previous lemma (since $t$ is Taylor) to see that $a+\RR$ is a binary absorbing subalgebra of $\bB$. Thus if $a \not\in \bC$, then either $a + \RR$ is a proper binary absorbing subalgebra of $\bB$, or else
\[
|t(c_1, ..., a, ..., c_k) + \RR| > |a+\RR|.\qedhere
\]
\end{proof}

\begin{rem} There is a variant of Theorem \ref{thm-center-implies-absorbing} which applies even if $\bB$ is not a Taylor algebra, from \cite{zhuk-strong}. Defining the concept of a projective subalgebra as in Definition \ref{defn-projective}, a slight modification of the proof of Theorem \ref{thm-center-implies-absorbing} shows that if $\RR \le_{sd} \bA \times \bB$ is subdirect and has a nontrivial left center $\bC$, then either $\bC$ absorbs $\bA$ or $\bB$ has a proper projective subalgebra $\PP < \bB$. Note that Lemma \ref{lem-projective-bin-absorbing} implies that any projective subalgebra of a Taylor algebra is automatically binary absorbing.
\end{rem}

Keeping the same setup, the left center $\bC$ has an additional nice property, which is much stronger than it looks.

\begin{thm}[Zhuk \cite{zhuk-dichotomy}]\label{thm-center-implies-central} Suppose $\bA, \bB$ are finite idempotent algebras. If $\RR \le_{sd} \bA \times \bB$ has a left center $\bC$ and $\bB$ has no proper binary absorbing subalgebras, then for any $a \in \bA$ we have
\[
a \not\in \bC \;\; \implies \;\; \begin{bmatrix}a\\ a\end{bmatrix} \not\in \Sg_{\bA^2}\left\{\begin{bmatrix}a\\ \bC\end{bmatrix}, \begin{bmatrix}\bC\\ \bC\end{bmatrix}, \begin{bmatrix}\bC\\ a\end{bmatrix}\right\}.
\]
\end{thm}
\begin{proof} Suppose otherwise. Then there are $i,j$ and $c_1, ..., c_i, c_j', ..., c_n' \in \bC$ with $j \le i+1$ and a term $t$ of arity $n$ such that
\[
\begin{bmatrix}a\\ a\end{bmatrix} = t\left(\begin{bmatrix}a\\ c_1\end{bmatrix}, \cdots, \begin{bmatrix}a\\ c_{j-1}\end{bmatrix}, \begin{bmatrix}c_j'\\ c_j\end{bmatrix}, \cdots, \begin{bmatrix}c_i'\\ c_i\end{bmatrix}, \begin{bmatrix}c_{i+1}'\\ a\end{bmatrix}, ..., \begin{bmatrix}c_n'\\ a\end{bmatrix}\right).
\]
Looking at the neighbors via $\RR$, we have
\[
\begin{bmatrix}a+\RR\\ a+\RR\end{bmatrix} \supseteq t\left(\begin{bmatrix} a+\RR & \cdots & a+\RR & \bB & \cdots & \bB & \bB & \cdots & \bB\\ \bB & \cdots & \bB & \bB & \cdots & \bB & a+\RR & \cdots & a+\RR \end{bmatrix}\right).
\]
Thus $a+\RR$ absorbs $\bB$ with respect to the binary term
\[
f(x,y) \coloneqq t(x,...,x,y,...,y)
\]
as long as the number of $x$s is between $j-1$ and $i$.
\end{proof}

We can combine the previous two results about left centers to define a new type of absorption. We won't need the full power of the previous result, and instead will use a slightly weaker property.

\begin{defn} We say that $\bC$ \emph{centrally absorbs} $\bA$, written $\bC \lhd_Z \bA$, if the following two properties hold:
\begin{itemize}
\item $\bC \lhd \bA$, and
\item for any $a \not\in \bC$, we have $\begin{bmatrix}a\\ a\end{bmatrix} \not\in \Sg_{\bA^2}\left\{\begin{bmatrix}a\\ \bC\end{bmatrix}, \begin{bmatrix}\bC\\ a\end{bmatrix}\right\}$.
\end{itemize}
\end{defn}

\begin{cor}\label{zhuk-center} Suppose $\bA, \bB$ are finite and idempotent. If $\RR \le_{sd} \bA \times \bB$ has left center $\bC$ and $\bB$ is Taylor and binary absorption free, then $\bC \lhd_Z \bA$.
\end{cor}
\begin{proof} By Theorem \ref{thm-center-implies-absorbing}, $\bC$ absorbs $\bA$, and then by Theorem \ref{thm-center-implies-central} the absorption is central.
\end{proof}

There is an unfortunate naming collision between the centers considered here, and the centers considered in commutator theory. Generally it should be clear from context which sort of center is meant. (I have proposed the alternate name \emph{stable absorption} instead of central absorption, but it seems unlikely to catch on.)

The key fact about central absorption that makes it so much more powerful than ordinary absorption is the following doubling trick due to Zhuk and Kozik.

\begin{lem}[Essential doubling trick \cite{zhuk-dichotomy}]\label{lem-essential-doubling} Suppose that $\RR \le \bA_0 \times \cdots \times \bA_{n+1}$ is $(\bC,\bB_1, ..., \bB_n, \bC')$-essential, with $\bC' \lhd_Z \bA_{n+1}$ and $\bA_{n+1}$ finite and idempotent. Then there is a relation
\[
\RR' \le \bA_0 \times \cdots \times \bA_n \times \bA_n \times \cdots \times \bA_0
\]
which is $(\bC,\bB_1, ..., \bB_n, \bB_n, ..., \bB_1, \bC)$-essential.
\end{lem}
\begin{proof} Suppose $\RR$ is chosen such that, subject to satisfying the assumptions of the lemma, the subalgebra $\bB' \le \bA_{n+1}$ defined by
\[
\bB' = \pi_{n+1}(\RR \cap \bC\times \bB_1 \times \cdots \times \bB_n \times \bA_{n+1})
\]
is as small as possible. Note that $\bB'$ is necessarily nonempty and disjoint from $\bC'$ if $\RR$ is $(\bC,\bB_1, ..., \bB_n, \bC')$-essential.

Since we may shrink $\RR$ to the subalgebra generated by any collection of tuples witnessing $\RR \cap (\bC\times\cdots \times \bA_i \times \cdots \times \bC') \ne \emptyset$ for all $i$ from $0$ to $n+1$, we see that
\[
b,b' \in \bB' \;\; \implies \;\; b' \in \Sg_{\bA_{n+1}}(\bC' \cup \{b\}).
\]
In particular, if we pick some $b \in \bB'$ and define the symmetric binary relation $\bS \le \bA_{n+1}\times \bA_{n+1}$ by
\[
\bS = \Sg_{\bA_{n+1}^2}\left\{\begin{bmatrix} b\\ \bC' \end{bmatrix}, \begin{bmatrix} \bC' \\ b \end{bmatrix}\right\},
\]
then $\pi_1(\bS) \supseteq \bB'$.

We now define the relation $\RR'$ by
\[
(x_0, ..., x_n, y_n, ..., y_0) \in \RR' \; \iff \; \exists x_{n+1},y_{n+1}\ (x_0, ..., x_{n+1}) \in \RR \wedge (x_{n+1},y_{n+1}) \in \bS \wedge (y_0, ..., y_{n+1}) \in \RR.
\]
To see that
\[
\RR' \cap \bC\times\bB_1\times \cdots \times \bA_i \times \cdots \times \bB_n \times \bB_n \times \cdots \times \bB_1 \times \bC \ne \emptyset
\]
for any $0 \le i \le n$, we choose $(x_0, ..., x_{n+1}) \in \RR \cap (\bC\times\cdots \times \bA_i \times \cdots \times \bC')$ and choose $(y_0, ..., y_{n+1}) \in \RR \cap \bC\times \bB_1 \times \cdots \times \bB_n \times \{b\}$, which is possible since $b \in \bB'$ and $\bC' \times \{b\} \subseteq \bS$. We can check that
\[
\RR' \cap \bC\times\bB_1\times \cdots \times \bB_n \times \bB_n \times \cdots \times \bA_i \times \cdots \times \bB_1 \times \bC \ne \emptyset
\]
for $0 \le i \le n$ similarly, by interchanging the roles of the $x_i$s and $y_i$s.

To finish, we just need to check that
\[
\RR' \cap \bC\times\bB_1\times \cdots \times \bB_n \times \bB_n \times \cdots \times \bB_1 \times \bC = \emptyset,
\]
or equivalently, that
\[
\bS \cap \bB'\times \bB' = \emptyset.
\]

So suppose for contradiction that there are $b',b'' \in \bB'$ with $(b',b'') \in \bS$. Since
\[
b \in \Sg(\bC' \cup \{b'\}) \subseteq \bB' - \bS,
\]
we see that there is some $b''' \in \bB'$ such that $(b,b''') \in \bS$. But then we have
\[
\{b\} + \bS \supseteq \Sg(\bC' \cup \{b'''\}) \supseteq \bB',
\]
so $(b,b) \in \bS$, contradicting our assumption that $\bC' \lhd_Z \bA_{n+1}$.
\end{proof}

\begin{cor}\label{center-ternary} If $\bC \lhd_Z \bA$ and $\bA$ is finite and idempotent, then $\bC$ absorbs $\bA$ with respect to some ternary term.
\end{cor}
\begin{proof} If $\bC$ does not absorb $\bA$ with respect to any ternary term, then by Theorem \ref{absorption-essential} there is some ternary $\bC$-essential relation $\RR \le \bA^3$. By repeatedly applying the doubling trick, we see that there exists some $\bC$-essential relation of arity $2+2^k$ for every $k \ge 0$, so $\bC$ can't absorb $\bA$ with respect to a term of any arity, contradicting the assumption $\bC \lhd_Z \bA$.
\end{proof}

\begin{cor}\label{cor-essential-sandwich} If $\bC_1 \lhd_Z \bA_1, \bB_2 \lhd \bA_2$, and $\bC_3 \lhd_Z \bA_3$ with $\bA_i$ finite and idempotent, then no $(\bC_1, \bB_2, \bC_3)$-essential relation can exist.
\end{cor}
\begin{proof} If a $(\bC_1, \bB_2, \bC_3)$-essential relation exists, then by repeatedly applying the doubling trick we can find $(\bC_1, \bB_2, ..., \bB_2, \bC_1)$-essential relations of arbitrarily high arity. By forcing the first and last coordinates to be in $\bC_1$ and existentially projecting, we see that there are $\bB_2$-essential relations of arbitrarily high arity, which contradicts the assumption $\bB_2 \lhd \bA_2$.
\end{proof}

\begin{cor} If $\bA_i$ are finite and idempotent, $\bC_i \lhd \bA_i$ for all $i$ and for all but at most one $i$ we have $\bC_i \lhd_Z \bA_i$, then for any relation $\RR \le \bA_1 \times \cdots \times \bA_n$ such that $\pi_{i,j}(\RR) \cap \bC_i \times \bC_j \ne \emptyset$ for all $i,j$, we have
\[
\RR \cap \bC_1 \times \cdots \times \bC_n \ne \emptyset.
\]
\end{cor}
\begin{proof} We show by induction on $|I|$ that for all $I \subseteq [n]$ we have
\[
\pi_I(\RR) \cap \prod_{i \in I} \bC_i \ne \emptyset.
\]
The base case $|I| \le 2$ is our assumption. For $|I| \ge 3$, pick $i,j,k \in I$ distinct. By the induction hypothesis, there are tuples $x_i, x_j, x_k \in \RR$ such that $\pi_{I\setminus\{i\}}(x_i) \in \prod_{i' \in I \setminus\{i\}} \bC_{i'}$, and similarly for $x_j,x_k$.

Now consider the subalgebra of $\bA_i\times \bA_j \times \bA_k$ generated by $\pi_{i,j,k}(x_i), \pi_{i,j,k}(x_j), \pi_{i,j,k}(x_k)$. Since this subalgebra can't be a $(\bC_i,\bC_j,\bC_k)$-essential relation (since at least two of $\bC_i,\bC_j,\bC_k$ are centrally absorbing and the third is absorbing), it must contain an element of $\bC_i\times\bC_j\times\bC_k$. Thus there is some $x \in \Sg\{x_i,x_j,x_k\}$ such that
\[
\pi_{i,j,k}(x) \in \bC_i\times\bC_j\times\bC_k,
\]
and this $x$ automatically satisfies
\[
\pi_{I\setminus\{i,j,k\}}(x) \in \prod_{i' \in I\setminus\{i,j,k\}} \bC_{i'}
\]
since each of $x_i, x_j, x_k$ do, which completes the inductive step.
\end{proof}

\begin{cor}\label{cor-central-ternary} If $\bA$ is finite and idempotent, then there is a ternary term $t \in \Clo_3(\bA)$ such that for all finite $\bB \in HSP(\bA)$ and each $\bC \lhd_Z \bB$, $\bC$ absorbs $\bB$ with respect to the term $t$.
\end{cor}
\begin{proof} For any finite collection of pairs $\bC_i \lhd_Z \bB_i \in HSP(\bA)$, we can apply the previous corollary to find a term $t \in \Clo_3(\bA)$ which simultaneously witnesses all $\bC_i \lhd \bB_i$. Since there are only finitely many ternary terms $t$ of $\bA$, some $t$ must work for all pairs $\bC \lhd_Z \bB \in HSP(\bA)$.
\end{proof}

\begin{rem} Since each case of these corollaries is a concrete statement about the elements of a finite collection of relations, by Theorem \ref{thm-finitely-related-jonsson-absorption} they should all still apply if we modify the definition of central absorption by replacing ``absorption'' by ``J\'onsson absorption''. In fact, the essential relations constructed through repeated application of Lemma \ref{lem-essential-doubling}, starting from a ternary $(\bC_1, \bB_2, \bC_3)$-essential relation, actually end up being definable by ``comb fomulas'' (using the language of Subsection \ref{ss-finitely-related-jonsson}). Furthermore, the assumptions of Corollary \ref{cor-essential-sandwich} can be weakened to just $\bB_2 \lhd_J \bA_2$ and $a \not\in \bC_i \implies (a,a) \not\in \Sg_{\bA_i^2}(\{a\}\times\bC_i \cup \bC_i\times\{a\})$.
\end{rem}

Central absorption turns out to be a good absorption concept (in the sense of the previous section), as long as we restrict ourselves to finite idempotent algebras. Unlike previous absorption concepts, in this case it is not so easy to see that $\lhd_Z$ is compatible with pp-formulas. For this, we need to consider the basic types of pp-formulas separately. The hardest case is the case of projections.

\begin{prop}\label{prop-central-abs-quotient} If $\bC \lhd_Z \bA$ with $\bA$ finite and idempotent, and if there is a surjective homomorphism $\pi : \bA \twoheadrightarrow \bB$, then $\pi(\bC) \lhd_Z \bB$.
\end{prop}
\begin{proof} Suppose there is some $b \in \bB\setminus\pi(\bC)$ such that $(b,b) \in \Sg(\pi(\bC)\times \{b\} \cup \{b\} \times \pi(\bC))$. Choose $a \in \pi^{-1}(b)$ such that the subalgebra $\Sg(\bC\cup \{a\})$ is as small as possible. Set
\[
\bS = \Sg_{\bA^2}\left\{\begin{bmatrix}a\\ \bC\end{bmatrix}, \begin{bmatrix} \bC \\ a\end{bmatrix}\right\}.
\]
By the choice of $b$, there exist $a',a'' \in \bA$ such that $(a',a'') \in \bS$ and $\pi(a') = \pi(a'') = b$. By the choice of $a$, we have $a \in \Sg(\bC \cup \{a''\})$. Thus we have
\[
\Sg\{a,a'\} + \bS \supseteq \Sg(\bC \cup \{a''\}) \supseteq \{a\},
\]
so there is some $a''' \in \Sg\{a,a'\}$ with $(a''',a) \in \bS$, and by idempotence we have $\pi(a''') = b$, so $a \in \Sg(\bC \cup \{a'''\})$. By idempotence we have $\{a\} \le \bA$, so
\[
\{a\} - \bS \supseteq \Sg(\bC \cup \{a'''\}) \supseteq \{a\},
\]
so $(a,a) \in \bS$, which contradicts the assumption $\bC \lhd_Z \bA$.
\end{proof}

\begin{prop} If $\bC \lhd_Z \bB \lhd_Z \bA$, then $\bC \lhd_Z \bA$. As a consequence, if $\bC_i \lhd_Z \bB_i \le \bA$, then $\bC_1 \cap \bC_2 \lhd_Z \bB_1 \cap \bB_2$.
\end{prop}
\begin{proof} Suppose there is some $a \in \bA$ such that $(a,a) \in \Sg(\bC\times \{a\} \cup \{a\}\times \bC)$. Since $\bC \le \bB$ and $\bB \lhd_Z \bA$, we must have $a \in \bB$. Then since $\bC \lhd_Z \bB$, we must have $a \in \bC$. Thus $\bC \lhd_Z \bA$.

For the second statement, note that $\bC_2 \lhd_Z \bB_2$ implies $\bC_1 \cap \bC_2 \lhd_Z \bC_1 \cap \bB_2$ and $\bC_1 \lhd_Z \bB_1$ implies $\bC_1\cap \bB_2 \lhd_Z \bB_1\cap\bB_2$.
\end{proof}

\begin{prop} If $\bC_1 \lhd_Z \bA_1$, then $\bC_1 \times \bA_2 \lhd_Z \bA_1 \times \bA_2$.
\end{prop}

Putting these three results together, we see that central absorption is a good absorption concept.

\begin{prop} The absorption concept $\lhd_Z$, restricted to finite idempotent algebras, is compatible with pp-formulas, is transitively closed, and transfers connectivity.
\end{prop}

\begin{rem}\label{rem-bin-abs} Annoyingly, binary absorption fails to be transitively closed or compatible with pp-formulas (the intersection of two binary absorbing subalgebras might not be binary absorbing). However, if we restrict ourselves to finite idempotent algebras which are \emph{prepared}, that is, such that $(b,b) \in \Sg\{(a,b),(b,a)\}$ implies that $\{a,b\}$ is a semilattice subalgebra with absorbing element $b$, then binary absorption becomes compatible with pp-formulas and transitively closed (see Proposition \ref{prepared-bin-abs}).
\end{rem}

In some cases central absorption implies binary absorption. To describe a criterion for when this happens, we will exploit partial semilattice operations.

\begin{prop}\label{prop-central-closed} Suppose that $\bC \lhd_Z \bA$ and that $s$ is any partial semilattice operation. Then $s(\bC,\bA) \subseteq \bC$.
\end{prop}
\begin{proof} Suppose $c \in \bC$ and $a \in \bA$, and let $b = s(c,a)$. Then $s(c,b) = s(b,c) = b$ by the defining property of partial semilattice operations, so $(b,b) \in \Sg(\{b\}\times \bC \cup \bC \times \{b\})$. Thus by the definition of central absorption, we have $b \in \bC$, that is, $s(c,a) \in \bC$.
\end{proof}

\begin{prop}\label{bin-central-criterion} Suppose that $\bC \lhd_Z \bA$ in a finite idempotent algebra $\bA$, or just that $s(\bC,\bA) \subseteq \bC$ for all partial semilattice terms $s \in \Clo(\bA)$. Then the following are equivalent:
\begin{itemize}
\item[(a)] $\bC$ binary absorbs $\bA$,
\item[(b)] for all $a \in \bA\setminus \bC$ and all $c \in \bC$, the subalgebra $\Sg\{a,c\}$ has a proper binary absorbing subalgebra,
\item[(c)] for all $a \in \bA$ and all $c \in \bC$, there is a sequence of elements $a = a_0, a_1, ..., a_n \in \Sg\{a,c\}$ with $a_n \in \bC$ such that $(a_i,a_i) \in \Sg\{(a_{i-1},a_i),(a_i,a_{i-1})\}$ for all $i$.
\end{itemize}
If $\bA$ is prepared, then the third condition is equivalent to the assumption that for all $a$ and for all $c \in \bC$, the subalgebra $\Sg\{a,c\}$ contains a directed path from $a$ to $\bC$.
\end{prop}
\begin{proof} To see that (a) implies (b), note that $\bC \lhd_{bin} \bA$ implies that $\bC \cap \Sg\{a,c\} \lhd_{bin} \Sg\{a,c\}$. To see that (b) implies (c), we induct on the size of $\Sg\{a,c\}$. Let $\bB$ be a proper binary absorbing subalgebra of $\Sg\{a,c\}$, and let $s$ be a partial semilattice term that witnesses this absorption (such an $s$ exists by Proposition \ref{bin-abs-semi}). Then for any $b \in \bB$ we have $s(a,b) \in \bB$, and if we take $a_1 = s(a,b)$ then $(a_1,a_1) \in \Sg\{(a,a_1),(a_1,a)\}$. Let $c_1 = s(c,b)$, then $c_1 \in \bB \cap \bC$, and so $\Sg\{a_1,c_1\} \subseteq \bB < \Sg\{a,c\}$, so by the inductive hypothesis we can complete this to a sequence $a_1, ..., a_n \in \Sg\{a_1,c_1\}$ as in (c).

Now suppose that (c) holds. For each $a,c$ with $c \in \bC$, we will construct a binary function $f_{ac}$ such that $f_{ac}(a,c) \in \bC$ and $f_{ac}(\bC,\bA) \subseteq \bC$. Then by cyclically composing the functions $f_{ac}$ together, we can produce a binary term which absorbs $\bC$. To construct $f_{ac}$, we pick a sequence of partial semilattice terms $s_i$ such that $s_i(a_{i-1},a_i) = a_i$ as well as binary terms $t_i$ such that $t_i(a,c) = a_i$. We set
\[
f_{ac}(x,y) \coloneqq s_n(\cdots s_2(s_1(x,t_1(x,y)),t_2(x,y))\cdots,t_n(x,y)).
\]
Then we have
\[
f_{ac}(a,c) = s_n(\cdots s_2(s_1(a,a_1),a_2)\cdots,a_n) = a_n \in \bC
\]
and
\[
f_{ac}(\bC,\bA) \subseteq s_n(\cdots s_2(s_1(\bC,\bA),\bA)\cdots,\bA) \subseteq \bC,
\]
as required.
\end{proof}



\section{Absorption Theorem and Loop Lemma: Binary relations in Taylor algebras}

In this section we'll go over two of the main results from Barto and Kozik's paper \cite{cyclic} about absorption, known as the ``absorption theorem'' and the ``loop lemma''. The first of these results can be used to constrain the possible subdirect binary relations in simple absorption free algebras, while the second result makes no direct mention of absorption, but combines the theory of absorbing subalgebras with an elementary argument in the absorption free case to give a criterion for a subdirect binary relation to intersect the diagonal.

The loop lemma was originally introduced in order to settle a special case of the dichotomy problem, where the template structure $\fA$ consists of a set together with a single subdirect binary relation (considered as a directed graph). As a bonus, the loop lemma easily implies the existence of a Taylor term operation of a special form, known as a \emph{Siggers} operation, named after the first person to notice that such special Taylor terms exist in the finite case \cite{siggers-original} (this result was quickly refined, after the initial discovery: see \cite{optimal-taylor} for the paper which introduced the $4$-ary operations which are now commonly known as Siggers operations).

Here is a strong form of the absorption theorem, stated in terms of Zhuk's centers.

\begin{thm}[Absorption Theorem \cite{cyclic}]\label{absorption-theorem} If $\RR \le_{sd} \bA \times \bB$ is a subdirect binary relation and $\bA,\bB$ are finite idempotent Taylor algebras, and if $\RR$ is \emph{linked}, then either
\begin{itemize}
\item $\RR = \bA \times \bB$,
\item $\bA$ has a proper binary absorbing or centrally absorbing subalgebra, or
\item $\bB$ has a proper subalgebra which is both binary absorbing and centrally absorbing.
\end{itemize}
\end{thm}

The absorption theorem can be viewed as a strengthening of Zhuk's results about central relations: as we will see, it actually follows from Zhuk's result by applying a few simple tricks. First we will bootstrap to the case of a subdirect relation $\RR$ such that $\RR \circ \RR^- = \bA\times \bA$.

\begin{lem} If $\RR \le_{sd} \bA \times \bB$ is a subdirect binary relation and $\bA,\bB$ are finite idempotent Taylor algebras, and if $\RR \circ \RR^- = \bA \times \bA$, then either
\begin{itemize}
\item there is some $b \in \bB$ such that $\bA \times \{b\} \subseteq \RR$,
\item $\bA$ has a proper binary absorbing subalgebra, or
\item every element of $\bA$ is contained in a proper centrally absorbing subalgebra.
\end{itemize}
\end{lem}
\begin{proof} Suppose that there is no $b \in \bB$ with $\bA \times \{b\} \subseteq \RR$ and that $\bA$ is binary absorption free, and choose any $a \in \bA$. Choose a sequence of subalgebras $\{a\} + \RR = \bD_0 \ge \bD_1 \ge \cdots \ge \bD_n$ such that each $\bD_{i+1}$ is a proper binary absorbing subalgebra of $\bD_i$ and such that $\bD_n$ has no proper binary absorbing subalgebras. We will first show that $\bD_n - \RR = \bA$, and then we will apply Zhuk's result (Corollary \ref{zhuk-center}) to the binary relation $\RR \cap (\bA \times \bD_n)$.

We will show that $\bD_i - \RR = \bA$ for each $i$, by induction on $i$. Note that $\bD_0 - \RR = \{a\} + \RR - \RR = \bA$ by the assumption $\RR\circ \RR^- = \bA\times \bA$. For the inductive step, note that since $\bD_{i+1} \lhd_{bin} \bD_i$, we have
\[
\bD_{i+1} - \RR \lhd_{bin} \bD_i - \RR = \bA,
\]
so we must have $\bD_{i+1} - \RR = \bA$ since $\bA$ has no proper binary absorbing subalgebra.

If we set $\RR' = \RR \cap (\bA \times \bD_n)$, then we have
\[
\{a\} \times \bD_n \subseteq \RR' \le_{sd} \bA \times \bD_n.
\]
Thus the left center $\bC$ of $\RR'$ contains $a$. Since $\bD_n$ is binary absorption free, we see that $\bC$ centrally absorbs $\bA$ by Corollary \ref{zhuk-center}. If $\bC = \bA$, then $\bA\times\bD_n \subseteq\RR$, contradicting the assumption that there is no $b \in \bB$ with $\bA \times \{b\} \subseteq \RR$.
\end{proof}

In the case where $\bA$ has no proper binary absorbing or centrally absorbing subalgebra and $\RR$ has a nontrivial right center, we will use the criterion developed in Proposition \ref{bin-central-criterion} to show that the right center of $\RR$ must actually be a binary absorbing subalgebra of $\bB$.

\begin{lem} If $\RR \le_{sd} \bA \times \bB$ is a subdirect binary relation and $\bA,\bB$ are finite idempotent Taylor algebras, and if there is some $b \in \bB$ such that $\bA \times \{b\} \subseteq \RR$, then either
\begin{itemize}
\item $\RR = \bA \times \bB$,
\item $\bA$ has a proper binary absorbing or centrally absorbing subalgebra, or
\item the right center of $\RR$ is a proper binary absorbing subalgebra of $\bB$.
\end{itemize}
\end{lem}
\begin{proof} Let $\bC \le \bB$ be the right center of $\RR$. By Corollary \ref{zhuk-center}, if $\bA$ has no proper binary absorbing subalgebra then we have $\bC \lhd_Z \bB$. If $\bC$ is not a binary absorbing subalgebra of $\bB$, then by Proposition \ref{bin-central-criterion} there must be some $b \in \bB\setminus \bC$ and $c \in \bC$ such that $\Sg\{b,c\}$ has no proper binary absorbing subalgebra.

Since $\RR$ is subdirect, there is some $a \in \bA$ such that $(a,b) \in \RR$. Since $c$ is in the right center of $\RR$, we also have $\bA \times \{c\} \subseteq \RR$. Thus if we set $\RR' = \RR \cap (\bA \times \Sg\{b,c\})$, then we have
\[
\{a\} \times \Sg\{b,c\} \subseteq \RR' \le_{sd} \bA \times \Sg\{b,c\}
\]
Since $b$ is \emph{not} in the right center of $\RR$, the left center of $\RR'$ is a proper subalgebra of $\bA$. Then since $\Sg\{b,c\}$ has no proper binary absorbing subalgebra, Corollary \ref{zhuk-center} shows that the left center of $\RR'$ is a proper centrally absorbing subalgebra of $\bA$.
\end{proof}

\begin{proof}[Proof of the Absorption Theorem] Let $\bS = \RR \circ \RR^- \le_{sd} \bA \times \bA$. If $\bS = \bA \times \bA$, we may apply the lemmas to see that either $\bA$ has a proper binary absorbing or centrally absorbing subalgebra, or that $\bB$ has a proper subalgebra which is both binary absorbing and centrally absorbing. Otherwise, by the fact that $\RR$ is linked and the finiteness of $\bA$ there must be some minimal $k > 1$ such that $\bS^{\circ k} = \bA\times \bA$. Then we can apply the lemmas to $\bS^{\circ(k-1)}$ to see that $\bA$ must have either a binary absorbing or centrally absorbing subalgebra.
\end{proof}

\begin{cor}\label{absorption-free-binary} Let $\bA, \bB$ be finite idempotent Taylor algebras with no proper binary or centrally absorbing subalgebras. If $\bB$ is simple, then every subdirect binary relation $\RR \le_{sd} \bA\times \bB$ is either the full relation or the graph of a surjective homomorphism $\bA \twoheadrightarrow \bB$.
\end{cor}
\begin{proof} Since $\bB$ is simple, the linking congruence of $\RR$ on $\bB$ is either trivial or is full. If the linking congruence of $\RR$ on $\bB$ is trivial, then $\RR$ must be the graph of a surjective homomorphism $\bA \twoheadrightarrow \bB$. Otherwise, $\RR$ is linked, so we can apply the Absorption Theorem \ref{absorption-theorem} to see that $\RR = \bA \times \bB$.
\end{proof}

Next we switch our focus to subdirect relations $\RR \le_{sd} \bA\times\bA$. In this case, it is often appropriate to think of $\RR$ as a digraph on the vertex set $\bA$, and we can ask questions about whether $\RR$ (viewed as a digraph) is weakly connected, strongly connected, whether it contains any loops, etc. To be precise, the associated digraph is the relational structure $\fR = (A,R)$, where $A$ is the underlying set of $\bA$ and $R \subseteq A\times A$ is the underlying set of $\RR$ (often I abuse notation and write $\fR = (\bA,\RR)$ instead of explicitly replacing $\bA,\RR$ with their underlying sets).

\begin{rem} Note that if $\RR \le \bA\times \bA$ and $\bS \le \bB\times \bB$ are subpowers of $\bA, \bB$, then a homomorphism $\RR \rightarrow \bS$ and a homomorphism $\fR \rightarrow \fS$ of the associated digraphs $\fR = (\bA,\RR), \fS = (\bB,\bS)$ are completely different things! The first is a homomorphism of algebraic structures, and doesn't depend on how $\RR,\bS$ are represented as collections of ordered pairs of elements in $\bA$ or $\bB$ (but does depend on how the algebraic operations behave). The second is a digraph homomorphism, which ignores the algebraic structure, and is completely determined by a map $A \rightarrow B$ of the underlying sets of $\bA,\bB$ which is compatible with the digraph structures $\RR,\bS$.
\end{rem}

In the context of digraphs, the case of a \emph{subdirect} relation $\RR \le_{sd} \bA \times \bA$ is actually rather special. The assumption $\pi_1(\RR) = \bA$ means that every vertex of the digraph $\fR$ has outdegree at least one, and the assumption $\pi_2(\RR) = \bA$ means that every vertex of $\RR$ has indegree at least one.

\begin{defn} A digraph $\fD = (V,E)$ is called \emph{smooth} if every vertex of $\fD$ has indegree at least one and outdegree at least one. Note that this is equivalent to the relation $E \subseteq V\times V$ being subdirect.
\end{defn}

If a digraph is not smooth, it is often desirable to find a smooth digraph within it. The natural thing to do is to simply prune all of the vertice with indegree $0$ or outdegree $0$. Unfortunately, after this pruning step we may find ourselves with more vertices that need to be pruned, and so on - possibly ending up with no vertices at all! For instance, this actually occurs if our initial digraph is a finite directed path. Additionally, it may not be clear that these pruning operations are compatible with the algebraic structures which we started with. Luckily, there is a standard way to describe the result of this pruning process via a primitive positive formula, as well as a simple criterion for when the pruned digraph will be nonempty.

\begin{prop} If $\fD = (V,E)$ is a digraph, then the largest smooth digraph $\fD_{sm}$ which is contained in $\fD$ is exactly the set of vertices $v$ of $\fD$ such that there exists a bi-infinite directed walk through $v$. If $\fD$ is finite, with $n$ vertices, then the vertex set of $\fD_{sm}$ may be defined by the pp-formula
\[
v \in \fD_{sm} \iff \exists v_{-n}, ..., v_n\ (v_0 = v) \wedge \bigwedge_{-n \le i < n} (v_i,v_{i+1}) \in E.
\]
The set $\fD_{sm}$ will be nonempty iff $\fD$ contains a directed cycle (or a bi-infinite directed path, in the infinite case).
\end{prop}

\begin{defn} If $\fD$ is a digraph and $\fD_{sm}$ is defined as in the previous proposition, then we call $\fD_{sm}$ the \emph{smooth part} of the digraph $\fD$.
\end{defn}

Note that the smooth part of a digraph may contain vertices which are not themselves part of any directed cycles: it may also contain intermediate vertices along directed paths connecting two directed cycles. In fact, the smooth part of a digraph enjoys the following convexity property.

\begin{prop} If $\fD$ is a digraph and $a,b$ are in the smooth part of $\fD$, then every vertex of $\fD$ which can be found along any directed path from $a$ to $b$ is also contained in the smooth part of $\fD$.
\end{prop}

One reason for introducing this terminology is that it lets us easily state results such as the following one.

\begin{prop} If $\bS \lhd \RR$ and $\RR,\bS \le \bA\times \bA$ correspond to digraphs $\fR, \fS$ with vertex set $\bA$, and if $\fR$ is smooth, then the smooth part $\fS_{sm}$ of the digraph $\fS$ has vertex set equal to an absorbing subalgebra of $\bA$, which will be nonempty as long as $\fS$ contains some directed cycle.
\end{prop}

Of course, we will often abuse notation a little further, and talk about the ``smooth part of the digraph $\bS$'' as long as this does not seem likely to cause confusion. It will be convenient to have the following criterion for the existence of a directed cycle contained in a subalgebra $\bB \le \bA$.

\begin{prop} If $\RR \le \bA \times \bA$, and if $\bB \le \bA$ is finite and satisfies either $\bB \subseteq \bB + \RR$ or $\bB \subseteq \bB-\RR$, then the restriction $\RR \cap (\bB \times \bB)$ of $\RR$ to $\bB$ has nonempty smooth part.
\end{prop}
\begin{proof} Suppose that $\bB \subseteq \bB - \RR$. Then every vertex in $\bB$ has an edge leaving it which lands in $\bB$, so we can find an arbitrarily long directed walk of $\RR$ which is entirely contained in $\bB$. Since $\bB$ is finite, this implies that there is some directed cycle which is entirely contained in $\bB$.
\end{proof}

As a warmup to the full loop lemma, we will first focus on the special case where the relation $\RR$ is linked. This special case is usually enough to handle most applications.

\begin{lem}[Loop Lemma, linked case]\label{loop-linked} Suppose that $\bA$ is a finite Taylor algebra and that $\RR \le_{sd} \bA\times \bA$ is a linked subdirect relation. Then $\RR$ contains a loop, that is, $\RR \cap \Delta_{\bA} \ne \emptyset$.
\end{lem}
\begin{proof} We prove this by induction on $|\bA|$. We may assume that $\bA$ is idempotent without loss of generality. If $\RR \ne \bA\times \bA$, then $\bA$ must have some proper absorbing subalgebra $\bB \lhd \bA$ by the Absorption Theorem \ref{absorption-theorem}. If we define a sequence of absorbing subalgebras $\bB = \bB_0, \bB_1, ...$ of $\bA$ by $\bB_{i+1} = \bB_i + \RR$ for $i$ even and $\bB_{i+1} = \bB_i - \RR$ for $i$ odd, then since $\RR$ is linked and $\bA$ is finite there must be some $i$ such that $\bB_{i+1} = \bA$ but $\bB_i \ne \bA$. Since this $\bB_i$ satisfies $\bB_i \subseteq \bA = \bB_{i+1}$, we see that either $\bB_i \subseteq \bB_i + \RR$ or $\bB_i \subseteq \bB_i - \RR$, so by the previous proposition the relation $\RR \cap (\bB_i \times \bB_i)$ has a nonempty smooth part $\bB_{sm} \lhd \bB_i$, with edge set $\bS = \RR \cap (\bB_{sm} \times \bB_{sm})$.

Since $\bB_{sm} \lhd \bA$, we have $\bS \lhd \RR$. Since $\bS$ is smooth, we can transfer the linkedness of $\RR$ to $\bS$ using Theorem \ref{absorbing-linked}, to see that $\bS$ must also be linked. By the inductive hypothesis applied to $\bB_{sm}$, we see that $\bS$ must contain a loop, and this loop will also be contained in $\RR$ since $\bS \le \RR$.
\end{proof}

To state the full loop lemma, we need another digraph concept.

\begin{defn} The \emph{algebraic length} of a weakly connected digraph $\fD$ is the least common multiple of all integers $k$ such that there is a digraph homomorphism from $\fD$ to a directed cycle of length $k$.
\end{defn}

\begin{prop} The algebraic length of a weakly connected digraph $\fD = (V,E)$ is the greatest common divisor of all integers $k$ such that there exist $v \in V$ and $k_1, k_2, ..., k_m \in \NN$ such that
\[
v \in \{v\} + k_1E - k_2E + \cdots \pm k_mE
\]
and $k = k_1 - k_2 + \cdots \pm k_m$.

Furthermore, there exists a digraph homomorphism from $\fD$ to a directed cycle $\fC$ iff the algebraic length of $\fD$ is a multiple of the length of the cycle $\fC$.
\end{prop}

\begin{prop} If $\fD = (V,E)$ is a smooth, weakly connected digraph of algebraic length $k$, then the digraph $\fD^{\circ m} = (V,E^{\circ m})$ has $\gcd(k,m)$ weakly connected components, and each weakly connected component of $\fD^{\circ m}$ has algebraic length $\frac{k}{\gcd(k,m)}$.
\end{prop}

\begin{prop}\label{prop-alg-length-linked} If $\fD = (V,E)$ is smooth and weakly connected, then $\fD$ has algebraic length $1$ if and only if there is some $m \ge 0$ such that the relation $E^{\circ m}$ is linked.
\end{prop}

\begin{cor} If $\fD = (V,E)$ is smooth and has a weakly connected component $C \subseteq V$ of algebraic length $1$, and if $v \in C$, then the set $C$ can be defined by a primitive positive formula using the singleton unary relation $\{v\}$ and the binary relation $E$.
\end{cor}

With these preliminaries out of the way, we can finally state the full version of the loop lemma for finite Taylor algebras.

\begin{thm}[Loop Lemma \cite{cyclic}]\label{loop-lemma} If $\bA$ is a finite Taylor algebra and $\RR \le_{sd} \bA\times \bA$ corresponds to a smooth digraph $\fR = (\bA,\RR)$ which has a weakly connected component of algebraic length $1$, then $\RR$ has a loop, i.e. $\RR \cap \Delta_{\bA} \ne \emptyset$.
\end{thm}
\begin{proof} We prove this by induction on $|\bA|$. We may assume that $\bA$ is idempotent without loss of generality. We may also assume that $\fR$ is weakly connected by restricting to a weakly connected component of algebraic length $1$ (which forms a subalgebra of $\bA$ by the results above). Let $m$ be minimal such that $\RR^{\circ m}$ is linked. We split into cases based on whether $\RR^{\circ m} = \bA \times \bA$ or not.

If $\RR^{\circ m} \ne \bA\times \bA$, then by the Absorption Theorem \ref{absorption-theorem} we see that $\bA$ must have some proper absorbing subalgebra. By a similar argument to the linked case (Lemma \ref{loop-linked}), we see that there is some proper absorbing $\bB \lhd \bA$ such that $\bS = \RR \cap (\bB\times \bB)$ is subdirect in $\bB \times \bB$. Then since $\bS^{\circ m} \lhd \RR^{\circ m}$, we can apply Theorem \ref{absorbing-linked} to see that $\bS^{\circ m}$ is linked, so the smooth digraph $\fS = (\bA,\bS)$ has algebraic length $1$ and $\bS$ has a loop by the inductive hypothesis.

If $\RR^{\circ m} = \bA\times \bA$, then we let $\bB$ be any linked component of $\RR^{\circ (m-1)}$ on the first coordinate (note that $\RR^{\circ (m-1)}$ is not linked by the choice of $m$, so $\bB$ is a proper subalgebra of $\bA$). First we will show that $\bB \subseteq \bB-\RR$. To see this, let $b \in \bB$ be arbitrary, pick any $c \in b + \RR^{\circ (m-1)}$. Then since $\RR^{\circ m} = \bA \times \bA$, we have
\[
c \in b+\RR^{\circ m},
\]
and if we let $d$ be the first element along a directed path of length $m$ from $b$ to $c$, then we have
\[
d \in (b+\RR) \cap (c - \RR^{\circ (m-1)}) \subseteq (b+\RR) \cap (b + \RR^{\circ (m-1)} - \RR^{\circ (m-1)}) \subseteq (b+\RR) \cap \bB.
\]
Thus $b \in \bB - \RR$, and since $b$ was an arbitrary element of $\bB$ we see that $\bB \subseteq \bB-\RR$. Thus the smooth part $\bB_{sm}$ of $\RR\cap (\bB\times \bB)$ is nonempty.

To finish, we just need to check that the smooth digraph corresponding to $\bS = \RR \cap (\bB_{sm} \times \bB_{sm})$ has algebraic length $1$. For this, we pick any $(b,c) \in \bS^{\circ (m-1)}$, and pick any directed path $b = b_1, ..., b_m = c$ of length $m-1$ with all $b_i \in \bB_{sm}$. Since $(b,c) \in \RR^{\circ m}$, we may also find directed path $b = c_0, ..., c_m = c$ from $b$ to $c$ of length $m$ in $\fR$. We will show that every $c_i$ along this path is actually in $\bB_{sm}$. For this, we just note that $b_i,c_i$ are in the same linked component of $\RR^{\circ (m-1)}$ for each $i \ge 1$ (since $b_i$ and $c_i$ can both reach $c$ in exactly $m-i$ steps), so each $c_i$ is at least in $\bB$, and then since each $c_i$ is along a directed path between two vertices of $\bB_{sm}$ we see that each $c_i$ belongs to the smooth part $\bB_{sm}$ as well. Thus $b \in b + \bS^{\circ (m-1)} - \bS^{\circ m}$, so $\bS$ has algebraic length $1$ and we may apply the inductive hypothesis to see that $\bS$ contains a loop.
\end{proof}

\begin{cor}[Siggers term \cite{siggers-original}, \cite{optimal-taylor}]\label{siggers-term} If $\bA$ is a finite Taylor algebra, then $\bA$ has a $4$-ary idempotent term $t$ which satisfies the identity
\[
t(x,x,y,z) \approx t(y,z,z,x).
\]
\end{cor}
\begin{proof} Assume without loss of generality that $\bA$ is idempotent. Let $\bF = \cF_{\bA}(x,y,z)$ be the free algebra on three generators in the variety generated by $\bA$. Let $\RR$ be the binary relation
\[
\RR = \Sg_{\bF^2}\left\{\begin{bmatrix} x\\ y\end{bmatrix}, \begin{bmatrix} x\\ z\end{bmatrix}, \begin{bmatrix} y\\ z\end{bmatrix}, \begin{bmatrix} z\\ x\end{bmatrix}\right\}.
\]
Then $\RR$ is clearly subdirect, and the generating set of $\RR$ forms the binary relation on $\{x,y,z\}$ pictured below, as both a bipartite graph and as a digraph.
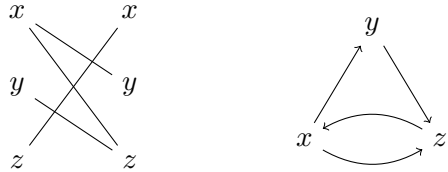
\begin{center}
\begin{tabular}{ccc}
\begin{tikzpicture}[scale=1,baseline=0.5cm]
  \node (x1) at (-1,2) {$x$};
  \node (y1) at (-1,1) {$y$};
  \node (z1) at (-1,0) {$z$};
  \node (x2) at (0.5,2) {$x$};
  \node (y2) at (0.5,1) {$y$};
  \node (z2) at (0.5,0) {$z$};
  \draw (x1) edge (y2) (y1) edge (z2);
  \draw (z1) edge (x2) (x1) edge (z2);
\end{tikzpicture} & \hspace{1cm} &
\begin{tikzpicture}[scale=1.5,baseline=0.2cm]
  \node (x) at (-0.6,0) {$x$};
  \node (y) at (0,1) {$y$};
  \node (z) at (0.6,0) {$z$};
  \draw [->] (x) edge (y) (y) edge (z);
  \draw [->] (z) edge[bend right] (x) (x) edge[bend right] (z);
\end{tikzpicture}
\end{tabular}
\end{center}
This digraph is smooth, strongly connected (in fact, it has $x + \RR^{\circ 3} = \bF$ and $\RR^{\circ 5} = \bF \times \bF$), and has algebraic length $1$ (since $x \in x + \RR^{\circ 2} - \RR^{\circ 1}$), so we can apply the Loop Lemma to see that $\RR$ contains some loop $(f,f)$ (we are using here the fact that $\bF \le \bA^{\bA^3}$ is finite and Taylor). Then since $(f,f) \in \RR$, there must be some $4$-ary term $t$ such that
\[
t\left(\begin{bmatrix} x\\ y\end{bmatrix}, \begin{bmatrix} x\\ z\end{bmatrix}, \begin{bmatrix} y\\ z\end{bmatrix}, \begin{bmatrix} z\\ x\end{bmatrix}\right) = \begin{bmatrix} f\\ f\end{bmatrix},
\]
and this $t$ then satisfies the identity
\[
t(x,x,y,z) = f = t(y,z,z,x).\qedhere
\]
\end{proof}

\begin{rem} Suppose that $t$ is a Siggers term, i.e. that $t(x,x,y,z) \approx t(y,z,z,x)$. If we substitute $y=z$ into the Siggers identity and rename variables, we see that
\[
t(y,y,x,x) \approx t(x,x,x,y),
\]
and if we substitute $x=y$ into the Siggers identity and rename variables, then we get
\[
t(x,x,x,y) \approx t(x,y,y,x).
\]
Thus there is some binary term $f(x,y)$ such that
\[
t\left(\begin{bmatrix}y & y & x & x\\ x & y & y & x\\ x & x & x & y\end{bmatrix}\right) \approx \begin{bmatrix}f(x,y)\\f(x,y)\\f(x,y)\end{bmatrix}.
\]
If we reorder the first and second inputs to $t$, the left hand side exactly becomes the left hand side of the equation for the $3$-edge term. If $f(x,y)$ was equal to $x$, then $t$ would become a $3$-edge term (up to reordering inputs).

If $f(x,y)$ was instead equal to $y$, then $p(x,y,z) = t(x,x,y,z)$ would become a Mal'cev term, which is even better than a $3$-edge term. However, if we allow for the possibility of semilattice subalgebras, then $f(x,y)$ must act as the semilattice operation on any two-element semilattice subalgebra, and of course in this case there couldn't possibly be any cube term operation of any arity. For this reason, the system of equations satisfied by $t$ above are often summarized by calling such a $t$ a ``weak $3$-edge term''.

The fact that a Siggers term looks suspiciously similar to a $3$-edge term is more than a coincidence: Theorem 6.2 of \cite{minimal-taylor} shows that every finite Taylor algebra either has a $3$-edge term or has some pair of elements $a \ne b$ such that $(b,b) \in \Sg\{(a,b),(b,a)\}$.
\end{rem}



\section{Finite abelian Taylor algebras are affine, and Zhuk's four cases}

First we recall the definition of an abelian algebra.

\begin{defn} An algebraic structure $\bA$ is called \emph{abelian} if there is a congruence $\Theta$ on $\bA\times \bA$ such that the diagonal $\Delta_\bA = \{(a,a) \mid a \in \bA\}$ is one of the congruence classes of $\Theta$.
\end{defn}

The reader might be skeptical about how often such a congruence $\Theta$ actually shows up. After all, such a congruence is most naturally viewed as a $4$-ary relation on $\bA$, and for the most part we have only been able to prove interesting structural results about binary relations so far. The next result illustrates the most common situation which leads to the existence of such a congruence.

\begin{prop}\label{ternary-abelian} Suppose that $\RR \le_{sd} \bA \times \bA \times \bA$ has the property that for each $a \in \bA$, and for each permutation $(i,j,k)$ of $(1,2,3)$, the binary relation
\[
\pi_{ij}(x\in\RR \wedge x_k = a)
\]
is the graph of an automorphism of $\bA$. Then $\bA$ is abelian.
\end{prop}
\begin{proof} Note that the assumption on $\RR$ can be rephrased as saying that if we fix any pair of coordinates of a tuple in $\RR$, then the last coordinate is uniquely determined. Therefore $\RR$ can be viewed as the graph of a homomorphism
\[
m: \bA\times \bA \twoheadrightarrow \bA
\]
such that the preimage $m^{-1}(a)$ is the graph of an automorphism of $\bA$ for every $a \in \bA$ (equivalently, $m$ is the multiplication of some quasigroup which commutes with the operations of $\bA$). In other words, every congruence class of the kernel $\ker m \in \Con(\bA \times \bA)$ is the graph of an automorphism of $\bA$. Twisting $\ker m$ by one of these automorphisms yields a congruence $\Theta \in \Con(\bA\times \bA)$ such that one of its congruence classes is the graph of the identity permutation of $\bA$.
\end{proof}

The proof we give in this section - following \cite{pointing-no-absorption} - of the fact that finite abelian Taylor algebras are affine breaks into three steps:
\begin{itemize}
\item every finite abelian algebra is (hereditarily) absorption free,

\item every finite, idempotent, Taylor, hereditarily absorption free algebra is Mal'cev, and

\item every abelian Mal'cev algebra is affine.
\end{itemize}
We have already completed the third step in Section \ref{s-abelian-malcev}, Theorem \ref{abelian-malcev}. We will complete the remaining steps in reverse order as well.

\begin{defn} We say that an algebra $\bA$ is \emph{hereditarily absorption free} if every subalgebra of $\bA$ is absorption free, that is, if $\bC \lhd \bB \le \bA$ implies that $\bC = \bB$ or $\bC = \emptyset$.
\end{defn}

\begin{prop} Suppose $\bA, \bB$ are idempotent and hereditarily absorption free. Then $\bA \times \bB$ is also hereditarily absorption free.
\end{prop}
\begin{proof} Suppose that $\bS \lhd \RR \le \bA \times \bB$, with $\bS \ne \emptyset$. Then since $\pi_1(\bS) \lhd \pi_1(\RR) \le \bA$ and $\bA$ is hereditarily absorption free, we see that $\pi_1(\bS) = \pi_1(\RR)$. Thus for every $a \in \pi_1(\RR)$ we have $a + \bS \ne \emptyset$, and since $\bA$ is idempotent, we have
\[
a+\bS \lhd a+\RR \le \bB.
\]
Then since $\bB$ is hereditarily absorption free, we see that $a+\bS = a+\RR$. Since $a$ was an arbitrary element of $\pi_1(\RR)$, we have $\bS = \RR$.
\end{proof}

\begin{thm}[HAF implies Mal'cev \cite{pointing-no-absorption}]\label{haf-malcev} If $\bA$ is finite, idempotent, Taylor, and hereditarily absorption free, then $\bA$ is Mal'cev.
\end{thm}
\begin{proof} By repeatedly applying the previous proposition, we see that the free algebra on two generators $\bF = \cF_{\bA}(x,y) \le \bA^{\bA^2}$ is absorption free. Consider the binary relation $\RR \le_{sd} \bF\times \bF$ defined by
\[
\RR = \Sg_{\bF^2}\left\{\begin{bmatrix} x\\ y\end{bmatrix}, \begin{bmatrix} x\\ x\end{bmatrix}, \begin{bmatrix} y\\ x\end{bmatrix}\right\}.
\]
Then $x+\RR \supseteq \Sg_\bF\{x,y\} = \bF$, so $x$ is contained in the left center of $\RR$. Thus by the Absorption Theorem \ref{absorption-theorem} (or just Zhuk's result Corollary \ref{zhuk-center}) we must have $\RR = \bF \times \bF$, and in particular $(y,y) \in \RR$. Thus there is some ternary term $p$ such that
\[
p\left(\begin{bmatrix} x\\ y\end{bmatrix}, \begin{bmatrix} x\\ x\end{bmatrix}, \begin{bmatrix} y\\ x\end{bmatrix}\right) = \begin{bmatrix} y\\ y\end{bmatrix}.\qedhere
\]
\end{proof}

To finish the proof that finite abelian Taylor algebras are affine, we just need to check that every abelian algebra is absorption free. Note that every subalgebra of an abelian algebra is also abelian, so this will imply that abelian algebras are \emph{hereditarily} absorption free as well. Additionally, every reduct of an abelian algebra is also abelian (since taking reducts can only increase the congruence lattice), so we see that the idempotent reduct of a finite abelian Taylor algebra will also be hereditarily absorption free, allowing us to apply the previous result to it.

It is not so easy to see how to use abelianness to rule out absorption. As a warmup, we will show that abelian algebras can't have any near-unanimity terms: this will give us the hint about how to show that finite abelian algebras are absorption free.

\begin{prop} If an algebra $\bA$ is abelian and has at least two elements, then $\bA$ does not have a near-unanimity term.
\end{prop}
\begin{proof} Let $\Theta \in \Con(\bA \times \bA)$ be a congruence with the diagonal $\Delta_{\bA}$ as a congruence class. Suppose for contradiction that $t$ is a near-unanimity term operation of minimal arity $n$, and note that $n$ must be at least $3$ since $\bA$ has at least two elements. Let $a,b$ be any pair of elements of $\bA$. Then we have
\[
t\left(\begin{bmatrix} a & b & b & \cdots & b\\ b & b & b & \cdots & b\end{bmatrix}\right) = \begin{bmatrix} b\\ b\end{bmatrix} \in \Delta_{\bA}.
\]
Since the second column of inputs to $t$ is $(b,b) \in \Delta_\bA$, we can replace it with any other element of $\Delta_\bA$ without changing the the result modulo $\Theta$. Thus we have
\[
t\left(\begin{bmatrix} a & a & b & \cdots & b\\ b & a & b & \cdots & b\end{bmatrix}\right) \equiv_\Theta \begin{bmatrix} b\\ b\end{bmatrix} \in \Delta_{\bA}.
\]
Since $t(b,a,b,...,b) = b$, we see that we must have
\[
t\left(\begin{bmatrix} a & a & b & \cdots & b\\ b & a & b & \cdots & b\end{bmatrix}\right) = \begin{bmatrix} b\\ b\end{bmatrix}.
\]
Since $a,b$ were arbitrary elements of $\bA$, we see that
\[
t(y,y,x,...,x) \approx x,
\]
so the term $t(x,x,y_2,...,y_{n-1})$ is a near-unanimity term operaion of arity $n-1$, contradicting the choice of $t$.
\end{proof}

In order to mimic this argument to rule out absorption, we will need to assume finiteness of $\bA$ and apply an iteration argument. (To see that the finiteness assumption is really needed, consider the infinite idempotent quasiaffine Taylor algebra $([0,1], \frac{x+y}{2})$, which has the open interval $(0,1)$ as a binary absorbing subalgebra.)

\begin{thm}[Abelian implies HAF \cite{pointing-no-absorption}]\label{abelian-haf} If a finite algebra $\bA$ is abelian, then it is absorption free.
\end{thm}
\begin{proof} Let $\Theta \in \Con(\bA \times \bA)$ be a congruence with the diagonal $\Delta_{\bA}$ as a congruence class. Suppose for contradiction that $\bB \lhd \bA$ is nonempty and proper, and let $t$ be a term of minimal arity $n$ among those which absorb $\bB$. Note that $n \ge 2$ since $\bB$ is a proper subalgebra of $\bA$. Now iterate $t$ on its first argument, i.e. define a sequence of terms $t_i$ with $t_1 = t$ and
\[
t_{i+1}(x,y_1, ..., y_{n-1}) \coloneqq t(t_i(x,y_1,...,y_{n-1}),y_1,...,y_{n-1}).
\]
By induction on $i$, each $t_i$ absorbs $\bB$. Since $\bA$ is finite, there is some $i$ such that $t_i = t_{2i}$, set $t_\infty = t_i$. Then we have
\[
t_\infty(t_\infty(x,y_1, ..., y_{n-1}),y_1, ..., y_{n-1}) \approx t_\infty(x,y_1, ..., y_{n-1}),
\]
and $t_\infty$ absorbs $\bB$.

Now we argue as in the near-unanimity case: let $a \in \bA$ and $b_1,b_2,...,b_{n-1} \in \bB$, and set
\[
b = t_\infty(a,b_1,b_2,...,b_{n-1}) \in \bB.
\]
Then we have
\[
t_\infty\left(\begin{bmatrix} a & b_1 & b_2 & \cdots & b_{n-1}\\ b & b_1 & b_2 & \cdots & b_{n-1}\end{bmatrix}\right) = \begin{bmatrix} b\\ b\end{bmatrix} \in \Delta_\bA,
\]
so since $(b_1,b_1) \equiv_\Theta (a,a)$, we have
\[
t_\infty\left(\begin{bmatrix} a & a & b_2 & \cdots & b_{n-1}\\ b & a & b_2 & \cdots & b_{n-1}\end{bmatrix}\right) \equiv_\Theta \begin{bmatrix} b\\ b\end{bmatrix} \in \Delta_{\bA}.
\]
Thus since $\Delta_\bA$ is a congruence class of $\Theta$ and $\bB$ absorbs $\bA$ with respect to $t_\infty$, we have
\[
t_{\infty}(a,a,b_2,...,b_{n-1}) = t_\infty(b,a,b_2,...,b_{n-1}) \in \bB.
\]
Since $a$ was an arbitrary element of $\bA$ and $b_2, ..., b_{n-1}$ were arbitrary elements of $\bB$, we see that the term
\[
t_\infty(x,x,y_2,...,y_{n-1})
\]
absorbs $\bB$ and has arity $n-1$, contradicting the choice of $t$.
\end{proof}

Now we can put all the pieces together and get our main result.

\begin{thm}[Fundamental Theorem of Abelian Algebras, finite Taylor case \cite{hobby-mckenzie}, \cite{pointing-no-absorption}, \cite{stronkowski-embedding}, \cite{zhuk-key}]\label{taylor-abelian} If $\bA$ is a finite abelian Taylor algebra, then $\bA$ is affine.
\end{thm}
\begin{proof} Let $\bA^{id}$ be the idempotent reduct of $\bA$, note that $\bA^{id}$ is still abelian and Taylor (since Taylor terms are idempotent by definition). Then every subalgebra of $\bA^{id}$ is also abelian, so by Theorem \ref{abelian-haf} $\bA^{id}$ is hereditarily absorption free. Since $\bA^{id}$ is finite, idempotent, Taylor, and hereditarily absorption free it has a Mal'cev term $p$ by Theorem \ref{haf-malcev}. Then $p$ is also a Mal'cev term operation of $\bA$, so we can apply Theorem \ref{abelian-malcev} to see that $\bA$ is affine.
\end{proof}

\begin{rem} It is not hard to generalize Theorem \ref{abelian-haf} to show that if a finite algebra $\bA$ is solvable, then $\bA$ is hereditarily absorption free (this follows from the fact that every solvable idempotent algebra $\bA$ has a congruence $\theta$ such that $\bA/\theta$ is abelian and every congruence class of $\theta$ is solvable). Thus finite solvable Taylor algebras are also Mal'cev by Theorem \ref{haf-malcev}.
\end{rem}

Now we can apply the fundamental theorem of abelian algebras to further constrain relations on absorption free algebras.

\begin{thm}[Zhuk \cite{zhuk-dichotomy}]\label{absorption-free} Suppose that $\bA$ is finite, simple, idempotent, Taylor, has no binary or centrally absorbing subalgebras, and is \emph{not} affine. Then every subdirect relation $\RR \le_{sd} \bA^n$ is the intersection of its binary projections, each of which is either a full relation or the graph of an automorphism of $\bA$.
\end{thm}
\begin{proof} We call a subdirect relation $\RR \le \bA^n$ \emph{irredundant} if no $\pi_{ij}(\RR)$ is the graph of an automorphism of $\bA$. We will prove by induction on $n$ that every irredundant subdirect relation on $\bA$ is the full relation.

The base cases of the induction are the cases $n = 1,2,3$. The case $n = 1$ is trivial (a unary subdirect relation must be full). The case $n=2$ follows from the Absorption Theorem \ref{absorption-theorem}, since every subdirect binary relation on $\bA$ is either the graph of an automorphism of $\bA$, or is linked (since $\bA$ is simple) and therefore is equal to the full relation (since $\bA$ has no binary or centrally absorbing subalgebras).

For the case $n=3$, note by the $n=2$ case both $\pi_{13}(\RR)$ and $\pi_{23}(\RR)$ must be full relations, so for any $a \in \bA$ the binary relation
\[
\RR^a \coloneqq \pi_{12}(\RR \cap (\bA^2\times\{a\}))
\]
is subdirect. Then by the $n=2$ case again, we see that $\RR^a$ is either the graph of an automorphism or is equal to $\bA^2$. If there is any $a \in \bA$ such that $\RR^a = \bA^2$, then $a$ is contained in the right center of $\RR$, considered as a binary relation on $(\bA^2)\times \bA$, so $\RR = \bA^3$ by the Absorption Theorem \ref{absorption-theorem} (or just Corollary \ref{zhuk-center}) and the fact that $\bA^2$ has no proper binary or centrally absorbing subalgebras by Proposition \ref{prop-absorption-in-powers}. Otherwise every $\RR^a$ is the graph of an automorphism, and a similar argument applies if we permute the coordinates of $\RR$, so we may apply Proposition \ref{ternary-abelian} to see that $\bA$ is abelian. But then by the fundamental theorem of abelian algebras \ref{taylor-abelian} we see that $\bA$ is affine, which contradicts our assumptions.

For the induction step, assume that $n > 3$. Then for every pair of distinct $i,j \le n-1$, the ternary relation $\pi_{ijn}(\RR)$ is full by the $n=3$ case, so for every $a \in \bA$, the binary relation
\[
\pi_{ij}(\RR \cap (\bA^{n-1}\times \{a\}))
\]
is the full relation $\bA^2$. Thus the relation
\[
\RR^a \coloneqq \pi_{[n-1]}(\RR \cap (\bA^{n-1}\times \{a\}))
\]
is irredundant, so by the inductive hypothesis, $\RR^a$ is the full relation $\bA^{n-1}$ for every $a \in \bA$. In other words, $\RR$ is the full relation $\bA^n$.
\end{proof}

Since the conclusion of Theorem \ref{absorption-free} is actually much stronger than merely being polynomially complete, we will give it a special name.

\begin{defn}\label{defn-subdirectly-complete} We say that an algebra $\bA$ is \emph{subdirectly simple} if every subdirect relation $\RR \le_{sd} \bA^n$ is the intersection of its binary projections, each of which is either a full relation or the graph of an automorphism of $\bA$.
\end{defn}

\begin{prop} Every subdirectly simple finite algebra is polynomially complete.
\end{prop}

\begin{rem} Every simple algebra with a Pixley term is subdirectly simple, in addition to being polynomially complete (see Theorem \ref{pixley-poly}).
\end{rem}

\begin{cor}[Zhuk's four cases \cite{zhuk-dichotomy}]\label{zhuk-four-cases} If $\bA$ is a nontrivial finite idempotent Taylor algebra, then at least one of the following is true.
\begin{itemize}
\item $\bA$ has a proper binary absorbing subalgebra,
\item $\bA$ has a proper centrally absorbing subalgebra,
\item $\bA$ has a nontrivial affine quotient, or
\item $\bA$ has a nontrivial subdirectly simple quotient.
\end{itemize}
\end{cor}
\begin{proof} Let $\theta \in \Con(\bA)$ be a maximal congruence on $\bA$, so $\bA/\theta$ is simple. If $\bA/\theta$ has a proper binary or centrally absorbing subalgebra $\bB$, then the preimage of $\bB$ under the projection $\bA \twoheadrightarrow \bA/\theta$ is a proper binary or centrally absorbing subalgebra of $\bA$. Otherwise, Theorem \ref{absorption-free} shows that if $\bA/\theta$ is not affine, then it is subdirectly simple.
\end{proof}

\begin{rem} For the sake of proving Theorem \ref{absorption-free} and Corollary \ref{zhuk-four-cases}, we only need to show that if $\bA$ is a finite idempotent Taylor algebra with a ternary relation $\RR \le_{sd} \bA^3$ as in Proposition \ref{ternary-abelian}, then $\bA$ is affine. It's possible to give a direct argument for this, as follows.

First, we reinterpret $\RR$ as the graph of a quasigroup operation $\cdot : \bA \times \bA \rightarrow \bA$. Using this quasigroup operation $\cdot$, we can define a Mal'cev operation $p : \bA^3 \rightarrow \bA$ which is centralized by the clone of $\bA$, such that $p$ is invertible in its first and last variables. We then pick any element $0 \in \bA$, and define the binary operation $m : \bA^2 \rightarrow \bA$ by $m(x,y) \coloneqq p(x,0,y)$. Then we have $m(x,0) = p(x,0,0) = x$ and $m(0,x) = p(0,0,x) = x$ for all $x \in \bA$, so we can apply the variant of the Eckmann-Hilton principle from Remark \ref{gen-eckmann-hilton} to see that $m$ must be commutative and associative. This $m$ will also be cancellative by construction, so by the finiteness of $\bA$ we see that $m$ defines an abelian group structure on $\bA$, which shows that $\bA$ is quasiaffine. One then needs to check that any finite Taylor algebra which is quasiaffine has a Mal'cev polynomial to finish the argument.
\end{rem}



\section{Bounded width: affine-free CSPs are solved by cycle-consistency}

Really, the title of this section should be referring to $pq$-consistency (see Definition \ref{defn-pq-consistent}), but I wanted to keep the table of contents understandable. We have already shown in Theorem \ref{absorbing-pq} that if we have a $pq$-consistent instance of a CSP, then we can reduce some of the domains to find a $pq$-consistent instance in which every domain is absorption free. In this section, we will show that if every domain is absorption free and affine free, then we can reduce the instance further while preserving $pq$-consistency.

\begin{defn} We say that a finite idempotent algebra $\bA$ is \emph{affine-free} if no quotient of any subalgebra of $\bA$ is affine.
\end{defn}

The argument strategy is very similar to the argument in the case of strongly connected algebras. We already have most of the pieces.
\begin{itemize}
\item If a binary subdirect relation $\RR \le_{sd} \bA \times \bB$ is linked and $\bA,\bB$ are absorption free and Taylor, then $\RR = \bA\times \bB$ by the Absorption Theorem \ref{absorption-theorem}.

\item If a binary relation $\RR \le \bA \times \bA$ absorbs the diagonal $\Delta_\bA$ and $\bA$ is absorption free, then $\Delta_\bA \subseteq \RR$ by Theorem \ref{absorbing-diagonal}.

\item If a binary subdirect relation $\RR \le_{sd} \bA \times \bA$ is linked and $\bA$ is Taylor, then $\RR \cap \Delta_\bA \ne \emptyset$ by the linked case of the Loop Lemma \ref{loop-linked}.

\item If $\bA$ is simple, idempotent, Taylor, absorption free, and not affine, then $\bA$ is subdirectly simple, by Theorem \ref{absorption-free}.
\end{itemize}
The missing ingredient is an analogue of Theorem \ref{strong-ternary}.

\begin{thm}\label{absorption-free-ternary} Suppose $\RR \le_{sd} \bA \times \bB \times \bC$ is subdirect, $\bA$ has no proper binary or centrally absorbing subalgebra and no affine quotient, $\pi_{23}(\RR)$ has no proper binary absorbing subalgebra, $\pi_{12}(\RR) = \bA \times \bB$, $\pi_{13}(\RR) = \bA\times \bC$, and $\bA,\bB,\bC$ are finite idempotent Taylor algebras. Then $\RR = \bA \times \pi_{23}(\RR)$.
\end{thm}

Note that by the Absorption Theorem \ref{absorption-theorem}, we just need to prove that if we consider $\RR$ as a subdirect binary relation $\RR \le_{sd} \bA \times \pi_{23}(\RR)$, then $\RR$ is linked. If not, then the linking congruence of $\RR$ on $\bA$ is contained in some maximal congruence $\theta \in \Con(\bA)$, and if we replace $\RR$ by the quotient $\RR/\theta \le_{sd} \bA/\theta \times \bB \times \bC$, then we have a smaller counterexample to Theorem \ref{absorption-free-ternary} such that $\bA$ is simple. So we just need to rule out the case where $\bA$ is simple and $\RR$ is the graph of a homomorphism $f : \pi_{23}(\RR) \twoheadrightarrow \bA$. For this, we will use a consequence of the linked case of the Loop Lemma \ref{loop-linked}.

\begin{lem} If $\bA, \bB$ are finite Taylor algebras, $\RR, \bS \le_{sd} \bA \times \bB$, and the linked components of $\RR$ on $\bA \sqcup \bB$ contain the corresponding linked components of $\bS$, then $\RR \cap \bS \ne \emptyset$.
\end{lem}
\begin{proof} Let $\bA' \le \bA, \bB' \le \bB$ be corresponding linked components of $\RR$, with $\RR \cap (\bA'\times\bB') \ne \emptyset$. By replacing $\bA,\bB$ with $\bA',\bB'$ and shrinking $\RR,\bS$, we may assume without loss of generality that $\RR$ is linked. Then $\RR \circ \bS^{-} \le_{sd} \bA \times \bA$ is also linked, so by the linked case of the Loop Lemma \ref{loop-linked}, there is some $a \in \bA$ such that $(a,a) \in \RR \circ \bS^-$. By the definition of $\RR \circ \bS^-$, this means that there is some $b \in \bB$ such that $(a,b) \in \RR$ and $(b,a) \in \bS^-$, so $(a,b) \in \RR\cap \bS$.
\end{proof}

\begin{proof}[Proof of Theorem \ref{absorption-free-ternary}] Write $\bS = \pi_{23}(\RR) \le_{sd} \bB \times \bC$. Assume for the sake of contradiction that $\bA$ is simple and that $\RR$ is the graph of a homomorphism $f : \bS \twoheadrightarrow \bA$. Note that by the idempotence of $\bA$, for each $a \in \bA$ the set $f^{-1}(a) \subseteq \bS$ is a subalgebra of $\bS$, and let $\bS_a \coloneqq f^{-1}(a)$. The assumptions $\pi_{12}(\RR) = \bA \times \bB, \pi_{13}(\RR) = \bA \times \bC$ are equivalent to each $\bS_a = f^{-1}(a)$ being a subdirect relation on $\bB \times \bC$.

If we can show that there are $a \ne a' \in \bA$ such that $\bS_a, \bS_{a'} \le_{sd} \bB \times \bC$ have the same linked components on $\bB \sqcup \bC$, then we can apply the lemma to see that $f^{-1}(a) \cap f^{-1}(a') = \bS_a \cap \bS_{a'} \ne \emptyset$, which will give us a contradiction. To accomplish this, we will show that each $\bS_a$ has the same linked components on $\bB \sqcup \bC$ as $\bS$. In fact, we will show that for every $a \in \bA$, we have $\bS \subseteq \bS_a\circ \bS_a^{-} \circ \bS_a$.

Let $(b,c)$ be any element of $\bS$. Define a subalgebra $\bX_{bc} \le \bA\times \bA\times \bA$ by
\[
\bX_{bc} \coloneqq \left\{\begin{bmatrix}x\\y\\z\end{bmatrix} \Bigg|\ \exists b' \in \bB, c' \in \bC\text{ s.t. }\begin{bmatrix}x\\b\\c'\end{bmatrix} \in \RR \wedge \begin{bmatrix}y\\b'\\c'\end{bmatrix} \in \RR \wedge \begin{bmatrix}z\\b'\\c\end{bmatrix} \in \RR\right\}.
\]
Equivalently, we have
\[
\begin{bmatrix}x\\y\\z\end{bmatrix} \in \bX_{bc} \iff \begin{bmatrix}b\\c\end{bmatrix} \in \bS_x \circ \bS_y^- \circ \bS_z.
\]
Since each $\bS_a$ is subdirect, we have $(b,b) \in \bS_a\circ \bS_a^-$ and $(c,c) \in \bS_a^- \circ \bS_a$. Thus for each $a \in \bA$, we have
\[
\begin{bmatrix}a\\a\\f(b,c)\end{bmatrix}, \begin{bmatrix}f(b,c)\\a\\a\end{bmatrix} \in \bX_{bc},
\]
so $\bX_{bc}$ is subdirect in $\bA^3$, and for each $i \ne j \le 3$ the projection $\pi_{ij}(\bX_{bc})$ is not the graph of an automorphism of $\bA$. Thus by Theorem \ref{absorption-free}, we see that $\bX_{bc} = \bA^3$, so in particular we have $(a,a,a) \in \bX_{bc}$ for all $a \in \bA$. Since this holds for every $(b,c) \in \bS$, we see that $\bS \subseteq \bS_a\circ \bS_a^{-} \circ \bS_a$ for all $a \in \bA$, so each $\bS_a$ has the same linked components on $\bB \sqcup \bC$ as $\bS$, which completes the contradiction.
\end{proof}

\begin{cor}\label{cor-abs-free-arc} If $\bA_1, ..., \bA_n$ are finite idempotent Taylor algebras with no proper binary or centrally absorbing subalgebras such that all but at most two of the $\bA_i$s have no affine quotients, and if $\RR \le_{sd} \bA_1 \times \cdots \times \bA_n$ is a subdirect relation such that each $\pi_{ij}(\RR)$ is full, then $\RR = \bA_1 \times \cdots \times \bA_n$.
\end{cor}

\begin{cor}\label{cor-abs-free-pq} If $\bA_1, ..., \bA_n$ are finite idempotent Taylor algebras with no proper binary or centrally absorbing subalgebras and no affine quotients, and if $\RR \le_{sd} \bA_1 \times \cdots \times \bA_n \times \bA_1$ is a subdirect relation such that $\Delta_{\bA_1} \subseteq \pi_{1,n+1}(\RR)$ and $\pi_{ij}(\RR)$ is full for all pairs $(i,j)$ other than $(1,n+1)$, then $\RR$ contains every tuple whose first and last coordinates are the same.
\end{cor}
\begin{proof} Suppose first that $\RR$ has a proper binary or centrally absorbing subalgebra $\RR'$. Note that each $\pi_{ij}(\RR')$ with $(i,j) \ne (1,n+1)$ is full since the $\bA_i$ have no binary or centrally absorbing subalgebras. Additionally, $\pi_{1,n+1}(\RR')$ absorbs $\Delta_{\bA_1}$, so by Theorem \ref{loop-linked} we see that $\Delta_{\bA_1} \subseteq \pi_{1,n+1}(\RR')$ as well. Thus we may replace $\RR$ by $\RR'$, until we eventually reach a situation where $\RR$ has no proper binary or central absorption. In particular, we may assume that $\pi_{1,n+1}(\RR)$ has no proper binary or centrally absorbing subalgebras.

For any $2 \le i \le n$, we may apply Theorem \ref{absorption-free-ternary} to $\pi_{i,1,n+1}(\RR)$ to see that $\pi_{i,1,n+1}(\RR) = \bA_i \times \pi_{1,n+1}(\RR)$. Now consider $\RR$ as an $n$-ary relation
\[
\RR \le_{sd} \pi_{1,n+1}(\RR) \times \bA_2 \times \cdots \times \bA_n,
\]
and apply the previous corollary to see that $\RR = \pi_{1,n+1}(\RR) \times \bA_2 \times \cdots \times \bA_n$. In particular, since we have $\Delta_{\bA_1} \subseteq \pi_{1,n+1}(\RR)$, we see that $\RR$ contains every tuple whose first and last coordinates are equal.
\end{proof}

Now that we've gathered up all the necessary ingredients, we argue as in the case of strongly connected algebras. We start by picking some variable $x$ with $|\bA_x| > 1$, pick a maximal congruence $\theta_x \in \Con(\bA_x)$, pick a congruence class $\bA_x' \le \bA_x$ of $\theta_x$. Then we refer back to Definition \ref{defn-proper-red} to define the ``proper'' variables $y$ to be the variables such that there exists a path $p$ from $y$ to $x$ such that
\[
\PP_p/\theta_x \le_{sd} \bA_y \times \bA_x/\theta_x
\]
is the graph of a homomorphism $\iota_y : \bA_y \twoheadrightarrow \bA_x/\theta_x$, and define $\theta_y$ to be the kernel of $\iota_y$ and $\bA_y'$ to be $\iota_y^{-1}(\bA_x')$.

As in the case of strongly connected algebras, we need to check that the homomorphism $\iota_y$ does not depend on the choice of path $p$. This time, we will check this using $pq$-consistency instead of cycle-consistency.

\begin{lem} Suppose that the instance $\fX$ is $pq$-consistent, and that $x, \theta_x$ are chosen as above. Suppose that $y$ is a proper variable, and that $p,q$ are two paths from $y$ to $x$ such that $\PP_p/\theta_x, \PP_q/\theta_x$ are the graphs of homomorphisms $\iota_p, \iota_q: \bA_y \twoheadrightarrow \bA_x/\theta_x$. Then $\iota_p = \iota_q$.
\end{lem}
\begin{proof} Consider the cycles $p-q$ and $q-p$ from $y$ to $y$, then by the definition of $pq$-consistency (Definition \ref{defn-pq-consistent}) we see that there must be some $j \ge 0$ such that for all $a \in \bA_y$, we have
\[
a \in \{a\} + j(p-q + q-p) + p-q.
\]
For any $b \in \bA_x$, we have $b/\theta_x - p + p = b/\theta_x$ and $b/\theta_x - q + q = b/\theta_x$ by the assumptions on $\PP_p, \PP_q$, so we see that
\[
\{a\} + j(p-q + q-p) + p-q \subseteq \iota_p(a) - q = \iota_q^{-1}(\iota_p(a)),
\]
so we must have $\iota_q(a) = \iota_p(a)$.
\end{proof}

As a consequence, we have the following analogue of Lemma \ref{ancestral-proper-path}.

\begin{lem}\label{proper-path} Suppose that the instance $\fX$ is $pq$-consistent, and that each domain has no proper binary or centrally absorbing subalgebra. Suppose $p$ is a path from $y$ to a proper variable $z$. Then one of the following is true:
\begin{itemize}
\item $\PP_p/\theta_z = \bA_y \times \bA_z/\theta_z$, or

\item $y$ is also proper, and $\PP_p/(\theta_y\times\theta_z)$ is the graph of an isomorphism $\iota_p : \bA_y/\theta_y \xrightarrow{\sim} \bA_z/\theta_z$ such that $\iota_y = \iota_z \circ \iota_p$.
\end{itemize}
\end{lem}
\begin{proof} Since $\bA_z/\theta_z$ is simple, the linking congruence of $\PP_p/\theta_z$ must either be trivial or full. If the linking congruence of $\PP_p/\theta_z$ is full, then by the Absorption Theorem \ref{absorption-theorem} we see that $\PP_p/\theta_z = \bA_y \times \bA_z/\theta_z$. Otherwise, $\PP_p/\theta_z$ is the graph of a homomorphism from $\bA_y$ to $\bA_z/\theta_z$, so then by joining the path $p$ with a path from $z$ to $x$ we see that $y$ is proper and $\iota_y = \iota_z \circ \iota_p$.
\end{proof}

To finish, we just need to show that restricting each proper variable's domain $\bA_x$ to $\bA_x'$ gives us a $pq$-consistent instance $\fX'$. To see that $\fX'$ is arc-consistent, we apply Corollary \ref{cor-abs-free-arc} as in the proof of Lemma \ref{ancestral-red-arc}. To see that $\fX'$ is $pq$-consistent, we apply Corollary \ref{cor-abs-free-pq} as in the proof of Lemma \ref{ancestral-red-cycle}. We have proven our main result.

\begin{thm}[Kozik \cite{pq-consistency}]\label{affine-free-pq} If $\fX$ is a $pq$-consistent instance of a CSP such that every domain is finite, idempotent, Taylor, and affine-free, then $\fX$ has a solution.
\end{thm}

As a curiously roundabout consequence, we see that we can't build an affine (or even abelian) algebra out of affine-free algebras.

\begin{cor}[A special case of Lemma \ref{strictly-simple-hs}]\label{affine-free-variety} If $\bA_1, ..., \bA_n$ are finite, idempotent, Taylor, and affine-free, then the variety $\cV(\bA_1, ..., \bA_n)$ which they generate does not contain any nontrivial abelian algebras.
\end{cor}
\begin{proof} Since the variety $\cV(\bA_1, ..., \bA_n)$ is finitely generated, it is locally finite, so any nontrivial abelian algebra in this variety must contain a finite abelian algebra $\bB$ with $|\bB| > 1$. Since $\bB$ is finite, Taylor, and abelian, we see that $\bB$ is affine by Theorem \ref{taylor-abelian}. But then $\bB$ is a subquotient of some finite product of $\bA_i$s, so $\CSP(\prod_i \bA_i^k)$ fails to have bounded width for some finite $k$, which contradicts the fact that $\CSP(\bA_1, ..., \bA_n)$ is solved by $pq$-consistency.
\end{proof}

Using commutator theory, we have the following consequence (see Corollary \ref{cor-sd-meet}).

\begin{cor} If $\bA$ is a finite idempotent algebra, then $\bA$ is Taylor and affine-free if and only if the variety $\cV(\bA)$ is congruence meet-semidistributive.
\end{cor}

Using the language of pp-constructability (see Definition \ref{defn-pp-construct}), we can rephrase Theorem \ref{affine-free-pq} as follows.

\begin{cor} A relational structure $\fA$ with a finite domain has $\CSP(\fA)$ solved by $pq$-consistency if and only if $\fA$ does not pp-construct any of the relational structures $(\ZZ/p,\{1\},x+y=z)$, $p$ prime.
\end{cor}
\begin{proof} Since $\fA$ pp-constructs its rigid core and vice-versa, we may assume without loss of generality that $\fA$ is a rigid core. Then the associated algebra $\bA$ is idempotent, so $\bA$ is Taylor if and only if there is any relational structure which $\fA$ does not pp-construct. To finish, we need to check that if $\bA$ is not affine-free, then $\fA$ pp-constructs some $(\ZZ/p,\{1\},x+y=z)$. Since restricting to a subalgebra of $\bA$ and taking a quotient can both be accomplished by pp-constructions, we may suppose that $\bA$ is affine and nontrivial.

If $\bA$ is affine, then by definition $\bA$ is polynomially equivalent to some module $\bM$. If $\bA$ is also idempotent, then the relation $x+y=z$ is preserved by $\bA$, as are all singleton unary relations, so $\fA$ pp-constructs the relational structure $(\bM,x+y=z)^{rig}$ (the superscript is shorthand for throwing in all unary singleton relations). Since $\bM$ is finite, some element of $\bM$ must have prime order, say order $p$. Then the set of all elements of $\bM$ with order $p$ is pp-definable, so we may suppose without loss of generality that every nonzero element of $\bM$ has order exactly $p$. As an abelian group we then have $\bM \cong (\ZZ/p)^k$ for some $k$. Letting $c$ be any nonzero element of $\bM$, we then see that $(\bM,\{c\},x+y=z)$ is homomorphically equivalent to $(\ZZ/p,\{1\},x+y=z)$.
\end{proof}

\subsection{Weak Prague instances}

The original proofs of the bounded width conjecture (i.e., that affine-free CSPs have bounded width) didn't use the concepts of $pq$-consistency or cycle-consistency. Bulatov's argument \cite{bulatov-bounded} used $(2,3)$-consistency, and leveraged a local structure theory of bounded width algebras in terms of two element semilattice and majority subalgebras. The early arguments due to Barto and Kozik \cite{barto}, \cite{local-consistency} used simpler algebraic ingredients, but used a more complicated consistency condition satisfied by instances called \emph{Prague instances}, which were then simplified to \emph{weak Prague instances}. We won't go over the original Prague instance concept until later, but weak Prague instances have a nice definition.

\begin{defn}\label{defn-weak-prague} An instance $\fX$ of a CSP with variable domains $\bA_x$ is called a \emph{weak Prague instance} if it satisfies the following three conditions.
\begin{itemize}
\item[(P1)] The instance $\fX$ is arc-consistent, that is, each constraint relation $\RR \le \prod_{x_i} \bA_{x_i}$ is subdirect.
\item[(P2)] For every variable $x$, every set $A \subseteq \bA_x$, and every cycle $p$ from $x$ to $x$, we have the implication
\[
A + p = A \;\; \implies \;\; A - p = A.
\]
\item[(P3)] For every variable $x$, every set $A \subseteq \bA_x$, and every pair of cycles $p,q$ from $x$ to $x$, we have the implication
\[
A + p + q = A \;\; \implies \;\; A + p = A.
\]
\end{itemize}
\end{defn}

We can understand what condition (P2) says about an individual cycle $p$ in terms of the digraph associated to the binary relation $\PP_p \le_{sd} \bA_x \times \bA_x$.

\begin{prop}\label{prop-p2} A subdirect binary relation $\PP \le_{sd} \bA \times \bA$ on a finite algebra $\bA$ satisfies the implication
\[
A + \PP = A \implies A - \PP = A
\]
for all $A \subseteq \bA$ if and only if the digraph $\fP = (\bA,\PP)$ satisfies one of the following equivalent conditions:
\begin{itemize}
\item every weakly connected component of $\fP$ is strongly connected,
\item every edge of $\fP$ is contained in a directed cycle of $\fP$,
\item there is some $k \ge 0$ such that $\PP^- \subseteq \PP^{\circ k}$.
\end{itemize}
\end{prop}

Based on Proposition \ref{prop-p2}, a more memorable name for consistency condition (P2) might be \emph{reversibility}. An alternative form of condition (P2) is given in \cite{sdp}.

\begin{prop}[Barto, Kozik \cite{sdp}]\label{prop-p2*} If an instance satisfies condition (P1), then (P2) is equivalent to the following condition.
\begin{itemize}
\item[(P2*)] For all variables $x$, sets $A \subseteq \bA_x$, and cycles $p$ from $x$ to $x$ such that $A + p = A$, if $p_1$ is the first step of the cycle $p$, then we have $A + p_1 - p_1 = A$.
\end{itemize}
Note that $A + p_1 - p_1 = A$ if and only if $A$ is a union of linked components of $p_1$.
\end{prop}
\begin{proof} It's easy to see that (P1) and (P2) imply (P2*), so we'll focus on proving the more difficult implication: that (P2*) implies (P2). Suppose that $A + p = A$, and write $p = p_1 + p_2 + \cdots + p_k$, where each $p_i$ has length one. By the assumption $A + p = A$, we have
\[
(A + p_1 + \cdots + p_i) + (p_{i+1} + \cdots + p_k + p_1 + \cdots + p_i) = (A + p) + p_1 + \cdots + p_i = A + p_1 + \cdots + p_i,
\]
so we can apply (P2*) to see that
\[
(A + p_1 + \cdots + p_i) + p_{i+1} - p_{i+1} = A + p_1 + \cdots + p_i.
\]
Thus we have
\begin{align*}
A - p &= (A + p) - p\\
&= A + p_1 + \cdots + p_{k-1} + p_k - p_k - p_{k-1} - \cdots - p_1\\
&= A + p_1 + \cdots + p_{k-1} - p_{k-1} - \cdots - p_1\\
&= \cdots\\
&= A + p_1 - p_1 = A.\qedhere
\end{align*}
\end{proof}

Conditions (P1) and (P2) are closely related to the basic linear relaxation of a CSP, from subsection \ref{subsection-lp}.

\begin{thm}\label{lp-p1-p2} If $\fX$ is an instance of a CSP such that the basic linear relaxation of $\fX$ has a solution assigning probability vectors $p_C$ to each constraint $C$ of $\fX$ and probability vectors $p_x$ to each variable $x$, then the instance $\fX'$ obtained by restricting each constraint relation of $\fX$ to the support of the corresponding probability distribution $p_C$ (and similarly for the variable domains) satisfies conditions (P1) and (P2).
\end{thm}
\begin{proof} Assume for simplicity that $\fX = \fX'$, that is, that all of the probability vectors have full support. The compatibility of the probability vectors $p_C$ with the probability vectors on the variable domains ensures that $\fX$ is arc-consistent, so (P1) is satisfied. For (P2), it is easier to check condition (P2*) from Proposition \ref{prop-p2*}. We attach to each set $A \subseteq \bA_x$ a probability $P(A)$, given by
\[
P(A) = \sum_{a \in A} p_{x,a}.
\]
Now consider any step $p_1$ from a variable $x$ to an adjacent variable $y$ within a constraint $C$. Let $\PP \subseteq \bA_x \times \bA_y$ be the binary projection of the corresponding constraint relation onto $x$ and $y$, and let $p_\PP$ be the corresponding marginal distribution of $p_C$. Then we have
\[
P(A + \PP) = \sum_{b \in A+\PP} p_{y,b} \ge \sum_{b \in A+\PP} \sum_{a \in A} p_{\PP,(a,b)} = \sum_{a \in A} p_{x,a} = P(A),
\]
with equality when $A + \PP - \PP = A$. Thus if $A + p = A$, then we have
\[
P(A) \le P(A+p_1) \le P(A+p) = P(A),
\]
so $P(A+p_1) = P(A)$, and thus we have $A + p_1 - p_1 = A$.
\end{proof}

In fact, Theorem \ref{lp-p1-p2} has a converse when we restrict our attention to a single cycle at a time.

\begin{thm}\label{lp-p1-p2-converse} If $\fX$ is an instance of a CSP such that the associated hypergraph of variables and relations consists of a single cycle, then $\fX$ has properties (P1) and (P2) if and only if the basic linear relaxation of $\fX$ has a solution such that for each constraint $C$ of $\fX$, the support of the corresponding probability distribution $p_C$ is exactly equal to the relation corresponding to $C$.
\end{thm}
\begin{proof} Let $v_1, ..., v_n$ be the variables of $\fX$ which occur in two constraints, in the order in which they appear around the cycle, and let the constraints $C_1, ..., C_n$ be numbered such that $v_i$ and $v_{i+1}$ are variables of $C_i$ for each $i$ (here we interpret the subscripts $i, i+1$ modulo $n$, so $v_{n+1} = v_1$).

Consider the following directed graph on the set of pairs $(i,a)$ where $i \in \ZZ/n\ZZ$ and $a \in \bA_{v_i}$: for every element $r$ of the relation corresponding to constraint $C_i$, we make a directed edge from $(i,\pi_{v_i}(r))$ to $(i+1,\pi_{v_{i+1}}(r))$. Then conditions (P1) and (P2) guarantee that every edge of this digraph is contained in a directed cycle. Choose some finite set of directed cycles $\mathcal{C}$ of this digraph which covers each edge at least once. Then for each constraint $C_i$, we let $p_{C_i}$ be the probability distribution defined by first choosing a cycle from $\mathcal{C}$ uniformly at random, and then choosing uniformly among the elements $r$ of the relation corresponding to the constraint $C_i$ such that the edge from $(i,\pi_{v_i}(r))$ to $(i+1,\pi_{v_{i+1}}(r))$ is contained in our chosen cycle.
\end{proof}

\begin{ex} When there is more than one cycle, properties (P1) and (P2) do not necessarily imply that the basic linear relaxation has a solution, even if all relations are binary. Consider the following instance, with two variables $x,y$ taking values in the domain $A_x = A_y = \{a,b,c\}$:
\begin{align*}
\begin{bmatrix} x\\ x \end{bmatrix} &\in \Big\{\begin{bmatrix} a\\ b \end{bmatrix}, \begin{bmatrix} b\\ a \end{bmatrix}, \begin{bmatrix} c\\ c \end{bmatrix}\Big\} \;\; \wedge\\
\begin{bmatrix} x\\ x \end{bmatrix} &\in \Big\{\begin{bmatrix} a\\ c \end{bmatrix}, \begin{bmatrix} b\\ c \end{bmatrix}, \begin{bmatrix} c\\ a \end{bmatrix}, \begin{bmatrix} c\\ b \end{bmatrix}\Big\} \;\; \wedge\\
\begin{bmatrix} x\\ y \end{bmatrix} &\in \Big\{\begin{bmatrix} a\\ a \end{bmatrix}, \begin{bmatrix} b\\ b \end{bmatrix}, \begin{bmatrix} b\\ c \end{bmatrix}, \begin{bmatrix} c\\ b \end{bmatrix}, \begin{bmatrix} c\\ c \end{bmatrix}\Big\} \;\; \wedge\\
\begin{bmatrix} y\\ y \end{bmatrix} &\in \Big\{\begin{bmatrix} a\\ b \end{bmatrix}, \begin{bmatrix} b\\ c \end{bmatrix}, \begin{bmatrix} c\\ a \end{bmatrix}\Big\}.
\end{align*}
First we check that the basic linear relaxation to this instance has no solutions. Suppose for a contradiction that the basic linear relaxation had a solution $p$. The first constraint implies that $p_{x,a} = p_{x,b}$, and the second constraint implies that $p_{x,a} + p_{x,b} = p_{x,c}$, so we must have $p_{x,a} = p_{x,b} = 1/4$ and $p_{x,c} = 1/2$. Similarly, the last constraint implies that $p_{y,a} = p_{y,b} = p_{y,c} = 1/3$. But the constraint connecting $x$ to $y$ implies that $p_{x,a} = p_{y,a}$, which is a contradiction as $1/4 \ne 1/3$.

Now we check that this instance satisfies (P2), using Proposition \ref{prop-p2*}. We only need to check condition (P2) for singletons, since the domains have size $3$ and for $A \subseteq A_x$ and for any cycle $p$ from $x$ to $x$, $A + p = A$ is equivalent to $(A_x \setminus A) - p = (A_x \setminus A)$ as long as the instance is arc-consistent. In order to check (P2) for singletons, first we check that there is no path through the instance which takes $\{a,c\}, \{b,c\} \subseteq A_x$ or $\{a,b\}, \{a,c\}, \{b,c\} \subseteq A_y$ to any singleton, and we check that there is no path from $x$ to $x$ which takes $\{c\} \subseteq A_x$ to either of $\{a\}, \{b\} \subseteq A_x$. Then we check that for any path $p_1$ of length $1$ such that $A + p_1$ doesn't contain one of the bad sets which can never reach a singleton, either $A$ is a union of linked components of $p_1$ or else $A$ is one of $\{a\}, \{b\} \subseteq A_x$, $p_1$ corresponds to the second constraint, and $A + p_1 = \{c\} \subseteq A_x$.
\end{ex}

The condition (P3) can be rephrased to look slightly more similar to the condition for $pq$-consistency.

\begin{prop} An instance $\fX$ with finite variable domains $\bA_x$ satisfies condition (P3) if and only if it satisfies the following condition.
\begin{itemize}
\item[(P3*)] For all variables $x$, for all pairs of cycles $p,q$ from $x$ to $x$, and for all $a \in \bA_x$, there is some $j \ge 0$ such that
\[
\{a\} + j(p+q) = \{a\} + j(p+q) + p = \{a\} + j(p+q) + p + q.
\]
\end{itemize}
\end{prop}
\begin{proof} First we show that (P3) implies (P3*). For this, note that if we define a sequence of subsets $A_i \subseteq \bA_x$ by $A_i = \{a\} + i(p+q)$, then by the finiteness of $\bA_x$ there must be some $j,k$ with $k > 0$ such that $A_j = A_{j+k}$. But then (P3) implies that $A_j + p = A_j$ and similarly that $(A_j + p) + q = A_j + p$.

For the reverse direction, let $A \subseteq \bA_x$ satisfy $A+p+q = A$. Then by the finiteness of $A$ we can find $j$ sufficiently large such that for each $a \in A$ we have $\{a\} + j(p+q) = \{a\} + j(p+q) + p$. For this choice of $j$, we then have
\[
A = A + j(p+q) = A + j(p+q) + p = A + p.\qedhere
\]
\end{proof}

There is also a natural way to certify that a given instance satisfies condition (P3), following a similar philosophy to the method we used to find absorbing reductions of cycle consistent majority CSPs.

\begin{prop}\label{prop-prec-p3} An instance $\fX$ satisfies condition (P3) at a variable $x$ if and only if there is a partial order $\preceq$ on the power set $\cP(\bA_x)$, such that for every cycle $p$ from $x$ to $x$ and every $A \subseteq \bA_x$, we have
\[
A \preceq A+p.
\]

The instance $\fX$ satisfies (P3) everywhere if and only if there is a quasiorder $\preceq$ on the set of ordered pairs $(x,A)$ with $A \subseteq \bA_x$, such that for each binary projection $\RR_{ij} \le \bA_x \times \bA_y$ of any constraint relation of $\fX$ and for each $A \subseteq \bA_x$, we have
\[
(x,A) \preceq (y,A+\RR_{ij}),
\]
and such that for each $x$, the restriction of $\preceq$ to $\{x\}\times \cP(\bA_x)$ defines a partial order on $\cP(\bA_x)$.
\end{prop}

Weak Prague instances are closely related to $pq$-consistent instances, but they are not quite the same.

\begin{thm}\label{weak-prague-pq} Every weak Prague instance is $pq$-consistent.
\end{thm}
\begin{proof} Suppose $\fX$ is a weak Prague instance, that $x$ is a variable of $\fX$, that $p,q$ are cycles from $x$ to $x$, and that $a \in \bA_x$. We need to check that there is some $j \ge 0$ such that
\[
a \in \{a\} + j(p+q) + p.
\]

Since $\bA_x$ is finite, there must be some $j > 0$ such that
\[
\{a\} + j(p+q) = \{a\} + 2j(p+q).
\]
Let $A = \{a\} + j(p+q)$ be the common value of both sides of the above equation (note that if $\bA_x$ is idempotent, then $A$ will actually be a subalgebra of $\bA_x$). Then by (P2) we have
\[
A = A + j(p+q) \implies A = A - j(p+q),
\]
so
\[
a \in \{a\} + j(p+q) - j(p+q) = A - j(p+q) = A.
\]
Additionally, by (P3) we have
\[
A = A + p + (q + (j-1)(p+q)) \implies A = A + p,
\]
so
\[
a \in A = A + p = \{a\} + j(p+q) + p.\qedhere
\]
\end{proof}

\begin{ex} Here we give an example of a $pq$-consistent instance (in fact, even a singleton arc-consistent instance!) which is not a weak Prague instance. Consider the instance of 2-SAT with just one variable $x$, domain $\bA_x = (\{0,1\}, \maj)$, and a binary constraint relation $\RR \le_{sd} \bA_x\times \bA_x$ imposed on $(x,x)$ given by $\RR = \{(0,0),(0,1),(1,1)\}$ (that is, $\RR$ is the binary relation $\le$).

Since $\Delta_{\{0,1\}} \subseteq \RR$, we see that this instance is $pq$-consistent. However, this instance does not satisfy property (P2) of a weak Prague instance: we have
\[
\{1\} + \RR = \{1\},
\]
but
\[
\{1\} - \RR = \{0,1\} \ne \{1\}.
\]
Alternatively, we can check that (P2) is not satisfied by noting that the digraph $(\{0,1\},\le)$ is weakly connected but not strongly conected.
\end{ex}

Although not every $pq$-consistent instance satisfies (P2), we at least have the following implication.

\begin{thm}[Kozik \cite{pq-consistency}]\label{pq-p1-p3} Every $pq$-consistent instance satisfies conditions (P1) and (P3).
\end{thm}
\begin{proof} Suppose $\fX$ is a $pq$-consistent instance, that $x$ is a variable of $\fX$, that $p,q$ are cycles from $x$ to $x$, and that $A \subseteq \bA_x$ satisfies
\[
A + p + q = A.
\]
By $pq$-consistency, there is some $j \ge 0$ such that
\[
A \subseteq A + j(p+q) + p = A + p.
\]
Similarly, from
\[
(A+p) + q+p = A+p,
\]
we see that
\[
A+p \subseteq (A+p) + q = A.
\]
Thus we have $A = A+p$.
\end{proof}

\begin{ex} There is an example of an instance which satisfies (P1) and (P3), but which is not $pq$-consistent. As in the previous example, this instance will have just a single variable $x$ and a single binary constraint $\RR \le_{sd} \bA_x \times \bA_x$. We take the algebra $\bA_x$ to be the three-element dual discriminator algebra $(\{0,1,2\},d(x,y,z))$ from Example \ref{ex-dual-discriminator}. The binary relation $\RR$ is the 0/1/all constraint displayed below.
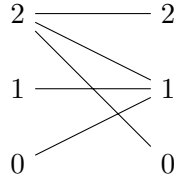
\begin{center}
\begin{tikzpicture}[scale=1]
  \node (a0) at (0,0) {$0$};
  \node (a1) at (0,1) {$1$};
  \node (a2) at (0,2) {$2$};
  \node (b0) at (2,0) {$0$};
  \node (b1) at (2,1) {$1$};
  \node (b2) at (2,2) {$2$};
  \draw (a0) -- (b1) (a1) -- (b1) (a2) -- (b0) (a2) -- (b1) (a2) -- (b2);
\end{tikzpicture}
\end{center}

To see that this is not $pq$-consistent, note that there is no $j$ such that $(0,0) \in \RR^{\circ j}$. To see that this instance satisfies condition (P3), we use the following total ordering on $\cP(\{0,1,2\})$:
\[
\emptyset \preceq \{0\} \preceq \{0,1\} \preceq \{1\} \preceq \{0,2\} \preceq \{2\} \preceq \{1,2\} \preceq \{0,1,2\}.
\]
\end{ex}

We can use weak Prague instances to see that there is a sense in which the linear programming relaxation \emph{almost} solves general CSPs of bounded width.

\begin{defn} We say that a probability distribution $\mu$ on a finite set $A$ is \emph{in general position} if we have $\mu(S) \ne \mu(T)$ for every pair of disjoint subsets $S,T \subseteq A$ with $\mu(S), \mu(T) \ne 0$. We say that a solution to the linear relaxation of an instance $\fX$ is in general position if the probability distribution which it assigns to each variable domain is in general position.
\end{defn}

\begin{prop} If there is a solution to the linear programming relaxation of $\fX$ which is in general position, then the instance we get by restricting each variable domain and relation to the support of this solution is a weak Prague instance (and is therefore $pq$-consistent as well).
\end{prop}
\begin{proof} We just need to verify condition (P3). Suppose that the solution to the linear relaxation assigns each variable $x$ to the probability distribution $\mu_x$ on the variable domain $A_x$. If $S \subseteq A_x$ is contained in the support of $\mu_x$ and $p,q$ are cycles from $x$ to $x$ such that $S + p + q = S$, then we have
\[
\mu_x(S) \le \mu_x(S+p) \le \mu_x(S+p+q) = \mu_x(S),
\]
so $\mu_x(S) = \mu_x(S+p)$. Thus we have
\[
\mu_x(S \setminus (S+p)) = \mu_x((S+p) \setminus S),
\]
so the assumption that $\mu_x$ is in general position implies that $S = S+p$.
\end{proof}

Libor Barto has raised the following question.

\begin{prob}\label{prob-p1-p3} Is it true that every instance of an affine-free CSP which satisfies conditions (P1) and (P3) has a solution?
\end{prob}

We will solve this problem later in these notes.

\section{Terms for bounded width and the meta-problem}

In this section we'll prove the existence of nice ternary terms characterizing bounded width algebras, which were first conjectured to exist by Jovanovi\'c \cite{jovanovic-terms} and later proved to exist using a Ramsey argument and the fact that bounded width CSPs are solved by $(2,3)$-consistency \cite{optimal-maltsev}. Using $pq$-consistency instead of $(2,3)$-consistency, it is possible to prove the existence of these terms directly, as noted by Kozik \cite{pq-consistency}. These nice ternary terms will allow us to efficiently solve the \emph{meta-problem} for bounded width CSPs: given a core relational structure $\fA$ as input, determine whether $\CSP(\fA)$ has bounded width.

\begin{thm}[Height $1$ identities for bounded width \cite{jovanovic-terms}, \cite{optimal-maltsev}, \cite{pq-consistency}]\label{bounded-width-terms} Suppose $\fA$ is a relational structure on a finite domain. Then $\CSP(\fA)$ has bounded relational width iff there are ternary polymorphisms $f,g \in \Pol_3(\fA)$ satisfying the height $1$ identities
\[
g(x,x,y) \approx g(x,y,x) \approx g(y,x,x) \approx f(x,x,y) \approx f(x,y,x) \approx f(x,y,y).
\]
In this case, every $pq$-consistent instance of $\CSP(\fA)$ has a solution.
\end{thm}

The identities in the statement of Theorem \ref{bounded-width-terms} may be interpreted as follows. If the common values $c(x,y)$ of $g(x,x,y)$, etc. are all equal to $x$, then $g$ is a majority function, and $f$ behaves as if it is first projection. If instead we have $c(x,y) = x \vee y$, then $f,g$ both behave as if they are the three-element semilattice operation $x \vee y \vee z$. Finally, if $c(x,y) = y$, then $f$ is a Pixley operation, so $f(x,f(x,y,z),z)$ is a majority operation, and additionally Theorem \ref{pixley-poly} applies.

Since having bounded relational width is preserved by homomorphic equivalence, we may reduce proving Theorem \ref{bounded-width-terms} to the special case where $\fA$ is a core, and then we can use Theorem \ref{rigid-core-reduction} to reduce to the case of a rigid core, so that the associated algebra $\bA$ is idempotent. Since any idempotent algebra $\bA$ such that $\CSP(\bA)$ has bounded width must be Taylor and affine-free, we see from Theorem \ref{affine-free-pq} that $\CSP(\bA)$ is solved by $pq$-consistency. Furthermore, by Corollary \ref{affine-free-variety} we see that the free algebra $\bF = \cF_{\bA}(x,y) \le \bA^{\bA^2}$ is also affine-free, so $\CSP(\bF)$ is also solved by $pq$-consistency. The plan is to construct a $pq$-consistent instance of $\CSP(\bF)$ which encodes the existence of such ternary terms $f,g$, but before we do this we need a basic result about taking closures under algebraic operations.

\begin{defn} Suppose that $\fX$ is an instance of a CSP such that every variable domain is contained in $\bA$, but possibly the variable domains and the relations of $\fX$ are not closed under the operations of $\bA$. Define $\Sg_\bA(\fX)$ to be the instance of $\CSP(\bA)$ where every variable domain and every relation of $\fX$ is replaced by the subalgebra it generates.
\end{defn}

\begin{prop}\label{pq-algebraic-closure} If $\fX$ is a $pq$-consistent instance as above, then $\Sg_\bA(\fX)$ is also $pq$-consistent.
\end{prop}
\begin{proof} For arc-consistency, let $R \subseteq \bA^n$ be any relation, and note that $\Sg_\bA(\pi_1(R)) = \pi_1(\Sg_\bA(R))$. For paths, let $R,S \subseteq \bA\times \bA$ be any binary relations, then we have $\Sg_\bA(R\circ S) \subseteq \Sg_\bA(R)\circ \Sg_\bA(S)$. For cycles interacting well with the diagonal, note that for any $B \subseteq \bA$ we have $\Sg_{\bA^2}(\Delta_B) = \Delta_{\Sg_\bA(B)}$.
\end{proof}

We have a similar result for weak Prague instances (Definition \ref{defn-weak-prague}), which we won't actually need.

\begin{prop} If $\fX$ is a weak Prague instance as above, then $\Sg_\bA(\fX)$ is also a weak Prague instance.
\end{prop}
\begin{proof} That $\Sg_\bA(\fX)$ satisfies (P1) and (P3) follows from the fact that $\fX$ is a $pq$-consistent instance (Theorem \ref{weak-prague-pq}), which implies that $\Sg_\bA(\fX)$ is also $pq$-consistent by the previous proposition, and this in turn implies that $\Sg_\bA(\fX)$ satisfies (P1) and (P3) (Theorem \ref{pq-p1-p3}). To check that $\Sg_\bA(\fX)$ satisfies (P2), we use Proposition \ref{prop-p2}: note that if $P \subseteq \bA\times \bA$ satisfies $P^- \subseteq P^{\circ k}$, then $\Sg_\bA(P)^- = \Sg_\bA(P^-) \subseteq \Sg_\bA(P^{\circ k}) \subseteq \Sg_\bA(P)^{\circ k}$.
\end{proof}

\begin{lem} Suppose $\fX$ is an instance of a CSP over the two-element domain $\{x,y\}$ with no unary relations, such that every binary projection $\pi_{i,j}(R)$ of every relation $R$ is subdirect in $\{x,y\}^2$ and has $(x,x) \in \pi_{i,j}(R)$. Then $\fX$ is $pq$-consistent.
\end{lem}
\begin{proof} The assumptions on $\fX$ directly imply that $\fX$ is arc-consistent. Now consider any pair of cycles $p,q$ from a variable $v$ of $\fX$ to itself. Note that the collection of binary relations on $\{x,y\}$ which are subdirect and contain $(x,x)$ is closed under composition and reversal, so $\PP_p,\PP_q$ are both subdirect and contain $(x,x)$. We just need to show that there is some $j$ such that $y \in \{y\} + j(p+q) + p$.

If $(y,y) \in \PP_p$, then we may take $j = 0$. Otherwise, we must have $\PP_p = \{(x,x), (x,y), (y,x)\}$, and since $(x,x) \in \PP_q$ this implies that $\PP_p \circ \PP_q \circ \PP_p = \{x,y\}^2$, so we may take $j = 1$. 
\end{proof}

\begin{proof}[Proof of Theorem \ref{bounded-width-terms}] First we prove the existence of such terms in any finite idempotent Taylor affine-free algebra $\bA$. Consider the ternary relations $R,S \subseteq \{x,y\}^3$ given by
\[
R = \left\{\begin{bmatrix} x\\ x\\ y\end{bmatrix}, \begin{bmatrix} x\\ y\\ x\end{bmatrix}, \begin{bmatrix} y\\ x\\ x\end{bmatrix}\right\}
\]
and
\[
S = \left\{\begin{bmatrix} x\\ x\\ x\end{bmatrix}, \begin{bmatrix} x\\ y\\ y\end{bmatrix}, \begin{bmatrix} y\\ x\\ y\end{bmatrix}\right\}.
\]
It's easy to check that each binary projection of $R$ and $S$ is subdirect in $\{x,y\}^2$ and contains $(x,x)$. Now consider the CSP instance $\fX$ with just a single variable $v$, and the apply the constraints $R$ and $S$ to the triple $(v,v,v)$ (if this makes you uncomfortable, you can instead use several different variables and impose equality constraints between them). By the lemma, $\fX$ is a $pq$-consistent instance.

If we let $\bF = \cF_{\bA}(x,y) \le \bA^{\bA^2}$, then we may consider $\{x,y\}$ to be a subset of $\bF$, and apply the proposition to see that $\Sg_\bF(\fX)$ is also $pq$-consistent. Since $\bF$ is finite, idempotent, Taylor, and affine-free, we can apply Theorem \ref{affine-free-pq} to see that $\Sg_\bF(\fX)$ has a solution. Suppose that this solution assigns the variable $v$ to the value $c \in \bF$. Then we have
\[
\begin{bmatrix} c\\ c\\ c\end{bmatrix} \in \Sg_\bF(R) \cap \Sg_\bF(S) = \Sg_\bF\left\{\begin{bmatrix} x & x & y\\ x & y & x\\ y & x & x\end{bmatrix}\right\} \cap \Sg_\bF\left\{\begin{bmatrix} x & x & y\\ x & y & x\\ x & y & y\end{bmatrix}\right\}.
\]
Thus there are ternary terms $f,g$ of $\bA$ such that
\[
g\left(\begin{bmatrix} x & x & y\\ x & y & x\\ y & x & x\end{bmatrix}\right) = \begin{bmatrix} c\\ c\\ c\end{bmatrix} = f\left(\begin{bmatrix} x & x & y\\ x & y & x\\ x & y & y\end{bmatrix}\right),
\]
and these $f,g$ satisfy the required identities.

For the converse direction, we will suppose that such terms $f,g$ exist for some idempotent algebra $\bA$, and prove that $\bA$ is Taylor and affine-free. It's easy to see that $\bA$ must be Taylor, since the identities satisfied by $g$ can't be satisfied by any projection. Since any identities which hold in $\bA$ also hold in any subquotient of $\bA$, we may suppose for contradiction that $\bA$ is a nontrivial idempotent affine algebra. Then $\bA$ is polynomially equivalent to some module $\bM$ over some ring $\RR$, and we may write
\[
g(x,y,z) \approx \alpha x + \beta y + \gamma z
\]
for some $\alpha,\beta,\gamma \in \RR$ with $\alpha+\beta+\gamma = 1$. Plugging in $x = 0$ to the identities
\[
g(x,x,y) \approx g(x,y,x) \approx g(y,x,x)
\]
gives $\alpha y \approx \beta y \approx \gamma y$, so
\[
g(x,y,z) \approx \alpha(x+y+z)
\]
and $3\alpha x \approx x$. Then if we plug in $x = 0$ to the identities
\[
2\alpha x + \alpha y \approx f(x,x,y) \approx f(x,y,x) \approx f(x,y,y),
\]
we see that $\alpha y \approx 2\alpha y$, so $\alpha y \approx 0$. Multiplying by $3$, we get $y \approx 0$, so in fact the algebra $\bA$ must consist of just the single element $0$, a contradiction.
\end{proof}

The proof technique of Theorem \ref{bounded-width-terms} can be used to produce many further terms which mimic the monotone self-dual functions found in the clone of a two-element majority algebra.

\begin{thm} Suppose $\CSP(\fA)$ has bounded relational width and $\fA$ is finite. Then there is a binary polymorphism $c(x,y)$, and an infinite family of polymorphisms $h_n^{\cF} \in \Pol_n(\fA)$ indexed by the collection of maximal intersecting families $\cF$ of subsets of $[n]$, such that for each set $S \in \cF$ with $S \ne [n]$, if we define $v_i^S$ by
\[
v_i^S = \begin{cases} x & i \in S,\\ y & i \not\in S,\end{cases}
\]
we have the identity
\[
h_n^\cF(v_1^S, ..., v_n^S) \approx c(x,y).
\]
\end{thm}

Now we show how we can use the ternary terms $f,g$ from Theorem \ref{bounded-width-terms} to solve the meta-problem.

\begin{thm} Suppose we are given a finite relational structure $\fA = (A, R_1, ..., R_n)$, where each relation $R_i$ has arity $m_i$ and is described by explicitly listing out its tuples, and suppose that we are promised that $\fA$ is core. Then we can determine whether $\CSP(\fA)$ has bounded width in polynomial time, and in the case where $\CSP(\fA)$ has bounded width, we can explicitly find ternary functions $f,g \in \Pol_3(\fA)$ as in Theorem \ref{bounded-width-terms}.
\end{thm}
\begin{proof} We will define an instance $\fX$ of $\CSP(\fA)$ such that every solution to $\fX$ corresponds to a pair of terms $f,g$ as in Theorem \ref{bounded-width-terms}. The instance $\fX$ will have two sets of $|A^3|$ variables, one variable for each value $f(a,b,c)$ for $a,b,c \in A$ and one variable for each value $g(a,b,c)$ for $a,b,c \in A$.

The relations of $\fX$ will do two jobs: they will ensure that $f, g \in \Pol_3(\fA)$, and they will ensure that $f,g$ satisfy the required identities. To ensure that $f \in \Pol_3(\fA)$, we consider every three tuples $a, b, c \in R_i$ (note that each of $a,b,c$ is an $m_i$-tuple of values in $A$), and we impose the constraint
\[
\begin{bmatrix} f(a_1,b_1,c_1)\\ f(a_2,b_2,c_2)\\ \vdots\\ f(a_{m_i}, b_{m_i}, c_{m_i})\end{bmatrix} \in R_i
\]
for each such tuple. The number of such constraints we need to impose to ensure that $f \in \Pol_3(\fA)$ is then
\[
\sum_i |R_i|^3,
\]
which is at most cubic in the size of the description of $\fA$. We ensure that $g \in \Pol_3(\fA)$ with a similar collection of constraints.

To enforce the required identities between $f,g$, for every pair $a,b \in A$, we impose the equality constraints
\[
g(a,a,b) = g(a,b,a) = g(b,a,a) = f(a,a,b) = f(a,b,a) = f(a,b,b).
\]
This requires a total of $5|A|^2$ equality constraints. Thus, the instance $\fX$ has overall size at most cubic in the size of the description of $\fA$.

In order to solve $\fX$, we view it as an instance of $\CSP(\fA^{rig})$, where $\fA^{rig}$ is the rigid core obtained from $\fA$ by adding a singleton unary relation $\{a\}$ for each element $a \in A$. Note that if $\fA$ is a core, then $\fA$ has bounded width iff $\fA^{rig}$ has bounded width (since each pp-constructs the other). We now attempt to solve the instance $\fX$ by using the cycle-consistency algorithm, as follows. For each variable $v$ of $\fX$, we go through the values $a \in A$ in order, and temporarily modify $\fX$ by adding the extra constraint $v \in \{a\}$ to make an instance $\fX_{v=a}$. Then we reduce $\fX_{v=a}$ until it either becomes cycle-consistent or until we reach a contradiction. If there is any $a \in A$ such that $\fX_{v=a}$ becomes cycle-consistent, then we replace $\fX$ by $\fX_{v=a}$ and move on to the next variable. If every choice of $a \in A$ leads to $\fX_{v=a}$ reaching a contradiction, then we give up and report that $\CSP(\fA)$ does not have bounded width.

If the procedure ends without us giving up, then we have found $f,g$ as in Theorem \ref{bounded-width-terms} and these terms prove that $\CSP(\fA)$ has bounded width. Conversely, if $\CSP(\fA)$ has bounded width, then the original instance $\fX$ has a solution, and each time we replace $\fX$ by $\fX_{v=a}$, the fact that $\fX_{v=a}$ can be reduced to a cycle-consistent instance implies that it has a solution, so the whole procedure will end by successfully finding a pair of functions $f,g$. Of course, if $\CSP(\fA)$ does not have bounded width, then we will fail to find a solution to $\fX$.
\end{proof}


A simple iteration argument allows us to give a criterion for bounded width involving just one ternary term and a binary term derived from it - however, the identities involved will not have height $1$, so these new terms are unsuitable for the application to the meta-problem.

\begin{thm}\label{bounded-width-term} A finite relational structure $\fA$ has bounded relational width if and only if it has a ternary polymorphism $g \in \Pol_3(\fA)$ such that, if $f$ is the binary term $f(x,y) \coloneqq g(x,x,y)$, we have
\[
g(x,x,y) \approx g(x,y,x) \approx g(y,x,x) \approx f(x,y) \approx f(f(x,y),f(y,x)) \approx f(f(x,y), f(x,y)).
\]
\end{thm}
\begin{proof} Suppose first that $\fA$ has bounded relational width, and let $f_3, g_3 \in \Pol_3(\fA)$ be terms as in Theorem \ref{bounded-width-terms}. By an iteration argument applied to the unary operation $x \mapsto g_3(x,x,x)$, we may assume without loss of generality that we have
\[
g_3(x,y,z) = h \circ g_3(x,y,z),
\]
where $h(x) \coloneqq g_3(x,x,x)$. Define a sequence of terms $g^i$ by $g^1 \coloneqq g_3$ and
\[
g^{i+1}(x,y,z) \coloneqq g^i(f_3(x,y,z),f_3(y,z,x),f_3(z,x,y)).
\]
Define binary terms $f^i$ by $f^i(x,y) \coloneqq g^i(x,x,y)$. Then we have
\[
f^1(x,y) \approx g_3(x,x,y) \approx f_3(x,x,y) \approx f_3(x,y,x) \approx f_3(x,y,y),
\]
and for each $i$ we have
\begin{align*}
f^{i+1}(x,y) &\approx g^{i+1}(x,x,y) \approx g^i(f_3(x,x,y), f_3(x,y,x), f_3(y,x,x))\\
&\approx g^i(f^1(x,y),f^1(x,y),f^1(y,x)) \approx f^i(f^1(x,y),f^1(y,x)).
\end{align*}
Thus the sequence $f^i(x,y)$ is generated by iterating the map $(x,y) \mapsto (f^1(x,y),f^1(y,x))$. Since $\fA$ is finite, there is some $N$ such that $g^N \approx g^{2N}$ and $f^N \approx f^{2N}$. Take $f \coloneqq f^N$ and $g \coloneqq g^N$ to finish the construction.

Now suppose that $f,g$ satisfy the assumed identities. Let $e$ be the unary operation $e(x) \coloneqq f(x,x) = g(x,x,x)$. The identity
\[
f(x,y) \approx f(f(x,y),f(x,y)) = e(f(x,y))
\]
implies that
\[
e(e(x)) \approx e(x),
\]
so $\fA$ is homomorphically equivalent to $e(\fA)$, and the restrictions of $f, e \circ g$ to $e(\fA)$ are idempotent. Let $\bA_e$ be the idempotent algebra $(e(\fA), f|_{e(\fA)}, e \circ g|_{e(\fA)})$. We will show that $\bA_e$ is Taylor and affine-free.

That $\bA_e$ is Taylor follows from the identity
\[
e \circ g(x,x,y) \approx e \circ g(x,y,x) \approx e \circ g(y,x,x).
\]
For the sake of contradiction, assume that $\bB \in HSP(\bA_e)$ is a nontrivial affine algebra. Then we can write
\[
e\circ g(x,y,z) \approx \alpha(x+y+z)
\]
on $\bB$, for some $\alpha$ with $3\alpha x \approx x$. Then we have
\[
f(x,y) \approx 2\alpha x + \alpha y,
\]
so
\[
f(f(x,y),f(y,x)) \approx 2\alpha(2\alpha x + \alpha y) + \alpha(2\alpha y + \alpha x) \approx 5\alpha^2 x + 4\alpha^2 y.
\]
Equating these and setting $y$ to $0$, we see that $2\alpha x \approx 5\alpha^2 x$. Multiplying by $9$ and using $3\alpha x \approx x$, we get $6x \approx 5x$, so $x \approx 0$ on $\bB$, a contradiction.
\end{proof}

The identities satisfied by the term $g$ of Theorem \ref{bounded-width-term} have the following nice consequence.

\begin{prop} Suppose that $g$ is a ternary term as in Theorem \ref{bounded-width-term}, and that $f$ is the associated binary term. Then for any $a,b$, either $f(a,b) = f(b,a)$, or the set $\{f(a,b), f(b,a)\}$ is closed under $g$, and $(\{f(a,b),f(b,a)\},g)$ is isomorphic to a two-element majority algebra.
\end{prop}

For small examples of bounded width algebras $\bA$ which do not contain large majority subalgebras, most of the structure of a bounded width algebra seems to be controlled by the binary term $f$ from Theorem \ref{bounded-width-term}, with the exact values of the ternary term $g$ only playing an important role on the majority subalgebras. I have also conjectured a very strong refinement of Theorem \ref{bounded-width-term}, which would give a much more explicit structure theory for bounded width algebras.

\begin{conj} A finite idempotent algebra $\bA$ has bounded relational width if and only if it has a ternary term $m$ and an associated binary term $s(x,y) \coloneqq m(x,x,y)$, which satisfy the identities
\[
m(x,x,y) \approx m(x,y,x) \approx m(y,x,x) \approx s(x,y)
\]
and
\[
s(x,s(x,y)) \approx s(s(x,y),x) \approx s(x,y).
\]
\end{conj}



\section{Stable subalgebras, and even weaker consistency for bounded width}

In this section we will introduce a new concept, which is similar to absorption but which is targeted at subdirect relations rather than arbitrary relations. This allows us to unify the treatment of centrally absorbing subalgebras with congruence classes of polynomially complete absorption free quotients, eliminating most of the casework we need to deal with. We will demonstrate the usefulness of this concept by solving the (P1)-(P3) problem (Problem \ref{prob-p1-p3}). The approach used in this section is based on Zhuk's theory of ``strong subalgebras'' \cite{zhuk-strong}.

Rather than directly defining stable subalgebras, we will give an axiomatic description of what we want from a concept of ``stability''.

\begin{defn}\label{defn-stability} Suppose that $\cV$ is a pseudovariety of finite algebras. We say that a binary relation $\Yleft$ on $\cV$ is a \emph{stability concept} (or just a \emph{stability}) on $\cV$ if $\Yleft$ satisfies the following axioms.
\begin{itemize}
\item[(Subalgebra)] If $\bB \Yleft \bA$, then $\bB \le \bA$.
\item[(Transitivity)] If $\bC \Yleft \bB \Yleft \bA$, then $\bC \Yleft \bA$.
\item[(Intersection)] If $\bB, \bC \Yleft \bA$ and $\bB \cap \bC \ne \emptyset$, then $\bB \cap \bC \Yleft \bB$.
\item[(Propagation)] If $f : \bA \twoheadrightarrow \bB$ is a surjective homomorphism, then
\begin{itemize}
\item[(Pushforward)] if $\bC \Yleft \bA$, then $f(\bC) \Yleft \bB$, and
\item[(Pullback)] if $\bD \Yleft \bB$, then $f^{-1}(\bD) \Yleft \bA$.
\end{itemize}
\item[(Helly)] If $\bB, \bC, \bD \Yleft \bA$ are such that $\bB \cap \bC \ne \emptyset$, $\bC \cap \bD \ne \emptyset$, and $\bB \cap \bD \ne \emptyset$, then $\bB \cap \bC \cap \bD \ne \emptyset$.
\item[(Ubiquity)] If $\bA \in \cV$ has $|\bA| > 1$, then either
\begin{itemize}
\item there is some $\bB \Yleft \bA$ such that $\bB \ne \bA, \emptyset$, or
\item there is some proper congruence $\theta \in \Con(\bA)$ such that $\bA/\theta$ is an affine algebra.
\end{itemize}
\end{itemize}
Given a stability concept $\Yleft$ on $\cV$, we say that $\bB$ is a \emph{stable subalgebra} of $\bA$ if $\bB \Yleft \bA$. We say that an element $a \in \bA$ is a \emph{stable element} if $\{a\} \Yleft \bA$.
\end{defn}

The axioms of a stability concept imply apparently stronger versions of themselves.

\begin{prop} If $\Yleft$ is a binary relation on $\cV$ which satisfies the propagation axiom from Definition \ref{defn-stability}, then for any subdirect relation $\RR \le_{sd} \bA \times \bB$ in $\cV$, we have
\[
\bC \Yleft \bA \;\; \implies \;\; \bC + \RR \Yleft \bB.
\]
\end{prop}
\begin{proof} Let $\pi_1, \pi_2$ be the surjective projection maps from $\RR$ to $\bA$ and $\bB$, respectively. Then $\pi_1^{-1}(\bC) \Yleft \bA$ by the pullback part of the propagation axiom, so $\bC + \RR = \pi_2(\pi_1^{-1}(\bC)) \Yleft \bB$ by the pushforward part of the propagation axiom.
\end{proof}

\begin{prop} If $\Yleft$ satisfies the Helly axiom and the intersection axiom from Definition \ref{defn-stability}, then for any $n$ and any $\bB_1, ... \bB_n \Yleft \bA$ such that $\bB_i \cap \bB_j \ne \emptyset$ for all $i, j \in [n]$, we have $\bigcap_{i\in [n]} \bB_i \ne \emptyset$.
\end{prop}
\begin{proof} We induct on $n$ - the base case $n = 3$ is the Helly axiom. For $n > 3$, set $\bA' = \bB_n$ and $\bB_i' = \bB_i \cap \bB_n$ for $i < n$, then by the Helly axiom we have
\[
\bB_i' \cap \bB_j' = \bB_i \cap \bB_j \cap \bB_n \ne \emptyset
\]
for all $i,j < n$, and by the intersection axiom we have $\bB_i' = \bB_i \cap \bB_n \Yleft \bB_n = \bA'$ for all $i < n$, so we can apply the induction hypothesis to $\bA'$ to see that
\[
\bigcap_{i\in [n]} \bB_i = \bigcap_{i < n} \bB_i' \ne \emptyset.\qedhere
\]
\end{proof}

\begin{prop}\label{prop-stable-essential} If $\Yleft$ satisfies the Helly, intersection, and propagation axioms from Definition \ref{defn-stability}, then for any subdirect relation $\RR \le_{sd} \bA_1 \times \cdots \times \bA_n$ in $\cV$, if $\bB_i \Yleft \bA_i$ for each $i$ and
\[
\pi_{ij}(\RR) \cap (\bB_i \times \bB_j) \ne \emptyset
\]
for all $i,j \in [n]$, then we have $\RR \cap \prod_{i \in [n]} \bB_i \ne \emptyset$.

If $\Yleft$ additionally satisfies the transitivity axiom, then we also have
\[
\RR \cap \Big(\prod_{i \in [n]} \bB_i\Big) \Yleft \RR.
\]
\end{prop}
\begin{proof} For each $i$, the pullback part of the propagation axiom implies that $\pi_i^{-1}(\bB_i) \Yleft \RR$, so the previous proposition implies that $\bigcap_{i\in [n]} \pi_i^{-1}(\bB_i) \ne \emptyset$.
\end{proof}

\begin{prop}\label{prop-stable-element} If $\Yleft$ satisfies the transitivity and ubiquity axioms from Definition \ref{defn-stability}, then for any idempotent bounded-width algebra $\bA \in \cV$ there is some $a \in \bA$ such that $\{a\} \Yleft \bA$.
\end{prop}

The precise choice of stability concept doesn't matter to us - we can use the following fact as a black box.

\begin{thm}\label{thm-stability} If $\cV$ is an affine-free pseudovariety of finite idempotent Taylor algebras, then there is at least one stability concept $\Yleft$ on $\cV$.
\end{thm}

Before we prove Theorem \ref{thm-stability}, we will apply it to prove that a weaker form of consistency suffices for bounded width CSPs.

\begin{defn} An arc-consistent instance $\fX$ of a CSP, with variable domains $\bA_x$, is called \emph{weakly consistent} if it satisfies
\begin{itemize}
\item[(W)] for all nonempty subsets $A \subseteq \bA_x$ and cycles $p,q$ from $x$ to $x$, we have
\[
A + p + q = A \;\; \implies \;\; A \cap (A + p) \ne \emptyset.
\]
\end{itemize}
\end{defn}

Weak consistency is clearly implied by properties (P1) and (P3) from Definition \ref{defn-weak-prague}. We can also rephrase weak consistency to make it look more similar to $pq$-consistency.

\begin{prop}\label{prop-weakly-consistent-cycle} An arc-consistent instance $\fX$ with finite variable domains $\bA_x$ is weakly consistent if and only if it satisfies
\begin{itemize}
\item[(W')] for all $a \in \bA_x$ and cycles $p,q$ from $x$ to $x$, there exist $j,k \ge 0$ such that
\[
a \in \{a\} + j(p+q) + p - k(p+q).
\]
\end{itemize}
In fact, if $\fX$ is weakly consistent then for each $x$ and each pair of cycles $p,q$ we can find some $j \ge 0$ such that $\Delta_{\bA_x} \subseteq \PP_{j(p+q) + p - j(p+q)}$.
\end{prop}
\begin{proof} First we prove that (W) implies (W'). By finiteness we can pick some $j \ge 1$ such that $\{a\} + j(p+q) = \{a\} + 2j(p+q)$ for all $a \in \bA_x$. Setting $A = \{a\} + j(p+q)$ and $q' = (j-1)(q+p) + q$, we have
\[
A + p + q' = A + j(p+q) = \{a\} + 2j(p+q) = A,
\]
so $A \cap A + p \ne \emptyset$, that is,
\[
(\{a\} + j(p+q)) \cap (a + j(p+q) + p) \ne \emptyset.
\]
Since this is true for all $a \in \bA_x$, we have $\Delta_{\bA_x} \subseteq \PP_{j(p+q) + p - j(p+q)}$.

Next we show that (W') implies (W). Suppose that $A \subseteq \bA_x$ satisfies $A + p + q = A$, and pick any element $a \in A$. If $j,k \ge 0$ are such that $a \in \{a\} + j(p+q) + p - k(p+q)$, then we have
\[
A \cap (A + p) \supseteq (\{a\} + k(p+q)) \cap (\{a\} + j(p+q) + p) \ne \emptyset.\qedhere
\]
\end{proof}

Our argument for showing that weakly consistent instances of bounded-width CSPs have solutions will follow the same general strategy as the argument for $pq$-consistent instances. First we will show that we can find an arc-consistent reduction where the reduced variable domains are all stable subalgebras of the original variable domains, and then we will try to show that any arc-consistent stable reduction is also weakly consistent. Unfortunately, we run into a snag: it is not clear that every arc-consistent stable reduction will really remain weakly consistent. To get around this, we introduce a still weaker condition which will make the strategy work.

\begin{defn} An arc-consistent instance $\fX$ of a CSP, with variable domains $\bA_x$ contained in a variety $\cV$ with a stability concept $\Yleft$, is called \emph{stably consistent} if it satisfies
\begin{itemize}
\item[(S)] for all nonempty stable subalgebras $\bB \Yleft \bA_x$ and cycles $p,q$ from $x$ to $x$, we have
\[
\bB + p + q = \bB \;\; \implies \;\; \bB \cap (\bB + p) \ne \emptyset.
\]
\end{itemize}
\end{defn}

If all of the variable domains are finite affine-free algebras, then stable consistency is equivalent to the following:
\begin{itemize}
\item[(S')] for all stable elements $\{a\} \Yleft \bA_x$ and cycles $p,q$ from $x$ to $x$, there exist $j,k \ge 0$ such that
\[
a \in \{a\} + j(p+q) + p - k(p+q).
\]
\end{itemize}

\begin{lem} If $\fX$ is stably consistent and the variable domains $\bA_x$ are affine free and are not all singletons, then there is some arc-consistent reduction $\fX'$ of $\fX$ such that every variable domain of $\fX'$ is a stable subalgebra of the corresponding variable domain in $\fX$, and such that at least one variable domain shrinks.
\end{lem}
\begin{proof} Consider the directed graph with vertices given by pairs $(x,\bB)$ such that $\bB \Yleft \bA_x$ and $\bB \ne \bA_x, \emptyset$, with an edge from $(x,\bB)$ to $(y,\bB + p)$ for each path $p$ from $x$ to $y$ such that $\bB + p \ne \bA_y$. Note that by the propagation axiom and the assumption that $\fX$ is arc-consistent, we always have
\[
\bB \Yleft \bA_x \;\; \implies \bB + p \Yleft \bA_y.
\]
Let $\cS$ be a maximal strongly connected component of this directed graph.

{\bf Claim:} For any $(x,\bB) \in \cS$ and any cycle $p$ from $x$ to $x$, we have $\bB \cap (\bB + p) \ne \emptyset$.

{\bf Proof of claim:} If $(x, \bB+p) \not\in \cS$, then by maximality of $\cS$ we must have $\bB + p = \bA_x$, so $\bB \cap (\bB+p) = \bB \ne \emptyset$. Otherwise, if $(x, \bB+p) \in \cS$, then there must be some cycle $q$ from $x$ to $x$ such that $\bB + p + q =\bB$. Then condition (S) from the definition of stable consistency implies that $\bB \cap (\bB + p) \ne \emptyset$.

Now define the \emph{universal cover} $\fT$ of $\fX$ to be the instance whose underlying constraint hypergraph is an infinite tree with a surjective map $\pi : \fT \twoheadrightarrow \fX$ on the sets of variables, such that for every path $p$ from $x$ to $y$ in $\fX$ and every preimage $u$ of $x$ in $\fT$ there is a unique lift of the path $p$ to $\fT$ which starts at $u$. The fact that $\fX$ is arc-consistent is equivalent to the fact that the solution set to the infinite instance $\fT$ is a subdirect relation of infinite arity.

Then for every pair of variables $u,v$ of $\fT$ and $\bB, \bC$ such that $(\pi(u), \bB), (\pi(v), \bC) \in \cS$, there is a unique non-backtracking path $p$ from $u$ to $v$ in $\fT$, and by the claim we have $(\bB + p) \cap \bC \ne \emptyset$. Applying Proposition \ref{prop-stable-essential}, we see that the set of solutions to $\fT$ such that $u \in \bB$ whenever $(\pi(u), \bB) \in \cS$ is a nonempty, stable subalgebra of the solution set to $\fT$ (well, every finite subinstance of $\fT$ has this property). Applying the pushforward part of the propagation axiom to this solution set, we obtain an arc-consistent stable reduction $\fX'$ of the instance $\fX$.
\end{proof}

\begin{lem} If $\fX$ is stably consistent and $\fX'$ is an arc-consistent stable reduction of $\fX$, then $\fX'$ is also stably consistent.
\end{lem}
\begin{proof} Let the variable domains of $\fX$ and $\fX'$ be $\bA_x, \bA_x'$, respectively. It's enough to show that if $\bB \Yleft \bA_x'$ is a nonempty stable subalgebra, and if $p$ is any cycle from $x$ to $x$ in $\fX$ such that $\bB \cap (\bB + p) \ne \emptyset$, then $\bB \cap (\bB + p')$ is also nonempty, where $p'$ is the corresponding path in $\fX'$. For this, consider the path instance $\mathbf{P}$ we get by unrolling the path $p$ in $\fX$, and let
\[
\RR \le_{sd} \bA_{x_0} \times \cdots \times \bA_{x_n}
\]
be the solution set to $\mathbf{P}$, with $x_0 = x_n = x$. Now apply Proposition \ref{prop-stable-essential} to $\RR$, setting $\bB_{x_0} = \bB_{x_n} = \bB$ and $\bB_{x_i} = \bA_{x_i}'$ for $i \ne 0,n$.
\end{proof}

Putting these results together, we see that weak consistency implies that a solution exists in bounded width CSPs.

\begin{thm}\label{thm-weakly-consistent} If $\fX$ is a weakly consistent instance of $\CSP(\bA_1, ..., \bA_n)$, where the $\bA_i$ are finite bounded width algebras, then $\fX$ has a solution. In fact, in this case $\fX$ has a \emph{stable} solution (i.e. a solution in which each variable $x$ is assigned to a stable element of its variable domain $\bA_x$).
\end{thm}

Of course, this all hinged on the existence of a stability concept: we still need to prove Theorem \ref{thm-stability}. Our construction of a stability concept won't be particularly elegant, but it will get the job done. First we will show that we can restrict to the case where every binary absorbing algebra is centrally absorbing.

\begin{prop} Suppose that $\cV$ is a pseudovariety of finite idempotent algebras, and that $\cV'$ is an affine-free reduct of $\cV$. If $\Yleft$ is a stability concept on $\cV'$, then the restriction of $\Yleft$ to $\cV$ is also a stability concept.
\end{prop}
\begin{proof} The only nontrivial axiom to check is ubiquity. For this, note that since $\cV'$ is affine-free, we can apply Proposition \ref{prop-stable-element} to see that every $\bA \in \cV'$ has a stable element $\{a\} \Yleft \bA$. Since $\cV$ is idempotent, we also have $\{a\} \le \bA$ in $\cV$.
\end{proof}

\begin{prop} A pseudovariety $\cV$ of finite algebras has a stability concept iff every finitely generated subvariety of $\cV$ has a stability concept.
\end{prop}
\begin{proof} This is an application of K\"onig's Lemma: for each algebra $\bA$ in $\cV$ with $|\bA| > 1$ which has no affine quotients, we need to choose at least one proper subalgebra $\bB$ to be a stable subalgebra of $\bA$. Since there are only a finite number of choices for $\bB$ for each finite $\bA$, if we can make a consistent set of choices for every finite collection of algebras $\bA_1, ..., \bA_n$, then there exists a globally consistent set of choices.
\end{proof}

\begin{defn} We say that a pseudovariety is \emph{strongly prepared} if every binary absorbing algebra $\bB \lhd_{bin} \bA \in \cV$ is also strongly absorbing.
\end{defn}

\begin{prop} Every locally finite variety of idempotent bounded width algebras has a strongly prepared, bounded width reduct.
\end{prop}
\begin{proof} This follows by repeatedly applying Proposition \ref{bin-abs-semi} and Proposition \ref{semilattice-preparation}, and using the fact that the set of two-variable terms can only be shrunk finitely many times if the free algebra on two generators is finite.
\end{proof}

To finish the proof of Theorem \ref{thm-stability}, we only need to construct stability concepts on pseudovarieties of idempotent, strongly prepared Taylor algebras.

\begin{defn}\label{defn-strong-stable} If $\cV$ is a pseudovariety of finite idempotent strongly prepared Taylor algebras, then we say that $\bB \Yleft \bA$ if there is a sequence of subalgebras
\[
\bA = \bA_0 \ge \bA_1 \ge \cdots \ge \bA_n = \bB
\]
such that for each $i$, one of the following is true:
\begin{itemize}
\item $\bA_{i+1}$ contains a strongly absorbing subalgebra of $\bA_i$,
\item $\bA_{i+1}$ centrally absorbs $\bA_i$, or
\item there is a congruence $\theta$ on $\bA_i$ such that $\bA_i/\theta$ is polynomially complete, binary absorption-free and central absorption-free, and such that $\bA_{i+1}$ is a congruence class of $\theta$.
\end{itemize}
We say that $\bB$ is \emph{stable in one step} if we can take $n = 1$ in the above. Following \cite{zhuk-strong}, if the third bullet point above holds for $\bB$, then we say that $\bB$ is a \emph{PC subalgebra} of $\bA$ with \emph{PC congruence} $\theta$.
\end{defn}

Algebras which are stable in one step are almost the same thing as what Zhuk calls \emph{strong} subalgebras in \cite{zhuk-strong} - the only difference is in how we handle the case of binary absorption.

\begin{thm} If $\cV$ is a pseudovariety of finite idempotent strongly prepared Taylor algebras, then the binary relation $\Yleft$ on $\cV$ from Definition \ref{defn-strong-stable} is a stability concept.
\end{thm}
\begin{proof} We just need to verify the axioms for $\Yleft$. Obviously $\Yleft$ satisfies the subalgebra axiom, and by Corollary \ref{zhuk-four-cases} and the assumption that $\cV$ is strongly prepared $\Yleft$ satisfies ubiquity as well. Transitivity holds for $\Yleft$ by construction. The remaining axioms are intersection, propagation, and the Helly property.

To verify the intersection axiom, suppose that $\bB, \bC \Yleft \bA$: we need to check that $\bB \cap \bC \Yleft \bB$. Inducting on $|\bA| + |\bB| + |\bC|$, we see that it's enough to prove this when $\bB$ and $\bC$ are both stable in one step. For this, we divide into four cases: either $\bC \lhd_Z \bA$ (1), $\bC$ contains a strongly absorbing subalgebra of $\bA$ (2), $\bB$ contains a centrally absorbing subalgebra of $\bA$ and $\bC$ is a PC subalgebra (3), or each of $\bB, \bC$ is a PC subalgebra (4). Case (1) follows from $\bB \cap \bC \lhd_Z \bB$. For case (2), we need a claim which we will also use elsewhere.

{\bf Claim:} If $\bS \lhd_{str} \bA$ and $\bB \Yleft \bA$, then $\bB \cap \bS \ne \emptyset$ and $\bB \cap \bS \lhd_{str} \bB$.

{\bf Proof of claim:} We just need to check this when $\bB$ is stable in one step. If $\bB$ contains a strongly absorbing subalgebra of $\bA$, we can apply Proposition \ref{prepared-bin-abs}. If $\bB \lhd_Z \bA$, then we can apply Proposition \ref{prepared-bin-abs} and Proposition \ref{prop-central-closed}. If $\bB$ is a PC subalgebra with PC congruence $\theta$, then $\bS/\theta \lhd_{str} \bA/\theta$ implies that $\bS/\theta = \bA/\theta$, so $\bS$ meets every congruence class of $\theta$, and in particular $\bB \cap \bS \ne \emptyset$.

{\bf Case (2) for intersection axiom:} Let $\bC'$ be a subalgebra of $\bC$ such that $\bC' \lhd_{str} \bA$. Then $\bC' \cap \bB \ne \emptyset$ and $\bB \cap \bC' \lhd_{str} \bB$, so $\bB \cap \bC$ contains a strongly absorbing subalgebra of $\bB$.

{\bf Case (3) for intersection axiom:} Let $\bB'$ be a subalgebra of $\bB$ such that $\bB' \lhd_Z \bA$, and let $\theta$ be the PC congruence for $\bC$. Then $\bB'/\theta \lhd_Z \bA/\theta$ by Proposition \ref{prop-central-abs-quotient}, so since $\bA/\theta$ is central absorption-free, we must have $\bB'/\theta = \bA/\theta$. Since $\bB' \le \bB \le \bA$, we must have $\bB/\theta = \bA/\theta$ as well, so $\bC$ is a PC subalgebra of $\bB$ with PC congruence $\theta|_\bB$.

{\bf Case (4) for intersection axiom:} Suppose that $\bB$ has PC congruence $\theta$ and $\bC$ has PC congruence $\psi$. We can consider $\bA/(\theta \wedge \psi)$ as a binary subdirect relation:
\[
\bA/(\theta \wedge \psi) \le_{sd} (\bA/\theta) \times (\bA/\psi).
\]
Since $\bA/\theta, \bA/\psi$ are both simple (all polynomially complete algebras are simple), this binary relation is either linked or is the graph of an isomorphism. If $\bA/(\theta \wedge \psi)$ is the graph of an isomorphism, then we must have $\theta = \psi$, so if $\bB \cap \bC \ne \emptyset$ then $\bB = \bC$. Otherwise, if $\bA/(\theta \wedge \psi)$ is linked, then by the Absorption Theorem \ref{absorption-theorem} we must have
\[
\bA/(\theta \wedge \psi) = (\bA/\theta) \times (\bA/\psi),
\]
since neither $\bA/\theta$ nor $\bA/\psi$ has a binary or centrally absorbing subalgebras. Thus $\bB/\psi = \bA/\psi$, so $\bB \cap \bC$ is a PC subalgebra of $\bB$ with PC congruence $\psi|_\bB$.

That completes the proof of the intersection axiom. For the propagation axiom, the pullback part is almost immediate from the definition. For the pushforward part of the propagation axiom, we can suppose that $\bA \twoheadrightarrow \bA/\theta$ and that $\bC$ is stable in one step in $\bA$. If $\bC \supseteq \bC' \lhd_{str} \bA$, then $\bC'/\theta \lhd_{str} \bA/\theta$. If $\bC \lhd_Z \bA$, then by Proposition \ref{prop-central-abs-quotient} we have $\bC/\theta \lhd_Z \bA/\theta$. The tricky case is the case where $\bC$ is a PC subalgebra of $\bA$ with PC congruence $\psi$. In this case, we consider $\bA/(\theta \wedge \psi)$ as a binary subdirect relation:
\[
\bA/(\theta \wedge \psi) \le_{sd} (\bA/\theta) \times (\bA/\psi).
\]
Since $\bA/\psi$ is simple, this binary relation is either the graph of a homomorphism $(\bA/\theta) \twoheadrightarrow (\bA/\psi)$ or is linked. If it is the graph of a homomorphism, then we must have $\theta \le \psi$, so $\bC/\theta$ is a PC subalgebra of $\bA/\theta$ with PC congruence $\psi/\theta$. Now suppose that it is linked, and let $\bS$ be a minimal strongly absorbing subalgebra of $\bA/\theta$, which exists by Proposition \ref{prepared-bin-abs}. Then by Theorem \ref{absorbing-linked}, the binary relation
\[
\bA/(\theta \wedge \psi) \cap (\bS \times (\bA/\psi)) \le_{sd} \bS \times (\bA/\psi)
\]
is also linked, and then by the Absorption Theorem \ref{absorption-theorem} it must be the full relation, since $\bS$ has no binary absorbing subalgebras and $\bA/\psi$ has no binary/centrally absorbing subalgebras. This means that $\bC/\theta$ contains $\bS$, so $\bC/\theta$ is a stable subalgebra of $\bA/\theta$.

To finish, we just need to verify the Helly axiom, which states that if $\bB,\bC,\bD \Yleft \bA$ and each pair has a nonempty intersection, then $\bB \cap \bC \cap \bD \ne \emptyset$. We induct on
\[
|\bA| + |\bA\setminus \bB| + |\bA \setminus \bC| + |\bA \setminus \bD|.
\]
Suppose that one of $\bB,\bC,\bD$ is not stable in one step, say $\bB$. Then there is some $\bA'$ with $\bB \Yleft \bA' \Yleft \bA$ such that $|\bA \setminus \bA'| < |\bA \setminus \bB|$ and $|\bA'| < |\bA|$. By the induction hypothesis, we have $\bA' \cap \bC \cap \bD \ne \emptyset$. Set $\bC' = \bA' \cap \bC$, $\bD' = \bA' \cap \bD$, and $\bB' = \bB$. Then by the intersection property we have $\bB', \bC', \bD' \Yleft \bA'$, and we have
\[
\bC' \cap \bD' = \bA' \cap \bC \cap \bD \ne \emptyset,
\]
while
\[
\bB' \cap \bC' = \bB \cap \bC \ne \emptyset
\]
and similarly $\bB' \cap \bD' \ne \emptyset$. Applying the induction hypothesis again, we see that
\[
\bB \cap \bC \cap \bD = \bB' \cap \bC' \cap \bD' \ne \emptyset.
\]
So we only need to check the Helly property when each of $\bB, \bC, \bD$ is stable in one step.

Suppose that one of $\bB, \bC, \bD$ contains a strongly absorbing subalgebra $\bS \lhd_{str} \bA$, say $\bB \supseteq \bS$. Then $\bS \cap \bC$ is a nonempty strongly absorbing subalgebra of $\bC$ by the earlier claim, and applying that claim again together with the fact that $\bC \cap \bD \Yleft \bC$, we see that $\bS \cap \bC \cap \bD \ne \emptyset$. So we may assume that each of $\bB, \bC, \bD$ is either centrally absorbing or is a PC subalgebra.

If all three of $\bB, \bC, \bD$ are centrally absorbing, then by Corollary \ref{cor-central-ternary} there is a ternary term $t$ such that each of $\bB, \bC, \bD$ absorbs $\bA$ with respect to $t$. If we pick $x \in \bB \cap \bC, y \in \bC \cap \bD, z \in \bB \cap \bD$, then we must have $t(x,y,z) \in \bB \cap \bC \cap \bD$.

If two of $\bB,\bC,\bD$ are centrally absorbing, suppose that $\bB, \bC$ are centrally absorbing and $\bD$ is a PC subalgebra with PC congruence $\theta$. Then $\bB \cap \bC \lhd_Z \bA$, so we must have $(\bB\cap \bC)/\theta \lhd_Z \bA/\theta$, so $(\bB \cap \bC)/\theta = \bA/\theta$. Thus $\bB \cap \bC$ intersects $\bD$.

If one of $\bB,\bC,\bD$ is centrally absorbing, we may suppose that $\bB \lhd_Z \bA$ and that $\bC, \bD$ are PC subalgebras with PC congruences $\theta, \psi$. As before, we must have $\bB/\theta = \bA/\theta$ and $\bB/\psi = \bA/\psi$. Considering $\bA/(\theta \wedge \psi)$ as a subdirect binary relation on $(\bA/\theta)\times (\bA/\psi)$, we see that it must be linked if $\bC \ne \bD$, since both of $\bA/\theta, \bA/\psi$ are simple. Then by Theorem \ref{absorbing-linked}, the binary relation
\[
\bB/(\theta \wedge \psi) \le_{sd} (\bA/\theta) \times (\bA/\psi)
\]
is also linked, so by the Absorption Theorem \ref{absorption-theorem} it must be the full relation. Thus we have $\bB \cap \bC \cap \bD \ne \emptyset$.

Finally, suppose that all three of $\bB, \bC, \bD$ are PC subalgebras, with PC congruences $\theta, \psi, \eta$. Then we can think of $\bA/(\theta\wedge \psi \wedge \eta)$ as a ternary subdirect relation:
\[
\bA/(\theta\wedge \psi \wedge \eta) \le_{sd} (\bA/\theta) \times (\bA/\psi) \times (\bA/\eta).
\]
If no two of $\bB, \bC, \bD$ are equal to each other, then every binary projection of this relation is linked, so every binary projection is full by the Absorption Theorem \ref{absorption-theorem}. If this ternary relation is the full relation, then we have $\bB \cap \bC \cap \bD \ne \emptyset$. Otherwise, pick some $x \in \bA/\theta$ such that the binary relation
\[
\pi_{23}(\bA/(\theta\wedge \psi \wedge \eta) \cap (\{x\} \times (\bA/\psi) \times (\bA/\eta))) \le_{sd} (\bA/\psi) \times (\bA/\eta)
\]
is not the full relation. Applying the Absorption Theorem \ref{absorption-theorem} again, we see that the binary relation above must be the graph of an isomorphism between $\bA/\psi$ and $\bA/\eta$. A similar argument shows that $\bA/\theta$ is isomorphic to $\bA/\psi$. Thus we can apply Theorem \ref{absorption-free} to show that all three of $\bA/\theta, \bA/\psi, \bA/\eta$ are affine, which contradicts the assumption that they are polynomially complete.
\end{proof}

\begin{prob} Is there a less ad-hoc stability concept?
\end{prob}

\subsection{Ramsey-theoretic upgrade: vague solutions imply solvability}

Unsatisfyingly, even weak consistency is too demanding to directly prove the existence of a $4$-ary Siggers term satisfying the identity $t(x,x,y,z) \approx t(y,z,z,x)$. Using Ramsey's theorem, we can cure this particular defect. The material in this subsection is based on the theory developed in \cite{brady-chromatic}.

In this subsection, we are mainly concerned with the following question: given an instance $\fX$ whose variable domains and relations are not assumed to be closed under the basic operations of our bounded width algebra, under what circumstances can we guarantee that $\Sg(\fX)$ has a solution? If the variable domains of $\fX$ consist of generating sets for free algebras, then this question is equivalent to asking which systems of height $1$ identities can be solved in every finite bounded width algebra. We will show that if the instance $\fX$ has a ``vague'' solution, then $\Sg(\fX)$ is guaranteed to have a solution in any finite bounded width algebra which contains the variable domains of $\fX$.

\begin{defn}\label{defn-vague} Let $\cP_\emptyset(A)$ be the set of non-empty subsets of a set $A$. A \emph{vague element} of $A$ is defined to be a total quasiorder $\preceq$ on $\cP_\emptyset(A)$ such that there is no pair of disjoint subsets $S, T \subset A$ with $S \sim T$, where $S \sim T$ means that $S \preceq T$ and $T \preceq S$.

If $R \subseteq_{sd} A_1 \times \cdots \times A_n$ is subdirect, then we say that a collection of vague elements $\preceq_i$ of the $A_i$s \emph{vaguely satisfies} the relation $R$ if there exists a total quasiorder $\preceq_R$ on the disjoint union
\[
\cP_\emptyset(A_1) \sqcup \cdots \sqcup \cP_\emptyset(A_n)
\]
such that the restriction of $\preceq_R$ to $\cP_\emptyset(A_i)$ is $\preceq_i$ for each $i$, and such that for each $i,j \in [n]$ and each nonempty $S \subseteq A_i$, we have
\[
S \preceq_R S+\pi_{ij}(R).
\]
If $\fX$ is an arc-consistent instance with variable domains $A_x$, then a collection of vague elements $\preceq_x$ of the $A_x$s is a \emph{vague solution} to $\fX$ if it vaguely satisfies every relation of $\fX$.
\end{defn}

Total quasiorders are also known as \emph{preference relations} - so a vague element of $A$ is a preference relation on the nonempty subsets of $A$ which our hypothetical element might live in, which avoids being caught out as incoherent by requiring that any pair of equally preferable subsets has a nonempty intersection. Vague solutions are closely connected to weak consistency.

\begin{prop} Every weakly consistent instance has a vague solution. In fact, there is always a vague solution where each vague element $\preceq_x$ extends the inclusion order $\subseteq$ on $\cP_\emptyset(A_x)$ (possibly identifying some subsets with each other as well).
\end{prop}
\begin{proof} Suppose $\fX$ is a weakly consistent instance with variable domains $A_x$, and define a quasiorder $\preceq_0$ on
\[
\bigsqcup_x \cP_\emptyset(A_x)
\]
by
\[
S \preceq_0 S+p
\]
for every path $p$ in $\fX$. Let $\preceq$ be any extension of $\preceq_0$ to a total quasiorder which does not identify any pair of sets which were not already identified by $\preceq_0$. Then if we take $\preceq_x$ to be the restriction of $\preceq$ to $\cP_\emptyset(A_x)$ for each variable $x$, we see that each $\preceq_x$ is a vague element of $A_x$ (since $\fX$ is weakly consistent) and that the collection of vague elements $\preceq_x$ is a vague solution to $\fX$.

For the second claim, we note that if we add additional tuples to any relation of a weakly consistent instance then it remains weakly consistent by Proposition \ref{prop-weakly-consistent-cycle}. If we add extra equality relations from a variable to itself then clearly the instance remains weakly consistent as well, so we see that we may add in the binary relation
\[
\Delta_{A_x} \cup \{(a,b)\} \subseteq_{sd} A_x \times A_x
\]
for any $a,b \in A_x$ without causing the instance to stop being weakly consistent. If we add all such binary relations in, then we see that the quasiorder $\preceq_0$ has $(x,S) \preceq_0 (x,T)$ for any $S \subseteq T$, so the same will be true in the extension $\preceq$ of $\preceq_0$.
\end{proof}

The converse isn't true - there are instances which have vague solutions, but which are not weakly consistent.

\begin{ex} Consider the Siggers instance, which has one variable $u$ with domain $A_u = \{x,y,z\}$, and one binary relation $R \subseteq_{sd} A_u \times A_u$ given by
\[
R = \left\{\begin{bmatrix}x\\y\end{bmatrix},\begin{bmatrix}x\\z\end{bmatrix},\begin{bmatrix}y\\z\end{bmatrix},\begin{bmatrix}z\\x\end{bmatrix}\right\}.
\]
This instance is \emph{not} weakly consistent: we have
\[
\{z\} + R = \{x\}, \;\; \{x\} - R = \{z\}
\]
and $\{x\} \cap \{z\} = \emptyset$. Nevertheless, the Siggers instance does have a vague solution: take $\preceq_u$ to be the total order
\[
\{y\} \prec_u \{x\} \prec_u \{z\} \prec_u \{x,y\} \prec_u \{y,z\} \prec_u \{x,z\} \prec_u \{x,y,z\}.
\]
To see that $(\preceq_u,\preceq_u)$ vaguely satisfies $R$, use the following total quasiorder $\preceq_R$ on $[2] \times \cP_\emptyset(A_u)$:
\begin{align*}
(1,\{y\}) &\prec_R (2,\{y\}) \prec_R (1,\{x\}) \prec_R (2,\{x\}) \sim_R (1,\{z\}) \prec_R (2,\{z\})\\
&\prec_R (2,\{x,y\}) \prec_R (1,\{x,y\}) \sim_R (2,\{y,z\}) \prec_R (1,\{y,z\}) \prec_R (2,\{x,z\}) \prec_R (1,\{x,z\})\\
&\prec_R (1,\{x,y,z\}) \sim_R (2,\{x,y,z\}).
\end{align*}
\end{ex}

In order to show that vaguely solvable instances have solutions (after taking the closure by algebraic operations), we will construct a very large weakly consistent instance with many copies of each variable and relation from the original instance. In fact, by K\"onig's Lemma we may even allow ourselves to build an infinitely large weakly consistent instance (although we will only truly need a finite portion of it). Roughly speaking, we will take the variables and relations of this instance to be indexed by subsets of $\NN$ of certain fixed sizes, which will put us in a position to apply Ramsey's theorem for hypergraphs.

\begin{defn} If $f : S \rightarrow \NN$ is any function, then we define the \emph{associated total quasiorder} $\preceq_f$ on $S$ by
\[
a \preceq_f b \;\; \iff \;\; f(a) \le f(b).
\]
If $\preceq$ is a vague element of $A$, then we say that a function $f : \cP_\emptyset(A) \rightarrow \NN$ is \emph{compatible} with $\preceq$ if $\preceq_f = \preceq$.

More generally, if $R \subseteq_{sd} A_1 \times \cdots \times A_n$ is a subdirect relation with a vague solution $(\preceq_1, ..., \preceq_n)$, then we say that a function
\[
f : \bigsqcup_i \cP_\emptyset(A_i) \rightarrow \NN
\]
is \emph{compatible} with $R$ and $(\preceq_1, ..., \preceq_n)$ if $\preceq_f$ can be used as the total quasiorder $\preceq_R$ from the definition of vague satisfaction.
\end{defn}

\begin{defn} If $\fX$ is an arc-consistent instance with a vague solution given by vague elements $\preceq_x$ of the variable domains $A_x$, then we define the \emph{associated weakly consistent instance} $\fX^*$ as follows:
\begin{itemize}
\item for each variable $x$ of $\fX$ and each $f : \cP_\emptyset(A_x) \rightarrow \NN$ which is compatible with $\preceq_x$, we have a variable $(x,f)$ of $\fX^*$, and
\item for each constraint relation $R \subseteq_{sd} A_{x_1} \times \cdots \times A_{x_n}$ of $\fX$ and each $f : \bigsqcup_i \cP_\emptyset(A_{x_i}) \rightarrow \NN$ which is compatible with $R$ and $(\preceq_{x_1}, ..., \preceq_{x_n})$, we impose
\[
((x_1,f|_{\cP_\emptyset(A_{x_1})}), ..., (x_n,f|_{\cP_\emptyset(A_{x_n})})) \in R
\]
as a constraint in $\fX^*$.
\end{itemize}
\end{defn}

\begin{prop} If $\fX$ has a vague solution, then the associated weakly consistent instance $\fX^*$ is indeed weakly consistent.
\end{prop}
\begin{proof} By the construction of $\fX^*$, if there is a path $p$ from $(x,f)$ to $(y,g)$ in $\fX^*$ and if $S \subseteq A_x$, then we have
\[
f(S) \le g(S+p).
\]
Thus if we have $S + p + q = S$ and $p,q$ are cycles from $(x,f)$ to $(x,f)$, then
\[
f(S) \le f(S+p) \le f(S+p+q) = f(S)
\]
implies that $f(S) = f(S+p)$. Then since $f$ is compatible with $\preceq_x$ we must have $S \sim_x S+p$, so $S \cap (S+p) \ne \emptyset$ since $\preceq_x$ is a vague element.
\end{proof}

\begin{prop} If $\fX^*$ is weakly consistent and its variable domains are contained in a finite bounded width algebra, then $\Sg(\fX^*)$ has a stable solution.
\end{prop}
\begin{proof} By Proposition \ref{prop-weakly-consistent-cycle}, an instance is weakly consistent if and only if certain cycles take every element of the variable domain back to themselves, so if $\fX^*$ is weakly consistent then so is $\Sg(\fX^*)$. By K\"onig's Lemma, in order to check that $\Sg(\fX^*)$ has a solution we just need to check that every finite subinstance of $\Sg(\fX^*)$ has a stable solution, and for this we can apply the main result of the previous section.
\end{proof}

\begin{thm}\label{thm-vague-solution} If $\fX$ is an arc-consistent instance which has a vague solution, and if the variable domains of $\fX$ are contained in a finite bounded width algebra, then $\Sg(\fX)$ has a stable solution.
\end{thm}
\begin{proof} Let $\fX^*$ be the associated weakly consistent instance. By the previous proposition, $\Sg(\fX^*)$ has a stable solution. Fix one particular stable solution to $\Sg(\fX^*)$.

To finish, we imagine ``coloring'' the compatible functions $f : \cP_\emptyset(A_x) \rightarrow \NN$ by the values that the variables $(x,f)$ are assigned to in our solution of $\Sg(\fX^*)$. Since the variable domains in the instance $\Sg(\fX^*)$ are finite, we only have finitely many colors to choose from, so Ramsey's theorem for hypergraphs implies that there is an infinite subset $S \subseteq \NN$ such that each compatible function $f : \cP_\emptyset(A_x) \rightarrow S$ has the same color. Iterating this for each variable $x$ of the instance $\fX$, we finally find an infinite subset $U \subseteq \NN$ such that for each variable $x$ of $\fX$ and each compatible $f : \cP_\emptyset(A_x) \rightarrow U$, the value assigned to $(x,f)$ in our solution to $\Sg(\fX^*)$ only depends on $x$ and does not depend on $f$.

We claim that assigning the variable $x$ of $\fX$ to the value assigned to any such $(x,f)$ (with $f : \cP_\emptyset(A_x) \rightarrow U$ compatible with $\preceq_x$) in our solution to $\Sg(\fX^*)$ solves the instance $\Sg(\fX)$. To see this, let $R \subseteq A_{x_1} \times \cdots \times A_{x_n}$ be any constraint relation of $\fX$. Then since $(\preceq_{x_1}, ..., \preceq_{x_n})$ vaguely satisfies $R$, there is some total quasiorder $\preceq_R$ as in the definition of vague satisfaction. Since $U \subseteq \NN$ is infinite (or just sufficiently large), we can then find a function $f : \bigsqcup_i \cP_\emptyset(A_{x_i}) \rightarrow U$ such that $\preceq_f = \preceq_R$, and for this $f$ we see that
\[
((x_1,f|_{\cP_\emptyset(A_{x_1})}), ..., (x_n,f|_{\cP_\emptyset(A_{x_n})})) \in \Sg(R)
\]
is a constraint of $\Sg(\fX^*)$. Therefore, the tuple of values assigned to $(x_1, ..., x_n)$ by this procedure satisfies the constraint $\Sg(R)$ of the instance $\Sg(\fX)$.
\end{proof}

\begin{rem} The same Ramsey argument can be used to show that if an instance has a vague solution, then it has a vague solution where each vague element $\preceq_x$ extends the inclusion order $\subseteq$ on $\cP_\emptyset(A_x)$. A refinement of this argument shows that we can also assume that our vague elements $\preceq_x$ have the following property: whenever $S \preceq_x T$, we also have $A_x\setminus T \preceq_x A_x \setminus S$.
\end{rem}

We can use stability to upgrade this result further: we don't need to find a vague solution to the full instance $\fX$, it's enough to find a vague solution to just the binary part of $\fX$.

\begin{defn} If $\fX$ is an instance, then we define the \emph{binary part} of $\fX$ to be the instance $\fX^{bin}$ given by replacing each $k$-ary constraint relation $R \subseteq A_{x_1} \times \cdots \times A_{x_k}$ by the collection of binary relations $\pi_{ij}(R) \subseteq A_{x_i} \times A_{x_j}$ for $i,j \in [k]$.
\end{defn}

\begin{cor}\label{cor-vague-bin} If $\fX$ is an arc-consistent instance such that $\fX^{bin}$ has a vague solution, and if the variable domains of $\fX$ are contained in a finite bounded width algebra, then $\Sg(\fX)$ has a stable solution.
\end{cor}
\begin{proof} By Proposition \ref{prop-stable-essential} and arc-consistency, any stable solution to $\Sg(\fX^{bin})$ is also a stable solution to $\Sg(\fX)$.
\end{proof}

\begin{ex} There is an arc-consistent instance $\fX$ such that $\fX^{bin}$ has a vague solution but $\fX$ does not: take the $5$-ary relation $R \subseteq_{sd} \{a,b,c,d\}^5$ given by
\[
R = \left\{\begin{bmatrix}a\\a\\b\\c\\d\end{bmatrix},\begin{bmatrix}a\\b\\c\\d\\a\end{bmatrix},\begin{bmatrix}b\\c\\d\\a\\a\end{bmatrix},\begin{bmatrix}c\\d\\a\\a\\b\end{bmatrix},\begin{bmatrix}d\\a\\a\\b\\c\end{bmatrix}\right\},
\]
and let $\fX$ be the instance with a single variable $x$ which is supposed to satisfy the constraint $(x,x,x,x,x) \in R$. To see that $\fX$ has no vague solution, it's enough to consider the relative order of the singleton sets $\{b\}$ and $\{c\}$. On the other hand, $\fX^{bin}$ has a vague solution where all five vague elements are given by
\[
\{d\} \prec \{c\} \prec \{b\} \prec \{c,d\} \prec \{b,d\} \prec \{a\} \prec \{b,c\} \prec \{a,d\} \sim \cdots \sim \{a,b,c,d\}.
\]
The fact that examples like this exist seems to be a hint that vague solutions are nowhere near to being the last word on bounded width CSPs.
\end{ex}

\begin{prob} Does every finite bounded width algebra have a $7$-ary term operation $t$ satisfying
\begin{align*}
t(a,a,b,c,d,e,f) &\approx t(f,a,a,b,c,d,e) \approx t(e,f,a,a,b,c,d)?
\end{align*}
\end{prob}

\begin{prob} If the sizes of the variable domains in $\fX$ are bounded by a constant (perhaps just $3$), how hard is it to determine whether $\fX^{bin}$ has a vague solution?
\end{prob}

We can at least make it a bit simpler for ourselves to check that a particular assignment of vague elements really gives us a vague solution.

\begin{prop} If $R \subseteq_{sd} A_x \times A_y$ and $\preceq_x, \preceq_y$ are vague elements of $A_x, A_y$, respectively, then the pair $(\preceq_x,\preceq_y)$ vaguely satisfies $R$ iff the following implication holds for all $S \subseteq A_x, T \subseteq A_y$:
\[
S + R \preceq_y T \;\; \wedge \;\; T - R \preceq_x S \;\; \implies \;\; S + R \sim_y T \;\; \wedge \;\; T - R \sim_x S.
\]
\end{prop}
\begin{proof} Consider the digraph on $\cP_\emptyset(A_x) \sqcup \cP_\emptyset(A_y)$ with an edge from $S \subseteq A_x$ to $T \subseteq A_y$ whenever $S + R \preceq_y T$, and similarly with an edge from $T \subseteq A_y$ to $S \subseteq A_x$ whenever $T - R \preceq_x S$. We just need to check that for every cycle $(S_1, T_1, S_2, T_2, ..., S_k, T_k)$ of this digraph, all $S_i$s are related by $\sim_x$ and all $T_j$s are related by $\sim_y$. For this, suppose that $S_i$ is $\preceq_x$-minimal and that $T_j$ is $\preceq_y$-minimal. Then from
\[
S_j + R \preceq_y T_j \preceq_y T_{i-1}
\]
and
\[
T_{i-1} - R \preceq_x S_i \preceq_x S_j,
\]
we see that $S_i \sim_x S_j$ and $T_{i-1} \sim_y T_j$, so $S_j$ is $\preceq_x$-minimal and $T_{i-1}$ is $\preceq_y$-minimal. Applying the same reasoning to $S_j, T_{i-1}$, we see that $S_{i-1}$ is $\preceq_x$-minimal and $T_{j-1}$ is $\preceq_y$-minimal. Continuing in this fashion, we see that all of the $S$s are related by $\sim_x$ and all of the $T$s are related by $\sim_y$, which is what we had to check.
\end{proof}

As our first illustration of the theory, we can show the existence of a Siggers term which satisfies a strong collection of additional identities.

\begin{prop} A finite algebra $\bA$ has bounded relational width if and only if it has a $4$-ary term $t$ which satisfies the identities
\[
t(x,x,y,z) \approx t(y,z,z,x) \approx t(z,x,y,x)
\]
and
\[
t(x,y,x,z) \approx t(x,z,y,x) \approx t(y,z,x,x)
\]
simultaneously.
\end{prop}
\begin{proof} It is easy to see that the identities satisfied by $t$ imply that the ternary terms $f, g$ defined by
\[
g(x,y,z) \coloneqq t(x,x,y,z), \;\;\; f(x,y,z) \coloneqq t(x,y,x,z)
\]
satisfy the equations
\[
g(x,x,y) \approx g(x,y,x) \approx g(y,x,x) \approx f(x,x,y) \approx f(x,y,x) \approx f(x,y,y).
\]
from Theorem \ref{bounded-width-terms}, so if such a $t$ exists then $\bA$ has bounded relational width.

Now suppose that $\bA$ has bounded relational width. Let $R$ be the following $6$-ary relation on $\{x,y,z\}$:
\[
R = \left\{\begin{bmatrix}x\\y\\z\\x\\x\\y\end{bmatrix},\begin{bmatrix}x\\z\\x\\y\\z\\z\end{bmatrix},\begin{bmatrix}y\\z\\y\\x\\y\\x\end{bmatrix},\begin{bmatrix}z\\x\\x\\z\\x\\x\end{bmatrix}\right\}.
\]
Consider the following instance $\fX$ on the six variables $a,b,c,d,e,f$:
\[
(a,b,c,d,e,f) \in R \;\; \wedge \;\; a = b = c \;\; \wedge \;\; d = e = f.
\]
If we consider the domain $\{x,y,z\}$ as a subset of the free algebra $\cF_\bA(x,y,z) \le \bA^{\bA^3}$ in the natural way, then we just need to show that $\Sg(\fX)$ has a solution. Since $\fX$ is arc-consistent, by Corollary \ref{cor-vague-bin} we just need to find a vague solution to the binary part $\fX^{bin}$. We assign the variables $a,b,c$ to the vague element $\preceq_g$ given by
\[
\{y\} \prec_g \{z\} \prec_g \{x\} \prec_g \{y,z\} \prec_g \{x,y\} \prec_g \{x,z\} \prec_g \{x,y,z\},
\]
and we assign the variables $d,e,f$ to the vague element $\preceq_f$ given by
\[
\{y\} \prec_f \{z\} \sim_f \{y,z\} \prec_f \{x\} \sim_f \{x,y\} \prec_f \{x,z\} \prec_f \{x,y,z\}.
\]
The reader may check that this assignment vaguely satisfies every binary projection of the relation $R$.
\end{proof}

As another illustration of the theory, we will show how the Loop Lemma \ref{loop-lemma} can be proved for finite bounded width algebras by constructing suitable vague solutions.

\begin{thm} If $R \subseteq A_x \times A_x$ is a smooth, weakly connected digraph of algebraic length $1$, then the instance $\fX$ which consists of only the variable $x$ and the constraint $(x,x) \in R$ has a vague solution. As a consequence, the instance $\Sg(\fX)$ has a stable solution in any finite bounded width algebra.
\end{thm}
\begin{proof} We will attempt to find a function $f : [2] \times \cP_\emptyset(A_x) \rightarrow \QQ$ such that for each proper nonempty $S \subset A_x$ we have
\[
|f(1,S) - f(2,S)| = 1,
\]
along with
\[
f(1,S) \le f(2,S+R)
\]
and
\[
f(2,S) \le f(1,S-R).
\]
To this end, we define a weighted directed graph $\cG$ with vertices corresponding to proper nonempty subsets $S \subset A_x$, and with an edge of weight $+1$ from $S$ to $S+R$ and an edge of weight $-1$ from $S$ to $S-R$ for each such $S$ (assuming $S+R, S-R \ne A_x$). We will handle each strongly connected component of $\cG$ separately.

We call a directed cycle of $\cG$ \emph{positive} if the sum of the weights along the cycle is strictly greater than $0$, and we define negative cycles similarly. For each positive directed cycle of $\cG$ from a vertex $S \subset A_x$ to $S$, there is a corresponding cycle $p$ of the instance $\fX$ which has strictly more $+\RR$ steps than $-\RR$ steps, with $S+p = S$, and we call such a cycle $p$ ``positive'' as well.

{\bf Claim.} No strongly connected component of $\cG$ contains both a positive directed cycle and a negative directed cycle.

{\bf Proof of claim.} Suppose otherwise. Then we can find a vertex $S \subset A_x$ of $\cG$, a positive cycle $p$, and a negative cycle $q$, such that
\[
S = S+p = S+q.
\]
We may assume without loss of generality that the total weights of $p$ and $q$ are opposite to each other, so $p+q$ has total weight $0$. Then we have
\[
S + R^{\circ j} - R^{\circ j} \subseteq S + jp + jq = S
\]
for all $j \ge 0$, so $S$ must be a union of linked components of $R^{\circ j}$ for all $j$. This contradicts Proposition \ref{prop-alg-length-linked}: some $R^{\circ j}$ must be linked if $R$ has algebraic length $1$.

Now suppose that $\cC$ is a strongly connected component of $\cG$ which does not contain any positive directed cycles. We will define the restriction of $f$ to $\cC$ such that
\[
f(2,S) = f(1,S) - 1
\]
for all $S \in \cC$. To do this, we pick any $S_0 \in \cC$ and any constant $c_\cC$, and define $f(1,T)$ to be $c_\cC$ plus the maximum total weight of any directed path from $S_0$ to $T$, for all $T \in \cC$. That this maximum total weight is well-defined follows from the fact that $\cC$ does not contain any positive directed cycles together with the finiteness of $\cC$. This definition is easily seen to satisfy
\begin{align*}
f(1,T) &\le f(2,T+R) = f(1,T+R)-1,\\
f(1,T)-1 &= f(2,T) \le f(1,T-R),
\end{align*}
so long as $T+R, T-R$ are in $\cC$. Additionally, if $f(1,T) = f(1,U)$ for some $T,U \in \cC$, then there is some $k$ such that $S_0$ has paths of total weight $k$ to each of $T$ and $U$ - in this case, we see that
\[
S + kR \subseteq T \cap U,
\]
so $T \cap U \ne \emptyset$.

We handle strongly connected components which do not have any negative directed cycles similarly, looking at the negative of the minimum total weight instead of the maximum total weight, and taking $f(2,S) = f(1,S) + 1$ on such components.

To finish, we pick any total order on the strongly connected components of $\cG$ which extends the reachability order, and we choose the constants $c_{\cC}$ for the various connected components $\cC$ according to this order, with sufficient distance between them that there is no interaction between the various strongly connected components of $\cG$.
\end{proof}

\section{Semidefinite Programming robustly solves bounded width CSPs}

In this section we finally touch on a difficult topic: trying to maximize the number of satisfied constraints in a CSP instance which has no perfect solution. We consider only a very special case of this problem here: the problem of trying to approximately solve a CSP when we are promised that there exists a way to satisfy all but a tiny fraction of the constraints. This problem was considered by Guruswami and Zhou in \cite{robust-horn-gap}.

\begin{defn}\label{defn-robust} We say that $\CSP(\fA)$ is \emph{robustly solvable} if there is a function $f : [0,1] \rightarrow [0,1]$ such that
\[
\lim_{\epsilon \rightarrow 0} f(\epsilon) = 0,
\]
and a polynomial time algorithm that takes as input an instance $\fX$ of $\CSP(\fA)$, and outputs an assignment to the variables of $\fX$ such that if it is possible to satisfy a $1-\epsilon$ fraction of the constraints of $\fX$, then the assignment found by the algorithm satisfies at least a $1-f(\epsilon)$ fraction of the constraints of $\fX$.
\end{defn}

Before we dive into our main topic, we first give evidence that certain CSPs are \emph{not} robustly solvable. We won't prove the next result here.

\begin{thm}[H\r{a}stad \cite{hastad-optimal}] Let $\fA$ be the affine CSP template with domain $\bA$, where $\bA$ is the idempotent reduct of any finite abelian group, with relations given by $\RR_c = \{(x,y,z) \mid x+y+z = c\} \le_{sd} \bA^3$ for every possible $c \in \bA$.

Then for every fixed $\epsilon > 0$, it is NP-hard to solve the following problem: given an instance $\fX$ of $\CSP(\fA)$ such that there exists an assignment satisfying at least a $1-\epsilon$ fraction of the constraints, find an assignment which satisfies at least a $\frac{1}{|\bA|} + \epsilon$ fraction of the constraints.
\end{thm}

Note that for the affine CSP defined above, randomly guessing values for variables will produce an assignment which satisfies a $\frac{1}{|\bA|}$ fraction of the constraints, on average. So H\r{a}stad's result tells us that it's NP-hard to find any improvement on randomly guessing, for affine CSPs which are not perfectly solvable.

\begin{cor} If $\CSP(\fA)$ is robustly solvable and $P \ne NP$, then $\fA$ must be affine-free (and therefore $\fA$ has bounded width).
\end{cor}

The best known approach to approximately solving CSPs, based on semidefinite programming, was laid out in Raghavendra's thesis \cite{raghavendra-thesis} (see \cite{raghavendra-optimal} for a short overview of the results). Under the Unique Games Conjecture, Raghavendra proved that this approach is actually optimal. The strategy is as follows.

As in the linear programming relaxation of a CSP, we imagine that we are looking for a probability distribution over solutions to the CSP. We do not give a full description of this unknown probability distribution: we only describe the marginal distribution over assignments to tuples of variables belonging to constraints of the CSP, as well as the marginal distribution over assignments to each pair of variables in the CSP. We impose compatibility conditions between the marginal distributions over each tuples of variables $(v_1, ..., v_m)$ belonging to some constraint and the marginal distribution over each pair $(v_i,v_j)$ for $i,j \le m$.

So far all the conditions given can be described by a system of linear inequalities. The semidefinite aspect comes from the following observation: every covariance matrix of any collection of random variables must be positive semidefinite.

To be more concrete, for each pair of variables $x,y$ and each pair of values $a \in \bA_x, b \in \bA_y$, we have some variable $p_{(x,a),(y,b)}$ between $0$ and $1$, describing the probability that $x$ is assigned the value $a$ and $y$ is assigned the value $b$. We create a matrix $M_p$ with rows and columns indexed by ordered pairs $(x,a)$ with $a \in \bA_x$, and fill the $(x,a),(y,b)$ entry with $p_{(x,a),(y,b)}$ (I like to imagine $M_p$ as a block matrix, with each block of rows or columns corresponding to a particular variable $x$). Then the matrix $M_p$ must be positive semidefinite if these probabilities come from an actual probability distribution.

Before defining everything formally, we give an example.

\begin{ex} Consider the following instance of 2-SAT: we have three variables $x,y,z$, and each pair of variables has a $\ne$ constraint imposed between them. This instance has no perfect solution, but the standard linear programming relaxation is incapable of noticing this. Let's see how the semidefinite relaxation does.

The matrix $M_p$ has six rows and six columns, corresponding to the pairs $(x,0),(x,1),(y,0)$, $(y,1),(z,0),(z,1)$, in that order. If $M_p$ comes from a probability distribution over perfect solutions to this instance of 2-SAT, then it must have the following shape:
\[
M_p = \begin{pmatrix}\begin{matrix} * & 0\\ 0 & *\end{matrix} & \rvline & \begin{matrix} 0 & *\\ * & 0\end{matrix} & \rvline & \begin{matrix} 0 & *\\ * & 0\end{matrix}\\
\hline \begin{matrix} 0 & *\\ * & 0\end{matrix} & \rvline & \begin{matrix} * & 0\\ 0 & *\end{matrix} & \rvline & \begin{matrix} 0 & *\\ * & 0\end{matrix}\\
\hline \begin{matrix} 0 & *\\ * & 0\end{matrix} & \rvline & \begin{matrix} 0 & *\\ * & 0\end{matrix} & \rvline & \begin{matrix} * & 0\\ 0 & *\end{matrix}\end{pmatrix}.
\]
Additionally, the entries in each block of $M_p$ must sum to $1$ (and be $\ge 0$), and for each fixed row or column of $M_p$, the sum of the entries in the intersection of the row/column with any block must only depend on the row/column. Putting these linear constraints together, we quickly see that every nonzero entry of $M_p$ must actually be equal to $\frac{1}{2}$. So far, this is exactly what the linear programming relaxation will guess.

The matrix $M_p$ found above, with all nonzero entries equal to $\frac{1}{2}$, is \emph{not} positive semidefinite. To see this, note that we have
\[
\begin{pmatrix}1 & -1 & \rvline & 1 & -1 & \rvline & 1 & -1\end{pmatrix}
\begin{pmatrix}\begin{matrix} \tfrac{1}{2} & 0\\ 0 & \tfrac{1}{2}\end{matrix} & \rvline & \begin{matrix} 0 & \tfrac{1}{2}\\ \tfrac{1}{2} & 0\end{matrix} & \rvline & \begin{matrix} 0 & \tfrac{1}{2}\\ \tfrac{1}{2} & 0\end{matrix}\\
\hline \begin{matrix} 0 & \tfrac{1}{2}\\ \tfrac{1}{2} & 0\end{matrix} & \rvline & \begin{matrix} \tfrac{1}{2} & 0\\ 0 & \tfrac{1}{2}\end{matrix} & \rvline & \begin{matrix} 0 & \tfrac{1}{2}\\ \tfrac{1}{2} & 0\end{matrix}\\
\hline \begin{matrix} 0 & \tfrac{1}{2}\\ \tfrac{1}{2} & 0\end{matrix} & \rvline & \begin{matrix} 0 & \tfrac{1}{2}\\ \tfrac{1}{2} & 0\end{matrix} & \rvline & \begin{matrix} \tfrac{1}{2} & 0\\ 0 & \tfrac{1}{2}\end{matrix}\end{pmatrix}
\begin{pmatrix} 1\\ -1\\ \hline 1\\ -1\\ \hline 1\\ -1\end{pmatrix} = -3 < 0.
\]
So the semidefinite relaxation of the problem can detect that we can't perfectly solve this instance of 2-SAT.

Now suppose that we give up on finding a perfect solution, and instead look for an approximate solution. This means that some of the entries of $M_p$ which were required to be $0$ before are instead required to be \emph{small}. One choice of $M_p$ that works is
\[
M_p = \frac{1}{8}\begin{pmatrix}\begin{matrix} 4 & 0\\ 0 & 4\end{matrix} & \rvline & \begin{matrix} 1 & 3\\ 3 & 1\end{matrix} & \rvline & \begin{matrix} 1 & 3\\ 3 & 1\end{matrix}\\
\hline \begin{matrix} 1 & 3\\ 3 & 1\end{matrix} & \rvline & \begin{matrix} 4 & 0\\ 0 & 4\end{matrix} & \rvline & \begin{matrix} 1 & 3\\ 3 & 1\end{matrix}\\
\hline \begin{matrix} 1 & 3\\ 3 & 1\end{matrix} & \rvline & \begin{matrix} 1 & 3\\ 3 & 1\end{matrix} & \rvline & \begin{matrix} 4 & 0\\ 0 & 4\end{matrix}\end{pmatrix},
\]
which the reader may verify is positive semidefinite. This seems to satisfy each particular constraint with a probability of $\frac{3}{4}$. So the semidefinite relaxation thinks it might be possible to satisfy a $\frac{3}{4}$ fraction of the constraints. An easy brute force search reveals that the best we can do in reality is to satisfy a $\frac{2}{3}$ fraction of the constraints.
\end{ex}

There is one further step we will take to analyze the semidefinite relaxation, based on a standard fact from linear algebra about positive semidefinite matrices.

\begin{prop} If $M$ is an $n\times n$ positive semidefinite matrix, then there is a collection of vectors $x_1, ..., x_n \in \RR^n$ such that $M_{ij} = x_i\cdot x_j$ for all $i,j \le n$. Such a collection of vectors $x_1, ..., x_n$ can be computed from $M$ in polynomial time.
\end{prop}
\begin{proof} Perhaps the simplest approach is to compute a Cholesky decomposition of $M$, writing $M = LL^T$ for some lower triangular matrix $L$. The columns of $L^T$ can then be used as the vectors $x_1, ..., x_n$.
\end{proof}

\begin{ex} The matrix $M_p$ from the end of the previous example is positive semidefinite, so there should exist vectors $x_0, x_1, y_0, y_1, z_0, z_1 \in \RR^6$ whose matrix of dot products is equal to $M_p$. Since $M_p$ has rank $3$, we should even be able to find such vectors in $\RR^3$. One particularly satisfying choice of vectors that works is
\[
x_0 = \frac{1}{\sqrt{24}}\begin{bmatrix}\sqrt{2} - \sqrt{3}\\ \sqrt{2}\\ \sqrt{2}+\sqrt{3}\end{bmatrix}, x_1 = \frac{1}{\sqrt{24}}\begin{bmatrix}\sqrt{2} + \sqrt{3}\\ \sqrt{2}\\ \sqrt{2}-\sqrt{3}\end{bmatrix}, y_0 = \frac{1}{\sqrt{24}}\begin{bmatrix}\sqrt{2}\\ \sqrt{2}+\sqrt{3}\\ \sqrt{2} - \sqrt{3}\end{bmatrix}, y_1 = \frac{1}{\sqrt{24}}\begin{bmatrix}\sqrt{2}\\ \sqrt{2}-\sqrt{3}\\ \sqrt{2} + \sqrt{3}\end{bmatrix}, 
\]
with $z_0, z_1$ defined similarly by cyclically shifting $y_0, y_1$.
\end{ex}

Now we can give the definition of the basic semidefinite relaxation of a CSP (this is the LC relaxation from \cite{raghavendra-thesis}).

\begin{defn} Given an instance $\fX$ of a CSP, with variable domains $\bA_v$ and constraints $C$ imposing relations $\RR_C \le \prod_{i \le m_C} \bA_{v_{C,i}}$ on the variables $v_{C,1}, ..., v_{C,m}$, the \emph{basic semidefinite relaxation} of $\fX$ is the following optimization problem. We wish to find a system of ``probabilities'' $p_{C,r}$ for $r \in \prod_{i \le m_C} \bA_{v_{C,i}}$, such that
\[
\sum_{r} p_{C,r} = 1
\]
for each constraint $C$ and
\[
p_{C,r} \ge 0
\]
for each $C,r$, and to find vectors
\[
x_a \in \RR^N
\]
for each variable $x$ and value $a \in \bA_x$, where $N = \sum_x |\bA_x|$ is the number of pairs $(x,a)$, such that for each $C$ and each pair of variables $x = v_{C,i}, y = v_{C,j}$ involved in the constraint $C$, we satisfy the compatibility condition
\[
x_a \cdot y_b = \sum_{r_i = a, r_j = b} p_{C,r}.
\]

For each pair of variables $x,y$ of $\fX$ which occur together in some constraint $C$, any solution to the basic semidefinite relaxation will automatically have the following properties:
\begin{itemize}
\item For all $a \ne b \in \bA_x$, we have $x_a \cdot x_b = 0$.
\item We have $\sum_{a \in \bA_x} \|x_a\|^2 = \|\sum_{a \in \bA_x} x_a\|^2 = 1$.
\item For all $a \in \bA_x, b \in \bA_y$, we have $x_a \cdot y_b \ge 0$.
\item We have $\sum_{a \in \bA_x, b \in \bA_y} x_a \cdot y_b = 1$.
\end{itemize}

Our goal is to find such a system of probabilities $p_{C,r}$ and vectors $x_a$ such that the quantity
\[
\tfrac{1}{\#C}\sum_C \sum_{r \in \RR_C} p_{C,r}
\]
is maximized. The maximum possible value of that sum is called the \emph{value} of the semidefinite relaxation. If the value of the semidefinite relaxation is equal to $1$, then we say that the system of probabilities $p_{C,r}$ and vectors $x_a$ \emph{perfectly solves} the semidefinite relaxation.
\end{defn}

Note that the constraints we make on the vectors $x_a$ only involve their dot products, and that they are always linear equalities/inequalities in terms of these dot products. We can deduce from these constraints a result which involves the vectors directly.

\begin{prop} Suppose that a system of probabilities $p_{C,r}$ and vectors $x_a$ are as in the basic semidefinite relaxation of a CSP instance $\fX$. Then for any variables $x,y$ of $\fX$ which occur together in some constraint, we have
\[
\sum_{a \in \bA_x} x_a = \sum_{b \in \bA_y} y_b.
\]
\end{prop}
\begin{proof} Let $x_{\bA_x} = \sum_{a \in \bA_x} x_a$ and similarly let $y_{\bA_y} = \sum_{b \in \bA_y} y_b$. Then we have $\|x_{\bA_x}\|^2 = \|y_{\bA_y}\|^2 = x_{\bA_x}\cdot y_{\bA_y} = 1$, so by the equality case of Cauchy-Schwarz we must have $x_{\bA_x} = y_{\bA_y}$.
\end{proof}

A useful generalization of the notation for the vectors $x_a$ was used heavily in \cite{sdp}.

\begin{defn} Suppose we are in the setup of the basic semidefinite relaxation of a CSP instance $\fX$. If $x$ is a variable of $\fX$ and $A \subseteq \bA_x$, then we define the vector $x_A$ by
\[
x_A = \sum_{a \in A} x_a.
\]
\end{defn}

Note that since the vectors $x_a$ are pairwise orthogonal for a fixed variable $x$, we have
\[
\|x_A\|^2 = \sum_{a \in A} \|x_a\|^2 = x_A \cdot x_{\bA_x}.
\]
Additionally, for each pair of variables $x,y$, we have
\[
x_A \cdot y_B = \sum_{a \in A, b \in B} x_a\cdot y_b
\]
by the distributive law.

Before explaining how to use the basic semidefinite relaxation to robustly solve affine-free CSPs, we will first show that if an instance $\fX$ of an affine-free CSP has a perfect solution to its basic semidefinite relaxation, then in fact $\fX$ has a solution. The main idea is to prove that the set of values $a \in \bA_x$ such that the vectors $x_a$ are nonzero can be used to restrict the instance to get a weak Prague instance, which will then be $pq$-consistent by Theorem \ref{weak-prague-pq}. The crucial computation is analyzing what happens to the vector $x_A$ when we take a single step along a path.

\begin{lem}\label{sdp-weak-prague} Suppose we are in the setup of the basic semidefinite relaxation of a CSP instance $\fX$. Let $x,y$ be variables of $\fX$ which occur together in some constraint, and define a binary relation $P \subseteq \bA_x \times \bA_y$ by
\[
(a,b) \in P \;\; \iff \;\; x_a \cdot y_b > 0.
\]
Then for any set $A \subseteq \bA_x$, we have
\[
\|x_A\|^2 \le \|y_{A+P}\|^2,
\]
with equality only when $x_A = y_{A+P}$. In fact, we have
\[
x_A \cdot (y_{A+P} - x_A) = 0,
\]
that is, $y_{A+P}$ is the sum of $x_A$ with a vector which is perpendicular to $x_A$. Furthermore, we have $x_A = y_{A+P}$ if and only if $A + P - P \subseteq A$.
\end{lem}
\begin{proof} Before diving into the algebraic details of the proof, it may be helpful to note that the dot products $x_a \cdot y_b$ define a probability distribution $\mu$ supported on $P \subseteq \bA_x \times \bA_y$ such that for $A \subseteq \bA_x, B \subseteq \bA_y$ we have $\mathbb{P}_\mu[A\times B] = x_A \cdot y_B$, and such that the marginal distributions $\mu_x, \mu_y$ on $\bA_x, \bA_y$ have probabilities given by $\mathbb{P}_{\mu_x}[A] = \|x_A\|^2$, $\mathbb{P}_{\mu_y}[B] = \|y_B\|^2$. It's possible to argue purely in terms of the probability distribution $\mu$, but the proof we give below won't directly refer to $\mu$ at all.

We have $x_A\cdot x_A = x_A \cdot x_{\bA_x}$, and by the definition of $P$ we have
\[
x_A \cdot y_{A+P} = \sum_{a\in A, b \in A+P} x_a \cdot y_b = \sum_{a\in A, b \in \bA_y} x_a \cdot y_b = x_A \cdot y_{\bA_y}.
\]
Since $x_{\bA_x} = y_{\bA_y}$, we have
\[
x_A \cdot x_A = x_A \cdot x_{\bA_x} = x_A \cdot y_{\bA_y} = x_A \cdot y_{A+P},
\]
so $x_A$ is orthogonal to $y_{A+P} - x_A$.

From the orthogonality of $x_A$ with $y_{A+P} - x_A$, we have
\[
\|y_{A+P}\|^2 = \|x_A + (y_{A+P} - x_A)\|^2 = \|x_A\|^2 + \|y_{A+P} - x_A\|^2 \ge \|x_A\|^2,
\]
with equality exactly when $y_{A+P} = x_A$.

For the last statement, note that we have $\|x_A\|^2 \le \|y_{A+P}\|^2 \le \|x_{A+P-P}\|^2$, so we just need to check the implication $x_A = y_{A+P} \implies x_{A+P-P} = x_A$. Under the assumption $x_A = y_{A+P}$, we have
\[
x_A \cdot y_{A+P} = \|y_{A+P}\|^2 = y_{\bA_y} \cdot y_{A+P} = x_{\bA_x} \cdot y_{A+P}.
\]
Suppose for contradiction that there was some $a \in (A+P-P) \setminus A$. Then by the definition of $A+P-P$ there would be some $b \in A+P$ such that $(a,b) \in P$, that is, such that $x_a \cdot y_b > 0$. But then we would have
\[
x_{\bA_x} \cdot y_{A+P} \ge x_a \cdot y_b + x_A \cdot y_{A+P} > x_A \cdot y_{A+P},
\]
which contradicts the assumption $x_A = y_{A+P}$.
\end{proof}

\begin{thm} Suppose $\fX$ is an instance of an affine-free CSP, and that there is a system of probabilities $p_{C,r}$ and vectors $x_a$ which perfectly solves the basic semidefinite relaxation of $\fX$. Then $\fX$ has a solution.
\end{thm}
\begin{proof} We define a restriction $\fX'$ of $\fX$ by restricting each relation $\RR_C$ to the support $R_C'$ of the marginal distribution $p_{C,r}$, and restricting each variable domain $\bA_x$ to the set $A_x'$ of $a \in \bA_x$ such that $\|x_a\|^2 \ne 0$. Note that each $R_C'$ will be contained in the original relation $\RR_C$ if we have a perfect solution to the basic semidefinite relaxation. Additionally, for each $C$ and each pair of variables $x = v_{C,i}, y = v_{C,j}$, the binary projection $\pi_{i,j}(R_C')$ will be equal to the set of ordered pairs $(a,b) \in A_x' \times A_y'$ such that $x_a \cdot y_b \ne 0$, by the compatibility between the probabilities $p_{C,r}$ and the vectors $x_a$.

We will check that $\fX'$ is a weak Prague instance (see Definition \ref{defn-weak-prague}). Arc-consistency (aka condition (P1)) of $\fX'$ follows from the compatibility between the probabilities and the vectors. To check (P2) and (P3), we use Lemma \ref{sdp-weak-prague}. Let $A \subseteq A_x'$ and let $p$ be a cycle from $x$ to $x$ in the instance $\fX'$, with $p_1$ from $x$ to $y$ the first step of the cycle $p$. If we have
\[
A+p = A,
\]
then by Lemma \ref{sdp-weak-prague} we have
\[
\|x_A\|^2 \le \|y_{A+p_1}\|^2 \le \|x_{A+p}\|^2 = \|x_A\|^2,
\]
so by the equality case of Lemma \ref{sdp-weak-prague} we must have
\[
x_A = y_{A+p_1}.
\]
Thus for any $a' \not\in A$, we must have $x_{a'} \cdot y_{A+p_1} = x_{a'} \cdot x_A = 0$, so we have
\[
A+p_1-p_1 = A.
\]
Thus by Proposition \ref{prop-p2*} we see that $\fX'$ satisfies condition (P2). We could have also checked condition (P2) using only the system of probabilities $p_{C,r}$, without mentioning the vectors $x_a$, by Theorem \ref{lp-p1-p2}.

To check (P3), let $A \subseteq A_x'$, and let $p,q$ be cycles from $x$ to $x$ in the instance $\fX'$, with
\[
A+p+q = A.
\]
Then by Lemma \ref{sdp-weak-prague} we have
\[
\|x_A\|^2 \le \|x_{A+p}\|^2 \le \|x_{A+p+q}\|^2 = \|x_A\|^2,
\]
so by the equality case of Lemma \ref{sdp-weak-prague} we must have
\[
x_A = x_{A+p}.
\]
In particular, we must have $A = A+p$, which proves (P3).

To finish, we note that $\fX'$ is $pq$-consistent by Theorem \ref{weak-prague-pq}, so $\Sg(\fX')$ is also $pq$-consistent by Proposition \ref{pq-algebraic-closure}, and $\Sg(\fX')$ is a restriction of the instance $\fX$ since $\fX = \Sg(\fX)$. Thus $\Sg(\fX')$ has a solution by Theorem \ref{affine-free-pq} and its corollaries, which is also a solution to the original instance $\fX$.
\end{proof}

In order to extend the previous result to an algorithm for \emph{robustly} solving affine-free CSPs, we need to find some approximate analogue of Lemma \ref{sdp-weak-prague}. The plan is to start by arguing as in Theorem \ref{lp-robust}, using the probabilities $p_{C,r}$ to produce an arc-consistent instance with variable domains $A_x'$ and constraint relations $R'$, such that each tuple $r \in R'$ has probability above some (random) threshold $\theta$. Then we will randomly cut the unit ball of possible values for the vectors $x_A$ into finitely many pieces, so that a version of Lemma \ref{sdp-weak-prague} holds when we forget what the exact values of the vectors $x_A$ are, and instead only keep track of which piece of the ball $x_A$ is contained in.

The proof becomes slightly simpler if we use the concept of weak consistency, from the previous section, instead of condition (P3). The simplification we get by aiming for weak consistency instead of (P3) is that when we have
\[
A + p + q = A,
\]
we only have to ensure that
\[
x_A \cdot x_{A+p} > 0 \;\; (\implies A \cap (A+p) \ne \emptyset),
\]
rather than needing to ensure that $x_A = x_{A+p}$. This means we only need to make sure that we chop the ball into fine enough pieces to separate each pair $x_A, x_B$ with $x_A \cdot x_B = 0$, instead of needing to separate each pair $x_A, x_B$ with $x_A \ne x_B$.

In order to produce our weakly consistent instance, we will define a quasiorder $\preceq$ (with an associated strict partial order $\prec$ and equivalence relation $\sim$) on the ball. We will first choose a (randomized) sequence of radii $r_1, r_2, ...$, and a (random) collection of hyperplanes $\mathcal{H}_i$ such that for every variable $x$, any pair of vectors $x_A, x_B$ which are orthogonal are (with high probability) separated from each other by at least one of the hyperplanes $\mathcal{H}_i$. The plan is to define $\preceq$ by
\[
u \prec v \;\; \iff \;\; \exists i \text{ s.t. } \|u\| < r_i \le \|v\|
\]
and
\[
u \sim v \;\; \iff \;\; (\forall i \; \|u\| < r_i \iff \|v\| < r_i) \wedge (\forall j \; u, v \text{ are on the same side of }\mathcal{H}_j).
\]
The plan is to throw away any constraint relation $R'$ which is incompatible with the quasiorder $\preceq$, where we say that $R'$ is incompatible with $\preceq$ if there are variables $x,y$ and some $A \subseteq A_x'$ such that
\[
x_A \not\preceq y_{A + \pi_{xy}(R')}.
\]
We will choose the radii $r_i$ to guarantee that
\[
A + \pi_{xy}(R') - \pi_{xy}(R') \ne A \;\; \implies \;\; x_A \prec y_{A + \pi_{xy}(R')},
\]
by taking $r_{i+1}^2 < r_i^2 + \theta$.
On the other hand, if
\[
A + \pi_{xy}(R') - \pi_{xy}(R') = A,
\]
then $x_A$ will be very close to $y_{A + \pi_{xy}(R')}$. In this case, we will ensure that the $r_i$ are spaced widely enough to make it unlikely that
\[
y_{A + \pi_{xy}(R')} \prec x_A.
\]
We also need to rule out the possibility that $x_A$ and $y_{A + \pi_{xy}(R')}$ are separated by some hyperplane $\mathcal{H}_i$. The chance that a particular random hyperplane separates $x_A$ and $y_{A + \pi_{xy}(R')}$ is proportional to the angle between $x_A$ and $y_{A + \pi_{xy}(R')}$, which will be low as long as $x_A$ is sufficiently close to $y_{A + \pi_{xy}(R')}$. All we have left to do is to carefully work out the details.

\begin{thm}[Slight refinement of \cite{sdp}]\label{refined-sdp} If $\Gamma$ is a finite constraint language and $\fA = (A, \Gamma)$ has bounded relational width, then the basic SDP relaxation can be used to \emph{robustly} solve $\CSP(\fA)$.

More precisely, if we are given an instance $\fX$ such that the basic SDP relaxation is $1 - \epsilon$ satisfiable, then we can efficiently find a solution which satisfies a
\[
1 - O\Big(\frac{\log(\log(\log(1/\epsilon)))}{\log(1/\epsilon)}\Big)
\]
fraction of the constraints.
\end{thm}
\begin{proof} Suppose we have a system of probabilities $p_{C,r}$ and vectors $x_a$ which solves the basic SDP relaxation of our instance $\fX$ with value at least $1-\epsilon$, that is, such that
\[
\tfrac{1}{\#C}\sum_C \sum_{r \in \RR_C} p_{C,r} \ge 1 - \epsilon.
\]
As an initial simplification to the problem, we will preemptively give up on any constraint $C$ such that
\[
\sum_{r \in \RR_C} p_{C,r} < 1 - \sqrt{\epsilon}.
\]
This gives up on at most a $\sqrt{\epsilon} = o(1/\log(1/\epsilon))$ fraction of the constraints of $\fX$, and allows us to focus on solving the problem in the special case where every individual constraint $C$ satisfies
\begin{equation}
\sum_{r \in \RR_C} p_{C,r} \ge 1 - \sqrt{\epsilon}.\label{sdp-cons-ineq}
\end{equation}
The advantage of this step is that from here on, we can look for a randomized algorithm which has a high probability of satisfying each constraint relation in isolation.

Let $N = \log(1/\sqrt{\epsilon})$. We will make a series of randomized choices in order to produce a weakly consistent instance $\fX'$, and after each choice we will give up on some of the constraints of our original instance $\fX$. At each step, we just need to confirm that for each constraint $C$ which satisfies \eqref{sdp-cons-ineq}, the chance of giving up on the constraint $C$ is at most
\[
O\Big(\frac{\log(\log(N))}{N}\Big).
\]

First we will try to produce an arc-consistent instance. Choose a threshold $\theta = \exp(-t) \in [\sqrt{\epsilon}, 1]$ by choosing $t = \log(1/\theta)$ uniformly at random from $[0,N]$. For each constraint $C$ with corresponding constraint relation $\RR_C$, we define the reduced relation $R_C'$ by
\[
R_C' = \{r \in \RR_C \mid p_{C,r} \ge 2\theta\}.
\]
For each variable $x$, we define the reduced variable domain $A_x'$ by
\[
A_x' = \{a \in \bA_x \mid \|x_a\|^2 \ge \theta\}.
\]
Our reduced instance $\fX'$ will have variable domains $A_x'$ and constraint relations $R_C'$ for all of the constraints $C$ which we do not choose to give up on by the end of our randomized procedure.
In order to ensure arc-consistency of $\fX'$, we preemptively give up on any constraint $C$ which does not satisfy
\begin{equation}
\sum_{r \not\in R_C'} p_{C,r} < \theta.\label{sdp-cons-arc}
\end{equation}

{\bf Claim 1.} If a constraint $C$ with $R_C' \subseteq A_{x_1}' \times \cdots \times A_{x_m}'$ satisfies \eqref{sdp-cons-arc}, then we have
\[
\pi_{x_i}(R_C') = A_{x_i}'
\]
for each $i \in \{1, ..., m\}$.

{\bf Proof of Claim 1.} Let $x = x_i$. First, note that for any $r \in R_C'$, if we set $a = \pi_x(r)$, then we have $\|x_a\|^2 \ge p_{C,r} \ge \theta$, so $a \in A_x'$. Thus we have $\pi_x(R_C') \subseteq A_x'$.

Conversely, if $a \not\in \pi_x(R_C')$ then we have
\[
\|x_a\|^2 = \sum_{\pi_x(r) = a} p_{C,r} \le \sum_{r \not\in R_C'} p_{C,r} < \theta
\]
by \eqref{sdp-cons-arc}, so $a \not\in \pi_x(R_C') \implies a \not\in A_x'$. Thus $\pi_x(R_C') \supseteq A_x'$, so $\pi_x(R_C') = A_x'$.

{\bf Claim 2.} If $\theta = \exp(-t)$ and $t$ is chosen uniformly at random from $[0, N]$, then
\[
\EE\Big[\frac{1}{\theta} \sum_{r \not\in R_C'} p_{C,r}\Big] = O(1/N),
\]
and the implied constant only depends on the number of tuples in the constraint $\RR_C$.
As a consequence, the probability that \eqref{sdp-cons-arc} fails to hold for any given constraint $C$ is $O(1/N)$.

{\bf Proof of Claim 2.} We have
\begin{align*}
\EE\Big[\frac{1}{\theta} \sum_{r \not\in R_C'} p_{C,r}\Big] &= \frac{1}{N} \int_{t = 0}^N \frac{1}{\theta} \sum_{r \not\in R_C'} p_{C,r}\ dt\\
&= \frac{1}{N} \sum_{r \not\in \RR_C} p_{C,r} \int_{t=0}^N \frac{1}{\theta}\ dt + \frac{1}{N} \sum_{r \in \RR_C} p_{C,r} \int_{t=0}^{\max(N,\log(2/p_{C,r}))} \frac{1}{\theta}\ dt\\
&\le \frac{1}{N} \int_{t=0}^N \frac{\sqrt{\epsilon}}{\theta}\ dt + \frac{1}{N} \sum_{r \in \RR_C} \int_{t=0}^{\log(2/p_{C,r})} \frac{p_{C,r}}{\theta}\ dt\\
&< \frac{2|\RR_C| + 1}{N}.
\end{align*}

Next we slice the ball $\mathcal{B}$ containing the $x_A$s into shells in order to ensure that the consistency condition (P2) is satisfied (note that (P2) is not required in the definition of weak consistency - however, enforcing it makes the rest of the proof simpler). To do this, we choose a second threshold $\theta_2$ uniformly at random from the interval $[0,\theta]$, and define radii $r_i$ by
\[
r_i^2 = i\theta + \theta_2.
\]
We define the strict partial order $\prec$ on the ball $\mathcal{B}$ by
\begin{align*}
u \prec v \;\; &\iff \;\; \exists i \text{ s.t. } \|u\| < r_i \le \|v\|\\
&\iff \;\; \Big\lfloor \frac{\|u\|^2 - \theta_2}{\theta} \Big\rfloor < \Big\lfloor \frac{\|v\|^2 - \theta_2}{\theta} \Big\rfloor.
\end{align*}
In order to ensure that (P2) is satisfied, we will preemptively give up on any constraint $C$ such that there is a pair of variables $x,y$ involved in the constraint $C$ which do not satisfy
\begin{equation}
\forall A \subseteq A_x', \;\; y_{A + \pi_{xy}(R_C')} \not\prec x_A.\label{sdp-cons-shells}
\end{equation}
By the next claim, we will only need to check \eqref{sdp-cons-shells} in the special case where $A + \pi_{xy}(R_C') - \pi_{xy}(R_C') = A$.

{\bf Claim 3.} If $A \subseteq A_x'$ has $A + \pi_{xy}(R_C') - \pi_{xy}(R_C') \ne A$, and if the constraint $C$ satisfies \eqref{sdp-cons-arc}, then we automatically have
\[
x_A \prec y_{A + \pi_{xy}(R_C')}.
\]

{\bf Proof of Claim 3.} We just need to prove that
\[
\|y_{A + \pi_{xy}(R_C')}\|^2 \ge \|x_A\|^2 + \theta.
\]
If $A + \pi_{xy}(R_C') - \pi_{xy}(R_C') \ne A$, then there must be some $a \not\in A$ and some $b \in A + \pi_{xy}(R_C')$ with $(a,b) \in \pi_{xy}(R_C')$. Picking $r \in R_C'$ with $\pi_{xy}(r) = (a,b)$, we see that
\[
x_a \cdot y_b \ge p_{C,r} \ge 2\theta
\]
by the definition of $R_C'$. Since
\[
\|y_{A + \pi_{xy}(R_C')}\|^2 = x_{\bA_x} \cdot y_{A + \pi_{xy}(R_C')} \ge x_A \cdot y_{A + \pi_{xy}(R_C')} + x_a \cdot y_b,
\]
we just need to check that
\[
x_A \cdot y_{A + \pi_{xy}(R_C')} \ge \|x_A\|^2 - \theta.
\]
For this, we just note that
\begin{align*}
\|x_A\|^2 - x_A \cdot y_{A + \pi_{xy}(R_C')} &= x_A \cdot y_{\bA_y} - x_A \cdot y_{A + \pi_{xy}(R_C')}\\
&= \sum_{\substack{\pi_x(r) \in A\\ \pi_y(r) \not\in A + \pi_{xy}(R_C')}} p_{C,r}\\
&\le \sum_{r \not\in R_C'} p_{C,r} < \theta.
\end{align*}

In order to put a bound on the probability that \eqref{sdp-cons-shells} fails to hold when $A + \pi_{xy}(R_C') - \pi_{xy}(R_C') = A$, we need to show that the expected value of
\[
\frac{\Big|\|x_A\|^2 - \|y_{A + \pi_{xy}(R_C')}\|^2\Big|}{\theta}
\]
is small when $\theta = \exp(-t)$ and $t$ is chosen uniformly at random from $[0, N]$. In this case we have
\begin{align*}
\Big|\|x_A\|^2 - \|y_{A + \pi_{xy}(R_C')}\|^2\Big| &\le \max\Big(\|x_A\|^2 - x_A \cdot y_{A + \pi_{xy}(R_C')}, \; \|y_{A+\pi_{xy}(R_C')}\|^2 - x_A \cdot y_{A + \pi_{xy}(R_C')}\Big)\\
&\le \sum_{r \not\in R_C'} p_{C,r},
\end{align*}
so by Claim 2 this is $O(1/N)$ on average.

So far, everything we have done could have been phrased in terms of the linear relaxation rather than the semidefinite relaxation. The next step, where we enforce weak consistency, is the step where we finally use the full power of the semidefinite relaxation. Set
\[
M = \lceil\log_2(N)\rceil,
\]
and independently pick $M$ uniformly random hyperplanes $\cH_1, ..., \cH_M$. Now define the equivalence relation $\sim$ on $\cB$ by
\[
u \sim v \;\; \iff \;\; (\forall i \; \|u\| < r_i \iff \|v\| < r_i) \wedge (\forall j \; u, v \text{ are on the same side of }\mathcal{H}_j).
\]
In order to guarantee weak consistency, we need to preemptively give up on any variable $x$ (along with any constraint involving $x$) which does not satisfy
\begin{equation}
\forall A, B \subseteq A_x', \; A,B \ne \emptyset, \;\;\;  x_A \sim x_B \;\; \implies \;\; A \cap B \ne \emptyset \label{sdp-vars-weak}
\end{equation}
and we need to preemptively give up on every constraint $C$ such that there is a pair of variables $x,y$ involved in the constraint $C$ which do not satisfy
\begin{equation}
\forall A \subseteq A_x', \;\; A + \pi_{xy}(R_C') - \pi_{xy}(R_C') = A \;\; \implies \;\; y_{A + \pi_{xy}(R_C')} \sim x_A.\label{sdp-cons-weak}
\end{equation}

{\bf Claim 4.} The chance that any particular variable $x$ fails to satisfy \eqref{sdp-vars-weak} is at most $\frac{3^{|A_x'|}}{2^M} = O(1/N)$. As a consequence, the chance that we give up on any particular constraint $C$ due to \eqref{sdp-vars-weak} is also $O(1/N)$ by the union bound.

{\bf Proof of Claim 4.} There are less than $3^{|A_x'|}$ pairs of disjoint, nonempty subsets $A,B$ of $A_x'$. For any particular pair $A,B$, the chance that any particular random hyperplane $\cH_i$ separates $x_A$ from $x_B$ is exactly $1/2$, since $x_A \cdot x_B = 0$. Thus the chance that none of the $M$ independent random hyperplanes $\cH_1, ..., \cH_M$ separate $x_A$ from $x_B$ is $1/2^M \le 1/N$.

To finish the argument, we need to find an upper bound on the probability that \eqref{sdp-cons-weak} is violated. For a given $A \subseteq A_x'$ with $A + \pi_{xy}(R_C') - \pi_{xy}(R_C') = A$, this probability depends mainly on the angle $\alpha$ between $x_A$ and $y_{A + \pi_{xy}(R_C')}$. From
\[
\max\Big(\|x_A\|^2 - x_A \cdot y_{A + \pi_{xy}(R_C')}, \; \|y_{A+\pi_{xy}(R_C')}\|^2 - x_A \cdot y_{A + \pi_{xy}(R_C')}\Big) \le \sum_{r \not\in R_C'} p_{C,r}
\]
we get
\[
x_A \cdot y_{A + \pi_{xy}(R_C')} \ge \|x_A\| \|y_{A + \pi_{xy}(R_C')}\| - \sum_{r \not\in R_C'} p_{C,r}.
\]
Using the fact that $\|x_A\|^2, \|y_{A + \pi_{xy}(R_C')}\|^2 \ge \theta$, we get
\[
\cos(\alpha) \ge 1 - \frac{1}{\theta}\sum_{r \not\in R_C'} p_{C,r},
\]
so
\[
\alpha^2 = O\Big(\frac{1}{\theta}\sum_{r \not\in R_C'} p_{C,r}\Big).
\]
The chance that at least one of the $M$ hyperplanes $\cH_i$ separates $x_A$ from $y_{A + \pi_{xy}(R_C')}$ is then given by
\begin{align*}
1 - \Big(1 - \frac{\alpha}{\pi}\Big)^M &\le \max\Big(1, \frac{M\alpha}{\pi}\Big)\\
&\ll \max\Big(1, M\Big(\frac{1}{\theta}\sum_{r \not\in R_C'} p_{C,r}\Big)^{1/2}\Big).
\end{align*}
To finish the proof, we just need to verify one final claim.

{\bf Claim 5.} If $\theta = \exp(-t)$ and $t$ is chosen uniformly at random from $[0, N]$, then
\[
\EE\Big[\max\Big(1, M\Big(\frac{1}{\theta}\sum_{r \not\in R_C'} p_{C,r}\Big)^{1/2}\Big)\Big] = O\Big(\frac{\log(M)}{N}\Big).
\]

{\bf Proof of Claim 5.} The argument is similar to the proof of Claim 2, just slightly more involved. First, we use the inequality $\sqrt{a+b} \le \sqrt{a} + \sqrt{b}$ to get the bound
\[
\max\Big(1, M\Big(\frac{1}{\theta}\sum_{r \not\in R_C'} p_{C,r}\Big)^{1/2}\Big) \le \max\Big(1, M\sqrt{\frac{\sqrt{\epsilon}}{\theta}}\Big) + \sum_{\substack{r \in \RR_C\\ p_{C,r} < \theta}} \max\Big(1, M\sqrt{\frac{p_{C,r}}{\theta}}\Big).
\]
We then bound the contribution of each summand individually. Each $\max$ takes the value $1$ for a range of values of $t = \log(1/\theta)$ of length $2\log(M)$, and from then on takes exponentially decaying values, so the total expectation ends up being bounded by
\[
\frac{2(\log(M) + 1)(|\RR_C| + 1)}{N} = O\Big(\frac{\log(M)}{N}\Big).
\]

Thus, after preemptively giving up on at most a
\[
O\Big(\frac{\log(M)}{N}\Big) = O\Big(\frac{\log(\log(N))}{N}\Big) = O\Big(\frac{\log(\log(\log(1/\epsilon)))}{\log(1/\epsilon)}\Big)
\]
fraction of the constraints, we finally manage to construct a weakly consistent instance $\fX'$. Then $\Sg(\fX') \subseteq \fX$ will also be weakly consistent, so it will have a solution by Theorem \ref{thm-weakly-consistent}.
\end{proof}

The bound given in Theorem \ref{refined-sdp} is slightly better than the bound from \cite{sdp}, which only guaranteed that we could find a solution of quality
\[
1 - O\Big(\frac{\log(\log(1/\epsilon))}{\log(1/\epsilon)}\Big).
\]
We can improve this further, getting rid of the unsightly $\log(\log(\log(1/\epsilon)))$ factor entirely, if we try to construct a vague solution (to the binary part) instead of aiming for weak consistency.

\begin{thm}[Perfected form of \cite{sdp}]\label{perfected-sdp} Suppose $\Gamma$ is a finite constraint language such that $\fA = (A, \Gamma)$ has bounded relational width. If we are given an instance $\fX$ such that the basic SDP relaxation is $1 - \epsilon$ satisfiable, then we can efficiently find a solution which satisfies a
\[
1 - O\Big(\frac{1}{\log(1/\epsilon)}\Big)
\]
fraction of the constraints. In other words, the function $f$ from Definition \ref{defn-robust} satisfies $f(\epsilon) = O(1/\log(1/\epsilon))$.
\end{thm}
\begin{proof} We argue as in the proof of Theorem \ref{refined-sdp}, using the same notation, up to the point where we obtained an arc-consistent instance $\fX'$ which satisfied (P2). From here on, we modify the argument: instead of trying to define an equivalence relation $\sim$ on the shell of the ball $\cB$ between radius $r_i$ and $r_{i+1}$, we extend $\prec$ to a (nearly) total order. We do this by picking a uniformly random unit vector $U$, and setting
\[
\forall \|u\|, \|v\| \in (r_i, r_{i+1}], \;\; u \prec v \iff U\cdot(u-v) < 0.
\]
Now we use the (almost) total order $\prec$ to define a vague element $\prec_x$ for each variable $x$, by setting
\[
A \prec_x B \;\; \iff \;\; x_A \prec x_B.
\]
Note that with probability $1$, each $\prec_x$ will be a total order on $\cP_\emptyset(A_x')$, so $\prec_x$ is indeed a vague element of $A_x'$.

We will now preemptively give up on any constraint relation $C$ if there is any pair of variables $x,y$ such that the vague elements $\prec_x, \prec_y$ fail to vaguely satisfy the binary relation $\pi_{xy}(R_C')$. We need to find an upper bound on the probability that we will give up on any particular constraint $C$ due to the pair of variables $x,y$. Set $R' = \pi_{xy}(R_C')$ to simplify the notation.

We begin trying to construct a quasiorder $\preceq_{R'}$ extending $\prec_x, \prec_y$ on $\cP_\emptyset(A_x') \sqcup \cP_\emptyset(A_y')$ by setting
\[
x_A \prec_{R'} y_B \text{ when } \exists i \; \|x_A\| < r_i \le \|y_B\|,
\]
and similarly with the roles of $x$ and $y$ reversed. This ensures that
\[
A \ne A + R' - R' \;\; \implies \;\; x_A \prec_{R'} y_{A + R'},
\]
and similarly with the roles of $x$ and $y$ reversed. So far we encounter no problems with the construction of $\preceq_{R'}$.

We will also need to require that
\[
A = A + R' - R' \;\; \implies \;\; x_A \sim_{R'} y_{A + R'},
\]
and this is the step which could potentially cause an issue. Note that if we haven't already given up on the constraint $C$ while ensuring that (P2) holds, then in this case we have the guarantee that $x_A$ and $y_{A + R'}$ are contained in the same shell of the ball $\cB$.

Since the different shells of $\cB$ don't interact in any of the orderings $\prec_x, \prec_y, \preceq_{R'}$, we can focus on just one particular shell of $\cB$ corresponding to a pair of adjacent radii $r_i, r_{i+1}$. Within such a shell, all $\preceq_{R'}$ does is identify certain vectors $x_A$ with corresponding vectors $y_{A+R'}$. The only way we could run into trouble is if there was some pair $A,B \subseteq A_x'$ with
\begin{align*}
A &\sim_{R'} A + R',\\
B &\sim_{R'} B + R',\\
A &\prec_x B,
\end{align*}
but
\[
B + R' \prec_y A + R'.
\]
In other words, we only need to preemptively give up on the constraint $C$ due to the pair $x,y$ if we can find such a pair of sets $A,B \subseteq A_x'$ with
\[
U\cdot (x_A - x_B) < 0 < U \cdot (y_{A+R'} - y_{B+R'}).
\]
For a randomly chosen $U$, the probability of this occuring is proportional to the angle between the vector $x_A - x_B$ and the vector $y_{A+R'} - y_{B+R'}$.

Since $A, B$ are different subsets of $A_x'$, at least one coordinate of $x_A - x_B$ must have absolute value at least $\sqrt{\theta}$, so we have $\|x_A - x_B\|^2 \ge \theta$, and similarly $\|y_{A+R'} - y_{B+R'}\|^2 \ge \theta$. As in the proof of Theorem \ref{refined-sdp}, we have
\[
\|x_A - y_{A+R'}\|^2, \|x_B - y_{B+R'}\|^2 \le \sum_{r \not\in R_C'} p_{C,r},
\]
so the angle between $x_A - x_B$ and $y_{A+R'} - y_{B+R'}$ is
\[
O\Big(\frac{1}{\sqrt{\theta}}\sqrt{\sum_{r \not\in R_C'} p_{C,r}}\Big).
\]
By Claim 5 of Theorem \ref{refined-sdp}, the average value of this upper bound is at most
\[
O(1/N) = O(1/\log(1/\epsilon)),
\]
so we only need to give up on $O(1/N)$ of our constraints in order to get an arc-consistent instance $\fX'$ whose binary part has a vague solution. Then by Corollary \ref{cor-vague-bin} the instance $\Sg(\fX') \subseteq \fX$ has a solution, which finishes the proof.
\end{proof}

\begin{rem} The algorithm from Theorem \ref{perfected-sdp} can be derandomized without too much effort using the method of conditional expectations. After throwing away constraints which are violated by more than a $\sqrt{\epsilon}$ fraction, we pick $\theta$ in $[\sqrt{\epsilon}, 1]$ such that
\[
\frac{1}{\#C} \sum_C \frac{1}{\theta}\sum_{r \not\in R_C'} p_{C,r} + \frac{1}{\#C} \sum_C \frac{1}{\sqrt{\theta}}\sqrt{\sum_{r \not\in R_C'} p_{C,r}}
\]
is minimized - the minimum is guaranteed to be $O(1/N)$, and we only have to examine values of $\theta$ which are equal to some $p_{C,r}/2$. Then we use $\theta$ to do the first bit of rounding to get to an arc-consistent instance. Next we pick $\theta_2 \in [0,\theta]$ to minimize the number of problems we run into while ensuring that (P2) is satisfied - we only need to try values of $\theta_2$ which are just above or just below a remainder of some $\|x_A\|^2$ modulo $\theta$, so this can be done efficiently.

For the final step, we want to pick a unit vector $U$ which is dual to a hyperplane which separates as few pairs of vectors $x_A - x_B$, $y_{A+R'} - y_{B+R'}$ with $A + R' - R' = A, B + R' - R' = B$ as possible. This step is trickier, but in Appendix C of \cite{sdp} they cite the paper \cite{karnin2012explicit}, which shows that you can find a $U$ which is worse by at most $O(1/N)$ compared to what you would get with a uniformly random choice, using a deterministic algorithm which runs in time $\|\fX\|^{O(1)}$ times $2^{\log(N)^2} = o(1/\epsilon)$.
\end{rem}

We may naturally wonder whether $f(\epsilon) = O(\frac{1}{\log(1/\epsilon)})$ is really the best possible asymptotic one can get in Theorem \ref{perfected-sdp}. The next example shows that at least for HORN-SAT, it is impossible to improve the asymptotic.

\begin{ex}[Guruswami and Zhou \cite{robust-horn-gap}]\label{ex-simul-ind} For every $n$, consider the following ``simultaneous induction'' instance of HORN-SAT, on the $2n+2$ variables $p_0, ..., p_n, q_0, ..., q_n$:
\begin{align*}
&(p_0 = 1) \wedge (q_0 = 1)\\
&\wedge (p_0 \wedge q_0 \implies p_1) \wedge (p_0 \wedge q_0 \implies q_1)\\
&\wedge\ \cdots\\
&\wedge (p_{n-1} \wedge q_{n-1} \implies p_n) \wedge (p_{n-1} \wedge q_{n-1} \implies q_n)\\
&\wedge (p_n = 0) \wedge (q_n = 0).
\end{align*}
This instance has $2n+4$ constraints, and it is possible to satisfy at most $2n+3$ of them. However, the basic SDP relaxation of this instance thinks it can satisfy a $1 - (3/4)^n$ fraction of the constraints!

Checking that the SDP relaxation has such a solution (for all $n$) is fairly tricky. One simplification which helps quite a bit is to note that the variables $p_i, q_i$ only interact with $p_{i\pm 1}, q_{i \pm 1}$, so we just need to analyze the set of solutions to the basic SDP relaxation for the four-variable HORN-SAT instance
\[
(x \wedge y \implies z) \wedge (x \wedge y \implies w)
\]
in order to understand the possible behavior in the general case where we string together many of these ``simultaneous induction'' steps.

A second simplification is that since the overall vector 
\[
x_{\{0,1\}} = y_{\{0,1\}} = z_{\{0,1\}} =  w_{\{0,1\}}
\]
is a fixed unit vector, and since $x_0$ and $x_1$ are perpendicular with
\[
x_0 + x_1 = x_{\{0,1\}},
\]
etc., we only have to describe the dot products between the four vectors $x_0, y_0, z_0, w_0$ to determine the whole configuration. (In particular, if our SDP relaxation of this instance has a solution, then it has a solution where all of the vectors live on a $4$-dimensional sphere of diameter $1$ in $\RR^5$, which passes through the origin and is centered at $\frac{1}{2}x_{\{0,1\}}$.)

The plan is to make the probabilities that $z,w$ are equal to $0$ grow to be a constant factor larger than the probabilities that $x,y$ are equal to $0$, while keeping the correlation between $z$ and $w$ under control. To this end, we claim that the SDP relaxation of this four-variable instance has a solution which satisfies
\begin{align*}
\|x_0\|^2 = \|y_0\|^2 &= 9\epsilon,\\
x_0 \cdot y_0 &= 3\epsilon,\\
\|z_0\|^2 = \|w_0\|^2 &= 12\epsilon,\\
z_0 \cdot w_0 &= 4\epsilon,\\
x_0 \cdot z_0 = ... = y_0 \cdot w_0 &= 6\epsilon
\end{align*}
for any $0 \le \epsilon \le 1/18$. Note that this way we have
\[
\frac{\|z_0\|^2}{\|x_0\|^2} = \frac{\|w_0\|^2}{\|y_0\|^2} = \frac{z_0 \cdot w_0}{x_0 \cdot y_0} = \frac{4}{3}
\]
and
\[
\frac{z_0\cdot w_0}{\|z_0\|\|w_0\|} = \frac{x_0\cdot y_0}{\|x_0\|\|y_0\|} = \frac{1}{3},
\]
so we can glue the solutions to many units like this together, as long as we can show that each unit is a valid solution to the SDP relaxation on its own.
For this, we check that the matrix of dot products between our hypothetical vectors $x_{\{0,1\}}, x_0, y_0, z_0, w_0$ below is positive semidefinite:
\[
\begin{bmatrix} 1 & 9\epsilon & 9\epsilon & 12\epsilon & 12\epsilon\\
9\epsilon & 9\epsilon & 3\epsilon & 6\epsilon & 6\epsilon\\
9\epsilon & 3\epsilon & 9\epsilon & 6\epsilon & 6\epsilon\\
12\epsilon & 6\epsilon & 6\epsilon & 12\epsilon & 4\epsilon\\
12\epsilon & 6\epsilon & 6\epsilon & 4\epsilon & 12\epsilon
\end{bmatrix} \stackrel{?}{\succeq} 0,
\]
and we use the probability distribution
\begin{align*}
\PP[(x,y,z) = (1,1,1)] &= 1-15\epsilon,\\
\PP[(x,y,z) = (1,0,0)] &= 6\epsilon,\\
\PP[(x,y,z) = (0,1,0)] &= 6\epsilon,\\
\PP[(x,y,z) = (0,0,1)] &= 3\epsilon
\end{align*}
over the set of solutions $(x,y,z)$ to the constraint $x \wedge y \implies z$ (and similarly for the constraint $x \wedge y \implies w$).

In order to check that the $5$ by $5$ matrix above is always positive semidefinite, the simplest method is to replace the $1$ in the upper left corner by $18\epsilon$ (which is at most $1$ by assumption) and divide out the $\epsilon$s, to get a fixed matrix that doesn't depend on $\epsilon$. In order to check that the resulting matrix is positive semidefinite, we can make our lives easier by rewriting it as the sum
\[
\begin{bmatrix} 18 & 9 & 9 & 12 & 12\\
9 & 9 & 3 & 6 & 6\\
9 & 3 & 9 & 6 & 6\\
12 & 6 & 6 & 12 & 4\\
12 & 6 & 6 & 4 & 12
\end{bmatrix} = \begin{bmatrix} 18 & 9 & 9 & 12 & 12\\
9 & 6 & 6 & 6 & 6\\
9 & 6 & 6 & 6 & 6\\
12 & 6 & 6 & 8 & 8\\
12 & 6 & 6 & 8 & 8
\end{bmatrix} + \begin{bmatrix} 0 & 0 & 0 & 0 & 0\\
0 & 3 & -3 & 0 & 0\\
0 & -3 & 3 & 0 & 0\\
0 & 0 & 0 & 4 & -4\\
0 & 0 & 0 & -4 & 4
\end{bmatrix}.
\]
Since the second summand is clearly positive semidefinite, we just need to check that the first summand is positive semidefinite, which boils down to checking positive semidefiniteness of the submatrix determined by the first, second, and fourth rows and columns. We can check this by directly computing the Cholesky decomposition:
\[
\begin{bmatrix} 18 & 9 & 12\\
9 & 6 & 6\\
12 & 6 & 8
\end{bmatrix} = \begin{bmatrix} 1 & 0\\
1/2 & 1\\
2/3 & 0
\end{bmatrix} \begin{bmatrix} 18 & 0\\
0 & 3/2
\end{bmatrix} \begin{bmatrix} 1 & 1/2 & 2/3\\
0 & 1 & 0
\end{bmatrix}.
\]
Examples of vectors $x_{\{0,1\}}, x_0, y_0, z_0, w_0 \in \mathbb{R}^5$ with the desired dot products can be extracted from this computation:
\[
x_{\{0,1\}} = \begin{bmatrix} \sqrt{1-18\epsilon}\\ \sqrt{18\epsilon}\\ 0\\ 0\\ 0 \end{bmatrix},\;\; x_0 = \begin{bmatrix} 0\\ \sqrt{9\epsilon/2}\\ \sqrt{3\epsilon/2}\\ \sqrt{3\epsilon}\\ 0 \end{bmatrix},\;\; y_0 = \begin{bmatrix} 0\\ \sqrt{9\epsilon/2}\\ \sqrt{3\epsilon/2}\\ -\sqrt{3\epsilon}\\ 0 \end{bmatrix},\;\; z_0 = \begin{bmatrix} 0\\ \sqrt{8\epsilon}\\ 0\\ 0\\ \sqrt{4\epsilon} \end{bmatrix},\;\; w_0 = \begin{bmatrix} 0\\ \sqrt{8\epsilon}\\ 0\\ 0\\ -\sqrt{4\epsilon} \end{bmatrix}.
\]
When we join copies of this configuration with varying values of $\epsilon$ together, they will need to be rotated appropriately so that the $x_0, y_0$ from the next unit become the $z_0, w_0$ from the previous unit, and the $x_{\{0,1\}}$ from the next unit is the same as the $x_{\{0,1\}}$ from the previous unit.

By joining together many copies of the unit described above, we can find a solution to the SDP relaxation where the only constraints which aren't satisfied exactly are the constraints $p_0 = 1, q_0 = 1$ - instead, these will be satisfied with the exponentially small error $\frac{3}{5}\cdot(\frac{3}{4})^{n-1}$. We arrange things so that at the second-to-last step, we have
\[
\PP[p_{n-1} = 0] = \PP[q_{n-1} = 0] = 3/5
\]
and
\[
\PP[p_{n-1} = q_{n-1} = 0] = 1/5,
\]
at which point we can no longer continue to use the unit described above (since $\epsilon$ exceeds $1/18$ at this point). For the last step, we join this with an honest probability distribution over the set of solutions to the $4$-variable instance we get by restricting to $p_{n-1}, q_{n-1}, p_n, q_n$:
\begin{align*}
\PP[(p_{n-1}, q_{n-1}, p_n, q_n) = (0,0,0,0)] &= 1/5,\\
\PP[(p_{n-1}, q_{n-1}, p_n, q_n) = (1,0,0,0)] &= 2/5,\\
\PP[(p_{n-1}, q_{n-1}, p_n, q_n) = (0,1,0,0)] &= 2/5.
\end{align*}

The construction described here was optimized for readability, rather than for obtaining the best constants, and can be improved. A slight modification to our basic unit, which instead has $\frac{z_0\cdot w_0}{\|z_0\|\|w_0\|} = \frac{x_0\cdot y_0}{\|x_0\|\|y_0\|} = \frac{1}{4}$ and $\frac{\|z_0\|^2}{\|x_0\|^2} = \frac{\|w_0\|^2}{\|y_0\|^2} = \frac{z_0 \cdot w_0}{x_0 \cdot y_0} = \frac{3}{2}$, can be used in all but the last few steps of the construction, leading to a solution to the SDP relaxation which satisfies a $1 - (2/3)^n$ fraction of the constraints. With some more fiddly work, we can improve the base of the exponential from $2/3$ to $\phi^2/4 \approx 0.6545...$, where $\phi = \frac{\sqrt{5}+1}{2} \approx 1.618...$ is the golden ratio - I don't know whether it's possible to do any better than that.
\end{ex}

As a consequence of this example, we can't hope for improved asymptotics for any relational structure $\fA$ which can pp-construct HORN-SAT unless the Unique Games Conjecture is false (by the main result of Raghavendra's thesis \cite{raghavendra-thesis}). On the other hand, it is certainly possible to get better asymptotics for 2-SAT! Under the Unique Games Conjecture, the best possible function $f$ as in Definition \ref{defn-robust} for 2-SAT is given by $f(\epsilon) \sim \sqrt{\epsilon}$ (see \cite{dalmau-robust-framework} for references and discussion).

\begin{prob}[From \cite{dalmau-robust-framework}] Suppose that a finite relational structure $\fA$ has bounded width and does not pp-construct HORN-SAT. Is it necessarily the case that there is some $k$ such that $\CSP(\fA)$ can be robustly solved with the function $f$ from Definition \ref{defn-robust} satisfying $f(\epsilon) = O(\epsilon^{1/k})$?
\end{prob}

In \cite{robust-near-unanimity}, the authors show that for every $\bA$ with a near-unanimity term, there is some $k$ such that $\CSP(\fA)$ can be robustly solved (via the basic SDP relaxation) with the function $f$ from Definition \ref{defn-robust} satisfying $f(\epsilon) = O(\epsilon^{1/k})$. The proof involves a consistency condition which doesn't rely on arc-consistency.

\begin{rem} The first part of the proof of Theorem \ref{refined-sdp} shows that given a solution to the LP relaxation of an instance $\fX$ with value $1-\epsilon$, we can find a subinstance $\fX'$ which satisfies (P1) and (P2) after giving up on an $O(1/\log(1/\epsilon))$ fraction of the constraints. A positive answer to the following problem would then prove that every CSP which is solved by its LP relaxation is also \emph{robustly} solved by its LP relaxation.
\end{rem}

\begin{prob}[Conjectured in \cite{mui-symmetric}]\label{prob-p1-p2} Suppose that $\CSP(\bA)$ is solved by its LP relaxation. Is it necessarily the case that every instance $\fX$ of $\CSP(\bA)$ which satisfies the consistency conditions (P1) and (P2) has a solution?
\end{prob}

We will solve this problem in the next section.


\section{Linear Programming, rounding rules, and reversibility}

The goal of this section is to show that every CSP which is solved by the basic linear programming relaxation is also \emph{robustly} solved by linear programming. This will follow once we prove that instances of such CSPs which satisfy the consistency conditions (P1) and (P2) always have solutions. Since we will be referring to these conditions frequently, we will give them a name.

\begin{defn} An instance is called \emph{reversible} if it satisfies the consistency conditions (P1) and (P2) from Definition \ref{defn-weak-prague}.
\end{defn}

Recall that by Theorems \ref{lp-p1-p2} and \ref{lp-p1-p2-converse}, an instance is reversible if and only if every cycle in the instance has a linear programming solution of full support (on both the variables and the relations). However, we will need to have a characterization of reversibility which is more similar to to the characterization of condition (P3) from Proposition \ref{prop-prec-p3}. Since both the variable domains and the relation domains are equally relevant to the basic linear programming relaxation, this turns out to be most convenient if we adopt an alternative perspective on CSPs which puts relations and variables on the same footing.

Mainly, the idea is to change the way we visualize a constraint relation. In order to visualize a ternary relation $\RR \le \bA_x \times \bA_y \times \bA_z$, we will use the following diagram:
\begin{center}
\begin{tikzpicture}[scale=1.3]
  \node (R) at (0,0) {$\RR$};
  \node (ax) at (0,1) {$\bA_{x}$};
  \node (ay) at (-1,-0.7) {$\bA_{y}$};
  \node (az) at (1,-0.7) {$\bA_{z}$.};
  \draw (R) [->] to ["$\pi_x$"'] (ax);
  \draw (R) [->] to ["$\pi_y$"'] (ay);
  \draw (R) [->] to ["$\pi_z$"] (az);
\end{tikzpicture}
\end{center}
In order to check that an assignment $x \mapsto a \in \bA_x, y \mapsto b \in \bA_y, z \mapsto c \in \bA_z$ satisfies the corresponding constraint, we need to find an element $r \in \RR$ such that
\[
\pi_x(r) = a, \;\; \pi_y(r) = b, \;\; \pi_z(r) = c.
\]
Note that it isn't necessary for $\RR$ to actually be a subalgebra of $\bA_x \times \bA_y \times \bA_z$: all that matters is that we know the three maps $\pi_x, \pi_y, \pi_z$ from $\RR$ to $\bA_x, \bA_y, \bA_z$. (Of course, we also gain no extra generality by allowing $\RR$ to not be a subpower, since we could always replace $\RR$ by its image in $\bA_x \times \bA_y \times \bA_z$.)

By stitching these pictures together, we can visualize any instance $\fX$ of the multisorted CSP $\CSP(\bA_1, ..., \bA_n)$ as a diagram $\mathcal{D}_\fX$ in the concrete category $HSP_{fin}(\bA_1, ..., \bA_n)$. A solution of the instance $\fX$ corresponds to an assignment from each object $\bA$ of the diagram $\mathcal{D}_\fX$ to an element $a_\bA \in \bA$, such that for every homomorphism $\pi : \RR \rightarrow \bA$ which occurs in the diagram $\mathcal{D}_\fX$, the elements $a_\RR \in \RR$ and $a_\bA \in \bA$ are related by
\[
\pi(a_\RR) = a_\bA.
\]
To a category-theorist, a solution to the instance $\fX$ is just another name for an element of the ``inverse limit'' of the diagram $\cD_\fX$.

Conversely, to any diagram $\cD$ in a concrete category $\cV$, we can associate an instance $\fX_\cD$ of a multisorted CSP. The instance $\fX_\cD$ will have one variable $x_\bA$ for each object $\bA$ of the diagram $\cD$, and for each morphism $\pi : \RR \rightarrow \bA$ of the diagram $\cD$, there will be a binary constraint
\[
\pi(x_\RR) = x_\bA.
\]
In this way, we can transform any instance $\fX$ of any multisorted CSP into an instance $\fX'$ where every constraint relation of $\fX'$ is binary, and is furthermore the graph of a homomorphism, at the cost of introducing a new variable (with a large domain) of $\fX'$ for each of the constraint relations of the original instance $\fX$. This is what makes this perspective well-suited to studying the basic linear programming relaxation: in the basic linear programming relaxation, we have to introduce these new variables just to describe the relaxation, so we may as well work these new variables into our theory from the very beginning.

\begin{defn} We say that an instance $\fX$ of a multisorted CSP is \emph{diagrammatic} if every constraint relation of $\fX$ is binary and is the graph of a homomorphism.
\end{defn}

As it turns out, reversibility has a nice characterization in terms of orderings when we restrict to diagrammatic instances.

\begin{prop} If an instance $\fX$ is arc-consistent and diagrammatic, then $\fX$ is reversible iff there is a quasiorder $\preceq$ on the collection of pairs $(x,B)$ with $B \subseteq \bA_x$ with the following properties:
\begin{itemize}
\item[(a)] for any $B \subset C \subseteq \bA_x$, we have the strict inequality $(x,B) \prec (x,C)$, and
\item[(b)] for any constraint of $\fX$ corresponding to a homomorphism $\pi : \bA_x \twoheadrightarrow \bA_y$, and for any $B \subseteq \bA_y$, $(y,B)$ and $(x,\pi^{-1}(B))$ are in the same equivalence class of $\preceq$.
\end{itemize}
\end{prop}
\begin{proof} First suppose that $\fX$ is reversible. Define a quasiorder $\preceq$ on the collection of pairs $(x,B)$ with $B \subseteq \bA_x$ by setting $(x,B) \preceq (y,C)$ iff there is a path $p$ from $x$ to $y$ in $\fX$ such that
\[
B + p \subseteq C.
\]
This will automatically satisfy (b).

In order to ensure that (a) is satisfied, we just need to show that there is no path $p$ from $x$ to $x$ in $\fX$ such that
\[
C + p \subset C
\]
for some $C \subseteq \bA_x$. Suppose for the sake of contradiction that there is such a $p$ and $C$. Then we have
\[
C + (k+1)p = (C+p) + kp \subseteq C + kp
\]
for all $k \ge 0$, so there must be some first $k$ such that
\[
C + (k+1)p = C + kp.
\]
Since $C + p \ne C$, we have $k \ge 1$, so if we define $B = C + kp$ then we have $B + p = B$, but
\[
B - p = C + kp - p \supseteq C + (k-1)p \supset B,
\]
where the strict containment follows from the minimality of $k$. However, $B - p \ne B$ contradicts the assumption that $X$ is reversible, so such a $p$ and $C$ do not exist.

To finish, we need to check that if such a quasiorder $\preceq$ can be found, then $\fX$ is reversible. For this, we will use Proposition \ref{prop-p2*}. So we may suppose that $B \subseteq \bA_x$ and that $B + p = B$ for a path $p$ from $x$ to $x$ with first step $p_1$ from $x$ to $x_1$, and we just need to prove that $B + p_1 - p_1 = B$. If $p_1$ is the graph of a homomorphism $\pi_1$ from $\bA_{x_1}$ to $\bA_x$, then we have
\[
B + p_1 - p_1 = \pi_1^{-1}(B) - p_1 = \pi_1(\pi_1^{-1}(B)) = B
\]
by the assumption that $\fX$ is arc-consistent. The interesting case is the case where $p_1$ is the graph of a homomorphism $\pi_1$ from $\bA_x$ to $\bA_{x_1}$: in this case, we have
\[
B + p_1 - p_1 = \pi_1(B) - p_1 = \pi_1^{-1}(\pi_1(B)) \supseteq B.
\]
Suppose for the sake of contradiction that $\pi_1^{-1}(\pi_1(B))$ is a strict superset of $B$. Suppose that
\[
p = p_1 + p_2 + \cdots + p_k
\]
where each $p_i$ is a single step, and define $B_i$ by
\[
B_i = B + p_1 + \cdots + p_i,
\]
so that $B = B_0 = B_k$. If $x_i$ is the $i$th variable along the path $p$, then we claim that for each $i$ we have
\[
(x_i, B_i) \preceq (x_{i+1}, B_{i+1}).
\]
If $p_{i+1}$ is the graph of a homomorphism $\pi_{i+1}$ from $\bA_{x_{i+1}}$ to $\bA_{x_i}$ then $B_{i+1} = \pi_{i+1}^{-1}(B_i)$, so $(x_i, B_i)$ and $(x_{i+1}, B_{i+1})$ must be in the same equivalence class of $\preceq$ by (b). Otherwise, $p_{i+1}$ is the graph of a homomorphism $\pi_{i+1}$ from $\bA_{x_i}$ to $\bA_{x_{i+1}}$ and $B_{i+1} = \pi_{i+1}(B_i)$, so
\[
B_i \subseteq \pi_{i+1}^{-1}(B_{i+1}) \;\; \implies \;\; (x_i, B_i) \preceq (x_i, \pi_{i+1}^{-1}(B_{i+1})) \preceq (x_{i+1}, B_{i+1})
\]
by (a) and (b). Thus we have
\[
(x, \pi_1^{-1}(B_1)) \preceq (x_1, B_1) \preceq \cdots \preceq (x_k, B_k) = (x,B),
\]
and this contradicts (a) if $B$ is a strict subset of $\pi_1^{-1}(B_1) = \pi_1^{-1}(\pi_1(B))$.
\end{proof}

\begin{cor}\label{cor-diagram-reversible} If an instance $\fX$ is diagrammatic and reversible, then for each variable domain $\bA_x$ we can find a total quasiorder $\preceq_x$ on $\cP(\bA_x)$ with the following properties:
\begin{itemize}
\item[(a)] for any $B \subset C \subseteq \bA_x$, we have the strict inequality $B \prec_x C$, and
\item[(b)] for any constraint of $X$ corresponding to a homomorphism $\pi : \bA_x \twoheadrightarrow \bA_y$, and for any $B, C \subseteq \bA_y$ we have
\[
B \preceq_y C \;\;\; \iff \;\;\; \pi^{-1}(B) \preceq_x \pi^{-1}(C).
\]
\end{itemize}
\end{cor}
\begin{proof} Extend the quasiorder $\preceq$ from the previous proposition to a total quasiorder $\preceq'$ with the same equivalence classes as $\preceq$, and then restrict $\preceq'$ to $\cP(\bA_x)$ to define $\preceq_x$.
\end{proof}

\begin{rem} With a little bit more work, we can also assume that the total quasiorders $\preceq_x$ satisfy
\[
B \preceq_x C \;\;\; \iff \;\;\; \bA_x \setminus C \preceq_x \bA_x \setminus B.
\]
\end{rem}

\begin{defn} A \emph{strictly monotone preference} on a set $A$ is defined to be a total quasiorder $\preceq$ on $\cP(A)$ such that
\[
B \subset C \subseteq A \;\;\; \implies \;\;\; B \prec C.
\]
Note the strict inequalities!
\end{defn}

Strictly monotone preferences are closely related to vague elements - the main difference is that a strictly monotone preference on $A$ may allow disjoint subsets $B,C$ of $A$ to be in the same equivalence class. The plan is to combine the fact that vague solutions can be converted to solutions if the variable domains have bounded width with the assumption that we have symmetric operations of every arity, by replacing each strictly monotone preference with a weighted average of vague elements by randomly breaking ties in a way that is compatible with homomorphisms. To make this precise, it is convenient to use some of the language of category theory.

If we restrict our attention to diagrammatic CSPs, then ``relaxations'' can be viewed as functors $F$ which take the one-element set to itself. Given a diagram $\cD$ in the category of sets, we can apply the functor $F$ to get a new diagram $F(\cD)$. The new diagram $F(\cD)$ will correspond to an easier CSP than the original diagram $\cD$, since any solution to $\cD$ corresponds to a \emph{cone} over $\cD$, that is, a diagram $\cD'$ consisting of $\cD$ together with a new one-element set with a collection of compatible maps to all of the objects in $\cD$, and $F$ maps any such cone $\cD'$ to a cone $F(\cD')$ over the relaxed diagram $F(\cD)$.

To make things a bit more concrete, let's see how we can describe the basic linear programming relaxation as a functor $\Delta$. The functor $\Delta$ transforms a set by
\[
\Delta(\{1, ..., n\}) = \Big\{(p_1, ..., p_n) \in \RR^n \mid \sum_i p_i = 1 \text{ and } p_i \ge 0 \text{ for all }i\Big\},
\]
and for a function $f : [n] \rightarrow [m]$, we have
\[
\Delta(f) : (p_1, ..., p_n) \mapsto \Big(\sum_{i \in f^{-1}(1)} p_i, ..., \sum_{i \in f^{-1}(m)} p_i\Big).
\]
Let's see what this does if we are given the unsolvable 2-SAT instance
\[
x,y \in \{0,1\} \wedge x = y \wedge x \ne y.
\]
This instance corresponds to a diagram with two sets $\bA_x, \bA_y$ of size two, and two homomorphisms between them, given by $y = x$ and $y = 1-x$. When we apply $\Delta$, we get a diagram with two sets $\Delta(\bA_x), \Delta(\bA_y)$ and two homomorphisms between them:
\[
\Delta(\bA_x) = \{(1-p_x, p_x) \mid 0 \le p_x \le 1\} \;\; \substack{\longrightarrow\\ \longrightarrow} \;\; \Delta(\bA_y) = \{(1-p_y, p_y) \mid 0 \le p_y \le 1\},
\]
where the two homomorphisms correspond to the constraints
\begin{align*}
(1-p_y, p_y) &= (1-p_x, p_x),\\
(1-p_y, p_y) &= (p_x, 1-p_x).
\end{align*}
While the original 2-SAT instance has no solution, the relaxed instance has the solution $p_x = p_y = 1/2$.

For another concrete example of a relaxation considered as a functor, we have arc-consistency. Arc-consistency corresponds to the functor $\cP_\emptyset$, which takes a set $A$ to the set of nonempty subsets $S$ of $A$, and acts on $f : A \rightarrow B$ by
\[
\cP_\emptyset(f) : S \mapsto \{f(a) \mid a \in S\}.
\]
We can make the claim that arc-consistency is more relaxed than the basic linear programming relaxation precise by exhibiting a natural transformation from $\Delta(A)$ to $\cP_\emptyset(A)$. This natural transformation takes a probability distribution to its \emph{support}:
\[
\supp : (p_a)_{a \in A} \mapsto \{a \mid p_a > 0\}.
\]

In order to show that a relaxation (corresponding to a functor $F$) solves $\CSP(\fA)$, where $\fA$ is a (possibly multisorted) relational structure, we need a \emph{rounding rule}. A rounding rule is a homomorphism of relational structures
\[
r : F(\fA) \rightarrow \fA,
\]
where we make sense of $F(\fA)$ by noticing that a relational structure corresponds in a natural way to a diagram, and we can apply functors to diagrams. As long as $F$ maps the one point set to itself, there will automatically be a homomorphism
\[
\fA \rightarrow F(\fA),
\]
so what a rounding rule is really doing is proving that $\fA$ and the relaxation $F(\fA)$ are homomorphically equivalent - which is exactly what it means for the relaxation $F$ to solve $\CSP(\fA)$, if you unwind the definitions (an element of the inverse limit of the diagram $F(\cD_\fX)$ is the same thing as a homomorphism $\fX \rightarrow F(\fA)$). We've already seen that $\CSP(\fA)$ is solved by arc-consistency iff there is a homomorphism $\cP_\emptyset(\fA) \rightarrow \fA$ (Theorem \ref{thm-set-polymorphism}) - now we can view that as a special case of a more general principle. Similarly, if you squint a bit, then Theorem \ref{lp-robust} can be viewed as saying that $\CSP(\fA)$ is solved by the basic linear programming relaxation iff there is a homomorphism from a simplification of $\Delta(\fA)$ (where we only consider probability vectors with rational entries) to $\fA$.

If we view an entire pseudovariety $\cV = HSP_{fin}(\bA_1, ..., \bA_n)$ as a giant multisorted relational structure, then a rounding rule for $F$ on $\cV$ is a natural transformation from $F$, considered as a functor from $\cV$ to the category of sets, to the functor which takes each algebra in $\cV$ to its underlying set. In other words, a rounding rule for $F$ on $\cV$ gives us a map
\[
r_\bA : F(A) \rightarrow A
\]
for each $\bA \in \cV$, where $A$ is the underlying set of $\bA$. These maps need to be compatible with homomorphisms $\pi : \bA \rightarrow \bB$ in the following sense: for any $x \in F(A)$, we must have
\[
\pi(r_\bA(x)) = r_\bB(F(\pi)(x)) \in B.
\]

In order to describe a rounding rule for $F$ on $\cV$ more compactly, it helps to note that we only need to know how the rounding rule behaves on free algebras, since every algebra in $\cV$ is a quotient of a free algebra. If $\cF_\cV(S)$ is the free algebra in $\cV$ with generators corresponding to the elements of $S$, then the map
\[
S \stackrel{i_s}{\hookrightarrow} \cF_\cV(S)
\]
which takes each element of $S$ to the corresponding generator of $\cF_\cV(S)$ turns into a map
\[
F(i_S) : F(S) \rightarrow F(\cF_\cV(S)),
\]
and I claim that all we need to know in order to describe a rounding rule is the result of composing $F(i_S)$ with the map
\[
r_{\cF_\cV(S)} : F(\cF_\cV(S)) \rightarrow \cF_\cV(S).
\]
Why is this the case? Well, for any algebra $\bA \in \cV$, we have maps
\[
A \stackrel{i_A}{\hookrightarrow} \cF_\cV(A) \stackrel{\pi_\bA}{\twoheadrightarrow} \bA,
\]
where $i_A$ is just a map of sets while $\pi_\bA$ is the homomorphism sending each generator of $\cF_\cV(A)$ to the corresponding element of $\bA$. Applying $F$, we get
\[
F(A) \rightarrow F(\cF_\cV(A)) \rightarrow F(\bA),
\]
and the composition of the two maps is the identity since $F$ is a functor. Then for any $x \in F(A)$, we have
\begin{align*}
r_\bA(x) &= r_\bA(F(\pi_\bA \circ i_A)(x))\\
&= r_\bA(F(\pi_\bA)(F(i_A)(x)))\\
&= \pi_\bA(r_{\cF_\cV(A)}(F(i_A)(x))),
\end{align*}
so
\[
r_\bA = \pi_\bA \circ r_{\cF_\cV(A)} \circ F(i_A),
\]
and we see that $r_\bA$ is determined by $r_{\cF_\cV(A)} \circ F(i_A) : F(A) \rightarrow \cF_\cV(A)$, as promised.

\begin{defn} If $F$ is a functor from finite sets to sets, and if $\cV$ is a variety, then a \emph{rounding rule} for $F$ on $\cV$ is defined to be a family of maps
\[
r_S : F(S) \rightarrow \cF_\cV(S)
\]
such that for every map $\pi : S \rightarrow T$ and every $x \in F(S)$ we have
\[
\cF_\cV(\pi)(r_S(x)) = r_T(F(\pi)(x)),
\]
where $\cF_\cV(\pi) : \cF_\cV(S) \rightarrow \cF_\cV(T)$ is the corresponding map of free algebras. In other words, $r$ is a natural transformation (also known as a \emph{minion homomorphism}) from $F$ to $\cF_\cV$.
\end{defn}

Another way of understanding this is that for every variety $\cV$, the free algebra functor $\cF_\cV$ is itself a type of relaxation which can be applied to diagrams, and for diagrams in $\cV$ it doesn't introduce new solutions to problems which didn't already have solutions. So to show that a relaxation $F$ is safe for $\cV$, we just need to show that it is no stronger than the relaxation $\cF_\cV$, and we do this by describing a recipe which turns elements of $F(S)$ into term operations of arity $S$.

The point of all of this was to introduce and justify the language we will use to repackage our results about vague solutions into something we can use to resolve our questions about how reversibility relates to linear programming.

\begin{defn} The \emph{vague minion} is defined to be the functor $\cW$ which takes a finite set $A$ to
\[
\cW(A) = \{(S, \preceq) \mid S \in \cP_\emptyset(A) \text{ and }\preceq\text{ is a vague element of }S\text{ which extends the inclusion order }\subseteq\},
\]
and which takes a map $\pi : A \rightarrow B$ to
\[
\cW(\pi) : (S, \preceq) \mapsto (\pi(S), \preceq^\pi),
\]
where $\preceq^\pi$ is the vague element of $\pi(S)$ defined by
\[
U \preceq^\pi V \;\; \iff \;\; \pi^{-1}(U) \preceq \pi^{-1}(V)
\]
for $U, V \subseteq \pi(S)$. (See Definition \ref{defn-vague} for the definition of vague elements.)
\end{defn}

\begin{thm}\label{thm-vague-minion} If $\cV$ is a locally finite bounded width variety, then there is a rounding rule for the vague minion $\cW$ on $\cV$, that is, there is a minion homomorphism
\[
r^\cW : \cW \rightarrow \cF_\cV.
\]
\end{thm}
\begin{proof} With an eye towards applying Theorem \ref{thm-vague-solution}, we define an (infinite) diagrammatic instance $\fX$ with $\Sg(\fX)$ an instance of $\CSP(\cV)$ as follows:
\begin{itemize}
\item for every finite set $S$ and for every vague element $\preceq$ on $S$ which extends the inclusion order $\subseteq$, introduce a variable $x_{(S,\preceq)}$ with variable domain $A_{x_{(S,\preceq)}} = S$ considered as a subset of the algebra $\bA_{x_{(S,\preceq)}} = \cF_\cV(S)$,
\item for every surjective map $\pi : S \twoheadrightarrow T$ and every vague element $\preceq$ of $S$ extending $\subseteq$, introduce the constraint
\[
x_{(T, \preceq^\pi)} = \pi(x_{(S, \preceq)}).
\]
\end{itemize}
By construction, the instance $\fX$ is arc-consistent, and we claim that it has a vague solution given by sending $x_{(S,\preceq)}$ to the vague element $\preceq$. For this, we just need to check that the pair of vague elements $(\preceq, \preceq^\pi)$ vaguely satisfies the binary relation
\[
\RR_\pi = \{(s, \pi(s)) \mid s \in S\},
\]
where vague satisfaction is defined in Definition \ref{defn-vague}. This is the step where it is important to restrict to vague elements which extend the inclusion order!

We define a total quasiorder $\preceq_{\RR_\pi}$ on $\cP_\emptyset(S) \sqcup \cP_\emptyset(\pi(S))$ by putting each element $U \in \cP_\emptyset(\pi(S))$ into the same equivalence class as $\pi^{-1}(U) \in \cP_\emptyset(S)$, and ordering $\cP_\emptyset(S)$ via $\preceq$. By definition, $\preceq_{\RR_\pi}$ restricts to $\preceq$ on $\cP_\emptyset(S)$ and restricts to $\preceq^\pi$ on $\cP_\emptyset(\pi(S))$, and for any $U \in \cP_\emptyset(\pi(S))$ we have
\[
U \preceq_{\RR_\pi} U - \RR_\pi = \pi^{-1}(U).
\]
Finally, for $V \in \cP_\emptyset(S)$, we have
\[
V \subseteq \pi^{-1}(\pi(V)),
\]
so since we only consider vague elements $\preceq$ which extend the inclusion order $\subseteq$ we have
\[
V \preceq_{\RR_\pi} \pi^{-1}(\pi(V)) \preceq_{\RR_\pi} \pi(V) = V + \RR_\pi.
\]
Therefore, the pair $(\preceq, \preceq^\pi)$ vaguely satisfies the binary relation $\RR_\pi$.

Now we can apply Theorem \ref{thm-vague-solution} to see that $\Sg(\fX)$ has a solution $r$ taking each variable $x_{(S,\preceq)}$ to a stable element of $\cF_\cV(S)$, such that for every surjective map $\pi : S \twoheadrightarrow T$ we have
\[
(r(x_{(S,\preceq)}), r(x_{(T,\preceq^\pi)})) \in \Sg_{\cF_\cV(S)\times \cF_\cV(T)}(\RR_\pi) = \{(f,\cF_\cV(\pi)(f)) \mid f \in \cF_\cV(S)\},
\]
that is,
\[
r(x_{(T,\preceq^\pi)}) = \cF_\cV(\pi)(r(x_{(S,\preceq)})).\qedhere
\]
\end{proof}

\begin{rem} The basic semidefinite relaxation considered in the previous section doesn't obviously fit into the functor framework, since it treats relations differently from the way it treats variable domains. However, if we restrict to diagrammatic instances, then it can be described as the following functor. Let $\cH$ be the infinite-dimensional inner product space
\[
\cH = \Big\{v : \NN \rightarrow \RR \mid \sum_i v_i^2 < \infty\Big\},
\]
and define the semidefinite relaxation functor $\cS$ on sets $A$ by
\[
\cS(A) = \Big\{v : A \rightarrow \cH \mid \sum_{a \in A} \|v(a)\|^2 = 1 \; \wedge \; \forall a \ne b, v(a) \cdot v(b) = 0\Big\},
\]
and on maps $\pi : A \rightarrow B$ by
\[
\cS(\pi) : v \mapsto \Big(b \mapsto \sum_{\pi(a) = b} v(a)\Big).
\]
There is a minion homomorphism
\begin{align*}
\cS &\rightarrow \cW,\\
v \in \cS(A) &\mapsto (\supp(v), \preceq_v) \in \cW(A),
\end{align*}
where $\supp(v)$ is the set of $a \in A$ such that $\|v(a)\|^2 > 0$, and $\preceq_v$ is the vague element of $\supp(v)$ given by
\[
U \prec_v V \;\; \iff \;\; \|v(U)\|^2 < \|v(V)\|^2 \; \vee \; \big(\|v(U)\|^2 = \|v(V)\|^2 \; \wedge \; v(U) <_{lex} v(V)\big),
\]
where $v(U) = \sum_{u \in U} v(u)$ and $<_{lex}$ is the total ordering on $\cH$ which compares vectors based on the first coordinate which differs between them. Composing this with the minion homomorphism
\[
r^\cW : \cW \rightarrow \cF_\cV
\]
from Theorem \ref{thm-vague-minion}, we get a rounding rule for (exact) solutions to the (diagrammatic) basic semidefinite relaxation in any locally finite variety $\cV$ which has bounded width. (The restriction to diagrammatic instances turns out to be inessential: it's a good exercise to check that if the basic SDP relaxation of an instance $\fX$ has a solution, then that solution can always be extended to a solution to $\cS(\cD_\fX)$.)
\end{rem}

Now that the preliminaries are out of the way, we can define the relaxation which is connected to reversibility.

\begin{defn} The \emph{reversible minion} is defined to be the functor $\cM$ which takes a finite set $A$ to
\[
\cM(A) = \{(S, \preceq) \mid S \in \cP_\emptyset(A) \text{ and }\preceq\text{ is a strictly monotone preference on }S\},
\]
and which takes a map $\pi : A \rightarrow B$ to
\[
\cM(\pi) : (S, \preceq) \mapsto (\pi(S), \preceq^\pi),
\]
where $\preceq^\pi$ is the strictly monotone preference on $\pi(S)$ defined by
\[
U \preceq^\pi V \;\; \iff \;\; \pi^{-1}(U) \preceq \pi^{-1}(V)
\]
for $U, V \subseteq \pi(S)$.
\end{defn}

Note that there is a natural transformation
\[
\operatorname{pref} : \Delta \rightarrow \cM
\]
from the basic linear programming relaxation to the reversible minion, given by
\[
\operatorname{pref}(p) = (\supp(p), \preceq_p),
\]
where $\preceq_p$ is the strictly monotone preference on $\supp(p)$ defined by
\[
U \preceq_p V \;\; \iff \;\; \sum_{u \in U} p_u \le \sum_{v \in V} p_v.
\]

\begin{thm}\label{thm-reversible-minion} If $\cV$ is a locally finite variety which has symmetric term operations of every arity, then there is a rounding rule for the reversible minion $\cM$ on $\cV$, that is, there is a minion homomorphism
\[
r^\cM : \cM \rightarrow \cF_\cV.
\]
\end{thm}
\begin{proof} Since no finite affine algebra can have symmetric term operations of all arities, $\cV$ must have bounded width by Theorem \ref{affine-free-pq}, so we can apply Theorem \ref{thm-vague-minion} to see that there is a rounding rule for the vague minion on $\cV$:
\[
r^\cW : \cW \rightarrow \cF_\cV.
\]
Additionally, by Theorem \ref{lp-robust} there is a rounding rule for the basic linear programming relaxation $\Delta$ on $\cV$:
\[
r^\Delta : \Delta \rightarrow \cF_\cV.
\]
We will apply these rounding rules as black-boxes, by finding a minion homomorphism
\[
\mu : \cM \rightarrow \Delta \circ \cW
\]
and composing it with the rounding rule
\[
r^\Delta \circ r^\cW : \Delta \circ \cW \rightarrow \cF_\cV \circ \cF_\cV
\]
and the fact that we can compose term operations
\[
\circ : \cF_\cV \circ \cF_\cV \rightarrow \cF_\cV.
\]

The map $\mu$ is rather simple: for any strictly monotone preference $\preceq$ on $S$, we take every equivalence class $E$ of $\preceq$, and break ties by choosing a uniformly random total ordering on $E$, independently of the total orderings we pick on the other equivalence classes of $\preceq$.

Alternatively, we can describe $\mu$ as follows. Say that a total order $\preceq_!$ on $\cP_\emptyset(S)$ is compatible with $\preceq$ if
\[
U \preceq_! V \;\; \implies \;\; U \preceq V.
\]
Then $\mu(S,\preceq)$ is defined to be the uniform distribution on the set of total orders $\preceq_!$ which are compatible with $\preceq$.

We need to check that $\mu(S,\preceq) \in \Delta(\cW(S))$, that is, that every compatible total order $\preceq_!$ is actually a vague element of $S$ which extends the inclusion order $\subseteq$. That $\preceq_!$ is a vague element follows immediately from the fact that it is a total order, so we just need to check that it extends the inclusion order. For this, note that since $\preceq$ was \emph{strictly} monotone, we have
\[
U \subset V \;\; \implies \;\; U \prec V \;\; \implies \;\; V \not\preceq_! U \;\; \implies \;\; U \prec_! V.
\]

To finish the proof we need to check that $\mu$ is actually a natural transformation, that is, that for any surjective map $\pi : S \twoheadrightarrow T$ we have
\[
(\Delta\circ\cW)(\pi)(\mu(S,\preceq)) = \mu(\cM(\pi)(T,\preceq)).
\]
The left hand side is the probability distribution we get by first picking uniformly among compatible total orders $\preceq_!$ on $\cP_\emptyset(S)$ and then applying $\pi$ to get a total order $(\preceq_!)^\pi$ on $\cP_\emptyset(T)$. The right hand side is the probability distribution we get by first applying $\pi$ to get a strictly monotone preference $\preceq^\pi$ on $T$, and then picking uniformly among compatible total orders $(\preceq^\pi)_!$ on $\cP_\emptyset(T)$.

To check that these two random processes have the same distributions, it's convenient use yet another way to describe the random tie-breaking process: independently for each set $U \subseteq S$, we pick a uniformly random real number $\epsilon_U$ in the range $[0,1]$, and we define $\preceq_!$ by
\[
U \preceq_! V \;\; \iff \;\; U \prec V \vee (U \sim V \wedge \epsilon_U \le \epsilon_V),
\]
where $\sim$ is the equivalence relation associated to $\preceq$. Then if we define $\epsilon_U^\pi$ for $U \subseteq T$ by
\[
\epsilon_U^\pi = \epsilon_{\pi^{-1}(U)},
\]
then the real numbers $\epsilon_U^\pi$ will also be independent (since $U \ne V \implies \pi^{-1}(U) \ne \pi^{-1}(V)$, as $\pi$ is surjective onto $T$) and uniformly distributed in $[0,1]$, and for this coupling between the random numbers we choose for $T$ and the ones we choose for $S$ we have $(\preceq_!)^\pi = (\preceq^\pi)_!$.
\end{proof}

\begin{rem} We could have defined the vague minion and the reversible minion to only involve self-dual quasiorders $\preceq$. Theorem \ref{thm-reversible-minion} would remain true with the revised definitions, with a slight modification to the map $\mu$ which occured in the proof.
\end{rem}

\begin{cor} If $\cV$ is a locally finite variety such that $\CSP(\cV)$ is solved by the basic linear programming relaxation, then every reversible instance of $\CSP(\cV)$ has a solution.
\end{cor}

\begin{cor} Suppose that $\Gamma$ is a finite constraint language such that $\CSP(\Gamma)$ is solved by the basic linear programming relaxation. If we are given an instance $\fX$ of $\CSP(\Gamma)$ such that the basic linear programing relaxation of $\fX$ is $1 - \epsilon$ satisfiable, then we can efficiently find a solution which satisfies a
\[
1 - O\Big(\frac{1}{\log(1/\epsilon)}\Big)
\]
fraction of the constraints of $\fX$.
\end{cor}

\begin{cor}\label{cor-lp-decidable} If $\fA$ is a finite relational structure, then $\CSP(\fA)$ is solved by the basic linear programming relaxation iff there is a homomorphism
\[
\cM(\fA) \rightarrow \fA.
\]
In particular, if $\fA$ has a finite relational signature, then we can decide whether $\CSP(\fA)$ is solved by the basic linear programming relaxation in a finite amount of time.
\end{cor}


While Corollary \ref{cor-lp-decidable} resolves the decidability of the meta-problem for finitely related structures, there are natural examples like Example \ref{lp-not-width-1} of algebraic structures $\bA$ with finitely many basic operations which have symmetric term operations of all arity, but which are not finitely related.

\begin{prob}\label{prob-lp-algebraic} Given a finite algebraic structure $\bA$ as input, can we decide whether or not $\bA$ has symmetric term operations of all arities?
\end{prob}

Define the relaxation $\Delta_n$ to be the subfunctor of $\Delta$ consisting of probability distributions $p$ such that $p_s$ is a multiple of $1/n$ for each $s$. Then it's easy to check that $\fA$ has a symmetric $n$-ary polymorphism iff there is a homomorphism
\[
\Delta_n(\fA) \rightarrow \fA.
\]
Additionally, if we define $\Delta_n^k$ to be the functor
\[
\Delta_n^k : S \mapsto S^n/\!\sim_k,
\]
where $\sim_k$ is the equivalence relation which identifies two sequences $a,b \in S^n$ as long as $b$ is a permutation of $a$ and at most $k$ distinct elements of $S$ occur in $a$, then a homomorphism
\[
\Delta_n^k(\fA) \rightarrow \fA
\]
corresponds to an $n$-ary polymorphism of $\fA$ which is symmetric on all tuples where at most $k$ distinct entries show up. So one approach to resolving Problem \ref{prob-lp-algebraic} is to try to find a minion homomorphism
\[
\bigcup_m \Delta_m^k \rightarrow \Delta_n \circ \cW
\]
for $k \ge |A|$, with $n$ depending only on $k$. Somewhat surprisingly, this approach actually works for $k = 2$ and $k = 3$!

\begin{defn} Call an operation \emph{$k$-symmetric} if it is symmetric on all tuples where at most $k$ distinct values show up.
\end{defn}

\begin{prop} If a locally finite variety $\cV$ has bounded width and has a binary symmetric term operation $s_2(x,y)$, then $\cV$ has $2$-symmetric term operations of every arity.
\end{prop}
\begin{proof} For every $n$, let $\preceq_0$ be the graded lexicographic total ordering on $\cP([n])$, given by
\[
U \prec_0 V \;\; \iff \;\; |U| < |V| \vee \big(|U| = |V| \; \wedge \; \exists k \in V \setminus U \text{ s.t. } [k-1] \cap U = [k-1] \cap V\big),
\]
and let $\preceq_1$ be the graded reverse lexicographic total order, given by
\[
U \prec_1 V \;\; \iff \;\; |U| < |V| \vee \big(|U| = |V| \; \wedge \; \exists k \in U \setminus V \text{ s.t. } [k-1] \cap U = [k-1] \cap V\big).
\]
By Theorem \ref{thm-vague-minion} there is a rounding rule
\[
r^\cW : \cW \rightarrow \cF_\cV.
\]
for the vague minion on $\cV$, and plugging in $\preceq_0, \preceq_1$ we get $n$-ary term operations
\[
r^\cW([n],\preceq_0), r^\cW([n],\preceq_1) \in \cF_\cV(n).
\]
Let $t$ be the $n$-ary term operation we get by applying $s_2$ to these two term operations, i.e.
\[
t = s_2\big(r^\cW([n],\preceq_0), r^\cW([n],\preceq_1)\big) \in \cF_\cV(n).
\]
We claim that $t$ is $2$-symmetric.

For any set $U \subseteq [n]$, we can plug in $x$s and $y$s into $t$ with the $x$s occuring at the positions corresponding to the elements of $U$ to get a binary term operation $t_U \in \cF_\cV(x,y)$. This $t_U$ is given by
\[
t_U = \cF_\cV(\pi_U)(t),
\]
where $\pi_U : [n] \rightarrow \{x,y\}$ is given by
\[
\pi_U(i) = \begin{cases} x & i \in U,\\ y & i \not\in U.\end{cases}
\]
We need to check that the value of $t_U$ only depends on $|U|$. We may as well assume that $|U| \ne 0, n$. Then since $r^\cW$ is a minion homomorphism, it satisfies
\begin{align*}
\cF_\cV(\pi_U)\big(r^\cW([n],\preceq_i)\big) &= r^\cW\big(\cW(\pi_U)([n],\preceq_i)\big)\\
&= r^\cW(\{x,y\},\preceq_i^{\pi_U}).
\end{align*}
Since $\preceq_i^{\pi_U}$ is a total ordering on $\cP(\{x,y\})$ which extends the inclusion order, it is completely determined by how it orders $\{x\}$ and $\{y\}$, which is determined by
\[
\{x\} \preceq_i^{\pi_U} \{y\} \;\; \iff \;\; \pi_U^{-1}(\{x\}) \preceq_i \pi_U^{-1}(\{y\}) \;\; \iff \;\; U \preceq_i [n]\setminus U.
\]
If $|U| \ne n/2$, then whether or not $U \preceq_i [n]\setminus U$ only depends on whether or not $|U| < n/2$, so if we let $s_1(x) = s_2(x,x)$ then we get
\begin{align*}
0 < |U| < n/2 \;\;\; &\implies \;\;\; t_U = s_1\Big(r^\cW\big(\{x,y\}, \emptyset \prec \{x\} \prec \{y\} \prec \{x,y\}\big)\Big),\\
n/2 < |U| < n \;\;\; &\implies \;\;\; t_U = s_1\Big(r^\cW\big(\{x,y\}, \emptyset \prec \{y\} \prec \{x\} \prec \{x,y\}\big)\Big).
\end{align*}
If $|U| = n/2$, then the definitions of $\preceq_0, \preceq_1$ imply that
\[
U \preceq_0 [n]\setminus U \;\;\; \iff \;\;\; [n]\setminus U \preceq_1 U,
\]
so the symmetry of $s_2$ implies that we have
\[
|U| = n/2 \;\;\; \implies \;\;\; t_U = s_2\Big(r^\cW\big(\{x,y\}, \emptyset \prec \{x\} \prec \{y\} \prec \{x,y\}\big), r^\cW\big(\{x,y\}, \emptyset \prec \{y\} \prec \{x\} \prec \{x,y\}\big)\Big).\qedhere
\]
\end{proof}

The combinatorial core of the construction for $3$-symmetric operations is encapsulated in the next lemma, which is a variation of a construction found by Marcin Kozik (personal communication).

\begin{lem}\label{lem-three-orders} For any set $S$, we can find three total orderings $\preceq_0, \preceq_1, \preceq_2$ on $\cP(S)$, such that for every partition of $S$ into three disjoint nonempty sets $U, V, W$ with
\[
W \prec_0 V \prec_0 U,
\]
the restrictions of $\preceq_1, \preceq_2$ to $\{U,V,W\}$ are given by
\begin{align*}
U &\prec_1 W \prec_1 V,\\
V &\prec_2 U \prec_2 W.
\end{align*}
In particular, the three restrictions $\preceq_i|_{\{U,V,W\}}$ together with their three reversals $\succeq_i|_{\{U,V,W\}}$ make up all six total orderings of the set $\{U,V,W\}$, each occuring once.
\end{lem}
\begin{proof} Let $<$ be a well-ordering of the set $S$. Define $\preceq_0$ to be the lexicographic ordering of $\cP(S)$ with respect to $<$, i.e.
\[
U \prec_0 V \;\; \iff \;\; \exists k \in V \setminus U \text{ s.t. } \forall i < k \; (i \in U \iff i \in V).
\]
Suppose (for the sake of keeping the notation under control) that $0$ is the $<$-minimal element in $S$. Then we define $\preceq_1$ to be the ordering given by
\[
U \prec_1 V \;\; \iff \;\; 0 \in U\setminus V \; \vee \; \big((0 \in U \iff 0 \in V) \; \wedge \; U \prec_0 V\big).
\]
Finally, introducing the notation
\[
U^\# = \begin{cases}U & 0 \not\in U,\\ S\setminus U & 0 \in U,\end{cases}
\]
we define $\preceq_2$ to be the ordering given by
\begin{align*}
U \prec_2 V \;\; \iff \;\; &\big(\min(U^\#) < \min(V^\#)\big) \; \vee \; \big(\min(U^\#) = \min(V^\#) \; \wedge \; 0 \in V\setminus U\big)\\
& \; \vee \; \big(\min(U^\#) = \min(V^\#) \; \wedge \; (0 \in U \iff 0 \in V) \; \wedge \; V^\# \prec_0 U^\#\big).
\end{align*}

Now suppose that $S$ is the disjoint union of three nonempty sets $U,V,W$. Exactly one of $U,V,W$ contains $0$, suppose without loss of generality that $0 \in U$. Additionally, since $U \ne S$, exactly one of $V,W$ contains $\min(S\setminus U)$. Suppose without loss of generality that we have $\min(S\setminus U) \in V$, in which case
\[
\min(V) = \min(S\setminus U) < \min(W).
\]
Then the three orderings $\preceq_i$ restricted to $\{U,V,W\}$ are given by
\begin{align*}
W &\prec_0 V \prec_0 U,\\
U &\prec_1 W \prec_1 V,\\
V &\prec_2 U \prec_2 W.\qedhere
\end{align*}
\end{proof}

\begin{thm}\label{thm-3-symmetric} If a locally finite variety $\cV$ has bounded width and has a binary symmetric term operation $s_2(x,y)$ and a ternary cyclic term operation $c_3(x,y,z)$, then $\cV$ has $3$-symmetric term operations of every arity.
\end{thm}
\begin{proof} Fix an arity $n$, and let $\preceq_0, \preceq_1, \preceq_2$ be the three total orderings on $\cP([n])$ from Lemma \ref{lem-three-orders}. For each $i$, define $\preceq_i^+$ and $\preceq_i^-$ by
\[
\preceq_i^+ \; = \; \preceq_i, \;\;\; \preceq_i^- \; = \;  \succeq_i.
\]
For any $i \in \{0,1,2\}$ and any $a,b \in \{+,-\}$, we will define the graded, self-dual total ordering $\preceq_{iab}$ by
\begin{align*}
U \preceq_{iab} V \;\;\; \iff \;\;\; |U| < |V| &\; \vee \; \big(|U| = |V| < n/2 \; \wedge \; U \preceq_i^a V\big)\\
& \; \vee \; \big(|U| = |V| = n/2 \; \wedge \; U \preceq_0^b V\big)\\
& \; \vee \; \big(|U| = |V| > n/2 \; \wedge \; ([n]\setminus V) \preceq_i^a ([n]\setminus U)\big).
\end{align*}
By Theorem \ref{thm-vague-minion} there is a rounding rule
\[
r^\cW : \cW \rightarrow \cF_\cV.
\]
for the vague minion on $\cV$. Define $n$-ary term operations $r_{iab} \in \cF_\cV(n)$ by
\[
r_{iab} = r^\cW([n], \preceq_{iab}),
\]
and define the $n$-ary term operation $t$ by
\begin{align*}
t = s_2(&s_2(c_3(r_{0++}, r_{1++}, r_{2++}), c_3(r_{2-+}, r_{1-+}, r_{0-+})),\\
&s_2(c_3(r_{0+-}, r_{1+-}, r_{2+-}), c_3(r_{2--}, r_{1--}, r_{0--}))).
\end{align*}
We claim that $t$ is $3$-symmetric.

If we apply $t$ to a tuple consisting of variables among $u,v,w$, and if we let $U,V,W \subseteq [n]$ be the sets of inputs assigned to $u,v,w$ respectively, then the result is the ternary term operation
\[
\cF_\cV(\pi_{UVW})(t) \in \cF_\cV(u,v,w),
\]
where $\pi_{UVW} : [n] \rightarrow \{u,v,w\}$ is the map
\[
\pi_{UVW} : i \mapsto \begin{cases}u & i \in U,\\ v & i \in V,\\ w & i \in W.\end{cases}
\]
Since $r^\cW$ is a minion homomorphism, we have
\[
\cF_\cV(\pi_{UVW})(r_{iab}) = r^\cW(\{u,v,w\}, \prec_{iab}^{\pi_{UVW}}),
\]
so $\cF_\cV(\pi_{UVW})(r_{iab})$ is determined by the restriction of $\preceq_{iab}$ to the eight sets of the form $\pi_{UVW}^{-1}(S)$ for $S \subseteq \{u,v,w\}$. Since $\preceq_{iab}$ is monotone and self-dual, so is $\preceq_{iab}^{\pi_{UVW}}$. We will need some notation for the monotone self-dual total orders on $\cP(\{u,v,w\})$.

Up to permuting the variables $u,v,w$, there are just two monotone self-dual total orderings on $\cP(\{u,v,w\})$, which we will name $\preceq_f$ and $\preceq_g$:
\[
\emptyset \preceq_f \{w\} \preceq_f \{v\} \preceq_f \{v,w\} \preceq_f \{u\} \preceq_f \{u,w\} \preceq_f \{u,v\} \preceq_f \{u,v,w\}
\]
and
\[
\emptyset \preceq_g \{w\} \preceq_g \{v\} \preceq_g \{u\} \preceq_g \{v,w\} \preceq_g \{u,w\} \preceq_g \{u,v\} \preceq_g \{u,v,w\}.
\]
Note that $\preceq_f$ is just the lexicographic order, while $\preceq_g$ is the graded lexicographic order. Define ternary term operations $f,g$ by
\[
f(u,v,w) = r^\cW(\{u,v,w\}, \preceq_f), \;\;\; g(u,v,w) = r^\cW(\{u,v,w\}, \preceq_g).
\]
As an aside, these ternary terms $f,g$ will automatically satisfy the identities
\[
f(x,x,y) \approx f(x,y,x) \approx f(x,y,y) \approx g(x,x,y) \approx g(x,y,x) \approx g(y,x,x).
\]

Writing
\[
t_{UVW} = \cF_\cV(\pi_{UVW})(t) \in \cF_\cV(u,v,w),
\]
our goal is to show that the ternary operation $t_{UVW}$ is completely determined by the sizes $|U|, |V|, |W|$. We may as well assume that
\[
|W| \le |V| \le |U|.
\]
There are several cases to consider. To verify that we have covered every case, I suggest drawing a triangle with vertices labeled $U,V,W$, drawing the midpoints of each side, and drawing the line segments connecting the midpoints to the vertices and to each other - each case corresponds to one or more of the smaller triangles, line segments, or points in the resulting diagram.

{\bf Case 1.} Suppose that $|U|, |V|, |W|$ and $|V|+|W|$ are all different from each other. In this case, the orderings $\preceq_{iab}^{\pi_{UVW}}$ do not depend on the choice of $i,a,b$, and if we let $s_1(x) = s_2(x,x)$ and $c_1(x) = c_3(x,x,x)$ then we get
\[
|W| < |V| < |W| + |V| < |U| \;\; \implies \;\; t_{UVW} = s_1(s_1(c_1(f(u,v,w)))
\]
\[
|W| < |V| < |U| < |W| + |V| \;\; \implies \;\; t_{UVW} = s_1(s_1(c_1(g(u,v,w))).
\]

{\bf Case 2.} Suppose that $|U| = |V| + |W| = n/2$ and $|V| \ne |W| \ne 0$. In this case, the orderings $\preceq_{iab}^{\pi_{UVW}}$ only differ in how they order $\{u\}$ with $\{v,w\}$, with the ordering determined by whether or not we have $V \cup W \preceq_0^b U$. Since
\[
V \cup W \prec_0^+ U \;\; \iff \;\; U \prec_0^- V \cup W,
\]
by the symmetry of $s_2$ we have
\[
|W| < |V| < |U| = |W| + |V| \;\; \implies \;\; t_{UVW} = s_2(s_1(c_1(f(u,v,w))), s_1(c_1(g(u,v,w)))).
\]

{\bf Case 3.} Suppose that $|U| = |V| + |W| = n/2$ and $|V| = |W| = n/4$. In this case, the orderings $\preceq_{iab}^{\pi_{UVW}}$ differ in how they order $\{u\}$ with $\{v,w\}$ as well as in how they order $\{v\}$ with $\{w\}$. The ordering between $\{u\}$ and $\{v,w\}$ is once again determined by $b$, but the ordering between $\{v\}$ and $\{w\}$ is a bit more interesting. By Lemma \ref{lem-three-orders}, either one out of three of the $i$s will have $W \preceq_i V$ and two out of three of the $i$s will have $V \preceq_i^- W$, or vice-versa, so we get
\begin{align*}
|W| = |V| < |U| =\; &|W| + |V| \;\; \implies\\
t_{UVW} = s_1(&s_2(c_3(f(u,v,w), f(u,v,w), f(u,w,v)), c_3(f(u,v,w), f(u,w,v), f(u,w,v))),\\
&s_2(c_3(g(u,v,w), g(u,v,w), g(u,w,v)), c_3(g(u,v,w), g(u,w,v), g(u,w,v)))).
\end{align*}

{\bf Case 4.} Suppose that $|V| = |W| \ne 0$ and $|U| \ne |V| + |W|$. By a similar argument to the previous case, we get
\begin{align*}
|W| = |V| < |U| <\; &|W| + |V| \;\; \implies\\
t_{UVW} = s_1(&s_2(c_3(g(u,v,w), g(u,v,w), g(u,w,v)), c_3(g(u,v,w), g(u,w,v), g(u,w,v))))
\end{align*}
and
\begin{align*}
|W| = |V| < |W| +\; &|V| < |U| \;\; \implies\\
t_{UVW} = s_1(&s_2(c_3(f(u,v,w), f(u,v,w), f(u,w,v)), c_3(f(u,v,w), f(u,w,v), f(u,w,v)))).
\end{align*}

{\bf Case 5.} Suppose that $|U| = |V| \ne |W| \ne 0$. In this case we get
\begin{align*}
|W| < |V| = |U| <\; &|W| + |V| \;\; \implies\\
t_{UVW} = s_1(&s_2(c_3(g(u,v,w), g(u,v,w), g(v,u,w)), c_3(g(u,v,w), g(v,u,w), g(v,u,w)))).
\end{align*}

{\bf Case 6.} Suppose that $|U| = |V| = |W| = n/3$. In this case, $b$ becomes irrelevant and Lemma \ref{lem-three-orders} shows that for $a$ fixed, the restrictions of the orders $\preceq_{iab}^{\pi_{UVW}}$ to $\{u\},\{v\},\{w\}$ are cyclic permutations of each other, so we get
\begin{align*}
|W| = |V| =\; &|U| \;\; \implies\\
t_{UVW} = s_1(&s_2(c_3(g(u,v,w), g(v,w,u), g(w,u,v)), c_3(g(v,u,w), g(u,w,v), g(w,v,u)))).
\end{align*}

{\bf Case 7.} Suppose that $|W| = 0$. If $|V| = 0$ as well, then we get
\[
0 = |W| = |V| < |U| = n \;\; \implies \;\; t_{UVW} = s_1(s_1(c_1(g(u,u,u)))),
\]
and if $0 < |V| < |U|$ then we get
\[
0 = |W| < |V| < |U| \;\; \implies \;\; t_{UVW} = s_1(s_1(c_1(g(u,u,v)))),
\]
while if $|U| = |V| = n/2$ then we get
\[
0 = |W| < |V| = |U| \;\; \implies \;\; t_{UVW} = s_2(s_1(c_1(g(u,u,v))), s_1(c_1(g(u,v,v)))).\qedhere
\]
\end{proof}

We can take things a little bit further, but how much further is unclear. We'll start with a construction of $4$-ary symmetric operations from \cite{symmetric-polymorphisms}, the existence of which is closely connected to the fact that the symmetric group $S_4$ on four elements is solvable.

\begin{prop}[Lemma 4 of \cite{symmetric-polymorphisms}] If $s_2(x,y)$ is a binary symmetric operation and $c_3(x,y,z)$ is a ternary cyclic operation, then
\[
s_3(x,y,z) \coloneqq s_2(c_3(x,y,z),c_3(z,y,x))
\]
is a ternary symmetric operation, and if we define $4$-ary operations $t, s_4$ by
\begin{align*}
t(x,y,z,w) &\coloneqq s_2(s_2(x,y), s_2(z,w)),\\
s_4(x,y,z,w) &\coloneqq s_3(t(x,y,z,w), t(x,z,w,y), t(x,w,y,z)),
\end{align*}
then $s_4$ is a $4$-ary symmetric operation.
\end{prop}

\begin{prop} If a locally finite variety $\cV$ has bounded width and has a binary symmetric term operation $s_2(x,y)$ and a ternary cyclic term operation $c_3(x,y,z)$, then $\cV$ has a $4$-symmetric term operation of arity $5$.
\end{prop}
\begin{proof} For every pair of permutations $\sigma, \tau \in S_5$, we define a graded self-dual total quasiorder $\preceq_{\sigma,\tau}$ on $\cP([5])$ as follows. First, as usual if $|U| < |V|$ then we have $U \prec_{\sigma,\tau} V$. Next, if $U$ and $V$ are singletons, then we define $\preceq_{\sigma,\tau}$ by
\[
\{u\} \preceq_{\sigma,\tau} \{v\} \;\;\; \iff \;\;\; \sigma(u) \le \sigma(v).
\]
For sets of size $2$, we define $\preceq_{\sigma,\tau}$ by
\[
\{u_1,u_2\} \preceq_{\sigma,\tau} \{v_1,v_2\} \;\;\; \iff \;\;\; \tau(\min(\sigma(u_1), \sigma(u_2))) \le \tau(\min(\sigma(v_1), \sigma(v_2))),
\]
so that the associated equivalence relation $\sim_{\sigma,\tau}$ on sets of size $2$ is given by
\[
\{u_1,u_2\} \sim_{\sigma,\tau} \{v_1,v_2\} \;\;\; \iff \;\;\; \min(\sigma(u_1),\sigma(u_2)) = \min(\sigma(v_1),\sigma(v_2)).
\]
The ordering on sets of size $3$ and $4$ is then defined to be the dual to the ordering on the sets of size $2$ and $1$:
\[
U \preceq_{\sigma,\tau} V \;\;\; \iff \;\;\; [5]\setminus V \preceq_{\sigma,\tau} [5]\setminus U.
\]
Then $\preceq_{\sigma,\tau}$ is a vague element of $[5]$ which extends the inclusion order $\subseteq$ for each pair $\sigma,\tau \in S_5$, so by Theorem \ref{thm-vague-minion} we can find associated $5$-ary term operations
\[
r_{\sigma,\tau} \coloneqq r^\cW([5], \preceq_{\sigma,\tau}) \in \cF_\cV(5).
\]
Letting $\iota \in S_5$ be the order-reversing permutation
\[
\iota : i \mapsto 6 - i,
\]
we define
\[
r_{\sigma,\{\tau,\iota\tau\}} \coloneqq s_2(r_{\sigma,\tau}, r_{\sigma,\iota\tau}).
\]
Next, by Theorem \ref{thm-3-symmetric} we can find in $\cV$ a $3$-symmetric term operation $f_{60}$ of arity $60$, and we define
\[
r_\sigma \coloneqq f_{60}(r_{\sigma,\{\tau_1,\iota\tau_1\}}, ..., r_{\sigma,\{\tau_{60},\iota\tau_{60}\}})
\]
for any ordering $\tau_i$ of a system of right coset representatives for the subgroup $\langle \iota \rangle < S_5$. To finish the construction, we define
\begin{align*}
r_{\{\sigma,\iota\sigma\}} &\coloneqq s_2(r_\sigma, r_{\iota\sigma}),\\
t_5 &\coloneqq f_{60}(r_{\{\sigma_1,\iota\sigma_1\}}, ..., r_{\{\sigma_{60},\iota\sigma_{60}\}}).
\end{align*}
In order to check that $t_5$ is a $4$-symmetric operation, we need to check that as long as two of the inputs to $t_5$ are equal, each application of $f_{60}$ has at most $3$ distinct inputs - this is left as an exercise to the reader.
\end{proof}

Similar techniques can be used to construct $4$-symmetric term operations of arities $6$ and $7$, but I haven't found a way to construct a $4$-symmetric term operation of arity $8$ with this approach.


\chapter{Finite Taylor Algebras}\label{chapter-taylor}

\section{Cyclic terms}

In this section we will prove that every finite Taylor algebra has a cyclic term.

\begin{defn} An $m$-ary operation $c$ is called \emph{cyclic} if it satisfies the identity
\[
c(x_1, x_2, ..., x_m) \approx c(x_2, ..., x_m, x_1).
\]
\end{defn}

Cyclic terms were first proved to exist for finite congruence modular algebras \cite{congruence-modular-cyclic}, and most of the basic facts about cyclic terms are developed in that paper. This was extended to finite congruence join-semidistributive algebras in \cite{sd-join-cyclic}, and then finally to all finite Taylor algebras in \cite{cyclic}. We'll start by showing that we only care about cyclic terms of prime arity.

\begin{prop}[Multiplicative property of cyclic terms \cite{congruence-modular-cyclic}] A variety $\cV$ has a cyclic term $c_{mn}$ of arity $mn$ if and only if $\cV$ has cyclic terms $c_m, c_n$ of arity $m$ and $n$, respectively.
\end{prop}
\begin{proof} Suppose first that $c_{mn}$ is a cyclic term operation of arity $mn$. Then we can define a cyclic term operation of arity $m$ by plugging in
\[
c_{mn}(\underbrace{x_1, ..., x_1}_{n}, \underbrace{x_2, ..., x_2}_{n}, ..., \underbrace{x_m, ..., x_m}_{n}),
\]
and we can define a cyclic term operation of arity $n$ similarly.

Conversely, suppose that $c_m, c_n$ are cyclic terms of arity $m$ and $n$. We define a cyclic term operation of arity $mn$ by renumbering the inputs of the star composition $c_n*c_m$:
\[
c_n\begin{pmatrix}c_m(x_1, & x_{n+1}, & ..., & x_{(m-1)n+1}),\\
c_m(x_2, & x_{n+2}, & ..., & x_{(m-1)n+2}),\\
\vdots & \vdots & \ddots & \vdots\\
c_m(x_n, & x_{2n}, & ..., & x_{mn})
\end{pmatrix} \approx
c_n\begin{pmatrix} c_m(x_2, & x_{n+2}, & ..., & x_{(m-1)n+2}),\\
\vdots & \vdots & \ddots & \vdots\\
c_m(x_n, & x_{2n}, & ..., & x_{mn}),\\
c_m(x_{n+1}, & ..., & x_{(m-1)n+1}, & x_1)
\end{pmatrix}.\qedhere
\]
\end{proof}

\begin{cor} A variety $\cV$ has a cyclic term operation of arity $m$ if and only if $\cV$ has a cyclic term operation of arity $p$ for every prime $p$ which divides $m$.
\end{cor}

Next we will describe the main obstruction to the existence of a cyclic term operation of a given arity.

\begin{prop}[Semantic meaning of cyclic terms \cite{congruence-modular-cyclic}] Suppose that $\cV$ is a variety. Then for any $m \in \NN$, the following are equivalent.
\begin{itemize}
\item[(a)] $\cV$ has no cyclic term operation of arity $m$.
\item[(b)] There is some $\bA \in \cV$ and an automorphism $\sigma \in \Aut(\bA)$ such that $\sigma^m = 1$ and $\sigma$ has no fixed point.
\end{itemize}
\end{prop}
\begin{proof} We start by showing that (b) implies (a). Suppose that $\bA, \sigma$ are as in (b), and suppose for contradiction that $\bA$ has some cyclic term $c_m$ of arity $m$. Let $a$ be any element of $\bA$, and define $a_i$ by $a_i = \sigma^i(a)$. Then we have
\[
c_m\begin{pmatrix}a_1, & a_2, & ..., & a_m\\ a_2, & a_3, & ..., & a_1\end{pmatrix} \in \sigma,
\]
so $c_m(a_1, ..., a_m)$ is a fixed point of $\sigma$, contradicting the assumption in (b).

Now suppose that (a) holds. Let $\bA = \cF_\cV(x_1, ..., x_m)$ be the free algebra on $m$ generators, and let $\sigma$ be the automorphism of $\bA$ defined by cyclically permuting the generators $x_1, ..., x_m$. Then a fixed point of $\sigma$ is precisely the same thing as a cyclic term operation of $\cV$ of arity $m$, so if $\cV$ has no cyclic term operation of arity $m$, then $\sigma$ has no fixed points (and satisfies $\sigma^m = 1$).
\end{proof}

For finite algebras, we can give a local criterion for the existence of a cyclic term.

\begin{prop}[Local criterion for cyclic terms \cite{congruence-modular-cyclic}]\label{cyclic-local-crit} If $\bA$ is a finite algebra, then $\bA$ has an $m$-ary cyclic term if and only if it is the case that for all $a_1, ..., a_m \in \bA$, there exists some $m$-ary term $t$ such that
\[
t(a_1, a_2, ..., a_m) = t(a_2, ..., a_m, a_1) = \cdots = t(a_m, a_1, ..., a_{m-1}).
\]
\end{prop}
\begin{proof} Say that an $m$-ary term $t$ is cyclic for a tuple $(a_1, ..., a_m)$ if it satisfies the displayed equation from the statement of the proposition. Let $c$ be an $m$-ary term which is cyclic for a maximal set of tuples (we are using finiteness of $\bA$ here). Suppose for contradiction that $c$ is not cyclic, and let $a = (a_1, ..., a_m)$ be any tuple such that $c$ is not cyclic for $a$.

Define a tuple $a' = (a_1', ..., a_m')$ by
\[
a_i' = c(a_i, a_{i+1}, ..., a_{i-1}),
\]
with indices taken modulo $m$. By assumption, there is some $m$-ary term $t$ which is cyclic for $a'$. But then the $m$-ary term
\[
t(c(x_1, x_2 ..., x_m), c(x_2, ..., x_m, x_1), ..., c(x_m, x_1, ..., x_{m-1}))
\]
is cyclic for $a$, and is also cyclic for every tuple which $c$ was cyclic for, contradicting the maximality assumption on $c$.
\end{proof}

For the sake of checking the local condition of Proposition \ref{cyclic-local-crit} for a particular tuple $a_1, ..., a_m$, the natural approach is to compute the $m$-ary relation
\[
\Sg_{\bA^m}\left\{\begin{bmatrix} a_1 & a_2 & \cdots & a_m\\ a_2 & a_3 & \cdots & a_1\\ \vdots & \vdots & \iddots & \vdots\\ a_m & a_1 & \cdots & a_{m-1}\end{bmatrix}\right\},
\]
and to check if it contains any constant tuples. This relation is invariant under cyclically permuting its coordinates, which leads us to make the following definition.

\begin{defn} A relation $\RR \le \bA^m$ is called \emph{cyclic} if $\RR$ is invariant under cyclically permuting its coordinates, that is,
\[
(a_1, a_2, ..., a_m) \in \RR \;\; \iff \;\; (a_2, ..., a_m, a_1) \in \RR.
\]
\end{defn}

\begin{cor}[Relational criterion for cyclic terms \cite{congruence-modular-cyclic}] If $\bA$ is a finite algebra, then $\bA$ has a cyclic term operation of arity $m$ if and only if every $m$-ary cyclic relation $\RR \le \bA^m$ contains a constant tuple.
\end{cor}

Now we are finally ready to prove one of the main results of \cite{cyclic}, which states that every finite Taylor algebra $\bA$ has cyclic terms of every prime arity $p > |\bA|$. In fact, we will prove a stronger version of this result due to Zhuk (currently unpublished).

\begin{thm}[Finite Taylor algebras have cyclic terms \cite{cyclic}, refined by Zhuk]\label{cyclic-terms} Suppose $\bA$ is a finite idempotent Taylor algebra and that $p$ is prime. Then one of the following is true:
\begin{itemize}
\item[(a)] either $\bA$ has a cyclic term operation of arity $p$, or
\item[(b)] there is some $\bB \in HS(\bA)$ and some automorphism $\sigma \in \Aut(\bB)$ such that $\sigma^p = 1$ and $\sigma$ has no fixed points.
\end{itemize}
In particular, if $p > |\bA|$ then $\bA$ has a cyclic term operation of arity $p$.
\end{thm}
\begin{proof} We prove this by induction on $|\bA|$. Suppose that there is no $\bB \in HS(\bA), \sigma \in \Aut(\bB)$ as in (b), and let $\RR \le \bA^p$ be any $p$-ary cyclic relation. It's enough to show that $\RR$ contains a constant tuple.

If $\pi_1(\RR) \ne \bA$, then since $\RR$ is cyclic we have $\RR \le \pi_1(\RR)^p$, so we can apply the induction hypothesis to the algebra $\pi_1(\RR)$ to see that $\RR$ has a constant tuple. Thus we may assume without loss of generality that $\RR$ is subdirect in $\bA^p$.

If $\bA$ has a nontrivial congruence $\theta \in \Con(\bA)$, then $\RR/\theta^p \le (\bA/\theta)^p$ is a cyclic relation on $\bA/\theta$, so by the induction hypothesis applied to $\bA/\theta$ there is some congruence class $a/\theta$ such that $\RR \cap (a/\theta)^p \ne \emptyset$. Setting $\RR' = \RR \cap (a/\theta)^p$, we see that $\RR'$ is a cyclic relation on $a/\theta$, so by the induction hypothesis applied to $a/\theta$ we see that $\RR'$ (and therefore also $\RR$) has a constant tuple. Thus we may assume that $\bA$ is simple.

If any $\pi_{ij}(\RR)$ is the graph of an automorphism $\sigma$ of $\bA$, then since $\RR$ is cyclic, we see that $\pi_{j,2j-i}(\RR)$ is also the graph of $\sigma$, and similarly so is $\pi_{2j-i,3j-2i}(\RR)$, etc., so
\[
\pi_{ii}(\RR) = \pi_{ij}(\RR) \circ \pi_{j,2j-i}(\RR) \circ \cdots \circ \pi_{2i-j,i}(\RR)
\]
is the graph of $\sigma^p$, which implies $\sigma^p = 1$. Since $p$ is prime, we see that in fact every $\pi_{kl}(\RR)$ is the graph of some power of the automorphism $\sigma$. In this case we see that $\RR$ has a constant tuple if and only if $\sigma$ has a fixed point. Thus we may assume without loss of generality that every $\pi_{ij}(\RR)$ is linked.

By Zhuk's four cases (Corollary \ref{zhuk-four-cases}), we see that $\bA$ is either affine, subdirectly simple, or has a proper ternary absorbing subalgebra.

If $\bA$ is affine, with underlying abelian group $(A, +, -, 0)$, then since $x-y+z$ is a term operation of $\bA$ (by the definition of an affine algebra), we see that $k_1x_1 + \cdots + k_mx_m$ is a term operation of $\bA$ for all $k_1, ..., k_m \in \ZZ$ such that $k_1 + \cdots + k_m \equiv 1 \pmod{|\bA|}$. In particular, if $p \nmid |\bA|$, then
\[
p^{-1}(x_1 + \cdots + x_p)
\]
is a $p$-ary cyclic term operation of $\bA$. On the other hand, if $p \mid |\bA|$, then by elementary group theory there must be some element $c \in \bA$ of order $p$, and then by the idempotence of $\bA$ the relation
\[
\{(x,y) \mid x = y + c\}
\]
is a subalgebra of $\bA^2$, and it is then the graph of an automorphism $\sigma$ of $\bA$ which has order $p$ and has no fixed points.

If $\bA$ is subdirectly simple, then since $\RR \le_{sd} \bA^p$ is subdirect and every $\pi_{ij}(\RR)$ is linked, we must have $\RR = \bA^p$. In this case $\RR$ contains \emph{every} constant tuple.

If $\bA$ has a proper ternary absorbing subalgebra, then we define a directed graph $\fD$ whose vertices are proper ternary absorbing subalgebras $\bB \lhd_3 \bA$, and with a directed edge $(\bB,\bC)$ whenever there is some $i \ne j$ with $\bB + \pi_{ij}(\RR) \subseteq \bC \lhd_3 \bA$.

{\bf Claim:} The digraph $\fD$ has no directed cycles.

{\bf Proof of claim:} Note first that since $\RR$ is cyclic we have
\[
\pi_{ij}(\RR)^- \subseteq \pi_{ij}(\RR)^{\circ (p-1)},
\]
so if $\bB + \pi_{ij}(\RR) \subseteq \bB$ then we must have
\[
\bB + \pi_{ij}(\RR) - \pi_{ij}(\RR) \subseteq \bB + \pi_{ij}(\RR)^{\circ p} \subseteq \bB,
\]
so
\[
\bB + \pi_{ij}(\RR) - \pi_{ij}(\RR) = \bB,
\]
so $\bB$ is a union of linked components of $\pi_{ij}(\RR)$. Since $\pi_{ij}(\RR)$ is linked and $\bB$ is proper, this is impossible. Thus $\fD$ has no directed cycles of length $1$. Since $\RR$ is cyclic, we also have
\[
\pi_{ij}(\RR)\circ \pi_{kl}(\RR) \supseteq \pi_{i+k,j+l}(\RR),
\]
so if $\fD$ has a directed cycle, then $\fD$ must have a directed cycle of length at most $2$, say
\[
\bB + \pi_{ij}(\RR) + \pi_{kl}(\RR) \subseteq \bB.
\]
If $i+k \ne j+l$ this gives us a directed cycle of length $1$, while if $i+k = j+l$ then we have $\bB + \pi_{ij}(\RR) - \pi_{ij}(\RR) \subseteq \bB$, so once again $\bB$ must be a union of linked components of $\pi_{ij}(\RR)$, which is impossible. The claim is proved.

Since the digraph $\fD$ is finite, nonempty, and has no directed cycles, there must be a proper ternary absorbing subalgebra $\bB \lhd_3 \bA$ such that
\[
\bB + \pi_{ij}(\RR) = \bA
\]
for all $i \ne j$. In particular, we see that $\pi_{ij}(\RR) \cap \bB^2 \ne \emptyset$ for all $i,j$. Since $\bB$ is ternary absorbing, this implies that in fact $\RR \cap \bB^p \ne \emptyset$ by Corollary \ref{cor-absorption-essential}. Setting $\RR' = \RR \cap \bB^p$, we can apply the induction hypothesis to $\bB$ to see that $\RR'$ contains a constant tuple. Thus $\RR$ contains a constant tuple, and we are done.
\end{proof}

\begin{cor} If $\bA$ is a finite Taylor algebra and $m$ has no prime factors $p$ which are less than or equal to $|\bA|$, then $\bA$ has an idempotent $m$-ary cyclic term.
\end{cor}

\begin{ex} Let $\bA_n$ be the dual discriminator algebra from Example \ref{ex-dual-discriminator} on a domain of size $n$. Then every subset of $\bA_n$ is a subalgebra with full automorphism group, so $\bA_n$ does not have cyclic terms of any arity between $2$ and $n$. By the previous results, we see that $\bA_n$ has a cyclic term operation of arity $m$ if and only if $m$ has no prime factors which are less than or equal to $n$.
\end{ex}

\begin{cor}[Siggers term from cyclic term] If $\bA$ is a finite Taylor algebra, then $\bA$ has an idempotent $4$-ary Siggers term $t$, satisfying the identity $t(x,x,y,z) \approx t(y,z,z,x)$.
\end{cor}
\begin{proof} Let $c$ be an idempotent $m$-ary cyclic term for some $m > 1$. Then there are numbers $a,b \in \NN$ such that $2a+3b = m$, and we may define $t$ by
\[
t(x,y,z,w) \coloneqq c(\underbrace{x,...,x}_{b}, \underbrace{y,...,y}_{a}, \underbrace{z,...,z}_{b}, \underbrace{w,...,w}_{a+b}).\qedhere
\]
\end{proof}

\begin{cor}[Daisy chain terms]\label{daisy-chain-terms} If $\bA$ is a finite Taylor algebra, then there are idempotent terms $w_i(x,y,z)$ for $i \in \ZZ$ such that for all $i$ we have
\[
w_i(x,x,y) \approx w_i(y,x,x) \approx w_{i-1}(x,y,x),
\]
and the sequence of terms $w_i$ is periodic with some finite period.
\end{cor}
\begin{proof} Choose $p$ to be an extremely huge prime, let $c$ be an idempotent $p$-ary cyclic term, and let $a = \lfloor\frac{p}{3}\rfloor$. Define a long sequence of numbers $a_0, a_1, ...$ by $a_0 = a$ and
\[
a_{i+1} = p - 2a_i,
\]
stopping as soon as we hit the first $a_i$ with $a_i > \frac{p}{2}$. Define terms $w_i'$ by
\[
w_i'(x,y,z) \coloneqq c(\underbrace{x,...,x}_{a_i}, \underbrace{y,...,y}_{a_{i+1}}, \underbrace{z,...,z}_{a_i}).
\]
Since $c$ is cyclic, these $w_i'$s will satisfy the identities
\[
w_i'(x,x,y) \approx w_i'(y,x,x) \approx w_{i-1}'(x,y,x).
\]
If $p$ is large enough, then there must be some $j < k$ such that $w_j' = w_k'$. Then we define $w_i$ by picking some $i' \in [j,k]$ such that $i \equiv i' \pmod{k-j}$ and setting $w_i = w_{i'}'$.
\end{proof}

\begin{cor}\label{consecutive-daisy} A finite algebra $\bA$ is Taylor if and only if it has a pair of idempotent ternary terms $p,q$ satisfying the identities
\begin{align*}
p(x,x,y) &\approx p(y,x,x),\\
q(x,x,y) &\approx q(y,x,x) \approx p(x,y,x).
\end{align*}
\end{cor}
\begin{proof} To see that such $p,q$ must exist in a Taylor algebra, we can take $p,q$ to be any pair of consecutive daisy chain terms from the previous corollary. To see that any such $p,q$ define Taylor terms, note that if $p$ is a projection then $p$ must be second projection, but in this case $q$ must be a Mal'cev term.
\end{proof}



\section{Minimal Taylor clones}

Since our main aim in these notes is to understand the most general CSPs which can be solved in polynomial time, it makes sense to study (core) relational structures $\fA$ such that $\CSP(\fA)$ is in P, but such that adding any additional relations to $\fA$ makes the problem NP-complete. According to the CSP dichotomy theorem of Bulatov \cite{bulatov-dichotomy} and Zhuk \cite{zhuk-dichotomy}, these maximal relational structures correspond under the $\Inv-\Pol$ Galois correspondence to \emph{minimal} Taylor clones.

\begin{defn} A clone $\cO$ on a finite domain is called a \emph{minimal Taylor clone} if $\cO$ is Taylor and every proper subclone of $\cO$ is not Taylor. A finite algebra $\bA$ is called a \emph{minimal Taylor algebra} if $\Clo(\bA)$ is a minimal Taylor clone.
\end{defn}

At first it may not be clear that minimal Taylor clones even exist: perhaps every Taylor clone contains a proper Taylor subclone, with the relevant Taylor operations having higher and higher arity. We can rule this out by using the existence of a Siggers term (Corollary \ref{siggers-term}).

\begin{prop} Every Taylor clone on a finite domain contains a minimal Taylor clone.
\end{prop}
\begin{proof} By Corollary \ref{siggers-term}, every Taylor clone contains a $4$-ary Siggers operation $t$ satisfying the identity $t(x,x,y,z) \approx t(y,z,z,x)$. Since any such $t$ is Taylor, and since there are only finitely many $4$-ary operations on a given finite domain, at least one of the $4$-ary Siggers operations $t \in \cO$ generates a minimal Taylor clone.
\end{proof}

Since every minimal Taylor clone is generated by a single $4$-ary operation, we see that the number of minimal Taylor clones on a domain of size $n$ is at most $n^{n^4}$. We can get a much better upper bound on the number of minimal Taylor clones by using the daisy chain terms from the previous section.

\begin{prop} The number of minimal Taylor clones on a domain of size $n$ is at most $n^{2n^3}$.
\end{prop}
\begin{proof} By Corollary \ref{consecutive-daisy}, every minimal Taylor clone $\cO$ contains a pair of ternary idempotent operations $p,q$ satisfying the identities
\begin{align*}
p(x,x,y) &\approx p(y,x,x),\\
q(x,x,y) &\approx q(y,x,x) \approx p(x,y,x).
\end{align*}
Since $\langle p,q\rangle$ generates a Taylor clone, we must have $\cO = \langle p,q\rangle$. Since the number of ordered pairs of ternary operations $p,q$ on a domain of size $n$ is $n^{2n^3}$, we see that the number of minimal Taylor clones is at most $n^{2n^3}$.
\end{proof}

\begin{rem}
The paper \cite{minimal-taylor} showed that the upper bound $n^{2n^3}$ can be reduced to $n^{n^3}$, by showing that every minimal Taylor clone is generated by a \emph{single} ternary operation. On a domain of size $2$, it is easy to check that every minimal Taylor algebra is term equivalent to either a semilattice, a majority algebra, or to the idempotent reduct of $\ZZ/2$. On a domain of size $3$, there turn out to be a total of $24$ minimal Taylor algebras, up to term equivalence and isomorphism.

Unfortunately, the number of minimal Taylor algebras grows quite rapidly as the size of the domain increases: even if we only consider majority algebras, it turns out that the number of minimal majority algebras (up to term-equivalence) such that every three-element subset is a subalgebra is $7^{\binom{n}{3}}$, and identifying isomorphic algebras can only reduce this by a factor of at most $n!$, which makes little difference to the asymptotics.
\end{rem}

The key fact that makes the theory of minimal Taylor algebras work is the following result, which essentially says that anything that ``looks like'' it ``could be'' a subalgebra or quotient of a minimal Taylor algebra actually \emph{is} a subalgebra or quotient, and is also minimal Taylor as well.

\begin{thm}\label{minimal-taylor-subalg} If $\bA$ is a minimal Taylor algebra and $\bB \in HSP_{fin}(\bA)$, then $\bB$ is also a minimal Taylor algebra.

In fact, if $S \subseteq \bB$ is a subset of $\bB$ (not assumed to be a subalgebra), $t \in \Clo(\bA)$ is any term operation of $\bA$, and $\theta$ is an equivalence relation on $S$ such that
\begin{itemize}
\item the set $S$ is closed under $t$,
\item every equivalence class of $\theta$ is a subalgebra of $\bB$,
\item the equivalence relation $\theta$ is a congruence of the algebraic structure $(S,t)$, and
\item the quotient $(S,t)/\theta$ is a Taylor algebra,
\end{itemize}
then in fact the following must all be true:
\begin{itemize}
\item the set $S$ is actually the underlying set of a subalgebra $\bS$ of $\bB$ (i.e. $S$ is a subuniverse of $\bB$),
\item the equivalence relation $\theta$ is actually a congruence on the subalgebra $\bS$, and
\item the restriction of every term operation of $\bA$ to the quotient $\bS/\theta$ is in the clone generated by the restriction of $t$ to $(S,t)/\theta$.
\end{itemize}
Note that taking $\theta$ to be the trivial equivalence relation $0_S$ is always allowed, since every minimal Taylor algebra is automatically idempotent.
\end{thm}
\begin{proof} Let $p$ be any prime such that $p > |\bA|$ and $p > |(S,t)/\theta|$. By Theorem \ref{cyclic-terms}, there is a $p$-ary cyclic term $c \in \Clo(\bA)$, as well as a $p$-ary term $u \in \Clo(t)$ such that the restriction of $u$ to $(S,t)/\theta$ is cyclic. Define a $p$-ary term $c'$ by
\[
c'(x_1, ..., x_p) \coloneqq c(u(x_1, ..., x_p), u(x_2, ..., x_p, x_1), ..., u(x_p, x_1, ..., x_{p-1})).
\]
Then since $c$ is cyclic, $c'$ will automatically be cyclic as well. Since $\bA$ is assumed to be minimal Taylor, we must have $\Clo(\bA) = \langle c'\rangle$.

Suppose that $x_1, ..., x_p \in S$. Then since $u \in \Clo(t)$ preserves $S$ and acts cyclically on $(S,t)/\theta$, we must have
\[
u(x_1, ..., x_p) \equiv_\theta u(x_2, ..., x_p, x_1) \equiv_\theta \cdots \equiv_\theta u(x_p, x_1, ..., x_{p-1}) \in S,
\]
and since equivalence classes of $\theta$ were assumed to be subalgebras of $\bB$, we have
\[
c'(x_1, ..., x_p) \equiv_\theta u(x_1, ..., x_p) \in S.
\]
Thus $c'$ preserves $S$ as well as the equivalence relation $\theta$, and the restriction of $c'$ to $(S,t)/\theta$ is the same as the restriction of $u$ to $(S,t)/\theta$. Since $c'$ generates $\Clo(\bA)$, this finishes the proof.
\end{proof}

An immediate consequence of Theorem \ref{minimal-taylor-subalg} is that minimal Taylor algebras are \emph{prepared} in the sense of Definition \ref{defn-prepared}.

\begin{prop}\label{prop-minimal-prepared} If $\bA$ is a minimal Taylor algebra, then $a, b \in \bA$ have
\[
\begin{bmatrix} b\\ b\end{bmatrix} \in \Sg_{\bA^2}\left\{\begin{bmatrix} a\\ b\end{bmatrix}, \begin{bmatrix} b\\ a\end{bmatrix}\right\}
\]
if and only if $\{a,b\}$ is a semilattice subalgebra of $\bA$ with absorbing element $b$.
\end{prop}
\begin{proof} If $(b,b) \in \Sg\{(a,b),(b,a)\}$, then there must be some binary term $t$ such that $t(a,b) = t(b,a) = b$. By idempotence, we automatically have $t(a,a) = a$ and $t(b,b) = b$, so the set $S = \{a,b\}$ is closed under $t$ and $(S,t)$ is a two-element semilattice. Thus we can apply Theorem \ref{minimal-taylor-subalg} to see that $\{a,b\}$ must be a subalgebra of $\bA$, and that the restriction of every term operation of $\bA$ to $\{a,b\}$ is in the clone generated by the restriction of $t$ to $\{a,b\}$.
\end{proof}

Similarly, we can recognize two-element majority subalgebras and $\ZZ/2^{\aff}$ subalgebras. To simplify the statements of these results, it is convenient to assume the existence of an order two automorphism.

\begin{prop}\label{prop-recognition-with-automorphism} If $\bA$ is a minimal Taylor algebra and $a, b \in \bA$ are such that $\Sg_{\bA^2}\left\{\begin{bmatrix} a\\ b\end{bmatrix}, \begin{bmatrix} b\\ a\end{bmatrix}\right\}$ is the graph of an automorphism of order two, then
\begin{itemize}
\item we have $\begin{bmatrix} a\\ a\\ a\end{bmatrix} \in \Sg_{\bA^3}\left\{\begin{bmatrix} a\\ a\\ b\end{bmatrix}, \begin{bmatrix} a\\ b\\ a\end{bmatrix}, \begin{bmatrix} b\\ a\\ a\end{bmatrix}\right\}$ iff $\{a,b\}$ is a majority subalgebra of $\bA$, and
\item we have $\begin{bmatrix} b\\ b\\ b\end{bmatrix} \in \Sg_{\bA^3}\left\{\begin{bmatrix} a\\ a\\ b\end{bmatrix}, \begin{bmatrix} a\\ b\\ a\end{bmatrix}, \begin{bmatrix} b\\ a\\ a\end{bmatrix}\right\}$ iff $\{a,b\}$ is a $\ZZ/2^{\aff}$ subalgebra of $\bA$.
\end{itemize}
\end{prop}

\begin{cor}\label{not-sd-complete} If a minimal Taylor algebra $\bA$ is generated by two elements $a,b$, then $\bA$ is \emph{not} subdirectly simple. As a consequence, either $\bA$ has an affine quotient or $\bA$ has a proper ternary absorbing subalgebra.
\end{cor}
\begin{proof} Suppose for contradiction that $\bA$ is subdirectly simple. Define a subdirect binary relation $\bS \le_{sd} \bA^2$ by
\[
\bS = \Sg_{\bA^2}\left\{\begin{bmatrix} a\\ b\end{bmatrix}, \begin{bmatrix} b\\ a\end{bmatrix}\right\}.
\]
If $(a,a)$ or $(b,b)$ is in $\bS$, then $\{a,b\}$ must be a two-element semilattice, which is not subdirectly simple. Otherwise, $\bS$ must be the graph of an automorphism of order two by our assumption that $\bA$ is subdirectly simple. Now define a subdirect ternary relation $\RR \le_{sd} \bA^3$ by
\[
\RR = \Sg_{\bA^3}\left\{\begin{bmatrix} a\\ a\\ b\end{bmatrix}, \begin{bmatrix} a\\ b\\ a\end{bmatrix}, \begin{bmatrix} b\\ a\\ a\end{bmatrix}\right\}.
\]
Since no $\pi_{ij}(\RR)$ can be the graph of an automorphism, we see that we must have $\RR = \bA^3$ by our assumption that $\bA$ is subdirectly simple. Thus we have $(a,a,a) \in \RR$, so $\{a,b\}$ must be a two-element majority algebra, which is not subdirectly simple. This contradiction proves that $\bA$ must not be subdirectly simple.

For the last claim, we recall Zhuk's four cases (Corollary \ref{zhuk-four-cases}), and note that both binary absorption and central absorption imply ternary absorption.
\end{proof}

\begin{prob} Suppose that a minimal Taylor algebra $\bA$ is generated by two elements. Is it possible for $\bA$ to be polynomially complete?
\end{prob}

The general recognition theorem for two-element majority subalgebras is as follows.

\begin{prop} If $\bA$ is a minimal Taylor algebra, then $a, b \in \bA$ have
\[
\begin{bmatrix} a & b\\ a & b\\ a & b\end{bmatrix} \in \Sg_{\bA^{3\times 2}}\left\{\begin{bmatrix} a & b\\ a & b\\ b & a\end{bmatrix}, \begin{bmatrix} a & b\\ b & a\\ a & b\end{bmatrix}, \begin{bmatrix} b & a\\ a & b\\ a & b\end{bmatrix}\right\}
\]
if and only if $\{a,b\}$ is a majority subalgebra of $\bA$.
\end{prop}

It is also easy to recognize copies of the free semilattice on two generators.

\begin{prop} If $\bA$ is a minimal Taylor algebra and $a,b,c \in \bA$ satisfy $a \rightarrow c$, $b \rightarrow c$ (i.e. $\{a,c\}$ and $\{b,c\}$ are semilattice subalgebras of $\bA$ with absorbing element $c$), then we have
\[
\begin{bmatrix} c\\ c\end{bmatrix} \in \Sg_{\bA^2}\left\{\begin{bmatrix} a\\ b\end{bmatrix}, \begin{bmatrix} b\\ a\end{bmatrix}\right\}
\]
if and only if $\{a,b,c\}$ is isomorphic to the free semilattice on two generators.
\end{prop}

We can also characterize binary absorbing subalgebras of minimal Taylor algebras, and show that they are always automatically strongly absorbing (and therefore are automatically centrally absorbing as well).

\begin{prop}\label{prop-minimal-binary-strong} Suppose that $\bA$ is a minimal Taylor algebra, and that $\bB \lhd_{bin} \bA$ is a binary absorbing subalgebra of $\bA$. Then the following must hold.
\begin{itemize}
\item[(a)] $\bB$ is a strongly absorbing subalgebra of $\bA$, that is, any term $f \in \Clo(\bA)$ which depends on its first input satisfies $f(\bB, \bA, ..., \bA) \subseteq \bB$.
\item[(b)] There is an equivalence relation $\theta_\bB \in \Con(\bA)$ such that $\bB$ is a congruence class of $\theta_\bB$, and all other congruence classes of $\theta_\bB$ are singletons.
\item[(c)] For every $a \not\in \bB$, $\bB \cup \{a\}$ is a subalgebra of $\bA$, and $(\bB \cup \{a\})/\theta_\bB$ is a two-element semilattice with absorbing element $\bB/\theta_\bB$.
\item[(d)] For every $a \not\in \bB$, there is some $b \in \bB$ such that $\{a,b\}$ is a two-element semilattice with absorbing element $b$.
\item[(e)] For every $a,b \not\in \bB$ such that $\Sg_\bA\{a,b\} \cap \bB \ne \emptyset$, $\Sg_{\bA}\{a,b\}/\theta_\bB$ is isomorphic to the free semilattice on two generators.
\item[(f)] For every $a_1, ..., a_k \not\in \bB$ such that $\Sg_\bA\{a_i,a_j\} \cap \bB \ne \emptyset$ for all $i \ne j$, $\Sg_{\bA}\{a_1, ..., a_k\}/\theta_\bB$ is isomorphic to a semilattice of size $k+1$.
\end{itemize}
In particular, if $\bA$ is generated by two elements and $\bB$ is a proper binary absorbing subalgebra, then $\bA/\theta_\bB$ is either a two-element semilattice, or is isomorphic to the free semilattice on two generators.
\end{prop}

For the sake of concretely writing down minimal Taylor algebras, we should pick convenient terms. My preference is to write them down in terms of the daisy chain terms from Corollary \ref{daisy-chain-terms}.

\begin{defn} We say that a sequence of idempotent ternary terms $w_i$, defined for all $i \in \ZZ$, is a sequence of \emph{daisy chain terms} if it satisfies the following properties:
\begin{itemize}
\item the sequence $w_i$ is purely periodic in $i$ with some finite period, and
\item for all $i \in \ZZ$, we have $w_i(x,x,y) \approx w_i(y,x,x) \approx w_{i-1}(x,y,x)$.
\end{itemize}
\end{defn}

It is useful to work out all possible sequences of daisy chain terms in our three basic examples of minimal Taylor algebras: semilattices, majority algebras, and affine algebras.

\begin{prop} If $\bA = (A,\vee)$ is a semilattice, then any sequence of daisy chain terms of $\bA$ must have
\[
w_i(x,y,z) \approx x \vee y \vee z
\]
for all $i \in \ZZ$.
\end{prop}
\begin{proof} It's enough to show that $w_i(x,x,y) \approx w_i(x,y,x) \approx w_i(y,x,x) \approx x \vee y$ for all $i$. Note that since $w_i(x,x,y) \approx w_i(y,x,x)$, we can't have $w_i(x,x,y) = y$, since semilattices have no Mal'cev terms. Additionally, if we had $w_i(x,x,y) = w_i(y,x,x) = x$, then $w_i$ could not depend on its first or last coordinates, so we would have $w_{i+1}(x,x,y) \approx w_{i+1}(y,x,x) \approx w_i(x,y,x) = y$, which again contradicts the fact that semilattices have no Mal'cev terms.

Since the only binary terms of a semilattice are $x, y,$ and $x \vee y$, we see by process of elimination that we must have $w_i(x,x,y) \approx w_i(y,x,x) \approx x \vee y$, and similar reasoning shows that $w_i(x,y,x) \approx w_{i+1}(x,x,y) \approx x \vee y$, so we are done.
\end{proof}

\begin{prop} If $\bA = (A,m)$ is a majority algebra, then in any sequence of daisy chain terms of $\bA$, each $w_i$ must be a majority term, that is, we have
\[
w_i(x,x,y) \approx w_i(x,y,x) \approx w_i(y,x,x) \approx x
\]
for all $i \in \ZZ$.
\end{prop}
\begin{proof} Note that every ternary term operation of a majority algebra is either a projection or a majority term (as is easily checked by induction on the construction of the term in terms of the majority operation $m$). If some $w_i$ is a projection, then the identity $w_i(x,x,y) \approx w_i(y,x,x)$ implies that it must be second projection, but then the identity $w_{i+1}(x,x,y) \approx w_{i+1}(y,x,x) \approx w_i(x,y,x) = y$ implies that $w_{i+1}$ is a Mal'cev term, which is impossible. Thus each $w_i$ must be a majority term.
\end{proof}

For the affine case, we have the following simplification in the setting of minimal Taylor algebras.

\begin{prop} If $\bA$ is minimal Taylor and affine, then there is an abelian group structure on the underlying set $A$ such that $\bA$ is term equivalent to $(A,x-y+z)$.
\end{prop}
\begin{proof} Every affine algebra has the ternary function $x - y + z$ as a term, by Proposition \ref{affine-malcev}. Since the ternary operation $x-y+z$ is Mal'cev, it generates a Taylor clone, so a minimal Taylor algebra is affine if and only if its clone is generated by $x-y+z$.
\end{proof}

Because of this result, we don't need to think about the general case of a module over a (possibly noncommutative) ring if we are only interested in minimal Taylor algebras: we only need to think about algebras of the form $(A, x-y+z)$ for $A$ an abelian group. By the classification of finite abelian groups, we can write such an algebra as a product of cyclic factors of prime power order. Recall that the \emph{exponent} of a group is the least number $n$ such that every cyclic subgroup has order dividing $n$.

\begin{prop}\label{affine-daisy} If $\bA = (A,x-y+z)$ is an affine algebra such that the abelian group $A$ has exponent $n$, then for any sequence of daisy chain terms $w_i$ of $\bA$, there is a sequence of elements $a_i \in \ZZ/n$ such that
\[
w_i(x,y,z) \approx a_ix + (1-2a_i)y + a_iz
\]
and
\[
a_{i+1} \equiv 1-2a_i \pmod{n}
\]
for all $i \in \ZZ$.
\end{prop}
\begin{proof} Every $m$-ary term $t \in \Clo(x-y+z)$ can be written in the form
\[
t(x_1, ..., x_m) \approx k_1x_1 + \cdots + k_mx_m,
\]
for some $k_i \in \ZZ$ satisfying
\[
k_1 + \cdots + k_m = 1.
\]
Of course, only the congruence classes of the values of the coefficients $k_i$ modulo $n$ matter, and the set of $m$-ary terms $t \in \Clo(\bA)$ is in bijection with the set of tuples of $k_i \in \ZZ/n$ such that $k_1 + \cdots + k_m \equiv 1 \pmod{n}$.

Thus we can write
\[
w_i(x,y,z) \approx a_ix + b_iy + c_iz
\]
for some $a_i,b_i,c_i \in \ZZ/n$ such that $a_i + b_i + c_i \equiv 1 \pmod{n}$. The identity $w_i(x,x,y) \approx w(y,x,x)$ then implies that $a_i \equiv c_i$, so $b_i \equiv 1 - 2a_i$, while the identity $w_{i+1}(y,x,x) \approx w_i(x,y,x)$ implies that $a_{i+1} \equiv b_i \equiv 1 - 2a_i$.
\end{proof}

\begin{prop} If $\bA = (A,x-y+z)$ is an affine algebra such that $|A|$ is a power of $2$, then any sequence of daisy chain terms of $\bA$ must have
\[
w_i(x,y,z) \approx \frac{x+y+z}{3}
\]
for all $i \in \ZZ$. In particular, if the abelian group $A$ has exponent $2$, then each $w_i$ is the Mal'cev operation $x-y+z \approx x+y+z$.
\end{prop}
\begin{proof} Suppose the exponent of $A$ is $2^k$. Then if we let $a_i \in \ZZ/2^k$ be the sequence from Proposition \ref{affine-daisy}, we see from $a_{i+1} \equiv 1 - 2a_i \pmod{2^k}$ that we have
\[
a_{i+1} - 1/3 \equiv -2(a_i - 1/3) \pmod{2^k}
\]
for all $i \in \ZZ$, so in fact we must have
\[
a_i - 1/3 \equiv 0 \pmod{2^k}
\]
for all $i \in \ZZ$.
\end{proof}

\begin{prop} If $\bA = (A,x-y+z)$ is an affine algebra such that the abelian group $A$ has exponent $3^k$, then any sequence of daisy chain terms of $\bA$ must have period equal to $3^k$, and there must be some $i \in \ZZ$ such that
\begin{align*}
w_{i-1}(x,y,z) &\approx \frac{x+z}{2},\\
w_i(x,y,z) &\approx y,\\
w_{i+1}(x,y,z) &\approx x-y+z.
\end{align*}
\end{prop}
\begin{proof} We just need to show that the map $a \mapsto 1-2a \pmod{3^k}$ defines a cyclic permutation of $\ZZ/3^k$ for all $k \ge 0$. To see this, it's enough to check that the cycle containing $0$ has length exactly $3^k$.

Lifting to the integers, if we define a sequence $a_i \in \ZZ$ by $a_0 = 0$ and $a_{i+1} = 1-2a_i$, then we can solve the recurrence to obtain
\[
a_i = \frac{1 - (-2)^i}{3}.
\]
Then we see that $a_i \equiv 0 \pmod{3^k}$ if and only if $3^{k+1} \mid 1 - (-2)^i$. By induction on $k$ we may assume that $3^{k-1}$ divides $i$, and by the binomial theorem, we have
\[
1 - (-2)^i = 1 - (1-3)^i = 3i - 9\binom{i}{2} + 27\binom{i}{3} -+ \cdots \equiv 3i - 0 + 0 -+ \cdots \pmod{3^{k+1}},
\]
so $3^{k+1} \mid 1 - (-2)^i$ if and only if $3^k$ divides $i$.
\end{proof}

By going back to the original construction of the daisy chain terms from a huge cyclic term, we can simplify the situation slightly for affine algebras of odd order.

\begin{prop} If $\bA$ is a minimal Taylor algebra, then it is possible to choose a sequence of daisy chain terms $w_i$ of $\bA$ such that for every affine $\bB \in HSP_{fin}(\bA)$ of odd order, the restriction of $w_1$ to $\bB$ is the Mal'cev operation $x-y+z$.
\end{prop}
\begin{proof} Since there are only finitely many ternary terms $w_1 \in \Clo(\bA)$, it's enough to prove that for every finite $k$ we can find a $w_1$ that is part of a sequence of daisy chain terms of $\bA$, such that $w_1$ restricts to the Mal'cev operation $x-y+z$ on every affine $\bB \in HSP(\bA)$ such that $|\bB|$ is odd and $|\bB| \le k$.

Note that for any large prime $p$, the restriction of a $p$-ary cyclic term $c$ to $\bB$ must be given by
\[
c(x_1, ..., x_p) = \frac{x_1 + \cdots + x_p}{p}.
\]
Thus, in the construction of the terms $w_i'$ from Corollary \ref{daisy-chain-terms} where we plugged in
\[
w_i'(x,y,z) \coloneqq c(\underbrace{x,...,x}_{a_i}, \underbrace{y,...,y}_{p-2a_i}, \underbrace{z,...,z}_{a_i}),
\]
we will have
\[
w_i'(x,y,z) = \frac{a_ix + (p-2a_i)y + a_iz}{p}
\]
on $\bB$. So as long as we choose $a_1$ such that $a_1 \equiv p \pmod{k!}$ and $a_1 \approx \frac{p}{3}$ (which is possible as long as we take $p$ much larger than $k!$), we will have $w_1'(x,y,z) = x-y+z$ on every affine algebra $\bB \in HSP(\bA)$ of size at most $k$. For $|\bB|$ odd, the restriction of the sequence of terms $w_i'$ to $\bB$ will be purely periodic, so the final sequence of daisy chain terms constructed will have $w_1 = w_1'$ on such $\bB$.
\end{proof}

\begin{rem} A similar argument shows that we can instead choose daisy chain terms $w_i$ such that $w_i(x,y,z) \approx \frac{x+y+z}{3}$ for all $i$ on every affine algebra $\bB \in HSP(\bA)$ such that $|\bB|$ is not a multiple of $3$. In fact, for any profinite integer $a \in \hat{\ZZ} = \varprojlim \ZZ/n$ such that $a \equiv \frac{1}{3} \pmod{2^k}$ for all $k$, we can choose daisy chain terms such that $w_0(x,y,z) \approx ax + (1-2a)y + az$ on every affine algebra $\bB \in HSP(\bA)$.
\end{rem}

We can also limit the collection of affine algebras which may show up in $HSP(\bA)$.

\begin{prop} If $\bA$ is a minimal Taylor algebra and $\bB \in HSP(\bA)$ is affine, then the exponent $n$ of $\bB$ is finite, with $n \le |\bA|^{|\bA|^2}$, and every prime $p$ which divides $n$ is bounded by $|\bA|$.
\end{prop}
\begin{proof} First we show that every $n$ such that $\ZZ/n^{\aff} \in HSP(\bA)$ has $n \le |\bA|^{|\bA|^2}$. To see this, note that $\ZZ/n^{\aff}$ is generated by two elements (to be more specific, it is generated by $0$ and $1$), so if $\ZZ/n^{\aff} \in HSP(\bA)$ then $\ZZ/n^{\aff}$ must be a quotient of the free algebra on two generators $\cF_\bA(x,y) \le \bA^{\bA^2}$.

Next, note that if $p$ is prime and $\ZZ/p^{\aff} \in HSP(\bA)$, then $\bA$ can't have any cyclic term operation of arity $p$, since $\ZZ/p^{\aff}$ has an automorphism of order $p$ with no fixed points. Thus by Theorem \ref{cyclic-terms} there is no prime $p > |\bA|$ such that $\ZZ/p^{\aff} \in HSP(\bA)$.
\end{proof}

We will end this section by characterizing Zhuk's centrally absorbing subalgebras in the case of minimal Taylor algebras, and using them to naturally produce majority subquotients of minimal Taylor algebras. First, we need a quick detour to show that (up to term equivalence) we can put any minimal Taylor algebra into a minimal Taylor variety which also contains a majority algebra.

\begin{prop} If $\bA$, $\bB$ are minimal Taylor algebras, then there are minimal Taylor algebras $\bA', \bB'$ such that
\begin{itemize}
\item $\bA$ is term equivalent to $\bA'$ and $\bB$ is term equivalent to $\bB'$,
\item $\bA'$ and $\bB'$ have the same signature, and
\item $\bA' \times \bB'$ is minimal Taylor.
\end{itemize}
\end{prop}
\begin{proof} We may assume without loss of generality that the basic operations of $\bA$ and $\bB$ are a sequence of daisy chain terms $w_i(x,y,z)$ by Corollary \ref{daisy-chain-terms}. In this case, $\bA$ and $\bB$ already have the same signature, and by the daisy chain identities we see that $\bA \times \bB$ is Taylor (alternatively, we can apply Lemma \ref{strictly-simple-hs} to see directly that if $\bA, \bB$ are idempotent Taylor algebras with the same signature, then $\bA \times \bB$ is also Taylor). Pick any minimal Taylor reduct $\bA' \times \bB'$ of $\bA \times \bB$ to complete the proof.
\end{proof}

\begin{thm}\label{thm-minimal-central} Suppose that $\bA$ is minimal Taylor, and that there is a majority algebra $\bM$ on the domain $\{0,1\}$ (with the same signature as $\bA$) such that $\bA \times \bM$ is minimal Taylor. For any $\bC \le \bA$, the following are equivalent:
\begin{itemize}
\item[(a)] $\bC$ is a ternary absorbing subalgebra of $\bA$,
\item[(b)] for every prime $p > |\bA|$ there is a $p$-ary cyclic term $c$ of $\bA \times \bM$ such that whenever $\#\{i \mid x_i \in \bC\} > \frac{p}{2}$, we have
\[
c(x_1, ..., x_p) \in \bC,
\]
and furthermore the restriction of $c$ to $\bM$ is the $p$-ary majority operation,
\item[(c)] the binary relation $\RR \subseteq \bA \times \bM$ given by
\[
\RR = (\bA \times \{0\}) \cup (\bC \times \{0,1\})
\]
is a subalgebra of $\bA \times \bM$,
\item[(d)] $\bC$ centrally absorbs $\bA$,
\item[(e)] every daisy chain term $w_i(x,y,z)$ witnesses the fact that $\bC$ ternary absorbs $\bA$.
\end{itemize}
\end{thm}
\begin{proof} For (a) implies (b): let $t$ be a ternary term which witnesses $\bC \lhd \bA$. If $m$ is a ternary term operation of $\bA$ which acts as majority on $\bM$, then the ternary term
\[
t'(x,y,z) \coloneqq m(t(x,y,z), t(y,z,x), t(z,x,y))
\]
also witnesses the absorption $\bC \lhd \bA$, and the restriction of $t'$ to $\bM$ is the majority operation. Now let $p > |\bA|$ be prime, and let $u \in \Clo(t')$ be any term such that the restriction of $u$ to $\bM$ is a $p$-ary majority operation. Any such $u$ must have the property that whenever $\#\{i \mid x_i \in \bC\} > \frac{p}{2}$, we have
\[
u(x_1, ..., x_p) \in \bC.
\]
Now let $c'$ be any $p$-ary cyclic term operation of $\bA$, and define $c$ by
\[
c(x_1, ..., x_p) \coloneqq c'(u(x_1, ..., x_p), u(x_2, ..., x_p, x_1), ..., u(x_p, x_1, ..., x_{p-1})).
\]

For (b) implies (c), note that the cyclic term $c$ must generate the clone of $\bA$, so it's enough to check that the relation $\RR$ is preserved by $c$, which is easy to prove directly.

That (c) implies (d) follows from Zhuk's Corollary \ref{zhuk-center}, since the left center of $\RR$ is $\bC$ and the majority algebra $\bM$ is binary absorption free.

That (c) implies (e) follows from a direct computation: we have
\[
\begin{bmatrix}w_i(\bC,\bA,\bC)\\ 1\end{bmatrix} = w_i\left(\begin{bmatrix}\bC\\ 1\end{bmatrix}, \begin{bmatrix}\bA\\ 0\end{bmatrix}, \begin{bmatrix}\bC\\ 1\end{bmatrix}\right) \subseteq \RR,
\]
so $w_i(\bC,\bA,\bC) \subseteq \bC$, and similarly $w_i(\bA,\bC,\bC), w_i(\bC,\bC,\bA) \subseteq \bC$.

That (d) implies (a) follows from Zhuk's Corollary \ref{center-ternary}, while (e) implies (a) is immediate.
\end{proof}

\begin{cor}\label{cor-minimal-central-majority} Suppose that $\bA$ is minimal Taylor and that $\bC, \bD \lhd_Z \bA$ are two ternary absorbing subalgebras of $\bA$. Then $\bC \cup \bD$ is a subalgebra of $\bA$.

If $\bC \cap \bD = \emptyset$, then the equivalence relation $\theta$ on $\bC \cup \bD$ with parts $\bC$ and $\bD$ is a congruence on $\bC \cup \bD$, and $(\bC \cup \bD)/\theta$ is isomorphic to the two element majority algebra.
\end{cor}
\begin{proof} Note that since $\bA$ is minimal Taylor, the clone of $\bA$ is generated by any pair of consecutive daisy chain terms, so we just need to check that $\bC \cup \bD$ is closed under each daisy chain term $w_i$. For any $a,b,c \in \bC\cup \bD$, we either have at least two of $a,b,c$ in $\bC$ or at least two of $a,b,c$ in $\bD$, so the fact that $w_i$ witnesses both $\bC \lhd \bA$ and $\bD \lhd \bA$ implies that $w_i(a,b,c) \in \bC \cup \bD$.

If $\bC \cap \bD = \emptyset$, then the fact that $w_i$ witnesses both $\bC \lhd \bA$ and $\bD \lhd \bA$ implies that $w_i$ is compatible with $\theta$, and that the restriction of $w_i$ to the two-element algebra $(\bC \cup \bD)/\theta$ is the majority operation.
\end{proof}

The following conjecture, if true, would wrap everything up quite neatly.

\begin{conj} If $\bA$ is a minimal Taylor algebra which is generated by two elements $a,b \in \bA$, then at least one of the following is true:
\begin{itemize}
\item there is a congruence $\theta \in \Con(\bA)$ such that $\bA/\theta$ is an affine algebra of prime order,
\item there is a congruence $\theta \in \Con(\bA)$ such that $\bA/\theta$ is a two element semilattice,
\item there is a congruence $\theta \in \Con(\bA)$ such that $\bA/\theta$ is a two element majority algebra, or
\item there are proper ternary absorbing subalgebras $\bC, \bD \lhd_Z \bA$ such that $a \in \bC, b\in \bD$, $\bC \cup \bD = \bA$, and $\bC \cap \bD \ne \emptyset$.
\end{itemize}
\end{conj}

Partial progress towards this conjecture was made in \cite{minimal-taylor}: if $\bA$ is minimal Taylor, generated by two elements, and has no affine quotient, then at least one of the two generators must be contained in a proper ternary absorbing subalgebra. Embarrassingly, we don't even know the answer to the following basic question.

\begin{prob} Is there any minimal Taylor algebra which is simple, is generated by two elements, has size at least $3$, and is not affine?
\end{prob}

A brute force search found no examples of size $3$ or $4$ - but this search was carried out by hand and never written up, so I may have made a mistake.



\section{Bulatov's colored graph}

In Bulatov's approach to the CSP dichotomy conjecture \cite{bulatov-dichotomy}, the theory of absorbing subalgebras isn't used. Instead, Bulatov introduces a colored graph in \cite{colored-graph-prelim} and \cite{colored-graph}, and uses connectivity properties of this graph to analyze finite Taylor algebras.

\begin{defn} Suppose $\bA$ is a finite idempotent algebra, and $a, b$ are any pair of distinct elements of $\bA$.
\begin{itemize}
\item We say that $(a,b)$ is a \emph{semilattice edge} if there is a binary term $t$ such that $t(a,b) = t(b,a) = b$.
\item We say that $\{a,b\}$ is a \emph{weak majority edge} if there is a congruence $\theta$ on $\Sg\{a,b\}$ and a ternary term $m$ such that $\{a/\theta,b/\theta\}$ is closed under $m$ and $(\{a/\theta,b/\theta\},m)$ is a two-element majority algebra.
\item We say that $\{a,b\}$ is a \emph{weak affine edge} if there is a congruence $\theta$ on $\Sg\{a,b\}$ and a term $p$ such that $(\Sg\{a,b\}/\theta, p)$ is an affine algebra.
\end{itemize}
We drop the modifier ``weak'' on an edge if $\theta$ is a maximal congruence on $\Sg\{a,b\}$, and for any $a', b' \in \Sg\{a,b\}$ such that $a' \equiv_\theta a$ and $b' \equiv_\theta b$, we have $\Sg\{a,b\} = \Sg\{a',b'\}$. Note that semilattice edges are directed, while majority and affine edges are undirected.
\end{defn}

Note that a semilattice edge might not be a subalgebra, and similarly the set $a/\theta \cup b/\theta$ might not be a subalgebra if $\{a,b\}$ is a majority edge. If $\bA$ is a \emph{minimal} Taylor algebra, however, then Theorem \ref{minimal-taylor-subalg} shows that $a/\theta \cup b/\theta$ is a subalgebra if $(a,b)$ is a weak majority edge, and similarly for semilattice edges. In \cite{bulatov-dichotomy}, Bulatov calls an algebra \emph{sm-smooth} if this special case of Theorem \ref{minimal-taylor-subalg} applies to it.

We could have also defined ``weak semilattice edges'' in a similar way to the way we defined weak majority edges, but this is unnecessary.

\begin{prop}\label{prop-weak-semilattice-edge} If there is a congruence $\theta$ on $\Sg\{a,b\}$ and an idempotent binary term $t$ such that $t(a,b) \equiv_\theta t(b,a) \equiv_\theta b$, then there is some $b' \in \Sg\{a,b\}$ such that $b' \equiv_\theta b$ and a partial semilattice term $s \in \Clo(t)$ such that $s(a,b') = s(b',a) = b'$.
\end{prop}

Bulatov defines his colored graph by coloring the semilattice edges red, coloring the majority edges yellow, and coloring the affine edges blue (I don't know why these particular colors were chosen).

\begin{defn} We say that a finite idempotent algebra $\bA$ has a \emph{hereditarily connected colored graph} if for all $\bB \le \bA$, the colored graph of $\bB$ is connected (ignoring the directions on the semilattice edges).
\end{defn}

For the purposes of checking if an algebra is hereditarily connected, weak edges are interchangeable with edges by the following result.

\begin{prop} Let $\bA$ be a finite idempotent algebra. If the colored graph of weak edges of $\bA$ is hereditarily connected, then the colored graph of edges of $\bA$ is also hereditarily connected.
\end{prop}
\begin{proof} Suppose that $\{a,b\}$ is a weak edge of $\bA$, with corresponding congruence $\theta$. We will prove by induction on $|\Sg\{a,b\}|$ that $a$ and $b$ are connected in the colored graph of edges of $\Sg\{a,b\}$. We may enlarge $\theta$ to a maximal congruence on $\Sg\{a,b\}$ without loss of generality, since any congruence of $\Sg\{a,b\}$ which identifies $a$ and $b$ is the full congruence. If $(a,b)$ is not an edge, then we may pick $a' \in a/\theta$ and $b' \in b/\theta$ such that $\Sg\{a',b'\}$ is strictly smaller than $\Sg\{a,b\}$, and by the inductive hypothesis we see that $a',b'$ are connected in the colored graph of edges of $\Sg\{a,b\}$. Since $\Sg\{a,a'\} \subseteq a/\theta$ and $\Sg\{b,b'\} \subseteq b/\theta$, we also see by the inductive hypothesis that $a$ is connected to $a'$ and $b$ is connected to $b'$ in the colored graph of edges of $\Sg\{a,b\}$.
\end{proof}

Algebras with hereditarily connected colored graphs are closed under the usual algebraic operations.

\begin{prop}\label{hereditarily-connected-product} If $\bA, \bB$ are finite idempotent algebras (of the same signature) with hereditarily connected colored graphs, then so is $\bA \times \bB$. More generally, if $\bA$ is a finite idempotent algebra and $\theta \in \Con(\bA)$ is such that $\bA/\theta$ is hereditarily connected and every congruence class of $\theta$ is also hereditarily connected, then $\bA$ is hereditarily connected.
\end{prop}
\begin{proof} We prove the more general statement. Let $a,b \in \bA$ be any pair of elements. We will show that $a$ and $b$ are connected by weak edges in $\Sg\{a,b\}$. If $a/\theta = b/\theta$, then $\Sg\{a,b\}$ is contained in a congruence class of $\theta$, so $a,b$ are connected by edges of $\Sg\{a,b\}$. Otherwise, since $a/\theta, b/\theta$ are connected by edges in $\Sg\{a,b\}/\theta$, we can find a sequence of elements $a = a_0, a_1, ..., a_n = b$ such that $(a_i/\theta,a_{i+1}/\theta)$ is an edge of $\bA/\theta$ for all $i$. Then each $(a_i,a_{i+1})$ will be a weak edge of $\Sg\{a,b\}$, with the corresponding congruence containing $\theta$.
\end{proof}

\begin{prop} If $\bA$ is a finite idempotent algebra with a hereditarily connected colored graph and $\theta \in \Con(\bA)$, then $\bA/\theta$ also has a hereditarily connected colored graph.
\end{prop}
\begin{proof} We just need to show that if $(a,b)$ is an edge of $\bA$ with $a/\theta \ne b/\theta$, then $a/\theta$ is connected to $b/\theta$ within the subalgebra they generate. We will induct on the size of $|\Sg\{a,b\}|$. If $(a,b)$ is a semilattice edge, then $(a/\theta,b/\theta)$ will automatically be a semilattice edge as well. Otherwise, let $\eta$ be the maximal congruence on $\Sg\{a,b\}$ corresponding to the edge $(a,b)$.

Since every congruence class of $\eta$ is a proper subalgebra of $\Sg\{a,b\}$, we see by induction that if $c \equiv_\eta d$, then $c/\theta$ and $d/\theta$ are connected in the subalgebra they generate, which is contained in $\Sg\{a/\theta,b/\theta\}$. Thus if the restriction of $\theta$ to $\Sg\{a,b\}$ is not contained in the maximal congruence $\eta$, then $a/\theta$ must be connected to $b/\theta$ in the subalgebra $\Sg\{a/\theta,b/\theta\}$. Otherwise, if the restriction of $\theta$ to $\Sg\{a,b\}$ \emph{is} contained in $\eta$, then $(a/\theta,b/\theta)$ is an edge, with witnessing congruence $\eta/\theta$.
\end{proof}

\begin{cor} If $\bA$ is a finite idempotent algebra with a hereditarily connected colored graph, then $\bA$ is Taylor.
\end{cor}

Bulatov's main result in \cite{colored-graph-prelim} and \cite{colored-graph} is that the converse to the above corollary holds. Since Bulatov didn't have the theory of absorbing subalgebras available to him, he proved this by using tame congruence theory. We will give a different proof, using a pair of consecutive daisy chain terms (whose existence followed from the existence of a cyclic term), and the fact that abelian Taylor algebras are affine.

\begin{thm}[Bulatov \cite{colored-graph-prelim}, \cite{colored-graph}]\label{colored-graph} A finite idempotent algebra $\bA$ is Taylor if and only if it has a hereditarily connected colored graph.
\end{thm}

We will prove Theorem \ref{colored-graph} by induction on $|\bA|$. A minimal counterexample $\bA$ must be simple by Proposition \ref{hereditarily-connected-product}, and every proper subalgebra of a minimal counterexample must have a hereditarily connected colored graph.

\begin{defn} If $\bA$ is any algebra, we define the \emph{hypergraph of proper subalgebras} of $\bA$ to be the hypergraph with vertex set equal to the underlying set of $\bA$, and with a hyperedge $\bB$ for every proper subalgebra $\bB \le \bA$. We say that $\bA$ is \emph{disconnected} if the hypergraph of proper subalgebras of $\bA$ is not connected.

We define the connected component equivalence relation $\sim_\bA$ (or just $\sim$ if $\bA$ is clear from context) on $\bA$ by $a \sim_\bA b$ if $a$ is connected to $b$ by a sequence of proper subalgebras of $\bA$ (note that in general, $\sim_\bA$ will \emph{not} be a congruence).
\end{defn}

\begin{prop} If $\bA$ is a disconnected algebra, then for any $a \not\sim_\bA b$ we have $\Sg\{a,b\} = \bA$.
\end{prop}

\begin{prop}\label{colored-claim-0} Suppose that $\bA$ is finite, idempotent, simple, and disconnected. For any binary relation $\RR \le \bA \times \bA$ with $\pi_2(\RR) = \bA$, either $\RR$ is the graph of an automorphism of $\bA$ or there is some $a \in \bA$ such that $\{a\} \times \bA \subseteq \RR$.
\end{prop}
\begin{proof} If $\RR$ is not the graph of an automorphism, then the linking congruence must be nontrivial, hence full (since $\bA$ is simple). Thus there is some $a \in \bA$ such that $(a,b),(a,c) \in \bA$ for some pair of elements $b,c$ with $b \not\sim c$, and from $\Sg_{\bA}\{b,c\} = \bA$ and idempotence, we see that $\{a\}\times \bA \subseteq \RR$.
\end{proof}

\begin{prop}\label{colored-claim-1} Suppose that $\bA$ is finite, idempotent, simple, and disconnected, and that $a \not\sim_\bA b$ are such that neither $(a,b)$ nor $(b,a)$ are semilattice edges. Then the binary relation
\[
\bS_{ab} \coloneqq \Sg_{\bA^2}\left\{\begin{bmatrix} a\\ b\end{bmatrix}, \begin{bmatrix} b\\ a\end{bmatrix}\right\}
\]
is the graph of an automorphism of order two which interchanges $a$ and $b$.
\end{prop}
\begin{proof} Assume not. Then by Proposition \ref{colored-claim-0}, there is some $c \in \bA$ with $\{c\}\times \bA \subseteq \bS_{ab}$. Since $a,b$ are in different connected components of $\bA$, at least one of them is in a different component than $c$, so we may suppose that $a$ and $c$ are in different connected components of $\bA$ without loss of generality. Then since $(b,a), (b,c) \in \bS_{ab}$ and $b \in \Sg_{\bA}\{a,c\} = \bA$, we have $(b,b) \in \bS_{ab}$, so $(a,b)$ is a semilattice edge.
\end{proof}

\begin{defn} Define the equivalence relation $\sim^s_\bA$ by $a \sim^s_\bA b$ if $a$ can be connected to $b$ by a chain of proper subalgebras and semilattice edges.
\end{defn}

\begin{cor}\label{colored-claim-1-cor} If $\bA$ is finite, idempotent, and simple, and if $\sim^s_\bA$ is not the full equivalence relation $\bA \times \bA$, then $\Aut(\bA)$ acts transitively on $\bA$.
\end{cor}
\begin{proof} For any pair $a,b \in \bA$, either $a \not\sim^s_\bA b$, or there is some $c \in \bA$ such that $a \not\sim^s_\bA c$ and $c \not\sim^s_\bA b$.
\end{proof}

\begin{prop}\label{colored-claim-2} Suppose that $\bA$ is finite, idempotent, and simple. For any $a \not\sim^s_\bA b$, if the ternary relation
\[
\RR_{ab} \coloneqq \Sg_{\bA^3}\left\{\begin{bmatrix} b\\ a\\ a\end{bmatrix}, \begin{bmatrix} a\\ b\\ a\end{bmatrix}, \begin{bmatrix} a\\ a\\ b\end{bmatrix}\right\}
\]
contains $(a,a,a)$ or $(b,a,b)$, then $\{a,b\}$ is a majority edge.
\end{prop}
\begin{proof} Suppose that there is some ternary term $t$ witnessing the presence of one of those tuples in $\RR_{ab}$. By Proposition \ref{colored-claim-1}, we see that $\{a,b\}$ is closed under $t$, and the restriction of $t$ to $\{a,b\}$ is either a majority operation or a Pixley operation. Either way, the ternary term $t(x,t(x,y,z),z)$ acts like a majority operation on $\{a,b\}$.
\end{proof}

\begin{proof}[Proof of Theorem \ref{colored-graph}] We only need to show that if $\bA$ is a finite idempotent Taylor algebra, then $\bA$ has a connected colored graph. Suppose that $\bA$ is a counterexample of minimal size, and note that $\bA$ must be simple by Proposition \ref{hereditarily-connected-product}.

Our aim is to show that for any $a \not\sim^s_\bA b$ such that $\{a,b\}$ is not a majority edge, the ternary relation
\[
\RR_{ab} \coloneqq \Sg_{\bA^3}\left\{\begin{bmatrix} b\\ a\\ a\end{bmatrix}, \begin{bmatrix} a\\ b\\ a\end{bmatrix}, \begin{bmatrix} a\\ a\\ b\end{bmatrix}\right\}
\]
must satisfy the conditions of Proposition \ref{ternary-abelian}, which will imply that $\bA$ is affine. Note that $a \not\sim_\bA b$ implies that $\RR_{ab}$ is subdirect.

The first step to proving that $\RR_{ab}$ satisfies the conditions of Proposition \ref{ternary-abelian} is showing that $\pi_{12}(\RR_{ab}) = \bA \times \bA$. Since $\pi_{12}(\RR_{ab})$ contains $(a,a), (a,b)$, and $(b,a)$, we need to check that $(b,b) \in \pi_{12}(\RR_{ab})$. So it is natural to study the set of tuples in $\RR_{ab}$ such that two of the coordinates are equal.

The next claim is the main place where we will use the fact that $\bA$ is Taylor.

{\bf Claim 1:} If we define $\bD_{ab} \le \bA\times\bA$ to be the set of pairs $(c,d)$ such that $(c,d,c) \in \RR_{ab}$, then $\pi_1(\bD_{ab}) \cap \pi_2(\bD_{ab}) \ne \emptyset$.

{\bf Proof of Claim 1:} Let $p,q$ be consecutive daisy chain terms, i.e. $p$ and $q$ are ternary terms satisfying the identities
\begin{align*}
p(x,x,y) &\approx p(y,x,x),\\
q(x,x,y) &\approx q(y,x,x) \approx p(x,y,x).
\end{align*}
If we set $c = p(a,a,b), d = p(a,b,a) = q(a,a,b)$, and $e = q(a,b,a)$, then we have
\[
p\left(\begin{bmatrix} b & a & a\\ a & b & a\\ a & a & b\end{bmatrix}\right) = \begin{bmatrix} c\\ d\\ c\end{bmatrix}
\]
and
\[
q\left(\begin{bmatrix} b & a & a\\ a & b & a\\ a & a & b\end{bmatrix}\right) = \begin{bmatrix} d\\ e\\ d\end{bmatrix},
\]
so $(c,d), (d,e) \in \bD_{ab}$, and $d \in \pi_1(\bD_{ab}) \cap \pi_2(\bD_{ab})$.

{\bf Claim 2:} The binary relation $\bD_{ab}$ from Claim 1 has $\pi_1(\bD_{ab}) = \bA$.

{\bf Proof of Claim 2:} Suppose not. First consider the case where neither $\pi_1(\bD_{ab}), \pi_2(\bD_{ab})$ are equal to $\bA$. Then if $c \in \pi_1(\bD_{ab}) \cap \pi_2(\bD_{ab})$, we see that both $\{a,c\}$ and $\{b,c\}$ are contained in proper subalgebras of $\bA$, so $a \sim c \sim b$, contradicting the assumption $a \not\sim b$.

Suppose now that $\pi_1(\bD_{ab}) \ne \bA$ but $\pi_2(\bD_{ab}) = \bA$. Then by Proposition \ref{colored-claim-0}, there is some $c \in \bA$ such that $\{c\} \times \bA \subseteq \bD_{ab}$. By Corollary \ref{colored-claim-1-cor}, there is an automorphism $\sigma \in \Aut(\bA)$ with $\sigma(a) = c$, and from $(\sigma(a),\sigma(b)) \in \{c\}\times\bA \subseteq \bD_{ab}$, we see that $(\sigma(a),\sigma(b),\sigma(a)) \in \RR_{ab}$, so in fact $\sigma(\RR_{ab}) \subseteq \RR_{ab}$, and so $\sigma(\bD_{ab}) = \bD_{ab}$. Thus $(a,a) = \sigma^{-1}(c,c) \in \bD_{ab}$, so by Proposition \ref{colored-claim-2} this contradicts the assumption that $\{a,b\}$ is not a majority edge.

{\bf Claim 3:} We have $\pi_{1,2}(\RR_{ab}) = \bA\times \bA$.

{\bf Proof of Claim 3:} By Claim 2, there is some $c$ such that $(b,c) \in \bD_{ab}$. Thus $(a,a),(a,b),(b,a),(b,b) \in \pi_{1,2}(\RR_{ab})$, and these four elements generate $\bA^2$.

{\bf Claim 4:} For any $c$, the binary relation $\RR_{ab}^c \le \bA^2$ defined as the set of pairs $(d,e)$ such that $(c,d,e) \in \RR_{ab}$ is the graph of an automorphism of order two.

{\bf Proof of Claim 4:} Suppose not. Note that by Claim 3, the relation $\RR_{ab}^c$ is subdirect. Thus if $\RR_{ab}^c$ is not the graph of an automorphism, then by Proposition \ref{colored-claim-0} there is some $d$ such that $\{d\}\times \bA \subseteq \RR_{ab}^c$.

First suppose that $c = a$, so $(a,b) \in \RR_{ab}^a$. Then either $a \not\sim d$ or $b \not\sim d$. If $a \not\sim d$, then from $(a,b),(d,b) \in \RR_{ab}^a$ we see that $(b,b) \in \RR_{ab}^a$, contradicting Proposition \ref{colored-claim-2}. Similarly if $b \not\sim d$, then from $(a,b),(a,d) \in \RR_{ab}^a$ we see that $(a,a) \in \RR_{ab}^a$, contradicting Proposition \ref{colored-claim-2}.

Now suppose that $c \ne a$. There is some $e \ne a$ such that $d \not\sim e$ (if not, then $\sim$ has just two equivalence classes, which are interchanged by the automorphism $\bS_{ab}$, and which both have size $1$, reducing us to the case $|\bA| = 2$). Then $\bS_{de}$ is the graph of an automorphism of order two which interchanges $d$ and $e$, and from $(d,e), (e,d) \in \RR_{ab}^c$ we see $\bS_{de} \subseteq \RR_{ab}^c$. Then since $e \ne a$ there is some $f \ne d$ such that $(a,f) \in \bS_{de} \subseteq \RR_{ab}^c$. Then from $(a,d), (a,f) \in \RR_{ab}^c$ we see that $(c,d), (c,f) \in \RR_{ab}^a$, contradicting the fact that $\RR_{ab}^a$ is the graph of an automorphism of order two.

To finish the proof, note that Claim 4 shows that the relation $\RR_{ab}$ satisfies the assumptions of Proposition \ref{ternary-abelian}, so $\bA$ must be abelian. Thus we can apply Theorem \ref{taylor-abelian} to see that $\bA$ must be affine.
\end{proof}

\section{Conservative Taylor algebras}\label{s-conservative}

Bulatov's colored graph was originally inspired by the study of conservative Taylor algebras. These algebras are easy to classify, and they are a great toy case for testing conjectures about general Taylor algebras.

\begin{defn} A $k$-ary operation $f : A^k \rightarrow A$ is \emph{conservative} if for all $a_1, ..., a_k \in A$ we have
\[
f(a_1, ..., a_k) \in \{a_1, ..., a_k\}.
\]
An algebraic structure $\bA$ is called \emph{conservative} if every basic operation of $\bA$ is conservative.
\end{defn}

Note that conservative algebras are automatically idempotent.

\begin{prop} An algebraic structure $\bA$ is conservative if and only if every subset $S \subseteq \bA$ is actually a subalgebra of $\bA$.
\end{prop}

On the relational side, we define conservative relational structures as follows.

\begin{defn} A relational clone $\Gamma$ on a domain $A$ is called \emph{conservative} if every unary relation $U \subseteq A$ is an element of $\Gamma$, i.e. $\cP(A) \subseteq \Gamma$. A relational structure $\fA$ is called \emph{conservative} if every unary relation $U$ can be primitively positively defined using the basic relations of $\fA$.
\end{defn}

\begin{prop} If a relational structure $\fA$ and an algebraic structure $\bA$ are related by the $\Inv-\Pol$ Galois correspondence, then $\fA$ is conservative if and only if $\bA$ is conservative.
\end{prop}

If we are handed a relational structure, then the next result can be useful to decrease the amount of work needed to verify that it is conservative.

\begin{prop} A relational structure $\fA$ with finite underlying set $A$ is conservative if and only if, for every $a \in A$, the unary relation $A \setminus \{a\}$ is primitively positively definable from the basic relations of $\fA$.
\end{prop}

\begin{ex} A natural example of a conservative CSP template (on an infinite domain) is the \emph{list-coloring} problem for graphs: the domain $A$ is an infinite set, and the relations consist of the binary $\ne$ relation and the collection of all possible subsets $U \subseteq A$ as unary relations.
\end{ex}

\begin{ex} A conservative $2$-semilattice is called a \emph{tournament}. The rock-paper-scissors algebra is probably the most famous example of a tournament which is not totally ordered.
\end{ex}

\begin{ex} If an affine CSP is conservative, then the domain must have size two: the only conservative affine algebra is $\ZZ/2^{\aff}$, up to term equivalence.
\end{ex}

Sometimes we will want to use the following refinement of the concept of conservative algebras.

\begin{defn} We say that a relational clone $\Gamma$ is $k$-\emph{conservative} if every unary relation $U \subseteq A$ with size $|U| \le k$ is an element of $\Gamma$, and we define $k$-conservative clones, algebras, and relational structures similarly.
\end{defn}

\begin{ex} An algebra is idempotent iff it is $1$-conservative.
\end{ex}

\begin{ex} The $k$-\emph{list-coloring} problem for graphs corresponds to the relational structure with infinite domain $A$, and relations consisting of the binary $\ne$ relation and the collection of all possible subsets $U \subseteq A$ with $|U| \le k$. This problem is equivalent to 2-SAT for $k = 2$, and is NP-hard for $k \ge 3$.
\end{ex}

\begin{ex} The only $2$-conservative affine algebras are $(\ZZ/2^{\aff})^k$, up to term equivalence and isomorphism.
\end{ex}

For the sake of understanding CSPs, we would like to focus on minimal Taylor algebras. The next result shows that we can reduce the study of conservative Taylor algebras to conservative minimal Taylor algebras without losing anything essential.

\begin{prop} Every reduct of a conservative algebra is also conservative. In particular, every conservative Taylor clone contains a minimal Taylor clone which is also conservative.
\end{prop}

Since every minimal Taylor clone can be generated by a pair of ternary terms (for instance, we can take a pair of consecutive daisy chain terms), we only have to focus on understanding conservative Taylor algebras of size $3$.

\begin{prop}\label{3-conservative} A minimal Taylor algebra is conservative if and only if it is $3$-conservative.
\end{prop}

In fact, looking carefully at how a daisy chain term must act on a conservative Taylor algebra, we have the following simplification.

\begin{prop}\label{conservative-ternary} If $\bA$ is a $2$-conservative minimal Taylor algebra and $w_i$ is any daisy chain term for $\bA$, then we have
\[
w_i(x,x,y) \approx w_i(x,y,x) \approx w_i(y,x,x),
\]
so in fact every $w_i$ is a ternary weak near-unanimity operation. The binary function
\[
f(x,y) \coloneqq w_i(x,x,y)
\]
is independent of $i$, and completely determines the colored graph of $\bA$. In addition, the binary function
\[
s(x,y) \coloneqq f(x,f(y,x))
\]
is a partial semilattice term operation of $\bA$.
\end{prop}
\begin{proof} Note that every pair of distinct elements $a,b \in \bA$ must form an edge of the colored graph of $\bA$ if $\bA$ is a $2$-conservative Taylor algebra. By our analysis of daisy chain terms on the basic two-element minimal Taylor algebras, we see that:
\begin{itemize}
\item if $(a,b)$ is a semilattice edge, then $w_i(x,x,y) = w_i(x,y,x) = w_i(y,x,x) = x \vee y$ for $x,y \in \{a,b\}$,
\item if $\{a,b\}$ is a majority edge, then $w_i(x,x,y) = w_i(x,y,x) = w_i(y,x,x) = x$ for $x,y \in \{a,b\}$,
\item if $\{a,b\}$ is a $\ZZ/2^{\aff}$ edge, then $w_i(x,x,y) = w_i(x,y,x) = w_i(y,x,x) = y$ for $x,y \in \{a,b\}$.
\end{itemize}
Thus we can tell what sort of edge $\{a,b\}$ is (as well as how it is directed, in case it is a semilattice edge) by examining the restriction of $f(x,y)$ to the set $\{a,b\}$. The claim about $s(x,y)$ follows easily by considering each of the three possible types of edge individually.
\end{proof}

Thus, from now on we imagine that all conservative minimal Taylor algebras live in a variety $\cV$ having just one ternary basic operation $w$, which satisfies the weak near-unanimity identity
\[
w(x,x,y) \approx w(x,y,x) \approx w(y,x,x).
\]

\begin{prop}\label{prop-conservative-free} The free algebra $\cF_\cV(x,y)$ on two generators in the variety $\cV$ generated by conservative minimal Taylor algebras has size $6$: its elements are $x, y, f(x,y), f(y,x), s(x,y), s(y,x)$ (defined as in the previous proposition). The colored graph of the algebra $\cF_\cV(x,y)$ is as follows.
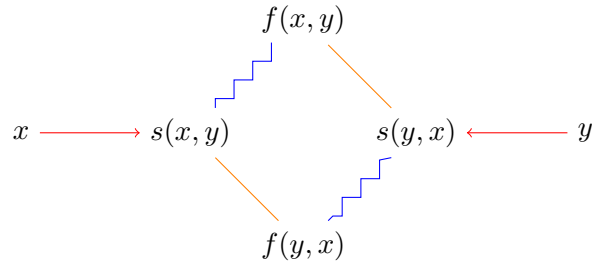
\begin{center}
\begin{tikzpicture}[scale=1.5,baseline=0.5cm]
  \node (x) at (-2.5,0) {$x$};
  \node (y) at (2.5,0) {$y$};
  \node (sx) at (-1,0) {$s(x,y)$};
  \node (sy) at (1,0) {$s(y,x)$};
  \node (fx) at (0,1) {$f(x,y)$};
  \node (fy) at (0,-1) {$f(y,x)$};
  \draw [->, color=red] (x) edge (sx) (y) edge (sy);
  \draw [color=orange] (sx) -- (fy) (sy) -- (fx);
  \draw [decorate, decoration={zigzag}, color=blue] (sx) to (fx) (sy) to (fy);
\end{tikzpicture}
\end{center}
Here the semilattice edges are directed, the majority edges are straight and undirected, and the $\ZZ/2^{\aff}$ edges are drawn as zigzags. This algebra has $\{f(x,y), f(y,x), s(x,y), s(y,x)\}$ as a binary absorbing subalgebra, corresponding to a semilattice quotient, and has $\{x,f(x,y),s(x,y)\}$ and $\{y,f(y,x),s(y,x)\}$ as ternary absorbing subalgebras, corresponding to a majority quotient.
\end{prop}

In order to understand conservative minimal Taylor algebras, Proposition \ref{3-conservative} implies that it's most important to understand the conservative algebras of size $3$. Additionally, Proposition \ref{conservative-ternary} implies that we just need to figure out which ternary weak near-unanimity operations on a three element set generate \emph{minimal} Taylor clones. We get a further simplification by dividing into cases based on whether or not there is a ternary cyclic term. In the case where there is no cyclic term, the following result is useful.

\begin{thm}\label{ternary-iteration} If $w$ is a ternary weak near-unanimity term operation of a finite algebra $\bA$, then there is a ternary weak near-unanimity term $g \in \Clo(w)$ which also satisfies the identity
\[
g(g(x,y,z),g(y,z,x),g(z,x,y)) \approx g(x,y,z).
\]
If $|\bA| = 3$ and $\bA$ has no ternary cyclic term, then any such $g$ satisfies $g(x,y,z) = x$ whenever $x,y,z$ are all different.
\end{thm}
\begin{proof} Let $\gamma : \bA^3 \rightarrow \bA^3$ be the map given by
\[
\gamma\left(\begin{bmatrix}x\\ y\\ z\end{bmatrix}\right) \coloneqq \begin{bmatrix}w(x,y,z)\\ w(y,z,x)\\ w(z,x,y)\end{bmatrix}.
\]
Then since $\bA^3$ is finite, there is some $k$ such that $\gamma^{\circ 2k} = \gamma^{\circ k}$. If we define $g$ by
\[
\gamma^{\circ k}\left(\begin{bmatrix}x\\ y\\ z\end{bmatrix}\right) = \begin{bmatrix}g(x,y,z)\\ g(y,z,x)\\ g(z,x,y)\end{bmatrix},
\]
then $\gamma^{\circ k} \circ \gamma^{\circ k} = \gamma^{\circ k}$ implies that $g$ satisfies the identity
\[
g(g(x,y,z),g(y,z,x),g(z,x,y)) \approx g(x,y,z).
\]
Note that since $w$ is a weak near-unanimity term, if any two of $x,y,z$ are equal, then $\gamma(x,y,z)$ is a constant tuple, and then by idempotence so is $\gamma^{\circ k}(x,y,z)$. Therefore $g$ is also a weak near-unanimity operation, and $\gamma^{\circ k}(x,y,z)$ can only avoid being a constant tuple if $x,y,z$ are all different.

If $\bA$ has no ternary cyclic term and has underlying set $\{a,b,c\}$, then by Theorem \ref{cyclic-terms} $\bA$ must have an automorphism of order three, so $\gamma(a,b,c)$ must be one of $(a,b,c), (b,c,a)$, or $(c,a,b)$, and in each case we have $\gamma^{\circ k}(a,b,c) = (a,b,c)$. Similarly, we must also have $\gamma^{\circ k}(a,c,b) = (a,c,b)$, so we have $g(x,y,z) = x$ whenever $x,y,z$ are all different.
\end{proof}

\begin{thm}\label{thm-minimal-no-cyclic} If a minimal Taylor algebra has size $3$ and has no ternary cyclic term, then (after renaming elements) it is term equivalent to one of the following four algebras:
\begin{itemize}
\item the affine algebra $\ZZ/3^{\aff}$,
\item the rock-paper-scissors algebra $\bA_{rps}$ from Section \ref{s-rps},
\item the three element dual discriminator algebra from Example \ref{ex-dual-discriminator}, or
\item the three element simple nonabelian Mal'cev algebra from Example \ref{ex-simple-nonab-malcev}.
\end{itemize}
All but the first are conservative, all but the second have a full automorphism group, and the first two have binary cyclic terms.
\end{thm}
\begin{proof} By Theorem \ref{cyclic-terms}, if a minimal Taylor algebra $\bA$ with underlying set $\{a,b,c\}$ has size $3$ and has no ternary cyclic term, then $\bA$ must have an automorphism of order three with no fixed points, so the permutation $(a\ b\ c)$ is in $\Aut(\bA)$. By Theorem \ref{colored-graph}, either $\bA$ is affine - in which case it must be term-equivalent to $\ZZ/3^{\aff}$ - or $\bA$ has some proper subalgebra of size $2$ (since $\bA$ has an edge $(a,b)$, and either $\Sg\{a,b\}$ has size $2$, or $\Sg\{a,b\} = \bA$ has a proper quotient, and one of the congruence classes is a subalgebra of size $2$). Since $(a\ b\ c) \in \Aut(\bA)$, if any $2$-element subset of $\bA$ is a subalgebra, then \emph{every} $2$-element subset of $\bA$ is a subalgebra, and all three $2$-element subalgebras of $\bA$ are isomorphic to $\{a,b\}$.

If $\{a,b\}$ is a semilattice, then any binary operation $s$ that acts like the semilattice term on $\{a,b\}$ has $(\{a,b,c\},s)$ isomorphic to the rock-paper-scissors algebra. If $\{a,b\}$ is a majority algebra and $g$ is a ternary weak near-unanimity operation as in the previous theorem, then $g$ is a majority operation which acts as first projection whenever all three of its inputs are distinct, so $(\{a,b,c\},g)$ is isomorphic to the three-element dual discriminator algebra. If $\{a,b\}$ is an affine algebra and $g$ is a ternary weak near-unanimity operation as in the previous theorem, then $g$ is a Mal'cev operation which acts as a minority operation whenever two of its inputs are equal, and which acts as first projection whenever all three of its inputs are distinct, so $(\{a,b,c\},g)$ is isomorphic to the three element simple nonabelian Mal'cev algebra from Example \ref{ex-simple-nonab-malcev}.
\end{proof}

In most of the remaining cases, the colored graph already does not have any automorphisms of order $3$. In these cases, it turns out to be relatively easy to pick out a specific ternary cyclic operation which is determined by the colored graph alone. In fact, we have the following slightly stronger statement.

\begin{thm} Suppose that a minimal Taylor algebra $\bA$ has the following properties:
\begin{itemize}
\item $\bA$ is $2$-conservative, that is, for all $a,b \in \bA$ the subset $\{a,b\}$ is a subalgebra of $\bA$,
\item the colored graph of $\bA$ does not contain any majority triangles, and
\item the colored graph of $\bA$ does not contain any affine triangles.
\end{itemize}
Then $\bA$ is conservative, and $\Clo(\bA)$ is determined by the colored graph of $\bA$.
\end{thm}
\begin{proof} Let $w$ be any daisy chain term for $\bA$. Define a map $\gamma : \bA^3 \rightarrow \bA^3$ as in Theorem \ref{ternary-iteration}. We will make sure to only apply $\gamma$ to triples where some pair of coordinates are equal, since the values $\gamma$ takes on such triples is completely determined by the colored graph of $\bA$ by Proposition \ref{conservative-ternary}. Define binary terms $f,s$ as in Proposition \ref{conservative-ternary}, and note that $f$ and $s$ are uniquely determined by the colored graph of $\bA$. Define maps $\alpha_f, \beta_f: \bA^3 \rightarrow \bA^3$ by
\[
\alpha_f\left(\begin{bmatrix}x\\ y\\ z\end{bmatrix}\right) \coloneqq \begin{bmatrix}f(x,y)\\ f(y,z)\\ f(z,x)\end{bmatrix}
\]
and
\[
\beta_f\left(\begin{bmatrix}x\\ y\\ z\end{bmatrix}\right) \coloneqq \begin{bmatrix}f(x,z)\\ f(y,x)\\ f(z,y)\end{bmatrix},
\]
and define maps $\alpha_s, \beta_s: \bA^3 \rightarrow \bA^3$ similarly, with $f$ replaced by $s$. Note that $\alpha_f, \beta_f$, etc. each have the property that if the input has two coordinates the same, then so does the output. As long as $a,b,c$ do not form a majority triangle, an affine triangle, or a rock-paper-scissors subalgebra, then the triple
\[
\alpha_f\circ\beta_f\circ\alpha_s\circ\beta_s\left(\begin{bmatrix}x\\ y\\ z\end{bmatrix}\right)
\]
has two of its three coordinates equal (to check this, consider the case where $\{a,b,c\}$ contains at least one semilattice edge separately from the case where it contains only majority and affine edges). Thus the ternary term
\[
t \coloneqq \pi_1\circ\gamma\circ\alpha_f\circ\beta_f\circ\alpha_s\circ\beta_s : \bA^3 \rightarrow \bA
\]
is cyclic on every such triple. Since we assumed that $\bA$ has no majority triangles or affine triangles, the only possible triples of $\bA$ such that the value of $t$ is not uniquely determined by the colored graph of $\bA$ are the rock-paper-scissors subsets of $\bA$, which are necessarily subalgebras of $\bA$ by Theorem \ref{minimal-taylor-subalg}. If we iterate $t$ as in Theorem \ref{ternary-iteration}, then the resulting ternary function $g$ has its values on rock-paper-scissors subalgebras fixed as well, so all of the values of $g$ are determined purely by the colored graph of $\bA$. Furthermore, this $g$ is conservative and generates a Taylor clone, so $\Clo(\bA) = \Clo(g)$ and $\bA$ is conservative.
\end{proof}

Finally, we need to understand the case of a majority triangle or affine triangle $\{a,b,c\}$ with a cyclic term. In these cases, it is helpful to keep track of the subalgebra
\[
\pi_1\left(\Sg_{\bA^2}\left\{\begin{bmatrix} a\\ b\end{bmatrix}, \begin{bmatrix} b\\ c\end{bmatrix}, \begin{bmatrix} c\\ a\end{bmatrix}\right\} \cap \Delta_\bA\right) \le \bA,
\]
since the set of possible outputs of a cyclic term applied to $(a,b,c)$ must be contained in this subalgebra. This subalgebra is an invariant of $\Clo(\bA)$, and it shrinks when $\Clo(\bA)$ shrinks.

\begin{defn} If $\bA$ is a three element minimal Taylor algebra with underlying set $\{a,b,c\}$, then we will say that an element $x$ of $\bA$ is \emph{circled} if $(x,x) \in \Sg\{(a,b),(b,c),(c,a)\}$. Note that the set of circled elements of $\bA$ does not depend on the ordering of $a,b,c$.
\end{defn}

\begin{thm} Suppose that $\bA$ is a conservative three element minimal Taylor algebra with a ternary cyclic term $g$, such that either all three of the edges of $\bA$ are majority or all three are affine. Then (after renaming elements) $\bA$ is term equivalent to one of the following three algebras:
\begin{itemize}
\item the three element solvable nonabelian Mal'cev algebra from Example \ref{ex-solvable-nonab-malcev}, with $*$ as the unique circled element,
\item the three element median algebra $\{0,1,2\}$, with the median element $1$ as the unique circled element, or
\item the three element minimal majority algebra $(\{a,b,c\}, m)$, with $m$ a cyclic majority operation such that $m(a,b,c) = b$ and $m(a,c,b) = c$, with $\{b,c\}$ as the set of circled elements.
\end{itemize}
In particular, every conservative three element minimal Taylor algebra is determined up to term equivalence by its colored graph and set of circled elements.

Furthermore, in any conservative minimal Taylor algebra, we can choose a ternary operation $g$ as in Theorem \ref{ternary-iteration} such that if we take $g$ as the basic operation, then every three element majority subalgebra with two circled elements is isomorphic to the third algebra listed above (not just term-equivalent).
\end{thm}
\begin{proof} Let $g$ be a ternary cyclic term for $\bA$, and suppose that $\bA$ has underlying set $\{a,b,c\}$. Once we know the types of the edges of $\bA$, we only need to know the values of $g(a,b,c)$ and $g(a,c,b)$ to completely determine $g$. For each choice of edges, we have two cases: either $g(a,b,c) = g(a,c,b)$, or $g(a,b,c) \ne g(a,c,b)$. This gives us four cases total.

First consider the case where all three edges of $\bA$ are affine (so $g$ is Mal'cev), and $g(a,b,c) \ne g(a,c,b)$. Without loss of generality, we may assume that $g(a,b,c) = b$ and $g(a,c,b) = c$. We will show that this case does not occur, by constructing a term $w$ which generates a strictly smaller Taylor clone. Note that the order two permutation which swaps $b$ and $c$ is an automorphism of $(\{a,b,c\},g)$. Then if we define the ternary operation $t$ by
\[
t(x,y,z) \coloneqq g(x,g(x,y,z),g(x,g(x,y,z),g(x,z,y))),
\]
then $t$ is also Mal'cev and satisfies
\[
t(a,b,c) = a, \;\; t(b,c,a) = c, \;\; t(c,a,b) = c,
\]
so if we define the ternary operation $w$ by
\[
w(x,y,z) \coloneqq g(t(x,y,z),t(y,z,x),t(z,x,y)),
\]
then $w$ is a symmetric Mal'cev operation, with $w(a,b,c) = w(a,c,b) = a$. Then $w$ generates a strictly smaller Taylor clone, since $w$ preserves the equivalence relation with equivalence classes $\{a\}$ and $\{b,c\}$, while $g$ does not. Thus this case does not occur.

In the remaining three cases, we get the three algebras described in the theorem statement. We need to check that these three algebras are really \emph{minimal} Taylor. Note that in each case, there is a nontrivial congruence on $\bA$ with quotient of size two and congruence classes of size at most two, so every Taylor reduct of $\bA$ is forced to have a cyclic ternary term. We will show that the clone of each of these algebras contains only one or two ternary cyclic operations $w$. Note that the only values of $w(x,y,z)$ which are not determined by the types of the edges are the ones where $x,y,z$ are all distinct.

In the case of the solvable Mal'cev algebra from Example \ref{ex-solvable-nonab-malcev} with underlying set $\{0,1,*\}$, the congruence with congruence classes $\{*\}, \{0,1\}$ forces the value of $w(0,1,*)$ to be $*$, and similarly for other permutations of the inputs. Thus there is only one ternary cyclic operation $w$ in the clone.

In the case of the three element median algebra $\{0,1,2\}$, the congruences corresponding to the partitions $\{0,1\},\{2\}$ and $\{0\},\{1,2\}$ force the value of $w(0,1,2)$ to be in $\{0,1\} \cap \{1,2\} = \{1\}$, and similarly for other permutations of the inputs. Thus there is only one ternary cyclic operation $w$ in the clone.

In the last case, the congruence corresponding to the partition $\{a\},\{b,c\}$ forces the value of $w(a,b,c)$ to be either $b$ or $c$. Additionally, the order two automorphism which interchanges $b$ and $c$ forces us to have
\[
w(a,b,c) = b \;\; \iff \;\; w(a,c,b) = c.
\]
Thus we either have $w(x,y,z) \approx m(x,y,z)$, or $w(x,y,z) \approx m(x,z,y)$, so there are exactly two ternary cyclic operations $w$ in the clone.

For the last statement, suppose that we have a minimal conservative algebra $\bA$, with several majority subalgebras with two circled elements. Let $g$ be any ternary operation as in Theorem \ref{ternary-iteration}. By the last case above, the restriction of $g$ to any of these majority subalgebras either acts like $m(x,y,z)$ or like $m(x,z,y)$. Suppose for contradiction that two of these subalgebras are not isomorphic, i.e., that $g$ acts as $m(x,y,z)$ on one and acts as $m(x,z,y)$ on the other. We will produce a ternary weak near-unanimity term $w$ which acts like $m(x,y,z)$ on both, which will generate a proper Taylor subclone. To this end, we define a ternary term $t$ by
\[
t(x,y,z) \coloneqq g(x,g(x,y,z),g(x,z,y)),
\]
and define $w$ by
\[
w(x,y,z) \coloneqq g(t(x,y,z),t(y,z,x),t(z,x,y)).
\]
Then $w$ is cyclic on any subalgebra of $\bA$ where $g$ is cyclic, so in particular $w$ is a weak near-unanimity operation. Note that if $\{a,b,c\}$ is a majority subalgebra of $\bA$ with $\{b,c\}$ as the set of circled elements, then regardless of whether the restriction of $g$ to $\{a,b,c\}$ is $m(x,y,z)$ or $m(x,z,y)$, we always have
\[
t(a,b,c) = b, \;\; t(b,c,a) = b, \;\; t(c,a,b) = c,
\]
so $w(a,b,c) = b$, and so the restriction of $w$ to $\{a,b,c\}$ is $m(x,y,z)$.
\end{proof}

Putting the proofs of the above theorems together, we get a procedure which puts the basic ternary weak near-unanimity operation $g$ of any minimal conservative Taylor algebra into a standard form, such that the restriction of $g$ to any three element subalgebra is completely determined by the edge types and the set of circled elements. In particular, we can exactly count the number of conservative minimal Taylor clones of a given size.

\begin{cor} The number of conservative minimal Taylor clones on a set of size $n$ is exactly
\[
\sum_{3\text{-edge-colorings of }K_n} 2^{\#(\text{\textcolor{red}{semilattice}})}4^{\Delta(\text{\textcolor{blue}{affine}})}7^{\Delta(\text{\textcolor{orange}{majority}})} = (1+o(1))\cdot 7^{\binom{n}{3}},
\]
where $\Delta(c)$ is the number of monochromatic triangles of color $c$. In particular, for large $n$ almost all conservative minimal Taylor clones are clones of majority algebras.
\end{cor}

If we only want to know the number of conservative minimal Taylor algebras of a given size up to term equivalence and \emph{isomorphism}, then we can use the Burnside counting theorem, together with the fact that the automorphism group of a conservative minimal Taylor algebra is determined by its colored graph and the choices of circled vertices on its three-element majority and Mal'cev subalgebras in the obvious way. The number of conservative minimal Taylor clones on domains of sizes $2,3,4,5$ is listed below (note: these were computed by hand, so there might be some mistakes).
\begin{center}
\begin{tabular}{c|c|c}
Domain size & \# up to term equiv. & \# up to term equiv. and iso. \\
\hline
2 & 4 & 3\\
3 & 73 & 19\\
4 & 9829 & 520\\
5 & 320668024 & 2686891
\end{tabular}
\end{center}

\subsection{Classification of three-element minimal Taylor algebras}\label{ss-minimal-three}

As it turns out, up to term equivalence and isomorphism there are just $24$ minimal Taylor algebras on a set of size $3$. Of these, $19$ are conservative, and the remaining $5$ are easy to describe. One of the most obvious non-conservative minimal Taylor algebras of size $3$ is the affine algebra $\ZZ/3^{\aff}$. Three more are subdirect products of two-element minimal Taylor algebras: specifically, the free semilattice on two generators (which is a subdirect product of two two-element semilattices), the subdirect product of a two-element semilattice and a two-element majority algebra, and the subdirect product of a two-element semilattice and $\ZZ/2^{\aff}$ (all three of these are quotients of the free algebra from Proposition \ref{prop-conservative-free}). There is no three-element subdirect product of a two-element majority algebra and $\ZZ/2^{\aff}$, but the final example is nearly this: it is the algebra from Example \ref{ex-few-subpowers}, which has a $3$-edge term, a $\ZZ/2^{\aff}$ quotient, and a two-element centrally absorbing algebra. In this subsection we will prove that this is the complete list of minimal Taylor algebras on a three-element set.

By Theorem \ref{thm-minimal-no-cyclic}, we only have to classify the minimal Taylor algebras of size $3$ which have a ternary cyclic term $g$, and since we have already classified the conservative ones, we just need to classify those which are generated by two elements. By Bulatov's Theorem \ref{colored-graph}, each of these algebras has a connected colored graph. Our first step will be to show that for minimal Taylor algebras of size $3$, every edge of Bulatov's colored graph is a two-element subalgebra.

\begin{prop} If $\bA$ is a minimal Taylor algebra of size $3$ other than $\ZZ/3^{\aff}$, then the graph on $\bA$ with edges given by the two-element subalgebras of $\bA$ is connected - in other words, $\bA$ has at least two subalgebras of size two.
\end{prop}
\begin{proof} Suppose for contradiction that the graph of two-element subalgebras of $\bA$ is disconnected, and suppose the underlying set of $\bA$ is $\{a,b,c\}$. By Theorem \ref{colored-graph}, Bulatov's colored graph of $\bA$ must be connected, so if $\bA$ is not $\ZZ/3^{\aff}$ then there must be some nontrivial congruence $\theta \in \Con(\bA)$ such that $\bA/\theta$ is either $\ZZ/2^{\aff}$ or the two-element majority algebra (it can't be the two-element semilattice by Proposition \ref{prop-weak-semilattice-edge}). Suppose without loss of generality that $\{b,c\}$ is the congruence class of $\theta$ which has size $2$, so that $\{b,c\}$ is a two-element subalgebra of $\bA$. By our assumption that the graph of two-element subalgebras is disconnected, neither $\{a,b\}$ nor $\{a,c\}$ can be a subalgebra of $\bA$. Let $g$ be a ternary cyclic operation on $\bA$, which exists by Theorem \ref{thm-minimal-no-cyclic} (or directly from Theorem \ref{cyclic-terms}).

Suppose first that $\bA/\theta$ is $\ZZ/2^{\aff}$. Then since any ternary cyclic term operation of $\ZZ/2^{\aff}$ acts as the minority operation, we must have
\[
g(a,b,b) = g(a,c,c) = a.
\]
Since $\{a,b\}$ and $\{a,c\}$ are not closed under $g$, we must then have
\[
g(a,a,b) = c, \;\;\; g(a,a,c) = b.
\]
Define a ternary term $t$ by
\[
t(x,y,z) = g(g(x,x,y), g(y,y,z), g(z,z,x)).
\]
Then $t$ is also cyclic, and we have
\[
t(a,a,b) = g(a,c,a) = b,
\]
so $\{a,b\}$ is closed under $t$, and similarly $\{a,c\}$ is also closed under $t$. Thus $t$ generates a strictly smaller Taylor clone, contradicting the assumption that $\bA$ is minimal Taylor.

Now suppose that $\bA/\theta$ is the two-element majority algebra. Then since any ternary cyclic term acts on a majority algebra as a majority operation, we must have
\[
g(a,a,b) = g(a,a,c) = a,
\]
and
\[
g(a,b,c), g(a,c,b) \in \{b,c\}.
\]
Since $\{a,b\}$ and $\{a,c\}$ are not closed under $g$, we must have
\[
g(a,b,b) = c, \;\;\; g(a,c,c) = b.
\]
Assume without loss of generality that $g(a,b,c) = b$, otherwise swap $b,c$ and reorder the last two arguments to $g$. Define a ternary term $t$ by
\[
t(x,y,z) = g(g(x,x,y), g(y,y,z), g(z,z,x)).
\]
Then $t$ is also cyclic, and we have
\[
t(a,b,b) = g(a,b,c) = b,
\]
so $\{a,b\}$ is closed under $t$. As before, this contradicts the assumption that $\bA$ is minimal Taylor.
\end{proof}

From here on, we just need to classify the minimal Taylor algebras on the set $\{a,b,c\}$ such that $a$ and $b$ generate the algebra, while $\{a,c\}$ and $\{b,c\}$ are two-element subalgebras. First we handle the case where one of these subalgebras is a two-element semilattice.

\begin{prop} If $\bA$ is a minimal Taylor algebra with underlying set $\{a,b,c\}$ such that $\bA$ is generated by $a$ and $b$, then $\bA$ does not have a semilattice subalgebra of the form $c \rightarrow a$.
\end{prop}
\begin{proof} Suppose for contradiction that $\bA$ is generated by $a$ and $b$, but that $c \rightarrow a$. Let $g$ be a ternary cyclic term operation of $\bA$, which exists by Theorem \ref{thm-minimal-no-cyclic}.

First suppose that we do not also have $c \rightarrow b$. We will initially attempt to prove that $(a,a) \in \Sg_{\bA^2}\{(a,b),(b,a)\}$, so that by Proposition \ref{prop-minimal-prepared} $\{a,b\}$ will be a two-element semilattice subalgebra with $b \rightarrow a$, which will contradict the assumption that $a$ and $b$ generate $\bA$. Let $s \in \Clo(g)$ be a partial semilattice term operation of $\bA$ with $s(a,c) = s(c,a) = a$. Since by assumption we do not have $a \rightarrow b$, we have $s(a,b) = a$ as well, and since we assumed that $\{b,c\}$ is a subalgebra and that we do not have $c \rightarrow b$, we have $s(c,b) = c$. Define $\bS$ by
\[
\bS = \Sg_{\bA^2}\Big\{\begin{bmatrix}a\\ b\end{bmatrix},\begin{bmatrix}b\\ a\end{bmatrix}\Big\}.
\]
If we have $(a,c) \in \bS$, then by symmetry we also have $(c,a) \in \bS$, so
\[
s\left(\begin{bmatrix}a\\ c\end{bmatrix},\begin{bmatrix}c\\ a\end{bmatrix}\right) = \begin{bmatrix}a\\ a\end{bmatrix}
\]
is in $\bS$ as well. If we have $(c,c) \in \bS$, then
\[
s\left(\begin{bmatrix}c\\ c\end{bmatrix},\begin{bmatrix}a\\ b\end{bmatrix}\right) = \begin{bmatrix}a\\ c\end{bmatrix},
\]
so $(a,c) \in \bS$, and then we have $(a,a) \in \bS$ as before. Since $c \in \Sg_\bA\{a,b\} = \pi_2(\bS)$, if neither $(a,c)$ nor $(c,c)$ are in $\bS$, then we must have $(b,c) \in \bS$. So the only way to avoid a contradiction in this case is for $\bS$ to be given by
\[
\bS = \Big\{\begin{bmatrix}a\\ b\end{bmatrix},\begin{bmatrix}b\\ a\end{bmatrix},\begin{bmatrix}b\\ c\end{bmatrix},\begin{bmatrix}c\\ b\end{bmatrix}\Big\}.
\]
Then the linking congruence $\theta$ of $\bS$ corresponds to the partition $\{a,c\}, \{b\}$, and $\bA/\theta$ is isomorphic to the two-element subalgebra $\{b,c\}$. Since $(b,b), (c,c) \not\in \bS$, $\bA/\theta \cong \{b,c\}$ is either a majority algebra or is $\ZZ/2^{\aff}$. Define a ternary term $t$ by
\[
t(x,y,z) = g(s(g(x,y,z),x), s(g(x,y,z),y), s(g(x,y,z),z)).
\]
Then $t$ is also cyclic, and we will show that $\{a,b\}$ is closed under $t$, contradicting the assumption that $\bA$ is minimal Taylor. If $\bA/\theta$ is $\ZZ/2^{\aff}$, then we have
\[
t(a,a,b) = g(s(b,a), s(b,a), s(b,b)) = g(b,b,b) = b,
\]
and since $\{a,b\}$ is not closed under $g$ we have
\[
g(a,b,b) = c,
\]
so
\[
t(a,b,b) = g(s(c,a), s(c,b), s(c,b)) = g(a, c, c) = a.
\]
If $\bA/\theta$ is a majority algebra, then we have
\[
t(a,b,b) = g(s(b,a), s(b,b), s(b,b)) = g(b,b,b) = b,
\]
and since $\{a,b\}$ is not closed under $g$ we have
\[
g(a,a,b) = c,
\]
so
\[
t(a,a,b) = g(s(c,a), s(c,a), s(c,b)) = g(a,a,c) = a.
\]
Either way, $\{a,b\}$ is closed under $t$, which gives us a contradiction.

Now suppose that we have both $c \rightarrow a$ and $c \rightarrow b$. Consider the binary relation $\bS = \Sg_{\bA^2}\{(a,b), (b,a)\}$ as we did before - if either $(a,c)$ or $(b,c)$ are in $\bS$, then we easily see that one of $(a,a), (b,b)$ is in $\bS$, which together with Proposition \ref{prop-minimal-prepared} would contradict the assumption that $\bA$ is generated by $a$ and $b$. Thus the only way to immediately avoid a contradiction is for $\bS$ to be given by
\[
\bS = \Big\{\begin{bmatrix}a\\ b\end{bmatrix},\begin{bmatrix}b\\ a\end{bmatrix},\begin{bmatrix}c\\ c\end{bmatrix}\Big\}.
\]
In particular, there must be some binary term $f \in \Clo(g)$ such that $f(a,b) = f(b,a) = c$, and since $\bA$ is prepared by Proposition \ref{prop-minimal-prepared}, this $f$ must be the commutative binary operation described in the following table.
\begin{center}
\begin{tabular}{c|ccc} $f$ & $a$ & $b$ & $c$\\ \hline $a$ & $a$ & $c$ & $a$\\ $b$ & $c$ & $b$ & $b$\\ $c$ & $a$ & $b$ & $c$\end{tabular}
\end{center}
Then $(\{a,b,c\},f)$ is isomorphic to the algebra $(\{-,0,+\},s_2)$ from Example \ref{lp-not-width-1}, with the isomorphism given by $a \mapsto +, b \mapsto -, c \mapsto 0$. This algebra is not minimal Taylor: $\Clo(s_2)$ properly contains the clone of the conservative bounded width algebra $(\{-,0,+\}, g)$ described in Example \ref{ex-slippery}. Explicitly, consider the ternary term $t$ on $\bA$ given by
\[
t(x,y,z) = f(f(f(x,y),f(y,z)), f(x,z)).
\]
It is easy to check that this $t$ is symmetric, and that $\{a,b\}$ is closed under $t$, contradicting the assumption that $\bA$ is minimal Taylor. In fact, even this operation $t$ does not generate a minimal Taylor clone (this claim is left as an exercise).
\end{proof}

\begin{prop} If $\bA$ is a minimal Taylor algebra with underlying set $\{a,b,c\}$ such that $\bA$ is generated by $a$ and $b$, and if $\bA$ has a semilattice subalgebra of the form $a \rightarrow c$, then $\bA$ is isomorphic to a subdirect product of its two-element subalgebras $\{a,c\}$ and $\{b,c\}$.
\end{prop}
\begin{proof} By assumption, we do not have $b \rightarrow a$ or $c \rightarrow a$, so $s(\{b,c\},\bA) \subseteq \{b,c\}$ for every partial semilattice term $s \in \Clo(\bA)$. Since $a \rightarrow c$ and $c \in \Sg_\bA\{a,b\}$, we can apply Proposition \ref{bin-central-criterion} to see that $\{b,c\}$ is a binary absorbing subalgebra of $\bA$. Then by Proposition \ref{prop-minimal-binary-strong}, we see that there is a congruence $\theta$ on $\bA$ corresponding to the partition $\{a\}, \{b,c\}$, such that $\bA/\theta \cong \{a,c\}$ is a two-element semilattice. To finish, we just need to show that the equivalence relation $\psi$ on $\bA$ corresponding to the partition $\{a,c\}, \{b\}$ is also a congruence of $\bA$.

If $b \rightarrow c$, then the same argument as in the last paragraph shows that $\{a,c\}$ is a binary absorbing subalgebra of $\bA$ and that $\psi$ is a congruence of $\bA$ - in this case, $\bA$ is the free semilattice on two generators. Additionally, the previous proposition shows that we can't have $c \rightarrow b$. Thus the only remaining cases are the case where $\{b,c\}$ is a majority algebra and the case where $\{b,c\}$ is $\ZZ/2^{\aff}$, and we assume from here on that we are in one of these two cases.

Consider the binary relation $\bS = \Sg_{\bA^2}\{(a,b),(b,a)\}$. Letting $s$ be a partial semilattice operation of $\bA$ with $s(a,c) = c$, we must have $s(a,b) \equiv s(a,c) = c \pmod{\theta}$, so since we do not have $a \rightarrow b$ we must have $s(a,b) = c$. Then by our assumption that we do not have $b \rightarrow c$, we have
\[
s\left(\begin{bmatrix} a\\ b\end{bmatrix}, \begin{bmatrix} b\\ a\end{bmatrix}\right) = \begin{bmatrix} c\\ b\end{bmatrix},
\]
so
\[
\bS \supseteq \Big\{\begin{bmatrix} a\\ b\end{bmatrix},\begin{bmatrix} b\\ a\end{bmatrix},\begin{bmatrix} b\\ c\end{bmatrix},\begin{bmatrix} c\\ b\end{bmatrix}\Big\}.
\]
If this containment is an equality, then the linking congruence of $\bS$ is $\psi$, which would prove that $\psi$ is a congruence of $\bA$ and finish the proof. Otherwise, since $(a,a), (b,b)$ are not contained in $\bS$ by Proposition \ref{prop-minimal-prepared}, at least one of $(a,c), (c,a)$ or $(c,c)$ must be an element of $\bS$. Since $a \rightarrow c$ and $\bS$ is symmetric, in each case we see that $(c,c) \in \bS$. From here on we assume that $(c,c) \in \bS$.

From $(b,c), (c,c), (c,b) \in \bS$ but $(b,b) \not\in \bS$, we see that $\{b,c\}$ can't be $\ZZ/2^{\aff}$ (since the parallelogram property fails for a binary relation on $\{b,c\}$). We are left with the case where $\{b,c\}$ is a majority algebra. Let $g$ be a cyclic ternary term on $\bA$, which exists by Theorem \ref{thm-minimal-no-cyclic}, and let $f$ be a binary term with $f(a,b) = f(b,a) = c$, which must exist if $(c,c) \in \bS$. Define a ternary term $t$ by
\[
t(x,y,z) = g(f(x,y), f(y,z), f(z,x)).
\]
Then $t$ is also cyclic, and we have
\[
t(a,a,b) = g(a,c,c) = c, \;\;\; t(a,b,b) = g(c,b,c) = c, \;\;\; t(a,b,c) \in g(c,\{b,c\},c) = \{c\}.
\]
This completely determines $t$, and shows that $\{a,c\}$ is a ternary absorbing subalgebra of $\bA$ with absorbing operation $t$. Therefore $\{a,c\}$ is a centrally absorbing subalgebra of $\bA$ by Theorem \ref{thm-minimal-central}. If we then define the cyclic ternary term $u$ by
\[
u(x,y,z) = t(s(x,y), s(y,z), s(z,x)),
\]
then it is easy to check that $\{b\}$ is a ternary absorbing subalgebra of $\bA$ with absorbing operation $u$, so by Theorem \ref{thm-minimal-central} and Corollary \ref{cor-minimal-central-majority} the equivalence relation $\psi$ is a congruence of $\bA$ (note that this actually contradicts the assumption $(c,c) \in \bS$).
\end{proof}

So far we have classified every minimal Taylor algebra of size $3$ which contains at least one semilattice subalgebra. The remaining cases are the cases where our algebra $\bA$ has two two-element subalgebras, each of which is either a majority algebra or a copy of $\ZZ/2^{\aff}$.

\begin{prop} If $\bA$ is a minimal Taylor algebra with underlying set $\{a,b,c\}$ which is generated by $a$ and $b$, then at least one of $\{a,c\}, \{b,c\}$ is not a $\ZZ/2^{\aff}$-subalgebra of $\bA$.
\end{prop}
\begin{proof} Suppose for contradiction that $\{a,c\},\{b,c\}$ are both $\ZZ/2^{\aff}$-subalgebras of $\bA$. By Corollary \ref{not-sd-complete}, $\bA$ either has a proper absorbing subalgebra or an affine quotient. If $\bA$ has a proper absorbing subalgebra, then one of $\{a,c\}, \{b,c\}$ has a proper absorbing subalgebra, which is impossible if they are both copies of $\ZZ/2^{\aff}$. Therefore there is a congruence $\theta$ of $\bA$ such that $\bA/\theta$ is affine, and we can assume without loss of generality that $\theta$ corresponds to the partition $\{a\}, \{b,c\}$ of $\bA$.

Let $g$ be a cyclic ternary term on $\bA$, which exists by Theorem \ref{thm-minimal-no-cyclic}. Then since $\bA/\theta \cong \{a,c\}$ is a copy of $\ZZ/2^{\aff}$ we must have
\[
g(a,\{b,c\},\{b,c\}) \in \{a\}, \;\;\; g(a,a,c) = c,
\]
and since $\{a,b\}$ is not closed under $g$ we must have
\[
g(a,a,b) = c
\]
as well. This, together with the fact that $\{b,c\}$ is also a copy of $\ZZ/2^{\aff}$, completely determines $g$. Define a ternary term $t$ by
\[
t(x,y,z) = g(x,g(x,y,y),g(x,y,z)),
\]
and note that
\[
t\left(\begin{bmatrix}a\\a\\b\end{bmatrix},\begin{bmatrix}a\\b\\a\end{bmatrix},\begin{bmatrix}b\\a\\a\end{bmatrix}\right) = g\left(\begin{bmatrix}a\\a\\b\end{bmatrix},\begin{bmatrix}a\\a\\c\end{bmatrix},\begin{bmatrix}c\\c\\c\end{bmatrix}\right) = \begin{bmatrix}c\\c\\b\end{bmatrix}.
\]
Thus the ternary term $u$ defined by
\[
u(x,y,z) = g(t(x,y,z), t(y,z,x), t(z,x,y))
\]
is a cyclic term with
\[
u(a,a,b) = g(c,c,b) = b,
\]
so $\{a,b\}$ is closed under $u$, which contradicts the assumption that $\bA$ is a minimal Taylor algebra.
\end{proof}

\begin{prop} If $\bA$ is a minimal Taylor algebra with underlying set $\{a,b,c\}$ which is generated by $a$ and $b$, then at least one of $\{a,c\}, \{b,c\}$ is not a majority subalgebra of $\bA$.
\end{prop}
\begin{proof} Suppose for contradiction that $\{a,c\},\{b,c\}$ are both majority subalgebras of $\bA$. Then no subquotient of $\bA$ can be affine, so by Theorem \ref{affine-free-pq} $\CSP(\bA)$ has bounded width, and so by Theorem \ref{bounded-width-term} $\bA$ has a binary term $f$ and a ternary term $g$ satisfying the identities
\[
g(x,x,y) \approx g(x,y,x) \approx g(y,x,x) \approx f(x,y) \approx f(f(x,y),f(y,x)).
\]
We may assume without loss of generality that this $g$ is also cyclic, by Theorem \ref{thm-minimal-no-cyclic} and the cyclic composition trick. Since $\{a,b\}$ is not closed under $g$, we may assume without loss of generality that $f(b,a) = c$.

First we show that $f(a,b) \ne b$. If $f(a,b) = b$, then we have $f(a,f(a,b)) = f(a,b) = b$ and
\[
f(b,f(b,a)) = f(f(a,b),f(b,a)) = f(a,b) = b,
\]
so $(b,b) \in \Sg_{\bA^2}\{(a,b),(b,a)\}$, and then by Proposition \ref{prop-minimal-prepared} we have $a \rightarrow b$, which contradicts the assumption that $\bA$ is generated by $a$ and $b$.

Now suppose that $f(a,b) = f(b,a) = c$. Let $t$ be the ternary term defined by
\[
t(x,y,z) = g(f(x,f(y,z)),f(y,f(z,x)),f(z,f(x,y))).
\]
Then $t$ is cyclic, and we have
\[
t(a,a,b) = g(f(a,c),f(a,c),f(b,a)) = g(a,a,c) = a
\]
and similarly $t(a,b,b) = b$, so $\{a,b\}$ is closed under $t$, contradicting the assumption that $\bA$ is a minimal Taylor algebra.

Finally, suppose that $f(a,b) = a, f(b,a) = c$. Let $h$ be the ternary term defined by
\[
h(x,y,z) = g(f(x,y),f(y,x),g(x,y,z)),
\]
let $i$ be the ternary term defined by
\[
i(x,y,z) = g(x,h(x,y,z),h(x,z,y)),
\]
and let $t$ be the cyclic ternary term defined by
\[
t(x,y,z) = g(i(x,y,z),i(y,z,x),i(z,x,y)).
\]
Then we have $h(a,a,b) = h(a,b,a) = h(b,a,a) = a$, so
\[
t(a,a,b) = g(g(a,a,a),g(a,a,a),g(b,a,a)) = g(a,a,a) = a,
\]
and $h(a,b,b) = h(b,a,b) = c, h(b,b,a) = b$, so
\[
t(a,b,b) = g(g(a,c,c),g(b,b,c),g(b,c,b)) = g(c,b,b) = b.
\]
Thus $\{a,b\}$ is closed under $t$, contradicting the assumption that $\bA$ is a minimal Taylor algebra.
\end{proof}

\begin{prop} Suppose that $\bA$ is a minimal Taylor algebra with underlying set $\{a,b,c\}$ which is generated by $a$ and $b$, that $\{a,c\}$ is a $\ZZ/2^{\aff}$-subalgebra, and that $\{b,c\}$ is a majority subalgebra of $\bA$. Then $\bA$ is the same as the algebra described in Example \ref{ex-few-subpowers}, up to isomorphism and term equivalence. In particular, $\{a,c\}$ is a centrally absorbing subalgebra of $\bA$, and $\bA$ has a congruence $\theta$ corresponding to the partition $\{a\},\{b,c\}$ such that $\bA/\theta \cong \{a,c\}$ is $\ZZ/2^{\aff}$.
\end{prop}
\begin{proof} Let $\bS = \Sg_{\bA^2}\{(a,b),(b,a)\}$, and let $g$ be a ternary cyclic term operation of $\bA$, which exists by Theorem \ref{thm-minimal-no-cyclic}. By Proposition \ref{prop-minimal-prepared}, we have $(a,a),(b,b) \not\in \bS$. If $(c,c) \in \bS$, then we have
\[
\begin{bmatrix}a\\c\end{bmatrix} = g\left(\begin{bmatrix}a\\b\end{bmatrix},\begin{bmatrix}c\\c\end{bmatrix},\begin{bmatrix}c\\c\end{bmatrix}\right) \in \bS,
\]
and then that
\[
\begin{bmatrix}a\\a\end{bmatrix} = g\left(\begin{bmatrix}a\\c\end{bmatrix},\begin{bmatrix}c\\a\end{bmatrix},\begin{bmatrix}c\\c\end{bmatrix}\right) \in \bS,
\]
which is a contradiction. Thus we have $(c,c) \not\in \bS$, so $\bS \cap \Delta_\bA = \emptyset$. If both $(a,c), (b,c) \in \bS$, then we have
\[
g\left(\begin{bmatrix}a\\b\end{bmatrix},\begin{bmatrix}b\\c\end{bmatrix},\begin{bmatrix}c\\a\end{bmatrix}\right) \in \bS \cap \Delta_\bA,
\]
which we just showed is impossible. Since $\pi_2(\bS) = \Sg_\bA\{a,b\} = \bA$, exactly one of $(a,c), (b,c)$ is in $\bS$.

Suppose first that $(b,c) \in \bS$. Then the linking congruence $\theta$ of $\bS$ corresponds to the partition $\{a,c\},\{b\}$ of $\bA$, and $\bA/\theta \cong \{b,c\}$ is a majority algebra. This implies that
\[
g(a,b,b) = g(b,b,c) = b, \; g(b,c,c) = c, \; g(a,c,c) = a.
\]
Since $\{a,b\}$ isn't closed under $g$, we must have
\[
g(a,a,b) = c.
\]
Let $f(x,y) = g(x,x,y)$, and define a ternary term $t$ by
\[
t(x,y,z) = g(f(g(x,y,z),f(x,y)),f(g(x,y,z),f(y,z)),f(g(x,y,z),f(z,x))).
\]
Then $t$ is cyclic,
\[
t(a,a,b) = g(f(c,a),f(c,c),f(c,b)) = g(a,c,c) = a,
\]
and
\[
t(a,b,b) = g(f(b,c),f(b,b),f(b,b)) = g(b,b,b) = b.
\]
Thus $\{a,b\}$ is closed under $t$, contradicting the assumption that $\bA$ is a minimal Taylor algebra.

Now suppose that $(a,c) \in \bS$. Then the linking congruence $\theta$ of $\bS$ corresponds to the partition $\{a\},\{b,c\}$ of $\bA$, and $\bA/\theta \cong \{a,c\}$ is $\ZZ/2^{\aff}$. This implies that
\[
g(a,\{b,c\},\{b,c\}) = \{a\}, \;\; g(a,a,c) = c, \; g(b,b,c) = b, \; g(b,c,c) = c.
\]
Since $\{a,b\}$ is not closed under $g$, we must have
\[
g(a,a,b) = c.
\]
This completely determines $g$, and we see that $\{a,c\}$ is a ternary absorbing subalgebra of $\bA$ with absorbing operation $g$, so by Theorem \ref{thm-minimal-central} $\{a,c\}$ is a centrally absorbing subalgebra of $\bA$. After swapping $a$ and $c$, $g$ is exactly the same operation as the one described in Example \ref{ex-few-subpowers}. To see that $\bA$ is really minimal Taylor in this case, note that any ternary cyclic $g' \in \Clo(g)$ must also satisfy all of the above identities, including $g'(a,a,b) = c$ since $\{a,c\}$ is centrally absorbing.
\end{proof}

Putting all of the pieces together, we have completed the classification of minimal Taylor algebras on a three-element domain.

\begin{thm} If $\bA$ is a minimal Taylor algebra on a set of size $3$, then up to term equivalence and isomorphism $\bA$ is one of the following $24$ algebras:
\begin{itemize}
\item one of the $19$ conservative minimal Taylor algebras classified in the previous section,
\item the affine algebra $\ZZ/3^{\aff} = (\ZZ/3, x-y+z)$,
\item the free semilattice on two generators,
\item the three-element subdirect product of $(\{0,1\},x\vee y\vee z)$ with $(\{0,1\}, \maj(x,y,z))$,
\item the three-element subdirect product of $(\{0,1\},x\vee y\vee x)$ with $\ZZ/2^{\aff} = (\ZZ/2, x+y+z)$,
\item the three-element algebra from Example \ref{ex-few-subpowers}, which has a $3$-edge term, a two-element centrally absorbing subalgebra, and a $\ZZ/2^{\aff}$ quotient.
\end{itemize}
\end{thm}

\begin{prob} For each one of the $24$ minimal Taylor algebras on a set of size $3$, find a generating set of relations for the corresponding relational clone. Are they all finitely related?
\end{prob}

So far, generating sets of relations have been found for $18$ of the $24$ minimal Taylor algebras on a set of size $3$. In a big surprise, Maximilian Hadek recently discovered that the conservative minimal Taylor algebra which has two majority edges and one semilattice edge is \emph{not} finitely related - so far, no explicit (infinite) generating set of relations is known for that algebra. (I think these results haven't been published yet.)

\section{The strands of an unlinked CSP instance, and a safe recursive strategy}

Generally speaking, in order to guarantee a polynomial running time for solving CSPs we attempt to avoid recursion. There is a form of recursion which can be safely applied, however: we can recursively solve polynomially many subproblems as long as the size of \emph{every} variable's domain is strictly reduced in each subproblem. The resulting algorithms will then have the property that the exponent in the running time will depend on the size of the largest domain of any variable. This approach seems to have been introduced with the solutions to the CSP dichotomy conjecture for conservative algebras by Bulatov \cite{bulatov-conservative}, \cite{bulatov-conservative-rerevisited} and Barto \cite{barto-conservative-revisited}, as well as Mikl\'os Mar\'oti's ``Tree on top of Maltsev'' algorithm \cite{tree-on-top-of-malcev} (which used this sort of recursive strategy in a \emph{very} different way from what we will consider in this section).

The challenge now is to find situations where we can usefully reduce to subproblems in which every single variable domain is reduced. The prototypical example of how this may occur is when a CSP instance is \emph{unlinked}.

\begin{defn} Let $\fX$ be an instance of a CSP, with variable domains $\bA_x$ for each variable $x$. We say that $\fX$ is \emph{unlinked} at a variable $x$ of $\fX$ if there are some $a,b \in \bA_x$ such that there are \emph{no} cycles $p$ of $\fX$ from $x$ to $x$ with $b \in \{a\} + p$. We say that $\fX$ is \emph{unlinked} if it is unlinked at every variable.

If an instance $\fX$ is unlinked at a variable $x$, then we define the \emph{linking relation} $\theta_\fX$ of $\fX$ at $x$ as the equivalence relation on $\bA_x$ defined by $(a,b) \in \theta_\fX$ iff there is some cycle $p$ from $x$ to $x$ such that $b \in \{a\} + p$.
\end{defn}

\begin{prop} If $\fX$ is cycle-consistent and unlinked at $x$, then the linking relation $\theta_\fX$ is a congruence of $\bA_x$.
\end{prop}
\begin{proof} We just need to check that $\theta_\fX$ is closed under unary polynomial operations of $\bA_x$. So suppose that $(a,b) \in \theta_\fX$, $c_1, ..., c_n \in \bA_x$, and $f$ is some $n+1$-ary polymorphism of $\bA_x$. Since $(a,b) \in \theta_\fX$, there is some cycle $p$ from $x$ to $x$ such that $b \in \{a\}+p$. By cycle-consistency, we also have $c_i \in \{c_i\} + p$ for each $i$, so if $\PP_p \le \bA_x^2$ is the binary relation corresponding to the cycle $p$, then we have
\[
f\left(\begin{bmatrix}a\\b\end{bmatrix},\begin{bmatrix}c_1\\c_1\end{bmatrix},\cdots,\begin{bmatrix}c_n\\c_n\end{bmatrix}\right) \in \PP_p.
\]
Thus we have $(f(a, c_1, ..., c_n), f(b, c_1, ..., c_n)) \in \theta_\fX$, which completes the proof.
\end{proof}

As a consequence, one way to discover unlinked subinstances of a cycle-consistent instance $\fX$ is to go through each maximal congruence $\theta$ on some variable domain $\bA_x$, and to greedily build up a subinstance $\fX'$ which is as large as possible subject to the condition $\theta_{\fX'} \le \theta$. The resulting subinstance $\fX'$ only depends on $\theta$ and not on the choices we make during the greedy construction of $\fX'$, so long as $\fX$ is cycle-consistent: if there are two paths $p,q$ from $x$ to $y$ such that $p-p$ and $q-q$ each have linking congruences contained in $\theta$, then cycle-consistency applied to $-q+p$ shows that the linking congruence of $p-q$ is also contained in $\theta$. This is the approach Zhuk used in his proof of the general CSP dichotomy conjecture \cite{zhuk-dichotomy}.

\begin{prop} Suppose that we have a multisorted CSP template $\CSP(\bA_1, ..., \bA_n)$ such that the CSP template $\CSP(\bB_1, ..., \bB_m)$ can be solved in polynomial time, where $\bB_1, ..., \bB_m$ is the collection of all algebras $\bB$ which are isomorphic to a proper subalgebra of some $\bA_i$ (note that $\max |\bB_i| < \max |\bA_j|$). Then we can solve any unlinked instance $\fX$ of $\CSP(\bA_1, ..., \bA_n)$ in polynomial time.
\end{prop}
\begin{proof} We assume that every pair of variables of $\fX$ can be connected by some path without loss of generality, and we shrink $\fX$ by enforcing cycle-consistency (if the shrunk instance is no longer unlinked, then every single variable domain has been reduced, and we can solve the instance). Pick any variable $x$ of $\fX$, and let $\theta_\fX$ be the linking congruence on $\bA_x$. Then for each congruence class $a/\theta_\fX \in \bA_x/\theta_\fX$, if we restrict the possible values of $x$ to $a/\theta$ and enforce arc-consistency, then every other variable $y$ of $\fX$ has its domain restricted to some subset of a congruence class of $\bA_y/\theta_\fX$ (where here we interpret $\theta_\fX$ as the linking congruence of $\fX$ on $\bA_y$). Since the value of $x$ must be in \emph{some} congruence class of $\bA_x/\theta_\fX$, and since there are only a constant number of such congruence classes to check, we can solve $\fX$ by solving a constant number of instances of CSPs with templates of the form $\CSP(\bB_1, ..., \bB_m)$, where for each variable $y$ the algebra $\bB_i$ is contained in an congruence class of $\bA_y/\theta_\fX$.
\end{proof}

Of course, an entire instance being unlinked is fairly rare. Additionally, the assumption of cycle-consistency is probably unnecessarily strong. The approach used in the algorithms for the conservative CSP dichotomy is a variant of the strategy of looking for unlinked subinstances, starting from the idea of looking for \emph{any} useful way of properly restricting the domains of some subset of the variables. This will naturally lead to different consistency principles (which are likely to also be unecessarily strong).

The most general thing we could do along these lines is the following. For each variable $x$, and for each element $a \in \bA_x$, restrict the domain of $x$ to the singleton $\{a\}$, and run arc-consistency (or cycle-consistency, etc.). Restrict our attention to the set $Y$ of variables $y$ whose domain has shrunk as a result of this restriction (replacing every relation which involves variables outside of $Y$ by its projection onto the variables contained in $Y$), and solve the resulting instance recursively to see if we can rule out the element $a \in \bA_x$ as a possible value for $x$.

Although this scheme is readily implemented, it is hard to algebraically control what happens as a result. Instead, we will consider the digraph of implications between restrictions on the individual domains, and ask under which conditions this has a nice structure.

\begin{defn} If $\fX$ is an instance of a multisorted CSP with variable domains $\bA_x$, then we define the \emph{implication digraph} to be the directed graph on $\fX$ where the vertices are pairs $(x, \bB)$ such that $x$ is a variable and $\bB \le \bA_x$, and where we have a directed edge from $(x,\bB)$ to $(y,\bC)$ if there is a path $p$ of length $1$ from $x$ to $y$ in $\fX$ such that $\bB + p = \bC$.

We write $(x,\bB) \preceq (y,\bC)$ if there is a path $p$ of any length from $x$ to $y$ such that $\bB + p = \bC$. The resulting quasiorder is the \emph{implication qoset}.

Let $\mathcal{E}$ be a subdigraph of the implication digraph, and consider $\cE$ as a subqoset of the implication qoset by taking its transitive closure. A \emph{strand} of $\mathcal{E}$ is just an equivalence class of $\mathcal{E}$. Often we might take $\cE$ to be the subdigraph of pairs $(x,\bB)$ such that $\bB$ is a proper subalgebra of $\bA_x$ (or a proper absorbing subalgebra, etc.). We say that $\cS$ is a \emph{maximal} strand of $\mathcal{E}$ if it is a maximal equivalence class of the qoset $\cE$. If $\cE$ is not specified, then a strand can be any strongly connected subset of the implication digraph.

A strand is called \emph{absorbing} if for all $(x,\bB) \in \cS$ we have $\bB \lhd \bA_x$. Note that as long as $\fX$ is arc-consistent, if this occurs for \emph{some} $(x,\bB) \in \cS$ then it occurs for \emph{all} $(x,\bB) \in \cS$.

We define the \emph{partial restriction} of $\fX$ to the strand $\cS$ by reducing the domain of each variable $x$ to $\bigcap_{(x,\bB)\in\cS} \bB$ if some $(x,\bB) \in \cS$, and leaving the domain of $x$ unchanged otherwise.
\end{defn}

\begin{prop}\label{prop-cycle-strand-finding} If $\fX$ is cycle-consistent, $\cE$ is a subdigraph of the implication digraph, and we are given some $(x,\bB) \in \cE$, then we can find a maximal strand of $\cE$ in polynomial time.
\end{prop}
\begin{proof} We start from any element of $\cE$ and keep following single-step paths until we stabilize at a maximal strand. To see that this takes only polynomially many steps, note that if $(x,\bB) \preceq (x,\bC)$, then cycle-consistency implies that $\bB \le \bC$, so if $\bB \ne \bC$ then $|\bC| \ge |\bB|+1$.
\end{proof}

In order to restrict an instance $\fX$ to a strand $\cS$ in a way that guarantees that the restricted instance is arc-consistent, it seems like we should require that the strand $\cS$ includes at most one pair $(x,\bB)$ for each variable $x \in \fX$. This will be guaranteed as long as the instance $\fX$ is sufficiently consistent.

\begin{prop} If $\fX$ is an arc-consistent instance, then each strand $\cS$ of $\fX$ includes at most one pair $(x,\bB)$ for each variable $x$ as long as $\fX$ satisfies property (P3) from Definiion \ref{defn-weak-prague}:
\[
\bB + p + q = \bB \;\; \implies \;\; \bB + p = \bB
\]
for all pairs of cycles $p,q$ from $x$ to $x$ and all $\bB \le \bA_x$. In particular, this always occurs if $\fX$ is $pq$-consistent.
\end{prop}

\begin{defn}\label{defn-full-restriction} Suppose that $\fX$ is an instance, and let $\cS$ be any strand of $\fX$ such that for each variable $x$, there is at most one $\bB_x$ such that $(x,\bB_x) \in \cS$. We define the \emph{full restriction} of $\fX$ to $\cS$ to be the instance $\fX'$ whose variable set is the set of variables $x$ such that some $(x,\bB_x) \in \cS$, with each variable domain reduced to $\bB_x$, where for each relation of $\fX$ we project it onto the set of variables in $\fX'$ and restrict each of its coordinates to $\bB_x$.
\end{defn}

In order to mimic the situation of an unlinked instance as well as possible, we might also like to have the property that for each variable $x$ such that some $(x,\bB)$ is in the strand $\cS$, we have parallel strands $\cS'$ which partition the domain $\bA_x$ into disjoint subalgebras $\bB'$ with each $(x,\bB') \in \cS'$.

\begin{prop} If $\fX$ is an arc-consistent instance, then for each strand $\cS$, each $(x,\bB)$, and each $a \not\in \bB$, there is a strand $\cS'$, involving the same set of variables as $\cS$, such that $(x,\bB') \in \cS'$ for some $\bB' \le \bA_x$ with $a \in \bB'$ and $\bB \cap \bB' = \emptyset$, as long as $\fX$ also satisfies property (P2) from Definition \ref{defn-weak-prague}:
\[
\bB + p = \bB \;\; \implies \;\; \bB - p = \bB
\]
for all cycles $p$ from $x$ to $x$ and all $\bB \le \bA_x$.

In this case, a more precise statement is true: for each $(x,\bB) \in \cS$ there is a congruence $\theta_\cS$ on $\bA_x$ such that $\bB$ is a union of congruence classes of $\theta_\cS$, and each congruence class $\bB'$ of $\theta_\cS$ is contained in a parallel strand $\cS'$.
\end{prop}
\begin{proof} Consider the set $\mathcal{C}$ of all cycles $p$ from $x$ to $x$ such that $\bB + p = \bB$. Property (P2) is the guarantee that for each $p \in \mathcal{C}$, $\bB$ will be a union of congruence classes of the linking congruence of $\PP_p$. Define $\theta_\cS$ to be the join of the collection of all linking congruences of $\PP_p$ for $p \in \mathcal{C}$. Then $\bB$ is still a union of congruence classes of $\theta_\cS$.

For any congruence class $\bB'$ of $\theta_\cS$, we define the parallel strand $\cS'$ as follows: for every way of splitting some $p \in \mathcal{C}$ as $p = q + r$, where $q$ is a path from $x$ to a variable $y$, we include $(y,\bB'+q)$ in $\cS'$. Since $\bB'$ is a congruence class of $\theta_\cS$, $\bB'$ will be a union of congruence classes of the linking congruence of $\PP_p$, so $\bB' + p - p = \bB'$, so $(\bB'+q) + r-p = \bB'$, and we see that $(x,\bB')$ and $(y,\bB'+q)$ are indeed part of the same equivalence class of the implication qoset of $\fX$.
\end{proof}

To guarantee that the restriction of our instance $\fX$ to any strand is arc-consistent, we need a still stronger consistency principle. This consistency principle was used to give one of the original proofs of the fact that $(2,3)$-consistency was strong enough to ensure satisfiability of any CSP with bounded width \cite{barto}. I will define it a bit differently here, but will show that the definition given here is equivalent to the original definition by following an argument from \cite{barto-conservative-revisited}.

\begin{defn} An arc-consistent instance $\fX$ is called a \emph{Prague instance} if for every $(x,\bB), (y,\bC)$ in the same equivalence class of the implication qoset, and for every path $p$ of length $1$ from $x$ to $y$, we have $\bB + p = \bC$.
\end{defn}

\begin{prop} If $\fX$ is a Prague instance, then $\fX$ satisfies conditions (P2) and (P3) of Definition \ref{defn-weak-prague}, so $\fX$ is a weak Prague instance.
\end{prop}

\begin{prop}[Barto and Kozik \cite{barto-kozik-bounded-width}, \cite{barto-conservative-revisited}]\label{prop-prague} If $\fX$ is arc-consistent, then the following are equivalent:
\begin{itemize}
\item[(a)] for every $(x,\bB), (y,\bC)$ in the same equivalence class of the implication qoset, and for every path $p$ of length $1$ from $x$ to $y$, we have $\bB + p \subseteq \bC$.
\item[(b)] $\fX$ is a Prague instance,
\item[(c)] for every variable $x$ and every pair of paths $p,q$ from $x$ to $x$ such that the variables involved in $p$ are a subset of the variables involved in $q$, if $b \in \{a\} + p$ then there is some $j \ge 0$ such that $b \in \{a\} + jq$,
\item[(d)] for every variable $x$ and every pair of paths $p,q$ from $x$ to $x$ such that the variables involved in $p$ are a subset of the variables involved in $q$, then for every sufficiently large $j$ we have $\{a\} + p \subseteq \{a\} + jq$.
\end{itemize}
\end{prop}
\begin{proof} For the implication from (a) to (b), note that from $\bB + p \subseteq \bC$ and $\bC - p \subseteq \bB$ we can immediately conclude that $\bB + p = \bC$.

For (b) $\implies$ (c), let $p,q$ be as in (c). Find a $j \ge 1$ such that $\{a\} + jq = \{a\} + 2jq$, write $\bB = \{a\} + jq$, and let $\cS$ be the strand of the implication digraph containing $(x, \bB)$. Since every Prague instance satisfies property (P3), for each variable $y$ which shows up in the cycle $q$ there is a unique $\bC$ such that $(y,\bC) \in \cS$. Since every step of the cycle $p$ goes between variables involved in $\cS$ and $\fX$ is a Prague instance, we have $\bB + p = \bB$, and for the same reason we have $\bB - jq = \bB$. From
\[
a \in \{a\} + jq - jq = \bB - jq = \bB
\]
and $b \in \{a\} + p$, we get
\[
b \in \{a\} + p \subseteq \bB + p = \bB = \{a\} + jq.
\]

For (c) $\implies$ (d), it's enough to prove that (c) implies $a \in \{a\} + jq$ for every sufficiently large $j$. Since $\fX$ is arc-consistent, there is some $b$ such that $b \in \{a\} + q$, and since $a \in \{b\} - q$ (c) implies there is some $j \ge 0$ such that $a \in \{b\} + jq$. Since $b \in \{a\} + q$, we have $a \in \{a\} + (j+1)q$. Now set $r = (j+1)q$, and by applying (c) again we see that there is some $k \ge 0$ such that
\[
a \in \{b\} + kr \subseteq \{a\} + q + k(j+1)q = \{a\} + (k(j+1) + 1)q.
\]
Since $j+1$ and $k(j+1) + 1$ are relatively prime positive integers, every sufficiently large number can be written as a positive combination of $j+1$ and $k(j+1)+1$, which proves (d).

For (d) $\implies$ (a), let $q$ be a path from $y$ to $x$ such that $\bC + q = \bB$, and let $r$ be a path from $x$ to $y$ such that $\bB + r = \bC$. Then for every $j$ we have $\bC + j(q+r) = \bC$. Since the cycle $q+p$ from $y$ to $y$ only involves a subset of the variables involved in $q+r$, by part (d) we have $\bC + (q+p) \subseteq \bC + j(q+r)$ for every sufficiently large $j$, so
\[
\bB + p = \bC + q+p \subseteq \bC + j(q+r) = \bC.\qedhere
\]
\end{proof}

\begin{prop}\label{prop-prague-reduction} If $\fX$ is a Prague instance, and if $\cS$ is a strand of $\fX$, then the full restriction of $\fX$ to $\cS$ (which restricts both the set of variables and the domains) will be arc-consistent - in fact, it will also be a Prague instance.

If $\cS$ is also a \emph{maximal} strand in the qoset of pairs $(x,\bB)$ such that $\bB < \bA_x$, then the partial restriction of $\fX$ to $\cS$ (which restricts the domains but not the set of variables) is also arc-consistent. If $\cS$ is additionally absorbing, then the partial restriction of $\fX$ to $\cS$ will even be a Prague instance.
\end{prop}
\begin{proof} The statement about the full restriction follows from the fact that for any tuple $r$ of any constraint relation $\RR \le_{sd} \bA_{x_1} \times \cdots \times \bA_{x_n}$, if $\pi_{x_i}(r) \in \bB_i$ and $(x_i,\bB_i) \in \cS$, then we necessarily have $\pi_{x_j}(r) \in \bB_j$ for any $(x_j, \bB_j) \in \cS$, by the definition of a Prague instance.

For the statement about the partial restriction, suppose that some variable $x_i$ does not occur in the strand $\cS$, and that the relation $\RR$ involves a variable $x_j$ with $(x_j, \bB_j) \in \cS$. Then by the maximality of $\cS$, we must have $\bB_j + \pi_{x_jx_i}(\RR) = \bA_{x_i}$, so for each $a \in \bA_{x_i}$, there is a tuple $r \in \RR$ with $\pi_{x_i}(r) = a$ and $\pi_{x_j}(r) \in \bB_j$, and then by the first part we have $\pi_{x_k}(r) \in \bB_k$ for each $(x_k, \bB_k) \in \cS$.

Now suppose that $\cS$ is absorbing. By the equivalence between parts (b) and (c) of Proposition \ref{prop-prague}, it's enough to show that for every cycle $p$ from $x$ to $x$ in $\fX$, if we let $p'$ be the corresponding path in the partially restricted instance, then for any $a,b \in \bB_x$ with $(x,\bB_x) \in \cS$ and
\[
b \in \{a\}+p,
\]
there is some $j \ge 1$ such that
\[
b \in \{a\} + jp'.
\]
For this, we let $\RR = \PP_p \le \bA_x\times \bA_x$ and $\bS = \PP_{p'} \le \bB_x \times \bB_x$, and we aim to apply Corollary \ref{cor-absorbing-directed-path}.

Since $\cS$ is absorbing, we have $\bS \lhd \RR$. Since the restriction of $\fX$ to $\cS$ is arc-consistent, $\bS$ is subdirect in $\bB_x \times \bB_x$. Thus, there is some directed cycle $\mathcal{C}_a$ of $\bS$ which can be reached from $a$ in $\bS$, and there is some directed cycle $\mathcal{C}_b$ of $\bS$ with a directed path from $\mathcal{C}_b$ to $b$ in $\bS$. Since $\fX$ is a Prague instance and $\mathcal{C}_a, \mathcal{C}_b$ are contained in the same weakly connected component of $\RR$, there is some directed path from $\mathcal{C}_a$ to $\mathcal{C}_b$ in $\RR$ by Proposition \ref{prop-p2}. Now we can apply Corollary \ref{cor-absorbing-directed-path} to see that there is a directed path from $\mathcal{C}_a$ to $\mathcal{C}_b$ in $\bS$, so there is a directed path from $a$ to $b$ in $\bS$, and we are done.
\end{proof}

In practice, the local consistency algorithm used to reduce general instances to Prague instances or to cycle-consistent instances actually produces an instance with an even greater level of consistency.

\begin{defn} An instance $\fX$ is $(l,k)$-\emph{minimal}, for $k \ge l$, if
\begin{itemize}
\item every set of at most $k$ variables of $\fX$ is in the scope of some constraint, and
\item for any set $S$ of at most $l$ variables and any pair of constraints $C_1, C_2$ of $\fX$ whose scopes contain $S$, the existential projections of $C_1$ and $C_2$ to the variables in $S$ are the same.
\end{itemize}
\end{defn}

\begin{prop} If $\fX$ is $(2,3)$-minimal, then $\fX$ is a Prague instance, and $\fX$ is cycle-consistent.
\end{prop}
\begin{proof} We leave cycle-consistency to the reader, and will prove that $\fX$ is a Prague instance. By the equivalence of (a) and (b) in Proposition \ref{prop-prague}, it's enough to show that for $(x,\bB) \preceq (y,\bC)$ and $p$ a $1$ step path from $x$ to $y$, we always have $\bB + p \subseteq \bC$.

Let $q$ be the path from $x$ to $y$ with $\bB + q = \bC$, and suppose that the variables which occur along the path $q$ are $x = x_0, x_1, ..., x_n = y$, and let $q = q_1 + \cdots + q_n$ be the decomposition of $q$ into single step paths. For each $i$, let $\bB_i = \bB + q_1 + \cdots + q_i$, and let $p_i$ be a single step path from $x_i$ to $y$ (which exists by $(2,3)$-minimality of $\fX$).

We prove by induction on $i$ that $\bB + p \subseteq \bB_i + p_i$. For the inductive step, we need to show that $\bB_i + p_i \subseteq \bB_{i+1} + p_{i+1}$, that is,
\[
\bB_i + p_i \subseteq \bB_i + q_{i+1} + p_{i+1}.
\]
This follows from the fact that there is some constraint whose scope contains the three variables $x_i, x_{i+1}$, and $y$, together with the fact that the two-variable projections of this constraint onto pairs of variables from $\{x_i,x_{i+1},y\}$ agree with the binary relations $\PP_{p_i}, \PP_{q_{i+1}}, \PP_{p_{i+1}}$, by $(2,3)$-minimality.
\end{proof}

The relationships between the various types of consistency introduced so far are summarized in the Hasse diagram below (to see that singleton arc-consistency is not implied by $(2,3)$-minimality, consider the four-variable instance of 1-IN-3 SAT where every group of three variables is required to satisfy the 1-IN-3 constraint).

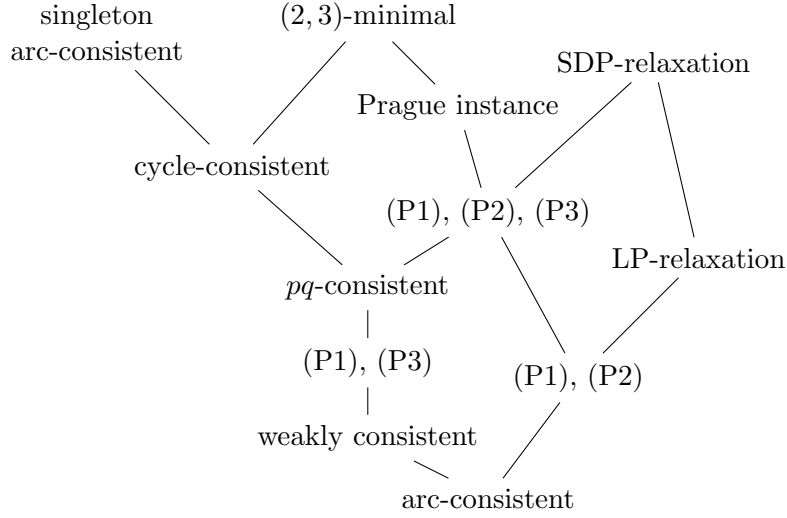
\begin{center}
\begin{tikzpicture}[scale=2]
  \node (1) at (0,1) {$(2,3)$-minimal};
  \node (2) at (0.6,0.4) {Prague instance};
  \node (3) [style={align=center}] at (0.8,-0.25) {weak Prague\\ (P1), (P2), (P3)};
  \node (4) at (0,-0.8) {$pq$-consistent};
  \node (5) at (-0.9,0) {cycle-consistent};
  \node (5a) [style={align=center}] at (-1.8,0.9) {singleton\\ arc-consistent};
  \node (6) at (0,-1.3) {(P1), (P3)};
  \node (7) at (1.9,0.7) {SDP-relaxation};
  \node (8) at (2.2,-0.6) {LP-relaxation};
  \node (9) [style={align=center}] at (1.4,-1.4) {reversible\\ (P1), (P2)};
  \node (10) [style={align=center}] at (0.8,-2.4) {arc-consistent\\ (P1)};
  \node (11) at (0, -1.8) {weakly consistent};
  \draw (1) -- (2) -- (3) -- (4) -- (5) -- (1);
  \draw (4) -- (6);
  \draw (5) -- (5a);
  \draw (3) -- (7) -- (8) -- (9) -- (3);
  \draw (9) -- (10) -- (11) -- (6);
\end{tikzpicture}
\end{center}

When our instance consists of just a single relation $\RR$ (with no repeated variables), all of these consistency conditions become equivalent to subdirectness of the relation $\RR$. Studying this very special case is what leads to the main algebraic ingredient we will need for the conservative CSP dichotomy.

\begin{defn}\label{defn-relation-strand} If $\RR \le_{sd} \bA_1 \times \cdots \times \bA_n$, and if $\bB_i \le \bA_i$ for all $i$, then an $(\RR,\bB)$-\emph{strand} is an equivalence class of the quasiorder $\preceq$ on $[n]$ which is the transitive closure of
\[
\bB_i + \pi_{ij}(\RR) = \bB_j \;\; \implies \;\; i \preceq j.
\]
The equivalence classes of this quasiorder are the same as the equivalence classes of the more permissive quasiorder
\[
i \preceq' j \;\; \iff \;\; \bB_i + \pi_{ij}(\RR) \subseteq \bB_j,
\]
and these quasiorders are the same in the special case where $\bB_i$ is a minimal absorbing subalgebra of $\bA_i$ for each $i$.
\end{defn}

Using either the fact that the instance consisting of just $\RR$ is cycle-consistent, or the fact that it satisfies (P1) and (P2), we have the following result.

\begin{prop}\label{prop-relation-strand} If $\RR \le_{sd} \bA_1 \times \cdots \times \bA_n$, and if $\bB_i \le \bA_i$ for all $i$, then for each $(\RR,\bB)$-strand $S \subseteq [n]$ there is a congruence $\theta_i \in \Con(\bA_i)$ for each $i \in S$ such that each $\bB_i$ is a union of congruence classes of $\theta_i$, and for any $i,j \in S$, the binary relation
\[
\pi_{ij}(\RR)/(\theta_i\times\theta_j) \le_{sd} \bA_i/\theta_i \times \bA_j/\theta_j
\]
is the graph of an isomorphism.
\end{prop}

The algebraic input needed for conservative CSPs is to show that if we take the $\bB_i$s to be minimal absorbing subalgebras of the $\bA_i$s such that $\RR \cap (\bB_1 \times \cdots \times \bB_n) \ne \emptyset$, then the $(\RR,\bB)$-strands do not interact with each other. This algebraic miracle is a special property of conservative Taylor algebras - it doesn't seem to hold in general.

\section{The rectangularity theorem for conservative Taylor algebras}

There are two versions of the Rectangularity Theorem for conservative algebras: one uses the theory of absorbing subalgebras and is proved in \cite{barto-conservative-revisited}, while the other uses Bulatov's theory of affine-semilattice components (shortened to \emph{as-components}) of the colored graph and is proved in \cite{bulatov-conservative-rerevisited}. Generally speaking, the affine-semilattice components of conservative Taylor algebras (and of subdirect products of conservative Taylor algebras) tend to behave like minimal absorbing subalgebras. We will mainly focus on the absorbing subalgebra approach.

First we will state both versions of the Rectangularity Theorem, before diving into the proofs.

\begin{thm}[Rectangularity Theorem for conservative Taylor algebras, absorbing version \cite{barto-conservative-revisited}]\label{thm-rectangularity-conservative-absorbing} If $\RR \le_{sd} \bA_1 \times \cdots \times \bA_n$ is a subdirect product of conservative Taylor algebras $\bA_i$, and if a system of minimal absorbing subalgebras $\bB_i \llhd \bA_i$ for each $i$ satisfies
\[
\RR \cap (\bB_1 \times \cdots \times \bB_n) \ne \emptyset,
\]
then
\[
\RR \cap (\bB_1 \times \cdots \times \bB_n) = \prod_{S\text{ an }(\RR,\bB)\text{-strand}} \Big(\pi_S(\RR) \cap \prod_{i\in S}\bB_i\Big).
\]
\end{thm}

The choice of absorption concept used in the Rectangularity Theorem is fairly arbitrary - we could use J\'onsson absorption or central absorption (as long as the algebras in question are minimal Taylor) instead, and it would still be true. We will use $\lhd_X$ in this section to refer to any choice of an absorption concept as in Section \ref{s-absorbing-arc-consistent} for which the Absorption Theorem \ref{absorption-theorem} applies (so if we take $\lhd_X$ to be central absorption, then we need to assume that we are in a context where binary absorption implies central absorption, such as the context of minimal Taylor algebras).

\begin{defn} If $\bA$ is a subdirect product of conservative minimal Taylor algebras, then we define a quasiorder on the elements of $\bA$ by $a \preceq_{as} b$ if there is a sequence $a = a_0, a_1, ..., a_n = b$ such that for each $i$ either $\{a_i,a_{i+1}\}$ is a $\ZZ/2^{\aff}$-subalgebra or $a_i \rightarrow a_{i+1}$.

We say that a subset $B \subseteq \bA$ is an \emph{as-component} of $\bA$ if $B$ is a maximal equivalence class of the quasiorder $\preceq_{as}$. We say that $\bA$ is \emph{strongly as-connected} if $\bA$ consists of a single equivalence class of the quasiorder $\preceq_{as}$.
\end{defn}

\begin{thm}[Rectangularity Theorem for conservative Taylor algebras, as-component version \cite{bulatov-conservative-rerevisited}]\label{thm-rectangularity-conservative-as} If $\RR \le_{sd} \bA_1 \times \cdots \times \bA_n$ is a subdirect product of conservative minimal Taylor algebras $\bA_i$, and if a system of as-components $\bB_i \subseteq \bA_i$ for each $i$ satisfies
\[
\RR \cap (\bB_1 \times \cdots \times \bB_n) \ne \emptyset,
\]
then
\[
\RR \cap (\bB_1 \times \cdots \times \bB_n) = \prod_{S\text{ an }(\RR,\bB)\text{-strand}} \Big(\pi_S(\RR) \cap \prod_{i\in S}\bB_i\Big).
\]
\end{thm}

We will frequently use the following consequence of the Absorption Theorem \ref{absorption-theorem} throughout the proof. Note that every strongly as-connected algebra is certainly absorption-free, so it applies in that context as well (it's a good exercise to give a direct proof of the analogue for strongly as-connected algebras, along the lines of Theorem \ref{strong-binary} - if you can't solve it, see \cite{bulatov-conservative-rerevisited}).

\begin{prop}\label{prop-absorption-free-linking} Let $\RR \le_{sd} \bA_1 \times \bA_2$ be a subdirect product of absorption-free idempotent algebras, and let $\theta_1$ be a maximal congruence on $\bA_1$. Then either $\theta_1$ contains the linking congruence of $\RR$ on $\bA_1$, in which case there is a maximal congruence $\theta_2$ on $\bA_2$ such that
\[
\RR/(\theta_1\times\theta_2) \le \bA_1/\theta_1 \times \bA_2/\theta_2
\]
is the graph of an isomorphism from $\bA_1/\theta_1 \xrightarrow{\sim} \bA_2/\theta_2$, or else we have
\[
(a/\theta_1) + \RR = \bA_2
\]
for each congruence class $a/\theta_1$ of $\theta_1$.
\end{prop}
\begin{proof} This is a restatement of Corollary \ref{absorption-free-binary}.
\end{proof}

The next lemma is one of the key places where conservativity is really used in the argument. Its analogue for as-components of conservative minimal Taylor algebras (replacing ``absorption-free'' with ``strongly as-connected'') is an easy exercise.

\begin{lem}[Barto \cite{barto-conservative-revisited}]\label{lem-conservative-subalgebra-absorption-free} If $\bA$ is an absorption-free conservative algebra, $\theta$ is a proper congruence on $\bA$, and $\bB \le \bA$ is any subalgebra such that $\bB$ has at least one element from each congruence class of $\theta$, then $\bB$ is also absorption-free.
\end{lem}
\begin{proof} Suppose for contradiction that $\bC \lhd_X \bB$, with $\bC \ne \bB$. Pick some $b \in \bB\setminus\bC$ such that $\bC \not\subseteq b/\theta$ (this is possible as long as $\theta$ has at least two congruence classes), and let $\bB'$ be $(\bB \setminus (b/\theta))\cup\{b\}$, that is, $\bB'$ is the subalgebra of $\bB$ formed by removing every element which is congruent to $b$ other than $b$ itself (that $\bB'$ is a subalgebra follows from the fact that $\bB$ is conservative). Then since $\lhd_X$ is compatible with pp-formulas, if we set $\bC' = \bC \cap \bB'$, we have $\bC' \lhd_X \bB'$. Applying compatibility with pp-formulas again, we have
\[
\bC'/\theta \lhd_X \bB'/\theta = \bA/\theta,
\]
and by the construction of $\bC'$ we see that $\bC'/\theta \ne \bA/\theta$, since $b/\theta \not\in \bC'/\theta$. Therefore $\bA$ has a proper absorbing subalgebra (i.e. the preimage of $\bC'/\theta$ under the quotient homomorphism $\bA \twoheadrightarrow \bA/\theta$), which is a contradiction.
\end{proof}

The next result is also specific to conservative algebras (see Example \ref{ex-subdirect-2-semi-not-absorption-free} for a counterexample in the 2-semilattice case). An analogue for strongly as-connected algebras can be proved using the same argument (see \cite{bulatov-conservative-rerevisited}).

\begin{thm}[Barto \cite{barto-conservative-revisited}]\label{thm-conservative-subdirect-absorption-free} If $\RR \le_{sd} \bA_1 \times \cdots \times \bA_n$ is a subdirect product of absorption-free conservative algebras, then $\RR$ is also absorption-free.
\end{thm}
\begin{proof} We induct on $n$ and on the sizes of the $\bA_i$s. Suppose that $\bS \lhd_X \RR$ - we aim to prove that $\bS = \RR$. By compatibility with pp-formulas we have $\pi_i(\bS) \lhd_X \pi_i(\RR) = \bA_i$ for each $i$, so since the $\bA_i$s are absorption-free $\bS$ is also a subdirect product of the $\bA_i$s.

Let $\theta_1$ be a maximal congruence on $\bA_1$. By Proposition \ref{prop-absorption-free-linking}, for each $i$ either
\begin{itemize}
\item there is a maximal congruence $\theta_i \in \Con(\bA_i)$ such that $\pi_{1i}(\RR)/\theta_1\times\theta_i$ is the graph of an isomorphism between $\bA_1/\theta_1$ and $\bA_i/\theta_i$, or
\item we have $a/\theta_1 + \pi_{1i}(\RR) = \bA_i$ for all congruence classes $a/\theta_1$ of $\theta_1$.
\end{itemize}
Rearrange the coordinates so that $\bA_1, ..., \bA_k$ are in the first case, with corresponding maximal congruences $\theta_i$, while $\bA_{k+1}, ..., \bA_n$ are in the second case. Define a congruence $\theta$ on the product by
\[
\theta = \theta_1 \times \cdots \times \theta_k \times 0_{\bA_{k+1}} \times \cdots \times 0_{\bA_n} \in \Con(\bA_1 \times \cdots \times \bA_n).
\]

First suppose that $\theta$ is the trivial congruence (i.e. each $\theta_i = 0_{\bA_i}$). In this case each $\bA_i$ with $1 < i \le k$ is a redundant coordinate (since $\pi_{1i}(\RR)$ is the graph of an isomorphism between $\bA_1$ and $\bA_i$), so we can assume without loss of generality that $k = 1$. Let $a$ be any element of $\bA$, and consider the relations $\bS_a, \RR_a$ given by
\[
\RR_a = \{(x_2, ..., x_n) \mid (a,x_2, ..., x_n) \in \RR\} \le \bA_2 \times \cdots \times \bA_n,
\]
with $\bS_a$ defined similarly. By compatibility with pp-formulas, we have $\bS_a \lhd_X \RR_a$, and since $a + \pi_{1i}(\RR) = \bA_i$ for each $i \ge 2$, $\RR_a$ is a subdirect product of $\bA_2, ..., \bA_n$. By the induction hypothesis, we then have $\bS_a = \RR_a$, and since this is true for every $a \in \bA$, we have $\bS = \RR$.

Now suppose that $\theta$ is nontrivial - suppose without loss of generality that $\theta_1$ is a nontrivial congruence of $\bA_1$, with $(a,a') \in \theta_1$ for some $a \ne a'$. Let $\bB = \bA \setminus \{a\}$, and let $\bB' = \bA\setminus\{a'\}$. By Lemma \ref{lem-conservative-subalgebra-absorption-free}, each of $\bB, \bB'$ is absorption-free. Define $\RR_\bB$ by
\[
\RR_\bB = \{(x_2, ..., x_n) \mid \exists b \in \bB \;\; (b,x_2, ..., x_n) \in \RR\} \le \bA_2 \times \cdots \times \bA_n,
\]
and similarly define $\bS_\bB$, $\RR_{\bB'}, \bS_{\bB'}$. We have $\bS_\bB \lhd_X \RR_\bB$ and $\bS_{\bB'} \lhd_X \RR_{\bB'}$ by compatibility with pp-formulas. For each $i \le k$, we have
\[
\pi_i(\RR_\bB)/\theta_i = \pi_i(\RR_\bB/\theta) = \bB/\theta_1 + \pi_{1i}(\RR/\theta) = \bA_1/\theta_1 + \pi_{1i}(\RR/\theta) = \bA_i/\theta_i,
\]
so by Lemma \ref{lem-conservative-subalgebra-absorption-free}, $\pi_i(\RR_\bB)$ is absorption-free for each $i \le k$. For $i > k$, we have
\[
\pi_i(\RR_\bB) = \bB + \pi_{1i}(\RR) = \bA_i,
\]
since $\bB$ contains at least one full congruence class of $\theta_1$. Thus for each $i \in [n]$, $\pi_i(\RR_\bB)$ is absorption-free, so since $|\bB| < |\bA_1|$ we can apply the induction hypothesis to see that $\bS_\bB = \RR_\bB$. A similar argument shows that $\bS_{\bB'} = \RR_{\bB'}$, and since $\bB \cup \bB' = \bA_1$ we see that $\bS = \RR$.
\end{proof}

In order to get a foothold into the Rectangularity Theorem, we start with the case of binary relations. We will only assume that one of the algebras involved is conservative - this way we can later apply the result with the other algebra equal to a larger subdirect product of conservative algebras. Once again, an analogous argument works for as-components (we need a version of Theorem \ref{strong-binary} for as-components to push the argument through, see \cite{bulatov-conservative-rerevisited}).

\begin{lem}[Barto \cite{barto-conservative-revisited}] Suppose $\RR \le_{sd} \bA_1 \times \bA_2$, where $\bA_1$ is conservative, $\bA_2$ is idempotent, and both are finite Taylor algebras. Suppose further that $\bB_i \llhd_X \bA_i$ for $i=1,2$, that $\RR \cap (\bB_1 \times \bB_2) \ne \emptyset$, and that there is some $(a,b) \in \RR$ with $a \in \bA_1\setminus\bB_1$ and $b \in \bB_2$. Then $\bB_1 \times \bB_2 \subseteq \RR$.
\end{lem}
\begin{proof} Since $(\bB_1 + \RR) \cap \bB_2$ is a nonempty absorbing subalgebra of $\bB_2$ and $\bB_2$ is absorption-free, we have $\bB_1 + \RR \supseteq \bB_2$, and similarly $\bB_1 \subseteq \bB_2 - \RR$. Thus $\RR \cap (\bB_1 \times \bB_2)$ is subdirect in $\bB_1 \times \bB_2$. Additionally, $\bB_1 \cup \{a\}$ is a subalgebra of $\bA_1$ since $\bA_1$ is conservative, so by the assumption $a \in \bB_2 - \RR$ we can assume without loss of generality that $\bA_1 = \bB_1 \cup \{a\}$ and $\bA_2 = \bB_2$.

If $\RR$ is linked, then by Theorem \ref{absorbing-linked} so is $\RR \cap (\bB_1 \times \bB_2)$, and then by the Absorption Theorem \ref{absorption-theorem} we have $\bB_1 \times \bB_2 \subseteq \RR$. Otherwise, the linking congruence of $\RR$ is a proper congruence $\theta_1$ on $\bA_1$. If $c \in (\{b\} - \RR) \cap \bB_1$, then we have $(a,c) \in \theta_1$, so $\bA_1/\theta_1 = \bB_1/\theta_1$. Let
\[
\bA_1' = (\bA_1 \setminus (a/\theta_1)) \cup \{a\},
\]
then by compatibility with pp-formulas we have
\[
\bB_1 \setminus (a/\theta_1) = \bB_1 \cap \bA_1' \lhd_X \bA_1',
\]
so
\[
(\bB_1/\theta_1) \setminus (a/\theta_1) \lhd_X \bA_1'/\theta_1 = \bB_1/\theta_1,
\]
which implies that $\bB_1 \setminus (a/\theta_1) \lhd_X \bB_1$, which is a contradiction.
\end{proof}

Now we bootstrap our way up.

\begin{lem}[Barto \cite{barto-conservative-revisited}] Suppose $\RR \le_{sd} \bA_1 \times \cdots \times \bA_n \times \bA_{n+1}$, where $\bA_1, ..., \bA_n$ are conservative, $\bA_{n+1}$ is idempotent, and each $\bA_i$ is a finite Taylor algebra. Suppose that we have $\bB_i \llhd_X \bA_i$ for all $i \in [n+1]$, that
\[
\RR \cap (\bB_1 \times \cdots \times \bB_n \times \bB_{n+1}) \ne \emptyset,
\]
that $[n]$ is an $(\RR,\bB)$-strand, and that there is some $(a_1, ..., a_n, b_{n+1}) \in \RR$ such that $a_i \in \bA_i\setminus \bB_i$ for $i \in [n]$ while $b_{n+1} \in \bB_{n+1}$. Then we have
\[
\RR \cap (\bB_1 \times \cdots \times \bB_n \times \bB_{n+1}) = (\pi_{[n]}(\RR) \cap (\bB_1 \times \cdots \times \bB_n)) \times \bB_{n+1}.
\]
\end{lem}
\begin{proof} We induct on $n$ and on the sizes of the $\bA_i$s, and note that we have already proved the case $n=1$ in the previous lemma. By Proposition \ref{prop-relation-strand} we can find maximal congruences $\theta_i \in \Con(\bA_i)$ for $i \in [n]$ such that
\[
\pi_{ij}(\RR)/(\theta_i \times \theta_j) \le \bA_i/\theta_i \times \bA_j/\theta_j
\]
is the graph of an isomorphism for all $i,j \in [n]$. If any $\theta_i$ is trivial (i.e. if $\bA_i$ is simple for some $i \in [n]$) then the $i$th coordinate of $\RR$ is redundant, so we can apply the induction hypothesis. Otherwise, we have $\theta_1 \ne 0_{\bA_1}$.

Let $b_1$ be any element of $\bB_1$. Let $\bA_1'$ be any proper subalgebra of $\bA_1$ such that $b_1, a_1 \in \bA_1'$ and such that $\bA_1'/\theta_1 = \bA_1/\theta_1$, and let $\bB_1' = \bB_1 \cap \bA_1'$. For each $i \in [n+1]$ define $\bA_i', \bB_i'$ by
\[
\bA_i' = \bA_1' + \pi_{1i}(\RR), \;\;\; \bB_i' = \bB_i \cap \bA_i'.
\]
Then since $\pi_{1i}(\RR)/(\theta_1 \times \theta_i)$ is the graph of an isomorphism for each $i \in [n]$, we have $\bA_i'/\theta_i = \bA_i/\theta_i$ and $\bB_i' = \bB_1' + \pi_{1i}(\RR)$ for all $i \in [n]$. Then since each $\bB_i'/\theta_i = \bB_i/\theta_i$, we can apply Lemma \ref{lem-conservative-subalgebra-absorption-free} to see that $\bB_i'$ is absorption-free for each $i \in [n]$. Additionally, by the previous lemma (i.e., the $n=1$ case) we have
\[
\bB_1 \times \bB_{n+1} \subseteq \pi_{1,n+1}(\RR),
\]
so $\bB_{n+1}' = \bB_{n+1}$. Thus if we set
\[
\RR' = \RR \cap (\bA_1' \times \cdots \times \bA_n' \times \bA_{n+1}'),
\]
then $(a_1, ..., a_n, b_{n+1}) \in \RR'$ and we can apply the induction hypothesis to $\RR'$ to see that
\[
\RR' \cap (\bB_1' \times \cdots \times \bB_n' \times \bB_{n+1}) = (\pi_{[n]}(\RR') \cap (\bB_1' \times \cdots \times \bB_n')) \times \bB_{n+1}.
\]
In particular, any tuple in $(\pi_{[n]}(\RR) \cap \prod_{i\in [n]} \bB_i) \times \bB_{n+1}$ such that the first coordinate is $b_1$ is contained in $\RR'$, and therefore is also contained in $\RR$. Since $b_1$ was an arbitrary element of $\bB_1$, we are done.
\end{proof}

\begin{proof}[Proof of the Rectangularity Theorem \ref{thm-rectangularity-conservative-absorbing}, following \cite{barto-conservative-revisited}] We induct on $n$, the number of algebras occuring in the product. Suppose for the sake of contradiction that $\RR$ is a counterexample, i.e. that there is some
\[
b = (b_1, ..., b_n) \in \prod_{S\text{ an }(\RR,\bB)\text{-strand}} \Big(\pi_S(\RR) \cap \prod_{i\in S}\bB_i\Big)
\]
such that $b \not\in \RR$. By the induction hypothesis, we have $\pi_{[n]\setminus\{i\}}(b) \in \pi_{[n]\setminus\{i\}}(\RR)$ for each $i \in [n]$.

Consider the quasiorder $\preceq$ on $[n]$ from Definition \ref{defn-relation-strand} defined by
\[
i \preceq j \;\; \iff \;\; \bB_i + \pi_{ij}(\RR) = \bB_j,
\]
and suppose without loss of generality that $[k]$ is a $\preceq$-minimal $(\RR,\bB)$-strand for some $k \le n$. If there is any tuple
\[
(a_1, ..., a_k, b_{k+1}', ..., b_n') \in \RR
\]
such that $a_i \in \bA_i\setminus \bB_i$ for $i \in [k]$ and $b_j' \in \bB_j$ for $j > k$, then we can apply the previous lemma to the situation
\[
\RR \le_{sd} \bA_1 \times \cdots \times \bA_k \times \pi_{[k+1,n]}(\RR)
\]
with
\[
\pi_{[k+1,n]}(\RR) \cap \prod_{j > k} \bB_j \; \llhd_X \; \pi_{[k+1,n]}(\RR),
\]
by compatibility with pp-formulas and Theorem \ref{thm-conservative-subdirect-absorption-free}, to finish the proof. To arrange for this situation, we consider the relation
\begin{align*}
\RR' &= \RR \cap \Big(\Big(\prod_{i\in [k]} (\bA_i\setminus \bB_i)\cup\{b_i\}\Big) \times \prod_{j>k} \bA_j \Big),\\
&= \RR \cap \Big(\Big(\{(b_1, ..., b_k)\} \cup \prod_{i\in [k]} (\bA_i\setminus \bB_i)\Big) \times \prod_{j>k} \bA_j \Big),
\end{align*}
set $\bA_i' = \pi_i(\RR')$, and set $\bB_i' = \bA_i' \cap \bB_i$ for each $i$. By the induction hypothesis applied to $\pi_{[k]\cup\{j\}}(\RR)$, we see that $\bB_j' = \bB_j$ for each $j > k$, so we have $\bB_i' \llhd_X \bA_i'$ for all $i \le n$. We will apply the induction hypothesis to $\pi_{[k+1,n]}(\RR')$, but first we need to check that $\pi_{[k+1,n]}(\RR')$ has more than one strand.

Let $S \subset [n]$ be some $(\RR,\bB)$-strand which is disjoint from $[k]$. By the assumption that $[k]$ was a $\preceq$-minimal $(\RR,\bB)$-strand, there must be some $(c_1, ..., c_n) \in \RR$ such that $c_i \not\in \bB_i$ for $i \in [k]$ and $c_j \in \bB_j$ for $j \in S$. Suppose without loss of generality that there is some $m \ge k$ such that $c_j \in \bB_j$ iff $j \in [m+1,n]$. If $m = k$, then we take $(a_1, ..., a_k, b_{k+1}', ..., b_n') = c$ to finish. Otherwise, since $c$ is an element of $\RR'$, we see that every strand of $\pi_{[k+1,n]}(\RR')$ is either contained in $[k+1,m]$ or contained in $[m+1,n]$. Thus, by the induction hypothesis we have
\[
\pi_{[k+1,n]}(\RR') \cap \prod_{j > k} \bB_j = \Big(\pi_{[k+1,m]}(\RR') \cap \prod_{j \in [k+1,m]} \bB_j \Big) \times \Big(\pi_{[m+1,n]}(\RR') \cap \prod_{j \in [m+1,n]} \bB_j \Big).
\]
Since we have $\pi_{[k+1,m]}(b) \in \pi_{[k+1,m]}(\RR')$ and $\pi_{[m+1,n]}(b) \in \pi_{[m+1,n]}(\RR')$ by the induction hypothesis applied to $\pi_{[1,m]}(\RR)$ and $\pi_{[k]\cup [m+1,n]}(\RR)$, we see that $\pi_{[k+1,n]}(b) \in \pi_{[k+1,n]}(\RR')$. Thus either $b \in \RR$, or there are some $a_i \in \bA_i \setminus \bB_i$ for $i \in [k]$ such that
\[
(a_1, ..., a_k, b_{k+1}, ..., b_n) \in \RR,
\]
which allows us to apply the previous lemma to finish the proof.
\end{proof}


\section{The algorithm for conservative CSPs}\label{s-alg-conservative}

In this section we present Barto's simple algorithm from \cite{barto-conservative-revisited}. Bulatov's algorithm from \cite{bulatov-conservative-rerevisited} is similar in spirit, but it relies on ideas from Mar\'oti's ``Tree on top of Maltsev'' algorithm \cite{tree-on-top-of-malcev} which we haven't covered yet.

The main idea of Barto's algorithm for conservative CSPs is to try to reduce to the case where all edges of the colored graphs occuring in each of the variable domains $\bA_x$ are affine. In this case, any daisy chain term will be a Mal'cev term for each variable domain, and we can solve the problem by using the algorithm for CSPs with a Mal'cev polymorphism. In Bulatov's algorithm from \cite{bulatov-conservative-rerevisited}, the main idea is to reduce to the case where there are no semilattice edges instead, in which case any daisy chain term will be a ternary generalized majority-minority operation (that a ternary generalized majority-minority operation exists in this case also follows from Theorem \ref{thm-swap-gmm}).

In order to accomplish this, we aim to show that if any semilattice or majority edge occurs in any variable domain of an instance $\fX$, then we can reduce some variable domain by solving an instance where every variable domain has been strictly decreased. We assume that our instance is $(2,3)$-minimal (or perhaps just that it is a Prague instance), and we consider the subdigraph $\cE$ of the implication digraph consisting of pairs $(x,\bB)$ such that that $\bB \le \bA_x$ is an algebra with at least one proper absorbing subalgebra. The digraph $\cE$ will be nonempty as long as any algebra $\bA_x$ has any non-affine edge.

We pick any maximal strand $\cS$ of the digraph $\cE$, and we note that the full restriction (see Definition \ref{defn-full-restriction}) of our instance $\fX$ to the strand $\cS$ is then a Prague instance in which every single domain has a proper absorbing subalgebra by Proposition \ref{prop-prague-reduction}. We can then repeatedly apply Proposition \ref{prop-prague-reduction} (or, alternatively, we could apply Kozik's \cite{pq-consistency} result from Section \ref{s-absorbing-arc-consistent}, using the fact that every Prague instance is $pq$-consistent), to find an arc-consistent absorbing reduction $\fX'$ of the full restriction of $\fX$ to the strand $\cS$, such that every variable domain in $\fX'$ is absorption free - and as a consequence, such that each variable domain in $\fX'$ is a proper subagebra of the corresponding variable domain in the original instance $\fX$. We can then apply the following consequence of the Rectangularity Theorem \ref{thm-rectangularity-conservative-absorbing}.

\begin{thm}[Barto \cite{barto-conservative-revisited}]\label{thm-rectangularity-application} Suppose that $\fX$ is a Prague instance such that each variable domain $\bA_x$ is a conservative Taylor algebra. Let $\cE$ be the subdigraph of the implication digraph which consists of pairs $(x,\bB)$ such that $\bB \le \bA_x$ and $\bB$ has a proper absorbing subalgebra, and let $\cS$ be any maximal strand of $\cE$.

Suppose that for each $(x,\bB_x) \in \cS$ we choose a minimal absorbing subalgebra $\bC_x \llhd \bB_x$, such that the system of variable domains $\bC_x$ defines an arc-consistent reduction of the full restriction of $\fX$ to the strand $\cS$. Then one of the following is true:
\begin{itemize}
\item the instance $\fX$ has no solutions with any variable $x$ assigned to any value in $\bB_x$, for any $(x,\bB_x) \in \cS$,
\item the instance $\fX'$ which we get by restricting to the variables in $\cS$ and by restricting each variable domain to the corresponding $\bC_x$ has a solution, and every solution of $\fX'$ extends to a solution of $\fX$, or
\item for some strand $\mathcal{T}$ of the subdigraph $\cC$ of the implication digraph consisting of the pairs $(x,\bC_x)$ for $x$ occuring in $\cS$, the full restriction of $\fX$ to $\mathcal{T}$ has no solutions.
\end{itemize}
\end{thm}
\begin{proof} Suppose that $a \in \prod_x \bA_x$ is a solution to the instance $\fX$ such that $a_x \in \bB_x$ for some $(x,\bB_x) \in \cS$, and suppose that for each strand $\mathcal{T}$ of the subqoset $\cC$ there is a solution $c^\mathcal{T} \in \prod_{(x,\bC_x) \in \mathcal{T}} \bC_x$ of the full restriction of $\fX$ to the strand $\mathcal{T}$. We construct a tuple $b \in \prod_x \bA_x$ by stitching these solutions together:
\[
b_x = \begin{cases} a_x & x\text{ does not occur in the strand }\cS,\\ c^{\mathcal{T}}_x & x\text{ occurs in the strand }\mathcal{T}\text{ of the subqoset }\cC.\end{cases}
\]
We claim that the tuple $b$ is also a solution to the instance $\fX$. For this, we can focus our attention on any particular constraint relation
\[
\RR \le_{sd} \bA_{x_1} \times \cdots \times \bA_{x_n}
\]
of the instance $\fX$. Suppose without loss of generality that the variables $x_1, ..., x_k$ occur in the strand $\cS$ and that the variables $x_{k+1}, ..., x_n$ do not, and suppose that $k \ge 1$. For each $i \in [n]$, define $\bB_i$ by
\[
\bB_i = \bB_{x_1} + \pi_{1i}(\RR).
\]
Since $\fX$ is a Prague instance, we have $\bB_i = \bB_{x_i}$ for all $i \le k$, and since $\cS$ is a maximal strand of $\cE$ each $\bB_j$ with $j \ge k+1$ is absorption-free. Define the relation $\RR_\bB \le \RR$ by
\[
\RR_\bB = \RR \cap \Big(\prod_i \bB_i\Big) \le_{sd} \bB_1 \times \cdots \times \bB_n,
\]
where subdirectness follows directly from the definition of the $\bB_i$s. Set $\bC_i = \bC_{x_i}$ for $i \le k$, and set $\bC_j = \bB_j$ for $j \ge k+1$, so we have $\bC_i \llhd \bB_i$ for all $i \in [n]$. Then since $\fX$ is a Prague instance, the $(\RR_\bB, \bC)$-strands are given by $[k+1,n]$ and by the intersections of the strands of $\cC$ to $\{x_1, ..., x_k\}$, and we have
\[
\RR_\bB \cap \prod_i \bC_i \ne \emptyset
\]
since the system of variable domains $\bC_x$ defines an arc-consistent reduction of the full restriction of $\fX$ to the strand $\cS$. Then the Rectangularity Theorem \ref{thm-rectangularity-conservative-absorbing} says that
\[
\RR \supseteq \prod_{T\text{ an }(\RR_\bB,\bC)\text{-strand}} \Big(\pi_T(\RR_\bB) \cap \prod_{i\in T}\bC_i\Big),
\]
so the tuple $b$ satisfies the constraint relation $\RR$.
\end{proof}

Using this result, we get the following algorithm for solving conservative CSPs.

\begin{algorithm}
\caption{Algorithm for solving an instance $\fX$ of a CSP with conservative Taylor variable domains $\bA_x$, from \cite{barto-conservative-revisited}.}\label{alg-conservative}
\begin{algorithmic}[1]
\State Run a local consistency algorithm until $\fX$ is a cycle-consistent Prague instance.
\State Let $\cE$ be the subdigraph of the implication digraph consisting of pairs $(x,\bB)$ such that $\bB \le \bA_x$ and $\bB$ has a proper absorbing subalgebra.
\If{$\cE$ is non-empty}
\State Let $\cS$ be any maximal strand of $\cE$. \Comment{Proposition \ref{prop-cycle-strand-finding}.}
\State Define $\bC_x = \bB_x$ for each $(x,\bB_x) \in \cS$.
\While{some $\bC_x$ is not absorption-free}
\State Let $\fX_\bC$ be the Prague instance we get by restricting to variables in $\cS$ and restricting each variable domain to $\bC_x$. \Comment{Proposition \ref{prop-prague-reduction}}
\State Pick any maximal strand $\mathcal{T}$ of the subqoset of the implication qoset of $\fX_\bC$ consisting of $(x,\bC')$ such that $\bC' \lhd \bC_x$ and $\bC' \ne \bC_x$.
\State Set $\bC_x \gets \bC'$ for each $(x,\bC') \in \mathcal{T}$.
\EndWhile
\State Let $\cC$ be the subqoset consisting of $(x,\bC_x)$ for $x$ occuring in the strand $\cS$.
\ForAll{strands $\mathcal{T}$ of $\cC$}
\State Let $\fX_\mathcal{T}$ be the full restriction of $\fX$ to the strand $\mathcal{T}$.
\State Solve the instance $\fX_\mathcal{T}$ recursively. \Comment{$\bC_x < \bB_x$ for all $(x,\bC_x) \in \mathcal{T}$.}
\If{$\fX_\mathcal{T}$ has no solutions}
\State Set $\bA_x \gets \bA_x \setminus \bC_x$ for each $(x,\bC_x) \in \mathcal{T}$.
\State \textbf{go to} Step 1.
\Else
\State Let $c^\mathcal{T}$ be a solution to the instance $\fX_\mathcal{T}$.
\EndIf
\EndFor
\State Set $\bA_x \gets (\bA_x \setminus \bB_x) \cup \{c^\mathcal{T}_x\}$ for all $(x,\bB_x) \in \cS$, where $\mathcal{T}$ is the strand of $\cC$ which contains $(x,\bC_x)$. \Comment{Theorem \ref{thm-rectangularity-application}}
\State \textbf{go to} Step 1.
\EndIf
\State Solve $\fX$ by using the algorithm for CSPs with a Mal'cev polymorphism. \Comment{Section \ref{s-malcev-algorithm}}
\end{algorithmic}
\end{algorithm}

\begin{thm} Algorithm \ref{alg-conservative} correctly solves every instance $\fX$ of any multisorted CSP where each variable domain is a conservative Taylor algebra. If each variable domain has size at most $k$, then Algorithm \ref{alg-conservative} runs in time $\|\fX\|^{O(k)}$, where $\|\fX\|$ is (up to a logarithmic factor) the number of bits needed to describe the instance $\fX$.
\end{thm}

\begin{cor} The CSP dichotomy conjecture is true for all CSP templates on a domain of size at most $3$.
\end{cor}
\begin{proof} By the classification of minimal Taylor algebras of size $3$ from Subsection \ref{ss-minimal-three}, every minimal Taylor algebra of size at most $3$ is either a subdirect product of conservative Taylor algebras, or has a $3$-edge term.
\end{proof}

It seems plausible that a much more careful analysis of the algorithm for conservative CSPs might show that it runs in time $\|\fX\|^{O(1)}$, regardless of the sizes of the variable domains.

\begin{prob}\label{conservative-uniform} Given as input an instance $\fX$ of any CSP together with a conservative ternary weak near-unanimity polymorphism which preserves the relations of $\fX$, can we solve the instance $\fX$ in time polynomial in $\|\fX\|$?
\end{prob}

The method we have been using to encode relations up to this point has been to explicitly list out the tuples contained in the relation. An alternate way of describing constraint relations on the domain $\{0,1\}$ via ``extension oracles'' was introduced in \cite{partial-poly-seth}, and this way of describing constraint relations seems to generalize naturally to conservative CSPs. I will use the phrase ``restriction oracle'' instead of ``extension oracle'' for the generalization I have in mind.

\begin{defn} A \emph{restriction oracle} $\mathcal{O}_R$ for a relation $R \subseteq A_1 \times \cdots \times A_n$ is defined as a black box function which takes as input a tuple of subsets $B_i \subseteq A_i$, and returns ``true'' if and only if we have
\[
R \cap (B_1 \times \cdots \times B_n) \ne \emptyset.
\]
\end{defn}

The idea behind a restriction oracle is that it is the bare minimum which is needed to be able to run the (generalized) arc-consistency algorithm. It's easy to see how we could use restriction oracle descriptions of constraint relations to establish cycle-consistency (or even singleton arc-consistency), but it is not clear if it is possible to use restriction oracles to establish $(2,3)$-minimality, or even to reduce to a subinstance which satisfies condition (P2).

\begin{ex} A concrete example of a high-arity relation which has an efficient restriction oracle is the \emph{all-different} relation $\bigwedge_{i \ne j} x_i \ne x_j$. This relation occurs naturally in Sudoku and its generalizations. For any sets $B_1, ..., B_n$, we can determine whether or not
\[
\{x \mid \forall i \ne j, x_i \ne x_j\} \cap (B_1 \times \cdots \times B_n) \ne \emptyset
\]
as follows. We start by drawing a bipartite graph with parts $A = \{1, ..., n\}$ and $B = \bigcup_i B_i$, with an edge from $i \in A$ to $b \in B$ exactly when $b \in B_i$. Then we use the standard augmenting path algorithm to find a maximum matching in this graph - if there is a matching of size $n$, then the edges of this matching can be viewed as an assignment from variables $x_i$ to values in $B_i$ which are all different.
\end{ex}

\begin{prob} Consider the problem where we are given an instance $\fX$ of a CSP together with a conservative ternary weak near-unanimity polymorphism which is promised to preserve the relations of $\fX$, but instead of having explicit descriptions of the constraint relations, the constraint relations are given to us implicitly in terms of restriction oracles. Is there an algorithm which determines whether $\fX$ has a solution and makes only polynomially many calls to the restriction oracles which describe the constraint relations?
\end{prob}

When we leave the context of conservative CSPs, restriction oracles become a less natural concept. The trouble is that it's only natural to call the restriction oracle when the sets $B_i$ are subalgebras of the variable domains. The next example shows how this can become an issue for the algebra $\ZZ/3^{\aff}$.

\begin{ex} We can efficiently describe high-arity relations $R$ on $\ZZ/3^{\aff}$ by writing down systems of linear equations. If we could convert a description of $R$ as the solution set of a system of linear equations into an efficient restriction oracle $\mathcal{O}_R$, however, then we would be able to solve 1-IN-3 SAT. To see this, note that for $x,y,z \in \ZZ/3^{\aff}$ we have
\[
(x,y,z) \in \{(0,0,1),(0,1,0),(1,0,0)\} \;\; \iff \;\; x,y,z \in \{0,1\} \; \wedge \; x+y+z \equiv 1 \pmod{3}.
\]
\end{ex}

\section{The meta-problem for conservative CSP templates}

In this section we will go over Carbonnel's solution to the meta-problem for conservative CSPs from their thesis \cite{carbonnel-thesis}. Recall that in the meta-problem, we are given a CSP template as a relational structure $\fA = (A, \Gamma)$ (which we usually assume to be a core), and we wish to either prove that $\CSP(\fA)$ is NP-complete or to find a Taylor polymorphism of $\fA$, in time polynomial in the total size $\|\fA\|$ of the description of $\fA$, which we define as
\[
\|\fA\| \coloneqq \sum_{R \in \Gamma} \arity(R)|R|.
\]
In the meta-problem for conservative CSPs, we restrict our attention to CSP templates where $\Gamma$ contains the unary relation $A\setminus \{a\}$ for each $a \in A$ (in particular, any such $\fA$ is automatically a rigid core). Note that by our classification of conservative minimal Taylor algebras, if a conservative CSP template $\fA$ has a Taylor polymorphism, then it has a ternary weak near-unanimity polymorphism, i.e. an idempotent ternary polymorphism $w$ satisfying the identities
\[
w(x,x,y) \approx w(x,y,x) \approx w(y,x,x).
\]
Furthermore, we can consider the case of $3$-conservative CSP templates without any additional difficulty, since any $3$-conservative Taylor algebra has a conservative Taylor reduct (by a $3$-conservative template, we mean a CSP template such that $\Gamma$ contains every unary relation of size at most $3$). These facts were not known at the time that Carbonnel wrote their thesis, and by using them we can make Carbonnel's algorithm more concrete.

In Carbonnel's thesis \cite{carbonnel-thesis}, the strategy for solving the meta-problem was described as being similar to a treasure hunt (or perhaps a puzzle hunt): we have a sequence of locked boxes, and a single key which opens the first box, such that each box contains the key to opening the next box. Here, the key is a metaphor for the Taylor polymorphism - once we know a Taylor polymorphism, we can use it to solve instances of our CSP. More specifically, the key is a metaphor for a \emph{partial description} of a Taylor polymorphism. In order to see how a partial description of a Taylor polymorphism can make sense, we first describe a useful rephrasing of the meta-problem for a $3$-conservative template in terms of a meta-problem for a multisorted CSP template.

\begin{defn} Suppose that $\fA = (A,\Gamma)$ is a $3$-conservative CSP template. We define the \emph{associated multisorted template} $\fA_3$ to have a sort for each subset of $A$ of size at most $3$, with two types of relations:
\begin{itemize}
\item for each relation $R \in \Gamma$ of arity $m$, and for every triple of elements $u,v,w \in \RR$ (not necessarily distinct), we have a multisorted relation
\[
R \cap (\{u_1,v_1,w_1\} \times \cdots \times \{u_n,v_n,w_n\}) \subseteq \{u_1,v_1,w_1\} \times \cdots \times \{u_n,v_n,w_n\},
\]
\item for every $a,b,c \in A$ (not necessarily distinct), we have the binary inclusion relation
\[
\{(a,a), (b,b)\} \subseteq \{a,b\} \times \{a,b,c\}.
\]
\end{itemize}
\end{defn}

Note that the size of the associated multisorted template $\fA_3$ is bounded by
\[
\|\fA_3\| \le \sum_{R \in \Gamma} \sum_{u,v,w \in R} \arity(R)|R| + \sum_{a,b,c \in A} 2|\{a,b\}| \le \|\fA\|^4 + 4|A|^3.
\]
A polymorphism of a multisorted relational structure is defined to be an operation with a different interpretation on each sort of the structure, such that applying the operation componentwise preserves each multisorted relation.

\begin{prop}\label{prop-conservative-multi-ternary} There is a bijection between conservative ternary polymorphisms of $\fA$ and ternary polymorphisms of the associated multisorted structure $\fA_3$ which preserves height $1$ identities.
\end{prop}
\begin{proof} There is an obvious way to convert any conservative ternary polymorphism of $\fA$ into a ternary polymorphism of $\fA_3$. Conversely, any ternary polymorphism $f$ of $\fA_3$ can be stitched together into a (necessarily conservative) ternary polymorphism $\tilde{f}$ of $\fA$: the fact that $f$ preserves the binary inclusion relations guarantees that the values of $\tilde{f}$ are well-defined on two-element sets, and for every $u,v,w \in R \in \Gamma$ the fact that $f$ preserves the multisorted relation
\[
R \cap (\{u_1,v_1,w_1\} \times \cdots \times \{u_n,v_n,w_n\})
\]
guarantees that $\tilde{f}(u,v,w) \in R$.
\end{proof}

So we have reduced the problem of determining whether a given relational structure $\fA$ has a conservative Taylor polymorphism to the problem of determining whether the associated multisorted relational structure $\fA_3$ has a ternary weak near-unanimity polymorphism. Now we have a way to understand what a partial description of a ternary weak near-unanimity polymorphism on $\fA_3$ should be.

\begin{defn} Suppose that $\fA = (A,\Gamma)$ is a $3$-conservative CSP template, and let $\mathcal{U} \subseteq \mathcal{P}(A)$ be a collection nonempty subsets of $A$ each of size at most $3$ which is closed under taking nonempty subsets. We define the multisorted template $\fA_\mathcal{U}$ to be the template whose sorts are exactly the elements of $\mathcal{U}$, such that for each $m$-ary relation $R$ of $\fA_3$ with at least one coordinate of a sort in $\mathcal{U}$, if $S \subseteq [m]$ is the set of coordinates of $R$ which have a sort in $\mathcal{U}$, then the relation $\pi_S(R)$ is a relation of $\fA_\cU$.
\end{defn}

Later we will want to fix certain coordinates of various relations of $\fA_3$ to have certain values, before projecting to the coordinates with sorts in $\mathcal{U}$. To reassure ourselves that this will not cause unexpected problems, we have the following result.

\begin{prop} Suppose that $\fA = (A,\Gamma)$ is a $3$-conservative CSP template. Let $\mathcal{U} \subseteq \cV \subseteq \mathcal{P}(A)$ be collections of nonempty subsets of $A$ each of size at most $3$ which are closed under taking nonempty subsets. Suppose that $f$ is a ternary polymorphism of the multisorted structure $\fA_\cU$, and that $R$ is an $m$-ary relation of $\fA_{\cV}$. Let $S \subseteq [m]$ be the set of coordinates of the relation $R$ with sorts from $\mathcal{U}$. Then for any $y \in \pi_{[m]\setminus S}(R)$, the relation
\[
R_y \coloneqq \{x \mid \exists r \in R\text{ s.t. }\pi_S(r) = x, \pi_{[m]\setminus S}(r) = y\} \subseteq \pi_S(R)
\]
is preserved by the polymorphism $f$.
\end{prop}
\begin{proof} We just need to check that for every three tuples $u,v,w \in R_y$, we have $f(u,v,w) \in R_y$. If $R$ is one of the binary inclusion relations, then this follows from the fact that $f$ is conservative on the sorts in $\cU$ (which folows from the fact that $\cU$ is closed under taking nonempty subsets). Otherwise, $R$ originally came from some relation $\tilde{R} \in \Gamma$. Then there are lifts $\tilde{u}, \tilde{v}, \tilde{w} \in \tilde{R}$ such that projecting to the coordinates of $\tilde{R}$ corresponding to $S$ gives us $u,v,w$, and such that projecting to the coordinates of $\tilde{R}$ corresponding to $[m]\setminus S$ gives us $y$ in each case. Then there is a relation $R'$ of $\fA_\cU$ given by restricting the $i$th coordinate of $\tilde{R}$ to the set $\{\tilde{u}_i, \tilde{v}_i, \tilde{w}_i\}$ for all $i$ and projecting to the coordinates corresponding to $S$, and we have
\[
f(u,v,w) \in R' \subseteq R_y
\]
since $f$ preserves $R'$.
\end{proof}

In order to find ternary weak near-unanimity polymorphisms of the structure $\fA_\cU$, we use the idea of solving an \emph{indicator instance} (we used this idea once already to solve the meta-problem for bounded width templates - this idea appears to have shown up for the first time in \cite{indicator-instance}).

\begin{prop}\label{prop-conservative-indicator} If $\fB = (\cU, \Gamma)$ is a multisorted relational structure with sorts $\cU$, then $\fB$ has a ternary weak near-unanimity polymorphism $f$ iff the following instance of $\CSP(\fB)$ has a solution:
\begin{align*}
\exists f=\{f_U \in U^{U^3}\}_{U \in \cU}\text{ s.t. } &\bigwedge_{U \in \cU} \Big(\bigwedge_{a \in U} f_U(a,a,a) \in \{a\} \wedge \bigwedge_{a,b \in U} f_U(a,a,b) = f_U(a,b,a) = f_U(b,a,a)\Big)\\
&\wedge \bigwedge_{R \in \Gamma} \bigwedge_{u,v,w \in R} f(u,v,w) \in R.
\end{align*}
\end{prop}

Now we can describe the treasure hunt algorithm: we iteratively build ternary weak near-unanimity polymorphisms $f_\cU$ of the multisorted structures $\fA_\cU$ for progressively larger collections $\cU$ of subsets of $\bA$ of size at most $3$. In each step, we add a single new subset $V$ of $A$ to $\cU$ to produce a larger collection $\cV$.

In order to solve the indicator instance for $\fA_\cV$, we brute force over all possible ternary weak near-unanimity operations on $V$, and for each one, we check if it extends to a solution to the indicator instance, using the fact that all of the remaining variables have sorts in $\cU$. In fact, we don't need to consider \emph{all} weak near-unanimity operations on $V$ - we only need to check the $73$ specific operations from the classification of conservative minimal Taylor algebras of size $3$ (if $V$ has size $2$, we only need to consider $4$ possible operations).

\begin{algorithm}
\caption{Treasure hunt algorithm for solving the meta-problem for a $3$-conservative relational structure $\fA = (A, \Gamma)$, from \cite{carbonnel-thesis}.}\label{alg-conservative-treasure-hunt}
\begin{algorithmic}[1]
\State Set $n \gets \binom{|A|}{3} + \binom{|A|}{2}$.
\State Pick a sequence $\cU_0 \subseteq \cU_1 \subseteq \cU_2 \subseteq \cdots \subseteq \cU_{n} = \{U \subseteq A \mid |U| \in [3]\}$ such that $\cU_0$ is the set of singleton subsets of $A$, each $\cU_{i+1}$ contains exactly one more subset of $A$ then $\cU_i$, and each $\cU_i$ is closed under taking nonempty subsets.
\State Let $f_{\cU_0}$ be the unique ternary polymorphism of $\fA_{\cU_0}$, given by $f_{\{a\}}(a,a,a) = a$ for all $a \in A$.
\For{$i \in [n]$}
\State Let $\fX$ be the indicator instance from Proposition \ref{prop-conservative-indicator} for the multisorted structure $\fA_{\cU_i}$.
\State Let $V$ be the new set in $\cU_i\setminus \cU_{i-1}$.
\ForAll{ternary weak near-unanimity operations $g_V$ on $V$}
\State Let $\fX'$ be the instance we get from $\fX$ by replacing each variable $f_V(a,b,c)$ of $\fX$ with sort $V$ by the constant $g_V(a,b,c)$.
\State Solve the instance $\fX'$ by using Algorithm \ref{alg-conservative} with the Taylor operation $f_{\cU_{i-1}}$.
\If{$\fX'$ has a solution}
\State Let $f_{\cU_i}$ be any solution to $\fX'$.
\EndIf
\EndFor
\If{$f_{\cU_i}$ hasn't been defined}
\State \Return ``NP-complete''.
\EndIf
\EndFor
\State Stitch $f_{\cU_n}$ into a ternary polymorphism $f$ of $\fA$ using Proposition \ref{prop-conservative-multi-ternary}.
\State \Return $f$.
\end{algorithmic}
\end{algorithm}

\begin{thm} If $\fA$ is a $3$-conservative relational structure, then Algorithm \ref{alg-conservative-treasure-hunt} runs in time polynomial in $\|\fA\|$, and either correctly determines that $\CSP(\fA)$ is NP-complete or produces a ternary weak near-unanimity polymorphism of $\fA$.
\end{thm}
\begin{proof} Given what we have already proved, the only thing left to check is that the algorithm runs in polynomial time. The only step which looks dangerous is the step where we apply Algorithm \ref{alg-conservative}, since the degree of the polynomial in the running time of that algorithm depends on the size of the largest sort which shows up as a variable domain. However, we only ever apply Algorithm \ref{alg-conservative} to multisorted structures where every sort has size at most $3$.
\end{proof}

\begin{rem}\label{rem-conservative-multisorted} Algorithm \ref{alg-conservative-treasure-hunt} easily generalizes to multisorted CSP templates. Since an arbitrary instance of an unstructured CSP can be thought of as an instance of a multisorted CSP where each variable has a different sort, with only those relations which actually show up in the instance, we can efficiently check whether there is any possible way to interpret each variable domain as a conservative Taylor algebra such that the relations are compatible with the algebraic structure. If we can impose such an algebraic structure, we can apply Algorithm \ref{alg-conservative} as long as the variable domains are not too large. If we can solve Problem \ref{conservative-uniform}, then we may not even need the restriction on the sizes of the variable domains!

\end{rem}

\section{Mar\'oti's ``Maltsev on top'': combining few subpowers and bounded width}

This section is a bit of a breather: we will go over Mar\'oti's algorithm from \cite{malcev-on-top}, which gives us a simple way to combine the few subpowers algorithm with bounded width reasoning that works in the case where we can mostly disentangle the few subpowers parts of our algebras from the bounded width parts. The nice thing about this algorithm is that we don't need any new algebraic ideas to understand it, and it gives us a hint that there may be a simple algorithm which can handle the case of general Taylor algebras.

The specific assumption that Mar\'oti uses in \cite{malcev-on-top} is that each variable domain $\bA_x$ is an idempotent algebra with a special congruence $\theta_x \in \Con(\bA_x)$ such that:
\begin{itemize}
\item $\bA_x/\theta_x$ has few subpowers, and
\item the congruence classes $a_x/\theta_x \le \bA_x$ have bounded width, for all $a_x \in \bA_x$.
\end{itemize}
We will say that such a pair $\bA_x,\theta_x$ has the \emph{Mal'cev-on-top} property to abbreviate this. Note that by Corollary \ref{cor-few-subpowers-product} and Corollary \ref{affine-free-variety}, this assumption implies that the multisorted CSP with variable domains $\bA_x/\theta_x$ has few subpowers, and that for any $a_x \in \bA_x$ the multisorted CSP with variable domains $a_x/\theta_x$ has bounded width (equivalently, the products $\bA = \prod_x \bA_x$ and $\theta = \prod_x \theta_x$ have the Mal'cev-on-top property).

In the few subpowers setting it's natural to consider global constraints (i.e. constraints which involve all of the variables at once), specified by their compact representations, but in the bounded width setting it is more natural to consider constraints with low arity. To get the best of both worlds, Mar\'oti introduced a new type of constraint, which he called a ``Maltsev constraint'' (we will call them ``M-constraints'').

\begin{defn}[Mar\'oti \cite{malcev-on-top}] If $X$ is the set of all variables, and if each variable domain $\bA_x$ has a congruence $\theta_x$ with the Mal'cev-on-top property, then an \emph{M-constraint} on $I \subseteq X$ is a subalgebra
\[
\RR \le \prod_{x \in I} \bA_x \times \prod_{y \in X\setminus I} \bA_y/\theta_y.
\]
We will say that the \emph{M-arity} of an M-constraint on $I$ is $|I|$. By a \emph{compact representation} of $\RR$, we mean a collection of compact representations of the relations
\[
a_I + \RR = \big\{b \in \prod_{y \in X\setminus I} \bA_y/\theta_y \mid (a_I,b) \in \RR\big\},
\]
for all $a_I \in \prod_{x \in I} \bA_x$.
\end{defn}

The way we will actually think about M-constraints on $I$ is reversed from how their compact representations are defined: we think of an M-constraint $\RR$ as a way of associating an $|I|$-ary relation $b - \RR \le \prod_{x \in I} \bA_x$ to every tuple of congruence classes $b \in \prod_{y \in X\setminus I} \bA_y/\theta_y$.

We might attempt to be a little bit more ambitious: let's consider the more general case of relations described by oracles which are just barely powerful enough to run the arc-consistency algorithm.

\begin{defn} A \emph{subalgebra restriction oracle} $\mathcal{O}_\RR$ for a relation $\RR \le \bA_1 \times \cdots \times \bA_n$ is defined as a black box function which takes as input a tuple of subalgebras $\bB_i \le \bA_i$, and returns ``true'' if and only if we have
\[
\RR \cap (\bB_1 \times \cdots \times \bB_n) \ne \emptyset.
\]
\end{defn}

\begin{prop} Let $\RR \le \prod_x \bA_x$, and for $I \subseteq X$ consider the M-relation on $I$ given by
\[
\RR/\theta_{X\setminus I} \le \bA_I \times \bA_{X \setminus I}/\theta_{X\setminus I},
\]
where we have made the abbreviations $\bA_I = \prod_{x \in I} \bA_x$, and similarly for $\bA_{X \setminus I}$ and $\theta_{X \setminus I}$. If $\bA_{X \setminus I}/\theta_{X\setminus I}$ has a $k$-edge term, and if we are given a subalgebra restriction oracle $\mathcal{O}_\RR$ for $\RR$, then we can compute a compact representation for $\RR/\theta_{X\setminus I}$ using at most
\[
\prod_{x \in I} |\bA_x| \times \Big(\sum_{y \in X \setminus I} |\bA_y/\theta_y|\Big)^k
\]
calls to the oracle $\mathcal{O}_\RR$.
\end{prop}
\begin{proof} Left as an exercise for now.
\end{proof}

Unfortunately, if we don't restrict ourselves to M-constraints (of known M-arity) then subalgebra restriction oracles seem to give us slightly too little information to run Mar\'oti's algorithm. So we make the following definition instead.

\begin{defn} Suppose that $\bA_i, \theta_i$ have the Mal'cev-on-top property. An \emph{M-restriction oracle} $\cO^\cM_\RR$ for a relation $\RR \le \bA_1 \times \cdots \times \bA_n$ is defined as a black box function which takes as input a tuple of compact representations of M-unary M-constraints
\[
\bB_i \le \bA_i \times \prod_{j \ne i} \bA_j/\theta_j,
\]
and returns a compact representation of
\[
\big(\RR \cap \bB_1 \cap \cdots \cap \bB_n\big)/\prod_i \theta_i.
\]
\end{defn}

M-restriction oracles are clearly at least as powerful as subalgebra restriction oracles. Furthermore, if we want to have any hope of solving CSPs involving $\RR$ along with M-unary M-constraints, then there had better be an efficient way to implement the M-restriction oracle $\cO^\cM_\RR$. In particular, if Mar\'oti's approach is going to work at all then there ought to be a way to efficiently implement the M-restriction oracle $\cO^\cM_\RR$ when $\RR$ happens to be an M-constraint!

We start with the case where $\RR = \bA_1 \times \cdots \times \bA_n$.

\begin{prop} Suppose that $\theta_i \in \Con(\bA_i)$ for all $i$, and that for each $i$ we have a relation
\[
\bB_i \le \bA_i \times \prod_{j \ne i} \bA_j/\theta_j.
\]
Then we have
\[
\big(\bB_1 \cap \cdots \cap \bB_n\big)/\prod_i \theta_i = (\bB_1/\theta_1) \cap \cdots \cap (\bB_n/\theta_n).
\]
\end{prop}
\begin{proof} It's clear that the left hand side is contained in the right hand side. For the other containment, let $\bar{b}$ be any tuple of congruence classes in $\bigcap_i \bB_i/\theta_i$. For each $i$, since $\bar{b} \in \bB_i/\theta_i$, there is some $b_i \in \bA_i$ such that $b_i/\theta_i = \pi_i(\bar{b})$ and
\[
(b_i, \pi_{[n]\setminus\{i\}}\bar{b}) \in \bB_i.
\]
Taking $b = (b_1, ..., b_n)$, we have $b/\prod_{j \ne i} \theta_j = (b_i, \pi_{[n]\setminus\{i\}}\bar{b}) \in \bB_i$ for each $i$, so
\[
b \in \bB_1 \cap \cdots \cap \bB_n,
\]
and $b / \prod_i \theta_i = \bar{b}$.
\end{proof}

Combining this with Algorithm \ref{alg-intersect} for intersecting compact representations of relations in algebras with few subpowers, we get the following result.

\begin{prop} Suppose that we are given a $k$-edge term $e$ for the algebras $\bA_x/\theta_x$, along with a compact representation of an M-constraint
\[
\RR \le \bA_I \times \bA_{X\setminus I}/\theta_{X \setminus I}
\]
and a tuple of compact representations of M-unary M-constraints
\[
\bB_x \le \bA_x \times \bA_{X\setminus \{x\}}/\theta_{X\setminus\{x\}}.
\]
Then we can compute a compact representation of
\[
\big(\RR \cap \bigcap_{x \in X} \bB_x\big)/\theta_X
\]
in time polynomial in the sum of the sizes of the compact representations of $\RR$ and the $\bB_x$s.
\end{prop}

\begin{prob} Given black box access to the subalgebra restriction oracle $\cO_\RR$ for an unknown relation $\RR \le \bA_1 \times \cdots \times \bA_n$ of high arity, not necessarily an M-constraint, together with a $k$-edge term $e$ for the algebras $\bA_i/\theta_i$, can we implement the M-restriction oracle $\cO^\cM_\RR$ using only polynomially many calls to $\cO_\RR$?
\end{prob}

Let's table this question for now, and move on to describing Mar\'oti's algorithm, assuming that we have been provided with M-restriction oracles for each constraint relation. We start with an easier case than the general case: let's assume that for each $a_x \in \bA_x$, the congruence class $a_x/\theta_x$ has totally symmetric polymorphisms of all arities, so that $\CSP(a_x/\theta_x)$ is solved by arc-consistency. Note that a minor variation of the proof of Proposition \ref{prop-totally-symmetric-preparation} implies that in this case, the multisorted problem $\CSP(\{a_x/\theta_x\})$ is also solved by arc-consistency.

\begin{defn} Suppose that $\bA_i, \theta_i$ have the Mal'cev-on-top property. We say that an instance of $\CSP(\{\bA_i\})$ is \emph{M-arc-consistent} if for each constraint relation $\RR_i \le \bA_X = \prod_x \bA_x$ and each tuple of congruences
\[
a/\theta_X \in \RR_i/\theta_X,
\]
the instance we get by restricting each variable domain $\bA_x$ to the congruence class $a_x/\theta_x$ is arc-consistent.
\end{defn}

The idea is to think of our instance of $\CSP(\{\bA_i\})$ as a giant bundle of restricted instances of $\CSP(\{a_x/\theta_x\})$, one for each tuple $a \in \bA_X$, and to try to run the arc-consistency algorithm on each restricted instance simultaneously. The M-restriction oracles give us the exact tool needed to accomplish this.

\begin{algorithm}
\caption{M-arc-consistency algorithm, given $\theta_i \in \Con(\bA_i)$ such that $\bA_i, \theta_i$ have the Mal'cev-on-top property and given M-restriction oracles $\cO^\cM_{\RR_i}$ for each constraint relation $\RR_i \le \bA_X = \prod_x \bA_x$.}\label{alg-M-arc-consistency}
\begin{algorithmic}[1]
\State Set $\bB_x \gets \bA_x \times \bA_{X\setminus\{x\}}/\theta_{X\setminus\{x\}}$ for each variable $x$. \Comment{$\bB_x$ is an M-unary M-constraint on $\{x\}$.}
\State Compute a compact representation for $\bB_x$ by sending each $a_x \in \bA_x$ to a compact representation of $\bA_{X\setminus\{x\}} / \theta_{X\setminus\{x\}}$.
\State For any set $S \subseteq X$, we use $\bB_S$ as an abbreviation for $\bigcap_{x\in S} \bB_x \le \bA_S \times \bA_{X\setminus S}/\theta_{X\setminus S}$.
\Repeat
\ForAll{constraint relations $\RR_i \le \bA_X$ and choices of $x \in X$}
\ForAll{$a \in \bA_x$}
\State Set $\bC_a \gets (\{a\} \times \bA_{X\setminus\{x\}}/\theta_{X\setminus\{x\}}) \cap \bB_x$. \Comment{$\bC_a$ is an M-unary M-constraint on $\{x\}$.}
\State Compute a compact representation for $\bC_a$ by sending $a$ to the known compact representation of $\pi_{X\setminus\{x\}}(a + \bB_x)$, and sending every $a' \ne a$ to the empty set.
\State Set $\bD_{a} \gets (\RR_i \cap \bB_{X\setminus \{x\}} \cap \bC_a) / \theta_X$.
\State Call the M-restriction oracle $\cO^\cM_{\RR_i}$ to compute a compact representation of $\bD_a$.
\EndFor
\State Set $\bB_x \gets \bigcup_{a \in \bA_x} \{a\} \times \pi_{X\setminus \{x\}}(\bD_a) = (\RR_i \cap \bB_X) / \theta_{X\setminus \{x\}}$.
\State Update the compact representation for $\bB_x$ by sending each $a \in \bA_x$ to the known compact representation of $\pi_{X\setminus \{x\}}(\bD_a)$.
\EndFor
\Until{no $\bB_x$ changes.}
\State If any $\bB_x = \emptyset$, there is no solution.
\State \Return a compact representation of $\bB_X/\theta_X$. \Comment{Each $\bB_x/\theta_x$ is equal to $\bB_X/\theta_X$.}
\end{algorithmic}
\end{algorithm}

\begin{prop} Suppose that $\bA_i, \theta_i$ have the Mal'cev-on-top property, and that $\fX$ is an instance of $\CSP(\{\bA_i\})$. Then Algorithm \ref{alg-M-arc-consistency} runs in polynomial time, and at the end of Algorithm \ref{alg-M-arc-consistency}, if we define a new instance $\fX'$ by restricting each variable domain to
\[
\bA_x' \coloneqq \pi_x(\bB_x)
\]
and restricting each constraint relation $\RR_i \le \bA_X$ of $\fX$ to
\[
\RR_i' \coloneqq \RR_i \cap \bB_X \le \bA_X',
\]
then the instance $\fX'$ is M-arc-consistent. Furthermore, for any tuple of congruence classes $a/\theta_X \in \bA_X/\theta_X$ we have
\[
a/\theta_X \in \bB_X/\theta_X \;\; \iff \;\; \text{running arc-consistency on } \fX_{a/\theta_X} \text{ doesn't prove that no solution exists},
\]
where $\fX_{a/\theta_X}$ is the instance we get by restricting each variable domain $\bA_x$ to the congruence class $a_x/\theta_x$.
\end{prop}
\begin{proof} Left as another exercise for now.
\end{proof}

To generalize from the case where each $a_x/\theta_x$ has width 1 to bounded width, we will use the fact that finite-domain CSPs have bounded width iff they are solved by \emph{singleton} arc-consistency (since singleton arc-consistency implies cycle consistency).

\begin{defn} Suppose that $\bA_i, \theta_i$ have the Mal'cev-on-top property. We say that an instance of $\CSP(\{\bA_i\})$ is \emph{singleton M-arc-consistent} if for each variable $x \in X$, each value $a_x \in \bA_x$, each constraint relation $\RR_i \le \bA_X$, and each tuple of congruences
\[
a_{X\setminus\{x\}}/\theta_{X\setminus\{x\}} \in (a_x + \RR_i)/\theta_{X\setminus\{x\}},
\]
running arc-consistency on the restricted instance we get by replacing $\bA_x$ by $\{a_x\}$ and replacing each $\bA_y$ with $a_y/\theta_y$ for $y \ne x$ doesn't prove that the restricted instance has no solution.
\end{defn}

The algorithm for singleton M-arc-consistency is nearly the same as the algorithm for M-arc-consistency, just with calls to the M-restriction oracles replaced by calls to the M-arc-consistency algorithm.

\begin{algorithm}
\caption{Singleton M-arc-consistency algorithm, given $\theta_i \in \Con(\bA_i)$ such that $\bA_i, \theta_i$ have the Mal'cev-on-top property and given M-restriction oracles $\cO^\cM_{\RR_i}$ for each constraint relation $\RR_i \le \bA_X = \prod_x \bA_x$.}\label{alg-singleton-M-arc-consistency}
\begin{algorithmic}[1]
\State Set $\bB_x \gets \bA_x \times \bA_{X\setminus\{x\}}/\theta_{X\setminus\{x\}}$ for each variable $x$. \Comment{$\bB_x$ is an M-unary M-constraint on $\{x\}$.}
\State Compute a compact representation for $\bB_x$ by sending each $a_x \in \bA_x$ to a compact representation of $\bA_{X\setminus\{x\}} / \theta_{X\setminus\{x\}}$.
\State For any set $S \subseteq X$, we use $\bB_S$ as an abbreviation for $\bigcap_{x\in S} \bB_x \le \bA_S \times \bA_{X\setminus S}/\theta_{X\setminus S}$.
\Repeat
\ForAll{variables $x \in X$}
\ForAll{$a \in \bA_x$}
\State Set $\bC_a \gets (\{a\} \times \bA_{X\setminus\{x\}}/\theta_{X\setminus\{x\}}) \cap \bB_x$. \Comment{$\bC_a$ is an M-unary M-constraint on $\{x\}$.}
\State Compute a compact representation for $\bC_a$ by sending $a$ to the known compact representation of $\pi_{X\setminus\{x\}}(\bB_x)$, and sending every $a' \ne a$ to the empty set.
\State Set $\RR_i' \gets \RR_i \cap \bB_{X\setminus \{x\}} \cap \bC_a$. \Comment{We can easily convert an M-restriction oracle for $\RR_i$ to an M-restriction oracle for $\RR_i'$.}
\State Let $\fX_a'$ be the restricted instance with each constraint relation $\RR_i$ replaced by $\RR_i'$.
\State Let $\bD_a \le \bA_X/\theta_X$ be the output of the M-arc-consistency algorithm (Algorithm \ref{alg-M-arc-consistency}) applied to the restricted instance $\fX_a'$, and remember the compact representation of $\bD_a$.
\EndFor
\State Set $\bB_x \gets \bigcup_{a \in \bA_x} \{a\} \times \pi_{X\setminus \{x\}}(\bD_a)$.
\State Update the compact representation for $\bB_x$ by sending each $a \in \bA_x$ to the known compact representation of $\pi_{X\setminus \{x\}}(\bD_a)$.
\EndFor
\Until{no $\bB_x$ changes.}
\State If any $\bB_x = \emptyset$, there is no solution.
\State Pick any tuple of congruences $a/\theta_X \in \bB_X/\theta_X$.
\State Let $\fX_{a/\theta_X}$ be the restricted instance we get by replacing each variable domain $\bA_x$ by the congruence class $a_x/\theta_x$.
\State Run your favorite bounded width algorithm to find a solution to $\fX_{a/\theta_X}$.
\end{algorithmic}
\end{algorithm}

\begin{prop} Suppose that $\bA_i, \theta_i$ have the Mal'cev-on-top property, and that $\fX$ is an instance of $\CSP(\{\bA_i\})$. Then Algorithm \ref{alg-singleton-M-arc-consistency} runs in polynomial time, and if it doesn't prove that there is no solution to $\fX$, then it will find a solution to $\fX$.
\end{prop}
\begin{proof} Another exercise. This time, we need to be careful to check that $\bB_x$ is actually a subalgebra of $\bA_X$ at each step of the algorithm, which follows from a description of singleton arc-consistency in terms of solutions to an instance defined by starting from a tree-like universal cover and identifying all of the covers of a particular variable to each other.
\end{proof}

\section{Mar\'oti's reduction: taking advantage of semilattice quotients}

Consider classifying Taylor algebras by which types of Bulatov's edges (i.e., semilattice, majority, and affine) occur in them. If only majority and semilattice edges occur, then our algebra has bounded width by Theorem \ref{affine-free-pq}, and we can handle it. If only majority and affine edges occur, then our algebra has few subpowers by Theorem \ref{thm-cube-term-blockers}, the argument of Corollary \ref{cor-cyclic-cube-term}, and Proposition \ref{prop-weak-semilattice-edge}, and we can handle that case too. How about the case where we only have affine and semilattice edges?

To make the problem simpler, we can start by considering the case where we have an algebra $\bA$ with a congruence $\theta$ such that $\bA/\theta$ only has edges of one type, while the congruence classes of $\theta$ only have edges of the other type. If $\bA/\theta$ is affine and the congruence classes of $\theta$ are semilattices (or more generally, ancestral algebras), then Mar\'oti's Mal'cev-on-top algorithm from the previous section handles the problem.

To make progress, we need to consider the case where $\bA/\theta$ is a semilattice and the congruence classes of $\theta$ are affine - so far, we have no algorithms to handle this case (unless $\bA$ happens to be conservative). A concrete example of such an algebra is the $4$-element algebra $\bA = (\{0,1,2,*\},\cdot)$ with $\cdot$ given by
\begin{center}
\begin{tabular}{c|cccc} $\cdot$ & $0$ & $1$ & $2$ & $*$\\ \hline $0$ & $0$ & $2$ & $1$ & $*$\\ $1$ & $2$ & $1$ & $0$ & $*$\\ $2$ & $1$ & $0$ & $2$ & $*$\\ $*$ & $*$ & $*$ & $*$ & $*$\end{tabular},
\end{center}
which showed up in Subsection \ref{ss-commutators} as an example of a Taylor algebra where the term condition commutator has
\[
0_\bA = [\theta, 1_\bA] \ne [1_\bA, \theta] = \theta.
\]

Intuitively, however, the case where $\bA/\theta$ is a semilattice should be the \emph{easy} case! The absorbing element of $\bA/\theta$ corresponds to a binary absorbing subalgebra $\bB \lhd_{bin} \bA$, and in ``many'' cases, we can assume that variables with domain equal to $\bA$ ``might as well'' have their domain replaced by $\bB$, unless there is some ``clear reason'' for the variable to not be contained in $\bB$ (as an exercise, use this intuition to find an algorithm which solves $\CSP(\bA)$ for the four-element algebra $\bA$ above). Binary absorption feels in many ways like having a unary operation which crushes $\bA$ into $\bB$, but which is somehow only available in certain contexts.

Mar\'oti's reduction from \cite{tree-on-top-of-malcev} makes the vague intuition above precise, by actually \emph{building} a bespoke unary operation out of a binary polymorphism satisfying the identity $t(x,t(x,y)) \approx t(x,y)$. The catch is that what he builds is not a polymorphism of our original relational structure, but is instead a polymorphism of a related multisorted relational structure which depends on the instance we are trying to solve.

Once we start thinking of CSPs in a multisorted context, we find that there is no reason to avoid giving each variable of an instance its own private variable domain, which is shared by no other variable (this idea came up previously in Remark \ref{rem-conservative-multisorted}). If the (possibly multisorted) relational structure $\fA$ was the original template for our CSP, then by giving each variable of an instance $\fX$ its own private variable domain we end up taking each sort in $\fA$ and copying it once for each variable of $\fX$ with that sort, and similarly for the relations. This ends up producing a new multisorted relational structure, with a new signature, for each instance $\fX$. Let's make up some notation for it.

\begin{defn} If $\fX, \fA$ are (possibly multisorted) relational structures of the same signature, then we define a new multisorted relational structure
\[
\coprod_\fX \fA
\]
with a new signature (depending on $\fX$) as follows:
\begin{itemize}
\item for each sort $X_i$ of $\fX$, and for each element $x \in X_i$, there is a sort $A_x$ of $\coprod_\fX \fA$, which is a copy of the corresponding sort $A_i$ of $\fA$, and
\item for each relation
\[
C_j \subseteq X_{i_1} \times \cdots \times X_{i_m}
\]
of $\fX$, and for each tuple $c = (x_1, ..., x_m) \in C_j$, there is a relation
\[
R_c \subseteq A_{x_1} \times \cdots \times A_{x_m}
\]
of $\coprod_\fX \fA$ which is a copy of the corresponding relation $R_j \subseteq A_{i_1} \times \cdots \times A_{i_m}$ of $\fA$.
\end{itemize}
\end{defn}

\begin{rem} Readers familiar with category theory may find it useful to relate $\coprod_\fX \fA$ to a category-theoretic construction known as the \emph{Grothendieck construction}. To make sense of this, we first note that a multisorted relational structure $\fX$ can be thought of as a bipartite diagram in the category of sets, where each map in the diagram is one of the projections $\pi_k$ from a relation $C_j$ of $\fX$ to a sort $X_{i_k}$ of $\fX$. The bipartite index category $\cS$ of this diagram corresponds directly to the common (multisorted) signature of $\fX$ and $\fA$. The ``category of elements'' of the diagram $F_\fX : \cS \rightarrow \operatorname{Sets}$ gives us a bipartite category $\int_\cS F_\fX$ which can be interpreted as the signature of $\coprod_\fX \fA$. The diagram
\[
F_{\coprod_\fX \fA} : \int_\cS F_\fX \rightarrow \operatorname{Sets}
\]
corresponding to the structure $\coprod_\fX \fA$ is formed by composing the natural functor $\pi : \int_\cS F_\fX \rightarrow \cS$ with the functor $F_\fA : \cS \rightarrow \operatorname{Sets}$ corresponding to $\fA$.
\end{rem}

So the main idea is to replace our original CSP template $\fA$ by the instance-dependent CSP template $\coprod_\fX \fA$, in the hope that it will have extra unary polymorphisms for us to take advantage of. First, however, we have to explain how to express the instance $\fX$ in terms of the new template. We need another definition.

\begin{defn} For a given (possibly multisorted) relational signature, we define the \emph{terminal structure} $*$ to have one element of each sort, and to have a single tuple in each relation.
\end{defn}

Note that for any relational structure $\fA$, the most natural definition of the structure $\fA^0$ is just
\[
\fA^0 \cong *,
\]
and a homomorphism $* \rightarrow \fA$ is the same as a constant (i.e. $0$-ary) polymorphism of $\fA$. Since there is always a unique homomorphism $\fA \rightarrow *$, $\fA$ has a constant polymorphism iff $*$ is the core of $\fA$ - in which case $\CSP(\fA)$ is a trivial problem.

\begin{prop} If $\fX, \fA$ are two relational structures of the same signature, then a homomorphism
\[
\fX \rightarrow \fA
\]
is the same thing as a homomorphism
\[
\coprod_\fX * \rightarrow \coprod_\fX \fA,
\]
which is the same thing as a constant (i.e., $0$-ary) polymorphism of $\coprod_\fX \fA$.
\end{prop}
\begin{proof} Just unwind the definitions.
\end{proof}

Now that we've reinterpreted our original instance $\fX$ of $\CSP(\fA)$ as the new instance $\coprod_\fX *$ of the bespoke $\CSP(\coprod_\fX \fA)$, we want to understand the polymorphisms of $\coprod_\fX \fA$. \begin{itemize}
\item First of all, every polymorphism of $\fA$ extends to a polymorphism of $\coprod_\fX \fA$.
\item Second, \emph{as long as the instance $\fX$ of $\CSP(\fA)$ had a solution}, the new structure $\coprod_\fX \fA$ has a constant polymorphism.
\end{itemize}
So the idea is to look for \emph{unary polynomials} of the original algebraic structure $\bA$ (i.e., unary operations built out of constants and polymorphisms of $\fA$) and hope that they become unary polymorphisms of $\coprod_\fX \fA$. If we manage to find such a unary polynomial, then we can use it to shrink some of the variable domains, thereby making progress towards solving the instance.

The simplest way to build a nontrivial unary polynomial of $\bA$ is to start with a binary polymorphism $t(x,y)$, and to plug a constant $c$ into one of the arguments of $t$ - we'll use the second argument here - making the unary polynomial
\[
t_c(x) \coloneqq t(x,c).
\]
The trouble is that we may need to choose a different constant $c$ for each sort of the new template $\coprod_\fX \fA$, and at first it seems that searching for such a collection of constants is just as hard as solving the original instance of the CSP. The trick is that if $t$ acts like a semilattice, then the values of the unary polynomial $t_c$ mostly won't depend on the particular choice of $c$, giving us more flexibility than we had in the original problem. If all goes according to plan, this will let us reduce the search for a unary polymorphism of $\coprod_\fX \fA$ to an instance of a simpler CSP with smaller domains, which we can solve recursively.

More concretely, we are searching for an unknown unary polymorphism
\[
\varphi : \coprod_\fX \fA \rightarrow \coprod_\fX \fA.
\]
Rather than searching for just any polymorphism - after all, the identity map $x \mapsto x$ is not going to be any help to us - we are searching for $\varphi$ satisfying
\[
\varphi(a) = t_c(a) = t(a, c)
\]
for each $a \in A_x$, with $c \in A_x$ not depending on $a$ (but possibly depending on $x$). So we need to shrink the domain of $\varphi(a)$ to $\{t(a,c) \mid c \in A_x\} = t^a(A_x)$, where $t^a$ is the unary polynomial defined by
\[
t^a(x) \coloneqq t(a,x).
\]
Unfortunately, $t^a(A_x)$ is probably not closed under the polymorphisms $f(x_1, ..., x_k)$ of $\fA$ - instead, it will be closed under operations of the form
\[
(x_1, ..., x_k) \mapsto t^a(f(x_1, ..., x_k)).
\]
So in order to find our hypothetical unary polymorphism $\varphi$ which behaves like the unary polynomials $t_c$, we will end up needing to solve an auxiliary problem involving relational structures that have been crushed down by various unary polynomials of the form $t^a$.

Since it will be important for us to actually \emph{find} a solution to the auxiliary problem, we will need to require the operations of the form $t^a \circ f : t^a(A_x)^k \rightarrow t^a(A_x)$ to be idempotent when $f : A_x^k \rightarrow A_x$ was idempotent - and for this, we need to require that $t$ satisfies
\[
t^a(t^a(x)) = t^a(x)
\]
for all $a$s and all $x$s. We won't need to make such an assumption about the operations $t_c$ or $\varphi$, however, since once we have constructed $\varphi$ we will be free to iterate it to our heart's content to produce a unary polymorphism $\varphi^\infty$ satisfying $\varphi^\infty \circ \varphi^\infty = \varphi^\infty$.

We still need to ensure that $\varphi^{A_x} = t_c$ for some $c \in A_x$. Mar\'oti's approach in \cite{tree-on-top-of-malcev} is to relax this requirement slightly: he just requires $\varphi^{A_x}$ to be contained in the multisorted relation
\[
\Sg_{\prod_a t^a(\bA_x)} \{t_c \mid c \in A_x\} \le \prod_{a \in A_x} t^a(\bA_x).
\]
This will allow $\varphi^{A_x}$ to have the form
\[
\varphi : a \mapsto t(a, f(t(a, c_1), ..., t(a, c_k)))
\]
for any $k$-tuple of elements $c_1, ..., c_k \in A_x$ and any $k$-ary polymorphism $f$ of $\fA$ - but more generally, it also allows $\varphi^{A_x}$ to have the form
\[
\varphi : a \mapsto t(a, f(\varphi_1(a), ..., \varphi_k(a)))
\]
where $\varphi_i$ are any previously constructed operations contained within this multisorted relation. The approach we will follow in these notes will be even simpler - we will just directly search for $\varphi$ such that $\varphi^{A_{x_*}}$ is not surjective for some particular $x_*$.

Time to iron out all of the details. The careful reader might have noticed that I cheated above, by referring to the algebraic structure $t^a(\bA_x)$ without defining it - let's start by addressing that.

\begin{defn}[Following Mar\'oti \cite{tree-on-top-of-malcev}] If $\bA = (A, \{f_i\})$ is an algebraic structure and $e : A \rightarrow A$ is any unary operation satisfying $e(e(a)) = e(a)$ for all $a \in A$, we define the \emph{retract} $e(\bA)$ of $\bA$ with respect to $e$ by
\[
e(\bA) \coloneqq (e(A), \{e \circ f_i\}),
\]
considered as an algebraic structure with the same signature as $\bA$ via $f_i^{e(\bA)} = e \circ f_i|_{e(\bA)}$.
\end{defn}

Note that unlike several other constructions which shrink an algebra with an idempotent unary operation, the retract $e(\bA)$ has the same signature as the algebra $\bA$, so it makes sense to consider multisorted CSPs which involve both $\bA$ and $e(\bA)$ as variable domains simultaneously. The price we pay for this ability is that identities satisfied by term operations of $\bA$ might not hold for the corresponding term operations of $e(\bA)$: the only identities which are definitely preserved are identities between the basic operations of height at most 1.

\begin{prop} If the basic operations $f_i, f_j$ of $\bA$ satisfy an identity of height at most one, i.e.
\[
f_i(x_1, ..., x_k) \approx f_j(y_1, ..., y_l)
\]
for some sets of (possibly repeated) variables $x_1, ..., x_k, y_1, ..., y_l$ or
\[
f_i(x_1, ..., x_k) \approx x,
\]
for some set of (possibly repeated) variables $x_1, ..., x_k, x$, then the basic operations $e \circ f_i, e \circ f_j$ of $e(\bA)$ satisfy the same identity (as long as $e(e(x)) \approx e(x)$, for the second type of identity).
\end{prop}

So to take advantage of retracts, we had better make sure that the basic operations of $\bA$ include a Taylor operation, or at least a collection of idempotent operations satisfying a system of height 1 identities that can't be satisfied by projections. Corollary \ref{consecutive-daisy} shows that any finite Taylor algebra always has a suitable pair of idempotent ternary term operations, and adding them to our collection of basic operations will ensure that $e(\bA)$ is also Taylor.

Recall that we were hoping to construct some interesting unary polymorphisms $\varphi : \coprod_\fX \fA \rightarrow \coprod_\fX \fA$. It turns out that we can turn the search for such polymorphisms into another CSP in a natural way.

\begin{prop} For any multisorted relational structures $\fX, \fA$ sharing a common signature, there are natural bijections between the set of unary polymorphisms
\[
\coprod_\fX \fA \rightarrow \coprod_\fX \fA,
\]
the set of homomorphisms
\[
\fA \times \fX \rightarrow \fA,
\]
and the set of constant polymorphisms
\[
\coprod_{\fA \times \fX} * \rightarrow \coprod_{\fA \times \fX} \fA.
\]
\end{prop}
\begin{proof} Again, just unwind the definitions.
\end{proof}

Now, here's the nice thing: we already know an obvious example of a unary polymorphism $\coprod_\fX \fA \rightarrow \coprod_\fX \fA$, namely, the identity map $\pi_1$, which happens to correspond to the obvious homomorphism
\[
\pi_1 : \fA \times \fX \rightarrow \fA.
\]
This gives us a corresponding obvious constant polymorphism of the structure $\coprod_{\fA \times \fX} \fA$, which we will name $a_{\pi_1}$.

Now recall that we have been assuming that $\fA$ has a binary polymorphism $t$ satisfying the identity
\[
t(x,t(x,y)) \approx t(x,y).
\]
This polymorphism extends to a polymorphism of the structure $\coprod_{\fA \times \fX} \fA$, and we can now combine it with the constant polymorphism $a_{\pi_1}$ to produce the unary polymorphism
\[
t^{a_{\pi_1}} : x \mapsto t(a_{\pi_1}, x)
\]
of $\coprod_{\fA \times \fX} \fA$. By construction, we have
\[
t^{a_{\pi_1}}(t^{a_{\pi_1}}(x)) \approx t^{a_{\pi_1}}(x),
\]
so we can use this unary polymorphism to crush $\coprod_{\fA \times \fX} \fA$ into the smaller structure $t^{a_{\pi_1}}(\coprod_{\fA \times \fX} \fA)$, which should be easier to handle. This leads to the following definition.

\begin{defn}[Mar\'oti \cite{tree-on-top-of-malcev}, slightly simplified] Let $\fX$ be an instance of $\CSP(\{\bA_i\})$, and let $t$ be a binary term operation of $\{\bA_i\}$ satisfying the identity $t(x,t(x,y)) \approx t(x,y)$. The \emph{decomposition} $t(\fX)$ is defined to be the instance of
\[
\CSP(\{t^a(\bA_i) \mid a \in \bA_i\}),
\]
where $t^a(x) \coloneqq t(a,x)$, whose variables and relations are defined as follows:
\begin{itemize}
\item for each variable $x$ of $\fX$ with variable domain $\bA_x$, and for each $a \in \bA_x$, we have a variable $(a,x)$ of $t(\fX)$ with variable domain $t^a(\bA_x)$, where $t^a(\bA_x)$ is the retract of $\bA_x$ with respect to $t^a$,
\item for each constraint tuple $c = (x_1, ..., x_m)$ with constraint relation $\RR \le \bA_{x_1} \times \cdots \times \bA_{x_m}$ of $\fX$, and for each tuple $r = (a_1, ..., a_m) \in \RR$, we have a constraint tuple
\[
(r,c) = ((a_1, x_1), ..., (a_m, x_m))
\]
with corresponding constraint relation given by
\[
t^r(\RR) \le t^{a_1}(\bA_{x_1}) \times \cdots \times t^{a_m}(\bA_{x_m}),
\]
where $t^r(\RR)$ is the retract of $\RR$ with respect to $t^r$.
\end{itemize}
The reader should take a moment to check that the constraint relation $t^r(\RR)$ really is a subalgebra of $t^{a_1}(\bA_{x_1}) \times \cdots \times t^{a_m}(\bA_{x_m})$.
\end{defn}

\begin{prop}[Mar\'oti \cite{tree-on-top-of-malcev}]\label{prop-tree-on-top-unary-poly} If the instance $\fX$ of $\CSP(\{\bA_i\})$ has a solution given by $x \mapsto c_x \in \bA_x$, then for any binary term operation $t$ of $\{\bA_i\}$ the decomposition $t(\fX)$ has a solution given by the mapping
\[
(a,x) \mapsto t(a,c_x).
\]
\end{prop}
\begin{proof} We need to check that for every constraint tuple $(x_1, ..., x_m)$ of $\fX$ with constraint relation $\RR \le \bA_{x_1} \times \cdots \times \bA_{x_m}$, and for each tuple $r = (a_1, ..., a_m) \in \RR$, we have
\[
(t(a_1, c_{x_1}), ..., t(a_m,c_{x_m})) \in t^r(\RR).
\]
But the left hand side is exactly $t^r((c_{x_1}, ..., c_{x_m}))$, and this is an element of $t^r(\RR)$ since $(c_{x_1}, ..., c_{x_m}) \in \RR$, by the assumption that $x \mapsto c_x$ was a solution to the original instance $\fX$.
\end{proof}

\begin{prop}[Mar\'oti \cite{tree-on-top-of-malcev}]\label{prop-tree-on-top-partial-solution} Suppose that $\cA = \{\bA_2, ...\}$ is a finite collection of similar finite idempotent algebras such that $\CSP(\cA)$ can be solved in polynomial time, and suppose that $\bA_1$ is a finite idempotent algebra of the same signature with a binary term operation $t$ satisfying the identity $t(x,t(x,y)) \approx t(x,y)$ such that
\[
a \in \bA \in \{\bA_1\} \cup \cA \;\;\; \implies \;\;\; t^a(\bA) \in \cA.
\]
Then for any instance $\fX$ of $\CSP(\{\bA_1\} \cup \cA)$, we can check whether any given partial assignment to the variables of $t(\fX)$ extends to a solution to the decomposition $t(\fX)$ in polynomial time.
\end{prop}
\begin{proof} The assumptions imply that the decomposition $t(\fX)$ is an instance of $\CSP(\cA)$, so we just need to check that it is at most polynomially larger than the original instance $\fX$. The number of variables of $t(\fX)$ is at most the product of the number of variables of $\fX$ with the maximum number of elements occuring in any of the finitely many finite domains $\bA \in \cA$, and we have $|t^a(\bA_x)| \le |\bA_x|$ for each one.

Similarly, the number of constraint relations of $t(\fX)$ is at most the product of the number of constraints of $\fX$ with the maximum number of tuples occuring in any constraint relation $\RR$ of the instance, and we have $|t^r(\RR)| \le |\RR|$ in each case.
\end{proof}

\begin{prop}[Mar\'oti \cite{tree-on-top-of-malcev}] Let $\fX$ be an instance of $\CSP(\fA)$, and let $t$ be any binary polymorphism of $\fA$. If the decomposition $t(\fX)$ has a solution given by
\[
(a,x) \mapsto \varphi^{A_x}(a),
\]
then $\varphi$ defines a homomorphism of multisorted relational structures
\[
\varphi : \coprod_\fX \fA \rightarrow \coprod_\fX \fA.
\]
\end{prop}
\begin{proof} Recall that the relations of $\coprod_\fX \fA$ were indexed by constraint tuples $c = (x_1, ..., x_m) \in C_j \subseteq X_{i_1} \times \cdots \times X_{i_m}$ of $\fX$, and that if the corresponding constraint relation in $\fA$ is given by $R_j \subseteq A_{i_1} \times \cdots \times A_{i_m}$, then the relation of $\coprod_\fX \fA$ associated with $c$ was given by
\[
R_c \subseteq A_{x_1} \times \cdots \times A_{x_m},
\]
which was a copy of $R_j$. We need to check that $\varphi$ sends each tuple $r = (a_1, ..., a_m) \in R_c$ to a tuple in $R_c$, that is, that
\[
(\varphi^{A_{x_1}}(a_1), ..., \varphi^{A_{x_m}}(a_m)) \stackrel{?}{\in} R_c = R_j.
\]
But the constraint of $t(\fX)$ corresponding to the pair $(r,c)$ requires that
\[
(\varphi^{A_{x_1}}(a_1), ..., \varphi^{A_{x_m}}(a_m)) \in t^r(R_j),
\]
and we have $t^r(R_j) \subseteq R_j$ since $r \in R_j$ and $t$ preserves $R_j$.
\end{proof}

Of course, the decomposition $t(\fX)$ might have useless solutions which fail to shrink the variable domains. The key to Mar\'oti's argument from \cite{tree-on-top-of-malcev} is that we only need a \emph{single} variable domain to be shrunk by some $\varphi^{A_x}$ in order to make progress. So, for each variable $x_*$ of our original instance $\fX$, we search for a solution $(a,x) \mapsto \varphi^{A_x}(a)$ to the decomposition $t(\fX)$ such that the unary operation
\[
\varphi^{A_{x_*}} : A_{x_*} \rightarrow A_{x_*}
\]
satisfies
\[
\varphi^{A_{x_*}}(a) = \varphi^{A_{x_*}}(b)
\]
for some pair $a \ne b \in A_{x_*}$, brute-forcing over all choices of $x_*$, $a \ne b \in A_{x_*}$, and all possible common values of $\varphi^{A_{x_*}}(a),\varphi^{A_{x_*}}(b)$ (using Proposition \ref{prop-tree-on-top-partial-solution} to test each possibility in turn). If we succeed, then we can use the unary polymorphism $\varphi^\infty : \coprod_\fX \fA \rightarrow \coprod_\fX \fA$ to produce a new instance where $A_{x_*}$ is mapped into
\[
(\varphi^\infty)^{A_{x_*}}(A_{x_*}) \subseteq \varphi^{A_{x_*}}(A_{x_*}) \subset A_{x_*},
\]
successfully shrinking at least one variable domain of our original instance $\fX$. Additionally, we can replace all of the relevant algebraic structures $\bA_x$ with the retracts $(\varphi^\infty)^{A_x}(\bA_x)$, which will remain Taylor as long as the basic operations of the algebras $\bA_x$ include, say, a pair of ternary operations as in Corollary \ref{consecutive-daisy}.

The bad case is the case where the decomposition $t(\fX)$ has solutions, but there is no variable $x_*$ such that the decomposition $t(\fX)$ has a solution with $\varphi^{A_{x_*}}$ which is not a permutation. By Proposition \ref{prop-tree-on-top-unary-poly} every solution $x \mapsto c_x \in \bA_x$ of the original instance $\fX$ produces a solution to the decomposition $t(\fX)$ where $\varphi$ is given by
\[
\varphi^{A_x}(a) = t(a, c_x),
\]
so this implies that for each variable $x_*$, and for each solution $x \mapsto c_x \in \bA_x$ of the original instance $\fX$, the map $a \mapsto t(a, c_{x_*})$ is a permutation. But this means that we can replace \emph{every single variable domain} $\bA_x$ by the subalgebra $\bC_x$ which is generated by the collection of elements $c \in \bA_x$ such that the map $a \mapsto t(a, c)$ is a permutation, without losing any solutions to the instance $\fX$ - and if we have $\bC_x \ne \bA_x$ for any $x$, then we have once again made progress towards solving $\fX$.

We have proved Mar\'oti's main result.

\begin{defn}[Mar\'oti \cite{tree-on-top-of-malcev}]\label{defn-tree-on-top-eliminated} If $\bA$ is a finite idempotent algebra, then we say that $\bA$ \emph{can be eliminated} if $\bA$ has a binary term operation $t$ satisfying the identity
\[
t(x,t(x,y)) \approx t(x,y),
\]
such that
\begin{itemize}
\item for each $a \in \bA$ there is some $x \in \bA$ such that $t(a, x) \ne x$, and
\item $\bA$ is \emph{not} generated by $\{c \in \bA \mid t_c : x \mapsto t(x,c) \text{ is a permutation}\}$.
\end{itemize}
\end{defn}

\begin{thm}[Mar\'oti \cite{tree-on-top-of-malcev}] Suppose that $\cA = \{\bA_2, ...\}$ is a finite collection of similar finite idempotent algebras such that $\CSP(\cA)$ can be solved in polynomial time, and suppose that $\bA_1$ is a finite idempotent algebra of the same signature such that
\[
\bA \in \{\bA_1\} \cup \cA,\; \varphi(\bA) \ne \bA,\; \varphi\circ\varphi = \varphi \;\;\; \implies \;\;\; \varphi(\bA) \in \cA,
\]
or in other words, such that $\{\bA_1\} \cup \cA$ is closed under taking retracts.

Suppose additionally that $\bA_1$ can be eliminated, in the sense of Definition \ref{defn-tree-on-top-eliminated}. Then $\CSP(\{\bA_1\} \cup \cA)$ can be solved in polynomial time.
\end{thm}

Mar\'oti went on to prove that in many cases, algebras $\bA$ with a semilattice quotient $\bA/\theta$ can be eliminated. Recall that an element $a$ of a semilattice with basic operation $s$ is called \emph{neutral} if we have $s(a,x) = x$ for all $x$. We extend this definition to partial semilattices (where partial semilattices are defined as in Definition \ref{defn-partial-semilattice}).

\begin{defn} We say that an element $a$ is a \emph{neutral element} of a partial semilattice operation $s$ on $A$ if we have $s(a,x) = x$ for all $x \in A$.
\end{defn}

\begin{prop}[Mar\'oti \cite{tree-on-top-of-malcev}, slight variation] If $\bA$ is a finite idempotent algebra with a partial semilattice term operation $s$ and a quotient $\bA/\theta$ such that $(A/\theta, s)$ is a semilattice with at least two elements, and if $s$ has no neutral element on $\bA$, then $\bA$ can be eliminated.
\end{prop}
\begin{proof} By the definition of a partial semilattice, $s$ satisfies the identity $s(x,s(x,y)) \approx s(x,y)$. Since $s$ has no neutral element by assumption, for every $a \in \bA$ there is some $x \in \bA$ with
\[
s(a,x) \ne x.
\]

There are now two cases, based on whether or not the semilattice $(A/\theta, s)$ has a neutral element. If the semilattice $(A/\theta, s)$ has no neutral element, then for every $c \in \bA$, there is some $a \in \bA$ such that
\[
s(a,c) \equiv_\theta s(c,a) \not\equiv_\theta a,
\]
and for this $a$ we see that $a/\theta$ is not contained in the image of $s_c : x \mapsto s(x,c)$. Thus in this case there is \emph{no} $c \in \bA$ such that $s_c$ is a permutation.

In the remaining case, suppose that the congruence class $a/\theta$ is a neutral element of the semilattice $(A/\theta, s)$. For any $c \not\in a/\theta$, we then have
\[
s(a/\theta, c/\theta) \subseteq c/\theta
\]
since $a/\theta$ was neutral, and
\[
s(c/\theta, c/\theta) \subseteq c/\theta
\]
by idempotence, so $s_c : x \mapsto s(x,c)$ can't be a permutation for any $c \not\in a/\theta$. Since $a/\theta$ is a subalgebra of $\bA$, it can't generate $\bA$, so $\bA$ can be eliminated.
\end{proof}

On the other hand, if we have a prepared algebra $\bA$ (as in Definition \ref{defn-prepared}) with a partial semilattice term $s$ which has a neutral element $a$ on $\bA$, then every binary term operation $t$ of $\bA$ which depends on both of its arguments will have $t(a,x) = x$ for all $x \in \bA$. So the previous result is essentially the best result of its kind when it comes to algebras with semilattices as quotients.

\begin{ex} The simplest algebra which the results we have proved so far still can't handle is the 4-element algebra $\bA = (\{0,1,2,*\}, \cdot)$, where $\cdot$ is the commutative binary operation given by
\begin{center}
\begin{tabular}{c|cccc} $\cdot$ & $0$ & $1$ & $2$ & $*$\\ \hline $0$ & $0$ & $2$ & $1$ & $0$\\ $1$ & $2$ & $1$ & $0$ & $1$\\ $2$ & $1$ & $0$ & $2$ & $2$\\ $*$ & $0$ & $1$ & $2$ & $*$\end{tabular},
\end{center}
which has $*$ as a neutral element. It's a good exercise to find an algorithm which solves $\CSP(\bA)$ in polynomial time.
\end{ex}

We mentioned at the beginning of Section \ref{s-alg-conservative} that Bulatov's algorithm for conservative CSPs from \cite{bulatov-conservative-rerevisited} was based on Mar\'oti's reduction. We can now go over the main ideas of this algorithm.

\begin{prop} If $\bA$ is an algebra and a unary operation $e : A \rightarrow A$ with $e \circ e = e$ has $e(A)$ closed under the basic operations of $\bA$, then $e(\bA)$ is the subalgebra of $\bA$ with underlying set $e(A)$. In particular, if $\bA$ is a conservative algebra, then every retract of $\bA$ is a subalgebra of $\bA$.
\end{prop}
\begin{proof} Let $f$ be a $k$-ary basic operation of $\bA$. Then for any $a_1, ..., a_k \in e(A)$, we have
\[
f(a_1, ..., a_k) = e(b)
\]
for some $b \in A$ by the assumption that $e(A)$ is closed under $f$, so
\[
e\circ f(a_1, ..., a_k) = e(e(b)) = e(b) = f(a_1, ..., a_k).
\]
Thus we have $f^{e(\bA)} = f|_{e(\bA)}$ for each basic operation $f$ of $\bA$, so $e(\bA)$ is exactly the subalgebra of $\bA$ with underlying set $e(A)$.
\end{proof}

\begin{prop}[Bulatov \cite{bulatov-conservative-rerevisited}] If $\bA$ is a finite conservative algebra with a partial semilattice term operation $s$ such that $s$ has no neutral element on $\bA$ and such that there is at least one pair of elements $a \ne b \in \bA$ with $s(a,b) = b$, then $\bA$ can be eliminated.
\end{prop}
\begin{proof} By the definition of a partial semilattice, $s$ satisfies the identity $s(x,s(x,y)) \approx s(x,y)$. Since $s$ has no neutral element by assumption, we just need to check that $\bA$ is not generated by the collection of elements $c \in \bA$ such that $s_c : x \mapsto s(x,c)$ is a permutation. For the pair $a \ne b$ with $s(a,b) = b$ we have
\[
s_b(a) = s_b(b) = b,
\]
so $s_b$ is not a permutation, and since $\bA$ is conservative $\bA \setminus \{b\}$ is a proper subalgebra of $\bA$.
\end{proof}

Recall that by Theorem \ref{thm-swap-gmm}, any conservative Taylor algebra with no semilattice edges has a ternary generalized majority-minority operation, and can therefore be solved by the few subpowers algorithm. It's easy to check that any finite collection of conservative algebras with no semilattice edges also has a ternary operation whose restriction to each conservative algebra is gmm - and that such an operation can be algorithmically constructed in polynomial time - so if we could eliminate every conservative algebra which had a semilattice edge then we could simply apply the few subpowers algorithm to finish. The only difficulty we face is the case where some of our variable domains have neutral elements.

Note that if $\bA$ is conservative and $s$ is a partial semilattice term operation with a neutral element $a$ on $\bA$, then we have
\[
\bA\setminus\{a\} \lhd_{bin} \bA,
\]
with $s$ witnessing the binary absorption. But then we are in position to apply Theorem \ref{thm-rectangularity-application}, with $\cS$ taken to be the collection of all \emph{full} variable domains which have a proper absorbing subalgebra! Thus we can make progress in this case as well.

\begin{ex} Here is a 4-element example which can't be handled by the results proved so far, but which \emph{doesn't} have a neutral element. Of course, it will neither be conservative nor have a semilattice quotient. We take $\bA = (\{a,b,c,*\}, g)$, where $g$ is the ternary operation which is determined by the following properties:
\begin{itemize}
\item $g$ is a symmetric function of its variables,
\item $*$ is a strongly absorbing element for $g$, i.e. $g(*, \bA, \bA) = \{*\}$,
\item each of the sets $\{a,b\}, \{a,c\}, \{b,c\}$ forms a two-element $\ZZ/2^{\aff}$ subalgebra of $\bA$, with $g$ acting as the minority operation on each of these pairs, and
\item $g(a,b,c) = *$.
\end{itemize}
Note that $\bA$ is $2$-conservative (i.e., every two-element subset of $\bA$ is closed under $g$), and it is not too hard to check that $\bA$ is also minimal Taylor. The reader should be able to find a polynomial time algorithm which solves $\CSP(\bA)$ without much difficulty.
\end{ex}

\bibliographystyle{plain}
\bibliography{all}

\begin{thebibliography}{100}

\bibitem{few-subpowers-finitely-related}
Erhard Aichinger, Peter Mayr, and Ralph McKenzie.
\newblock On the number of finite algebraic structures.
\newblock {\em Journal of the European Mathematical Society},
  016(8):1673--1686, 2014.

\bibitem{pac-membership}
D.~Angluin and M.~Kharitonov.
\newblock When won't membership queries help?
\newblock {\em Journal of Computer and System Sciences}, 50(2):336 -- 355,
  1995.

\bibitem{angluin-learning}
Dana Angluin.
\newblock Queries and concept learning.
\newblock {\em Machine Learning}, 2(4):319--342, Apr 1988.

\bibitem{near-subgroups-aschbacher}
Michael Aschbacher.
\newblock Near subgroups of finite groups.
\newblock {\em J. Group Theory}, 1(2):113--129, 1998.

\bibitem{pebble-game-width}
Albert Atserias and V\'ictor Dalmau.
\newblock A combinatorial characterization of resolution width.
\newblock {\em Journal of Computer and System Sciences}, 74(3):323 -- 334,
  2008.
\newblock Computational Complexity 2003.

\bibitem{2+epsilon}
Per Austrin, Venkatesan Guruswami, and Johan H{\aa}stad.
\newblock (2+$\varepsilon$)-{S}at is {NP}-hard.
\newblock {\em SIAM Journal on Computing}, 46(5):1554--1573, 2017.

\bibitem{baker-pixley}
Kirby~A. Baker and Alden~F. Pixley.
\newblock Polynomial interpolation and the {C}hinese remainder theorem for
  algebraic systems.
\newblock {\em Math. Z.}, 143(2):165--174, 1975.

\bibitem{barto-conservative-revisited}
L.~{Barto}.
\newblock The dichotomy for conservative constraint satisfaction problems
  revisited.
\newblock In {\em 2011 IEEE 26th Annual Symposium on Logic in Computer
  Science}, pages 301--310, 2011.

\bibitem{pcsp-coloring-full}
L.~{Barto}, J.~{Bul{\'{\i}}n}, A.~{Krokhin}, and J.~{Opr{\v s}al}.
\newblock {Algebraic approach to promise constraint satisfaction}.
\newblock {\em arXiv e-prints}, November 2018.

\bibitem{barto-kozik-bounded-width}
L.~{Barto} and M.~{Kozik}.
\newblock Constraint satisfaction problems of bounded width.
\newblock In {\em 2009 50th Annual IEEE Symposium on Foundations of Computer
  Science}, pages 595--603, 2009.

\bibitem{near-unanimity-congruence-distributive}
Libor Barto.
\newblock Finitely related algebras in congruence distributive varieties have
  near unanimity terms.
\newblock {\em Canadian Journal of Mathematics}, 65(1):3--21, 2013.

\bibitem{barto}
Libor Barto.
\newblock The collapse of the bounded width hierarchy.
\newblock {\em Journal of Logic and Computation}, 2014.

\bibitem{barto-valeriote-conjecture}
Libor Barto.
\newblock Finitely related algebras in congruence modular varieties have few
  subpowers.
\newblock {\em Journal of the European Mathematical Society}, 20(6):1439--1471,
  2018.

\bibitem{minimal-taylor}
Libor Barto, Zarathustra Brady, Andrei Bulatov, Marcin Kozik, and Dmitriy Zhuk.
\newblock Unifying the three algebraic approaches to the {CSP} via minimal
  {T}aylor algebras.
\newblock {\em TheoretiCS}, 3, 2024.

\bibitem{deciding-absorption-relational}
Libor Barto and Jakub Bul{\'i}n.
\newblock Deciding absorption in relational structures.
\newblock {\em Algebra universalis}, 78(1):3--18, Sep 2017.

\bibitem{near-unanimity-minimal}
Libor Barto and Ond{\v{r}}ej Draganov.
\newblock The minimal arity of near unanimity polymorphisms.
\newblock {\em Mathematica Slovaca}, 69(2):297--310, 2019.

\bibitem{deciding-absorption}
Libor Barto and Alexandr Kazda.
\newblock Deciding absorption.
\newblock {\em International Journal of Algebra and Computation},
  26(05):1033--1060, 2016.

\bibitem{sd-join-cyclic}
Libor Barto and Marcin Kozik.
\newblock Cyclic terms for {SD$\vee$} varieties revisited.
\newblock {\em Algebra universalis}, 64(1):137--142, 2010.

\bibitem{cyclic}
Libor Barto and Marcin Kozik.
\newblock Absorbing subalgebras, cyclic terms, and the constraint satisfaction
  problem.
\newblock {\em Log. Methods Comput. Sci.}, 8(1):1:07, 27, 2012.

\bibitem{sdp}
Libor Barto and Marcin Kozik.
\newblock Robust satisfiability of constraint satisfaction problems.
\newblock In {\em Proceedings of the Forty-fourth Annual ACM Symposium on
  Theory of Computing}, STOC '12, pages 931--940, New York, NY, USA, 2012. ACM.

\bibitem{local-consistency}
Libor Barto and Marcin Kozik.
\newblock Constraint satisfaction problems solvable by local consistency
  methods.
\newblock {\em J. ACM}, 61(1):Art. 3, 19, 2014.

\bibitem{congruence-modular-cyclic}
Libor Barto, Marcin Kozik, Mikl{\'o}s Mar{\'o}ti, Ralph McKenzie, and Todd
  Niven.
\newblock Congruence modularity implies cyclic terms for finite algebras.
\newblock {\em Algebra universalis}, 61(3):365--380, 2009.

\bibitem{pointing-no-absorption}
Libor Barto, Marcin Kozik, and David Stanovsk{\'y}.
\newblock Mal'tsev conditions, lack of absorption, and solvability.
\newblock {\em Algebra universalis}, 74(1):185--206, Sep 2015.

\bibitem{nu-pathwidth}
Libor Barto, Marcin Kozik, and Ross Willard.
\newblock Near unanimity constraints have bounded pathwidth duality.
\newblock In {\em 2012 27th Annual IEEE Symposium on Logic in Computer
  Science}, pages 125--134. IEEE, 2012.

\bibitem{barto-reflections}
Libor Barto, Jakub Opr{\v{s}}al, and Michael Pinsker.
\newblock The wonderland of reflections.
\newblock {\em Israel Journal of Mathematics}, 223(1):363--398, 2018.

\bibitem{topology-irrelevant}
Libor Barto and Michael Pinsker.
\newblock Topology is irrelevant (in the infinite domain dichotomy conjecture
  for constraint satisfaction problems).
\newblock {\em Preprint}, 2018.

\bibitem{few-subpowers}
Joel Berman, Pawe{\l} Idziak, Petar Markovi{\'c}, Ralph McKenzie, Matthew
  Valeriote, and Ross Willard.
\newblock Varieties with few subalgebras of powers.
\newblock {\em Transactions of the American Mathematical Society},
  362(3):1445--1473, 2010.

\bibitem{singleton-arc-consistency-quadratic}
Christian Bessiere, St{\'e}phane Cardon, Romuald Debruyne, and Christophe
  Lecoutre.
\newblock Efficient algorithms for singleton arc consistency.
\newblock {\em Constraints}, 16:25--53, 2011.

\bibitem{birkhoff}
Garrett Birkhoff.
\newblock On the structure of abstract algebras.
\newblock In {\em Mathematical proceedings of the Cambridge philosophical
  society}, volume~31, pages 433--454. Cambridge University Press, 1935.

\bibitem{birkhoff-lattice}
Garrett Birkhoff.
\newblock {\em Lattice theory}, volume~25.
\newblock American Mathematical Soc., 1940.

\bibitem{birkhoff-subdirect}
Garrett Birkhoff.
\newblock Subdirect unions in universal algebra.
\newblock {\em Bull. Amer. Math. Soc.}, 50:764--768, 1944.

\bibitem{median-distributive}
Garrett Birkhoff and Stephen~A Kiss.
\newblock A ternary operation in distributive lattices.
\newblock {\em Bulletin of the American Mathematical Society}, 53(8):749--752,
  1947.

\bibitem{vc-dim-learnability}
Anselm Blumer, A.~Ehrenfeucht, David Haussler, and Manfred~K. Warmuth.
\newblock Learnability and the {V}apnik-{C}hervonenkis {D}imension.
\newblock {\em J. ACM}, 36(4):929--965, October 1989.

\bibitem{bodirsky-thesis}
Manuel Bodirsky.
\newblock Complexity classification in infinite-domain constraint satisfaction.
\newblock {\em arXiv preprint arXiv:1201.0856}, 2012.

\bibitem{small-orbit-growth}
Manuel Bodirsky and Bertalan Bodor.
\newblock Structures with small orbit growth.
\newblock {\em arXiv preprint arXiv:1810.05657}, 2018.

\bibitem{csp-comparing-solution-sets}
Simone Bova, Hubie Chen, and Matthew Valeriote.
\newblock Generic expression hardness results for primitive positive formula
  comparison.
\newblock {\em Information and Computation}, 222:108--120, 2013.
\newblock 38th International Colloquium on Automata, Languages and Programming
  (ICALP 2011).

\bibitem{bowditch-median}
Brian~H Bowditch.
\newblock Median algebras, 2022.

\bibitem{brady-chromatic}
Zarathustra Brady.
\newblock Chromatic numbers of directed hypergraphs with no ``bad'' cycles.
\newblock {\em arXiv preprint arXiv:1806.00783}, 2018.

\bibitem{mui-symmetric}
Zarathustra Brady and Holden Mui.
\newblock Symmetric operations on domains of size at most 4.
\newblock {\em arXiv preprint arXiv:2102.07329}, 2021.

\bibitem{pcsp-symmetric-boolean}
Joshua Brakensiek and Venkatesan Guruswami.
\newblock Promise constraint satisfaction: Algebraic structure and a symmetric
  boolean dichotomy.
\newblock {\em arXiv preprint arXiv:1704.01937}, 2017.

\bibitem{colored-graph-prelim}
A.~A. {Bulatov}.
\newblock A graph of a relational structure and constraint satisfaction
  problems.
\newblock In {\em Proceedings of the 19th Annual IEEE Symposium on Logic in
  Computer Science, 2004.}, pages 448--457, July 2004.

\bibitem{learnable-gmm}
Andrei Bulatov, Hubie Chen, and V{\'i}ctor Dalmau.
\newblock Learnability of relatively quantified generalized formulas.
\newblock In Shoham Ben-David, John Case, and Akira Maruoka, editors, {\em
  Algorithmic Learning Theory}, pages 365--379, Berlin, Heidelberg, 2004.
  Springer Berlin Heidelberg.

\bibitem{bulatov-dalmau-malcev}
Andrei Bulatov and V{\'\i}ctor Dalmau.
\newblock A simple algorithm for mal'tsev constraints.
\newblock {\em SIAM Journal on Computing}, 36(1):16--27, 2006.

\bibitem{subpower-residually-small}
Andrei Bulatov, Peter Mayr, and {\'A}gnes Szendrei.
\newblock The subpower membership problem for finite algebras with cube terms.
\newblock {\em arXiv preprint arXiv:1803.08019}, 2018.

\bibitem{2-semilattice}
Andrei~A. Bulatov.
\newblock Combinatorial problems raised from 2-semilattices.
\newblock {\em J. Algebra}, 298(2):321--339, 2006.

\bibitem{bulatov-bounded}
Andrei~A. Bulatov.
\newblock Bounded relational width.
\newblock {\em manuscript.
  \url{http://www.cs.sfu.ca/~abulatov/papers/relwidth.pdf}}, 2009.

\bibitem{bulatov-conservative}
Andrei~A. Bulatov.
\newblock Complexity of conservative constraint satisfaction problems.
\newblock {\em ACM Trans. Comput. Logic}, 12(4), July 2011.

\bibitem{bulatov-conservative-rerevisited}
Andrei~A. Bulatov.
\newblock Conservative constraint satisfaction re-revisited.
\newblock {\em Journal of Computer and System Sciences}, 82(2):347--356, 2016.

\bibitem{colored-graph}
Andrei~A. Bulatov.
\newblock Graphs of finite algebras, edges, and connectivity.
\newblock {\em CoRR}, abs/1601.07403, 2016.

\bibitem{bulatov-dichotomy}
Andrei~A Bulatov.
\newblock A dichotomy theorem for nonuniform {CSP}s.
\newblock In {\em 2017 IEEE 58th Annual Symposium on Foundations of Computer
  Science (FOCS)}, pages 319--330. IEEE, 2017.

\bibitem{bulatov-jeavons-varieties}
Andrei~A. Bulatov and Peter~G. Jeavons.
\newblock Algebraic structures in combinatorial problems.
\newblock Technical Report MATH-AL-4-2001, Technische universit{\"a}t
  {D}resden, Dresden, {G}ermany, 2001.

\bibitem{carbonnel-thesis}
Cl{\'e}ment Carbonnel.
\newblock {\em Harnessing tractability in constraint satisfaction problems}.
\newblock PhD thesis, Institut National Polytechnique de Toulouse, 2016.

\bibitem{symmetric-polymorphisms}
Catarina Carvalho and Andrei Krokhin.
\newblock On algebras with many symmetric operations.
\newblock {\em International Journal of Algebra and Computation},
  26(05):1019--1031, 2016.

\bibitem{sectionally-pseudocomplemented-lattices}
Ivan Chajda and Sándor Radeleczki.
\newblock On varieties defined by pseudocomplemented nondistributive lattices.
\newblock {\em Publicationes Mathematicae}, 63, 11 2003.

\bibitem{chen-few-subpowers}
Hubie Chen.
\newblock The expressive rate of constraints.
\newblock {\em Annals of Mathematics and Artificial Intelligence},
  44(4):341--352, Aug 2005.

\bibitem{arc}
Hubie Chen, Victor Dalmau, and Berit Gru{\ss}ien.
\newblock Arc consistency and friends.
\newblock {\em Journal of Logic and Computation}, 23(1):87--108, 2013.

\bibitem{membership-query-hardness}
Hubie Chen and Matthew Valeriote.
\newblock Learnability of solutions to conjunctive queries.
\newblock {\em Journal of Machine Learning Research}, 20(67):1--28, 2019.

\bibitem{dalmau-thesis}
V~Dalmau.
\newblock {\em Computational complexity of problems over generalized formulas,
  2000}.
\newblock PhD thesis, PhD thesis, Universitat Polit{\'e}cnica de Catalunya.

\bibitem{dalmau-gmm}
Victor Dalmau.
\newblock Generalized majority-minority operations are tractable.
\newblock In {\em 20th Annual IEEE Symposium on Logic in Computer Science
  (LICS'05)}, pages 438--447. IEEE, 2005.

\bibitem{no-pure-width-2}
V\'ictor Dalmau.
\newblock There are no pure relational width 2 constraint satisfaction
  problems.
\newblock {\em Information Processing Letters}, 109(4):213 -- 218, 2009.

\bibitem{robust-near-unanimity}
V{\'\i}ctor Dalmau, Marcin Kozik, Andrei Krokhin, Konstantin Makarychev, Yury
  Makarychev, and Jakub Opr{\v{s}}al.
\newblock Robust algorithms with polynomial loss for near-unanimity {CSP}s.
\newblock In {\em Proceedings of the Twenty-Eighth Annual ACM-SIAM Symposium on
  Discrete Algorithms}, pages 340--357. SIAM, 2017.

\bibitem{dalmau-robust-framework}
V{\'\i}ctor Dalmau and Andrei Krokhin.
\newblock Robust satisfiability for {CSP}s: {H}ardness and algorithmic results.
\newblock {\em ACM Transactions on Computation Theory (TOCT)}, 5(4):1--25,
  2013.

\bibitem{dalmau-approximation}
V{\'\i}ctor Dalmau, Andrei Krokhin, and Rajsekar Manokaran.
\newblock Towards a characterization of constant-factor approximable
  finite-valued {CSP}s.
\newblock {\em Journal of Computer and System Sciences}, 97:14 -- 27, 2018.

\bibitem{dalmau-width-1}
V{\'\i}ctor Dalmau and Justin Pearson.
\newblock Closure functions and width 1 problems.
\newblock In {\em International Conference on Principles and Practice of
  Constraint Programming}, pages 159--173. Springer, 1999.

\bibitem{dechter}
Rina Dechter.
\newblock From local to global consistency.
\newblock {\em Artificial intelligence}, 55(1):87--107, 1992.

\bibitem{dedekind-modular}
Richard Dedekind.
\newblock {\"U}ber die von drei moduln erzeugte dualgruppe.
\newblock {\em Mathematische Annalen}, 53(3):371--403, 1900.

\bibitem{deutsch-quantum-turing-machine}
David Deutsch.
\newblock Quantum theory, the {C}hurch--{T}uring principle and the universal
  quantum computer.
\newblock {\em Proceedings of the Royal Society of London. A. Mathematical and
  Physical Sciences}, 400(1818):97--117, 1985.

\bibitem{hypergraph-promise-hardness}
Irit Dinur, Oded Regev, and Clifford Smyth.
\newblock The hardness of 3-uniform hypergraph coloring.
\newblock {\em Combinatorica}, 25(5):519--535, 2005.

\bibitem{hornsat-linear}
William~F. Dowling and Jean~H. Gallier.
\newblock Linear-time algorithms for testing the satisfiability of
  propositional horn formulae.
\newblock {\em The Journal of Logic Programming}, 1(3):267--284, 1984.

\bibitem{eckmann-hilton}
Beno Eckmann and Peter~J Hilton.
\newblock Group-like structures in general categories {I} multiplications and
  comultiplications.
\newblock {\em Mathematische Annalen}, 145(3):227--255, 1962.

\bibitem{pseudovarieties-ultimate}
Samuel Eilenberg and Marcel~P Sch{\"u}tzenb{\'e}rger.
\newblock {\em On pseudovarieties}.
\newblock IRIA. Laboratoire de Recherche en Informatique et Automatique, 1975.

\bibitem{near-subgroups-feder}
Tomas Feder.
\newblock Constraint satisfaction on finite groups with near subgroups.
\newblock In {\em Electronic Colloquium on Computational Complexity (ECCC),
  TR05-005}, 2005.

\bibitem{feder-vardi}
Tom{\'a}s Feder and Moshe~Y Vardi.
\newblock The computational structure of monotone monadic {SNP} and constraint
  satisfaction: A study through {D}atalog and group theory.
\newblock {\em SIAM Journal on Computing}, 28(1):57--104, 1998.

\bibitem{commutator-theory}
Ralph Freese and Ralph McKenzie.
\newblock {\em Commutator theory for congruence modular varieties}, volume 125.
\newblock CUP Archive, 1987.

\bibitem{group-algorithms}
Merrick Furst, John Hopcroft, and Eugene Luks.
\newblock Polynomial-time algorithms for permutation groups.
\newblock In {\em 21st Annual Symposium on Foundations of Computer Science
  (sfcs 1980)}, pages 36--41. IEEE, 1980.

\bibitem{geiger-galois}
David Geiger.
\newblock Closed systems of functions and predicates.
\newblock {\em Pacific journal of mathematics}, 27(1):95--100, 1968.

\bibitem{goldreich-majority}
Oded Goldreich.
\newblock Valiant’s polynomial-size monotone formula for majority, 2011.

\bibitem{semimodular-jordan-holder}
G~Gr{\"a}tzer and JB~Nation.
\newblock Prime intervals and maximal chains in finite dimensional semimodular
  lattices.
\newblock 2010.

\bibitem{treewidth-homomorphism}
Martin Grohe.
\newblock The complexity of homomorphism and constraint satisfaction problems
  seen from the other side.
\newblock {\em Journal of the ACM (JACM)}, 54(1):1, 2007.

\bibitem{gumm-geometric}
Heinz~Peter Gumm.
\newblock {\em Geometrical methods in congruence modular algebras}, volume 286.
\newblock American Mathematical Soc., 1983.

\bibitem{robust-horn-gap}
Venkatesan Guruswami and Yuan Zhou.
\newblock Tight bounds on the approximability of almost-satisfiable {H}orn
  {SAT} and {E}xact {H}itting {S}et.
\newblock In {\em Proceedings of the twenty-second annual ACM-SIAM symposium on
  Discrete algorithms}, pages 1574--1589. Society for Industrial and Applied
  Mathematics, 2011.

\bibitem{sorting-word-RAM-deterministic}
Yijie Han.
\newblock Deterministic sorting in {$O(n \log \log n)$} time and linear space.
\newblock In {\em Proceedings of the thiry-fourth annual ACM symposium on
  Theory of computing}, pages 602--608, 2002.

\bibitem{h-coloring}
Pavol Hell and Jaroslav Ne\v{s}et\v{r}il.
\newblock On the complexity of {$H$}-coloring.
\newblock {\em J. Combin. Theory Ser. B}, 48(1):92--110, 1990.

\bibitem{core-graph}
Pavol Hell and Jaroslav Ne\v{s}et\v{r}il.
\newblock The core of a graph.
\newblock {\em Discrete Mathematics}, 109(1):117 -- 126, 1992.

\bibitem{free-modular-undecidable}
Christian Herrmann.
\newblock On the word problem for the modular lattice with four free
  generators.
\newblock {\em Mathematische Annalen}, 265(4):513--527, 1983.

\bibitem{higmans-lemma}
Graham Higman.
\newblock Ordering by divisibility in abstract algebras.
\newblock {\em Proceedings of the London Mathematical Society},
  s3-2(1):326--336, 1952.

\bibitem{hobby-mckenzie}
David Hobby and Ralph McKenzie.
\newblock {\em The structure of finite algebras}, volume~76 of {\em
  Contemporary Mathematics}.
\newblock American Mathematical Society, Providence, RI, 1988.

\bibitem{hodges-model}
Wilfrid Hodges.
\newblock {\em Model theory}.
\newblock Cambridge University Press, 1993.

\bibitem{hodges-shorter}
Wilfrid Hodges.
\newblock {\em A shorter model theory}.
\newblock Cambridge university press, 1997.

\bibitem{horowitz-thesis}
Jonah Horowitz.
\newblock {\em Results on the Computational Complexity of Linear Idempotent
  Mal'cev Conditions}.
\newblock PhD thesis, 2012.

\bibitem{horowitz-complexity-malcev-conditions}
Jonah Horowitz.
\newblock Computational complexity of various {M}al'cev conditions.
\newblock {\em International Journal of Algebra and Computation},
  23(06):1521--1531, 2013.

\bibitem{hastad-optimal}
Johan H\r{a}stad.
\newblock Some optimal inapproximability results.
\newblock {\em J. ACM}, 48(4):798–859, July 2001.

\bibitem{few-subpowers-algorithm}
Pawe{\l} Idziak, Petar Markovi{\'c}, Ralph McKenzie, Matthew Valeriote, and
  Ross Willard.
\newblock Tractability and learnability arising from algebras with few
  subpowers.
\newblock {\em SIAM Journal on Computing}, 39(7):3023--3037, 2010.

\bibitem{P-vs-BPP-exponential-circuits}
Russell Impagliazzo and Avi Wigderson.
\newblock {P = BPP} if {E} requires exponential circuits: {D}erandomizing the
  {XOR} lemma.
\newblock In {\em Proceedings of the twenty-ninth annual ACM symposium on
  Theory of computing}, pages 220--229, 1997.

\bibitem{jeavons}
Peter Jeavons.
\newblock On the algebraic structure of combinatorial problems.
\newblock {\em Theoretical Computer Science}, 200(1-2):185--204, 1998.

\bibitem{indicator-instance}
Peter Jeavons, David Cohen, and Marc Gyssens.
\newblock Closure properties of constraints.
\newblock {\em J. ACM}, 44(4):527–548, July 1997.

\bibitem{affine-quandles}
P\v{r}emysl Jedli\v{c}ka, Agata Pilitowska, David Stanovsk\'y, and Anna
  Zamojska-Dzienio.
\newblock Subquandles of affine quandles.
\newblock {\em Journal of Algebra}, 510:259 -- 288, 2018.

\bibitem{jonsson-distributive}
Bjarni J{\'o}nsson.
\newblock Algebras whose congruence lattices are distributive.
\newblock {\em Mathematica Scandinavica}, pages 110--121, 1968.

\bibitem{jovanovic-terms}
Jelena Jovanovi{\'c}.
\newblock On terms describing omitting unary and affine types.
\newblock {\em Filomat}, 27(1):183--199, 2013.

\bibitem{optimal-maltsev}
Jelena Jovanovi{\'{c}}, Petar Markovi{\'{c}}, Ralph McKenzie, and Matthew
  Moore.
\newblock Optimal strong mal'cev conditions for congruence
  meet-semidistributivity in locally finite varieties.
\newblock {\em Algebra universalis}, pages 1--21, 2016.

\bibitem{karnin2012explicit}
Zohar~S Karnin, Yuval Rabani, and Amir Shpilka.
\newblock Explicit dimension reduction and its applications.
\newblock {\em SIAM Journal on Computing}, 41(1):219--249, 2012.

\bibitem{directed-gumm}
Alexandr Kazda, Marcin Kozik, Ralph McKenzie, and Matthew Moore.
\newblock {\em Absorption and directed J{\'o}nsson terms}, pages 203--220.
\newblock Springer International Publishing, Cham, 2018.

\bibitem{deciding-minority}
Alexandr Kazda, Jakub Opršal, Matt Valeriote, and Dmitriy Zhuk.
\newblock Deciding the existence of minority terms.
\newblock {\em Canadian Mathematical Bulletin}, 63(3):577–591, 2020.

\bibitem{deciding-some-malcev-conditions}
Alexandr Kazda and Matt Valeriote.
\newblock Deciding some {M}altsev conditions in finite idempotent algebras.
\newblock {\em The Journal of Symbolic Logic}, 85(2):539–562, 2020.

\bibitem{cube-terms-chipped-cubes}
Alexandr Kazda and Dmitriy Zhuk.
\newblock Existence of cube terms in finite algebras.
\newblock {\em Algebra universalis}, 82:1--29, 2021.

\bibitem{optimal-taylor}
Keith Kearnes, Petar Markovi{\'c}, and Ralph McKenzie.
\newblock Optimal strong {M}al'cev conditions for omitting type 1 in locally
  finite varieties.
\newblock {\em Algebra Universalis}, 72(1):91--100, 2014.

\bibitem{kearnes-quasi-affine}
Keith~A Kearnes.
\newblock A quasi-affine representation.
\newblock {\em International Journal of Algebra and Computation}, 5:673--702,
  1995.

\bibitem{kearnes-difference}
Keith~A Kearnes.
\newblock Varieties with a difference term.
\newblock {\em Journal of Algebra}, 177(3):926--960, 1995.

\bibitem{kearnes-simple}
Keith~A. Kearnes.
\newblock Idempotent simple algebras.
\newblock In {\em Logic and algebra ({P}ontignano, 1994)}, volume 180 of {\em
  Lecture Notes in Pure and Appl. Math.}, pages 529--572. Dekker, New York,
  1996.

\bibitem{kearnes-taylor-affine}
Keith~A Kearnes and {\'A}gnes Szendrei.
\newblock The relationship between two commutators.
\newblock {\em International Journal of Algebra and Computation},
  8(04):497--531, 1998.

\bibitem{parallelogram-terms}
Keith~A Kearnes and {\'A}gnes Szendrei.
\newblock Clones of algebras with parallelogram terms.
\newblock {\em International Journal of Algebra and Computation},
  22(01):1250005, 2012.

\bibitem{cube-terms-crosses}
Keith~A Kearnes and {\'A}gnes Szendrei.
\newblock Cube term blockers without finiteness.
\newblock {\em Algebra universalis}, 78(4):437--459, 2017.

\bibitem{strongly-abelian-hamiltonian}
E.~Kiss and M.~Valeriote.
\newblock Strongly abelian varieties and the hamiltonian property.
\newblock {\em Canadian Journal of Mathematics}, 43(2):331–346, 1991.

\bibitem{vcsp-algorithm}
Vladimir Kolmogorov, Andrei Krokhin, and Michal Rolinek.
\newblock The complexity of general-valued {CSP}s.
\newblock {\em SIAM Journal on Computing}, 46(3):1087--1110, 2017.

\bibitem{majority-sorting-networks}
Alexander Kozachinskiy and Vladimir Podolskii.
\newblock Multiparty {K}archmer-{W}igderson games and threshold circuits.
\newblock {\em arXiv preprint arXiv:2002.07444}, 2020.

\bibitem{kozik-subpower-exptime}
Marcin Kozik.
\newblock A finite set of functions with an {EXPTIME}-complete composition
  problem.
\newblock {\em Theoretical Computer Science}, 407(1):330 -- 341, 2008.

\bibitem{slac}
Marcin Kozik.
\newblock Weaker consistency notions for all the {CSP}s of bounded width.
\newblock {\em CoRR}, abs/1605.00565, 2016.

\bibitem{pq-consistency}
Marcin Kozik.
\newblock Solving {CSP}s using weak local consistency.
\newblock 2018.

\bibitem{maltsev}
Marcin Kozik, Andrei Krokhin, Matt Valeriote, and Ross Willard.
\newblock Characterizations of several {M}altsev conditions.
\newblock {\em Algebra Universalis}, 73(3-4):205--224, 2015.

\bibitem{vcsp-hardness}
Marcin Kozik and Joanna Ochremiak.
\newblock Algebraic properties of valued constraint satisfaction problem.
\newblock In {\em International Colloquium on Automata, Languages, and
  Programming}, pages 846--858. Springer, 2015.

\bibitem{lp-width-1}
Gabor Kun, Ryan O'Donnell, Suguru Tamaki, Yuichi Yoshida, and Yuan Zhou.
\newblock Linear programming, width-1 {CSP}s, and robust satisfaction.
\newblock In {\em Proceedings of the 3rd Innovations in Theoretical Computer
  Science Conference}, pages 484--495. ACM, 2012.

\bibitem{ladner}
Richard~E Ladner.
\newblock On the structure of polynomial time reducibility.
\newblock {\em Journal of the ACM (JACM)}, 22(1):155--171, 1975.

\bibitem{partial-poly-seth}
Victor Lagerkvist and Magnus Wahlstr{\"o}m.
\newblock The ({C}oarse) {F}ine-{G}rained {S}tructure of {NP}-{H}ard {SAT} and
  {CSP} problems.
\newblock {\em ACM Transactions on Computation Theory (TOCT)}, 14(1):1--54,
  2021.

\bibitem{zadori-near-unanimity-graphs}
Benoit Larose, Cynthia Loten, and L{\'a}szl{\'o} Z{\'a}dori.
\newblock A polynomial-time algorithm for near-unanimity graphs.
\newblock {\em Journal of Algorithms}, 55(2):177--191, 2005.

\bibitem{ability-to-count}
Benoit Larose, Matt Valeriote, and L{\'a}szl{\'o} Z{\'a}dori.
\newblock Omitting types, bounded width and the ability to count.
\newblock {\em Internat. J. Algebra Comput.}, 19(5):647--668, 2009.

\bibitem{lau-clone-theory}
Dietlinde Lau.
\newblock {\em Function algebras on finite sets: Basic course on many-valued
  logic and clone theory}.
\newblock Springer Science \& Business Media, 2006.

\bibitem{lipparini-difference}
Paolo Lipparini.
\newblock Difference terms and commutators.

\bibitem{littlestone-online-learning}
Nick Littlestone.
\newblock Learning quickly when irrelevant attributes abound: A new
  linear-threshold algorithm.
\newblock {\em Machine Learning}, 2(4):285--318, Apr 1988.

\bibitem{cube-term-blockers}
Petar Markovi{\'c}, Mikl{\'o}s Mar{\'o}ti, and Ralph McKenzie.
\newblock Finitely related clones and algebras with cube terms.
\newblock {\em Order}, 29(2):345--359, 2012.

\bibitem{near-unanimity-maroti}
Mikl{\'o}s Mar{\'o}ti.
\newblock The existence of a near-unanimity term in a finite algebra is
  decidable.
\newblock {\em The Journal of Symbolic Logic}, 74(3):1001--1014, 2009.

\bibitem{malcev-on-top}
Mikl{\'o}s Mar{\'o}ti.
\newblock Maltsev on top.
\newblock {\em Manuscript, available at
  \url{http://www.math.u-szeged.hu/~mmaroti/pdf/200x\%20Maltsev\%20on\%20top.pdf}},
  2011.

\bibitem{tree-on-top-of-malcev}
Mikl{\'o}s Mar{\'o}ti.
\newblock Tree on top of {M}altsev.
\newblock {\em Manuscript, available at
  \url{http://www.math.u-szeged.hu/~mmaroti/pdf/200x\%20Tree\%20on\%20top\%20of\%20Maltsev.pdf}},
  2011.

\bibitem{subpower-supernilpotent}
Peter Mayr.
\newblock The subpower membership problem for mal'cev algebras.
\newblock {\em International Journal of Algebra and Computation},
  22(07):1250075, 2012.

\bibitem{finite-forbidden-lattices}
Ralph McKenzie.
\newblock Finite forbidden lattices.
\newblock In Ralph~S. Freese and Octavio~C. Garcia, editors, {\em Universal
  Algebra and Lattice Theory}, pages 176--205, Berlin, Heidelberg, 1983.
  Springer Berlin Heidelberg.

\bibitem{commutator-uses}
Ralph McKenzie and John Snow.
\newblock Congruence modular varieties: commutator theory and its uses.
\newblock In {\em Structural theory of automata, semigroups, and universal
  algebra}, pages 273--329. Springer, 2005.

\bibitem{symmetric-terms-simple-groups}
Sebastian Meyer and Florian Starke.
\newblock Finite simple groups in the primitive positive constructability
  poset.
\newblock {\em arXiv preprint arXiv:2409.06487}, 2024.

\bibitem{maximal-intersecting}
Aaron Meyerowitz.
\newblock Maximal intersecting families.
\newblock {\em European Journal of Combinatorics}, 16(5):491 -- 501, 1995.

\bibitem{finitely-related-undecidable}
Matthew Moore.
\newblock Finite degree clones are undecidable.
\newblock {\em Theoretical Computer Science}, 796:237--271, 2019.

\bibitem{quantum-computation}
Michael~A Nielsen and Isaac~L Chuang.
\newblock {\em Quantum computation and quantum information}.
\newblock Cambridge university press, 2010.

\bibitem{olsak-weak}
Miroslav Ol{\v{s}}{\'a}k.
\newblock The weakest nontrivial idempotent equations.
\newblock {\em Bulletin of the London Mathematical Society}, 49(6):1028--1047,
  2017.

\bibitem{residually-small-groups}
A~Ju Ol'{\v{s}}anski{\u\i}.
\newblock Varieties of finitely approximable groups.
\newblock {\em Mathematics of the USSR-Izvestiya}, 3(4):867, 1969.

\bibitem{commutator-notes}
Peter Ouwehand.
\newblock Commutator theory and abelian algebras.
\newblock {\em arXiv preprint arXiv:1309.0662}, 2013.

\bibitem{palfy-permutational}
P{\'e}ter~P{\'a}l P{\'a}lfy.
\newblock Unary polynomials in algebras, {I}.
\newblock {\em Algebra Universalis}, 18(3):262--273, 1984.

\bibitem{palfy-pudlak}
P{\'e}ter~P{\'a}l P{\'a}lfy and Pavel Pudl{\'a}k.
\newblock Congruence lattices of finite algebras and intervals in subgroup
  lattices of finite groups.
\newblock {\em Algebra Universalis}, 11(1):22--27, 1980.

\bibitem{abelian-malcev-affine}
H.~Peter~Gumm.
\newblock Algebras in permutable varieties: Geometrical properties of affine
  algebras.
\newblock {\em algebra universalis}, 9(1):8--34, Dec 1979.

\bibitem{pinsker-rosenberg}
Michael Pinsker.
\newblock {\em Rosenberg's characterization of maximal clones}.
\newblock na, 2002.

\bibitem{Pixley-term}
Alden~F Pixley.
\newblock Distributivity and permutability of congruence relations in
  equational classes of algebras.
\newblock {\em Proceedings of the American Mathematical Society},
  14(1):105--109, 1963.

\bibitem{plonka-p-cyclic}
J.~P{\l}onka.
\newblock On {$k$}-cyclic groupoids.
\newblock {\em Math. Japon.}, 30(3):371--382, 1985.

\bibitem{poschel-galois}
Reinhard P{\"o}schel.
\newblock A general galois theory for operations and relations and concrete
  characterization of related algebraic structures.
\newblock 1980.

\bibitem{post-lattice}
Emil~L Post.
\newblock {\em The Two-Valued Iterative Systems of Mathematical Logic}.
\newblock Princeton University Press, 1942.

\bibitem{quasi-affine-quackenbush}
Robert~W. Quackenbush.
\newblock Quasi-affine algebras.
\newblock {\em algebra universalis}, 20(3):318--327, Oct 1985.

\bibitem{raghavendra-optimal}
Prasad Raghavendra.
\newblock Optimal algorithms and inapproximability results for every {CSP}?
\newblock In {\em Proceedings of the fortieth annual ACM symposium on Theory of
  computing}, pages 245--254. ACM, 2008.

\bibitem{raghavendra-thesis}
Prasad Raghavendra.
\newblock {\em Approximating {NP}-hard problems: efficient algorithms and their
  limits}.
\newblock University of Washington, 2009.

\bibitem{pseudovarieties-implicit}
Jan Reiterman.
\newblock The {B}irkhoff theorem for finite algebras.
\newblock {\em Algebra universalis}, 14(1):1--10, 1982.

\bibitem{median-poc}
Martin Roller.
\newblock Poc sets, median algebras and group actions.
\newblock {\em arXiv preprint arXiv:1607.07747}, 2016.

\bibitem{romov-multisorted-galois}
B.~A. Romov.
\newblock On the lattice of subalgebras of direct products of post algebras of
  finite degree.
\newblock {\em Mathematical Models of Complex Systems, IK AN UkrSSR, Kiev},
  65(1):156--168, 1973.

\bibitem{rosenberg-completeness}
Ivo Rosenberg.
\newblock {\em {\"U}ber die funktionale Vollst{\"a}ndigkeit in den mehrwertigen
  Logiken}.
\newblock Academia, 1970.

\bibitem{salomaa-essential}
Arto Salomaa.
\newblock {\em On essential variables of functions, especially in the algebra
  of logic}.
\newblock Suomalainen Tiedeakatemia, 1963.

\bibitem{schaefer}
Thomas~J. Schaefer.
\newblock The complexity of satisfiability problems.
\newblock In {\em Conference {R}ecord of the {T}enth {A}nnual {ACM} {S}ymposium
  on {T}heory of {C}omputing ({S}an {D}iego, {C}alif., 1978)}, pages 216--226.
  ACM, New York, 1978.

\bibitem{subpower-hardness}
Jeff Shriner.
\newblock Hardness results for the subpower membership problem.
\newblock {\em International Journal of Algebra and Computation},
  28(05):719--732, 2018.

\bibitem{siggers-original}
Mark~H Siggers.
\newblock A strong mal’cev condition for locally finite varieties omitting
  the unary type.
\newblock {\em Algebra universalis}, 64(1-2):15--20, 2010.

\bibitem{stronkowski-embedding}
Micha{\l} Stronkowski and David Stanovsk{\'y}.
\newblock Embedding general algebras into modules.
\newblock {\em Proceedings of the American Mathematical Society},
  138(8):2687--2699, 2010.

\bibitem{semiprojection-lemma}
S.~{\'S}wierczkowski.
\newblock Algebras which are independently generated by every n elements.
\newblock {\em Fundamenta Mathematicae}, 49:93--104, 1960.

\bibitem{taylor-varieties}
Walter Taylor.
\newblock Varieties obeying homotopy laws.
\newblock {\em Canadian Journal of Mathematics}, 29(3):498--527, 1977.

\bibitem{turing-computable}
Alan~Mathison Turing.
\newblock On computable numbers, with an application to the
  {E}ntscheidungsproblem.
\newblock {\em J. of Math}, 58(345-363):5, 1936.

\bibitem{deciding-n-permutable}
M.~Valeriote and R.~Willard.
\newblock {Idempotent $n$-permutable varieties}.
\newblock {\em Bulletin of the London Mathematical Society}, 46(4):870--880, 06
  2014.

\bibitem{valiant-pac}
L.~G. Valiant.
\newblock A theory of the learnable.
\newblock In {\em Proceedings of the Sixteenth Annual ACM Symposium on Theory
  of Computing}, STOC '84, pages 436--445, New York, NY, USA, 1984. ACM.

\bibitem{valiant-majority}
L.G Valiant.
\newblock Short monotone formulae for the majority function.
\newblock {\em Journal of Algorithms}, 5(3):363 -- 366, 1984.

\bibitem{sparse-systems}
Douglas Wiedemann.
\newblock Solving sparse linear equations over finite fields.
\newblock {\em IEEE transactions on information theory}, 32(1):54--62, 1986.

\bibitem{witt-wedderburn-little}
Ernst Witt.
\newblock {\"U}ber die kommutativit{\"a}t endlicher schiefk{\"o}rper.
\newblock In {\em Abhandlungen aus dem Mathematischen Seminar der
  Universit{\"a}t Hamburg}, volume~8, pages 413--413. Springer, 1931.

\bibitem{uncountable-clones}
Yu~I Yanov and AA~Muchnik.
\newblock On the existence of k-valued closed classes that do not have a basis.
\newblock In {\em Soviet Acad. Sci. Dokl}, volume 127, pages 144--146, 1959.

\bibitem{zhuk-selfdual}
Dmitriy Zhuk.
\newblock The lattice of all clones of self-dual functions in three-valued
  logic.
\newblock {\em Journal of Multiple-Valued Logic \& Soft Computing}, 24, 2015.

\bibitem{zhuk-dichotomy}
Dmitriy Zhuk.
\newblock A proof of {CSP} dichotomy conjecture.
\newblock In {\em 2017 IEEE 58th Annual Symposium on Foundations of Computer
  Science (FOCS)}, pages 331--342. IEEE, 2017.

\bibitem{zhuk-strong}
Dmitriy Zhuk.
\newblock Strong subalgebras and the constraint satisfaction problem.
\newblock {\em J. Multiple Valued Log. Soft Comput.}, 36(4-5):455--504, 2021.

\bibitem{near-unanimity-zhuk}
Dmitriy~N Zhuk.
\newblock The existence of a near-unanimity function is decidable.
\newblock {\em Algebra universalis}, 71(1):31--54, 2014.

\bibitem{zhuk-key}
Dmitriy~N Zhuk.
\newblock Key (critical) relations preserved by a weak near-unanimity function.
\newblock {\em Algebra universalis}, 77(2):191--235, 2017.

\bibitem{uncountable-polynomial-clones}
I.~Ágoston, J.~Demetrovics, and L.~Hannák.
\newblock On the number of clones containing all constants (a problem of {R.
  McKenzie}).
\newblock In L.~Szabó and Á. Szendrei, editors, {\em Lectures in Universal
  Algebra}, Colloquia Mathematica Societatis Janos Bolyai, pages 21--25.
  North-Holland, Amsterdam, 1986.

\end{thebibliography}


\begin{appendices}

\chapter{Commutator theory in congruence modular varieties}\label{a-commutator}

Before diving into commutator theory, we'll review of some of the theory of modular lattices. The theory really begins with the observation that in any module, the lattice of submodules is always \emph{ranked} (so long as there are no infinite chains of submodules). In fact, not only is this lattice ranked, but also every (finite) \emph{sublattice} of the lattice of submodules is ranked as well. So it is natural to study lattices which have this property.

\begin{defn} The \emph{length} of a finite chain is the number of elements in the chain minus $1$. The \emph{length} of a poset is the supremum of the lengths of all of its chains.
\end{defn}

\begin{defn} A poset satisfies the \emph{Jordan-Dedekind chain condition} if for any $a \le b$, any two maximal chains from $a$ to $b$ have equal length.
\end{defn}

The simplest situation to consider is the situation where some element $a$ has two distinct covers $b,c$. Then $a = b \wedge c$, and we may start by considering sublattices of the interval $\llbracket a, b \vee c \rrbracket$. The claim is that in this scenario, if we want every sublattice of the interval $\llbracket a, b \vee c \rrbracket$ to be ranked, then we need $b\vee c$ to cover both $b$ and $c$ (so the interval $\llbracket a, b \vee c \rrbracket$ must have length two). If $b \vee c$ does \emph{not} cover $c$, say $c < d < b \vee c$ for some $d$, then we have a problem: the sublattice generated by $b,c,d$ is a copy of the pentagon lattice $\cN_5$, which is not ranked. The only hard part of verifying this is checking that $b \wedge d = a$, but this follows from $a \le b\wedge d \le b$ and $b \not\le d$.

\begin{defn} A poset is called \emph{upper semimodular} if whenever an element $a$ has two distinct covers $b,c$, there is some element $d$ which covers both $b$ and $c$.
\end{defn}

Surprisingly, it turns out that any upper semimodular poset which has no infinite chains satisfies the Jordan-Dedekind chain condition. Note that every chain is contained in a maximal chain (by Zorn's Lemma).




\begin{prop} If $a$ is any element of an upper semimodular poset which has no infinite chains, then any two maximal chains starting at $a$ (going upwards) have the same length.
\end{prop}
\begin{proof} Let $a < a_1 < \cdots$ and $a < a_1' < \cdots$ be two maximal chains starting from $a$ of lengths $m, n$, and induct on $\min(m,n)$. We may assume without loss of generality that $m \le n$. By upper semimodularity, there is some element $a_2''$ which covers both $a_1$ and $a_1'$. Pick some maximal chain $a_2'' < a_3'' < \cdots$ starting from $a_2''$. Then the maximal chains $a_1 < a_2 < \cdots$ and $a_1 < a_2'' < \cdots$ must both have length $m-1$ by the induction hypothesis. Since the maximal chain $a_1' < a_2'' < \cdots$ then also has length $m-1$, we can apply the induction hypothesis to see that the maximal chain $a_1' < a_2' < \cdots$ has length $m-1$ as well, so $m = n$.
\end{proof}

\begin{cor}[Birkhoff \cite{birkhoff-lattice}] An upper semimodular poset which has no infinite chains satisfies the Jordan-Dedekind chain condition.
\end{cor}
\begin{proof} If $a \le b$, then we can pick some fixed maximal chain $b < b_1 < \cdots$ starting from $b$. By appending it to any two maximal chains from $a$ to $b$ of different lengths, we obtain two maximal chains starting from $a$ which have different lengths, contradicting the previous proposition.
\end{proof}

On any poset of finite length which satisfies the Jordan-Dedekind chain condition and has upper or lower bounds, we can define a \emph{height function} $h$ such that whenever $a$ is covered by $b$, we have $h(b) = h(a)+1$.

\begin{prop}[Birkhoff \cite{birkhoff-lattice}] A ranked lattice of finite length is upper semimodular if and only if its height function satisfies the inequality
\[
h(x) + h(y) \ge h(x \vee y) + h(x \wedge y).
\]
\end{prop}
\begin{proof} The inequality clearly implies upper semimodularity. Now suppose our lattice is upper semimodular, and pick maximal chains
\begin{align*}
x\wedge y &= x_0 < x_1 < \cdots < x_m = x,\\
x\wedge y &= y_0 < y_1 < \cdots < y_n = y.
\end{align*}
We claim that for each $i,j$, $x_i \vee y_j$ is either covered by or equal to $x_{i+1} \vee y_j$ and $x_i \vee y_{j+1}$. We can prove this by induction on $i,j$: if it's true for $i,j$, then by upper semimodularity $x_{i+1} \vee y_{j+1}$ will either cover or be equal to both of $x_{i+1} \vee y_j$ and $x_i \vee y_{j+1}$.

Thus, the sequence
\[
x = x \vee y_0 \le x \vee y_1 \le \cdots \le x \vee y_n = x \vee y
\]
has every adjacent pair either equal or a cover, so
\[
h(x \vee y) - h(x) \le h(y) - h(x \wedge y).\qedhere
\]
\end{proof}

There is also a corresponding notion of lower semimodularity, and a dual version of the above result. Putting them together, we get the following.

\begin{thm}[Birkhoff \cite{birkhoff-lattice}] A lattice of finite length is modular iff it satisfies the Jordan-Dedekind chain condition and its height function satisfies
\[
h(x) + h(y) = h(x\vee y) + h(x \wedge y).
\]
\end{thm}
\begin{proof} Since modular implies both upper and lower semimodular, it implies the chain condition and the condition on the height function. For the other direction, suppose that we have a ranked lattice whose height function satisfies the given condition.

Suppose for contradiction that there is a sublattice isomorphic to the pentagon $\cN_5$ (recall from the discussion around Definition \ref{modular-defn} that a lattice is modular iff it doesn't have $\cN_5$ as a sublattice). Suppose this sublattice is generated by $a,b,c$, with $b < c$ and $a \wedge b = a \wedge c$, $a \vee b = a \vee c$. Then we have
\[
h(a) + h(b) = h(a\vee b) + h(a \wedge b) = h(a\vee c) + h(a \wedge c) = h(a) + h(c),
\]
so $h(b) = h(c)$, contradicting $b < c$.
\end{proof}

The next result can be viewed as a strengthening of the fact that a modular lattice is both upper and lower semimodular.

\begin{thm}[Diamond Isomorphism Theorem]\label{diamond-isom-lattice} If $a,b$ are elements of a modular lattice, then the maps $\phi : \llbracket a, a\vee b\rrbracket \rightarrow \llbracket a\wedge b, b\rrbracket$ and $\varphi: \llbracket a\wedge b, b\rrbracket \rightarrow \llbracket a, a \vee b\rrbracket$ given by
\[
\phi : x \mapsto x\wedge b \text{ and } \varphi : y \mapsto y\vee a
\]
are lattice isomorphisms.
\end{thm}
\begin{proof} First we check that $\phi, \varphi$ are inverse to each other. By the modular law, for $x \in \llbracket a, a \vee b \rrbracket$ we have
\[
\varphi(\phi(x)) = (x \wedge b) \vee a = x \wedge (b \vee a) = x,
\]
and for $y \in \llbracket a \wedge b, b \rrbracket$ we have
\[
\phi(\varphi(y)) = (y \vee a) \wedge b = y \vee (a \wedge b) = y.
\]
It is clear that $\phi$ respects meets and that $\varphi$ respects joins, so from the fact that they are inverse to each other we see that they are both lattice isomorphisms.
\end{proof}

\begin{defn} If $a, b$ are elements of a lattice, then we say that the intervals $\llbracket a, a \vee b \rrbracket$ and $\llbracket a \wedge b, b\rrbracket$ are \emph{perspective} to each other, and we abbreviate this with either the notation
\[
\llbracket a, a \vee b \rrbracket \searrow \llbracket a \wedge b, b \rrbracket
\]
or the notation
\[
\llbracket a \wedge b, b \rrbracket \nearrow \llbracket a, a \vee b\rrbracket.
\]
If two intervals in a lattice can be connected by a chain of perspectivities, then we say that they are \emph{projective} to each other.
\end{defn}

The fact that all maximal chains in a finite length semimodular lattice have the same length can be strengthened to a lattice version of the Jordan-H\"older Theorem.

\begin{thm}[Jordan-H\"older for semimodular lattices \cite{semimodular-jordan-holder}]\label{thm-jordan-holder} Suppose we have two maximal chains
\begin{align*}
0 &= a_0 < a_1 < \cdots < a_n = 1,\\
0 &= b_0 < b_1 < \cdots < b_n = 1
\end{align*}
in an upper semimodular lattice of finite length. Then there is a permutation $\sigma \in S_n$ such that each $\llbracket a_{i-1}, a_i\rrbracket$ is projective in two steps (going $\nearrow$, $\searrow$) to $\llbracket b_{\sigma(i)-1}, b_{\sigma(i)} \rrbracket$.
\end{thm}
\begin{proof} We induct on the length $n$. If $a_1 = b_1$ then we can apply the inductive hypothesis. Otherwise, for each $i$, let $c_i = a_1 \vee b_i$. If $k$ is maximal such that $a_1 \not\le b_k$, then
\[
a_1 = c_0 < c_1 < \cdots < c_k = c_{k+1} < \cdots < c_n = 1
\]
where the strict inequalities up to $c_k$ follow from upper semimodularity, and in the portion after $c_{k+1}$ we have $c_j = b_j$.

Applying the induction hypothesis, we get a bijection $\sigma' : [n]\setminus\{1\} \rightarrow [n]\setminus\{k+1\}$ such that each $\llbracket a_{i-1}, a_i\rrbracket$ is projective going $\nearrow$, $\searrow$ to $\llbracket c_{\sigma'(i)-1}, c_{\sigma'(i)} \rrbracket$. Since $\llbracket c_{\sigma'(i)-1}, c_{\sigma'(i)} \rrbracket \searrow \llbracket b_{\sigma'(i)-1}, b_{\sigma'(i)} \rrbracket$, and since $\llbracket a_0, a_1 \rrbracket \nearrow \llbracket b_k, b_{k+1}\rrbracket$, we can take $\sigma$ to be the extension of $\sigma'$ given by setting $\sigma(1) = k+1$.
\end{proof}

To relate this to the usual Jordan-H\"older Theorem, we have to consider the lattice of \emph{subnormal subgroups} of a group. A subgroup $\bM \le \bG$ is called subnormal if there is a finite chain of subgroups connecting it to $\bG$, such that each is a normal subgroup of the next.

\begin{prop} A subgroup $\bM \le \bG$ is subnormal iff the sequence of groups $\bG = \bG_0 \rhd \bG_1 \rhd \cdots$ defined by taking $\bG_{i+1}$ to be the normal closure of $\bM$ inside $\bG_i$ eventually reaches $\bM$. As a consequence, the intersection of two subnormal subgroups is also subnormal.
\end{prop}

\begin{prop} If $\bG$ is a group of finite composition length, then the collection of subnormal subgroups of $\bG$ forms a lower semimodular lattice. If $\llbracket \bN_1, \bM_1 \rrbracket$, $\llbracket \bN_2, \bM_2 \rrbracket$ are $\searrow, \nearrow$ projective covers in this lattice, then $\bM_1/\bN_1 \cong \bM_2/\bN_2$.
\end{prop}

Note that the modular law is equivalent to the following identity, which recovers the usual modular law in the case $a \le b$ by replacing $a \wedge b$ with $a$:
\[
(a \wedge b) \vee (c \wedge b) \approx ((a \wedge b) \vee c) \wedge b.
\]
Thus modular lattices form a variety of lattices. We finish our review of modular lattices by mentioning a famous result of Dedekind.

\begin{prop}[Dedekind \cite{dedekind-modular}] The free modular lattice on $3$ generators is finite, with exactly $28$ elements and length $8$. It is isomorphic to a subdirect product of six copies of the two-element lattice and a single copy of the diamond lattice $\cM_3$.

In particular, one can test whether a given $3$-variable lattice identity is a consequence of modularity in finite time, by testing whether it holds on $\cM_3$.
\end{prop}

A corresponding result for $4$ generators does not exist: the free modular lattice on $4$ generators is infinite. To see this, note that if you start with four generic points on the projective plane and repeatedly generate new points and lines, the resulting set of points and lines you obtain is infinite. Determining whether a $4$-variable lattice identity follows from the modular law is undecidable in general \cite{free-modular-undecidable}.

\section{The Shifting Lemma and the Day terms}\label{s-shifting}

We will follow Freese and McKenzie \cite{commutator-theory}, with some arguments taken from Gumm \cite{gumm-geometric} and some from \cite{commutator-notes}. The starting point for proving things in congruence modular varieties is the Shifting Lemma (this is the main place in the theory where the modular law is actually used).

\begin{lem}[Shifting Lemma] If $\bA$ is congruence modular, $x,y,z,w \in \bA$ and $\alpha,\beta,\gamma \in \Con(\bA)$ with $\alpha \wedge \beta \le \gamma$ and $x \equiv_\alpha y, z\equiv_\alpha w, x \equiv_\beta z, y\equiv_\beta w$, then $z \equiv_\gamma w \implies x \equiv_\gamma y$.
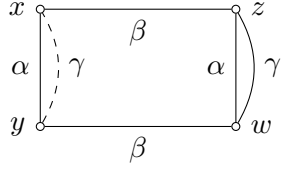
\begin{center}
\begin{tikzpicture}[scale=1.3]
\node[circle, minimum width=3pt, draw, inner sep=0pt, label=left:{$x$}] (a) at (0,1.2){};
\node[circle, minimum width=3pt, draw, inner sep=0pt, label=right:{$z$}] (c) at (2.0,1.2){};
\node[circle, minimum width=3pt, draw, inner sep=0pt, label=left:{$y$}] (b) at (0,0){};
\node[circle, minimum width=3pt, draw, inner sep=0pt, label=right:{$w$}] (d) at (2.0,0){};
\draw (a) to ["$\beta$"'] (c) to ["$\alpha$"'] (d) to ["$\beta$"] (b) to ["$\alpha$"] (a);
\draw [bend left] (c) to ["$\gamma$"] (d);
\draw [bend right, dashed] (b) to ["$\gamma$"'] (a);
\end{tikzpicture}
\end{center}
\end{lem}
\begin{proof} We have $(x,y) \in \alpha \wedge (\beta \circ (\alpha \wedge \gamma) \circ \beta) \subseteq \alpha \wedge (\beta \vee (\alpha \wedge \gamma))$. Since $\alpha \wedge \gamma \le \alpha$, we can apply the modular law to get $\alpha \wedge (\beta \vee (\alpha \wedge \gamma)) = (\alpha \wedge \beta) \vee (\alpha \wedge \gamma)$, and this is contained in $\gamma$ by the assumption $\alpha \wedge \beta \le \gamma$, so $(x,y) \in \gamma$.
\end{proof}

\begin{cor}[Day terms]\label{day-terms} In any congruence modular variety $\cV$, if $\cF_\cV(x,y,z,w)$ is the free algebra on four generators, and if we let $\theta_{a,b}$ be the congruence generated by identifying $a,b$, then there are quaternary terms $m_0, ..., m_n \in \cF_\cV(x,y,z,w)$ such that
\begin{align*}
m_0 &= x,\\
m_i\ &(\theta_{x,y}\vee \theta_{z,w})\wedge (\theta_{x,z}\vee \theta_{y,w})\ m_{i+1}\text{ for $i$ even,}\\
m_i\ &\theta_{z,w}\ m_{i+1}\text{ for $i$ odd,}\\
m_n &= y.
\end{align*}
In other words, the $m_i$ satisfy the following system of identities:
\begin{align*}
m_0(x,y,z,w) &\approx x,\\
m_i(x,x,z,z) &\approx x\text{ for all $i$,}\\
m_i(x,y,x,y) &\approx m_{i+1}(x,y,x,y)\text{ for $i$ even,}\\
m_i(x,y,z,z) &\approx m_{i+1}(x,y,z,z)\text{ for $i$ odd,}\\
m_n(x,y,z,w) &\approx y.
\end{align*}
\end{cor}
\begin{proof} Apply the Shifting Lemma with $\alpha = \theta_{x,y}\vee \theta_{z,w}, \beta = \theta_{x,z}\vee \theta_{y,w}$, and $\gamma = (\alpha \wedge \beta)\vee \theta_{z,w}$ to see that $(x,y) \in (\alpha \wedge \beta)\vee \theta_{z,w} = \bigcup_n ((\alpha \wedge \beta)\circ \theta_{z,w})^{\circ n}$.
\end{proof}

\begin{lem}\label{day-congruence} Let $\bA$ be an algebra with Day terms $m_0, ..., m_n$, $\theta \in \Con(\bA)$, and $a,b,c,d \in \bA$ with $(c,d) \in \theta$. Then $(a,b) \in \theta$ iff for all $i \le n$ we have $m_i(a,b,a,b) \equiv_\theta m_i(a,b,c,d)$.
\end{lem}
\begin{proof} If $(a,b) \in \theta$, then for each $i$ we have $m_i(a,b,a,b) \equiv_\theta m_i(a,a,a,a) = a$ and $m_i(a,b,c,d) \equiv_\theta m_i(a,a,c,c) = a$. For the converse direction, we will show that if $c \equiv_\theta d$ and $m_i(a,b,a,b) \equiv_\theta m_i(a,b,c,d)$ for all $i$, then $m_i(a,b,c,d) \equiv_\theta m_{i+1}(a,b,c,d)$ for all $i$, and then we can conclude $a = m_0(a,b,c,d) \equiv_\theta m_n(a,b,c,d) = b$.

For $i$ even, we use $m_i(a,b,a,b) = m_{i+1}(a,b,a,b)$ together with the assumed congruences relating $m_i(a,b,a,b)$ to $m_i(a,b,c,d)$, while for $i$ odd we use $m_i(a,b,c,c) = m_{i+1}(a,b,c,c)$ together with $c \equiv_\theta d$.
\end{proof}

The existence of Day terms implies a result slightly stronger than the Shifting Lemma, called the Shifting Principle.

\begin{lem}[The Shifting Principle] If $\bA$ has Day terms $m_0, ..., m_n$, then $\bA$ satisfies the Shifting Principle: if $x,y,z,w \in \bA$ and $\alpha,\gamma \in \Con(\bA)$ and $\Lambda \le \bA^2$ is a reflexive relation preserved by the $m_i$ with $\alpha \cap \Lambda \subseteq \gamma$ and $x \equiv_\alpha y, z\equiv_\alpha w, (x,z)\in \Lambda, (y,w)\in \Lambda$, then $z \equiv_\gamma w \implies x \equiv_\gamma y$.
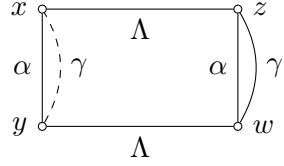
\begin{center}
\begin{tikzpicture}[scale=1.3]
\node[circle, minimum width=3pt, draw, inner sep=0pt, label=left:{$x$}] (a) at (0,1.2){};
\node[circle, minimum width=3pt, draw, inner sep=0pt, label=right:{$z$}] (c) at (2.0,1.2){};
\node[circle, minimum width=3pt, draw, inner sep=0pt, label=left:{$y$}] (b) at (0,0){};
\node[circle, minimum width=3pt, draw, inner sep=0pt, label=right:{$w$}] (d) at (2.0,0){};
\draw (a) to ["$\Lambda$"'] (c) to ["$\alpha$"'] (d) to ["$\Lambda$"] (b) to ["$\alpha$"] (a);
\draw [bend left] (c) to ["$\gamma$"] (d);
\draw [bend right, dashed] (b) to ["$\gamma$"'] (a);
\end{tikzpicture}
\end{center}
\end{lem}
\begin{proof} By Lemma \ref{day-congruence}, it's enough to show that $m_i(x,y,x,y) \equiv_\gamma m_i(x,y,z,w)$ for each $i$. Since $\Lambda$ is preserved by the $m_i$ and is reflexive, we have
\[
\begin{bmatrix} m_i(x,y,x,y)\\ m_i(x,y,z,w)\end{bmatrix} = m_i\left(\begin{bmatrix} x\\ x\end{bmatrix}, \begin{bmatrix} y\\ y\end{bmatrix}, \begin{bmatrix} x\\ z\end{bmatrix}, \begin{bmatrix} y\\ w\end{bmatrix}\right) \in \Lambda,
\]
while $m_i(x,y,x,y) \equiv_\alpha m_i(x,y,z,w)$ by Lemma \ref{day-congruence}, so $(m_i(x,y,x,y),m_i(x,y,z,w)) \in \alpha \cap \Lambda \subseteq \gamma$.
\end{proof}

\begin{lem} If the Shifting Principle holds for an algebra $\bA$ in the special case where $\alpha \ge \gamma$, then $\bA$ is congruence modular.
\end{lem}
\begin{proof} Suppose that $\alpha, \beta, \gamma \in \Con(\bA)$ with $\alpha \ge \gamma \ge \alpha \wedge \beta$, then to verify congruence modularity we just need to check that $\alpha \wedge (\beta \vee \gamma) \le \gamma$, as this rules out the existence of a sublattice of $\Con(\bA)$ isomorphic to the pentagon $\cN_5$. Defining reflexive, symmetric relations $\Lambda_i$ by $\Lambda_i = \beta\circ (\gamma\circ\beta)^{\circ i}$, we see that we just need to prove that $\alpha \cap \Lambda_i \subseteq \gamma$ for each $i$.

We will prove this by induction on $i$: note that the base case $i = 0$ is trivial, since $\Lambda_0 = \beta$. For the inductive step, we apply the Shifting Principle to $\alpha, \Lambda_i$, and $\gamma$ see that if $\alpha \cap \Lambda_i \subseteq \gamma$, then
\[
\alpha \cap \Lambda_{2i+1} = \alpha \cap (\Lambda_i \circ \gamma \circ \Lambda_i) = \alpha \cap (\Lambda_i \circ (\alpha\wedge \gamma) \circ \Lambda_i) \subseteq \gamma.\qedhere
\]
\end{proof}

\begin{cor} A variety is congruence modular iff it has Day terms.
\end{cor}

\begin{ex} If $p(x,y,z)$ is a Mal'cev term, then we can take
\begin{align*}
m_0(x,y,z,w) &= x,\\
m_1(x,y,z,w) &= p(z,w,y),\\
m_2(x,y,z,w) &= y
\end{align*}
as a sequence of Day terms. Rather than laboriously checking the Day identities, it is easier to verify that this sequence of terms can be used in the Shifting Lemma setup to show that $x \equiv_\gamma y$. We have $x\ (\alpha \wedge \beta)\ p(z,w,y)\ \gamma\ y$, so from $\alpha \wedge \beta \le \gamma$ we get $x\ \gamma\ y$.
\end{ex}

\begin{ex} If $g(x,y,z)$ is a majority term, then we can take
\begin{align*}
m_0(x,y,z,w) &= x,\\
m_1(x,y,z,w) &= g(x,y,z),\\
m_2(x,y,z,w) &= g(x,y,w),\\
m_3(x,y,z,w) &= y
\end{align*}
as a sequence of Day terms. Again, in the Shifting Lemma setup, we have $x\ (\alpha\wedge\beta)\ g(x,y,z)\ \gamma\ g(x,y,w)\ (\alpha\wedge\beta)\ y$, so $x\ \gamma\ y$.
\end{ex}

The next corollary gives us a large class of examples of congruence modular varieties, generalizing groups and rings.

\begin{defn} An algebra is \emph{congruence regular} if every congruence on $\bA$ is uniquely determined by any of its congruence classes.
\end{defn}

\begin{cor}[Gumm \cite{gumm-geometric}] If every subalgebra of $\bA^2$ is congruence regular, then $\bA$ is congruence modular.
\end{cor}
\begin{proof} We just need to verify the Shifting Principle for $\bA$. Let $\Lambda \le \bA^2$ be reflexive, let $\alpha \ge \gamma$ be congruences on $\bA$ with $\alpha \cap \Lambda \subseteq \gamma$, and consider the congruences $\alpha \times \gamma$ and $\gamma \times \gamma$ restricted to $\Lambda$. We will show that for any $a \in \bA$, the congruence classes containing $(a,a)$ in these restrictions are equal, so congruence regularity will imply that $\alpha\times\gamma|_\Lambda = \gamma\times\gamma|_\Lambda$, which is the Shifting Principle for $\alpha,\Lambda,\gamma$.

So suppose that $(a,a) \equiv_{\alpha\times\gamma} (b,c) \in \Lambda$. Then $(b,c) \in \alpha \circ \gamma = \alpha$, so $(b,c) \in \alpha \cap \Lambda \subseteq \gamma$, and this implies that $(a,b) \in \gamma \circ \gamma = \gamma$. Thus $(a,a) \equiv_{\gamma\times\gamma} (b,c)$ as well, and we are done.
\end{proof}

To finish this section, we will prove one of Gumm's ``geometric'' results on congruence modular varieties, which generalizes the result used to prove associativity of the loop operation in the case of abelian Mal'cev algebras.

\begin{lem}[The Cube Lemma \cite{gumm-geometric}]\label{mod-cube-lemma} Suppose every subalgebra of $\bA^2$ satisfies the Shifting Lemma. If $\alpha, \beta, \gamma \in \Con(\bA)$ with $\gamma \ge \alpha \wedge \beta$, and if $a,b,c,d,a',b',c',d' \in \bA$ with $(a,b),(c,d),(a',b'),(c',d') \in \alpha$, $(a,d),(b,c),(a',d'),(b',c') \in \beta$, and $(a,a'),(b,b'),(c,c') \in \gamma$, then $(d,d') \in \gamma$.
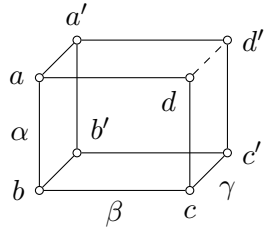
\begin{center}
\begin{tikzpicture}[xscale=0.5,yscale=0.5]
\node[circle, minimum width=3pt, draw, inner sep=0pt, label=left:{$b$}] (11) at (1,1){};
\node[circle, minimum width=3pt, draw, inner sep=0pt, label=above right:{$b'$}] (22) at (2,2){};
\node[circle, minimum width=3pt, draw, inner sep=0pt, label=left:{$a$}] (13) at (1,4){};
\node[circle, minimum width=3pt, draw, inner sep=0pt, label=above:{$a'$}] (24) at (2,5){};
\node[circle, minimum width=3pt, draw, inner sep=0pt, label=below:{$c$}] (31) at (5,1){};
\node[circle, minimum width=3pt, draw, inner sep=0pt, label=right:{$c'$}] (42) at (6,2){};
\node[circle, minimum width=3pt, draw, inner sep=0pt, label=below left:{$d$}] (33) at (5,4){};
\node[circle, minimum width=3pt, draw, inner sep=0pt, label=right:{$d'$}] (44) at (6,5){};
\draw (11) to ["$\beta$"'] (31);
\draw (22) -- (42);
\draw (13) -- (33);
\draw (24) -- (44);
\draw (11) to ["$\alpha$"] (13);
\draw (22) -- (24);
\draw (31) -- (33);
\draw (42) -- (44);
\draw (11) -- (22);
\draw (13) -- (24);
\draw (31) to ["$\gamma$"'] (42);
\draw[dashed] (33) -- (44);
\end{tikzpicture}
\end{center}
\end{lem}
\begin{proof} We apply the Shifting Lemma to the algebra $\beta \le \bA^2$, and the congruences $\gamma \times 1_\bA|_\beta, \alpha\times\alpha|_\beta$, and $\gamma \times \gamma|_\beta$. By the Shifting Lemma applied to $\alpha, \beta, \gamma$, we have $(\alpha\times\alpha|_\beta) \wedge (\gamma \times 1_\bA|_\beta) \le \gamma \times \gamma|_\beta$, so the Shifting Lemma applies to $\gamma \times 1_\bA|_\beta, \alpha\times\alpha|_\beta, \gamma \times \gamma|_\beta$.

Thus, from $((b,c),(b',c')) \in \gamma \times \gamma|_\beta$, $((a,d),(a',d')) \in \gamma \times 1_\bA|_\beta$, and $((b,c),(a,d)), ((b',c'),(a',d')) \in \alpha\times\alpha|_\beta$, the Shifting Lemma allows us to conclude that $((a,d),(a',d')) \in \gamma \times \gamma|_\beta$, so $(d,d') \in \gamma$.
\end{proof}

\section{The modular commutator}

First we go over a slick proof of the main properties of the commutator, using the Day terms to construct an explicit set of generators $X(\alpha,\beta)$ for the congruence $[\alpha,\beta]$ - however, as the approach feels somewhat ad-hoc, we will also prove these properties via a different approach based on the Shifting Lemma applied to congruences (as in the proof of the Cube Lemma). The definition of $X(\alpha,\beta)$ is based on the algebra of matrices $\bM(\alpha,\beta)$ used to visualize the term condition (Definition \ref{commutator-matrix}) and Lemma \ref{day-congruence}.

\begin{defn} Suppose $\bA$ has Day terms $m_0, ..., m_n$. For $\alpha, \beta \in \Con(\bA)$, we define $X(\alpha,\beta)$ to be the set of pairs $(m_i(a,b,a,b),m_i(a,b,c,d))$ for $\begin{bmatrix} a & c\\ b & d\end{bmatrix} \in \bM(\alpha,\beta)$ and $i \le n$.
\end{defn}

\begin{ex} If we have a Mal'cev term $p(x,y,z)$ and take $m_0 = x, m_1 = p(z,w,y), m_2 = y$ as our sequence of Day terms, then $X(\alpha,\beta)$ is the set of pairs $(a,p(c,d,b))$ for $\begin{bmatrix} a & c\\ b & d\end{bmatrix} \in \bM(\alpha,\beta)$.
\end{ex}

\begin{ex} If we have a majority term $g(x,y,z)$ and take $m_0 = x, m_1 = g(x,y,z), m_2 = g(x,y,w), m_3 = y$ as our sequence of Day terms, then $X(\alpha,\beta)$ is the set of pairs $(a,g(a,b,c))$ for $\begin{bmatrix} a & c\\ b & d\end{bmatrix} \in \bM(\alpha,\beta)$. For $(a,b) \in \alpha\wedge\beta$, we have
\[
g\left(\begin{bmatrix} a & a\\ b & b\end{bmatrix}, \begin{bmatrix} a & b\\ a & b\end{bmatrix}, \begin{bmatrix} b & b\\ b & b\end{bmatrix}\right) = \begin{bmatrix} a & b\\ b & b\end{bmatrix} \in \bM(\alpha,\beta),
\]
so $(a,g(a,b,b)) = (a,b) \in X(\alpha,\beta)$. Thus $X(\alpha,\beta) = \alpha\wedge\beta$ for majority algebras.
\end{ex}

\begin{thm}[Commutator via Day terms]\label{day-commutator} If $\bA$ has Day terms $m_0, ..., m_n$ and $\alpha,\beta,\delta \in \Con(\bA)$, then the following are equivalent.
\begin{itemize}
\item[(i)] $X(\alpha,\beta) \subseteq \delta$,

\item[(ii)] $X(\beta,\alpha) \subseteq \delta$,

\item[(iii)] $C(\alpha,\beta;\delta)$ holds,

\item[(iv)] $C(\beta,\alpha;\delta)$ holds,

\item[(v)] $[\alpha,\beta] \le \delta$.
\end{itemize}
\end{thm}
\begin{proof} It's enough to show (iii) $\implies$ (i) $\implies$ (iv). For (iii) $\implies$ (i), suppose that $\begin{bmatrix} a & c\\ b & d\end{bmatrix} \in \bM(\alpha,\beta)$, then
\[
m_i\left(\begin{bmatrix} a & a\\ a & a\end{bmatrix}, \begin{bmatrix} a & a\\ b & b\end{bmatrix}, \begin{bmatrix} a & c\\ a & c\end{bmatrix}, \begin{bmatrix} a & c\\ b & d\end{bmatrix}\right) = \begin{bmatrix} a & a\\ m_i(a,b,a,b) & m_i(a,b,c,d)\end{bmatrix} \in \bM(\alpha,\beta),
\]
so $C(\alpha,\beta;\delta)$ implies that we have $(m_i(a,b,a,b),m_i(a,b,c,d)) \in \delta$.

For (i) $\implies$ (iv), we apply Lemma \ref{day-congruence} to see that if $\begin{bmatrix} a & c\\ b & d\end{bmatrix} \in \bM(\alpha,\beta)$, $(c,d) \in \delta$, and $X(\alpha,\beta) \subseteq \delta$, then we must have $(a,b) \in \delta$ as well, so $C(\beta,\alpha;\delta)$ holds.
\end{proof}

Now we can finally prove some useful properties of commutators.

\begin{prop}\label{mod-commutator} If $\bA$ is contained in a congruence modular variety, then for congruences on $\bA$ we have
\begin{itemize}
\item[(a)] $[\alpha,\beta] = [\beta,\alpha]$,

\item[(b)] $[\alpha \wedge \gamma, \beta] \le [\alpha,\beta]\wedge \gamma$,

\item[(c)] $[\bigvee_i \alpha_i, \beta] = \bigvee_i[\alpha_i,\beta]$,

\item[(d)] if $f : \bA \twoheadrightarrow \bB$ is surjective, then $f([\alpha,\beta]\vee \ker f) = [f(\alpha \vee \ker f), f(\beta \vee \ker f)]$,

\item[(e)] if $\bB \le \bA$, then $[\alpha|_\bB,\beta|_\bB] \le [\alpha,\beta]|_\bB$,

\item[(f)] if $\bA = \prod_{i \in I} \bA_i$, then $[\bigoplus_i \alpha_i, \bigoplus_i \beta_i] = \bigoplus_i[\alpha_i,\beta_i]$, where $\bigoplus_i \alpha_i$ is the set of pairs $(a,b)$ in $\prod_i \alpha_i$ such that for all but finitely many $i$ we have $a_i = b_i$,

\item[(g)] if $\bA = \prod_{i \in I} \bA_i$, then $[\prod_i \alpha_i, \prod_i \beta_i] \le \prod_i [\alpha_i, \beta_i]$.
\end{itemize}
\end{prop}
\begin{proof} Part (a) follows from Theorem \ref{day-commutator}, part (b) follows from Proposition \ref{gen-commutator}(d), and part (e) is Proposition \ref{gen-commutator}(g). For part (c), Theorem \ref{day-commutator} shows that $C(\alpha_j,\beta;\bigvee_i [\alpha_i,\beta])$ holds for each $j$, so we can use Proposition \ref{gen-commutator}(e) to see that $[\bigvee_i \alpha_i, \beta] \le \bigvee_i[\alpha_i,\beta]$, while the other inequality follows from monotonicity of the commutator.

For part (d), note that part (c) implies $[\alpha,\beta]\vee \ker f = [\alpha \vee \ker f, \beta \vee \ker f] \vee \ker f$, so we may assume that $\alpha,\beta \ge \ker f$ without loss of generality. By Theorem \ref{day-commutator}, $[\alpha,\beta]\vee \ker f$ is the congruence generated by $X(\alpha,\beta) \cup \ker f$, and $[f(\alpha),f(\beta)]$ is the congruence generated by $X(f(\alpha),f(\beta)) = f(X(\alpha,\beta))$, so $f([\alpha,\beta]\vee \ker f) = [f(\alpha),f(\beta)]$.

Parts (f) and (g) follow directly from Theorem \ref{day-commutator}, but they can also be proved using only parts (a) - (d) (left as an exercise to the reader).
\end{proof}

Proposition \ref{mod-commutator}(d) tells us that we can compute commutators on quotients of $\bA$ directly in $\bA$. Since $\Con(\bA/\pi)$ is naturally isomorphic to the interval $\llbracket \pi, 1_\bA\rrbracket$ in $\Con(\bA)$, computing commutators on $\bA/\pi$ is equivalent to computing \emph{relative commutators} on $\bA$. Recall that if $\alpha,\beta \ge \pi$, then their \emph{relative commutator} $[\alpha,\beta]_\pi$ is defined to be the least $\delta \ge \pi$ which satisfies the term condition $C(\alpha,\beta;\delta)$.

\begin{cor} If $\alpha,\beta \ge \pi$ are congruences in a congruence modular variety, then their relative commutator is given by the formula $[\alpha,\beta]_\pi = [\alpha,\beta]\vee \pi$.
\end{cor}

\begin{thm}[Diamond Isomorphism Theorem for relative commutators]\label{diamond-isom} If $\bA$ is in a congruence modular variety and $\alpha, \beta \in \Con(\bA)$, then the maps $\phi : \llbracket \alpha,\alpha\vee \beta\rrbracket \rightarrow \llbracket \alpha\wedge \beta, \beta\rrbracket$ and $\varphi: \llbracket\alpha\wedge \beta, \beta\rrbracket \rightarrow \llbracket\alpha,\alpha\vee \beta\rrbracket$ given by
\[
\phi : x \mapsto x\wedge \beta\ \text{ and }\ \varphi : y \mapsto y\vee \alpha
\]
are lattice isomorphisms which respect the relative commutators $[\cdot,\cdot]_\alpha, [\cdot,\cdot]_{\alpha\wedge\beta}$.

Furthermore, in this case we have the equality of relative centralizers $(\alpha:\alpha\vee\beta) = (\alpha\wedge\beta:\beta)$.
\end{thm}
\begin{proof} By Theorem \ref{diamond-isom-lattice}, $\phi$ and $\varphi$ are lattice isomorphisms. 
If $\gamma,\delta \ge \alpha\wedge \beta$, then from $[\gamma \vee \alpha, \delta \vee \alpha] \le [\gamma,\delta] \vee \alpha$ we have
\[
\varphi([\gamma,\delta]_{\alpha\wedge \beta}) = [\gamma,\delta]_{\alpha\wedge \beta} \vee \alpha = [\gamma,\delta]\vee \alpha = [\gamma \vee \alpha, \delta \vee \alpha] \vee \alpha = [\varphi(\gamma),\varphi(\delta)]_\alpha.
\]
If $\gamma,\delta \in \llbracket\alpha, \alpha\vee\beta\rrbracket$, then from $\gamma = \phi(\gamma) \vee \alpha, \delta = \phi(\delta) \vee \alpha$ and $[\phi(\gamma),\phi(\delta)] \le \beta$ we have
\begin{align*}
\phi([\gamma,\delta]_\alpha) &= [\gamma,\delta]_\alpha \wedge \beta = [\phi(\gamma)\vee \alpha, \phi(\delta) \vee \alpha]_\alpha \wedge \beta\\
&= ([\phi(\gamma),\phi(\delta)]\vee \alpha) \wedge \beta = [\phi(\gamma),\phi(\delta)]\vee (\alpha \wedge \beta) = [\phi(\gamma),\phi(\delta)]_{\alpha\wedge \beta}.
\end{align*}

For the last statement, note that
\[
[\delta,\alpha\vee\beta] \le \alpha \iff [\delta,\alpha]\vee[\delta,\beta] \le \alpha \iff [\delta,\beta] \le \alpha \iff [\delta,\beta] \le \alpha \wedge \beta,
\]
so $\delta \le (\alpha : \alpha\vee \beta) \iff \delta \le (\alpha \wedge \beta : \beta)$.
\end{proof}


For the second approach to the commutator, we will follow Gumm \cite{gumm-geometric} and take the transitive closure of $\bM(\alpha,\beta)$ to produce a congruence on $\alpha$, considered as a subalgebra of $\bA^2$.

\begin{defn} For $\alpha,\beta \in \Con(\bA)$, if we consider $\alpha \le \bA^2$ as an algebra of column vectors then we can treat $\bM(\alpha,\beta)$ (from Definition \ref{commutator-matrix}) as a binary relation on $\alpha$, so $\left(\begin{bmatrix} a\\ b\end{bmatrix},\begin{bmatrix} c \\ d\end{bmatrix}\right) \in \bM(\alpha,\beta)$ means that $\begin{bmatrix} a & c\\ b & d\end{bmatrix} \in \bM(\alpha,\beta)$. We define $\Delta_\alpha^\beta$ to be the transitive closure of this binary relation on $\alpha$.
\end{defn}

Note that $\Delta_\alpha^\beta$ is the least congruence on $\alpha$ which contains the binary relation $\beta\times\beta|_{\Delta_\bA}$. When $\bA$ has a Mal'cev term, $\Delta_\alpha^\beta$ simplifies to $\bM(\alpha,\beta)$.

\begin{prop} If $\bA$ has a Mal'cev polynomial $p$ and $\alpha,\beta \in \Con(\bA)$, then $\Delta_\alpha^\beta = \bM(\alpha,\beta)$.
\end{prop}
\begin{proof} We just need to check that $\bM(\alpha,\beta)$ is transitively closed, so supposed that $\begin{bmatrix} a & c\\ b & d\end{bmatrix}, \begin{bmatrix} c & e\\ d & f\end{bmatrix} \in \bM(\alpha,\beta)$. Then we have
\[
\begin{bmatrix} a & e\\ b & f\end{bmatrix} = p\left(\begin{bmatrix} a & c\\ b & d\end{bmatrix}, \begin{bmatrix} c & c\\ d & d\end{bmatrix}, \begin{bmatrix} c & e\\ d & f\end{bmatrix}\right) \in \bM(\alpha,\beta).\qedhere
\]
\end{proof}

\begin{thm}\label{shifting-commutator} Suppose that the Shifting Lemma holds for every subalgebra of $\bA^2$. Then for $x,y \in \bA$ and $\alpha,\beta \in \Con(\bA)$, the following are equivalent:
\begin{itemize}
\item[(a)] $(x,y) \in [\beta,\alpha]$,

\item[(b)] $\begin{bmatrix} x & y\\ y & y\end{bmatrix} \in \Delta_\alpha^\beta$,

\item[(c)] there exists $a \in \bA$ such that $\begin{bmatrix} x & a\\ y & a\end{bmatrix} \in \Delta_\alpha^\beta$,

\item[(d)] there exists $b \in \bA$ such that $\begin{bmatrix} x & y\\ b & b\end{bmatrix} \in \Delta_\alpha^\beta$.
\end{itemize}
\end{thm}
\begin{proof} That (b) implies (c), (d) are clear, and (c) implies (a) directly from the term condition $C(\beta,\alpha;[\beta,\alpha])$ and the definition of $\Delta_{\alpha}^{\beta}$. That (c) implies (b) follows from the fact that $(a,y) \in \beta \implies \begin{bmatrix} a & y\\ a & y\end{bmatrix} \in \bM(\alpha,\beta)$, and since $\Delta_{\alpha}^{\beta}$ is the transitive closure of $\bM(\alpha,\beta)$ we have $\begin{bmatrix} x & y\\ y & y\end{bmatrix} \in \Delta_{\alpha}^{\beta}\circ \bM(\alpha,\beta) = \Delta_{\alpha}^{\beta}$.

For (d) $\implies$ (b) we apply the Shifting Lemma to the algebra $\alpha \le_{sd} \bA\times \bA$, the congruences $\ker \pi_1, \ker \pi_2, \Delta_{\alpha}^{\beta} \in \Con(\alpha)$, and the elements $\begin{bmatrix} x\\ b\end{bmatrix},\begin{bmatrix} y\\ b\end{bmatrix},\begin{bmatrix} x\\ y\end{bmatrix},\begin{bmatrix} y\\ y\end{bmatrix} \in \alpha$ (that $x \equiv_\alpha y$ follows from $x \equiv_\alpha b \equiv_\alpha y$).
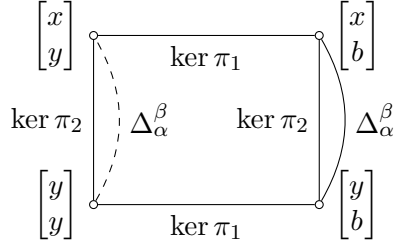
\begin{center}
\begin{tikzpicture}[scale=1.5]
\node[circle, minimum width=3pt, draw, inner sep=0pt, label=left:{$\begin{bmatrix} x\\ y\end{bmatrix}$}] (a) at (0,1.5){};
\node[circle, minimum width=3pt, draw, inner sep=0pt, label=right:{$\begin{bmatrix} x\\ b\end{bmatrix}$}] (c) at (2.0,1.5){};
\node[circle, minimum width=3pt, draw, inner sep=0pt, label=left:{$\begin{bmatrix} y\\ y\end{bmatrix}$}] (b) at (0,0){};
\node[circle, minimum width=3pt, draw, inner sep=0pt, label=right:{$\begin{bmatrix} y\\ b\end{bmatrix}$}] (d) at (2.0,0){};
\draw (a) to ["$\ker \pi_1$"'] (c) to ["$\ker \pi_2$"'] (d) to ["$\ker \pi_1$"] (b) to ["$\ker \pi_2$"] (a);
\draw [bend left] (c) to ["$\Delta_{\alpha}^{\beta}$"] (d);
\draw [bend right, dashed] (b) to ["$\Delta_{\alpha}^{\beta}$"'] (a);
\end{tikzpicture}
\end{center}

For (a) $\implies$ (b), we will show that the relation $\Theta$ defined by $(x,y) \in \Theta \iff \begin{bmatrix} x & y\\ y & y\end{bmatrix} \in \Delta_{\alpha}^{\beta}$ is a congruence which satisfies $C(\beta,\alpha;\Theta)$, which will show that $[\beta,\alpha] \le \Theta$. That $\Theta$ is reflexive is obvious, that it is symmetric follows from the equivalence of (b) with (c) or (d). If $(x,y), (y,z) \in \Theta$, then from $\begin{bmatrix} x & y\\ y & y\end{bmatrix}, \begin{bmatrix} y & z\\ y & y\end{bmatrix} \in \Delta_{\alpha}^{\beta}$ and the fact that $\Delta_{\alpha}^{\beta}$ is transitively closed, we get $\begin{bmatrix} x & z\\ y & y\end{bmatrix} \in \Delta_{\alpha}^{\beta}$, so $(x,z) \in \Theta$ by the equivalence of (b) and (d).

To finish, we just need to show that $\Theta$ satisfies $C(\beta,\alpha;\Theta)$, that is, if $\begin{bmatrix} a & c\\ b & d\end{bmatrix} \in \bM(\alpha,\beta)$ with $(c,d) \in \Theta$, then $(a,b) \in \Theta$. But if $(c,d) \in \Theta$, then $\begin{bmatrix} c & d\\ d & d\end{bmatrix} \in \Delta_{\alpha}^{\beta}$, so since $\Delta_{\alpha}^{\beta}$ is the transitive closure of $\bM(\alpha,\beta)$, we see that $\begin{bmatrix} a & d\\ b & d\end{bmatrix} \in \Delta_{\alpha}^{\beta}$, so by the equivalence of (b) with (c) we have $(a,b) \in \Theta$.
\end{proof}

\begin{cor} If every subalgebra of $\bA^2$ satisfies the Shifting Lemma, then for $\alpha, \beta_i \in \Con(\bA)$ we have $[\bigvee_i \beta_i, \alpha] = \bigvee_i[\beta_i,\alpha]$.
\end{cor}
\begin{proof} We have $\Delta_\alpha^{\bigvee_i \beta_i} = \bigvee_i \Delta_\alpha^{\beta_i}$, so $(x,y) \in [\bigvee_i \beta_i, \alpha]$ iff $\begin{bmatrix} x & z\\ y & z\end{bmatrix} \in \bigvee_i \Delta_\alpha^{\beta_i}$ for some $z$. So there must be a sequence $(x_i,y_i) \in \alpha$, $j_i$, with $\begin{bmatrix} x_i & x_{i+1}\\ y_i & y_{i+1}\end{bmatrix} \in \Delta_\alpha^{\beta_{j_i}}$ and $(x,y) = (x_n,y_n)$, $x_0 = y_0$.

We show by induction on $i$ that $\begin{bmatrix} x_i & y_i\\ y_i & y_i\end{bmatrix} \in \bigvee_{j} (\Delta_\alpha^{\beta_{j}} \wedge \ker \pi_2)$, this will show that $(x_i,y_i) \in \bigvee_j [\beta_j,\alpha]$. For the inductive step, we apply the Shifting Lemma to $\alpha \le \bA^2$ with the congruences $\ker \pi_2, \Delta_\alpha^{\beta_{j_i}}, \bigvee_{j} (\Delta_\alpha^{\beta_{j}} \wedge \ker \pi_2)$.
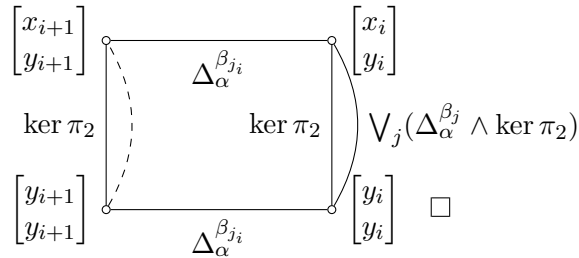
\begin{center}
\begin{tikzpicture}[scale=1.5]
\node[circle, minimum width=3pt, draw, inner sep=0pt, label=left:{$\begin{bmatrix} x_{i+1}\\ y_{i+1}\end{bmatrix}$}] (a) at (0,1.5){};
\node[circle, minimum width=3pt, draw, inner sep=0pt, label=right:{$\begin{bmatrix} x_i\\ y_i\end{bmatrix}$}] (c) at (2.0,1.5){};
\node[circle, minimum width=3pt, draw, inner sep=0pt, label=left:{$\begin{bmatrix} y_{i+1}\\ y_{i+1}\end{bmatrix}$}] (b) at (0,0){};
\node[circle, minimum width=3pt, draw, inner sep=0pt, label=right:{$\begin{bmatrix} y_i\\ y_i\end{bmatrix}$\qedhere}] (d) at (2.0,0){};
\draw (a) to ["$\Delta_\alpha^{\beta_{j_i}}$"'] (c) to ["$\ker \pi_2$"'] (d) to ["$\Delta_\alpha^{\beta_{j_i}}$"] (b) to ["$\ker \pi_2$"] (a);
\draw [bend left] (c) to ["$\bigvee_{j} (\Delta_\alpha^{\beta_{j}} \wedge \ker \pi_2)$"] (d);
\draw [bend right, dashed] (b) to (a);
\end{tikzpicture}
\end{center}
\end{proof}

\begin{cor} If every subalgebra of $\bA^2$ satisfies the Shifting Lemma, then for $f : \bA \twoheadrightarrow \bB$ surjective and $\alpha, \beta \ge \ker f$, we have $f([\beta,\alpha]\vee \ker f) = [f(\beta),f(\alpha)]$.
\end{cor}
\begin{proof} The hard direction is to check that if $(f(x), f(y)) \in [f(\beta),f(\alpha)]$, then $(x,y) \in [\beta,\alpha]\vee \ker f$. In this case we have $\begin{bmatrix} x & y\\ y & y\end{bmatrix} \in \Delta_\alpha^\beta \vee (\ker f\times \ker f|_\alpha)$. Using a similar argument to the previous corollary, we can show that this implies $\begin{bmatrix} x & y\\ y & y\end{bmatrix} \in (\Delta_\alpha^\beta \wedge \ker \pi_2) \vee (\ker f\times \ker f|_\alpha)$ by repeatedly applying the Shifting Lemma on $\alpha$. Thus we have $(x,y) \in [\beta,\alpha] \vee \ker f$ by Theorem \ref{shifting-commutator}.
\end{proof}

To prove the symmetry of the commutator, we will actually prove a stronger statement: $\Delta_\alpha^\beta$ is in fact the \emph{transpose} of $\Delta_\beta^\alpha$. In particular, if we view $\Delta_\alpha^\beta$ as a binary relation on \emph{row} vectors in $\beta$, then $\Delta_\alpha^\beta$ will be transitively closed (which is far from obvious from the definition!).

\begin{thm}\label{commutator-transpose} Suppose that the Shifting Lemma holds for every subalgebra of $\bA^4$. If $\overline{\Delta}_\alpha^\beta$ denotes the set of transposes of matrices from $\Delta_\alpha^\beta$, then $\overline{\Delta}_\alpha^\beta$ is transitively closed as a binary relation on $\beta$ and we have $\overline{\Delta}_\alpha^\beta = \Delta_\beta^\alpha$. In particular, we have $[\alpha,\beta] = [\beta,\alpha]$.
\end{thm}
\begin{proof} It's enough to prove that $\overline{\Delta}_\alpha^\beta$ is transitively closed as a binary relation on $\beta$, as we will then have $\Delta_\beta^\alpha = \bigcup_n \bM(\beta,\alpha)^{\circ n} \subseteq \overline{\Delta}_\alpha^\beta$, and a symmetric argument with $\alpha,\beta$ swapped will show that $\Delta_\alpha^\beta \subseteq \overline{\Delta}_\beta^\alpha$, so $\Delta_\beta^\alpha \subseteq \overline{\Delta}_\alpha^\beta \subseteq \Delta_\beta^\alpha$.

Suppose that $\begin{bmatrix} a & b\\ c & d\end{bmatrix}, \begin{bmatrix} c & d\\ e & f\end{bmatrix} \in \Delta_{\alpha}^\beta$. To finish, we just need to show that $\begin{bmatrix} a & b\\ e & f\end{bmatrix} \in \Delta_{\alpha}^\beta$. This follows from the following application of the Cube Lemma (Lemma \ref{mod-cube-lemma}) applied to the congruences $\ker \pi_1, \ker \pi_2, \Delta_\alpha^\beta$ on $\alpha$.
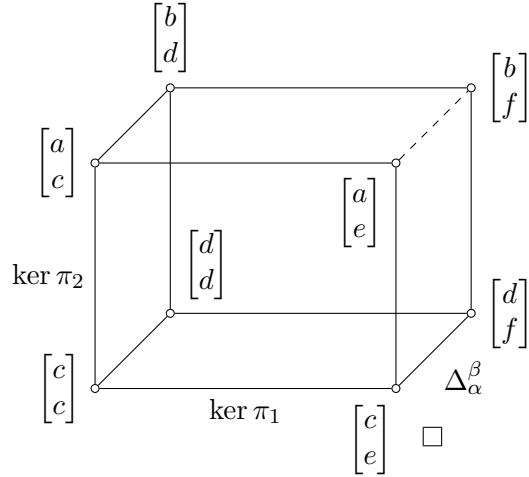
\begin{center}
\begin{tikzpicture}[xscale=1.0,yscale=1.0]
\node[circle, minimum width=3pt, draw, inner sep=0pt, label=left:{$\begin{bmatrix} c\\ c\end{bmatrix}$}] (11) at (1,1){};
\node[circle, minimum width=3pt, draw, inner sep=0pt, label=above right:{$\begin{bmatrix} d\\ d\end{bmatrix}$}] (22) at (2,2){};
\node[circle, minimum width=3pt, draw, inner sep=0pt, label=left:{$\begin{bmatrix} a\\ c\end{bmatrix}$}] (13) at (1,4){};
\node[circle, minimum width=3pt, draw, inner sep=0pt, label=above:{$\begin{bmatrix} b\\ d\end{bmatrix}$}] (24) at (2,5){};
\node[circle, minimum width=3pt, draw, inner sep=0pt, label=below:{$\begin{bmatrix} c\\ e\end{bmatrix}$\qedhere}] (31) at (5,1){};
\node[circle, minimum width=3pt, draw, inner sep=0pt, label=right:{$\begin{bmatrix} d\\ f\end{bmatrix}$}] (42) at (6,2){};
\node[circle, minimum width=3pt, draw, inner sep=0pt, label=below left:{$\begin{bmatrix} a\\ e\end{bmatrix}$}] (33) at (5,4){};
\node[circle, minimum width=3pt, draw, inner sep=0pt, label=right:{$\begin{bmatrix} b\\ f\end{bmatrix}$}] (44) at (6,5){};
\draw (11) to ["$\ker \pi_1$"'] (31);
\draw (22) -- (42);
\draw (13) -- (33);
\draw (24) -- (44);
\draw (11) to ["$\ker \pi_2$"] (13);
\draw (22) -- (24);
\draw (31) -- (33);
\draw (42) -- (44);
\draw (11) -- (22);
\draw (13) -- (24);
\draw (31) to ["$\Delta_\alpha^\beta$"'] (42);
\draw[dashed] (33) -- (44);
\end{tikzpicture}
\end{center}
\end{proof}

\begin{thm} In a congruence modular variety, any alternative commutator $[\cdot,\cdot]'$ which satisfies $[\alpha,\beta]' \le \alpha \wedge \beta$ and $f([\alpha,\beta]'\vee \ker f) = [f(\alpha),f(\beta)]'$ for $f$ surjective and $\alpha,\beta \ge \ker f$ has $[\alpha,\beta]' \le [\alpha,\beta]$ for all $\alpha, \beta$.
\end{thm}
\begin{proof} Consider congruences on $\alpha \le_{sd} \bA\times \bA$. We have $[\Delta_\alpha^\beta,\ker \pi_2]' \le \Delta_\alpha^\beta\wedge \ker \pi_2 \le \pi_1^{-1}[\beta,\alpha]$ by Theorem \ref{shifting-commutator}. Also, $\alpha = \pi_1(\ker \pi_2 \vee \ker \pi_1), \beta = \pi_1(\Delta_\alpha^\beta \vee \ker \pi_1)$, so
\[
[\beta,\alpha]' = [\pi_1(\Delta_\alpha^\beta \vee \ker \pi_1), \pi_1(\ker \pi_2 \vee \ker \pi_1)]' = \pi_1([\Delta_\alpha^\beta, \ker \pi_2]' \vee \ker \pi_1) \le [\beta,\alpha].\qedhere
\]
\end{proof}

\section{The Gumm difference term}

In this section we prove that congruence modular varieties have a ternary term $p$, called a \emph{Gumm difference term}, which acts like a Mal'cev operation on all abelian algebras. This will imply that abelian algebras in congruence modular varieties are affine.

\begin{thm}[Gumm difference term]\label{gumm-difference} For any variety with Day terms $m_0, ..., m_n$, there is a ternary term $p$ satisfying the following two properties:
\begin{itemize}
\item[(i)] $p$ satisfies the identity $p(y,y,x) \approx x$, and

\item[(ii)] for any $(x,y) \in \theta$, $\theta$ any congruence, we have $p(x,y,y)\ [\theta,\theta]\ x$.
\end{itemize}
Furthermore, in a congruence modular variety, a ternary term $p$ satisfies (i) and (ii) iff it satisfies the following property:
\begin{itemize}
\item[(iii)] for any congruences $\alpha,\beta,\gamma$ with $\alpha \wedge \beta \le \gamma$, the implication in the following picture holds.
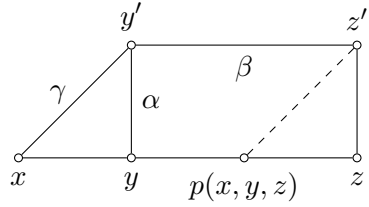
\begin{center}
\begin{tikzpicture}[scale=1.5]
\node[circle, minimum width=3pt, draw, inner sep=0pt, label=below:{$p(x,y,z)$}] (d) at (2,0){};
\node[circle, minimum width=3pt, draw, inner sep=0pt, label=below:$z$] (z) at (3,0){};
\node[circle, minimum width=3pt, draw, inner sep=0pt, label=below:$x$] (x) at (0,0){}; 
\node[circle, minimum width=3pt, draw, inner sep=0pt, label=below:$y$] (y) at (1,0){};
\node[circle, minimum width=3pt, draw, inner sep=0pt, label=above:{$z'$}] (z') at (3,1){};
\node[circle, minimum width=3pt, draw, inner sep=0pt, label=above:{$y'$}] (y') at (1,1){};
\node[circle, minimum width=3pt, inner sep=0pt, label=above:{$\gamma$}] (g) at (0.35,0.35){};
\draw (x) -- (y) -- (d) -- (z);
\draw (z') to ["$\beta$"] (y');
\draw (z) to (z');
\draw (y) to ["$\alpha$"'] (y');
\draw (x) to (y');
\draw[dashed] (d) -- (z');
\end{tikzpicture}
\end{center}
\end{itemize}
Finally, if a variety has a term $p$ which satisfies (iii), then it is congruence modular.
\end{thm}
\begin{proof} Recall the identities satisfied by Day terms:
\begin{align*}
m_0(x,y,z,w) &\approx x,\\
m_i(x,x,z,z) &\approx x\text{ for all $i$,}\\
m_i(x,y,x,y) &\approx m_{i+1}(x,y,x,y)\text{ for $i$ even,}\\
m_i(x,y,z,z) &\approx m_{i+1}(x,y,z,z)\text{ for $i$ odd,}\\
m_n(x,y,z,w) &\approx y.
\end{align*}
We inductively define a sequence of ternary terms $q_i(x,y,z)$ by $q_0(x,y,z) = z$, and
\[
q_{i+1}(x,y,z) = \begin{cases} m_{i+1}(q_{i}(x,y,z),q_{i}(x,y,z),y,x) & i\text{ odd,}\\ m_{i+1}(q_{i}(x,y,z),q_{i}(x,y,z),x,y) & i\text{ even,}\end{cases}
\]
and we set $p(x,y,z) = q_n(x,y,z)$.

To see that (i) holds, we just check inductively that $q_i(y,y,x) \approx x$:
\[
q_{i+1}(y,y,x) = m_{i+1}(q_i(y,y,x),q_i(y,y,x),y,y) \approx m_{i+1}(x,x,y,y) \approx x.
\]

For (ii), we will inductively check that
\[
q_i(y,x,x)\ [\theta,\theta]\ \begin{cases} m_i(x,y,x,y) & i\text{ even,}\\ m_i(x,y,x,x) & i\text{ odd.}\end{cases}
\]
Taking $i = n$, this will give us $p(y,x,x)\ [\theta,\theta]\ m_n(x,y,x,?) = y$.

The base case is easy: $q_0(y,x,x) = x = m_0(x,y,x,y)$. For the inductive step, we divide into cases based on whether $i$ is even or odd.

If $i$ is even, then the induction hypothesis gives
\[
q_{i+1}(y,x,x) = m_{i+1}(q_i(y,x,x),q_i(y,x,x),y,x)\ [\theta,\theta]\ m_{i+1}(m_i(x,y,x,y),m_i(x,y,x,y),y,x).
\]
Using the term condition $C(\theta,\theta;[\theta,\theta])$, from
\begin{align*}
m_{i+1}(m_i(x,y,x,y),m_i(x,y,x,y),y,\boxed{y}) &= m_i(x,y,x,y) = m_{i+1}(x,y,x,y)\\
&= m_{i+1}(m_i(x,x,x,x),m_i(y,y,y,y),x,\boxed{y}),
\end{align*}
we conclude
\begin{align*}
m_{i+1}(m_i(x,y,x,y),m_i(x,y,x,y),y,\boxed{x})\ [\theta,\theta]\ &m_{i+1}(m_i(x,x,x,x),m_i(y,y,y,y),x,\boxed{x})\\
=\ &m_{i+1}(x,y,x,x),
\end{align*}
so $q_{i+1}(y,x,x)\ [\theta,\theta]\ m_{i+1}(x,y,x,x)$.

When $i$ is odd, the proof is very similar. Inductively, we have
\[
q_{i+1}(y,x,x) = m_{i+1}(q_i(y,x,x),q_i(y,x,x),x,y)\ [\theta,\theta]\ m_{i+1}(m_i(x,y,x,x),m_i(x,y,x,x),x,y).
\]
Using the term condition $C(\theta,\theta;[\theta,\theta])$, from
\begin{align*}
m_{i+1}(m_i(x,y,x,x),m_i(x,y,x,x),x,\boxed{x}) &= m_i(x,y,x,x) = m_{i+1}(x,y,x,x)\\
&= m_{i+1}(m_i(x,x,x,x),m_i(y,y,y,y),x,\boxed{x}),
\end{align*}
we conclude
\begin{align*}
m_{i+1}(m_i(x,y,x,x),m_i(x,y,x,x),x,\boxed{y})\ [\theta,\theta]\ &m_{i+1}(m_i(x,x,x,x),m_i(y,y,y,y),x,\boxed{y})\\
=\ &m_{i+1}(x,y,x,y),
\end{align*}
so $q_{i+1}(y,x,x)\ [\theta,\theta]\ m_{i+1}(x,y,x,y)$. This conclude the proof of (ii).

Now we show that (i) and (ii) imply (iii). Suppose we have the configuration
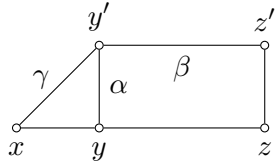
\begin{center}
\begin{tikzpicture}[scale=1.1]
\node[circle, minimum width=3pt, draw, inner sep=0pt, label=below:$z$] (z) at (3,0){};
\node[circle, minimum width=3pt, draw, inner sep=0pt, label=below:$x$] (x) at (0,0){}; 
\node[circle, minimum width=3pt, draw, inner sep=0pt, label=below:$y$] (y) at (1,0){};
\node[circle, minimum width=3pt, draw, inner sep=0pt, label=above:{$z'$}] (z') at (3,1){};
\node[circle, minimum width=3pt, draw, inner sep=0pt, label=above:{$y'$}] (y') at (1,1){};
\node[circle, minimum width=3pt, inner sep=0pt, label=above:{$\gamma$}] (g) at (0.3,0.3){};
\draw (x) -- (y) -- (z);
\draw (z') to ["$\beta$"] (y');
\draw (z) to (z');
\draw (y) to ["$\alpha$"'] (y');
\draw (x) to (y');
\end{tikzpicture}
\end{center}
with $\gamma \ge \alpha \wedge \beta$. From $x \equiv_\beta y \equiv_\beta z$, we have $p(x,y,z) \equiv_\beta z$. Additionally, we have $p(x,y,z) \equiv_\gamma p(y',y,z)$, so we just need to prove that $p(y',y,z) \equiv_\gamma z'$ to finish.

We have $p(y',y,z) \equiv_\alpha p(y',y',z') = z'$, and $p(y',y,z) \equiv_\beta p(z',z,z)$. From $(z,z') \in \alpha \wedge (\beta \vee \gamma)$, we have
\[
p(z',z,z)\ [\alpha\wedge (\beta\vee \gamma),\alpha \wedge (\beta \vee \gamma)]\ z'.
\]
The commutator above is bounded by
\[
[\alpha\wedge (\beta\vee \gamma),\alpha \wedge (\beta \vee \gamma)] \le [\alpha, \beta \vee \gamma] = [\alpha,\beta] \vee [\alpha,\gamma] \le (\alpha \wedge \beta) \vee (\alpha \wedge \gamma) = \alpha \wedge \gamma.
\]
Thus, we have $(p(y',y,z), z') \in \alpha \wedge (\beta \vee (\alpha \wedge \gamma))$, and by the modular law this is $(\alpha \wedge \beta) \vee (\alpha \wedge \gamma) = \alpha \wedge \gamma \le \gamma$.

Finally, assume that $p$ is a term which satisfies (iii). Taking $x = y = y', z = z', \alpha = \gamma = 0_\bA, \beta = 1_\bA$, we get $p(y,y,z) = z$, which is (i). Taking $x = y$ and using $p(y,y,z) = z$, we see that (iii) implies the Shifting Lemma in every algebra, so our variety is congruence modular.

To prove that (iii) implies (ii), suppose $(x,y) \in \theta$, and consider the congruences $\ker \pi_1, \ker \pi_2, \Delta_\theta^\theta$ on $\theta$. Applying (iii) in the picture
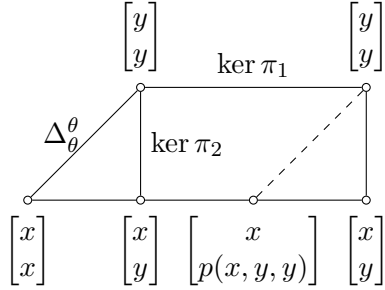
\begin{center}
\begin{tikzpicture}[scale=1.5]
\node[circle, minimum width=3pt, draw, inner sep=0pt, label=below:{$\begin{bmatrix} x \\ p(x,y,y)\end{bmatrix}$}] (d) at (2,0){};
\node[circle, minimum width=3pt, draw, inner sep=0pt, label=below:{$\begin{bmatrix} x \\ y\end{bmatrix}$}] (x) at (3,0){};
\node[circle, minimum width=3pt, draw, inner sep=0pt, label=below:{$\begin{bmatrix} x \\ x\end{bmatrix}$}] (z) at (0,0){}; 
\node[circle, minimum width=3pt, draw, inner sep=0pt, label=below:{$\begin{bmatrix} x \\ y\end{bmatrix}$}] (y) at (1,0){};
\node[circle, minimum width=3pt, draw, inner sep=0pt, label=above:{$\begin{bmatrix} y \\ y\end{bmatrix}$}] (x') at (3,1){};
\node[circle, minimum width=3pt, draw, inner sep=0pt, label=above:{$\begin{bmatrix} y \\ y\end{bmatrix}$}] (y') at (1,1){};
\node[circle, minimum width=3pt, inner sep=0pt, label=above:{$\Delta_\theta^\theta$}] (g) at (0.3,0.3){};
\draw (z) -- (y) -- (d) -- (x);
\draw (x') to ["$\ker \pi_1$"'] (y');
\draw (x) to (x');
\draw (y) to ["$\ker \pi_2$"'] (y');
\draw (z) to (y');
\draw[dashed] (d) -- (x');
\end{tikzpicture}
\end{center}
we see that $\begin{bmatrix} x & y\\ p(x,y,y) & y\end{bmatrix} \in \Delta_\theta^\theta$, so by Theorem \ref{shifting-commutator} we have $p(x,y,y)\ [\theta,\theta]\ x$.
\end{proof}

\begin{cor}[Factor Permutability] If $\bA = \bA_1\times \bA_2$ is contained in a congruence modular variety, then the factor congruences $\ker \pi_1,\ker \pi_2$ permute with every congruence $\gamma \in \Con(\bA)$.
\end{cor}
\begin{proof} A pair of congruences $\alpha,\beta \in \Con(\bA)$ correspond to a pair of factor congruences iff they satisfy $\alpha \wedge \beta = 0_\bA$ and $\alpha \circ \beta = 1_\bA$. Thus, if $x\ \gamma\ y'\ \alpha\ z'$, then by $\alpha \circ \beta = 1_\bA$ we can find $y, z \equiv_\alpha x$ with $(y,y'), (z,z') \in \beta$. Then from $\gamma \ge 0_\bA = \alpha \wedge \beta$ we can use property (iii) of a difference term to see that $x\ \alpha\ p(x,y,z)\ \gamma\ z'$, so $(x,z') \in \gamma \circ \alpha \implies (x,z') \in \alpha \circ \gamma$.
\end{proof}

\begin{cor} Any abelian algebra which is contained in a congruence modular variety is affine.
\end{cor}

\begin{cor} A nontrivial algebra $\bA$ in a congruence modular variety is abelian iff there is some $\bB \le_{sd} \bA\times \bA$ such that $\cM_3$ is a $0,1$-sublattice of $\Con(\bB)$.
\end{cor}
\begin{proof} If $\bA$ is abelian, then it is affine and we can take $\bB = \bA\times \bA$. For the other direction, it suffices to prove that $\bB$ is abelian if $\cM_3$ is a $0,1$-sublattice of $\Con(\bB)$, since then $\bB$ is affine and $\bA$ is a quotient of $\bB$, so $\bA$ is also affine.

Let $\alpha, \beta, \gamma \in \Con(\bB)$ generate a copy of $\cM_3$ which is a $0,1$-sublattice. Then
\[
[1,1] = [\alpha \vee \beta, \alpha \vee \gamma] = [\alpha,\alpha]\vee [\alpha,\gamma]\vee [\beta,\alpha]\vee [\beta,\gamma] \le \alpha \vee (\beta \wedge \gamma) = \alpha.
\]
Similarly we have $[1,1] \le \beta$, so $[1,1] \le \alpha \wedge \beta = 0$.
\end{proof}

By plugging a difference term into itself, we can strengthen property (ii) of a Gumm difference term, to get terms which act as Mal'cev operations on solvable algebras.

\begin{defn} For any congruence $\alpha$, define $[\alpha]^n$ inductively by $[\alpha]^0 = \alpha, [\alpha]^{n+1} = [[\alpha]^n,[\alpha]^n]$.
\end{defn}

\begin{prop} If $p$ is a Gumm difference term, and if we define terms $p_n$ inductively by $p_0 = p$ and
\[
p_{n+1}(x,y,z) = p_n(x,p_n(x,y,y),p_n(x,y,z)),
\]
then each $p_n$ is also a Gumm difference term, and for any $(x,y) \in \theta$ we have $p_n(x,y,y)\ [\theta]^{2^n}\ x$.
\end{prop}
\begin{proof} Inductively, we have
\[
p_{n+1}(y,y,x) = p_n(y,p_n(y,y,y),p_n(y,y,x)) = p_n(y,y,x) = x,
\]
and from $(x,p_n(x,y,y)) \in [\theta]^{2^n}$, we have
\[
p_{n+1}(x,y,y) = p_n(x,p_n(x,y,y),p_n(x,y,y))\ [[\theta]^{2^n}]^{2^n}\ x.\qedhere
\]
\end{proof}

\begin{cor}\label{cor-modular-solvable-malcev} Any solvable algebra in a congruence modular variety is Mal'cev.
\end{cor}

The last result of this section is useful for understanding the center of an algebra in terms of the difference term.

\begin{thm}\label{difference-commutator} Suppose $p$ is a Gumm difference term for a congruence modular variety and $\alpha \ge \beta$. Then $\Delta_\beta^\alpha$ is given by
\[
\begin{bmatrix} x & w\\ y & z\end{bmatrix} \in \Delta_\beta^\alpha \iff (p(x,y,z)\ [\alpha,\beta]\ w) \wedge (x\ \beta\ y\ \alpha\ z).
\]
\end{thm}
\begin{proof} If $\begin{bmatrix} x & w\\ y & z\end{bmatrix} \in \Delta_\beta^\alpha$ then clearly $(x\ \beta\ y\ \alpha\ z)$, and from
\[
p\left(\begin{bmatrix} x & x\\ x & x\end{bmatrix}, \begin{bmatrix} x & x\\ y & y\end{bmatrix}, \begin{bmatrix} x & w\\ y & z\end{bmatrix}\right) = \begin{bmatrix} x & w\\ p(x,y,y) & p(x,y,z)\end{bmatrix} \in \Delta_\beta^\alpha,
\]
we see that from $p(x,y,y)\ [\beta,\beta]\ x$, $[\beta,\beta] \le [\alpha,\beta]$, and the term condition for $[\alpha,\beta]$ we have $w\ [\alpha,\beta]\ p(x,y,z)$.

For the other direction, if $x\ \beta\ y\ \alpha\ z$ then from $\alpha \ge \beta$ we have $(x,y),(x,z)\in\alpha$, so we can apply the defining property (iii) of the difference term to congruences on $\alpha$ to see the implication in the following picture.
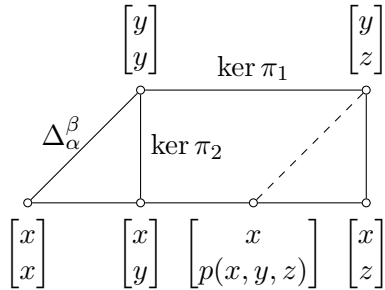
\begin{center}
\begin{tikzpicture}[scale=1.5]
\node[circle, minimum width=3pt, draw, inner sep=0pt, label=below:{$\begin{bmatrix} x \\ p(x,y,z)\end{bmatrix}$}] (d) at (2,0){};
\node[circle, minimum width=3pt, draw, inner sep=0pt, label=below:{$\begin{bmatrix} x \\ z\end{bmatrix}$}] (x) at (3,0){};
\node[circle, minimum width=3pt, draw, inner sep=0pt, label=below:{$\begin{bmatrix} x \\ x\end{bmatrix}$}] (z) at (0,0){}; 
\node[circle, minimum width=3pt, draw, inner sep=0pt, label=below:{$\begin{bmatrix} x \\ y\end{bmatrix}$}] (y) at (1,0){};
\node[circle, minimum width=3pt, draw, inner sep=0pt, label=above:{$\begin{bmatrix} y \\ z\end{bmatrix}$}] (x') at (3,1){};
\node[circle, minimum width=3pt, draw, inner sep=0pt, label=above:{$\begin{bmatrix} y \\ y\end{bmatrix}$}] (y') at (1,1){};
\node[circle, minimum width=3pt, inner sep=0pt, label=above:{$\Delta_\alpha^\beta$}] (g) at (0.3,0.35){};
\draw (z) -- (y) -- (d) -- (x);
\draw (x') to ["$\ker \pi_1$"'] (y');
\draw (x) to (x');
\draw (y) to ["$\ker \pi_2$"'] (y');
\draw (z) to (y');
\draw[dashed] (d) -- (x');
\end{tikzpicture}
\end{center}
Taking transposes, we have $\begin{bmatrix} x & p(x,y,z)\\ y & z\end{bmatrix} \in \Delta_\beta^\alpha$ by Theorem \ref{commutator-transpose}. By Theorems \ref{shifting-commutator} and \ref{commutator-transpose}, if $p(x,y,z)\ [\alpha,\beta]\ w$ then $\begin{bmatrix} p(x,y,z) & w\\ z & z\end{bmatrix} \in \Delta_\beta^\alpha$, so $\begin{bmatrix} x & w\\ y & z\end{bmatrix} \in \Delta_\beta^\alpha$ by the fact that $\Delta_\beta^\alpha$ is transitively closed.
\end{proof}

\begin{cor}\label{difference-graph} If $\alpha \ge \beta$ and $[\alpha,\beta] = 0$, then the restriction of the graph of $p(x,y,z)$ to triples with $x\ \beta\ y\ \alpha\ z$ is preserved by all polynomial operations of $\bA$.
\end{cor}

Using this, it's possible to show that if $\bA$ has center $\zeta$, then we can write $\bA$ as an extension of the quotient $\bA/\zeta$ by the abelian algebra $\zeta/\Delta_\zeta^1$ after making a choice of a section $s : \bA/\zeta \rightarrow \bA$, with each $n$-ary basic operation $f$ inducing a map $t : (\bA/\zeta)^n \rightarrow \zeta/\Delta_\zeta^1$ so that the action of $f$ on $\bA$ can be decomposed as $(x,y) \mapsto (f^{\zeta/\Delta_\zeta^1}(x) + t(y), f^{\bA/\zeta}(y))$. If $\bA$ is idempotent, then we can simplify this description slightly by noting that in this case, $\zeta/\Delta_\zeta^1$ is isomorphic to any congruence class of $\zeta$.

As a consequence of the decomposition of an algebra via its center, nilpotent algebras in congruence modular varieties turn out to be very well-behaved (e.g. they are always Mal'cev and they have regular congruences), and after selecting an element to serve as the identity, one can define an associated nilpotent loop. See Chapter 7 of Freese and McKenzie \cite{commutator-theory} for details.

\section{(Directed) J\'onsson and Gumm terms}

First we give J\'onsson's \cite{jonsson-distributive} characterization of congruence distributive varieties.

\begin{defn} A variety $\cV$ is congruence distributive if for every $\bA \in \cV$, $\Con(\bA)$ is a distributive lattice, that is, if the inequality
\[
\alpha \wedge (\beta \vee \gamma) \le (\alpha \wedge \beta) \vee (\alpha \wedge \gamma)
\]
holds for all $\alpha, \beta, \gamma \in \Con(\bA)$.
\end{defn}

The prototypical modular lattice which is \emph{not} distributive is the lattice $\cM_3$, as the next proposition shows.

\begin{prop}[Birkhoff \cite{birkhoff-lattice}]\label{distributive-m3} In any modular lattice, if $a,b,c$ do not satisfy the distributive law, and if we define elements $d,e,f$ by
\[
d = (b \wedge c) \vee (a \wedge (b \vee c)) = ((b\wedge c) \vee a)\wedge (b\vee c),
\]
with $e,f$ defined by cyclic permutations of the variables $a,b,c$ in the above formula, then $d,e,f$ generate a sublattice isomorphic to the diamond lattice $\cM_3$.
\end{prop}
\begin{proof} Using the modular law, we can check the formulas
\[
d \wedge e = e \wedge f = f \wedge d = (a\wedge b)\vee (b\wedge c) \vee (c \wedge a)
\]
and
\[
d \vee e = e \vee f = f \vee d = (a\vee b)\wedge (b\vee c) \wedge (c \vee a).
\]
If any two of $d,e,f$ are equal, then so are the two displayed expressions, and if we take the wedge of both with $a$ we get
\[
a \wedge ((a\vee b)\wedge (b\vee c) \wedge (c \vee a)) = a \wedge (b \vee c)
\]
and (using the modular law again)
\[
a\wedge ((a\wedge b)\vee (b\wedge c) \vee (c \wedge a)) = (a \wedge b) \vee (a \wedge c).\qedhere
\]
\end{proof}

\begin{prop} A variety is congruence distributive iff it is congruence modular and none of its algebras has a nontrivial abelian congruence. In particular, the commutator is given by $[\alpha,\beta] = \alpha\wedge \beta$ and $p(x,y,z) = z$ is a Gumm difference term in any congruence distributive variety.
\end{prop}
\begin{proof} If $\alpha$ is an abelian congruence, then $\ker \pi_1, \ker \pi_2, \Delta_\alpha^\alpha \in \Con(\alpha)$ generate a sublattice isomorphic to $\cM_3$, with top element $\alpha \times \alpha|_\alpha$. The other direction follows from Proposition \ref{sd-meet-commutator}, since $\cM_3$ does not satisfy the meet-semidistributive law SD($\wedge$).
\end{proof}

\begin{ex} The variety of unital rings is not congruence distributive, even though it is congruence modular (in fact, congruence permutable, since it has a Mal'cev term $x-y+z$) and contains no nontrivial abelian algebras (any such algebra would have $x\cdot y = 0$ for all $x,y$, and plugging in $y = 1$ would give $x = 0$ for all $x$). The reason for this is that the congruence on the ring $\ZZ/p^2$ corresponding to the ideal $(p)$ is abelian, but no congruence class of this ideal forms a unital subring of $\ZZ/p^2$.
\end{ex}

\begin{thm}[J\'onsson terms]\label{jonsson-terms} A variety is congruence distributive iff it has ternary terms $q_0, ..., q_n$ satisfying the system of identities
\begin{align*}
q_0(x,y,z) &\approx x,\\
q_i(x,y,x) &\approx x\text{ for all }i,\\
q_i(x,y,y) &\approx q_{i+1}(x,y,y)\text{ for }i\text{ odd,}\\
q_i(x,x,y) &\approx q_{i+1}(x,x,y)\text{ for }i\text{ even,}\\
q_n(x,y,z) &\approx z.
\end{align*}
\end{thm}
\begin{proof} Consider the congruences $\theta_{x,y}, \theta_{y,z}, \theta_{x,z}$ corresponding to identifying pairs of variables on the free algebra $\cF(x,y,z)$ in a congruence distributive variety. From $(x,z) \in \theta_{x,z}\wedge (\theta_{x,y} \vee \theta_{y,z})$ and distributivity, we have
\[
x\ (\theta_{x,z} \wedge \theta_{x,y}) \vee (\theta_{x,z} \wedge \theta_{y,z})\ z.
\]
Thus there exist $q_0, ..., q_n \in \cF(x,y,z)$ with $q_0 = x$, $q_i\ (\theta_{x,z} \wedge \theta_{x,y})\ q_{i+1}$ for $i$ even, $q_i\ (\theta_{x,z} \wedge \theta_{y,z})\ q_{i+1}$ for $i$ odd, and $q_n(x,y,z) = z$. In particular, we have $q_i\ \theta_{x,z}\ x$ for all $i$ by induction on $i$. Thus $q_0, ..., q_n$ satisfy the desired system of identities.

For the converse, suppose that $\alpha,\beta,\gamma$ are congruences on any algebra and that $(a,c) \in \alpha \wedge (\beta \vee \gamma)$. We need to show that $(a,c) \in (\alpha \wedge \beta) \vee (\alpha \wedge \gamma)$.

From $(a,c) \in \beta\vee\gamma$, there is a sequence $b_0, ..., b_m$ with $a = b_0$, $b_j\ \beta\cup\gamma\ b_{j+1}$ for all $j$, and $b_m = c$. Since $q_i(a,b_j,c)\ \alpha\ q_i(a,b_j,a) = a$ for all $i,j$, we then have
\[
q_i(a,b_j,c)\ (\alpha \wedge \beta) \cup (\alpha \wedge \gamma)\ q_i(a,b_{j+1},c)
\]
for each $i,j$, so $q_i(a,a,c)\ (\alpha \wedge \beta) \vee (\alpha \wedge \gamma)\ q_i(a,c,c)$ for all $i$. Stringing these together with the identities relating $q_i$ to $q_{i+1}$, we see that $a = q_0(a,c,c)\ (\alpha \wedge \beta) \vee (\alpha \wedge \gamma)\ q_n(a,a,c) = c$.
\end{proof}

\begin{ex} The variety of lattices is congruence distributive. For the J\'onsson terms, we may take $n = 2$ and $q_1(x,y,z)$ to be the majority term $(x\wedge y)\vee (y\wedge z) \vee (z \wedge x)$. More generally, any variety with a near-unanimity term is congruence distributive.
\end{ex}

We now prove a permutability result which is directly related to the fact that every congruence modular variety has a sequence of ternary terms known as \emph{Gumm terms}, which look like J\'onsson terms ``glued to'' a Mal'cev term.

\begin{thm}\label{commutator-permute} If $\alpha, \beta$ are any two congruences in a congruence modular variety, then
\[
\alpha \circ \beta \subseteq [\alpha,\alpha] \circ \beta \circ \alpha.
\]
\end{thm}
\begin{proof} If $(a,c) \in \alpha \circ \beta$, then there is some $b$ with $a\ \alpha\ b\ \beta\ c$. Applying the Gumm difference term $p$, we have
\[
a\ [\alpha,\alpha]\ p(a,b,b)\ \beta\ p(a,b,c)\ \alpha\ p(b,b,c) = c.\qedhere
\]
\end{proof}

\begin{cor} If $\alpha, \beta, \gamma$ are congruences in a congruence modular variety, then
\[
(\alpha \circ \beta) \cap \gamma \subseteq ((\alpha \wedge \beta)\vee (\alpha \wedge \gamma))\circ \beta \circ \alpha.
\]
\end{cor}
\begin{proof} We have $(\alpha \circ \beta) \cap \gamma = ((\alpha \wedge (\beta\vee \gamma))\circ \beta) \cap \gamma$, and
\[
[\alpha \wedge (\beta\vee \gamma),\alpha \wedge (\beta \vee \gamma)] \le [\alpha, \beta \vee \gamma] = [\alpha, \beta] \vee [\alpha, \gamma] \le (\alpha \wedge \beta) \vee (\alpha \wedge \gamma).
\]
Thus, by the previous theorem we have
\[
(\alpha \circ \beta) \cap \gamma \subseteq (\alpha \wedge (\beta\vee \gamma))\circ \beta \subseteq ((\alpha \wedge \beta)\vee (\alpha \wedge \gamma))\circ \beta \circ \alpha.\qedhere
\]
\end{proof}

A very similar argument shows that
\[
(\alpha \circ \beta) \cap \gamma \subseteq ((\alpha \wedge \gamma)\vee (\beta \wedge \gamma))\circ \beta \circ \alpha,
\]
which we will use to prove the following result (the corollary above could also be used to prove it, but there is an extra step of reordering the variables if we do it that way). Note that this containment can be viewed as a combination of a distributivity result with a permutability result.

\begin{thm}[Gumm terms]\label{gumm-terms} A variety is congruence modular iff it has ternary terms $q_0, ..., q_n, p$ satisfying the system of identities
\begin{align*}
q_0(x,y,z) &\approx x,\\
q_i(x,y,x) &\approx x\text{ for all }i,\\
q_i(x,y,y) &\approx q_{i+1}(x,y,y)\text{ for }i\text{ odd,}\\
q_i(x,x,y) &\approx q_{i+1}(x,x,y)\text{ for }i\text{ even,}\\
q_n(x,y,y) &\approx p(x,y,y),\\
p(x,x,y) &\approx y.
\end{align*}
Furthermore, a ternary term $p$ is a Gumm difference term iff there exist terms $q_0, ..., q_n$ satisfying the above system of identities.
\end{thm}
\begin{proof} Consider the congruences $\theta_{x,y}, \theta_{y,z}, \theta_{x,z}$ corresponding to identifying pairs of variables on the free algebra $\cF(x,y,z)$ in a congruence distributive variety. From $(x,z) \in \theta_{x,z}\wedge (\theta_{x,y} \vee \theta_{y,z})$ and
\[
[\theta_{x,z}\wedge (\theta_{x,y} \vee \theta_{y,z}),\theta_{x,z}\wedge (\theta_{x,y} \vee \theta_{y,z})] \le (\theta_{x,z} \wedge \theta_{x,y}) \vee (\theta_{x,z} \wedge \theta_{y,z}),
\]
which is proved as in the previous corollary, we see that for any Gumm difference term $p$ we have
\[
x\ (\theta_{x,z} \wedge \theta_{x,y}) \vee (\theta_{x,z} \wedge \theta_{y,z})\ p(x,z,z).
\]
Thus there exist $q_0, ..., q_n \in \cF(x,y,z)$ with $q_0 = x$, $q_i\ (\theta_{x,z} \wedge \theta_{x,y})\ q_{i+1}$ for $i$ even, $q_i\ (\theta_{x,z} \wedge \theta_{y,z})\ q_{i+1}$ for $i$ odd, and $q_n(x,y,z) = p(x,z,z)$. Therefore $q_0, ..., q_n, p$ satisfy the desired system of identities.

To see that Gumm terms imply congruence modularity, we just need to show that they imply the existence of Day terms. If we assume without loss of generality that $n$ is odd and take
\begin{align*}
m_0(x,y,z,w) &= x,\\
m_{2i-1}(x,y,z,w) &= q_i(x,w,y)\text{ for }i\text{ even},\\
m_{2i}(x,y,z,w) &= q_i(x,z,y)\text{ for }i\text{ even},\\
m_{2i-1}(x,y,z,w) &= q_i(x,z,y)\text{ for }i\text{ odd},\\
m_{2i}(x,y,z,w) &= q_i(x,w,y)\text{ for }i\text{ odd},\\
m_{2n+1}(x,y,z,w) &= p(z,w,y),\\
m_{2n+2}(x,y,z,w) &= y,
\end{align*}
then we have $m_i(x,x,z,z) \approx x$ for all $i$, $m_i(x,y,x,y) \approx m_{i+1}(x,y,x,y)$ for $i$ even, and $m_i(x,y,z,z) \approx m_{i+1}(x,y,z,z)$ for $i$ odd, so $m_0, ..., m_{2n+2}$ are Day terms.

To show that any such $p$ is a Gumm difference term, we just need to show that if $(x,y) \in \theta$, then $p(x,y,y)\ [\theta,\theta]\ x$. We will show by induction that $q_i(x,y,y)\ [\theta,\theta]\ x$ for all $i$. For the inductive step, we just need to show that for all $i$, we have $q_i(x,y,y)\ [\theta,\theta]\ q_i(x,x,y)$. This follows from the term condition for $[\theta,\theta]$:
\[
q_i(x,y,\boxed{x}) = q_i(x,x,\boxed{x}) \implies q_i(x,y,\boxed{y})\ [\theta,\theta]\ q_i(x,x,\boxed{y}).\qedhere
\]
\end{proof}

The need to constantly divide into cases for even vs. odd $i$ can be eliminated by the main result of \cite{directed-gumm}, which establishes the existence of \emph{directed} J\'onsson and Gumm terms. The idea behind the directed variants is that if we have idempotent ternary terms $f,g$ which satisfy
\[
f(x,y,y) \approx g(x,x,y),
\]
then they can also be indirectly connected by a ternary term $h$ which satisfies $h(x,y,x) \approx x$ and joins $f,g$ by $f\ \theta_{y,z}\ h\ \theta_{x,y}\ g$, that is,
\begin{align*}
f(x,y,y) &\approx h(x,y,y),\\
h(x,x,y) &\approx g(x,x,y).
\end{align*}
In fact, we can just take $h(x,y,z) = f(x,z,z)$: then we will have $h(x,y,y) = h(x,x,y) = f(x,y,y) = g(x,x,y)$, and $h(x,y,x) = f(x,x,x) = x$. The goal of the directed J\'onsson and Gumm terms is to cut out the middleman $h$, to obtain a substantially stronger system of identities.

Another reason to prefer the directed equations $f_i(x,y,y) \approx f_{i+1}(x,x,y)$ is that they have a clearer connection to higher arity terms, especially near-unanimity terms. Suppose that $\phi$ is an $n$-ary operation, and define terms $f_i$ by
\[
f_i(x,y,z) = \phi(x,...,x,y,z,...,z),
\]
where the lone $y$ occurs in the $i$-th position from the right (so there are $i-1$ $z$s). Then the $f_i$ will automatically satisfy
\[
f_i(x,y,y) = \phi(x,...,x,y,y,...,y) = f_{i+1}(x,x,y),
\]
and if $\phi$ is idempotent they will satisfy $f_1(x,x,y) \approx x$ and $f_n(x,y,y) \approx y$. Finally, $\phi$ will be a near-unanimity term iff each $f_i$ satisfies $f_i(x,y,x) \approx x$.

\begin{thm}[Directed Gumm terms \cite{directed-gumm}]\label{directed-gumm-terms} A variety is congruence modular iff it has ternary terms $f_1, ..., f_m, p$ with
\begin{align*}
f_1(x,x,y) &\approx x,\\
f_i(x,y,x) &\approx x\text{ for all }i,\\
f_i(x,y,y) &\approx f_{i+1}(x,x,y)\text{ for all }i,\\
f_m(x,y,y) &\approx p(x,y,y),\\
p(x,x,y) &\approx y,
\end{align*}
and if the variety is congruence distributive then we can take $f_m(x,y,y) \approx y$ (directed J\'onsson terms).
\end{thm}
\begin{proof} Assume without loss of generality that our variety is idempotent. Suppose that there are Gumm terms $q_1, ..., q_{2k+1}, p_1$ with
\begin{align*}
q_1(x,x,y) &\approx x,\\
q_i(x,y,x) &\approx x\text{ for all }i,\\
q_{2i-1}(x,y,y) &\approx q_{2i}(x,y,y)\text{ for all }i,\\
q_{2i}(x,x,y) &\approx q_{2i+1}(x,x,y)\text{ for all }i,\\
q_{2k+1}(x,y,y) &\approx p_1(x,y,y),\\
p_1(x,x,y) &\approx y.
\end{align*}
Let $\cF$ be the free algebra on $x,y$. Let $\rightsquigarrow$ be the transitive closure of the binary relation on $\cF$ generated by $x\rightsquigarrow x, x\rightsquigarrow y, y\rightsquigarrow y$, so binary terms $a(x,y),b(x,y)$ have $a \rightsquigarrow b$ iff there is a sequence of ternary terms $t_i$ with $t_1(x,x,y) = a(x,y)$, $t_i(x,y,y) = t_{i+1}(x,x,y)$, and $t_n(x,y,y) = b(x,y)$.

Additionally, let $\rightarrow$ be the relation on $\cF$ with $a \rightarrow b$ iff there is a sequence of ternary terms $t_i$ with $t_1(x,x,y) = a(x,y)$, $t_i(x,y,y) = t_{i+1}(x,x,y)$, $t_n(x,y,y) = b(x,y)$, and additionally $t_i(x,y,x) = x$ for all $i$. Then for any ternary term $q$ satisfying $q(x,y,x) = x$, we have
\[
q(\rightarrow,\rightsquigarrow,\rightarrow) \subseteq \rightarrow.
\]

For any binary term $a(x,y)$, we define $a^n(x,y)$ recursively by $a^0(x,y) = y, a^1(x,y) = a(x,y)$, and
\[
a^{n+1}(x,y) = a(x,a^n(x,y))
\]
for each $n$.

Setting $b_k(x,y) = q_{2k+1}(x,y,y) = p_1(x,y,y)$, our goal will be to prove that
\[
\exists b \in \cF\;\;\; x \rightarrow b_k^{2^k}(b(x,y),b_k^{2^k-1}(x,y)).
\]
It will then be easy to construct a ternary term $p$ with $p(x,y,y) = b_k^{2^k}(b,b_k^{2^k-1})$ and $p(x,x,y) = y$, by recursively plugging $p_1$ into itself in a similar way to the way we constructed Mal'cev terms on solvable algebras.

{\bf Claim 1:} If $a \rightsquigarrow b$ and $c(x,y) \rightarrow d(x,y)$, then $c(a,b) \rightarrow d(a,b)$.

{\bf Proof of Claim 1:} We just have to check this in the case where $c \rightarrow d$ in one step. So suppose that $t(x,x,y) = c(x,y), t(x,y,y) = d(x,y), t(x,y,x) = x$. Then
\[
\begin{bmatrix} c(a,b)\\ d(a,b)\end{bmatrix} = t\left(\begin{bmatrix} a\\ a\end{bmatrix}, \begin{bmatrix} a\\ b\end{bmatrix}, \begin{bmatrix} b\\ b\end{bmatrix}\right) \in t(\rightarrow,\rightsquigarrow,\rightarrow) \subseteq \rightarrow.
\]

{\bf Claim 1.5:} If $a \rightsquigarrow b$ and $c(x,y) \leftarrow d(x,y)$, then $c(b,a) \rightarrow d(b,a)$.

{\bf Proof of Claim 1.5:} This follows from Claim 1 and the fact that $c(x,y) \leftarrow d(x,y) \iff c(y,x) \rightarrow d(y,x)$.

{\bf Claim 2:} If $a \rightarrow b$, then $a^n \rightarrow b^n$ for every $n$.

{\bf Proof of Claim 2:} Induct on $n$. For the inductive step, we have
\[
a^{n+1}(x,y) = a(x,a^n(x,y)) \rightarrow b(x,a^n(x,y)) \rightarrow b(x,b^n(x,y)) = b^{n+1}(x,y),
\]
where the first $\rightarrow$ follows from $x = a^n(x,x) \rightsquigarrow a^n(x,y)$ and Claim 1, while the second $\rightarrow$ follows from the fact that $\rightarrow$ is preserved by $b$ and the inductive hypothesis.

The sequence of Gumm terms $q_1, ..., q_{2k+1}$ gives us a $k$-\emph{fence}:
\[
x = a_0 \rightarrow b_0 \leftarrow a_1 \rightarrow b_1 \leftarrow a_2 \rightarrow \cdots \leftarrow a_k \rightarrow b_k,
\]
where $a_i(x,y) = q_{2i+1}(x,x,y) = q_{2i}(x,x,y)$, $b_i(x,y) = q_{2i+1}(x,y,y) = q_{2i+2}(x,y,y)$. Our strategy will be to use Claims 1 and 1.5 to iteratively reduce the length of the fence.

{\bf Claim 3:} If $x \rightarrow b \leftarrow a \rightarrow c$ is a $1$-fence, then $x \rightarrow b_k(b,c(b,c))$.

{\bf Proof of Claim 3:} We define a sequence of terms $d_i$ by $d_0 = x$ and
\[
d_{i+1} = b(d_i,a),
\]
and define terms $e_i$ by
\[
e_i = a(d_i,a).
\]
We claim that for each $i$ we have
\begin{itemize}
\item $d_i \rightsquigarrow a$, $d_i \rightarrow d_{i+1}$, $d_i\rightsquigarrow e_i$, $d_i \rightarrow b$,
\item $e_i\rightarrow d_{i+1}$, $e_i \rightarrow e_{i+1}$, $e_i \rightarrow c(b,c)$.
\end{itemize}
\begin{center}
\begin{tikzpicture}[scale=1]
  \node (d0) at (0,1) {$x=d_0$};
  \node (d1) at (2,1) {$d_1$};
  \node (d2) at (4,1) {$d_2$};
  \node (b) at (6,1) {$b$};
  \node (e0) at (1,0) {$e_0$};
  \node (e1) at (3,0) {$e_1$};
  \node (e2) at (5,0) {$e_2$};
  \node (cbc) at (7,0) {$c(b,c)$};
  \draw [->, line join=round, decorate, decoration={zigzag, segment length=4, amplitude=.9,post=lineto, post length=2pt}] (d0) -- (e0);
  \draw [->, line join=round, decorate, decoration={zigzag, segment length=4, amplitude=.9,post=lineto, post length=2pt}] (d1) -- (e1);
  \draw [->, line join=round, decorate, decoration={zigzag, segment length=4, amplitude=.9,post=lineto, post length=2pt}] (d2) -- (e2);
  \draw [->] (d0) edge (d1) (d1) edge (d2) (d2) edge (b);
  \draw [->] (e0) edge (e1) (e1) edge (e2) (e2) edge (cbc);
  \draw [->] (e0) edge (d1) (e1) edge (d2);
\end{tikzpicture}
\end{center}
To see this, note first that $d_0 = x \rightsquigarrow a$, so by induction on $i$ we have $d_{i+1} = b(d_i,a) \rightsquigarrow b(a,a) = a$ for each $i$. So from $x \rightarrow b \leftarrow a$ we get $d_i \rightarrow b(d_i,a) \leftarrow a(d_i,a)$ by Claim 1, that is, $d_i \rightarrow d_{i+1} \leftarrow e_i$ for each $i$.

Then we have $e_i = a(d_i,a) \rightarrow a(d_{i+1},a) = e_{i+1}$ for each $i$, and $d_i = a(d_i,d_i) \rightsquigarrow a(d_i,a) = e_i$ for each $i$. This finishes up all of the arrows other than the rightmost two in the picture.

For $d_i \rightarrow b$, note that $d_0 = x \rightarrow b$ by assumption, and $d_{i+1} = b(d_i,a) \rightarrow b(b,b) = b$ inductively. Finally, for each $i$ we have
\[
e_i = a(d_i,a) \rightarrow c(d_i,a) \rightarrow c(b,c),
\]
where the first arrow follows from Claim 1.

Now we can use all these arrows to see that
\[
x = d_0 = a_0(d_0,e_0) \rightarrow b_0(d_0,e_0) \rightarrow b_0(d_1,e_0) \rightarrow a_1(d_1,e_0) \rightarrow a_1(d_1,e_1) \rightarrow b_1(d_1,e_1) \rightarrow \cdots,
\]
where we have used Claim 1 and Claim 1.5 several times. Chaining these together, we get
\[
x \rightarrow b_k(d_k,e_k) \rightarrow b_k(b,c(b,c)).
\]
This completes the proof of Claim 3.

{\bf Claim 4:} For each $i < k$, there is a $k-i$-fence
\[
x \rightarrow b_{0,i} \leftarrow a_{1,i} \rightarrow b_{1,i} \leftarrow a_{2,i} \rightarrow \cdots \leftarrow a_{k-i,i} \rightarrow b_{k-i,i} = b_k^{2^{i+1}-1}.
\]

{\bf Proof of Claim 4:} We prove this by induction on $i$. The base case $i=0$ comes from the Gumm terms. Suppose it is known for $i$, then by Claim 3 we have
\[
x \rightarrow b_k(b_{0,i},b_{1,i}(b_{0,i},b_{1,i})),
\]
and from $b_{0,i} \leftarrow x$ we have
\[
b_k(b_{0,i},b_{1,i}(b_{0,i},b_{1,i})) \leftarrow b_k(x,b_{1,i}(x,b_{1,i})) = b_k(x,b_{1,i}^2).
\]
By Claim 2, we have $b_{1,i}^2 \leftarrow a_{2,i}^2 \rightarrow b_{2,i}^2 \leftarrow \cdots$, so if we take
\[
b_{0,i+1} = b_k(b_{0,i},b_{1,i}(b_{0,i},b_{1,i}))
\]
and
\[
a_{j,i+1} = b_k(x,a_{j+1,i}^2), \;\;\; b_{j,i+1} = b_k(x,b_{j+1,i}^2),
\]
we get
\[
x \rightarrow b_{0,i+1} \leftarrow a_{1,i+1} \rightarrow b_{1,i+1} \leftarrow a_{2,i+1} \rightarrow \cdots \leftarrow a_{k-i-1,i+1} \rightarrow b_{k-i-1,i+1},
\]
and
\[
b_{k-i-1,i+1} = b_k(x,b_{k-i,i}^2) = b_k(x,(b_k^{2^{i+1}-1})^2) = b_k^{2^{i+2}-1}.
\]
This completes the proof of Claim 4.

By Claim 4 applied with $i = k-1$, we get a $1$-fence
\[
x \rightarrow b_{0,k-1} \leftarrow a_{1,k-1} \rightarrow b_{1,k-1} = b_k^{2^k-1}.
\]
Applying Claim 3, we get
\[
x \rightarrow b_k(b_{0,k-1},b_k^{2^k-1}(b_{0,k-1},b_k^{2^k-1})) = b_k^{2^k}(b_{0,k-1},b_k^{2^k-1}).
\]
Letting $b = b_{0,k-1}$, we see that we have succeeded in showing that $x \rightarrow b_k^{2^k}(b,b_k^{2^k-1})$. Thus there exist ternary terms $f_i$ with
\begin{align*}
f_1(x,x,y) &\approx x,\\
f_i(x,y,x) &\approx x\text{ for all }i,\\
f_i(x,y,y) &\approx f_{i+1}(x,x,y)\text{ for all }i,\\
f_m(x,y,y) &\approx b_k^{2^k}(b(x,y),b_k^{2^k-1}(x,y))).
\end{align*}
Note that if $b_k(x,y) = y$, then we also have $b_k^{2^k}(b(x,y),b_k^{2^k-1}(x,y))) = y$, so the above becomes a sequence of directed J\'onsson terms.

To finish, we just need to construct $p$ with $p(x,y,y) = b_k^{2^k}(b(x,y),b_k^{2^k-1}(x,y)))$ and $p(x,x,y) = y$. Recall that $p_1$ satisfied $p_1(x,y,y) = b_k(x,y)$ and $p_1(x,x,y) = y$. We construct terms $p_i$ inductively. For $2 \le i+1 < 2^k$, we set
\[
p_{i+1}(x,y,z) = p_1(x,p_i(x,y,y),p_i(x,y,z)),
\]
and for $2^k \le i+1$, we set
\[
p_{i+1}(x,y,z) = p_1(b(x,y),p_i(x,y,y),p_i(x,y,z)),
\]
and finally we set $p(x,y,z) = p_{2^{k+1}-1}(x,y,z)$.
\end{proof}


\section{Subdirectly irreducible algebras, ultraproducts, and residually small varieties}\label{s-subdirectly-irred}

In this section, we go over the proof of an extension of J\'onsson's Lemma \cite{jonsson-distributive}, which shows that subdirectly irreducible algebras in a finitely generated congruence distributive variety have bounded size, to the congruence modular case. The key technical tool is the concept of an ultraproduct, and the fact that any ultrapower of a finite algebra $\bA$ is isomorphic to $\bA$.

Before we discuss ultraproducts, we first review some basic results about subdirect representations of algebras due to Birkhoff \cite{birkhoff-subdirect}. The following result is elementary.

\begin{prop} If $\bA \le_{sd} \prod_{i \in I} \bA_i$ is a subdirect product, then $\bigwedge_{i \in I} \ker \pi_i = 0_\bA$. In particular, if no $\pi_i$ is an isomorphism then the congruence $0_\bA$ can be written as a meet of some family of nontrivial congruences.

Conversely, if $0_\bA$ can be written as a meet of congruences $\alpha_i \in \Con(\bA)$ for $i \in I$, then $\bA \le_{sd} \prod_{i \in I} \bA/\alpha_i$.
\end{prop}

\begin{defn} An algebraic structure $\bA$ is \emph{subdirectly irreducible} if every way of writing $\bA$ as a subdirect product $\bA \le_{sd} \prod_{i \in I} \bA_i$ has at least one coordinate $i$ such that the projection map $\pi_i : \bA \rightarrow \bA_i$ is an isomorphism. The least nontrivial congruence on a subdirectly irreducible algebra is called its \emph{monolith}.
\end{defn}

The preceeding proposition can now be rephrased as saying that $\bA$ is subdirectly irreducible iff $0_\bA$ is \emph{meet-irreducible}.

\begin{defn} An element $\alpha$ of a complete lattice $\cL$ is \emph{meet-irreducible} if for any set of elements $\alpha_i \in \cL$ with $\bigwedge_{i \in I} \alpha_i = \alpha$, some $\alpha_i$ is equal to $\alpha$. In this case, we define the \emph{cover} of $\alpha$, written $\alpha^*$, to be the least element of $\cL$ with $\alpha < \alpha^*$.
\end{defn}

In particular, the monolith of a subdirectly irreducible algebra is the cover $0_\bA^*$ of $0_\bA$.

\begin{thm}[Birkhoff's Subdirect Representation Theorem] Any algebraic structure $\bA$ can be represented as a subdirect product of subdirectly irreducible algebras.
\end{thm}
\begin{proof} For any $a\ne b \in \bA$, Zorn's Lemma implies that there is a maximal congruence $\theta_{a,b}'$ such that $(a,b) \not\in \theta_{a,b}'$. Any such $\theta_{a,b}'$ is necessarily meet-irreducible, since any congruence which properly contains $\theta_{a,b}'$ necessarily contains $(a,b)$, and therefore contains the congruence generated by $\theta_{a,b}'$ and the pair $(a,b)$.

Since we clearly have $0_\bA = \bigwedge_{a \ne b} \theta_{a,b}'$, we have the subdirect representation $\bA \le_{sd} \prod_{a \ne b} \bA/\theta_{a,b}'$.
\end{proof}

Birkhoff's subdirect representation theorem has a purely lattice-theoretic generalization to \emph{algebraic} lattices.

\begin{defn}\label{defn-algebraic-lattice} An element $\alpha$ of a complete lattice is called \emph{compact} if for any family $\alpha_i$ such that $\alpha \le \bigvee_{i \in I} \alpha_i$, there is some finite subset $\{i_1, ..., i_k\} \subseteq I$ such that $\alpha \le \alpha_{i_1} \vee \cdots \vee \alpha_{i_k}$. A complete lattice is called \emph{algebraic} if every element can be written as a join of compact elements.
\end{defn}

Every congruence lattice $\Con(\bA)$ is an algebraic lattice, since for any $a,b \in \bA$ the congruence $\theta_{a,b}$ generated by $(a,b)$ is compact, and every congruence is a join of such congruences.

\begin{prop}\label{meet-irreducible-rep} Let $\cL$ be an algebraic lattice. Then every element $\alpha$ of $\cL$ can be written as a meet of some family of meet-irreducible elements of $\cL$.
\end{prop}
\begin{proof} Let $\theta$ be any compact element of $\cL$ with $\alpha \not\ge \theta$. By Zorn's Lemma and the compactness of $\theta$, there is some $\theta' \ge \alpha$ which is maximal such that $\theta' \not\ge \theta$, and this $\theta'$ is necessarily meet-irreducible with cover $\theta' \vee \theta$. Then $\bigwedge_{\theta \not\le \alpha} \theta'$ is $\ge \alpha$, and is not $\ge$ any compact element $\theta$ with $\alpha \not\ge \theta$, so it must be equal to $\alpha$.
\end{proof}

\begin{cor}\label{meet-irreducible} If $\alpha < \beta$ in an algebraic lattice, then there is a meet-irreducible $\gamma$ such that $\gamma \ge \alpha$ but $\gamma \not\ge \beta$.
\end{cor}

Now we can briefly discuss ultrafilters and ultraproducts before moving on to the main result of this section.

\begin{defn} If $I$ is a set, then a collection of subsets $\cU \subseteq \cP(I)$ is a \emph{filter} if $\cU$ does not contain $\emptyset$, $U,V \in \cU \implies U\cap V \in \cU$, and $U \subseteq V, U \in \cU \implies V \in \cU$. We say that $\cU$ is an \emph{ultrafilter} if additionally for every $U \subseteq I$, one of $U, I\setminus U$ is in $\cU$.
\end{defn}

\begin{prop} Any filter is contained in an ultrafilter.
\end{prop}
\begin{proof} We apply Zorn's Lemma to see that any filter is contained in a maximal filter. To finish, we just need to show that any maximal filter is an ultrafilter. Suppose that $U, I\setminus U \not\in \cU$, and let $\cU'$ be the collection of $V \subseteq I$ such that $V\cup U \in \cU$. Then $\cU'$ is a filter which strictly contains $\cU$.
\end{proof}

\begin{defn} If $\bA_i$ is a collection of structures which share a common signature $\sigma$ and are indexed by $i \in I$, and if $\cU$ is an ultrafilter on $I$, then we define the \emph{ultraproduct} $\prod_i \bA_i/\cU$ to be the quotient of $\prod_i \bA_i$ by the congruence defined by
\[
a \equiv_\cU b \iff \{i \mid a_i = b_i\} \in \cU.
\]
That $\equiv_\cU$ is compatible with functions $f \in \sigma$ follows from the fact that $\cU$ is a filter. If $R \in \sigma$ is an $m$-ary relation, then $R$ is interpreted on $\prod_i \bA_i/\cU$ by
\[
R(a^1/\cU, ..., a^m/\cU) \iff \{i \mid R(a^1_i, ..., a^m_i)\} \in \cU.
\]
If all the $\bA_i$ are isomorphic to $\bA$, then we call $\bA^I/\cU$ an \emph{ultrapower} of $\bA$.
\end{defn}

Note that in terms of the congruence lattice $\Con(\prod_i \bA_i)$, the congruence $\equiv_\cU$ is equal to the join
\[
\bigvee_{U \in \cU} \ker \pi_U,
\]
where $\pi_U : \prod_{i\in I} \bA_i \rightarrow \prod_{i \in U} \bA_i$ is projection onto the coordinates in $U$. That this join is equal to the union $\bigcup_{U \in \cU} \ker \pi_U$ follows from the fact that $\cU$ is a filter.

\begin{prop} If $\cU$ is an ultrafilter on $I$ and $U_1, ..., U_k$ partition $I$ into $k$ disjoint parts, then exactly one of $U_1, ..., U_k$ is in $\cU$.
\end{prop}

\begin{cor}\label{ultraproduct-finite} If $|\bA_i| \le n$ for all $i \in I$, then $|\prod_i \bA_i/\cU| \le n$ as well. If each $\bA_i$ is finite and only finitely many isomorphism classes occur among the $\bA_i$, then $\prod_i \bA_i/\cU$ is isomorphic to some $\bA_i$.
\end{cor}

In fact, much more is true about ultraproducts, and the corollary above also follows from the following result from model theory.

\begin{thm}[{\L}o\'s's Theorem] Let $\varphi(x_1, ..., x_n)$ be any first order formula in the signature $\sigma$ with parameters $x_1, ..., x_n$, then for any $a^1, ..., a^n \in \prod_{i \in I} \bA_i$ and any ultrafilter $\cU$ on $I$, we have
\[
\prod_i \bA_i/\cU \models \varphi(a^1/\cU, ..., a^n/\cU) \iff \{i \mid \bA_i \models \varphi(a^1_i, ..., a^n_i)\} \in \cU.
\]
\end{thm}
\begin{proof} If $\varphi$ is atomic, then this follows directly from the definitions. Otherwise, $\varphi$ can be built up from atomic formulas via $\neg, \wedge, \exists$, and we can induct on the structure of $\varphi$: for $\neg$, we use the ultrafilter property that exactly one of $U, I\setminus U$ is in $\cU$ for each $U$, for $\wedge$ we use the filter property that intersections of sets in $\cU$ are in $\cU$, and for $\exists$ we just need the fact that supersets of sets in $\cU$ are in $\cU$.
\end{proof}

Now for the main result. We extend Birkhoff's $H,S,P$ notation by the operation $P_u$, where $P_u(\{\bA_i\})$ is the collection of ultraproducts of the $\bA_i$s. Recall that for $\beta$ a congruence, the centralizer $(0:\beta)$ of $\beta$ is defined as the largest $\alpha$ such that $[\alpha,\beta] = 0$, and more generally $(\delta:\beta)$ is defined as the largest $\alpha$ such that $[\alpha,\beta] \le \delta$.

\begin{thm}\label{jonsson-modular} Let $\{\bA_i\}$ be a family of algebras, and let $\cV = \cV(\{\bA_i\})$ be the variety they generate. If $\cV$ is congruence modular, $\bB \in \cV$ is subdirectly irreducible, and $\alpha = (0_\bB : 0_\bB^*)$ is the centralizer of the monolith $0_\bB^*$ of $\bB$, then $\bB/\alpha$ is a homomorphic image of a subalgebra of an ultraproduct of the $\bA_i$s, that is, $\bB/\alpha \in HSP_u(\{\bA_i\})$.
\end{thm}
\begin{proof} (From \cite{commutator-theory}, where a stronger statement is proved.) By Birkhoff's HSP Theorem, we can write $\bB = \bC/\theta$ for $\bC \le \prod_i \bA_i$. Then $\bB$ will be subdirectly irreducible iff $\theta$ is meet-irreducible in $\Con(\bC)$, so $\theta$ will have a cover $\theta^*$. The preimage $\varphi$ of $\alpha$ under $\bC \rightarrow \bB$ is the largest congruence on $\bC$ such that $[\varphi,\theta^*] \le \theta$ (i.e. $\varphi = (\theta:\theta^*)$), and we have $\bB/\alpha = \bC/\varphi$.

The main step of the proof is the following {\bf claim:} if $\beta \wedge \gamma \le \theta$ but $\gamma \not\le \theta$, then $\beta \le \varphi$.

{\bf Proof of claim:} We have
\[
[\beta,\theta^*] \le [\beta,\gamma\vee\theta] = [\beta,\gamma]\vee [\beta,\theta] \le (\beta \wedge \gamma) \vee \theta = \theta,
\]
so $\beta \le \varphi$ by $\varphi = (\theta:\theta^*)$.

Using the claim, we can now argue as follows: let $\cF$ be a maximal filter such that $U \in \cF$ implies $\ker \pi_U \le \theta$, and let $\cU$ be any ultrafilter which extends $\cF$. Then for any $U \in \cU$, we were unable to adjoin its complement to $\cF$, so there is some $V \in \cF$ such that $\ker \pi_{V\setminus U} \not\le \theta$. Then
\[
\ker \pi_U \wedge \ker \pi_{V\setminus U} = \ker \pi_{U\cup V} \le \ker \pi_V \le \theta,
\]
so by the claim we have $\ker \pi_U \le \varphi$. Thus the congruence $\bigvee_{U \in \cU} \ker \pi_U$ corresponding to $\cU$ is also $\le \varphi$, and we see that $\bB/\alpha = \bC/\varphi$ is a quotient of $\bC/\cU \le \prod_i \bA_i/\cU$.
\end{proof}

\begin{cor}[J\'onsson's Lemma \cite{jonsson-distributive}] Let $\{\bA_i\}$ be a family of algebras, and let $\cV = \cV(\{\bA_i\})$ be the variety they generate. If $\cV$ is congruence distributive and $\bB \in \cV$ is subdirectly irreducible, then $\bB \in HSP_u(\{\bA_i\})$. In particular, if $\{\bA_i\}$ is a finite set of finite algebras, then $\bB \in HS(\{\bA_i\})$.
\end{cor}

\begin{cor} For any two finite subdirectly irreducible algebras $\bA, \bB$ with the same signature which generate congruence distributive varieties, we have $\bA \cong \bB$ iff the set of identities that hold in $\bA$ is the same as the set of identities that hold in $\bB$.
\end{cor}

\begin{ex} Consider the variety of distributive lattices, and the two-element lattice $(\{0,1\},\max,\min)$. It is easy to see that every identity that holds in the two-element lattice is implied by the lattice axioms together with distributivity (since these allow us to put every term into conjunctive normal form), so the variety of distributive lattices is generated by the two-element lattice.

By J\'onsson's Lemma, the only subdirectly irreducible distributive lattice is the two-element lattice itself, so we see that in fact every distributive lattice is a sublattice of $\{0,1\}^I$ for some index set $I$, that is, every distributive lattice is a sublattice of the lattice of subsets of some set $I$.
\end{ex}

In order to understand subdirectly irreducible algebras in congruence modular varieties, we need to combine the above results with an understanding of subdirectly irreducible modules over rings.

\begin{prop} Let $\bG,\bM$ be abelian groups and let $\RR$ be a finite subgroup of $\Hom(\bG,\bM)$, such that there is a nonzero element $a \in \bM$ so that for all $x \in \bG$ there is an $r \in \RR$ with $rx = a$. Then $|\bG|$ is a prime power dividing $|\RR|$.
\end{prop}
\begin{proof} First we show that $\bG$ is finite, following \cite{commutator-theory}. Let $r_1, ..., r_k$ be the nontrivial elements of $\RR$.

We will show by induction on $k$ that $|\bG| \le (k+1)!$. For the base case, if $k = 0$ then $\bG$ can have no nonzero elements, so $|\bG| = 1 = (k+1)!$. For the inductive step, note that by the pigeonhole principle there is some $r_i$ such that at least $\frac{|\bG|-1}{k}$ elements are mapped to $a$ by $r_i$, so $|\ker r_i| \ge \frac{|\bG|-1}{k}$ (this is the ordinary group theoretic kernel), and every nonzero element of $\ker r_i$ can be mapped to $a$ by some $r_j$ with $j \ne i$, so $|\ker r_i| \le k!$ by the induction hypothesis. Thus $|\bG| \le 1 + k\cdot k! \le (k+1)!$.

Now that we know that $\bG$ is finite, we know that every element of $\bG$ has finite order, so some element $x$ has order $p$ for some prime $p$. Then there is some $r \in \RR$ with $rx = a$, so $a$ must also have order $p$. From this argument, we see that every element of $\bG$ must have order a power of $p$, so $|\bG|$ is also a power of $p$.

We may assume without loss of generality that $\bM$ is generated by the image of $\bG$ under all elements of $\RR$, so in particular that $\bM$ is finite. Then there exists an element $\pi \in \hat{\bM} = \Hom(\bM,\QQ/\ZZ)$ such that $\pi(a) \ne 0$.

Define a linear map $\phi : \RR \rightarrow \hat{\bG} = \Hom(\bG,\QQ/\ZZ)$ by $\phi : r \mapsto \phi_r$, where $\phi_r$ is the linear map $\phi_r: x \mapsto \pi(rx)$. Then $\phi$ must be surjective, or else the image will be a proper subgroup of $\hat{\bG}$ and so there will be some nonzero $x \in \bG$ with $\phi_r(x) = 0$ for all $r \in \RR$, which implies $rx \ne a$ for all $r$. Thus $|\bG| = |\hat{\bG}|$ divides $|\RR|$.
\end{proof}

\begin{cor} Let $\RR$ be a finite ring, and let $\bM$ be a subdirectly irreducible module over $\RR$. Then $|\bM|$ is a prime power dividing $|\RR|$.
\end{cor}
\begin{proof} If $\bM$ is subdirectly irreducible, then it has a least nontrivial submodule $\bN$, which is generated by some nonzero element $a \in \bN$. Then for each nonzero $x \in \bM$ we have $\bN \le \RR x$, so there is some $r \in \RR$ with $rx = a$. Thus we can apply the previous proposition with $\bG = \bM$.
\end{proof}

Now we can use this result to bound the sizes of subdirectly irreducible algebras in congruence modular varieties in the special case where the centralizer of the monolith is abelian.

\begin{thm}\label{subdirect-ab-prime-power} Suppose that $\bB \in \cV$ is subdirectly irreducible, and $\cV$ is locally finite and congruence modular. If $\alpha \in \Con(\bB)$ is abelian and $|\bB/\alpha| = k$, then every congruence class of $\alpha$ has size a prime power bounded by $|\cF_\cV(k+1)|$.
\end{thm}
\begin{proof} (Adapted from \cite{commutator-uses}.) Assume $\alpha$ is nontrivial, so $0_\bB^* \le \alpha$. Let $p$ be a Gumm difference term. By Corollary \ref{difference-graph}, the restriction of the graph of $p$ to the blocks of $\alpha$ is preserved by every operation of $\bB$. Choose elements $0 \ne a$ with $(0,a) \in 0_\bB^*$, and note that $0,a$ are in the same congruence block of $\alpha$.

Pick constants $c_0, ..., c_{k-1}$ with $c_0 = 0$ such that each congruence class of $\alpha$ contains some $c_i$. We will treat each congruence class $c_i/\alpha$ of $\alpha$ as an abelian group with zero element $c_i$, addition given by $x +_i y = p(x,c_i,y)$, and subtraction given by $x -_i y = p(x,y,c_i)$.

Suppose that $x \ne y$ with $(x,y) \in \alpha$. Then since $0_\bB^*$ is the least nontrivial congruence, the pair $(0,a)$ must be in the congruence generated by $(x,y)$, so there must be a chain of unary polynomials $f_i$ such that $0 = f_0(x)$, $f_i(y) = f_{i+1}(x)$, and $f_m(y) = a$. Note that this implies that $f_i(x), f_i(y)$ are all in the congruence class $0/\alpha$. Thus, it makes sense to define a unary polynomial $f$ such that
\[
f(z) = f_0(z) +_0 f_1(z) -_0 f_1(x) +_0 \cdots +_0 f_m(z) -_0 f_m(x)
\]
for $z$ in the congruence class $x/\alpha$. One explicit way to construct such an $f$ is given by
\[
f(z) = p(p(\cdots p(p(f_0(z),f_1(x),f_1(z)), f_2(x), f_2(z)), \cdots), f_m(x), f_m(z)).
\]
It's easy to check that we have $f(x) = 0$ and $f(y) = a$. Since $f$ preserves the graph of $p$ restricted to congruence classes of $\alpha$, if $x,y \in c_i/\alpha$ then we have $f(x -_i y) -_0 f(c_i) = a$, and the unary polynomial $z \mapsto f(z) -_0 f(c_i)$ defines a linear map in $\Hom(c_i/\alpha, c_0/\alpha)$.

To finish, we just need to bound the size of the subgroup $\RR_{i,0}$ of linear maps in $\Hom(c_i/\alpha, c_0/\alpha)$ which can be defined by unary polynomials $f$. Suppose that $f(z) = t(z,b_1, ..., b_m)$ for some term $t$ and constants $b_1, ..., b_m \in \bB$, such that $f(c_i) = c_0$. For each $b_i$, we choose $j_i$ such that $b_i \in c_{j_i}/\alpha$. Define a unary polynomial $f'$ by
\[
f'(z) = t(z,c_{j_1}, ..., c_{j_m}) -_0 t(c_i, c_{j_1}, ..., c_{j_m}).
\]
Then for $z \in c_i/\alpha$, we have $f'(z) \in c_0/\alpha$, and since $t$ preserves the graph of $p$ restricted to congruence classes of $\alpha$, we have $f'(z) = f(z)$ for $z \in c_i/\alpha$ (alternatively, we could prove this by the term condition for $[\alpha,\alpha] = 0_\bB$). Thus every element of $\Hom(c_i/\alpha,c_0/\alpha)$ which can be defined by a unary polynomial can also be defined by a polynomial $f'$ which has the form $f'(z) = t'(z,c_0,...,c_{k-1})$ for some $k+1$-ary term $t'$, so
\[
|\RR_{i,0}| \le |\cF_\cV(k+1)|.
\]
Applying the previous proposition, we see that $|c_i/\alpha|$ is a prime power dividing $|\RR_{i,0}|$.
\end{proof}

\begin{cor}\label{subdirect-bound-ab} If $|\bA| = m$ is finite and $\cV(\bA)$ is congruence modular, and if $\bB \in \cV(\bA)$ is subdirectly irreducible with $(0_\bB:0_\bB^*)$ abelian, then $|\bB| \le m\cdot m^{m^{m+1}}$.
\end{cor}

\begin{defn} A variety $\cV$ is called \emph{residually small} if there is a cardinal $\kappa$ such that every subdirectly irreducible algebra $\bB \in \cV$ has $|\bB| < \kappa$, and \emph{residually finite} if every subdirectly irreducible algebra in $\cV$ is finite. An algebra $\bA$ is called residually small if the variety $\cV(\bA)$ generated by $\bA$ is residually small.
\end{defn}

First we show that if a locally finite variety contains an infinite subdirectly irreducible algebra, then it contains infinitely many distinct finite subdirectly irreducible algebras.

\begin{thm}\label{ultraproduct-fin-gen} If $\bB$ is subdirectly irreducible, then $\bB$ is a subalgebra of an ultraproduct of a family of finitely generated subdirectly irreducible algebras in $HS(\bB)$.
\end{thm}
\begin{proof} (From \cite{commutator-theory}.) Let the monolith $0_\bB^*$ of $\bB$ be generated (as a congruence) by the pair $(a,b)$. Let $I$ be the family of finitely generated subalgebras $\bS \le \bB$ with $a,b \in \bS$, and for each $\bS \in I$, pick a congruence $\alpha_\bS$ on $\bS$ which is maximal among all congruences which do not contain $(a,b)$. Then each $\bS/\alpha_\bS$ is subdirectly irreducible, since every congruence which properly contains $\alpha_\bS$ contains $(a,b)$.

Let $\cU$ be an ultrafilter on $I$ such that the set $U_S = \{\bS \mid S \subseteq \bS\}$ is in $\cU$ for every finite $S \subseteq \bB$. Such an ultrafilter exists since for any $S_1, S_2$ we have $U_{S_1} \cap U_{S_2} = U_{S_1\cup S_2}$, and for $S$ finite $U_S$ is nonempty since it contains $\Sg_\bB(S \cup \{a,b\})$.

Define a map $\varphi : \bB \rightarrow (\prod_{\bS \in I} \bS/\alpha_\bS)/\cU$ as the ultraproduct of the family of maps $\varphi_\bS$ given by $\varphi_\bS(x) = x/\alpha_\bS$ for $x \in \bS$ and $\varphi_\bS(x) = a/\alpha_\bS$ for $x \not\in \bS$. Then for $x_1, ..., x_k \in \bB$ and $t$ a $k$-ary term of $\bB$, we have
\[
\{\bS \mid \varphi_\bS(t(x_1, ..., x_k)) = t(\varphi_\bS(x_1), ..., \varphi_\bS(x_k))\} \in \cU,
\]
since it contains $U_{\{x_1, ..., x_k\}}$. Thus $\varphi$ is a homomorphism. To see that it is injective, just note that $\varphi_\bS(a) \ne \varphi_\bS(b)$ for all $\bS \in I$.
\end{proof}

\begin{cor}\label{residually-infinite} If a locally finite variety contains an infinite subdirectly irreducible algebra, then it contains arbitrarily large finite subdirectly irreducible algebras.
\end{cor}

It turns out that finite residually small algebras can be understood in terms of a commutator condition. We say that an algebra $\bA$ satisfies a commutator identity \emph{hereditarily} if every congruence lattice of every subalgebra of $\bA$ satisfies the identity.

\begin{prop} The commutator identity $[\alpha \wedge \beta, \beta] = \alpha \wedge [\beta,\beta]$ is equivalent to the implication $\alpha \le [\beta,\beta] \implies [\alpha,\beta] = \alpha$.
\end{prop}
\begin{proof} The implication clearly follows from the identity. For the other direction, we apply the implication to $\alpha \wedge [\beta,\beta] \le [\beta,\beta]$ to see that
\[
\alpha \wedge [\beta,\beta] = [\alpha \wedge [\beta,\beta], \beta] \le [\alpha\wedge \beta,\beta] \le \alpha \wedge [\beta,\beta].\qedhere
\]
\end{proof}

\begin{prop} If $\bA$ is in a congruence modular variety and satisfies the commutator identity $[\alpha \wedge \beta,\beta] = \alpha \wedge [\beta,\beta]$ hereditarily, then so does every quotient $\bB$ of $\bA$.
\end{prop}
\begin{proof} Suppose $\bB = \bA/\gamma$ and $\alpha,\beta \in \Con(\bA)$ with $\alpha,\beta \ge \gamma$. We need to check that if $\alpha \le [\beta,\beta]_\gamma$, then $\alpha = [\alpha,\beta]_\gamma$. By the modular law, if $\alpha \le [\beta,\beta]\vee \gamma$ then
\[
\alpha = \alpha \wedge ([\beta,\beta]\vee \gamma) = (\alpha \wedge [\beta,\beta]) \vee \gamma = [\alpha,\beta] \vee \gamma = [\alpha,\beta]_\gamma.\qedhere
\]
\end{proof}

\begin{prop} If $\bA_1,\bA_2$ are in a congruence modular variety and satisfy the commutator identity $[\alpha\wedge \beta, \beta] = \alpha\wedge [\beta,\beta]$ hereditarily, then so does their product $\bA_1 \times \bA_2$.
\end{prop}
\begin{proof} Let $\bB \le \bA_1\times \bA_2$, we can assume without loss of generality that this inclusion is subdirect by replacing the $\bA_i$ with $\pi_i(\bB)$. Suppose $\alpha,\beta \in \Con(\bB)$ with $\alpha \le [\beta,\beta]$, we will show that $[\alpha,\beta] = \alpha$. We have
\[
\alpha \vee \ker\pi_1 \le [\beta\vee\ker\pi_1, \beta\vee\ker\pi_1]_{\ker \pi_1},
\]
so from the assumption on $\bA_1$ we get
\[
\alpha \vee \ker \pi_1 = [\alpha \vee \ker\pi_1, \beta \vee \ker \pi_1]_{\ker \pi_1} = [\alpha,\beta]\vee \ker \pi_1.
\]
Thus by the modular law and $[\alpha,\beta] \le \alpha$, we have
\[
\alpha = \alpha\wedge (\ker \pi_1 \vee [\alpha,\beta]) = (\alpha \wedge \ker \pi_1) \vee [\alpha,\beta].
\]
Similarly, we have $\alpha = (\alpha \wedge \ker \pi_2) \vee [\alpha,\beta]$. Since $\alpha \wedge \ker \pi_2 \le \alpha \le [\beta,\beta]$, we may apply the same reasoning to $\alpha \wedge \ker \pi_2$ to see that
\[
\alpha \wedge \ker \pi_2 = (\alpha \wedge \ker \pi_2 \wedge \ker \pi_1) \vee [\alpha \wedge \ker \pi_2, \beta],
\]
so $\alpha \wedge \ker \pi_2 \le [\alpha,\beta]$, so
\[
\alpha = (\alpha \wedge \ker \pi_2) \vee [\alpha,\beta] = [\alpha,\beta].\qedhere
\]
\end{proof}

\begin{thm} If $|\bA| = m$ is finite and $\cV(\bA)$ is congruence modular, and if $\bA$ satisfies the commutator identity $[\alpha \wedge \beta,\beta] = \alpha \wedge [\beta,\beta]$ hereditarily, then every subdirectly irreducible algebra $\bB \in \cV(\bA)$ has $|\bB| \le m\cdot m^{m^{m+1}}$.
\end{thm}
\begin{proof} By Corollary \ref{residually-infinite}, we just need to check the bound in the case where $\bB$ is finite. In this case, we have $\bB \in HSP_{fin}(\bA)$, so $\bB$ satisfies the commutator identity $[\alpha \wedge \beta,\beta] = \alpha \wedge [\beta,\beta]$ by the previous propositions. Let $0_\bB^*$ be the monolith of $\bB$, and let $\alpha = (0_\bB:0_\bB^*)$ be its centralizer.

We claim that $\alpha$ is abelian. To see this, note that from $[\alpha,0_\bB^*] = 0_\bB$ we have
\[
0_\bB = [0_\bB^*\wedge \alpha,\alpha] = 0_\bB^*\wedge [\alpha,\alpha],
\]
so $[\alpha,\alpha] = 0_\bB$. Now we can apply Corollary \ref{subdirect-bound-ab} to see that $|\bB| \le m\cdot m^{m^{m+1}}$.
\end{proof}

\begin{ex} The symmetric group $S_3$ on three letters is residually small, since it satisfies the commutator identity $[\alpha \wedge \beta, \beta] = \alpha \wedge [\beta,\beta]$ hereditarily: the only interesting case to check is that $[A_3,S_3] = A_3$, where $A_3$ is the alternating group on three letters. We have $HS(S_3) = \{1, \ZZ/2, \ZZ/3, S_3\}$, and all three nontrivial elements are subdirectly irreducible.

The general theory shows that every subdirectly irreducible $\bG \in \cV(S_3)$ has an abelian normal subgroup $\bN$ with $\bG/\bN \in HS(S_3)$, with $|\bN|$ a prime power bounded by $|\cF_{\cV(S_3)}(|\bG/\bN|+1)| \le 6^{6^7}$. Since $\bN \in \cV(S_3)$ and every element of $S_3$ has order dividing $6$, $\bN$ has exponent $2$ or $3$. From here it is not too hard to check that the only nontrivial subdirectly irreducible algebras in $\cV(S_3)$ are $\ZZ/2,\ZZ/3, S_3$, and all three of these are subgroups of $S_3$. Thus every group in $\cV(S_3)$ is a subgroup of a power of $S_3$.
\end{ex}

\begin{prop} If $\bA$ is contained in a congruence modular variety but does not satisfy the commutator identity $[\alpha \wedge \beta,\beta] = \alpha \wedge [\beta,\beta]$ hereditarily, then there is some subdirectly irreducible $\bB \in HS(\bA)$ such that the centralizer of the monolith of $\bB$ is not abelian.
\end{prop}
\begin{proof} Suppose that $\bA$ fails to satisfy the commutator identity. In this case there must be $\alpha, \beta \in \Con(\bA)$ with $\alpha \le [\beta,\beta]$ and $[\alpha,\beta] < \alpha$. Let $\theta$ be a meet-irreducible congruence such that $\theta \ge [\alpha,\beta]$ but $\theta \not\ge \alpha$, and let $\theta^*$ be its cover. Then
\[
\theta^* \le \alpha\vee\theta \le [\beta,\beta]\vee \theta \le [\beta\vee\theta,\beta\vee\theta]_\theta
\]
and
\[
[\theta^*,\beta\vee\theta]_\theta \le [\alpha\vee\theta,\beta\vee\theta]_\theta = [\alpha,\beta]\vee\theta = \theta,
\]
so if we take $\bB$ to be $\bA/\theta$, then the monolith of $\bB$ is $\theta^*/\theta$, and $\beta\vee\theta/\theta$ is contained in the centralizer of the monolith of $\bB$ but is not abelian.
\end{proof}

\begin{thm} If $\bA$ is contained in a congruence modular variety but does not satisfy the commutator identity $[\alpha \wedge \beta,\beta] = \alpha \wedge [\beta,\beta]$, then $\cV(\bA)$ is not residually small. In fact, for every cardinal $\kappa$, $\cV(\bA)$ contains a subdirectly irreducible algebra whose congruence lattice has size at least $\kappa$.
\end{thm}
\begin{proof} (From \cite{commutator-theory}.) By the proposition, we can reduce to the case where $\bA$ is subdirectly irreducible and the centralizer $\beta$ of the monolith $0^*$ is not abelian.

Consider $\beta$ as a subalgebra of $\bA^2$, and $\Delta_\beta^{0^*}$ as a congruence on $\beta$. From $[\beta,0_\bA^*] = 0_\bA$ we have $\Delta_\beta^{0^*} \wedge \ker \pi_1 = \Delta_\beta^{0^*} \wedge \ker \pi_2 = 0_\beta$, and from the definition of $\Delta_\beta^{0^*}$ we have $\Delta_\beta^{0^*} \vee \ker \pi_i = \pi_i^{-1}(0^*)$. Set $0^*_i = \pi_i^{-1}(0^*)$ and $\theta = 0^*_1 \wedge 0^*_2, \theta_i = 0^*_i \wedge \ker \pi_{\{1,2\}\setminus \{i\}}$, then (after several applications of the modular law - don't worry about the details just yet) we have the following sublattice in $\Con(\beta)$.
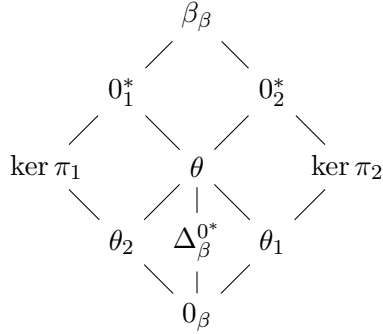
\begin{center}
\begin{tikzpicture}[scale=1.0]
  \node (1) at (0,2) {$\beta_{\beta}$};
  \node (01) at (-1,1) {$0^*_1$};
  \node (02) at (1,1) {$0^*_2$};
  \node (p1) at (-2,0) {$\ker \pi_1$};
  \node (c) at (0,0) {$\theta$};
  \node (p2) at (2,0) {$\ker \pi_2$};
  \node (t1) at (-1,-1) {$\theta_2$};
  \node (t2) at (1,-1) {$\theta_1$};
  \node (d) at (0,-1) {$\Delta_\beta^{0^*}$};
  \node (0) at (0,-2) {$0_\beta$};
  \draw (1) -- (01) -- (p1) -- (t1) -- (0) -- (t2) -- (p2) -- (02) -- (1);
  \draw (01) -- (c) -- (02);
  \draw (t1) -- (c) -- (t2);
  \draw (c) -- (d) -- (0);
\end{tikzpicture}
\end{center}
In the picture, we see that $\Delta_\beta^{0^*}$ appears to be meet-irreducible in $\Con(\beta)$, and the interval $\llbracket \Delta_\beta^{0^*}, \beta_\beta\rrbracket$ contains the incomparable elements $0^*_1, 0^*_2$. If $\Delta_\beta^{0^*}$ isn't meet-irreducible, we can still try to find a meet-irreducible congruence $\lambda$ on $\beta$ which is above $\Delta_\beta^{0^*}$ but not above $\theta$, and then $\beta/\lambda$ should give us a subdirectly irreducible algebra whose congruence lattice contains two distinct elements coming from $\ker \pi_1 \vee \lambda$ and $\ker \pi_2 \vee \lambda$ (that neither of these is equal to $\beta_\beta$ will come from the assumption that $\beta$ is not abelian). This is the basic idea behind the general construction, but we will need to scale up by considering higher dimensional analogues of $\beta \le \bA^2$.

Let $\kappa$ be any cardinal, considered as the set of all ordinals below $\kappa$. Define $\bB \le \bA^\kappa$ by
\[
\bB = \{a \in \bA^\kappa \mid a_i \equiv_\beta a_j\ \forall i,j \in \kappa\}.
\]
Then $\bB$ has a natural map to $\bA/\beta$, and we call the kernel of this map $\beta_\bB$. Inside $\Con(\bB)$, we have $\ker \pi_i \vee \ker \pi_j = \beta_\bB$ for all $i \ne j \in \kappa$. The strategy is to construct a congruence $\lambda$ on $\bB$ such that $\bB/\lambda$ is subdirectly irreducible and $\ker \pi_i \vee \lambda \not\ge \beta_\bB$ for all $i$, which will guarantee that the congruences $\ker \pi_i \vee \lambda/\lambda \in \Con(\bB/\lambda)$ are pairwise distinct. The congruence $\lambda$ will be constructed by first constructing congruences $\Delta, \theta$ with $\Delta < \theta$ and $\theta \vee \ker \pi_i \not\ge \beta_\bB$.

We need a congruence on $\bB$ generalizing $\Delta_\beta^{0^*}$ on $\beta$. We define $\Delta_i$ by
\[
(a,b) \in \Delta_i \iff \begin{bmatrix} a_0 & b_0\\ a_i & b_i\end{bmatrix} \in \Delta_\beta^{0^*} \wedge (a_j = b_j\ \forall j \ne 0,i),
\]
and define $\Delta$ by
\[
\Delta = \bigvee_{0 < i < \kappa} \Delta_i.
\]
We also define congruences $\theta_i$ by
\[
(a,b) \in \theta_i \iff (a_i,b_i) \in 0_\bA^* \wedge (a_j = b_j\ \forall j \ne i),
\]
and define $\theta$ by
\[
\theta = \bigvee_{i < \kappa} \theta_i.
\]

We need to check some basic properties of these congruences, to see that they behave as in the picture of $\Con(\beta)$. First, we check that $\theta_0 \le \theta_i \vee \Delta_i$ for all $i$. Letting $\pi_{i'}$ be the projection onto all coordinates other than $i$, then it's easy to check that $\theta_0 \le \ker \pi_{i'} \vee \Delta$ by reasoning about just the two coordinates $0,i$ and keeping all other coordinates fixed:
\[
\begin{bmatrix} a_0\\ a_i\end{bmatrix}\ \ker \pi_{i'}\ \begin{bmatrix} a_0\\ a_0\end{bmatrix}\ \Delta_i\ \begin{bmatrix} b_0\\ b_0\end{bmatrix}\ \ker \pi_{i'}\ \begin{bmatrix} b_0\\ b_i\end{bmatrix}.
\]
Then by the modular law, if we let $0^*_i = \pi_i^{-1}(0^*_\bA)$ and note that $\Delta_i \le 0^*_i$, we get
\[
\theta_0 = \theta_0 \wedge 0^*_i \le (\Delta_i \vee \ker \pi_{i'}) \wedge 0^*_i = \Delta_i \vee (\ker \pi_{i'} \wedge 0^*_i) = \Delta_i \vee \theta_i.
\]
Similarly, we get $\theta_i \le \theta_0 \vee \Delta_i$ for all $i$.

Next, for each $i$ we have $\Delta \vee \theta_i = \theta$: for each $j \in \kappa$, we have
\[
\Delta \vee \theta_i = \Delta \vee \Delta_i \vee \theta_i \ge \Delta \vee \theta_0 \ge \Delta_j \vee \theta_0 \ge \theta_j,
\]
so $\Delta \vee \theta_i \ge \bigvee_{j \in \kappa} \theta_j = \theta$, while the other containment follows from $\Delta_i \le \theta_0 \vee \theta_i$ for all $i$.

We now check that $\Delta \ne \theta$. It's enough to check that $\theta_0 \not\le \Delta$, since $\Delta \vee \theta_0 = \theta$. Note first that $\theta_0$ is compact, since $0^*_\bA$ is compact. Thus we just need to check that $\theta_0 \not\le \bigvee_{j \le n} \Delta_{i_j}$ for all $i_1, ..., i_n$. In fact, we can assume that $i_1, ..., i_n$ are $1, ..., n$ by a symmetry argument.

We will show by induction on $n$ that $\theta_0\wedge (\Delta_1 \vee \cdots \vee \Delta_n) = 0_\bB$ for all $n$. The base case follows from the fact that $[\beta,0^*_\bA] = 0_\bA \implies \ker \pi_2 \wedge \Delta_\beta^{0^*} = 0_\beta$ in $\Con(\beta)$, which in turn implies $\theta_0 \wedge \Delta_1 = 0_\bB$. For the inductive step, we argue as follows:
\begin{align*}
\theta_0\wedge (\Delta_1 \vee \cdots \vee \Delta_n) &= \theta_0\wedge (\theta_0\vee \theta_n)\wedge (\Delta_1 \vee \cdots \vee \Delta_n)\\
&= \theta_0 \wedge (((\theta_0\vee \theta_n) \wedge (\Delta_1 \vee \cdots \vee \Delta_{n-1})) \vee \Delta_n)\\
&= \theta_0 \wedge ((\theta_0 \wedge (\Delta_1 \vee \cdots \vee \Delta_{n-1})) \vee \Delta_n)\\
&= \theta_0 \wedge \Delta_n = 0_\bB,
\end{align*}
where the second equality used the modular law and the fact that $\Delta_n \le \theta_0\vee \theta_n$, the third equality used the fact that $\theta_n$ is independent of everything that happens on the coordinates $0, ..., n-1$, and the last two equalities used the inductive hypothesis.

We have shown that $\Delta < \theta$. We can now apply Corollary \ref{meet-irreducible} to see that there is some meet-irreducible congruence $\lambda$ with $\lambda \ge \Delta$ but $\lambda \not\ge \theta$. To finish, we just need to check that $\lambda \vee \ker \pi_i \not\ge \beta_\bB$. To see this, note that $\lambda \not\ge \theta_i$, since otherwise we would have $\lambda \ge \Delta \vee \theta_i = \theta$, a contradiction. Since $\theta_i$ is the minimal nonzero element of the interval $\llbracket 0_\bB, \ker \pi_{i'}\rrbracket$, this means that $\lambda \wedge \ker \pi_{i'} = 0_\bB$. Thus if (for contradiction) $\lambda \vee \ker \pi_i \ge \beta_\bB$, then we would have
\[
[\beta_\bB, \beta_\bB] \le [\ker \pi_{i'} \vee \ker \pi_i, \lambda \vee \ker \pi_i] \le (\lambda \wedge \ker \pi_{i'}) \vee \ker \pi_i = \ker \pi_i,
\]
and applying $\pi_i$ we would get $[\beta,\beta] = 0_\bA$, a contradiction to the assumption that $\beta$ was not abelian.

Putting it all together, we have a meet-irreducible congruence $\lambda$ such that $\lambda \vee \ker \pi_i \not\ge \beta_\bB$ for each $i$, but $\ker \pi_i \vee \ker \pi_j \ge \beta_\bB$ for all $i \ne j$. Thus $\bB/\lambda$ is subdirectly irreducible, and the congruences $\ker \pi_i \vee \lambda/\lambda$ are mutually distinct elements of $\Con(\bB/\lambda)$.
\end{proof}

\begin{cor}\label{residual-crit} Let $\cV$ be a finitely generated congruence modular variety. Then the following are equivalent:
\begin{itemize}
\item $\cV$ is residually small,

\item every algebra in $\cV$ satisfies the commutator identity $[\alpha \wedge \beta, \beta] = \alpha \wedge [\beta,\beta]$,

\item $\cV$ is generated by a finite algebra $\bA$ such that for every subdirectly irreducible $\bB \in HS(\bA)$, the centralizer of the monolith of $\bB$ is abelian,

\item $\cV$ is generated by a finite algebra which satisfies the commutator identity $[\alpha \wedge \beta, \beta] = \alpha \wedge [\beta,\beta]$ hereditarily,

\item $\cV$ has a finite bound on the size of its subdirectly irreducible elements.
\end{itemize}
\end{cor}

\begin{cor} If $\bA$ is in a congruence modular variety and has size $|\bA| \le 3$, then $\cV(\bA)$ is residually small.
\end{cor}
\begin{proof} Suppose for contradiction that $\bA$ is subdirectly irreducible with a monolith $0_\bA^*$ whose centralizer $(0_\bA : 0_\bA^*)$ is not abelian. Then since $\Con(\bA)$ has height at most $2$, we necessarily have
\[
0_\bA < 0_\bA^* < (0_\bA : 0_\bA^*) = 1_\bA.
\]
Thus $|\bA| = 3$, and we may name the elements of $\bA$ as $a,b,c$, such that $0_\bA^*$ corresponds to the partition $\{a,b\},\{c\}$ of $\bA$. Letting $p(x,y,z)$ be a Gumm difference term, we see from Theorem \ref{difference-commutator} that
\[
\begin{bmatrix} a & p(a,b,c)\\ b & c\end{bmatrix} \in \Delta_{0_\bA^*}^{1_\bA}.
\]
Modulo $0_\bA^*$, we have $p(a,b,c) \equiv_{0_\bA^*} p(a,a,c) = c$, so we must have $p(a,b,c) = c$. Then by Theorem \ref{shifting-commutator} we have $(a,b) \in [1_\bA, 0_\bA^*]$, which contradicts $(0_\bA : 0_\bA^*) = 1_\bA$.
\end{proof}

\begin{prop}\label{residual-nilpotent} If $\bA$ satisfies the commutator identity $[\alpha\wedge \beta,\beta] = \alpha \wedge [\beta,\beta]$, then every nilpotent congruence on $\bA$ is abelian.
\end{prop}
\begin{proof} The commutator identity implies that
\[
[[\alpha,\alpha],\alpha] = [[\alpha,\alpha] \wedge \alpha,\alpha] = [\alpha,\alpha]\wedge[\alpha,\alpha] = [\alpha,\alpha].\qedhere
\]
\end{proof}

\begin{prop}[Ol'{\v{s}}anski{\u\i} \cite{residually-small-groups}]\label{sylow-abelian} If all the Sylow subgroups of a finite group $\bG$ are abelian, then the center $Z(\bG)$ and the commutator subgroup $[\bG,\bG]$ intersect trivially, that is, $Z(\bG) \wedge [\bG,\bG] = 0_\bG$.
\end{prop}
\begin{proof} Fix a Sylow subgroup $\bS$ of $\bG$, and consider the transfer map $\bG \rightarrow \bS/[\bS,\bS]$. Recall that the transfer homomorphism from a finite group to the abelianization of a subgroup is defined by making a choice of coset representatives $x_i$ with $\bG = \bigcup_i x_i\bS$, and sending $g \in \bG$ to $\prod_i s_i/[\bS,\bS]$, where for each $i$, $s_i \in \bS$ is given by $gx_i = x_js_i$ for some $j$. Since $\bS$ is assumed to be abelian, this gives us a homomorphism from $\bG$ to $\bS$.

Now consider any $g \in Z(\bG) \cap \bS$. The transfer homomorphism sends $g$ to $\prod_i g = g^{[\bG:\bS]}$ since $gx_i = x_ig$ for each $i$, and if $g \ne 1$ then $g^{[\bG:\bS]} \ne 1$ as well since $[\bG:\bS]$ is relatively prime to the order of $g$. Thus there is a map from $\bG$ to an abelian group such that $g$ is not in the kernel, so $g \not\in [\bG,\bG]$. Since every nontrivial element of $Z(\bG) \wedge [\bG,\bG]$ has a power which has prime order and is therefore contained in a Sylow subgroup of $\bG$, we must have $Z(\bG) \wedge [\bG,\bG] = 0_\bG$ to avoid a contradiction.
\end{proof}

\begin{cor}[Ol'{\v{s}}anski{\u\i} \cite{residually-small-groups}] A finite group is residually small iff all of its Sylow subgroups are abelian.
\end{cor}
\begin{proof} By Proposition \ref{residual-nilpotent}, all nilpotent subgroups of a finite residually small group must be abelian, so in particular the Sylow subgroups must be abelian since all $p$-groups are nilpotent.

For the other direction, note that for any $\bB \in HS(\bA)$, the Sylow subgroups of $\bB$ are quotients of subgroups of the Sylow subgroups of $\bA$ by the Sylow theorems. Thus we just have to check that if the Sylow subgroups of a subdirectly irreducible group are abelian, then the centralizer $\bC$ of its  monolith $0^*$ is abelian.

Note that if $\bC$ centralizes $0^*$, then $0^* \le Z(\bC)$. By Proposition \ref{sylow-abelian}, we have $Z(\bC) \wedge [\bC,\bC] = 0$, so $0^* \wedge [\bC,\bC] = 0$, which implies that $[\bC,\bC] = 0$.
\end{proof}

\subsection{Similarity}\label{ss-similarity}

Even if a finitely generated congruence modular variety is not residually small, we can still classify its subdirectly irreducible algebras by using the concept of \emph{similarity} from Freese and McKenzie \cite{commutator-theory}. We will use a different definition of similarity than their definition, but which they prove to be equivalent.

\begin{defn} We say that subdirectly irreducible algebras $\bA, \bB$ in a congruence modular variety $\cV$ are \emph{similar} if there exists an algebra $\bC \in \cV$ with congruences $\alpha, \beta, \gamma, \delta \in \bC$ such that $\bC/\alpha \cong \bA$, $\bC/\beta \cong \bB$, and
\[
\llbracket \alpha, \alpha^* \rrbracket \searrow \llbracket \gamma, \delta \rrbracket \nearrow \llbracket \beta, \beta^* \rrbracket.
\]
If furthermore $\bC \le_{sd} \bA\times \bB$ and $\alpha, \beta$ are the kernels of the projections to $\bA,\bB$, then we say that $\bC$ is the \emph{graph of a similarity} from $\bA$ to $\bB$.
\end{defn}

\begin{prop}\label{similarity-graph} If $\bA, \bB$ are similar, then there is a witnessing algebra $\bC \le_{sd} \bA\times\bB$ which is the graph of a similarity from $\bA$ to $\bB$. If $\alpha,\beta$ are the kernels of the projections to $\bA,\bB$, then $(\alpha:\alpha^*) = (\beta:\beta^*)$ and $\bC/(\alpha:\alpha^*)$ is the graph of an isomorphism
\[
\bA/(0_\bA:0_\bA^*) \xrightarrow{\sim} \bB/(0_\bB:0_\bB^*).
\]
If $\bA, \bB$ are similar but not isomorphic, then they must both have abelian monoliths.
\end{prop}
\begin{proof} For the first statement, let $\bC \in \cV$ and $\alpha,\beta,\gamma,\delta \in \Con(\bC)$ be as in the definition of similarity. It's enough to show that we have
\[
\llbracket \alpha, \alpha^* \rrbracket \searrow \llbracket \alpha\wedge\beta, (\alpha\wedge\beta)\vee \delta \rrbracket \nearrow \llbracket \beta, \beta^* \rrbracket,
\]
since then we can replace $\bC$ by $\bC/(\alpha\wedge\beta)$, which is a subdirect product of $\bC/\alpha \cong \bA$ and $\bC/\beta \cong \bB$. We have
\[
\alpha \vee ((\alpha\wedge\beta)\vee \delta) = \alpha\vee\delta = \alpha^*,
\]
and by the modular law and the fact that $\gamma \le \alpha \wedge \beta$, we have
\[
\alpha \wedge ((\alpha\wedge\beta)\vee \delta) = (\alpha \wedge \delta) \vee (\alpha \wedge \beta) = \gamma \vee (\alpha \wedge \beta) = (\alpha \wedge \beta),
\]
so $\llbracket \alpha, \alpha^* \rrbracket \searrow \llbracket \alpha\wedge\beta, (\alpha\wedge\beta)\vee \delta \rrbracket$, and the other perspectivity follows by a symmetric argument.

The remaining statements follow from the Diamond Isomorphism Theorem \ref{diamond-isom}: if $\llbracket \alpha, \alpha^* \rrbracket \searrow \llbracket \gamma, \delta \rrbracket \nearrow \llbracket \beta, \beta^* \rrbracket$, then $(\alpha:\alpha^*) = (\gamma:\delta) = (\beta:\beta^*)$, so
\[
\bA/(0_\bA:0_\bA^*) \cong \bC/(\alpha:\alpha^*) = \bC/(\beta:\beta^*) \cong \bB/(0_\bB:0_\bB^*),
\]
and
\[
[\alpha^*,\alpha^*]_\alpha = \alpha \iff [\delta,\delta]_\gamma = \gamma \iff [\beta^*,\beta^*]_\beta = \beta,
\]
so $0_\bA^*$ is abelian iff $0_\bB^*$ is abelian, and if neither is abelian then $\alpha = (\alpha:\alpha^*) = (\beta:\beta^*) = \beta$ and $\bA \cong \bC/\alpha = \bC/\beta \cong \bB$.
\end{proof}

\begin{prop} If $\bA, \bB$ are similar such that $\sigma$ is the corresponding isomorphism
\[
\sigma: \bA/(0_\bA:0_\bA^*) \xrightarrow{\sim} \bB/(0_\bB:0_\bB^*),
\]
then they are similar via the algebra $\RR = \{(x,y) \in \bA \times \bB \mid \sigma(x/(0_\bA:0_\bA^*)) = y/(0_\bB:0_\bB^*)\}$.
\end{prop}
\begin{proof} Suppose $\bC \le \RR$ is the graph of a similarity from $\bA$ to $\bB$, with
\[
\llbracket \ker \pi_1, (\ker \pi_1)^* \rrbracket \searrow \llbracket 0_\bC, \delta \rrbracket \nearrow \llbracket \ker \pi_2, (\ker \pi_2)^* \rrbracket
\]
in $\Con(\bC)$. We may assume that $\bA, \bB$ have abelian monoliths, so $[\delta,\delta] = 0_\bC$ by the Diamond Isomorphism Theorem \ref{diamond-isom}. Then by Theorem \ref{commutator-permute}, $\delta$ permutes with all congruences in $\Con(\bC)$, so in particular $(\ker \pi_1)^* = \delta \circ \ker \pi_1$. In other words, for any $(a,b) \in \bC$ and any $a' \in a/0_\bA^*$, there exists a $b'$ such that
\[
\begin{bmatrix} a\\ b\end{bmatrix}\ \delta\ \begin{bmatrix} a'\\ b'\end{bmatrix}.
\]
In fact, this $b'$ is uniquely determined by $a,b,a'$, since $\delta \wedge \ker \pi_1 = 0_\bC$. Additionally, we must have $b' \in b/0_\bB^*$, since $\delta \le (\ker \pi_2)^*$.

Now we can extend $\delta$ to a congruence $\delta_\RR \in \Con(\RR)$ as follows. For $(a,b), (a',b') \in \RR$ with $a\ 0_\bA^*\ a'$ and $b\ 0_\bB^*\ b'$, we pick any $(u,v) \in \bC$ with $u\ (0_\bA:0_\bA^*)\ a$ and write
\[
\begin{bmatrix} a & a'\\ b & b'\end{bmatrix} \in \delta_\RR\ \iff\ \begin{bmatrix} p(a,a',u) & u\\ p(b,b',v) & v\end{bmatrix} \in \delta,
\]
where $p$ is a Gumm difference term. Note that by Corollary \ref{difference-graph}, this choice of $\delta_\RR$ is preserved by the operations of $\bA$ so long as it is well-defined. To check that this is in independent of the choice of $(u,v) \in \bC$, suppose $(u',v') \in \bC$ with $u'\ (0_\bA:0_\bA^*)\ a$, and apply Corollary \ref{difference-graph} again to see that
\[
p\left(\begin{bmatrix} p(a,a',u) & u\\ p(b,b',v) & v\end{bmatrix}, \begin{bmatrix} p(a,a,u) & u\\ p(b,b,v) & v\end{bmatrix}, \begin{bmatrix} p(a,a,u') & u'\\ p(b,b,v') & v'\end{bmatrix}\right) = \begin{bmatrix} p(a,a',u') & u'\\ p(b,b',v') & v'\end{bmatrix},
\]
where we have used $0_\bA^*, 0_\bB^*$ abelian to see that $p(a',a,a) = a'$ and $p(b',b,b) = b'$.

We need to check that $\delta_\RR$ is a congruence on $\RR$. It clearly contains the equality relation on $\RR$. For symmetry and transitivity, note that
\[
p\left(\begin{bmatrix} p(a,a',u)\\ p(b,b',v)\end{bmatrix}, \begin{bmatrix} p(a'',a',u)\\ p(b'',b',v)\end{bmatrix}, \begin{bmatrix} u\\ v\end{bmatrix}\right) = p\left(\begin{bmatrix} p(a,a',u)\\ p(b,b',v)\end{bmatrix}, \begin{bmatrix} p(a'',a',u)\\ p(b'',b',v)\end{bmatrix}, \begin{bmatrix} p(a'',a'',u)\\ p(b'',b'',v)\end{bmatrix}\right) = \begin{bmatrix} p(a,a'',u)\\ p(b,b'',v)\end{bmatrix}.
\]

Finally, we need to check that $\delta_\RR \wedge \ker \pi_1 = 0_\RR$ and $\delta_\RR \vee \ker \pi_1 = (\ker \pi_1)^*$. That $\delta_\RR \wedge \ker \pi_1 = 0_\RR$ follows from the fact that if we pick $u$ such that $(u,b') \in \bC$, then
\[
\begin{bmatrix} a & a\\ b & b'\end{bmatrix} \in \delta_\RR\ \iff\ \begin{bmatrix} p(a,a,u) & u\\ p(b,b',b') & b'\end{bmatrix} = \begin{bmatrix} u & u\\ b & b'\end{bmatrix} \in \delta,
\]
and so this can only occur when $b = b'$ since $\delta \wedge \ker \pi_1 = 0_\bC$ (by assumption). That $\delta_\RR \vee \ker \pi_1 = (\ker \pi_1)^*$ follows from $\delta \subseteq \delta_\RR \subseteq (\ker \pi_1)^*$ and $\delta \not\subseteq \ker \pi_1$.
\end{proof}

\begin{cor}\label{similar-detail} A similarity from $\bA$ to $\bB$ can be described by the following data: an isomorphism
\[
\sigma: \bA/(0_\bA:0_\bA^*) \xrightarrow{\sim} \bB/(0_\bB:0_\bB^*)
\]
together with a congruence $\delta \in \Con(\RR)$, where $\RR = \{(x,y) \in \bA \times \bB \mid \sigma(x/(0_\bA:0_\bA^*)) = y/(0_\bB:0_\bB^*)\}$, such that for every $(a,b) \in \RR$ and every $a' \in a/0_{\bA}^*$, there exists a unique $b' \in b/0_\bB^*$ such that
\[
\begin{bmatrix} a & a'\\ b & b'\end{bmatrix} \in \delta.
\]
In particular, if $\bA, \bB$ are idempotent, then for any $(a,b) \in \RR$ the congruence classes $a/0_\bA^*$ and $b/0_\bB^*$ are isomorphic to each other.
\end{cor}

\begin{cor} Similarity is an equivalence relation on subdirectly irreducible algebras.
\end{cor}
\begin{proof} Suppose we have similarities from $\bA$ to $\bB$ and from $\bB$ to $\bC$, described by isomorphisms
\[
\bA/(0_\bA:0_\bA^*) \xrightarrow{\sigma} \bB/(0_\bB:0_\bB^*) \xrightarrow{\sigma'} \bC/(0_\bC:0_\bC^*)
\]
and congruences $\delta, \delta'$. We define a congruence $\delta \circ \delta'$ by
\[
\begin{bmatrix} a & a'\\ c & c'\end{bmatrix} \in \delta \circ \delta'\ \iff\ \exists (b,b') \in 0_\bB^*\ \left(\begin{bmatrix} a & a'\\ b & b'\end{bmatrix} \in \delta\right) \wedge \left(\begin{bmatrix} b & b'\\ c & c'\end{bmatrix} \in \delta'\right).
\]
We need to check that for each $a,c,a'$ there exists a unique $c'$ satisfying the above. Existence is easy: for each $b$, we can fill in a unique $b'$ to satisfy $\delta$, and then there is a unique $c'$ which satisfies $\delta'$. We just need to show that the choice of $b$ doesn't affect the final $c'$ we get. Suppose that instead of $b$ we had picked $v$. Then the claim is that if we leave $a,a',c,c'$ unchanged and replace $b$ by $v$ and $b'$ by $p(b',b,v)$, we get another valid solution. For $\delta$, this follows from
\[
p\left(\begin{bmatrix} a & a'\\ b & b'\end{bmatrix}, \begin{bmatrix} a & a\\ b & b\end{bmatrix}, \begin{bmatrix} a & a\\ v & v\end{bmatrix}\right) = \begin{bmatrix} a & a'\\ v & p(b',b,v)\end{bmatrix},
\]
and it follows for $\delta'$ similarly.
\end{proof}

We will show that every subdirectly irreducible algebra $\bA$ with abelian monolith is similar to a subdirectly irreducible algebra $D(\bA)$ such that the monolith of $D(\bA)$ is equal to its own centralizer. The size of the algebra $D(\bA)$ can then be bounded using Theorem \ref{jonsson-modular} and the following proposition.

\begin{prop}\label{monolith-bound} If $\bB \in \cV(\bA)$ is subdirectly irreducible, $\bA$ is finite, and $\cV(\bA)$ is congruence modular, then every congruence class of $0_\bB^*$ has size at most $|\bA|$.
\end{prop}
\begin{proof} By Theorem \ref{ultraproduct-fin-gen} and Corollary \ref{ultraproduct-finite}, we may assume without loss of generality that $\bB$ is finite. By Theorem \ref{jonsson-modular}, we may also assume that $0_\bB^*$ is abelian. Take $m$ minimal such that there exists $\bC \le \bA^m$ and $\theta \in \Con(\bC)$ with $\bB \cong \bC/\theta$, so $[\theta^*,\theta^*] \le \theta$.

Let $\pi_{1'}$ be the projection onto all but the first coordinate, then by the minimality of $m$ we have $\ker \pi_{1'} \not\le \theta$. Thus we have
\[
\llbracket \theta, \theta^* \rrbracket \searrow \llbracket \theta \wedge \ker \pi_{1'}, \theta^*\wedge \ker\pi_{1'} \rrbracket.
\]
By Theorem \ref{commutator-permute}, the congruences $\theta$ and $\theta^*\wedge \ker \pi_{1'}$ permute. Thus for every congruence class $C^*$ of $\theta^*$ containing some $c \in \bC$, the size of $C^*/\theta$ is equal to the size of $C'/(\theta\wedge \ker\pi_{1'})$, where $C'$ is the congruence class of $\theta^*\wedge\ker\pi_{1'}$ containing $c$. But $|C'/(\theta\wedge \ker\pi_{1'})| \le |\bC/\ker \pi_{1'}| = |\bA|$, so every congruence class of $0_\bB^*$ has size bounded by $|\bA|$.
\end{proof}

\begin{defn} Suppose $\bA$ is a subdirectly irreducible algebra in a congruence modular variety. If $0_\bA^*$ is nonabelian, define $D(\bA)$ to be $\bA$. Otherwise, consider $0_\bA^*$ as a subalgebra of $\bA^2$ and $\Delta_{0_\bA^*}^{(0:0^*)}$ as a congruence on $0_\bA^*$, and define $D(\bA) = 0_\bA^*/\Delta_{0_\bA^*}^{(0:0^*)}$.
\end{defn}

Recall that by Theorem \ref{difference-commutator}, if $0_\bA^*$ is abelian and $p$ is a Gumm difference term, then $(0_\bA:0_\bA^*) \ge 0_\bA^*$ and $[(0_\bA:0_\bA^*),0_\bA^*] = 0_\bA$, so we have
\[
\begin{bmatrix} x & w\\ y & z\end{bmatrix} \in \Delta_{0_\bA^*}^{(0:0^*)} \iff (p(x,y,z) = w) \wedge (x\ \equiv_{0_\bA^*}\ y\ \equiv_{(0:0^*)}\ z).
\]
In this case, the subalgebra $\{(x,x)/\Delta_{0_\bA^*}^{(0:0^*)}\} \le D(\bA)$ meets every congruence class of $(0_\bA:0_\bA^*)_{D(\bA)}$ (that is, the congruence $(0_\bA:0_\bA^*)$ considered as a congruence on $D(\bA)$) exactly once, and is isomorphic to $\bA/(0_\bA:0_\bA^*)$.

\begin{prop} If $\bA$ is a subdirectly irreducible algebra in a congruence modular variety with an abelian monolith, then $D(\bA)$ is subdirectly irreducible with monolith $(0_\bA:0_\bA^*)_{D(\bA)}$, and $\bA, D(\bA)$ are similar via the algebra $0_\bA^*$ and the congruences $\ker \pi_1, \Delta_{0_\bA^*}^{(0:0^*)} \in \Con(0_\bA^*)$. Furthermore, the monolith $(0_\bA:0_\bA^*)_{D(\bA)}$ of $D(\bA)$ is its own centralizer.
\end{prop}
\begin{proof} Note that $\ker \pi_1$ is covered by $\ker \pi_1 \vee \ker \pi_2$, since $\pi_1(\ker \pi_1 \vee \ker \pi_2) = 0_\bA^*$. First we check that in $\Con(0_\bA^*)$ we have the perspectivities
\[
\llbracket \ker \pi_1, \ker \pi_1 \vee \ker \pi_2\rrbracket \searrow \llbracket 0_{0_\bA^*}, \ker \pi_2 \rrbracket \nearrow \llbracket \Delta_{0_\bA^*}^{(0:0^*)}, (0_\bA : 0_{\bA}^*)_{0_\bA^*}\rrbracket.
\]
The hardest step here is checking that $\ker \pi_2 \vee \Delta_{0_\bA^*}^{(0:0^*)} = (0_\bA : 0_{\bA}^*)_{0_\bA^*}$: if $(x,y), (w,z) \in 0_\bA^*$ with $(y,z) \in (0_\bA : 0_{\bA}^*)$, then we have
\[
\begin{bmatrix} x \\ y \end{bmatrix}\ \Delta_{0_\bA^*}^{(0:0^*)}\ \begin{bmatrix} p(x,y,z)\\ z\end{bmatrix}\ \ker \pi_2\ \begin{bmatrix} w\\ z\end{bmatrix}.
\]
To see that $\ker \pi_2 \wedge \Delta_{0_\bA^*}^{(0:0^*)} = 0_{0_\bA^*}$, note that by Theorem \ref{shifting-commutator} the inequality $\ker \pi_2 \wedge \Delta_{0_\bA^*}^{(0:0^*)} \le \ker \pi_1$ is equivalent to $[(0_\bA:0_\bA^*),0_\bA^*] = 0_\bA$.

Next we show that $(0_\bA : 0_{\bA}^*)_{0_\bA^*}$ is the unique cover of $\Delta_{0_\bA^*}^{(0:0^*)}$ in $\Con(0_\bA^*)$. Note first that $(0_\bA : 0_{\bA}^*)_{0_\bA^*}$ is \emph{a} cover of $\Delta_{0_\bA^*}^{(0:0^*)}$, since the interval $\llbracket \Delta_{0_\bA^*}^{(0:0^*)}, (0_\bA : 0_{\bA}^*)_{0_\bA^*}\rrbracket$ is isomorphic to $\llbracket \ker \pi_1, \ker \pi_1 \vee \ker \pi_2\rrbracket \cong \llbracket 0_\bA, 0_\bA^* \rrbracket$ by the Diamond Isomorphism Theorem \ref{diamond-isom}.

Suppose that $\psi$ is any congruence in $\Con(0_\bA^*)$ with $\psi > \Delta_{0_\bA^*}^{(0:0^*)}$. If $\psi \ge \ker \pi_2$, then $\psi \ge \Delta_{0_\bA^*}^{(0:0^*)} \vee \ker \pi_2 = (0_\bA : 0_{\bA}^*)_{0_\bA^*}$, and we are done. Otherwise, since $\ker \pi_2$ is a cover of $0_{0_\bA^*}$, we must have $\psi \wedge \ker \pi_2 = 0_{0_\bA^*}$. Then we have
\[
[\psi \vee \ker \pi_1, \ker \pi_2 \vee \ker \pi_1]_{\ker \pi_1} \le [\psi, \ker \pi_2] \vee \ker \pi_1 \le (\psi \wedge \ker \pi_2) \vee \ker \pi_1 = \ker \pi_1.
\]
Applying $\pi_1$ to both sides, we see that $\pi_1(\psi \vee \ker \pi_1) \le (0_\bA : 0_\bA^*)$, so $\psi \vee \ker \pi_1 \le (0_\bA : 0_{\bA}^*)_{0_\bA^*}$. Thus $\psi \in \llbracket \Delta_{0_\bA^*}^{(0:0^*)}, (0_\bA : 0_{\bA}^*)_{0_\bA^*}\rrbracket$, so again we must have $\psi = (0_\bA : 0_{\bA}^*)_{0_\bA^*}$. We have finished showing that $D(\bA)$ is subdirectly irreducible.

To see that the monolith $(0_\bA:0_\bA^*)_{D(\bA)}$ of $D(\bA)$ is its own centralizer, note that by the Diamond Isomorphism Theorem \ref{diamond-isom} we have
\[
(\Delta_{0_\bA^*}^{(0:0^*)} : (0_\bA : 0_{\bA}^*)_{0_\bA^*}) = (\ker \pi_1 : \ker \pi_1 \vee \ker \pi_2) = \pi_1^{-1}((0_\bA : 0_\bA^*)) = (0_\bA : 0_{\bA}^*)_{0_\bA^*}.\qedhere
\]
\end{proof}

\begin{prop} If $\bA, \bB$ are subdirectly irreducible algebras in a congruence modular variety, then $\bA$ is similar to $\bB$ iff $D(\bA) \cong D(\bB)$.
\end{prop}
\begin{proof} Since similarity is an equivalence relation, we may as well replace $\bA, \bB$ by $D(\bA), D(\bB)$. Thus we just need to prove that if $\bA, \bB$ have monoliths equal to their own centralizers, and have subalgebras $X_\bA, X_\bB$ which intersect their monoliths transversely, then they are similar iff they are isomorphic.

Let $\sigma : \bA/0_\bA^* \rightarrow \bB/0_\bB^*$ be the isomorphism and $\delta \in \Con(\RR)$, where $\RR = \{(x,y) \in \bA \times \bB \mid \sigma(x/(0_\bA:0_\bA^*)) = y/(0_\bB:0_\bB^*)\}$, be the data describing a similarity from $\bA$ to $\bB$. Then $\sigma$ induces an isomorphism $\sigma_X : X_\bA \rightarrow X_\bB$, and the graph of $\sigma_X$ is a subalgebra of $\RR$. Let $\bS$ be the subalgebra of $(a,b) \in \RR$ such that $(a,b)$ is congruent to some element of $\sigma_X$ modulo $\delta$. Then $\bS$ must be the graph of an isomorphism from $\bA$ to $\bB$.
\end{proof}

\begin{thm} If $\bB \in \cV(\bA)$ is subdirectly irreducible, $\bA$ is finite, and $\cV(\bA)$ is congruence modular, then $\bB$ is similar to a subdirectly irreducible algebra in $HS(\bA)$.
\end{thm}
\begin{proof} We may as well replace $\bB$ by $D(\bB)$, so assume without loss of generality that the monolith of $\bB$ is either nonabelian or equal to its own centralizer. If the monolith of $\bB$ is nonabelian, then $\bB \in HS(\bA)$ by Theorem \ref{jonsson-modular}, so we just need to handle the case where $0_\bB^* = (0_\bB : 0_\bB^*)$. In this case, Theorem \ref{jonsson-modular} implies that $\bB/0_\bB^* \in HS(\bA)$, so by Proposition \ref{monolith-bound} we have $|\bB| \le |\bA|^2 < \infty$.

Since $\bB$ is finite, we can write $\bB = \RR/\theta$ for some $\RR \le \bA^n$ and $\theta \in \Con(\RR)$. Then we can write $\RR$ as a subdirect product $\RR \le_{sd} \bA_1 \times \cdots \times \bA_m$ of finitely many subdirectly irreducible algebras $\bA_i \in HS(\bA)$. We assume that the $\bA_i$ are chosen such that none of them can be replaced by a subdirect product of some number of proper quotients of $\bA_i$ while still keeping the isomorphism $\RR/\theta \cong \bB$.

Then for any $i$, we must have $\theta \wedge \ker \pi_{[m]\setminus\{i\}} = 0_\RR$: if not, we could replace $\bA_i$ with a subdirect representation of $\RR/(\ker \pi_i \vee (\theta \wedge \ker \pi_{[m]\setminus\{i\}}))$, since by the modular law we have
\[
\ker \pi_{[m]\setminus \{i\}} \wedge (\ker \pi_i \vee (\theta \wedge \ker \pi_{[m]\setminus\{i\}})) = (\ker \pi_{[m]\setminus \{i\}} \wedge \ker \pi_i) \vee (\theta \wedge \ker \pi_{[m]\setminus\{i\}}) \le \theta.
\]
Since $\ker \pi_{[m]\setminus \{i\}} \ne 0_\RR$, we have $\theta \vee \ker \pi_{[m]\setminus \{i\}} \ge \theta^*$, so $\theta^* \wedge \ker \pi_{[m]\setminus \{i\}}$ is a cover of $0_\RR$, and we have
\[
\llbracket \theta, \theta^* \rrbracket \searrow \llbracket 0_\RR, \theta^* \wedge \ker \pi_{[m]\setminus \{i\}} \rrbracket \nearrow \llbracket \ker \pi_i, (\ker \pi_i)^* \rrbracket,
\]
so $\bB = \RR/\theta$ is similar to $\bA_i = \RR/\ker \pi_i$.
\end{proof}

\begin{ex} Let's work out what $D(\bG)$ is when $\bG$ is a subdirectly irreducible group. Let $\bM \lhd \bG$ be the normal subgroup corresponding to the monolith $0_\bG^*$, and let $\bN = C_{\bG}(\bM) \lhd \bG$ be the normal subgroup corresponding to the centralizer $(0_\bG : 0_\bG^*)$. First off, what is the group structure on the congruence $0_\bG^*$?

By definition, we have
\[
0_\bG^* = \{(x,y) \in \bG^2 \mid x^{-1}y \in \bM\}.
\]
We have a natural exact sequence of groups
\[
0 \rightarrow \bM \hookrightarrow 0_\bG^* \twoheadrightarrow \bG \rightarrow 0,
\]
where the inclusion is the map $m \mapsto (1,m)$ and the quotient map is the first projection $\pi_1$. The quotient $0_\bG^* \twoheadrightarrow \bG$ has a section $\Delta : \bG \hookrightarrow 0_\bG^*$ given by $g \mapsto (g,g)$. Thus we can write $0_\bG^*$ as a semidirect product
\[
0_\bG^* \cong \bM \rtimes \bG,
\]
where the action of $\bG$ on $\bM$ is the standard conjugation action.

How about the congruence $\Delta_{0_\bG^*}^{(0:0^*)} \in \Con(0_\bG^*)$? By Theorem \ref{difference-commutator}, we have
\[
\begin{bmatrix} x & w\\ y & z\end{bmatrix} \in \Delta_{0_\bG^*}^{(0:0^*)} \iff (xy^{-1}z = w) \wedge (x\ \equiv_{\bM}\ y\ \equiv_{\bN}\ z).
\]
Since this is a congruence on a group, we just need to understand the congruence class of the identity, so we plug in $x = y = 1$ and ask what values $(w,z)$ can take. We find that $\Delta_{0_\bG^*}^{(0:0^*)}$ corresponds to the normal subgroup
\[
\{(n,n) \mid n \in \bN\},
\]
so under the isomorphism $0_\bG^* \cong \bM \rtimes \bG$ it corresponds to $\bN$, considered as a subgroup of $\bG$. Thus we have
\[
D(\bG) = 0_\bG^*/\Delta_{0_\bG^*}^{(0:0^*)} \cong (\bM \rtimes \bG)/\bN \cong \bM \rtimes (\bG/\bN).
\]
That any of this makes sense follows from $\bN = C_{\bG}(\bM)$. We see that $\bM$ is the normal subgroup corresponding to the monolith of $D(\bG)$, that $\bM$ is equal to its own centralizer in $D(\bG)$, and that the natural map $\bG/\bN \hookrightarrow D(\bG)$ has image transverse to the monolith, and induces an isomorphism
\[
\sigma : \bG/\bN \xrightarrow{\sim} D(\bG)/\bM.
\]
To complete the description of the similarity from $\bG$ to $D(\bG)$, we let $\RR$ be the fiber product of $\bG$ and $D(\bG)$ over $\bG/\bN$, and define the congruence $\delta \in \Con(\RR)$ as the $4$-ary relation
\[
\begin{bmatrix} a & a'\\ b & b'\end{bmatrix} \in \delta\ \iff\ \begin{bmatrix} a\\ b \end{bmatrix}, \begin{bmatrix} a'\\ b'\end{bmatrix} \in \RR\ \wedge\ a^{-1}a' = b^{-1}b' \in \bM.
\]
That $\delta$ is closed under multiplication must be checked - it follows from the fact that $\bN$ centralizes $\bM$, and the fact that for any $a,b,a',b'$ satisfying the above conditions all of $a,b,a',b'$ must necessarily map to the same element of $\bG/\bN$.

What are the possible values for $D(\bG)$, assuming the monolith is abelian? Note that if we consider $\bM$ as a module via the $\bG/\bN$ action, then it must be a \emph{simple} module, since if it has any nontrivial submodule $\bM'$, then $\bM'$ will be a smaller normal subgroup of $\bG$. Thus the general situation is that $\bM$ is some simple module over the ring $\ZZ[\bG/\bN]$ (where $\bG/\bN$ acts faithfully on $\bM$), and $D(\bG) \cong \bM \rtimes (\bG/\bN)$.
\end{ex}

\begin{ex} If we take $\bG = S_3$ in the above, we find that $D(S_3) \cong \ZZ/3 \rtimes \ZZ/2 \cong S_3$. The $4$-ary relation $\delta \le S_3^{2\times 2}$ corresponding to the trivial similarity from $S_3$ to itself is given by
\[
\begin{bmatrix} a & a'\\ b & b'\end{bmatrix} \in \delta\ \iff\ s(a) = s(b) = s(a') = s(b')\ \wedge\ a^{-1}a' = b^{-1}b',
\]
where $s : S_3 \rightarrow \{\pm 1\}$ is the sign homomorphism.

We can think of the relation $\delta$ as having two ``strands'' corresponding to the two possible signs of permutations, and if we restrict to either strand then $\delta$ becomes an affine relation over $\ZZ/3$. The fact that we can multiply elements of $\delta$ which come from different strands and still get an element of $\delta$ is worth thinking about.

Now suppose that $\bG$ is some other subdirectly irreducible group such that $D(\bG) \cong S_3$, with monolith corresponding to $\bM \lhd \bG$ and $\bN = C_G(\bM)$. Then since $\bG$ is similar to $S_3$, we must have $\bM \cong \ZZ/3$ and $\bG/\bN \cong \ZZ/2$ by Corollary \ref{similar-detail}, with $\bG/\bN$ acting on $\bM$ by negation since $D(\bG) \cong \bM \rtimes (\bG/\bN) \cong S_3$. If the action of $\bG/\bN$ on $\bN$ is given by an involution $\tau$, then for any $n \in \bN \setminus \{1\}$ we must have $\bM$ contained in the normal subgroup of $\bN$ generated by $n, n^\tau$.

In particular, if $\bN$ is abelian then we see that $n + n^\tau, n - n^\tau \in \bM$ for all $n \in \bN$, and additionally in this case $\bN$ must have prime power order by Theorem \ref{subdirect-ab-prime-power}. Thus if $\bN$ is abelian then we must actually have $\bN = \bM$, and $\bG \cong S_3$.
\end{ex}

\begin{ex} If we take $\bG = Q_8 = \{\pm 1, \pm i, \pm j, \pm k\}$ the quaternion group with $i^2 = j^2 = k^2 = ijk = -1$, then the monolith is equal to the center, corresponding to the normal subgroup $\{\pm 1\}$, and the centralizer of the monolith is the full congruence $1_{Q_8}$. Thus
\[
D(Q_8) \cong \{\pm 1\} \cong \ZZ/2.
\]
The relation $\delta \le (Q_8 \times \ZZ/2)^2$ is then given by
\[
\begin{bmatrix} a & a'\\ b & b'\end{bmatrix} \in \delta\ \iff\ a' = (-1)^{b + b'}a.
\]
This relation closely resembles an affine relation over $\ZZ/2$.
\end{ex}



\chapter{Tame Congruence Theory}\label{a-tct}

Tame congruence theory was introduced by Hobby and McKenzie \cite{hobby-mckenzie} in order to answer questions about congruence lattices of finite algebras. Since every congruence contains the diagonal, every congruence is automatically invariant under the \emph{polynomial} clone of our algebra, and in fact the congruence lattice of an algebra is completely determined by the collection of \emph{unary} polynomials of the algebra.

\begin{prop}\label{prop-quasiorder-poly} An equivalence relation $\theta$ on an algebra $\bA$ is a congruence of $\bA$ iff $\theta$ is preserved by every unary polynomial operation of $\bA$. More generally, a quasiorder $\preceq$ on $\bA$ is a subalgebra of $\bA^2$ if and only if $\preceq$ is preserved by every unary polynomial of $\bA$.
\end{prop}
\begin{proof} Recall that a quasiorder is just a binary relation which contains the diagonal and is transitively closed, so any quasiorder which is preserved by every basic operation of $\bA$ will also be preserved by any unary polynomial of $\bA$.

Conversely, suppose that the quasiorder $\preceq$ is closed under all unary polynomials of $\bA$. Let $t$ be any $k$-ary term of $\bA$, and suppose that $a_i \preceq b_i$ for $i \in [k]$. Then for each $i$, we have
\[
t(b_1, ..., b_{i-1}, a_i, a_{i+1}, ..., a_k) \preceq t(b_1, ..., b_{i-1}, b_i, a_{i+1}, ..., a_k),
\]
since $\preceq$ is closed under the unary polynomial
\[
x \mapsto t(b_1, ..., b_{i-1}, x, a_{i+1}, ..., a_k).
\]
Since $\preceq$ is transitively closed, we can string these inequalities together to show that
\begin{align*}
t(a_1, a_2, ..., a_k) &\preceq t(b_1, a_2, ..., a_k)\\
&\preceq t(b_1, b_2, ..., a_k)\\
&\preceq \cdots\\
&\preceq t(b_1, b_2, ..., b_k).\qedhere
\end{align*}
\end{proof}

\begin{cor}\label{cor-quasiorder-poly} Suppose $\bA$ is an algebra and let $\Pol_1(\bA)$ be the set of unary polynomials of $\bA$. If $R$ is any quasiorder of $\bA$, then among all quasiorders on $\bA$ which are also subalgebras of $\bA^2$, the minimal compatible quasiorder containing $R$ is the transitive closure of the set
\[
\{(f(a), f(b)) \mid (a,b) \in R \text{ and } f \in \Pol_1(\bA)\},
\]
and the maximal compatible quasiorder contained in $R$ is given by
\[
\{(a,b) \mid \forall f \in \Pol_1(\bA), (f(a), f(b)) \in R\}.
\]
\end{cor}

So tame congruence theory is really about the how the identities satisfied in a (usually locally finite) variety affect the behavior of unary polynomials, and how this in turn affects the behavior of congruences. Because of the important role of polynomial operations in tame congruence theory, we will use $\Pol_n(\bA)$ to represent the set of $n$-ary polynomial operations of $\bA$ throughout this appendix (hopefully this doesn't cause any confusion with the notation for the polymorphism clone of a relational structure).

The material in this appendix is mostly taken from Hobby and McKenzie's wonderful book \cite{hobby-mckenzie} (some of it is my solutions to various exercises from their book).

\section{Shrinking algebras with unary polynomials, minimal sets, and traces}

As soon as you have a unary operation $\varphi$ on a finite set, the natural thing to do with it is to iterate it until we get the compositionally idempotent operation
\[
\varphi^{\infty} \coloneqq \lim_{n\rightarrow \infty} \varphi^{\circ n!}.
\]
This gives us a large collection of (compositionally) idempotent unary polynomial operations on any finite algebra. For unary operations, I will drop the qualifier ``compositionally'' on ``idempotent'', since idempotent unary operations in the usual sense are not very interesting.

\begin{defn} For any algebra $\bA$, we define $E(\bA)$ to be the set of unary polynomials $e \in \Pol_1(\bA)$ such that $e \circ e = e$. Elements of $E(\bA)$ might be called the \emph{idempotents} or \emph{projections} of $\bA$.
\end{defn}

Recall from Section \ref{s-partial-semi} that for any idempotent $e \in E(\bA)$, the clone of restrictions to $e(\bA)$ of polynomial operations of $\bA$ which preserve $e(\bA)$ is essentially the same as the clone of operations of the form
\[
(x_1, ..., x_n) \mapsto e(f(e(x_1), ..., e(x_n))),
\]
for $f \in \Pol_n(\bA)$ (strictly speaking, the $e$s on the inside are not really necessary, they are only there to stop us from caring about how $f$ behaves outside the set $e(\bA)$). This gives us a rich enough source of polynomial operations which preserve $e(\bA)$ to make it worth studying the restricted clone and introducing notation for it.

\begin{defn} If $\bA$ is an algebra and $U \subseteq \bA$, then we define the restriction $\Pol(\bA)|_U$ to be the set of restrictions $f|_U$ of polynomial operations $f \in \Pol(\bA)$ which preserve the subset $U$, and we define the \emph{induced algebra} $\bA|_U$ to be $(U, \Pol(\bA)|_U)$ (up to term equivalence).
\end{defn}

Restrictions are related to the congruence lattice using the following result.

\begin{lem}[P\'alfy and Pudl\'ak \cite{palfy-pudlak}, \cite{hobby-mckenzie}]\label{lem-idempotent-surjective-lattice} For any idempotent $e \in E(\bA)$, if we set $U = e(\bA)$, then the map taking the congruence $\theta \in \Con(\bA)$ to $e(\theta) = \theta|_U$ defines a surjective lattice homomorphism:
\[
\theta \mapsto \theta|_U : \Con(\bA) \twoheadrightarrow \Con(\bA|_U).
\]
More generally, if $N \subseteq U$, then:
\begin{itemize}
\item $\bA|_N = (\bA|_U)|_N$,
\item if $N$ is a union of $\theta|_U$ congruence classes, then the map $\alpha \mapsto \alpha|_N$ defines a lattice homomorphism from the interval $\llbracket 0_\bA, \theta \rrbracket$ of $\Con(\bA)$ to the interval $\llbracket 0_N, \theta|_N \rrbracket$ of $\Con(\bA|_N)$, and
\item if $N$ is equal to a congruence class of $\theta|_U$, then the map $\alpha \mapsto \alpha|_N$ is a surjective lattice homomorphism $\llbracket 0_\bA, \theta \rrbracket \twoheadrightarrow \Con(\bA|_N)$.
\end{itemize}
\end{lem}
\begin{proof} (Following \cite{hobby-mckenzie}) First, we will prove the statements about the map $\theta \mapsto \theta|_U$. For any $\alpha \in \Con(\bA|_U)$, we define the equivalence relation $\hat{\alpha}$ on $\bA$ to be the largest congruence of $\bA$ which is contained in the equivalence relation $e^{-1}(\alpha)$, which is given by
\[
\hat{\alpha} \coloneqq \{(x,y) \mid \forall f \in \Pol_1(\bA), (e(f(x)), e(f(y))) \in \alpha\}
\]
by Corollary \ref{cor-quasiorder-poly}. Then for any $\alpha \in \Con(\bA|_U)$ and $\theta \in \Con(\bA)$ we have $\hat{\alpha} \in \Con(\bA)$, $\hat{\alpha}|_U = \alpha$, and
\[
\theta|_U \le \alpha \; \iff \; \theta \le \hat{\alpha}.
\]
Since restriction obviously preserves meets of congruences, we just need to check that it preserves joins. For this, let $\alpha = \theta_1|_U \vee \theta_2|_U$, and note that
\[
\theta_i|_U \le \alpha \implies \theta_i \le \hat{\alpha} \implies \theta_1 \vee \theta_2 \le \hat{\alpha} \implies (\theta_1\vee\theta_2)|_U \le \alpha,
\]
while the inequality $\alpha \le (\theta_1\vee\theta_2)|_U$ is obvious.

To see that $\bA|_N = (\bA|_U)|_N$, note that if $f \in \Pol(\bA)$ preserves $N$, then $e\circ f$ preserves $U$ and $(e \circ f)|_N = f|_N$. Since we have $\alpha|_N = (\alpha|_U)|_N$ for $\alpha \in \Con(\bA)$, to prove the remaining claims we just need to think about the restriction map $\llbracket 0_U, \theta|_U \rrbracket \rightarrow \llbracket 0_N, \theta|_N \rrbracket$.

If $N$ is a union of congruence classes of $\theta|_U$, then for any $\theta_1, \theta_2 \le \theta|_U$ neither $\theta_i$ connects any element inside $N$ to any element outside $N$, so $(\theta_1 \vee \theta_2)|_N = \theta_1|_N \vee \theta_2|_N$. Thus the map $\llbracket 0_U, \theta|_U \rrbracket \rightarrow \llbracket 0_N, \theta|_N \rrbracket$ is a lattice homomorphism.

To finish, we need to show that if $N$ is a congruence class of $\theta|_U$ then this map is surjective. For this, we extend a congruence $\alpha \in \Con(\bA|_N)$ to the largest congruence $\check{\alpha}$ on $\bA|_U$ which is contained in
\[
\alpha \cup (U \setminus N)^2,
\]
which is given by
\[
\check{\alpha} \coloneqq \{(x,y) \mid \forall f \in \Pol_1(\bA|_U), f(x) \in N\text{ or }f(y)\in N \implies (f(x),f(y)) \in \alpha\}
\]
by Corollary \ref{cor-quasiorder-poly}. Since every $f \in \Pol_1(\bA|_U)$ preserves $\theta|_U$, if $f(x) \in N$ for any $x \in N$ then $f$ must preserve $N$, so we have $\check{\alpha}|_N = \alpha$.
\end{proof}

\begin{defn}\label{defn-restrict-order} Write $\bB \preceq_| \bA$ if there is some idempotent $e \in E(\bA)$, congruence $\theta \in \Con(\bA)$, and $a \in e(\bA)$, such that if we define
\[
U = e(\bA)
\]
and
\[
N = U \cap (a/\theta),
\]
then $\bB$ is polynomially equivalent to $\bA|_N$. (Note that in general $\bB$ will have a different signature than $\bA$.)
\end{defn}

\begin{prop}\label{prop-restrict-transitive} The relation $\preceq_|$ is transitively closed on finite algebras: if $\bC \preceq_| \bB \preceq_| \bA$ and $\bA$ is finite, then $\bC \preceq_| \bA$.
\end{prop}
\begin{proof} Suppose that $\bB = \bA|_N$ for $N = U \cap (a/\theta)$, $U = e(\bA)$, $e \in E(\bA)$. Additionally, suppose that $\bC = \bB|_{N'}$, for $N' = U' \cap (a'/\alpha)$, $\alpha \in \Con(\bB)$, $U' = e'(\bB)$, $e' \in E(\bB)$.

Since $\alpha$ is a congruence on $\bB$, by Lemma \ref{lem-idempotent-surjective-lattice} there is some congruence $\bar{\alpha} \in \llbracket 0_\bA, \theta \rrbracket$ such that $\bar{\alpha}|_N = \alpha$. Additionally, $e'$ is the restriction of some unary polynomial $\hat{e}'$ of $\bA$ to $N$, and by composing $e \circ \hat{e}'$ and iterating it, we get $\bar{e}' \in E(\bA)$ such that
\[
\bar{e}'(\bA) \subseteq U \;\; \text{ and } \;\; \bar{e}'(N) = U'.
\]
Thus we have
\[
N' \subseteq \bar{e}'(\bA) \cap (a'/\bar{\alpha}) \subseteq \bar{e}'(U \cap (a'/\bar{\alpha})) \cap (a'/\bar{\alpha}) \subseteq \bar{e}'(N) \cap (a'/\bar{\alpha}) \subseteq N'.
\]
To finish, we need to check that for every polynomial $f \in \Pol(\bA)$ which preserves $N'$, there is a polynomial $\bar{f} \in \Pol(\bA|_N)$ such that $\bar{f}|_{N'} = f|_{N'}$. If we take $\bar{f} = e \circ f$, then $\bar{f}$ automatically preserves $U$, and since $f$ preserves $N'$, we have
\[
\bar{f}(a', ..., a') \in a'/\bar{\alpha} \subseteq a/\theta,
\]
so $\bar{f}$ also preserves $a/\theta$, and therefore $\bar{f}$ preserves $U \cap (a/\theta) = N$.
\end{proof}


\begin{prop}\label{prop-restrict-variety} If $\bA$ is finite and $\bB \preceq_| \bA$ is such that every constant of $\bB$ is a term operation of $\bB$, and if $\bD \in HSP_{fin}(\bB)$, then there is some $\bC \in HSP_{fin}(\bA)$ such that $\bD \preceq_| \bC$.
\end{prop}
\begin{proof} Suppose that $\bB = \bA|_B$ with $B = e(\bA) \cap (a/\alpha)$, where $e \in E(\bA)$, $\alpha \in \Con(\bA)$, and $a \in e(\bA)$. We handle quotients and subpowers separately - for quotients, we will not need to assume that every constant of $\bB$ is a term operation of $\bB$. In fact, if $\bD = \bB/\theta$, we just choose $\bar{\theta} \in \llbracket 0_\bA, \alpha \rrbracket$ such that $\bar{\theta}|_{\bB} = \theta$ using Lemma \ref{lem-idempotent-surjective-lattice}, and take $\bC = \bA/\bar{\theta}$.

Now suppose that $\bD \le \bB^n$. Note that since every constant operation of $\bB$ is a term operaion of $\bB$, every $\bD \le \bB^n$ must contain the diagonal $\bB^{(n)} = \{b^{(n)} \mid b \in \bB\}$, where $b^{(n)} = (b, b, ..., b)$. Let $\bC \le \bA^n$ be given by
\[
\bC = \Sg_{\bA^n}(\bA^{(n)} \cup \bD).
\]
Note that $\bC$ is exactly the closure of $\bD$ under coordinatewise application of polynomials of $\bA$. Define $e^{(n)} \in \Pol(\bA^n)$ by replacing each constant $c$ in the definition of $e$ by $c^{(n)}$. Since each $c^{(n)}$ is also an element of $\bC$, we see that $e^{(n)}$ is also a polynomial of $\bC$, which acts like $e$ on each coordinate, so $e^{(n)} \in E(\bC)$. We need to check that if we set
\[
D = e^{(n)}(\bC) \cap (a^{(n)}/\alpha^n),
\]
then $D$ is the underlying set of $\bD$ and $\bC|_D$ is polynomially equivalent to $\bD$. This follows from the facts that $e^{(n)}(\bC)$ is the closure of $\bD$ under coordinatewise application of polynomials of $\bA$ which have been composed with $e$, and that a polynomial $f \in \Pol(\bA)$ has $f^{(n)}(\bD, ..., \bD) \cap (a^{(n)}/\alpha^n) \ne \emptyset$ iff $e \circ f$ preserves $B$.
\end{proof}

In order to get any use out of this to study an interval $\llbracket \alpha, \beta \rrbracket$ of the congruence lattice, we need to find idempotents $e$ such that $e(\alpha) \ne e(\beta)$, or equivalently such that $e(\beta) \not\subseteq \alpha$.

\begin{defn} If $\alpha < \beta \in \Con(\bA)$, then we define $U_\bA(\alpha, \beta)$ to be the collection of sets given by
\[
U_\bA(\alpha, \beta) = \{f(\bA) \mid f \in \Pol_1(\bA) \text{ and } f(\beta) \not\subseteq \alpha\}.
\]
We define $M_\bA(\alpha, \beta)$ to be the collection of minimal sets in $U_\bA(\alpha, \beta)$, and we call the sets in $M_\bA(\alpha, \beta)$ the $(\alpha, \beta)$-\emph{minimal sets} of $\bA$.
\end{defn}

\begin{prop}\label{prop-prime-tame} If $\beta$ is a cover of $\alpha$ in $\Con(\bA)$ and $\bA$ is finite, then for each $(\alpha,\beta)$-minimal set $U \in M_\bA(\alpha,\beta)$, there is some $e \in E(\bA)$ such that $U = e(\bA)$.
\end{prop}
\begin{proof} Pick any $g \in \Pol_1(\bA)$ such that $g(\bA) = U$ and $g(\beta) \not\subseteq \alpha$. Then since $\beta$ covers $\alpha$ and $g(\beta) \subseteq \beta$, the congruence generated by $g(\beta) \cup \alpha$ must be $\beta$. Thus for any $(x,y) \in \beta$, there must be some $h_i \in \Pol_1(\bA)$ and $(u_i,v_i) \in g(\beta) \cup \alpha$ such that $x = h_1(u_1)$, $h_i(v_i) = h_{i+1}(u_{i+1})$, and $h_n(v_n) = y$ for some $n$.

If we choose $(x,y) \in \beta$ such that $(g(x),g(y)) \not\in \alpha$, we see that there must be some $i$ such that $(g(h_i(u_i)), g(h_i(v_i))) \not\in \alpha$, and for this $i$ we must have $(u_i,v_i) \in g(\beta)$. Setting $f = g\circ h_i$, we see that
\[
f(g(\beta)) \not\subseteq \alpha
\]
and
\[
f(g(\bA)) \subseteq U.
\]
Since $U$ is $(\alpha,\beta)$-minimal, we must have $f(g(\bA)) = U$, so $f(U) = U$. Iterating $f$ gives us $e = f^{\infty} \in E(\bA)$ with $e(\bA) = U$.
\end{proof}

We will also want to prove similar results for certain other pairs of congruences $\alpha < \beta$. Precisely stating what we are going for requires a bit more work.

\begin{defn} A $0,1$-\emph{lattice} is defined to be a lattice with constants $0$ and $1$ which satisfy $0 \le x \le 1$ for all $x$. Note that every interval $\llbracket \alpha, \beta \rrbracket$ in a lattice can be regarded as a $0,1$-lattice, with $0$ interpreted as $\alpha$ and $1$ interpreted as $\beta$.

A $0,1$-\emph{separating homomorphism} is a lattice homomorphism such that $f^{-1}(f(0)) = \{0\}$ and $f^{-1}(f(1)) = \{1\}$.
\end{defn}

\begin{defn}\label{defn-tame} A \emph{congruence quotient} is defined to be an ordered pair of congruences $(\alpha,\beta)$ such that $\alpha < \beta$. A congruence quotient $(\alpha,\beta)$ is called \emph{prime} if $\beta$ covers $\alpha$.

A congruence quotient $(\alpha,\beta)$ is called \emph{tame} if there is some $U \in M_\bA(\alpha,\beta)$ and some $e \in E(\bA)$ such that $e(\bA) = U$, and such that the restriction homomorphism
\[
\llbracket \alpha, \beta \rrbracket \twoheadrightarrow \llbracket \alpha|_U, \beta|_U \rrbracket
\]
is a $0,1$-separating homomorphism, that is, for $\alpha \le \gamma \le \beta$ we have $\gamma|_U = \alpha|_U \implies \gamma = \alpha$ and $\gamma|_U = \beta|_U \implies \gamma = \beta$. An algebra $\bA$ is called \emph{tame} if $(0_\bA,1_\bA)$ is tame.
\end{defn}

So far we have shown that every prime quotient on a finite algebra is tame. There is a more general lattice theoretic condition that implies tameness, but first we should try to see what being tame is good for. The first important result is that all of the minimal sets for a tame quotient look the same as each other.

\begin{defn} If $U,V \subseteq \bA$, then we say that $U,V$ are \emph{polynomially isomorphic} in $\bA$ if there are unary polynomials $f,g \in \Pol_1(\bA)$ such that
\[
f(U) = V, \; g(V) = U, \; g\circ f|_U = \operatorname{id}_U, \; f\circ g|_V = \operatorname{id}_V.
\]
In this case we write $f : U \simeq V$.
\end{defn}

\begin{prop} If $f : U \simeq V$ in $\bA$, then $f|_U$ defines an isomorphism from $\bA|_U$ to $\bA|_V$ (up to term equivalence). Furthermore, for any $\theta \in \Con(\bA)$ we have $f|_U(\theta|_U) = \theta|_V$.
\end{prop}

\begin{prop} If $U,V \subseteq \bA$ and $\bA$ is finite, then $U \simeq V$ iff there are $f, g \in \Pol_1(\bA)$ such that $f(U) = V$ and $g(V) = U$. If additionally we have $f(\bA) = V$, then there is some idempotent $e \in E(\bA)$ such that $e(\bA) = U$.
\end{prop}
\begin{proof} Take $g' = (g \circ f)^{\infty - 1}\circ g$ and $e = g' \circ f = (g \circ f)^{\infty}$, so $e \in E(\bA)$. Then $g'(V) = U$, $g'\circ f|_U = \operatorname{id}_U$, $f\circ g'|_V = \operatorname{id}_V$, and if $f(\bA) = V$ then $e(\bA) = g'(f(\bA)) = g'(V) = U$.
\end{proof}

\begin{thm}[Minimal sets for tame quotients \cite{hobby-mckenzie}]\label{thm-minimal-sets} If $(\alpha,\beta)$ is a tame congruence quotient on a finite algebra $\bA$, then all of the following are true.
\begin{itemize}
\item[(a)] For all $U,V \in M_\bA(\alpha,\beta)$, we have $U \simeq V$ in $\bA$.

\item[(b)] For all $U \in M_\bA(\alpha,\beta)$, there is some $e \in E(\bA)$ such that $e(\bA) = U$, and the restriction homomorphism $\llbracket \alpha, \beta \rrbracket \twoheadrightarrow \llbracket \alpha|_U, \beta|_U \rrbracket$ is a $0,1$-separating homomorphism.

\item[(c)] For all $U \in M_\bA(\alpha, \beta)$ and $(x,y) \in \beta\setminus\alpha$, there is some $f \in \Pol_1(\bA)$ such that $f(\bA) = U$ and $(f(x),f(y)) \not\in \alpha$.

\item[(d)] For all $U \in M_\bA(\alpha, \beta)$, $\beta$ is the transitive closure of $\alpha \cup \bigcup_{g \in \Pol_1(\bA)} g(\beta|_U)$.

\item[(e)] For all $U \in M_\bA(\alpha, \beta)$ and $f \in \Pol_1(\bA)$ such that $f(\beta|_U) \not\subseteq \alpha$, we have $f(U) \in M_\bA(\alpha, \beta)$ and $f : U \simeq f(U)$.

\item[(f)] For any $f \in \Pol_1(\bA)$ such that $f(\beta) \not\subseteq \alpha$, there is some $U \in M_\bA(\alpha,\beta)$ such that $f : U \simeq f(U)$.

\item[(g)] For any $(x,y) \in \beta\setminus\alpha$, there is some $U \in M_\bA(\alpha, \beta)$ and $e \in E(\bA)$ such that $e(\bA) = U$ and $(e(x), e(y)) \not\in \alpha$.
\end{itemize}
\end{thm}
\begin{proof} (Following \cite{hobby-mckenzie}) By the definition of tameness, there is \emph{some} $U \in M_\bA(\alpha, \beta)$ that satisfies (b). We will first show that (c) and (d) hold for this $U$, and then use this to prove (a), which will imply that (b) is true in general. Then we will use these to prove (e), (f), and (g).

Suppose that (b) holds for $U$. To prove (c), let $\gamma \in \llbracket \alpha, \beta \rrbracket$ be the congruence generated by $\alpha$ and $(x,y)$. Then since $\gamma \ne \alpha$ we must have $\gamma|_U \ne \alpha|_U$ (since restriction is $0,1$-separating). Since $\gamma|_U = e(\gamma)$, this means that $\gamma \not\subseteq e^{-1}(\alpha)$, so there must be some $g \in \Pol_1(\bA)$ such that $(g(x), g(y)) \not\in e^{-1}(\alpha)$. Taking $f = e\circ g$ proves (c).

To see that (b) implies (d), let $\gamma \in \llbracket \alpha, \beta \rrbracket$ be the transitive closure of $\alpha \cup \bigcup_{g \in \Pol_1(\bA)} g(\beta|_U)$. Then $\gamma|_U = \beta|_U$, so we must have $\gamma = \beta$ (since restriction is $0,1$-separating).

Now suppose that $U,V \in M_\bA(\alpha, \beta)$ and that $V$ satisfies (b), (c), (d). Since $U \in M_\bA(\alpha,\beta)$, there must be some $h \in \Pol_1(\bA)$ and $(x,y) \in h(\beta) \setminus \alpha$ with $h(\bA) = U$, so by (c) applied to $V$ there is some $f \in \Pol_1(\bA)$ such that $f(\bA) = V$ and $(f(x), f(y)) \not\in \alpha$. Then from $f(h(\bA)) \subseteq V$ and $f(h(\beta)) \not\subseteq \alpha$ we have $f(h(\bA)) = V$ by $(\alpha,\beta)$-minimality of $V$, so $f(U) = V$. Next, by (d) applied to $V$ we see that there must be some $g \in \Pol_1(\bA)$ such that $g(\beta|_V) \not\subseteq h^{-1}(\alpha)$. Then if $e \in E(\bA)$ has $e(\bA) = V$, we see that
\[
h(g(e(\beta))) \not\subseteq \alpha
\]
and $h(g(e(\bA))) \subseteq U$, so since $U$ is $(\alpha,\beta)$-minimal we see that $h(g(V)) = U$. Now we can apply the previous proposition to see that $f : U \simeq V$ and that $U$ satisfies (b) as well.

To prove (e), we use (b) to see that there is some $e \in E(\bA)$ with $e(\bA) = U$, and note that $f(e(\beta)) = f(\beta|_U) \not\subseteq \alpha$, so $f(U) = f(e(\bA))$ must contain some minimal $V \in M_\bA(\alpha,\beta)$. By (a) we see that $|V| = |U|$, so we must in fact have $f(U) = V$ and $f : U \simeq V$.

To prove (f), we apply (d) to any $V \in M_\bA(\alpha, \beta)$ to see that there is some $g \in \Pol_1(\bA)$ such that $g(\beta|_V) \not\subseteq f^{-1}(\alpha)$. Then by applying (e) twice we see that we can take $U = g(V)$.

To prove (g), we apply (c) to any $V \in M_\bA(\alpha, \beta)$ to see that there is some $f \in \Pol_1(\bA)$ such that $f(\bA) = V$ and $(f(x), f(y)) \not\in \alpha$. By (f), there is some $U \in M_\bA(\alpha, \beta)$ such that $f : U \simeq f(U)$, and since $f(U) \subseteq f(\bA) = V$, we must have $f(U) = V$ by $(\alpha,\beta)$-minimality. Thus there is some $g \in \Pol_1(\bA)$ such that $g\circ f|_U = \operatorname{id}_U$, and we can take $e = g \circ f$.
\end{proof}

\begin{cor} If $(\alpha,\beta)$ is a tame congruence quotient on a finite algebra $\bA$ and $U \in M_\bA(\alpha,\beta)$, then every unary polynomial $f$ of $\bA|_U$ is either a permutation of $U$ or has $f(\beta|_U) \subseteq \alpha|_U$.
\end{cor}

If we restrict to a congruence class of $\beta|_U$, we get an even stronger result.

\begin{defn} If $(\alpha,\beta)$ is a tame congruence quotient on $\bA$ and $U \in M_\bA(\alpha,\beta)$, then a set $N \subseteq U$ is called an $(\alpha,\beta)$-\emph{trace} in $U$ if $N$ is a congruence class of $\beta|_U$ which is not also a congruence class of $\alpha|_U$.

We define the \emph{body} of the $(\alpha,\beta)$-minimal set $U$ to be the union of the $(\alpha,\beta)$-traces, and we define the \emph{tail} of $U$ to be the set of congruence classes of $\beta|_U$ which are also congruence classes of $\alpha|_U$.
\end{defn}

Since $\beta|_U \not\subseteq \alpha|_U$ by the definition of an $(\alpha,\beta)$-minimal set, we see that every $(\alpha,\beta)$-minimal set $U$ has a nonempty body, i.e. there is at least one $(\alpha,\beta)$-trace $N$ in $U$.

\begin{cor}\label{cor-trace-permutational} If $(\alpha,\beta)$ is a tame congruence quotient on a finite algebra $\bA$ and $N$ is an $(\alpha,\beta)$-trace, then every unary polynomial $f$ of $\bA|_N$ is either a permutation, or has $f(N)$ contained in some congruence class of $\alpha|_N$. In particular, every unary polynomial of $\bA|_N/\alpha|_N$ is either a permutation or is constant.
\end{cor}

\begin{defn} We say that an algebra is \emph{permutational} or \emph{minimal} if every unary polynomial is either a permutation or is constant.

More generally, if $(\alpha,\beta)$ is a congruence quotient on $\bA$, we say that $\bA$ is $(\alpha,\beta)$-\emph{minimal} if every unary polynomial $f$ of $\bA$ is either a permutation or has $f(\beta) \subseteq \alpha$. Note that in this case, $(\alpha,\beta)$ is necessarily tame, with $M_\bA(\alpha,\beta) = \{A\}$.
\end{defn}

The general strategy will be to understand an algebra by first understanding the structure of the traces, and then reconstructing the algebra from the traces. When we apply unary polynomials of $\bA$ to $(\alpha,\beta)$-traces, the result will often also be an $(\alpha,\beta)$-trace.

\begin{cor}\label{cor-trace-iso} If $(\alpha,\beta)$ is a tame congruence quotient on a finite algebra $\bA$ and $N$ is an $(\alpha,\beta)$-trace, then for every unary polynomial $f$ of $\bA$ either $f(N)$ is contained in some congruence class of $\alpha$, or $f(N)$ is another $(\alpha,\beta)$-trace and $f : N \simeq f(N)$.
\end{cor}
\begin{proof} This follows directly from Theorem \ref{thm-minimal-sets}(e).
\end{proof}

For the purpose of stitching together the traces to reconstruct the algebra, we have another consequence of Theorem \ref{thm-minimal-sets}.

\begin{cor}\label{cor-trace-closure} If $(\alpha,\beta)$ is a tame congruence quotient on a finite algebra $\bA$, then $\beta$ is the transitive closure of
\[
\alpha \cup \{N^2 \mid N \text{ is an }(\alpha,\beta)\text{-trace}\}.
\]
\end{cor}
\begin{proof} This follows from Theorem \ref{thm-minimal-sets}(d) and the previous corollary, once we note that for each $(\alpha,\beta)$-minimal set $U$, every congruence class of $\beta|_U$ is either a congruence class of $\alpha|_U$ or is an $(\alpha,\beta)$-trace.
\end{proof}

\begin{prop}\label{prop-prime-poly-iso} If $(\alpha,\beta)$ is a prime congruence quotient of a finite algebra $\bA$ (i.e., if $\beta$ is a cover of $\alpha$), then all of the $(\alpha,\beta)$-traces are polynomially isomorphic in $\bA$.
\end{prop}
\begin{proof} Let $U$ be any $(\alpha,\beta)$-minimal set. By Theorem \ref{thm-minimal-sets}(a) it's enough to show that any pair of $(\alpha,\beta)$-traces $N$, $K$ contained in $U$ are polynomially isomorphic in $\bA|_U$. By Lemma \ref{lem-idempotent-surjective-lattice}, the restriction homomorphism
\[
\llbracket \alpha, \beta \rrbracket \rightarrow \llbracket \alpha|_U, \beta|_U \rrbracket
\]
is surjective, so $\beta|_U$ covers $\alpha|_U$ in $\bA|_U$. Then since $N$ is not contained in a congruence class of $\alpha$, the congruence of $\bA|_U$ generated by $\alpha|_U \cup N^2$ must be $\beta|_U$. In particular, there must be some unary polynomial $f \in \Pol_1(\bA|_U)$ such that $f(N) \subseteq K$ and $f(N)$ is not contained in any congruence class of $\alpha$. Then $f$ must be a permutation of $U$, and so we have $f : N \simeq K$.
\end{proof}

Tameness can also be derived from some of its consequences.

\begin{prop} If $\alpha < \beta \in \Con(\bA)$ and there is some finite $(\alpha,\beta)$-minimal set $U$ which satisfies Theorem \ref{thm-minimal-sets}(c) and \ref{thm-minimal-sets}(d), then $(\alpha,\beta)$ is tame.
\end{prop}
\begin{proof} If we make a digraph on $(x,y) \in \beta|_U\setminus\alpha|_U$ with an edge from $(x,y)$ to $(u,v)$ whenever there is some $f \in \Pol_1(\bA)$ with $f(\bA) = U$ and $f(x) = u, f(y) = v$, then Theorem \ref{thm-minimal-sets}(c) for $U$ implies that every vertex in this digraph has outdegree at least one, so there must be a directed cycle.

By composing the unary polymorphisms corresponding to the edges of this directed cycle, we find an $f \in \Pol_1(\bA)$ such that $f(\bA) \subseteq U$ and $f(x) = x, f(y) = y$ for some $(x,y) \in \beta|_U\setminus\alpha|_U$. Taking $e = f^\infty$, we see that $e \in E(\bA)$ with $e(\bA) \subseteq U$ and $e(\beta) \not\subseteq \alpha$. Finally, Theorem \ref{thm-minimal-sets}(c) and \ref{thm-minimal-sets}(d) for $U$ directly imply that the restriction homomorphism $\llbracket \alpha, \beta \rrbracket \twoheadrightarrow \llbracket \alpha|_U, \beta|_U \rrbracket$ is $0,1$-separating.
\end{proof}

Many of the previous results simplify if $\alpha = 0_\bA$. To relate the general case to that situation, we need to know how tameness behaves when we pass to a quotient.

\begin{prop}\label{prop-tame-quotient} If $\delta \le \alpha < \beta$ are congruences on a finite algebra $\bA$, then $(\alpha,\beta)$ is tame on $\bA$ iff $(\alpha/\delta, \beta/\delta)$ is tame on $\bA/\delta$. If $(\alpha,\beta)$ is tame, then we have
\[
M_{\bA/\delta}(\alpha/\delta,\beta/\delta) = \{U/\delta \mid U \in M_\bA(\alpha,\beta)\}.
\]
In particular, the $(\alpha/\delta,\beta/\delta)$-traces are exactly the quotients of the $(\alpha,\beta)$-traces by $\delta$.
\end{prop}
\begin{proof} (Following \cite{hobby-mckenzie}) For any unary polynomial $f \in \Pol_1(\bA)$, we have $f(\beta) \subseteq \alpha$ iff $f(\beta/\delta) \subseteq \alpha/\delta$ in $\bA/\delta$, and for any $\gamma, \gamma' \in \llbracket \alpha, \beta \rrbracket$ we have $\gamma = \gamma' \iff \gamma/\delta = \gamma'/\delta$ and similarly for restrictions. The challenge is showing that the minimal sets correspond.

First we show that if $U \in M_\bA(\alpha,\beta)$ and there is an $e \in E(\bA)$ with $e(\bA) = U$, then $U/\delta \in M_{\bA/\delta}(\alpha/\delta,\beta/\delta)$. To see this, suppose that some $f \in \Pol_1(\bA)$ has $f(\beta) \not\subseteq \alpha$ and $f(\bA/\delta) \subseteq U/\delta$. Then $e(f(\bA)) \subseteq U$ and $e(f(\beta)) \not\subseteq \alpha$, so $e(f(\bA)) = U$, so $f(\bA/\delta) = U/\delta$. This shows that if $(\alpha,\beta)$ is tame then $(\alpha/\delta, \beta/\delta)$ is tame.

Now suppose $(\alpha/\delta, \beta/\delta)$ is tame, and let $U$ be any $(\alpha,\beta)$-minimal set. Pick $f \in \Pol_1(\bA)$ with $f(\bA) = U$ and $f(\beta) \not\subseteq \alpha$. By Theorem \ref{thm-minimal-sets}(f) we see that there is some
\[
V \in M_{\bA/\delta}(\alpha/\delta, \beta/\delta)
\]
such that $f : V \simeq f(V)$. Pick $g \in \Pol_1(\bA)$ such that $g : f(V) \rightarrow V$ inverts $f : V \rightarrow f(V)$, then by iterating $f\circ g$ we get an idempotent $e \in E(\bA)$ with $e(f(V)) = f(V)$. Thus $e(\beta) \not\subseteq \alpha$, and from
\[
e(\bA) \subseteq f(\bA) = U
\]
we get $e(\bA) = U$. Then by the previous paragraph we see that $U/\delta$ is $(\alpha/\delta,\beta/\delta)$-minimal, which allows us to conclude that the restriction homomorphism $\llbracket \alpha, \beta \rrbracket \twoheadrightarrow \llbracket \alpha|_U, \beta|_U \rrbracket$ is $0,1$-separating.

To finish, we need to show that any $V \in M_{\bA/\delta}(\alpha/\delta,\beta/\delta)$ is a quotient of an $(\alpha,\beta)$-minimal set when $(\alpha,\beta)$ is tame. Pick any $U \in M_\bA(\alpha,\beta)$, then since $U/\delta$ is $(\alpha/\delta,\beta/\delta)$-minimal we can apply Theorem \ref{thm-minimal-sets}(a) to see that there is an $f \in \Pol_1(\bA)$ with $f : U/\delta \simeq V$. Then $f(\beta|_U) \not\subseteq \alpha$, so by Theorem \ref{thm-minimal-sets}(e) we have $f(U) \in M_\bA(\alpha,\beta)$, and $V = f(U)/\delta$.
\end{proof}

\begin{prop} If $\alpha \le \gamma < \beta$ and $(\alpha,\beta), (\gamma,\beta)$ are both tame quotients on a finite algebra $\bA$, then $M_\bA(\alpha,\beta) = M_\bA(\gamma,\beta)$.
\end{prop}
\begin{proof} If $U \in M_\bA(\alpha,\beta)$, then $\beta|_U \not\subseteq \gamma|_U$ since the restriction map is $0,1$-separating (by Theorem \ref{thm-minimal-sets}(b)), so $U \in U_\bA(\gamma,\beta)$. If $f(\bA) \subseteq U$ and $f(\beta) \not\subseteq \gamma$, then we have $f(\beta) \not\subseteq \alpha$, so $f(\bA) = U$ by $(\alpha,\beta)$-minimality, so $U \in M_\bA(\gamma,\beta)$.

Conversely, if $V \in M_\bA(\gamma,\beta)$, then by Theorem \ref{thm-minimal-sets}(a) we have $V \simeq U$ for some $U \in M_\bA(\alpha,\beta)$, so $V \in M_\bA(\alpha,\beta)$ by Theorem \ref{thm-minimal-sets}(e).
\end{proof}

Recall that two intervals are perspective, written $\llbracket \alpha, \beta \rrbracket \searrow \llbracket \gamma, \delta \rrbracket$, if $\alpha \wedge \delta = \gamma$ and $\alpha \vee \delta = \beta$.

\begin{prop}\label{prop-minimal-sets-perspective} If $\llbracket \alpha, \beta \rrbracket \searrow \llbracket \gamma, \delta \rrbracket$ in $\Con(\bA)$, then $M_\bA(\alpha,\beta) = M_\bA(\gamma,\delta)$ and $M_\bA(\gamma,\beta) \subseteq M_\bA(\gamma,\alpha) \cup M_\bA(\gamma, \delta)$.
\end{prop}
\begin{proof} For any $f \in \Pol_1(\bA)$, we have
\[
f(\beta) = f(\alpha \vee \delta) \subseteq \alpha \iff f(\alpha) \cup f(\delta) \subseteq \alpha \iff f(\delta) \subseteq \alpha \iff f(\delta) \subseteq \alpha \cap \delta = \gamma,
\]
so $M_\bA(\alpha,\beta) = M_\bA(\gamma,\delta)$. Similarly, we have $f(\beta) \subseteq \gamma$ iff $f(\beta) \subseteq \alpha$ and $f(\beta) \subseteq \delta$, so
\[
U_\bA(\gamma,\beta) = U_\bA(\alpha,\beta) \cup U_\bA(\delta,\beta),
\]
so $M_\bA(\gamma,\beta) \subseteq M_\bA(\alpha,\beta) \cup M_\bA(\delta,\beta) = M_\bA(\gamma,\alpha) \cup M_\bA(\gamma, \delta)$.
\end{proof}

\begin{cor} If $\alpha,\beta$ are congruences on a finite algebra and the interval $\llbracket \alpha, \beta \rrbracket$ is isomorphic to the diamond lattice $\cM_n$ for some $n \ge 3$, then $(\alpha, \beta)$ is tame and every congruence quotient contained in $\llbracket \alpha, \beta \rrbracket$ has the same collection of minimal sets.
\end{cor}

\begin{prop} If $(\alpha,\beta)$ is a tame congruence quotient on a finite algebra $\bA$ and $U$ is an $(\alpha,\beta)$-minimal set, then for any $\gamma' < \delta' \in \llbracket \alpha|_U ,\beta|_U \rrbracket$, there are lifts $\gamma < \delta \in \llbracket \alpha, \beta \rrbracket$ such that $\gamma|_U = \gamma', \delta|_U = \delta'$, with $(\gamma,\delta)$ a tame quotient. For any such $\gamma, \delta$ we have $M_\bA(\gamma,\delta) = M_\bA(\alpha,\beta)$.
\end{prop}
\begin{proof} Since any unary polynomial which collapses $\beta$ into $\alpha$ necessarily collapses $\delta$ into $\gamma$ for any $\gamma,\delta \in \llbracket \alpha, \beta \rrbracket$, we see that $U$ is $(\gamma,\delta)$-minimal for any $\gamma < \delta$ which restrict to $\gamma', \delta'$, so we just need to ensure that the restriction homomorphism from $\llbracket \gamma, \delta \rrbracket$ is $0,1$-separating.

Since restriction to $U$ is a lattice homomorphism, we can take $\gamma$ to be maximal among congruences which restrict to $\gamma'$ and are $\le \beta$, and take $\delta$ to be minimal among congruences which restrict to $\delta'$ and are $\ge \gamma$. By Theorem \ref{thm-minimal-sets}(a) and \ref{thm-minimal-sets}(e) every $(\gamma,\delta)$-minimal set $V$ has $U \simeq V$ and so is also $(\alpha,\beta)$-minimal, and similarly every $(\alpha,\beta)$-minimal set is $(\gamma,\delta)$-minimal.
\end{proof}

Some basic examples to keep in mind follow.

\begin{ex} If $\bA$ is a finite lattice and $\alpha < \beta \in \Con(\bA)$, then the $(\alpha,\beta)$-minimal sets all have the form $\{a,b\}$ with $a < b$ and $(a,b) \in \beta\setminus \alpha$, and any pair $\{a,b\}$ which satisfies those conditions and additionally has $b$ covering $a$ is an $(\alpha,\beta)$-minimal set. Each such set with $b$ covering $a$ is the image of the idempotent unary polynomial
\[
x \mapsto (x \wedge b) \vee a,
\]
however, in order for the restriction homomorphism to be $0,1$-separating, $\beta$ must be a cover of $\alpha$. Thus the tame congruence quotients of $\bA$ are exactly the same as the prime quotients. Additionally, every $(\alpha,\beta)$-minimal set consists of just a single $(\alpha,\beta)$-trace (i.e., the minimal sets have no tails), and every trace is the image of some pair $\{a,b\}$ with $b$ covering $a$ under some unary polynomial of $\bA$, and is polynomially equivalent to a two element lattice.
\end{ex}

\begin{ex} If $\bA$ is a finite module over a ring $\RR$ (which we may assume acts faithfully on $\bA$ without loss of generality) and $\alpha < \beta \in \Con(\bA)$, then we can represent the congruences $\alpha,\beta$ by the submodules $\bM_\alpha = 0/\alpha$, $\bM_\beta = 0/\beta$ of $\bA$. Every unary polynomial $f$ of $\bA$ has the form $f : x \mapsto rx + c$, and we have
\[
f(\beta) \subseteq \alpha \iff r\bM_\beta \subseteq \bM_\alpha.
\]
The set of $r \in \RR$ such that $r \bM_\beta \subseteq \bM_\alpha$ is called the \emph{annihilator} of $\bM_\beta/\bM_\alpha$, and forms a (two-sided) ideal $\bI$ of $\RR$. Then $\bM_\beta/\bM_\alpha$ is a module over $\RR/\bI$, and $\RR/\bI$ acts faithfully on $\bM_\beta/\bM_\alpha$.

We will show that $(\alpha,\beta)$ is a tame congruence quotient if and only if $\bI$ is maximal among two-sided ideals of $\RR$, that is, iff $\RR/\bI$ is a simple ring. By the classification of finite simple rings, this is equivalent to proving that $\RR/\bI$ is isomorphic to a matrix ring $M_n(\bF_{p^k})$ over a finite field $\bF_{p^k}$ for some $n$ and some prime power $p^k$. Additionally, we will show that $\bM_\beta/\bM_\alpha$ is one of the modules $\bF_{p^k}^{n\times m}$ for some $m$, with the action of $\RR/\bI$ on $\bM_\beta/\bM_\alpha$ given by matrix multiplication. To prove this, it is simpler to think only about the module $\bM_\beta/\bM_\alpha$ - note that $\bM_\beta/\bM_\alpha$ is a tame algebra, with $(e\bA\cap \bM_\beta)/\bM_\alpha = e\bM_\beta/\bM_\alpha$ as a $(0,1)$-minimal set, for any $e \in \RR$ such that $e^2 = e$ and $e\bA$ is an $(\alpha,\beta)$-minimal set.
\end{ex}

\begin{prop} Suppose that $\bM$ is a finite module over a finite ring $\RR$ which acts faithfully on $\bM$, and suppose that $\bM$ is a tame algebra. Then $\RR$ is isomorphic to a matrix ring $M_n(\bF_{p^k})$ over a finite field $\bF_{p^k}$ for some $n$ and some prime power $p^k$, and $\bM$ is $\bF_{p^k}^{n\times m}$ for some $m$.
\end{prop}
\begin{proof}
By Corollary \ref{cor-trace-permutational}, every $(0_\bM,1_\bM)$-trace $N$ has $\bM|_N$ a permutational algebra. If $N = e\bM$ for some $e \in E(\bA)$, the restriction $\Pol(\bM)|_N$ consists of the linear functions with coefficients in the ring $e\RR e$, so we see that every nonzero element of the ring $e\RR e$ is invertible, that is, $e\RR e$ is a division ring. Since every finite division ring is a field by Wedderburn's little theorem, we have $e \RR e \cong \bF_{p^k}$ for some finite field $\bF_{p^k}$ (note, the rest of the argument still works over a division ring rather than a field).

We claim that we can pick $e_1, ..., e_n \in \RR$ idempotent such that $\sum_i e_i = 1$, $e_ie_j = e_je_i = 0$ for $i \ne j$, and each $e_i\bM$ is $(0_\bM,1_\bM)$-minimal. Suppose that $e_1, ..., e_n$ is a maximal collection of idempotents such that $e_ie_j = e_je_i = 0$ for $i \ne j$ and each $e_i\bM$ is $(0_\bM,1_\bM)$-minimal, and let $f = 1 - \sum_i e_i$ (that such a maximal set exists and is finite follows from the fact that $e_i\bM \cap e_j\bM = \emptyset$ for $i \ne j$, since $e_ie_j = 0$). Then we have $f^2 = f$, and if $f \ne 0$ then there must be some $(0_\bM,1_\bM)$-minimal set $U \subseteq f\bM$. If $e' \in E(\bM)$ has $e'\bM = U$, then we set $e_{n+1} = e'f$, and note that we have $e_{n+1}\bM = e'\bM = U$, and $e_{n+1}e_i = e_ie_{n+1} = 0$ for each $i \le n$, contradicting the maximalilty of the collection $e_1, ..., e_n$. Therefore we must have $f = 0$, that is, $\sum_i e_i = 1$. Note that every element $x \in \bM$ has a unique decomposition
\[
x = \sum_i e_ix_i,
\]
so as a group, $\bM$ is the direct sum of the $e_i\bM$s.

By Theorem \ref{thm-minimal-sets}(a), for each pair $i,j \le n$ we have $e_i\bM \simeq e_j\bM$. Pick $f_i : e_1\bM \simeq e_i\bM$ for each $i$, and pick inverses $g_i : e_i\bM \simeq e_1\bM$ to each $f_i$ (with $f_1 = g_1 = e_1$ for $i = 1$). Then for any $i,j$, define the matrix element $e_{ij}$ by
\[
e_{ij} = f_ig_je_j,
\]
and note that each $e_{ij}$ is an isomorphism $e_{ij} : e_j\bM \simeq e_i\bM$, with $e_ie_{ij} = e_{ij}$ and $e_{ij}e_j = e_{ij}$, with $e_{ij}e_{ji} = e_{ii} = e_i$, with $e_{ij}e_{jk} = e_{ik}$, and $e_{ij}e_{kl} = 0$ for $j \ne k$.
For each $r_1 \in e_1\RR e_1$ we can additionally define the corresponding scalar $r \in \RR$ by
\[
r_1 \in e_1\RR e_1 \mapsto r = \sum_i e_{i1}r_1e_{1i},
\]
and we identify the set of such scalars $r$ with $\bF_{p^k}$, noting that $re_{ij} = e_{ij}r$ for all $r \in \bF_{p^k}$ and $i,j \le n$, and that the multiplication in $\bF_{p^k}$ is the same as the multiplication in $e_1\RR e_1$. We claim that every element $m \in \RR$ can be written uniquely in the form
\[
m = \sum_{i,j} r_{i,j}e_{ij}
\]
for some $r_{i,j} \in \bF_{p^k}$. To prove this, note that since $\sum_i e_i = 1$, we have
\[
m = \sum_{i,j} e_ime_j,
\]
and each $e_ime_j$ defines a map $e_j\bM \rightarrow e_i\bM$. If we define the element $r_{i,j} \in \bF_{p^k}$ by
\[
r_{i,j} \coloneqq \sum_k e_{ki}me_{jk},
\]
then we have
\[
r_{i,j}e_{ij} = \sum_k e_{ki}me_{jk}e_{ij} = e_{ii}me_{ji}e_{ij} = e_ime_j,
\]
so $m = \sum_{i,j} r_{i,j}e_{ij}$. For the uniqueness, note that $\sum_{i,j} r_{i,j}e_{ij} = 0$ implies that each $r_{i,j}e_{ij} = 0$, and if $r_{i,j} \ne 0$ then $r_{i,j}$ is invertible, since $\bF_{p^k}$ is a field. Thus we have an explicit isomorphism $M_n(\bF_{p^k}) \cong \RR$. Finally, $e_1\bM$ is a vector space of some dimension $m$ over $\bF_{p^k}$, and since $\bM$ is the direct sum of $n$ copies of $e_1\bM$ as a vector space over $\bF_{p^k}$ we have $\bM \cong \bF_{p^k}^{n\times m}$, with the action of $\RR$ on $\bM$ corresponding to matrix multiplication.
\end{proof}

For the sake of completeness, we include Witt's proof of Wedderburn's little theorem here.

\begin{thm}[Wedderburn's little theorem]\label{thm-wedderburn-little} If $\RR$ is a finite division ring, then $\RR$ is a field.
\end{thm}
\begin{proof} (Following Witt \cite{witt-wedderburn-little}) Let $\RR^\times$ be the group of nonzero elements of $\RR$, and let $Z(\RR^\times)$ be the center of the group $\RR^\times$. Then $\bF = Z(\RR) = Z(\RR^\times) \cup \{0\}$ is a finite field, of some prime power order $q = p^k$. Since $\RR$ is an $\bF$-algebra, $\RR$ is in particular a vector space over $\bF$ of some dimension $n$, so $|\RR| = q^n$ and $|\RR^\times| = q^n-1$.

We consider the conjugation action of $\RR^\times$ on itself: if $x \in \RR \setminus \bF$, then the centralizer $C_\RR(x) = \{r \in \RR \mid rx = xr\}$ is a proper $\bF$-subalgebra of $\RR$, so $|C_\RR(x)| = q^k$ for some $k < n$, and since $\RR$ can be thought of as a module over the division ring $C_\RR(x)$, we have $k \mid n$. Then the conjugacy class of $x$ in $\RR^\times$ has size $\frac{q^n-1}{q^k-1}$, so we have
\[
q^n-1 = |\RR^\times| = |Z(\RR^\times)| + \sum_{\text{conj. classes of }\RR\setminus\bF} \frac{q^n-1}{q^{k_i}-1} = q-1 + \sum_{\text{conj. classes of }\RR\setminus\bF} \frac{q^n-1}{q^{k_i}-1}.
\]
If we let $\Phi_n(x)$ be the $n$th cyclotomic polynomial, then we have $\Phi_n(q) \mid \frac{q^n-1}{q^{k_i}-1}$ for each conjugacy class, so $q-1$ must be a multiple of $\Phi_n(q)$. However, $|\Phi_n(q)|$ is the product of $|q-\zeta|$ over various $n$th roots of unity $\zeta$, so $|\Phi_n(q)| > q-1$ for $n > 1$, a contradiction.
\end{proof}

\subsection{Tight lattices produce tame quotients}

The purpose of this subsection is to give a purely lattice-theoretic criterion which we can use to prove that certain congruence quotients are tame. As we will see later, nontrivial occurences of this sort of sublattice imply the existence of abelian congruence quotients.

\begin{defn} Suppose $\cL$ is a lattice with a $0$ and a $1$. A lattice homomorphism $f : \cL \rightarrow \cL'$ is \emph{$0,1$-separating} if we have
\[
f^{-1}(f(0)) = \{0\}, \;\;\; f^{-1}(f(1)) = \{1\}.
\]
A lattice $\cL$ is \emph{$0,1$-simple} if it has a $0$ and a $1$ which are not equal to each other, and if every nonconstant lattice homomorphism $\cL \rightarrow \cL'$ is $0,1$-separating.

A \emph{meet endomorphism} of a lattice $\cL$ is a function $\mu : \cL \rightarrow \cL$ which preserves $\wedge$, i.e. such that
\[
\mu(x \wedge y) = \mu(x) \wedge \mu(y).
\]
A function $\mu : \cL \rightarrow \cL$ is called \emph{increasing} if
\[
\mu(x) \ge x
\]
for all $x \in \cL$, and is called \emph{strictly increasing} if
\[
\mu(x) > x
\]
for all $x \in \cL \setminus \{1\}$.

A lattice $\cL$ is called \emph{tight} if $\cL$ is $0,1$-simple and every strictly increasing meet endomorphism of $\cL$ is constant.
\end{defn}

\begin{thm}\label{thm-tight-tame} If $\bA$ is a finite algebra and if the interval $\llbracket \alpha, \beta \rrbracket$ of $\Con(\bA)$ is tight, then the congruence quotient $(\alpha, \beta)$ is tame: for every $U \in M_\bA(\alpha, \beta)$, there is an idempotent unary polynomial $e \in E(\bA)$ such that $e(\bA) = U$.
\end{thm}
\begin{proof} (Following \cite{hobby-mckenzie}) Since $\cL$ is $0,1$-simple, the restriction homomorphism will automatically be $0,1$-separating once we show that such an $e$ exists. It's enough to show that there is some $f \in \Pol_1(\bA)$ such that $f(\bA) = U$ and $f(U) = U$, since then we can iterate $f$ to produce $e$. To find such an $f$, we just need to find a pair $f,g \in \Pol_1(\bA)$ such that $f(\bA), g(\bA) \subseteq U$ and $f(g(\beta)) \not\subseteq \alpha$.

Let $K$ be the set of unary polynomials $f \in \Pol_1(\bA)$ such that $f(\bA) \subseteq U$. One way to check whether there is some $f \in K$ with $f(\beta) \not\subseteq \alpha$ is to try to find the largest congruence $\mu$ below $\beta$ such that $f(\mu) \subseteq \alpha$ for all $f \in K$, and then to check if $\mu = \beta$. This leads to defining the following mapping on congruences:
\[
\mu(\theta) \coloneqq \{(x,y) \in \beta \mid \forall f \in K,\; (f(x), f(y)) \in \theta\}.
\]
It's easy to see that $\mu(\theta)$ is automatically a congruence, that $\theta \le \mu(\theta)$, and that
\[
\mu(\theta_1 \wedge \theta_2) = \mu(\theta_1) \wedge \mu(\theta_2).
\]
Thus $\mu$ is an increasing meet endomorphism of $\llbracket\alpha, \beta\rrbracket$.

Since $U \in M_\bA(\alpha,\beta)$, there must be some $f \in K$ such that $f(\beta) \not\subseteq \alpha$, so
\[
\mu(\alpha) < \beta.
\]
Thus $\mu$ is not constant. By the assumption that $\llbracket \alpha, \beta \rrbracket$ is tight, $\mu$ must not be \emph{strictly} increasing, so there must be some $\theta < \beta$ such that
\[
\mu(\theta) = \theta.
\]
Thus we have
\[
\mu(\mu(\alpha)) \le \mu(\mu(\theta)) = \theta < \beta.
\]
The point is that $\mu\circ \mu$ is what we would get if we replaced $K$ by $K^2$ in the definition of $\mu$, that is,
\[
\mu(\mu(\alpha)) = \{(x,y) \in \beta \mid \forall f, g \in K,\; (f(g(x)), f(g(y))) \in \alpha\},
\]
so from $\mu(\mu(\alpha)) \ne \beta$ we conclude that there must be some $f,g \in K$ such that $f(g(\beta)) \not\subseteq \alpha$, and we are done.
\end{proof}

At first it may seem that the proof only needs us to require that $\llbracket \alpha, \beta \rrbracket$ has no nonconstant increasing meet endomorphisms $\mu$ such that $\mu \circ \mu$ is constant. However, this is actually equivalent to having no nonconstant strictly increasing meet endomorphisms: if $\mu$ is a nonconstant strictly increasing meet endomorphism, then there is some minimal $k > 1$ such that $\mu^{\circ k}(\alpha) = \beta$, and then $\mu^{\circ (k-1)}$ will be nonconstant but $\mu^{\circ (k-1)} \circ \mu^{\circ (k-1)}$ will be constant.

In the remainder of this subsection, we will give alternative lattice-theoretic characterizations of what it means for a finite lattice to tight. We start by examining what it means for a lattice to be $0,1$-simple.

\begin{prop} A lattice $\cL$ is $0,1$-simple iff there is a unique dual atom $\theta \prec 1_\cL \in \Con(\cL)$ and the associated map $\cL \twoheadrightarrow \cL/\theta$ is $0,1$-separating.
\end{prop}
\begin{proof} Call a congruence $\eta$ on $\cL$ $0,1$-separating if the quotient map $\cL \twoheadrightarrow \cL/\eta$ is $0,1$-separating, that is, if $0/\eta = \{0\}$ and $1/\eta = \{1\}$. Then any join of $0,1$-separating congruences is $0,1$-separating, so there is always a unique maximal $0,1$-separating congruence $\theta \in \Con(\cL)$. Then $\cL$ is $0,1$-simple iff all congruences $\eta < 1_\cL$ satisfy $\eta \le \theta$.
\end{proof}

\begin{prop} If a complete lattice $\cL$ is $0,1$-simple and $\theta$ is a proper congruence on $\cL$, then $\cL$ is tight iff $\cL/\theta$ is tight.
\end{prop}
\begin{proof} Let $f : \cL/\theta \rightarrow \cL$ be the meet homomorphism given by
\[
f(a/\theta) = \bigvee_{b \in a/\theta} b.
\]
Then $f$ is a section of the quotient map $\pi : \cL \rightarrow \cL/\theta$, i.e. $\pi \circ f$ is the identity on $\cL/\theta$, and furthermore $f \circ \pi$ is an increasing meet endomorphism of $\cL$ which maps $0$ to $0$ and $1$ to $1$. Then for any nonconstant strictly increasing meet endomorphism $\mu$ of $\cL$, the map
\[
\pi \circ \mu \circ f
\]
is a strictly increasing meet endomorphism of $\cL/\theta$ which sends $0/\theta$ to $\mu(0)/\theta \ne 1/\theta$, and similarly for any nonconstant strictly increasing meet endomorphism $\mu'$ of $\cL/\theta$, the map
\[
f \circ \mu' \circ \pi
\]
is a strictly increasing meet endomorphism of $\cL$ which sends $0$ to $f(\mu'(0/\theta)) \ne 1$.
\end{proof}

From this we see that we only need to characterize \emph{simple} tight lattices. Recall that a \emph{tolerance} on an algebraic structure is a compatible binary relation which is symmetric and which contains the diagonal, and that a tolerance is called \emph{connected} if its transitive closure is the full congruence.

\begin{prop} A simple lattice $\cL$ of finite length is tight iff it has no nontrivial tolerances. Equivalently, a $0,1$-simple lattice of finite length is tight iff it has no proper connected tolerances.
\end{prop}
\begin{proof} If $\bS \le_{sd} \cL \times \cL$ is a tolerance, then we can define a corresponding increasing meet endomorphism $\mu_\bS$ by
\[
\mu_\bS(a) = \bigvee_{(a,b) \in \bS} b.
\]
The tolerance $\bS$ is connected iff $\mu_\bS$ is strictly increasing: the largest element $0$ can be connected to via $\bS$ in $k$ steps is $\mu_\bS^{\circ k}(0)$.

Conversely, if $\mu$ is an increasing meet endomorphism, then we can define a corresponding tolerance $\bS_\mu$ by
\[
(a,b) \in \bS_\mu \;\; \iff \;\; (a \le \mu(b)) \wedge (b \le \mu(a)).
\]
In fact, the constructions $\bS \mapsto \mu_\bS$ and $\mu \mapsto \bS_\mu$ invert each other.
\end{proof}

We say that a lattice is \emph{order polynomially complete} if every monotone operation $\cL^n \rightarrow \cL$ is a polynomial of $\cL$.

\begin{prop} If $\cL$ is a finite lattice, then the following are equivalent:
\begin{itemize}
\item[(a)] $\cL$ is simple and tight,
\item[(b)] for any $a < b \in \cL$, there is a unary polynomial $f$ of $\cL$ such that $f(a) = 0$ and $f(b) = 1$,
\item[(c)] the only compatible binary relations on $\cL$ which contain the diagonal are the diagonal $\Delta_\cL$, the partial orders $\le_\cL$ and $\ge_\cL$, and the full relation $\cL^2$,
\item[(d)] $\cL$ is order polynomially complete.
\end{itemize}
\end{prop}
\begin{proof} For (a) $\implies$ (b), note that if $\cL$ has no nontrivial tolerances, then the tolerance generated by the diagonal and $\{(a,b),(b,a)\}$ must contain $(0,1)$, so there is some binary polynomial $g$ such that
\[
g\Big(\begin{bmatrix} a\\ b \end{bmatrix}, \begin{bmatrix} b\\ a \end{bmatrix}\Big) = \begin{bmatrix} 0\\ 1 \end{bmatrix}.
\]
Since $g$ is monotone, we must also have
\[
g\Big(\begin{bmatrix} a\\ b \end{bmatrix}, \begin{bmatrix} a\\ a \end{bmatrix}\Big) = \begin{bmatrix} 0\\ 1 \end{bmatrix},
\]
so we can take $f(x) = g(x,a)$.

For (b) $\implies$ (c), we just have to prove that the binary relation $\RR$ generated by the diagonal and $\{(a,b)\}$ contains $\le_{\cL}$ as long as $a \not\ge b$. Note that $a \not\ge b$ implies $a < a \vee b$, so by (b) there is some unary polynomial $f$ such that $f(a) = 0$ and $f(a \vee b) = 1$. Since
\[
(a, a\vee b) = (a,a) \vee (a,b) \in \RR,
\]
we see that $(0,1) \in \RR$. But then for any $c \le d$, we have
\[
(c,d) = ((c,c) \vee (0,1)) \wedge (d,d) \in \RR,
\]
so $\le_\cL$ is contained in $\RR$.

For (c) $\implies$ (d), note that the collection of $n$-ary polynomials of $\cL$ is equal to the sublattice $\RR \le \cL^{\cL^n}$ which is generated by the constant tuples and the projections $\pi_i : x \mapsto x_i$. Since every lattice has a majority term, we see that $\RR$ is equal to the intersection of its binary projections, each of which contains the diagonal. Applying (c), we see that $f \in \RR$ if and only if for every $x, y \in \cL^n$ such that
\[
\pi_i(x) = x_i \le \pi_i(y) = y_i
\]
for all $i$, we have $f(x) \le f(y)$.

For (d) $\implies$ (b), we check that the map $f : \cL \rightarrow \cL$ given by
\[
f : x \mapsto \begin{cases} 0 & x \le a,\\ 1 & x \not\le a \end{cases}
\]
is monotone. Finally, (c) $\implies$ (a) follows from the previous proposition.
\end{proof}

\begin{prop} If $\cL$ is a lattice of finite length, then the smallest connected tolerance on $\cL$ is generated by the diagonal and the pairs $(x,y)$ such that either $y$ covers $x$ or $x$ covers $y$.
\end{prop}
\begin{proof} Let $\bS$ be any connected tolerance of $\cL$. Suppose that $y$ covers $x$, and consider an increasing path $x = x_0 < x_1 < \cdots < x_n = 1$ from $x$ to $1$ through $\bS$. Then there must be some first $i$ such that $x_i \ge y$, and we see that
\[
(x,y) = (y,y) \wedge (x_{i-1}, x_i) \in \bS.\qedhere
\]
\end{proof}

\begin{prop} If the join of the atoms of a $0,1$-simple lattice $\cL$ of finite length is equal to $1$, or if the meet of the co-atoms is equal to $0$, then $\cL$ is tight.
\end{prop}
\begin{proof} We will check that in either case $\cL$ has no proper connected tolerances. Suppose that $\bS$ is a connected tolerance, and let $a$ be any atom of $\cL$. In order for $a$ to be connected to $0$ via $\bS$ in any number of steps, $a$ must be connected to something strictly less than $a$ via $\bS$ in one step, so we must have
\[
(0,a) \in \bS.
\]
Since this is true for all atoms of $\cL$, joining them together we see that $(0,1) \in \bS$ if the join of the atoms is $1$.
\end{proof}

\begin{prop} The lattice $\cL_\bM$ of subspaces of a finite-dimensional vector space $\bM$ over a field $\bF$ is always tight.
\end{prop}
\begin{proof} By the previous result, we just need to check that $\cL_\bM$ is in fact simple. If $\theta$ is a nontrivial congruence on $\cL_\bM$ which identifies subspaces $u \ne v$, then by taking meets with a one-dimensional subspace which is contained in one of $u,v$ but not the other, we see that $0$ is congruent to some atom $a = \Sg_\bM\{x\}$ of $\cL_\bM$.

Now let $b = \Sg_\bM\{y\}$ be any other atom, and note that $c = \Sg_\bM\{x+y\}$ is necessarily different from both $a$ and $b$. Then $0, a, b, c$, and $a \vee b = \Sg_{\bM}\{x,y\}$ form a sublattice of $\cL_\bM$ isomorphic to the diamond lattice $\cM_3$. Since $\cM_3$ is simple, we see that $(0,a) \in \theta$ implies $(0,b) \in \theta$ - so in fact, any nontrivial congruence $\theta$ on $\cL_\bM$ must contain every atom in $0/\theta$. Since the join of the atoms is the whole space, we have $(0,1) \in \theta$, so $\theta$ was not a proper congruence on $\cL_\bM$.
\end{proof}

\begin{prop} If $A$ is a finite set, then the lattice $\cL_A$ of equivalence relations on $A$ is tight.
\end{prop}
\begin{proof} The proof is very similar to the previous proof - this time we use the fact that for any distinct $x,y,z \in A$, the equivalence relations $0_A, \Cg_A\{(x,y)\}, \Cg_A\{(y,z)\}, \Cg_A\{(x,z)\}, \Cg_A\{(x,y),(y,z)\}$ form a sublattice of $\cL_A$ which is isomorphic to $\cM_3$.
\end{proof}

We say that a $0,1$-lattice $\cL$ is \emph{complemented} if for all $x \in \cL$ there is some $x' \in \cL$ such that
\[
x \vee x' = 1, \;\;\; x \wedge x' = 0.
\]
Such an $x'$ is called a \emph{complement} of $x$ (and in general there may be more than one complement). Both types of lattices just considered (the subspaces of a finite dimensional vector space and the equivalence relations on a finite set) are complemented.

\begin{prop} If a lattice $\cL$ of finite length is complemented, then the join of the atoms of $\cL$ is $1$ and the meet of the co-atoms is $0$.
\end{prop}
\begin{proof} Let $x$ be the join of the atoms of $\cL$, and suppose that $x'$ is a complement of $x$. Then since $x \wedge x' = 0$, $x'$ is not greater than any atom of $\cL$, so $x' = 0$. Thus we have $1 = x \vee x' = x$.
\end{proof}

\begin{prop} If $\cL \twoheadrightarrow \cL'$ is $0,1$-separating, then $\cL$ is complemented iff $\cL'$ is complemented, and the atoms of $\cL$ join to $1$ iff the atoms of $\cL'$ join to $1$.
\end{prop}

The theory of tight lattices simplifies dramatically when we restrict our attention to modular lattices.

\begin{prop} If $\cL$ is a modular lattice of finite length, then the map $\mu$ given by
\[
\mu : x \mapsto x \vee \bigvee\{y \mid x \prec y\},
\]
which takes $x$ to the join of the collection of covers of $x$, is a strictly increasing meet endomorphism.
\end{prop}
\begin{proof} First we check that $\mu$ is monotone, i.e. that $x \le z$ implies $\mu(x) \le \mu(z)$. If $x \le z$ and $x \prec y$, then modularity of $\cL$ implies that either $y \le z$ or $y \vee z$ is a cover of $z$. Thus we have
\[
y \le y \vee z \le \mu(z)
\]
for all $x \prec y$, so $\mu(x) \le \mu(z)$.

Define a dual map $\sigma$ by
\[
\sigma : x \mapsto x \wedge \bigwedge\{y \mid y \prec x\}.
\]
Note that $\sigma$ is also monotone (by a dual argument to the above). Our strategy is to prove that
\begin{equation}
x \le \mu(y) \;\; \iff \;\; \sigma(x) \le y.\tag{$*$}\label{polarity}
\end{equation}
If we prove \eqref{polarity}, then we will have
\[
x \le \mu(a \wedge b) \iff \sigma(x) \le a \wedge b \iff x \le \mu(a) \wedge \mu(b),
\]
which will prove that $\mu$ is a meet endomorphism.

By the monotonicity of $\sigma$, we just need to check that we have $\sigma(x) \le y$ when $x = \mu(y)$ in order to verify the forward direction of \eqref{polarity}. Let $y_1, ..., y_k$ be a minimal collection of covers of $y$ such that
\[
x = y \vee y_1 \vee \cdots \vee y_k.
\]
For each $i$, define $x_i$ by
\[
x_i = y \vee y_1 \vee \cdots \vee y_{i-1} \vee y_{i+1} \vee \cdots \vee y_k.
\]
By modularity of $\cL$ and the choice of $k$, we have $x_i \prec x$ for all $i$. It's now easy to prove by induction that for any $I \subset [k]$, we have
\[
x \wedge \bigwedge_{i \in I} x_i = y \vee \bigvee_{j \in [k]\setminus I} y_j,
\]
so
\[
\sigma(x) \le x \wedge \bigwedge_{i \in [k]} x_i = y.\qedhere
\]
\end{proof}

\begin{prop} If $\cL$ is a modular lattice of finite length, then $\cL$ is tight iff $\cL$ is simple and complemented.
\end{prop}
\begin{proof} We've already proven that if $\cL$ is simple and complemented then $\cL$ is tight, so suppose that $\cL$ is tight. Let $\mu$ be the strictly increasing meet endomorphism from the previous result. Since $\cL$ is tight, we must have $\mu(0) = 1$, so $1$ must be a join of atoms.

By the Jordan-H\"older Theorem \ref{thm-jordan-holder} and the fact that $1$ is a join of atoms, we see that any cover $x \prec y$ of $\cL$ is projective to $0 \prec a$ for some atom $a$. Thus any nontrivial congruence $\theta$ of $\cL$ which includes $(x,y)$ also includes $(0,a)$, so in order for $\cL$ to be $0,1$-simple $\cL$ must actually be simple.

To finish, we just need to check that $\cL$ is complemented. Letting $x$ be any element of $\cL$, pick a minimal set of atoms $a_1, ..., a_k$ such that
\[
x \vee a_1 \vee \cdots \vee a_k = 1,
\]
and let $x' = a_1 \vee \cdots \vee a_k$. We claim that $x \wedge x' = 0$. Suppose for contradiction that there is some atom $a'$ with
\[
a' \le x \wedge x'.
\]
Then since
\[
a' \vee a_1 \vee \cdots \vee a_k = x'
\]
is not a cover of $x'$, modularity of $\cL$ implies that there must be some $i$ such that
\[
a' \vee a_1 \vee \cdots \vee a_{i-1} = a_1 \vee \cdots \vee a_i.
\]
But then we can leave $a_i$ out of the list of atoms and we still have
\[
x \vee a_1 \vee \cdots \vee a_{i-1} \vee a_{i+1} \vee \cdots \vee a_k = 1,
\]
contradicting the choice of $k$.
\end{proof}

\section{P\'alfy's classification of finite permutational algebras: the five types}

In the last section we proved that if $(\alpha,\beta)$ is a tame congruence quotient of a finite algebra $\bA$, then for every $(\alpha,\beta)$-trace $N$ the restriction $\bA|_N / \alpha|_N$ is permutational, i.e. every unary polynomial of $\bA|_N$ is either a constant (modulo $\alpha|_N$) or a permutation. In \cite{palfy-permutational}, P\'alfy gave a complete classification of the finite permutational algebras (up to polynomial equivalence), which was one of the key ingredients needed for tame congruence theory.

The classification splits into two very different cases: algebras of size $2$, and algebras of size $\ge 3$. Since every unary operation on a set of size $2$ is either constant or is a permutation, the classification of permutational algebras on a set of size $2$ is the same as the classification of \emph{all} algebras on a set of size $2$, up to polynomial equivalence. There turn out to be exactly $7$ of these. On the other hand, the number of polynomial clones on any set of size $\ge 3$ is uncountable \cite{uncountable-polynomial-clones} - but as we will see, the permutational algebras on a set of size $\ge 3$ are all either unary or affine algebras, so they end up being much simpler than general algebras.

We start by giving some definitions in order to rule out the least interesting case - the case of unary operations only.

\begin{defn} An operation $f$ of arity $n$ \emph{depends on} its $i$th input if there is some tuple $a_1, ..., a_n$ and some $b_i$ such that
\[
f(a_1, ..., a_n) \ne f(a_1, ..., a_{i-1}, b_i, a_{i+1}, ..., a_n).
\]
An operation $f$ is \emph{essentially unary} if it only depends on one of its inputs - equivalently, $f$ is essentially unary if it can be written as the composition of a projection $\pi_i$ and a unary operation.
\end{defn}

If $f$ does not depend on its $i$th input, then we can express $f$ in terms of the function we get by replacing its $i$th input by some other input, such as its first input. So there is no need to ever think too deeply about functions which do not depend on all their inputs. In order to gain a foothold, it is helpful to start by considering the case of a binary operation which depends on all of its inputs - for this, we will replace one of the inputs of a higher arity polynomial with some constant to get a lower arity polynomial which also depends on all its inputs.

The next result is much stronger than what we will need: all we really need is the fact that if $f$ depends on at least two of its inputs, then there is a way to plug in constants for some subset of the inputs to $f$ to get a polynomial in two variables that depends on both of its inputs.

\begin{prop}[Salomaa \cite{salomaa-essential}]\label{prop-depend-inputs} If a polynomial $f$ of arity $n$ depends on all of its inputs, then it is possible to substitute a constant for one of its inputs to get a polynomial of arity $n-1$ which also depends on all of its inputs.

In fact, if $n \ge 2$, then it is possible to find at least two different inputs to $f$ where constants can be substituted to get polynomials depending on all $n-1$ of their inputs.
\end{prop}
\begin{proof} Following \cite{hobby-mckenzie}, we write $f[a,i]$ for the polynomial we get by substituting $a$ for the $i$th input of $f$. Suppose that for some $a$ and $i,j,k$, $f[a,i]$ does not depend on the $j$th input but does depend on the $k$th input. Then for \emph{every} $b$ we see that $f[b,j]$ depends on the $k$th input, by considering the case where we plug in $a$ in the $i$th input and $b$ in the $j$th input. Additionally, since $f$ depends on all its inputs there must be some $a'$ such that $f[a',i]$ depends on its $j$th input, so there must be some $b$ such that $f[b,j]$ depends on the $i$th coordinate (consider plugging in a tuple with an $a'$ in the $i$th input such that varying the $j$th input changes the value, and note that if we change $a'$ to $a$ in the $i$th position then varying the $j$th input no longer changes the value of $f$). Thus if $f[a,i]$ does not depend on its $j$th input, then there is some $b$ such that $f[b,j]$ depends on a strictly larger subset of its inputs than $f[a,i]$ does, which proves the first claim.

For the second claim, note that for each $i \le n$ and each $j \ne i$, there is some $a$ such that $f[a,j]$ depends on the $i$th input, as long as $f$ depends on its $i$th input. Then if we choose a pair $a,j$ such that $j\ne i$, $f[a,j]$ depends on the $i$th input, and $f[a,j]$ depends on as many inputs as possible subject to the previous constraints, then the argument of the previous paragraph shows that $f[a,j]$ must depend on all of its inputs.
\end{proof}

\begin{lem}\label{lem-permutational-quasigroup} Suppose that $f \in \Pol_2(\bA)$ is a binary polynomial of a finite permutational algebra $\bA$ which depends on both of its inputs, and suppose that $|\bA| \ge 3$. Then $f$ is a quasigroup operation, that is, every unary polynomial of the form $f(a,\cdot)$ or $f(\cdot, b)$ is a permutation.
\end{lem}
\begin{proof} Suppose for the sake of contradiction that $f$ depends on both of its inputs, but that there is some $a$ such that $f(a,\cdot)$ is constant, with $f(a,y) = e$ for all $y \in \bA$. Since $f$ depends on its second coordinate, there must be some $a' \ne a$ such that $f(a',\cdot)$ is a permutation, which implies that there is some $b \in \bA$ such that $f(a',b) = e$. Then since $f(a,b) = e = f(a',b)$, we must have $f(x,b) = e$ for all $x \in \bA$ as well.

For any $a' \ne a$, if $f(a',\cdot)$ is constant, then since $f(a',b) = e$, we must have $f(a',y) = e$ for all $y \in \bA$, and then for each $y$ from $f(a',y) = e = f(a,y)$ we conclude that $f(x,y) = e$ for all $x$, so $f$ is constant. This contradicts the assumption that $f$ depends on its inputs, so for all $a' \ne a$ the unary polynomial $f(a', \cdot)$ must be a permutation.

So far we have not used the fact that $|\bA| \ge 3$, and we have not fully exploited the fact that $\bA$ is finite and permutational. For this, we iterate $f$ on its second argument: define $f^1 = f$, and for each $n$ define $f^{n+1}(x,y)$ by
\[
f^{n+1}(x,y) = f(x,f^n(x,y)),
\]
and take $f^\infty(x,y) = \lim_{n \rightarrow \infty} f^{n!}(x,y)$, so
\[
f^\infty(x,y) = f^\infty(x,f^\infty(x,y)).
\]
Then $f^\infty(a,\cdot)$ is constant, while for $a' \ne a$ we have $f^\infty(a',y) = y$ for all $y \in \bA$ since each $f(a',\cdot)$ is a permutation. But then for any distinct $a',a'' \ne a$, we have $f^\infty(a',y) = y = f^\infty(a'',y)$, so $f^\infty(\cdot,y)$ must be constant for all $y$, and in particular $f^\infty(a,y) = y$ for all $y$, which contradicts the fact that $f^\infty(a,\cdot)$ is constant.
\end{proof}

\begin{cor}\label{cor-permutational-malcev} Suppose $\bA$ is a finite permutational algebra with $|\bA| \ge 3$, and suppose that some operation of $\bA$ is not essentially unary. Then $\bA$ has a Mal'cev polynomial $p(x,y,z)$.
\end{cor}
\begin{proof} This follows from the previous lemma and Proposition \ref{prop-quasigroup-malcev}.
\end{proof}

\begin{cor}\label{cor-permutational-depend} Suppose $\bA$ is a finite permutational algebra with $|\bA| \ge 3$, and suppose $f \in \Pol_n(\bA)$ has $f(a_1, ..., a_n) = f(a_1, ..., a_{i-1}, b_i, a_{i+1}, ..., a_n)$ for some $a_1, ..., a_n$ and some $b_i \ne a_i$. Then $f$ does not depend on its $i$th coordinate.
\end{cor}
\begin{proof} We will show that for any $j \ne i$, any $a'_j$, and any $b'_i$, we have
\[
f(a_1, ..., a_{j-1}, a'_j, a_{j+1}, ..., a_n) = f(a_1, ..., a_{i-1}, b'_i, a_{i+1}, ..., a_{j-1}, a'_j, a_{j+1}, ..., a_n).
\]
For this, we define a two-variable polynomial from $f$ by substituting the $k$th input of $f$ with $a_k$ for all $k \ne i,j$, and apply the previous lemma to this two variable polynomial to see that it can't depend on its $i$th input. Applying this repeatedly, we can mutate the tuple $a_1, ..., a_n$ into any tuple $a'_1, ..., a'_n$, so $f$ does not depend on its $i$th coordinate.
\end{proof}

\begin{lem}\label{lem-permutational-zero} Suppose $\bA$ is a finite permutational algebra with $|\bA| \ge 3$, $f,g \in \Pol_n(\bA)$, and suppose that for some $0 \in \bA$ we have $f(x_1, ..., x_n) = g(x_1, ..., x_n)$ for all $x_1, ..., x_n \in \bA$ such that all but one $x_i$ is $0$. Then $f = g$.
\end{lem}
\begin{proof} If every operation of $\bA$ is essentially unary, then this is obvious. Otherwise, let $p \in \Pol_3(\bA)$ be the Mal'cev operation from Corollary \ref{cor-permutational-malcev}. Then the polynomial
\[
h(x_1, ..., x_n) = p(f(x_1, ..., x_n), g(x_1, ..., x_n), 0)
\]
is $0$ whenever $f(x_1, ..., x_n) = g(x_1, ..., x_n)$, and since any Mal'cev operation must depend on its second input, we can apply Corollary \ref{cor-permutational-depend} to $p(x,y,z)$ to see that $h(x_1, ..., x_n) = 0$ if and only if $f(x_1, ..., x_n) = g(x_1, ..., x_n)$. For any input $i$, since we have
\[
h(0, ..., 0, x_i, 0, ..., 0) = 0
\]
for all $x_i$, we can apply Corollary \ref{cor-permutational-depend} to see that $h$ does not depend on any of its inputs, so $h$ must be constantly $0$, which implies that $f = g$.
\end{proof}

\begin{thm}[P\'alfy \cite{palfy-permutational}]\label{thm-permutational-big} Suppose $\bA$ is a finite permutational algebra with $|\bA| \ge 3$, and suppose that some operation of $\bA$ is not essentially unary. Then $\bA$ is affine, and in fact $\bA$ is polynomially equivalent to a vector space over a finite field.
\end{thm}
\begin{proof} Let $p(x,y,z) \in \Pol_3(\bA)$ be the Mal'cev operation from Corollary \ref{cor-permutational-malcev}, pick an element to call $0$ in $\bA$, and define a binary polynomial $+ \in \Pol_2(\bA)$ by
\[
x + y = p(x,0,y).
\]
Then since $p$ is Mal'cev, we have $0 + x = x = x + 0$ for all $x$, so by Lemma \ref{lem-permutational-zero} we have $x+y = y+x$ for all $x,y$. Similarly, from
\[
0+(0+x) = x = (0+0)+x, \;\; 0 + (x + 0) = x = (0 + x) + 0, \;\; x + (0 + 0) = x = (x + 0) + 0,
\]
we can apply Lemma \ref{lem-permutational-zero} to conclude that $x + (y + z) = (x + y) + z$ for all $x,y,z$. Since $x+0 = x = 0+x$ for all $x$, $+$ depends on both of its arguments, so by Lemma \ref{lem-permutational-quasigroup}, we see that every element $x \in \bA$ has an inverse $-x$ such that $x + (-x) = 0$. Thus $+$ defines an abelian group structure on $\bA$ with identity element $0$.

For any $f \in \Pol_1(\bA)$, if $f(0) = c$, then we can define $r \in \Pol_1(\bA)$ by $r(x) = f(x) - c$, so that $r(0) = 0$. We will show that any such $r$ distributes over addition: since
\[
r(x+0) = r(x) = r(x) + r(0), \;\; r(0+y) = r(y) = r(0) + r(y),
\]
we can apply Lemma \ref{lem-permutational-zero} to conclude that $r(x+y) = r(x) + r(y)$ for all $x,y$. Thus every unary polynomial $f \in \Pol_1(\bA)$ can be written in the form $f(x) = r(x) + c$, where $r$ distributes over addition and takes $0$ to $0$. Letting $\bF$ be the ring of $r \in \Pol_1(\bA)$ such that $r(0) = 0$, we see that every nonzero element of $\bF$ is invertible, so $\bF$ is a finite division ring, and therefore $\bF$ is a finite field by Wedderburn's little theorem \ref{thm-wedderburn-little}.

Now suppose $f \in \Pol_n(\bA)$ is any $n$-ary polynomial. Then if we define unary polynomials $r_i$ by
\[
r_i(x_i) = f(0, ..., 0, x_i, 0, ..., 0) - f(0, ..., 0),
\]
then each $r_i$ has $r_i(0) = 0$, so $r_i \in \bF$ for all $i$. If we define $g \in \Pol_n(\bA)$ by
\[
g(x_1, ..., x_n) = r_1(x_1) + \cdots + r_n(x_n) + f(0, ..., 0),
\]
then we can apply Lemma \ref{lem-permutational-zero} to see that $f = g$, so every operation of $\bA$ is linear over the finite field $\bF$.
\end{proof}

To complete the classification, we just need to classify the polynomial clones on the two-element set $\{0,1\}$. We use $\neg$ to denote the unary negation operation on $\{0,1\}$, $\oplus$ to denote xor, and of course $\wedge, \vee$ to denote and and or.

\begin{prop} If $\bA$ has underlying set $\{0,1\}$, then $\bA$ is polynomially complete iff $\neg, \wedge \in \Pol(\bA)$. Additionally, we have $\oplus \in \Pol(\bA) \implies \neg \in \Pol(\bA)$.
\end{prop}

\begin{lem}\label{lem-permutational-monotone} If $\bA$ has underlying set $\{0,1\}$ and if there is some $f \in \Pol(\bA)$ which is not monotone, then the unary negation $\neg$ is a polynomial of $\bA$.
\end{lem}
\begin{proof} If $f$ is not monotone, then there is some tuple $a_1, ..., a_n \in \{0,1\}$ and some $i$ such that
\[
f(a_1, ..., a_{i-1}, 0, a_{i+1}, ..., a_n) > f(a_1, ..., a_{i-1}, 1, a_{i+1}, ..., a_n).
\]
Then the left hand side of the displayed inequality must be $1$ and the right hand side must be $0$, so we have
\[
\neg(x) = f(a_1, ..., a_{i-1}, x, a_{i+1}, ..., a_n).\qedhere
\]
\end{proof}

\begin{lem}\label{lem-permutational-xor} If $\bA$ has underlying set $\{0,1\}$ and $\oplus \in \Pol(\bA)$, and if there is any $f \in \Pol(\bA)$ which is not affine-linear over $\ZZ/2$, then $\wedge \in \Pol(\bA)$, so $\bA$ is polynomially complete.
\end{lem}
\begin{proof} By xoring with an affine-linear function over $\ZZ/2$, we may assume without loss of generality that $f(x_1, ..., x_n) = 0$ whenever at most one $x_i$ is nonzero. Since $f$ is not identically $0$, there must be some $a_1, ..., a_n \in \{0,1\}$ with $f(x_1, ..., x_n) = 1$, and we may suppose that $\sum_i a_i$ is minimal. By our assumption on $f$, the number of nonzero $a_i$ must be at least $2$, so there is some pair of coordinates $i\ne j$ such that $a_i = a_j = 1$. If we decrease either $a_i$ or $a_j$, then by the minimality assumption $f$ becomes $0$, so we have
\[
x \wedge y = f(a_1, ..., a_{i-1}, x, a_{i+1}, ..., a_{j-1}, y, a_{j+1}, ..., a_n).\qedhere
\]
\end{proof}

\begin{prop} The polynomial clone of $(\{0,1\}, \wedge, \vee)$ is exactly the clone of all monotone functions.
\end{prop}
\begin{proof} We prove that every monotone function $f$ of arity $n$ is in the clone generated by $\wedge, \vee$ by induction on $n$:
\[
f(x_1, ..., x_n) = f(x_1, ..., x_{n-1}, 0) \vee (f(x_1, ..., x_{n-1}, 1) \wedge x_n).\qedhere
\]
\end{proof}

\begin{lem}\label{lem-permutational-or} If $\bA$ has underlying set $\{0,1\}$ and $\vee \in \Pol(\bA)$, and if $f \in \Pol(\bA)$ is monotone but is not contained in the clone generated by $\vee$, then $\wedge \in \Pol(\bA)$, so $\Pol(\bA)$ contains the clone of all monotone functions.
\end{lem}
\begin{proof} We may suppose without loss of generailty that $f$ depends on all of its inputs. If $f$ is not contained in the clone generated by $\vee$, then in particular $f(x_1, ..., x_n) \ne x_1 \vee \cdots \vee x_n$, so since $f$ is monotone there must be some input $i$ such that $f(0, ..., 0, 1, 0, ..., 0) = 0$. Since $f$ is monotone and depends on its $i$th input, there must be some $a_1, ..., a_n \in \{0,1\}$ such that
\[
f(a_1, ..., a_{i-1}, 0, a_{i+1}, ..., a_n) = 0, \;\;\; f(a_1, ..., a_{i-1}, 1, a_{i+1}, ..., a_n) = 1.
\]
Choose the $a_1, ..., a_n$ such that $\sum_j a_j$ is minimized subject to the displayed equations above. Then there is at least one $j$ such that $a_j = 1$, by the choice of $i$, and for this $j$ we have
\[
x \wedge y = f(a_1, ..., a_{i-1}, x, a_{i+1}, ..., a_{j-1}, y, a_{j+1}, ..., a_n).\qedhere
\]
\end{proof}

Putting everything together, we have the following classification of finite permutational algebras.

\begin{thm}\label{palfy-five-types} If $\bA$ is a finite permutational algebra, then up to isomorphism and polynomial equivalence, $\bA$ is one of the following:
\begin{itemize}
\item[(1)] a unary algebra, with the set of unary operations equal to a finite permutation group,

\item[(2)] a vector space over a finite field,

\item[(3)] the boolean algebra $(\{0,1\},\vee,\wedge,\neg)$,

\item[(4)] the lattice $(\{0,1\},\vee,\wedge)$, or

\item[(5)] the semilattice $(\{0,1\},\vee)$.
\end{itemize}
\end{thm}
\begin{proof} If every polynomial of $\bA$ is essentially unary, then we are in case (1). Otherwise, by Proposition \ref{prop-depend-inputs} there is some binary polynomial $f \in \Pol_2(\bA)$ which depends on both of its inputs. If $|\bA| \ge 3$, then Theorem \ref{thm-permutational-big} shows that we are in case (2). Otherwise, we assume that the underlying set of $\bA$ is $\{0,1\}$.

If $\oplus \in \Pol_2(\bA)$, then Lemma \ref{lem-permutational-xor} shows that we are either in case (2) or case (3). If $\oplus \not\in \Pol_2(\bA)$, then we must also have $\neg \not\in \Pol_2(\bA)$, since $\oplus$ is in the clone generated by $\neg$ and any binary operation $f$ which depends on both its inputs. Then by Lemma \ref{lem-permutational-monotone} every polynomial operation of $\bA$ is monotone, and $f$ is either $\vee$ or $\wedge$. By possibly swapping $0$ and $1$, we may assume without loss of generality that $f = \vee$. Then Lemma \ref{lem-permutational-or} shows that we are either in case (4) or (5).
\end{proof}

\begin{cor} If $\bA$ is a finite permutational algebra, then $\Pol(\bA)$ is generated by the binary polynomials of $\bA$.
\end{cor}

The previous corollary is a general feature of tame congruence theory: most concrete computations in tame congruence theory depend only on the set of binary polynomials.

\begin{defn} If $(\alpha,\beta)$ is a tame congruence quotient of a finite algebra $\bA$, and if $N$ is an $(\alpha,\beta)$-trace, then we say that $N$ has
\begin{itemize}
\item \emph{unary type}, or type \textbf{1}, if $\bA|_N/\alpha|_N$ is polynomially equivalent to a unary algebra,
\item \emph{affine type}, or type \textbf{2}, if $\bA|_N/\alpha|_N$ is polynomially equivalent to a vector space over a finite field,
\item \emph{boolean type}, or type \textbf{3}, if $\bA|_N/\alpha|_N$ is polynomially equivalent to a boolean algebra,
\item \emph{lattice type}, or type \textbf{4}, if $\bA|_N/\alpha|_N$ is polynomially equivalent to a lattice,
\item \emph{semilattice type}, or type \textbf{5}, if $\bA|_N/\alpha|_N$ is polynomially equivalent to a semilattice.
\end{itemize}
We say that the tame congruence quotient $(\alpha,\beta)$ has type \textbf{i} if there is any $(\alpha,\beta)$-trace with type \textbf{i}. As we will see in the next section, all of the traces of a tame congruence quotient $(\alpha,\beta)$ have the same type as each other.
\end{defn}

The numbering of the five types can be remembered with the following visual mnemonic lattice, where the ordering corresponds to the richness of the operations in the polynomial clone.

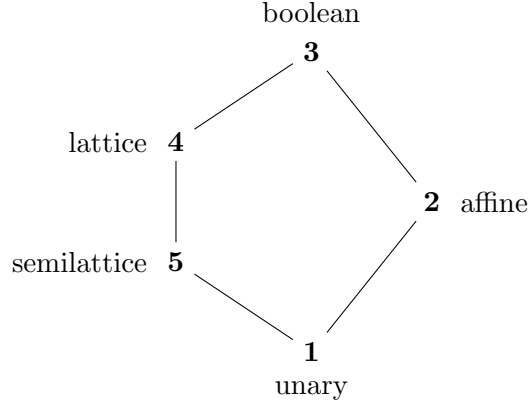
\begin{center}
\begin{tikzpicture}[scale=2]
  \node (1) [label=below:unary] at (0,-1) {\textbf{1}};
  \node (2) [label=right:affine] at (0.8,0) {\textbf{2}};
  \node (3) [label=above:boolean] at (0,1) {\textbf{3}};
  \node (4) [label=left:lattice] at (-0.9,0.4) {\textbf{4}};
  \node (5) [label=left:semilattice] at (-0.9,-0.4) {\textbf{5}};
  \draw (1) -- (2) -- (3) -- (4) -- (5) -- (1);
\end{tikzpicture}
\end{center}

\section{The structure of minimal sets}

So far we have classified the traces of tame congruence quotients, by classifying the permutational algebras. In order to classify the minimal sets of an algebra, we note that for any $(\alpha,\beta)$-minimal set $U$, the restriction $\bA|_U$ has the property that for each unary polynomial $f \in \Pol_1(\bA|_U)$, either $f$ is a permutation or $f(\beta|_U) \subseteq \alpha|_U$. Thus the restricted algebra $\bA|_U$ is $(\alpha|_U,\beta|_U)$-minimal:

\begin{defn} A finite algebra $\bA$ is called $(\alpha,\beta)$-\emph{minimal}, for $\alpha < \beta \in \Con(\bA)$, if for every unary polynomial $f \in \Pol_1(\bA)$, either $f$ is a permutation or $f(\beta) \subseteq \alpha$.
\end{defn}

By Proposition \ref{prop-tame-quotient}, an algebra $\bA$ is $(\alpha,\beta)$-minimal iff $\bA/\alpha$ is $(0_{\bA/\alpha}, \beta/\alpha)$-minimal, and for each $(\alpha,\beta)$-trace $N$ there is a corresponding $(0_{\bA/\alpha},\beta/\alpha)$-trace $N/\alpha$ of the same type, so we can often reduce to the case $\alpha = 0$ without loss of generality. We can simplify some of the arguments of \cite{hobby-mckenzie} about types \textbf{3}, \textbf{4}, and \textbf{5} by using the concept of a partial semilattice operation from Section \ref{s-partial-semi}.

\begin{defn} We say that an idempotent binary operation $s$ is a \emph{partial semilattice} if it satisfies the identity
\[
s(x,s(x,y)) \approx s(s(x,y),x) \approx s(x,y).
\]
Equivalently, $s$ is a partial semilattice if for all $x,y$, the set $\{x,s(x,y)\}$ is closed under $s$, and acts like a semilattice subalgebra with absorbing element $s(x,y)$ under $s$.

We write $a \rightarrow_s b$ if $s$ is a partial semilattice and $s(a,b) = b$.
\end{defn}

\begin{prop}\label{prop-tame-partial-semi} If $\bA$ is $(\alpha,\beta)$-minimal and has a trace $N$ of type \textbf{3}, \textbf{4}, or \textbf{5} (that is, of either boolean, lattice, or semilattice type), then $\bA$ has a partial semilattice polynomial $s \in \Pol_2(\bA)$ such that $N$ is closed under $s$, and such that $(N/\alpha, s)$ is a two-element semilattice.

Furthermore, if $N$ has type  \textbf{3} or \textbf{4} (i.e. boolean or lattice type), then there is another partial semilattice $s' \in \Pol_2(\bA)$ such that $(N, s, s')$ is a two-element lattice.

If $s$ is a partial semilattice term and $a,b \in N$ have $s(a,b) \not\equiv_\alpha a$, then $a \rightarrow_s x$ for all $x \in \bA$ and $a/\alpha = \{a\}$.
\end{prop}
\begin{proof} Let $t \in \Pol_2(\bA)$ be such that $N$ is closed under $t$ and $(N/\alpha, t)$ is a two-element semilattice. Since the unary polynomial $t(x,x)$ is not constant on $N/\alpha$, $(\alpha,\beta)$-minimality implies the unary polynomial $t(x,x)$ must be invertible, so we may assume without loss of generality that $t$ is idempotent. Then we may apply the semilattice iteration argument from Proposition \ref{prop-semilattice-iteration} to produce a partial semilattice polynomial $s \in \Clo(t)$ such that the restriction of $s$ and $t$ to $N/\alpha$ agree.

For the last statement, if $a,b \in N$ have $s(a,b) \not\equiv_\alpha s(a,a) = a$, then by $(\alpha,\beta)$-minimality the unary polynomial $x \mapsto s(a,x)$ must be a permutation, and since $s(a,s(a,x)) = s(a,x)$, it must be the identity, so $s(a,x) = x$ for all $x \in \bA$. If $a' \in a/\alpha$, then the same argument shows that $s(a',x) = x$ for all $x$, so $a \rightarrow_s a'$ and $a' \rightarrow_s a$, which is only possible if $a' = a$.
\end{proof}

Hobby and McKenzie \cite{hobby-mckenzie} like to think of their semilattices as meet-semilattices, so they call the partial semilattice polynomial $s$ from Proposition \ref{prop-tame-partial-semi} a \emph{pseudo-meet} operation. If the type is \textbf{3} or \textbf{4}, then they call the second partial semilattice operation $s'$ a \emph{pseudo-join} operation. If the type is \textbf{3} (i.e. boolean), then you can additionally find a unary polynomial $f$ which preserves $N$ and swaps the elements of $N$, and any such $f$ will be invertible. We can assume without loss of generality that this $f$ has even order, and we might call such an $f$ a \emph{pseudo-negation} operation.

\begin{prop}\label{prop-single-trace} If $\bA$ is $(\alpha,\beta)$-minimal and has at least two different $(\alpha,\beta)$-traces, then all of the $(\alpha,\beta)$-traces have type \textbf{1} or \textbf{2} (that is, they all have either unary or affine type).
\end{prop}
\begin{proof} We assume without loss of generality that $\alpha = 0_\bA$. Suppose that $N$ is a $(0_\bA,\beta)$-trace of type \textbf{3}, \textbf{4}, or \textbf{5} (that is, of either boolean, lattice, or semilattice type). Then $N$ has two elements, call them $a,b$, and by Proposition \ref{prop-tame-partial-semi} there is some partial semilattice polynomial $s \in \Pol_2(\bA)$ such that $N = \{a,b\}$ is closed under $s$ and such that $s$ acts as a semilattice operation on $\{a,b\}$, say with $s(a,b) = b$. By the second part of Proposition \ref{prop-tame-partial-semi}, we then have $s(a,x) = s(x,a) = x$ for all $x \in \bA$.

Now suppose, for the sake of a contradiction, that $K$ is a different $(0_\bA,\beta)$-trace. Since $s(a,b) = s(b,b) = b$, $(0_\bA,\beta)$-minimality implies that the unary polynomial $f : x \mapsto s(x,b)$ has $f(\beta) \subseteq 0_\bA$, so $s(K,b)$ is a singleton. Thus there must be some $c \in K$ such that $s(c,b) \ne c$. Since $s(c,a) = c \ne s(c,b)$, the unary polynomial $g : x \mapsto s(c,x)$ must be a permutation by $(0_\bA,\beta)$-minimality. However, we have $g(K) = s(c,K) \subseteq s(c,c)/\beta = c/\beta = K$ and $g(N) = s(c,N) \subseteq s(c,a)/\beta = c/\beta = K$, so $g$ can't be a permutation, which is a contradiction.
\end{proof}

\begin{prop}\label{prop-abelian-minimal-set} If $\bA$ is $(\alpha,\beta)$-minimal, then $\beta$ is abelian over $\alpha$ if and only if all of the $(\alpha,\beta)$-traces have type \textbf{1} or \textbf{2} (i.e., unary or affine type).
\end{prop}
\begin{proof} We assume without loss of generality that $\alpha = 0_\bA$. If some $(\alpha,\beta)$-trace $N$ has type \textbf{3}, \textbf{4}, or \textbf{5} (i.e., boolean, lattice, or semilattice), then there is a partial semilattice polynomial $s \in \Pol_2(\bA)$ which acts nontrivially on $N$ by Proposition \ref{prop-tame-partial-semi}. Since semilattices aren't abelian, $\bA|_N$ is not abelian, and therefore $\beta$ isn't abelian either (since $N$ is a congruence class of $\beta$).

Now suppose for contradiction that all the traces have type \textbf{1} or \textbf{2}, but that $\beta$ is not abelian. The plan is to transport the nonabelianness of $\beta$ into one of the $(0_\bA,\beta)$-traces to contradict the fact that traces of type \textbf{1} or \textbf{2} must be abelian. Recall that $\beta$ not being abelian means that there is some polynomial $g \in \Pol(\bA)$ and some $(u,v), (a_1,b_1), ..., (a_n,b_n) \in \beta$ such that
\[
g(u,a_1,...,a_n) = g(u,b_1,...,b_n)
\]
but
\[
g(v,a_1,...,a_n) \ne g(v,b_1,...,b_n),
\]
and we may assume without loss of generality that $n$ is minimal. By the minimality of $n$, we have $a_i \ne b_i$ for all $i$ (else we could just substitute $a_i$ for the $i$th argument), and we clearly have $u \ne v$, so there are $(0_\bA,\beta)$-traces $N_0, N_1, ..., N_n$ such that $u,v \in N_0$ and $a_i,b_i \in N_i$ for each $i$. Let $K$ be the $(0_\bA,\beta)$-trace which contains $g(u,a_1,...,a_n)$, then since $g$ is compatible with $\beta$ we have
\[
g(N_0,N_1, ..., N_n) \subseteq K.
\]
The restriction of $g$ to $N_0\times N_1 \times \cdots \times N_n$ must depend on all of its inputs by the minimality of $n$, so for each $i$ there are $c_j \in N_j$ such that the unary polynomial
\[
f_i : x \mapsto g(c_0,...,c_{i-1},x,c_{i+1},...,c_n)
\]
is not constant on $N_i$. Each such $f_i$ must be a permutation by $(0_\bA,\beta)$-minimality, so we have $f_i : N_i \simeq K$ for each $i$. Then the polynomial $h$ given by
\[
h(x_0,...,x_n) = g(f_0^{-1}(x_0), ..., f_n^{-1}(x_n))
\]
preserves $K$, and we have
\[
h(f_0(u), f_1(a_1), ..., f_n(a_n)) = h(f_0(u), f_1(b_1), ..., f_n(b_n))
\]
but
\[
h(f_0(v), f_1(a_1), ..., f_n(a_n)) \ne h(f_0(v), f_1(b_1), ..., f_n(b_n)).
\]
Thus $\bA|_K$ is not abelian, so $\bA|_K$ can't be polynomially equivalent to a unary or affine algebra, which is a contradiction.
\end{proof}

The next challenge is to construct a \emph{pseudo-Mal'cev} operation when the type is \textbf{2}, and to use it to prove that all of the $(\alpha,\beta)$-traces are isomorphic when at least one of them has type \textbf{2}.

\begin{lem}\label{lem-pseudo-malcev} If $\bA$ is $(\alpha,\beta)$-minimal and has an $(\alpha,\beta)$-trace $N$ of type \textbf{2}, then there is a ternary polynomial $p \in \Pol_3(\bA)$ such that, if $B$ is the union of all $(\alpha,\beta)$-traces (the ``body''), we have
\begin{itemize}
\item[(a)] $N$ is closed under $p$, and the restriction of $p(x,y,z)$ to $\bA|_N/\alpha$ is the Mal'cev operation $x - y + z$,
\item[(b)] $p$ is idempotent, that is $p(a,a,a) = a$ for all $a \in \bA$,
\item[(c)] for all $a \in \bA, b \in B$ we have $p(a,b,b) = a$, and
\item[(d)] for all $a \in \bA, b \in B$ we have $p(b,b,a) = a$.
\end{itemize}
\end{lem}
\begin{proof} (Following \cite{hobby-mckenzie}) We construct $p$ in stages, in each step getting a ternary polynomial which satisfies one more of (a), (b), (c), (d). To start, since $N$ has type \textbf{2}, there is a polynomial $f \in \Pol_3(\bA)$ satisfying (a). Next, since the restriction of the unary polynomial $g(x) = f(x,x,x)$ to $\bA|_N/\alpha$ is nonconstant, $(\alpha,\beta)$-minimality implies that $g(x)$ is a permutation, and since the restriction of $g$ to $\bA|_N/\alpha$ is the identity, the polynomial $h(x,y,z) = g^{-1}(f(x,y,z))$ satisfies (a) and (b).

{\bf Claim.} Suppose that $f$ is any polynomial satisfying (a) and (b). For any $b \in B$, the polynomials $x \mapsto f(x,b,b)$ and $x \mapsto f(b,b,x)$ are permutations.

{\bf Proof of claim.} Suppose not - suppose for contradiction that the unary polynomial $x \mapsto f(x,b,b)$ is not a permutation for some $b \in B$, and let $K$ be the $(\alpha,\beta)$-trace which contains $b$. Iterate $f$ on its first argument to get $f^\infty \in \Pol_3(\bA)$ such that
\[
f^\infty(f^\infty(x,y,z),y,z) = f^\infty(x,y,z)
\]
for all $x,y,z$, define a unary polynomial $g$ by
\[
g(x) = f^\infty(x,b,b),
\]
and note that if $x \mapsto f(x,b,b)$ is not a permutation, then $g$ is also not a permutation. By (a) and $(\alpha,\beta)$-minimality, $K$ can't be $N$: for any $a \in N$, the restriction of $x \mapsto f^\infty(x,a,a)$ to $\bA|_N/\alpha$ is the identity, so $(\alpha,\beta)$-minimality and the identity $f^\infty(f^\infty(x,a,a),a,a) = f^\infty(x,a,a)$ imply that
\[
f^\infty(x,a,a) = x
\]
for all $x \in \bA$ and $a \in N$. Since $g$ is not a permutation, $(\alpha,\beta)$-minimality implies that $g(K)$ is contained in a single $\alpha$-congruence class, so in particular there is some $c \in K$ such that $g(c) \not\equiv_\alpha c$. Then if we define the unary polynomial $h$ by $h(x) = f^\infty(c,x,x)$, we have
\[
h(c) = f^\infty(c,c,c) = c \not\equiv_\alpha g(c) = f^\infty(c,b,b) = h(b),
\]
so by $(\alpha,\beta)$-minimality $h$ is a permutation. But then
\[
h(c) = f^\infty(c,c,c) = c = f^\infty(c,a,a) = h(a)
\]
for any $a \in N$, so $h$ is not injective, a contradiction.

Now we use the claim to upgrade an $f$ satisfying (a) and (b) to one which also satisfies (c). Let $t(x,y) = f(x,y,y)$, and iterate $t$ on its first argument, to get $t^\infty \in \Pol_2(\bA)$ with $t^\infty(t^\infty(x,y),y) = t^\infty(x,y)$. By the claim, for any $b \in B$ we have $t^\infty(x,b) = x$. Now define $g \in \Pol_3(\bA)$ by
\[
g(x,y,z) = t^{\infty-1}(f(x,y,z),z).
\]
The restriction of $t$ to $\bA|_N/\alpha$ is just first projection, so $g$ satisfies (a), since $f$ is idempotent $g$ will be idempotent as well, and by construction we have $g(x,b,b) = t^\infty(x,b) = x$ for any $b \in B$.

Finally, we use the claim to upgrade an $f$ satisfying (a), (b), (c) to one which also satisfies (d). By swapping the first and third inputs to $f$, this is equivalent to upgrading an $f$ which satisfies (a), (b), (d) to one which also satisfies (c). We use the exact same construction for this as in the previous step - we just have to check that the resulting $g$ also satisfies (d): for $b \in B$, we have
\[
g(b,b,x) = t^{\infty-1}(f(b,b,x),x) = t^{\infty-1}(x,x) = x,
\]
where the last step follows from idempotence.
\end{proof}

\begin{defn} If $\bA$ is $(\alpha,\beta)$-minimal and has a trace of type \textbf{2} (i.e. affine type), and if $B$ is the union of the $(\alpha,\beta)$-traces, then we call any idempotent ternary polynomial $p \in \Pol_3(\bA)$ such that $p(a,b,b) = p(b,b,a) = a$ for all $a \in \bA$ and all $b \in B$ a \emph{pseudo-Mal'cev} polynomial for $\bA$.
\end{defn}


\begin{thm}\label{thm-pseudo-malcev} Suppose that $\bA$ is $(\alpha,\beta)$-minimal and has a trace of type \textbf{2} (i.e. affine type), let $B$ be the union of the $(\alpha,\beta)$-traces, and let $p$ be any pseudo-Mal'cev polynomial for $\bA$. Then
\begin{itemize}
\item for all $a,b \in B$, the unary polynomials $x \mapsto p(x,a,b), p(a,x,b), p(a,b,x)$ are all permutations,
\item $B$ is closed under $p$ and the restriction of $p$ to $B$ is Mal'cev, and
\item all of the $(\alpha,\beta)$-traces are polynomially isomorphic.
\end{itemize}
In particular, if one of the $(\alpha,\beta)$-traces has type \textbf{2} then they all do.
\end{thm}
\begin{proof} (Following \cite{hobby-mckenzie}) We assume without loss of generality that $\alpha = 0_\bA$. First we show that for $a,b \in B$, the unary polynomial $x \mapsto p(x,a,b)$ is a permutation iff $x \mapsto p(a,x,b)$ is. For this, let $c$ be any element of the $(0_\bA,\beta)$-trace $N = a/\beta$, and note that since $\beta$ is abelian (by Proposition \ref{prop-abelian-minimal-set}) and $(a,c) \in \beta$, we have
\[
p(c,a,b) = p(a,a,b) = b \;\;\; \iff \;\;\; p(a,c,b) = p(c,c,b) = b,
\]
so $x \mapsto p(x,a,b)$ is not a permutation iff $p(N,a,b) = \{b\}$, which happens iff $p(a,N,b) = \{b\}$, which happens iff $x \mapsto p(a,x,b)$ is not a permutation (and in fact, these all occur iff $p(N,N,b) = \{b\}$). Similarly, $x \mapsto p(a,x,b)$ is a permutation iff $x \mapsto p(a,b,x)$ is a permutation.

Now suppose for a contradiction that $x \mapsto p(a,x,b)$ and $x \mapsto p(a,b,x)$ are not permutations, and consider the unary poynomial $f(x) = p(a,p(a,x,b),x)$. If $x \in N = a/\beta$, then we have
\[
f(x) = p(a,p(a,a,b),x) = p(a,b,x) = p(a,b,a),
\]
so $f$ is not a permutation. If $x \in b/\beta$, then we have
\[
f(x) = p(a,p(a,b,b),x) = p(a,a,x) = x,
\]
so by $(0_\bA,\beta)$-minimality $f$ must be a permutation, which is a contradiction.

To see that $B$ is closed under $p$, let $a,b \in B$, then since the unary polynomial $x \mapsto p(a,b,x)$ is a permutation and therefore takes $(0_\bA,\beta)$-traces to $(0_\bA,\beta)$-traces, we see that it takes $B$ to $B$. Finally, if $N, K$ are two $(0_\bA,\beta)$-traces and $a \in N, b \in K$, then the unary polynomial $g(x) = p(a,x,b)$ takes $N$ to $K$ bijectively.
\end{proof}

Putting these results together, we have proved the main result of this section.

\begin{thm}\label{thm-same-type} If $(\alpha,\beta)$ is a tame congruence of a finite algebra $\bA$, then all of the $(\alpha,\beta)$-traces have the same type. If this type is not \textbf{1}, or if $(\alpha,\beta)$ is a prime quotient, then all of the $(\alpha,\beta)$-traces are polynomially isomorphic to each other. If the type is not \textbf{1} or \textbf{2}, then each $(\alpha,\beta)$-minimal set has just one trace, and if the type is \textbf{3} or \textbf{4} then every $(\alpha,\beta)$-trace has size two.
\end{thm}

In order to rule out type \textbf{1} in most cases, we introduce a stronger version of abelianness which is characteristic of unary algebras.

\begin{defn} If $\alpha \le \beta \in \Con(\bA)$, then $\beta$ is \emph{strongly abelian} over $\alpha$ if for all $f \in \Pol(\bA)$, all $(u,v) \in \beta$ and all tuples $x, y, z$ with $x_i \equiv y_i \equiv z_i \pmod{\beta}$, we have
\[
f(u,x_1, ..., x_n) \equiv_\alpha f(v, y_1, ..., y_n) \;\; \implies \;\; f(u, z_1, ..., z_n) \equiv_\alpha f(v, z_1, ..., z_n).
\]
An algebra $\bA$ is \emph{strongly abelian} if $1_\bA$ is strongly abelian over $0_\bA$.
\end{defn}

It's easy to see that every unary algebra is strongly abelian, while any group or semilattice is not strongly abelian. We can now characterize type \textbf{1} in terms of strong abelianness.

\begin{prop}\label{prop-strong-abelian-minimal-set} If $\bA$ is $(\alpha,\beta)$-minimal, then $\beta$ is strongly abelian over $\alpha$ iff all of the $(\alpha,\beta)$-traces have type \textbf{1} (i.e. unary type).
\end{prop}
\begin{proof} The proof is almost identical to the proof of Proposition \ref{prop-abelian-minimal-set}.
\end{proof}

\begin{prop} Suppose that $\bA$ has a Taylor polynomial and that $\alpha < \beta \in \Con(\bA)$. Then $\beta$ is not strongly abelian over $\alpha$. As a consequence, if $\bA$ is a finite Taylor algebra then no tame congruence quotient of $\bA$ has type \textbf{1}.
\end{prop}
\begin{proof} Suppose for contradiction that $\beta$ is strongly abelian over $\alpha$, and pick some $(a,b) \in \beta\setminus\alpha$. Let $t$ be a Taylor term, then for each input of $t$ we have an equation of the form
\[
t(?, ..., ?, a, ?, ..., ?) = t(?, ..., ?, b, ?, ..., ?),
\]
where each $?$ is either an $a$ or a $b$, so strong abelianness of $\beta$ over $\alpha$ implies that
\[
t(a,...,a,a,b,...,b) \equiv_\alpha t(a,...,a,b,b,...,b)
\]
at each input. Stringing these equations together and using the idempotence of $t$, we get
\[
a = t(a,...,a) \equiv_\alpha \cdots \equiv_\alpha t(b,...,b) = b,
\]
which contradicts the assumption $(a,b) \not\in \alpha$.

For the last statement, note that if $t$ is a Taylor polynomial for $\bA$ and $e \in E(\bA)$ has $e(\bA) = U$ for some $(\alpha,\beta)$-minimal set $U$, then $e(t(x_1, ..., x_n))$ is a Taylor polynomial for $\bA|_U$. In fact, by the idempotence of $t$, the restriction of $e\circ t$ to any $(\alpha,\beta)$-trace $N$ would be a Taylor polynomial for the unary algebra $\bA|_N/\alpha|_N$, which gives an even simpler contradiction.
\end{proof}

The following reformulation of strong abelianness from \cite{finite-forbidden-lattices} should give a more concrete idea of just how strong it is.

\begin{prop} An algebra $\bA$ is strongly abelian iff, for each $n$-ary polynomial $t$ of $\bA$, there are equivalence relations $R_i$ on $\bA$ such that
\[
t(a_1, ..., a_n) = t(b_1, ..., b_n) \;\; \iff \;\; \forall i \le n, \; (a_i, b_i) \in R_i.
\]
In particular, if $\bA$ is finite and strongly abelian then every polynomial of $\bA$ depends on at most $\log_2|\bA|$ of its inputs.
\end{prop}

\begin{cor} Every finite, idempotent, strongly abelian algebra can be written as a product of algebras where every operation is a projection.
\end{cor}
\begin{proof} Let $t$ be any $m$-ary operation of $\bA$, and let the equivalence relations $R_i$ be as in the previous proposition. If $t$ is idempotent, then $t$ is the graph of a bijection between $\prod_i \bA/R_i$ and $\bA$: the inverse map takes $a \in \bA$ to $(a/R_1, ..., a/R_m)$. If any $R_i$ is $0_\bA$, then $t$ must be the $i$th projection, otherwise each $\bA/R_i$ is smaller than $\bA$. To finish the proof, we just need to verify that each $R_i$ is a congruence of $\bA$. It's enough to prove this for $R_1$.

Let $s$ be any other operation of $\bA$, say of arity $n$, and consider the term
\[
t(s(y_1, ..., y_n), x_2, ..., x_m).
\]
Then by the previous proposition, there are equivalence relations $S_1, ..., S_n$ on $\bA$ such that
\[
t(s(a_1, ..., a_n), c_2, ..., c_m) = t(s(b_1, ..., b_n), d_2, ..., d_m)
\]
iff each $(a_i,b_i) \in S_i$ and each $(c_j,d_j) \in R_j$. Taking all $a_i$ to be equal to $a$ and all $b_i$ to be equal to $b$, by idempotence we see that $(a,b) \in R_1$ iff $(a,b) \in S_i$ for all $i$. In particular, we have $(a,b) \in R_1 \implies (a,b) \in S_i$, so
\[
\forall i \;\; (a_i,b_i) \in R_1 \;\; \implies \;\; \forall i \;\; (a_i,b_i) \in S_i \;\; \implies \;\; (s(a_1, ..., a_n), s(b_1, ..., b_n)) \in R_1.
\]
Since $s$ was arbitrary, $R_1$ is a congruence of $\bA$.
\end{proof}

\begin{ex} A \emph{rectangular band} is an idempotent semigroup which satisfies the identity
\[
xyx \approx x.
\]
This identity implies the apparently stronger identity
\[
xyz \approx xz,
\]
as follows:
\[
xyz \approx xy(zxz) \approx (xyzx)z \approx xz.
\]
Every rectangular band $\bA$ is strongly abelian, and is therefore isomorphic to a product of two semigroups $\bA_1, \bA_2$ such that $\cdot^{\bA_1} = \pi_1$ and $\cdot^{\bA_2} = \pi_2$. The multiplication on the rectangular band $\bA_1\times\bA_2$ is explicitly given by the rule
\[
(a,b)\cdot(c,d) = (a,d).
\]
\end{ex}

\begin{ex}[From \cite{strongly-abelian-hamiltonian}]\label{ex-strongly-abelian-not-quasiaffine} There is a $5$-element strongly abelian algebra $\bA = (\{a,b,c,d,e\},\cdot_1,\cdot_2)$ which is not quasiaffine (of course, this algebra is not idempotent). The basic binary operations of $\bA$ are given below.
\begin{center}
\begin{tabular}{cc}
\begin{tabular}{c | c c c c c} $\cdot_1$ & $a$ & $b$ & $c$ & $d$ & $e$\\ \hline $a$ & $a$ & $b$ & $a$ & $b$ & $b$\\ $b$ & $a$ & $b$ & $a$ & $b$ & $b$\\ $c$ & $c$ & $d$ & $c$ & $d$ & $d$\\ $d$ & $c$ & $d$ & $c$ & $d$ & $d$\\ $e$ & $c$ & $d$ & $c$ & $d$ & $d$ \end{tabular} &
\begin{tabular}{c | c c c c c} $\cdot_2$ & $a$ & $b$ & $c$ & $d$ & $e$\\ \hline $a$ & $a$ & $b$ & $a$ & $b$ & $b$\\ $b$ & $a$ & $b$ & $a$ & $b$ & $b$\\ $c$ & $c$ & $e$ & $c$ & $e$ & $e$\\ $d$ & $c$ & $e$ & $c$ & $e$ & $e$\\ $e$ & $c$ & $e$ & $c$ & $e$ & $e$ \end{tabular}
\end{tabular}
\end{center}
This algebra fails to be quasiaffine because it fails to satisfy the two term condition:
\[
a\cdot_1 a = a\cdot_2 a, \;\; a\cdot_1 b = a\cdot_2 b, \;\; c\cdot_1 a = c\cdot_2 a, \;\text{ but }\; c\cdot_1 b \ne c\cdot_2 b.
\]
To see that it is strongly abelian, note that for each $i,j$ we have the identities
\[
x\cdot_i (y\cdot_j z) \approx x\cdot_i z, \;\; (x\cdot_i y)\cdot_j z \approx x\cdot_j z,
\]
so every term operation of $\bA$ is one of $x, x \cdot_i x, x \cdot_i y$ for some $i \in \{1,2\}$, up to permuting its inputs.
\end{ex}

\begin{ex} A \emph{$p$-cyclic groupoid} is an idempotent binary operation $\cdot$ which satisfies the following identities:
\begin{align*}
x(yz) &\approx xy,\\
(xy)z &\approx (xz)y,\\
(...((x\underbrace{y)y)\cdots) y}_{p\text{ }y\text{s}} &\approx x.
\end{align*}
See \cite{plonka-p-cyclic} for the theory of $p$-cyclic groupoids for arbitrary primes $p$.

The free $2$-cyclic groupoid on two generators $a,b$ is isomorphic to the idempotent algebra $\bA = (\{a,b,c,d\}, \cdot)$ with basic operation $\cdot$ given below.
\begin{center}
\begin{tabular}{c | c c c c} $\cdot$ & $a$ & $b$ & $c$ & $d$\\ \hline $a$ & $a$ & $c$ & $a$ & $c$\\ $b$ & $d$ & $b$ & $d$ & $b$\\ $c$ & $c$ & $a$ & $c$ & $a$\\ $d$ & $b$ & $d$ & $b$ & $d$ \end{tabular}
\end{center}
This algebra has a congruence $\theta$ corresponding to the partition $\{a,c\}, \{b,d\}$, such that $\cdot$ on $\bA/\theta$ is first projection - in particular, $1_\bA$ is strongly abelian over $\theta$. Additionally, $\theta$ is strongly abelian over $0_\bA$, so $\bA$ is strongly solvable. The algebra $\bA$ is $(0_\bA,\theta)$-minimal, and the $(0_\bA,\theta)$-traces are $\{a,c\}$ and $\{b,d\}$. The reader may check that the $(0_\bA,\theta)$-traces $\{a,c\}$ and $\{b,d\}$ are \emph{not} polynomially isomorphic in $\bA$. The algebra $\bA|_{\{a,c\}}$ is polynomially equivalent to the unary algebra with the unary operation which swaps $a$ and $c$, corresponding to the polynomial $x \mapsto x\cdot b$.

The reader may check that $\bA$ is abelian (and even quasiaffine) as well. If we let $\eta$ be the congruence corresponding to the partition $\{a\}, \{c\}, \{b,d\}$, however, then we see that $\bA/\eta$ is \emph{not} abelian - so quotients of idempotent abelian algebras are not necessarily abelian. This is one of the senses in which type \textbf{1} can be pathological.

More generally, the free $p$-cyclic groupoid on $n$ generators is (up to isomorphism) the subalgebra of
\[
\big((\ZZ/p^2)^n,\; (x,y) \mapsto x + p(y-x)\big)
\]
generated by the basis vectors $(1,0,...,0), (0,1,...,0), ..., (0,...,0,1)$ - this algebra has $np^{n-1}$ elements. The free $p$-cyclic groupoid on $n$ generators is always quasiaffine and strongly solvable in $2$ steps via the congruence corresponding to reduction modulo $p$, and for $p,n \ge 2$ it always has a quotient which is \emph{not} abelian.
\end{ex}




\section{The abelian types: type \textbf{1} (unary) and \textbf{2} (affine)}

In case the reader has lost track, we briefly recap what we have done so far before moving on.
\begin{itemize}
\item For any $\alpha < \beta \in \Con(\bA)$, we defined $M_\bA(\alpha, \beta)$ to be the collection of minimal sets $U \subseteq \bA$ such that there is some unary polynomial $f \in \Pol_1(\bA)$ with $f(\bA) = U$ and $f(\beta) \not\subseteq \alpha$.

\item We showed that if $\beta$ is a cover of $\alpha$ in $\Con(\bA)$, then the congruence quotient $(\alpha, \beta)$ is automatically \emph{tame}, that is, for every $U \in M_\bA(\alpha, \beta)$ there is some idempotent unary polynomial $e \in E(\bA)$ with $e(\bA) = U$, and the restriction homomorphism $\llbracket \alpha, \beta \rrbracket \twoheadrightarrow \llbracket \alpha|_U, \beta|_U \rrbracket$ is a $0,1$-separating homomorphism (Proposition \ref{prop-prime-tame}).

\item We showed that if $(\alpha, \beta)$ is tame, then any two minimal sets $U, V \in M_\bA(\alpha, \beta)$ are polynomially isomorphic (Theorem \ref{thm-minimal-sets}).

\item We defined an $(\alpha, \beta)$-trace $N$ to be any congruence class of $\beta|_U$ which is not contained in a congruence class of $\alpha|_U$, for any minimal set $U \in M_\bA(\alpha, \beta)$.

\item We showed that if $(\alpha, \beta)$ is tame, then $\beta$ is the transitive closure of $\alpha \cup \{N^2 \mid N \text{ is an }(\alpha,\beta)\text{-trace}\}$ (Corollary \ref{cor-trace-closure}).

\item We showed that if $(\alpha, \beta)$ is tame, then each trace $N$ has $\bA|_N / \alpha|_N$ a permutational algebra (Corollary \ref{cor-trace-permutational}), and we classified the permutational algebras into five types (Theorem \ref{palfy-five-types}).

\item We showed that if $(\alpha, \beta)$ is tame and $U \in M_\bA(\alpha, \beta)$ is a minimal set, then all of the $(\alpha,\beta)$-traces $N \subseteq U$ have the same type, and if that type is not \textbf{1} (i.e. unary) or if $(\alpha, \beta)$ is prime, then in fact all of the $(\alpha, \beta)$-traces are polynomially isomorphic (Theorem \ref{thm-same-type}).

\item We showed that if $(\alpha, \beta)$ is tame with type \textbf{1} or \textbf{2} (i.e. unary or affine), then $\beta|_U$ is abelian over $\alpha|_U$ for any minimal set $U \in M_\bA(\alpha, \beta)$ (Proposition \ref{prop-abelian-minimal-set}), and if the type is \textbf{1} then $\beta|_U$ is strongly abelian over $\alpha|_U$ (Proposition \ref{prop-strong-abelian-minimal-set}).
\end{itemize}

We would like to have some results which don't directly reference minimal sets or traces. In this section, we will upgrade the results in the last bullet point to the claims that if $(\alpha, \beta)$ is tame with type \textbf{1} or \textbf{2}, then $\beta$ is abelian over $\alpha$, and if the type is \textbf{1} then $\beta$ is strongly abelian over $\alpha$. We start with the strongly abelian case, but the reader may prefer to read the next two results in the opposite order (or even to skip the strongly abelian case entirely, if they only care about Taylor algebras).

\begin{thm}\label{type-strong-abelian} If $\bA$ is a finite algebra and $(\alpha, \beta)$ is a tame congruence quotient, then $\beta$ is strongly abelian over $\alpha$ if and only if the type of $(\alpha, \beta)$ is \textbf{1} (i.e., unary).
\end{thm}
\begin{proof} (Following \cite{hobby-mckenzie}) We assume without loss of generality that $\alpha = 0_\bA$. If the type of $(0_\bA, \beta)$ is not \textbf{1}, then every $(0_\bA,\beta)$-trace $N$ has $\bA|_N$ not strongly abelian, so in this case $\beta$ definitely can't be strongly abelian. We just need to prove that if the type is \textbf{1} then $\beta$ is strongly abelian.

Suppose for contradiction that $\beta$ is \emph{not} strongly abelian, i.e. that there is some $f \in \Pol_n(\bA)$ and $a_i \equiv b_i \equiv c_i \pmod{\beta}$ such that
\[
f(a_1, ..., a_n) = f(b_1, ..., b_n)
\]
but
\[
f(a_1, c_2, ..., c_n) \ne f(b_1, c_2, ..., c_n).
\]
Since
\[
f(a_1, c_2, ..., c_n) \equiv f(b_1, c_2, ..., c_n) \pmod{\beta},
\]
by Theorem \ref{thm-minimal-sets}(c) or (g) there is some $e \in \Pol_1(\bA)$ and $U \in M_\bA(0_\bA, \beta)$ such that $e(\bA) = U$ and
\[
e(f(a_1, c_2, ..., c_n)) \ne e(f(b_1, c_2, ..., c_n)) \in U.
\]
Let $f'$ be the restriction of $e \circ f$ to $C = \prod_i c_i/\beta$. Then if we let $N$ be the $(0_\bA, \beta)$-trace in $U$ which contains $f'(c)$, we have
\[
f'(C) \subseteq U \cap f'(c)/\beta = N.
\]
Since
\[
f'(a_1, c_2, ..., c_n) \ne f'(b_1, c_2, ..., c_n),
\]
we see that $f'$ must depend on its first variable, and since
\[
f'(a) = e(f(a)) = e(f(b)) = f'(b),
\]
we see that $f'$ must depend on at least one other variable. By Salomaa's Proposition \ref{prop-depend-inputs}, there must be some binary polynomial $g \in \Pol_2(\bA)$ which we get by fixing some of the coordinates of $f'$ to constants which depends on both of its inputs (with the inputs restricted to the relevant $c_i/\beta$s).

In other words, we have a binary polynomial $g \in \Pol_2(\bA)$ and a pair of congruence classes $c_i/\beta, c_j/\beta$ such that
\[
g(c_i/\beta, c_j/\beta) \subseteq N,
\]
such that the restriction of $g$ to $C_{ij} = (c_i/\beta) \times (c_j/\beta)$ depends on both of its inputs. We will show that this already gives us a contradiction.

{\bf Claim.} If $N_1, N_2$ are a pair of $(0_\bA, \beta)$-traces such that $g(N_1, N_2) \subseteq N$, then the restriction of $g$ to $N_1 \times N_2$ depends on at most one of its arguments.

{\bf Proof of Claim.} Suppose not, for a contradiction. Then by plugging in constants to the first and second argument of $g$, we can apply Corollary \ref{cor-trace-iso} to see that $N_i \simeq N$. Thus we may assume without loss of generality that $N_1 = N_2 = N$. But in this case, $g$ preserves $N$, so $g|_N$ must be unary since $\bA|_N$ is a unary algebra. This contradiction proves the claim.

Now note that if $(0_\bA,\beta)$-traces $N_2, N_2'$ overlap, and if the restriction of $g$ to $N_1 \times N_2$ depends on its first input, then by the claim $g$ restricts to a nonconstant unary function of its first input on $N_1 \times N_2$, so the restriction of $g$ to $N_1 \times N_2'$ also depends on its first input, and is equal to the same nonconstant unary function of its first input.

Since every congruence class of $\beta$ is connected through $(0_\bA,\beta)$-traces by Corollary \ref{cor-trace-closure}, we see that if the restriction of $g$ to $C_{ij}$ depends on its first input, then there is some trace $N_i \subseteq c_i/\beta$ such that the restriction of $g$ to $N_1 \times (c_j/\beta)$ is a nonconstant unary function of its first input. Similarly, there is some trace $N_2 \subseteq c_j/\beta$ such that the restriction of $g$ to $(c_i/\beta) \times N_2$ is a nonconstant unary function of its second input. But then the restriction of $g$ to $N_1 \times N_2$ depends on both of its inputs, which is a contradiction.
\end{proof}

\begin{thm}\label{type-abelian} If $\bA$ is a finite algebra and $(\alpha, \beta)$ is a tame congruence quotient, then $\beta$ is abelian over $\alpha$ if and only if the type of $(\alpha, \beta)$ is \textbf{1} or \textbf{2} (i.e., unary or affine type).
\end{thm}
\begin{proof} (Following P\'alfy's argument from \cite{hobby-mckenzie}) We assume without loss of generality that $\alpha = 0_\bA$. If the type of $(0_\bA, \beta)$ is \textbf{3}, \textbf{4}, or \textbf{5}, then every $(0_\bA,\beta)$-trace $N$ has $\bA|_N$ nonabelian, so in this case $\beta$ definitely can't be abelian. We just need to prove that if the type is \textbf{1} or \textbf{2} then $\beta$ is abelian. We handled the case where the type is \textbf{1} in Theorem \ref{type-strong-abelian}, so from here on we assume that $(0_\bA, \beta)$ has type \textbf{2}.

Suppose for contradiction that $\beta$ is \emph{not} abelian, i.e. that there is some $f \in \Pol_{n+1}(\bA)$ and $a \equiv b \pmod{\beta}$, $c_i \equiv d_i \pmod{\beta}$, such that
\[
f(a, c_1, ..., c_n) = f(a, d_1, ..., d_n),
\]
but
\[
f(b, c_1, ..., c_n) \ne f(b, d_1, ..., d_n).
\]
Since every congruence class of $\beta$ is connected through $(0_\bA,\beta)$-traces by Corollary \ref{cor-trace-closure}, we may assume without loss of generality that $a,b$ are both contained in some $(0_\bA,\beta)$-trace $N$.

By Theorem \ref{thm-minimal-sets}(c) or (g), we see that there is some unary polynomial $e$ with $e(\bA) \in M_\bA(0_\bA,\beta)$ such that
\[
e(f(b, c_1, ..., c_n)) \ne e(f(b, d_1, ..., d_n)).
\]
For this choice of $e$, we see that we have
\[
e(f(b/\beta, c_1/\beta, ..., c_n/\beta)) \subseteq e(\bA) \cap e(f(b,c))/\beta = N'
\]
for some $(0_\bA,\beta)$-trace $N'$. Since the type of $(0_\bA,\beta)$ is not \textbf{1}, we can apply Theorem \ref{thm-same-type} to see that the traces $N$ and $N'$ are polynomially isomorphic. Thus we may assume without loss of generality that $N' = N$, and to simplify the notation we replace $f$ with $e \circ f$, so that we have
\[
f(b/\beta, c_1/\beta, ..., c_n/\beta) \subseteq N.
\]

The purpose of ensuring that the output of $f$ is in the same trace as the elements $a,b$ which we used in the first input is as follows. Suppose that we have traces $N_i \subseteq c_i/\beta$ for each $i$. Then by Theorem \ref{thm-same-type} there are unary polynomials $g_i \in \Pol_1(\bA)$ such that
\[
g_i : N \simeq N_i
\]
for each $i$. Then the function
\[
f(x,g_1(y_1), ..., g_n(y_n))
\]
has
\[
f(N,g_1(N), ..., g_n(N)) = f(N, N_1, ..., N_n) \subseteq N,
\]
so it preserves $N$. Since $\bA|_N$ is affine, we can fix once and for all a vector space structure on $\bA|_N$ with coefficients in some fixed finite field $\bF$. Then there are coefficients $r, r_1, ..., r_n, c \in \bF$ such that
\[
f(x,g_1(y_1), ..., g_n(y_n))|_N \approx rx + r_1y_1 + \cdots + r_ny_n + c.
\]
The coefficients $r_i$ depend on the choice of the maps $g_i : N \simeq N_i$, but the coefficient $r$ on $x$ does not - this is what we will exploit to complete the proof.

{\bf Claim.} Suppose that $N_i, N_i' \subseteq c_i/\beta$ are $(0_\bA,\beta)$-traces for each $i$, such that each $N_i$ overlaps with $N_i'$. If we choose unary polynomials $g_i, g_i' \in \Pol_1(\bA)$ with
\[
g_i : N \simeq N_i, \;\;\; g_i' : N \simeq N_i',
\]
and if we let $r, r', r_i, ..., r_i', c, c' \in \bF$ be the coefficients which satisfy
\begin{align*}
f(x,g_1(y_1), ..., g_n(y_n))|_N &\approx rx + r_1y_1 + \cdots + r_ny_n + c,\\
f(x,g_1'(y_1), ..., g_n'(y_n))|_N &\approx r'x + r_1'y_1 + \cdots + r_n'y_n + c',
\end{align*}
then $r = r'$.

{\bf Proof of Claim.} Since $N_i$ and $N_i'$ overlap, we can find $u_i, u_i' \in N$ such that
\[
g_i(u_i) = g_i'(u_i') \in N_i \cap N_i'.
\]
Pluggin in $u_i, u_i'$ for the $y_i$s, we get
\[
rx + \sum_i r_iu_i + c = f(x, g(u)) = f(x,g'(u')) = r'x + \sum_i r_i'u_i' + c'
\]
for all $x \in N$. Thus the difference $(r-r')x$ is a constant function on $N$, so $r = r'$, which proves the claim.

Since every congruence class of $\beta$ is connected through $(0_\bA,\beta)$-traces by Corollary \ref{cor-trace-closure}, we can apply the claim repeatedly to see that if we let $N_i, N_i' \subseteq c_i/\beta$ be any $(0_\bA,\beta)$-traces with $c_i \in N_i$ and $d_i \in N_i'$ (with $N_i, N_i'$ not necessarily overlapping any more), and if we define unary polynomials $g_i, g_i'$ and coefficients $r, r', r_i, ..., r_i', c, c' \in \bF$ as in the claim, then we must still have $r = r'$.

Let $u_i, u_i' \in N$ have $g_i(u_i) = c_i, g_i'(u_i') = d_i$. Then from
\[
f(a, c_1, ..., c_n) = f(a, d_1, ..., d_n)
\]
we conclude that
\[
ra + \sum_i r_iu_i + c = f(a, g(u)) = f(a, g'(u')) = ra + \sum_i r_i'u_i' + c'.
\]
Adding $r(b-a)$ to both sides, we see that
\[
f(b, g(u)) = rb + \sum_i r_iu_i + c = rb + \sum_i r_i'u_i' + c' = f(b, g'(u')),
\]
and this contradicts our assumption that
\[
f(b, c_1, ..., c_n) \ne f(b, d_1, ..., d_n),
\]
completing the proof.
\end{proof}

\begin{cor} If $\bA$ is a finite algebra and if the interval $\llbracket \alpha, \beta \rrbracket$ is a tight sublattice of $\Con(\bA)$ of size at least $3$, then $\beta$ is abelian over $\alpha$.

If additionally $\llbracket \alpha, \beta \rrbracket$ does not have a $0,1$-separating homomorphism onto the congruence lattice of a vector space, then $\beta$ is strongly abelian over $\alpha$.
\end{cor}

\section{The basic tolerance, and orderability}

The main idea of this section is to take a prime congruence quotient $(\alpha, \beta)$ in $\Con(\bA)$, and try to study the simplest binary relations $\RR$ with $\alpha \le \RR \le \beta$. Actually, we want to study this in a way that doesn't give a different answer for the pair $(\alpha, \beta)$ on $\bA$ from the answer it gives for the pair $(0_{\bA/\alpha}, \beta/\alpha)$ on $\bA/\alpha$. So we mainly focus on relations which are compatible with $\alpha$ in the following sense.

\begin{defn} If $\alpha$ is a congruence on $\bA$ and $\RR$ is a binary relation on $\bA$, then we say that $\RR$ is \emph{$\alpha$-closed} if whenever $(a,b) \in \RR$ and $a \equiv c \pmod{\alpha}, b \equiv d \pmod{\alpha}$, we also have $(c,d) \in \RR$.

We define the \emph{$\alpha$-closure} of $\RR$ to be the binary relation
\[
\alpha \circ \RR \circ \alpha.
\]
\end{defn}

\begin{prop} For $\alpha \in \Con(\bA)$ and $\RR \le \bA^2$, the $\alpha$-closure of $\RR$ is the smallest $\alpha$-closed relation on $\bA$ which contains $\RR$. There is a bijection between $\alpha$-closed binary relations on $\bA$ and binary relations on $\bA/\alpha$.
\end{prop}

So we can mainly focus on the case $\alpha = 0_\bA$. In this case, we are studying the simplest binary relations $\RR$ which contain the diagonal (and are contained in some atomic congruence $\beta$) - but relations which contain the diagonal are exactly the same as relations which are preserved by $\Pol(\bA)$, so we are really studying the simplest binary relations on the algebraic structure $(A,\Pol(\bA))$.

The simplest thing we can do is to take some pair $(a,b) \in \beta$ with $a \ne b$, and consider the binary relation generated by $\Delta_\bA \cup \{(a,b)\}$, where $\Delta_\bA = 0_\bA$ is the diagonal of $\bA$. This can be written down explicitly as
\[
\Sg_{\bA^2}(\Delta_\bA \cup \{(a,b)\}) = \{(f(a),f(b)) \mid f \in \Pol_1(\bA)\}.
\]
So we really just need to know what \emph{unary} polynomials do to the pair $(a,b)$. Now we can see how tame congruence theory will be helpful: each trace $N$ of $(0_\bA, \beta)$ contains an image of every pair $(a,b) \in \beta \setminus 0_\bA$ under some unary polynomial by Theorem \ref{thm-minimal-sets}(c) and Proposition \ref{prop-prime-poly-iso}.

We start by studying tolerances - recall that a tolerance on $\bA$ is just a symmetric reflexive relation which is compatible with the algebraic structure of $\bA$.

\begin{thm}\label{thm-basic-tolerance} If $(\alpha,\beta)$ is a prime congruence quotient of a finite algebra $\bA$ with type different from \textbf{1}, then there is a unique minimal $\alpha$-closed tolerance $\tau$ with
\[
\alpha \subsetneq \tau \subseteq \beta.
\]
This tolerance $\tau$ is the $\alpha$-closure of the relation
\[
\Sg_{\bA^2}\big(\Delta_\bA \cup \{N^2 \mid N \text{ is an $(\alpha,\beta)$-trace}\}\big).
\]
Furthermore, if the type is \textbf{2} or \textbf{3}, then $\tau$ is also minimal among reflexive $\alpha$-closed relations which properly contain $\alpha$ and are contained in $\beta$.
\end{thm}
\begin{proof} We can assume without loss of generality that $\alpha = 0_\bA$. We just need to prove that every nontrivial tolerance $\tau \subseteq \beta$ contains $N^2$ for every $(0_\bA,\beta)$-trace $N$.

Since $\tau$ is nontrivial, it must contain some $(a,b) \in \beta\setminus 0_\bA$, and by Theorem \ref{thm-minimal-sets}(c) and Proposition \ref{prop-prime-poly-iso} we can assume without loss of generality that $a,b \in N$, for any particular $(\alpha,\beta)$-trace $N$. Thus $\tau \cap N^2$ is a nontrivial tolerance on $\bA|_N$, and we just have to check that $\bA|_N$ has no nontrivial proper tolerances to finish the proof.

If the type is \textbf{3}, \textbf{4}, or \textbf{5}, then $|N| = 2$, so in this case we have $N = \{a,b\}$, and
\[
N^2 = \Delta_N \cup \{(a,b), (b,a)\} \subseteq \tau,
\]
since $(a,b) \in \tau$ and since $\tau$ is symmetric and reflexive.

If the type is \textbf{2} or \textbf{3}, then $\bA|_N$ is a Mal'cev algebra, so every reflexive relation on $\bA|_N$ is a congruence, and we see that $N^2 \subseteq \tau$ in these cases, even without the assumption that $\tau$ is symmetric.
\end{proof}

\begin{ex} If the type is equal to \textbf{1}, then there might not be a unique minimal $\alpha$-closed tolerance containing $\alpha$. Consider the unary algebra $\bA = (\ZZ/5, x \mapsto x+1 \pmod{5})$, which is simple and permutational. The minimal tolerances of $\bA$ are
\[
\tau_1 = \{(x,y) \mid x-y \in \{-1,0,1\} \pmod{5}\}
\]
and
\[
\tau_2 = \{(x,y) \mid x-y \in \{-2,0,2\} \pmod{5}\}.
\]
\end{ex}

\begin{defn} If $(\alpha,\beta)$ is a prime congruence quotient with type different from \textbf{1}, then we define the \emph{basic tolerance} of $(\alpha,\beta)$ to be the minimal $\alpha$-closed tolerance $\tau$ with $\alpha \subsetneq \tau \subseteq \beta$.
\end{defn}

If the type is \textbf{2} or \textbf{3} the situation simplifies - in these cases, the basic tolerance really is basic.

\begin{thm} If $(\alpha,\beta)$ is a prime congruence quotient of a finite algebra $\bA$ with type \textbf{2} or \textbf{3}, then the basic tolerance of $(\alpha,\beta)$ is just the $\alpha$-closure of
\[
\Delta_\bA \cup \{N^2 \mid N \text{ is an $(\alpha,\beta)$-trace}\},
\]
without needing to apply $\Sg_{\bA^2}$.
\end{thm}
\begin{proof} We can assume without loss of generality that $\alpha = 0_\bA$. By the argument of Theorem \ref{thm-basic-tolerance}, for any $(0_\bA,\beta)$-trace $N$ and for any $a \ne b \in N$, the basic tolerance $\tau$ is given by
\begin{align*}
\tau &= \Sg_{\bA^2}(\Delta_\bA \cup \{(a,b)\})\\
&= \{(f(a),f(b)) \mid f \in \Pol_1(\bA)\}.
\end{align*}
By Corollary \ref{cor-trace-iso}, for every $f \in \Pol_1(\bA)$ either $f(a) = f(b)$ or $f(N)$ is another $(0_\bA,\beta)$-trace and $f : N \simeq f(N)$. In other words, we have
\[
\tau = \Delta_\bA \cup \{f(N)^2 \mid f(N) \text{ is a $(0_\bA,\beta)$-trace}\}.\qedhere
\]
\end{proof}

If the type is \textbf{4} or \textbf{5}, then we can find smaller reflexive relations within the basic tolerance.

\begin{thm}\label{thm-basic-reflexive} If $(\alpha,\beta)$ is a prime congruence quotient of a finite algebra $\bA$ with type \textbf{4} or \textbf{5}, then there are exactly two minimal $\alpha$-closed reflexive relations $\rho_0, \rho_1$ which strictly contain $\alpha$ and are contained in $\beta$. These relations have the following properties:
\begin{itemize}
\item $\rho_1 = \rho_0^-$, that is, $\rho_1 = \{(y,x) \mid (x,y) \in \rho_0\}$,
\item $\rho_0 \cap \rho_1 = \alpha$,
\item $\rho_0 \cup \rho_1$ is the $\alpha$-closure of $\Delta_\bA \cup \{N^2 \mid N \text{ is an $(\alpha,\beta)$-trace}\}$,
\item the basic tolerance of $(\alpha,\beta)$ is the $\alpha$-closure of $\Sg_{\bA^2}(\rho_0 \cup \rho_1)$.
\end{itemize}
\end{thm}
\begin{proof} We can assume without loss of generality that $\alpha = 0_\bA$. By the argument of Theorem \ref{thm-basic-tolerance}, if $\rho$ is a nontrivial reflexive relation contained in $\beta$, then for any $(0_\bA,\beta)$-trace $N$ the restriction $\rho \cap N^2$ is a nontrivial reflexive relation on $\bA|_N$.

Since the type is \textbf{4} or \textbf{5}, $N$ has size $2$, say $N = \{a,b\}$. Then we see that $\rho$ must either contain $(a,b)$ or $(b,a)$, so the minimal $\alpha$-closed relations are
\[
\rho_0 = \Sg_{\bA^2}(\Delta_\bA \cup \{(a,b)\})
\]
and
\[
\rho_1 = \Sg_{\bA^2}(\Delta_\bA \cup \{(b,a)\}) = \rho_0^-.
\]
As in the previous argument, we have
\[
\rho_0 = \{(f(a), f(b)) \mid f \in \Pol_1\} \subseteq \Delta_\bA \cup \{f(N)^2 \mid f(N) \text{ is a $(0_\bA,\beta)$-trace}\},
\]
and
\[
\rho_0 \cup \rho_1 = \Delta_\bA \cup \{f(N)^2 \mid f(N) \text{ is a $(0_\bA,\beta)$-trace}\}.
\]
To finish the proof, we just need to check that $\rho_0 \cap \rho_1 = \Delta_\bA$, or equivalently that $(b,a) \not\in \rho_0$. To see this, note that if there was a unary polynomial $f$ such that
\[
(f(a), f(b)) = (b,a),
\]
then $f|_N$ would be a unary operation of $\bA|_N$ which swaps the elements of $N$. In this case, $\bA|_N$ would actually have type \textbf{3} (i.e. boolean type), contradicting our assumption that the type was \textbf{4} or \textbf{5} (i.e. lattice or semilattice type, respectively).
\end{proof}

The ``$\alpha$-antisymmetry'' of the relation $\rho_0$ is intriguing, and leads us to wonder if we can produce a nice quasiorder by taking the transitive closure of $\rho_0$ when the type is \textbf{4} or \textbf{5}.

\begin{defn} We say that a compatible binary relation $\zeta \le \bA^2$ is an \emph{$(\alpha,\beta)$-preorder} if
\begin{itemize}
\item $\zeta$ is a quasiorder on $\bA$,
\item $\zeta \cap \zeta^- = \alpha$, and
\item the transitive closure of $\zeta \cup \zeta^-$ is $\beta$.
\end{itemize}
Note that every compatible quasiorder $\zeta$ on $\bA$ is an $(\alpha,\beta)$-preorder for some pair of congruences $(\alpha,\beta)$, since the transitive closure of $\zeta \cup \zeta^-$ is exactly the linking congruence of $\zeta$.

We say that a congruence quotient $(\alpha,\beta)$ is \emph{orderable} if an $(\alpha,\beta)$-preorder exists.
\end{defn}

\begin{thm} If $(\alpha,\beta)$ is a tame congruence quotient with type different from \textbf{1}, then $(\alpha, \beta)$ is orderable if and only if the type of $(\alpha,\beta)$ is \textbf{4} or \textbf{5}.

In fact, if the type of $(\alpha,\beta)$ is \textbf{4} or \textbf{5}, then there are exactly two minimal $(\alpha,\beta)$-preorders $\zeta_0, \zeta_1$ and two maximal $(\alpha,\beta)$-preorders $\xi_0, \xi_1$ such that every $(\alpha,\beta)$-preorder $\eta$ satisfies
\[
\zeta_i \subseteq \eta \subseteq \xi_i
\]
for either $i = 0$ or $i = 1$.
\end{thm}
\begin{proof} (Following \cite{hobby-mckenzie}) Once again, we assume without loss of generality that $\alpha = 0_\bA$. Theorem \ref{thm-basic-tolerance} shows that if the type is \textbf{2} or \textbf{3} then every nontrivial reflexive relation $\eta$ contained in $\beta$ contains the basic tolerance, and therefore can't be a $(0_\bA, \beta)$-preorder.

Now suppose that the type is \textbf{4} or \textbf{5}. Let $\rho_0, \rho_1$ be the minimal nontrivial reflexive relations contained in $\beta$ from Theorem \ref{thm-basic-reflexive}. Note that any compatible relation $\eta$ which contains both $\rho_0$ and $\rho_1$ also contains the basic tolerance, and therefore can't be a $(0_\bA,\beta)$-quasiorder. Clearly we need to let $\zeta_i$ be the transitive closure of $\rho_i$, but the difficulty lies in verifying that $\zeta_0 \cap \rho_1 = \Delta_\bA$. To pull this off, we need to understand the maximal quasiorder $\xi_0 \subseteq \beta$ which satisfies $\xi_0 \cap \rho_1 = \Delta_\bA$.

Let $N = \{a,b\}$ be a $(0_\bA,\beta)$-trace, and suppose that $(a,b) \in \rho_0$. Then we define $\xi_0$ by
\[
\xi_0 = \{(x,y) \in \beta \mid \forall f \in \Pol_1(\bA) \text{ s.t. } f(x/\beta) \subseteq N \text{ and } f(x) = b, \text{ we also have } f(y) = b\}.
\]
Since $(x,y) \in \xi_0$ is defined in terms of an implication from $x$ to $y$, we see that $\xi_0$ is a quasiorder. By the definition of $\xi_0$ and the fact that $(0_\bA,\beta)$ is tame, we have
\[
(b,a) \not\in \xi_0,
\]
so $\xi_0 \cap \rho_1 = \Delta_\bA$, and $\xi_0$ is clearly maximal among quasiorders which are contained in $\beta$ and only meet $\rho_1$ at $\Delta_\bA$. Additionally, we have
\[
(a,b) \in \xi_0
\]
since there is no unary polynomial which swaps $a$ and $b$ if $(0_\bA,\beta)$ has type \textbf{4} or \textbf{5}, so $\zeta_0 \subseteq \xi_0$. By Proposition \ref{prop-quasiorder-poly}, to finish we just need to check that $\xi_0$ is closed under unary polynomials - but this follows directly from the definition of $\xi_0$.
\end{proof}

\begin{ex} If the type of $(\alpha,\beta)$ is \textbf{1}, then $(\alpha,\beta)$ can sometimes be orderable and sometimes not. To see this, consider the unary algebra $\bA_1 = (\{0,1\})$ with no operations, and the unary algebra $\bA_2 = (\{0,1\}, 1-x)$ with just a single operation which swaps the two elements. Then $(0_{\bA_1},1_{\bA_1})$ is orderable but $(0_{\bA_2}, 1_{\bA_2})$ is not.
\end{ex}

\section{Snags and (strong) solvability}\label{a-sec-snags}

Recall that Theorem \ref{type-abelian} says that a tame congruence quotient $(\alpha,\beta)$ of $\bA$ is abelian iff it has type \textbf{1} or \textbf{2}. Since the type of a tame congruence quotient is determined by the collection of binary polynomials $\Pol_2(\bA)$, we should be able to tell if $(\alpha,\beta)$ is abelian by examining $\Pol_2(\bA)$. We can make this more explicit by recalling that minimal sets of congruence quotients of type \textbf{3}, \textbf{4}, and \textbf{5} all have a binary pseudo-meet polynomial $s$, by Proposition \ref{prop-tame-partial-semi}. This naturally leads to the concept of a \emph{snag}.

\begin{defn} If $\bA$ is an algebra, then an ordered pair of elements $(a,b) \in \bA^2$ is a \emph{2-snag} if there is a binary polynomial $s \in \Pol_2(\bA)$ such that
\[
s(a,a) = s(a,b) = s(b,a) = a, \;\;\; s(b,b) = b.
\]
In other words, we require $(\{a,b\},s)$ to be a semilattice with $b \rightarrow_s a$. We write $\Sn_2(\bA) \subseteq \bA^2$ for the set of 2-snags of $\bA$.
\end{defn}

Similarly, Theorem \ref{type-strong-abelian} says that a tame congruence quotient $(\alpha,\beta)$ is strongly abelian iff it has type \textbf{1}. If we already know that $(\alpha,\beta)$ is abelian, then we can use the fact that the minimal sets for congruence quotients of type \textbf{2} all have a ternary pseudo-Mal'cev polynomial $p$, by Lemma \ref{lem-pseudo-malcev}. The trick to deal with this case is to pick an element $b$ in the body of one of the minimal sets (recall that the ``body'' of a minimal set is defined to be the union of the traces contained in it), and to examine the binary polynomial
\[
s(x,y) = p(x,b,y).
\]
The pseudo-Mal'cev property ensures that for any $a$ in the minimal set, we have $s(a,b) = p(a,b,b) = a$ and $s(b,a) = p(b,b,a) = a$, while $s(b,b) = p(b,b,b) = b$, so $s$ depends on both of its arguments in a way that can't occur in a strongly abelian algebra.

\begin{defn} If $\bA$ is an algebra, then an ordered pair of elements $(a,b) \in \bA^2$ is a \emph{1-snag} if there is a binary polynomial $s \in \Pol_2(\bA)$ such that
\[
s(a,b) = s(b,a) = a, \;\;\; s(b,b) = b.
\]
We write $\Sn_1(\bA)$ for the set of 1-snags of $\bA$, and note that $\Sn_2(\bA) \subseteq \Sn_1(\bA)$.
\end{defn}

\begin{thm} If $(\alpha,\beta)$ is a tame congruence quotient of a finite algebra $\bA$, then
\begin{itemize}
\item $\beta$ is abelian over $\alpha$ iff $\beta \cap \Sn_2(\bA) = \alpha \cap \Sn_2(\bA)$, and
\item $\beta$ is strongly abelian over $\alpha$ iff $\beta \cap \Sn_1(\bA) = \alpha \cap \Sn_1(\bA)$.
\end{itemize}
As a consequence, for \emph{any} $\alpha \le \beta \in \Con(\bA)$, we have
\begin{itemize}
\item $\beta$ is solvable over $\alpha$ iff $\beta \cap \Sn_2(\bA) = \alpha \cap \Sn_2(\bA)$, and
\item $\beta$ is strongly solvable over $\alpha$ iff $\beta \cap \Sn_1(\bA) = \alpha \cap \Sn_1(\bA)$.
\end{itemize}
\end{thm}

This motivates the definition of two equivalence relations $\stackrel{s}{\sim}$, $\stackrel{ss}{\sim}$ on $\Con(\bA)$.

\begin{defn} If $\bA$ is an algebra and $\alpha, \beta \in \Con(\bA)$, then we write
\[
\alpha \stackrel{s}{\sim} \beta
\]
when $\beta \cap \Sn_2(\bA) = \alpha \cap \Sn_2(\bA)$, and
\[
\alpha \stackrel{ss}{\sim} \beta
\]
when $\beta \cap \Sn_1(\bA) = \alpha \cap \Sn_1(\bA)$.
\end{defn}

\begin{thm}\label{thm-locally-solvable-congruence} If $\bA$ is a finite algebra, then each of the equivalence relations $\stackrel{s}{\sim}, \stackrel{ss}{\sim}$ defines a congruence on the lattice $\Con(\bA)$.

In particular, we have $\alpha \stackrel{s}{\sim} \beta$ iff $\alpha \vee \beta$ is solvable over $\alpha \wedge \beta$, and similarly $\alpha \stackrel{ss}{\sim} \beta$ iff $\alpha \vee \beta$ is strongly solvable over $\alpha \wedge \beta$.
\end{thm}
\begin{proof} That $\stackrel{s}{\sim}, \stackrel{ss}{\sim}$ are compatible with $\wedge$ is immediate from the definition, so we just have to prove that they are compatible with $\vee$. Note that the compatibility with $\wedge$ immediately implies that
\[
\alpha \stackrel{s}{\sim} \beta \;\;\; \iff \;\;\; \alpha \stackrel{s}{\sim} \alpha \wedge \beta \text{ and } \alpha \wedge \beta \stackrel{s}{\sim} \beta,
\]
so we just have to check that
\[
\alpha \stackrel{s}{\sim} \beta \text{ and } \gamma \stackrel{s}{\sim} \delta \;\;\; \implies \;\;\; \alpha \vee \gamma \stackrel{s}{\sim} \beta \vee \delta
\]
in the special case where $\alpha \le \beta$ and $\gamma \le \delta$ (and similarly for $\stackrel{ss}{\sim}$). In fact, we just have to check this in the special case where $\delta = \gamma$ and $(\alpha,\beta)$ is a prime congruence quotient, and we may as well assume further that $\alpha \le \gamma$. By taking $\gamma$ as large as possible among potential counterexamples, we see that we just need to prove the following claim.

{\bf Claim.} If $(\alpha, \beta)$ and $(\gamma, \eta)$ are tame congruence quotients such that
\[
\alpha \le \gamma < \eta \le \beta \vee \gamma,
\]
and if $\eta \setminus \gamma$ contains a snag, then $\beta \setminus \alpha$ contains a snag of the same type.

{\bf Proof of Claim.} Let $U \in M_\bA(\gamma, \eta)$ be a $(\gamma, \eta)$-minimal set. By Lemma \ref{lem-idempotent-surjective-lattice}, we have
\[
\eta|_U \subseteq \beta|_U \vee \gamma|_U.
\]
If $\eta \setminus \gamma$ contains a snag, then $(\gamma, \eta)$ fails to be abelian (or strongly abelian), so $\eta|_U \setminus \gamma|_U$ will also contain a snag $(a,b)$ of the same type, with $a,b$ contained in some $(\gamma, \eta)$-trace $N$. From $(a,b) \in \beta|_U \vee \gamma|_U$, we see that there must be some $(a', b') \in \beta|_U$ such that
\[
b' \in b/\gamma, \;\;\; a' \not\in b/\gamma.
\]
We clearly have $(a',b') \in \beta \setminus \gamma \subseteq \beta \setminus \alpha$, so we just need to check that $(a',b')$ is a snag of the same type as $(a,b)$.

If $(\gamma, \eta)$ has type \textbf{3}, \textbf{4}, or \textbf{5} (which must always occur if $(a,b)$ is a 2-snag), then by Proposition \ref{prop-tame-partial-semi} we see that $b' = b$, and that $\bA|_U$ has a partial semilattice polynomial $s$ such that $b \rightarrow_s x$ for all $x \in U$. In particular, $(a',b') = (a',b)$ is a 2-snag via $s$.

If $(\gamma, \eta)$ has type \textbf{2} (which may only occur if $(a,b)$ is a 1-snag), then by Lemma \ref{lem-pseudo-malcev} we see that $\bA|_U$ has a pseudo-Mal'cev operation $p$ which satisfies
\[
p(x,b',b') = p(b',b',x) = x
\]
for all $x \in U$, since $b' \in b/\gamma \cap U \subseteq N$ is contained in the body of $U$. In particular, defining the binary polynomial $s$ by
\[
s(x,y) = p(x,b',y),
\]
we see that $(a',b')$ is a 1-snag via $s$.
\end{proof}

The equivalence relations $\stackrel{s}{\sim}, \stackrel{ss}{\sim}$ still make sense on infinite algebras $\bA$, but they lose some of their meaning. We can still make use of them for \emph{locally finite} algebras - recall that an algebra $\bA$ is locally finite if every finitely generated subalgebra of $\bA$ is finite.

\begin{cor} If $\bA$ is a locally finite algebra, then each of the equivalence relations $\stackrel{s}{\sim}, \stackrel{ss}{\sim}$ defines a congruence on the lattice $\Con(\bA)$.

In particular, we have $\alpha \stackrel{s}{\sim} \beta$ iff $\alpha \vee \beta|_\bB$ is solvable over $\alpha \wedge \beta|_\bB$ for every finite subalgebra $\bB$ of $\bA$, and similarly for $\stackrel{ss}{\sim}$.
\end{cor}

It therefore makes sense to read $\stackrel{s}{\sim}$ as ``locally solvably equivalent'', and $\stackrel{ss}{\sim}$ as ``locally strongly solvably equivalent'' when studying locally finite algebras. In the infinite case, we may want to know slightly more than just the fact that $\stackrel{s}{\sim}, \stackrel{ss}{\sim}$ are congruences - we want to know if they are compatible with infinite meets and joins, for instance. Recall from Definition \ref{defn-algebraic-lattice} that a complete lattice is called \emph{algebraic} if every element can be written as a join of compact elements.

\begin{prop}\label{prop-solvable-quotient-algebraic} If $\bA$ is locally finite, then the congruences $\stackrel{s}{\sim}, \stackrel{ss}{\sim}$ are compatible with arbitrary meets and joins, and the lattices $\Con(\bA)/\!\stackrel{s}{\sim}, \Con(\bA)/\!\stackrel{ss}{\sim}$ are algebraic.
\end{prop}
\begin{proof} The only tricky claim to check is that the quotient lattices are algebraic. For this, we use the fact that each $\stackrel{s}{\sim}$-class $\alpha/\!\stackrel{s}{\sim}$ is determined by the intersection $\alpha \cap \Sn_2(\bA)$, and for any 2-snag $(a,b) \in \Sn_2(\bA)$, we can prove that $\Cg_{\bA}\{(a,b)\}/\!\stackrel{s}{\sim}$ is a compact element of $\Con(\bA)/\!\stackrel{s}{\sim}$. The argument for $\stackrel{ss}{\sim}$ is similar.
\end{proof}

We would like to claim that as long as we avoid type \textbf{1}, locally solvable algebras behave like Mal'cev algebras - i.e., that they are congruence modular. We can actually prove a much stronger claim about copies of the pentagon lattice $\cN_5$ in $\Con(\bA)$.

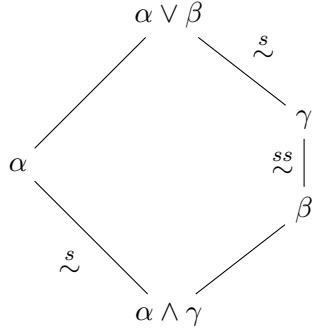
\begin{center}
\begin{tikzpicture}[scale=2]
  \node (ab) at (0,1) {$\alpha\vee\beta$};
  \node (a) at (-1,0) {$\alpha$};
  \node (ag) at (0,-1) {$\alpha\wedge\gamma$};
  \node (b) at (0.9,-0.3) {$\beta$};
  \node (g) at (0.9,0.3) {$\gamma$};
  \draw (g) to ["$\stackrel{s}{\sim}$"'] (ab) -- (a) to ["$\stackrel{s}{\sim}$"'] (ag) -- (b) to ["$\stackrel{ss}{\sim}$"] (g);
\end{tikzpicture}
\end{center}

\begin{thm} If $\bA$ is locally finite and $\alpha, \beta, \gamma \in \Con(\bA)$ form a copy of the pentagon lattice $\cN_5$ with
\[
\alpha \wedge \gamma \le \beta \le \gamma \le \alpha \vee \beta,
\]
then
\[
\alpha \stackrel{s}{\sim} \alpha \wedge \gamma \;\;\; \implies \;\;\; \beta \stackrel{ss}{\sim} \gamma.
\]
\end{thm}
\begin{proof} (Following \cite{hobby-mckenzie}) We may assume without loss of generality that $\bA$ is finite and that $(\beta, \gamma)$ is a prime congruence quotient. Since $\stackrel{s}{\sim}$ is a congruence on $\Con(\bA)$, we know that $(\beta,\gamma)$ must be abelian, so assume for the sake of contradiction that $(\beta, \gamma)$ has type \textbf{2}.

Let $U$ be a $(\beta, \gamma)$-minimal set, let $B$ be the body of $U$ (i.e., $B$ is the union of the $(\beta,\gamma)$-traces), and $T = U \setminus B$ the ``tail'' of $U$. Let $p$ be a pseudo-Mal'cev operation for $\bA|_U$. Then $(B, p)$ is a Mal'cev algebra by Theorem \ref{thm-pseudo-malcev}, so
\[
\gamma|_B \not\subseteq \alpha|_B \vee \beta|_B,
\]
since otherwise we would have a copy of the pentagon lattice $\cN_5$ in the congruence lattice of $(B, p)$, contradicting Proposition \ref{prop-malcev-modular}. Since $\gamma|_U \le \alpha|_U \wedge \beta|_U$ by Lemma \ref{lem-idempotent-surjective-lattice}, we see that we must have
\[
\alpha \cap (B\times T) \ne \emptyset.
\]
We will use this to show that $\alpha$ can't possibly be solvable over $\alpha \wedge \gamma$, which will give us our desired contradiction. We just need to prove the following claim.

{\bf Claim.} If $\delta < \theta$ is a pair of congruences such that
\[
\delta \cap (B \times T) = \emptyset \;\; \text{ and } \;\; \theta \cap (B \times T) \ne \emptyset,
\]
then $\theta$ is not abelian over $\delta$.

{\bf Proof of Claim.} Pick $(b,t) \in \theta \cap (B \times T)$, and assume for contradiction that $\theta$ is abelian over $\delta$. Since $p$ is pseudo-Mal'cev, we have
\[
p(b,b,t) = p(t,t,t) = t,
\]
so abelianness of $\theta$ over $\delta$ implies that
\[
b = p(b,b,b) \equiv p(t,t,b) \pmod{\delta},
\]
so
\[
p(t,t,b) \in B
\]
by the assumption $\delta \cap (B \times T) = \emptyset$. Let $a$ be any element in the $(\beta,\gamma)$-trace $U \cap b/\gamma$ which is not in $b/\beta$. Then we have
\[
p(t,a,b) \equiv p(t,b,b) = t \pmod{\gamma},
\]
and since $t$ is in the tail of $U$, we have $U \cap t/\gamma = U \cap t/\beta$, so
\[
p(t,a,b) \equiv t \pmod{\beta}.
\]
Then if we define the unary polynomial $f$ by
\[
f(x) = p(x,p(t,p(t,x,b),b),b),
\]
we have
\begin{align*}
f(b) &= p(b, p(t,t,b), b),\\
f(a) &\equiv p(a, p(t,t,b), b) \pmod{\beta},\\
f(t) &\equiv p(t, p(t, b, b), b) = p(t, t, b) \equiv b \pmod{\delta}.
\end{align*}
Since $p(t,t,b) \in B$, Theorem \ref{thm-pseudo-malcev} implies that $f(a) \not\equiv f(b) \pmod{\beta}$, so $f|_U$ is a permutation of $U$ by the $(\beta,\gamma)$-minimality of $U$. But then we must have $f(T) = T$, so
\[
(b, f(t)) \in \delta \cap (B \times T),
\]
which is a contradiction.
\end{proof}

\begin{cor} If $\bA$ is locally finite, then every equivalence class of $\stackrel{s}{\sim}\!/\!\stackrel{ss}{\sim}$ is a modular sublattice of $\Con(\bA)/\!\stackrel{ss}{\sim}$.
\end{cor}

\begin{defn} If $\cV$ is a locally finite variety, then we say that $\cV$ \emph{omits} type \textbf{i} if for every finite $\bA \in \cV$ and every tame congruence quotient $(\alpha,\beta)$ of $\bA$, the type of $(\alpha,\beta)$ is not \textbf{i}. We write $\typ(\cV)$ for the set of types which $\cV$ does not omit.
\end{defn}

\begin{prop} If $\cV$ is a locally finite variety, then the locally solvable algebras in $\cV$ form a subvariety $\cV_s$, and similarly the locally strongly solvable algebras in $\cV$ form a subvariety $\cV_{ss}$.
\end{prop}

\begin{cor} If $\cV$ is a locally finite variety which omits type \textbf{1}, then the subvariety $\cV_s$ of locally solvable algebras in $\cV$ has a Mal'cev term.
\end{cor}
\begin{proof} We've already shown that in this case $\cV_s$ is congruence modular, so by Corollary \ref{cor-modular-solvable-malcev} applied to the free algebra on two generators in $\cV_s$ (which is finite, and therefore solvable), the variety $\cV_s$ has a Mal'cev term.
\end{proof}

If we are only studying idempotent algebras, then this result is satisfying - but for general algebras, we need to be able to ``restrict to a congruence class'' if we want to get the most use out of this. Recall the partial order $\preceq_|$ from Definition \ref{defn-restrict-order}.

\begin{prop} If $\bA$ is finite and $\bB \preceq_| \bA$ is such that every constant of $\bB$ is a term operation of $\bB$, then $\typ(\cV(\bB)) \subseteq \typ(\cV(\bA))$.
\end{prop}
\begin{proof} If there is some finite $\bC \in \cV(\bB)$ which has a congruence quotient $(\alpha, \beta)$ of type \textbf{i}, then we pick an $(\alpha,\beta)$-trace $N$ and note that $\bC|_N / \alpha|_N$ is a permutational algebra (of type \textbf{i}) by Corollary \ref{cor-trace-permutational}. Then by Propositions \ref{prop-restrict-transitive} and \ref{prop-restrict-variety}, we see that there is some finite $\bA' \in \cV(\bA)$ such that $\bC|_N/\alpha|_N \preceq_| \bA'$.

Pick $e \in E(\bA'), \theta \in \Con(\bA')$, and $a \in e(\bA')$ such that the permutational algebra $\bC|_N/\alpha|_N$ is polynomially equivalent to $\bA'|_{N'}$, where $N'$ is given by
\[
N' = e(\bA') \cap (a / \theta).
\]
Pick $\eta \le \theta$ maximal such that $\eta|_{N'} \ne \theta|_{N'}$, and let $\theta' \le \theta$ be a cover of $\eta$. Then $(\eta, \theta')$ is tame by Proposition \ref{prop-prime-tame}, so to finish we just need to check that $N'$ is an $(\eta,\theta')$-trace.

Note that $e(a/\theta') \not\subseteq a/\eta$, so by Corollary \ref{cor-trace-closure} there must be some $(\eta,\theta')$-trace $N'' \subseteq a/\theta'$ and some $b,c \in N''$ such that $e(b)/\eta \ne e(c)/\eta$. Then by Corollary \ref{cor-trace-iso} $e(N'')$ is also an $(\eta,\theta')$-trace, and we have
\[
e(N'') \subseteq e(\bA') \cap (e(a)/\theta') = e(\bA') \cap (a/\theta) = N'.
\]
Since $e(N'')$ is an $(\eta,\theta')$-trace contained in $a/\theta'$, there is some $e' \in e(\bA')$ such that
\[
e(N'') = e'(\bA') \cap (a/\theta'),
\]
and we may assume without loss of generality that $e' = e \circ e'$. But then $e'|_{N'}$ is a polynomial of $\bA'|_{N'}$, which is permutational, so in fact we have $e(N'') = N'$, so $N'$ is an $(\eta,\theta')$-trace.
\end{proof}

Putting this together with the previous results, we can prove the existence of a Mal'cev-like term which behaves nicely on every locally solvable congruence.

\begin{thm}\label{thm-loc-solvable-malcev} If $\cV$ is a locally finite variety which omits type \textbf{1}, then $\cV$ has an idempotent ternary term $p$ such that for any $\bA \in \cV$ and any $a,b \in \bA$,
\[
\Cg_{\bA}\{(a,b)\} \stackrel{s}{\sim} 0_\bA \;\;\; \implies \;\;\; p(a,b,b) = p(b,b,a) = a.
\]
\end{thm}
\begin{proof} Let $\bF = \cF_{\cV}(x,y)$ be the (finite) free algebra on two generators in $\cV$. Define $\beta \in \Con(\bF)$ to be the congruence $\Cg_{\bF}\{(x,y)\}$, that is, the least congruence which identifies $x$ with $y$, so that $\bF/\beta \cong \cF_{\cV}(x)$. Then $x/\beta$ consists of all binary terms $t(x,y)$ of $\cV$ which satisfy $t(x,x) \approx x$, that is, $x/\beta$ corresponds exactly to the set of idempotent binary terms of $\cV$. Additionally, let $\alpha \in \Con(\bF)$ be minimal such that $\alpha \stackrel{s}{\sim} \beta$.

Taking $N = x/\beta$, we have $\bF|_N \preceq_| \bF$, so the variety generated by $\bF|_N$ omits type \textbf{1}. Additionally, $\bF|_N/\alpha|_N$ is solvable, so
\[
\typ(\cV(\bF|_N/\alpha|_N)) = \{\textbf{2}\}.
\]
In particular, we see that $\bF|_N/\alpha|_N$ has a ternary Mal'cev term $p_0$. By the definition of $\bF|_N$, $p_0$ is the restriction to $N$ of some polynomial $p_1$ of $\bF$ which preserves $N$. Since $\bF$ is generated by $x$ and $y$, we see that there is some $5$-ary term $t$ of $\cV$ such that
\[
p_1(u,v,w) = t(u,v,w,x,y)
\]
for all $u,v,w \in \bF$. Since $p_1$ preserves $N = x/\beta$, we have
\[
t(x,x,x,x,y) = p_1(x,x,x) \in x/\beta,
\]
so $t$ is idempotent. Define an idempotent ternary term $p$ of $\cV$ by
\[
p(u,v,w) = t(u,v,w,u,w).
\]
Then we have
\[
p(x,x,y) = t(x,x,y,x,y) = p_1(x,x,y) \equiv_\alpha y
\]
and
\[
p(x,y,y) = t(x,y,y,x,y) = p_1(x,y,y) \equiv_\alpha x.
\]

Now suppose that $\bA \in \cV$ and $a,b \in \bA$ have $\Cg_{\bA}\{(a,b)\} \stackrel{s}{\sim} 0_\bA$. Let $\pi : \bF \rightarrow \bA$ be the unique map with $\pi(x) = a, \pi(y) = b$. Then
\[
\pi^{-1}(\Cg_{\bA}\{(a,b)\}) \supseteq \beta,
\]
so we have $\beta \stackrel{s}{\sim} \ker \pi$, which implies that $\alpha \le \ker \pi$ by our choice of $\alpha$. In particular, we have
\[
p(x,x,y) \equiv_{\ker \pi} y, \;\;\; p(x,y,y) \equiv_{\ker \pi} x,
\]
so $p(a,a,b) = b$ and $p(a,b,b) = a$. Interchanging $a$ and $b$ in the argument gives $p(b,b,a) = a$ as well, so we are done.
\end{proof}

\begin{defn} An idempotent ternary term $p$ is called a \emph{weak difference term} for $\cV$ if for any $\bA \in \cV$, any $a,b \in \bA$, and any $\theta \in \Con(\bA)$ with $(a,b) \in \theta$, we have
\[
p(a,b,b) \equiv_{[\theta,\theta]} p(b,b,a) \equiv_{[\theta,\theta]} a.
\]
\end{defn}

\begin{cor} A locally finite variety omits type \textbf{1} iff it has a weak difference term. In particular, every locally finite Taylor variety has a weak difference term.
\end{cor}
\begin{proof} We always have $\theta \stackrel{s}{\sim} [\theta,\theta]$, so any term $p$ as in Theorem \ref{thm-loc-solvable-malcev} is automatically a weak difference term. Conversely, we need to show that if $\cV$ has a weak difference term $p$, then $\cV$ omits type \textbf{1}.

Suppose for contradiction that $\bA \in \cV$ has a tame congruence quotient $(\alpha,\beta)$ of type \textbf{1}. We may assume without loss of generality that $\alpha = 0_\bA$, in which case Theorem \ref{type-abelian} implies that $[\beta,\beta] = 0_\bA$. Letting $U = e(\bA)$ (with $e \in E(\bA)$) be a $(0_\bA, \beta)$-minimal set, we see that $e \circ p$ restricts to a Mal'cev operation on any $(0_\bA,\beta)$-trace $N$, contradicting the assumption that $\bA|_N$ is a unary algebra.
\end{proof}

\begin{rem} Theorem 4.8 of \cite{kearnes-taylor-affine} shows that any variety (not necessarily locally finite) with a weak difference term also has a Taylor term. In general, having a weak difference term is slightly stronger than having a Taylor term (an example of an algebra which has a Taylor term but no weak difference term is the algebra $(\RR, \frac{x+y}{2})$, which is abelian but has no Mal'cev term).
\end{rem}

\begin{cor} If $\cV$ is a locally finite variety which omits type \textbf{1}, $\bA \in \cV$, and $\alpha, \beta \in \Con(\bA)$, then
\[
\alpha \stackrel{s}{\sim} \beta \;\;\; \implies \;\;\; \alpha \vee \beta = \alpha \circ \beta = \beta \circ \alpha.
\]
\end{cor}
\begin{proof} By symmetry, we just need to check that $\alpha \circ \beta \subseteq \beta \circ \alpha$. Pick $p$ as in Theorem \ref{thm-loc-solvable-malcev}, and suppose that $x,y,z \in \bA$ satisfy
\[
x\ \alpha\ y\ \beta\ z.
\]
Then we have
\[
x\ \beta\ p(x,y,z)\ \alpha\ z,
\]
where
\[
p(x,y,z)/\beta = p(x,y,y)/\beta = x/\beta
\]
follows from $(\alpha \vee \beta)/\beta \stackrel{s}{\sim} 0_{\bA/\beta}$, and
\[
p(x,y,z)/\alpha = p(x,x,z)/\alpha = z/\alpha
\]
follows from $(\alpha \vee \beta)/\alpha \stackrel{s}{\sim} 0_{\bA/\alpha}$.
\end{proof}

\begin{cor} If $\cV$ is a locally finite variety which omits type \textbf{1}, $\bA \in \cV$, and $\alpha, \beta \in \Con(\bA)$, then
\[
\beta \stackrel{s}{\sim} \alpha \wedge \beta \;\;\; \implies \;\;\; \alpha \vee \beta = \alpha \circ \beta \circ \alpha.
\]
\end{cor}
\begin{proof} Note that the assumption is equivalent to $\alpha \vee \beta \stackrel{s}{\sim} \alpha$, or equivalently $(\alpha \vee \beta)/\alpha \stackrel{s}{\sim} 0_{\bA/\alpha}$. We just need to check that $\beta \circ \alpha \circ \beta \subseteq \alpha \circ \beta \circ \alpha$. Pick $p$ as in Theorem \ref{thm-loc-solvable-malcev}, and suppose that $w,x,y,z \in \bA$ satisfy
\[
w\ \beta\ x\ \alpha\ y\ \beta\ z.
\]
Then we have
\[
w\ \alpha\ p(w,y,y)\ \beta\ p(x,y,z)\ \alpha\ z,
\]
with the $\alpha$ congruences following as in the previous corollary, while the $\beta$ congruence follows directly from $w \equiv_\beta x$ and $y \equiv_\beta z$.
\end{proof}

\section{Pseudocomplements and semidistributivity}

In this section we start investigating the consequences of avoiding the abelian types on the congruence lattices of finite algebras. We start with some lattice-theoretic preliminaries.

\begin{defn} If $\cL$ is a lattice and $\alpha \le \beta \in \cL$, then we say that $\delta$ is the \emph{weak pseudocomplement of $\beta$ over $\alpha$} if $\delta$ is the greatest element of $\cL$ such that $\beta \wedge \delta = \alpha$, that is, if
\[
\beta \wedge \gamma = \alpha \;\;\; \iff \;\;\; \gamma \in \llbracket \alpha, \delta \rrbracket.
\]
A closely related concept is the relative pseudocomplement: for any $\alpha, \beta \in \cL$, $\delta$ is called the \emph{relative pseudocomplement} of $\beta$ with respect to $\alpha$ if
\[
\beta \wedge \gamma \le \alpha \;\;\; \iff \;\;\; \gamma \le \delta,
\]
and this is written in symbols as $\delta = \beta \to \alpha$ or $\delta = \beta \thinsupset \alpha$. If $\alpha = 0$, then a weak or relative pseudocomplement $\delta$ of $\beta$ over $0$ is just called a \emph{pseudocomplement} of $\beta$, and written in symbols as $\delta = \neg \beta$ or $\delta = \beta^*$.

Similarly, for $\alpha \le \beta \in \cL$ we say that $\delta$ is the \emph{dual weak pseudocomplement of $\alpha$ under $\beta$} if
\[
\alpha \vee \gamma = \beta \;\;\; \iff \;\;\; \gamma \in \llbracket \delta, \beta \rrbracket.
\]
Additionally, for any $\alpha, \beta \in \cL$, $\delta$ is called the \emph{dual relative pseudocomplement} of $\alpha$ with respect to $\beta$ if
\[
\alpha \vee \gamma \ge \beta \;\;\; \iff \;\;\; \gamma \ge \delta,
\]
and some authors write this in symbols as $\delta = \beta - \alpha$ or $\delta = \beta \backslash \alpha$.
\end{defn}

Of course, pseudocomplements don't always exist - for instance, the diamond lattice $\cM_3$ is not pseudocomplemented. Note that there can be a weak pseudocomplement of $\beta$ over $\alpha$ even if there is no relative pseudocomplement of $\beta$ with respect to $\alpha$ - this situation occurs in the pentagon lattice $\cN_5$. We at least have the following implication between the two concepts.

\begin{prop} If $\alpha \le \beta \in \cL$ and the relative pseudocomplement of $\beta$ with respect to $\alpha$ exists and is equal to $\delta$, then $\delta$ is also the weak pseudocomplement of $\beta$ over $\alpha$.
\end{prop}

To put these concepts in context, we recall the definition of a Heyting algebra, from intuitionistic logic.

\begin{defn} A \emph{Heyting algebra} is an algebraic structure $\cH = (H, \wedge, \vee, \thinsupset, 0, 1)$ such that $(H, \wedge, \vee, 0, 1)$ is a $0,1$-lattice and for every pair of elements $\alpha, \beta \in H$, $\beta \thinsupset \alpha$ is the relative pseudocomplement of $\beta$ with respect to $\alpha$.
\end{defn}

\begin{prop} A complete lattice is the lattice reduct of a Heyting algebra iff it satisfies the infinite distributive law
\begin{equation}
\alpha \wedge \Big(\bigvee_{\beta \in S} \beta\Big) = \bigvee_{\beta \in S} (\alpha \wedge \beta).\label{heyting-distributive}\tag{D$_\infty$($\wedge$)}
\end{equation}
\end{prop}
\begin{proof} First we check that any complete lattice $\cL$ which satisfies the infinite distributive law \eqref{heyting-distributive} can be expanded to a Heyting algebra. For $\alpha, \beta \in \cL$, the least possible value for the relative pseudocomplement $\beta \thinsupset \alpha$ is given by
\[
\beta \thinsupset \alpha = \bigvee_{\beta \wedge \gamma \le \alpha} \gamma.
\]
To check that this definition works, we just need to check that it actually satisfies $\beta \wedge (\beta \thinsupset \alpha) \le \alpha$, which follows from
\[
\beta \wedge \Big(\bigvee_{\beta \wedge \gamma \le \alpha} \gamma\Big) = \bigvee_{\beta \wedge \gamma \le \alpha} (\beta \wedge \gamma) \le \bigvee_{\beta \wedge \gamma \le \alpha} \alpha = \alpha,
\]
where the first equality is a special case of \eqref{heyting-distributive}.

Conversely, we need to check that any complete Heyting algebra satisfies the infinite distributive law \eqref{heyting-distributive}. For this, we argue as follows:
\begin{align*}
\alpha \wedge \Big(\bigvee_{\beta \in S} \beta\Big) \le \gamma &\iff \bigvee_{\beta \in S} \beta \le \alpha \thinsupset \gamma\\
&\iff \forall \beta \in S,\ \beta \le \alpha \thinsupset \gamma\\
&\iff \forall \beta \in S,\ \alpha \wedge \beta \le \gamma\\
&\iff \bigvee_{\beta \in S} (\alpha \wedge \beta) \le \gamma.\qedhere
\end{align*}
\end{proof}

Now we compare this to the relationship between weak pseudocomplements and semidistributivity.

\begin{defn} A lattice $\cL$ is \emph{meet-semidistributive}, written SD($\wedge$), if for all $\alpha, \beta, \gamma \in \cL$ we have
\[
\alpha \wedge \beta = \alpha \wedge \gamma \;\;\; \implies \;\;\; \alpha \wedge (\beta \vee \gamma) = \alpha \wedge \beta.
\]
Similarly, a lattice $\cL$ is \emph{join-semidistributive}, written SD($\vee$), if for all $\alpha, \beta, \gamma \in \cL$ we have
\[
\alpha \vee \beta = \alpha \vee \gamma \;\;\; \implies \;\;\; \alpha \vee (\beta \wedge \gamma) = \alpha \vee \beta.
\]
A lattice is called \emph{semidistributive} if it is both meet-semidistributive and join-semidistributive.
\end{defn}

Recall from Definition \ref{defn-algebraic-lattice} that a lattice is called \emph{algebraic} if it is complete and every element can be written as a join of compact elements.

\begin{prop}[\cite{sectionally-pseudocomplemented-lattices}]\label{prop-infinite-semidistributive-algebraic} An algebraic lattice $\cL$ is meet-semidistributive iff for all $\alpha \le \beta \in \cL$, there is a weak pseudocomplement of $\beta$ over $\alpha$. In this case, $\cL$ also satisfies the following infinite form of meet-semidistributivity:
\[
\forall i,j\  \alpha \wedge \beta_i = \alpha \wedge \beta_j \ \implies\ \forall i\ \alpha \wedge \Big(\bigvee_j \beta_j\Big) = \alpha \wedge \beta_i.\label{infinite-meet-semidistributive-law}\tag{SD$_\infty$($\wedge$)}
\]
\end{prop}
\begin{proof} First we prove that every meet-semidistributive lattice has weak pseudocomplements. Note that the sublattice $\llbracket \alpha, 1 \rrbracket$ of elements of $\cL$ which are above $\alpha$ also forms an algebraic lattice: if $\theta$ is compact in $\cL$, then $\alpha \vee \theta$ is compact as an element of $\llbracket \alpha, 1 \rrbracket$. Thus we may assume without loss of generality that $\alpha = 0$, in which case we just need to prove that every element $\beta$ has a pseudocomplement $\neg \beta$.

If $\alpha = 0$, then the least possible value for $\neg \beta$ is given by
\[
\neg \beta = \bigvee_{\beta \wedge \gamma = 0} \gamma.
\]
Note that meet-semidistributivity implies that every join of finitely many elements $\gamma_i$ satisfying $\beta \wedge \gamma_i = 0$ will satisfy
\[
\beta \wedge \Big(\bigvee_{i \le n} \gamma_i\Big) = 0.
\]
We reduce the infinite case to the finite case by using the algebraicity of the lattice $\cL$. Suppose for the sake of contradiction that $\beta \wedge (\neg \beta) \ne 0$, then since $\cL$ is algebraic there is some nonzero compact element $\theta$ of $\cL$ such that
\[
\theta \le \beta \wedge (\neg \beta) \le \bigvee_{\beta \wedge \gamma = 0} \gamma.
\]
Since $\theta$ is compact, there is a finite collection of $\gamma_i$ satisfying $\beta \wedge \gamma_i = 0$ such that $\theta \le \bigvee_{i \le n} \gamma_i$. But then we have
\[
\theta \le \beta \wedge \Big(\bigvee_{i \le n} \gamma_i\Big) = 0,
\]
contradicting the assumption that $\theta$ is nonzero.

Now suppose that for all $\alpha \le \beta \in \cL$ there is a weak pseudocomplement of $\beta$ over $\alpha$. We will prove the infinite form of the meet-semidistributivity property. Suppose that there is a family $\beta_i$ such that
\[
\alpha \wedge \beta_i = \gamma
\]
for all $i$. Let $\delta$ be a weak pseudocomplement of $\alpha$ over $\gamma$. Then by the definition of a weak pseudocomplement, we have $\beta_i \in \llbracket \gamma, \delta \rrbracket$ for all $i$, so
\[
\bigvee_i \beta_i \in \llbracket \gamma, \delta \rrbracket,
\]
which in turn implies that
\[
\alpha \wedge \Big(\bigvee_i \beta_i\Big) = \gamma.\qedhere
\]
\end{proof}

\begin{rem} A similar argument can be used to show that an algebraic lattice satisfies the finite distributive law if and only if it satisfies the infinite distributive law \eqref{heyting-distributive}. In particular, if a variety $\cV$ is congruence distributive, then for every $\bA \in \cV$ the congruence lattice $\Con(\bA)$ forms a Heyting algebra.
\end{rem}


For lattices of finite length, we can show that meet-semidistributivity is a consequence of the existence of weak pseudocomplements for covers $\alpha \prec \beta$.

\begin{prop}\label{prop-semidistributive-prime-pseudocomplement} If $\alpha \wedge \beta = \alpha \wedge \gamma = \delta$ but $\alpha \wedge (\beta \vee \gamma) \ne \delta$, then for any $\epsilon$ such that
\[
\delta \prec \epsilon \le \alpha \wedge (\beta \vee \gamma),
\]
there is no weak pseudocomplement of $\epsilon$ over $\delta$.
\end{prop}
\begin{proof} Suppose for the sake of contradiction that there was some weak pseudocomplement $\theta$ of $\epsilon$ over $\delta$. Then from
\[
\delta \le \epsilon \wedge \beta \le \alpha \wedge \beta = \delta
\]
we see that $\beta \le \theta$, and similarly $\gamma \le \theta$. But then $\beta \vee \gamma \le \theta$, so we have
\[
\epsilon \le \epsilon \wedge \alpha \wedge (\beta \vee \gamma) = \epsilon \wedge (\beta \vee \gamma) = \delta,
\]
contradicting the assumption $\delta \prec \epsilon$.
\end{proof}

With the lattice-theoretic preliminaries out of the way, our task is now to show that weak pseudocomplements exist when we avoid the abelian types. We will use the concept of the relative centralizer $(\alpha : \beta)$ from Definition \ref{defn-centralizer}.

\begin{prop}\label{prop-centralizer-pseudocomplement} If $(\alpha, \beta)$ is a nonabelian prime congruence quotient on $\bA$, then the relative centralizer $(\alpha : \beta)$ is the weak pseudocomplement of $\beta$ over $\alpha$ in $\Con(\bA)$. In particular, in this case the weak pseudocomplement of $\beta$ over $\alpha$ exists.

More generally, if $\alpha \le \beta$ then $(\alpha : \beta)$ is the weak pseudocomplement of $\beta$ over $\alpha$ if and only if
\[
\beta \wedge (\alpha : \beta) = \alpha,
\]
and this occurs if the lattice $\llbracket \alpha, \beta \rrbracket$ is atomic and every prime congruence quotient $(\alpha,\delta)$ with $\alpha \prec \delta \le \beta$ is nonabelian.
\end{prop}
\begin{proof} We may assume without loss of generality that $\alpha = 0_\bA$, so we just have to prove that if $\beta \in \Con(\bA)$ is a nonabelian atomic congruence, then the centralizer $(0_\bA : \beta)$ is the pseudocomplement of $\beta$. We just need to check that
\[
C(\gamma, \beta; 0_\bA) \;\; \iff \;\; \beta \wedge \gamma = 0_\bA.
\]
By Proposition \ref{gen-commutator}(b) we see that
\[
\beta \wedge \gamma = 0_\bA \;\; \implies \;\; C(\gamma, \beta; 0_\bA)
\]
without any assumptions on $\beta$. For the other direction, since $\beta$ is an atom we have
\[
\beta \wedge \gamma \ne 0_\bA \;\; \iff \;\; \gamma \ge \beta,
\]
and by Proposition \ref{gen-commutator}(c) we see that if $\beta$ is nonabelian and $\gamma \ge \beta$ then $C(\gamma, \beta; 0_\bA)$ can't be true.

For the more general statement, note that by the argument above every $\gamma$ which satisfies $\beta \wedge \gamma = \alpha$ also satisfies $\gamma \le (\alpha : \beta)$. If $\beta \wedge (\alpha : \beta) \ne \alpha$, then picking any $\delta$ with $\alpha \prec \delta \le \beta \wedge (\alpha : \beta)$ we see that $C(\delta, \beta; \alpha)$ holds, so $\delta$ is abelian over $\alpha$ by Proposition \ref{gen-commutator}(c).
\end{proof}

In \cite{hobby-mckenzie}, the following alternative tame congruence theoretic characterization of the weak pseudocomplement of $\beta$ over $\alpha$ is given, based on Proposition \ref{prop-tame-partial-semi}.

\begin{prop} Suppose $(\alpha,\beta)$ is a nonabelian prime congruence quotient on a finite algebra $\bA$. Let $U \in M_\bA(\alpha,\beta)$ be any $(\alpha,\beta)$-minimal set, and as in Proposition \ref{prop-tame-partial-semi} let $a \in U$ be an element of the unique $(\alpha,\beta)$-trace $N$, such that there is a partial semilattice polynomial $s \in \Pol(\bA|_U)$ with $s(a,x) = x$ for all $x \in U$.

Then the weak pseudocomplement of $\beta$ over $\alpha$ is equal to the largest congruence $\delta \in \Con(\bA)$ such that $\delta|_U$ has $\{a\}$ as a congruence class (and this $\delta$ exists).
\end{prop}
\begin{proof} To see that such a $\delta$ exists, we first let $\delta'$ be the largest congruence on $\bA|_U$ with $\{a\}$ as a congruence class, and then we apply Lemma \ref{lem-idempotent-surjective-lattice} to see that the restriction map $\theta \mapsto \theta|_U$ is a surjective homomorphism from $\Con(\bA)$ to $\Con(\bA|_U)$, and we take $\delta$ to be the join of all preimages of $\delta'$ under this map.

To see that $\delta$ is the weak pseudocomplement of $\beta$ over $\alpha$, first we note that Proposition \ref{prop-single-trace} implies that $N$ is the unique $(\alpha,\beta)$-trace contained in $U$, so since $\{a\}$ is a congruence class of $\delta|_U$ we must have $(\beta \wedge \delta)|_U \subseteq \alpha|_U$. Additionally, since $\{a\}$ is a congruence class of $\alpha|_U$ we must have $\alpha \le \beta \wedge \delta$. Then since the restriction map $\llbracket \alpha, \beta \rrbracket \twoheadrightarrow \llbracket \alpha|_U, \beta|_U \rrbracket$ is $0,1$-separating we see that we must in fact have $\beta \wedge \delta = \alpha$.

Additionally, for any $\alpha \le \gamma \not\le \delta$, $\{a\}$ is not a congruence class of $\gamma|_U$, so there is some $c \in U \setminus \{a\}$ such that $(a,c) \in \gamma$. If $c \in N \setminus \{a\}$, then $(a,c) \in \beta \setminus \alpha$, so $\beta \wedge \gamma \ne \alpha$. Otherwise, let $b$ be any element of $N \setminus \{a\}$, so we have $N/\alpha|_N = \{a,b\}/\alpha|_N$. Then we have
\[
c = s(a,c) \equiv_{\beta|_U} s(b,c),
\]
and since $c/\beta|_U = c/\alpha|_U$ (since $c \not\in N$ and $N$ is the unique $(\alpha,\beta)$-trace contained in $U$), we have
\[
a \equiv_\gamma c \equiv_{\alpha|_U} s(b,c) \equiv_\gamma s(b,a) = b.
\]
Thus, in this case we have $(a,b) \in (\beta \wedge \gamma) \setminus \alpha$, so $\beta \wedge \gamma \ne \alpha$. Either way, $\gamma \not\le \delta$ implies $\beta \wedge \gamma \ne \alpha$.
\end{proof}

There is a corresponding result for dual weak pseudocomplements, as long as we exclude both the abelian types and the semilattice type.

\begin{prop}\label{prop-dual-weak-pseudo} Suppose $(\alpha,\beta)$ is a prime congruence quotient of type \textbf{3} or \textbf{4} (i.e. boolean or lattice type) on a finite algebra $\bA$, and let $N$ be any $(\alpha,\beta)$-trace. Then the dual weak pseudocomplement of $\alpha$ under $\beta$ exists and is equal to $\Cg_\bA(N^2)$.
\end{prop}
\begin{proof} By Proposition \ref{prop-tame-partial-semi}, $|N| = 2$, and by Lemma \ref{lem-idempotent-surjective-lattice}, the restriction map $\llbracket 0_\bA, \beta \rrbracket \rightarrow \Con(\bA|_N)$ is a surjective lattice homomorphism. Since $|N| = 2$ we have $\alpha|_N = 0_{\bA|_N}$, so we have
\begin{align*}
\alpha \vee \gamma = \beta \;\; &\implies \;\; \gamma \in \llbracket 0_\bA, \beta \rrbracket \text{ and } 0_{\bA|_N} \vee \gamma|_N = 1_{\bA|_N}\\
&\implies \;\; \gamma \le \beta \text{ and } N^2 \subseteq \gamma\\
&\implies \;\; \gamma \le \beta \text{ and } \gamma \not\le \alpha\\
&\implies \;\; \alpha \vee \gamma = \beta,
\end{align*}
where the last implication follows from the fact that $(\alpha,\beta)$ is a prime congruence quotient.
\end{proof}

The fact that we had to exclude type \textbf{5} from the last result isn't just an artifact of the proof: if $\bA = (\{0, 1\}, \vee)$ is a two-element semilattice, then $\Con(\bA^2)$ is depicted in Example \ref{ex-con-semi-square}, and we can see that the congruence $\Theta = \Cg_{\bA^2}\{((0,1), (1,0))\}$ has no dual weak pseudocomplement under $1_{\bA^2}$. As an abstract lattice, $\Con(\bA^2)$ is isomorphic to the lattice pictured below, which is called $\cD_2$.
\begin{center}
\begin{tikzpicture}[scale=0.5]
  \node[circle, minimum width=3pt, draw, inner sep=0pt] (1) at (0,2) {};
  \node[circle, minimum width=3pt, draw, inner sep=0pt] (p1) at (-2,0) {};
  \node[circle, minimum width=3pt, draw, inner sep=0pt] (t) at (0,0.2) {};
  \node[circle, minimum width=3pt, draw, inner sep=0pt] (p2) at (2,0) {};
  \node[circle, minimum width=3pt, draw, inner sep=0pt] (t1) at (-1,-1) {};
  \node[circle, minimum width=3pt, draw, inner sep=0pt] (t2) at (1,-1) {};
  \node[circle, minimum width=3pt, draw, inner sep=0pt] (0) at (0,-2) {};
  \draw (1) -- (p1) -- (t1) -- (0) -- (t2) -- (p2) -- (1);
  \draw (1) -- (t);
  \draw (t1) -- (t) -- (t2);
\end{tikzpicture}
\end{center}
The occurence of the lattice $\cD_2$ in $\Con(\bA^2)$ is not restricted to this particular example - the next result from \cite{hobby-mckenzie} shows that something like this occurs whenever we have a prime congruence quotient of type \textbf{5}.

\begin{prop}[Theorem 5.27 of \cite{hobby-mckenzie}] Suppose $(\alpha,\beta)$ is a nonabelian prime quotient on a finite algebra $\bA$ and let $\RR \le \bA^2$ be the basic tolerance for $(\alpha, \beta)$. Consider the sublattice
\[
\cL = \llbracket (\alpha \times \alpha)|_\RR, (\beta \times \beta)|_\RR \rrbracket
\]
of $\Con(\RR)$. If $(\alpha,\beta)$ has type \textbf{3} or \textbf{4} then $\cL$ is isomorphic to the four-element diamond lattice $\cM_2$, and if $(\alpha,\beta)$ has type \textbf{5} then $\cL$ is isomorphic to the lattice $\cD_2$ depicted above.
\end{prop}
\begin{proof} We can assume without loss of generality that $\alpha = 0_\bA$. Let $N$ be a $(0_\bA, \beta)$-trace, then $|N| = 2$ and $\bA|_N$ is polynomially equivalent to either a boolean algebra, a lattice, or a semilattice according to the type of $(0_\bA,\beta)$. Suppose $N = \{a,b\}$ and pick $e \in E(\bA)$ such that
\[
N = e(\bA) \cap a/\beta.
\]
Additionally, if $(0_\bA, \beta)$ has type \textbf{5} then assume that $a$ is the neutral element of $\bA|_N$ and that $b$ is the absorbing element.

First we check that each congruence on $(\bA|_N)^2$ extends to a congruence on $\RR$ which is contained in $(\beta \times \beta)|_\RR$. By Theorem \ref{thm-basic-tolerance}, we have
\[
\RR = \Sg_{\bA^2}(\Delta_\bA \cup N^2).
\]
Defining the unary polynomial $e^{(2)}$ on $\RR$ as in the proof of Proposition \ref{prop-restrict-variety}, we see that
\[
N^2 = e^{(2)}(\RR) \cap (a,a)/(\beta \times \beta)|_\RR,
\]
and $(\bA|_N)^2$ is polynomially equivalent to $\RR|_{N^2}$ by the argument of Proposition \ref{prop-restrict-variety}, so Lemma \ref{lem-idempotent-surjective-lattice} shows that restriction to $N^2$ defines a surjective lattice homomorphism from $\llbracket 0_\RR, (\beta \times \beta)|_\RR \rrbracket$ to $\Con((\bA|_N)^2)$.

The main difficulty is to check that every congruence $\theta$ on $\RR$ which is contained in $(\beta \times \beta)|_\RR$ is equal to $\Cg_\RR(\theta|_{N^2})$ - this requires some tedious casework. It's helpful to note that since the transitive closure of the tolerance $\RR$ is $\beta$, the congruence $(\beta \times \beta)|_\RR$ is actually the linking congruence of $\RR \le_{sd} \bA \times \bA$. In other words, we have
\[
(\beta \times \beta)|_\RR = \ker \pi_1 \vee \ker \pi_2.
\]
First we will show that the containment
\[
\Cg_\RR(\ker \pi_1|_{N^2}) \subseteq \ker \pi_1 = (0_\bA \times 1_\bA)|_\RR
\]
is an equality. Consider any $(c,d) \in \RR$. Since $\RR$ is generated by $\Delta_\bA \cup N^2$, there is some binary polynomial $p \in \Pol_2(\bA)$ such that $p(a,b) = c, p(b,a) = d$. Then we have
\[
\begin{bmatrix} c\\ d\end{bmatrix} = p\Big(\begin{bmatrix} a\\ b\end{bmatrix}, \begin{bmatrix} b\\ a\end{bmatrix}\Big) \equiv p\Big(\begin{bmatrix} a\\ a\end{bmatrix}, \begin{bmatrix} b\\ b\end{bmatrix}\Big) = \begin{bmatrix} c\\ c\end{bmatrix} \pmod{\Cg_\RR(\ker \pi_1|_{N^2})}.
\]
Since this is true for any $(c,d) \in \RR$, we see that $\Cg_\RR(\ker \pi_1|_{N^2}) = \ker \pi_1$.

Now for any $\theta \le (\beta \times \beta)|_\RR$, if $\pi_1(\theta) \ne \beta$ then $\theta \subseteq \ker \pi_1$ since $\beta$ is atomic. If $\pi_1(\theta) = \beta$, then $(a,b) \in \pi_1(e^{(2)}(\theta)) = \pi_1(\theta|_{N^2})$, so we have the implication
\[
\theta|_{N^2} \subseteq \ker \pi_1|_{N^2} \;\;\; \implies \;\;\; \theta \subseteq \ker \pi_1 = \Cg_\RR(\ker \pi_1|_{N^2}).
\]
Together with $\Cg_\RR(N^2) \supseteq \ker \pi_1 \vee \ker \pi_2$, this shows that $\theta = \Cg_\RR(\theta|_{N^2})$ if $\theta|_{N^2}$ is one of $0_N \times 0_N, 0_N \times 1_N, 1_N \times 0_N, 1_N\times 1_N$. This handles the cases where $(0_\bA, \beta)$ has type \textbf{3} or \textbf{4} (lattices and boolean algebras are congruence distributive, so congruences on $(\bA|_N)^2$ are determined by their first and second projections in these cases), so from here on we may assume that $(0_\bA, \beta)$ has type \textbf{5} (i.e. semilattice type).

By the analysis of the congruences on the square of the two-element semilattice from Example \ref{ex-con-semi-square}, we have two remaining cases: either $\theta|_{N^2}$ is the congruence generated by $((a,b), (b,a))$, or (possibly after swapping coordinates) $\theta|_{N^2}$ is the congruence generated by $((b,a),(b,b))$. In the first case, $\theta|_{N^2}$ contains both $((b,a),(b,b))$ and $((a,b),(b,b))$, so for any $(c,d) \in \RR$, if we choose the binary polynomial $p$ satisfying $p(a,b) = c, p(b,a) = d$ as before, we get
\[
\begin{bmatrix} c\\ d\end{bmatrix} = p\Big(\begin{bmatrix} a\\ b\end{bmatrix}, \begin{bmatrix} b\\ a\end{bmatrix}\Big) \equiv p\Big(\begin{bmatrix} b\\ b\end{bmatrix}, \begin{bmatrix} b\\ b\end{bmatrix}\Big) \in \Delta_\bA \pmod{\Cg_\RR(\theta|_{N^2})}.
\]
Thus in this case, every element of $\RR$ is congruent modulo $\theta$ to a diagonal element, so $\theta$ is determined by its restriction to $\Delta_\bA \cong \bA$. Since $\theta|_{\Delta_\bA} \subseteq (\beta \times \beta)|_{\Delta_\bA}$ and $(a,a)$ is not congruent to $(b,b)$ modulo $\theta$, we see that $\theta|_{\Delta_\bA}$ is trivial, so every pair of elements of $\RR$ which are congruent modulo $\theta$ are congruent to the same diagonal element of $\Delta_\bA$ via the congruence $\Cg_\RR\{((a,b),(b,a))\}$.

To finish the proof, we consider the case where $\theta|_{N^2}$ is the congruence generated by $((b,a),(b,b))$, so $\theta \subsetneq \ker \pi_1$. Suppose that the pairs $(c,d_1), (c,d_2) \in \RR$ are congruent modulo $\theta$, and choose binary polynomials $p_1,p_2$ such that $p_i(a,b) = c$ and $p_i(b,a) = d_i$. Then we have
\[
\begin{bmatrix} c\\ d_i\end{bmatrix} = p_i\Big(\begin{bmatrix} a\\ b\end{bmatrix}, \begin{bmatrix} b\\ a\end{bmatrix}\Big) \equiv p_i\Big(\begin{bmatrix} a\\ b\end{bmatrix}, \begin{bmatrix} b\\ b\end{bmatrix}\Big) = \begin{bmatrix} c\\ p_i(b,b)\end{bmatrix} \pmod{\Cg_\RR(\theta|_{N^2})}.
\]
We claim that $p_1(b,b) = p_2(b,b)$. Suppose not, for the sake of contradiction. Since $p_1(b,b) \ne p_2(b,b) \in c/\beta$, we can apply Theorem \ref{thm-minimal-sets}(c) to see that there is some unary $f \in \Pol_1(\bA)$ such that $f(p_1(b,b)) \ne f(p_2(b,b))$ and $f(c/\beta) = N$. Suppose without loss of generality that $f(p_1(b,b)) = a$, and note that since $(a,a)$ is not congruent to $(a,b)$ modulo $\theta$ we must have $f(c) = b$. Then the unary polynomial $g(x) = f(p_1(x,b))$ preserves $N$ and satisfies
\[
g(b) = f(p_1(b,b)) = a, \;\;\; g(a) = f(p_1(a,b)) = f(c) = b.
\]
Then $g|_N$ is not monotone, which contradicts the assumption that $\bA|_N$ is polynomially equivalent to a semilattice, so we must have had $p_1(b,b) = p_2(b,b)$ after all. Therefore $(c,d_1)$ is congruent to $(c,d_2)$ modulo $\Cg_\RR(\theta|_{N^2})$, and since this is true for any $(c,d_1),(c,d_2)$ which are congruent modulo $\theta$, we are done.
\end{proof}

The fact that prime congruences of types \textbf{3} and \textbf{4} have dual weak pseudocomplements has a nice concrete consequence.

\begin{prop} If $\bB$ is a finite simple algebra of boolean or lattice type (i.e., if $(0_\bB,1_\bB)$ has type \textbf{3} or \textbf{4}), then for any finite collection of finite algebras $\bA_i$, if $\bB \in \cV(\bA_1, ..., \bA_n)$ then $\bB \in HS(\bA_i)$ for some $i$.
\end{prop}
\begin{proof} Since $\bB, \bA_i$ are finite, if $\bB \in \cV(\bA_1, ..., \bA_n)$ then $\bB \in HSP_{fin}(\bA_1, ..., \bA_n)$, so there is some $\RR \le \prod_i \bA_i^{k_i}$ and some congruence $\theta \in \Con(\RR)$ such that $\bB \cong \RR/\theta$. Assume for simplicity that the $k_i$ are all $1$, by repeating some of the $\bA_i$s if necessary. We just need to prove that $\ker \pi_i \le \theta$ for some $i$ to complete the proof, since then $\bB$ will be isomorphic to a quotient of $\pi_i(\RR) \le \bA_i$.

By Proposition \ref{prop-tame-quotient}, the prime quotient $(\theta,1_\RR)$ has the same type as $(0_\bB,1_\bB)$, so by Proposition \ref{prop-dual-weak-pseudo} we see that $\theta$ has a dual weak pseudocomplement $\delta$ under $1_\RR$. If every $i$ has $\ker \pi_i \not\le \theta$, then each $i$ has $\theta \vee \ker \pi_1 = 1_\RR$, in which case we must have $\delta \le \ker \pi_i$ for all $i$. But then we have $\delta = 0_\RR$, which contradicts $\theta \vee \delta = 1_\RR$.
\end{proof}

Even though prime congruences of type \textbf{5} might not have dual weak pseudocomplements in general, the fact that they always have weak pseudocomplements can be use to prove that they have dual weak pseudocomplements in some special cases.

\begin{prop} If $(\alpha,\beta)$ is a nonabelian prime quotient on a finite algebra $\bA$, and if there is some $\gamma$ such that $\alpha \vee \gamma = \beta$ and $\alpha \wedge \gamma = 0_\bA$, then $\alpha$ has a dual weak pseudocomplement under $\beta$.

More generally, any nonabelian prime quotient $(\alpha,\beta)$ of $\bA$ has the following property: for all $\gamma$ such that $\alpha \vee \gamma = \beta$, there is a least $\delta$ such that $\alpha \vee \delta = \beta$ and $\alpha \wedge \gamma \le \delta$.
\end{prop}
\begin{proof} The more general statement follows from the first statement by replacing $\bA$ by $\bA/(\alpha \wedge \gamma)$, so suppose that $\alpha \vee \gamma = \beta$ and $\alpha \wedge \gamma = 0_\bA$.

Let $\delta$ be any atom of the lattice $\llbracket 0_\bA, \gamma \rrbracket$. Then we have
\[
\alpha \wedge \delta \le \alpha \wedge \gamma = 0_\bA
\]
and
\[
\delta \not\le \alpha, \;\;\; \delta \le \gamma \le \beta \;\;\; \implies \;\;\; \alpha \vee \delta = \beta.
\]
Thus the prime congruence quotient $(0_\bA,\delta)$ is perspective to $(\alpha,\beta)$, so it must be nonabelian since $\stackrel{s}{\sim}$ is a congruence on $\Con(\bA)$ by Theorem \ref{thm-locally-solvable-congruence} (in fact $(0_\bA,\delta)$ has the same type as $(\alpha,\beta)$ by Proposition \ref{prop-minimal-sets-perspective}).

By Proposition \ref{prop-centralizer-pseudocomplement}, $\delta$ has a pseudocomplement $(0_\bA : \delta)$. Then for any $\theta$ such that $\delta \not\le \theta$ we have
\[
\delta \wedge \theta = 0_\bA,
\]
and together with $\alpha \wedge \delta = 0_\bA$ we see that
\[
\alpha \vee \theta \le (0_\bA : \delta).
\]
Since $\delta \le \beta$ we have $\beta \not\le (0_\bA : \delta)$, so we have proven that
\[
\delta \not\le \theta \;\;\; \implies \;\;\; \alpha \vee \theta \ne \beta.
\]
Thus $\delta$ is the dual weak pseudocomplement of $\alpha$ under $\beta$.
\end{proof}

Putting together the results we have shown so far, we can give a sufficient condition for intervals in $\Con(\bA)$ to be congruence semidistributive.

\begin{prop} If $\bA$ is a finite algebra and $\alpha \le \beta \in \Con(\bA)$, then
\begin{itemize}
\item if no prime congruence quotient $(\gamma,\delta)$ with $\alpha \le \gamma \prec \delta \le \beta$ has type \textbf{1} or \textbf{2}, then $\llbracket \alpha, \beta \rrbracket$ is meet-semidistributive, and

\item if no prime congruence quotient $(\gamma,\delta)$ with $\alpha \le \gamma \prec \delta \le \beta$ has type \textbf{1}, \textbf{2}, or \textbf{5}, then $\llbracket \alpha, \beta \rrbracket$ is join-semidistributive.
\end{itemize}
\end{prop}
\begin{proof} This follows from Proposition \ref{prop-semidistributive-prime-pseudocomplement}, Proposition \ref{prop-centralizer-pseudocomplement}, and Proposition \ref{prop-dual-weak-pseudo}.
\end{proof}

We can prove much stronger results by making use of the congruence $\stackrel{s}{\sim}$ on $\Con(\bA)$.

\begin{thm} If $\bA$ is locally finite, then $\Con(\bA)/\!\stackrel{s}{\sim}$ satisfies the infinite meet-semidistributivity law \eqref{infinite-meet-semidistributive-law}.
\end{thm}
\begin{proof} Suppose that $\alpha, \beta_i \in \Con(\bA)$ satisfy $\alpha \wedge \beta_i \stackrel{s}{\sim} \alpha \wedge \beta_j$ for all $i,j$. By Proposition \ref{prop-solvable-quotient-algebraic}, we may assume without loss of generality that each $\beta_i$ is the join of all the elements of its $\stackrel{s}{\sim}$-class, in which case we must actually have
\[
\alpha \wedge \beta_i = \alpha \wedge \beta_j
\]
for all $i,j$. Let $\delta$ be the common value of $\alpha \wedge \beta_i$. By Proposition \ref{gen-commutator}(b), we have
\[
\alpha \wedge \beta_i = \delta \;\;\; \implies \;\;\; C(\beta_i, \alpha; \delta)
\]
for all $i$, so by Proposition \ref{gen-commutator}(e) we have
\[
C\big(\bigvee_i \beta_i, \alpha; \delta\big).
\]
Then by Proposition \ref{gen-commutator}(c) $\alpha \wedge \big(\bigvee_i \beta_i\big)$ is abelian over $\delta$, which implies $\alpha \wedge \big(\bigvee_i \beta_i\big) \stackrel{s}{\sim} \delta$.
\end{proof}

\begin{thm} If $\bA$ is finite and a convex sublattice $\cL \le \Con(\bA)$ contains no prime congruence quotients of type \textbf{5}, then $\cL/\!\!\stackrel{s}{\sim}$ is semidistributive (i.e. both meet-semidistributive and join-semidistributive).
\end{thm}
\begin{proof} We've already shown that $\cL/\!\stackrel{s}{\sim}$ is meet-semidistributive, so we only need to check that it is join-semidistributive. Suppose that $\alpha, \beta, \gamma \in \cL$ satisfy $\alpha \vee \beta \stackrel{s}{\sim} \alpha \vee \gamma$. We can assume without loss of generality that $\alpha$ is minimal in $\alpha/\!\stackrel{s}{\sim} \cap \cL$, and similarly for $\beta$ and $\gamma$, in which case we actually have
\[
\alpha \vee \beta = \alpha \vee \gamma,
\]
and if we call the common value $\delta$ then $\delta$ is minimal in $\delta/\!\!\stackrel{s}{\sim} \cap \cL$. If we assume for the sake of contradiction that $\alpha \vee (\beta \wedge \gamma) \ne \delta$, then there is some prime congruence quotient $(\epsilon,\delta)$ such that
\[
\alpha \vee (\beta \wedge \gamma) \le \epsilon \prec \delta,
\]
and by the dual to Proposition \ref{prop-semidistributive-prime-pseudocomplement} there can't be any dual weak pseudocomplement to $\epsilon$ under $\delta$. Since $\epsilon \in \cL$ and $\delta$ is minimal in $\delta/\!\stackrel{s}{\sim} \cap \cL$, we see that $\epsilon \not\stackrel{s}{\sim} \delta$, and by our assumption on $\cL$ the type of $(\epsilon,\delta)$ must therefore be \textbf{3} or \textbf{4}, contradicting Proposition \ref{prop-dual-weak-pseudo}.
\end{proof}



\end{appendices}

\end{document}